# QUANTUM TOPOLOGY WITHOUT TOPOLOGY

DANIEL TUBBENHAUER

Abstract. These lecture notes cover 13 sessions and are presented as an e-print, intended to evolve over time.

Quantum invariants do more than distinguish topological objects; they build bridges between topology, algebra, number theory and quantum physics helping to transfer ideas, and stimulating mutual development. They also possess deep and intriguing connections to representation theory, particularly through representations of quantum groups.

These lecture notes aim to illustrate how categorical algebra provides a framework for studying both algebra and topology. Specifically, they demonstrate how quantum invariants emerge naturally from a mostly categorical perspective.

## Contents

These lecture notes are a draft. In particular, the notes might change in the future by correcting typos, adding additional material, or improving the exposition. If you find typos or errors, let me know, mentioning the version number:

**v1.99, June 10, 2025.**

References for further reading are, for example, [**BS11**], [**BK01a**], [**EGNO15**], [**HV19**] and [**TV17**]. There will be many more references in the main body of the text.

**Acknowledgment.** I would like to thank the participants of the class "Quantum invariants" at the Universität Zürich in the Spring 2020 for feedback and their patience with me being slow. Special thanks to Super Cao, Rea Dalipi, Liam Rogel, Davide Saccardo, Joel Schmitz, and Ulrich Thiel for carefully reading these notes, spotting many typos and mistakes. Also, I am in debt of the Corona crisis 2020, which started the project of writing these notes in the first place. I also thank ChatGPT for helpful comments.

## Introduction

Motivated by the Rosetta Stone, see Figure 1, here is the ***categorical Rosetta stone***.

| Category theory | Algebra | Topology | Physics | Logic |
|---|---|---|---|---|
| objects $\mathtt{X}$ | algebraic data $\mathtt{X}$ | manifold $\mathtt{X}$ | system $\mathtt{X}$ | proposition $\mathtt{X}$ |
| morphism $\mathrm{f}\colon \mathtt{X}\to\mathtt{Y}$ | relation $\mathrm{f}\colon \mathtt{X}\to\mathtt{Y}$ | cobordism $\mathrm{f}\colon \mathtt{X}\to\mathtt{Y}$ | process $\mathrm{f}\colon \mathtt{X}\to\mathtt{Y}$ | proof $\mathrm{f}\colon \mathtt{X}\to\mathtt{Y}$ |
| monoidal product $\mathtt{X}\otimes\mathtt{Y}$ | product data $\mathtt{X}\otimes\mathtt{Y}$ | disjoint union $\mathtt{X}\otimes\mathtt{Y}$ | joint systems $\mathtt{X}\otimes\mathtt{Y}$ | conjunction $\mathtt{X}\otimes\mathtt{Y}$ |
| monoidal product $\mathrm{f}\otimes\mathrm{g}$ | parallel relations $\mathrm{f}\otimes\mathrm{g}$ | disjoint union $\mathrm{f}\otimes\mathrm{g}$ | parallel process $\mathrm{f}\otimes\mathrm{g}$ | parallel proofs $\mathrm{f}\otimes\mathrm{g}$ |

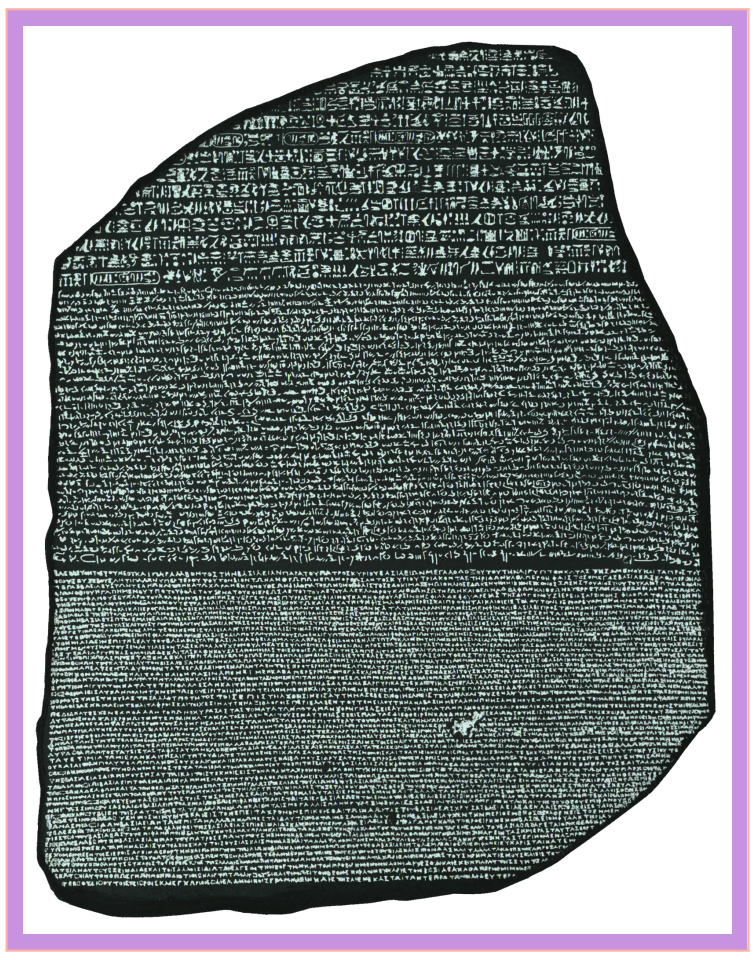

Figure 1. The Rosetta stone: the top and middle texts are in ancient Egyptian using hieroglyphic and Demotic scripts, respectively, while the bottom is in ancient Greek. The decree has only minor differences among the three versions, so the Rosetta stone became key to deciphering Egyptian hieroglyphs.

Picture from https://commons.wikimedia.org/wiki/File:Rosetta_Stone_BW.jpeg

In the 1980s, we have witnessed the birth of a fascinating new mathematical field, often called quantum algebra or quantum topology. The most spectacular achievements of this was to combine various fields of mathematics and mathematical physics such as the theory of monoidal categories, von Neumann algebras and subfactors, Hopf algebras, representations of semisimple Lie algebras, quantum field theories, the topology of knots, *etc.*, all centered around the so-called ***quantum invariants*** of links.

In these lecture notes, we focus our attention on the categorical aspects of the theory. Our goal is the construction and study of invariants of knots and links using techniques from categorical algebra only:

**Goal.** Use the left column of the categorical Rosetta stone to say something interesting about the others; especially with the focus on quantum invariants.

Summarized in a picture, the goal is to describe the categorical analog of:

Algebra: non-commutative structures
Topology: knots and links $\rightsquigarrow$

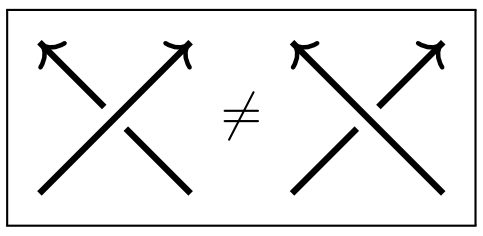

$\leftrightsquigarrow$ Physics: particles in $\mathbb{R}^3$
Logic: ribbon logic

## 1. Categories – definitions, examples and graphical calculus

The slogan for this first section is:

*Classical mathematics is based on sets, modern mathematics is based on categories.*

The meaning of “modern” is left to the reader.

### 1A. A word about conventions.

*Convention* 1A.1. Throughout, categories will be denoted by bold letters such as $\mathbf{C}$ or $\mathbf{D}$, objects by $\mathtt{X}$, $\mathtt{Y}$ *etc.* and morphisms by *e.g.* f, g. Moreover, functors are denoted by F, G *etc.*, while natural transformations are denoted by Greek letters such as $\alpha$. Further, for the sake of simplicity, we will write $\mathtt{X}\in\mathbf{C}$ for objects and $(\mathrm{f}\colon\mathtt{X}\to\mathtt{Y})\in\mathbf{C}$ (or just $\mathrm{f}\in\mathbf{C}$) for morphisms $\mathrm{f}\in\mathrm{Hom}_{\mathbf{C}}(\mathtt{X},\mathtt{Y})$, and also $\mathrm{gf}=\mathrm{g}\circ\mathrm{f}$ for composition, which is itself denoted by $\circ$. (Note our reading conventions from right to left, called ***operator notation***, which originates in the way we read functions.) When we write these, we assume that the expression makes sense, *e.g.* f’s target is g’s source. $\diamond$

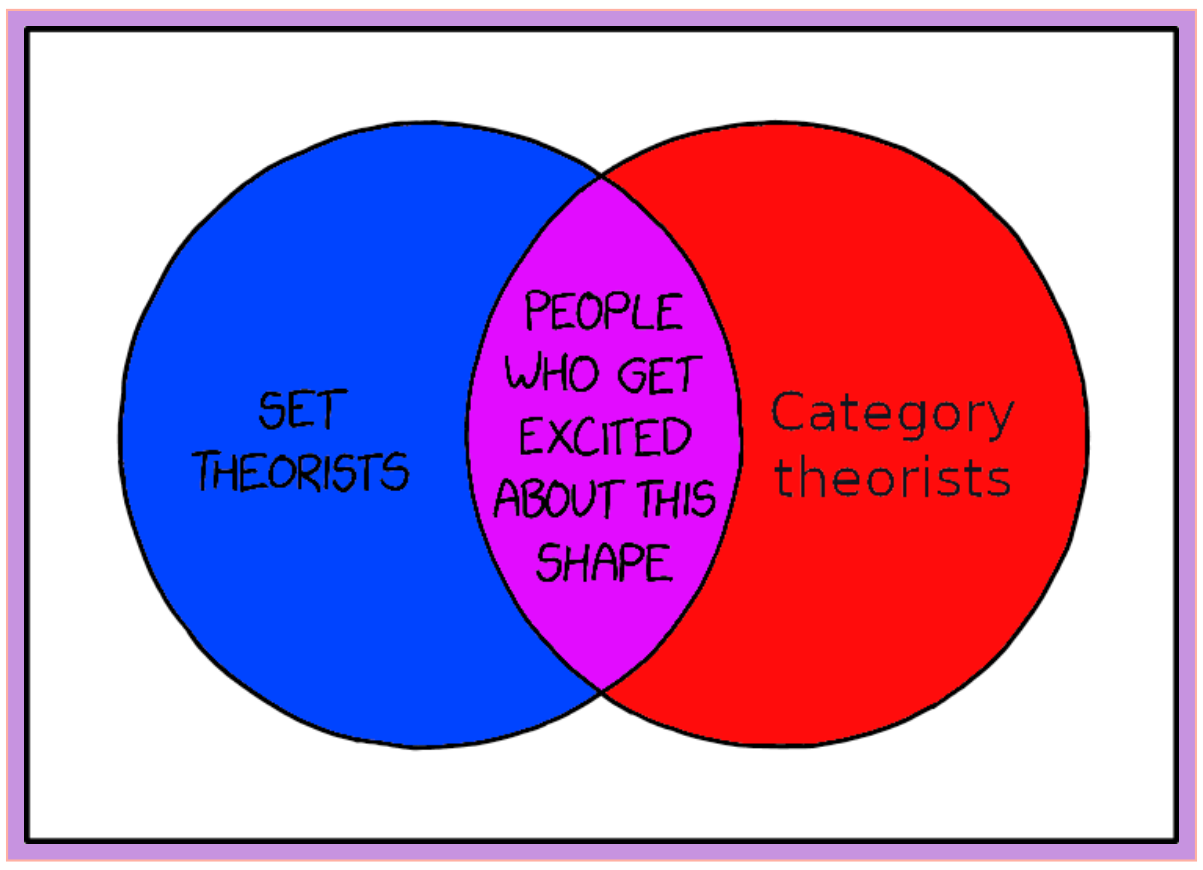

FIGURE 2. It does not matter what is modern and what is not modern. In the end we all love the same thing: Venn diagrams ☺

The picture is a variation of https://xkcd.com/2769/.

*Convention* 1A.2. There are some set theoretical issues with the definitions of some categories. For example, the objects of **Set** are all sets, but such a collection does not form a set. These issues are completely irrelevant for the purposes of these notes and are ignored throughout. ◇

*Convention* 1A.3. We will read any diagrammatic calculus from bottom to top, *cf.* Example 1B.9, and from right to left, *cf.* (2E-1). Moreover, the Feynman diagrams which we will use should be oriented, but we employ the convention that "No orientation on Feynman diagrams means upward oriented by default.". ◇

*Convention* 1A.4. $\Bbbk$ will always denote some field which we sometimes specialize to be *e.g.* of characteristic zero. If we need an algebraically closed field, we write $\mathbb{K}$, and a general associative and unital ring such as $\mathbb{Z}$ is denoted by $\mathbb{S}$. (Many constructions that we will see are stated over a field $\Bbbk$, but could also be formulated over $\mathbb{S}$. However, we find it easier to think about a field $\Bbbk$ and leave potential and easy generalizations to the reader.) ◇

1B. **Basics.** We begin at the beginning.

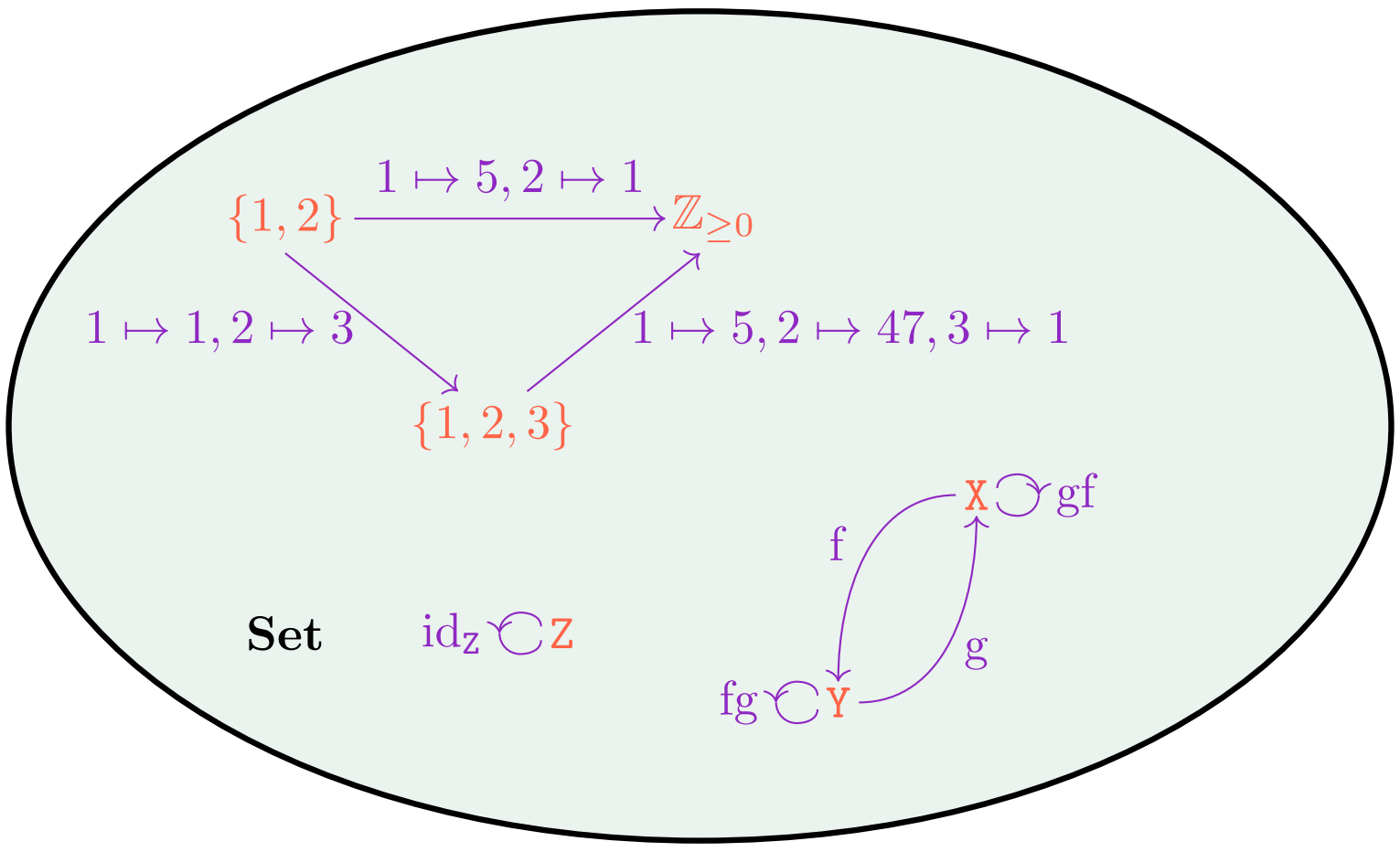

The collection of sets **Set** (category) contains sets (objects). It also contains maps (arrows or morphisms or ...), and we can compose maps in a nice way, *i.e.* compositions exist, we have identities and associativity. Formalizing this gives us the following.

**Definition 1B.1.** A ***category*** $\mathbf{C}$ consists of

- a collection of ***objects*** $\mathrm{Ob}(\mathbf{C})$;
- a set of ***morphisms*** $\mathrm{Hom}_{\mathbf{C}}(\mathtt{X},\mathtt{Y})$ for all $\mathtt{X},\mathtt{Y}\in\mathbf{C}$;

such that

(i) there exists a morphism $\mathrm{gf}\in\mathrm{Hom}_{\mathbf{C}}(\mathtt{X},\mathtt{Z})$ for all $\mathrm{f}\in\mathrm{Hom}_{\mathbf{C}}(\mathtt{X},\mathtt{Y})$ and $\mathrm{g}\in\mathrm{Hom}_{\mathbf{C}}(\mathtt{Y},\mathtt{Z})$;

(ii) there exists a morphism $\mathrm{id}_{\mathtt{X}}$ for all $\mathtt{X}\in\mathbf{C}$ satisfying $\mathrm{id}_{\mathtt{Y}}\mathrm{f}=\mathrm{f}=\mathrm{f}\mathrm{id}_{\mathtt{X}}$ for all $\mathrm{f}\in\mathrm{Hom}_{\mathbf{C}}(\mathtt{X},\mathtt{Y})$;

(iii) we have $\mathrm{h}(\mathrm{gf})=(\mathrm{hg})\mathrm{f}$ whenever this makes sense.

We also write $\mathrm{End}_{\mathbf{C}}(\mathtt{X})=\mathrm{Hom}_{\mathbf{C}}(\mathtt{X},\mathtt{X})$. ◇

The morphism gf is called the ***composition*** of g after f, while the morphism $\mathrm{id}_{\mathtt{X}}$ is called the ***identity*** on $\mathtt{X}$. The last condition in Definition 1B.1 is called ***associativity*** of morphism composition, as it is equivalent to associativity, and henceforth we omit all brackets. (To get us used to diagrammatic methods, we will sketch a proof why one can ignore brackets after the first few examples of categories.)

**Example 1B.2.** Categories generalize many familiar concepts.

(a) Categories generalize monoids: $\mathrm{End}_{\mathbf{C}}(\mathtt{X})$ is always a monoid. In some sense, a category is a bunch of monoids glued together: it is a ***monoidoid*** where the final "oid" is often used meaning "many objects".

(b) Categories generalize monoids in another way: given a monoid M, there is a category $\mathbf{M}$ with $\mathrm{Ob}(\mathbf{M}) = \{\bullet\}$ (a dummy) and $\mathrm{End}_{\mathbf{M}}(\bullet) = \mathrm{M}$, where the composition is the multiplication in M. The picture for $\mathrm{M} = \mathbb{Z}/4\mathbb{Z}$ being the cyclic group with four elements is

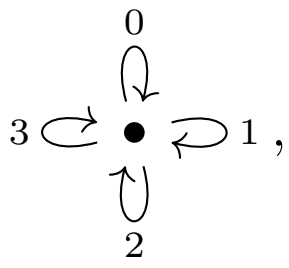

$$\text{gf} = (\text{f} + \text{g}) \bmod 4.$$

Here we write $\mathbb{Z}/4\mathbb{Z}$ additively.

(c) Categories generalize monoids in another way: There is a category **Mon** whose objects are monoids and whose morphisms are monoid homomorphisms.

(d) Categories generalize sets: there is a category **Set** whose objects are sets and whose morphisms are maps.

(e) Categories generalize vector spaces: There is a category $\mathbf{Vec}_{\Bbbk}$ whose objects are $\Bbbk$ vector spaces and whose morphisms are $\Bbbk$ linear maps. More general, the same construction gives the category of $\mathbb{S}$ modules also, abusing notation a bit, denoted by $\mathbf{Vec}_{\mathbb{S}}$.

(f) Categories generalize vector spaces in another way: There is a category $\mathbf{fdVec}_{\Bbbk}$ whose objects are finite-dimensional $\Bbbk$ vector spaces and whose morphisms are $\Bbbk$ linear maps.

With these examples in mind, one can think of categories as objects of discrete mathematics. ◇

*Remark* 1B.3. Note that categories are traditionally named after their objects, as *e.g.* **Set**, but the main players are actually the morphisms. Following this naming scheme, we will use several categories in examples with hopefully self-explanatory names, like **Group** (groups and group homomorphisms), **Ring** (rings and ring homomorphisms), or **Field** (fields and field homomorphisms). ◇

**Example 1B.4.** Later we often have categories which are denoted by **Mod**(A), which will be module categories of A. For now we observe that $\mathbf{Mod}(\mathbb{Z})$, the ***category of abelian groups***, whose objects are abelian groups (equivalently, $\mathbb{Z}$ modules $\mathbf{Vec}_{\mathbb{Z}}$) and whose morphisms are group homomorphisms, is a category. ◇

*Remark* 1B.5. Since the terminology is a bit confusing, let us spell it out: An XYZ module is nothing else than a vector space over XYZ. For example, $\mathbb{Z}$ modules $\mathbf{Mod}(\mathbb{Z})$ should be thought of as vector spaces over $\mathbb{Z}$, hence the notation $\mathbf{Vec}_{\mathbb{Z}}$. With this in mind, it is not difficult to show that every $\mathbb{Z}$ module gives to an abelian group and *vice versa*, justifying the terminology in Example 1B.4. ◇

**Example 1B.6.** It is helpful (and recommended), but formally not correct to think of morphisms as maps. For example, there is a category $\mathbf{A_3}$ having three objects and three non-identity morphisms arranged via

(1B-7)

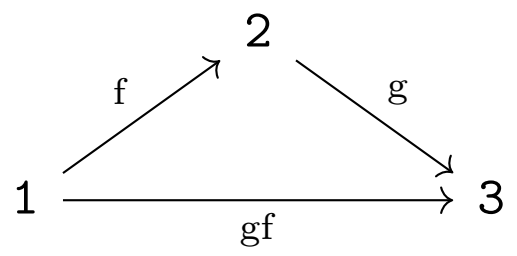

having the evident, and illustrated, composition rule. Thus, morphisms are more like "arrows" and not maps. ◇

*Remark* 1B.8. In (1B-7) we have seen the first ***commutative diagram***, which in general is a certain oriented graph, in these lecture notes. This is always to be understood that all ways composing along the various edges of the graph give the same result. In (1B-7) this is easy, as the commutative diagram is a triangle and there are only two paths to compare, which are equal by definition. However, things can get more complicated, of course, *cf.* Exercise 1K.1.

To give more familiar examples:

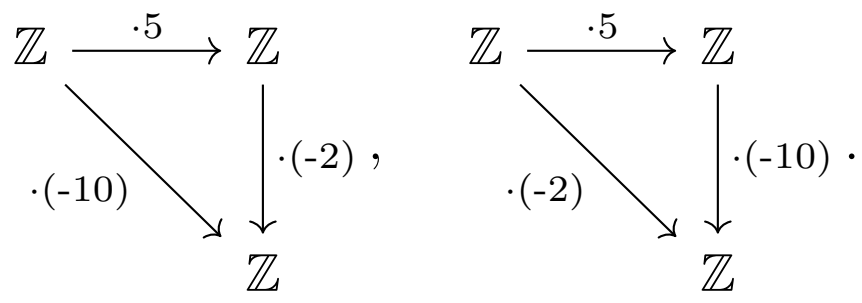

The left diagram is commutative, while the right is not. In fact, talking about biases, almost no diagram in a category is commutative, yet we predominantly illustrate only commutative diagrams. ◇

**Example 1B.9.** Very important for these lecture notes are the following examples. We will not define these formally, which is a bit painful, but rather stay with an informal, but handy, definition. (Later we will be able to give alternative and rigorous constructions.)

(a) The category **1Cob** of 1 ***dimensional cobordisms***. Its objects are 0 dimensional manifolds, a.k.a. points $\bullet^n = \bullet ... \bullet$ for $n \in \mathbb{Z}_{\geq 0}$, and its morphisms are 1 dimensional cobordisms between these, a.k.a. strands, illustrated as follows:

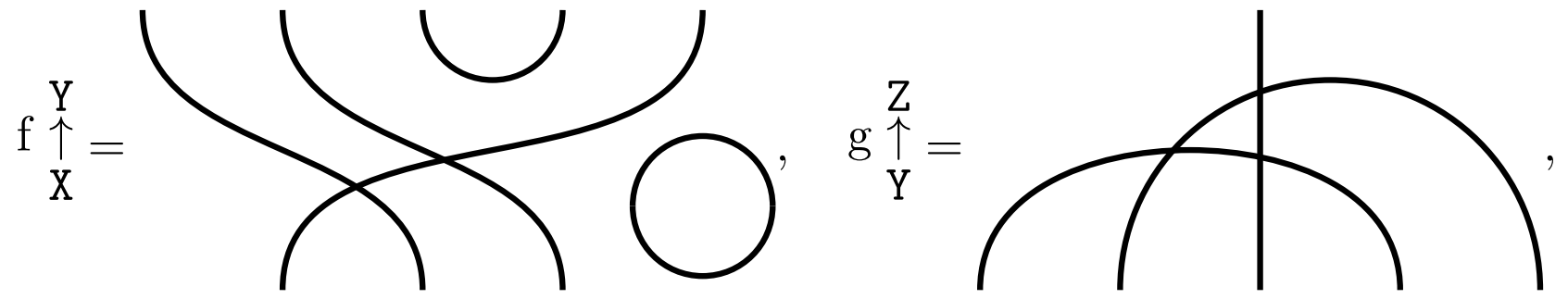

where $\mathtt{X} = \bullet\bullet\bullet = \bullet^3$, $\mathtt{Y} = \bullet\bullet\bullet\bullet\bullet = \bullet^5$ and $\mathtt{Z} = \bullet = \bullet^1$. Composition is stacking g on top of f:

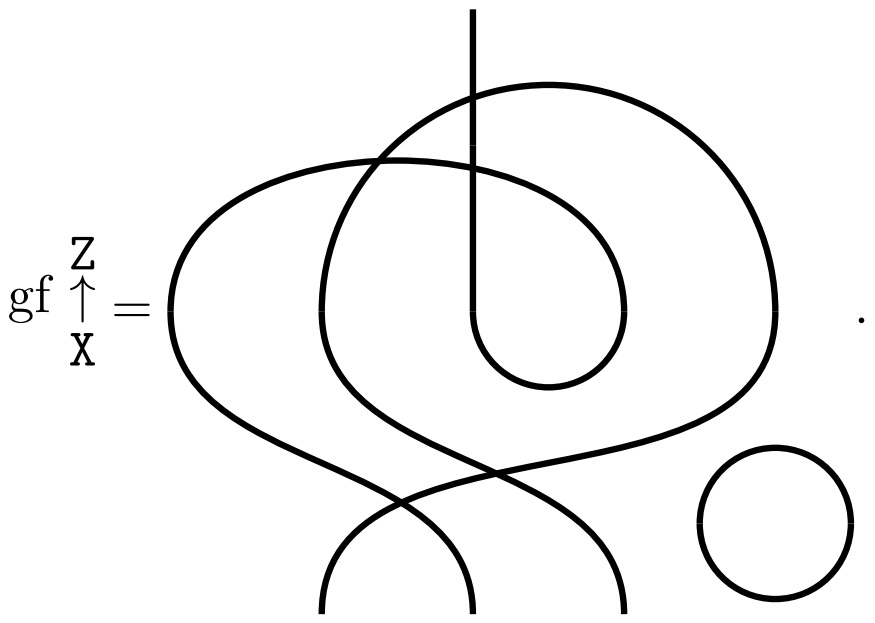

In a commuting diagram:

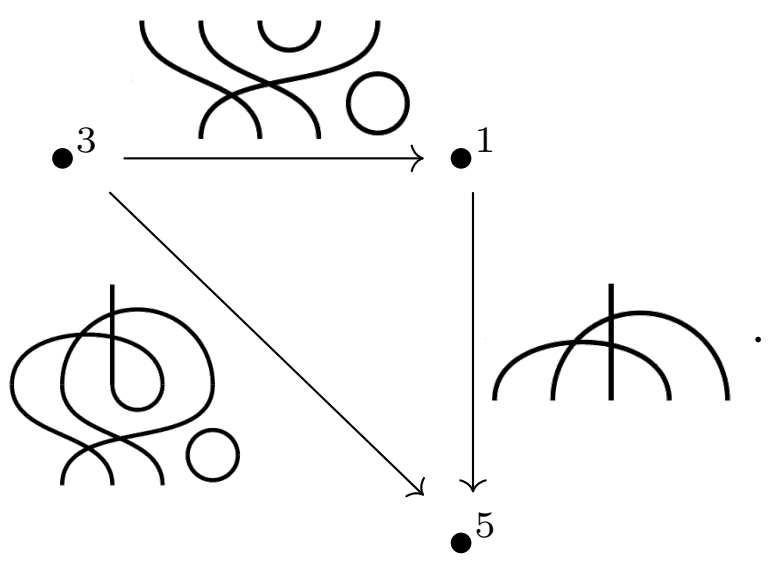

(b) The category **1Tan** of 1 ***dimensional tangles***. This is the same as **1Cob**, but now remembering some embedding into $\mathbb{R}^3$, illustrated as follows:

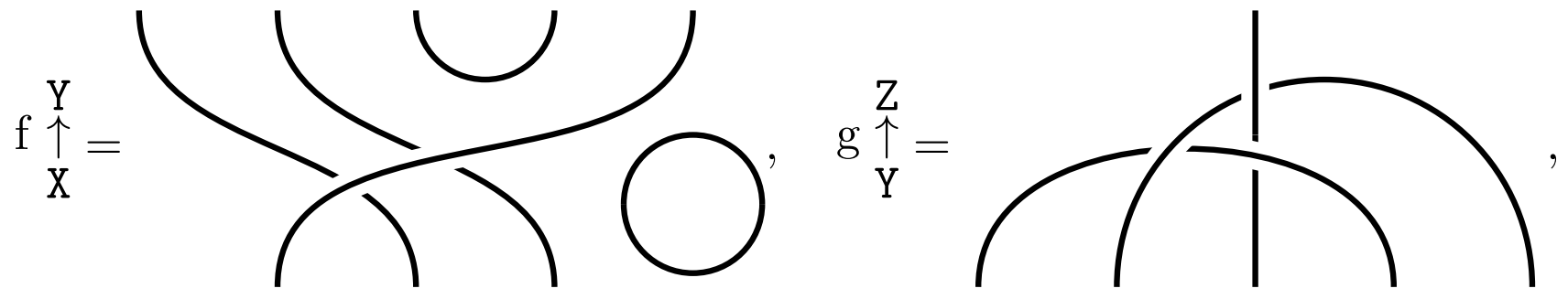

where $\mathtt{X} = \bullet\bullet\bullet$, $\mathtt{Y} = \bullet\bullet\bullet\bullet\bullet$ and $\mathtt{Z} = \bullet$. Composition is stacking g on top of f:

$\mathrm{gf}\,\overset{\mathtt{Z}}{\underset{\mathtt{X}}{\uparrow}} =$ 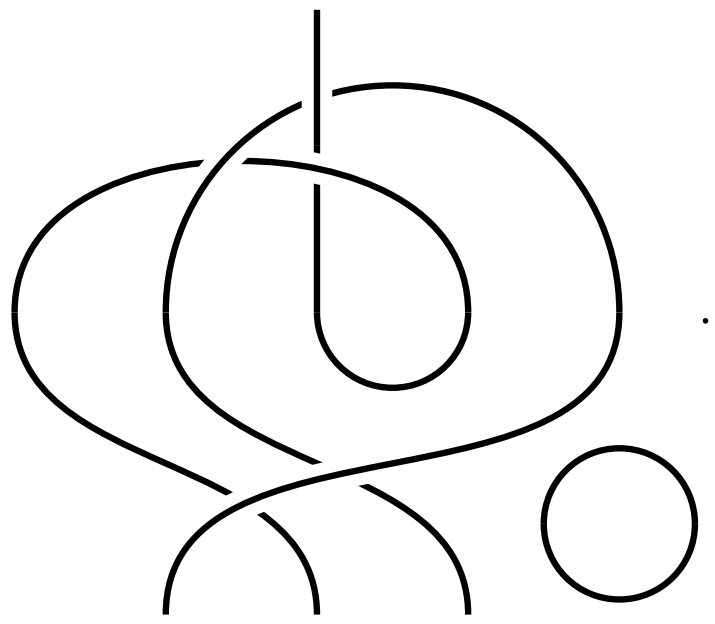 .

(c) The category **1State** of 1 ***dimensional states*** (sometimes called ***oriented tangles***), which is the category of particles moving in space with objects being particles and morphisms being worldlines. Said otherwise, it is the same as **1Tan**, but now remembering some orientation, illustrated as follows:

$\mathrm{f}\,\overset{\mathtt{Y}}{\underset{\mathtt{X}}{\uparrow}} =$ 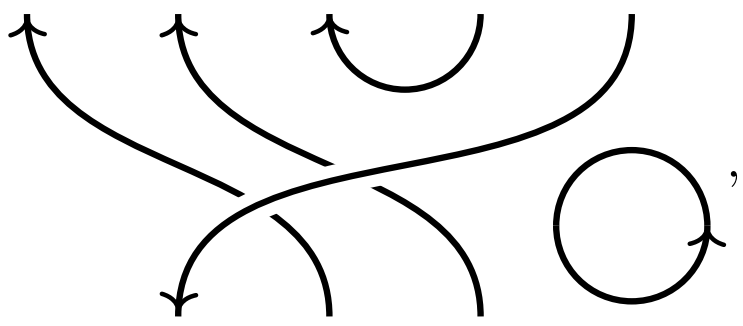 , $\mathrm{g}\,\overset{\mathtt{Z}}{\underset{\mathtt{Y}}{\uparrow}} =$ 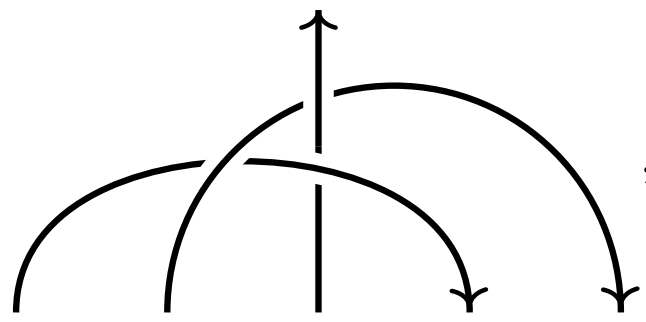 ,

where $\mathtt{X} = (\bullet^{\star})\bullet\bullet$, $\mathtt{Y} = \bullet\bullet\bullet(\bullet^{\star})\bullet^{\star}$ and $\mathtt{Z} = \bullet$, a notation which will become clear in later sections. Composition is stacking g on top of f:

$\mathrm{gf}\,\overset{\mathtt{Z}}{\underset{\mathtt{X}}{\uparrow}} =$ 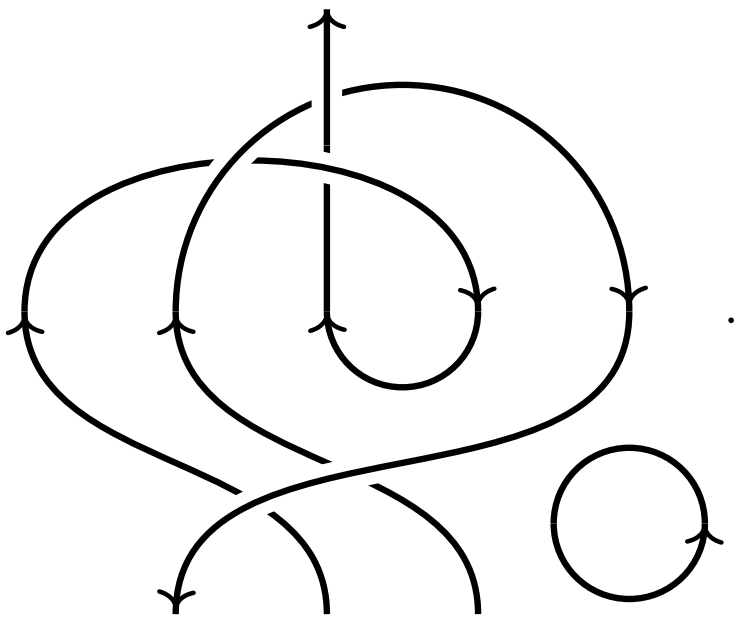 .

We will see many more examples of this flavor. ◇

Returning to associativity, let (A) be the statement "h(gf) = (hg)f" and let (B) be the statement "Same result regardless of how valid pairs of parentheses are inserted". Philosophically correct would be to use (B) as the definition and show that (A) and (B) are equivalent. However, traditionally (A) is used as the definition and then one secretly uses the following ***coherence theorem***.

**Lemma 1B.10.** *(A) ⇔ (B).*

*Proof.* Clearly, (B) implies (A). For the other way around, let us illustrate (A) and (B) using the following diagrams where every composition is replaced by a trivalent vertex:

(A): h(gf) = (hg)f ⇝ 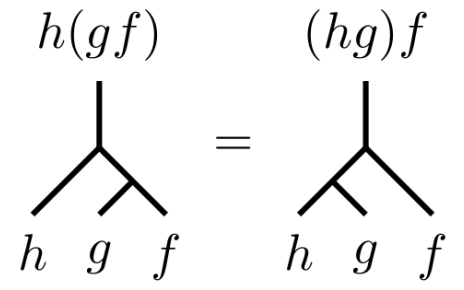

(B) ⇝ 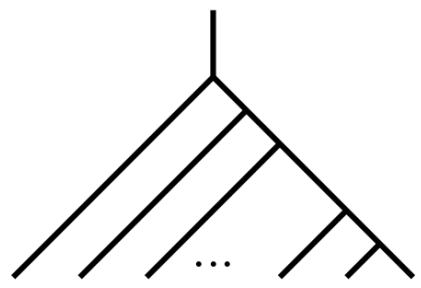 = 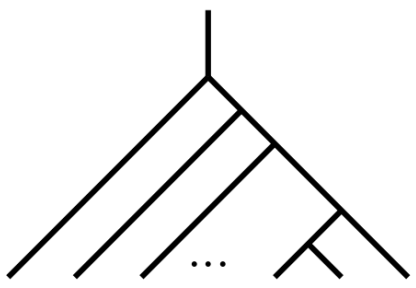 = ... .

In this notation, to show that (A) implies (B) in the case with three compositions boils down to drawing a graph where every vertex is one of the possible bracketing configurations and edges correspond to applications

of (A):

(1B-11)

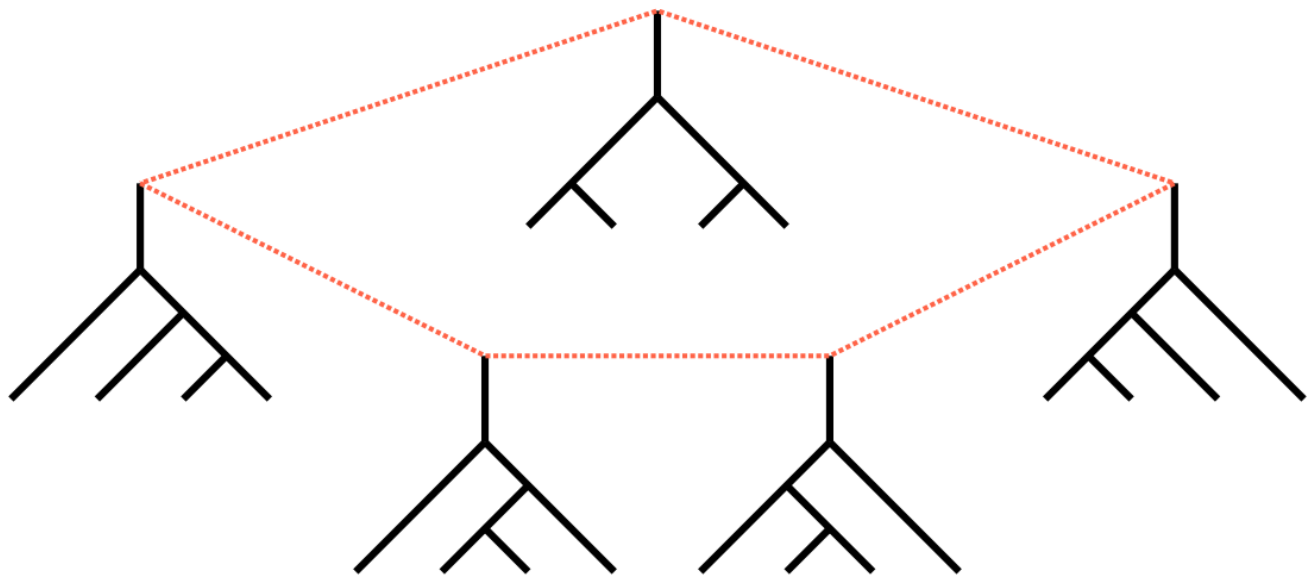

The task is to show that this graph is connected, which is immediate for three compositions by the above picture. The general case is a pleasant exercise. The geometric object that one encounters (above in Equation 1B-11 this is the dotted pentagon) is called ***associahedron***. □

**Definition 1B.12.** For any category $\mathbf{C}$, the ***pair category*** $\mathbf{C}\times\mathbf{C}$ is the category whose objects and morphisms are pairs of their corresponding types, *i.e.*

$$\mathrm{Ob}(\mathbf{C}\times\mathbf{C}) = \big\{(\mathtt{X},\mathtt{Y}) \mid \mathtt{X},\mathtt{Y}\in\mathbf{C}\big\},\quad \mathrm{Hom}_{\mathbf{C}\times\mathbf{C}}\big((\mathtt{X},\mathtt{Y}),(\mathtt{Z},\mathtt{A})\big) = \mathrm{Hom}_{\mathbf{C}}(\mathtt{X},\mathtt{Z})\times\mathrm{Hom}_{\mathbf{C}}(\mathtt{Y},\mathtt{A}),$$

and whose composition is defined componentwise. Iterating this process, we define $\mathbf{C}\times\mathbf{C}\times\ldots\times\mathbf{C}$. ◇

**Definition 1B.13.** For any category $\mathbf{C}$, the ***opposite category*** $\mathbf{C}^{op}$ is the category with the same objects and morphisms, but reversed composition:

(1B-14)

| | $\mathbf{C}$ | $\mathbf{C}^{op}$ |
|---|---|---|
| Reversed $\circ$? | No | Yes |

◇

We also write $f^{op}$ for opposite morphisms.

1C. **The duality principle.** In linear algebra, we quickly learn that statements for matrices have always analogs for the transposed matrices. In pictures:

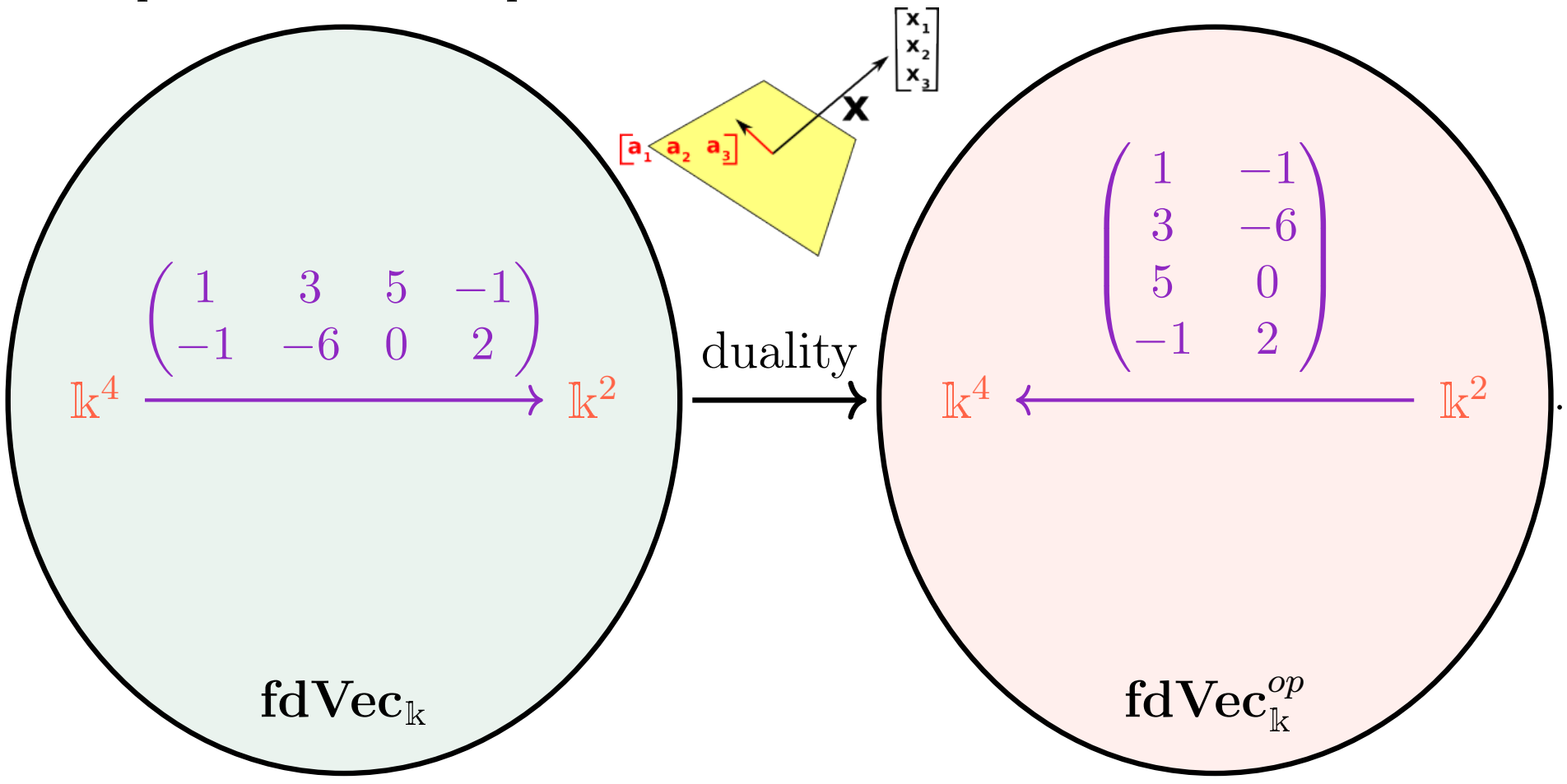

In **1Cob** the picture is

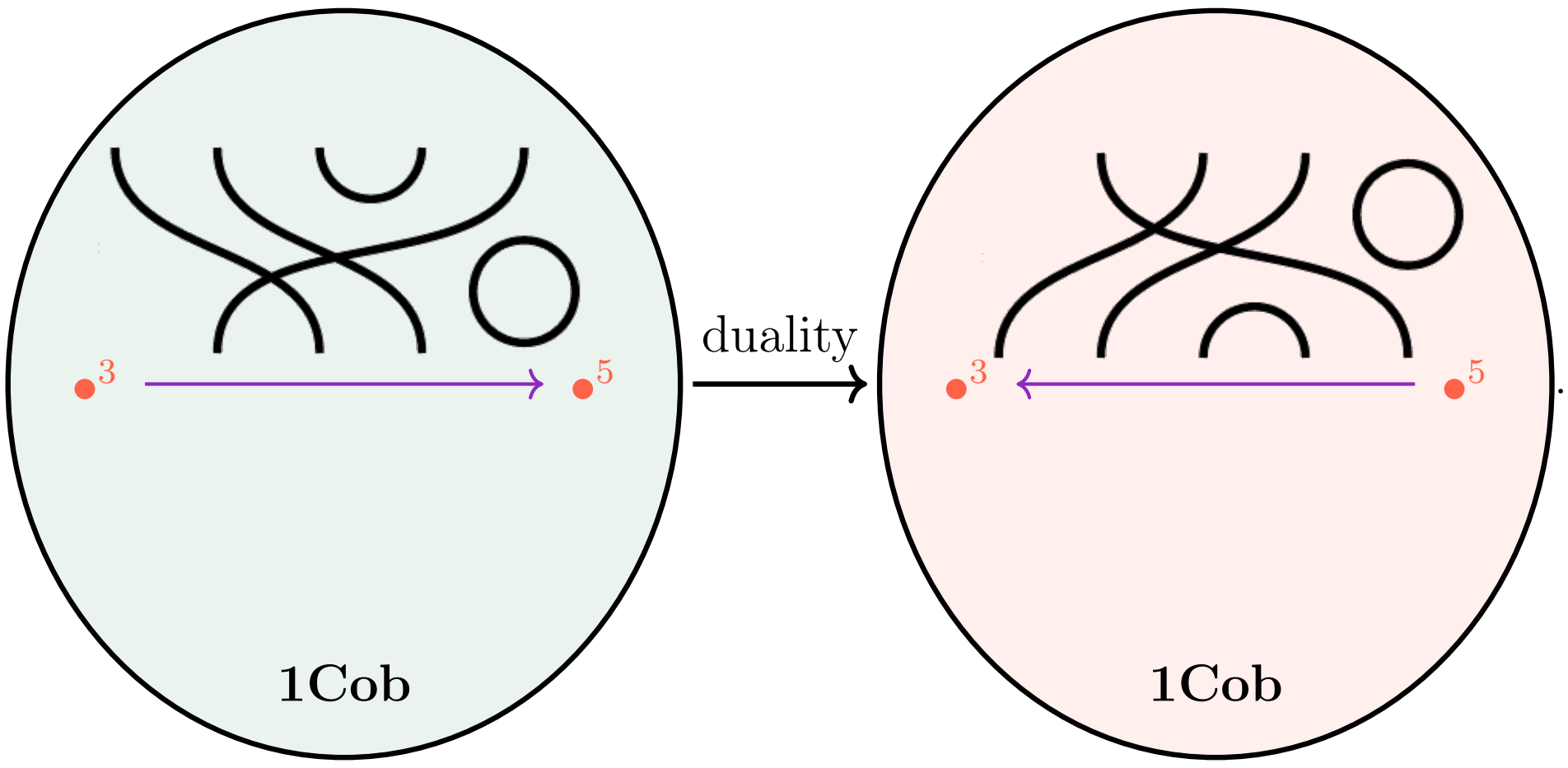

In both cases, the categories under study are essentially the same. And indeed, statements have ***dual/co*** statements, *e.g.*:

- $P_{\mathbf{C}}(\mathtt{X})$ is "$\forall\ \mathtt{Y} \in \mathbf{C}\ \exists!\ \mathrm{f}\colon \mathtt{X} \to \mathtt{Y} \in \mathbf{C}$".
- $P_{\mathbf{C}^{op}}(\mathtt{X})$ is "$\forall\ \mathtt{Y} \in \mathbf{C}^{op}\ \exists!\ \mathrm{f}\colon \mathtt{X} \to \mathtt{Y}$ in $\mathbf{C}^{op}$".
- $P^{op}_{\mathbf{C}}(\mathtt{X})$ is "$\forall\ \mathtt{Y} \in \mathbf{C}\ \exists!\ \mathrm{f}\colon \mathtt{X} \leftarrow \mathtt{Y}$ in $\mathbf{C}$".

**Theorem 1C.1.** *The* **Duality Principle** *is:*

$$\textit{Property } P \textit{ holds } \forall \textit{ categories} \Leftrightarrow \textit{property } P^{op} \textit{ holds } \forall \textit{ categories},$$

*or a variation of it.*

*Proof.* This is an experimental fact. (Theorem 1C.1 is a metatheorem not an honest theorem, so we allow this type of proof ☺.) □

The reason why we mention Theorem 1C.1 is that it allows us to be lazy: many statements have variations, and these variations are often neither stated nor proven, but still used.

1D. **Feynman diagrams.** We now discuss a convenient notation for categories, sometimes called ***Feynman*** (or ***Penrose*** or ***string*** or...) ***diagrams***, but we will also say *e.g.* ***diagrammatics***.

Given a category $\mathbf{C}$, we will denote objects $\mathtt{X} \in \mathbf{C}$ and morphisms $\mathrm{f} \in \mathbf{C}$ as

(1D-1)
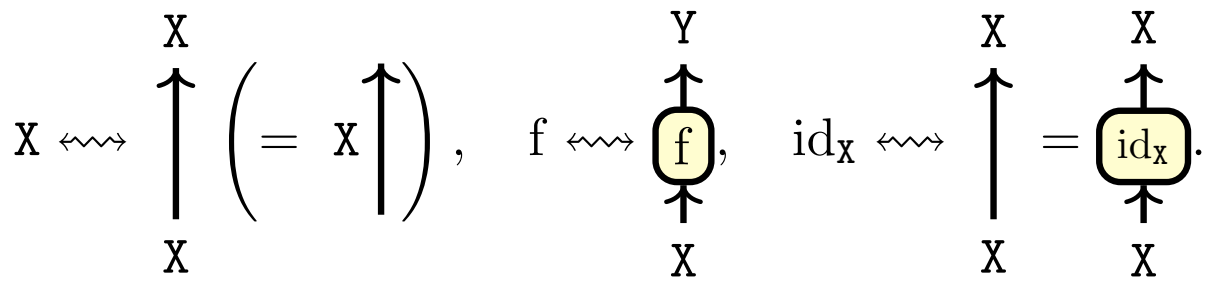

From now on, we use the convention from Convention 1A.3, meaning that we omit the orientations.

*Remark* 1D.2. This notation is "Poincaré dual" to the one $\mathrm{f}\colon \mathtt{X} \to \mathtt{Y}$ since, in diagrammatic notation, objects are strands and morphisms points, illustrated as coupons; see (1D-1). ◇

The composition is horizontal stacking, *i.e.*

(1D-3)
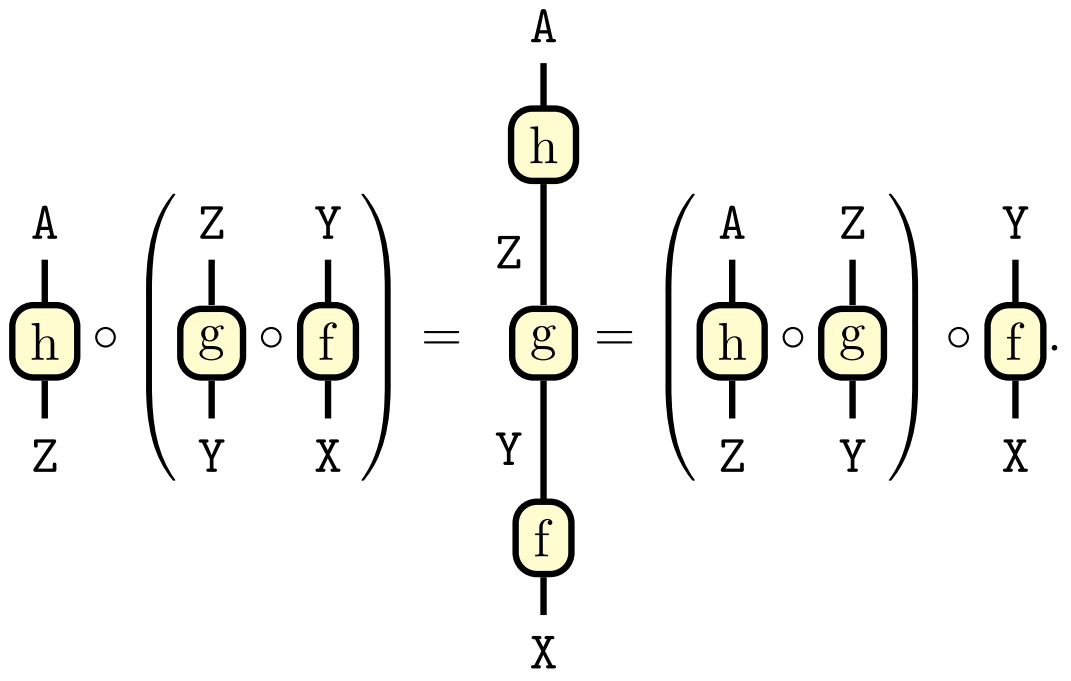

The formal rule of manipulation of these diagrams is:

(1D-4) "Two diagrams are equivalent if they are related by scaling."

The following is (almost) immediate.

**Theorem 1D.5.** *The graphical calculus is consistent; i.e. two morphisms are equal if and only if their diagrams are related by* (1D-4)*.*

*Proof.* Note that associativity is implicitly used as we have only one way to illustrate $\mathrm{h}(\mathrm{gf}) = (\mathrm{hg})\mathrm{f}$ as shown in (1D-3), while

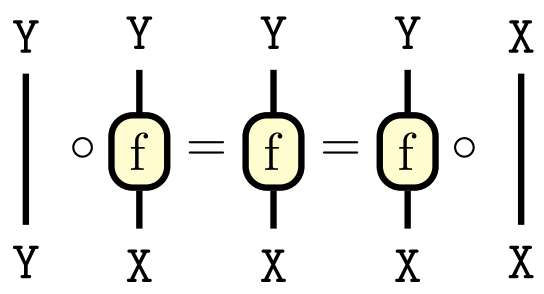

shows the identity axiom. □

*Remark* 1D.6. Later, with more structure at hand, these diagrams will turn out to be a (quite useful) 2 dimensional calculus. For now they are rather 1 dimensional, even when commutativity is imposed:

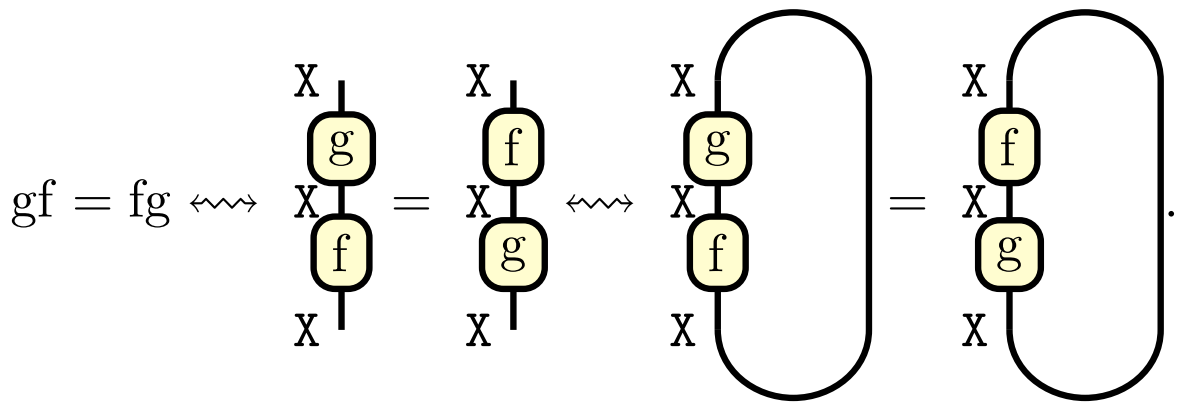

In this sense "classical mathematics lives on a line or circle". ◇

1E. **Maps between categories.** Forgetting is easy, and this is a fact that is true in the real world and mathematics at the same time. If the reader ever forgot their umbrella somewhere, then we already have a real-world example. For mathematics, since $\mathbb{k}$ vector spaces have underlying sets and $\mathbb{k}$ linear maps have underlying maps, we can forget that we have a $\mathbb{k}$ linear structure and get a map between categories (functor):

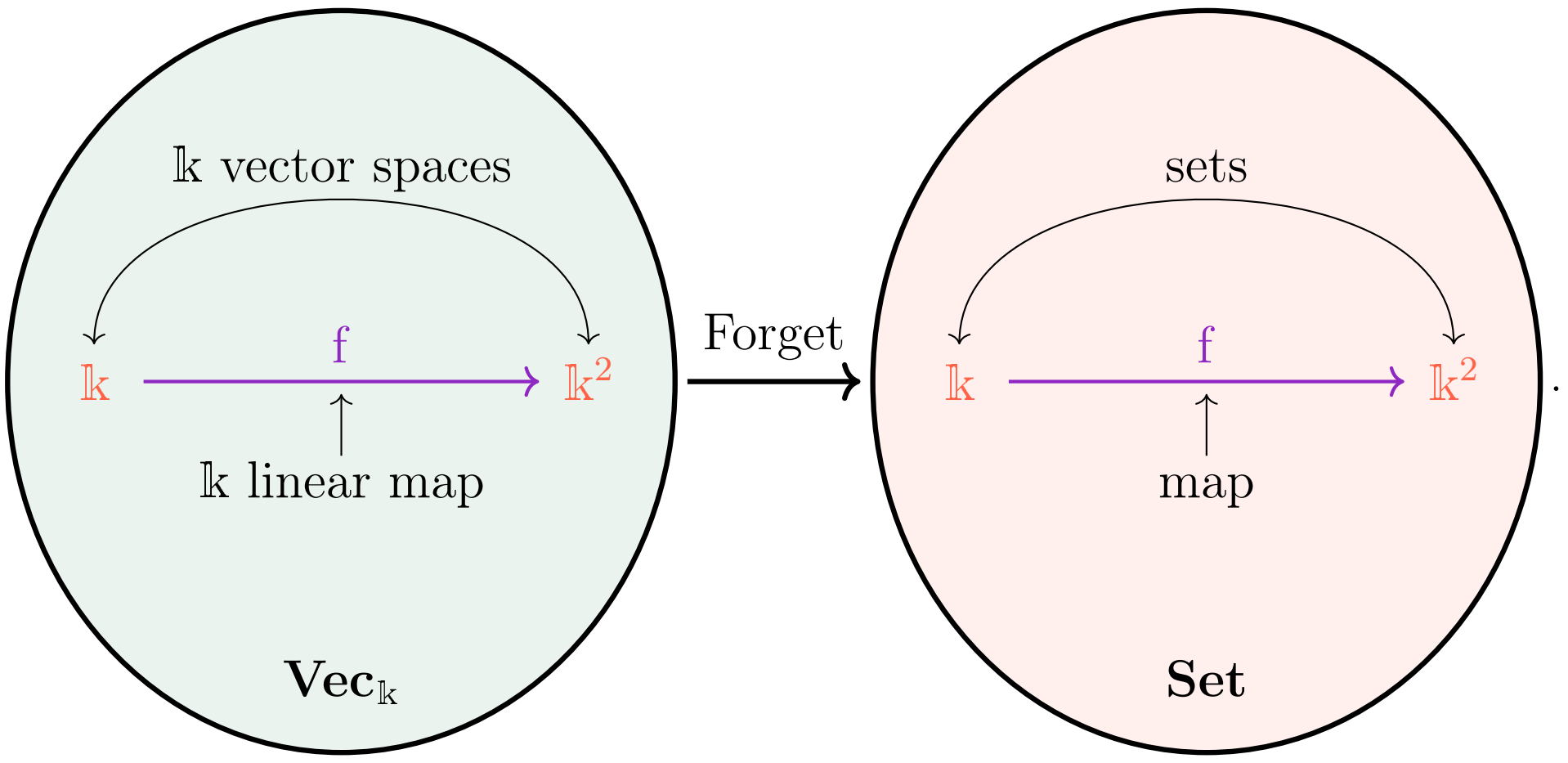

Formalizing this gives:

**Definition 1E.1.** A ***functor*** $\mathrm{F}\colon \mathbf{C} \to \mathbf{D}$ between categories $\mathbf{C}$ and $\mathbf{D}$ is a map sending

- $\mathtt{X} \in \mathbf{C}$ to an object $\mathrm{F}(\mathtt{X}) \in \mathbf{D}$;
- $(\mathrm{f}\colon \mathtt{X} \to \mathtt{Y}) \in \mathbf{C}$ to a morphism $\big(\mathrm{F}(\mathrm{f})\colon \mathrm{F}(\mathtt{X}) \to \mathrm{F}(\mathtt{Y})\big) \in \mathbf{D}$;

such that

(i) composition is preserved, *i.e.* $\mathrm{F}(\mathrm{gf}) = \mathrm{F}(\mathrm{g})\mathrm{F}(\mathrm{f})$;
(ii) identities are preserved, *i.e.* $\mathrm{F}(\mathrm{id}_{\mathtt{X}}) = \mathrm{id}_{\mathrm{F}(\mathtt{X})}$.

(As a historical side remark, in linguistics, functors are words that express grammatical relationships among other words; functors in category theory express relationships among categories.) ◇

Note that the two defining properties of categories, composition and identities, are preserved by functors. Hence, functors are maps between categories.

**Example 1E.2.** There is an ***identity functor*** $\mathrm{Id}_{\mathbf{C}}\colon \mathbf{C} \to \mathbf{C}$, sending each object and each morphism to themselves. ◇

A functor sends objects to objects and morphisms to morphisms in such a way that all relevant structures are preserved and can thus be seen as morphisms between categories. Note further that one can compose functors in the evident way (with the identity functors being identities), and the result is again a functor:

**Lemma 1E.3.** *If* F *and* G *are functors, then so is* GF. □

**Example 1E.4.** Hence, by Lemma 1E.3, we get the prototypical example of a category: **Cat**, the ***category of categories***, whose objects are categories and whose morphisms are functors.

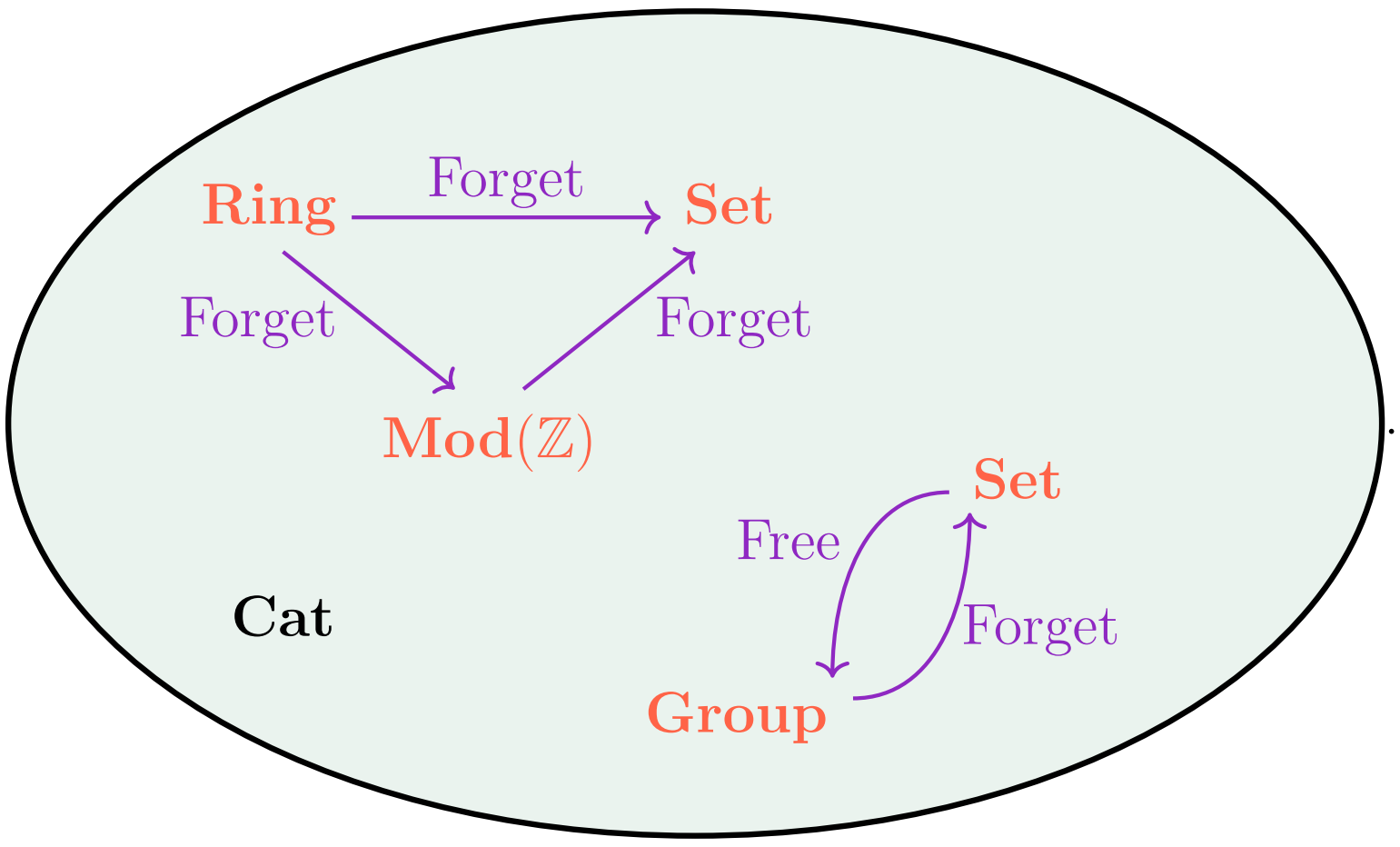

Well, there are set-theoretical issues with **Cat** (similar to what happens if one tries to define the set of sets), but let us ignore that. ◇

**Example 1E.5.** Functors generalize many familiar concepts.

(a) Functors generalize monoid maps: A functor on $\mathrm{End}_{\mathbf{C}}(\mathtt{X})$ is a homomorphism of monoids.

(b) Functors generalize monoid maps, take 2: A functor $\mathrm{F}\colon \mathbf{M}\to\mathbf{M}'$ between monoid categories $\mathbf{M}$ and $\mathbf{M}'$, as in Example 1B.2.(b), is a homomorphism of monoids.

(c) Functors generalize models: A functor $\mathrm{F}\colon \mathbf{M}\to\mathbf{Set}$ between a monoid category $\mathbf{M}$ and $\mathbf{Set}$ assigns a set $\mathrm{F}(\bullet)$ to $\bullet$ and an endomorphism $\mathrm{F}(f)$ of $\mathrm{F}(\bullet)$ to $f\in\mathbf{M}$, which can be seen as a concrete model of the underlying monoid M.

(d) Functors generalize representations: A functor $\mathrm{F}\colon \mathbf{M}\to\mathbf{Vec}_{\Bbbk}$ between a monoid category $\mathbf{M}$ and $\mathbf{Vec}_{\Bbbk}$ assigns a $\Bbbk$ vector space $\mathrm{F}(\bullet)$ to $\bullet$ and a $\Bbbk$ linear endomorphism $\mathrm{F}(f)$ of $\mathrm{F}(\bullet)$ to $f\in\mathbf{M}$, which can be seen as a representation of the underlying monoid M.

(e) Functors generalize forgetting: There is a functor $\mathrm{Forget}\colon \mathbf{Vec}_{\Bbbk}\to\mathbf{Set}$ that forgets the underlying $\Bbbk$ linear structure.

(f) Functors generalize free structures: There is a functor $\mathrm{Free}\colon \mathbf{Set}\to\mathbf{Vec}_{\Bbbk}$ for which $\mathrm{Free}(\mathtt{X})$ is the free $\Bbbk$ vector space with basis $\mathtt{X}$ and $\mathrm{Free}(f)$ is the $\Bbbk$ linear extension of $f$.

There are many more examples, as we shall see. ◇

Finally, note that any functor $\mathrm{F}\colon \mathbf{C}\to\mathbf{D}$ gives rise to a natural map

$$\mathrm{Hom}_{\mathbf{C}}(\mathtt{X},\mathtt{Y})\to\mathrm{Hom}_{\mathbf{D}}\big(\mathrm{F}(\mathtt{X}),\mathrm{F}(\mathtt{Y})\big),\ f\mapsto\mathrm{F}(f),$$

which we often use without further comment. In particular:

**Example 1E.6.** There are ***hom functors***:

$$\mathrm{Hom}_{\mathbf{C}}(\mathtt{X},_)\colon \mathbf{C}\to\mathbf{Set},\ \begin{cases}\mathtt{Y}\mapsto\mathrm{Hom}_{\mathbf{C}}(\mathtt{X},\mathtt{Y}),\\ f\mapsto(f\circ_),\end{cases}\qquad \mathrm{Hom}_{\mathbf{C}}(_,\mathtt{X})\colon \mathbf{C}^{op}\to\mathbf{Set},\ \begin{cases}\mathtt{Y}\mapsto\mathrm{Hom}_{\mathbf{C}}(\mathtt{Y},\mathtt{X}),\\ f\mapsto(_\circ f).\end{cases}$$

(For the reason to write $\mathbf{C}^{op}$, see Remark 1E.7.) ◇

*Remark* 1E.7. A functor $\mathrm{F}\colon \mathbf{C}^{op}\to\mathbf{D}$, such as $\mathrm{Hom}_{\mathbf{C}}(_,\mathtt{X})$, is sometimes seen as a ***contravariant functor*** $\mathrm{F}\colon \mathbf{C}\to\mathbf{D}$, meaning that $\mathrm{F}(gf)=\mathrm{F}(f)\mathrm{F}(g)$ holds instead of $\mathrm{F}(gf)=\mathrm{F}(g)\mathrm{F}(f)$. ◇

Let us use functors to revisit commuting diagrams.

**Definition 1E.8.** An ***abstract diagram*** is a directed graph $\mathbf{J}$, seen as a category (similarly as in Example 1B.6). Given a category $\mathbf{C}$, a ***diagram*** D ***of shape*** $\mathbf{J}$ ***in*** $\mathbf{C}$ is a functor

$$\mathrm{D}\colon \mathbf{J}\to\mathbf{C}.$$

It ***commutes*** if all directed paths in $\mathrm{D}(\mathbf{J})$ with the same start and endpoints lead to the same result in $\mathbf{C}$. ◇

**Example 1E.9.** Consider the triangle directed graph:

$$\mathbf{J} = \begin{array}{ccc} \bullet & \longrightarrow & \bullet \\ & \searrow & \downarrow \\ & & \bullet \end{array} .$$

This shape has many diagrams, e.g. $\mathrm{D}\colon \mathbf{J} \to \mathbf{Mod}(\mathbb{Z})$ or $\mathrm{E}\colon \mathbf{J} \to \mathbf{1Cob}$ given by

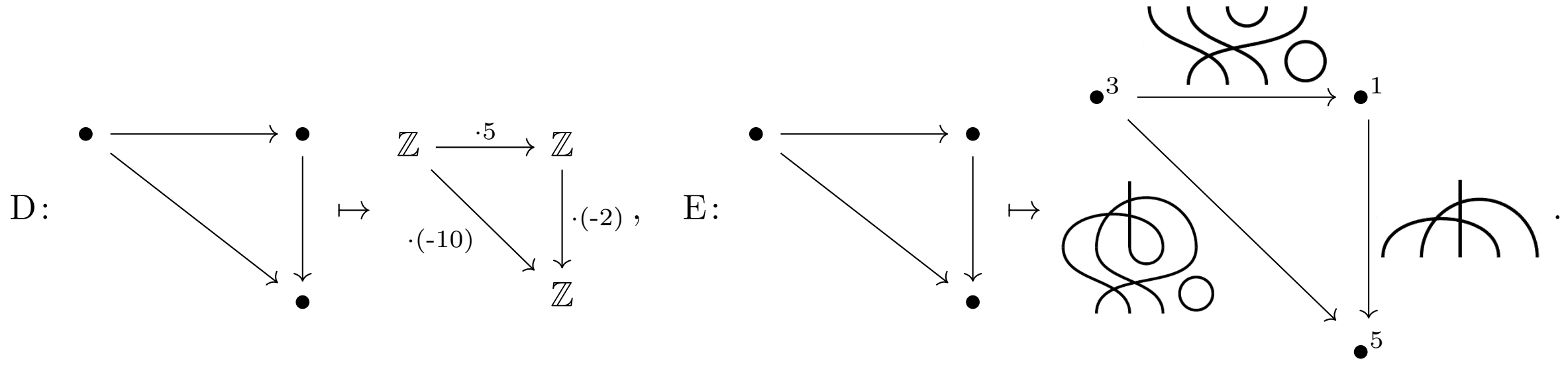

The objects and morphisms in $\mathbf{J}$ are largely irrelevant, only the shape (in the linguistic sense) plays a role. ◇

1F. **Maps between maps between categories.** A key fact that one learns in linear algebra is that $\mathtt{V} \in \mathbf{fdVec}_{\Bbbk}$ is isomorphic to its dual $\mathtt{V}^{\star}$, and an isomorphism is obtained by choosing a basis of $\mathtt{V}$ by sending the basis elements to their indicator functions. However, this isomorphism depends on the choice of basis. In contrast, $\mathtt{V}$ is also isomorphic to the double dual $\mathtt{V}^{\star\star}$, but in this case there is a basis free isomorphism given by $v \mapsto \mathrm{eval}_v$, where $\mathrm{eval}_v\colon \mathtt{V}^{\star} \to \Bbbk, \mathrm{eval}_v(f) = f(v)$ evaluates $f \in \mathtt{V}^{\star}$ at $v$.

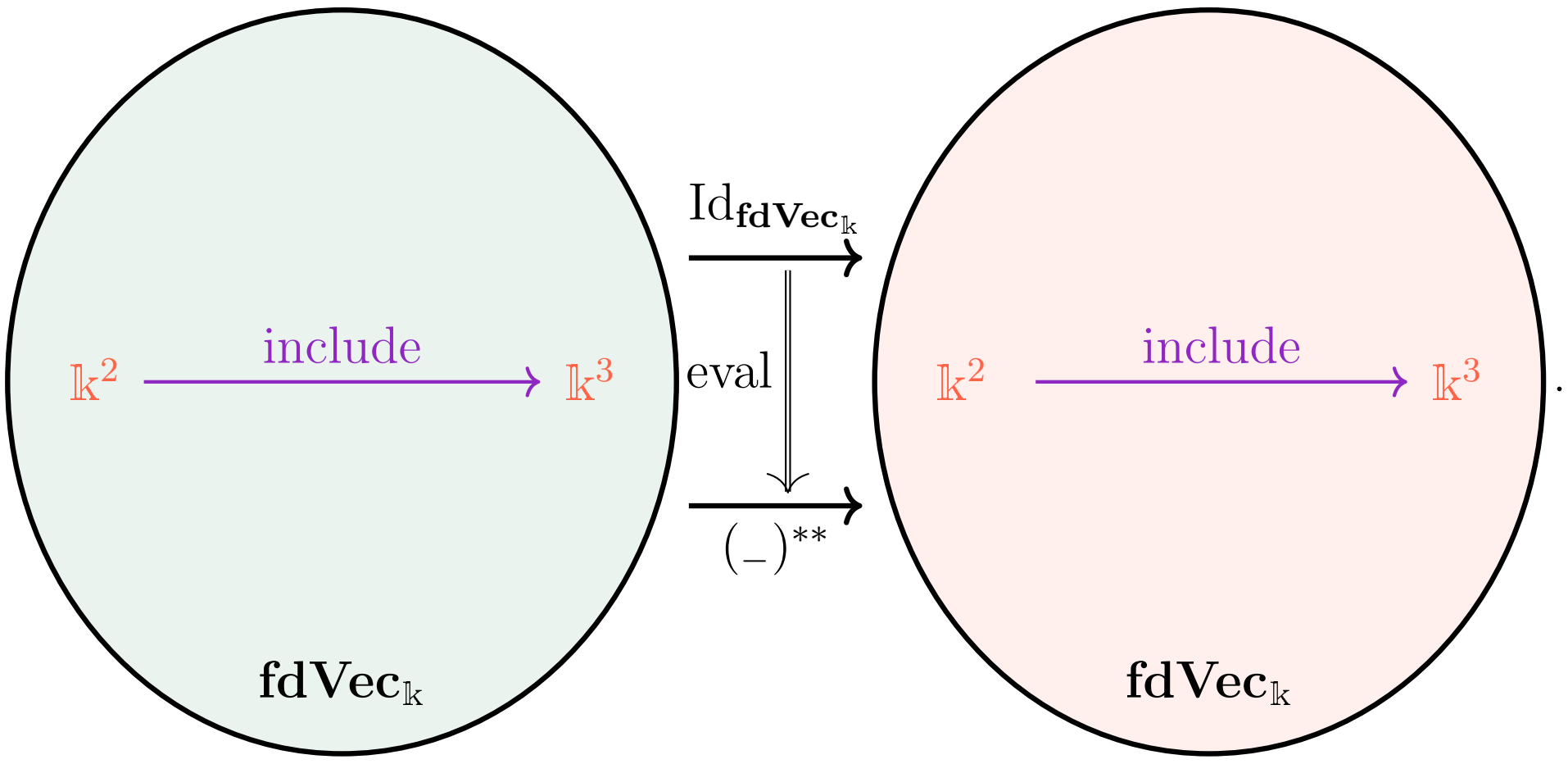

Generalizing this gives maps between functor:

**Definition 1F.1.** A ***natural transformation*** $\alpha\colon \mathrm{F} \Rightarrow \mathrm{G}$ between functors $\mathrm{F}, \mathrm{G}\colon \mathbf{C} \to \mathbf{D}$ is a collection of morphisms in $\mathbf{D}$

$$\{\alpha_{\mathtt{X}}\colon \mathrm{F}(\mathtt{X}) \to \mathrm{G}(\mathtt{X}) \mid \mathtt{X} \in \mathbf{C}\}$$

such that the following diagram commutes for all $\mathtt{X}, \mathtt{Y}, \mathrm{f} \in \mathbf{C}$:

(1F-2)
$$\begin{array}{ccc} \mathrm{F}(\mathtt{X}) & \xrightarrow{\mathrm{F}(\mathrm{f})} & \mathrm{F}(\mathtt{Y}) \\ {\scriptstyle\alpha_{\mathtt{X}}}\downarrow & & \downarrow{\scriptstyle\alpha_{\mathtt{Y}}} \\ \mathrm{G}(\mathtt{X}) & \xrightarrow[\mathrm{G}(\mathrm{f})]{} & \mathrm{G}(\mathtt{Y}) \end{array} .$$

Note that the functors $\mathrm{F}, \mathrm{G}\colon \mathbf{C} \to \mathbf{D}$ are evaluated. ◇

*Remark* 1F.3. The diagram in (1F-2) (which is the classical way of illustrating natural transformations, sometimes also called ***natural*** or ***naturality***) is of course the same as

$$\alpha \overset{\mathrm{G}}{\underset{\mathrm{F}}{\Uparrow}} \rightsquigarrow \begin{array}{ccc} \mathrm{G}(\mathtt{X}) & \xrightarrow{\mathrm{G}(\mathrm{f})} & \mathrm{G}(\mathtt{Y}) \\ {\scriptstyle\alpha_{\mathtt{X}}}\uparrow & & \uparrow{\scriptstyle\alpha_{\mathtt{Y}}} \\ \mathrm{F}(\mathtt{X}) & \xrightarrow[\mathrm{F}(\mathrm{f})]{} & \mathrm{F}(\mathtt{Y}) \end{array} ,$$

which, using our reading conventions, is saying that $\alpha$ can be seen as a morphism from F to G. ◇

There is a composition of natural transformations, called the ***vertical composition*** and denoted by $\circ$, of natural transformations given by

$$
\begin{array}{c} \mathrm{H} \\ \beta \Uparrow \\ \mathrm{G} \\ \alpha \Uparrow \\ \mathrm{F} \end{array}
\rightsquigarrow \beta\alpha_{\mathtt{X}}
\left(
\begin{array}{ccc}
\mathrm{H}(\mathtt{X}) & \xrightarrow{\mathrm{H}(\mathrm{f})} & \mathrm{H}(\mathtt{Y}) \\
\beta_{\mathtt{X}} \uparrow & & \uparrow \beta_{\mathtt{Y}} \\
\mathrm{G}(\mathtt{X}) & \xrightarrow[\mathrm{G}(\mathrm{f})]{} & \mathrm{G}(\mathtt{Y}) \\
\alpha_{\mathtt{X}} \uparrow & & \uparrow \alpha_{\mathtt{Y}} \\
\mathrm{F}(\mathtt{X}) & \xrightarrow[\mathrm{F}(\mathrm{f})]{} & \mathrm{F}(\mathtt{Y})
\end{array}
\right) \beta\alpha_{\mathtt{Y}} .
\tag{1F-4}
$$

**Example 1F.5.** There is an ***identity natural transformation*** $\mathrm{ID}_{\mathrm{F}} \colon \mathrm{F} \to \mathrm{F}$, $(\mathrm{ID}_{\mathrm{F}})_{\mathtt{X}} = \mathrm{id}_{\mathtt{X}}$. ◇

Clearly:

**Lemma 1F.6.** *If $\alpha$ and $\beta$ are natural transformations, then so is $\beta\alpha$.* □

**Example 1F.7.** By Lemma 1F.6, there is a category $\mathbf{Hom}(\mathbf{C}, \mathbf{D})$, the ***category of functors*** from $\mathbf{C}$ to $\mathbf{D}$. Its objects are all such functors, and its morphisms are natural transformations, with composition being vertical composition. A special case are ***endofunctors***, whose category we denote by $\mathbf{End}(\mathbf{C}) = \mathbf{Hom}(\mathbf{C}, \mathbf{C})$, which will play an important role. ◇

**Example 1F.8.** Natural transformations generalize intertwiners (a.k.a. homomorphisms of representations): given two representations $\mathrm{F}, \mathrm{G} \colon \mathbf{M} \to \mathbf{Vec}_{\Bbbk}$ as in Example 1E.5.(c), a natural transformation between them would provide a commuting diagram

$$
\begin{array}{ccc}
\mathrm{F}(\bullet) & \xrightarrow{\mathrm{F}(\mathrm{f})} & \mathrm{F}(\bullet) \\
\alpha_{\bullet} \downarrow & & \downarrow \alpha_{\bullet} \\
\mathrm{G}(\bullet) & \xrightarrow[\mathrm{G}(\mathrm{f})]{} & \mathrm{G}(\bullet)
\end{array}
\Rightarrow \mathrm{F}(\mathrm{f})\alpha_{\bullet} = \alpha_{\bullet}\mathrm{G}(\mathrm{f}),
$$

where $\mathrm{F}(\bullet)$ and $\mathrm{G}(\bullet)$ are the $\Bbbk$ vector spaces associated to the representations, and $\alpha_{\bullet} \colon \mathrm{F}(\bullet) \to \mathrm{G}(\bullet)$ is a $\Bbbk$ linear map between them. ◇

**Example 1F.9.** Having a monoid category $\mathbf{M}$, the category $\mathbf{Hom}(\mathbf{M}, \mathbf{Vec}_{\Bbbk})$ can be identified with all representations of the underlying monoid. ◇

1G. **2d Feynman diagrams.** We will come back to this later, but let us mention right away that the calculus of functors gets a 2 dimensional flavor via ***colored Feynman diagrams***.

That is, we color faces by categories, strings by functors, and use coupons to illustrate natural transformations. For example, say that the category $\mathbf{C}$ is colored orchid, and $\mathbf{D}$ is colored spinach, then, in operator notation:

$\mathrm{F} \colon \mathbf{C} \to \mathbf{D} \rightsquigarrow$ 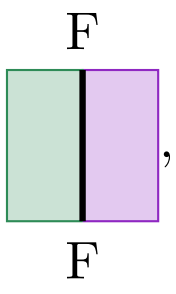, $\mathrm{G} \colon \mathbf{D} \to \mathbf{C} \rightsquigarrow$ 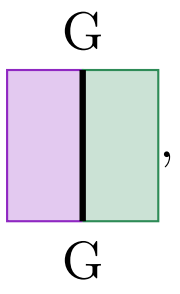,

$\epsilon \colon \mathrm{FG} \Rightarrow \mathrm{Id}_{\mathbf{D}} \rightsquigarrow$ 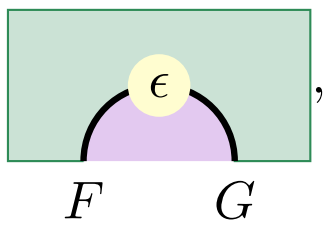, $\iota \colon \mathrm{Id}_{\mathbf{C}} \Rightarrow \mathrm{GF} \rightsquigarrow$ 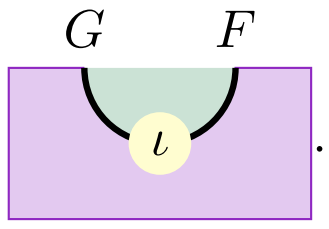.

Note that, by convention, we do not draw identity functors, *e.g.*

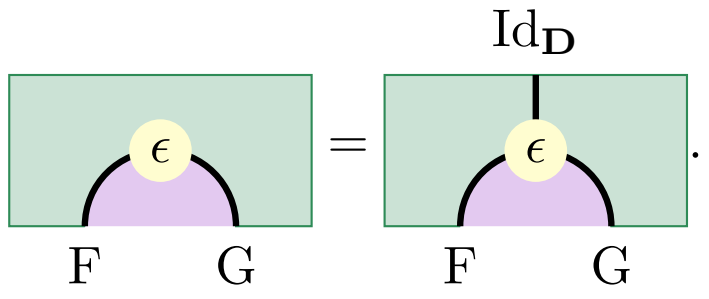

In this notation, compatibility conditions will look like, for example,

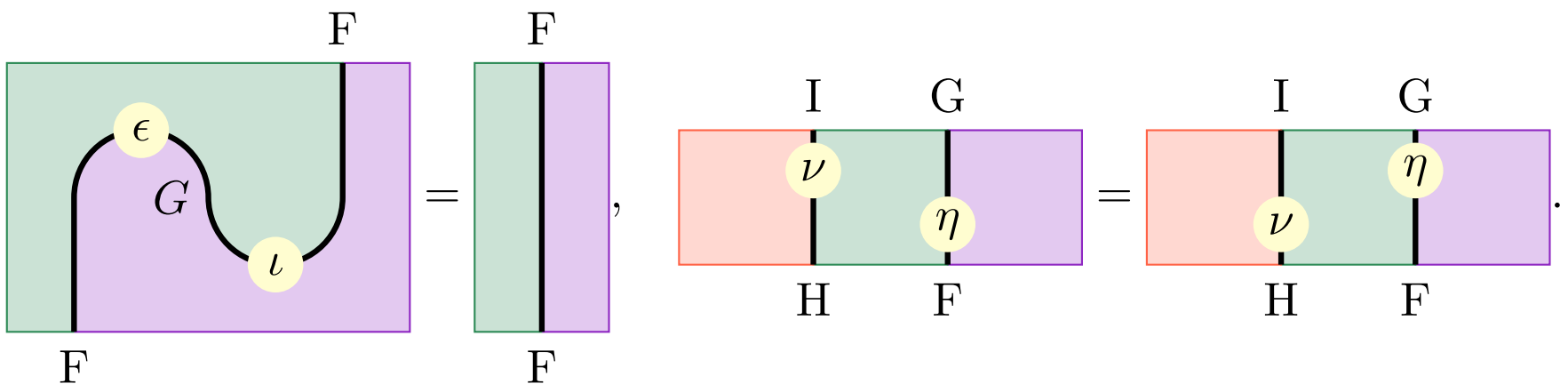

**Theorem 1G.1.** *The graphical calculus is consistent; i.e. two morphisms are equal if and only if their diagrams are related by* (1D-4)*.*

*Proof.* This can be proven, *mutatis mutandis*, as Theorem 1D.5. □

1H. **Some notions that we will need.** Up next, some category theoretical notions.

**Definition 1H.1.** Let $(\mathrm{f}\colon \mathtt{X} \to \mathtt{Y}) \in \mathbf{C}$.

(i) f is called an ***isomorphism*** if there exists a $(\mathrm{g}\colon \mathtt{Y} \to \mathtt{X}) \in \mathbf{C}$ such that

$$\mathrm{gf} = \mathrm{id}_{\mathtt{X}}, \quad \mathrm{fg} = \mathrm{id}_{\mathtt{Y}}.$$

(ii) f is called a ***monomorphism*** if it is left-cancellative, *i.e.*

$$(\mathrm{fh} = \mathrm{fi}) \Rightarrow (\mathrm{h} = \mathrm{i}) \text{ for all } \mathrm{h}, \mathrm{i} \in \mathbf{C}.$$

(iii) f is called an ***epimorphism*** if it is right-cancellative, *i.e.*

$$(\mathrm{hf} = \mathrm{if}) \Rightarrow (\mathrm{h} = \mathrm{i}) \text{ for all } \mathrm{h}, \mathrm{i} \in \mathbf{C}.$$

We also say f is ***monic*** (in (ii)) or ***epic*** (in (iii)). ◇

The following is the usual Yoga:

**Lemma 1H.2.** *If* $\mathrm{f} \in \mathbf{C}$ *is an isomorphism, then* $\mathrm{g} \in \mathbf{C}$ *as in Definition 1H.1.(i) is unique. Moreover, such an* f *is monic and epic.* □

Thus, we can just denote g as in Definition 1H.1.(i) as $\mathrm{f}^{-1}$ and call it ***inverse*** of f.

**Example 1H.3.** In a lot of categories, *e.g.* **Set** or $\mathbf{Vec}_{\Bbbk}$ the three notions in Definition 1H.1 correspond to bijective, injective and surjective morphisms, respectively. However, this is slightly misleading: all non-identity morphisms in $\mathbf{A_3}$, *cf.* Example 1B.6, are monic and epic, but none of these is an isomorphism, nor does being injective or surjective make sense. ◇

**Example 1H.4.** In **1Cob** the isomorphisms are precisely the diagrams without turnbacks, *e.g.*

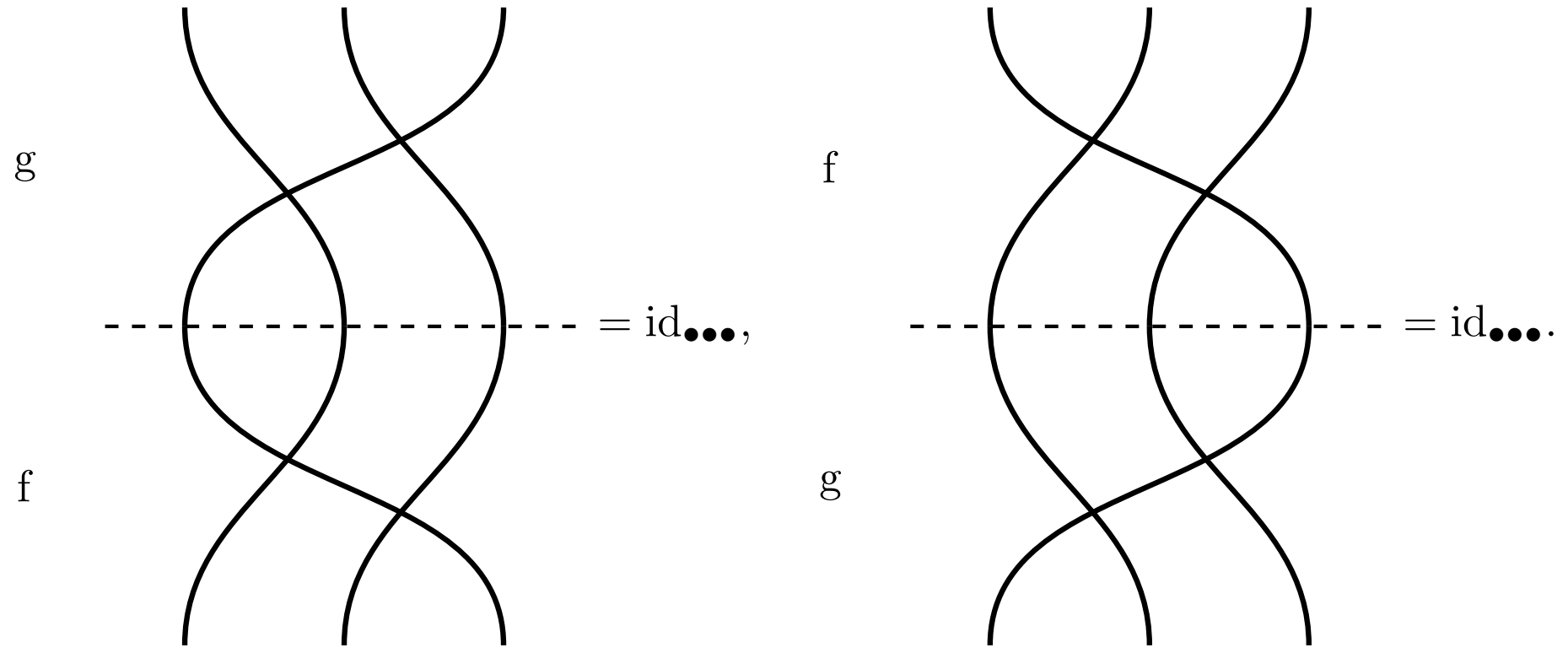

Note that the category **1Cob** is another example where monic and epic make sense, while injective and surjective have no meaning. ◇

**Definition 1H.5.** Let $\mathtt{X}, \mathtt{Y}, \mathtt{Z} \in \mathbf{C}$, and all morphisms are assumed to be in **C**.

(i) $\mathtt{X}$ and $\mathtt{Y}$ are called ***isomorphic***, denoted by $\mathtt{X} \cong \mathtt{Y}$, if there exists an isomorphism $\mathrm{f}\colon \mathtt{X} \to \mathtt{Y}$.

(ii) $\mathtt{X}$ is called a ***subobject*** of $\mathtt{Y}$, denoted by $\mathtt{X} \hookrightarrow \mathtt{Y}$, if there exists a monic morphism $\mathrm{f}\colon \mathtt{X} \to \mathtt{Y}$.

(iii) $\mathtt{Y}$ is called a ***quotient*** of $\mathtt{X}$, denoted by $\mathtt{X} \twoheadrightarrow \mathtt{Y}$, if there exists an epic morphism $\mathrm{f}\colon \mathtt{X} \to \mathtt{Y}$.

(iv) $\mathtt{X}$ is called a ***subquotient*** of $\mathtt{Z}$ if there exists $\mathtt{Y}$ and a sequence $\mathtt{X} \twoheadleftarrow \mathtt{Y} \hookrightarrow \mathtt{Z}$.

In the final case, $\mathtt{X}$ is a quotient of a subobject of $\mathtt{Z}$. ◇

Note that fixing an isomorphism $\mathrm{f}\colon \mathtt{X} \to \mathtt{Y}$ also gives us a unique isomorphism $\mathrm{f}^{-1}\colon \mathtt{Y} \to \mathtt{X}$, a fact which we will use silently throughout.

**Example 1H.6.** Note that *e.g.* being isomorphic depends on the category in which one is working. Explicitly, $\mathbb{Z}/4\mathbb{Z}$ and $\mathbb{Z}/2\mathbb{Z}\times\mathbb{Z}/2\mathbb{Z}$ are clearly isomorphic in **Set** (isomorphisms in **Set** are bijective maps, and these exist if and only if the sets are of the same size), but not in $\mathbf{Mod}(\mathbb{Z})$ since the corresponding morphisms in **Set** are not homomorphisms of abelian groups. ◇

**Example 1H.7.** For $\mathbf{C},\mathbf{D}\in\mathbf{Cat}$, by using Definition 1H.5.(a), we get the notions of two categories being isomorphic, denoted by $\mathbf{C}\cong\mathbf{D}$. ◇

**Example 1H.8.** For $\mathrm{F},\mathrm{G}\in\mathbf{Hom}(\mathbf{C},\mathbf{D})$, by using Definition 1H.5.(a), we get the notion of two functors being isomorphic. In particular, $\mathrm{F}\in\mathbf{Hom}(\mathbf{C},\mathbf{Set})$ is called ***representable*** (a particularly nice functor), if it is isomorphic to a hom functor as in Example 1E.6. ◇

As above, we will also use the notion $\cong$, $\hookrightarrow$ and $\twoheadrightarrow$ for the morphisms; *e.g.* $\mathrm{f}\colon \mathtt{X}\hookrightarrow\mathtt{Y}$ means that f is monic. The following is clear.

**Lemma 1H.9.** *The three notions $\cong$, $\hookrightarrow$ and $\twoheadrightarrow$ are reflexive and transitive, meaning e.g.*

$$\big(\mathtt{X}\hookrightarrow\mathtt{Y}\text{ and }\mathtt{Y}\hookrightarrow\mathtt{Z}\big)\Rightarrow(\mathtt{X}\hookrightarrow\mathtt{Z}),\ \text{for all }\mathtt{X},\mathtt{Y},\mathtt{Z}\in\mathbf{C},$$

*and $\cong$ is symmetric, thus, an equivalence relation.* □

**Definition 1H.10.** Let $\mathbf{C},\mathbf{D}\in\mathbf{Cat}$.

(i) $\mathbf{C}$ is called a ***subcategory*** of $\mathbf{D}$, denoted by $\mathbf{C}\subset\mathbf{D}$, if $\mathrm{Ob}(\mathbf{C})\subset\mathrm{Ob}(\mathbf{D})$, $\mathrm{Hom}_{\mathbf{C}}(\mathtt{X},\mathtt{Y})\subset\mathrm{Hom}_{\mathbf{D}}(\mathtt{X},\mathtt{Y})$ for all $\mathtt{X},\mathtt{Y}\in\mathbf{C}$, and $\mathrm{id}_{\mathtt{X}}\in\mathbf{C}$ for all $\mathtt{X}\in\mathbf{C}$.

(ii) Such a subcategory is called ***dense*** if for all $\mathtt{Y}\in\mathbf{D}$ there exists $\mathtt{X}\in\mathbf{C}$ such that $\mathtt{X}\cong\mathtt{Y}$.

(iii) Such a subcategory is called ***full*** if $\mathrm{Hom}_{\mathbf{C}}(\mathtt{X},\mathtt{Y})=\mathrm{Hom}_{\mathbf{D}}(\mathtt{X},\mathtt{Y})$ for all $\mathtt{X},\mathtt{Y}\in\mathbf{C}$.

Dense is also sometimes called ***essentially surjective***. ◇

**Example 1H.11.** We have $\mathbf{fdVec}_{\Bbbk}\subset\mathbf{Vec}_{\Bbbk}$, and $\mathbf{fdVec}_{\Bbbk}$ is full, but not dense, in $\mathbf{Vec}_{\Bbbk}$. ◇

Using Lemma 1H.9, we can define:

**Definition 1H.12.** Let $\overline{\mathrm{Ob}(\mathbf{C})/\cong}$ be a choice of representatives of $\mathrm{Ob}(\mathbf{C})/\cong$. Given a category $\mathbf{C}$, its ***skeleton*** $\mathbf{Sk}(\mathbf{C})$ is the full subcategory with objects $\overline{\mathrm{Ob}(\mathbf{C})/\cong}$. ◇

Formally, the skeleton depends on the choice of representatives. However, we can (and will) be sloppy and say that there is "the" skeleton:

**Lemma 1H.13.** *For any $\overline{\mathrm{Ob}(\mathbf{C})/\cong}$, the corresponding skeletons are isomorphic.* □

A category is called ***skeletal***, if its isomorphic to its skeleton.

**Example 1H.14.** The skeleton of $\mathbf{fdVec}_{\Bbbk}$ can be identified with $\mathbf{Mat}_{\Bbbk}$, *i.e.* $\mathbf{Sk}(\mathbf{fdVec}_{\Bbbk})\cong\mathbf{Mat}_{\Bbbk}$. Here $\mathbf{Mat}_{\Bbbk}$ is the ***category of matrices*** whose objects are natural numbers $m,n\in\mathbb{Z}_{\geq 0}$, and $\mathrm{Hom}_{\mathbf{Mat}_{\Bbbk}}(\mathtt{n},\mathtt{m})=\mathrm{Mat}_{m\times n}(\Bbbk)$, *i.e.* matrices with entries in $\Bbbk$, and $\mathbf{Mat}_{\Bbbk}$ is skeletal. ◇

**Definition 1H.15.** We let $\mathrm{K}_0(\mathbf{C})=\overline{\mathrm{Ob}(\mathbf{C})/\cong}$ , and call it the ***Grothendieck class of*** $\mathbf{C}$. Elements in $\mathrm{K}_0(\mathbf{C})$ are ***Grothendieck classes*** of $\mathtt{X}\in\mathbf{C}$ and are denoted by $[\mathtt{X}]$. ◇

*Remark* 1H.16. In the literature the notation $\mathrm{K}_0$ is also used to denote the Grothendieck group, where one (artificially) adds formal inverses. For us $\mathrm{K}_0$ is not a group, so we do not have formal inverses. ◇

We think of $\mathrm{K}_0(\mathbf{C})$ as capturing all the information about the objects of $\mathbf{C}$. For an arbitrary category, $\mathrm{K}_0(\mathbf{C})$ is just a set, but when $\mathbf{C}$ has more structure, then so does $\mathrm{K}_0(\mathbf{C})$.

**Example 1H.17.** We can identify $\mathrm{K}_0(\mathbf{fdVec}_{\Bbbk})\xrightarrow{\cong}\mathbb{Z}_{\geq 0}$ as sets, the map being $[\Bbbk^n]\mapsto n$, since any $\mathtt{X}\in\mathbf{fdVec}_{\Bbbk}$ is isomorphic to $\Bbbk^n$ for some $n\in\mathbb{Z}_{\geq 0}$. ◇

Note that $\mathbf{fdVec}_{\Bbbk}$ and $\mathbf{Mat}_{\Bbbk}$ in some sense contain the same information, but they are not isomorphic: $\mathbf{fdVec}_{\Bbbk}\ncong\mathbf{Mat}_{\Bbbk}$. This is due to the fact that isomorphisms of categories give bijections on objects. But recall that we "do not care about objects". So maybe asking

$$\mathrm{Id}_{\mathbf{C}}=\mathrm{GF},\quad \mathrm{FG}=\mathrm{Id}_{\mathbf{D}}$$

is a bit too much. This motivates the following "correct" notion of equivalence of categories:

**Definition 1H.18.** Let $\mathbf{C},\mathbf{D}\in\mathbf{Cat}$. The categories $\mathbf{C}$ and $\mathbf{D}$ are called ***equivalent***, denoted by $\mathbf{C}\simeq\mathbf{D}$, if there exists $\mathrm{F}\colon\mathbf{C}\to\mathbf{D}$ and $\mathrm{G}\colon\mathbf{D}\to\mathbf{C}$ such that

$$\mathrm{Id}_{\mathbf{C}}\cong\mathrm{GF},\quad \mathrm{FG}\cong\mathrm{Id}_{\mathbf{D}},$$

where $\cong$ is taken in $\mathbf{Hom}(\mathbf{C},\mathbf{D})$, *cf.* Example 1H.8 ◇

*Remark* 1H.19. In Feynman diagrammatics for $\mathbf{Hom}(\mathbf{C},\mathbf{D})$ there is a nice interpretation of equivalence. To this end, let us fix

$$\iota\colon \mathrm{Id}_{\mathbf{C}} \xrightarrow{\cong} \mathrm{GF}, \quad \varepsilon\colon \mathrm{FG} \xrightarrow{\cong} \mathrm{Id}_{\mathbf{D}},$$

sometimes also called ***unit*** and ***counit***. Then these can be pictured as caps and cups

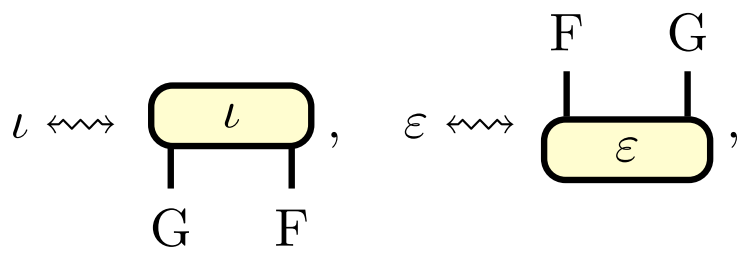

where we have not drawn strands for the identity functors, by the usual convention. Later, with more structure at hand, we will revisit such diagrams, which then become topological objects. ◇

Functors as in Definition 1H.18 are called ***equivalences*** and they are ***quasi-inverse*** to each other. Clearly, isomorphic categories are equivalent, but the converse is not true:

**Example 1H.20.** Any category $\mathbf{C}$ is equivalent to its skeleton, but not necessarily isomorphic. Explicitly, $\mathbf{fdVec}_{\Bbbk} \simeq \mathbf{Mat}_{\Bbbk}$, but $\mathbf{fdVec}_{\Bbbk} \not\cong \mathbf{Mat}_{\Bbbk}$. ◇

**Example 1H.21.** The category $\mathbf{fSet}$, which is the full subcategory of $\mathbf{Set}$ with objects being finite sets, is not skeletal. ◇

**Example 1H.22.** The category $\mathbf{fSet}^{\cong}$, which is the subcategory of $\mathbf{fSet}$ with the same objects, but only bijections, is not skeletal. Its skeleton is $\mathbf{Sym}^{top}$ which is the subcategory of $\mathbf{1Cob}$, *cf.* Example 1B.9.(a), with the same objects but only cobordisms without Morse points (a.k.a. no turnbacks). A typical diagram in $\mathbf{Sym}^{top}$ is a permutation diagram, *e.g.*

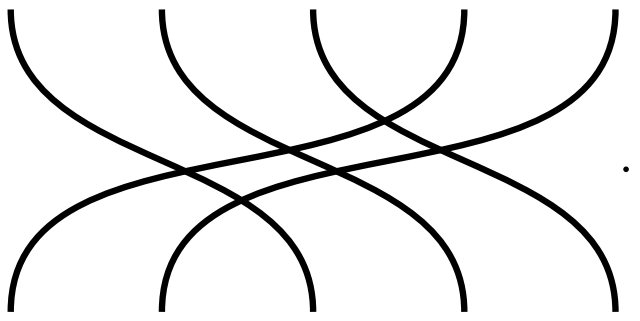

Note that $\mathrm{Hom}_{\mathbf{Sym}^{top}}(\bullet^{m},\bullet^{n}) = \emptyset$ unless $m = n$. ◇

Let us also note:

**Lemma 1H.23.** *Any functor* $\mathrm{F} \in \mathbf{Hom}(\mathbf{C},\mathbf{D})$ *induces a map*

$$K_0(\mathrm{F})\colon K_0(\mathbf{C}) \to K_0(\mathbf{D}),\ [\mathtt{X}] \mapsto [\mathrm{F}(\mathtt{X})].$$

*Further, if* F *is an equivalence, then* $K_0(\mathrm{F})$ *is an isomorphism.* □

If one wants to check whether two categories are equivalent one almost always uses:

**Proposition 1H.24.** *A functor* $\mathrm{F}\colon \mathbf{C} \to \mathbf{D}$ *is an equivalence if and only if*

- *it is* **dense** *(also called* **essentially surjective***), i.e.*

$$\text{for all } \mathtt{Y} \in \mathbf{D} \text{ there exist } \mathtt{X} \in \mathbf{C} \text{ such that } \mathrm{F}(\mathtt{X}) \cong \mathtt{Y};$$

- *it is* **faithful***, i.e.*

$$\mathrm{Hom}_{\mathbf{C}}(\mathtt{X},\mathtt{Y}) \hookrightarrow \mathrm{Hom}_{\mathbf{D}}\big(\mathrm{F}(\mathtt{X}),\mathrm{F}(\mathtt{Y})\big) \text{ for all } \mathtt{X},\mathtt{Y} \in \mathbf{C};$$

- *it is* **full***, i.e.*

$$\mathrm{Hom}_{\mathbf{C}}(\mathtt{X},\mathtt{Y}) \twoheadrightarrow \mathrm{Hom}_{\mathbf{D}}\big(\mathrm{F}(\mathtt{X}),\mathrm{F}(\mathtt{Y})\big) \text{ for all } \mathtt{X},\mathtt{Y} \in \mathbf{C}.$$

If a functor is full and faithful, then we also say its ***fully faithful.***

*Proof.* The proof is what is called ***diagram chasing***.

⇒. Let $(\mathrm{F},\mathrm{G},\iota,\varepsilon)$ as in Remark 1H.19 define the equivalence. By $\varepsilon_{\mathtt{X}}\colon \mathrm{FG}(\mathtt{Y}) \xrightarrow{\cong} \mathtt{Y}$ we see that F is dense. To see that F is faithful consider the commuting diagram

$$\begin{array}{ccc} \mathtt{X} & \xrightarrow[\iota_{\mathtt{X}}]{\cong} & \mathrm{GF}(\mathtt{X}) \\ {\scriptstyle \mathrm{f} \text{ or } \mathrm{g}}\downarrow & & \downarrow{\scriptstyle \mathrm{GF}(\mathrm{f}) \text{ or } \mathrm{GF}(\mathrm{g})} \\ \mathtt{X}' & \xrightarrow[\cong]{\iota_{\mathtt{X}'}} & \mathrm{GF}(\mathtt{X}') \end{array}.$$

Assuming that $\mathrm{GF}(\mathrm{f}) = \mathrm{GF}(\mathrm{g})$, by Exercise 1K.2, implies that $\mathrm{f} = \mathrm{g}$ which in turn implies that F is faithful. Very similar arguments, using again Exercise 1K.2, show that F is full.

$\Leftarrow$. Suppose that $\mathrm{F}$ is dense and fully faithful, so we need to construct the quadruple $(\mathrm{F},\mathrm{G},\iota,\varepsilon)$ as in Remark 1H.19. First, using density, we find an object $\mathrm{G}(\mathtt{Y})$ for all $\mathtt{Y}\in\mathbf{D}$ as well as an isomorphism $\varepsilon_{\mathtt{Y}}\colon \mathrm{FG}(\mathtt{Y})\xrightarrow{\cong}\mathtt{Y}$. Thus, for each $\mathrm{f}\colon\mathtt{Y}\to\mathtt{Y}'$ we find a unique solution $\mathrm{FG}(\mathrm{f})$ to make

$$\begin{array}{ccc}
\mathrm{FG}(\mathtt{Y}) & \xrightarrow[\varepsilon_{\mathtt{Y}}]{\cong} & \mathtt{Y} \\
{\scriptstyle \mathrm{FG}(\mathrm{f})}\downarrow & & \downarrow{\scriptstyle \mathrm{f}} \\
\mathrm{FG}(\mathtt{Y}') & \xrightarrow[\cong]{\varepsilon_{\mathtt{Y}'}} & \mathtt{Y}'
\end{array}$$

commutative, by Exercise 1K.2. Hence, fully faithfulness of $\mathrm{F}$ defines us $\mathrm{G}(\mathrm{f})$. Scrutiny of this construction actually show that $\mathrm{G}(\mathtt{Y})$ and $\mathrm{G}(\mathrm{f})$, and $\varepsilon_{\mathtt{Y}}$ assemble into a functor and a natural transformation, respectively. It remains to construct $\iota_{\mathtt{X}}$ (and prove that these give rise to a natural transformation), which can be done in a similar fashion. ☐

**Definition 1H.25.** A category $\mathbf{C}$ is called ***concrete*** if it admits a faithful functor, called its ***realization***, $\mathrm{R}\colon\mathbf{C}\to\mathbf{Set}$. ◇

**Example 1H.26.** The functor Forget, *cf.* Example 1E.5.(d), realizes $\mathbf{Vec}_{\Bbbk}$ as a concrete category. ◇

The following is arguably the most important statement in classical category theory, known as the ***Yoneda lemma***. We will not need it, and only give a reference for its proof, but any text on category theory without it feels like "missing something". So here we go:

**Theorem 1H.27.** *For any* $\mathrm{F}\in\mathbf{Hom}(\mathbf{C},\mathbf{Set})$ *and any* $\mathtt{X}\in\mathbf{C}$ *there is a bijection*

$$\mathrm{Hom}_{\mathbf{Hom}(\mathbf{C},\mathbf{Set})}\big(\mathrm{Hom}_{\mathbf{C}}(\mathtt{X},_),\mathrm{F}\big)\to\mathrm{F}(\mathtt{X}),\quad \big(\alpha\colon\mathrm{Hom}_{\mathbf{C}}(\mathtt{X},_)\Rightarrow\mathrm{F}\big)\mapsto\alpha_{\mathtt{X}}(\mathrm{id}_{\mathtt{X}}).$$

*Moreover, this correspondence is natural in both* $\mathrm{F}$ *and* $\mathtt{X}$.

*Proof.* Proofs tend to be a bit technical and longish. We do not need the Yoneda lemma much, so we refer to [ML98, Section III.2]. ☐

As a consequence, we have the ***Yoneda embedding(s)*** given by the ***Yoneda functor(s)***:

**Proposition 1H.28.** *Fix* $\mathbf{C}\in\mathbf{Cat}$. *We have fully faithful functors*

$$\begin{cases}\mathrm{Y}\colon\mathbf{C}\to\mathbf{Hom}(\mathbf{C}^{op},\mathbf{Set}),\\ \mathtt{X}\mapsto\mathrm{Hom}_{\mathbf{C}}(_,\mathtt{X}),(\mathrm{f}\colon\mathtt{X}\to\mathtt{Y})\mapsto\big(\mathrm{f}\circ_\colon\mathrm{Hom}_{\mathbf{C}}(_,\mathtt{X})\to\mathrm{Hom}_{\mathbf{C}}(_,\mathtt{Y}),\mathrm{g}\mapsto\mathrm{fg}\big),\end{cases}$$

$$\begin{cases}\mathrm{Y}^{op}\colon\mathbf{C}^{op}\to\mathbf{Hom}(\mathbf{C},\mathbf{Set}),\\ \mathtt{X}\mapsto\mathrm{Hom}_{\mathbf{C}}(\mathtt{X},_),(\mathrm{f}\colon\mathtt{X}\to\mathtt{Y})^{op}\mapsto\big(_\circ\mathrm{f}\colon\mathrm{Hom}_{\mathbf{C}}(\mathtt{Y},_)\to\mathrm{Hom}_{\mathbf{C}}(\mathtt{X},_),\mathrm{g}\mapsto\mathrm{gf}\big).\end{cases}$$

*Hence,* $\mathbf{C}$ *and* $\mathbf{C}^{op}$ *are full subcategories of* $\mathbf{Hom}(\mathbf{C}^{op},\mathbf{Set})$ *respectively of* $\mathbf{Hom}(\mathbf{C},\mathbf{Set})$.

What is this saying? It is a ***gold foil experiment***, see Figure 3, in mathematics. Idea (reinterpreted): "Shoot on $\mathtt{X}$ to learn more about $\mathtt{X}$".

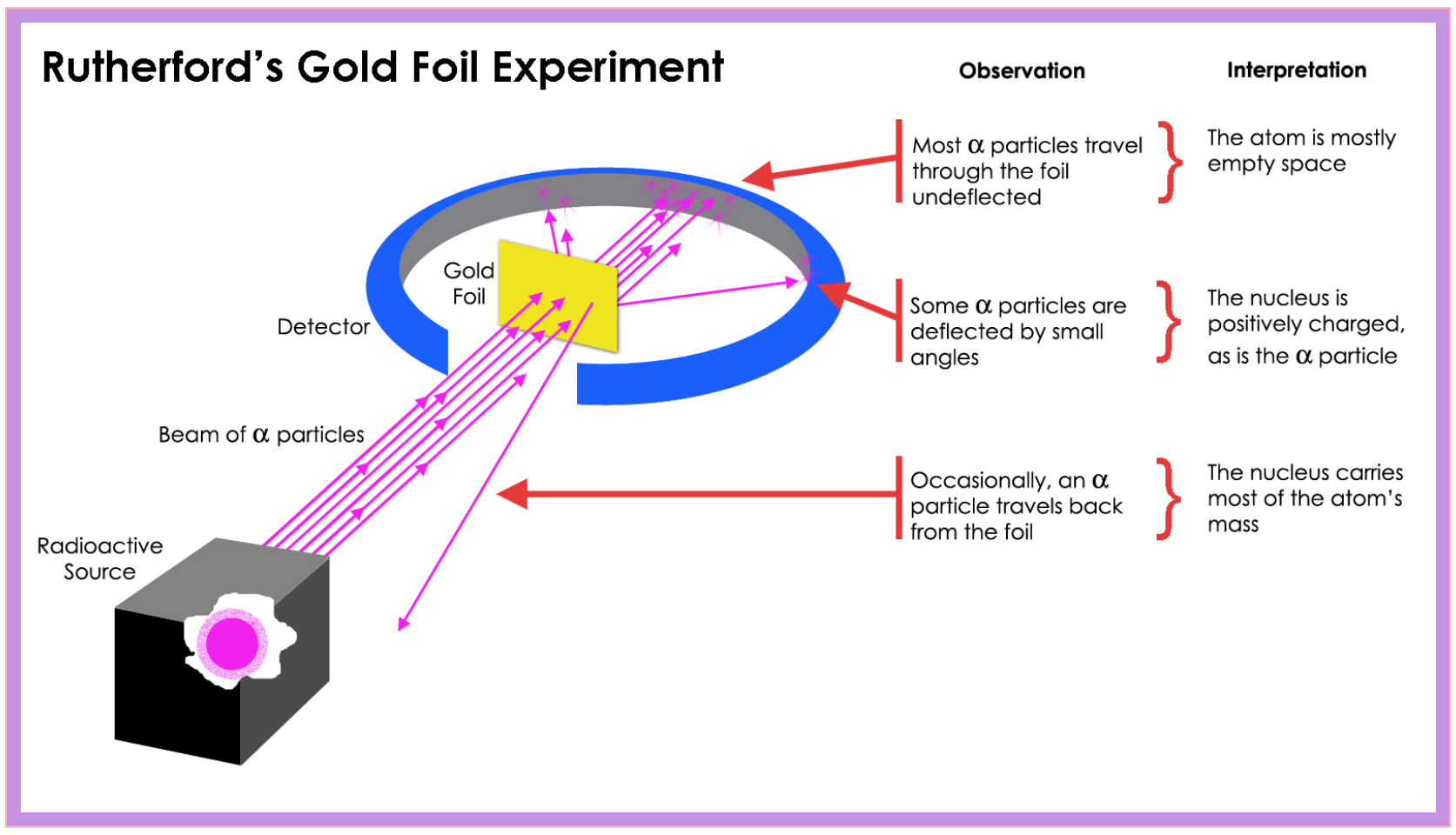

FIGURE 3. A picture of Geiger–Marsden (Rutherford's gold foil) landmark experiments from ∼1910. See, for example, https://en.wikipedia.org/wiki/Rutherford_scattering_experiments.

Picture from https://en.wikipedia.org/wiki/Rutherford_scattering_experiments

*Proof of Proposition 1H.28.* From the construction of the Yoneda functors we see that we have injections

$$\begin{aligned}\mathrm{Hom}_{\mathbf{C}}(\mathtt{X},\mathtt{Y}) &\hookrightarrow \mathrm{Hom}_{\mathbf{Hom}(\mathbf{C},\mathbf{Set})}\big(\mathrm{Hom}_{\mathbf{C}}(\mathtt{X},_),\mathrm{Hom}_{\mathbf{C}}(\mathtt{Y},_)\big),\\ \mathrm{Hom}_{\mathbf{C}}(\mathtt{X},\mathtt{Y}) &\hookrightarrow \mathrm{Hom}_{\mathbf{Hom}(\mathbf{C},\mathbf{Set})}\big(\mathrm{Hom}_{\mathbf{C}}(_,\mathtt{X}),\mathrm{Hom}_{\mathbf{C}}(_,\mathtt{Y})\big).\end{aligned} \tag{1H-29}$$

Furthermore, Theorem 1H.27 implies that every natural transformation between represented functors arises in this way, showing that (1H-29) are bijections. Comparing this to the second and third bullet points in Proposition 1H.24, which define the notion of being fully faithful, shows the claim. □

**Example 1H.30.** For the category $\mathbf{A_3}$ from Example 1B.6 the Yoneda functor $\mathrm{Y}^{op}$ associates

$$\mathrm{Y}(\mathtt{1}) = \left(\mathrm{Hom}_{\mathbf{A_3}}(\mathtt{1},_)\colon \mathbf{A_3}\to\mathbf{Set}, \begin{cases}\mathtt{1}\mapsto\{\mathrm{id}_1\},\mathtt{2}\mapsto\{\mathrm{f}\},\mathtt{3}\mapsto\{\mathrm{gf}\}\\ \mathrm{f}\mapsto(\mathrm{id}_1\mapsto\mathrm{f}),\mathrm{g}\mapsto(\mathrm{id}_2\mapsto\mathrm{g}),\mathrm{gf}\mapsto(\mathrm{id}_1\mapsto\mathrm{gf}),\end{cases}\right),$$

*etc.*, which identifies $(\mathbf{A_3})^{op}$ with the functors of the form $\mathrm{Hom}_{\mathbf{A_3}}(\mathtt{i},_)$ for $\mathtt{i}\in\{\mathtt{1},\mathtt{2},\mathtt{3}\}$. ◇

1I. **Adjoint functors.** Consider the ceiling and floor functions:

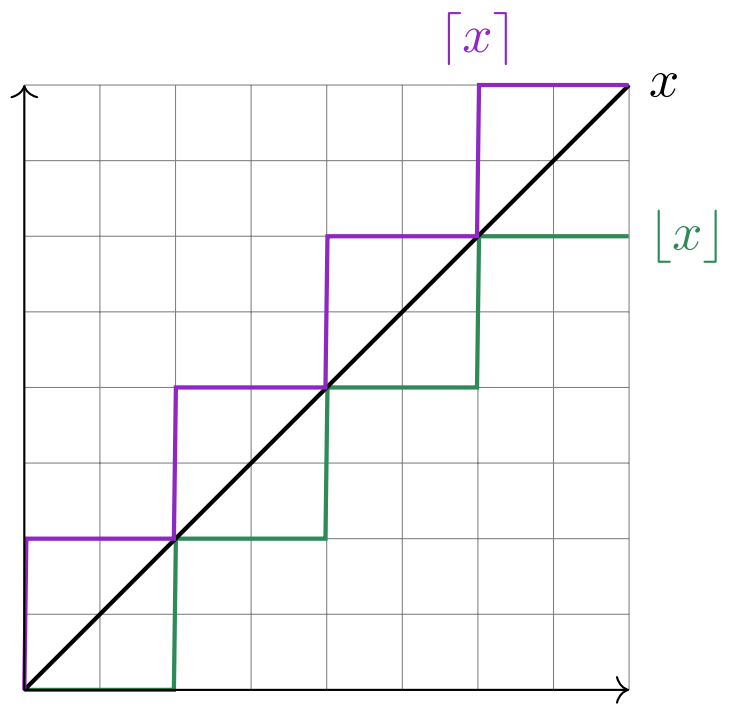

.

The inclusion $\iota\colon\mathbb{Z}\hookrightarrow\mathbb{R}$ is not invertible and $\mathbb{Z}\not\cong\mathbb{R}$. In fact, there is no invertible way to assign an integer to a real number. However, ceiling and floor serve as approximations of inverses, or ***pseudo inverses***. Adjoint functors are a way to category theoretically catch the notion of pseudo inverse.

**Definition 1I.1.** Two functors $(\mathrm{F},\mathrm{G})=(\mathrm{F}\colon\mathbf{C}\to\mathbf{D},\mathrm{G}\colon\mathbf{D}\to\mathbf{C})$ form an ***adjoint pair*** if:

**(a)** There exists a natural transformation $\alpha\colon\mathrm{Hom}_{\mathbf{D}}\big(\mathrm{F}(_),_\big)\Rightarrow\mathrm{Hom}_{\mathbf{C}}\big(_,\mathrm{G}(_)\big)$.

**(b)** For all $\mathtt{X},\mathtt{Y}$ there are isomorphism

$$\alpha_{\mathtt{X},\mathtt{Y}}\colon\mathrm{Hom}_{\mathbf{D}}\big(\mathrm{F}(\mathtt{X}),\mathtt{Y}\big)\xrightarrow{\cong}\mathrm{Hom}_{\mathbf{C}}\big(\mathtt{X},\mathrm{G}(\mathtt{Y})\big).$$

In this case F is the ***left adjoint*** of G, and G is the ***right adjoint*** of F ◇

A functor might not have left/right adjoints. But the usual yoga shows that we can say the instead of a left/right adjoint:

**Lemma 1I.2.** *If a left/right adjoint exists, then it is unique up to unique isomorphism.* □

**Example 1I.3.** Here is a list of important adjoint functors (F, G) with $\mathrm{F}\colon\mathbf{C}\rightleftarrows\mathbf{D}\colon\mathrm{G}$.

**(a)** The ***free-forget adjunction*** (F = Free, G = Forget) as in this table:

| **D** | $\mathbf{Vec}_{\Bbbk}$ | **Monoid** | **Group** | $\mathbf{Mod}(\mathbb{Z})$ | **Ring** | $\mathbf{Alg}_{\Bbbk}$ | **Cat** | **Field** |
|---|---|---|---|---|---|---|---|---|
| **C** | **Set** | **Set** | **Set** | **Set** | **Set** | **Set** | **Quiver** | **Domain** |

.

And indeed, Free is the pseudo inverse of Forget.

The rightmost free-forget adjunction also goes under the name "Field of fractions." However, the forgetful functor Forget: **Field** → **Set** does not have a left adjoint, and indeed there is no "free field."

**(b)** $\mathbf{C}=\mathbf{Vec}_{\mathbb{R}}$, $\mathbf{D}=\mathbf{Vec}_{\mathbb{C}}$, F=Scalar extension, G = Forget (one can vary $\mathbb{R}$, $\mathbb{C}$).

**(c)** **C** = **Ring**, **D** = **pRing**, F=Polynomial ring, G = Forget.

**(d)** The ***tensor-hom adjunction***, a.k.a. currying, given by

$$\mathrm{Hom}_{\mathbf{C}}(\mathtt{X}\otimes_{\mathbf{D}}\mathtt{Y},\mathtt{Z})\cong\mathrm{Hom}_{\mathbf{D}}\big(\mathtt{X},\mathrm{Hom}_{\mathbf{C}}(\mathtt{Y},\mathtt{Z})\big).$$

Here, tensoring is the left adjoint.

**(e)** **C** = **Group**, **D** = **Ring**, F=Group ring, G=Group of units.

**(f)** **C** = **Group**, **D** = $\mathbf{Mod}(\mathbb{Z})$, F=Abelianization, G=Include.

**(g)** Any equivalence is an adjoint pair. The converse is far from true, and a very convincing example that this is not supposed to be the case is Example 1I.4 below.

**(h)** And there are many more!

(We leave it to the reader to define the undefined notation.) ◇

**Example 1I.4.** Consider $\mathbb{R}$ as a category with $x \to y$ if $x \leq y$. Similarly for $\mathbb{Z} \subset \mathbb{R}$, and consider the inclusion $\iota\colon \mathbb{Z} \to \mathbb{R}$. Then we have

$$\lceil y \rceil \leq x \Leftrightarrow y \leq x, \quad x \leq \lfloor y \rfloor \Leftrightarrow x \leq y, \quad x \in \mathbb{Z}, y \in \mathbb{R},$$

which in turn gives adjoint pairs $(\lceil _ \rceil, \iota)$ and $(\iota, \lfloor _ \rfloor)$. ◇

**Example 1I.5.** Consider the category of topological spaces and continuous maps **Top** and:

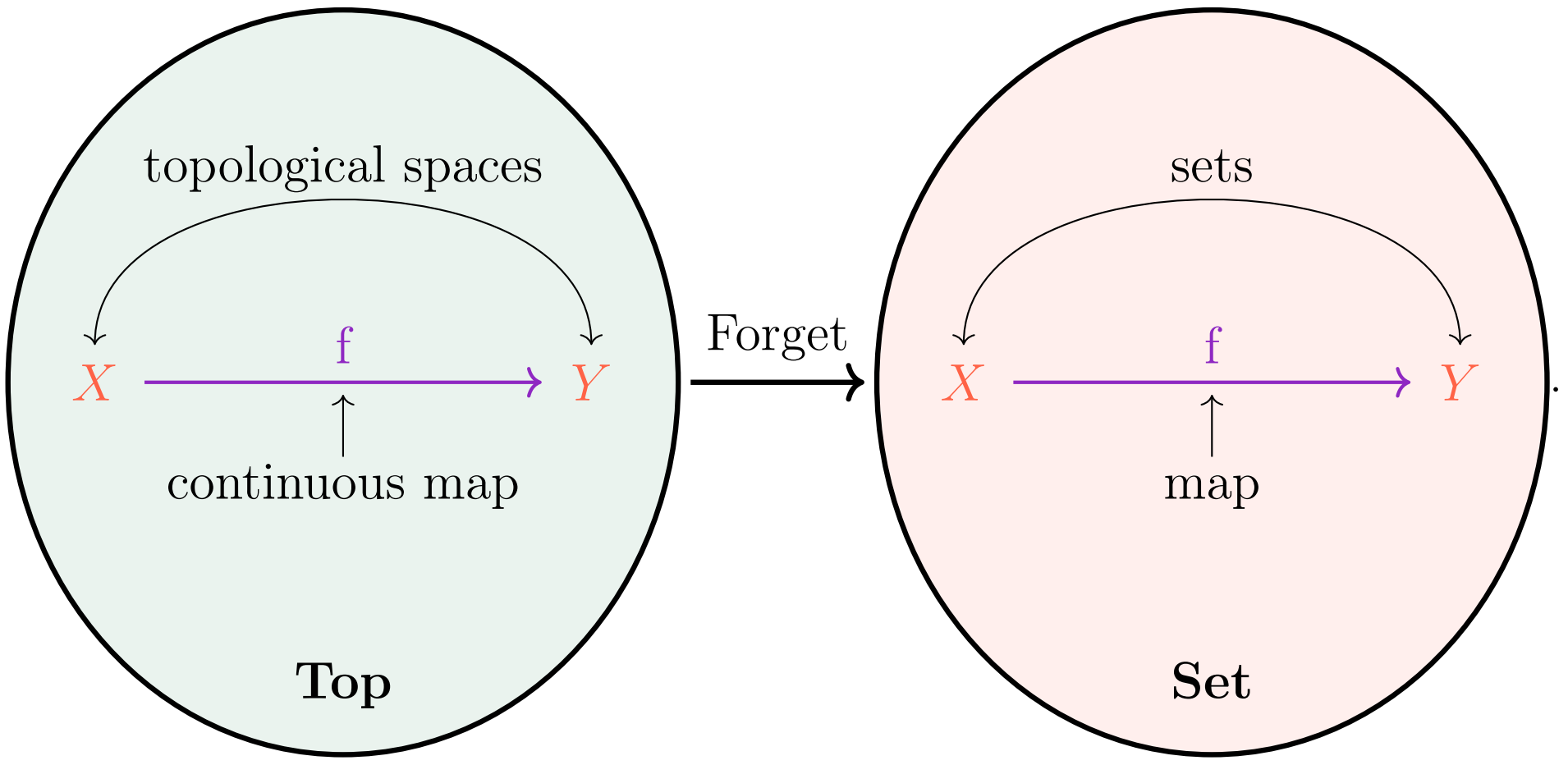

The functor Forget: **Top** $\to$ **Set** has a left and a right adjoint. The left adjoint uses the discrete topology (free topology), the right adjoint uses the indiscrete topology (cofree topology). ◇

**Example 1I.6.** Back to cobordisms:

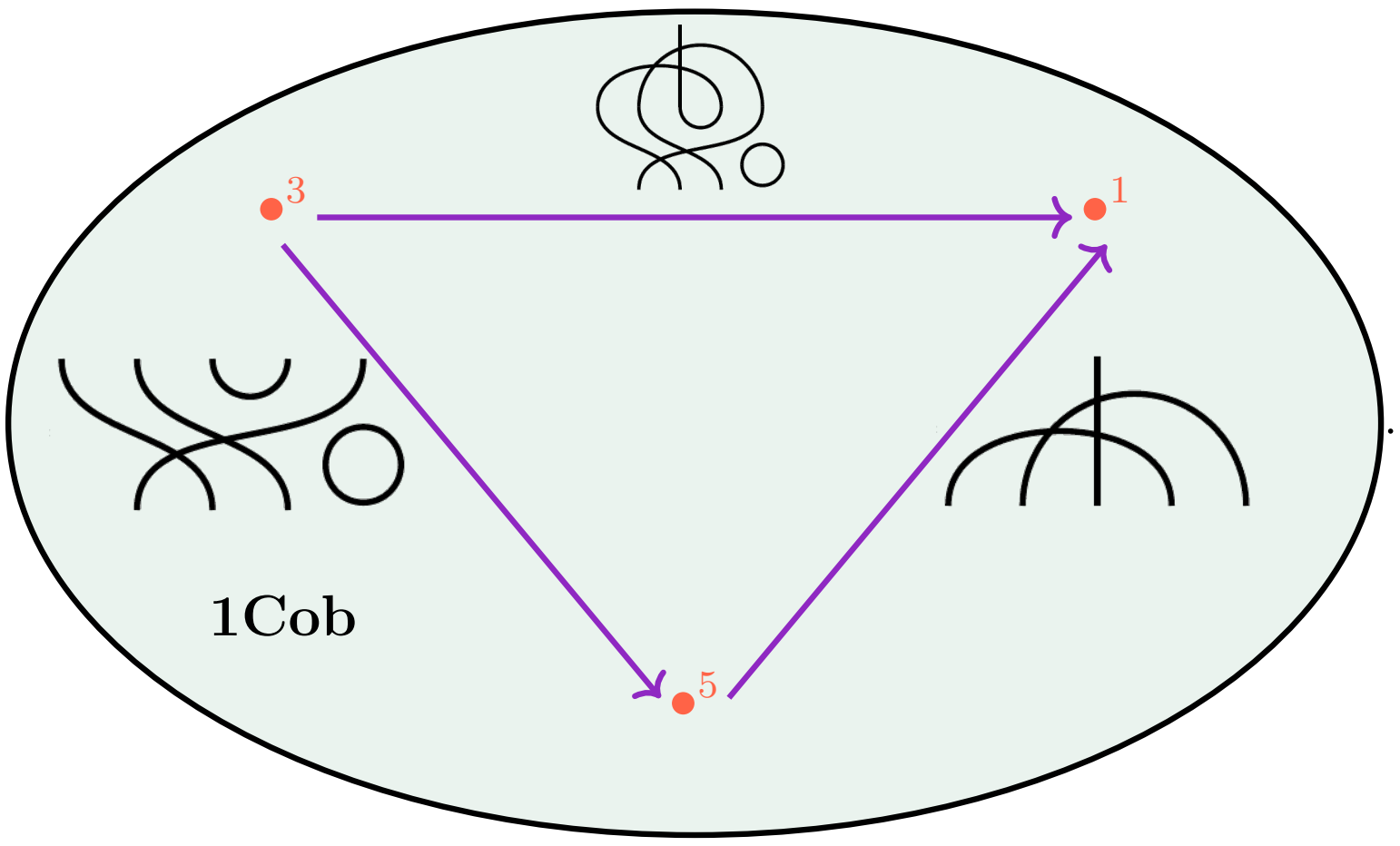

Let **PivSym** be the category of pivotal symmetric monoidal categories, with details to come later on. Then Forget: **PivSym** $\to$ **Set** has a left adjoint, and **1Cob** arises in this way. ◇

Let us now redefine adjunctions using colored Feynman diagrams.

**Definition 1I.7.** Two functors $(\mathrm{F}, \mathrm{G}) = (\mathrm{F}\colon \mathbf{C} \to \mathbf{D}, \mathrm{G}\colon \mathbf{D} \to \mathbf{C})$ form an ***adjoint pair*** if:

**(a)** We draw them as

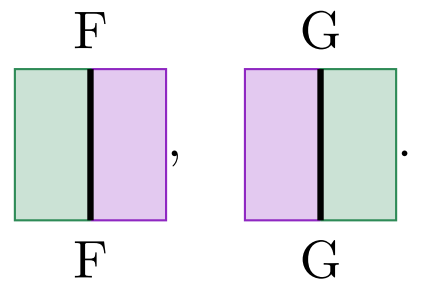

**(b)** We assume there exist natural transformations $\epsilon\colon \mathrm{FG} \Rightarrow \mathrm{Id}_{\mathbf{D}}$ and $\iota\colon \mathrm{Id}_{\mathbf{C}} \Rightarrow \mathrm{GF}$

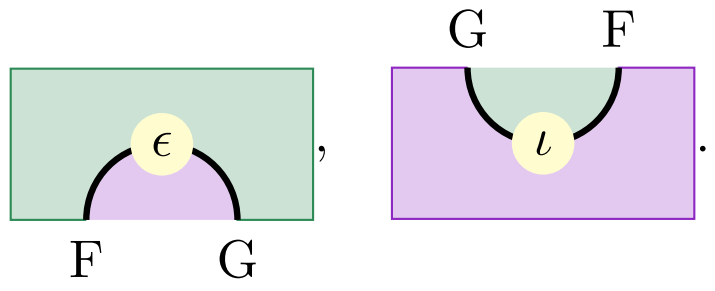

These are called ***counit*** and ***unit***, respectively.

**(c)** We assume the ***zigzag relations***, *i.e.*:

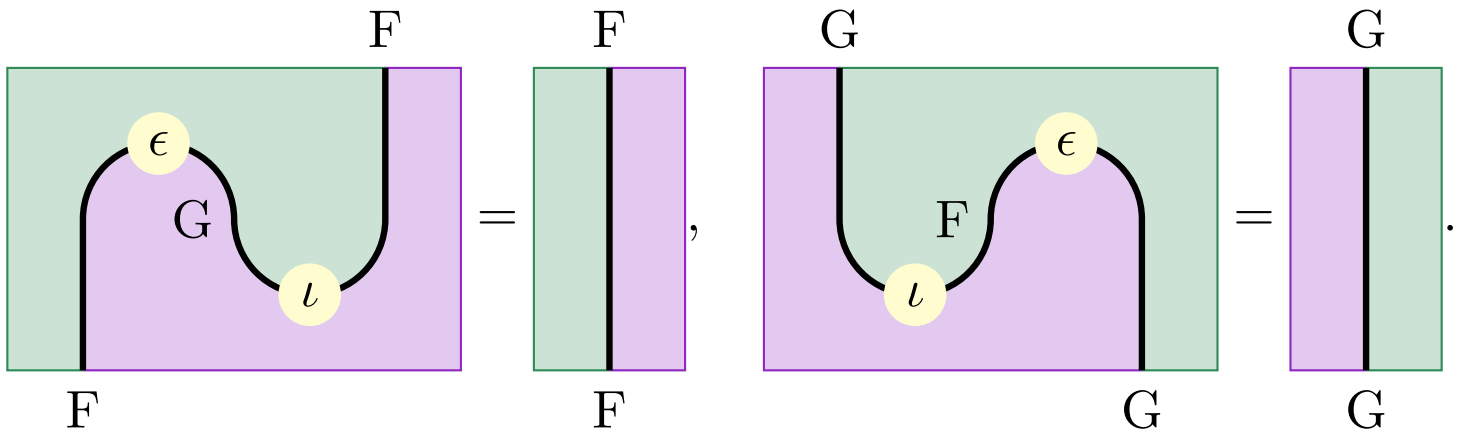

In this case F is the ***left adjoint*** of G, and G is the ***right adjoint*** of F ◇

**Proposition 1I.8.** *Definition 1I.1 and Definition 1I.7 are equivalent.*

*Proof. Definition 1I.1 implies Definition 1I.7.* Given the natural isomorphism $\alpha_{\mathtt{X},\mathtt{Y}}\colon \mathrm{Hom}_{\mathbf{D}}\big(\mathrm{F}(\mathtt{X}),\mathtt{Y}\big) \xrightarrow{\cong} \mathrm{Hom}_{\mathbf{C}}\big(\mathtt{X},\mathrm{G}(\mathtt{Y})\big)$, we can define the unit and counit as follows: For each $\mathtt{X}\in\mathbf{C}$, the morphism $\iota_{\mathtt{X}}\colon \mathtt{X}\to\mathrm{GF}(\mathtt{X})$ corresponds to the identity morphism $\mathrm{Id}_{\mathrm{F}(\mathtt{X})}$ under $\alpha_{\mathtt{X},\mathrm{F}(\mathtt{X})}$. Thus, the unit can be defined using the $\iota_{\mathtt{X}}$. Dually, for each $\mathtt{Y}\in\mathbf{D}$, the morphism $\epsilon_{\mathtt{Y}}\colon \mathrm{GF}(\mathtt{Y})\to\mathtt{Y}$ corresponds to the identity morphism $\mathrm{Id}_{\mathrm{G}(\mathtt{Y})}$ under $\alpha_{\mathrm{G}(\mathtt{Y}),\mathtt{Y}}$. Hence, the counit can be defined using the $\epsilon_{\mathtt{Y}}$. It is then easy to see that the zigzag relation holds.

*Definition 1I.7 implies Definition 1I.1.* We use

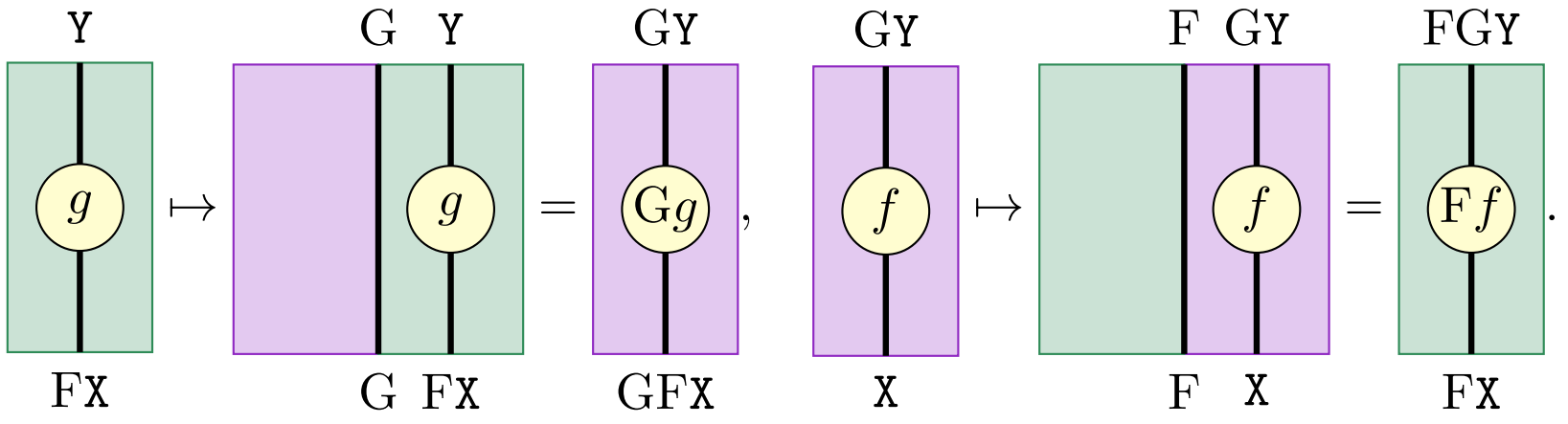

(We omitted several brackets to improve readability.) Scrutiny shows that these give the desired isomorphism $\mathrm{Hom}_{\mathbf{D}}\big(\mathrm{F}(\mathtt{X}),\mathtt{Y}\big) \xrightarrow{\cong} \mathrm{Hom}_{\mathbf{C}}\big(\mathtt{X},\mathrm{G}(\mathtt{Y})\big)$. □

1J. **Categorification? Sure!** Categorification is the process of replacing set theoretic or numerical concepts with higher level categorical structures to gain deeper insights. It involves replacing sets with categories, functions with functors, and equations with natural isomorphisms. Categorification has applications mostly in representation theory, topology, and mathematical physics, providing more refined invariants and structures. In some sense, categorification is the object and what we have categorified is its shadow.

For example, categorifying a number might yield a vector space, with its dimension recovering the original number. Or, from a birds-eye point of view, replace all of $\mathbb{Z}_{\geq 0}$ with the category $\mathbf{Vec}_{\Bbbk}$ as in Figure 4. It is hereby crucial that the categorification reflects and extends the set-theoretic or numerical concepts one starts with.

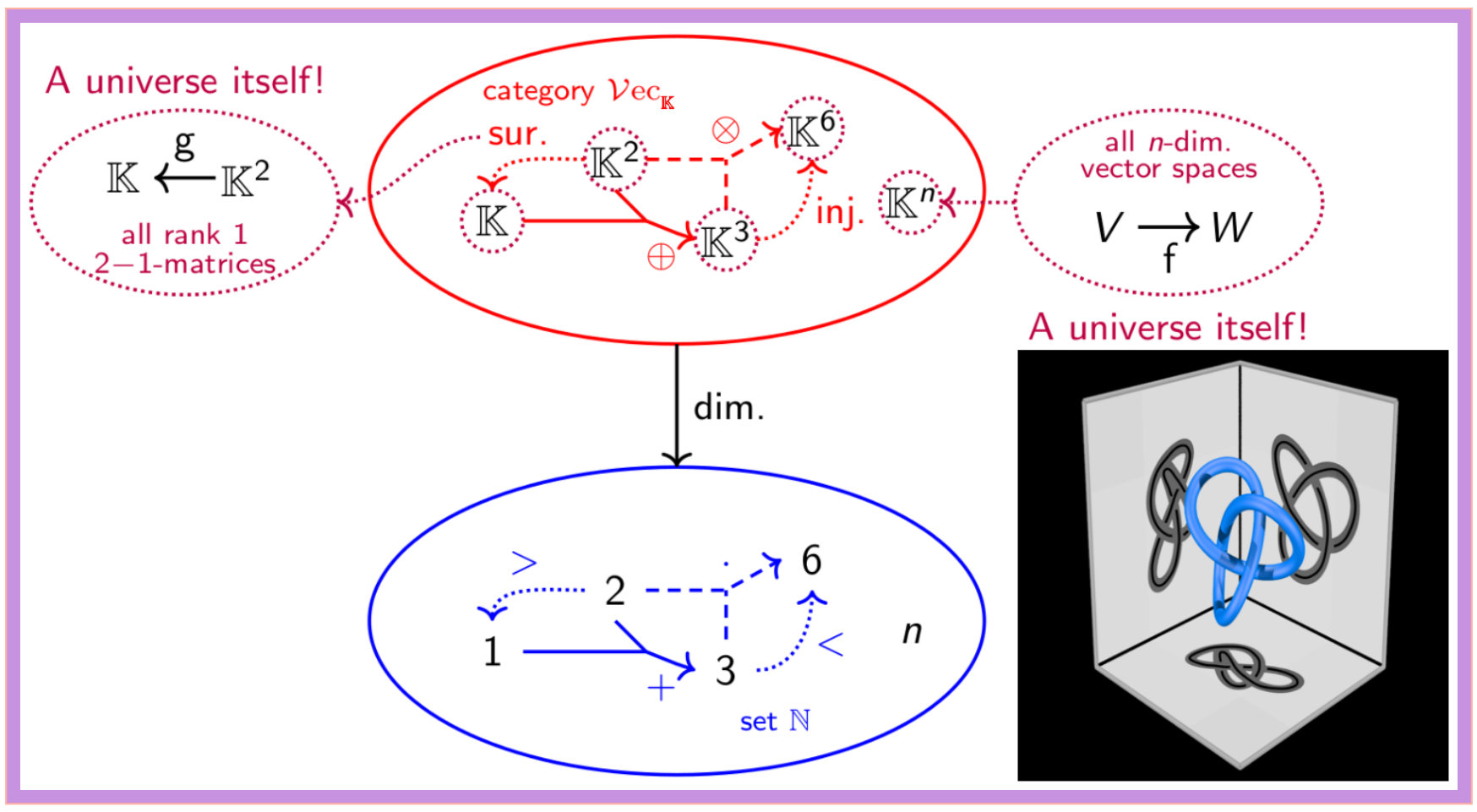

FIGURE 4. The category $\mathbf{Vec}_{\Bbbk}$ categorifies $\mathbb{Z}_{\geq 0}$ and a lot of its structures, like addition, multiplication, *etc.*

A goal of these notes is to categorify the following notions:

| Notion | Categorification |
|---|---|
| Set | Category |
| Monoid | Monoidal category |
| Involutive monoid | Pivotal category |
| Commutative monoid | Braided category |
| Abelian group | Additive/abelian category |
| Ring/algebra | Fiat/tensor category |
| Group ring | Fusion category |

So far we have only seen what a category is; the remaining notions will follow in later chapters.

## 1K. **Exercises.**

*Exercise* 1K.1. Consider the following diagram in some category.

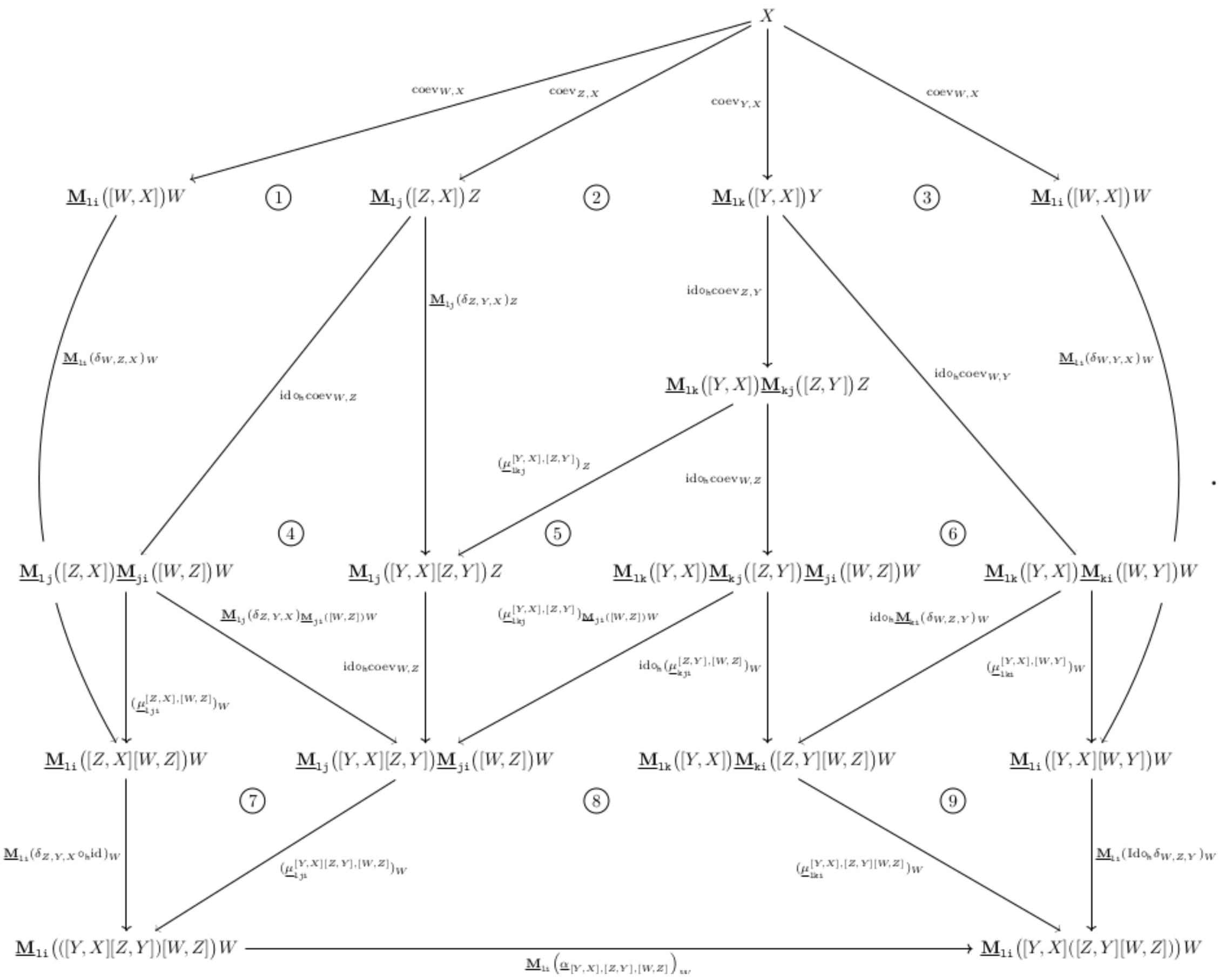

If all the numbered subdiagrams commute, does it follow that the diagram itself is commutative? ◇

*Exercise* 1K.2. Given $f\colon X \to Y$ and fixed isomorphisms $X \cong X'$ and $Y \cong Y'$, there exists a unique $f'\colon X' \to Y'$ such that any, or, equivalently, all, of the following diagrams commute:

$$\begin{CD} X @>{\cong}>> Y \\ @V{f}VV @VV{f'}V \\ X' @>>{\cong}> Y' \end{CD}\quad,\quad \begin{CD} X @<{\cong}<< Y \\ @V{f}VV @VV{f'}V \\ X' @>>{\cong}> Y' \end{CD}\quad,\quad \begin{CD} X @>{\cong}>> Y \\ @V{f}VV @VV{f'}V \\ X' @<<{\cong}< Y' \end{CD}\quad,\quad \begin{CD} X @<{\cong}<< Y \\ @V{f}VV @VV{f'}V \\ X' @<<{\cong}< Y' \end{CD}\quad.$$

Prove that statement. ◇

*Exercise* 1K.3. Consider the following statement: "In every concrete category **C** with realization R, a morphism $f \in \mathbf{C}$ is an isomorphism $\Leftrightarrow$ $R(f) \in \mathbf{Set}$ is an isomorphism.". Is this claim true or false? Is at least one of the two directions, meaning $\Rightarrow$ or $\Leftarrow$, correct? ◇

*Exercise* 1K.4. What is the skeleton of the category **fSet** from Example 1H.21? ◇

*Exercise* 1K.5. Let $F \in \mathbf{Hom}(\mathbf{C}, \mathbf{D})$ be an equivalence of categories. Show that $f \in \mathbf{C}$ is monic (or epic, or an isomorphism) if and only if $F(f) \in \mathbf{D}$ is monic (or epic, or an isomorphism). ◇

1L. **Additional topics.** Topics we have not touched on in these notes include:

**(a)** (Co)limits and, as special cases, (co)products *etc.*;

**(b)** Monads and algebras;

**(c)** Kan extensions.

These topics can be found *e.g.* in [**ML98**]. For the interested reader, these are covered on the "What is...category theory?" playlist on [**Tub21**].

## 2. Monoidal categories I – definitions, examples and graphical calculus

We have seen Feynman diagrams for categories, but they are one dimensional (with the exception of functor calculus, but we have not made this precise yet). So:

What are the right axioms to get a two dimensional diagrammatic calculus?

Why we want this is illustrated and explained in Figure 5.

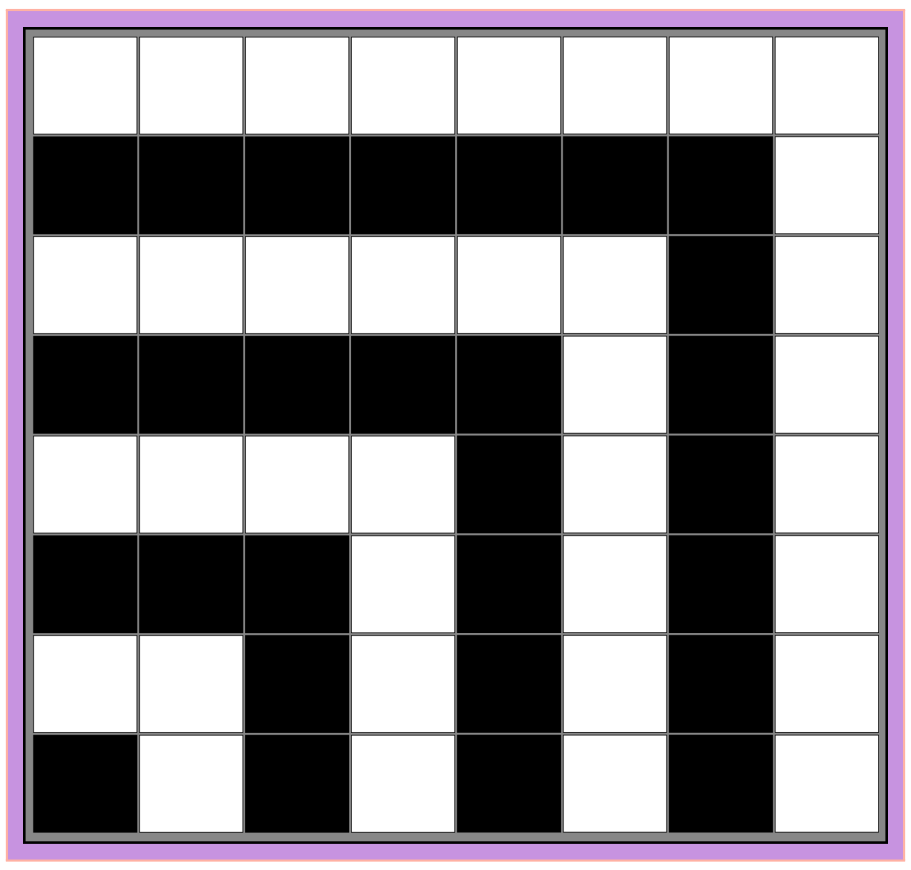

Figure 5. A proof without words: The sum $1+3+5+7+9+...+2n-1$ is the square $n^2$. Wikipedia (as in the link below) writes: "For a proof to be accepted by the mathematical community, it must logically show how the statement it aims to prove follows totally and inevitably from a set of assumptions. A proof without words might imply such an argument, but it does not make one directly, so it cannot take the place of a formal proof where one is required. Rather, mathematicians use proofs without words as illustrations and teaching aids for ideas that have already been proven formally." The point of diagrammatic algebra, as part of categorical algebra, is that this is not necessarily true: There can be rigorous proofs without words proof.

The picture is a variation of https://en.wikipedia.org/wiki/Proof_without_words.

2A. **Motivating example.** If one considers the vertical composition $\circ$ of natural transformations Equation 1F-4, say

$$\begin{array}{c} \mathrm{H}(\mathtt{X}) \\ {\scriptstyle\gamma_{\mathtt{X}}}\uparrow \\ \mathrm{G}(\mathtt{X}) \end{array} \circ \begin{array}{c} \mathrm{G}(\mathtt{X}) \\ {\scriptstyle\alpha_{\mathtt{X}}}\uparrow \\ \mathrm{F}(\mathtt{X}) \end{array} = \begin{array}{c} \mathrm{H}(\mathtt{X}) \\ {\scriptstyle\gamma_{\mathtt{X}}}\uparrow \\ \mathrm{G}(\mathtt{X}) \\ {\scriptstyle\alpha_{\mathtt{X}}}\uparrow \\ \mathrm{F}(\mathtt{X}) \end{array},$$

then it seems there should be a second, ***horizontal composition*** $\otimes$:

$$\begin{array}{ccc} \mathrm{I}(\mathtt{X}) & \xrightarrow{\mathrm{I}(\mathrm{f})} & \mathrm{I}(\mathtt{Y}) \\ {\scriptstyle\beta_{\mathtt{X}}}\uparrow & & \uparrow{\scriptstyle\beta_{\mathtt{Y}}} \\ \mathrm{H}(\mathtt{X}) & \xrightarrow[\mathrm{H}(\mathrm{f})]{} & \mathrm{H}(\mathtt{Y}) \end{array} \otimes \begin{array}{ccc} \mathrm{G}(\mathtt{X}) & \xrightarrow{\mathrm{G}(\mathrm{f})} & \mathrm{G}(\mathtt{Y}) \\ {\scriptstyle\alpha_{\mathtt{X}}}\uparrow & & \uparrow{\scriptstyle\alpha_{\mathtt{Y}}} \\ \mathrm{F}(\mathtt{X}) & \xrightarrow[\mathrm{F}(\mathrm{f})]{} & \mathrm{F}(\mathtt{Y}) \end{array} = \begin{array}{ccc} \mathrm{IG}(\mathtt{X}) & \xrightarrow{\mathrm{IG}(\mathrm{f})} & \mathrm{IG}(\mathtt{Y}) \\ {\scriptstyle(\beta\otimes\alpha)_{\mathtt{X}}}\uparrow & & \uparrow{\scriptstyle(\beta\otimes\alpha)_{\mathtt{Y}}} \\ \mathrm{HF}(\mathtt{X}) & \xrightarrow[\mathrm{HF}(\mathrm{f})]{} & \mathrm{HF}(\mathtt{Y}) \end{array}. \tag{2A-1}$$

The picture to have in mind in colored Feynman diagrams is:

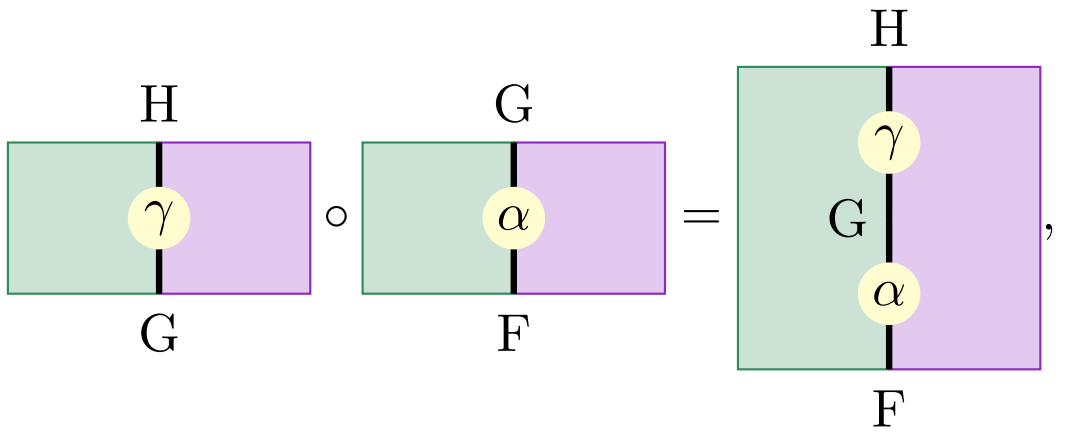

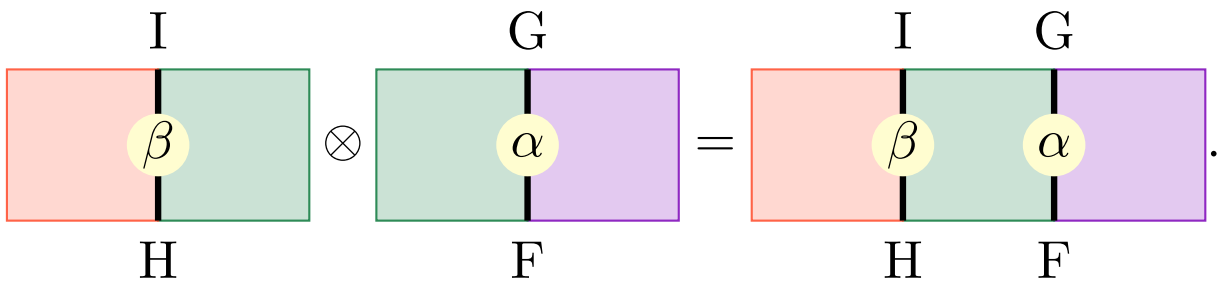

As we will see, there is indeed such a second composition.

2B. **A more down to earth motivating example.** Recall from Definition 1B.12 that we can form the pair category **Set**×**Set**. Note that we have a functor

$$\otimes\colon \mathbf{Set}\times\mathbf{Set}\to\mathbf{Set},\quad \otimes\big((\mathtt{X},\mathtt{Y})\big)=\mathtt{X}\otimes\mathtt{Y}=\mathtt{X}\times\mathtt{Y},\otimes\big((\mathrm{f},\mathrm{g})\big)=\mathrm{f}\otimes\mathrm{g}=\mathrm{f}\times\mathrm{g},$$

where we already use the usual standard notation, meaning writing *e.g.* $\mathtt{X}\otimes\mathtt{Y}$ instead of $\otimes\big((\mathtt{X},\mathtt{Y})\big)$, for these kinds of functors.

The functor $\otimes$ is actually a bit better: it is a ***bifunctor***. This mean that it satisfies an identity rule and the ***interchange law***, *i.e.*

$$\mathrm{id}_{\mathtt{X}}\otimes\mathrm{id}_{\mathtt{Y}}=\mathrm{id}_{\mathtt{X}\otimes\mathtt{Y}},\quad (\mathrm{gf})\otimes(\mathrm{kh})=(\mathrm{g}\otimes\mathrm{k})(\mathrm{f}\otimes\mathrm{h}).\tag{2B-1}$$

We observe the following:

- This is only weakly associative, *i.e.*

  $$\mathtt{X}\otimes(\mathtt{Y}\otimes\mathtt{Z})\neq(\mathtt{X}\otimes\mathtt{Y})\otimes\mathtt{Z},\ \text{but rather}\ \mathtt{X}\otimes(\mathtt{Y}\otimes\mathtt{Z})\cong(\mathtt{X}\otimes\mathtt{Y})\otimes\mathtt{Z},$$

  because the set $\mathtt{X}\otimes(\mathtt{Y}\otimes\mathtt{Z})$ contains elements of the form $\big(x,(y,z)\big)$, while $(\mathtt{X}\otimes\mathtt{Y})\otimes\mathtt{Z}$ contains elements of the form $\big((x,y),z\big)$.
- Similarly, this operation has $\mathbb{1}=\{\bullet\}$ as a unit, but it is again only a weak unit, meaning

  $$\mathbb{1}\otimes\mathtt{X}\neq\mathtt{X}\neq\mathtt{X}\otimes\mathbb{1},\ \text{but rather}\ \mathbb{1}\otimes\mathtt{X}\cong\mathtt{X}\cong\mathtt{X}\otimes\mathbb{1}.$$

The reader might object that the difference between, for example, $\big(x,(y,z)\big)$ and $\big((x,y),z\big)$ is rather silly. Indeed, in this case it is, but in general the symbol $\cong$ can be a nontrivial operation, so let us keep the difference around for now.

2C. **A word about conventions.** As we have seen in the example above, there are two operations for morphisms $\circ$ and $\otimes$, but only one $\otimes$ for objects. Recall, *cf.* Convention 1A.1, that we already abbreviate $\mathrm{gf}=\mathrm{g}\circ\mathrm{f}$, and we will do the same for objects:

*Convention* 2C.1. We will write $\mathtt{XY}=\mathtt{X}\otimes\mathtt{Y}$ for simplicity, and similarly we write $\mathtt{X}^k$ instead of $k\in\mathbb{Z}_{\geq 0}$ factors of the form $\mathtt{X}\otimes\ldots\otimes\mathtt{X}$. ◇

*Convention* 2C.2. Although monoidal categories, functor *etc.* usually consists of a choice of extra data, we will for brevity often just write *e.g.* **C** for a monoidal category. We also *e.g.* write "**C** is a monoidal category" when the choice of monoidal structure is clear from the context. ◇

*Convention* 2C.3. There will be several places where we have two or more monoidal categories with potentially different structures. However, in order not to overload the notation we will write *e.g.* $\mathbb{1}$ for all of them instead of for example $\mathbb{1}_{\mathbf{C}}$. ◇

2D. **Basics.** The definition of a monoidal category is a mouthful (but we will get rid most of the complication later in Theorem 2I.5; for now recall that $\cong$ could indicate something nontrivial):

**Definition 2D.1.** A ***monoidal category*** $(\mathbf{C},\otimes,\mathbb{1},\alpha,\lambda,\rho)$ consists of

- a category **C**;
- a bifunctor (*cf.* Equation 2B-1)

  $$\otimes\colon\mathbf{C}\times\mathbf{C}\to\mathbf{C},\quad\otimes\big((\mathtt{X},\mathtt{Y})\big)=\mathtt{XY},\otimes\big((\mathrm{f},\mathrm{g})\big)=\mathrm{f}\otimes\mathrm{g},$$

  called ***monoidal product***;
- a ***unit*** (object) $\mathbb{1}\in\mathbf{C}$;

- a collection of natural isomorphisms

$$\alpha_{\mathtt{X},\mathtt{Y},\mathtt{Z}}\colon \mathtt{X}(\mathtt{YZ}) \xrightarrow{\cong} (\mathtt{XY})\mathtt{Z}, \tag{2D-2}$$

  for all $\mathtt{X},\mathtt{Y},\mathtt{Z}\in\mathtt{C}$, called ***associators***;
- a collection of natural isomorphisms

$$\lambda_{\mathtt{X}}\colon \mathbb{1}\mathtt{X} \xrightarrow{\cong} \mathtt{X},\quad \rho_{\mathtt{X}}\colon \mathtt{X}\mathbb{1} \xrightarrow{\cong} \mathtt{X}, \tag{2D-3}$$

  for all $\mathbf{X}\in\mathbf{C}$, called ***left and right unitors***;

such that

(i) the ⬠ ***equality*** holds, *i.e.* we have commuting diagrams

$$\begin{array}{ccccc} & & \big((\mathtt{XY})\mathtt{Z}\big)\mathtt{A} & & \\ & \nearrow^{\alpha_{\mathtt{XY},\mathtt{Z},\mathtt{A}}} & & \nwarrow^{\alpha_{\mathtt{X},\mathtt{Y},\mathtt{Z}}\otimes\mathrm{id}_{\mathtt{A}}} & \\ (\mathtt{XY})(\mathtt{ZA}) & & & & \big(\mathtt{X}(\mathtt{YZ})\big)\mathtt{A} \\ \nwarrow^{\alpha_{\mathtt{X},\mathtt{Y},\mathtt{ZA}}} & & & & \nearrow^{\alpha_{\mathtt{X},\mathtt{YZ},\mathtt{A}}} \\ & \mathtt{X}\big(\mathtt{Y}(\mathtt{ZA})\big) & \xrightarrow[\mathrm{id}_{\mathtt{X}}\otimes\alpha_{\mathtt{Y},\mathtt{Z},\mathtt{A}}]{} & \mathtt{X}\big((\mathtt{YZ})\mathtt{A}\big) & \end{array},$$

for all $\mathtt{X},\mathtt{Y},\mathtt{Z},\mathtt{A}\in\mathbf{C}$.

(ii) the △ ***equality*** holds, *i.e.* we have commuting diagrams

$$\begin{array}{ccc} & \mathtt{XY} & \\ \nearrow^{\mathrm{id}_{\mathtt{X}}\otimes\lambda_{\mathtt{Y}}} & & \nwarrow^{\rho_{\mathtt{X}}\otimes\mathrm{id}_{\mathtt{Y}}} \\ \mathtt{X}(\mathbb{1}\mathtt{Y}) & \xrightarrow[\alpha_{\mathtt{X},\mathbb{1},\mathtt{Y}}]{} & (\mathtt{X}\mathbb{1})\mathtt{Y} \end{array},$$

for all $\mathtt{X},\mathtt{Y}\in\mathbf{C}$.

(The ⬠ equality can be compared with Equation 1B-11.) ◇

*Remark* 2D.4. There is a hidden □ ***equality***, coming from naturality,

$$\begin{array}{ccc} \mathtt{X}'(\mathtt{Y}'\mathtt{Z}') & \xrightarrow{\alpha_{\mathtt{X}',\mathtt{Y}',\mathtt{Z}'}} & (\mathtt{X}'\mathtt{Y}')\mathtt{Z}' \\ {\scriptstyle \mathrm{f}\otimes(\mathrm{g}\otimes\mathrm{h})}\uparrow & & \uparrow{\scriptstyle (\mathrm{f}\otimes\mathrm{g})\otimes\mathrm{h}} \\ \mathtt{X}(\mathtt{YZ}) & \xrightarrow[\alpha_{\mathtt{X},\mathtt{Y},\mathtt{Z}}]{} & (\mathtt{XY})\mathtt{Z} \end{array},$$

which holds for all for all $\mathtt{X},\mathtt{Y},\mathtt{Z}\in\mathbf{C}$ and all $\mathrm{f},\mathrm{g},\mathrm{h}\in\mathbf{C}$. ◇

**Definition 2D.5.** A monoidal category $\mathbf{C}$ is called ***strict*** if all associators and unitors are identities, and ***non-strict*** otherwise. ◇

**Example 2D.6.** Monoidal categories arise in the wild.

(a) As seen above, **Set** with $\otimes=\times$ and $\mathbb{1}=\{\bullet\}$ is a non-strict monoidal category.

(b) Similarly, $\mathbf{Vec}_{\Bbbk}$ or $\mathbf{fdVec}_{\Bbbk}$ with $\otimes=\otimes_{\Bbbk}$ and $\mathbb{1}=\Bbbk$ are non-strict monoidal categories.

(c) The skeletons of the three examples above, with the same monoidal structures, are strict monoidal categories.

For completeness, on $\mathbf{Mat}_{\Bbbk}$, the skeleton of $\mathbf{fdVec}_{\Bbbk}$, the monoidal product is, for example,

$$\begin{pmatrix}1&2\\3&4\end{pmatrix}\otimes\begin{pmatrix}5&6\end{pmatrix}=\begin{pmatrix}1\cdot\begin{pmatrix}5&6\end{pmatrix}&2\cdot\begin{pmatrix}5&6\end{pmatrix}\\3\cdot\begin{pmatrix}5&6\end{pmatrix}&4\cdot\begin{pmatrix}5&6\end{pmatrix}\end{pmatrix}=\begin{pmatrix}5&6&10&12\\15&18&20&24\end{pmatrix},$$

$$\begin{pmatrix}5&6\end{pmatrix}\otimes\begin{pmatrix}1&2\\3&4\end{pmatrix}=\begin{pmatrix}5\cdot\begin{pmatrix}1&2\\3&4\end{pmatrix}&6\cdot\begin{pmatrix}1&2\\3&4\end{pmatrix}\end{pmatrix}=\begin{pmatrix}5&10&6&12\\15&20&18&24\end{pmatrix}.$$

This is known as the ***Kronecker product*** of matrices. ◇

**Example 2D.7.** Monoidal structures on categories are far from being unique. For example, $\mathbf{Vec}_{\Bbbk}$ and $\mathbf{fdVec}_{\Bbbk}$ have another monoidal structure given by $\otimes=\oplus$ and $\mathbb{1}=\{0\}$, which is again non-strict. We will however always use the monoidal structures in Example 2D.6.(b) for vector spaces. ◇

**Example 2D.8.** Diagrammatic categories such as **1Cob**, **1Tan** and **1State**, *cf.* Example 1B.9, have (often) a monoidal structure given by $\otimes$ being juxtaposition (called "Ignore the $\otimes$ symbol."), *e.g.*

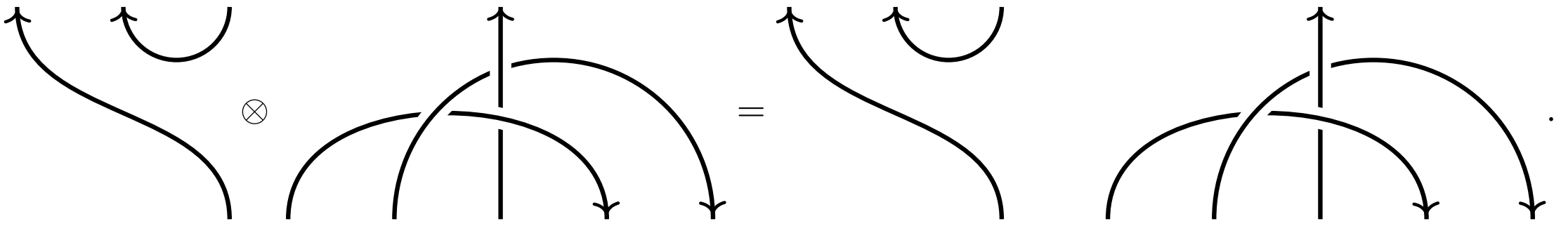

and $\mathbb{1}$ being the empty diagram. These monoidal structures are strict. $\diamond$

The following is in some sense the motivation for the name "monoidal category." Recall hereby the Grothendieck classes $K_0(\mathbf{C})$ of $\mathbf{C}$, see Definition 1H.15.

**Proposition 2D.9.** *For any monoidal category* $\mathbf{C}$ *its Grothendieck classes* $K_0(\mathbf{C})$ *form a monoid with multiplication and unit*

$$[\mathtt{X}][\mathtt{Y}] = [\mathtt{XY}], \quad 1 = [\mathbb{1}].$$

*Proof.* Directly from the definitions, *e.g.* the associator Equation 2D-2 and the unitors Equation 2D-3 descent to associativity and unitality on $K_0(\mathbf{C})$. $\square$

**Example 2D.10.** Coming back to Example 1H.17, $K_0(\mathbf{fdVec}_{\Bbbk}) \xrightarrow{\cong} \mathbb{Z}_{\geq 0}$ with $[\Bbbk^n] \mapsto n$ is an isomorphism of monoids. Note that the Kronecker product has no analog on $K_0(\mathbf{fdVec}_{\Bbbk})$. $\diamond$

Example 2D.8 gives important examples of strict monoidal categories, while crucial examples of non-strict monoidal categories are the monoidal incarnations of groups. These are very different from the ones we have, noting that every group is of course a monoid, seen in Example 1B.2.(a):

**Example 2D.11.** Let G be a group.

(a) The category $\mathbf{Vec}(\mathrm{G})$ is the category with $\mathrm{Ob}\big(\mathbf{Vec}(\mathrm{G})\big) = \mathrm{G}$, and whose morphisms are only identities. The monoidal product is $\mathtt{i} \otimes \mathtt{j} = \mathtt{ij}$, with $i, j, ij \in \mathrm{G}$. For example, if $\mathrm{G} = \mathbb{Z}/2\mathbb{Z} \times \mathbb{Z}/2\mathbb{Z}$ (written additively), then we have

$$\overset{\mathrm{id}_{(0,0)}}{\circlearrowright} \atop (0,0) \qquad \overset{\mathrm{id}_{(1,0)}}{\circlearrowright} \atop (1,0) \qquad \overset{\mathrm{id}_{(0,1)}}{\circlearrowright} \atop (0,1) \qquad \overset{\mathrm{id}_{(1,1)}}{\circlearrowright} \atop (1,1) \;, \quad (\mathtt{a},\mathtt{b}) \otimes (\mathtt{c},\mathtt{d}) = (\mathtt{a}+\mathtt{c},\mathtt{b}+\mathtt{d}),$$

$$\mathrm{Hom}_{\mathbf{Vec}(\mathrm{G})}(\mathtt{i},\mathtt{j}) = \begin{cases} \{\mathrm{id}_{\mathtt{i}}\} & \text{if } i = j, \\ \emptyset & \text{if } i \neq j. \end{cases}$$

Thus, as a category $\mathbf{Vec}(\mathrm{G})$ is rather boring and the point is the monoidal structure, which is strict, by construction.

(b) We also have the $\Bbbk$ linear version $\mathbf{Vec}_{\Bbbk}(\mathrm{G})$ of $\mathbf{Vec}(\mathrm{G})$. The only difference is that the endomorphisms are now given by scalars times the identities:

$$\mathrm{Hom}_{\mathbf{Vec}_{\Bbbk}(\mathrm{G})}(\mathtt{i},\mathtt{j}) \cong \begin{cases} \Bbbk & \text{if } i = j, \\ 0 & \text{if } i \neq j. \end{cases}$$

The monoidal category $\mathbf{Vec}_{\Bbbk}(\mathrm{G})$ is strict.

(c) Let $\omega \in Z^3(\mathrm{G}, \mathbb{C}^{\times} = \mathbb{C}/\mathbb{Z})$ be a 3 cocycle of G, see Remark 2D.13. Then we can define a monoidal category $\mathbf{Vec}^{\omega}_{\mathbb{C}}(\mathrm{G})$ exactly as above, but with associator and unitors

$$\alpha_{\mathtt{i},\mathtt{j},\mathtt{k}} = \omega(i,j,k)\mathrm{id}_{\mathtt{ijk}} \quad \lambda_{\mathtt{i}} = \omega(1,1,i)^{-1}\mathrm{id}_{\mathtt{i}}, \quad \rho_{\mathtt{i}} = \omega(i,1,1)\mathrm{id}_{\mathtt{i}}. \tag{2D-12}$$

Explicitly, for $\mathrm{G} = \mathbb{Z}/2\mathbb{Z}$ we have $H^3(\mathrm{G}, \mathbb{C}^{\times}) \cong \mathbb{Z}/2\mathbb{Z}$ and the nontrivial $\omega$ has $\omega(1,1,1) = -1$. Finally, note that $\mathbf{Vec}^{1}_{\mathbb{C}}(\mathrm{G}) = \mathbf{Vec}_{\mathbb{C}}(\mathrm{G})$, but for a nontrivial $\omega \in H^3(\mathrm{G}, \mathbb{C}^{\times})$ the monoidal category $\mathbf{Vec}^{\omega}_{\mathbb{C}}(\mathrm{G})$ is non-strict and skeletal at the same time.

In (c) we saw the first instance where $\cong$ is a nontrivial operation. $\diamond$

*Remark* 2D.13. For a group G, one can define a cohomology theory $H^*(\mathrm{G}, \mathtt{A})$, called ***group cohomology*** with values in $\mathtt{A}$. Here $\mathtt{A} \in \mathbf{Mod}(\mathbb{Z})$ is assumed to be a G-representation. See [Wei94, Chapter 6] for a reference.

For these lecture notes, all we need to know about group cohomology is summarized now for convenience. As usual these are constructed from a certain cochain complex and $H^i(\mathrm{G}, \mathrm{A}) = Z^i(\mathrm{G}, \mathrm{A})/B^i(\mathrm{G}, \mathrm{A})$, so $i$ cocycles modulo $i$ coboundaries. The 3 cocycles which are functions $\omega\colon \mathrm{G}\times\mathrm{G}\times\mathrm{G}\to\mathbb{C}^{\times}$ satisfying

$$\omega(j,k,l)\omega(i,jk,l)\omega(i,j,k)=\omega(ij,k,l)\omega(i,j,kl),$$

(2D-14) pictorially: 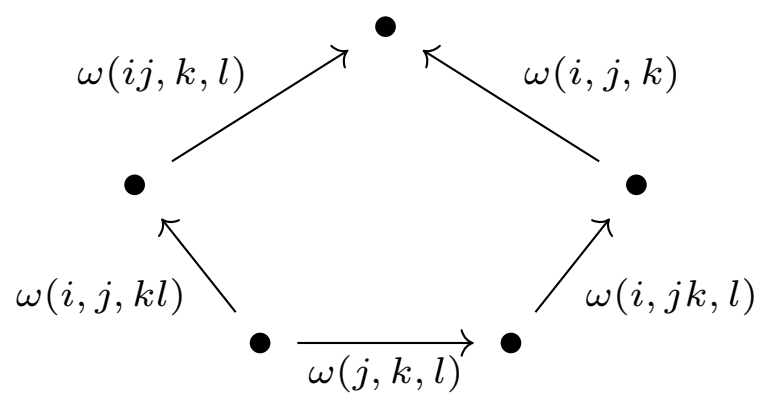

.

Comparing Equation 2D-14 and Definition 2D.1 shows that scaling as in Equation 2D-12 satisfies the △ and ⬠ equations. This is again similar to Equation 1B-11.

For completeness:

**(a)** For $\mathrm{G}=\mathbb{Z}/n\mathbb{Z}$ we have $H^3(\mathrm{G},\mathbb{C}^{\times})\cong\mathbb{Z}/n\mathbb{Z}$. In this case the nontrivial 3 cocycles $\omega_s$ are indexed by $s\in\mathbb{Z}/n\mathbb{Z}$ and are given by

$$\omega_s(i,j,k)=\zeta^{si(j+k-((j+k)\bmod n))/n}$$

where $\zeta$ is a primitive $n$th root of unity given by, for example, $\zeta=\exp(2\pi i/n)$.

**(b)** For a general abelian group one can use (a) and a Künneth-type-formula to compute the cohomology.

(Since we mix elements of $\mathbb{Z}/n\mathbb{Z}$ and $\mathbb{Z}$ in the exponent of $\zeta$ above, let us be completely explicit to avoid confusion. For $n=10$, $s=4$, $i=5$, $j=6$ and $k=7$ the above exponent is 20; for $n=2$, $s=1$, $i=1$, $j=1$ and $k=1$ the above exponent is 1 and we recover $\omega$ from Example 2D.11.(c).) ◇

*Remark* 2D.15. Note that for $\mathbf{Vec}(\mathrm{G})$ or $\mathbf{Vec}_{\mathbb{C}}(\mathrm{G})$ we can also allow monoids M instead of groups G, or work over rings $\mathbb{S}$, but for $\mathbf{Vec}_{\mathbb{S}}^{\omega}(\mathrm{M})$ one would need to be careful how to define it. ◇

A good question is whether we can "ignore" non-strict monoidal categories since working with associators and unitors is a bit messy. However, Example 2D.11.(c) suggests that one can not simply go to the skeleton, although this works for monoidal categories such as $\mathbf{fdVec}_{\mathbb{k}}$. We can only answer this question after we have a bit more technology at hand.

2E. **Feynman diagrams for monoidal categories.** Motivated by Example 2D.8, we get the following Feynman diagrammatics for strict monoidal categories. That is, given a strict monoidal category $\mathbf{C}$, we can depict $\otimes$ as juxtaposition and the unit as an empty diagram, *e.g.*

(2E-1) $\mathbb{1}\rightsquigarrow\emptyset$, $\mathtt{XY}\rightsquigarrow$ 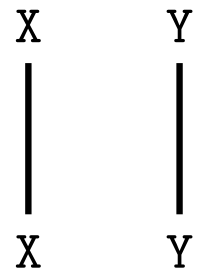

, $\mathrm{g}\otimes\mathrm{f}\rightsquigarrow$ 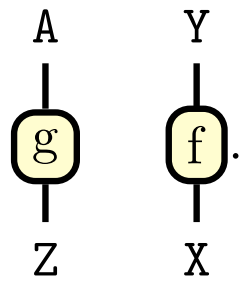

.

Note the cute fact that we do not need to be careful with the relative heights in Equation 2E-1 since the interchange law Equation 2B-1 implies that

(2E-2) $(\mathrm{id}_{\mathtt{A}}\otimes\mathrm{f})\circ(\mathrm{g}\otimes\mathrm{id}_{\mathtt{X}})\rightsquigarrow$ 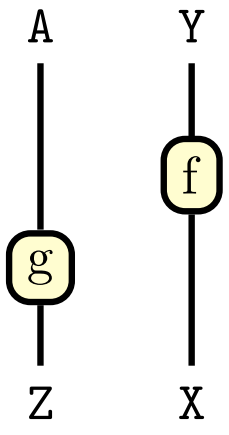

$=$ 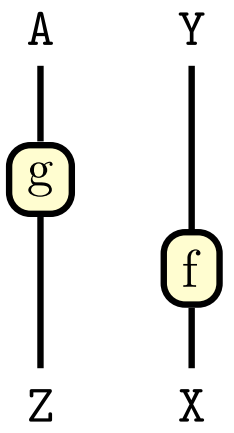

$\rightsquigarrow(\mathrm{g}\otimes\mathrm{id}_{\mathtt{Y}})\circ(\mathrm{id}_{\mathtt{Z}}\otimes\mathrm{f})$.

We can also illustrate morphisms with many $\otimes$ inputs nicely, *e.g.*

(2E-3) $\mathrm{f}\colon\mathtt{XYZ}\to\mathtt{AB}\rightsquigarrow$ 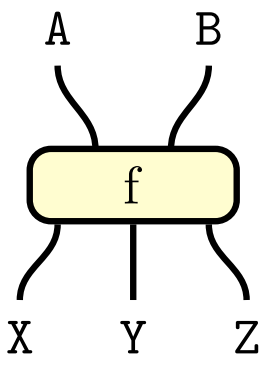

.

However, note that there are two drawbacks. First, diagrammatic calculus, by its very definition, is not suitable for non-strict monoidal categories. Second, although Equation 2E-2 looks promising, we do not have a two dimensional calculus yet as we are not allowed to change the upwards orientation of diagrams (recall Convention 1A.3), *e.g.*

(2E-4)

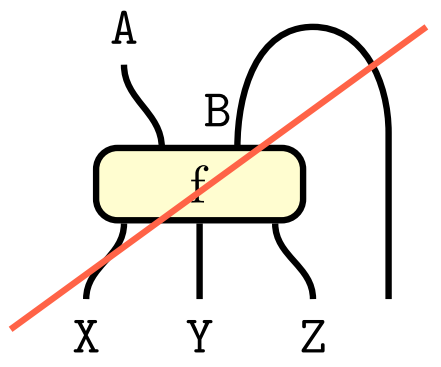

is not an allowed diagram.

*Remark* 2E.5. One should stress here that Equation 2E-4 and the text around it is not a contradiction to Example 2D.8: in that example the diagrams actually are just abbreviations for upwards oriented Feynman diagrams, *e.g.*

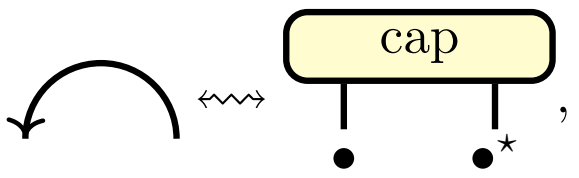

where • and •⋆ are the two generating objects of **1State**, as we will see. (Note that the unit is omitted from diagrams, *cf.* Equation 2E-1.) ◇

**Example 2E.6.** Every morphisms $\mathrm{f}\colon \mathbb{1}\to\mathbb{1}$ is presented by a ***floating diagram***:

$$\mathrm{f}\colon \mathbb{1}\to\mathbb{1} \rightsquigarrow \boxed{\mathrm{f}}.$$

This follows from our convention that $\mathbb{1}$ is diagrammatically presented by the empty diagram. ◇

We discuss how to incorporate non-strict monoidal categories below; the flaw of having only upward oriented diagrams will be taken care of in Section 4.

2F. **Coherence for monoids.** The following should be compared to Lemma 1B.10, and we are happy to repeat the arguments here since they are nice and important at the same time. For starters, let us compare two definitions of a monoid M, with Def1 being the one that you will usually find in written texts:

(2F-1)

| | | | | |
|---|---|---|---|---|
| Def1 | a set M | multiplication | unit | $h(gf) = (hg)f$ |
| Def2 | a set M | multiplication | unit | associativity |

where "associativity" means that all ways of using parentheses agree. Both definitions have their advantages: Def2 is arguably the correct definition, but Def1 is much more useful in practice and one only needs to check $h(fg) = (hg)f$ instead of infinitely many bracketings. So one would like to have the following, called ***coherence theorem for monoids***, which is rarely stated:

**Theorem 2F.2.** *The two definitions in Equation 2F-1 are equivalent.*

*Proof.* Clearly, Def2 implies Def1. To see that Def1 implies Def2, we argue diagrammatically. The condition $h(gf) = (hg)f$ can be pictured as

(2F-3)

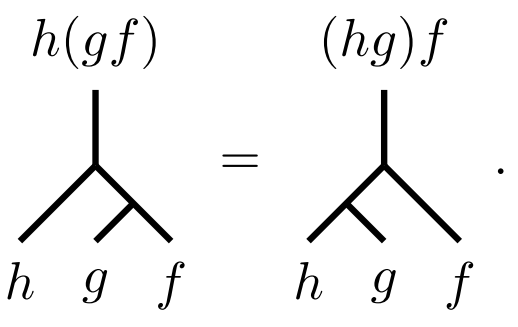

However, successively applying this equality gives

(2F-4)

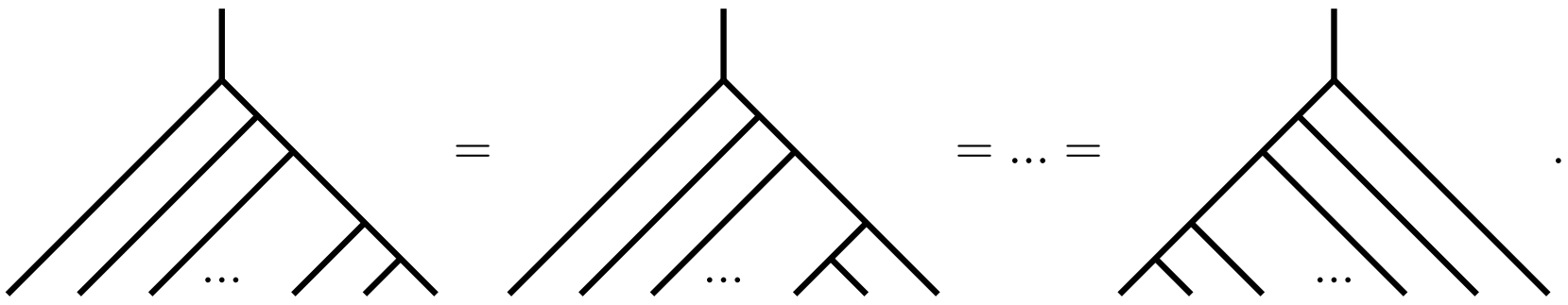

(Actually, these are not aligned, *cf.* Equation 2F-5.) Thus, all ways of putting parenthesis agree. □

The above can also be stated differently. Let $\mathrm{K}_n$ be the one dimensional CW complex (a.k.a. graph) obtained by adding an edge to the disjoint union of the graphs in Equation 2F-4 (with $n$ endpoints) for each application of Equation 2F-3, connecting the corresponding graphs. For example,

(2F-5) $K_4 =$

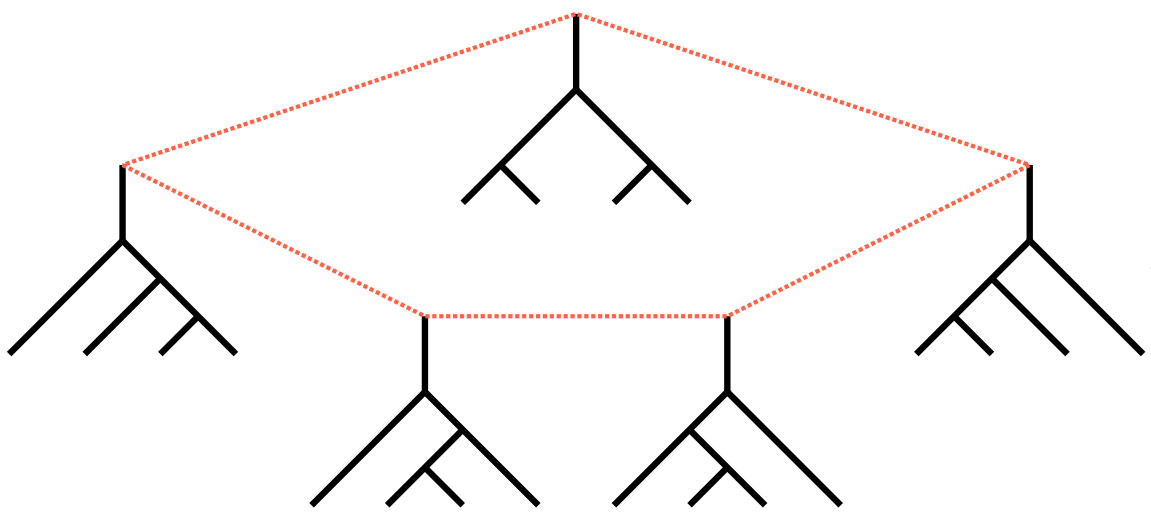

.

Then the above can be rephrased as $\pi_0(K_n) \cong 1$.

2G. **Coherence for monoidal categories.** With respect to the discussion about coherence for monoids, in particular, Equation 2F-1, here is Def2 for monoidal categories with Def1 being Definition 2D.1.

**Definition 2G.1.** A ***monoidal category*** $(\mathbf{C}, \otimes, \mathbb{1}, \alpha, \lambda, \rho)$ consists of

- a category $\mathbf{C}$;
- a bifunctor (*cf.* Equation 2B-1)
$$\otimes\colon \mathbf{C}\times\mathbf{C} \to \mathbf{C}, \quad \otimes\big((\mathtt{X},\mathtt{Y})\big) = \mathtt{XY}, \otimes\big((f,g)\big) = f \otimes g,$$
called ***monoidal product***;
- a ***unit*** (object) $\mathbb{1} \in \mathbf{C}$;
- a collection of natural isomorphisms
$$\alpha_{\mathtt{X},\mathtt{Y},\mathtt{Z}}\colon \mathtt{X} \otimes (\mathtt{Y} \otimes \mathtt{Z}) \xrightarrow{\cong} (\mathtt{X} \otimes \mathtt{Y}) \otimes \mathtt{Z},$$
for all $\mathtt{X},\mathtt{Y},\mathtt{Z} \in \mathtt{C}$, called ***associators***;
- a collection of natural isomorphisms
$$\lambda_{\mathtt{X}}\colon \mathbb{1}\mathtt{X} \xrightarrow{\cong} \mathtt{X}, \quad \rho_{\mathtt{X}}\colon \mathtt{X}\mathbb{1} \xrightarrow{\cong} \mathtt{X},$$
for all $\mathbf{X} \in \mathbf{C}$, called ***left and right unitors***;

such that "every formal diagram" made up of associators and unitors commutes. ◇

Similarly as for monoids, the difference between the two definitions is in this final sentence. We will not define what "every formal diagram" means precisely as this gets a bit technical. Moreover, we will only sketch a proof of the ***coherence theorem for monoidal categories*** (also known as ***Mac Lane's coherence theorem***), which is up next, for the very same reason.

**Theorem 2G.2.** *The two definitions Definition 2D.1 and Definition 2G.1 are equivalent.*

*Proof.* Let us sketch how this can be proven, following the exposition in [**Kap93**]. (A completely different proof is due to Mac Lane, see [**ML98**, Section VII.2].) Let us focus on associators, the idea of the proof with unitors is exactly the same.

The proof works by constructing certain polytopes $K_n$, sometimes called ***Stasheff polytopes*** or ***associahedron***. These are two dimensional analogs of the graphs we have seen in the proof of Theorem 2F.2, and constructed from the two relevant commuting diagrams, □ and ⬠ equations. For examples see Figure 6.

Then one needs to show that $\pi_1(K_n) \cong 1$. □

Note the analogy: In the one dimensional case (for monoids, categories *etc.*) one needs to assume that $K_3$ is "nice", and all other polytopes will then also be "nice". On the other hand, in the two dimensional case (for monoidal categories *etc.*) one needs to assume that $K_4$ is "nice".

2H. **Monoidal functors, natural transformations and equivalences.** First things first:

**Definition 2H.1.** A ***monoidal functor*** $(F, \xi, \xi_{\mathbb{1}})$ with $F \in \mathbf{Hom}(\mathbf{C}, \mathbf{D})$ consists of

- a functor F;
- a collection of natural isomorphisms
$$\xi_{\mathtt{X},\mathtt{Y}}\colon F(\mathtt{X})F(\mathtt{Y}) \xrightarrow{\cong} F(\mathtt{XY}),$$
for all $\mathtt{X},\mathtt{Y} \in \mathtt{C}$;
- a natural isomorphism
$$\xi_{\mathbb{1}}\colon \mathbb{1} \xrightarrow{\cong} F(\mathbb{1});$$

such that

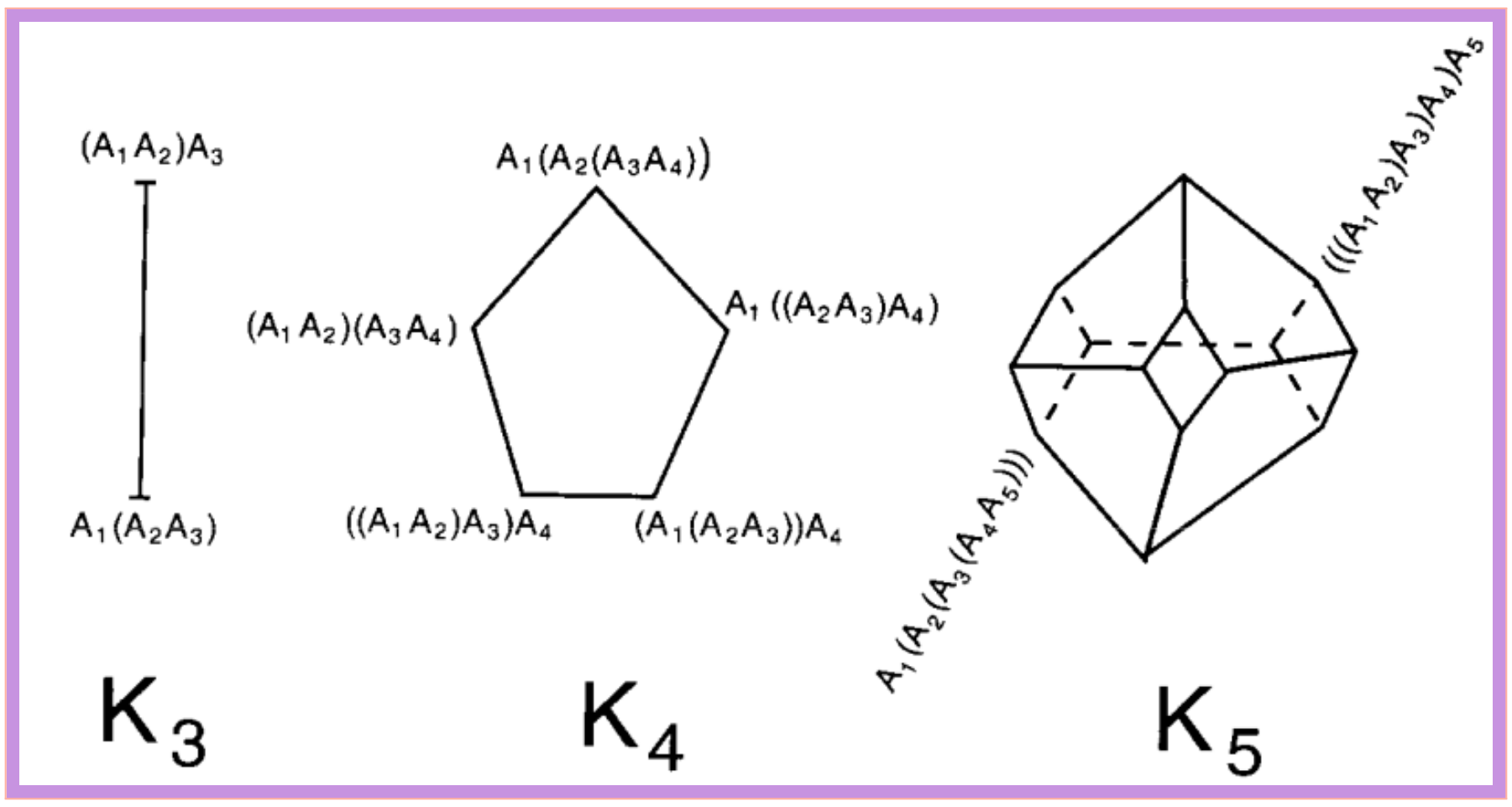

FIGURE 6. The first few Stasheff polytopes. The middle is, of course, again a pentagon.
Picture from [Kap93]

(i) the ⬡ ***equality*** holds, *i.e.* we have a commuting diagram

$$\begin{array}{ccccc} & \big(\mathtt{F(X)F(Y)}\big)\mathtt{F(Z)} & \xrightarrow{\xi_{\mathtt{X,Y}}\otimes\mathrm{id}_{\mathtt{F(Z)}}} & \mathtt{F(XY)F(Z)} & \\ {}^{\alpha_{\mathtt{F(X),F(Y),F(Z)}}}\nearrow & & & & \searrow^{\xi_{\mathtt{XY,Z}}} \\ \mathtt{F(X)}\big(\mathtt{F(Y)F(Z)}\big) & & & & \mathtt{F}\big((\mathtt{XY})\mathtt{Z}\big) \\ {}_{\mathrm{id}_{\mathtt{F(X)}}\otimes\xi_{\mathtt{Y,Z}}}\searrow & & & & \nearrow_{\mathtt{F}(\alpha_{\mathtt{X,Y,Z}})} \\ & \mathtt{F(X)F(YZ)} & \xrightarrow[\xi_{\mathtt{X,YZ}}]{} & \mathtt{F}\big(\mathtt{X(YZ)}\big) & \end{array},$$

for all $\mathtt{X,Y,Z}\in\mathtt{C}$;

(ii) a ***left*** and a ***right*** □ ***equation*** holds, *i.e.* we have commuting diagrams

$$\begin{array}{ccc} \mathbb{1}\mathtt{F(X)} & \xrightarrow{\xi_{\mathbb{1}}\otimes\mathrm{id}_{\mathtt{F(X)}}} & \mathtt{F}(\mathbb{1})\mathtt{F(X)} \\ {\scriptstyle\lambda_{\mathtt{F(X)}}}\downarrow & & \downarrow{\scriptstyle\xi_{\mathbb{1},\mathtt{X}}} \\ \mathtt{F(X)} & \xleftarrow[\mathtt{F}(\lambda_{\mathtt{X}})]{} & \mathtt{F}(\mathbb{1}\mathtt{X}) \end{array}, \quad \begin{array}{ccc} \mathtt{F(X)}\mathbb{1} & \xrightarrow{\mathrm{id}_{\mathtt{F(X)}}\otimes\xi_{\mathbb{1}}} & \mathtt{F(X)F}(\mathbb{1}) \\ {\scriptstyle\rho_{\mathtt{F(X)}}}\downarrow & & \downarrow{\scriptstyle\xi_{\mathtt{X},\mathbb{1}}} \\ \mathtt{F(X)} & \xleftarrow[\mathtt{F}(\rho_{\mathtt{X}})]{} & \mathtt{F}(\mathtt{X}\mathbb{1}) \end{array},$$

for all $\mathtt{X}\in\mathtt{C}$.

Such functors can also be called non-strict monoidal. ◇

**Definition 2H.2.** A ***monoidal natural transformation*** $\alpha\colon\mathrm{F}\Rightarrow\mathrm{G}$ between monoidal functors $\mathrm{F},\mathrm{G}\in\mathbf{Hom}(\mathbf{C},\mathbf{D})$ is a natural transformation such that

(i) for all $\mathrm{X},\mathrm{Y}\in\mathbf{C}$ there is a commuting diagram

$$\begin{array}{ccc} \mathtt{G(X)G(Y)} & \xrightarrow{\xi_{\mathtt{X,Y}}} & \mathtt{G(XY)} \\ {\scriptstyle\alpha_{\mathtt{X}}\otimes\alpha_{\mathtt{Y}}}\uparrow & & \uparrow{\scriptstyle\alpha_{\mathtt{XY}}} \\ \mathtt{F(X)F(Y)} & \xrightarrow[\xi_{\mathtt{X,Y}}]{} & \mathtt{F(XY)} \end{array};$$

(ii) there is a commuting diagram

$$\begin{array}{ccc} & \mathbb{1} & \\ {}^{\xi_{\mathbb{1}}}\swarrow & & \searrow^{\xi_{\mathbb{1}}} \\ \mathtt{F}(\mathbb{1}) & \xrightarrow[\alpha_{\mathbb{1}}]{} & \mathtt{G}(\mathbb{1}) \end{array}.$$

◇

**Lemma 2H.3.** *We have the following.*

*(i) If* F *and* G *are monoidal functors, then so is* GF.

*(ii) If* $\alpha$ *and* $\beta$ *are monoidal natural transformations, then so is* $\beta\alpha$. □

Thus, since the identity functor has an evident structure of a monoidal functor:

**Example 2H.4.** We get further examples of (plain) categories.

(a) There is a category $\mathbf{MCat}$, the ***category of monoidal categories***. Its objects are monoidal categories and its morphisms are monoidal functors.

(b) There is a category $\mathbf{Hom}_{\otimes}(\mathbf{C},\mathbf{D})$, the ***category of monoidal functors*** from $\mathbf{C}$ to $\mathbf{D}$. Its objects are monoidal functors and its morphisms are monoidal natural transformations, with vertical composition Equation 1F-4.

We elaborate on (b) in the case $\mathbf{C}=\mathbf{D}$ momentarily. $\diamond$

**Example 2H.5.** Given any category $\mathbf{C}$, the category $\mathbf{End}(\mathbf{C})$ of its endofunctors is a strict monoidal category:

- the composition $\circ$ is vertical composition of natural transformations Equation 1F-4;
- the monoidal product on objects is $\mathrm{G}\otimes\mathrm{F}=\mathrm{GF}$, *i.e.* composition of functors;
- the monoidal product on morphisms is $\beta\otimes\alpha=\beta\alpha$, *i.e.* horizontal composition of natural transformation Equation 2A-1.

Indeed, in colored Feynman diagrams fixing $\mathbf{C}$ corresponds to color all faces with, say, white and we obtain the calculus from Section 2E. $\diamond$

**Definition 2H.6.** $\mathbf{C},\mathbf{D}\in\mathbf{MCat}$ are called ***monoidally equivalent***, denoted by $\mathbf{C}\simeq_{\otimes}\mathbf{D}$, if there exists an equivalence $\mathrm{F}\in\mathbf{Hom}(\mathbf{C},\mathbf{D})$ which is additionally a monoidal functor. $\diamond$

**Example 2H.7.** Equivalent monoidal categories need not, but can be, monoidally equivalent:

(a) Recall that $\mathbf{fdVec}_{\Bbbk}\simeq\mathbf{Mat}_{\Bbbk}$. Together with the choice of monoidal structures being the usual tensor products, this is an monoidal equivalence $\mathbf{fdVec}_{\Bbbk}\simeq_{\otimes}\mathbf{Mat}_{\Bbbk}$.

(b) We have $\mathbf{Vec}_{\Bbbk}(\mathrm{G})\simeq\mathbf{Vec}_{\Bbbk}(\mathrm{G}')$ are equivalent as categories if and only if $\#\mathrm{G}=\#\mathrm{G}'$. However, $\mathbf{Vec}_{\Bbbk}(\mathrm{G})\simeq_{\otimes}\mathbf{Vec}_{\Bbbk}(\mathrm{G}')$ if and only if $\mathrm{G}=\mathrm{G}'$.

(c) Similarly, $\mathbf{Vec}_{\Bbbk}^{\omega}(\mathrm{G})\simeq\mathbf{Vec}_{\Bbbk}^{\omega'}(\mathrm{G}')$ holds always, *i.e.* regardless of the 3 cocycles. However, $\mathbf{Vec}_{\Bbbk}^{\omega}(\mathrm{G})$ and $\mathbf{Vec}_{\Bbbk}^{\omega'}(\mathrm{G})$ are rarely equivalent as monoidal categories. Explicitly, let $\omega$ be the nontrivial 3 cocycle of $\mathrm{G}=\mathbb{Z}/2\mathbb{Z}$. Then $\mathbf{Vec}_{\mathbb{C}}(\mathbb{Z}/2\mathbb{Z})\simeq\mathbf{Vec}_{\mathbb{C}}^{\omega}(\mathbb{Z}/2\mathbb{Z})$ but $\mathbf{Vec}_{\mathbb{C}}(\mathbb{Z}/2\mathbb{Z})\not\simeq_{\otimes}\mathbf{Vec}_{\mathbb{C}}^{\omega}(\mathbb{Z}/2\mathbb{Z})$.

Here $G$ denotes a finite group. $\diamond$

*Remark* 2H.8. More general as in Example 2H.7.(c), one can check that $\mathbf{Vec}_{\mathbb{C}}^{\omega}(\mathrm{G})\not\simeq_{\otimes}\mathbf{Vec}_{\mathbb{C}}^{\omega'}(\mathrm{G})$ unless $\omega$ and $\omega'$ are cohomologically equivalent, see *e.g.* [**EGNO15**, Proposition 2.6.1]. (The philosophy is that $H^{3}(\mathrm{G},\mathbb{C}^{\times})$ "measures" how much choice there is to twist the associativity constrain.) One can further show that $\mathbf{Vec}_{\mathbb{C}}^{\omega}(\mathrm{G})$ is only monoidally equivalent to a skeletal category if $\omega$ is cohomologically trivial. $\diamond$

Again, we have:

**Lemma 2H.9.** *Any functor* $\mathrm{F}\in\mathbf{Hom}_{\otimes}(\mathbf{C},\mathbf{D})$ *induces a monoid homomorphism*

$$K_{0}(\mathrm{F})\colon K_{0}(\mathbf{C})\to K_{0}(\mathbf{D}),\ [\mathtt{X}]\mapsto[\mathrm{F}(\mathtt{X})].$$

*Further, if* $\mathrm{F}$ *is an equivalence, then* $K_{0}(\mathrm{F})$ *is an isomorphism.* $\square$

We leave it to the reader to define monoidal analogs of notions which we have seen in Section 1 (whenever appropriate), *e.g.* what a ***monoidal subcategory*** is. We only mention here that there are now three opposite categories (four, if one takes $\mathbf{C}$ itself into account):

**Definition 2H.10.** For any monoidal category $\mathbf{C}$, we define three additional monoidal categories

(2H-11)

| | $\mathbf{C}$ | $\mathbf{C}^{op}$ | $\mathbf{C}^{co}$ | $\mathbf{C}^{coop}$ |
|---|---|---|---|---|
| Reversed $\circ$? | No | Yes | No | Yes |
| Reversed $\otimes$? | No | No | Yes | Yes |

.

Using $^{op}$ is called taking the ***opposite***, *cf.* Definition 1B.13, taking $^{co}$ is called taking the ***coopposite***, and $\mathbf{C}^{coop}$ is called the ***biopposite*** of $\mathbf{C}$. $\diamond$

2I. **Strict vs. non-strict.** Let us start the comparison of strict and non-strict monoidal categories with a crucial example of a strict monoidal category, very much in the spirit of Example 2H.5.

**Definition 2I.1.** Given a monoidal category $\mathbf{C}$, define the category of ***right*** $\mathbf{C}$ ***module endofunctors***, denoted by $\mathbf{End}_{\mathbf{C}}(\mathbf{C})$, via:

- the objects are pairs $(\mathrm{F},\rho)$ with $\mathrm{F}\in\mathbf{End}(\mathbf{C})$ and natural isomorphisms $\rho_{\mathtt{X},\mathtt{Y}}\colon \mathrm{F}(\mathtt{X})\mathtt{Y}\to\mathrm{F}(\mathtt{XY})$ such that we have a commuting diagram

$$
\text{(2I-2)}\qquad
\begin{array}{ccccc}
 & & \mathrm{F}\big((\mathtt{XY})\mathtt{Z}\big) & & \\
 & \nearrow^{\mathrm{F}(\alpha_{\mathtt{X},\mathtt{Y},\mathtt{Z}})} & & \nwarrow^{\rho_{\mathtt{XY},\mathtt{Z}}} & \\
\mathrm{F}\big(\mathtt{X}(\mathtt{YZ})\big) & & & & \mathrm{F}(\mathtt{XY})\mathtt{Z} \\
 & \nwarrow_{\rho_{\mathtt{X},\mathtt{YZ}}} & & \nearrow_{\rho_{\mathtt{X},\mathtt{Y}}\otimes\mathrm{id}_{\mathtt{Z}}} & \\
 & \mathrm{F}(\mathtt{X})(\mathtt{YZ}) & \xrightarrow[\alpha_{\mathrm{F}(\mathtt{X}),\mathtt{Y},\mathtt{Z}}]{} & (\mathrm{F}(\mathtt{X})\mathtt{Y})\mathtt{Z} &
\end{array}
$$

for all $\mathtt{X},\mathtt{Y},\mathtt{Z}\in\mathbf{C}$;

- the morphisms $\alpha\colon(\mathrm{F},\rho)\Rightarrow(\mathrm{G},\rho')$ are natural transformations $\alpha\colon\mathrm{F}\Rightarrow\mathrm{G}$ such that we have a commuting diagram

$$
\text{(2I-3)}\qquad
\begin{array}{ccc}
\mathrm{G}(\mathtt{X})\mathtt{Y} & \xrightarrow{\rho'_{\mathtt{X},\mathtt{Y}}} & \mathrm{G}(\mathtt{XY}) \\
{\scriptstyle \alpha_{\mathtt{X}}\otimes\mathrm{id}_{\mathtt{Y}}}\uparrow & & \uparrow{\scriptstyle \alpha_{\mathtt{XY}}} \\
\mathrm{F}(\mathtt{X})\mathtt{Y} & \xrightarrow[\rho_{\mathtt{X},\mathtt{Y}}]{} & \mathrm{F}(\mathtt{XY})
\end{array}
$$

for all $\mathtt{X},\mathtt{Y}\in\mathbf{C}$;

- the composition $\circ$ is vertical composition of natural transformations.

(Note the pentagon again.) $\diamond$

**Lemma 2I.4.** *For $\mathbf{End}_{\bullet\mathbf{C}}(\mathbf{C})$ as in Definition 2I.1 the rules*

- *$\otimes$ on objects is $(\mathrm{G},\rho')(\mathrm{F},\rho)=(\mathrm{GF},\rho'')$, where*

$$
\rho''_{\mathtt{X},\mathtt{Y}}=\ \big(\mathrm{GF}(\mathtt{X})\big)\mathtt{Y}\xrightarrow{\rho'_{\mathrm{F}(\mathtt{X}),\mathtt{Y}}}\mathrm{G}\big(\mathrm{F}(\mathtt{X})\mathtt{Y}\big)\xrightarrow{\mathrm{G}(\rho_{\mathtt{X},\mathtt{Y}})}\mathrm{GF}(\mathtt{XY})\ ;
$$

- *$\otimes$ on morphisms is horizontal composition of natural transformations;*

*define the structure of a strict monoidal category on $\mathbf{End}_{\bullet\mathbf{C}}(\mathbf{C})$ with $\mathbb{1}=\mathrm{Id}_{\mathbf{C}}$.*

*Proof.* All appearing structures use compositions, either of maps, functors or of natural transformations, which are associative by definition. Thus, the only calculation one needs to check is that $\beta\alpha$ satisfies Equation 2I-3 if $\alpha$ and $\beta$ do. This is straightforward. $\square$

Comparing the definitions of a monoidal category (in particular, the $\triangle$ and the $\pentagon$ equations) and of a strict monoidal category, the following seems to be surprising.

**Theorem 2I.5.** *For any monoidal category $\mathbf{C}$ there exists a strict monoidal category $\mathbf{C}^{st}$ which is monoidally equivalent to $\mathbf{C}$, i.e. $\mathbf{C}\simeq_{\otimes}\mathbf{C}^{st}$.*

The statement of Theorem 2I.5 is called ***strictification***, and $\mathbf{C}^{st}$ is called a ***strictification*** of $\mathbf{C}$. Strictification allows us to very often "ignore" that we have to worry about associators and unitors. For example, we get diagrammatics for any monoidal category by passing to $\mathbf{C}^{st}$.

*Proof.* The idea is as follows. As a matter of fact, every monoid M is isomorphic to the monoid $\mathrm{End}_{\bullet\mathrm{M}}(\mathrm{M})$ consisting of maps from M to itself commuting with the right multiplication of M; the isomorphism is given by left multiplication. We will prove the theorem by copying this fact, *i.e.* we will show that $\mathbf{C}^{st}$ can be chosen to be $\mathbf{End}_{\bullet\mathbf{C}}(\mathbf{C})$.

By Lemma 2I.4 we have a strict monoidal category $\mathbf{End}_{\bullet\mathbf{C}}(\mathbf{C})$, which has a left action functor

$$
\mathrm{L}\colon\mathbf{C}\to\mathbf{End}_{\bullet\mathbf{C}}(\mathbf{C}),\quad \mathrm{L}(\mathtt{X})=(\mathtt{X}\otimes{}_{-},\alpha^{-1}_{\mathtt{X},{}_{-},{}_{-}}),\mathrm{L}(\mathrm{f})=\mathrm{f}\otimes{}_{-}.
$$

Note that Equation 2I-2 for L is the $\pentagon$ equation.

The functor L is an equivalence of categories, which we verify using Proposition 1H.24.

- The functor L is dense since any $(\mathrm{F},\rho)$ is isomorphic to $\mathrm{L}\big(\mathrm{F}(\mathbb{1})\big)$.
- The functor L is faithful, since $\mathrm{L}(\mathrm{f})=\mathrm{L}(\mathrm{g})$ implies $\mathrm{f}\otimes\mathrm{id}_{\mathbb{1}}=\mathrm{g}\otimes\mathrm{id}_{\mathbb{1}}$, which in turn gives $\mathrm{f}=\mathrm{g}$, by naturality of the unitor $\rho$. That is, commutativity of

$$
\begin{array}{ccc}
\mathtt{X}\mathbb{1} & \xrightarrow{\rho_{\mathtt{X}}} & \mathtt{X} \\
{\scriptstyle \mathrm{f}\otimes\mathrm{id}_{\mathbb{1}}}\downarrow & & \downarrow{\scriptstyle \mathrm{f}\otimes\mathrm{id}_{\mathbb{1}}} \\
\mathtt{Y}\mathbb{1} & \xrightarrow[\rho_{\mathtt{Y}}]{} & \mathtt{Y}
\end{array}
,\qquad
\begin{array}{ccc}
\mathtt{X}\mathbb{1} & \xrightarrow{\rho_{\mathtt{X}}} & \mathtt{X} \\
{\scriptstyle \mathrm{g}\otimes\mathrm{id}_{\mathbb{1}}}\downarrow & & \downarrow{\scriptstyle \mathrm{g}\otimes\mathrm{id}_{\mathbb{1}}} \\
\mathtt{Y}\mathbb{1} & \xrightarrow[\rho_{\mathtt{Y}}]{} & \mathtt{Y}
\end{array}
,
$$

with bottom and top being isomorphisms and $\mathrm{f}\otimes\mathrm{id}_{\mathbb{1}}=\mathrm{g}\otimes\mathrm{id}_{\mathbb{1}}$ implies $\mathrm{f}=\mathrm{g}$.

- Given a morphism $\alpha \in \mathbf{End}_{\mathbf{C}}(\mathbf{C})$, define a morphism

$$f = \mathtt{X} \xrightarrow{\rho_{\mathtt{X}}^{-1}} \mathtt{X}\mathbb{1} \xrightarrow{\alpha_{\mathbb{1}}} \mathtt{Y}\mathbb{1} \xrightarrow{\rho_{\mathtt{Y}}} \mathtt{Y} .$$

  Direct verification shows that $L(f) = \alpha$, thus L is full.

Finally, we define the structure of a monoidal functor on L via defining

- $\xi_{\mathtt{X},\mathtt{Y}} \colon \big(\mathtt{X} \otimes (\mathtt{Y} \otimes _), \mathrm{id}_{\mathtt{X}} \otimes \alpha^{-1}_{\mathtt{Y},_,_}\alpha^{-1}_{\mathtt{X},\mathtt{Y}_,_}\big) \xrightarrow{\cong} \big((\mathtt{X} \otimes \mathtt{Y}) \otimes _, \alpha^{-1}_{\mathtt{XY},_,_}\big)$ to be the associator $\alpha_{\mathtt{X},\mathtt{Y},_}$;
- $\xi_{\mathbb{1}} \colon (\mathrm{Id}_{\mathbf{C}}, \mathrm{id}) \xrightarrow{\cong} (\mathbb{1} \otimes _, \alpha^{-1}_{\mathbb{1},_,_})$ to be given by the inverse of the left unitor $\lambda$.

One verifies that this satisfies the axioms in Definition 2H.1. □

*Remark* 2I.6. Alternatively, and historically they first proof, Theorem 2I.5 can be proven using Theorem 2G.2, see *e.g.* [**ML98**, Section XI.3]. ◇

**Example 2I.7.** In most cases, like **Set** or $\mathbf{Vec}_{\Bbbk}$, the strictification $\mathbf{C}^{st}$ of **C** is the skeleton. Although formally incorrect (as we will see in this example), we still encourage the reader to think of strictification as going to the skeleton.

However, whenever $\cong$ is no-trivial one should not expect $\mathbf{C}^{st}$ to be skeletal. For example, $\mathbf{Vec}^{\omega}_{\mathbb{C}}(G)$ admits a strict skeleton if and only if $\omega$ is trivial. We will revisit this example diagrammatically later in Example 3G.6, and we will construct a strictification of $\mathbf{Vec}^{\omega}_{\mathbb{C}}(G)$ with infinitely many objects.

https://mathoverflow.net/questions/228715 ◇

2J. **More graphical calculus.** Recall the rules for diagrammatics of strict monoidal categories, *i.e.* Equation 2E-1 and Equation 2E-3. The formal rule of manipulation of these diagrams is:

(2J-1) "Two diagrams are equivalent if they are related by scaling or by a planar isotopy keeping the upwards orientation." .

**Theorem 2J.2.** *The graphical calculus is consistent, i.e. two morphisms are equal if and only if their diagrams are related by Equation 2J-1.*

*Proof.* This basically boils down to Equation 2E-2. □

**Example 2J.3.** Note that the condition of keeping the upwards orientation is a bit strange. In fact, it is probably not needed and can be dropped. The condition of only allowing planar isotopies is however crucial and *e.g.*

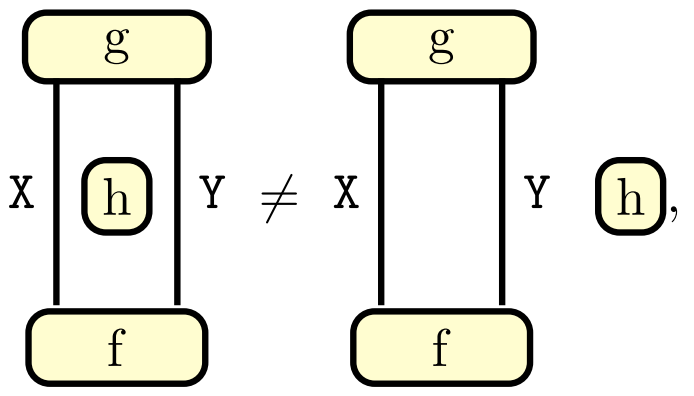

present different morphisms in general. ◇

Let us finish by showing the first hints why the diagrammatic calculus is very useful.

**Proposition 2J.4.** *For* $\mathbf{C} \in \mathbf{MCat}$ *the space* $\mathrm{End}_{\mathbf{C}}(\mathbb{1})$ *is a commutative monoid.*

*Proof.* The diagrammatic equation

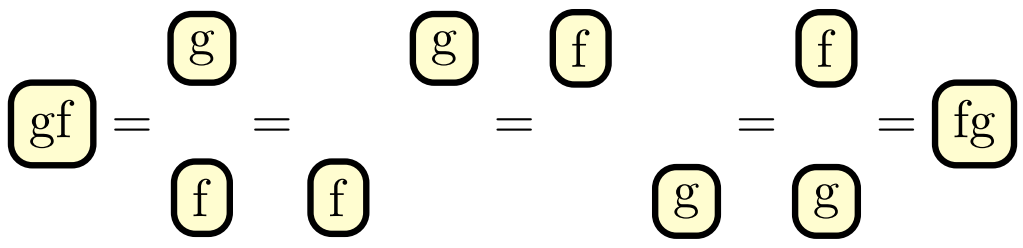

proves commutativity. □

**Proposition 2J.5.** *For* $\mathbf{C} \in \mathbf{MCat}$ *and for any* $\mathtt{X}, \mathtt{Y} \in \mathbf{C}$, *we have commuting actions*

$$\mathrm{End}_{\mathbf{C}}(\mathbb{1}) \circlearrowright \mathrm{Hom}_{\mathbf{C}}(\mathtt{X}, \mathtt{Y}) \circlearrowleft \mathrm{End}_{\mathbf{C}}(\mathbb{1}),$$

*given by*

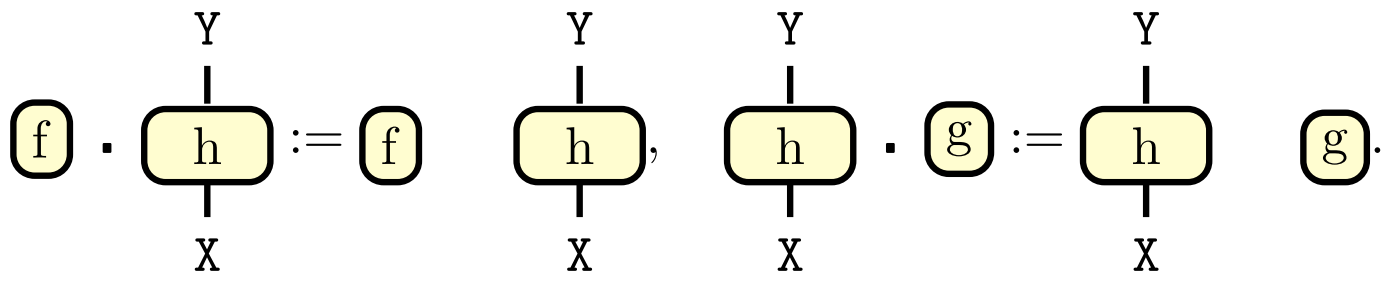

*Thus,* $\mathrm{Hom}_{\mathbf{C}}(\mathtt{X}, \mathtt{Y})$ *is an* $\mathrm{End}_{\mathbf{C}}(\mathbb{1})$*-bimodule.*

*Proof.* Associativity and unitality of the left action reads as

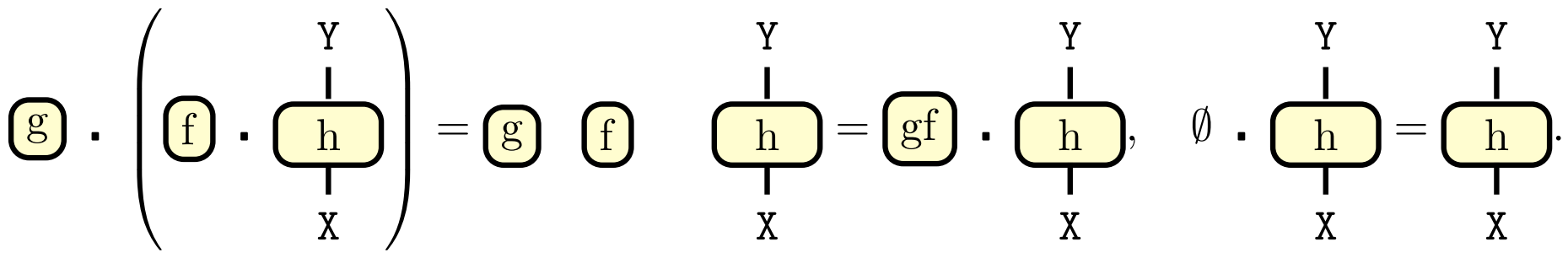

By reflecting the diagrams right to left, the same follows for the right action. Finally,

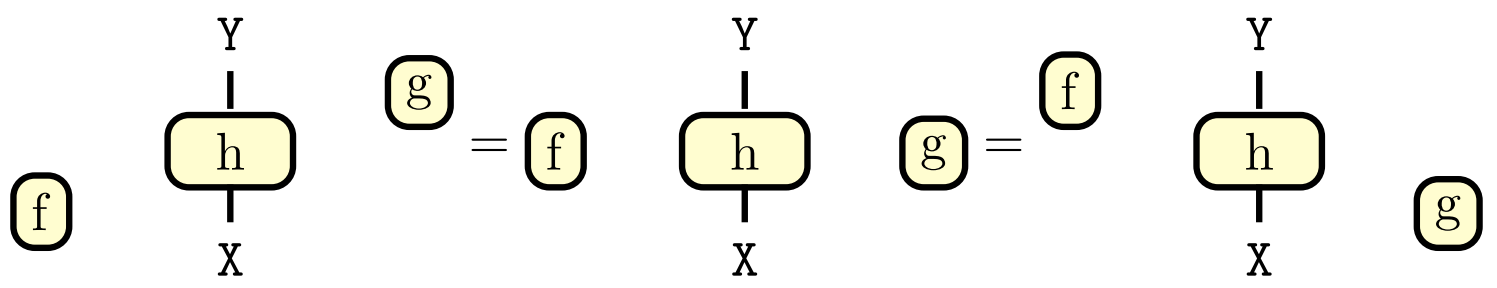

shows that the two actions commute. □

**Proposition 2J.6.** *The bimodule structure on* $\mathrm{Hom}_{\mathbf{C}}(\mathtt{X},\mathtt{Y})$ *from Proposition 2J.5 is compatible with* $\circ$ *and* $\otimes$.

*Proof.* This is Exercise 2K.5. □

## 2K. Exercises.

*Exercise* 2K.1. Explain explicitly what the four opposites from Equation 2H-11 are for the monoidal categories **1Cob**, **1Tan** and **1State**. ◇

*Exercise* 2K.2. Verify that $\mathbf{End}(\mathbf{C})$ is a strict monoidal category, *cf.* Example 2H.5. ◇

*Exercise* 2K.3. Show that $\mathbf{Vec}_{\mathbb{C}}(\mathbb{Z}/2\mathbb{Z}) \not\simeq_{\otimes} \mathbf{Vec}_{\mathbb{C}}^{\omega}(\mathbb{Z}/2\mathbb{Z})$ if $\omega$ is the nontrivial 3 cocycle. What happens for $\Bbbk = \mathbb{F}_2$ compared to $\Bbbk = \mathbb{C}$? ◇

*Exercise* 2K.4. Verify that $\mathbf{End}_{\mathbf{C}}(\mathbf{C})$ is a category, *cf.* Definition 2I.1. ◇

*Exercise* 2K.5. Prove Proposition 2J.6 diagrammatically. Hereby, compatibility means

$$\mathrm{f}\centerdot(\mathrm{jh}) = (\mathrm{f}\centerdot\mathrm{j})\mathrm{h} = \mathrm{j}(\mathrm{f}\centerdot\mathrm{h}), \quad \mathrm{f}\centerdot(\mathrm{h}\otimes\mathrm{j}) = (\mathrm{f}\centerdot\mathrm{h})\otimes\mathrm{j},$$

and *vice versa* for the right action. ◇

## 3. Monoidal categories II – more graphical calculus

The next question we will address, which will in particular give a rigorous, non-topological, construction of **1Cob**, **1Tan** and **1State**, is:

How to construct monoidal categories or algebraic objects using diagrammatic calculus?

The main player will be generators and relations presentations of monoidal categories, motivated the analog for groups, *cf.* Figure 7.

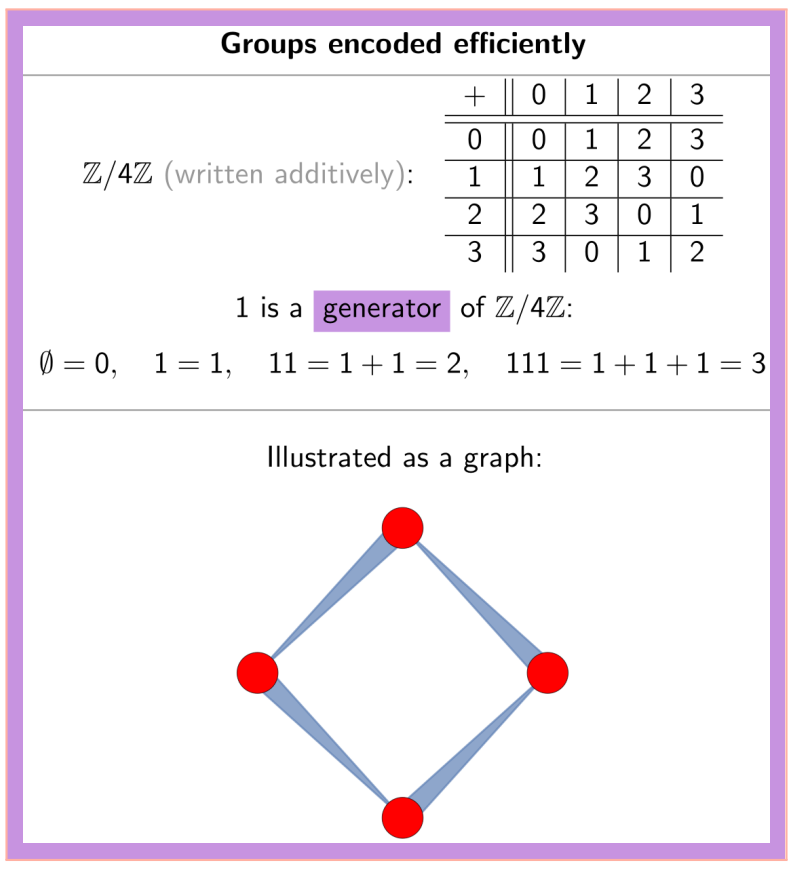

Figure 7. Generators and relations are well-studied for groups and an efficient way to encode the information.

Picture from the "What is...algebra?" playlist on [Tub21].

3A. **A word about conventions.** This section is all about the algebra of diagrams.

*Convention* 3A.1. Terminology that we will use several times is: "free as an ABC" means that no relation except the ones forced by "being an ABC" hold, and we will write "generated by XYZ" for short instead of "generated as an ABC by XYZ". Moreover, we sometimes do not define what "generated (as an ABC) by XYZ" means precisely as it will be clear from the context, see *e.g.* Example 3B.1.(b) for a non-defined phrase. ◇

*Convention* 3A.2. We usually simplify notation involving generators and relations as long as no confusion can arise. For example, all generators and relations will be elements of sets, but we omit the set brackets to make the notation less cumbersome. ◇

*Convention* 3A.3. All diagrammatics in this section are defined by generators and relations, in particular, not topologically. However, to simplify illustrations we draw diagrams sometimes in a topological fashion, and say some relations are mirrors of one another, *e.g.* the relations in Equation 3D-8 without mirrors are

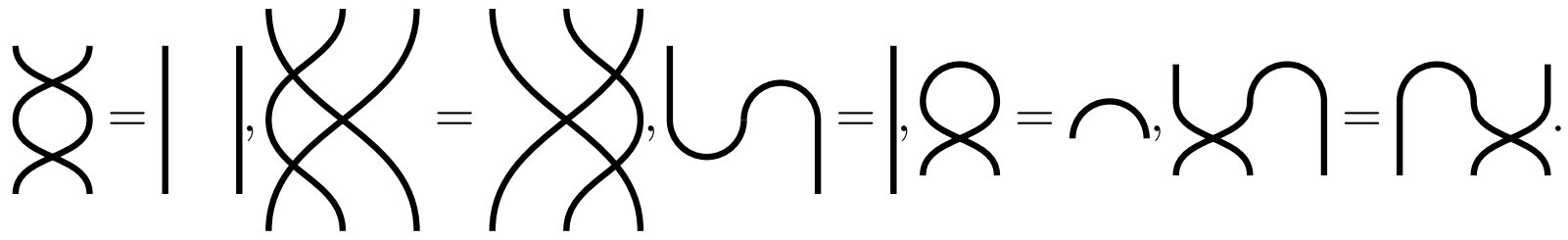

◇

*Convention* 3A.4. If certain notions only make sense under specific assumptions, then we tend to not to repeat these assumptions, *e.g.* we write "algebras" rather than "algebras in monoidal categories". ◇

3B. **Generator-relation presentation for monoids.** Recall the following constructions.

Given a set S, we obtain the ***free monoid generated (as a monoid) by*** S, denoted by $\langle \mathrm{S} \mid \emptyset\rangle$, by:

- elements are finite ***words*** $s_{i_r}...s_{i_1}$, where $s_{i_j} \in \mathrm{S}$ are the ***letters***, and $r \in \mathbb{Z}_{\geq 0}$;
- composition is concatenation of words;
- the unit is the empty word $\emptyset$;
- associativity is the only relation among words.

The elements of S are called ***generators*** (of $\langle \mathrm{S} \mid \emptyset\rangle$). The slogan is: 'free = minimal number of relations', and the only relation a monoid need is associativity.

Similarly, the ***free group generated (as a group) by*** S, also denoted by $\langle \mathrm{S} \mid \emptyset\rangle$, is defined *verbatim*, but having additional formal letters $s_i^{-1}$ satisfying $s_i s_i^{-1} = 1 = s_i^{-1} s_i$.

**Example 3B.1.** We stress that being free depends on the adjectives:

(a) The free monoid generated by $\mathrm{S} = \{\bullet\}$ is isomorphic to $\mathbb{N}$, while the free group generated by $\mathrm{S} = \{\bullet\}$ is isomorphic to $\mathbb{Z}$;

(b) The free commutative monoid generated by $\mathrm{S} = \{\bullet, *\}$ is isomorphic to $\mathbb{N}^2$, while the free monoid generated by $\mathrm{S} = \{\bullet, *\}$ is not isomorphic to $\mathbb{N}^2$.

(b) has the same analog for groups instead of monoids. ◇

Moreover, fix two sets S and R, where

$$\mathrm{R} \subset \langle \mathrm{S} \mid \emptyset\rangle \times \langle \mathrm{S} \mid \emptyset\rangle.$$

The elements of R will be written as $r = r'$ for $r, r' \in \langle \mathrm{S} \mid \emptyset\rangle$ and we call them ***relations***. We obtain the ***monoid generated by*** S ***with relations*** R, denoted by $\langle \mathrm{S} \mid \mathrm{R}\rangle$ as the quotient

$$\langle \mathrm{S} \mid \mathrm{R}\rangle = \langle \mathrm{S} \mid \emptyset\rangle/\mathrm{R},$$

meaning that two words in $\langle \mathrm{S} \mid \emptyset\rangle$ are equal in $\langle \mathrm{S} \mid \mathrm{R}\rangle$ if and only if they can be obtained from one another by using a finite number of relations from R. Said otherwise, we take the quotient by the congruence generated by the set R.

If $\mathrm{M} \cong \langle \mathrm{S} \mid \mathrm{R}\rangle$, then we say S | R give a ***generator-relation presentation of*** M.

**Example 3B.2.** Again, this depends on the adjectives:

(a) For $\mathrm{S} = \{\bullet\}$ and $\mathrm{R} = \{\bullet\bullet = 1\}$ we get $\langle \mathrm{S} \mid \mathrm{R}\rangle \cong \mathbb{Z}/2\mathbb{Z}$, regardless of whether we want to view this as being generated as a monoid or as a group.

(b) The symmetry group of the triangle, *i.e.* the dihedral group $\mathrm{I}_2(3)$ of order 6, has the generator-relation presentations

$$\mathrm{I}_2(3) \cong \langle a, b \mid a^2 = 1, b^3 = 1, aba = b^{-1}\rangle \cong \langle s, t \mid s^2 = 1, t^2 = 1, sts = tst\rangle,$$

where the middle expression is read to be as a group, while the right expression can be either as a monoid or a group.

It is often preferable to have presentations that avoid inverses. ◇

The set-theoretical issues of the following lemma are as usual ignored.

**Lemma 3B.3.** *Every monoid has a generator-relation presentation.*

*Proof.* Take $\mathrm{S} = \mathrm{M}$ and let R be given by all equations coming from multiplication in M. □

The presentation obtained via Lemma 3B.3 is, of course, useless. In general it is hard question whether one can find a good generator-relation presentation for a given monoid, group *etc.* But it is a good question, since we have the following evident, but useful, fact:

**Lemma 3B.4.** *To define a morphism* $f\colon \langle \mathrm{S} \mid \mathrm{R}\rangle \to \mathrm{M}$ *to any monoid* M *it suffices to*

- *specify* $f(s)$ *for* $s \in \mathrm{S}$*;*
- *check that* $f(r) = f(r')$ *for* $r = r' \in \mathrm{R}$. □

**Example 3B.5.** In the free case we have $\mathrm{R} = \emptyset$. Thus, every choice $f(s)$ for $s \in \mathrm{S}$ defines a morphism, regardless of M. ◇

3C. **Generator-relation presentation for monoidal categories.** We generalize the above:

**Definition 3C.1.** A set T is called a set of ***morphism generators*** if it consists of triples $(\mathrm{f}, \mathtt{X}, \mathtt{Y})$. Such a set is compatible with a set S if

$$\mathtt{X}, \mathtt{Y} \in \langle \mathrm{S} \mid \emptyset\rangle,$$

in which case we call S a set of ***object generators***. ◇

Of course we think of elements of S and T as being objects and morphisms $\mathrm{f}\colon \mathtt{X} \to \mathtt{Y}$, respectively.

**Definition 3C.2.** We define words as follows.

- An ***object word*** (in S) is a word in $\langle \mathrm{S} \mid \emptyset\rangle$, which can be concatenated as for monoids.
- A ***morphism word*** (in T) is defined recursively as follows. A morphism word of length 1 is either an element of T or of the form $(\mathrm{id}_{\mathtt{X}}, \mathtt{X}, \mathtt{X})$ for $\mathtt{X} \in \langle \mathrm{S} \mid \emptyset\rangle$. Suppose all morphism words of length $n \geq 1$ are already defined. A morphism word of length $n+1$ is either of the form
  - $(\mathrm{gf}, \mathtt{X}, \mathtt{Z})$ ($\circ$ concatenation),
  - $(\mathrm{f} \otimes \mathrm{h}, \mathtt{XA}, \mathtt{YB})$ ($\otimes$ concatenation),

  where $(\mathrm{f}, \mathtt{X}, \mathtt{Y})$, $(\mathrm{g}, \mathtt{Y}, \mathtt{Z})$ and $(\mathrm{h}, \mathtt{A}, \mathtt{B})$ are morphism words of length $n$. The two ways to create new words are also the two possible concatenations of morphism words.

We denote the collections of objects and morphism words by $\langle \mathrm{S} \mid \emptyset\rangle$ and $\langle \mathrm{T} \mid \emptyset\rangle_{\circ,\otimes}$. ◇

**Definition 3C.3.** Given sets S and T of object and morphism generators, we define the ***free strict monoidal category (monoidally generated) by*** S ***and*** T, denoted by $\boldsymbol{\langle \mathrm{S}, \mathrm{T} \mid \emptyset\rangle}$, as:

- the objects are $\langle \mathrm{S} \mid \emptyset\rangle$;
- the morphisms are $\langle \mathrm{T} \mid \emptyset\rangle_{\circ,\otimes}$;
- composition is $\circ$ concatenation of morphism words;
- the monoidal product is $\otimes$ concatenation of object respectively morphism words;
- the unit is $\mathbb{1} = 1 \in \langle \mathrm{S} \mid \emptyset\rangle$;
- the relations among object words are
$$\mathtt{X}(\mathtt{YZ}) = (\mathtt{XY})\mathtt{Z}, \quad \mathbb{1}\mathtt{X} = \mathtt{X} = \mathtt{X}\mathbb{1},$$
  where $\mathtt{X}, \mathtt{Y}, \mathtt{Z} \in \langle \mathrm{S} \mid \emptyset\rangle$;
- the relations among morphism words are
$$\mathrm{h}(\mathrm{gf}) = (\mathrm{hg})\mathrm{f}, \quad \mathrm{id}_{\mathtt{Y}}\mathrm{f} = \mathrm{f} = \mathrm{f}\mathrm{id}_{\mathtt{X}},$$
$$\mathrm{f} \otimes (\mathrm{g} \otimes \mathrm{h}) = (\mathrm{f} \otimes \mathrm{g}) \otimes \mathrm{h}, \quad \mathrm{id}_{\mathbb{1}} \otimes \mathrm{f} = \mathrm{f} = \mathrm{f} \otimes \mathrm{id}_{\mathbb{1}},$$
$$\mathrm{id}_{\mathtt{X}} \otimes \mathrm{id}_{\mathtt{Y}} = \mathrm{id}_{\mathtt{XY}}, \quad (\mathrm{gf}) \otimes (\mathrm{kh}) = (\mathrm{g} \otimes \mathrm{k})(\mathrm{f} \otimes \mathrm{h}),$$
  where $\mathtt{X}, \mathtt{Y} \in \langle \mathrm{S} \mid \emptyset\rangle$ and $\mathrm{f}, \mathrm{g}, \mathrm{h} \in \langle \mathrm{T} \mid \emptyset\rangle_{\circ,\otimes}$

◇

*Remark* 3C.4. The last two bullet points in Definition 3C.3 should be read as follows. The only relations among object words are those ensuring that $\langle \mathrm{S} \mid \emptyset\rangle$ is the free monoid generated by S. The only relations among morphism words are those ensuring that $\circ$ is a composition in a category and $\otimes$ is a bifunctor, *i.e.* Equation 2B-1, for a strict monoidal category. ◇

**Example 3C.5.** We have $\langle\{\bullet\},\emptyset \mid \emptyset\rangle \simeq_{\otimes} \mathbf{Vec}(\mathbb{Z}_{\geq 0})$, the latter being the evident adaption of Example 2D.11.(a) to the monoid $\mathbb{Z}_{\geq 0}$. ◇

As before we can choose

$$\mathrm{R} \subset \langle \mathrm{T} \mid \emptyset\rangle_{\circ,\otimes} \times \langle \mathrm{T} \mid \emptyset\rangle_{\circ,\otimes}.$$

The elements of R will be written as $\mathrm{r} = \mathrm{r}'$ for $\mathrm{r}, \mathrm{r}' \in \langle \mathrm{T} \mid \emptyset\rangle_{\circ,\otimes}$ and we call them ***relations***.

**Definition 3C.6.** We define the ***strict monoidal category generated by*** S ***and*** T ***with relations*** R, denoted by $\langle\mathbf{S},\mathbf{T} \mid \mathbf{R}\rangle$, as the quotient

$$\langle\mathbf{S},\mathbf{T} \mid \mathbf{R}\rangle = \langle\mathbf{S},\mathbf{T} \mid \emptyset\rangle/\mathrm{R},$$

meaning that two morphism words in $\langle\mathbf{S},\mathbf{T} \mid \emptyset\rangle$ are equal in $\langle\mathbf{S},\mathbf{T} \mid \mathbf{R}\rangle$ if and only if they can be obtained from one another by using a finite number of relations from R. ◇

The following is immediate.

**Lemma 3C.7.** *$\langle\mathbf{S},\mathbf{T} \mid \mathbf{R}\rangle$ with $\otimes$ being concatenation is a strict monoidal category.* □

**Definition 3C.8.** If $\mathbf{C} \simeq_{\otimes} \langle\mathbf{S},\mathbf{T} \mid \mathbf{R}\rangle$, then we say $\mathrm{S}, \mathrm{T} \mid \mathrm{R}$ give a ***generator-relation presentation of*** $\mathbf{C}$, where $\mathbf{C} \in \mathbf{MCat}$. ◇

The following two lemmas can be proven *verbatim* as for monoids, using beforehand Theorem 2I.5 for Lemma 3C.9 if the monoidal category of interest is not strict.

**Lemma 3C.9.** *Every monoidal category has a generator-relation presentation.* □

**Lemma 3C.10.** *To define a monoidal functor $\mathrm{F}\colon \langle\mathbf{S},\mathbf{T} \mid \mathbf{R}\rangle \to \mathbf{C}$ to a strict monoidal category $\mathbf{C} \in \mathbf{MCat}$ it suffices to*

- *specify* $\mathrm{F}(\mathtt{X})$ *for* $\mathtt{X} \in \mathrm{S}$;
- *specify* $\mathrm{F}(\mathrm{f})$ *for* $\mathrm{f} \in \mathrm{T}$;
- *check that* $\mathrm{F}(\mathrm{r}) = \mathrm{F}(\mathrm{r}')$ *for* $\mathrm{r} = \mathrm{r}' \in \mathrm{R}$. □

**Example 3C.11.** In the free case again any choice works, regardless of $\mathbf{C}$. ◇

3D. **Examples for generator-relation presentations.** Recall that we write $\mathtt{X}^k = \mathtt{X}...\mathtt{X}$.

*Remark* 3D.1. All diagrams we will use below are not topological objects, but rather formal symbols. However, as we will see, one should think of them as being topological objects, see also Convention 3A.3. ◇

**Example 3D.2.** The ***category of symmetric groups*** $\mathbf{Sym}$ can be defined as follows. We let $\mathbf{Sym} = \langle\mathbf{S},\mathbf{T} \mid \mathbf{R}\rangle$ with

(3D-3) $\mathrm{S}: \bullet,\quad \mathrm{T}:$ 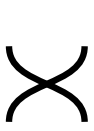 $: \bullet^2 \to \bullet^2,\quad \mathrm{R}:$ 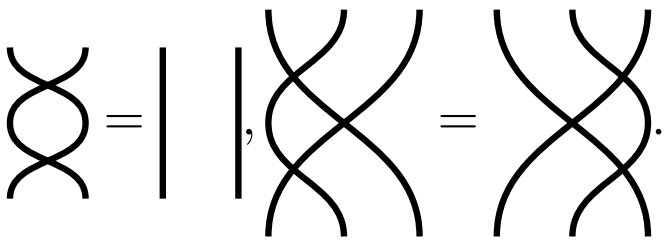

The monoid $\mathrm{End}_{\mathbf{Sym}}(\bullet^n)$ is called the ***symmetric group of*** $\{1,\dots,n\}$. ◇

*Remark* 3D.4. There are extra relations which are implicit, *e.g.*

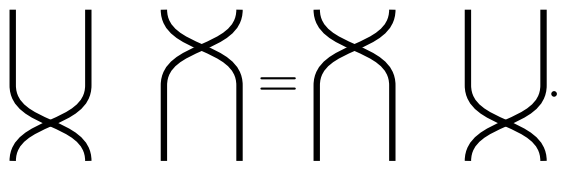

We do not need to add this relation to Equation 3D-3 since it follows from the interchange law, *cf.* Equation 2E-2. ◇

It is a (non-trivial) fact that $\mathbf{Sym} \simeq_{\otimes} \mathbf{Sym}^{top} \subset \mathbf{1Cob}$, see Example 1H.22, and the above can be seen as a purely algebraic construction of the (topologically defined) category $\mathbf{Sym}^{top}$.

Let $\mathrm{Aut}(\{1,\dots,n\})$ denote all bijections $\{1,\dots,n\} \to \{1,\dots,n\}$ with composition as the group structure (this is the usual symmetric group whose elements are permutations), and let $(i,i+1)$ be the simple transposition that swaps $i, i+1 \in \{1,\dots,n\}$.

**Lemma 3D.5.** *The map (the crossing is at the $i$ and $(i+1)$th strands)*

$\phi\colon \mathrm{End}_{\mathbf{Sym}}(\bullet^n) \to \mathrm{Aut}(\{1,\dots,n\}),$ 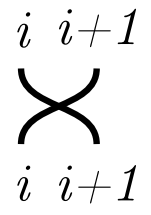 $\mapsto (i,i+1) \in \mathrm{Aut}(\{1,\dots,n\})$

*is a monoid isomorphism. In particular,* $\mathrm{End}_{\mathbf{Sym}}(\bullet^n)$ *is a group.*

*Proof.* By Lemma 3C.10 it suffices to observe that the defining relations of **Sym** hold in $\mathrm{Aut}(\{1,...,n\})$ to see that $\phi$ is a monoid map. The map $\phi$ is surjective since the simple transpositions generate $\mathrm{Aut}(\{1,...,n\})$. Injectivity is the most difficult part. To see this, let $\psi\colon \mathrm{Aut}(\{1,...,n\}) \to \mathrm{Aut}(\mathbb{C}[x_1,...,x_n])$ be the map that sends every permutation $\sigma$ to the operator on $\mathbb{C}[x_1,...,x_n]$ that permutes the variables according to $\sigma$. The map $\psi\circ\phi$ is injective, so $\phi$ is also injective. □

**Example 3D.6.** The ***(generic) Rumer–Teller–Weyl category*** **TL** (also known as the generic Temperley–Lieb category, hence the notation) is defined as follows. We let

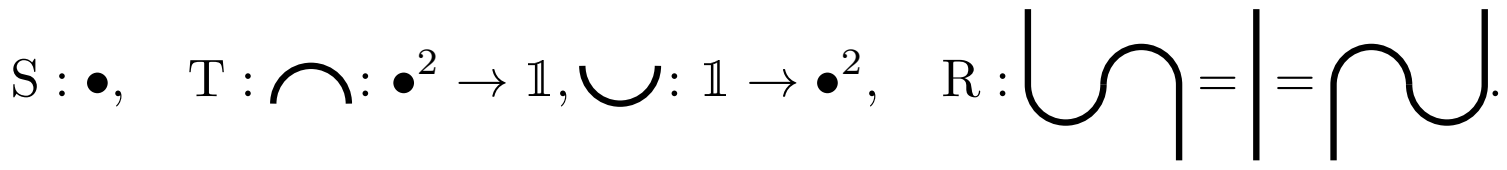

The monoid $\mathrm{End}_{\mathbf{TL}}(\bullet^n)$ is called the ***Temperley–Lieb monoid***. ◇

Again, its non-trivial, but visually clear, that **TL** is a monoidal subcategory of **1Cob**.

**Example 3D.7.** The ***(generic) Brauer category*** **Br** is defined as follows. We let

(3D-8)

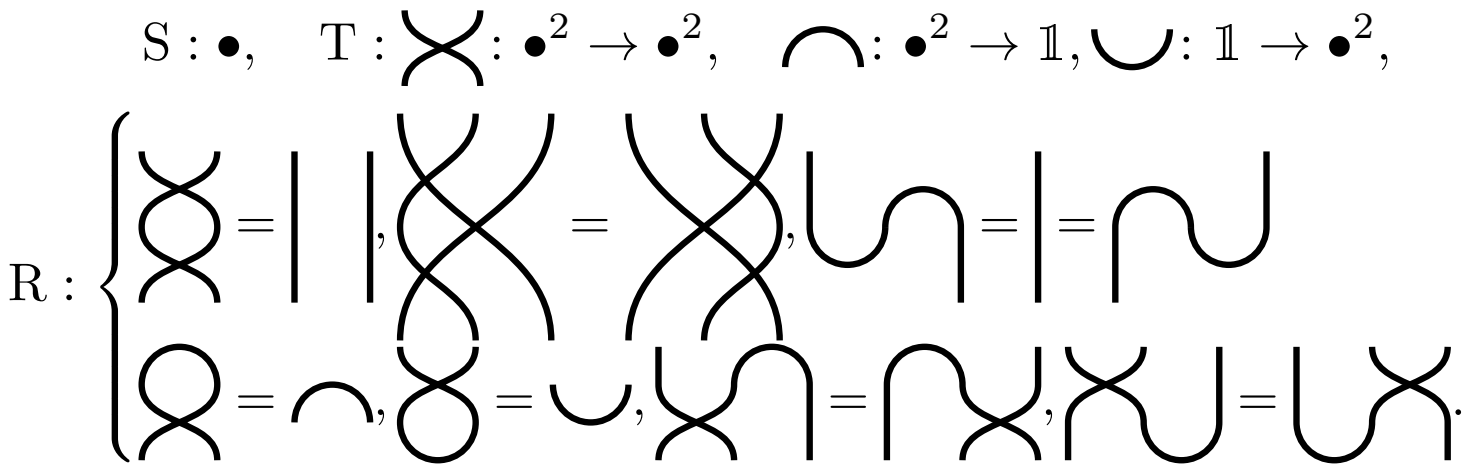

The monoid $\mathrm{End}_{\mathbf{Br}}(\bullet^n)$ is called the ***Brauer monoid***. ◇

**Example 3D.9.** With Lemma 3C.10 it is easy to define monoidal functors

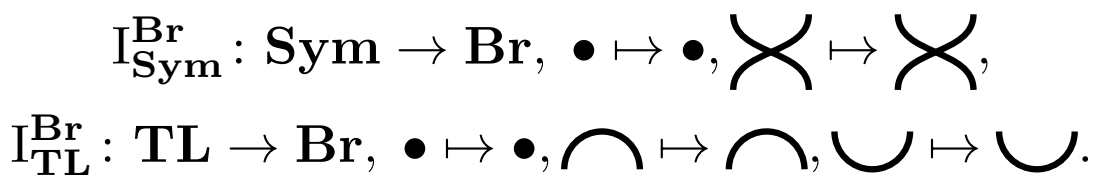

These are dense by construction, and with a bit more work one can show that they are faithful. Thus, **Sym** and **TL** are (non-full) monoidal subcategories of **Br**. ◇

The punchline is that $\mathbf{Br} \simeq_{\otimes} \mathbf{1Cob}$. Let us state this as a theorem, whose proof we will sketch, highlighting what is easy and what is non-trivial about this statement.

**Theorem 3D.10.** *There exists a monoidal functor*

(3D-11)

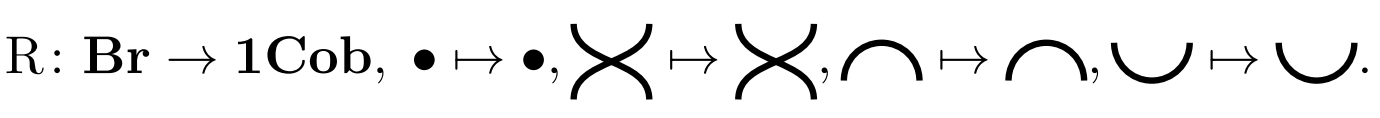

*The functor* R *is dense and fully faithful, thus,* $\mathbf{Br} \simeq_{\otimes} \mathbf{1Cob}$.

Crucial: the left diagrams in Equation 3D-11 are just algebraic symbols, while the right diagrams are just placeholder symbols for topological objects.

*Proof.* There are several things to check, namely:

(a) The functor R is a well-defined monoidal functor. Using Lemma 3C.10, this is just the observation that the Brauer relations Equation 3D-8 hold in **1Cob**, which is easy.

(b) That R is dense is clear.

(c) That R is full is not hard: Every 1 dimensional cobordism in **1Cob** has locally a Morse point, or not. But Morse points in this situation are just caps or cups. Moreover, taking the immersion (into the plane) of the cobordism into account, locally a 1 dimensional cobordism in **1Cob** is of the form

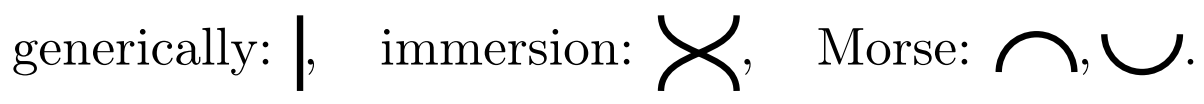

In particular, we have a Morse positioning of such cobordisms. Here is an example:

(3D-12) 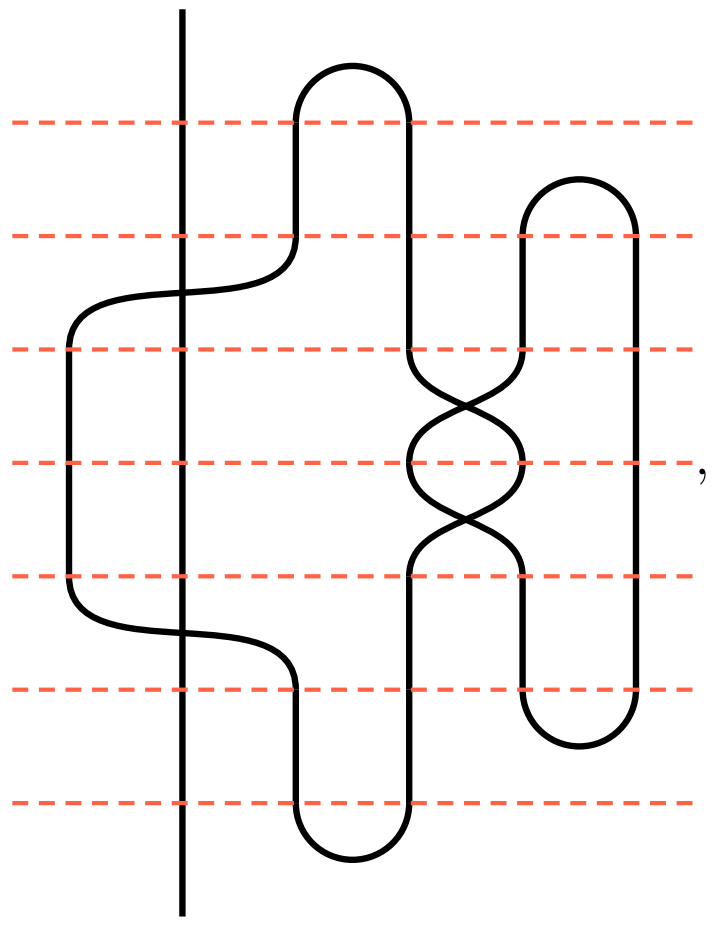 ,

where the horizontal and dashed lines indicate height levels. Said otherwise, crossings, caps and cups generate **1Cob**, so R is full as all generators of **1Cob** appear in its image. (Note that already here one would need to be precise what one means by a "1 dimensional cobordism". But this is not the hard part.)

(d) The proof that R is faithful is hard and painful, because one needs to show that the topologically defined **1Cob** has the Brauer relations Equation 3D-8 as generating relations. (See also Exercise 3H.2.)

The proof is complete. $\square$

*Remark* 3D.13. By the same reason as in (d) in the proof of Theorem 3D.10, it is hard to write down any functor $\mathbf{1Cob} \to \mathbf{Br}$. That is, the inverse of R is of course

$\mathrm{R}^{-1}\colon \mathbf{1Cob} \to \mathbf{Br}$, $\bullet \mapsto \bullet$, 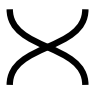 $\mapsto$ 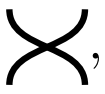, 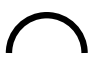 $\mapsto$ 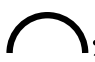, 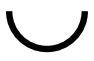 $\mapsto$ 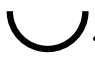.

But showing that this is well-defined boils down to (d). $\diamond$

There are several variations of the Brauer category **Br**, *e.g.* with orientations, which we will revisit later. For now we are brief:

**Example 3D.14.** The ***(generic) quantum Brauer category* qBr**, the ***(generic) oriented Brauer category* oBr** and the ***(generic) oriented quantum Brauer category* oqBr** are defined *verbatim* as the Brauer category, with a few differences:

- the adjective "oriented" means that one has two object generators $\bullet$ and $\bullet^{\star}$ and one has oriented caps and cups generators, *i.e.*

(3D-15) 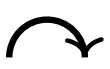$\colon \bullet\bullet^{\star} \to \mathbb{1}$, 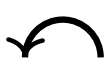$\colon (\bullet^{\star})\bullet \to \mathbb{1}$, 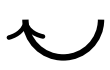$\colon \mathbb{1} \to \bullet\bullet^{\star}$, 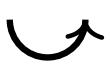$\colon \mathbb{1} \to (\bullet^{\star})\bullet$;

- the adjective "quantum" means that one distinguishes over- and undercrossings, *i.e.*

(3D-16) overcrossing: 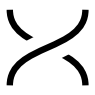, undercrossing: 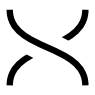.

All of these have analogs of Theorem 3D.10, *e.g.* $\mathbf{qBr} \simeq_{\otimes} \mathbf{1Tan}$, and monoid versions. $\diamond$

**Example 3D.17.** We also have **qSym**, the ***category of braids***, being the analog of **Sym** with crossing as in Equation 3D-16, but as a subcategory of **qBr**. Similarly, we also have **oTL**, the oriented version of **TL**, with oriented diagrams as in Equation 3D-15. $\diamond$

3E. **Algebras in monoidal categories.** Next, we aim to generalize the notion of an algebra.

**Definition 3E.1.** An ***algebra*** $\mathtt{A} = (\mathtt{A}, \mathrm{m}, \mathrm{i})$ in a monoidal category **C** consists of

- an object $\mathtt{A} \in \mathbf{C}$;
- a ***multiplication***, *i.e.* a morphism $\mathrm{m}\colon \mathtt{AA} \to \mathtt{A}$;
- a ***unit***, *i.e.* a morphism $\mathrm{i}\colon \mathbb{1} \to \mathtt{A}$;

such that

(i) we have a commuting diagram

(3E-2)
$$\begin{array}{ccccc} \mathtt{A(AA)} & & \xrightarrow{\alpha_{\mathtt{A},\mathtt{A},\mathtt{A}}} & & \mathtt{(AA)A} \\ \downarrow{\scriptstyle \mathrm{id}_{\mathtt{A}}\otimes\mathrm{m}} & & & & \downarrow{\scriptstyle \mathrm{m}\otimes\mathrm{id}_{\mathtt{A}}} \\ \mathtt{AA} & \xrightarrow[\mathrm{m}]{} & \mathtt{A} & \xleftarrow[\mathrm{m}]{} & \mathtt{AA} \end{array};$$

(ii) we have commuting diagrams

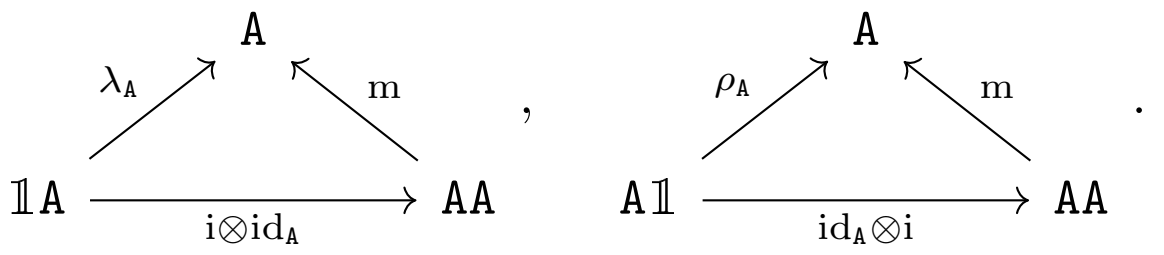

A ***coalgebra*** $\mathtt{C} = (\mathtt{C}, \mathrm{d}, \mathrm{e})$ in a monoidal category $\mathbf{C}$ is, by definition, an algebra in $\mathbf{C}^{op}$. ◇

*Remark* 3E.3. There is a bias and algebras seem to be much more beloved than coalgebras. However, Theorem 1C.1 applies and they are "essentially the same". ◇

Of course, the three commuting diagrams in Definition 3E.1 are associativity (recall that Equation 3E-2 implies honest associativity) and unitality, but in the context of (not necessary strict) monoidal categories. Similarly for coalgebras. The Feynman diagrams can be simplified:

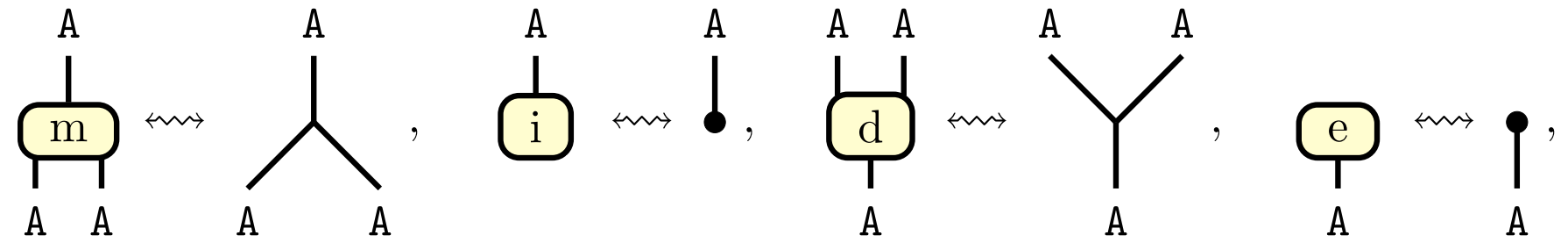

are the structure maps and

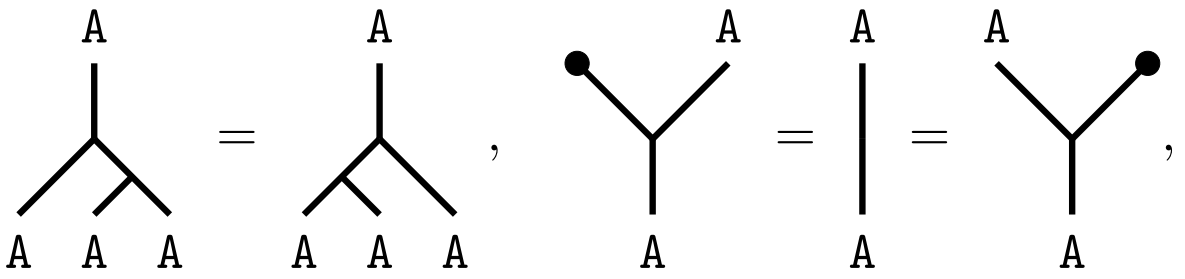

present associativity and counitality, respectively.

**Example 3E.4.** In any monoidal category $\mathbb{1}$ has the structure of a (co)algebra. ◇

**Example 3E.5.** Algebras and coalgebras generalize many notions:

(*a*) Algebras and coalgebras in **Set** are monoids and comonoids.

(*b*) Algebras and coalgebras in $\mathbf{Vec}_{\Bbbk}$ are algebras and coalgebras over $\Bbbk$.

In (b), the easiest algebra is "plain multiplication", that is $\mathbb{C}$ with $\mathrm{m}(a \otimes b) = ab$ being multiplication of complex numbers and the unit $\mathrm{i}(1) = 1$ being 1.

Coalgebras on the other hand tend to be generous: they distribute evenly. The easiest example of this type, again in (b), is $\mathbb{C}[x]$ with $\mathrm{d}(x^n) = \sum_{k=0}^{n} x^k \otimes x^{n-k}$, for $n \in \mathbb{Z}_{\geq 0}$ with $x^0 = 1$, and $\mathrm{e}(x^n) = 0$ for $n > 0$ and $\mathrm{e}(x^0) = 1$. ◇

**Example 3E.6.** Note that the final example in Example 3E.5 is also an algebra, using plain multiplication. The 2 dimensional quotient $\mathrm{D} = \mathbb{C}[x]/(x^2)$ inherits both structures, and is one of the most crucial examples: the so-called ***dual numbers***.

Relabeling $x$ to $\epsilon$ so that $\mathrm{D} = \mathbb{C}[\epsilon]/(\epsilon^2)$ gives is the historically correct notation: these were introduced in the 19th century and $\epsilon$ can be thought of as an infinitesimal, so that D mimics linear approximation as we "ignore" quadratic terms.

It comes up from at least three perspectives, justifying its importance:

**(a)** The dual numbers are the cohomology ring of the sphere (equivalently, of $\mathbb{P}^1$).

**(b)** The dual numbers are the exterior algebra of $\mathbb{C}$.

**(c)** The dual numbers is the only 2 dimensional non-semisimple algebra.

We will prove (c) later in Proposition 7E.4. ◇

**Example 3E.7.** The object $\bullet^2 \in \mathbf{Br}$ (the Brauer category, see Example 3D.7) is an algebra with structure maps

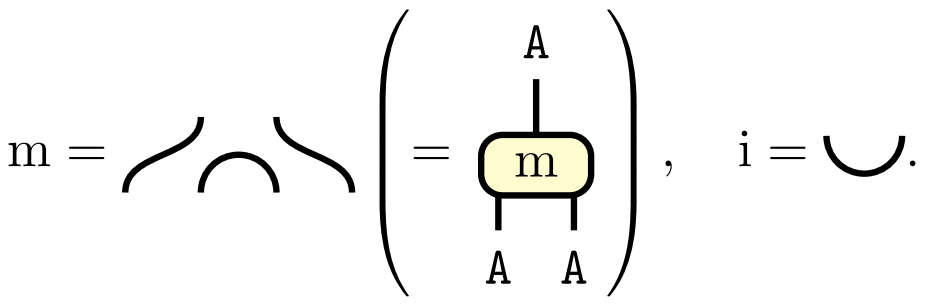

Associativity and unitality are topologically clear:

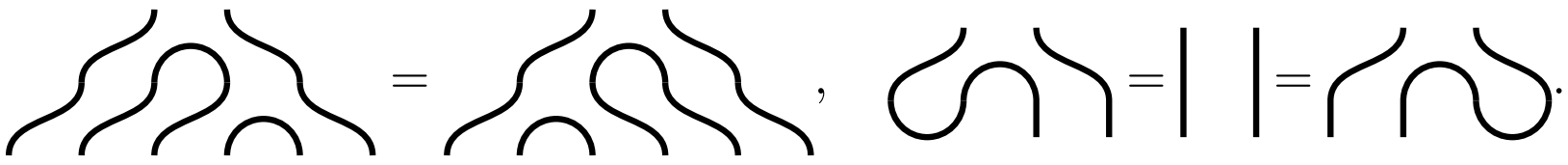

Similarly, the object $\bullet^2 \in \mathbf{Br}$ is also a coalgebra, by mirroring the diagrams. ◇

**Definition 3E.8.** A ***Frobenius algebra*** $\mathtt{A} = (\mathtt{A}, \mathrm{m}, \mathrm{i}, \mathrm{d}, \mathrm{e})$ in $\mathbf{C}$ is an algebra $(\mathtt{A}, \mathrm{m}, \mathrm{i})$ and a coalgebra $(\mathtt{A}, \mathrm{d}, \mathrm{e})$ in $\mathbf{C}$ satisfying a compatibility condition, *i.e.* we have commuting diagrams

(3E-9)
$$\begin{array}{ccccc} \mathtt{A}(\mathtt{AA}) & & \xrightarrow{\alpha_{\mathtt{A},\mathtt{A},\mathtt{A}}} & & (\mathtt{AA})\mathtt{A} \\ {\scriptstyle \mathrm{id}_{\mathtt{A}}\otimes \mathrm{d}}\uparrow & & & & \downarrow{\scriptstyle \mathrm{m}\otimes \mathrm{id}_{\mathtt{A}}} \\ \mathtt{AA} & \xrightarrow[\mathrm{m}]{} & \mathtt{A} & \xrightarrow[\mathrm{d}]{} & \mathtt{AA} \end{array}, \qquad \begin{array}{ccccc} (\mathtt{AA})\mathtt{A} & & \xrightarrow{\alpha^{-1}_{\mathtt{A},\mathtt{A},\mathtt{A}}} & & \mathtt{A}(\mathtt{AA}) \\ {\scriptstyle \mathrm{d}\otimes \mathrm{id}_{\mathtt{A}}}\uparrow & & & & \downarrow{\scriptstyle \mathrm{id}_{\mathtt{A}}\otimes \mathrm{m}} \\ \mathtt{AA} & \xrightarrow[\mathrm{m}]{} & \mathtt{A} & \xrightarrow[\mathrm{d}]{} & \mathtt{AA} \end{array}.$$

These are the Frobenius conditions. ◇

Diagrammatically Equation 3E-9 is

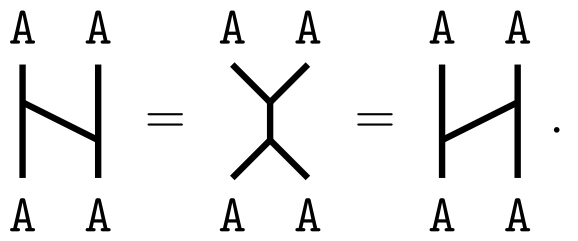

**Example 3E.10.** Frobenius algebras in $\mathbf{Vec}_{\Bbbk}$ are classical Frobenius algebras over $\Bbbk$. The dual numbers D with the structure morphisms as in Example 3E.6 is an example for $\Bbbk = \mathbb{C}$. ◇

**Lemma 3E.11.** *Let* $\mathtt{A}$ *be a Frobenius algebra in a strict monoidal category. Define*

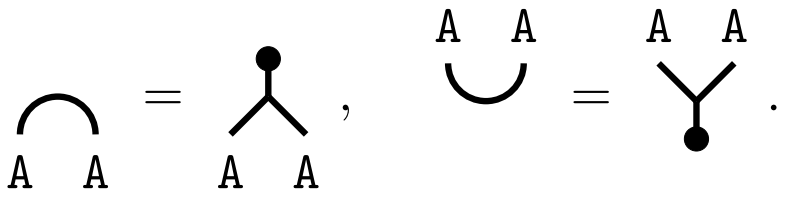

*Then the following hold, including mirrors:*

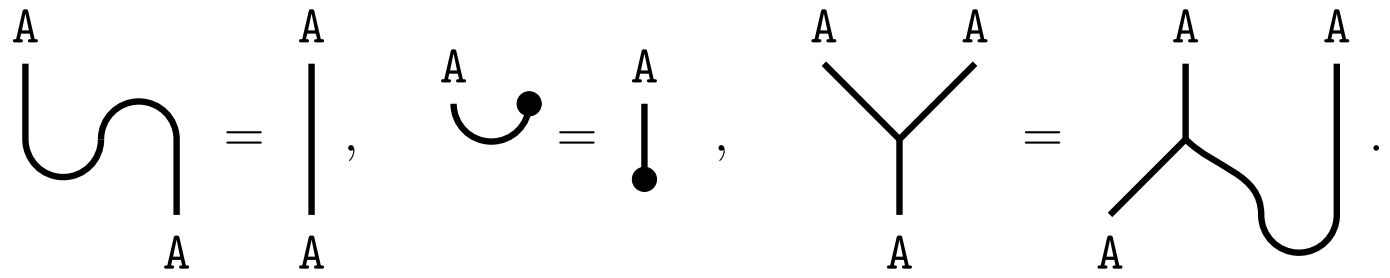

Thus, Frobenius algebras are topologically in nature since Lemma 3E.11 shows that the diagrams for Frobenius algebras satisfy all planar isotopies.

*Proof.* This is Exercise 3H.3 □

*Remark* 3E.12. There is of course also a non-strict version of Lemma 3E.11 which looks almost exactly the same. ◇

## 3F. Modules of algebras.

Arguably modules of algebras are more interesting than the algebras themselves. For example, replacing algebras by groups, abstract groups formalize the concept of symmetry via their actions on spaces (and such a space is then a module of the group). For example, the symmetric group $\mathrm{Aut}(\{1, 2, 3, 4\})$ arises in the wild as the symmetries of the tetrahedron in the same way as the number 3 arises in the wild as "three apples":

| | Abstract | Incarnation |
|---|---|---|
| Numbers | 3 | 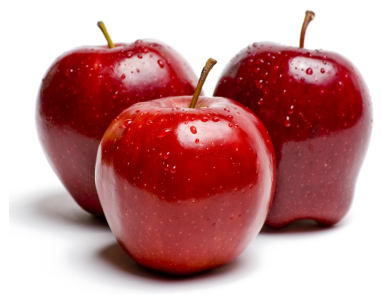 or... |
| Groups | $S_4 = \langle s, t, u \mid \text{some relations}\rangle$ | 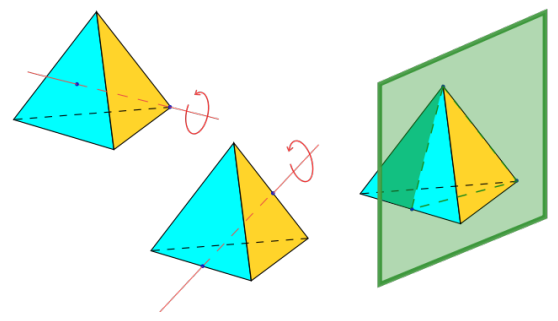 or... |

So:

**Definition 3F.1.** Let $\mathtt{A}$ be an algebra. A ***right*** $\mathtt{A}$ ***module*** $\mathtt{M} = (\mathtt{M}, \triangleleft)$ in $\mathbf{C}$ consists of

- an object $\mathtt{M} \in \mathbf{C}$;
- a ***right action***, *i.e.* a morphism $\triangleleft\colon \mathtt{MA} \to \mathtt{M}$;

such that

(i) we have a commuting diagram

$$\begin{array}{ccccc} \mathtt{M}(\mathtt{AA}) & & \xrightarrow{\alpha_{\mathtt{M},\mathtt{A},\mathtt{A}}} & & (\mathtt{MA})\mathtt{A} \\ {\scriptstyle \mathrm{id}_{\mathtt{M}}\otimes \mathrm{m}}\downarrow & & & & \downarrow{\scriptstyle (\centerdot_)\otimes \mathrm{id}_{\mathtt{A}}} \\ \mathtt{MA} & \xrightarrow[\centerdot_]{} & \mathtt{M} & \xleftarrow[\centerdot_]{} & \mathtt{MA} \end{array};$$

(ii) we have a commuting diagram

$$\begin{array}{ccc} & \mathtt{M} & \\ {\scriptstyle \rho_{\mathtt{M}}}\nearrow & & \nwarrow{\scriptstyle \centerdot_} \\ \mathtt{M}\mathbb{1} & \xrightarrow[\mathrm{id}_{\mathtt{M}}\otimes \mathrm{i}]{} & \mathtt{MA} \end{array}.$$

Here right ⇝ "A comes from the right". ◇

By duality, we also get the notions of ***right* A *comodules***, ***left* A *modules*** and ***left* A *comodules***.

**Example 3F.2.** In $\mathbf{Vec}_{\Bbbk}$, right A modules are the classical notion of algebra modules. Explicitly, for $\Bbbk = \mathbb{C}$ and cyclic groups (a.k.a. rotational symmetries), say $\mathbb{Z}/3\mathbb{Z} = \{g, h = g^2, 1 = g^3\}$, we can let $\mathtt{A} = \mathbb{C}[\mathbb{Z}/3\mathbb{Z}]$ and take $\mathtt{M} = \mathbb{C}^3$ with the following action. Call the basis vectors of $\mathbb{C}^3$ "north east", "north west" and "south". The action is then determined by

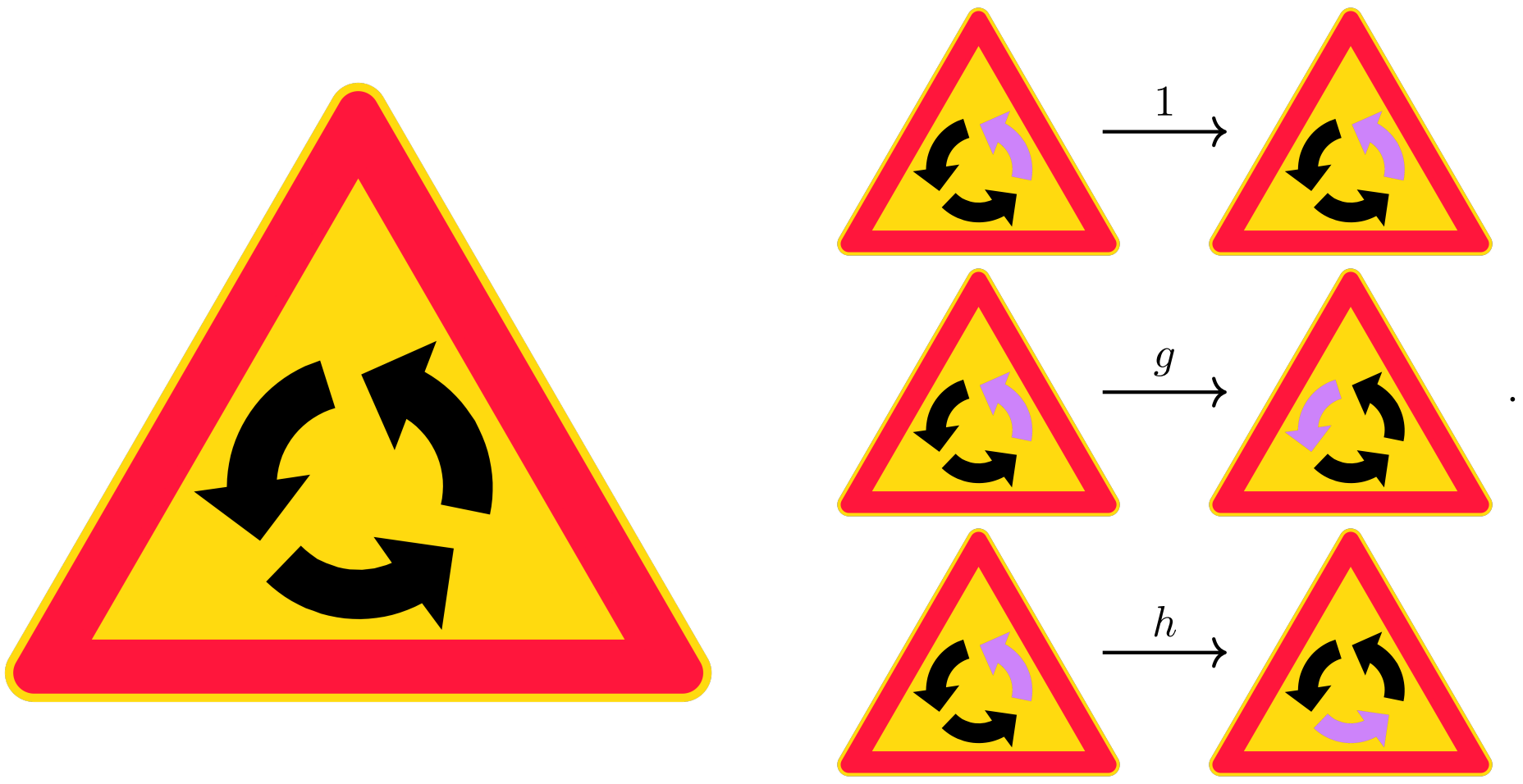

In this sense, right A modules in $\mathbf{Vec}_{\Bbbk}$ are symmetries in vector spaces, and Definition 3F.1 generalizes this to monoidal categories. ◇

**Definition 3F.3.** Let M and N be right A modules. A morphism $\mathrm{f}\colon \mathtt{M} \to \mathtt{N}$ is said to be A ***equivariant*** if it intertwines the right A action, *i.e.* $\mathrm{f}(\centerdot_) = (\centerdot_)\mathrm{f}$. ◇

In pictures these notions are again nice:

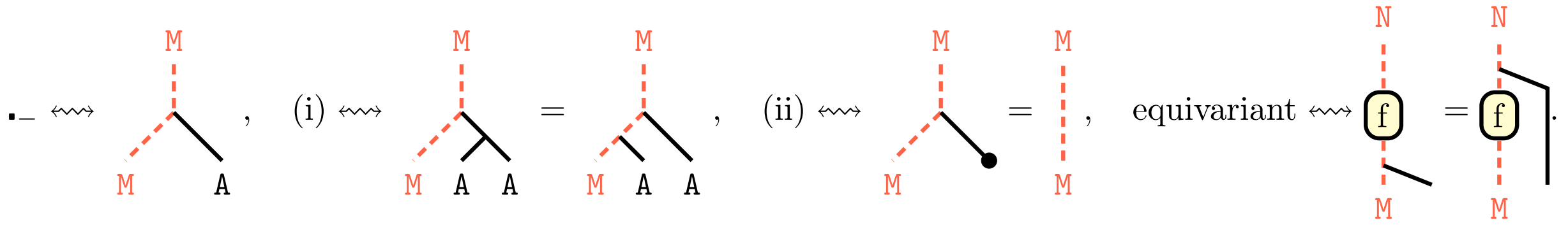

**Lemma 3F.4.** *The composition of* A *equivariant morphisms is* A *equivariant.*

*Proof.* The proof in diagrams is easy:

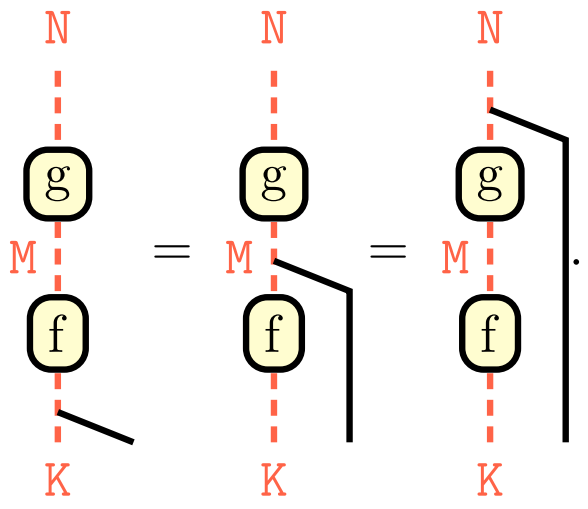

The non-strict version is thus also true. □

Since the identity is always A equivariant, we get another category:

**Example 3F.5.** We have a category $\mathbf{Mod}_{\mathbf{C}}(\mathtt{A})$, the ***category of right*** A ***modules***, whose objects are right A modules and morphisms are A equivariant morphisms. ◇

**Example 3F.6.** All of these notions generalize the classical notions of algebras, modules and their categories if we work in $\mathbf{Vec}_{\mathbb{k}}$. ◇

We leave it to the reader to write down the definitions of other classical notions from algebra in the categorical sense, see also Exercise 3H.4. Let us instead wrap-up this part with a diagrammatic proof generalizing a classical fact which is actually messy to prove classically.

**Proposition 3F.7.** *Let* A *be a Frobenius algebra in* **C**. *Then every right* A *module has a compatible structure of a right* A *comodule and vice versa.*

*Proof.* We can assume that **C** is strict. Let M be a right A module. Then we define the coaction $(_)^{co}$ via

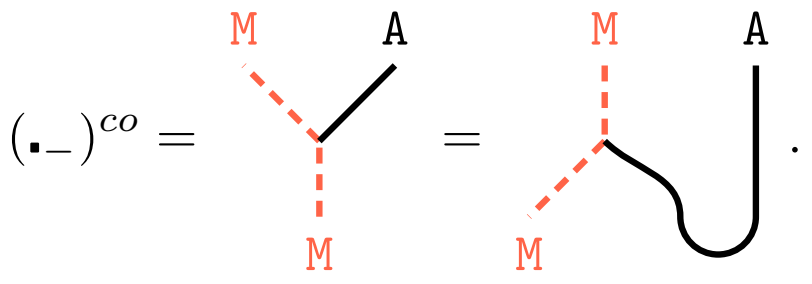

This defines a right A comodule since, by Lemma 3E.11, we have *e.g.*

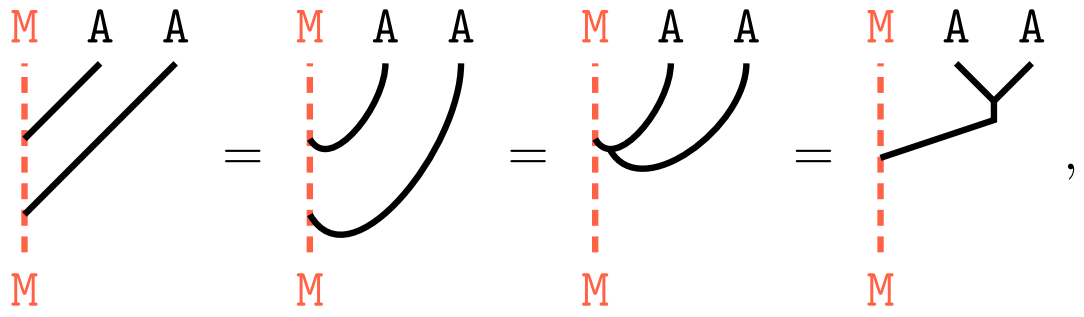

and unitality follows *mutatis mutandis*. By mirroring the diagrams we can get from comodules to modules. □

3G. **Web categories.** The Rumer–Teller–Weyl category from Example 3D.6 can be viewed as the free category generated by a non-degenerate bilinear form and coform.

We now describe categories generated by non-degenerate trilinear (co)forms, following the exposition in [**Tub23**]. Before reading Example 3G.2 we point out that, using isotopy, one can generate many new morphisms. For example,

(3G-1) 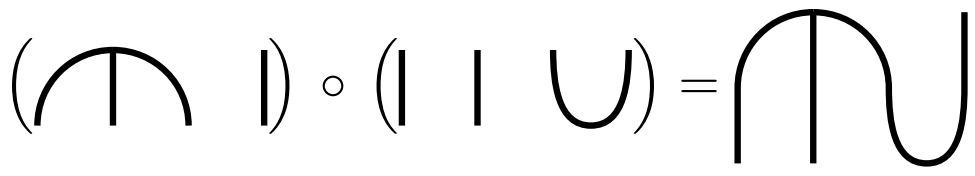 is isotopic to 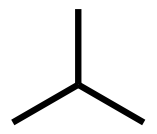 .

We use this to define picture as on the right-hand side in Equation 3G-1.

**Example 3G.2.** The ***(generic) web category*** **Web** (also known as the generic ***trivalent*** category or ***spider*** category or generic ***birdtracks*** category or ...) is defined as follows. We let

S : •, T : 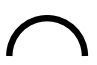: $\bullet\otimes\bullet\to\mathbb{1}$, 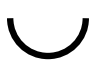: $\mathbb{1}\to\bullet\otimes\bullet$, 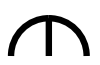: $\bullet\otimes\bullet\otimes\bullet\to\mathbb{1}$, 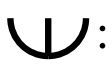: $\mathbb{1}\to\bullet\otimes\bullet\otimes\bullet$,

and the generating morphisms are called ***bilinear*** and ***trilinear (co)form***, respectively. Using Equation 3G-1 and similar expressions to define

R : 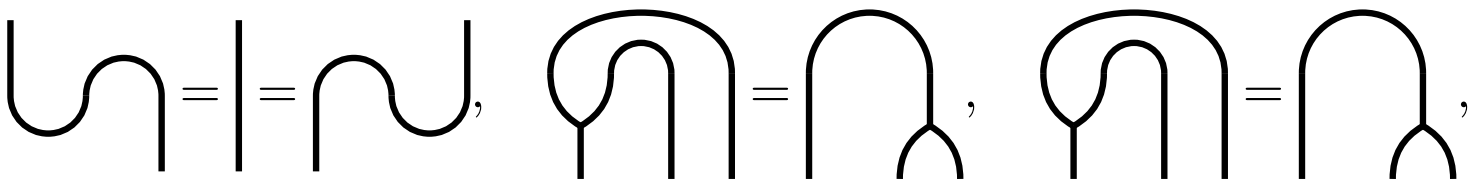

called ***isotopies*** (or ***non-degeneracy conditions***). ◇

*Remark* 3G.3. We will see several categories similar to **Web** and we will call them web categories and their morphisms webs. The picture, however, one should have in mind is more that of a soap foam, see Figure 8 and Exercise 3H.6. In all of these, using cups and caps, we can always use various rotations of the defining pictures, e.g. Equation 3G-1 could (and will be) used as an alternative generator.

Moreover, Theorem 3G.5 allows us to make even more simplifications, for example:

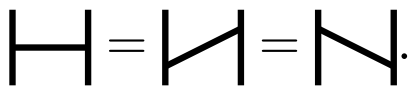

We will use this, often silently, throughout. ◇

Similarly to Example 1B.9 we also define:

**Example 3G.4.** The category **TGr** of embedded ***trivalent graphs***. Its objects are vertices $\bullet^n=\bullet\ldots\bullet$ for $n\in\mathbb{Z}_{\geq 0}$, and its morphisms are trivalent graphs (trivalent means every internal vertex has degree three) embedded in $\mathbb{R}^2$ between these, a.k.a. strands, illustrated as follows: ◇

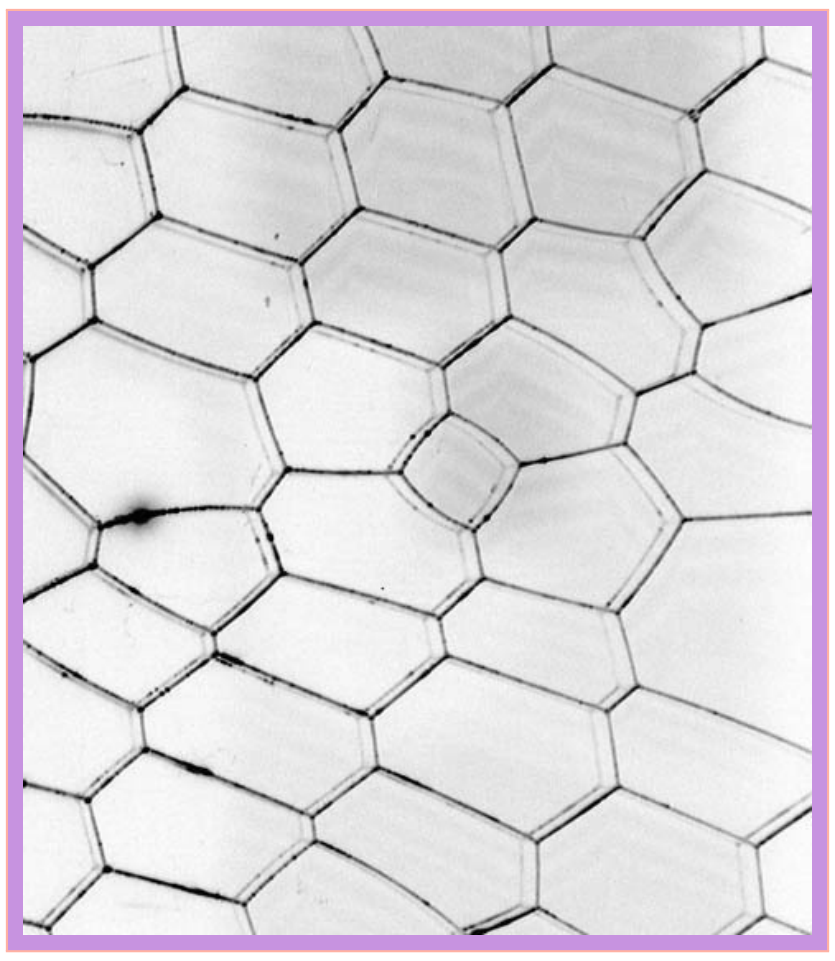

FIGURE 8. A cut through ("projection") a soap foam. Note the trivalent vertices, and the various faces (4,5,6,7 gons).

Picture from https://commons.wikimedia.org/wiki/File:2-dimensional_foam_(colors_inverted).jpg

**Theorem 3G.5.** *There exists a monoidal functor*

$$\mathrm{R}\colon \mathbf{Web} \to \mathbf{TGr},\ \bullet \mapsto \bullet,$$

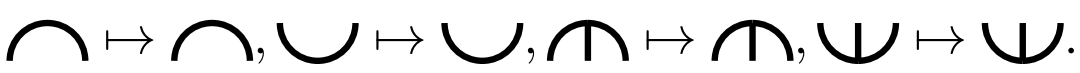

*The functor* R *is dense and fully faithful, thus,* $\mathbf{Web} \simeq_{\otimes} \mathbf{TGr}$.

*Proof (Sketch).* Very similar as the proof of Theorem 3D.10, there are several things to check, namely:

(a) Using Lemma 3C.10, that the functor R is a well-defined is easy.

(b) That R is dense is clear.

(c) That R is full is not hard: Every graph in **TGr** is embedded, so is locally of the form

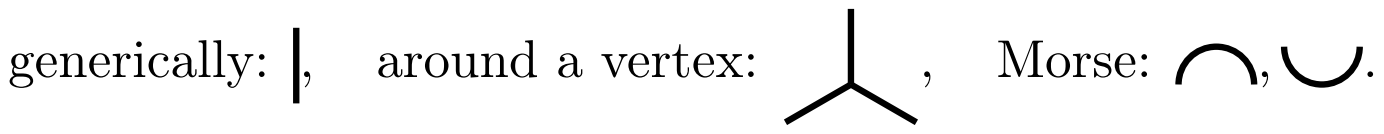

The same arguments as in the proof of Theorem 3D.10 apply.

(d) Proving faithfulness is painful. A crucial fact is the rotational invariance of the trivalent vertex, *i.e.*:

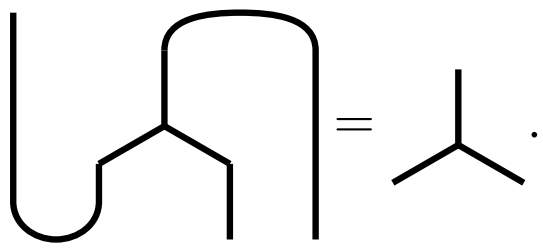

This follows by applying cups and caps to the isotopies in **Web** that contain a trivalent vertex. (This is Exercise 3H.7.)

The proof is complete. □

Here is another example in the same style.

**Example 3G.6.** We will now explain now to construct $\mathbf{Vec}_{\mathbb{C}}(G)$ diagrammatically. To explain it, we simplify notation below and an oriented strand labeled $g$ is the same as a strand labeled $g^{-1}$ with the opposite orientation. We use this convention to omit to display labels altogether: all labels can be recovered from the ones we display. Moreover, when we omit labels or orientations then we mean all suitable labels or orientations.

With these conventions in place, we let

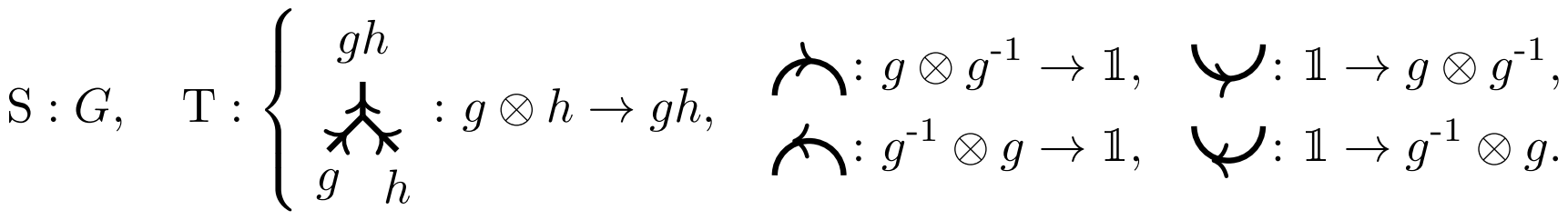

(3G-7)

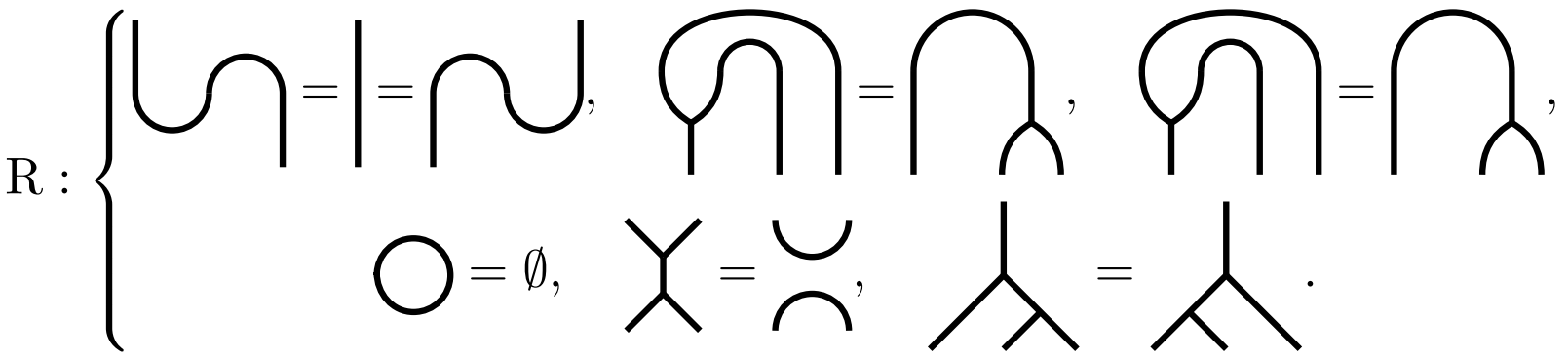

Here, similarly as in Equation 3G-1, we define various version of the trivalent vertex by rotating the defining picture. We denote the resulting category by $\mathbf{Vec}^{dia}_{\mathbb{C}}(G)$. Note that $\mathbf{Vec}^{dia}_{\mathbb{C}}(G)$ has infinitely many objects. ◇

**Theorem 3G.8.** *We have* $\mathbf{Vec}^{dia}_{\mathbb{C}}(G) \simeq_{\otimes} \mathbf{Vec}_{\mathbb{C}}(G)$.

*Proof.* Let us assume that $g = g^{-1}$ (thus, $G \cong (\mathbb{Z}/2\mathbb{Z})^k$) so that we can ignore orientations.

*Step 1, the number of isomorphism classes.* Note that, in $\mathbf{Vec}^{dia}_{\mathbb{C}}(G)$, we have an isomorphism $\mathrm{f} = \mathrm{f}(g,h)\colon g \otimes h \to gh$ and its inverse given by

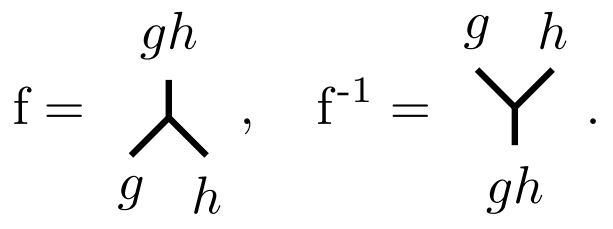

Indeed, the defining relations imply

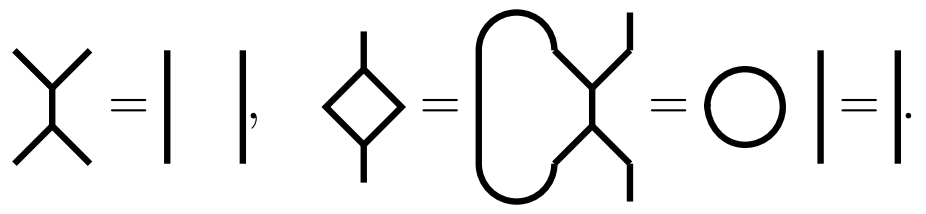

In particular, the category $\mathbf{Vec}^{dia}_{\mathbb{C}}(G)$ has exactly $\#G$ objects up to isomorphism.

*Step 2, the size of hom spaces.* Between any boundary sequence we can put at bottom and top, there is at most one diagram we can draw. To see this note that let us get rid of all internal edges of diagrams.

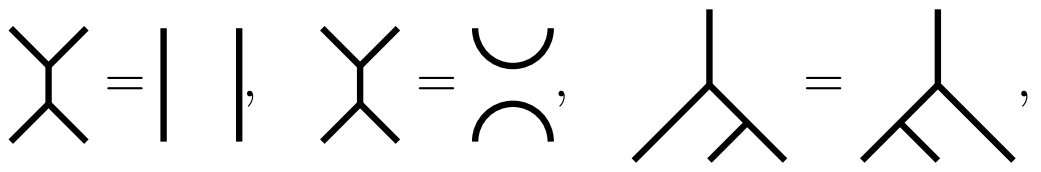

*Step 3, a monoidal functor.* Define $\mathrm{F}\colon \mathbf{Vec}_{\mathbb{C}}(G) \to \mathbf{Vec}^{dia}_{\mathbb{C}}(G)$ determined by sending $g \in \mathbf{Vec}_{\mathbb{C}}(G)$ to the object with the same name in $\mathbf{Vec}^{dia}_{\mathbb{C}}(G)$. Step 1 implies that $\mathrm{F}(gh) = \mathrm{F}(g \otimes h) \cong \mathrm{F}(g) \otimes \mathrm{F}(h)$. Step 1 also implies that F is dense, while Step 2 shows that F is full. Faithfulness is immediate.

The details are left as an exercise, see Exercise 3H.8. □

## 3H. Exercises.

*Exercise* 3H.1. Let $\mathrm{S}_n$ be the symmetric group of the set $\{1,...,n\}$ for $n > 2$. Show that

$$\mathrm{S}_n \cong \langle s_1,...,s_{n-1} \mid s_i^2 = 1, s_i s_{i\pm 1} s_i = s_{i\pm 1} s_i s_{i\pm 1}, s_i s_j = s_j s_i \text{ for } |i-j| > 1 \rangle$$

as groups. Deduce that $\mathrm{End}_{\mathbf{Sym}}(\bullet^n) \cong \mathrm{S}_n$. ◇

*Exercise* 3H.2. Recall the construction of Brauer category **Br** from Example 3D.7. Prove that the defining relations Equation 3D-8 imply that the following hold in **Br**:

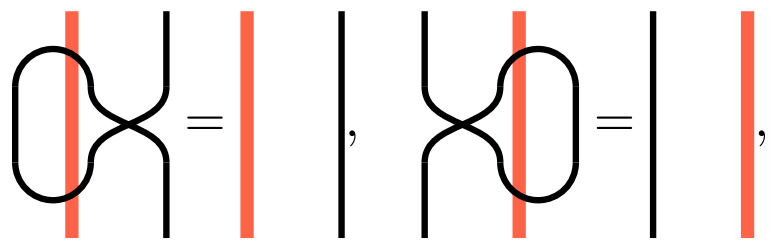

where the thick red strands represent an arbitrary number of straight strands. ◇

*Exercise* 3H.3. Prove Lemma 3E.11. ◇

*Exercise* 3H.4. Think about how to define right, left, bi(co)modules, their homomorphisms, subalgebras, ideals, submodules, *etc.* in the categorical setting, and choose your favorite notion and write down its categorical definition. ◇

*Exercise* 3H.5. Verify that $\bullet^2$ is a Frobenius algebra in **Br**, *cf.* Example 3E.7. ◇

*Exercise* 3H.6. Realize the picture in Figure 8 as a morphism in **Web**. ◇

*Exercise* 3H.7. In the proof sketch of Theorem 3G.5, verify that the rotational invariance of the trivalent vertex follows from the defining relations in **Web**. ◇

*Exercise* 3H.8. Finish the proof of Theorem 3G.8. ◇

# 4. Pivotal categories – definitions, examples and graphical calculus

Recall that Feynman diagrams for monoidal categories in general need to be upwards oriented, *i.e.* they do not have Morse points (a.k.a. cups and caps). So:

What kind of categories allow Morse points in their graphical calculus?

For us Morse points are one dimensional versions of "saddle points", *cf.* Figure 9.

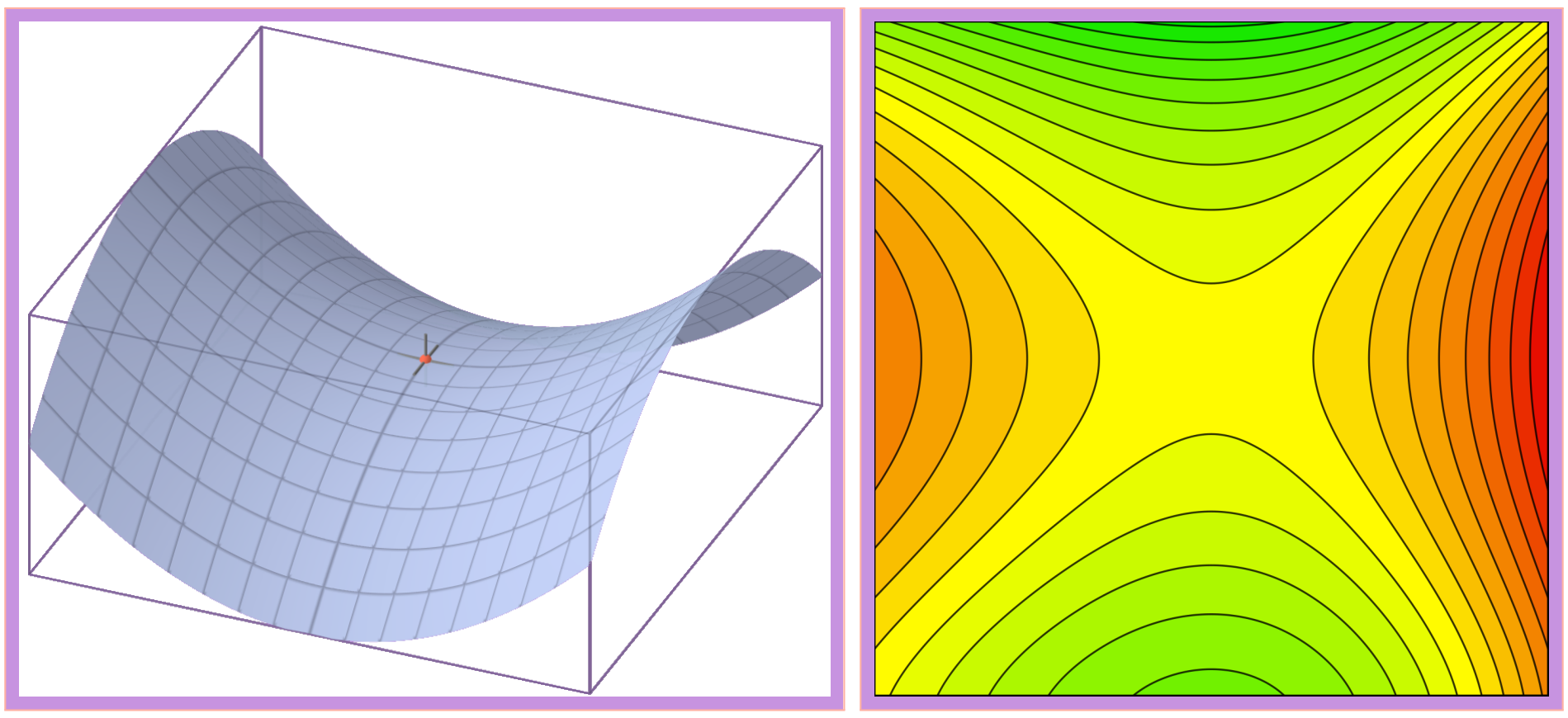

Figure 9. In a nutshell, Morse theory is the idea that most point on a manifold have nicely behaved neighborhoods and only a few critical points (called ***Morse points***) where the behavior changes drastically need to be studied. The pictures show a saddle point on a surface and its contour lines; this is a prototypical examples of a critical point where the behavior of a system changes.

Pictures from https://en.wikipedia.org/wiki/Morse_theory.

## 4A. A word about conventions.

This section is all about duals.

*Convention* 4A.1. We will use the symbol $\star$ for duality. Because it is always confusing, let us state right away that right duals will have their $\star$ on the right, and left duals on the left, *e.g.*

$$\text{object } \mathtt{X}, \quad \text{right dual } \mathtt{X}^{\star}, \quad \text{left dual } {}^{\star}\mathtt{X}.$$

If the left and the right dual agree, then we use the right dual $_^{\star}$ as the notation, and similar conventions for traces and dimensions. ◇

*Convention* 4A.2. Again, there will be several choices which we tend to omit when no confusion can arise. Moreover, whenever we write *e.g.* $\mathtt{X}^{\star}$ we implicitly assume existence of the right dual. ◇

*Convention* 4A.3. For pivotal categories we use the convention that strands labeled $\mathtt{X}$ are directed upwards, and those labeled with duals are oriented downward, see (4G-1). In particular, it suffices to label each strand once and in contrast to the general situation, *cf.* Convention 1A.3, we usually orient diagrams. ◇

*Convention* 4A.4. From now on we will use diagrammatics most of the time, and leave it to the reader to work out some of the non-strict versions of definitions and statements. For diagrams we use the terminology "taking mirrors" as before, but this also includes orientation reversals, *e.g.*

original: 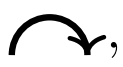, mirrors: 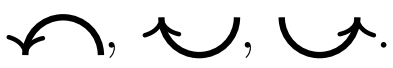.

(This is more general than using an actual mirror.) ◇

*Convention* 4A.5. If $\mathrm{End}_{\mathbf{C}}(\mathbb{1})$ is *e.g.* $\Bbbk$ , then we often identify the endomorphisms with actual elements, *e.g.* instead of "multiplication by $a \in \Bbbk$" we just write $a$. ◇

## 4B. Duality in monoidal categories.

Since duality is a powerful concept, we start with:

**Definition 4B.1.** A ***right dual*** $(\mathtt{X}^{\star}, \mathrm{ev}_{\mathtt{X}}, \mathrm{coev}_{\mathtt{X}})$ of $\mathtt{X} \in \mathbf{C}$ in a category $\mathbf{C} \in \mathbf{MCat}$ consists of

- an object $\mathtt{X}^{\star} \in \mathbf{C}$;
- a ***(right) evaluation*** $\mathrm{ev}_{\mathtt{X}}$ and a ***(right) coevaluation*** $\mathrm{coev}_{\mathtt{X}}$, *i.e.* morphisms

(4B-2) $\mathrm{ev}_{\mathtt{X}} \colon \mathtt{X}\mathtt{X}^{\star} \to \mathbb{1} \rightsquigarrow$ 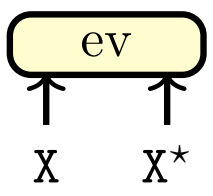

, $\mathrm{coev}_{\mathtt{X}} \colon : \mathbb{1} \to (\mathtt{X}^{\star})\mathtt{X} \rightsquigarrow$ 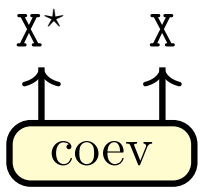

;

such that they are ***non-degenerate***, *i.e.*

(4B-3) 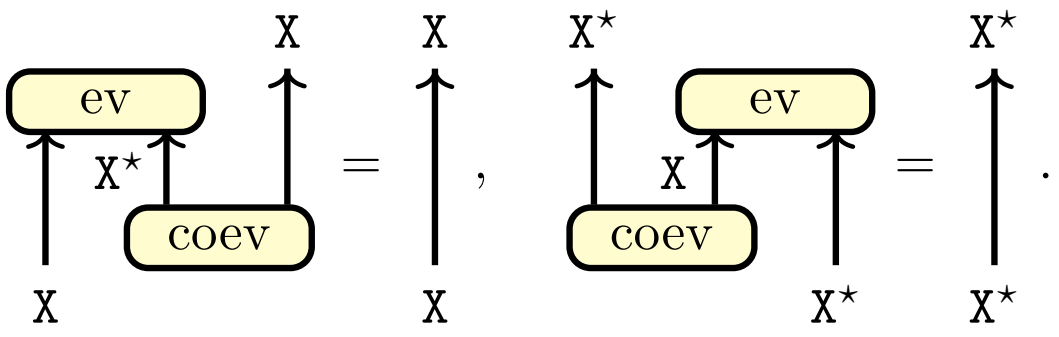

Similarly, a ***left dual*** $({}^{\star}\mathtt{X}, \mathrm{ev}^{\mathtt{X}}, \mathrm{coev}^{\mathtt{X}})$ of $\mathtt{X} \in \mathbf{C}$ in a category $\mathbf{C} \in \mathbf{MCat}$ consists of

- an object ${}^{\star}\mathtt{X} \in \mathbf{C}$;
- a ***(left) evaluation*** $\mathrm{ev}^{\mathtt{X}}$ and a ***(left) coevaluation*** $\mathrm{coev}^{\mathtt{X}}$, *i.e.* morphisms

(4B-4) $\mathrm{ev}^{\mathtt{X}} \colon {}^{\star}\mathtt{X}\mathtt{X} \to \mathbb{1} \rightsquigarrow$ 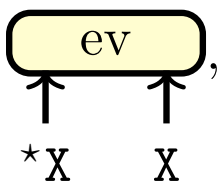, $\mathrm{coev}^{\mathtt{X}} \colon \colon \mathbb{1} \to \mathtt{X}({}^{\star}\mathtt{X}) \rightsquigarrow$ 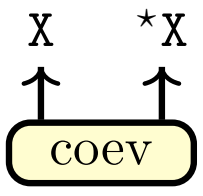;

such that they are ***non-degenerate***, *i.e.*

(4B-5) 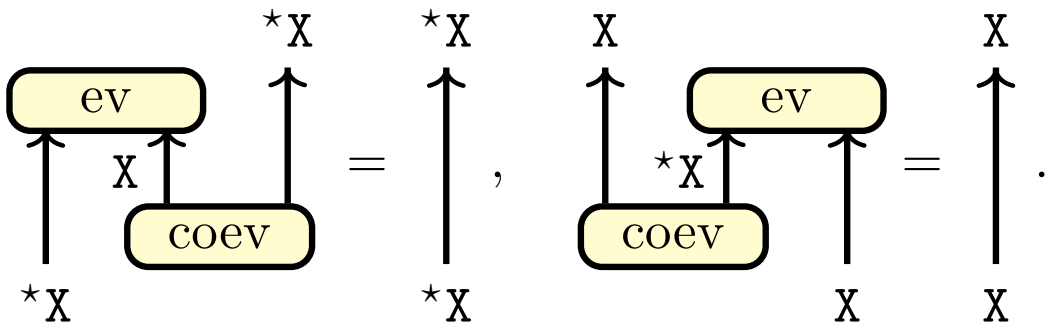

We call (4B-3) and (4B-5) the ***zigzag relations***. ◇

*Remark* 4B.6. Note that we do not distinguish the right and left (co)evaluation in coupons since the position of $\star$ will determine whether its right or left, *cf.* (4B-2) and (4B-4). ◇

The following justifies to say "the" right and left dual.

**Lemma 4B.7.** *Right and left duals, if they exist, are unique up to unique isomorphism.*

*Proof.* Let $\mathtt{X}^{\star}$ and $\overline{\mathtt{X}^{\star}}$ be two right duals of $\mathtt{X}$. These come with evaluation and coevaluation morphisms, $\mathrm{ev}_{\mathtt{X}}$ and $\mathrm{coev}_{\mathtt{X}}$, and $\overline{\mathrm{ev}}_{\mathtt{X}}$ and $\overline{\mathrm{coev}_{\mathtt{X}}}$, respectively. We use these to define two morphisms

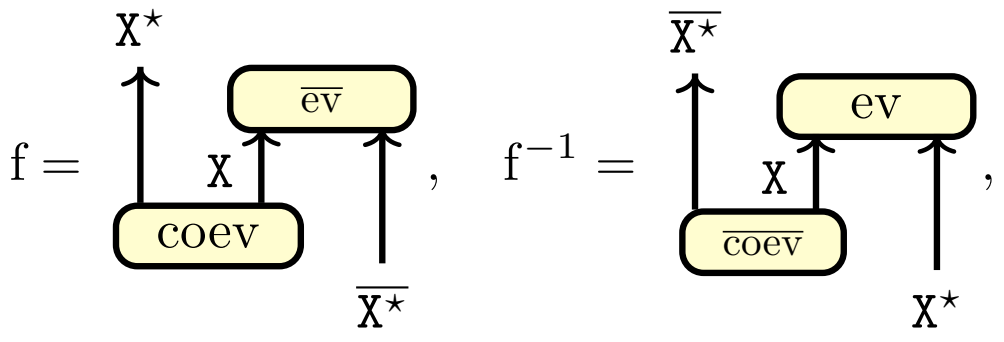

which are inverses by the zigzag relations, *e.g.*

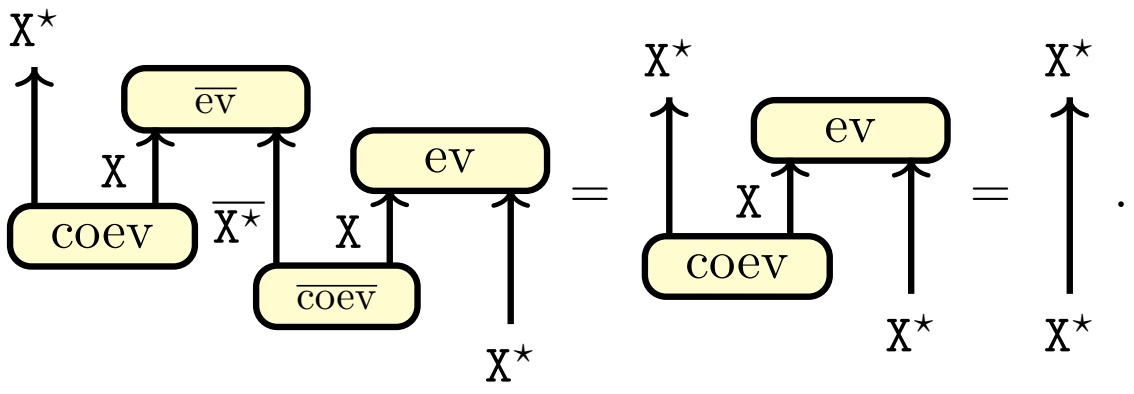

Moreover, it is easy to check that f is the only isomorphism which preserves the (co)evaluation. The proof for left duals is similar. □

**Lemma 4B.8.** *Fix* $\mathbf{C} \in \mathbf{MCat}$.

*(i) The monoidal unit is* **self-dual**, *meaning*

$$\mathbb{1} \cong \mathbb{1}^{\star} \cong {}^{\star}\mathbb{1}.$$

*(ii) For any* $\mathtt{X} \in \mathbf{C}$ *which has a right and a left dual we have*

$$(4B\text{-}9) \qquad {}^{\star}(\mathtt{X}^{\star}) \cong \mathtt{X} \cong ({}^{\star}\mathtt{X})^{\star}.$$

*(iii) If* $\mathtt{X} \in \mathbf{C}$ *has a right dual, then* $\mathtt{X} \in \mathbf{C}^{co}$ *has a left dual, and vice versa.*

*Proof.* (i). This follows since we can take the unitors as (co)evaluation morphisms.

(ii). The isomorphisms are similar to the ones in the proof of Lemma 4B.7, where we again use the zigzag relations to show that they invert one another.

(iii). Clear by comparing (4B-2) and (4B-4). □

Lemma 4B.8.(iii) is the first instance of what we call ***right-left symmetry***. It says in words that "Every statement about right duals has a left counterpart and *vice versa*.", see also Section 1C.

4C. **Some first examples of duals.** The following can be taken as an example or as our definition of adjoint functors:

**Example 4C.1.** In the monoidal category $\mathbf{End}(\mathbf{C})$, *cf.* Example 2H.5, the right dual $\mathrm{F}^{\star}$ of a functor is called its ***right adjoint***, while the left dual ${}^{\star}\mathrm{F}$ is called its ***left adjoint***. ◇

**Lemma 4C.2.** *Example 4C.1 is not lying (meaning duals in* $\mathbf{End}(\mathbf{C})$ *are adjoints).*

*Proof.* By comparing Definition 1I.7 and Definition 4B.1. □

Duals in general might not exist, *e.g.* not every functor has adjoints. A more down to earth example is:

**Example 4C.3.** Not every objects in $\mathbf{Vec}_{\mathbb{k}}$ has duals. However, if $\mathtt{X}$ is finite dimensional, then $\mathtt{X}^{\star}={}^{\star}\mathtt{X}$ is the vector space dual with the (co)evaluations being the usual maps, *e.g.*

$$\mathrm{ev}_{\mathtt{X}}\colon \mathtt{X}\mathtt{X}^{\star}\to\mathbb{k},\ (x,y^{\star})\mapsto y^{\star}(x),\quad \mathrm{coev}_{\mathtt{X}}\colon \mathbb{k}\to(\mathtt{X}^{\star})\mathtt{X},\ 1\mapsto\textstyle\sum_{i=1}^{n}x_i^{\star}\otimes x_i,$$

where $\{x_1,...,x_n\}$ and $\{x_1^{\star},...,x_n^{\star}\}$ are choices of dual bases of $\mathtt{X}$ and $\mathtt{X}^{\star}$. ◇

**Example 4C.4.** Most of the diagrammatic categories which we have seen have duals. For example, in $\mathbf{TL}$ or $\mathbf{Br}$ the generating object $\bullet$ is self-dual. More precisely,

$$\bullet=\bullet^{\star}={}^{\star}\bullet,\quad \mathrm{ev}_{\bullet}=\mathrm{ev}^{\bullet}=$$ 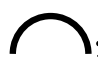, $\quad\mathrm{coev}_{\bullet}=\mathrm{coev}^{\bullet}=$ 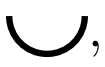,

where the evaluation and the coevaluation are cap and cup morphisms, as illustrated. In the algebraic model of $\mathbf{1State}$, the oriented quantum Brauer category $\mathbf{oqBr}$, we have $\bullet^{\star}={}^{\star}\bullet$ with (3D-15) being the four (co)evaluations. In other words, cups represent evaluations and caps coevaluations. ◇

Note that duals, if they exist, are unique, but the evaluation and coevaluation are not unique. In particular, they usually can be scaled if we are in a $\mathbb{k}$ linear setting. The crucial example where scaling will matter later on is:

**Example 4C.5.** Let us consider $\mathbf{Vec}_{\mathbb{C}}^{\omega}(\mathbb{Z}/2\mathbb{Z})$ for the non-trivial 3 cocycle $\omega$. In this category all objects are self-dual, *i.e.* $\mathtt{1}=\mathtt{1}^{\star}={}^{\star}\mathtt{1}$ and

$$\mathtt{1}\mathtt{1}^{\star}=\mathtt{0},\quad(\mathtt{1}^{\star})\mathtt{1}=\mathtt{0}.$$

But the object $\mathtt{1}$ admits several (co)evaluations, which we explain for the right duality, the left being similar by right-left symmetry. The (hidden) associativity constrains in (4B-3) are

$$\mathtt{1}\xrightarrow{\mathrm{id}_{1}\otimes\mathrm{coev}_{1}}\mathtt{1}(\mathtt{1}\mathtt{1})\xrightarrow{\alpha_{1,1,1}}(\mathtt{1}\mathtt{1})\mathtt{1}\xrightarrow{\mathrm{ev}_{1}\otimes\mathrm{id}_{1}}\mathtt{1}\ ,\qquad \mathtt{1}\xrightarrow{\mathrm{coev}^{1}\otimes\mathrm{id}_{1}}(\mathtt{1}\mathtt{1})\mathtt{1}\xrightarrow{\alpha_{1,1,1}^{-1}}\mathtt{1}(\mathtt{1}\mathtt{1})\xrightarrow{\mathrm{id}_{1}\otimes\mathrm{ev}^{1}}\mathtt{1}\ .$$

Thus, recalling that $\alpha_{1,1,1}$ gives a sign, whatever non-zero scalar $a\in\mathbb{C}^{*}$ we like to scale $\mathrm{ev}_{1}$ with, we then need to scale $\mathrm{coev}_{1}$ by $-a^{-1}$. The minus sign is the crucial part here: one can also scale the (co)evaluations for $\mathtt{0}$, but then only with $a$ and $a^{-1}$. ◇

Duality is actually a functor, as we will see next.

**Definition 4C.6.** For $(\mathrm{f}\colon\mathtt{X}\to\mathtt{Y})\in\mathbf{C}$, in a category $\mathbf{C}\in\mathbf{MCat}$, its ***right*** $\mathrm{f}^{\star}\colon\mathtt{Y}^{\star}\to\mathtt{X}^{\star}$ and ***left mate*** ${}^{\star}\mathrm{f}\colon{}^{\star}\mathtt{Y}\to{}^{\star}\mathtt{X}$ are defined as

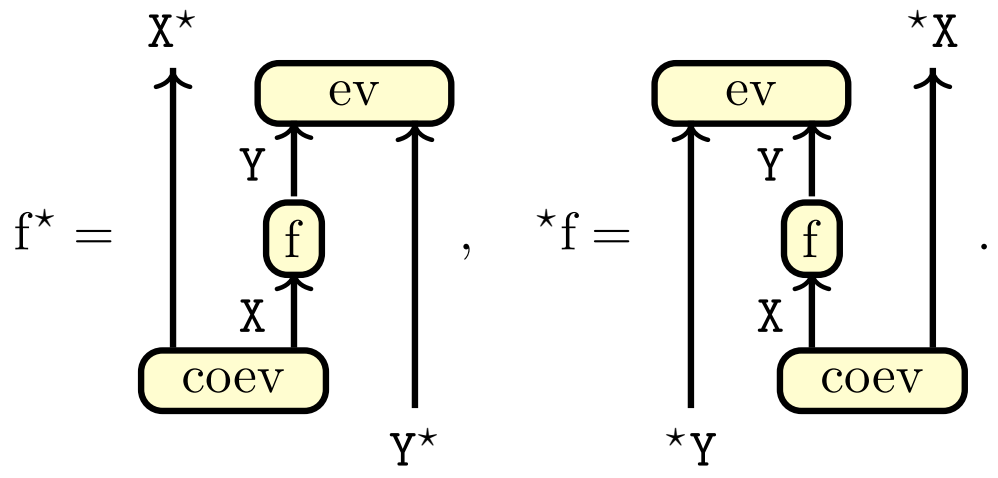

(Instead of mate the terminology ***dual morphism*** is used in other literature.) ◇

**Lemma 4C.7.** *Fix any* $\mathbf{C}\in\mathbf{MCat}$*. Then, for all* $\mathtt{X},\mathtt{Y},\mathrm{f},\mathrm{g}\in\mathbf{C}$*:*

(i) *We have* $(\mathrm{gf})^{\star}=(\mathrm{f}^{\star})(\mathrm{g}^{\star})$ *and* ${}^{\star}(\mathrm{gf})=({}^{\star}\mathrm{f})({}^{\star}\mathrm{g})$*.*

(ii) *We have* $(\mathtt{X}\mathtt{Y})^{\star}\cong(\mathtt{Y}^{\star})(\mathtt{X}^{\star})$ *and* ${}^{\star}(\mathtt{X}\mathtt{Y})\cong({}^{\star}\mathtt{Y})({}^{\star}\mathtt{X})$*.*

*Proof.* This is Exercise 4J.1. □

The most useful consequence of having duals in practice is:

**Theorem 4C.8.** *Let* $\mathtt{X},\mathtt{Y},\mathtt{Z}\in\mathbf{C}$ *be objects in any* $\mathbf{C}\in\mathbf{MCat}$*. Then we have*

$$\mathrm{Hom}_{\mathbf{C}}(\mathtt{XY},\mathtt{Z})\cong\mathrm{Hom}_{\mathbf{C}}\big(\mathtt{Y},(\mathtt{X}^{\star})\mathtt{Z}\big),\quad \mathrm{Hom}_{\mathbf{C}}(\mathtt{YX},\mathtt{Z})\cong\mathrm{Hom}_{\mathbf{C}}\big(\mathtt{Y},\mathtt{Z}({}^{\star}\mathtt{X})\big),$$
$$\mathrm{Hom}_{\mathbf{C}}(\mathtt{Y},\mathtt{ZX})\cong\mathrm{Hom}_{\mathbf{C}}(\mathtt{YX}^{\star},\mathtt{Z}),\quad \mathrm{Hom}_{\mathbf{C}}(\mathtt{Y},\mathtt{XZ})\cong\mathrm{Hom}_{\mathbf{C}}({}^{\star}\mathtt{XY},\mathtt{Z}).$$

*(Of course assuming that the corresponding duals exist for* $\mathtt{X}$*.)*

*Proof.* Let us construct isomorphisms for the first case, all others cases are similar. We define

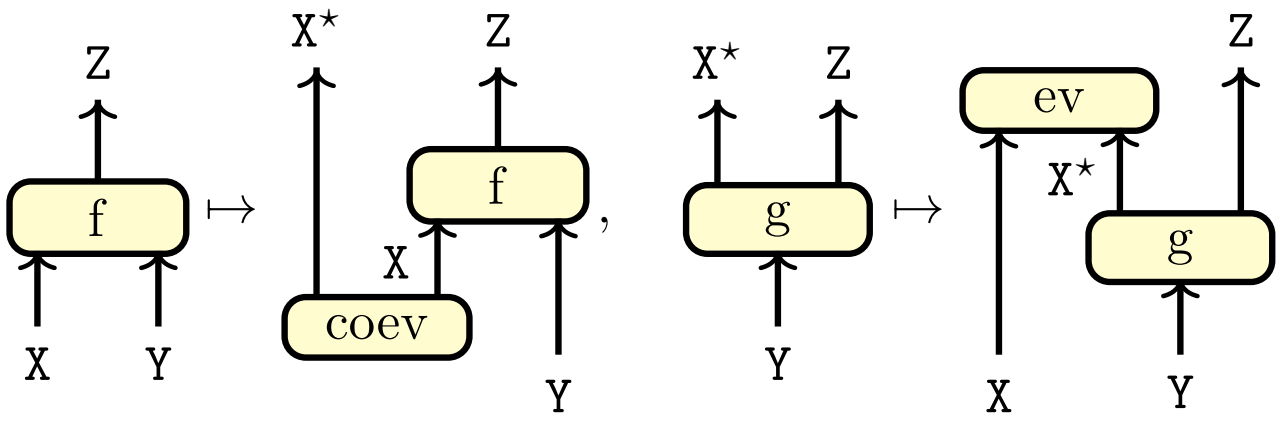

That these are inverses follows from the zigzag relations. □

**Proposition 4C.9.** *Let* $\mathrm{F}\in\mathbf{Hom}_{\otimes}(\mathbf{C},\mathbf{D})$*. If* $\mathtt{X}^{\star}$ *is a right dual of* $\mathtt{X}\in\mathbf{C}$*, then* $\mathrm{F}(\mathtt{X}^{\star})$ *is a right dual of* $\mathrm{F}(\mathtt{X})\in\mathbf{D}$*. Similarly for left duals.*

*Proof.* By right-left symmetry, it suffices to define

$$\mathrm{ev}_{\mathrm{F}(\mathtt{X})}\colon\ \mathrm{F}(\mathtt{X})\mathrm{F}(\mathtt{X}^{\star})\xrightarrow{\xi_{\mathtt{X},\mathtt{X}^{\star}}}\mathrm{F}(\mathtt{XX}^{\star})\xrightarrow{\mathrm{F}(\mathrm{ev}_{\mathtt{X}})}\mathrm{F}(\mathbb{1})\xrightarrow{\xi_{\mathbb{1}}^{-1}}\mathbb{1}\ ,$$

$$\mathrm{coev}_{\mathrm{F}(\mathtt{X})}\colon\ \mathbb{1}\xrightarrow{\xi_{\mathbb{1}}}\mathrm{F}(\mathbb{1})\xrightarrow{\mathrm{F}(\mathrm{coev}_{\mathtt{X}})}\mathrm{F}\big((\mathtt{X}^{\star})\mathtt{X}\big)\xrightarrow{\xi_{\mathtt{X}^{\star},\mathtt{X}}^{-1}}\mathrm{F}(\mathtt{X}^{\star})\mathrm{F}(\mathtt{X})\ .$$

These are the corresponding (co)evaluations, as a straightforward calculation verifies. □

## 4D. **Rigidity.** Recall that duals might not exists. This motivates:

**Definition 4D.1.** A category $\mathbf{C}\in\mathbf{MCat}$ is called ***rigid*** if every object has right and left duals. ◇

**Example 4D.2.** Several examples which we have seen are rigid.

(a) $\mathbf{fdVec}_{\mathbb{k}}$ is rigid, *cf.* Example 4C.3.

(b) $\mathbf{Vec}_{\mathbb{k}}^{\omega}(\mathrm{G})$ (for the duration, we will always use the $\mathbb{k}$ linear incarnation of $\mathbf{Vec}^{\omega}(\mathrm{G})$) is rigid with $\mathtt{g}^{-1}=\mathtt{g}^{\star}={}^{\star}\mathtt{g}$.

(c) The diagrammatic categories **TL** and **Br** are rigid with a self-dual generator $\bullet$.

(d) The diagrammatic categories **oTL** and **oBr** are rigid with $\bullet^{\star}={}^{\star}\bullet$.

In contrast, $\mathbf{Vec}_{\mathbb{k}}$ is not rigid. ◇

Let $\mathtt{X}^{\star\star}=(\mathtt{X}^{\star})^{\star}$ and ${}^{\star\star}\mathtt{X}={}^{\star}({}^{\star}\mathtt{X})$ denote the ***double duals***. Note that all the examples in Example 4D.2 satisfy

$$\mathtt{X}^{\star}\cong{}^{\star}\mathtt{X}\xRightarrow{\text{(4B-9)}}\mathtt{X}^{\star\star}\cong\mathtt{X}\cong{}^{\star\star}\mathtt{X}\ . \tag{4D-3}$$

This is not always true:

**Example 4D.4.** The free rigid category generated by one object $\bullet$ has

$$\bullet^{\star}\ncong{}^{\star}\bullet\Longrightarrow\bullet^{\star\star}\ncong\bullet\ncong{}^{\star\star}\bullet\ . \tag{4D-5}$$

The proof of (4D-5) this requires non-trivial arguments, *i.e.* by constructing models: examples where (4D-3) fails exist, but are not easy to construct. ◇

**Example 4D.6.** By Proposition 4C.9, we see that monoidal functors are already the correct morphisms between rigid categories, as long as we do not care for the choices of (co)evaluations. Thus, we get the ***category of rigid categories*** **RCat**. ◇

**Lemma 4D.7.** *If* $\mathbf{C}\in\mathbf{RCat}$*, then* $\mathbf{C}^{op},\mathbf{C}^{co},\mathbf{C}^{coop}\in\mathbf{RCat}$

*Proof.* Immediate by taking mirrors of diagrams. □

Lemma 4C.7 shows that we can define important functors between rigid categories:

**Definition 4D.8.** For $\mathbf{C}\in\mathbf{RCat}$ we define ***right*** and ***left duality functors***

$$_^{\star}\colon\mathbf{C}\to\mathbf{C}^{coop},\ \mathtt{X}\mapsto\mathtt{X}^{\star},\mathrm{f}\mapsto\mathrm{f}^{\star},\quad {}^{\star}_\colon\mathbf{C}^{coop}\to\mathbf{C},\ \mathtt{X}\mapsto{}^{\star}\mathtt{X},\mathrm{f}\mapsto{}^{\star}\mathrm{f}.$$

(Note the coop.) ◇

**Proposition 4D.9.** *For* $\mathbf{C} \in \mathbf{RCat}$ *we have equivalences*

$$_^{\star} \colon \mathbf{C} \xrightarrow{\simeq_{\otimes}} \mathbf{C}^{coop}, \quad {}^{\star}_ \colon \mathbf{C}^{coop} \xrightarrow{\simeq_{\otimes}} \mathbf{C}, \tag{4D-10}$$

$$^{\star}(_^{\star}) \cong \mathrm{Id}_{\mathbf{C}}, \quad ({}^{\star}_)^{\star} \cong \mathrm{Id}_{\mathbf{C}^{coop}}. \tag{4D-11}$$

*Proof.* As usual it suffices to prove (4D-10) for the right duality. By Lemma 4C.7 we see that $_^{\star}$ is a well-defined monoidal functor, while Theorem 4C.8 shows fully faithfulness of $_^{\star}$. Moreover, (4B-9) proves that $_^{\star}$ is dense, thus, an equivalence. The second part (4D-11) follows easily from (4D-3). □

*Remark* 4D.12. The only reason to define the right duality be a functor from $\mathbf{C}$ to $\mathbf{C}^{coop}$ and the left duality the other way around is to get a cleaner statement in (4D-11), but for the duration we rather have the left duality also defined to be from $\mathbf{C}$ to $\mathbf{C}^{coop}$. Furthermore, alternatively right and left dualities also give equivalences (either way) $\mathbf{C}^{op} \simeq_{\otimes} \mathbf{C}^{co}$. ◇

Immediate consequences of Proposition 4D.9 are:

**Proposition 4D.13.** *For* $\mathbf{C} \in \mathbf{RCat}$ *we have equivalences*

$$_^{\star\star} \colon \mathbf{C} \xrightarrow{\simeq_{\otimes}} \mathbf{C}, \quad {}^{\star\star}_ \colon \mathbf{C} \xrightarrow{\simeq_{\otimes}} \mathbf{C}. \tag{4D-14}$$

*Both equivalences can also be stated between* $\mathbf{C}^{op}$ *and* $\mathbf{C}^{op}$*,* $\mathbf{C}^{co}$ *and* $\mathbf{C}^{co}$*, or* $\mathbf{C}^{coop}$ *and* $\mathbf{C}^{coop}$*.* □

**Proposition 4D.15.** *For any rigid category* $\mathbf{C}$ *its Grothendieck classes* $\mathrm{K}_0(\mathbf{C})$ *form a monoid, with multiplication and unit as in Proposition 2D.9, and two homomorphisms*

$$[_^{\star}] \colon \mathrm{K}_0(\mathbf{C}) \to \mathrm{K}_0(\mathbf{C}^{coop}),\ [\mathtt{X}] \mapsto [\mathtt{X}^{\star}], \quad [{}^{\star}_] \colon \mathrm{K}_0(\mathbf{C}) \to \mathrm{K}_0(\mathbf{C}^{coop}),\ [\mathtt{X}] \mapsto [{}^{\star}\mathtt{X}].$$

*Moreover, they are inverse of one another.*

*Proof.* By Lemma 4C.7 we have

$$[(\mathtt{XY})^{\star}] = [\mathtt{Y}^{\star}\mathtt{X}^{\star}] = [\mathtt{Y}^{\star}][\mathtt{X}^{\star}]$$

and we get the left analog by right-left symmetry. They are inverses by (4B-9). □

**Example 4D.16.** On $\mathrm{K}_0(\mathbf{fdVec}_{\Bbbk})$, *cf.* Example 1H.17, the two homomorphisms $[_^{\star}]$ and $[{}^{\star}_]$ agree and are the identities. In contrast, the functor itself are only isomorphic to identity functors. ◇

Let us now take care of the choice of (co)evaluations:

**Definition 4D.17.** A functor $\mathrm{F} \in \mathbf{Hom}_{\otimes}(\mathbf{C},\mathbf{D})$ for $\mathbf{C},\mathbf{D} \in \mathbf{RCat}$ is called ***rigid*** if

$$\mathrm{F}(\mathrm{ev}_{\mathtt{X}}) = \mathrm{ev}_{\mathrm{F}(\mathtt{X})}, \quad \mathrm{F}(\mathrm{coev}_{\mathtt{X}}) = \mathrm{coev}_{\mathrm{F}(\mathtt{X})}, \quad \mathrm{F}(\mathrm{ev}^{\mathtt{X}}) = \mathrm{ev}^{\mathrm{F}(\mathtt{X})}, \quad \mathrm{F}(\mathrm{coev}^{\mathtt{X}}) = \mathrm{coev}^{\mathrm{F}(\mathtt{X})},$$

holds for all $\mathtt{X} \in \mathbf{C}$. ◇

The following lemma is immediate.

**Lemma 4D.18.** *The identity functor on a rigid category is rigid. Moreover, if* F *and* G *are rigid functors, then so is* GF*.* □

**Example 4D.19.** We get a (dense, in the monoidal sense, but non-full) subcategory $\mathbf{R^{+}Cat} \subset \mathbf{RCat}$, the ***category of rigid categories and rigid functors***. Also, we have a (non-dense, but full) subcategory $\mathbf{Hom}_{\star}(\mathbf{C},\mathbf{D}) \subset \mathbf{Hom}_{\otimes}(\mathbf{C},\mathbf{D})$, the ***category rigid functors***. ◇

**Definition 4D.20.** $\mathbf{C},\mathbf{D} \in \mathbf{RCat}$ are called ***equivalent as rigid categories***, denoted by $\mathbf{C} \simeq_{\star} \mathbf{D}$, if there exists an equivalence $\mathrm{F} \in \mathbf{Hom}_{\star}(\mathbf{C},\mathbf{D})$. ◇

**Example 4D.21.** Recall that $\mathbf{Vec}_{\mathbb{C}}^{\omega}(\mathbb{Z}/2\mathbb{Z})$ allowed several choices of (co)evaluations, some of which differ by signs. A monoidal functor does not take these choices into account, so they are all monoidally equivalent. However, Lemma 4H.10 below will show that not all of these choice give $\simeq_{\star}$ equivalent rigid categories. ◇

4E. **Categorical groups.** In some sense, see also Example 4D.2.(b) or Exercise 4J.3, rigid categories are like categorical versions of groups. Let us make this a bit more precise.

**Definition 4E.1.** Let $\mathbf{C} \in \mathbf{RCat}$. Then $\mathtt{X} \in \mathbf{C}$ is called ***invertible*** if $\mathrm{ev}_{\mathtt{X}} \colon \mathtt{XX}^{\star} \to \mathbb{1}$ and $\mathrm{coev}_{\mathtt{X}} \colon \mathbb{1} \to (\mathtt{X}^{\star})\mathtt{X}$ are isomorphisms. ◇

That Definition 4E.1 seems to favor right over left is a mirage:

**Lemma 4E.2.** *If* $\mathtt{X},\mathtt{Y} \in \mathbf{C}$ *are invertible, then:*

*(i)* *We have* $\mathtt{X}^{\star} \cong {}^{\star}\mathtt{X}$*.*

*(ii)* *The object* $\mathtt{X}^{\star}$ *is invertible.*

*(iii)* *The object* $\mathtt{XY}$ *is invertible.*

*Proof.* (i). Note that we have $\mathtt{X}\mathtt{X}^{\star} \cong \mathbb{1} \cong (\mathtt{X}^{\star})\mathtt{X}$ by invertibility of $\mathtt{X}$. Thus, taking duals we also have ${}^{\star}\mathtt{X}\mathtt{X} \cong \mathbb{1} \cong \mathtt{X}({}^{\star}\mathtt{X})$, which we can put together to get $\mathtt{X}^{\star} \cong (\mathtt{X}^{\star})\mathtt{X}({}^{\star}\mathtt{X}) \cong {}^{\star}\mathtt{X}$.

(ii). Clear by (4B-9).

(iii). This follows since $\mathrm{ev}_{\mathtt{X}\mathtt{Y}}$, respectively $\mathrm{coev}_{\mathtt{X}\mathtt{Y}}$, can be defined as compositions of $\mathrm{ev}_{\mathtt{X}}$ with $\mathrm{ev}_{\mathtt{Y}}$, respectively of $\mathrm{coev}_{\mathtt{X}}$ with $\mathrm{coev}_{\mathtt{Y}}$. □

**Example 4E.3.** Lemma 4E.2 says that we get a monoidal category $\mathbf{Inv}(\mathbf{C})$ as well as a group $\mathrm{Inv}(\mathbf{C}) = \mathrm{Inv}\big(\mathrm{K}_0(\mathbf{C})\big)$ of ***invertible objects***. ◇

**Definition 4E.4.** Let $\mathbf{C} \in \mathbf{RCat}$. Then $\mathbf{C}$ is called a ***categorical group*** if $\mathbf{Inv}(\mathbf{C}) = \mathbf{C}$. ◇

**Example 4E.5.** With respect to the examples in Example 4D.2:

(i) $\mathbf{Inv}(\mathbf{fdVec}_{\mathbb{k}})$ has, up to isomorphisms, only the object $\mathbb{k}$. Hence, $\mathrm{Inv}(\mathbf{fdVec}_{\mathbb{k}}) \cong 1$, which is the submonoid of invertible elements in $\mathrm{K}_0(\mathbf{fdVec}_{\mathbb{k}}) \cong \mathbb{Z}_{\geq 0}$.

(ii) For $\mathbf{Vec}_{\mathbb{k}}^{\omega}(\mathrm{G})$ one clearly has $\mathrm{Inv}\big(\mathbf{Vec}_{\mathbb{k}}^{\omega}(\mathrm{G})\big) \cong \mathrm{K}_0\big(\mathbf{Vec}_{\mathbb{k}}^{\omega}(\mathrm{G})\big) \cong \mathrm{G}$, and $\mathbf{Vec}_{\mathbb{k}}^{\omega}(\mathrm{G})$ is a categorical group.

(iii) For the diagrammatic categories à la Brauer one always has $\mathrm{Inv}(\mathbf{Br}) \cong 1$.

Note that the second example shows that any group can appear as $\mathrm{Inv}(\mathbf{C})$. ◇

4F. **Pivotality.** Note that Example 4D.4 shows that the equivalences from (4D-14) might not be trivial. In fact, they can be of infinite order. This motivates the following definition.

**Definition 4F.1.** A category $\mathbf{C} \in \mathbf{RCat}$ is called ***pivotal*** if $_^{\star} \cong_{\otimes} {}^{\star}_$. A ***pivotal structure*** on a pivotal category is a choice of an isomorphism $\pi\colon _^{\star} \xRightarrow{\cong_{\otimes}} {}^{\star}_$. ◇

In other words, in a pivotal category we have (4D-3). Thus:

**Proposition 4F.2.** *For any pivotal category* $\mathbf{C}$ *we have* $\mathrm{Id}_{\mathbf{C}} \cong_{\otimes} _^{\star\star} \cong_{\otimes} {}^{\star\star}_$, *and hence the functor* $_^{\star} \cong_{\otimes} {}^{\star}_$ *is of order two.* □

On the other hand, a pivotal structure on a pivotal category is a further choice of isomorphisms

$$\pi_{\mathtt{X}}\colon \mathtt{X}^{\star} \xrightarrow{\cong} {}^{\star}\mathtt{X}, \tag{4F-3}$$

natural in $\mathtt{X}$, satisfying $\pi_{\mathtt{X}\mathtt{Y}} = \pi_{\mathtt{X}} \otimes \pi_{\mathtt{Y}}$. Alternatively, a pivotal structure on a pivotal category is a further choice of isomorphisms

$$\pi_{\mathtt{X}}\colon \mathtt{X} \xrightarrow{\cong} \mathtt{X}^{\star\star}, \tag{4F-4}$$

satisfying exactly the same conditions.

*Remark* 4F.5. It is more natural to define a pivotal structure as isomorphisms identifying right and left duals, *i.e.* using (4F-3). However, in practice the choice of isomorphisms as in (4F-4) turns out to be more useful, and we will use both interchangeable. ◇

**Example 4F.6.** All examples in Example 4D.2 are pivotal. More precisely:

(a) $\mathbf{fdVec}_{\mathbb{k}}$ has a pivotal structure coming from the classical $\mathbb{k}$ vector space duality $\mathrm{V} \cong \mathrm{V}^{\star\star}$.

(b) For $\mathbf{Vec}_{\mathbb{k}}(\mathrm{G})$ one can choose the pivotal structure to be the identity.

(c) The diagrammatic categories à la Brauer usually have $\bullet^{\star} = {}^{\star}\bullet$ or even $\bullet = \bullet^{\star} = {}^{\star}\bullet$, which gives them an evident pivotal structure.

(For the first example recall the beginning of Section 1F.) ◇

**Example 4F.7.** Note the difference between being free:

(a) The free rigid category generated by one object $\bullet$, *cf.* Example 4D.4, is not pivotal.

(b) The free pivotal category generated by one object $\bullet$ is $\mathbf{oTL}$.

(c) The free pivotal category generated by one self-dual object $\bullet$ is $\mathbf{TL}$.

Essentially, self-dual means no orientation. ◇

In these notes we tend to omit to choose a pivotal structure. To be precise, we take the one in (4G-7) which only involves choices of (co)evaluations, so:

**Example 4F.8.** We also have the category $\mathbf{PCat} \subset \mathbf{R^{+}Cat}$, the ***category of pivotal categories***, whose morphisms are rigid functors. ◇

**Lemma 4F.9.** *In a pivotal category right and left mates are conjugate, i.e.* $\pi_{\mathtt{X}}\mathrm{f}^{\star} = {}^{\star}\mathrm{f}\pi_{\mathtt{Y}}$, *where* $\pi\colon _^{\star} \xRightarrow{\cong_{\otimes}} {}^{\star}_$ *is a choice of pivotal structure.*

*Proof.* The claim follows directly from $_^{\star} \cong_{\otimes} {}^{\star}_$ and its commuting diagram. □

**Definition 4F.10.** A category $\mathbf{C} \in \mathbf{PCat}$ is called ***strict***, if $_^{\star} = {}^{\star}_$ as functors. ◇

Thus, we can write $f^{\star}$ for the mate in case $\mathbf{C} \in \mathbf{PCat}$ is strict.

Similarly as in Theorem 2I.5 we have the pivotal strictification, which we will use in all diagrammatics:

**Theorem 4F.11.** *For any pivotal category* $\mathbf{C}$ *there exists a strict pivotal category* $\mathbf{C}^{st}$ *which is pivotal equivalent to* $\mathbf{C}$, *i.e.* $\mathbf{C} \simeq_{\star} \mathbf{C}^{st}$.

*Proof.* It is not hard, but also not trivial, to generalize the arguments in Theorem 2I.5 to pivotal categories by constructing an appropriate functor category, see *e.g.* [**NS07**, Theorem 2.2] Alternatively, this can be deduced from a version of the monoidal coherence theorem for pivotal categories similarly as the proof of Theorem 2I.5 can be deduced from Theorem 2G.2. (Such a pivotal coherence theorem is stated in [**BW99**, Theorem 1.9].) Details are omitted for brevity. □

4G. **Feynman diagrams for pivotal categories.** The diagrams we can draw for strict pivotal categories are now topological in nature, as well will see. The diagrammatic conventions are the ones for monoidal categories, see *e.g.* (2E-1), together with diagrammatic rules for duals:

(4G-1)

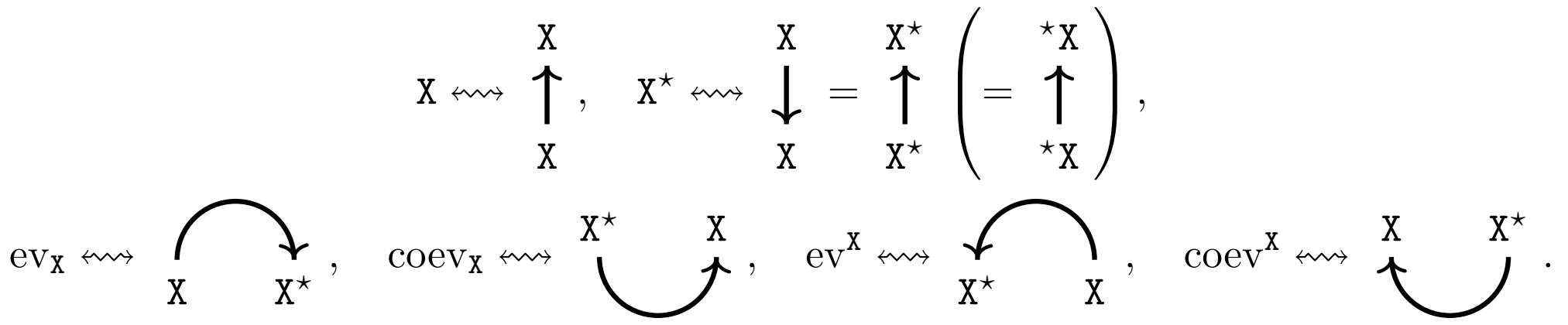

Note our reading conventions for duals, see also Convention 4A.3. The zigzag relations (4B-3) and (4B-5) in these diagrams are

(4G-2)

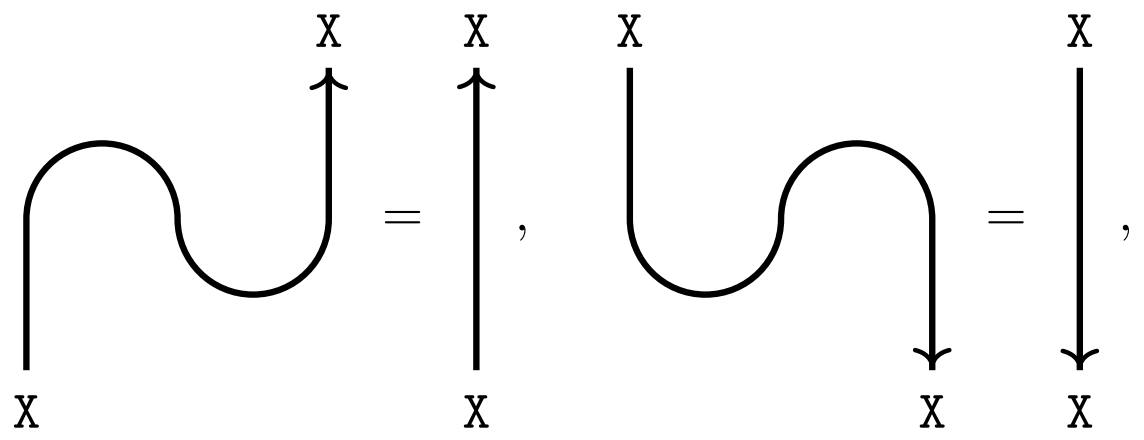

including mirrors, which imply:

(4G-3)

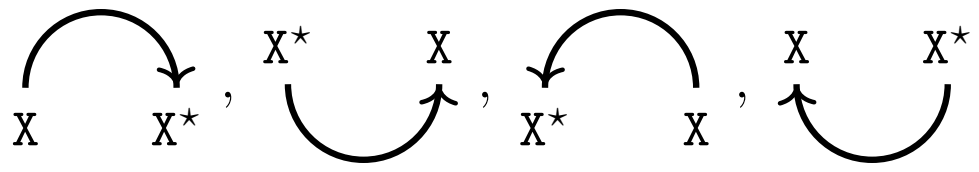

are invertible operations.

Let us prove some lemmas using this diagrammatics.

**Lemma 4G.4.** *For all* $f \in \mathbf{C}$, *where* $\mathbf{C} \in \mathbf{PCat}$, *we have*

(4G-5)

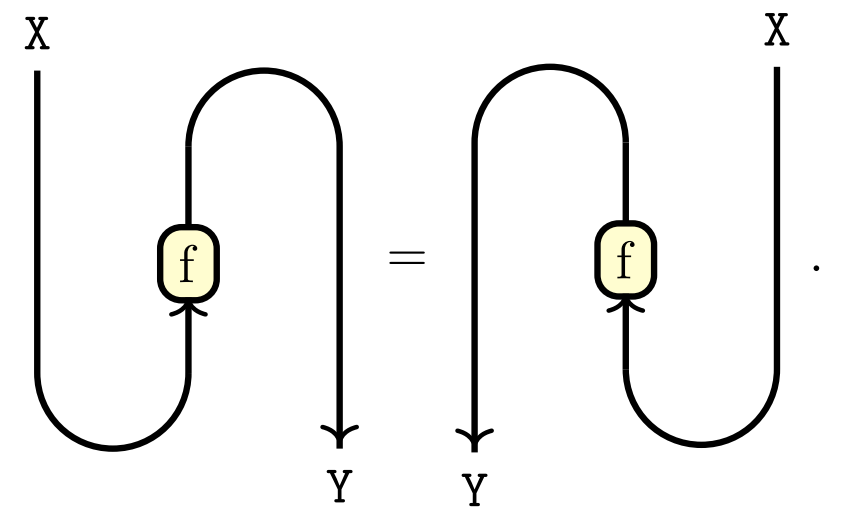

*Proof.* We calculate

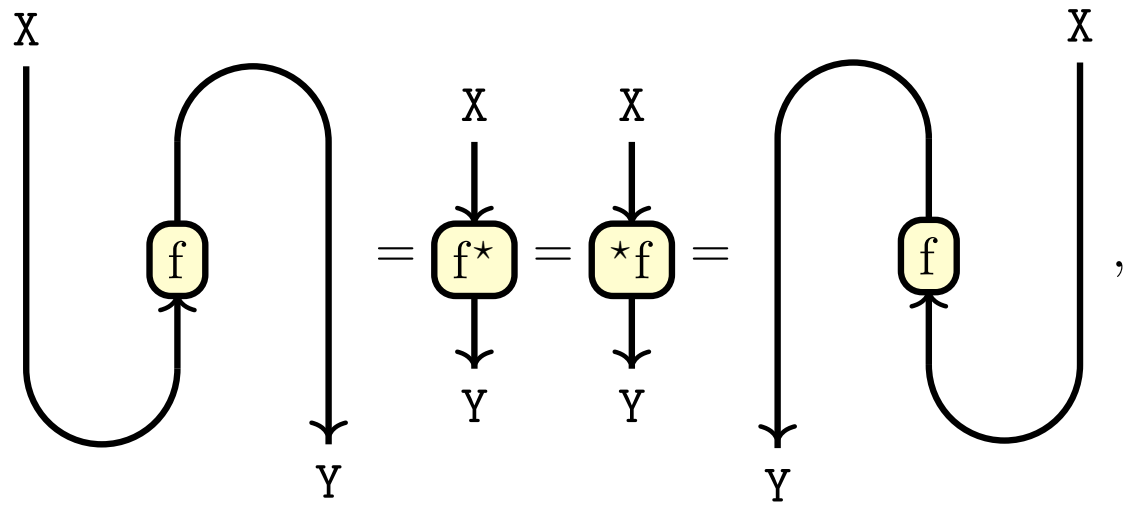

which is an application of Lemma 4F.9. □

**Lemma 4G.6.** *For all* f ∈ **C**, *where* **C** ∈ **PCat**, *we have*

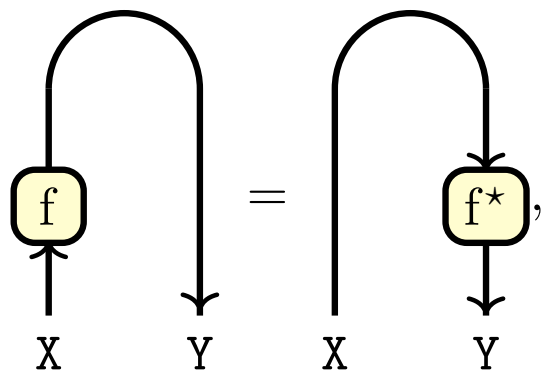

*including mirrors.*

These relations are called ***sliding***.

*Proof.* Using (4G-3) this is easy:

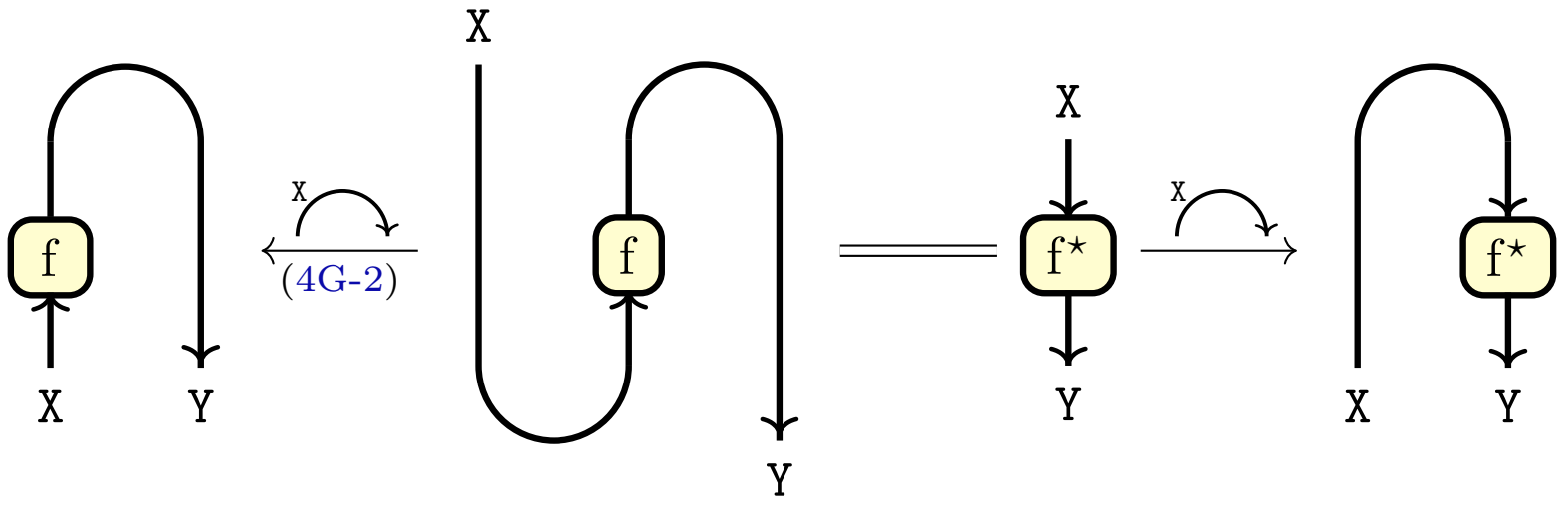

where we have used (4G-2). □

Recall that a pivotal structure was an additional choice of an isomorphism $\pi_X\colon X \xrightarrow{\cong} X^{\star\star}$. One such choice, sometimes called the ***canonical choice***, is

(4G-7) $\pi_X^{can}\colon X \to X^{\star\star}, \quad \pi_X =$

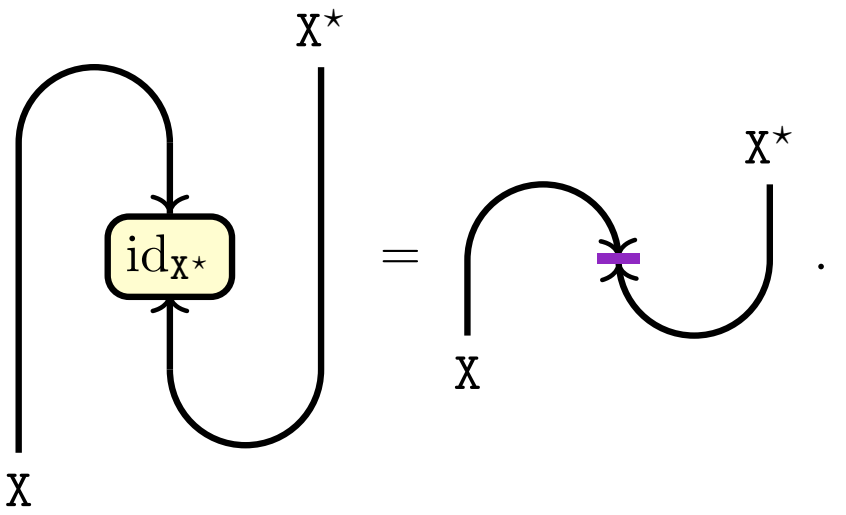

The colored marker is a shorthand notation for the corresponding identity morphism, which we also use below for different identities.

**Lemma 4G.8.** *For all* X ∈ **C**, *where* **C** ∈ **PCat**, *the morphism* $\pi_X^{can}$ *is invertible and*

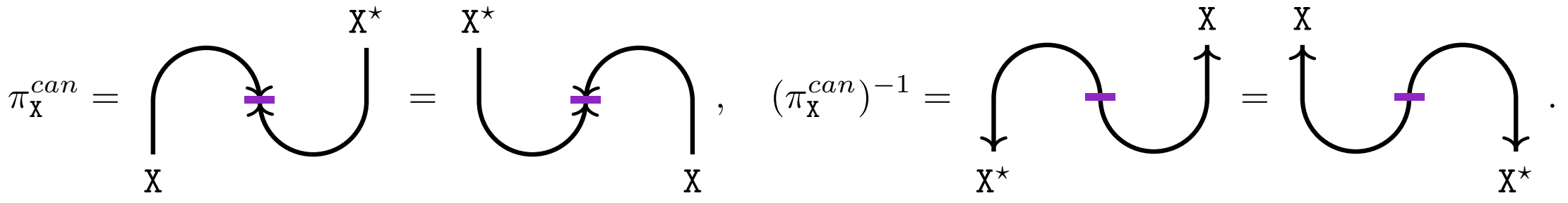

*Proof.* Note that, by definition, markers are identities and just turn orientations on diagrams around. Moreover, they are morphisms, so they slide. Hence, we have the diagrammatic equations

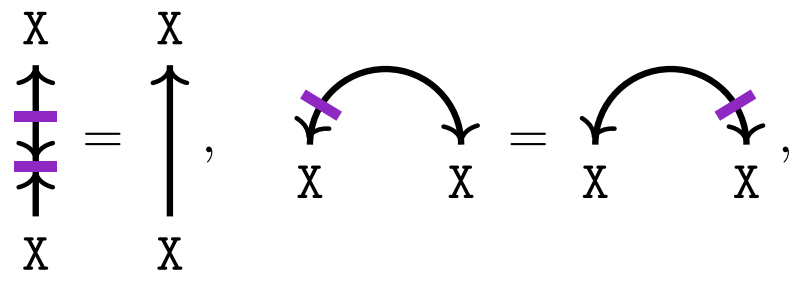

including mirrors. Now

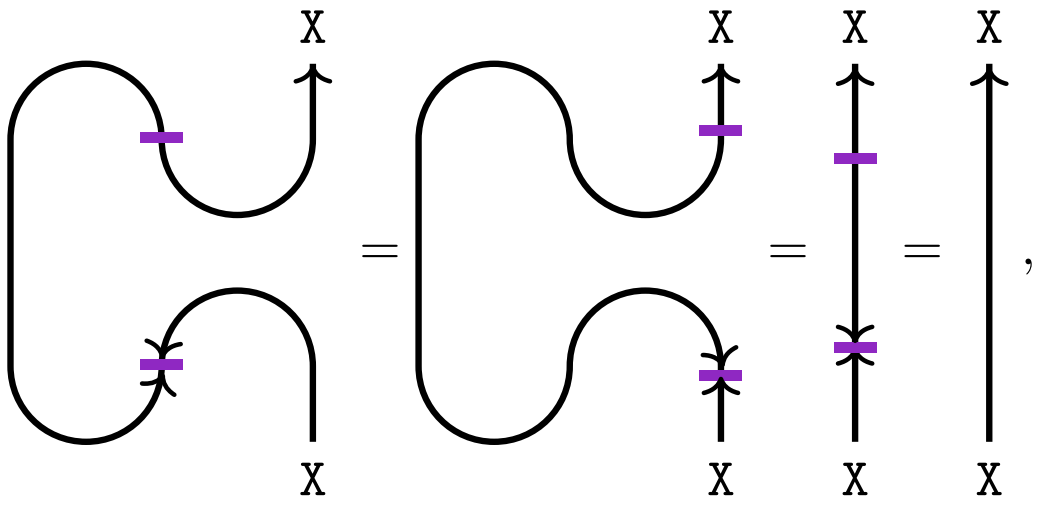

is one of the equalities we need to check; the others being similar. □

**Example 4G.9.** The canonical pivotal structure in examples is as follows.

(i) Using the choice and notation from Example 4C.3 for $\mathbf{fdVec}_{\Bbbk}$, we see that

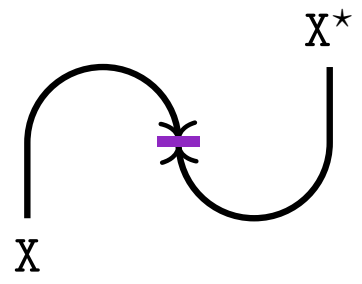

$$: \mathtt{X} \to \mathtt{X}^{\star\star},\ x \mapsto \textstyle\sum_{i=1}^{n} x \otimes x_i^{\star} \otimes x_i^{\star\star} \mapsto x^{\star\star},$$

which is independent of the choice of dual bases, *cf.* to the begin of Section 1F.

(ii) Recall from Example 4C.5 that for $\mathbf{Vec}_{\mathbb{C}}^{\omega}(\mathbb{Z}/2\mathbb{Z})$ the (co)evaluations are basically multiplication with $\pm 1$, giving the two possible choices

$$\pi_{\mathtt{1}}^{can}\colon \mathtt{1} \to \mathtt{1}, 1 \mapsto 1, \quad \pi_{\mathtt{1}}^{can}\colon \mathtt{1} \to \mathtt{1}, 1 \mapsto -1.$$

(Careful in the final example with the difference between the fonts.) ◇

The formal rule of manipulation of these diagrams is:

(4G-10) "Two diagrams are equivalent if they are related by scaling or by a planar isotopy." .

**Theorem 4G.11.** *The graphical calculus is consistent, i.e. two morphisms are equal if and only if their diagrams are related by* (4G-10)*.*

*Proof.* Note that it is crucial to have $\mathtt{X} \cong \mathtt{X}^{\star\star}$ which is key to have well-defined diagrammatics:

(4G-12)
$$\overset{\mathtt{X}^{\star\star}}{\underset{\mathtt{X}^{\star\star}}{\uparrow}} = \overset{\mathtt{X}^{\star}}{\underset{\mathtt{X}^{\star}}{\downarrow}} = \overset{\mathtt{X}}{\underset{\mathtt{X}}{\uparrow}},$$

and the isomorphism between the left and right sides in (4G-12) is the choice of pivotal structure, see *e.g.* (4G-7). Moreover, the zigzag relations in terms of diagrams (4G-2) and the identification of functors $_^{\star} = {}^{\star}_$, which gives (4G-5), ensure that one has all planar isotopies. □

4H. **Generalizing traces.** Let us continue with a motivating example.

**Example 4H.1.** Take $\mathbf{Mat}_{\Bbbk}$, the skeleton of $\mathbf{Vec}_{\Bbbk}$, which is pivotal with

$$\mathtt{n} = \mathtt{n}^{\star} = {}^{\star}\mathtt{n}, \quad \mathrm{ev}_{\mathtt{n}} = \mathrm{ev}^{\mathtt{n}}\colon \mathtt{nn} \to \mathtt{1},\ \mathrm{ev}_{\mathtt{n}} = \begin{pmatrix} e_1 & \cdots & e_n \end{pmatrix}, \quad \mathrm{coev}_{\mathtt{n}} = \mathrm{coev}^{\mathtt{n}}\colon \mathtt{1} \to \mathtt{nn},\ \mathrm{coev}_{\mathtt{n}} = \begin{pmatrix} e_1 \\ \vdots \\ e_n \end{pmatrix}.$$

Here $\{e_1, ..., e_n\}$ denotes the standard basis of $\Bbbk^n$ (which is secretly $\mathtt{n}$, of course) . Thus, given any $\mathrm{f} = (a_{ij})_{i,j=1,...,n} \in \mathrm{End}_{\mathbf{Mat}_{\Bbbk}}(\mathtt{n})$, we can calculate, keeping Convention 4A.5 in mind, that

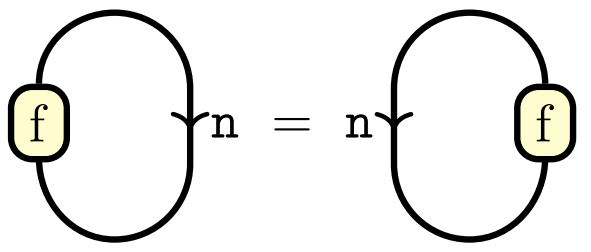

$$= \textstyle\sum_{i=1}^{n} a_{ii}.$$

This is the classical trace of the matrix f. Very explicitly, if $n = 2$ and $\mathrm{f} = \left(\begin{smallmatrix} a & b \\ c & d \end{smallmatrix}\right)$, then the calculation boils down to the matrix multiplication

$$\underbrace{\begin{pmatrix} 1 & 0 & 0 & 1 \end{pmatrix}}_{\mathrm{ev}_2} \underbrace{\begin{pmatrix} a & b & 0 & 0 \\ c & d & 0 & 0 \\ 0 & 0 & a & b \\ 0 & 0 & c & d \end{pmatrix}}_{\mathrm{f} \otimes \mathrm{id}_2} \underbrace{\begin{pmatrix} 1 \\ 0 \\ 0 \\ 1 \end{pmatrix}}_{\mathrm{coev}_2} = a + d.$$

Moreover, we get the dimension of $\mathtt{n}$ via

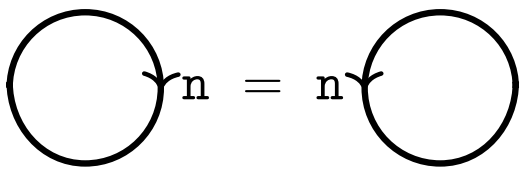

$$= n.$$

Indeed, we will think of circles as dimensions. ◇

**Definition 4H.2.** For $\mathrm{f} \in \mathrm{End}_{\mathbf{C}}(\mathtt{X})$, where $\mathbf{C} \in \mathbf{PCat}$, we define the ***right*** $\mathrm{tr}^{\mathbf{C}}(\mathrm{f})$ and ***left trace*** ${}^{\mathbf{C}}\mathrm{tr}(\mathrm{f})$ as the endomorphisms $\mathrm{tr}^{\mathbf{C}}(\mathrm{f}), {}^{\mathbf{C}}\mathrm{tr}(\mathrm{f}) \in \mathrm{End}_{\mathbf{C}}(\mathbb{1})$ given by

$\mathrm{tr}^{\mathbf{C}}(\mathrm{f}) =$ 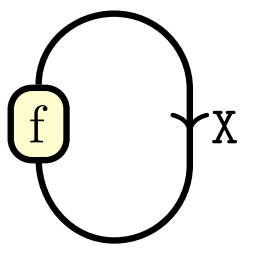

, ${}^{\mathbf{C}}\mathrm{tr}(\mathrm{f}) =$ 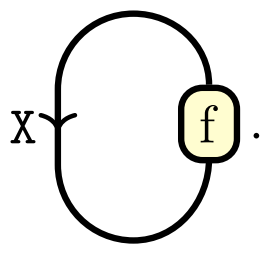

.

(Right trace ⇝ closing $\mathrm{f}\colon \mathtt{X} \to \mathtt{X}$ to the right, and *vice versa* for the left trace.) ◇

**Definition 4H.3.** For $\mathtt{X} \in \mathbf{C}$, where $\mathbf{C} \in \mathbf{PCat}$, we define the ***right*** $\dim^{\mathbf{C}}(\mathtt{X})$ and ***left dimension*** ${}^{\mathbf{C}}\dim(\mathtt{X})$ as the endomorphisms $\dim^{\mathbf{C}}(\mathtt{X}), {}^{\mathbf{C}}\dim(\mathtt{X}) \in \mathrm{End}_{\mathbf{C}}(\mathbb{1})$ given by

$$\dim^{\mathbf{C}}(\mathtt{X}) = \mathrm{tr}^{\mathbf{C}}(\mathrm{id}_{\mathtt{X}}) =$$ 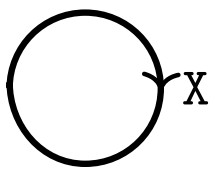 , $${}^{\mathbf{C}}\dim(\mathtt{X}) = {}^{\mathbf{C}}\mathrm{tr}(\mathrm{id}_{\mathtt{X}}) =$$ 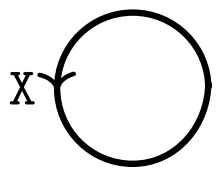 .

Again, the mantra is "circles = dimensions". ◇

Keeping track of right *versus* left is a bit tedious, so we hope for:

**Definition 4H.4.** A category $\mathbf{C} \in \mathbf{PCat}$ is called ***spherical*** if

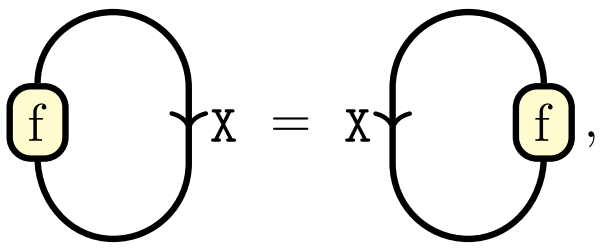

for $\mathtt{X} \in \mathbf{C}$ and all $\mathrm{f} \in \mathrm{End}_{\mathbf{C}}(\mathtt{X})$. ◇

*Remark* 4H.5. The name "spherical" comes from the idea that we can also see Feynman diagrams for endomorphisms as living on a sphere rather than being planar. Then

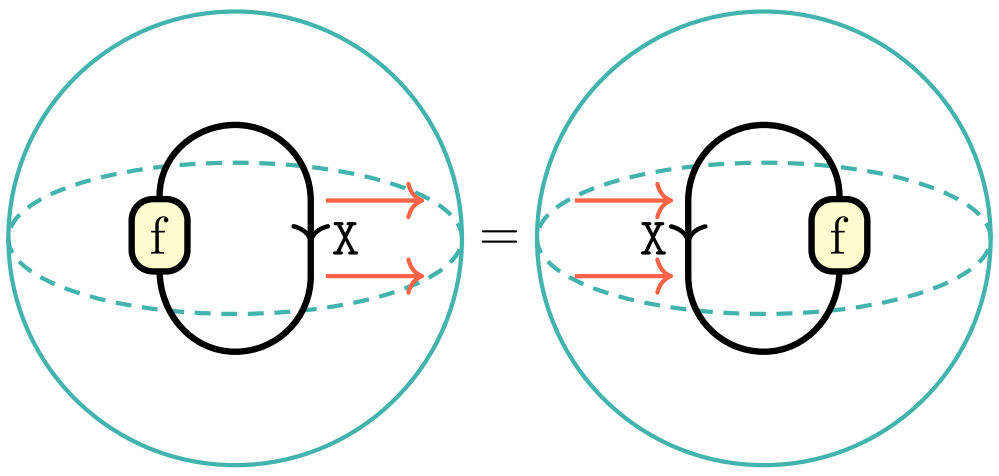

is just an isotopy which moves the strand around the sphere, a.k.a. the ***lasso move***. ◇

**Example 4H.6.** We have already seen in Example 4H.1 that traces and dimensions generalize traces and dimensions for matrices. Here are a few more examples.

(a) The category $\mathbf{fdVec}_{\mathbb{k}}$ with the standard (co)evaluations is spherical and traces and dimensions are the basis free definitions of the ones for $\mathbf{Mat}_{\mathbb{k}}$.

(b) The category $\mathbf{Vec}_{\mathbb{k}}^{\omega}(\mathrm{G})$ with the standard (co)evaluations is spherical and one has $\dim^{\mathbf{Vec}_{\mathbb{k}}^{\omega}(\mathrm{G})}(\mathtt{g}) = 1$ for all $\mathtt{g} \in \mathbf{Vec}_{\mathbb{k}}(\mathrm{G})$.

(c) The category $\mathbf{TL}$ with its generators being the (co)evaluations is spherical. The dimension of its generating object $\bullet$ is the morphism

$$\dim^{\mathbf{TL}}(\bullet) =$$ 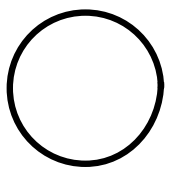 $$\in \mathrm{End}_{\mathbf{TL}}(\mathbb{1}).$$

Note that dimensions do not need to be integers, see the final example. ◇

As a warning, being spherical or not depends on choices:

**Example 4H.7.** For $\mathrm{G} = \mathbb{Z}/3\mathbb{Z}$, take $\zeta \in \mathbb{C}$ to be a complex primitive third root of unity, and let

$$d(\mathtt{i}) = \zeta^{i},\ i \in \{0,1,2\}.$$

Then there is a choice of (co)evaluations on $\mathbf{Vec}_{\mathbb{C}}(\mathbb{Z}/3\mathbb{Z})$ given by

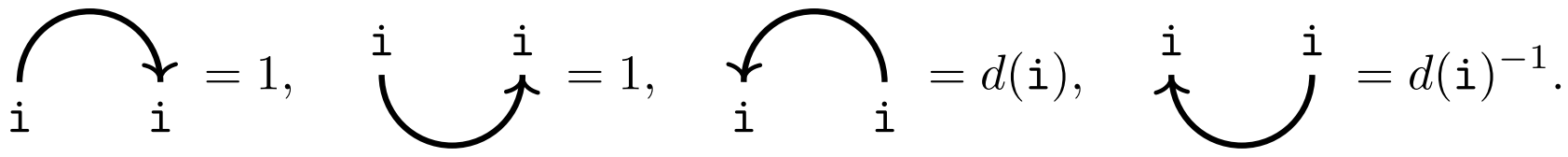

This gives

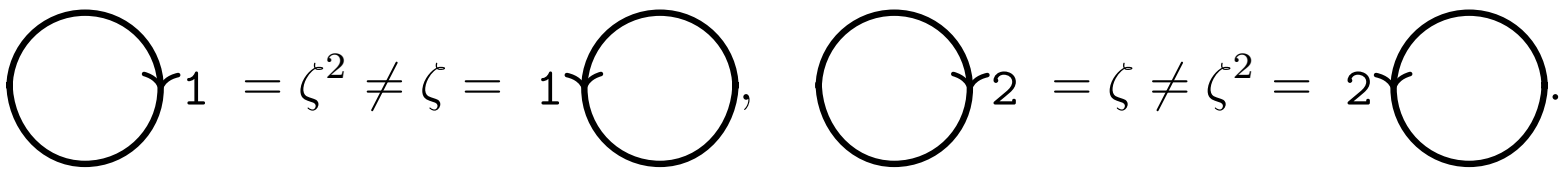

Thus, with this choice $\mathbf{Vec}_{\mathbb{C}}(\mathbb{Z}/3\mathbb{Z})$ is pivotal, but not spherical. ◇

*Remark* 4H.8. We have:

$$\text{rigid} \Leftarrow \text{pivotal} \Leftarrow \text{spherical},$$
$$\text{rigid} \not\Rightarrow \text{pivotal} \not\Rightarrow \text{spherical}.$$

This follows from definition and Example 4D.4 as well as Example 4H.7. ◇

We now discuss the generalization of the well-known properties of traces of matrices.

**Proposition 4H.9.** *For any* $\mathbf{C} \in \mathbf{PCat}$ *the following hold.*

*(i) We have*

$$\mathrm{tr}^{\mathbf{C}}(\mathrm{f}) = {}^{\mathbf{C}}\mathrm{tr}(\mathrm{f}) = \mathrm{f},$$

*for all* $\mathrm{f} \in \mathrm{End}_{\mathbf{C}}(\mathbb{1})$*. In particular,*

$$\dim^{\mathbf{C}}(\mathbb{1}) = {}^{\mathbf{C}}\dim(\mathbb{1}) = \mathrm{id}_{\mathbb{1}}.$$

*(ii) Traces are* $\mathrm{End}_{\mathbf{C}}(\mathbb{1})$*-linear, i.e.*

$$\mathrm{tr}^{\mathbf{C}}(\mathrm{f}\,.\,\mathrm{g}) = \mathrm{f}\,.\,\mathrm{tr}^{\mathbf{C}}(\mathrm{g}), \quad \mathrm{tr}^{\mathbf{C}}(\mathrm{g}\,.\,\mathrm{f}) = \mathrm{tr}^{\mathbf{C}}(\mathrm{g})\,.\,\mathrm{f}, \quad {}^{\mathbf{C}}\mathrm{tr}(\mathrm{f}\,.\,\mathrm{g}) = \mathrm{f}\,.\,{}^{\mathbf{C}}\mathrm{tr}(\mathrm{g}), \quad {}^{\mathbf{C}}\mathrm{tr}(\mathrm{g}\,.\,\mathrm{f}) = {}^{\mathbf{C}}\mathrm{tr}(\mathrm{g})\,.\,\mathrm{f},$$

*for all* $\mathrm{f} \in \mathrm{End}_{\mathbf{C}}(\mathbb{1})$ *and* $\mathrm{g} \in \mathrm{End}_{\mathbf{C}}(\mathtt{X})$*.*

*(iii) Traces are cyclic, i.e.*

$$\mathrm{tr}^{\mathbf{C}}(\mathrm{gf}) = \mathrm{tr}^{\mathbf{C}}(\mathrm{fg}), \quad {}^{\mathbf{C}}\mathrm{tr}(\mathrm{gf}) = {}^{\mathbf{C}}\mathrm{tr}(\mathrm{fg}),$$

*for all* $\mathrm{f} \in \mathrm{Hom}_{\mathbf{C}}(\mathtt{X},\mathtt{Y})$ *and* $\mathrm{g} \in \mathrm{Hom}_{\mathbf{C}}(\mathtt{Y},\mathtt{X})$*.*

*(iv) We have*

$$\mathrm{tr}^{\mathbf{C}}(\mathrm{f}) = {}^{\mathbf{C}}\mathrm{tr}(\mathrm{f}^{\star}), \quad {}^{\mathbf{C}}\mathrm{tr}(\mathrm{f}) = \mathrm{tr}^{\mathbf{C}}(\mathrm{f}^{\star})$$

*for all* $\mathrm{f} \in \mathrm{End}_{\mathbf{C}}(\mathtt{X})$*. In particular, for all* $\mathtt{X} \in \mathbf{C}$*, we get*

$$\dim^{\mathbf{C}}(\mathtt{X}) = {}^{\mathbf{C}}\dim(\mathtt{X}^{\star}) = \dim^{\mathbf{C}}(\mathtt{X}^{\star\star}), \quad {}^{\mathbf{C}}\dim(\mathtt{X}) = \dim^{\mathbf{C}}(\mathtt{X}^{\star}) = {}^{\mathbf{C}}\dim(\mathtt{X}^{\star\star}).$$

*Proof.* (i) and (ii). The short argument is that f is a floating bubble, *cf.* Proposition 2J.4.

(iii). By right-left symmetry, we only need to calculate

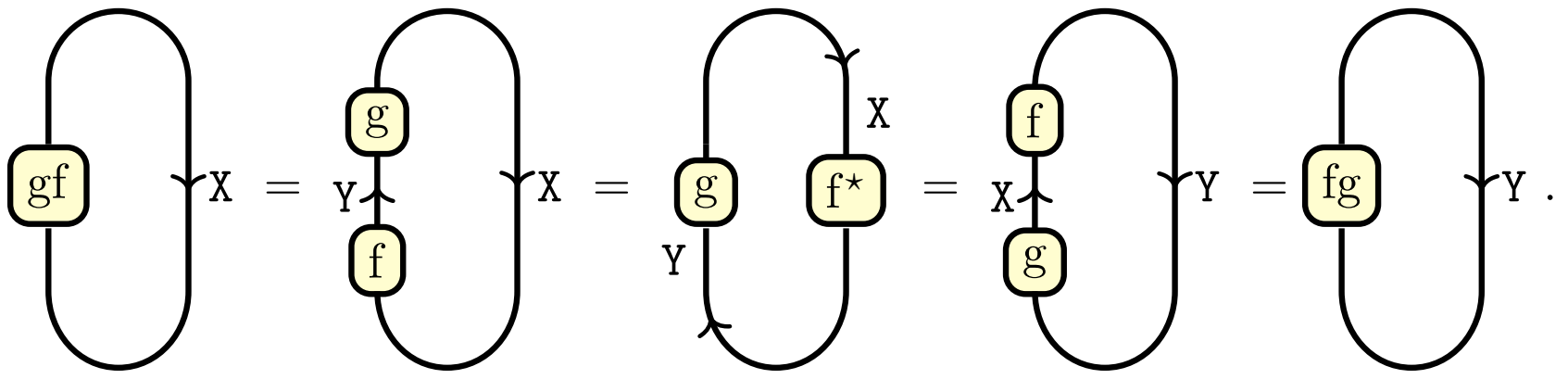

(iv). Sliding immediately gives

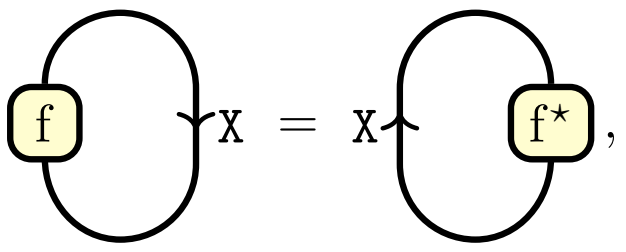

which proves the claim by right-left symmetry. □

**Lemma 4H.10.** *For* $\mathrm{F} \in \mathrm{Hom}_{\star}(\mathbf{C},\mathbf{D})$ *and all* $\mathtt{X} \in \mathbf{C}$ *and* $\mathrm{f} \in \mathrm{End}_{\mathbf{C}}(\mathtt{X})$ *we have*

$$\mathrm{tr}^{\mathbf{D}}\big(\mathrm{F}(\mathrm{f})\big) = \mathrm{F}\big(\mathrm{tr}^{\mathbf{C}}(\mathrm{f})\big), \quad {}^{\mathbf{D}}\mathrm{tr}\big(\mathrm{F}(\mathrm{f})\big) = \mathrm{F}\big({}^{\mathbf{C}}\mathrm{tr}(\mathrm{f})\big),$$
$$\dim^{\mathbf{D}}\big(\mathrm{F}(\mathtt{X})\big) = \mathrm{F}\big(\dim^{\mathbf{C}}(\mathtt{X})\big), \quad {}^{\mathbf{D}}\dim\big(\mathrm{F}(\mathtt{X})\big) = \mathrm{F}\big({}^{\mathbf{C}}\dim(\mathtt{X})\big).$$

*Proof.* Note that rigid functors preserve (co)evaluations. □

In words, rigid functors preserve traces and dimensions, which motivates:

**Definition 4H.11.** $\mathbf{C},\mathbf{D} \in \mathbf{PCat}$ are called ***equivalent as pivotal categories***, if they are equivalent as rigid categories. ◇

**Example 4H.12.** Back to Example 4D.21: there are sign choices for $\mathbf{Vec}_{\mathbb{C}}^{\omega}(\mathbb{Z}/2\mathbb{Z})$ such that:

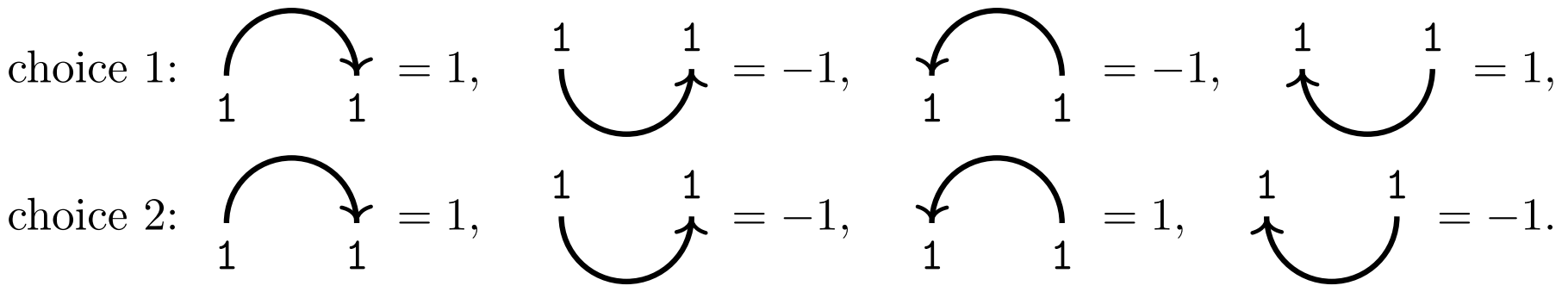

This gives

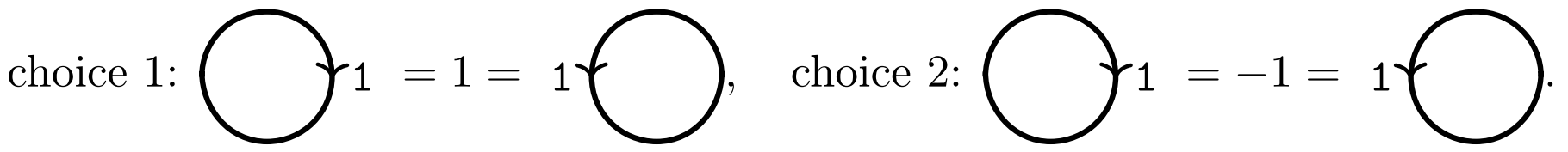

This shows, by Lemma 4H.10, that these choices do not give pivotal categories which are equivalent as pivotal categories. ◇

4I. **Algebras and coalgebras revisited.** We conclude with (co)algebras in rigid categories, whose modules have a right-left symmetry:

**Proposition 4I.1.** *Let* $\mathtt{A} \in \mathbf{C}$ *for* $\mathbf{C} \in \mathbf{RCat}$ *be an algebra.*

(i) *For every* $\mathtt{M} \in \mathbf{Mod}_{\mathbf{C}}(\mathtt{A})$ *its right dual* $\mathtt{M}^{\star}$ *has the structure of a left* $\mathtt{A}$ *module.*

(ii) *For every* $\mathtt{N} \in (\mathtt{A})\mathbf{Mod}_{\mathbf{C}}$ *its left dual* ${}^{\star}\mathtt{N}$ *has the structure of a right* $\mathtt{A}$ *module.*

*Similarly for coalgebras.*

*Proof.* By symmetry, it suffices to prove (i).

(i). We define a left action on $\mathtt{M}^{\star}$ via

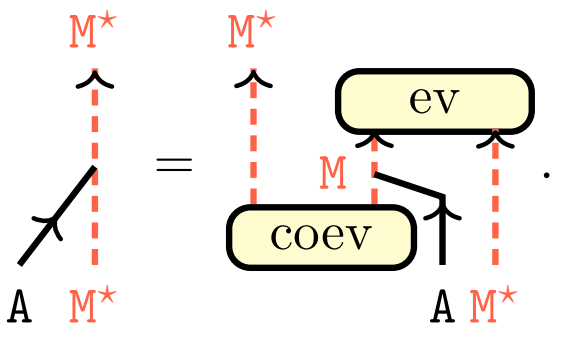

To show that this satisfies associativity and unitality is an easy zigzag argument. □

Thus:

**Proposition 4I.2.** *Let* $\mathtt{A} \in \mathbf{C}$ *for* $\mathbf{C} \in \mathbf{PCat}$ *be an algebra. Then, for every* $\mathtt{M} \in \mathbf{Mod}_{\mathbf{C}}(\mathtt{A})$ *its dual* $\mathtt{M}^{\star}$ *has the structure of a right and left* $\mathtt{A}$ *module. In particular,* $\mathtt{M}$ *itself has the structure of a left* $\mathtt{A}$ *module. Similarly for coalgebras.* □

4J. **Exercises.**

*Exercise* 4J.1. Prove Lemma 4C.7. ◇

*Exercise* 4J.2. Show that Theorem 4C.8 implies that the functors

$$\mathtt{X} \otimes _ : \mathbf{C} \to \mathbf{C}, \quad _ \otimes \mathtt{X} : \mathbf{C} \to \mathbf{C}$$

have duals (a.k.a. adjoints) given by

$$(\mathtt{X} \otimes _)^{\star} \cong \mathtt{X}^{\star} \otimes _, \quad {}^{\star}(\mathtt{X} \otimes _) \cong {}^{\star}\mathtt{X} \otimes _, \quad (_ \otimes \mathtt{X})^{\star} \cong _ \otimes {}^{\star}\mathtt{X}, \quad {}^{\star}(_ \otimes \mathtt{X}) \cong _ \otimes \mathtt{X}^{\star},$$

assuming the existence of duals of $\mathtt{X} \in \mathbf{C}$, where $\mathbf{C} \in \mathbf{MCat}$, of course. ◇

*Exercise* 4J.3. Show that $\mathbf{Vec}_{\Bbbk}(\mathrm{M}) \in \mathbf{RCat}$ if and only if M is a group. ◇

*Exercise* 4J.4. Show that $\mathbf{Vec}_{\Bbbk}$ is not rigid. ◇

*Exercise* 4J.5. Verify the claims in Example 4H.7. ◇

## 5. Braided categories – definitions, examples and graphical calculus

Recall that the difference between **1Cob** and **1Tan** was a choice of embedding. So how can we distinguish between these two categories using categorical algebra, *i.e.*:

What categorical framework can detect embeddings?

Studying the various methods of binding shoes, among other things, falls under the umbrella of knot theory, as shown in Figure 10. Following our usual approach, we aim to establish a categorical framework for studying knots.

5A. **A word about conventions.** This section is all about crossings.

*Convention* 5A.1. We will have over- and undercrossings, which are algebraic and not topological in nature. Our diagrammatic conventions for these are

overcrossing: $\beta = \beta^{+1} \rightsquigarrow$ 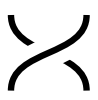, undercrossing: $\beta^{-1} \rightsquigarrow$ 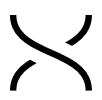.

These will come in various incarnations, *e.g.* with orientations, and our preferred choice will be to use overcrossings, and the undercrossings will be the inverses of the overcrossings. ◇

*Convention* 5A.2. We use the same conventions as in Convention 4A.2 for the various choices involved in the notions which we will see in this section. Moreover, and also as before, since we will use diagrammatics most of the time we usually omit the associators and unitors. ◇

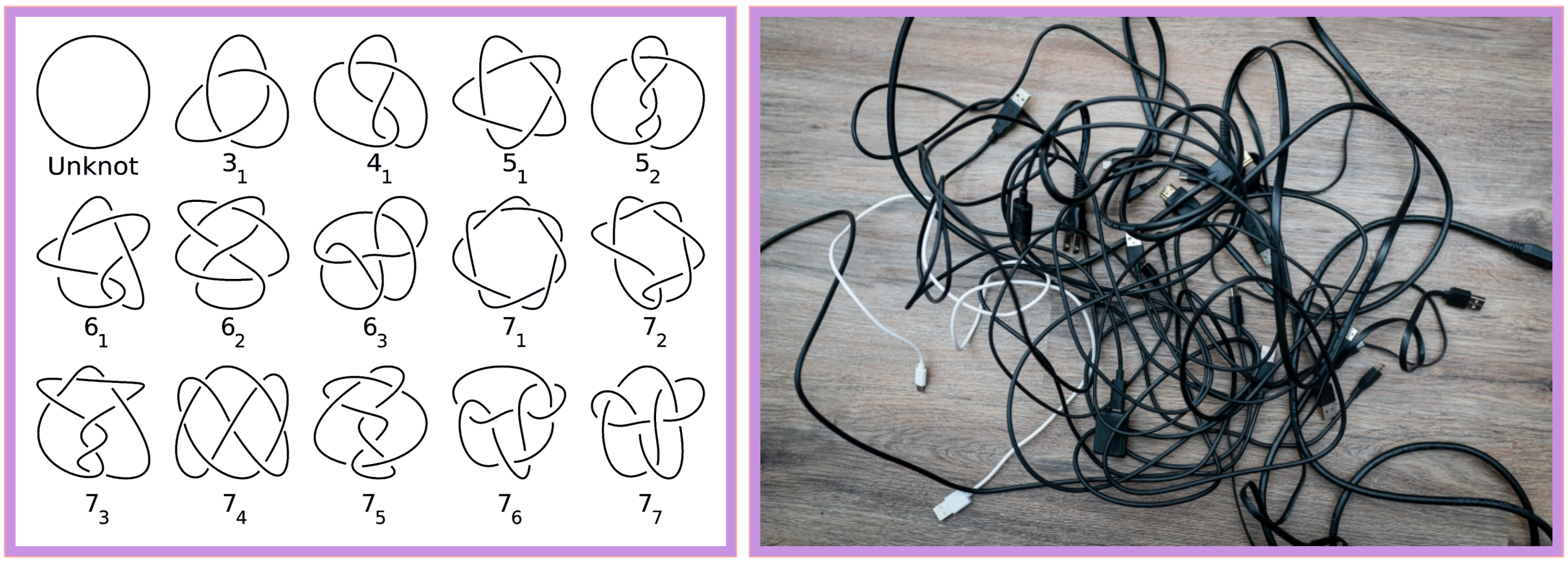

FIGURE 10. Left: Various embeddings of a string into three space. These are all the same object, say a rope, that only differ how they sit in space. Right: A real world knotting.
Pictures from https://en.wikipedia.org/wiki/Knot_theory and https://www.dtubbenhauer.com/lecture-geotop-2023.html.

*Convention* 5A.3. Our terminology "taking mirrors" includes crossing reversals of all displayed crossings, *e.g.*

original: 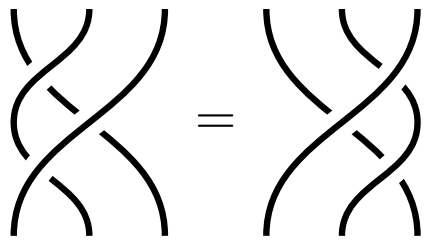, mirror: 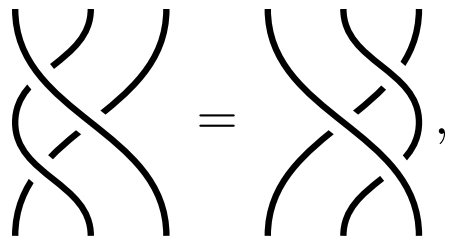, not a mirror: 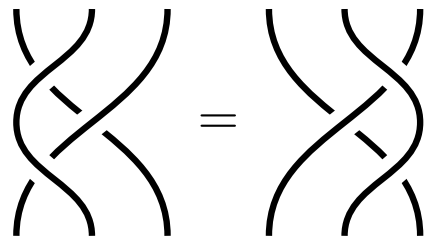.

All of these are valid relations, but the right equation is not a mirror of the left equation. ◇

5B. **Braided categories.** First, the main definition of this section:

**Definition 5B.1.** A ***braided category*** $(\mathbf{C},\beta)$ consists of

- a category $\mathbf{C}\in\mathbf{MCat}$;
- a collection of natural isomorphisms

(5B-2) $$\beta_{\mathtt{X},\mathtt{Y}}\colon \mathtt{XY}\xrightarrow{\cong}\mathtt{YX},$$

for all $\mathtt{X},\mathtt{Y}\in\mathtt{C}$, called ***braiding***;

such that

(i) the ***braided ⬡ equalities*** hold, *i.e.* we have commuting diagrams

(5B-3)
$$\begin{array}{c}\mathtt{X(YZ)}\xrightarrow{\alpha_{\mathtt{X},\mathtt{Y},\mathtt{Z}}}\mathtt{(XY)Z}\xrightarrow{\beta_{\mathtt{XY},\mathtt{Z}}}\mathtt{Z(XY)}\xrightarrow{\alpha_{\mathtt{Z},\mathtt{X},\mathtt{Y}}}\mathtt{(ZX)Y}\\ \mathtt{X(YZ)}\xrightarrow{\mathrm{id}_{\mathtt{X}}\otimes\beta_{\mathtt{Y},\mathtt{Z}}}\mathtt{X(ZY)}\xrightarrow{\alpha_{\mathtt{X},\mathtt{Z},\mathtt{Y}}}\mathtt{(XZ)Y}\xrightarrow{\beta_{\mathtt{X},\mathtt{Z}}\otimes\mathrm{id}_{\mathtt{Y}}}\mathtt{(ZX)Y}\end{array},$$
$$\begin{array}{c}\mathtt{(XY)Z}\xrightarrow{\alpha^{-1}_{\mathtt{X},\mathtt{Y},\mathtt{Z}}}\mathtt{X(YZ)}\xrightarrow{\beta_{\mathtt{X},\mathtt{YZ}}}\mathtt{(YZ)X}\xrightarrow{\alpha^{-1}_{\mathtt{Y},\mathtt{Z},\mathtt{X}}}\mathtt{Y(ZX)}\\ \mathtt{(XY)Z}\xrightarrow{\beta_{\mathtt{X},\mathtt{Y}}\otimes\mathrm{id}_{\mathtt{Z}}}\mathtt{(YX)Z}\xrightarrow{\alpha^{-1}_{\mathtt{Y},\mathtt{X},\mathtt{Z}}}\mathtt{Y(XZ)}\xrightarrow{\mathrm{id}_{\mathtt{Y}}\otimes\beta_{\mathtt{X},\mathtt{Z}}}\mathtt{Y(ZX)}\end{array},$$

for all $\mathtt{X},\mathtt{Y},\mathtt{Z}\in\mathtt{C}$.

◇

*Remark* 5B.4. Similarly as in Remark 2D.4, there is a hidden ***braided*** □ ***equality***:

(5B-5)
$$\begin{array}{ccc} \mathtt{ZX} & \xrightarrow{\mathrm{g}\otimes\mathrm{f}} & \mathtt{AY} \\ {\scriptstyle\beta_{\mathtt{X},\mathtt{Z}}}\uparrow & & \uparrow{\scriptstyle\beta_{\mathtt{Y},\mathtt{A}}} \\ \mathtt{XZ} & \xrightarrow[\mathrm{f}\otimes\mathrm{g}]{} & \mathtt{YA} \end{array},$$

which holds for all for all $\mathtt{X},\mathtt{Y},\mathtt{Z},\mathtt{A}\in\mathbf{C}$ and all $(\mathrm{f}\colon\mathtt{X}\to\mathtt{Y}),(\mathrm{g}\colon\mathtt{Z}\to\mathtt{A})\in\mathbf{C}$. ◇

**Lemma 5B.6.** *In any braided category* $\mathbf{C}$ *we have the* **Reidemeister 2 moves**, *i.e. for all* $\mathtt{X},\mathtt{Y}\in\mathbf{C}$ *there exist a natural isomorphism* $\beta_{\mathtt{Y},\mathtt{X}}^{-1}\colon\mathtt{YX}\xrightarrow{\cong}\mathtt{XY}$ *such that*

(5B-7)
$$\beta_{\mathtt{Y},\mathtt{X}}^{-1}\beta_{\mathtt{X},\mathtt{Y}}=\mathrm{id}_{\mathtt{XY}},\quad \beta_{\mathtt{X},\mathtt{Y}}\beta_{\mathtt{Y},\mathtt{X}}^{-1}=\mathrm{id}_{\mathtt{YX}}\quad\big(\Leftrightarrow\beta_{\mathtt{Y},\mathtt{X}}^{-1}\beta_{\mathtt{X},\mathtt{Y}}=\mathrm{id}_{\mathtt{XY}}=\beta_{\mathtt{Y},\mathtt{X}}\beta_{\mathtt{X},\mathtt{Y}}^{-1}\big)$$

*Proof.* (5B-2) implies that $\beta_{\mathtt{X},\mathtt{Y}}$ has an inverse, and we denote by $\beta_{\mathtt{Y},\mathtt{X}}^{-1}$. Note that $\beta_{\mathtt{Y},\mathtt{X}}^{-1}$ is natural in $\mathtt{X}$ and $\mathtt{Y}$ because $\beta_{\mathtt{X},\mathtt{Y}}$ is. □

**Lemma 5B.8.** *In any braided category* $\mathbf{C}$ *we have the* **Reidemeister 3 move**, *i.e.*

(5B-9)
$$\begin{array}{ccc} & (\mathtt{XY})\mathtt{Z}\xrightarrow{\beta_{\mathtt{X},\mathtt{Y}}\otimes\mathrm{id}_{\mathtt{Z}}}(\mathtt{YX})\mathtt{Z}\xrightarrow{\alpha_{\mathtt{Y},\mathtt{X},\mathtt{Z}}^{-1}}\mathtt{Y}(\mathtt{XZ})\xrightarrow{\mathrm{id}_{\mathtt{Y}}\otimes\beta_{\mathtt{X},\mathtt{Z}}}\mathtt{Y}(\mathtt{ZX})\xrightarrow{\alpha_{\mathtt{Y},\mathtt{Z},\mathtt{X}}}(\mathtt{YZ})\mathtt{X} & \\ {\scriptstyle\alpha_{\mathtt{X},\mathtt{Y},\mathtt{Z}}}\nearrow & & \searrow{\scriptstyle\beta_{\mathtt{Y},\mathtt{Z}}\otimes\mathrm{id}_{\mathtt{X}}} \\ \mathtt{X}(\mathtt{YZ}) & & (\mathtt{ZY})\mathtt{X} \\ {\scriptstyle\mathrm{id}_{\mathtt{X}}\otimes\beta_{\mathtt{Y},\mathtt{Z}}}\searrow & & \nearrow{\scriptstyle\alpha_{\mathtt{Z},\mathtt{Y},\mathtt{X}}} \\ & \mathtt{X}(\mathtt{ZY})\xrightarrow[\alpha_{\mathtt{X},\mathtt{Z},\mathtt{Y}}]{}(\mathtt{XZ})\mathtt{Y}\xrightarrow[\beta_{\mathtt{XZ}}\otimes\mathrm{id}_{\mathtt{Y}}]{}(\mathtt{ZX})\mathtt{Y}\xrightarrow[\alpha_{\mathtt{Z},\mathtt{X},\mathtt{Y}}^{-1}]{}\mathtt{Z}(\mathtt{XY})\xrightarrow[\mathrm{id}_{\mathtt{Z}}\otimes\beta_{\mathtt{X},\mathtt{Y}}]{}\mathtt{Z}(\mathtt{YX}) & \end{array}$$

*commutes for all* $\mathtt{X},\mathtt{Y},\mathtt{Z}\in\mathbf{C}$.

*Proof.* The diagrammatic proof is easy as we will see in Lemma 5C.5 below, so we are done after strictification, *cf.* Theorem 5D.7 stated below. □

**Definition 5B.10.** A braided category $\mathbf{C}$ is called ***symmetric***, if $\beta_{\mathtt{Y},\mathtt{X}}\beta_{\mathtt{X},\mathtt{Y}}=\mathrm{id}_{\mathtt{XY}}$ holds for all $\mathtt{X},\mathtt{Y}\in\mathbf{C}$. ◇

**Example 5B.11.** The motivating examples are the usual ones:

(a) The category **Set** with $\otimes=\times$ and $\mathbb{1}=\{\bullet\}$ can be endowed with a symmetric braiding by using the ***swap*** map
$$\tau_{\mathtt{X},\mathtt{Y}}\colon\mathtt{XY}\to\mathtt{YX},\ (x,y)\mapsto(y,x).$$

(b) Similarly, $\mathbf{Vec}_{\Bbbk}$ or $\mathbf{fdVec}_{\Bbbk}$ with $\otimes=\otimes_{\Bbbk}$ and $\mathbb{1}=\Bbbk$ can be endowed with a symmetric braiding by also using (the $\Bbbk$ linear incarnation of) the swap map.

Although technically speaking wrong, we invite the reader to think about a symmetric braiding as being some kind of swap map. ◇

As usual with choices, they are often not unique:

**Example 5B.12.** The category $\mathbf{Vec}_{\mathbb{C}}(\mathbb{Z}/2\mathbb{Z})$ with its standard monoidal structure can be endowed with two braidings (note that braidings $\beta_{\mathtt{i},\mathtt{j}}\in\mathrm{End}_{\mathbf{Vec}_{\mathbb{C}}(\mathbb{Z}/2\mathbb{Z})}(\mathtt{ij})\cong\mathbb{C}$ are scalars):

- by using the so-called ***standard braiding*** $\beta_{\mathtt{1},\mathtt{1}}^{st}=1$;
- by using the so-called ***super braiding*** $\beta_{\mathtt{1},\mathtt{1}}^{su}=-1$.

In both cases, the braiding involving $\mathtt{0}$ are trivial (meaning 1). ◇

**Proposition 5B.13.** *For any braided category* $\mathbf{C}$ *its Grothendieck classes* $\mathrm{K}_0(\mathbf{C})$ *form a commutative monoid with multiplication and unit*
$$[\mathtt{X}][\mathtt{Y}]=[\mathtt{XY}]=[\mathtt{Y}][\mathtt{X}],\quad 1=[\mathbb{1}].$$

*Proof.* The only thing we need to observe in addition to Proposition 2D.9 is that the existence of the braiding gives $[\mathtt{X}][\mathtt{Y}]=[\mathtt{Y}][\mathtt{X}]$ because we have (5B-2). □

**Example 5B.14.** The category $\mathbf{Vec}_{\Bbbk}(\mathrm{G})$ with its standard monoidal structure can only be braided if G is abelian. To see this observe that $\mathrm{K}_0\big(\mathbf{Vec}_{\Bbbk}(\mathrm{G})\big)\cong\mathrm{G}$ as groups, and thus Proposition 5B.13 applies. ◇ 5C. **Feynman diagrams for braided categories.** Example 5C.8 already suggests the following diagrammatic conventions. We take the ones for monoidal categories, see *e.g.* (2E-1), together with diagrammatic rules for braidings:

(5C-1)

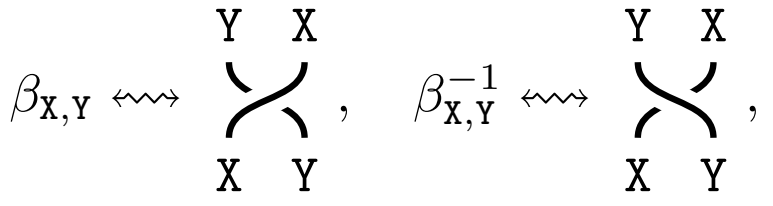

where (5B-2) implies that $\beta_{X,Y}$ has an inverse. Being inverses gives the ***Reidemeister* 2 *moves***, *i.e.* the diagrammatic analogs of (5B-7):

(5C-2)

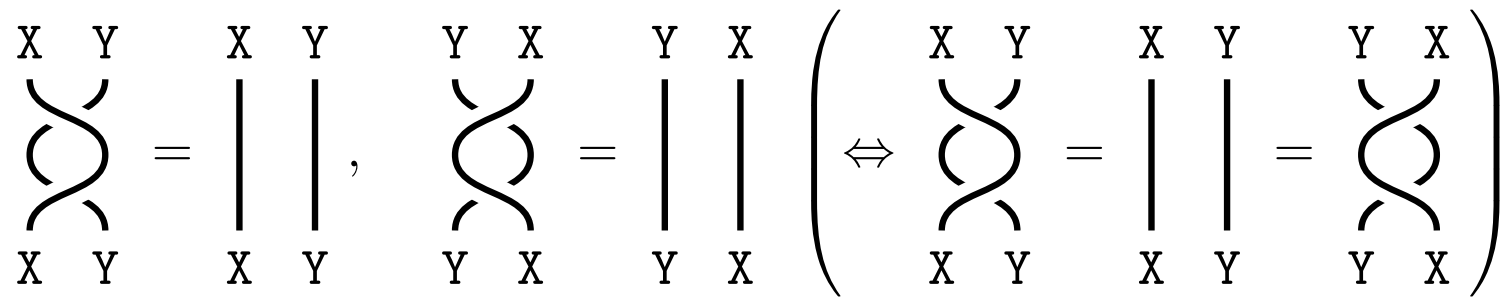

hold for all $X, Y \in \mathbf{C}$.

*Remark* 5C.3. We usually, following history, use the right-hand side in (5C-2) as the Reidemeister 2 moves. Further, beware that $\beta^{-1}_{Y,X}$ is the inverse of $\beta_{X,Y}$ and not $\beta^{-1}_{X,Y}$. ◇

The diagrammatic incarnations of the braided ⬡ and □ equalities in (5B-3) and (5B-5) are

(5C-4)

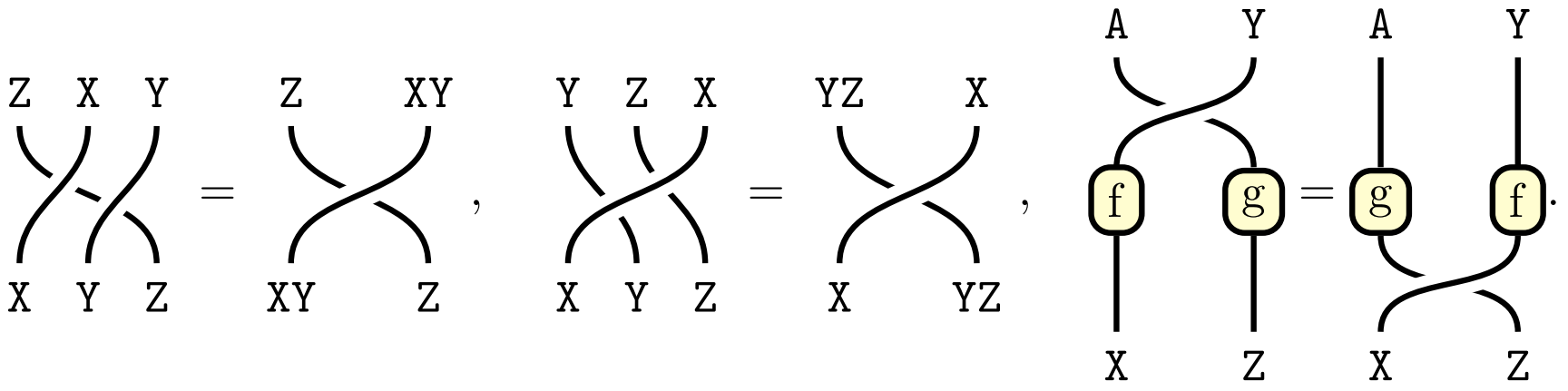

Similarly as in Lemma 4G.6, we call the right relation ***sliding***. We also have the ***Reidemeister* 3 *move***, *i.e.* the diagrammatic analog of (5B-9):

**Lemma 5C.5.** *In any braided category* $\mathbf{C}$ *we have the* **Reidemeister** 3 **move**, *i.e.*

(5C-6)

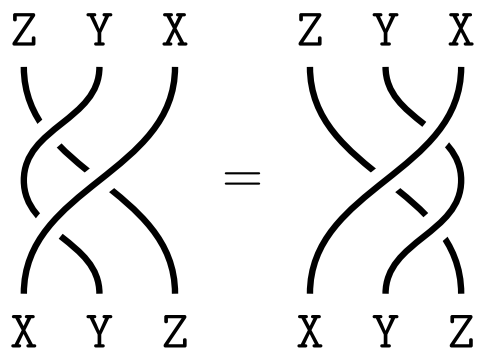

*holds for all* $X, Y, Z \in \mathbf{C}$.

*Proof.* We use the right equation in (5C-4) for a specific choice, *i.e.* we replace the labels and morphisms in (5C-4) as follows:

$$X \rightsquigarrow XY, Y \rightsquigarrow YX, Z \rightsquigarrow Z, A \rightsquigarrow Z,$$

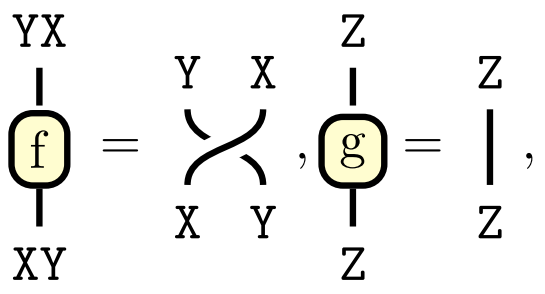

with the right-hand sides of all equations being the choices. This shows (5C-6). □

There are six versions of the Reidemeister 3 move (when fixing the boundary objects): each crossing can be an overcrossing or an undercrossing, but not all possibilities are topologically sensible. We list them now:

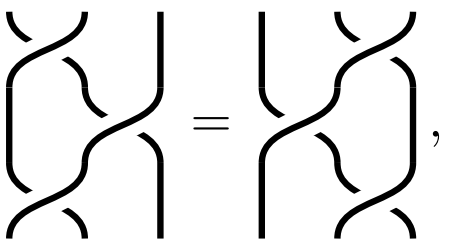

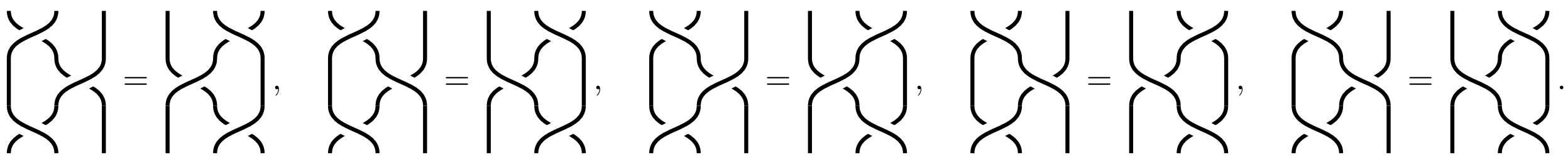

The ***non-defining Reidemeister* 3 *moves*** are the ones in the second row. These are in general valid, but not covered by Lemma 5C.5 so:

**Lemma 5C.7.** *In any braided category* **C** *we have the* **non-defining Reidemeister 3 moves** *hold for all* $\mathtt{X},\mathtt{Y},\mathtt{Z}\in\mathbf{C}$.

*Proof.* The fastest way to do this is algebraically. Write the defining relation as $\beta_1\beta_2\beta_1=\beta_2\beta_1\beta_2$ (the notation meaning that $\beta_i$ swaps strands $i$ and $i+1$) or, simplified, just $121=212$. These are all overcrossings, undercrossing will be indicated by $1^{-1}$ and $2^{-1}$. Now, to prove *e.g.* $1^{-1}21=212^{-1}$ we use simply compute:

$$(121=212)\Rightarrow(21=1^{-1}212)\Rightarrow(212^{-1}=1^{-1}21).$$

This uses $1^{-1}1=\mathrm{id}_{\mathtt{XYZ}}$ and $22^{-1}=\mathrm{id}_{\mathtt{XYZ}}$, see Lemma 5B.6. The other can be proven similarly. □

In the symmetric braided case we have $\beta_{\mathtt{Y},\mathtt{X}}\beta_{\mathtt{X},\mathtt{Y}}=\mathrm{id}_{\mathtt{XY}}$, which implies that

$$\beta_{\mathtt{X},\mathtt{Y}}=\beta_{\mathtt{X},\mathtt{Y}}^{-1}\rightsquigarrow$$ 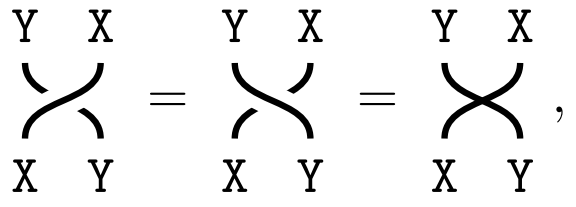

where the right-hand side is thus an appropriate shorthand notation. Hence, in symmetric braided categories the Reidemeister moves (5C-2) and (5C-6) then become

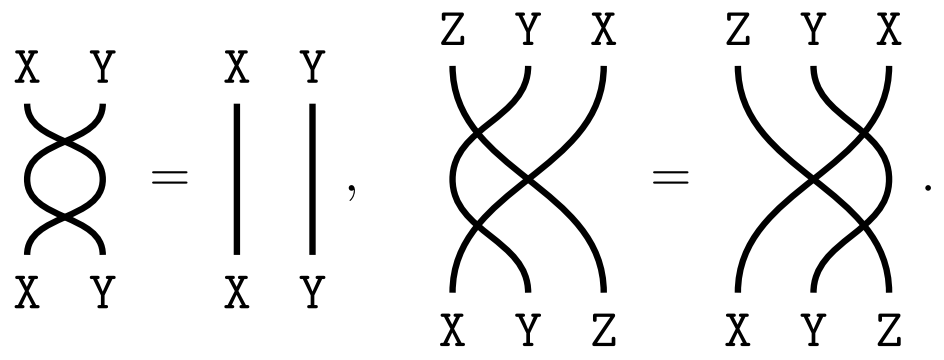

Note that, in this symmetric setting, there is no difference anymore between the defining Reidemeister 3 move and its five other forms.

**Example 5C.8.** Again, we have the notions of being "free as an XYZ":

(a) The free braided category generated by one object • is **qSym** from Example 3D.17. This category is important, so let us be completely explicit. We let **qSym** = ⟨**S**, **T** | **R**⟩ with

(5C-9) S : •, T : 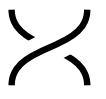 $:\bullet^2\to\bullet^2$, R : 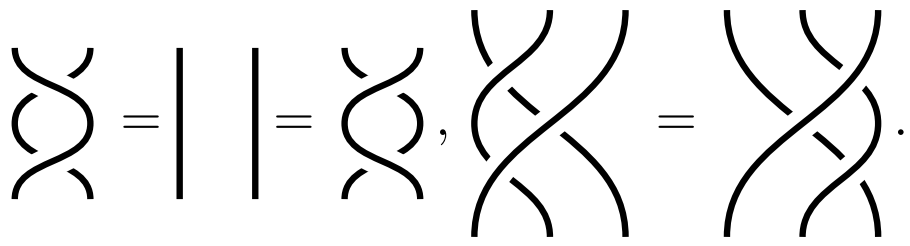

(We do not take mirrors.)

(b) The free symmetric braided category generated by one object • is **Sym** from Example 3D.2.

To see that this we observe that for one object (5C-2) and (5C-6) are equivalent to the braided ⬡ and □ equalities in (5B-3) and (5B-5). ◇

*Remark* 5C.10. Note that **qSym** has only overcrossings appearing in its definition. The undercrossings come into the game via invertibility. In particular,

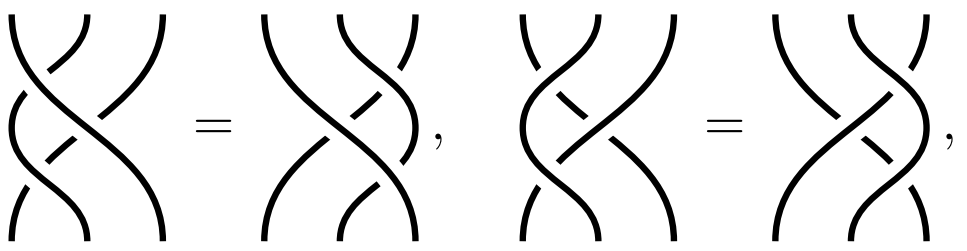

and all other versions of Reidemeister 3 moves are consequences and need not to be imposed, see Lemma 5C.7. ◇

The formal rule for braided Feynman diagrams is thus:

(5C-11) "Two diagrams are equivalent if they are related by scaling, by a planar isotopy, or braided ⬡ and □ equalities (5C-4)."

**Theorem 5C.12.** *The graphical calculus is consistent, i.e. two morphisms are equal if and only if their diagrams are related by* (5C-11).

*Proof.* The statement of the theorem just summarizes the discussion above: we have the Reidemeister 2 and 3 moves for strands, see (5C-2) and (5C-6), and we can slide coupons (5C-4). □

5D. **Braided functors.** As usual, we want the notion of functors between braided categories. To this end, recall that a monoidal functor $(\mathrm{F},\xi,\xi_{\mathbb{1}})$ was a functor with an additional choice of data, *cf.* Definition 2H.1. In contrast, being braided is a property:

**Definition 5D.1.** A functor $\mathrm{F}\in\mathbf{Hom}_{\otimes}(\mathbf{C},\mathbf{D})$ between braided categories is called ***braided*** if

$$\begin{array}{ccc} \mathrm{F}(\mathtt{X})\mathrm{F}(\mathtt{Y}) & \xrightarrow{\beta_{\mathrm{F}(\mathtt{X}),\mathrm{F}(\mathtt{Y})}} & \mathrm{F}(\mathtt{Y})\mathrm{F}(\mathtt{X}) \\ {\scriptstyle\xi_{\mathtt{X},\mathtt{Y}}}\downarrow & & \downarrow{\scriptstyle\xi_{\mathtt{Y},\mathtt{X}}} \\ \mathrm{F}(\mathtt{XY}) & \xrightarrow[\mathrm{F}(\beta_{\mathtt{X},\mathtt{Y}})]{} & \mathrm{F}(\mathtt{YX}) \end{array},$$

commutes for all $\mathtt{X},\mathtt{Y}\in\mathbf{C}$. ◇

We proceed as usual:

**Lemma 5D.2.** *The identity functor on a braided category is braided. Moreover, if* F *and* G *are braided functors, then so is* GF*.* □

**Example 5D.3.** We get the ***category of braided categories*** **BCat** and the ***category of braided functors*** $\mathbf{Hom}_{\beta}(\mathbf{C},\mathbf{D})$, whose natural transformation are monoidal natural transformations. ◇

**Definition 5D.4.** $\mathbf{C},\mathbf{D}\in\mathbf{BCat}$ are called ***equivalent as braided categories***, denoted by $\mathbf{C}\simeq_{\beta}\mathbf{D}$, if there exists an equivalence $\mathrm{F}\in\mathbf{Hom}_{\beta}(\mathbf{C},\mathbf{D})$. ◇

**Example 5D.5.** We also have the ***category of braided pivotal categories*** **BPCat** and the notion of equivalence for these is denoted by $\mathbf{C}\simeq_{\beta,\star}\mathbf{D}$. These equivalences use braided rigid functors which also form the ***category of braided rigid functors*** $\mathbf{Hom}_{\beta,\star}(\mathbf{C},\mathbf{D})$. ◇

**Definition 5D.6.** A category $(\mathbf{C},\beta)$ is called ***strict***, if it is strict as a monoidal category. ◇

As usual:

**Theorem 5D.7.** *For any braided category* $\mathbf{C}$ *there exists a strict braided category* $\mathbf{C}^{st}$ *which is braided equivalent to* $\mathbf{C}$*, i.e.* $\mathbf{C}\simeq_{\beta}\mathbf{C}^{st}$*.*

*Proof.* This is an almost immediate consequence of Theorem 2I.5, see [**JS93**, Theorem 2.5] for a detailed argument. □

5E. **Classifying braidings.** ***Classifying braidings***, meaning finding all possible braidings on $\mathbf{C}\in\mathbf{MCat}$ up to braided equivalence, is very difficult. So Theorem 5E.3 below is quite remarkable. Before we state it we need some preparation.

**Lemma 5E.1.** *Let* G *be abelian. Then the braidings on* $\mathbf{Vec}_{\Bbbk}^{\omega}(\mathrm{G})$ *(with its usual monoidal structure) are classified by twisted group homomorphisms* $\beta\colon\mathrm{G}\times\mathrm{G}\to\Bbbk^{*}$*, i.e. maps satisfying*

$$\begin{gathered}\omega(k,i,j)\beta(ij,k)\omega(i,j,k)=\beta(i,k)\omega(i,k,j)\beta(j,k)\\ \omega(j,k,i)^{-1}\beta(i,jk)\omega(i,j,k)^{-1}=\beta(i,k)\omega(j,i,k)^{-1}\beta(i,j).\end{gathered}\tag{5E-2}$$

*In particular, if* $\omega$ *is trivial, then braidings are classified by group homomorphisms* $\beta\colon\mathrm{G}\times\mathrm{G}\to\Bbbk^{*}$*.*

Such maps as in Lemma 5E.1 are also known as ***(twisted) bicharacters***. The above thus says that every braiding has an associated twisted bicharacter, which we denote by the same symbol.

*Proof.* By comparing (5B-3) and (5E-2) we see that each such $\beta$ can be used to define a braiding, and *vice versa.* □

A general philosophy, which we already have seen in Remark 2H.8, is that "Some cohomology theory should measure the obstruction of two braiding to be equivalent.". In fact, it is easy to see that functions satisfying (5E-2) for fixed $\omega$ form an abelian group $Z^{3}_{\omega}(\mathrm{G},\Bbbk^{\times})$ , which are the 3 cocyles of a cohomology group $H^{3}_{\omega}(\mathrm{G},\Bbbk^{\times})$, see *e.g.* [**EGNO15**, Section 8.4] for the definition. Indeed:

**Theorem 5E.3.** *Let* G *be abelian and fix a* 3 *cocycle* $\omega$*. Then* $\big(\mathbf{Vec}_{\Bbbk}^{\omega}(\mathrm{G}),\beta\big)\cong_{\beta}\big(\mathbf{Vec}_{\Bbbk}^{\omega}(\mathrm{G}),\beta'\big)$ *if and only if* $\beta$ *and* $\beta'$ *are cohomologically equivalent.*

*Proof.* In the end this is just a careful, but demanding, check of the involved definitions and commuting diagrams. Details are discussed in [**EGNO15**, Section 8.4]. □

**Example 5E.4.** Via Theorem 5E.3 we get the following, always using the standard monoidal structures.

(a) The $\mathrm{G}=1$ case of Theorem 5E.3 implies that $\mathbf{fdVec}_{\Bbbk}$ allows only one braiding if one fixes its standard monoidal structure since one can check that $H^{3}_{1}(1,\Bbbk^{\times})\cong 1$. See also Example 6F.7 later on.

(b) For $\mathrm{G} = \mathbb{Z}/2\mathbb{Z}$ and non-trivial $\omega$, has only two braidings:

$$\beta^{\pm}_{0,0} = \beta^{\pm}_{0,1} = \beta^{\pm}_{1,0} = 1, \quad \beta^{\pm}_{1,1} = \pm i \in \mathbb{C}.$$

The crucial calculation hereby is

$$\omega(1,1,1)\beta(11,1)\omega(1,1,1) = (-1)1(-1) = \beta(1,1)(-1)\beta(1,1) = \beta(1,1)\omega(1,1,1)\beta(1,1)$$
$$\Rightarrow \beta(1,1)^2 = -1.$$

It turns out that these are equivalent, *i.e.* $H^3_\omega(\mathbb{Z}/2\mathbb{Z}, \mathbb{C}^*) \cong 1$, and has only one (non-trivial) braiding up to equivalence.

(c) For $\mathrm{G} = \mathbb{Z}/2\mathbb{Z}$ and trivial $\omega$, we find precisely two possible solutions to (5E-2) and recover Example 5B.12.(a) and (b), since $H^3_1(\mathbb{Z}/2\mathbb{Z}, \mathbb{C}^\times) \cong \mathbb{Z}/2\mathbb{Z}$.

For groups of bigger order the necessary calculations to classify braidings are possible, but also much more difficult. ◇

5F. **The Reidemeister calculus.** Recall that Theorem 3D.10 identifies **1Cob** algebraically. We are now ready to state the analogs for **1Tan** and **1State**.

But first things first, let us be clear about the definitions of **qBr** and **oqBr**:

**Example 5F.1.** The ***(generic) quantum Brauer category*** **qBr** is the braided pivotal category generated by one self-dual object • with relations

(5F-2) R : 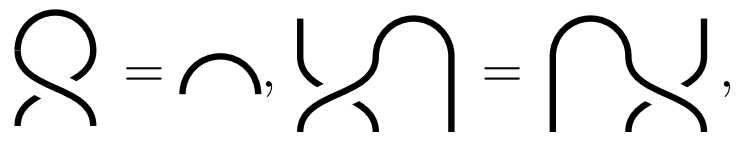

including mirrors, with structure maps

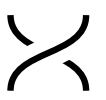: ••→••, 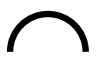: •• → $\mathbb{1}$, 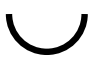: $\mathbb{1}$ → ••.

(Note that several relations come for free from "the braided pivotal category generated by".) ◇

**Example 5F.3.** The ***(generic) oriented quantum Brauer category*** **oqBr** is the braided pivotal category generated by one object • with relations

(5F-4) R : 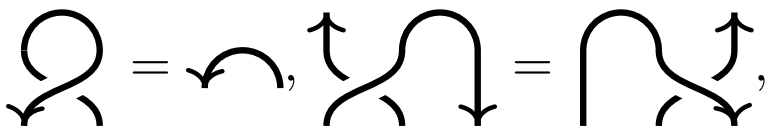

including mirrors, with structure maps

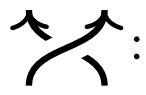: ••→••, 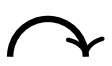: ••$^\star$ → $\mathbb{1}$, 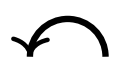: (•$^\star$)• → $\mathbb{1}$, 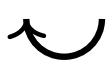: $\mathbb{1}$ → ••$^\star$, 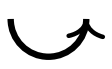: $\mathbb{1}$ → (•$^\star$)•.

(As in Example 5F.1, several relations are hidden in the definition.) ◇

The category **qBr** is also called ***BMW (Birman–Murakami–Wenzl) category*** (sometimes just ***BM category***) in the literature, while **oqBr** is sometimes called ***quantum walled Brauer category***.

*Remark* 5F.5. Note that (5F-2) and (5F-4) are not all defining relations as some are hidden in the phrase "generated as an XYZ". For example,

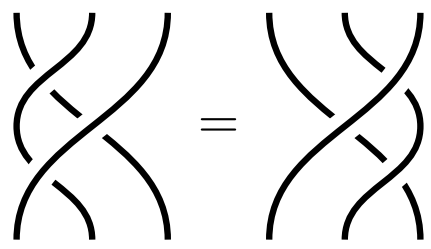

holds in both categories (with upward orientations for **oqBr**) and is part of being braided. ◇

Clearly, **1Tan** and **1State** are braided and pivotal with the evident structures. Recall also that we have the ***Reidemeister theorem***:

(5F-6) "Two (oriented) tangles in three space are isotopic if their projections (also known as tangle diagrams) are related by planar isotopies and Reidemeister moves 1-3, see (5F-7).

The ***topological Reidemeister*** 1***,*** 2 ***and*** 3 ***moves*** are all versions (not just mirrors) of

(5F-7) 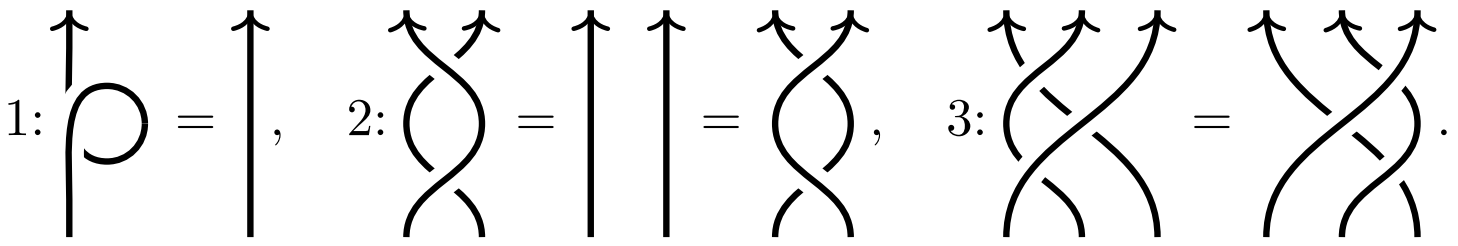

*Remark* 5F.8. Traditionally the topological Reidemeister 1 moves are usually illustrated vertically as in (5F-7), while their analogs in the Brauer calculus are traditionally sideways, see *e.g.* (5F-4), as it is a shorter composition of the generators. By (4G-3), these are the same data, and we will call the both the Reidemeister 1 moves. ◇

In our language, the categorical version of the Reidemeister theorem (5F-6) is:

**Theorem 5F.9.** *There exist braided rigid functors*

qR: **qBr** → **1Tan**, • ↦ •, 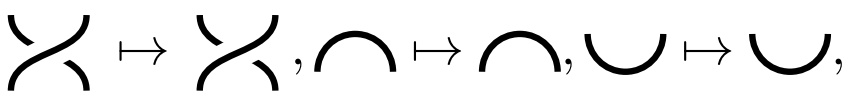

oqR: **oqBr** → **1State**, • ↦ •, 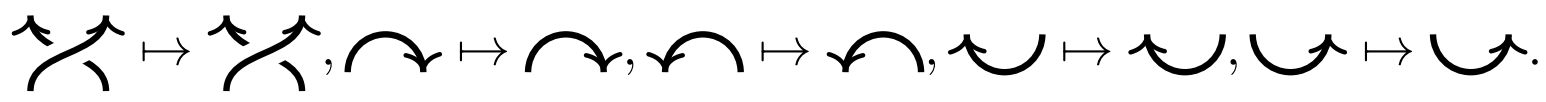

*Both functors are dense and fully faithful, thus,* **qBr** $\simeq_{\beta,\star}$ **1Tan** *and* **oqBr** $\simeq_{\beta,\star}$ **1State***.*

*Proof.* Let us sketch a proof. Exactly as for Theorem 3D.10, the main problem is to prove that these functors are faithful. To show faithfulness one needs to identify the generating relations of, say, **1Tan**. In words, one needs to identify what "isotopies of tangles" means for their projections.

To this end, the first step is to show that any tangle, appropriately defined, has a piecewise linear Morse presentation. A Morse presentation was already needed for the proof of fullness and has exactly the same meaning as in (3D-12), while piecewise linear basically is

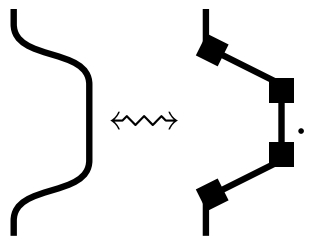

Here the ■ denotes the boundary points of the piecewise linear parts. (Note that all pictures in this proof are meant to represent topological objects.)

Then one needs to identify what isotopies are on these piecewise linear presentations and one gets the notion of △ ***equivalence*** $\sim_{\triangle}$ via △ ***moves***:

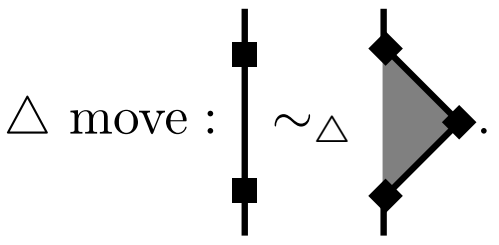

Here the triangle is not part of the link, but rather an illustration that no other strand is allowed to pass through it while performing the △ move. In words, two piecewise linear tangles are isotopic if and only if they are △ equivalent.

The first consequence of △ equivalence to notice is subdivision, *i.e.*

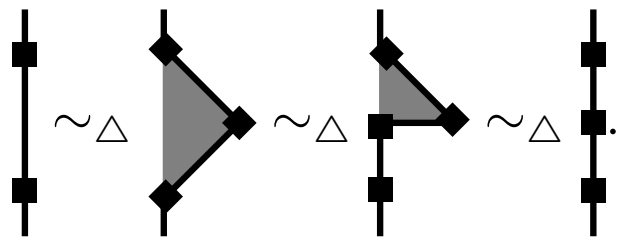

Thus, it remains to analyze how △ moves generically and locally project. In fact, the main upshot is that there are only finitely many possibilities one needs to check and one ends with precisely all possible versions (not just mirrors) of the Reidemeister moves *e.g.*:

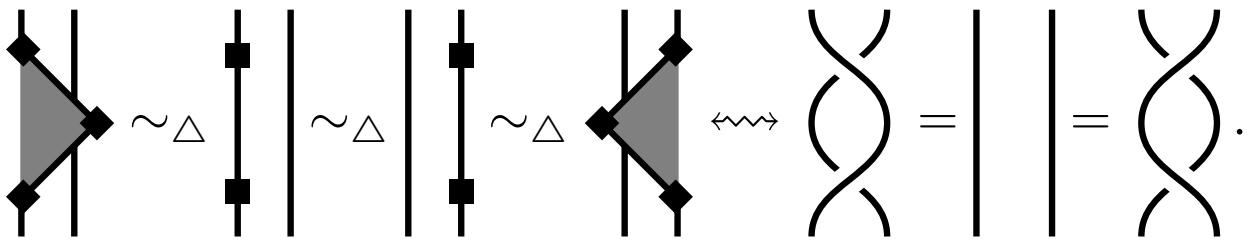

This established the Reidemeister theorem that two tangles are isotopic (in three space) if and only if their projections are related by planar isotopies and Reidemeister moves.

The final thing to check is that **oqBr** has enough relations to obtain all versions of the Reidemeister moves as well as all possible planar isotopies. Again, this is non-trivial as we *e.g.* imposed only certain types of Reidemeister relations such as only upwards oriented Reidemeister 2 moves. □

5G. **Minimal presentations of (oriented) quantum Brauer categories.** Without orientations we have Lemma 5C.7 which implies that we only need to impose one Reidemeister 3 move. Let us now state the analog for orientations.

**Lemma 5G.1.** *The following relations are sufficient to generate* **1State** *together with isotopies.*

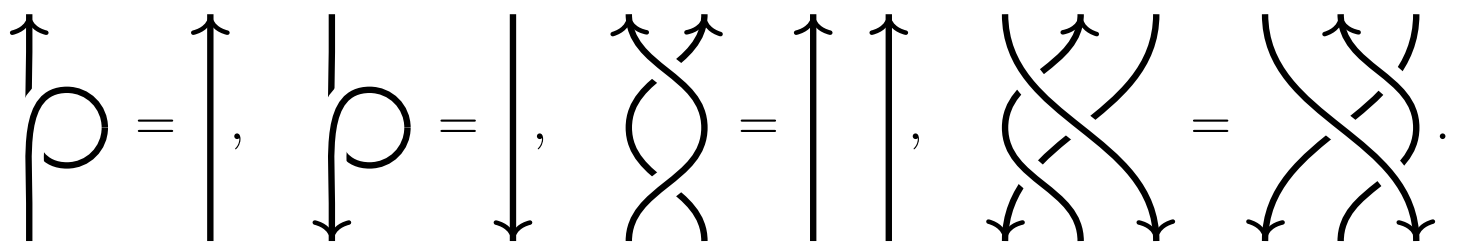

*Here we mean only the displayed versions, not mirrors or orientation reversals.*

*Proof.* [Pol10] does a good job, illustrating all the proofs carefully. □

To Do: correct form for Brauer

5H. **Twists.** Recall that the Reidemeister 2 and 3 moves (5C-2) and (5C-6) are consequences of the axioms of a braided category. For a braided pivotal category a good question would be whether the Reidemeister 1 moves as in (5F-7) follows from the combined axioms. Let us address this question.

**Definition 5H.1.** For $\mathtt{X} \in \mathbf{C}$ with $\mathbf{C} \in \mathbf{BPCat}$ the ***right*** $\mathrm{t}_{\mathtt{X}}$ and ***left*** $\mathrm{t}^{\mathtt{X}}$ ***twists*** are defined via

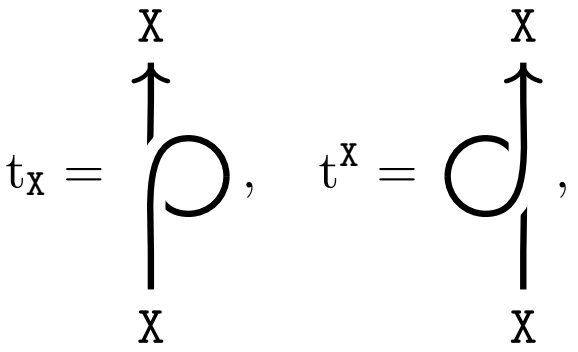

thus, the right twist has a cup and cap to the right, and *vice versa* for the left twist. ◇

There are three possible version of Reidemeister 1 moves. To explain them fix $\mathtt{X} \in \mathbf{C}$ for $\mathbf{C} \in \mathbf{BPCat}$. First, the ***(classical) Reidemeister*** 1 ***moves*** are

(5H-2)

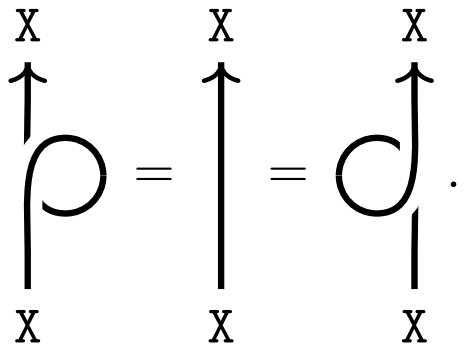

Second, the ***ribbon equation*** is

(5H-3)

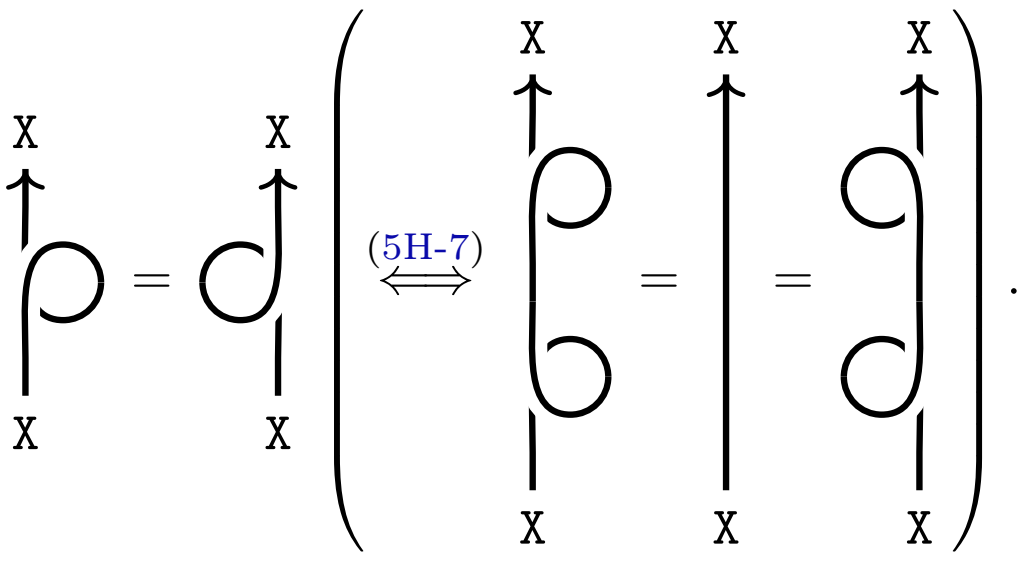

Finally, the ***framed Reidemeister*** 1 ***moves*** are

(5H-4)

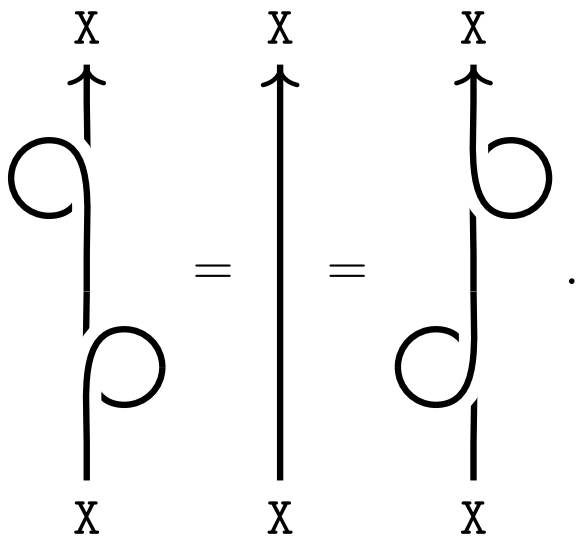

**Example 5H.5.** Clearly, (5H-2) ⇒ (5H-3) ⇒ (5H-4). But:

(a) In $\mathbf{fdVec}_{\mathbb{k}}$ with its standard monoidal structure, pairing and braiding we have

$$\mathrm{t}_{\mathtt{X}} \colon \mathtt{X} \to \mathtt{X},\ x_j \mapsto \textstyle\sum_{i=1}^{n} x_j \otimes x_i \otimes x_i^{\star} \mapsto \sum_{i=1}^{n} x_i \otimes x_j \otimes x_i^{\star} \mapsto x_j,$$
$$\mathrm{t}^{\mathtt{X}} \colon \mathtt{X} \to \mathtt{X},\ x_j \mapsto \textstyle\sum_{i=1}^{n} x_i^{\star} \otimes x_i \otimes x_j \mapsto \sum_{i=1}^{n} x_i^{\star} \otimes x_j \otimes x_i \mapsto x_j.$$

(Recall hereby our notation from Example 4C.3.) Thus, (5H-2) holds.

(b) For $\mathbf{Vec}_{\mathbb{C}}(\mathbb{Z}/2\mathbb{Z})$ we have discussed two pivotal structures in Example 4H.12. If we take the second pivotal strucutre therein together with the trivial braiding, then

$$\mathtt{t}_1\colon \mathtt{1}\to\mathtt{1},\ \mathtt{t}_1=-1,\quad \mathtt{t}^1\colon \mathtt{1}\to\mathtt{1},\ \mathtt{t}^1=-1.$$

Thus, (5H-4) and (5H-3) hold, but (5H-2) fails to hold.

(c) For $\mathbf{Vec}_{\mathbb{C}}(\mathbb{Z}/3\mathbb{Z})$ and $\zeta\in\mathbb{C}$ a primitive third root of unity we have discussed a pivotal structure in Example 4H.7. Taking this structure together with the trivial braiding we get

$$\mathtt{t}_1\colon \mathtt{1}\to\mathtt{1},\ \mathtt{t}_1=\zeta^2,\quad \mathtt{t}_2\colon \mathtt{2}\to\mathtt{2},\ \mathtt{t}_2=\zeta,$$
$$\mathtt{t}^1\colon \mathtt{1}\to\mathtt{1},\ \mathtt{t}^1=\zeta,\quad \mathtt{t}^2\colon \mathtt{2}\to\mathtt{2},\ \mathtt{t}^2=\zeta^2.$$

Thus, (5H-4) holds, but neither do (5H-2) or (5H-3).

In other words, (5H-2) $\not\Leftarrow$ (5H-3) $\not\Leftarrow$ (5H-4). ◇

**Lemma 5H.6.** *Fix* $\mathtt{X}\in\mathbf{C}$ *for* $\mathbf{C}\in\mathbf{BPCat}$.

*(i) The right and left twists are invertible with inverses*

(5H-7) $(\mathtt{t}_{\mathtt{X}})^{-1} =$ 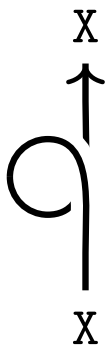 , $(\mathtt{t}^{\mathtt{X}})^{-1} =$ 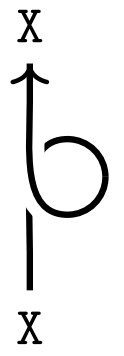 .

*(ii) We have* **sliding**, *i.e.*

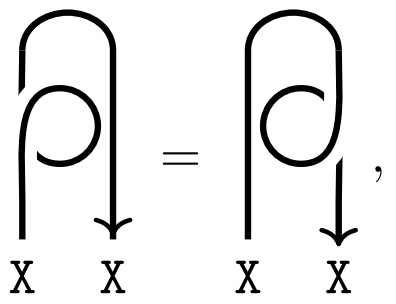

*including mirrors.*

*(iii) With respect to duality we have*

$(\mathtt{t}_{\mathtt{X}})^{\star} = \mathtt{t}^{\mathtt{X}^{\star}} =$ 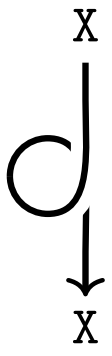 , $(\mathtt{t}^{\mathtt{X}})^{\star} = \mathtt{t}_{\mathtt{X}^{\star}} =$ 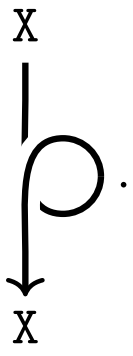 .

*(iv) All of the above maps are natural, i.e. they assemble into natural transformations.*

*Proof. (i).* See [TV17, Lemma 3.2].

*(ii).* Using (i) and Reidemeister moves we get

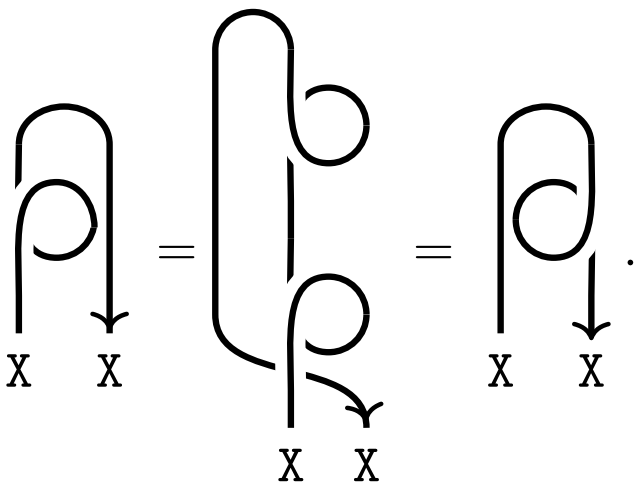

*(iii).* This is a direct application of sliding (ii) and zigzag.

*(iv).* The naturality of the twists is a direct consequence of the naturality of the braiding and the dualities. The proof is complete. □

Lemma 5H.6 immediately implies:

**Proposition 5H.8.** *For* $\mathbf{C}\in\mathbf{BPCat}$ *and all* $\mathtt{X}\in\mathbf{C}$ *the equation* (5H-4) *holds.* □

*Remark* 5H.9. We have

$$\text{Framed Reidemeister 1} \Leftarrow \text{ribbon equation} \Leftarrow \text{Reidemeister 1},$$
$$\text{framed Reidemeister 1} \not\Rightarrow \text{ribbon equation} \not\Rightarrow \text{Reidemeister 1}.$$

This follows by definition and Example 5H.5. ◇

**Example 5H.10.** The ***(generic) oriented framed quantum Brauer category*** **ofqBr** is the braided pivotal category generated by one object • with relations

(5H-11) R : 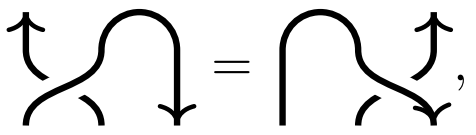

including mirrors, with structure maps

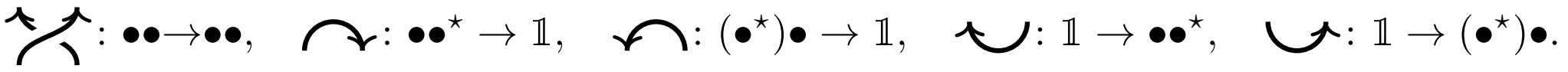

Note that Proposition 5H.8 implies that (5H-4) holds in **ofqBr**. In fact, a non-trivial argument shows that **ofqBr** is the free braided pivotal category generated by one object. ◇

**Example 5H.12.** Without further definition, there are of course also a non-oriented **fqBr**, non-quantum **ofBr** and only framed **fBr** versions of **ofBr**. ◇

5I. **Ribbon categories.** We have seen that the Reidemeister 1 moves (5H-2) are motivated from the topology of tangles in three space. As we will see, the framed Reidemeister 1 moves (5H-4) are also related to topology, but only in the form as in (5H-3). This motivates:

**Definition 5I.1.** A category $\mathbf{C} \in \mathbf{BPCat}$ is called ***ribbon*** if (5H-3) holds for all $\mathtt{X} \in \mathbf{C}$. ◇

**Example 5I.2.** The two categories in Example 5H.5.(a) and (b) are ribbon, while Example 5H.5.(c) is not. It will also follow from Lemma 5I.4 below that the non-spherical category from Example 5H.5.(c) cannot be ribbon. ◇

**Example 5I.3.** Note that being ribbon is a property and not a structure. So we can let the ***category of ribbon categories*** **RiCat** be the corresponding full subcategory of **BPCat**. ◇

**Lemma 5I.4.** *Let* $\mathbf{C} \in \mathbf{RiCat}$. *Then* $\mathbf{C}$ *is spherical.*

*Proof.* Using the Reidemeister calculus and sliding this is the calculation

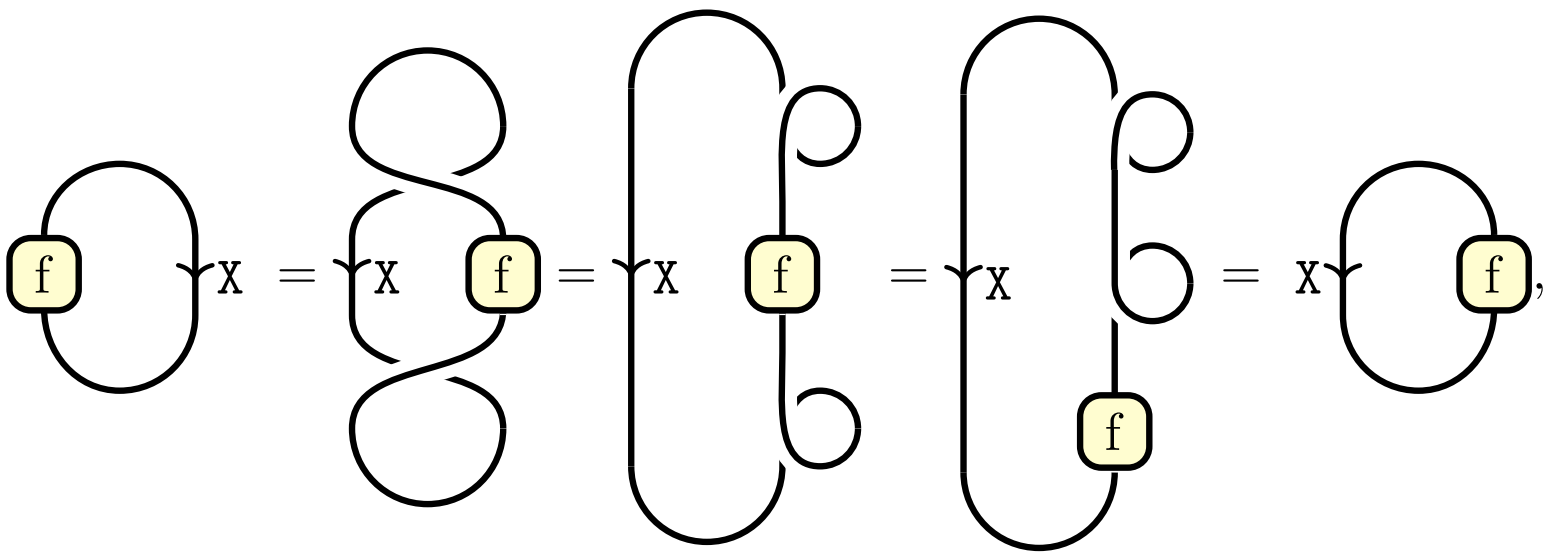

where the last step uses (5H-3). □

*Remark* 5I.5. The name ribbon comes from the following. If one takes a strip of paper (a thin and long strip works best) and performs the following

(5I-6) 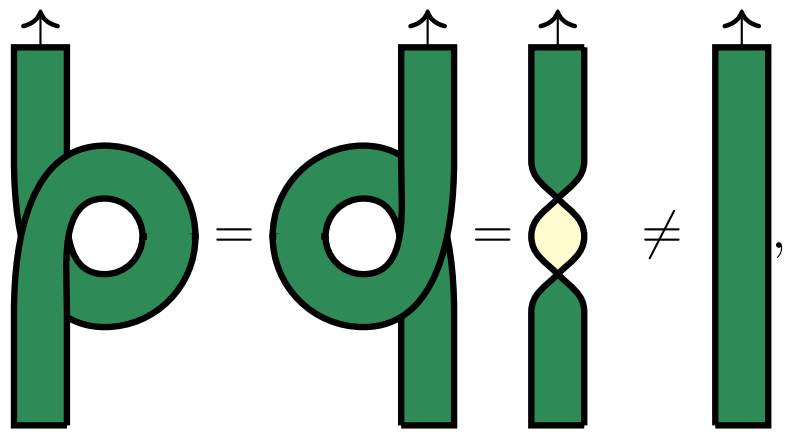

then we get the ribbon equation (5H-3). However, the paper strip is twisted, so (5H-2) does not hold. This motivates the definition of an important category in low-dimensional topology, called the ***category of ribbons*** (a.k.a. paper strips, see Figure 11) **1Ribbon**, which is, of course, a ribbon category. This category consists of oriented ribbons embedded in three space, *e.g.*

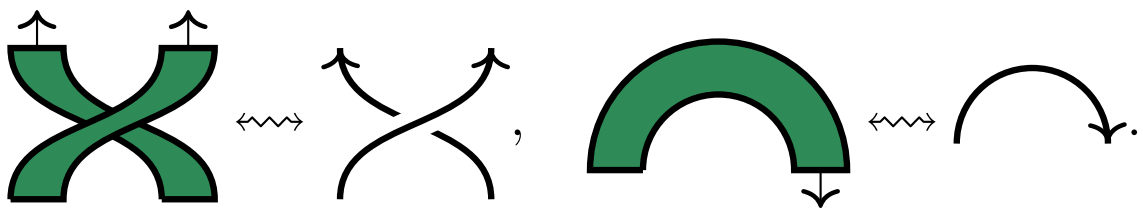

By making them arbitrary thin, these ribbons can be identified with the usual diagrams in the Reidemeister calculus with (5I-6) being the difference between ribbons, which have two sides, say green and white colored, and strings, which do not have any sides. ◇

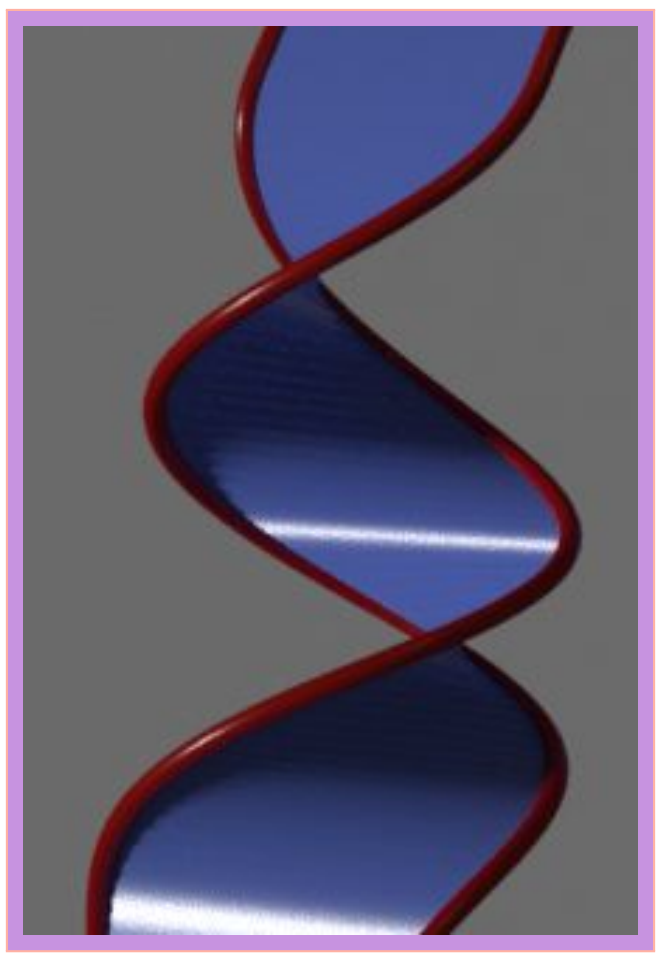

FIGURE 11. A twisted ribbon.

Picture from https://upload.wikimedia.org/wikipedia/commons/0/09/Ribbon_knot_8_20.jpg

**Example 5I.7.** The ***(generic) oriented ribbon quantum Brauer category*** **orqBr** is the braided pivotal category generated by one object • with relations

R : 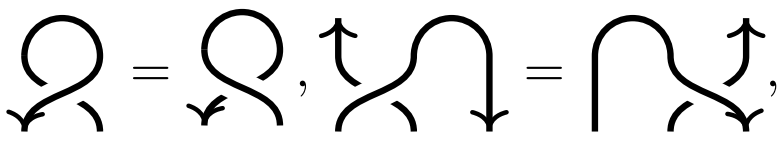

including mirrors. The structure maps are the usual ones, see *e.g.* Example 5H.10. ◇

The following is the point, but again non-trivial to prove.

**Theorem 5I.8.** *There exist a braided rigid functor*

$$\mathrm{orqR}\colon \mathbf{orqBr} \to \mathbf{1Ribbon}, \bullet \mapsto \bullet,$$

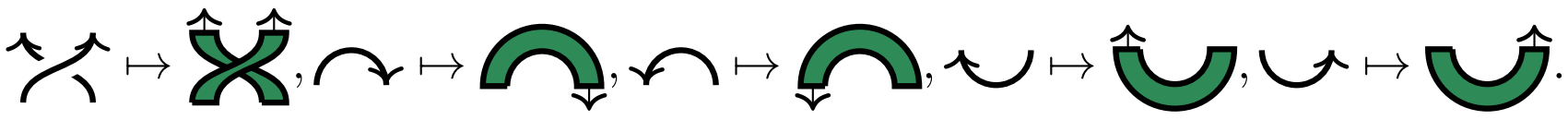

*The functor is dense and fully faithful, thus,* $\mathbf{orqBr} \simeq_{\beta,\star} \mathbf{1Ribbon}$.

*Proof.* A version of this theorem, which can be used to prove the formulation of it as above, is proven in [CP94, Section 5.3]. □

5J. **Algebras in braided categories.** Let us conclude this part with an continuation of Section 3E. A classical problem which we will address is to determine what condition on an algebra $\mathtt{A} \in \mathbf{Vec}_{\Bbbk}$ ensures that $\mathbf{Mod}_{\mathbf{C}}(\mathtt{A})$ is monoidal. Two known answers are:

- The category $\mathbf{Mod}_{\mathbf{C}}(\mathtt{A})$ is monoidal if $\mathtt{A}$ is commutative.
- The category $\mathbf{Mod}_{\mathbf{C}}(\mathtt{A})$ is monoidal if $\mathtt{A}$ is a bialgebra.

We will now discuss the categorical versions of these facts.

**Definition 5J.1.** A ***commutative algebra*** $\mathtt{A} = (\mathtt{A}, \mathrm{m}, \mathrm{i})$ in a braided category $\mathbf{C} \in \mathbf{BCat}$ is an algebra in $\mathbf{C}$ such that

(5J-2) 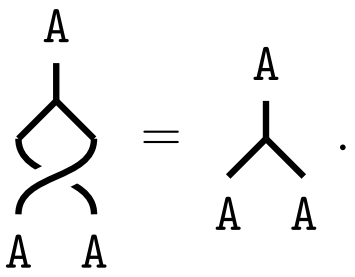

A ***cocommutative coalgebra*** $\mathtt{C} = (\mathtt{C}, \mathrm{d}, \mathrm{e})$ in a braided category $\mathbf{C}$ is, by definition, a commutative algebra in $\mathbf{C}^{op}$. ◇

By up-down symmetry, we can focus on algebras, since all constructions and statements for coalgebras are similar.

**Example 5J.3.** Definition 5J.1 generalizes the notions of (co)commutative (co)algebras, which can be recovered by taking $\mathbf{C} = \mathbf{Vec}_{\Bbbk}$. More generally, as we will see later, commutative algebras in $\big(\mathbf{Vec}_{\Bbbk}(\mathbb{Z}/2\mathbb{Z}), \beta^{su}_{1,1}\big)$, see Example 5B.12, are supercommutative algebras. ◇

A classical result is that for commutative algebras the notions of left, right and bimodules agree. Categorically this is also the case:

**Proposition 5J.4.** *Let* $\mathtt{A} \in \mathbf{C}$ *be a commutative algebra.*

*(i) Every* $\mathtt{M} \in \mathbf{Mod}_{\mathbf{C}}(\mathtt{A})$ *has the structure of a left* $\mathtt{A}$ *module.*

*(ii) Every* $\mathtt{N} \in (\mathtt{A})\mathbf{Mod}_{\mathbf{C}}$ *has the structure of a right* $\mathtt{A}$ *module.*

*(iii) We have equivalences of categories*

$$\mathbf{Mod}_{\mathbf{C}}(\mathtt{A}) \simeq (\mathtt{A})\mathbf{Mod}_{\mathbf{C}} \simeq (\mathtt{A})\mathbf{Mod}_{\mathbf{C}}(\mathtt{A}).$$

*Consequently,* $\mathbf{Mod}_{\mathbf{C}}(\mathtt{A})$ *is a monoidal category with* $\otimes = \otimes_{\mathtt{A}}$ *and* $\mathbb{1} = \mathtt{A}$*.*

*Proof. (i).* We can define a left action on $\mathtt{M}$ via

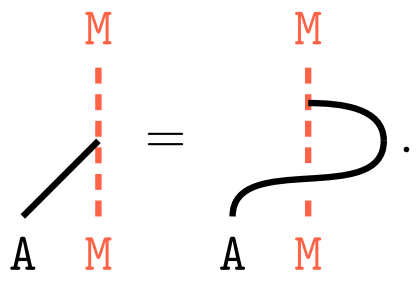

We now need to check associativity and unitality using the Reidemeister calculus as well a the diagrammatics for right actions:

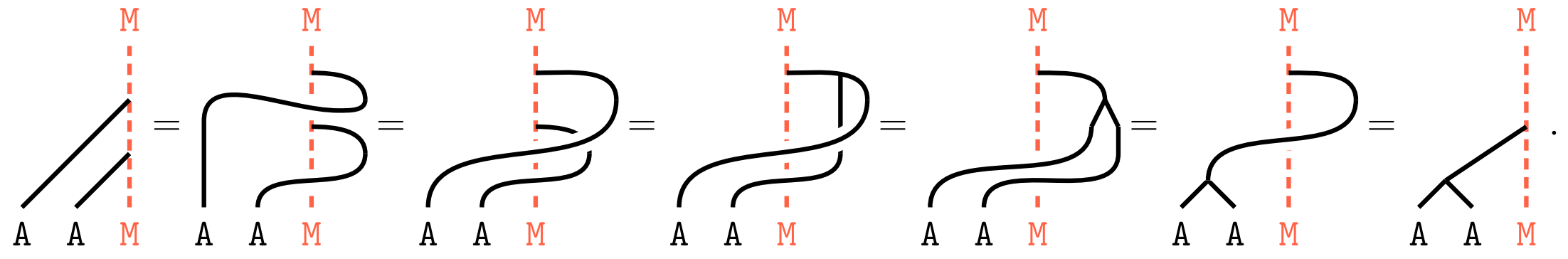

Note that the forth equality uses commutativity (5J-2). Similarly, we compute

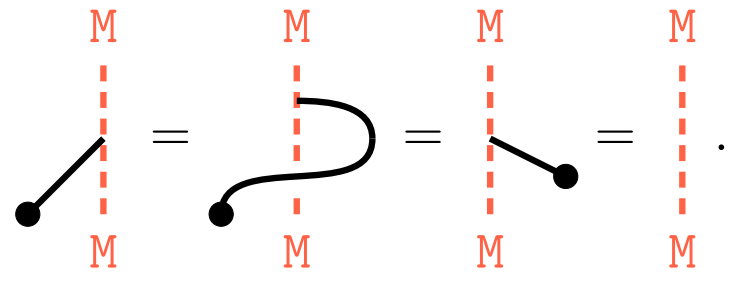

This shows (i).

*(ii).* By symmetry.

*(iii).* The following verifies that right and left action commute:

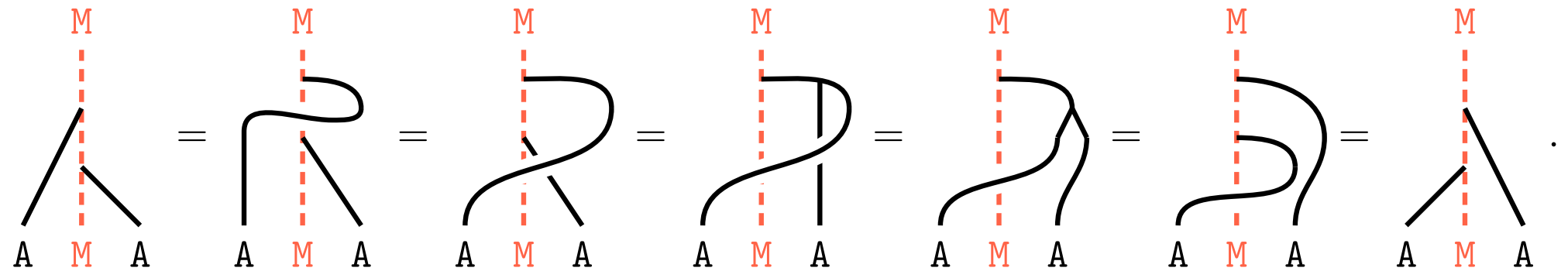

Thus, we can define a bimodule structure on $\mathtt{M} \in \mathbf{Mod}_{\mathbf{C}}(\mathtt{A})$ and also, by symmetry, on $\mathtt{N} \in (\mathtt{A})\mathbf{Mod}_{\mathbf{C}}$. Moreover, we can also match the equivariant morphisms, *e.g.*

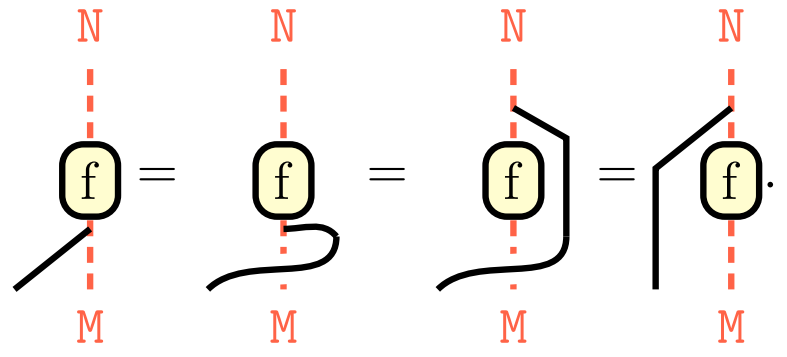

One then easily verifies that one gets the claimed equivalences of categories. ∎

We can also define the monoidal structure on $\mathbf{Mod}_{\mathbf{C}}(\mathtt{A})$ diagrammatically:

(5J-5)

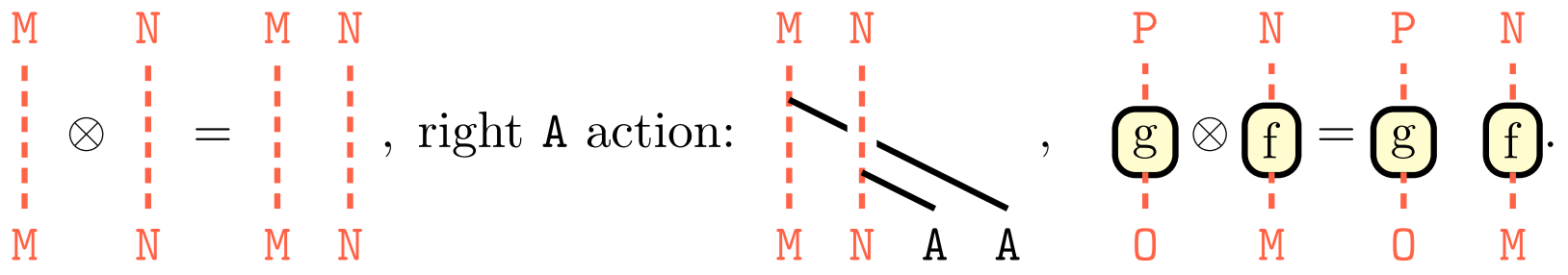

Thus, a good question would be whether $\mathbf{Mod}_{\mathbf{C}}(\mathtt{A})$ is also braided. This is not quote the case:

**Definition 5J.6.** For any commutative algebra $\mathtt{A} \in \mathbf{C}$ let $\mathbf{Mod}^{\beta}_{\mathbf{C}}(\mathtt{A}) \subset \mathbf{Mod}_{\mathbf{C}}(\mathtt{A})$ denote the full subcategory with objects satisfying

(5J-7)

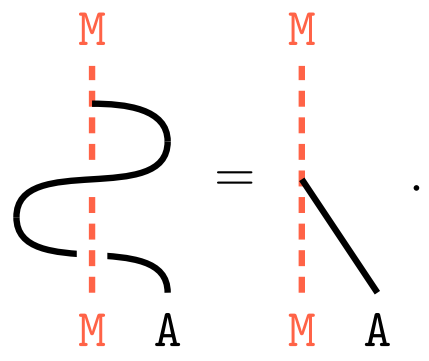

(Note that the equations needs $\mathtt{A}$ to be commutative.) ◇

The right $\mathtt{A}$ modules satisfying (5J-7) are also sometimes called ***braided*** for the following reason.

**Proposition 5J.8.** *For any commutative algebra $\mathtt{A} \in \mathbf{C}$ the category $\mathbf{Mod}^{\beta}_{\mathbf{C}}(\mathtt{A})$ is braided with braiding inherited from $\mathbf{C}$.*

*Proof.* This is Exercise 5M.5. □

**Example 5J.9.** If $\mathbf{C} \in \mathbf{BCat}$ is symmetric, then, using Reidemeister calculus, we see that (5J-7) becomes

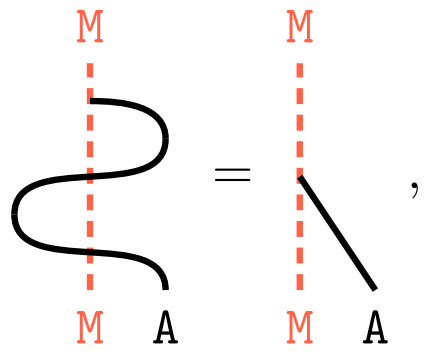

which holds for all objects, *i.e.* $\mathbf{Mod}^{\beta}_{\mathbf{C}}(\mathtt{A}) = \mathbf{Mod}_{\mathbf{C}}(\mathtt{A})$. Furthermore, we have $\mathbf{Mod}_{\mathbf{C}}(\mathtt{A}) \in \mathbf{BCat}$ is also symmetric. In particular, the module categories of any commutative algebra in $\mathbf{Vec}_{\Bbbk}$ are symmetric. This is a classic example. ◇

Another answers is that $\mathbf{Mod}_{\mathbf{C}}(\mathtt{A})$ is monoidal if $\mathtt{A}$ is a bialgebra. We will see that the same is true categorically. Since this is important let us be precise:

**Definition 5J.10.** A ***bialgebra*** $\mathtt{A} = (\mathtt{A}, \mathrm{m}, \mathrm{i}, \mathrm{d}, \mathrm{e})$ in a category $\mathbf{C} \in \mathbf{BCat}$ consist of

- an algebra $\mathtt{A} = (\mathtt{A}, \mathrm{m}, \mathrm{i}) \in \mathbf{C}$;
- a coalgebra $\mathtt{A} = (\mathtt{A}, \mathrm{d}, \mathrm{e}) \in \mathbf{C}$;

such that

(i) we have the unitality conditions

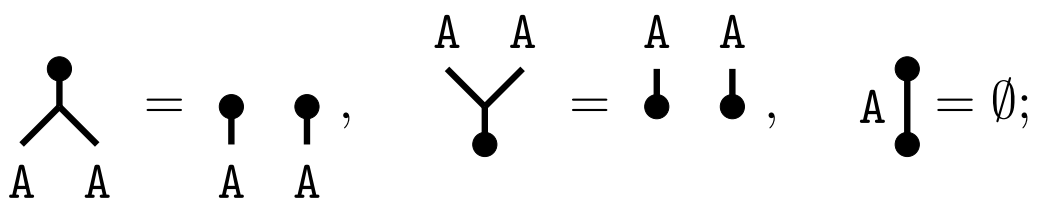

(ii) we have the compatibility condition

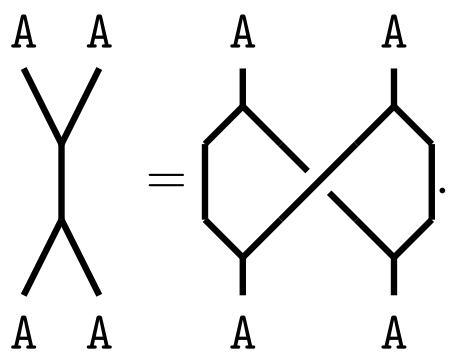

(As far as we know, the final equation does not have any name.) ◇

**Proposition 5J.11.** *For any bialgebra $\mathtt{A} \in \mathbf{C}$ the assignment as in (5J-5) but*

(5J-12) *right $\mathtt{A}$ action:*

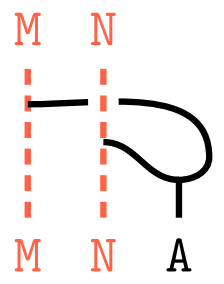

*defines a monoidal structure on* $\mathbf{Mod}_{\mathbf{C}}(\mathtt{A})$.

*Proof.* First note that (5J-12) is a well-defined right A action. For example,

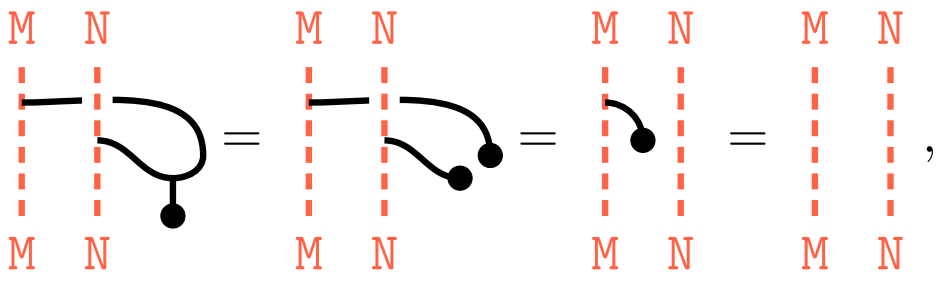

verifies unitality, where the first equality is using Definition 5J.10.(i). Moreover,

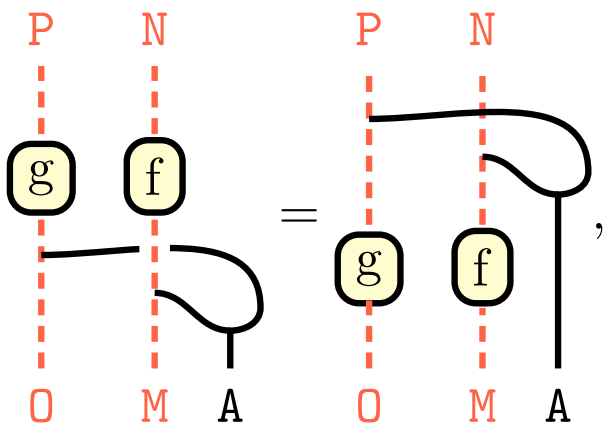

shows that A equivariant morphisms go to A equivariant morphisms under $\otimes$, showing that the assignment is well-defined, *i.e.* $\otimes$ stays within $\mathbf{Mod_C}(\mathtt{A})$. Verifying that this assembles into a monoidal structure works then exactly in the same way as for any horizontal juxtaposition. □

5K. **Hopf algebras in braided rigid categories.** To get our examples later on, we do some diagrammatics again.

**Definition 5K.1.** A ***pre Hopf algebra*** $\mathtt{A} = (\mathtt{A}, \mathrm{m}, \mathrm{i}, \mathrm{d}, \mathrm{e}, \mathrm{s})$ in a category $\mathbf{C} \in \mathbf{BCat}$ consist of

- a bialgebra $\mathtt{A} = (\mathtt{A}, \mathrm{m}, \mathrm{i}, \mathrm{d}, \mathrm{e}) \in \mathbf{C}$;
- an ***antipode*** $(\mathrm{s}\colon \mathtt{A} \to \mathtt{A}) \in \mathbf{C}$, illustrated by

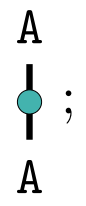

such that

(i) we have the antipode condition

(5K-2)

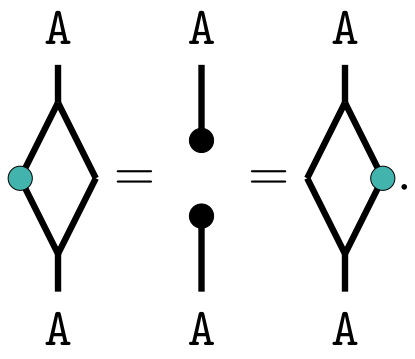

If s is invertible, then we call A a ***Hopf algebra***. ◇

Hopf algebras kind of generalize commutative algebras, *e.g.* compare (5J-2) to:

**Lemma 5K.3.** *For any pre Hopf algebra* $\mathtt{A} \in \mathbf{C}$ *for* $\mathbf{C} \in \mathbf{BCat}$ *we have* **sliding**, *i.e.*

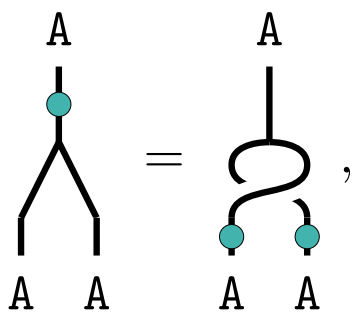

*including its horizontal mirror.*

*Proof.* Via the Reidemeister calculus and (5K-2), see *e.g.* [**Maj94**, Lemma 2.3]. □

Further, the following should be compared to Proposition 5J.4:

**Proposition 5K.4.** *Let* $\mathtt{A} \in \mathbf{C}$ *for* $\mathbf{C} \in \mathbf{BCat}$ *be a pre Hopf algebra.*

*(i) Every* $\mathtt{M} \in \mathbf{Mod_C}(\mathtt{A})$ *has the structure of a left* A *module.*

*(ii) Every* $\mathtt{N} \in (\mathtt{A})\mathbf{Mod_C}$ *has the structure of a right* A *module.*

*Proof.* Of course, (i) and (ii) are equivalent up to right-left symmetry, and it suffices to prove (i).

*(i).* We can define a left action on M via

(5K-5)

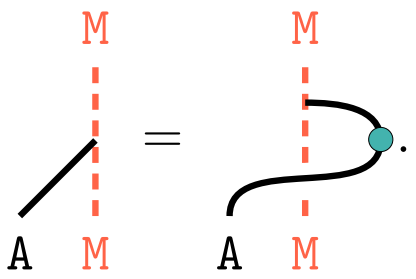

Using sliding Lemma 5K.3, we see that (5K-5) satisfies associativity:

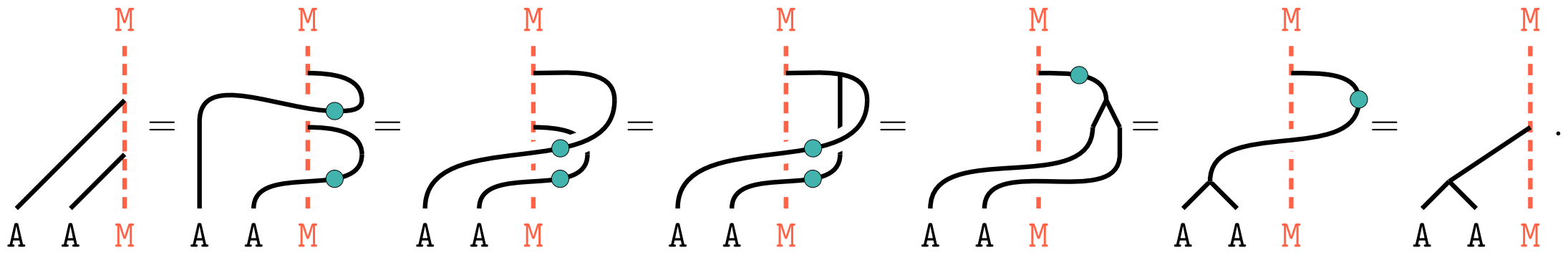

(Note that this is almost the same argument as in Proposition 5J.4.) The easier proof that (5K-5) also satisfies unitality is omitted. □

Let us say ***right rigid*** and ***left rigid*** in case only the right respectively left duals need to exist.

**Theorem 5K.6.** *For any pre Hopf algebra* A $\in$ **C** *where* **C** $\in$ **BRCat** *we have:*

*(i) The category* $\mathbf{Mod}_{\mathbf{C}}(\mathtt{A})$ *is right rigid with duality inherited from* **C**.

*(ii) The category* $(\mathtt{A})\mathbf{Mod}_{\mathbf{C}}$ *is left rigid with duality inherited from* **C**.

*(iii) If* A *is a Hopf algebra, then both,* $\mathbf{Mod}_{\mathbf{C}}(\mathtt{A})$ *and* $(\mathtt{A})\mathbf{Mod}_{\mathbf{C}}$*, are rigid with duality inherited from* **C**.

*Proof. (i).* First, Proposition 5J.11 shows that we get a monoidal structure, and by Proposition 4I.1 we then get that right duals, using this monoidal structure, have an action of A, but from the wrong side. However, by Proposition 5K.4 we can swap sides of the actions.

*(ii).* By (i) via symmetry.

*(iii).* To define a left action on the left dual $^{\star}$M we first denote the inverse of the antipode by

$s^{-1} \rightsquigarrow$ 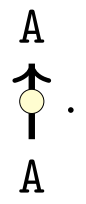

In order to define a right action on $^{\star}$M recall from (5K-5) that M has a left action. Thus, we can define a right action on $^{\star}$M via

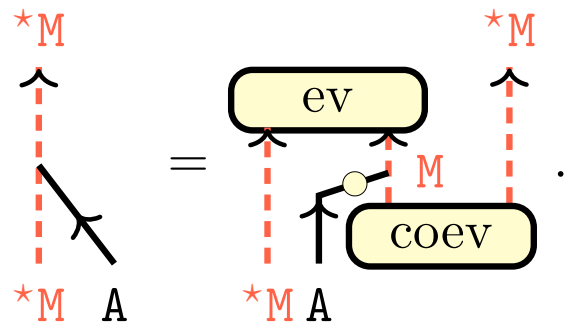

The invertibility comes into play since one needs *e.g.*

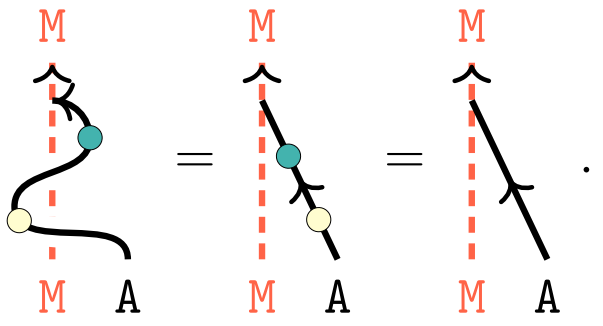

The second claim in (iii) follows again by symmetry. □

**Example 5K.7.** Again, this generalizes several notions:

(a) Hopf algebras in $\mathbf{Vec}_{\Bbbk}$ are classical Hopf algebras. A particular example is $\Bbbk[G]$.

(b) Hopf algebras in $\mathbf{Vec}_{\mathbb{C}}(\mathbb{Z}/2\mathbb{Z})$ with its super braiding could be called super Hopf algebras.

Thus, by (a) we recover the classical result that $\mathbf{fdMod}\big(\Bbbk[G]\big)$ is rigid. ◇

5L. **Summary of the interplay between topology and categorical algebra.** Some (the most important) Brauer categories we have seen are summarized in the following table.

| | monoidal | braided | pivotal | symmetric | self-dual • | Reidemeister 1 | topology |
|---|---|---|---|---|---|---|---|
| **Br** | **Y** | **Y** | **Y** | Y | Y | Y | **1Cob** |
| **qBr** | **Y** | **Y** | **Y** | N | Y | Y | **1Tan** |
| **oqBr** | **Y** | **Y** | **Y** | N | N | Y | **1State** |
| **orqBr** | **Y** | **Y** | **Y** | N | N | N | **1Ribbon** |

.

We leave it to the reader to fill in all the various versions using the adjectives "oriented", "quantum" and "ribbon". Let us use the placeholder _, which can be filled in with these adjective.

The point is that they are all equivalent to their topological incarnations while "free XYZ with properties ABC". Thus, we define:

A ***quantum invariant*** Q is a structure preserving functor

$$\mathrm{Q}\colon _\mathbf{Br} \to \mathbf{C},$$

where **C** is "a linear algebra like category".

The aim of the following lectures is to make precise what "a linear algebra like category", a.k.a. "a category where we can compute", might mean, the guiding example being $\mathbf{fdVec}_{\Bbbk}$.

5M. **Exercises.**

*Exercise* 5M.1. Prove that $\big(\mathbf{Vec}_{\mathbb{C}}(\mathbb{Z}/2\mathbb{Z}),\beta^{st}_{1,1}\big) \not\simeq_{\beta} \big(\mathbf{Vec}_{\mathbb{C}}(\mathbb{Z}/2\mathbb{Z}),\beta^{su}_{1,1}\big)$, where the standard and super braidings are defined in Example 5B.12. ◇

*Exercise* 5M.2. Verify the claims in Example 5C.8 and Remark 5C.10. ◇

*Exercise* 5M.3. For any $\mathbf{C} \in \mathbf{BCat}$ show that $\mathbf{C}^{op}, \mathbf{C}^{co}, \mathbf{C}^{coop} \in \mathbf{BCat}$ by defining braidings on them using the braiding of **C**. ◇

*Exercise* 5M.4. With respect to Remark 5F.5, write down all only implicitly stated relations for **qBr** and **oqBr**. What about *e.g.*

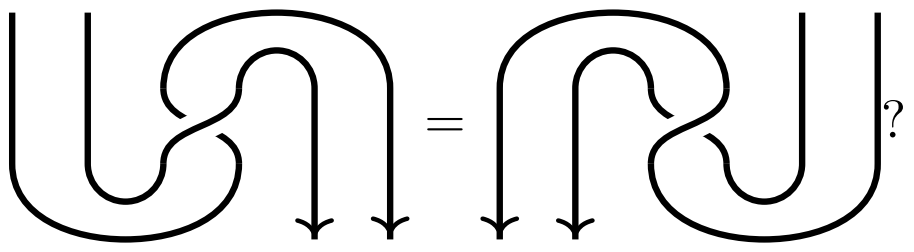

What about other equations that are topologically true (pick your favorite)? ◇

*Exercise* 5M.5. Prove Proposition 5J.8 and verify the missing claims in Proposition 5J.11. ◇

## 6. Additive, linear and abelian categories – definitions and examples

A topological invariant should be computable, *i.e.* within the realm of linear or homological algebra. So:

What is the analog of linear or homological algebra in a categorical language?

The author firmly believes in "Linear algebra might be the single most useful thing ever produced from mathematics.", *cf.* Figure 12, and generalizing linear algebra is the topic of this section ☺.

6A. **Conventions.** We keep all conventions from before, including all abbreviations that we used. Let us stress one:

*Convention* 6A.1. Similarly as in Convention 2C.3, we will write *e.g.* $\oplus$ instead of $\oplus_{\mathbf{C}}$. ◇

*Convention* 6A.2. We will see a lot of objects defined via some universal property. By a very general type of argument, which we will call a ***universality argument***, see Section 6D, such objects are unique up to unique isomorphism. Since these arguments are very similar, we usually omit the corresponding proofs. Moreover, these objects are usually objects together with extra data such as a morphism, but we tend to treat them as objects if no confusion can arise. ◇

*Convention* 6A.3. Recall from Section 3E that a $\Bbbk$ algebra is an algebra in $\mathbf{Vec}_{\Bbbk}$. Similarly, a ***ring*** for us is an algebra object in $\mathbf{Vec}_{\mathbb{Z}}$ (the category of abelian groups, see Example 1B.4), in particular, associative and unital. Moreover, throughout, $\mathbb{S}$ denotes a commutative ring, *i.e.* a commutative algebra in $\mathbf{Vec}_{\mathbb{Z}}$. ◇

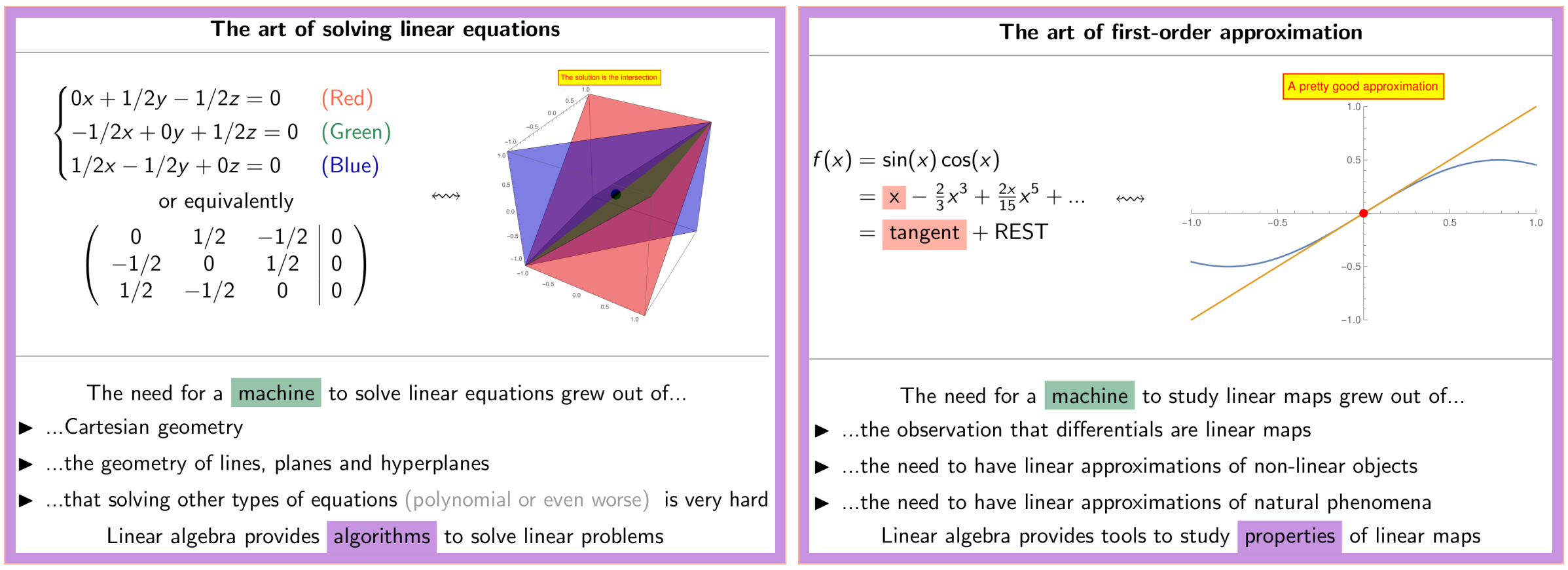

FIGURE 12. The essence of linear algebra.

Picture from the "What is...linear algebra?" playlist on [Tub21].

6B. **A motivating example.** As we have already seen, the "multiplication" $\otimes_{\mathbb{k}}$ of $\mathbf{Vec}_{\mathbb{k}}$ generalizes to the notion of monoidal categories. Let us now focus on the "addition" $\oplus$ of $\mathbf{Vec}_{\mathbb{k}}$.

- Firstly, we note that $\mathrm{Hom}_{\mathbf{Vec}_{\mathbb{k}}}(\mathtt{X},\mathtt{Y}) \in \mathbf{Vec}_{\mathbb{k}}$. Or in words, hom spaces between $\mathbb{k}$ vector spaces are, of course, $\mathbb{k}$ vector spaces again. In particular, they are abelian groups, meaning that we can add and subtract morphisms $f,g \in \mathrm{Hom}_{\mathbf{Vec}_{\mathbb{k}}}(\mathtt{X},\mathtt{Y})$ and the results $f \pm g$ are still in $\mathrm{Hom}_{\mathbf{Vec}_{\mathbb{k}}}(\mathtt{X},\mathtt{Y})$. There is also an additive unit, the zero map 0, additive inverses and composition is biadditive.

  *Remark* 6B.1. We can, of course, also use scalars from $\mathbb{k}$, and this property will below be called $\mathbb{k}$ ***linear***. However, for now only the condition of being an abelian group or, equivalently, $\mathbb{Z}$ linear is relevant for us. ◇

- Second, we have a ***zero object***, the zero $\mathbb{k}$ vector spaces $\mathtt{0}$, which satisfies:

  (6B-2) $$\text{For all } \mathtt{X} \in \mathbf{Vec}_{\mathbb{k}} \text{ there exist unique morphisms } 0\colon \mathtt{X} \to \mathtt{0},\ 0\colon \mathtt{0} \to \mathtt{X}.$$

  The morphisms in (6B-2) are called the ***zero morphisms*** and they are the zero maps.

- Finally, we consider the pair category (see Definition 1B.12) $\mathbf{Vec}_{\mathbb{k}} \times \mathbf{Vec}_{\mathbb{k}}$ and we have a bifunctor

  $$\oplus\colon \mathbf{Vec}_{\mathbb{k}} \times \mathbf{Vec}_{\mathbb{k}} \to \mathbf{Vec}_{\mathbb{k}}, \quad \oplus\big((\mathtt{X},\mathtt{Y})\big) = \mathtt{X} \oplus \mathtt{Y}, \oplus\big((f,g)\big) = f \oplus g,$$

  called the ***direct sum***, using again abbreviations of the form $\mathtt{X} \oplus \mathtt{Y}$ instead of $\oplus\big((\mathtt{X},\mathtt{Y})\big)$. We note that the object $\mathtt{X} \oplus \mathtt{Y}$ has a universal-type property, namely: there exist morphisms $i_{\mathtt{X}}, i_{\mathtt{Y}}, p_{\mathtt{X}}, p_{\mathtt{Y}} \in \mathbf{Vec}_{\mathbb{k}}$ such that

  (6B-3) $$\begin{array}{ccc} & \mathtt{X} \oplus \mathtt{Y} & \\ {\scriptstyle i_{\mathtt{X}}} \nearrow \swarrow {\scriptstyle p_{\mathtt{X}}} & & {\scriptstyle p_{\mathtt{Y}}} \searrow \nwarrow {\scriptstyle i_{\mathtt{Y}}} \\ \mathtt{X} & & \mathtt{Y} \end{array}, \quad p_{\mathtt{X}} i_{\mathtt{X}} = \mathrm{id}_{\mathtt{X}}, \quad p_{\mathtt{Y}} i_{\mathtt{Y}} = \mathrm{id}_{\mathtt{Y}}, \quad i_{\mathtt{X}} p_{\mathtt{X}} + i_{\mathtt{Y}} p_{\mathtt{Y}} = \mathrm{id}_{\mathtt{X} \oplus \mathtt{Y}}.$$

  The two morphisms $i_{\mathtt{X}}$ and $i_{\mathtt{Y}}$ are called ***inclusions***, the other two $p_{\mathtt{X}}$ and $p_{\mathtt{Y}}$ ***projections*** (of $\mathtt{X}$ and $\mathtt{Y}$, respectively).

6C. **An even more down to earth motivating example.** Let us be completely explicit and consider the situation of $\mathbf{Mat}_{\mathbb{k}}$. In this case the three observations above take the following form.

- Matrices can be added and this is bilinear with respect to multiplication $\circ$, *e.g.*

  $$\left(s \cdot \begin{pmatrix} a & b \\ c & d \end{pmatrix}\right) \circ \left(t \cdot \begin{pmatrix} e & f \\ g & h \end{pmatrix}\right) = st \cdot \left(\begin{pmatrix} a & b \\ c & d \end{pmatrix} \circ \begin{pmatrix} e & f \\ g & h \end{pmatrix}\right).$$

- There exists a zero $\mathtt{0}$ and a zero matrix $0 = (\,0\,)$.

- We can add numbers and there exist block matrices and corresponding inclusions and projections of blocks, *e.g.*

  $$\mathtt{2} \oplus \mathtt{2} = \mathtt{4}, \quad i_2^{\leftarrow} = \begin{pmatrix} 1 & 0 \\ 0 & 1 \\ 0 & 0 \\ 0 & 0 \end{pmatrix}, i_2^{\rightarrow} = \begin{pmatrix} 0 & 0 \\ 0 & 0 \\ 1 & 0 \\ 0 & 1 \end{pmatrix}, \quad p_2^{\leftarrow} = \begin{pmatrix} 1 & 0 & 0 & 0 \\ 0 & 1 & 0 & 0 \end{pmatrix}, p_2^{\rightarrow} = \begin{pmatrix} 0 & 0 & 1 & 0 \\ 0 & 0 & 0 & 1 \end{pmatrix}.$$

6D. **A brief reminder on universality.** Let $F \in \mathbf{Hom}(\mathbf{C},\mathbf{D})$ be a functor between the categories $\mathbf{C},\mathbf{D} \in \mathbf{Cat}$. Then a pair $\big(X, f\colon Y \to F(X)\big)$ in $\mathbf{D}$ satisfies a ***universal property for*** $Y$ ***and*** $F$ if for any $g\colon Y \to F(Z)$ there exist a unique $u\colon X \to Z$ making commutative. A ***universal property for*** $F$ ***and*** $Y$ is defined similarly, using the opposite categories $\mathbf{C}^{op}$ and $\mathbf{D}^{op}$ instead of $\mathbf{C}$ and $\mathbf{D}$.

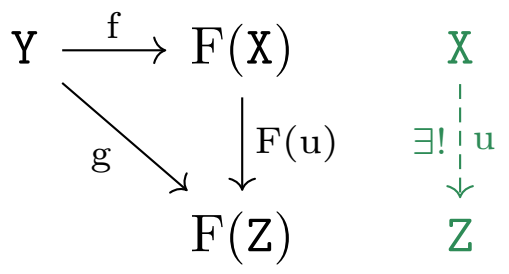

**Example 6D.1.** Let $\mathbf{D} = \mathbf{Set}$. Then the product $X_1 \times X_2$ comes with the two coordinate projections $p_1$ and $p_2$ and satisfies a universal property. This universal property will take place in $\mathbf{Set}\times\mathbf{Set}$ (the pair category of Definition 1B.12) can be formulated as follows. First, let $F\colon \mathbf{Set} \to \mathbf{Set}\times\mathbf{Set}$ be the diagonal functor. Then the pair $\big(X_1 \times X_2, (p_1, p_2)\big)$ satisfies a universal property from $F$ to $(X_1, X_2)$, *i.e.* we have

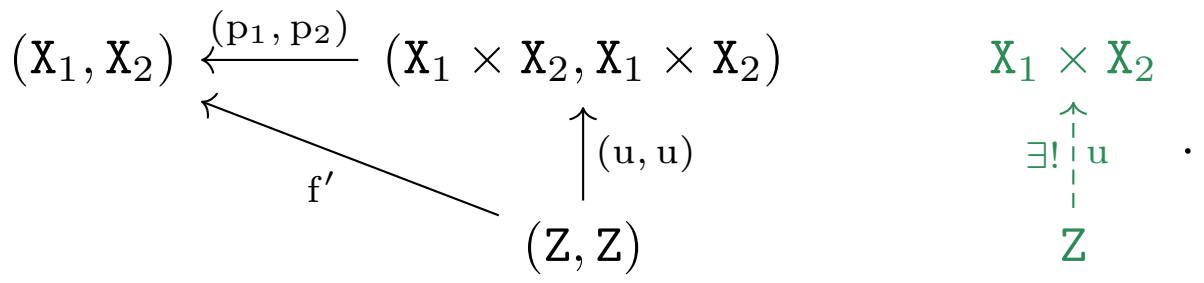

Here $X = X_1 \times X_2$ and $Y = (X_1, X_2)$. ◇

**Lemma 6D.2.** *A pair* $(X, f)$ *satisfying a universal property for* $Y$ *and* $F$*, if it exists, is unique up to unique isomorphism, i.e. if* $(X', f')$ *is another pair, then there exists a unique isomorphism* $h\colon X \xrightarrow{\cong} X'$ *such that* $f' = F(h)f$*. Similarly for pairs* $F$ *and* $Y$.

*Proof.* It follows by substituting $(X', f')$ into the definition for $(X, f)$ that $h$ exists uniquely, *i.e.*

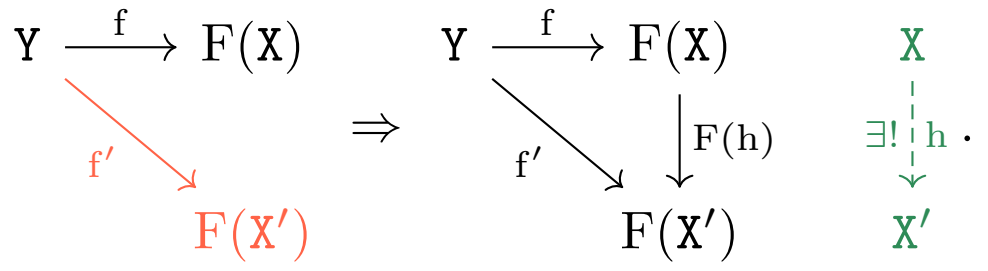

By symmetry, we also get a unique $h'\colon X' \to X$, which we paste together with $h$:

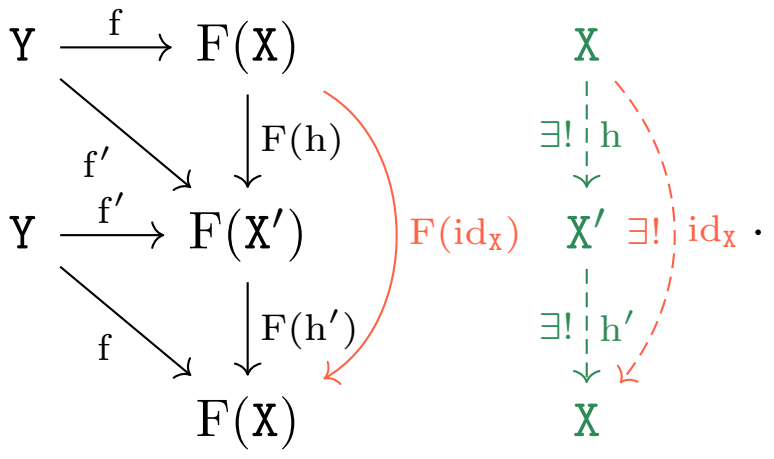

However, there already exists a morphisms $X \to X$ making the outer part above commutative, namely $\mathrm{id}_X$. Thus, $h'h = \mathrm{id}_X$ and, by symmetry, $hh' = \mathrm{id}_{X'}$. □

Notions defined by universal properties are unique if they exist, but do not have to exist in general. However, recall from Convention 3A.4 that we do not write "assuming that XYZ exists" below, and this will be an implicit assumption for *e.g.* statements such as those in Lemma 6G.6 to make sense.

6E. **Linear algebra in categories.** We generalize the situations of $\mathbf{Vec}_{\mathbb{k}}$ and $\mathbf{Mat}_{\mathbb{k}}$.

**Definition 6E.1.** A category $\mathbf{C} \in \mathbf{Cat}$ is called $\mathbb{S}$ ***linear*** if the space $\mathrm{Hom}_{\mathbf{C}}(X, Y)$ is an $\mathbb{S}$ module for all $X, Y \in \mathbf{Cat}$ and composition is $\mathbb{S}$ bilinear. ◇

**Definition 6E.2.** A category $\mathbf{C} \in \mathbf{Cat}$ is called ***additive*** if

- it is $\mathbb{Z}$ linear;
- there exists a ***zero object*** $0 \in \mathbf{Cat}$ (meaning an object satisfying (6B-2));
- for all $X, Y \in \mathbf{Cat}$ there exists an object $X \oplus Y$ called ***direct sum***, which satisfies the universal property in (6B-3).

A category satisfying only the first bullet point is called ***preadditive***. ◇

The following slogan should be compared to Equation 6F-5:

**Slogan.** $\mathbb{S}$ linear and additive categories generalize linear algebra.

**Example 6E.3.** The properties of being $\mathbb{S}$ linear and additive are parallel to each other: one asks for linearity of hom spaces, the other for existence of direct sums.

(a) The category **Set** is neither $\mathbb{S}$ linear nor additive.

(b) The categories $\mathbf{Vec}_{\Bbbk}$ and $\mathbf{fdVec}_{\Bbbk}$ are $\Bbbk$ linear additive.

(c) The categories of the form $\mathbf{Vec}_{\Bbbk}^{\omega}(\mathrm{G})$ are $\Bbbk$ linear, but not additive.

(d) The category $\mathbf{Vec}_{\mathbb{Z}}$ is additive, but not $\Bbbk$ linear (but, of course, $\mathbb{Z}$ linear).

Here, and throughout, we write $\Bbbk$ linear additive instead of $\Bbbk$ linear and additive. ◇

**Example 6E.4.** Additive categories need to be closed under $\oplus$. For example, the full subcategory $\mathbf{evenVec}_{\Bbbk} \subset \mathbf{fdVec}_{\Bbbk}$ of ***even dimensional*** $\Bbbk$ ***vector spaces*** is $\Bbbk$ linear additive, while the full subcategory $\mathbf{oddVec}_{\Bbbk} \subset \mathbf{fdVec}_{\Bbbk}$ of ***odd dimensional*** $\Bbbk$ ***vector spaces*** is only $\Bbbk$ linear. ◇

The following is (almost) immediate.

**Lemma 6E.5.** *Let* $\mathbf{C} \in \mathbf{Cat}$ *be additive. Then there exists a bifunctor* $\oplus\colon \mathbf{Cat}\times\mathbf{Cat} \to \mathbf{Cat}$ *called* **direct sum**. □

We can again say "the" direct sum, justified by the following lemma whose proof is a universality argument:

**Lemma 6E.6.** *Up to unique isomorphism,* $\mathtt{X} \oplus \mathtt{Y}$ *is the only object in* $\mathbf{C}$ *satisfying* (6B-3).

**Definition 6E.7.** An $\mathbb{S}$ ***linear functor*** $\mathrm{F} \in \mathbf{Hom}(\mathbf{C},\mathbf{D})$ between $\mathbb{S}$ linear categories is a functor such that the induced map

$$\mathbf{Hom}_{\mathtt{C}}(\mathtt{X},\mathtt{Y}) \to \mathbf{Hom}_{\mathtt{D}}\big(\mathrm{F}(\mathtt{X}),\mathrm{F}(\mathtt{Y})\big)$$

is $\mathbb{S}$ linear for all $\mathtt{X},\mathtt{Y} \in \mathtt{C}$. ◇

At first glance Definition 6E.7 looks like the "wrong definition" for additive categories, since it does not involve direct sums. However, the slogan is "linear implies additive":

**Lemma 6E.8.** *Let* $\mathrm{F} \in \mathbf{Hom}(\mathbf{C},\mathbf{D})$ *be* $\mathbb{Z}$ *linear, and let* $\mathbf{C}$ *and* $\mathbf{D}$ *be additive. Then there exists a natural isomorphism* $\mathrm{F}(\mathtt{X}\oplus\mathtt{Y}) \cong \mathrm{F}(\mathtt{X}) \oplus \mathrm{F}(\mathtt{Y})$.

*Proof.* Note that being $\mathbb{Z}$ linear gives us the equality

$$\mathrm{F}(\mathrm{f}+\mathrm{g}) = \mathrm{F}(\mathrm{f}) + \mathrm{F}(\mathrm{g}).$$

Thus, all the equations in (6B-3) are preserved by F, which implies that $\mathrm{F}(\mathtt{X}\oplus\mathtt{Y})$ is a direct sum of $\mathrm{F}(\mathtt{X})$ and $\mathrm{F}(\mathtt{Y})$, and the claim follows from Lemma 6E.6. □

The following is as usual not hard to see:

**Lemma 6E.9.** *The identity functor on an* $\mathbb{S}$ *linear category is* $\mathbb{S}$ *linear. Moreover, if* F *and* G *are* $\mathbb{S}$ *linear functors, then so is* GF.

**Example 6E.10.** By Lemma 6E.8, we get the ***categories of*** $\mathbb{S}$ ***linear categories*** $\mathbf{Cat}_{\mathbb{S}}$ and the ***category of additive categories*** $\mathbf{Cat}_{\oplus}$ at the same time, morphism being the appropriate linear functors. Also important is the category $\mathbf{Cat}_{\mathbb{S}\oplus}$ ***of*** $\mathbb{S}$ ***linear additive categories***, the morphisms being $\mathbb{S}$ linear functors, as well as the corresponding functor categories $\mathbf{Hom}_{\mathbb{S}}(\mathbf{C},\mathbf{D})$, $\mathbf{Hom}_{\oplus}(\mathbf{C},\mathbf{D})$ and $\mathbf{Hom}_{\mathbb{S}\oplus}(\mathbf{C},\mathbf{D})$. ◇

**Definition 6E.11.** $\mathbf{C},\mathbf{D} \in \mathbf{Cat}_{\mathbb{S}\oplus}$ are called ***equivalent as*** $\mathbb{S}$ ***linear additive categories***, denoted by $\mathbf{C} \simeq_{\mathbb{S}\oplus} \mathbf{D}$, if there exists an equivalence $\mathrm{F} \in \mathbf{Hom}_{\mathbb{S}\oplus}(\mathbf{C},\mathbf{D})$. Similarly in the $\mathbb{S}$ linear and additive setup, using the appropriate linear functors. ◇

6F. **The linear extension and the additive closure.** The following constructions allow us to perform linear algebra in almost all categories.

**Definition 6F.1.** The $\mathbb{S}$ ***linear extension*** of $\mathbf{C} \in \mathbf{Cat}$, denoted by $\mathbf{C}_{\mathbb{S}}$, is the category with $\mathrm{Ob}(\mathbf{C}_{\mathbb{S}}) = \mathrm{Ob}(\mathbf{C})$ and

$$\mathrm{Hom}_{\mathbf{C}_{\mathbb{S}}}(\mathtt{X},\mathtt{Y}) = \mathbb{S}\big\{\mathrm{Hom}_{\mathbf{C}}(\mathtt{X},\mathtt{Y})\big\},$$

and the composition being the evident $\mathbb{S}$ linear extension of the composition in $\mathbf{C}$. ◇

In other words, the hom spaces of $\mathbf{C}_{\mathbb{S}}$ are the free $\mathbb{S}$ modules with the basis set being the corresponding hom space in $\mathbf{C}$. In particular, if $\mathbb{S}$ is a field, then the hom spaces of $\mathbf{C}_{\mathbb{S}}$ are free $\mathbb{S}$ vector spaces.

**Example 6F.2.** To match our previous conventions, let $\mathbb{S} = \Bbbk$.

(a) We have $\mathbf{Set}_{\Bbbk} \simeq_{\Bbbk} \mathbf{Vec}_{\Bbbk}$.

(b) We have $\mathbf{Vec}(\mathrm{G})_{\Bbbk} \simeq_{\Bbbk} \mathbf{Vec}_{\Bbbk}(\mathrm{G})$.

Note that the notation, using the ground ring as a subscript for the $\mathbb{S}$ linear extension, is a bit annoying. For example, $(\mathbf{Vec}_{\Bbbk})_{\Bbbk}$ is not the same as $\mathbf{Vec}_{\Bbbk}$. ◇

**Example 6F.3.** For diagrammatic categories such as the Brauer category **Br**, taking the $\mathbb{S}$ linear extension amounts to taking formal sums of pictures with the same endpoints, *e.g.*

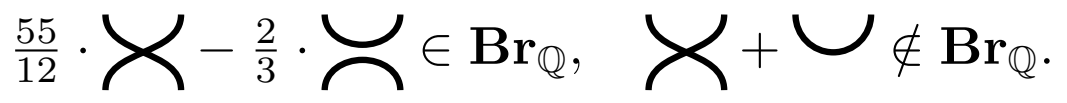

Since each diagram is a basis element, by definition, simplification of scalars is only allowed if the diagrams are the same. Moreover, composition is bilinear, meaning *e.g.*

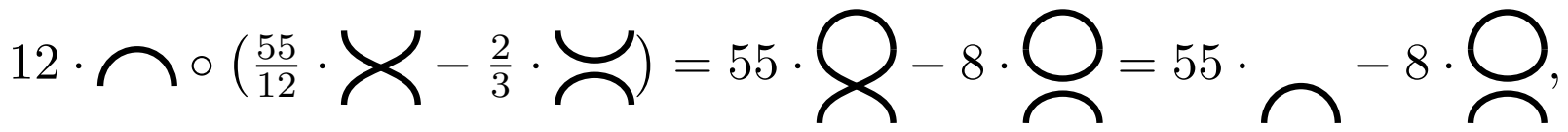

where we have used one of the Brauer relations (3D-8) in the last step. ◇

Recall our notation for free monoids from Section 3B.

**Definition 6F.4.** The ***additive closure*** of $\mathbf{C} \in \mathbf{Cat}_{\mathbb{S}}$, denoted by $\mathbf{C}_{\oplus}$, is the category with $\mathrm{Ob}(\mathbf{C}_{\oplus}) = \langle \mathrm{Ob}(\mathbf{C}) \mid \emptyset \rangle$ (the composition is written $\oplus$) and

$$\mathrm{Hom}_{\mathbf{C}_{\oplus}}(\mathtt{X}_1 \oplus \dots \oplus \mathtt{X}_c, \mathtt{Y}_1 \oplus \dots \oplus \mathtt{Y}_r) = \big\{ (\mathrm{f}_{ij})_{i=1,\dots,c}^{j=1,\dots,r} \mid \mathrm{f}_{ij} \in \mathrm{Hom}_{\mathbf{C}}(\mathtt{X}_j, \mathtt{Y}_i) \big\},$$

and the composition being matrix multiplication. ◇

In words, the objects of $\mathbf{C}_{\oplus}$ are formal (finite) direct sums of objects of $\mathbf{C}$ and the hom spaces of $\mathbf{C}_{\oplus}$ are matrices whose entries are morphisms of $\mathbf{C}$:

$$(\mathrm{f}_{ij})_{i=1,\dots,c}^{j=1,\dots,r} \leftrightsquigarrow \begin{pmatrix} \mathrm{f}_{11} & \dots & \mathrm{f}_{1c} \\ \vdots & \ddots & \vdots \\ \mathrm{f}_{r1} & \dots & \mathrm{f}_{rc} \end{pmatrix}. \tag{6F-5}$$

In particular, the following is clear, saying that the additive closure is a closure:

**Lemma 6F.6.** *For all* $\mathbf{C} \in \mathbf{Cat}_{\mathbb{S}}$ *we have* $\mathbf{C}_{\oplus} \simeq_{\oplus} (\mathbf{C}_{\oplus})_{\oplus}$. □

In contrast, the $\mathbb{S}$ linear extension is not a closure operation; see the end of Example 6F.2.

**Example 6F.7.** With respect to $\mathbf{Vec}(\mathrm{G})$ we have interesting examples:

(a) In $\mathbf{Vec}_{\Bbbk\oplus}(\mathrm{G})$, which is defined as first taking the $\Bbbk$ linear extension and then the additive closure of $\mathbf{Vec}(\mathrm{G})$, objects are formal direct sums of group elements, while morphisms are honest matrices with a restriction on entries coming from $\mathrm{Hom}_{\mathbf{Vec}_{\Bbbk}(\mathrm{G})}(\mathtt{i}, \mathtt{j}) = 0$ if $\mathtt{i} \neq \mathtt{j}$.

(b) The category $\mathbf{Vec}_{\Bbbk\oplus}(\mathrm{G})$ is called the ***category of*** G***-graded*** $\Bbbk$ ***vector spaces***. Important special cases are:
- For $\mathrm{G} = 1$ we have $\mathbf{Vec}_{\Bbbk\oplus}(1) \simeq \mathbf{fdVec}_{\Bbbk}$. In particular, Theorem 5E.3 implies that, after choosing the monoidal structure, $\mathbf{fdVec}_{\Bbbk}$ has only one structure of a braided category.
- For $\mathrm{G} = \mathbb{Z}/2\mathbb{Z}$ another common name is the ***category of (finite dimensional) super vector spaces***. This category has two braidings by Theorem 5E.3, the non-symmetric one is called the ***super braiding***.

On a side note, the term "super" is derived from Latin and means "beyond". This reflects the idea that super vector spaces go beyond standard vector spaces by including $\mathbb{Z}/2\mathbb{Z}$, which allows for the exploration of concepts in supersymmetry and other areas of mathematics and theoretical physics. ◇

As in Example 6F.7 the notation $\mathbb{S}\oplus$ means taking first the $\mathbb{S}$ linear extension, and then the additive closure.

**Example 6F.8.** For diagrammatic categories we get a diagram calculus of matrices. For example

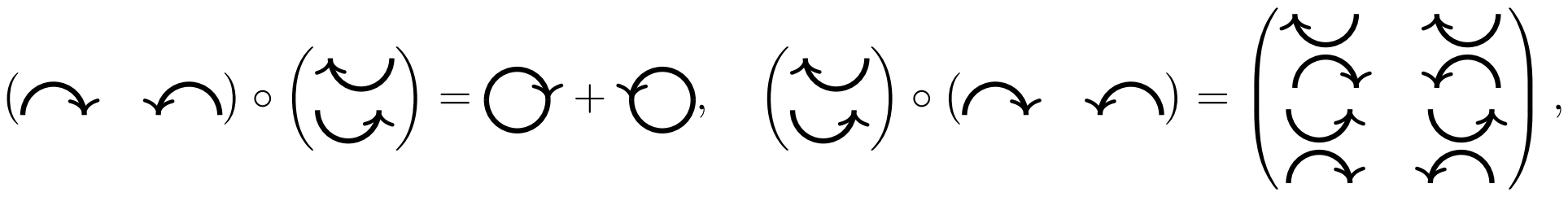

are computations in $\mathbf{oBr}_{\mathbb{S}\oplus}$. ◇

**Proposition 6F.9.** *We have the following.*

*(i) For any* $\mathbf{C} \in \mathbf{Cat}$ *we have* $\mathbf{C}_{\mathbb{S}} \in \mathbf{Cat}_{\mathbb{S}}$.

*(ii) There is a dense and faithful functor* $\mathrm{L}\colon \mathbf{C} \hookrightarrow \mathbf{C}_{\mathbb{S}}$ *given by* $\mathbb{S}$ *linearization.*

*(iii) If* $\mathbf{C} \in \mathbf{Cat}$ *is monoidal (or rigid or pivotal or braided etc.), then so is* $\mathbf{C}_{\mathbb{S}}$ *with its structure induced from* $\mathbf{C}$.

*Similarly for* $\mathbf{C}_{\oplus}$*, except that the corresponding functor in (ii) is fully faithful, but not dense.*

*Proof.* This is Exercise 6L.2. □

Here is the analog of Proposition 6F.10.

**Proposition 6F.10.** *Let* $\mathrm{F} \in \mathbf{Hom}(\mathbf{C}, \mathbf{D})$. *Then there exists a unique* $\mathrm{F}_{\mathbb{S}} \in \mathbf{Hom}_{\mathbb{S}}(\mathbf{C}_{\mathbb{S}}, \mathbf{D}_{\mathbb{S}})$ *such that we have a commuting diagram*

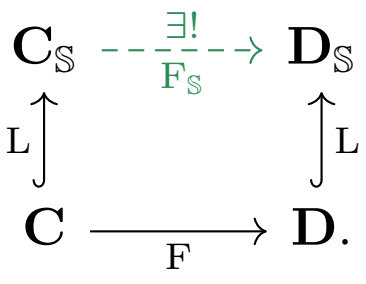

*Similarly for additive closures.*

*Proof.* The functor $\mathrm{F}_{\mathbb{S}}$ is the $\mathbb{S}$ linear extension of F. □

The following is not an honest theorem (clearly, the "we care about" is not well-defined ☺), but we state it nonetheless:

**Theorem 6F.11.** *All properties we care about behave nicely with* $\mathbb{S}$ *linear extensions and additive closures, e.g. if* F *is monoidal, then so is* $\mathrm{F}_{\mathbb{S}}$.

*Proof.* This is an experimental fact, taking Proposition 6F.9 and Proposition 6F.10 together. □

6G. **The first steps towards homological algebra in categories.** First, the generalization of a kernel:

**Definition 6G.1.** For a category $\mathbf{C} \in \mathbf{Cat}_{\oplus}$ and $\mathrm{f} \in \mathbf{C}$ we say $\mathtt{Ker}(\mathrm{f}) = \big(\mathtt{Ker}(\mathrm{f}), \mathrm{k} \colon \mathtt{Ker}(\mathrm{f}) \to \mathtt{X}\big)$ is a ***kernel of*** f if it has the universal property of the form

(6G-2)

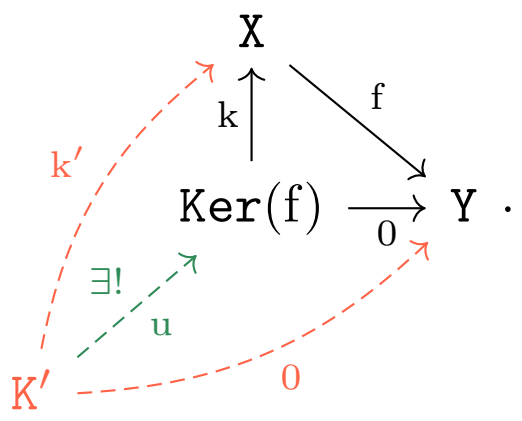

A ***cokernel of*** f, denoted by $\mathtt{Coker}(\mathrm{f}) = \big(\mathtt{Coker}(\mathrm{f}), \mathrm{c}\big)$, is a kernel of f in $\mathbf{C}^{op}$. ◇

Universality gives:

**Lemma 6G.3.** *Up to unique isomorphisms,* $\mathtt{Ker}(\mathrm{f})$ *is the only object in* $\mathbf{C}$ *satisfying* (6G-2). *Similarly for the cokernel.* □

**Example 6G.4.** As usual with universal objects, they might not exist:

(a) In $\mathbf{Vec}_{\Bbbk}$ and $\mathbf{fdVec}_{\Bbbk}$ (co)kernels exist and are the usual (co)kernels.

(b) The category $\mathbf{evenVec}_{\Bbbk}$ does not have (co)kernels in general, as the morphisms in $\mathbf{evenVec}_{\Bbbk}$ might be of odd rank.

(c) Diagrammatic categories such as $\mathbf{Br}_{\Bbbk\oplus}$ generally do not have (co)kernels.

We elaborate on (a) in Example 6G.5. ◇

The usual convention to identify kernels with their objects $\mathtt{Ker}(\mathrm{f})$ is a bit misleading, particularly in skeletal categories:

**Example 6G.5.** In $\mathbf{Mat}_{\mathbb{Q}}$ (co)kernels exist and can be described as follows. Take *e.g.* the matrix $\mathrm{f} = \left(\begin{smallmatrix} 1 & 2 \\ 2 & 4 \end{smallmatrix}\right)$. Then, for any $a \in \mathbb{Q}^{*}$, the pairs

$$\mathtt{Ker}(\mathrm{f}) = \Big(\mathtt{1}, \left(\begin{smallmatrix} 2a \\ -a \end{smallmatrix}\right)\Big), \quad \mathtt{Coker}(\mathrm{f}) = \Big(\mathtt{1}, \left(\begin{smallmatrix} 2a & -a \end{smallmatrix}\right)\Big),$$

are kernels and cokernels of f, respectively. The unique isomorphism u in (6G-2) for different scalars is the corresponding scaling map.

In particular, the interesting information is the vector $(2a, -a)$, and indeed:

$$\left(\begin{smallmatrix} 1 & 2 \\ 2 & 4 \end{smallmatrix}\right)\left(\begin{smallmatrix} 2a \\ -a \end{smallmatrix}\right) = \left(\begin{smallmatrix} 0 \\ 0 \end{smallmatrix}\right), \quad \left(\begin{smallmatrix} 2a & -a \end{smallmatrix}\right)\left(\begin{smallmatrix} 1 & 2 \\ 2 & 4 \end{smallmatrix}\right) = \left(\begin{smallmatrix} 0 & 0 \end{smallmatrix}\right),$$

so these vectors are what is non-categorically called kernel and cokernel. ◇

The following should remind the reader of a classical fact from linear algebra and gives us a good way to describe monic and epic morphisms, isomorphisms, subobjects, quotient objects *etc.* (Recall that these notions as defined in Section 1H.)

**Lemma 6G.6.** *Let* $\mathrm{f} \in \mathbf{C}$ *with* $\mathbf{C} \in \mathbf{Cat}_{\oplus}$.

(i) *We have* $\mathtt{Ker}\big(\mathtt{Ker}(f)\big) = 0$ *and* $\mathtt{Coker}\big(\mathtt{Coker}(f)\big) = 0$.

(ii) $\mathtt{Ker}(f) = 0$, *respectively* $\mathtt{Coker}(f) = 0$, *if and only if* $f$ *is monic, respectively epic.*

(iii) *The morphism* $k$ *of* $\mathtt{Ker}(f)$ *is monic, and the morphism* $c$ *of* $\mathtt{Coker}(f)$ *is epic.*

(iv) $\mathtt{Ker}(f) = 0 = \mathtt{Coker}(f) = 0$ *if and only if* $f$ *is monic and epic if and only if* $f$ *is an isomorphism.*

(v) *If* $f$ *is monic, then* $\mathtt{Y}/\mathtt{X} = \mathtt{Coker}(f)$ *is a quotient object of* $\mathtt{Y}$.

*Proof.* By symmetry, it suffices to prove the claims for the kernels.

(i). We write

$$\begin{array}{ccc} \mathtt{Ker}(f) & & \\ {\scriptstyle 0}\uparrow & \searrow^{k} & \\ 0 & \xrightarrow[0]{} & \mathtt{X} \end{array}$$

and observe that the zero object clearly satisfies the universal property of a kernel.

(ii). For $fh = fi$, we calculate $fh - fi = f(h - i) = 0$. Hence, letting $k' = h - i$, we see that $h - i = 0$ by (6G.1). Conversely, if $f$ is monic, then $fg = 0$ gives $g = 0$, which implies that the zero object is the kernel of $f$.

(iii). By combining (i) and (ii).

(iv). We already know by (ii) that the first two statements are equivalent, and one direction of the last statement is always true; see Lemma 1H.2. So, suppose that $f$ monic and epic. Then $f = \mathtt{Coker}\big(\mathtt{Ker}(f)\big) = \mathtt{Coker}(0)$, the latter being always an isomorphism.

(v). By the definition, we have a morphism $c\colon \mathtt{Y} \to \mathtt{Y}/\mathtt{X}$, which is epic by (iii). □

Note that we use the notation $\mathtt{Y}/\mathtt{X}$ in Lemma 6G.6.(v) for quotient objects, since there is a natural choice for epic morphisms in this case. The same notation will be used below, and we will also write $\mathtt{Y} \subset \mathtt{X}$ for subobjects, partially justified by Theorem 6H.6.

**Definition 6G.7.** An ***epic-monic factorization*** $(f, m, e)$ of $f \in \mathbf{C}$ with $\mathbf{C} \in \mathbf{Cat}_{\oplus}$ consists of

- a kernel $\big(\mathtt{ker}(f), k\big)$ and a cokernel $\big(\mathtt{coker}(f), c\big)$ for $f$;
- a kernel for $c$ and a cokernel for $k$;
- an object $\mathtt{I}$ and two morphisms $e\colon \mathtt{X} \to \mathtt{I}$ and $m\colon \mathtt{I} \to \mathtt{Y}$;

such that

(i) $f = me$;

(ii) we have $(\mathtt{I}, e) = \mathtt{coker}(k)$ and $(\mathtt{I}, m) = \mathtt{ker}(c)$ giving a sequence

(6G-8)
$$\mathtt{ker}(f) \xrightarrow{k} \underbrace{\mathtt{X} \xrightarrow{e} \mathtt{I} \xrightarrow{m} \mathtt{Y}}_{f} \xrightarrow{c} \mathtt{coker}(f)\ .$$

Note that it is not assumed that $e$ is epic and $m$ in monic. ◇

**Example 6G.9.** Let us spell out the epic-monic factorization in $\mathbf{Mat}_{\mathbb{Q}}$. Take the matrix

$$f = \begin{pmatrix} 1 & 2 & 3 \\ 4 & 5 & 6 \\ 7 & 8 & 9 \end{pmatrix}.$$

The rank of $f$ is 2, and we then set $\mathtt{I} = 2$. Now, we perform row reduction (rr) and get

$$\begin{pmatrix} 1 & 2 & 3 \\ 4 & 5 & 6 \\ 7 & 8 & 9 \end{pmatrix} \sim_{\mathrm{rr}} \begin{pmatrix} 1 & 0 & -1 \\ 0 & 1 & 2 \\ 0 & 0 & 0 \end{pmatrix} \xRightarrow[\text{rows}]{\text{remove zero}} e = \begin{pmatrix} 1 & 0 & -1 \\ 0 & 1 & 2 \end{pmatrix},$$

(this is the ***row reduced echelon form***) and the epic map $e$ is obtained by removing all zero rows. For $m$ we take the columns for the pivots of the row reduction of $f$:

$$\begin{pmatrix} 1 & 2 & 3 \\ 4 & 5 & 6 \\ 7 & 8 & 9 \end{pmatrix} \sim_{\mathrm{rr}} \begin{pmatrix} 1 & 0 & -1 \\ 0 & 1 & 2 \\ 0 & 0 & 0 \end{pmatrix} \xRightarrow[\text{1 and 2}]{\text{pivots}} m = \begin{pmatrix} 1 & 2 \\ 4 & 5 \\ 7 & 8 \end{pmatrix}.$$

And indeed, $f = me$. ◇

As in Example 6G.9, the reduced echelon form of the row is a good and highly recommended way to think about an epic-monic factorization. It, however, comes with one catch: in the process of row reduction of an integer valued matrix one divides by non-zero integers. Over a field this is not an issue, but for $\mathbb{S} = \mathbb{Z}$ one needs the following generalization:

**Example 6G.10.** Over $\mathbf{Mat}_{\mathbb{Z}}$ we need the ***Smith normal form***: Given an integer matrix $f$, one can find two matrices $e'$, $m'$ with determinant $\pm 1$ (thus, they are invertible in $\mathbf{Mat}_{\mathbb{Z}}$), as well as a diagonal matrix $d$ such that $f = m'de'$. Explicitly,

$$f = \begin{pmatrix} 1 & 2 & 3 \\ 4 & 5 & 6 \\ 7 & 8 & 9 \end{pmatrix} \Rightarrow m' = \begin{pmatrix} 1 & 0 & 0 \\ 4 & -1 & 0 \\ 7 & -2 & 1 \end{pmatrix}, \quad d = \begin{pmatrix} 1 & 0 & 0 \\ 0 & 3 & 0 \\ 0 & 0 & 0 \end{pmatrix}, \quad e' = \begin{pmatrix} 1 & 2 & 3 \\ 0 & 1 & 2 \\ 0 & 0 & 1 \end{pmatrix}.$$

Now we run the same procedure as for the row reduced echelon form and get:

$$\mathrm{m}=\begin{pmatrix}1&0\\4&-1\\7&-2\end{pmatrix},\quad \mathrm{e}=\begin{pmatrix}1&0\\0&3\\0&0\end{pmatrix}\begin{pmatrix}1&2&3\\0&1&2\end{pmatrix}=\begin{pmatrix}1&2&3\\0&3&6\end{pmatrix}.$$

In $\mathbf{Vec}_{\mathbb{Z}}$ (but not in $\mathbf{Mat}_{\mathbb{Z}}$) this gives as an epic-monic factorization with $\mathtt{I}=\mathbb{Z}\oplus 3\mathbb{Z}$. ◇

Using Lemma 6G.6, we get:

**Lemma 6G.11.** *In (6G-8), the morphism* e *is epic and the morphism* m *is monic.* ☐

**Definition 6G.12.** For $\mathrm{f}\in\mathbf{C}$ with $\mathbf{C}\in\mathbf{Cat}$ we say $\mathtt{Im}(\mathrm{f})=\big(\mathtt{Im}(\mathrm{f}),\mathrm{m}\colon\mathtt{Im}(\mathrm{f})\hookrightarrow\mathtt{Y}\big)$ (thus, we assume that m is monic) is an ***image of*** f if it has the universal property of the form

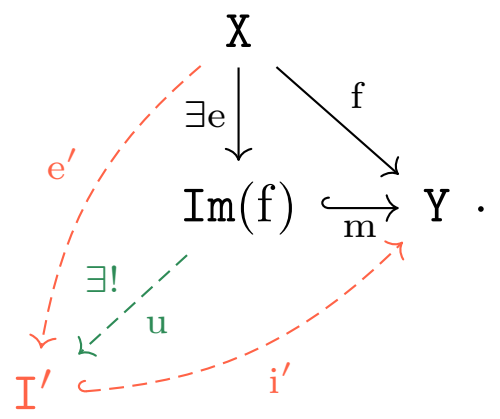

A ***coimage of*** f, denoted by $\mathtt{Coim}(\mathrm{f})=(\mathtt{Coim}(\mathrm{f}),\mathrm{e}\colon\mathtt{Coim}(\mathrm{f})\to\mathtt{Y})$, is an image of f in $\mathbf{C}^{op}$. ◇

**Example 6G.13.** The images in $\mathbf{Vec}_{\Bbbk}$ (in its various incarnations) are the classical images of morphisms. To be completely explicit, take $\mathrm{f}=\begin{pmatrix}1&2\\2&4\end{pmatrix}\in\mathbf{Mat}_{\mathbb{Q}}$. Then we can let

$$\mathtt{Im}(\mathrm{f})=\Big(1,\begin{pmatrix}1\\2\end{pmatrix}\Big),\quad \mathrm{e}=\begin{pmatrix}1&2\end{pmatrix},\quad \mathrm{f}=\begin{pmatrix}1\\2\end{pmatrix}\circ\begin{pmatrix}1&2\end{pmatrix}=\mathrm{me}.$$

An epic-monic factorization is also illustrated above: as in Example 6G.9, we compute that the row reduced echelon form of f is $\begin{pmatrix}1&2\\0&0\end{pmatrix}$, and the factorization follows. ◇

**Lemma 6G.14.** *For* $\mathrm{f}\in\mathbf{C}$ *with* $\mathbf{C}\in\mathbf{Cat}_{\oplus}$ *we have:*

*(i) For an image* $(\mathtt{Im}(\mathrm{f}),\mathrm{m})$ *the morphism* e *is unique and epic. Similarly for coimages.*

*(ii)* $\mathtt{Im}(\mathrm{f})\cong\mathtt{Ker}\big(\mathtt{Coker}(\mathrm{f})\big)$ *and* $\mathtt{Coim}(\mathrm{f})\cong\mathtt{Coker}\big(\mathtt{Ker}(\mathrm{f})\big)$.

*(iii) If* f *has an epic-monic factorization, then* $\mathtt{Im}(\mathrm{f})\cong\mathtt{I}$ *in* (6G-8).

*Proof.* This is Exercise 6L.3. ☐

Consequently, by universality:

**Lemma 6G.15.** *For* $\mathrm{f}\in\mathbf{C}$ *with* $\mathbf{C}\in\mathbf{Cat}_{\oplus}$ *with an epic-monic factorization and an image, this factorization is unique up to unique isomorphism.* ☐

*Remark* 6G.16. Since we will later almost always work $\Bbbk$ linearly, let us stress that, clearly, all the statements above have analogs for $\mathbf{C}\in\mathbf{Cat}_{\mathbb{S}\oplus}$ instead of $\mathbf{C}\in\mathbf{Cat}_{\oplus}$. Similarly, all statements below can be (appropriately) linearized. For example, Theorem 6H.6 holds *verbatim* with A then being a $\Bbbk$ algebra, called a ***presenting algebra*** of $\mathbf{C}$. ◇

## 6H. Abelian categories.

Most invariants of classic topology take place in the following type of categories.

**Definition 6H.1.** A category $\mathbf{C}\in\mathbf{Cat}_{\oplus}$ is called ***abelian*** if every morphism $\mathrm{f}\in\mathbf{C}$ has an epic-monic factorization. ◇

Keeping Example 6G.9 in mind, we can say:

**Slogan.** Abelian categories generalize the row reduced echelon form.

**Example 6H.2.** Here are the first examples:

(a) The names comes from the fact that $\mathbf{Vec}_{\mathbb{Z}}$, see Example 1B.4, is abelian. This follows from Example 6G.10.

(b) Of course, $\mathbf{Vec}_{\Bbbk}$ and $\mathbf{fdVec}_{\Bbbk}$ are abelian.

(c) The categories $\mathbf{Vec}^{\omega}_{\Bbbk\oplus}(\mathrm{G})$ are all abelian.

(d) Diagrammatic categories such as $\mathbf{Br}_{\Bbbk\oplus}$ are usually not abelian.

Abelian categories are in some sense "rare" as we will see in Theorem 6H.6. ◇

*Remark* 6H.3. In Section 6F we have seen linear and additive closures of categories, which allowed us to extend the power of linear algebra to basically any category. There are also several notions of an ***abelian envelope*** (one classical reference is [Fre66]). However, they always come with some form of catch: either they do not preserve structures one might care about, *e.g.* of being monoidal, or they do not exist in general. In other words: we can not naively abelianize our favorite categories. ◇

**Definition 6H.4.** Let A be a ring, which we view as an algebra in $\mathbf{Vec}_{\mathbb{Z}}$. Then ***category of right*** A ***modules*** is defined to be $\mathbf{Mod}(\mathrm{A}) = \mathbf{Mod}_{\mathbf{Vec}_{\mathbb{Z}}}(\mathtt{A})$, *cf.* Section 3F. ◇

The prototypical example of an abelian category is $\mathbf{Mod}(\mathrm{A})$ as will also be seen in Theorem 6H.6.

*Remark* 6H.5. Note that Definition 6H.1 implicitly assumes that kernels and cokernels exist. (This, by Lemma 6G.14, implies that images exist.) In fact, there are many equivalent definitions of abelian categories. However, all of these are meant to be intrinsic descriptions of the "definition" of the concrete abelian categories in Theorem 6H.6. ◇

The following, called ***Freyd–Mitchell theorem***, is the reason why all the above looks very familiar. Note that some of the involved notions will be defined later, but we want to have the theorem stated as soon as possible (a ring is finite dimensional if it is finite dimensional as a $\mathbb{Z}$ algebra):

**Theorem 6H.6.** *We have the following.*

(i) *For every abelian category* $\mathbf{C} \in \mathbf{Cat}_{\oplus}$ *there exist a ring* A *such that*

$$\mathbf{C} \xhookrightarrow{exact} \mathbf{Mod}(\mathrm{A}) ,$$

*i.e.* $\mathbf{C}$ *is equivalent, as an abelian category, to a full subcategory of* $\mathbf{Mod}(\mathrm{A})$.

(ii) *If* $\mathbf{C}$ *is additionally finite, then one can find a finite dimensional* A *such that*

$$\mathbf{C} \xrightarrow{\simeq_e} \mathbf{fdMod}(\mathrm{A}) ,$$

*i.e.* $\mathbf{C}$ *is equivalent, as an abelian category, to* $\mathbf{fdMod}(\mathrm{A})$.

*Proof.* We will sketch a proof later, for now see *e.g.* [Fre64, Theorem 7.34 and Exersice F]. □

The ring A in Theorem 6H.6 is called a ***presenting ring*** of $\mathbf{C}$.

*Remark* 6H.7. The psychologically useful statement in Theorem 6H.6 in words says that we can think of objects of an abelian category as being A modules and of the notions we have seen above, such as *e.g.* kernels, as being those of linear algebra. However, the statement has two drawbacks: neither is A unique nor is it easy to compute in practice. ◇

**Example 6H.8.** For $\mathbf{Vec}_{\Bbbk}$ one can let $\mathrm{A} = \Bbbk$ in Theorem 6H.6, but $\mathrm{A} = \mathrm{M}_{n\times n}(\Bbbk)$, the $\Bbbk$ algebra of $n \times n$ matrices with values in $\Bbbk$, works as well for any $n \in \mathbb{Z}_{\geq 0}$. (This is a special case of ***Morita equivalence***. Roughly, both, $\Bbbk$ and $\mathrm{M}_{n\times n}(\Bbbk)$ have only one simple module, which is either $\Bbbk$ or $\Bbbk^n$ with the evident action.) In this case $\mathbf{Vec}_{\Bbbk} \simeq_{\Bbbk\oplus} \mathbf{Mod}(\mathrm{A})$ for any such A, and these are equivalences of abelian categories. ◇

In order to define appropriate versions of functors between abelian categories (these are called ***exact*** in Theorem 6H.6) we need to understand abelian categories better.

6I. **Exact sequences and functors.** Among others, a 20th century way of doing linear algebra is homological algebra. In homological algebra one usually has a sequence of maps $\delta_i$ and defines homology groups $H_i(X) = \ker(\delta_i)/\mathrm{im}(\delta_{i+1})$ and these serve as a measurement of the complexity of $X$. Hence, one takes some kernel (whose elements are called cycles or cocycles) modulo some image (whose elements are called boundaries or coboundaries), the resulting space are abelian groups or vector spaces. Easy examples of homological algebra can already be found in school arithmetic: carrying is a cocycle [Isa02]. The historically first example is homology for topological spaces, *cf.* Figure 13. Although not obvious at first sight, homology is everywhere in mathematics, and many different homologies are known for topological spaces, groups, Lie groups, knots, graphs, point clouds, *etc.*

Note that kernel, image and quotient make sense in abelian categories. Thus:

**Slogan.** In abelian categories we can run homological algebra.

Recall that in homological algebra one always has certain sequences satisfying exactness properties. Here is the analog:

**Definition 6I.1.** A ***cohomologically written sequence***, or sequence for short, in $\mathbf{C} \in \mathbf{Cat}_{\oplus}$ is a collection of objects $\mathtt{X}_i$ and morphisms $\mathrm{f}_i \colon \mathtt{X}_i \to \mathtt{X}_{i+1}$ for $i \in \mathbb{Z}$. We write $(\mathtt{X}_i, \mathrm{f}_i)^{\bullet} \in \mathbf{C}$ for such sequences, with zero objects sometimes omitted. A ***homologically written sequence*** in $\mathbf{C} \in \mathbf{Cat}_{\oplus}$ is a cohomologically written sequence in $\mathbf{C}^{op}$. ◇

The usual way to illustrate this is

$$\dots \xrightarrow{\mathrm{f}_{i-2}} \mathtt{X}_{i-1} \xrightarrow{\mathrm{f}_{i-1}} \mathtt{X}_i \xrightarrow{\mathrm{f}_i} \mathtt{X}_{i+1} \xrightarrow{\mathrm{f}_{i+1}} \dots , \qquad \dots \xleftarrow{\mathrm{f}_{i-2}} \mathtt{X}_{i-1} \xleftarrow{\mathrm{f}_{i-1}} \mathtt{X}_i \xleftarrow{\mathrm{f}_i} \mathtt{X}_{i+1} \xleftarrow{\mathrm{f}_{i+1}} \dots .$$

By symmetry, we can focus on cohomologically written sequences from now on.

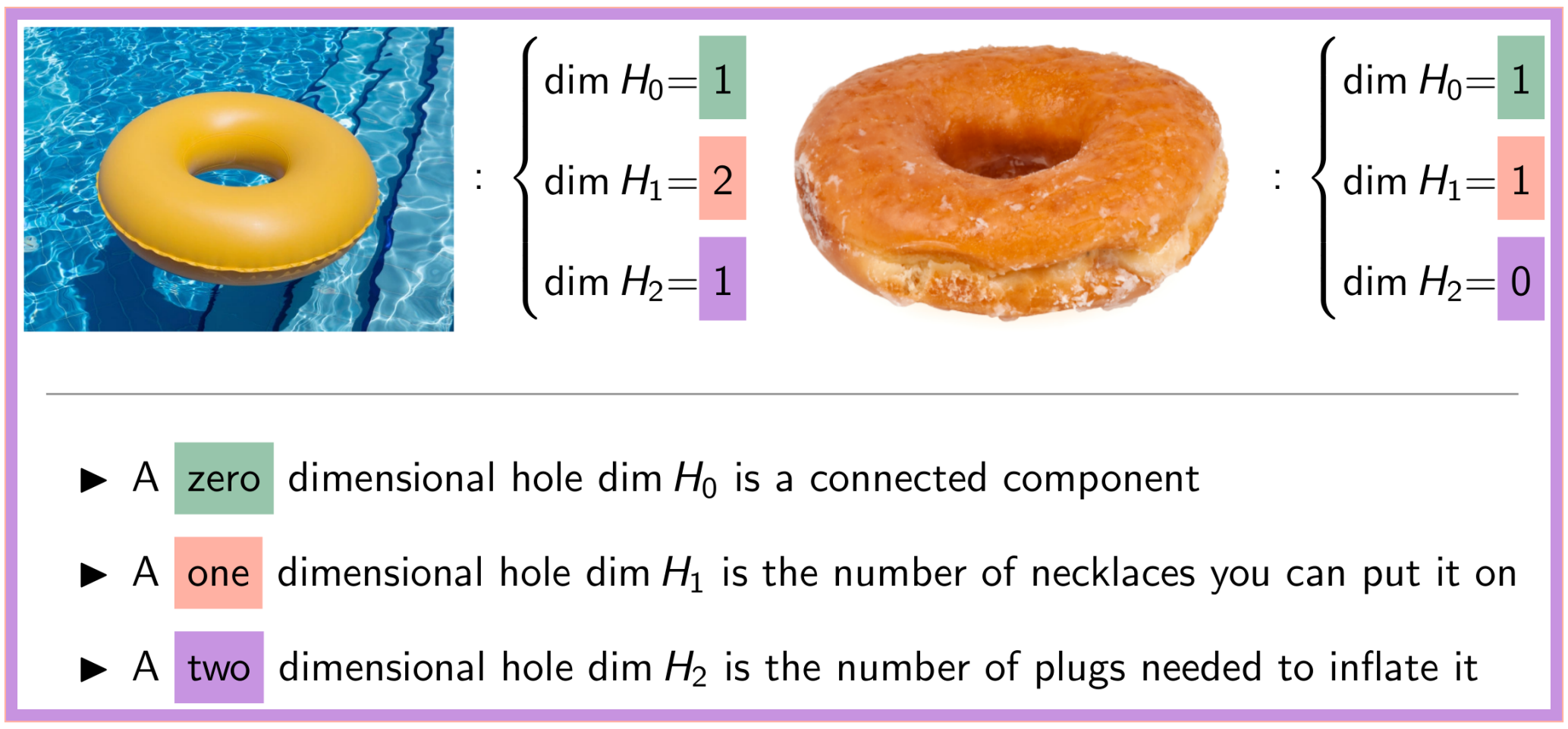

FIGURE 13. Left: homology of a torus = swim ring; right: homology of a solid torus = donut. Picture from the "What is...algebraic topology?" playlist from [Tub21].

**Definition 6I.2.** A sequence $(\mathtt{X}_i,\mathrm{f}_i)^{\bullet}\in\mathbf{C}$ is called ***exact in*** $i$ if $\mathtt{Ker}(\mathrm{f}_i)=\mathtt{Im}(\mathrm{f}_{i-1})$, and ***exact*** if its exact in $i$ for all $i\in\mathbb{Z}$. ◇

**Example 6I.3.** A so-called ***short exact sequence (SES)*** is an exact sequence

$$\mathtt{X}\xhookrightarrow{\mathrm{i}}\mathtt{Y}\xrightarrow{\mathrm{p}}\mathtt{Z}\;=\;\dots\xrightarrow{0}0\xrightarrow{0}0\xrightarrow{0}\mathtt{X}\xrightarrow{\mathrm{i}}\mathtt{Y}\xrightarrow{\mathrm{p}}\mathtt{Z}\xrightarrow{0}0\xrightarrow{0}0\xrightarrow{0}\dots,$$

where i is monic and p is epic by exactness and Lemma 6G.6. Note also that $\mathtt{Z}\cong\mathtt{Y}/\mathtt{X}$.

(a) To be completely explicit, here is a SES in $\mathbf{Mat}_{\mathbb{Q}}$:

$$\mathsf{2}\xhookrightarrow{\mathrm{i}=\left(\begin{smallmatrix}1&0\\0&1\\-1&-2\end{smallmatrix}\right)}\mathsf{3}\xrightarrow{\mathrm{p}=\left(\begin{smallmatrix}1&2&1\end{smallmatrix}\right)}\mathsf{1}\,,\quad \mathsf{2}\xleftarrow{\mathrm{p}'=\left(\begin{smallmatrix}1&0&0\\0&1&0\end{smallmatrix}\right)}\mathsf{3}\xleftarrow{\mathrm{i}'=\left(\begin{smallmatrix}0\\0\\1\end{smallmatrix}\right)}\mathsf{1}\,. \tag{6I-4}$$

The right sequence is a so-called ***splitting*** of the left sequence, meaning that

$$\mathrm{p}'\mathrm{i}=\mathrm{id}_2,\quad \mathrm{p}\mathrm{i}'=\mathrm{id}_1,\quad \mathrm{i}\mathrm{p}'+\mathrm{i}'\mathrm{p}=\mathrm{id}_3. \tag{6I-5}$$

By comparing (6I-5) to (6B-3) we thus, not surprisingly, get that $\mathsf{3}\cong\mathsf{2}\oplus\mathsf{1}$.

(b) A SES in $\mathbf{Vec}_{\mathbb{C}}$ is

$$\begin{aligned}&\mathbb{C}\{X\}\xhookrightarrow{X\mapsto X}\mathbb{C}[X]/(X^2)\cong\mathbb{C}\{1,X\}\xrightarrow{1\mapsto 1,\,X\mapsto 0}\mathbb{C}\{1\}\,,\\ &\mathbb{C}\{X\}\xleftarrow{1\mapsto 0,\,X\mapsto X}\mathbb{C}[X]/(X^2)\cong\mathbb{C}\{1,X\}\xleftarrow{1\mapsto 1}\mathbb{C}\{1\}\,.\end{aligned} \tag{6I-6}$$

The bottom sequence is a splitting of the top, and thus $\mathbb{C}[X]/(X^2)\cong\mathbb{C}\{1\}\oplus\mathbb{C}\{X\}$.

There are a few things to check: First, we have to make sure that the morphisms used are in the correct category. (This sounds obvious, but is crucial: a lot of categories we will see are subcategories of $\mathbf{Vec}_{\Bbbk}$, but not all $\Bbbk$ linear maps are, in general, in such subcategories.) Second, we need to make sure that the left morphism is monic and the right morphism epic. Third, we have to check $\mathtt{Ker}(\mathrm{g})=\mathtt{Im}(\mathrm{f})$. (All of this is easy to see for (6I-4) and (6I-6).) ◇

**Definition 6I.7.** A functor $\mathrm{F}\in\mathbf{Hom}_{\oplus}(\mathbf{C},\mathbf{D})$ is called ***exact*** if

$$\big(\,\mathtt{X}\xhookrightarrow{\mathrm{i}}\mathtt{Y}\xrightarrow{\mathrm{p}}\mathtt{Z}\;\;\text{SES}\,\big)\Rightarrow\big(\,\mathrm{F}(\mathtt{X})\xhookrightarrow{\mathrm{F(i)}}\mathrm{F}(\mathtt{Y})\xrightarrow{\mathrm{F(p)}}\mathrm{F}(\mathtt{Z})\;\;\text{SES}\,\big).$$

Alternatively, we say F preserves SES. ◇

As usual:

**Lemma 6I.8.** *The identity functor on an additive category is exact. Moreover, if* F *and* G *are exact functors, then so is* GF.

**Example 6I.9.** We have a (non-dense, but full) subcategory $\mathbf{Hom}_{e}(\mathbf{C},\mathbf{D})\subset\mathbf{Hom}_{\oplus}(\mathbf{C},\mathbf{D})$, the ***category exact functors***. ◇

Exact functors are the correct functors between abelian categories:

**Lemma 6I.10.** *Let* $\mathrm{F}\in\mathbf{Hom}_{e}(\mathbf{C},\mathbf{D})$ *be a functor between abelian categories. Then:*

*(i) If* $(\mathtt{X}_i,\mathrm{f}_i)^{\bullet}\in\mathbf{C}$ *is exact, then* $\big(\mathrm{F}(\mathtt{X}_i),\mathrm{F}(\mathrm{f}_i)\big)^{\bullet}\in\mathbf{D}$ *is also exact.*

*(ii) If* $(\mathtt{Ker}(f), k) \in \mathbf{C}$ *is a kernel, then* $(F(\mathtt{Ker}(f)), F(k)) \in \mathbf{D}$ *is also a kernel. Similarly for cokernels.*

*(iii) If* $(\mathtt{Im}(f), m) \in \mathbf{C}$ *is an image, then* $(F(\mathtt{Im}(f)), F(m)) \in \mathbf{D}$ *is also an image. Similarly for coimages.*

*(iv) If* f *is monic, respectively epic, then* F(f) *is monic, respectively epic.*

*(v) If*

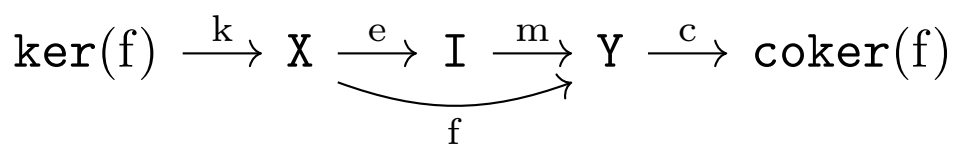

*is an epic-monic factorization in* **C***, then*

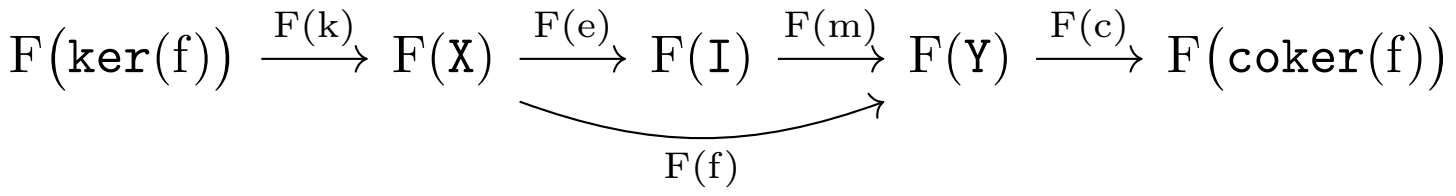

*is also an epic-monic factorization in* **D***.*

*Proof. (i).* Note that being an exact sequence implies that we have

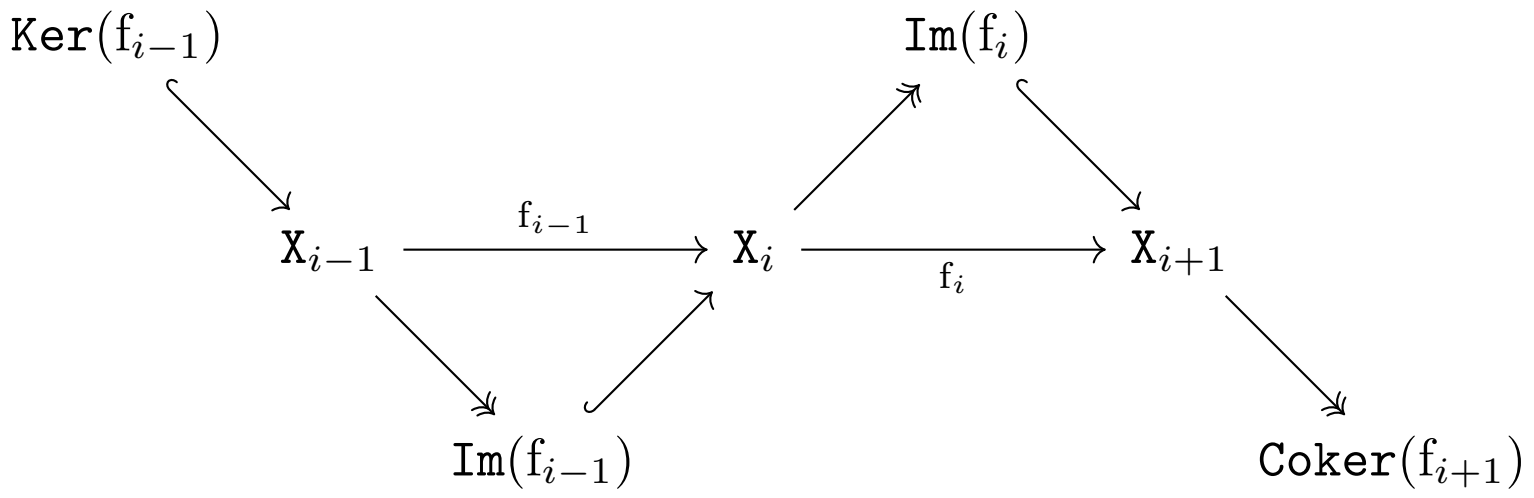

which commutes and has SES diagonals. Applying F yields a commuting diagram

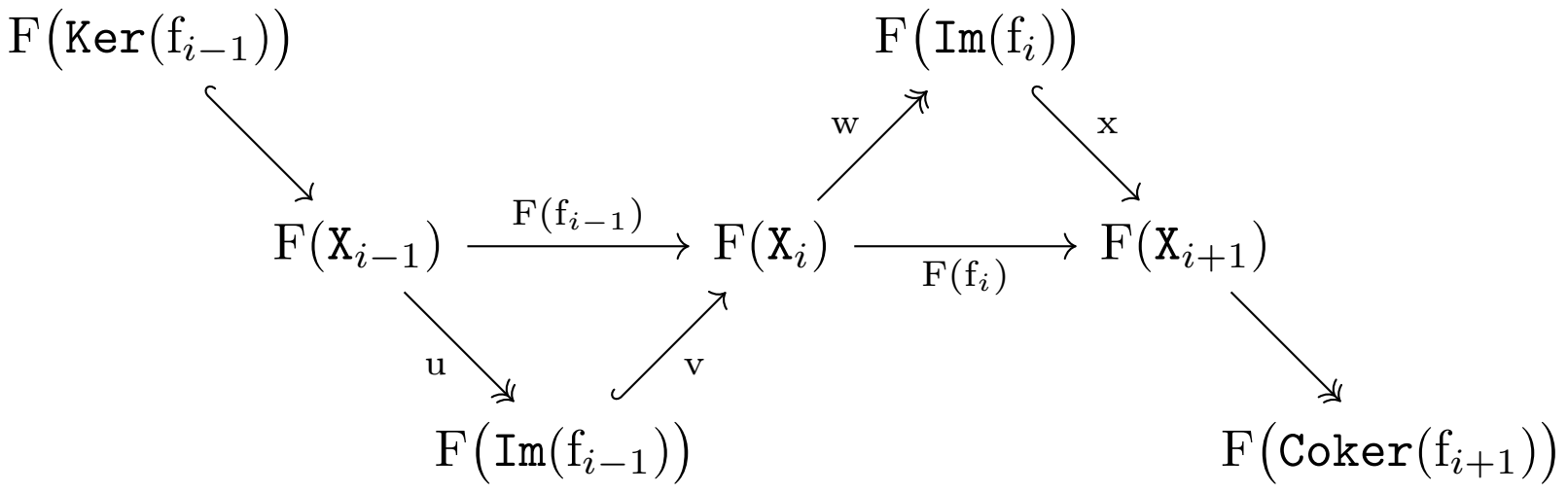

with still has SES diagonals. Then

$$\mathtt{Im}\big(F(f_{i-1})\big) = \mathtt{Im}(vu) = \mathtt{Im}(v) = \mathtt{Ker}(w) = \mathtt{Ker}(xw) = \mathtt{Ker}\big(F(f_i)\big),$$

shows the claim.

*(ii).* Kernels and cokernels are special cases of exact sequences.

*(iii).* We use (ii) and Lemma 6G.14.(iii).

*(iv).* By (ii) and Lemma 6G.6.(iii).

*(v).* Clear by the other statements. □

**Example 6I.11.** We get a (non-full) subcategory $\mathbf{Cat}_A \subset \mathbf{Cat}_\oplus$, the ***category of abelian categories*** with morphisms being exact functors. ◇

**Definition 6I.12.** $\mathbf{C}, \mathbf{D} \in \mathbf{Cat}_A$ are called ***equivalent as abelian categories***, denoted by $\mathbf{C} \simeq_e \mathbf{D}$, if there exists an equivalence $F \in \mathbf{Hom}_e(\mathbf{C}, \mathbf{D})$. ◇

Recall the hom functors from Example 1E.6, which will now again play a crucial role.

**Example 6I.13.** As we have seen, any $\mathbf{C} \in \mathbf{Cat}_A$ is equivalent as an abelian category to some full subcategory of $\mathbf{Mod}(A)$ for some appropriate ring A. The additive versions of the Yoneda embedding, *cf.* Proposition 1H.28, almost do the job:

$$Y\colon \mathbf{C} \to \mathbf{Hom}_\oplus(\mathbf{C}^{op}, \mathbf{Vec}_\mathbb{Z}), \quad Y^{op}\colon \mathbf{C}^{op} \to \mathbf{Hom}_\oplus(\mathbf{C}, \mathbf{Vec}_\mathbb{Z})$$

is fully faithful, and thus, **C** is equivalent as an additive category to a full subcategory of right **C** modules. (Here we think of $\mathbf{Hom}_\oplus(\mathbf{C}^{op}, \mathbf{Vec}_\mathbb{Z})$ as right **C** modules.) However, an SES is in general not preserved since

$$\begin{gathered}\mathrm{Hom}_\mathbf{C}(_, \mathtt{X}) \xrightarrow{\mathrm{Hom}_\mathbf{C}(_, i)} \mathrm{Hom}_\mathbf{C}(_, \mathtt{Y}) \xrightarrow{\mathrm{Hom}_\mathbf{C}(_, p)} \mathrm{Hom}_\mathbf{C}(_, \mathtt{Z}), \\ \mathrm{Hom}_\mathbf{C}(\mathtt{X}, _) \xrightarrow{\mathrm{Hom}_\mathbf{C}(i, _)} \mathrm{Hom}_\mathbf{C}(\mathtt{Y}, _) \xrightarrow{\mathrm{Hom}_\mathbf{C}(p, _)} \mathrm{Hom}_\mathbf{C}(\mathtt{Z}, _),\end{gathered} \tag{6I-14}$$

are not exact in the rightmost, respectively leftmost, position, even if one starts with an SES. ◇

So, the Yoneda embedding is not exact and does not prove Theorem 6H.6, at least not directly:

*Sketch of a proof of Theorem 6H.6.(i).* It is easy to see that $\mathrm{Hom}_{\mathbf{C}}(_,\mathtt{X})$ is an additive and ***right exact functor***, which means that it sends an SES to a sequence as in the first row of (6I-14), which is exact except to the far left. Such functors form a category $\mathbf{Hom}_{re}\big(\mathbf{C}^{\mathrm{op}},\mathbf{Vec}_{\mathbb{Z}}\big)$, and one shows the following (non-trivial) statements:

- The category $\mathbf{Hom}_{re}(\mathbf{C}^{\mathrm{op}},\mathbf{Vec}_{\mathbb{Z}})$ is abelian.
- The adjusted Yoneda embedding $\mathrm{Y}^{re}\colon \mathbf{C}\to\mathbf{Hom}_{re}\big(\mathbf{C}^{\mathrm{op}},\mathbf{Vec}_{\mathbb{Z}}\big)$ with $\mathtt{X}\mapsto\mathrm{Hom}_{\mathbf{C}}(_,\mathtt{X})$ is additive exact and fully faithful.
- There exists an object $\mathtt{I}\in\mathbf{Hom}_{re}\big(\mathbf{C}^{\mathrm{op}},\mathbf{Vec}_{\mathbb{Z}}\big)$ whose endomorphism ring
$$\mathrm{A}=\mathrm{End}_{\mathbf{Hom}_{re}(\mathbf{C}^{\mathrm{op}},\mathbf{Vec}_{\mathbb{Z}})^{op}}(\mathtt{I})$$
provides an abelian category $\mathbf{Mod}(\mathrm{A})$ equivalent to $\mathbf{Hom}_{re}(\mathbf{C}^{\mathrm{op}},\mathbf{Vec}_{\mathbb{Z}})$ as an abelian category.

The proof is complete. □

The projective, respectively, injective, objects correct the "failure" in (6I-14) (in fact, the object $\mathtt{I}$ from the above proof sketch is a certain nice injective object called a ***cogenerator***):

**Definition 6I.15.** Let $\mathbf{C}\in\mathbf{Cat}$. Then $\mathtt{P}\in\mathbf{C}$ is ***projective*** if it has the universal property of the form

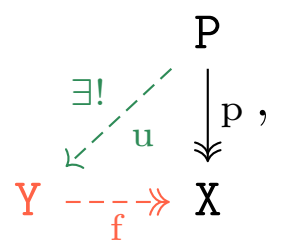

for any epic morphism f. Moreover, $\mathtt{I}\in\mathbf{C}$ is ***injective*** if it has the universal property of the form

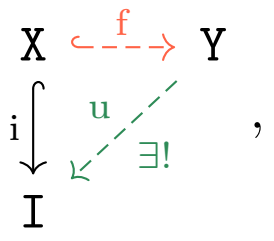

for any monic morphism f. ◇

We view projective and injective as dual in the sense of Section 1C. The following will be justified in Remark 6K.9. We park it here to demystify Definition 6I.15 a bit:

**Slogan.** Projective objects are the elemental summands of the algebra presenting an abelian category.

**Example 6I.16.** The projective objects in $\mathbf{Vec}_{\mathbb{Z}}$ are the free $\mathbb{Z}$ modules. Explicitly, the only finite generated projective objects in $\mathbf{Vec}_{\mathbb{Z}}$ are $\mathbb{Z}^n$ for some $n\in\mathbb{Z}_{\geq 0}$. ◇

The following are (almost) immediate.

**Lemma 6I.17.** *If $\mathbf{C}\in\mathbf{Cat}_A$, then $\mathtt{P}\in\mathbf{C}$ is projective if and only if*
$$\mathrm{Hom}_{\mathbf{C}}(\mathtt{P},_)\in\mathbf{Hom}_e\big(\mathbf{C},\mathbf{Vec}_{\mathbb{Z}}\big).$$
*Moreover, $\mathtt{I}\in\mathbf{C}$ is injective if and only if*
$$\mathrm{Hom}_{\mathbf{C}}(_,\mathtt{I})\in\mathbf{Hom}_e\big(\mathbf{C}^{op},\mathbf{Vec}_{\mathbb{Z}}\big).$$
*That is, the corresponding hom functors are exact.* □

**Lemma 6I.18.** *In additive categories, being projective is an additive property: two objects $\mathtt{P},\mathtt{P}'\in\mathbf{C}$ are projective if and only if $\mathtt{P}\oplus\mathtt{P}'$ is projective. Similarly for injective objects.* □

**Example 6I.19.** By Lemma 6I.18 we have two additive full subcategories $\mathbf{Proj}(\mathbf{C})$, the ***category of projective objects***, and $\mathbf{Inj}(\mathbf{C})$ the ***category of injective objects***, for all $\mathbf{C}\in\mathbf{Cat}_{\oplus}$. ◇

**Definition 6I.20.** Let $\mathbf{C}\in\mathbf{Cat}$ and $\mathtt{X}\in\mathbf{C}$. We say $\mathtt{P}(\mathtt{X})=\big(\mathtt{P}(\mathtt{X}),\mathrm{f}\colon\mathtt{P}(\mathtt{X})\twoheadrightarrow\mathtt{X}\big)$ is a ***projective cover of*** $\mathtt{X}$ if $\mathtt{P}(\mathtt{X})$ is projective and has the universal property of the form

(6I-21)
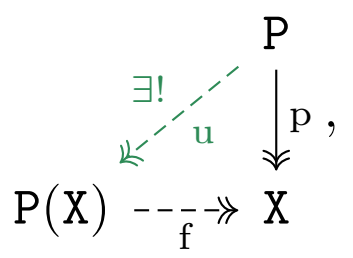

where $\mathtt{P}$ is projective.

An ***injective hull of*** f, denoted by $\mathtt{I}(\mathtt{X})=(\mathtt{I}(\mathtt{X}),\mathrm{i}\colon\mathtt{X}\hookrightarrow\mathtt{I})$, is a projective cover of $\mathtt{X}$ in $\mathbf{C}^{op}$. ◇

The philosophy is a bit that every object is a quotient of a projective object and a subobject of an injective object, and the projective cover and the injective envelope are the universal objects achieving that. Thus, not surprisingly:

**Lemma 6I.22.** *Up to unique isomorphisms,* $\mathtt{P}(\mathtt{X})$ *is the only object in* $\mathbf{C}$ *satisfying* (6I-21). *Similarly for the injective hull.* □

6J. **The "elements" of additive and abelian categories.** Recall that, in chemistry, an ***element*** is a substance that cannot be broken down into other substances by reactions. They are simple (the words "simple" is meant in the sense that they are "as simple as possible", and not meaning they are easy), and make up all matter.

There are (at least) two competing ways to define "elements": Either these are objects without substructure, called ***simple***. Or, these are objects that cannot be further decomposed, called ***indecomposable***. The analogy between chemistry and abelian categories goes even further, *cf.* Figure 14.

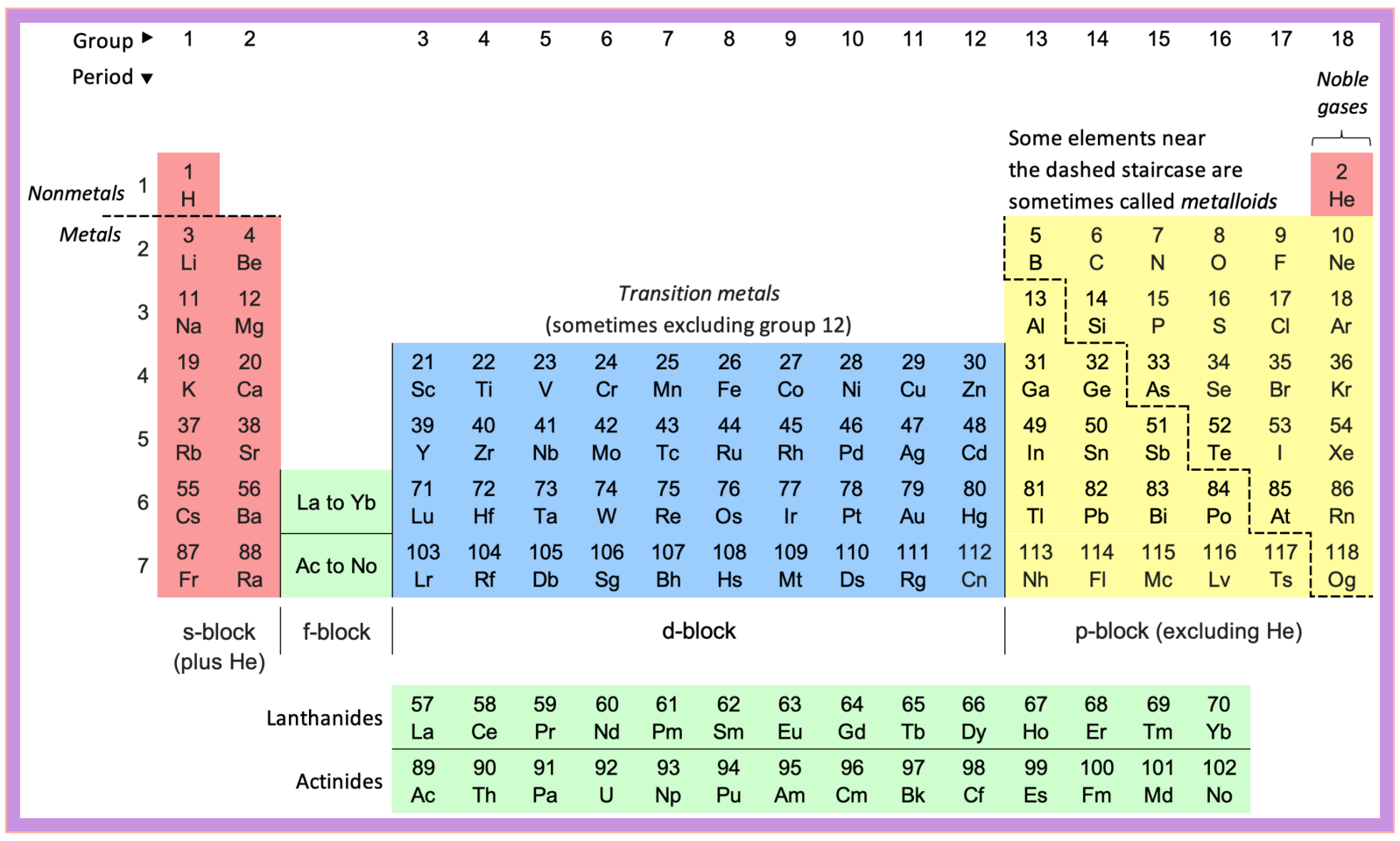

FIGURE 14. The simple/indecomposable objects play the role of elements in our story. A main goal in abelian categories is to finite the corresponding periodic table.

Picture from https://en.wikipedia.org/wiki/Periodic_table#/media/File:Colour_18-col_PT_with_labels.png

**Definition 6J.1.** Let $\mathbf{C} \in \mathbf{Cat}_{\oplus}$.

(i) A non-zero object $\mathtt{L} \in \mathbf{C}$ is called ***simple*** if

$$(\mathtt{X} \subset \mathtt{L}) \Rightarrow (\mathtt{X} \cong 0 \text{ or } \mathtt{X} \cong \mathtt{L}).$$

(ii) A non-zero object $\mathtt{Z} \in \mathbf{C}$ is called ***indecomposable*** if

$$\mathtt{Z} \cong \mathtt{X} \oplus \mathtt{Y} \Rightarrow (\mathtt{X} \cong 0 \text{ or } \mathtt{Y} \cong 0).$$

***Finding the periodic table of*** $\mathbf{C}$ means classifying either of these up to isomorphism. ◇

We also say a decomposition $\mathtt{X}' \cong \mathtt{X} \oplus \mathtt{Y}$ is ***non-trivial*** if neither $\mathtt{X}$ nor $\mathtt{Y}$ are zero. Similarly, a subobject $\mathtt{Y} \subset \mathtt{X}$ is ***non-trivial*** if it is neither 0 nor (isomorphic to) $\mathtt{X}$.

*Remark* 6J.2. Note that indecomposable means that an object has no non-trivial decomposition, while simple means that an object has no non-trivial subobjects. ◇

The following lemma is clear and enables us to define the set of simples $\mathrm{Si}(\mathbf{C}) \subset \mathrm{Ob}(\mathbf{C})/\cong$, respectively, indecomposables $\mathrm{In}(\mathbf{C}) \subset \mathrm{Ob}(\mathbf{C})/\cong$ (up to isomorphism). We, abusing notation, write *e.g.* $\mathtt{L} \in \mathrm{Si}(\mathbf{C})$ for simplicity.

**Lemma 6J.3.** *The properties of being indecomposable or simple are preserved under isomorphisms.* □

Note that being projective or injective is also preserved under isomorphism. Hence, by Lemma 6I.18 we also have the sets of projective indecomposables $\mathrm{Pi}(\mathbf{C}) \subset \mathrm{Ob}(\mathbf{C})/\cong$ and injective indecomposables $\mathrm{Ii}(\mathbf{C}) \subset \mathrm{Ob}(\mathbf{C})/\cong$, respectively.

**Lemma 6J.4.** *Every simple object* $\mathtt{L} \in \mathbf{C}$ *is indecomposable.*

*Proof.* Clearly, any non-trivial decomposition $\mathtt{L} \cong \mathtt{X} \oplus \mathtt{Y}$ gives non-trivial subobjects $\mathtt{X}$ and $\mathtt{Y}$. □

**Example 6J.5.** Being simple or indecomposable depends on the ambient category, compare (6J-6) and (6J-8):

(a) For $\mathbf{Vec}_{\Bbbk}$ it is easy to see that $\mathrm{Si}(\mathbf{Vec}_{\Bbbk}) = \mathrm{In}(\mathbf{Vec}_{\Bbbk}) = \mathrm{Pi}(\mathbf{Vec}_{\Bbbk}) = \mathrm{Ii}(\mathbf{Vec}_{\Bbbk}) = \{\Bbbk\}$. Thus, $\mathtt{A} = \mathbb{C}[X]/(X^2) \in \mathbf{Vec}_{\mathbb{C}}$ is neither simple nor indecomposable and we have

$$\mathtt{A} \cong \mathbb{C} \oplus \mathbb{C}, \quad (\text{in } \mathbf{Vec}_{\mathbb{C}}), \tag{6J-6}$$

*cf.* (6I-6). Moreover, every object in $\mathbf{Vec}_{\Bbbk}$ is projective and injective.

(b) Consider now $\mathrm{A} = \mathbb{C}[X]/(X^2)$ as a $\mathbb{C}$ algebra. Then A acts on itself by multiplication; thus A can be seen as an object $\mathtt{A}$ of (the $\mathbb{C}$ linear abelian category) $\mathbf{Mod}(\mathrm{A})$. The $\mathbb{C}$ algebra A also acts on $\mathbb{C}$ through evaluation, so we have two objects $\mathbb{C}, \mathtt{A} \in \mathbf{Mod}(\mathrm{A})$. Choose $\{1, X\}$ as a basis of $\mathtt{A}$. Looking at the action matrices on this basis gives

$$1 \mapsto \begin{pmatrix} \boxed{1} & 0 \\ 0 & \boxed{1} \end{pmatrix}, \quad X \mapsto \begin{pmatrix} \boxed{0} & 0 \\ 1 & \boxed{0} \end{pmatrix}, \quad \begin{matrix} \text{lower right block entry} \rightsquigarrow \mathbb{C}, \\ \text{upper left block entry} \rightsquigarrow \mathbb{C}. \end{matrix} \tag{6J-7}$$

This shows that $\mathbb{C}$ is a subobject of $\mathtt{A}$ (indicated in (6J-7)), and hence $\mathtt{A}$ is not simple. However, the very same action matrices show that the complement space $\mathbb{C}\{1\}$ is not a subobject (the entry 1 in the lower left ruins this). However, one easily sees that $\mathtt{A}$ is (projective injective) indecomposable and

$$\underbrace{0 \subset \mathbb{C}}_{\boxed{\mathbb{C}}} = \underbrace{\mathbb{C} \subset \mathtt{A}}_{\boxed{\mathbb{C}}} \quad \mathtt{A} \ncong \mathbb{C} \oplus \mathbb{C}, \quad (\text{in } \mathbf{Mod}(\mathrm{A})), \tag{6J-8}$$

with the right copy of $\mathbb{C}$ being the upper left block entry and the left copy of $\mathbb{C}$ being the lower right block entry in (6J-7).

The difference between $\mathbf{Vec}_{\Bbbk}$ and $\mathbf{Mod}(\mathrm{A})$ is that the maps defining the decomposition from (6J-6) are not A equivariant, *i.e.* they are not morphisms in $\mathbf{Mod}(\mathrm{A})$. Precisely, we still have

$$\boxed{\mathbb{C}} \cong \mathbb{C}\{X\} \xhookrightarrow{X \mapsto X} \mathbb{C}[X]/(X^2) \xrightarrow{1 \mapsto 1, X \mapsto 0} \mathbb{C}\{1\} \cong \boxed{\mathbb{C}} \quad \text{SES},$$

but it does not split in contrast to (6I-6). ◇

Thus, Lemma 6J.4 and Example 6J.5 give:

$$\text{indecomposable} \Leftarrow \text{simple},$$
$$\text{indecomposable} \not\Rightarrow \text{simple}.$$

The following is known as ***Schur's lemma*** (or at least (i) of it).

**Lemma 6J.9.** *Let* $\mathbf{C} \in \mathbf{Cat}_{\oplus}$.

*(i) If* $\mathbf{C}$ *has kernels and cokernels, then, for any* $\mathtt{L}, \mathtt{L}' \in \mathrm{Si}(\mathbf{C})$ *with* $\mathtt{L} \ncong \mathtt{L}'$*:*

$$\mathrm{End}_{\mathbf{C}}(\mathtt{L}) \textit{ is a division ring}, \quad \mathrm{Hom}_{\mathbf{C}}(\mathtt{L}, \mathtt{L}') \cong 0.$$

*(ii) For any* $\mathtt{Z} \in \mathrm{In}(\mathbf{C})$ *we have*

$$\mathrm{End}_{\mathbf{C}}(\mathtt{Z}) \textit{ is a local ring}.$$

*Proof.* This is Exercise 6L.5. □

Schur's lemma, part II:

**Lemma 6J.10.** *Let* $\mathbf{C} \in \mathbf{Cat}_{\mathbb{K}A}$ *with* $\mathbb{K}$ *being algebraically closed. Then, for any* $\mathtt{L}, \mathtt{L}' \in \mathrm{Si}(\mathbf{C})$ *with* $\mathtt{L} \ncong \mathtt{L}'$*, we have*

$$\mathrm{End}_{\mathbf{C}}(\mathtt{L}) \cong \mathbb{K}, \quad \mathrm{Hom}_{\mathbf{C}}(\mathtt{L}, \mathtt{L}') \cong 0. \tag{6J-11}$$

*Proof.* This is also Exercise 6L.5. □

**Example 6J.12.** With respect to Example 6J.5 we have

$$\mathrm{End}_{\mathbf{Vec}_{\mathbb{C}}}(\mathbb{C}) \cong \mathbb{C} \cong \mathrm{End}_{\mathbf{Mod}(\mathrm{A})}(\mathbb{C}),$$
$$\mathrm{End}_{\mathbf{Vec}_{\mathbb{C}}}\big(\mathbb{C}[X]/(X^2)\big) \cong \mathrm{Mat}_{2\times 2}(\mathbb{C}), \quad \mathrm{End}_{\mathbf{Mod}(\mathrm{A})}\big(\mathbb{C}[X]/(X^2)\big) \cong \mathbb{C}[X]/(X^2),$$

and the idempotents $\left(\begin{smallmatrix} 1 & 0 \\ 0 & 0 \end{smallmatrix}\right), \left(\begin{smallmatrix} 0 & 0 \\ 0 & 1 \end{smallmatrix}\right) \in \mathrm{Mat}_{2\times 2}(\mathbb{C})$ give the decomposition in (6J-6). ◇

The assumption in Lemma 6J.10 that $\mathbb{K}$ is algebraically closed is necessary. For an explicit example, recall Frobenius theorem (on real division algebras): the only finite dimensional real division algebras are $\mathbb{R}$, $\mathbb{C}$ and $\mathbb{H}$ of real dimensions 1, 2 and 4. Here $\mathbb{H}$ are the famous (but slightly useless) ***quaternions***, see Figure 15. All of these can appear as endomorphism rings when the ground field is $\mathbb{R}$ as the next example shows.

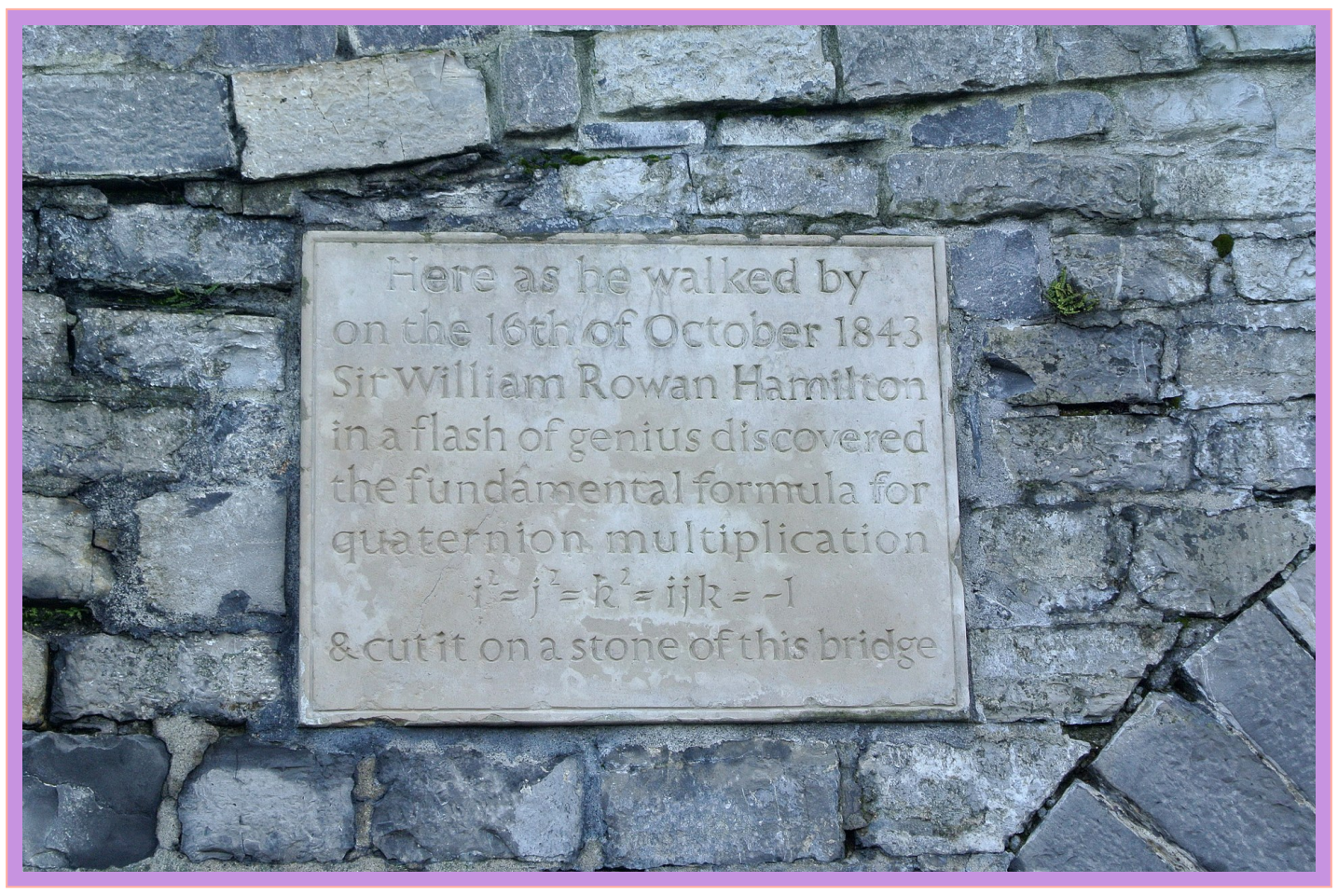

FIGURE 15. Quaternion plaque: "Here as they walked by on the 16th of October 1843 William Rowan Hamilton in a flash of genius discovered the fundamental formula for quaternion multiplication $i^2 = j^2 = k^2 = ijk = -1$ & cut it on a stone of this bridge"

Picture from https://en.wikipedia.org/wiki/Quaternion

**Example 6J.13.** Take $G = \mathbb{Z}/3\mathbb{Z} = \langle a | a^3 = 1 \rangle$ and $\mathrm{A} = \mathbb{C}G$ and $\mathrm{B} = \mathbb{R}G$. Fix $\zeta = \exp(2\pi i/3)$ (a primitive third root of unity). There are three one dimensional simple objects $\mathtt{L}_i \in \mathbf{Mod}(\mathrm{A})$, for $i \in \{0, 1, 2\}$, determined by $a \mapsto \zeta^i$. However, for B only $\mathtt{L}_0$ makes sense as $\zeta, \zeta^2 \notin \mathbb{R}$, but $\mathtt{L}_1 \oplus \mathtt{L}_2 \in \mathbf{Mod}(\mathrm{B})$ is simple. We have

$$\mathrm{End}_{\mathbf{Mod}(\mathrm{B})}(\mathtt{L}_1 \oplus \mathtt{L}_2) \cong \mathbb{C} \quad \text{(as a division } \mathbb{R} \text{ algebra)}.$$

In general, for $\mathrm{A} = \mathbb{C}\mathbb{Z}/n\mathbb{Z}$ and $\mathrm{B} = \mathbb{R}\mathbb{Z}/n\mathbb{Z}$ the picture is illustrated in Figure 16: all simple $\mathtt{L}_\zeta$ objects in $\mathbf{Mod}(\mathrm{A})$ are given by sending a generator of $\mathbb{Z}/n\mathbb{Z}$ to an $n$th root of unity $\zeta$. Say $\zeta = \exp(ik\Theta) \neq \pm 1$. for $\Theta = 2\pi/n$. Using

$$\begin{pmatrix} \exp(ik\Theta) & 0 \\ 0 & \overline{\exp(ik\Theta)} \end{pmatrix} \sim_{\mathbb{C}} \begin{pmatrix} \cos(k\Theta) & -\sin(k\Theta) \\ \sin(k\Theta) & \cos(k\Theta) \end{pmatrix}$$

one gets a real, rotation, action on the sum $\mathtt{L}_\zeta \oplus \mathtt{L}_{\overline{\zeta}}$, which is a simple object in $\mathbf{Mod}(\mathrm{B})$, and

$$\mathrm{End}_{\mathbf{Mod}(\mathrm{B})}(\mathtt{L}_\zeta \oplus \mathtt{L}_{\overline{\zeta}}) \cong \mathbb{C} \quad \text{(as a division } \mathbb{R} \text{ algebra)}.$$

Finding an example where the endomorphism ring is $\mathbb{H}$ requires a bit more work, but not too much. If one takes $G$ to be the quaternion group of order eight, then there is a four dimensional simple object in $\mathbf{Mod}(\mathbb{R}G)$ with endomorphism ring $\mathbb{H}$. But more on this later. ◇

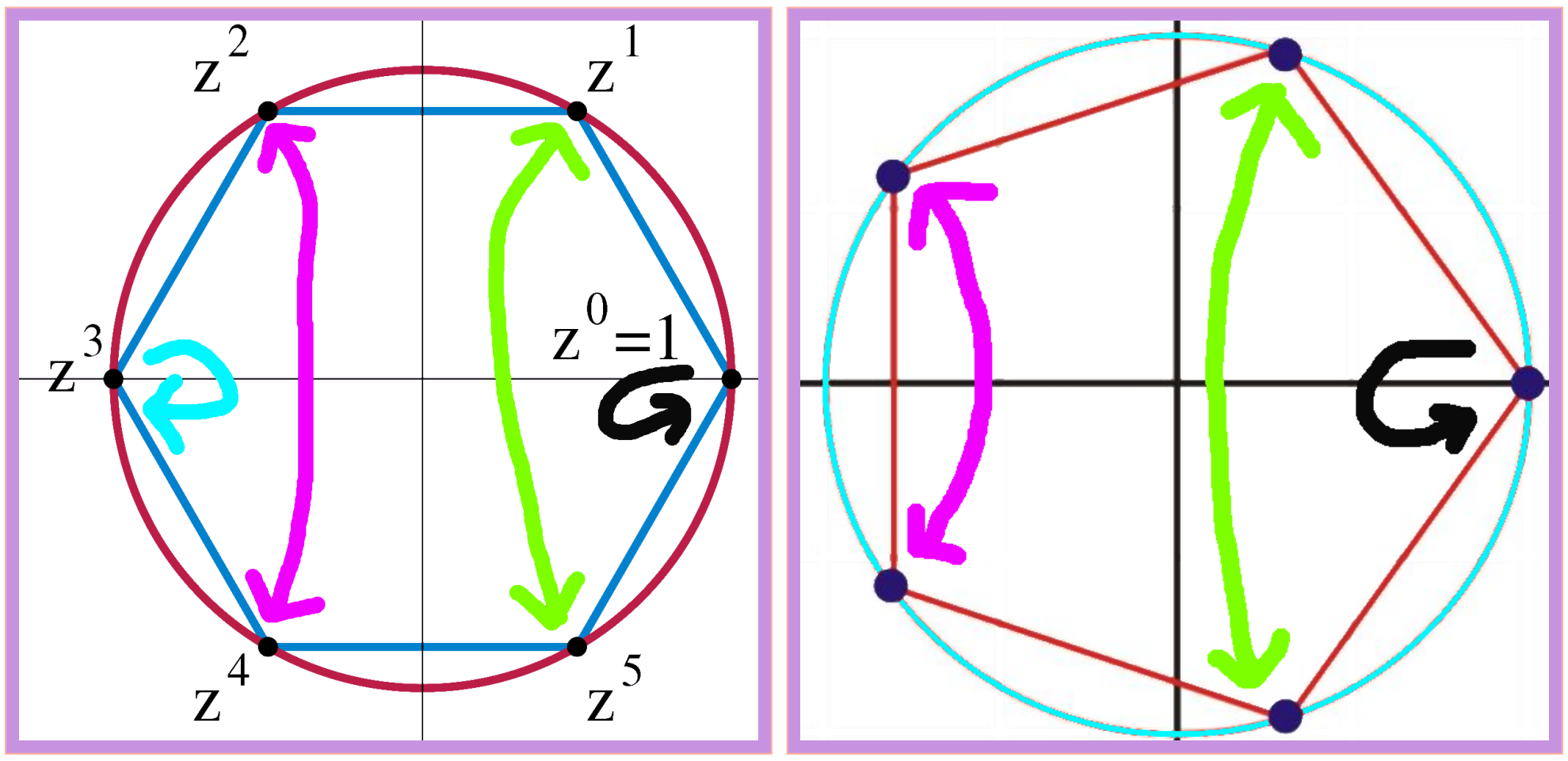

FIGURE 16. $n = 5$ and $n = 6$ in Example 6J.13 illustrated geometrically.

Picture from https://en.wikipedia.org/wiki/Quaternion

Note that any ("finite") $\mathtt{X} \in \mathbf{C}$ with $\mathbf{C} \in \mathbf{Cat}_\oplus$, by definition, decomposes additively into indecomposables. However, Example 6J.5 shows that it is too much to hope that $\mathtt{X}$ decomposes additively into simples. We rather need the analog of (6J-8):

**Definition 6J.14.** Assume that $\mathbf{C} \in \mathbf{Cat}_{\oplus}$ has kernels and cokernels. For a non-zero $\mathtt{X} \in \mathbf{C}$ a sequence of subobjects of the form

$$0 = \mathtt{X}_0 \subset \mathtt{X}_1 \subset \dots \subset \mathtt{X}_{n-1} \subset \mathtt{X}_n = \mathtt{X} \tag{6J-15}$$

is called a ***filtration by simples*** or a ***composition series*** if $\mathtt{X}_i/\mathtt{X}_{i-1} \cong \mathtt{L}_i \in \mathrm{Si}(\mathbf{C})$. A non-zero $\mathtt{X} \in \mathbf{C}$ is called of ***finite length*** if it has such a filtration, and in this case the appearing $\mathtt{L}_i$ are called the ***simple factors*** of $\mathtt{X}$. ◇

We stress that the main point in (6J-15) is that successive quotients are simple:

$$0 = \underbrace{\mathtt{X}_0 \subset \mathtt{X}_1}_{\mathtt{L}_1} \underbrace{\subset \dots \subset}_{\dots} \underbrace{\mathtt{X}_{n-1} \subset \mathtt{X}_n}_{\mathtt{L}_n} = \mathtt{X}.$$

**Definition 6J.16.** For a non-zero $\mathtt{X} \in \mathbf{C}$ with $\mathbf{C} \in \mathbf{Cat}_{\oplus}$ a decomposition of the form

$$\mathtt{X} \cong \mathtt{Z}_1 \oplus \dots \oplus \mathtt{Z}_n \tag{6J-17}$$

is called a ***decomposition by indecomposables*** if $\mathtt{Z}_i \in \mathrm{In}(\mathbf{C})$. A non-zero $\mathtt{X} \in \mathbf{C}$ is called of ***finite decomposition length*** if it has such a decomposition, and in this case the appearing $\mathtt{Z}_i$ are called the ***indecomposable summands*** of $\mathtt{X}$. ◇

**Example 6J.18.** In $\mathbf{Vec}_{\mathbb{k}}$ an object is of finite length if and only if it is finite dimensional. Moreover, for finite dimensional $\mathbb{k}$ vector spaces (6J-15) and (6J-17) agree. ◇

The following theorem is our justification for using the analogy to elements in chemistry, where Theorem 6J.19.(i) is known as the ***Jordan–Hölder theorem***, and (ii) as the ***Krull–Schmidt theorem***.

**Theorem 6J.19.** *Let $\mathbf{C} \in \mathbf{Cat}_{\oplus}$.*

*(i) Assume that $\mathbf{C} \in \mathbf{Cat}_{\oplus}$ has kernels and cokernels. Let $\mathtt{X} \in \mathbf{C}$ be of finite length. Then a filtration as in (6J-15) is unique up to reordering and isomorphisms of subobjects.*

*(ii) Let $\mathtt{X} \in \mathbf{C}$ be of finite decomposition length. Then a decomposition as in (6J-17) is unique up to reordering and isomorphisms of summands.*

In particular, we can define the following ***numerical invariants*** of such $\mathtt{X}$.

- The ***length*** $\ell(\mathtt{X})$ of $\mathtt{X}$ can be defined to be $n$ in (6J-15), and the ***decomposition length*** $d(\mathtt{X})$ of $\mathtt{X}$ can be defined to be $n$ in (6J-17).
- The ***multiplicities*** of $\mathtt{L}$ (simple) respectively $\mathtt{Z}$ (indecomposable) in $\mathtt{X}$ denoted by

$$[\mathtt{X} : \mathtt{L}] = \#\{i \mid \mathtt{X}_i/\mathtt{X}_{i-1} \cong \mathtt{L}\}, \quad (\mathtt{X} : \mathtt{Z}) = \#\{i \mid \mathtt{Z}_i \cong \mathtt{I}\}.$$

- The sets

$$\big\{(\mathtt{L}, m) \mid \mathtt{L} \text{ is a simple factor of } \mathtt{X} \text{ with multiplicity } m\big\},$$
$$\big\{(\mathtt{Z}, m) \mid \mathtt{Z} \text{ is an indecomposable summand of } \mathtt{X} \text{ with multiplicity } m\big\}.$$

*Remark* 6J.20. We will use the notion "numerical" quite often and this is to be understood as reducing notions from categorical algebra to "something easier" such as classical algebra, combinatorics, linear algebra *etc.* Thus, a "numerical invariant" for us is not necessarily a number, but simply something that is "easier" than the problem at hand. ◇

*Proof of Theorem 6J.19.* We only prove (i), the arguments for (ii) are similar. The proof works by induction over $n \geq 1$, with $n$ being the smallest possible length of a filtration by simples. For $n = 1$ there is nothing to show, since $\mathtt{X}$ is then itself simple. So assume that $n > 1$ and that we have two filtrations

$$0 = \mathtt{X}_0 \subset \mathtt{X}_1 \subset \dots \subset \mathtt{X}_{n-1} \subset \mathtt{X}_n = \mathtt{X}, \quad \text{simple factors } \mathtt{L}_1, \dots, \mathtt{L}_n,$$
$$0 = \mathtt{X}'_0 \subset \mathtt{X}'_1 \subset \dots \subset \mathtt{X}'_{n'-1} \subset \mathtt{X}'_{n'} = \mathtt{X}, \quad \text{simple factors } \mathtt{L}'_1, \dots, \mathtt{L}'_{n'},$$

with $n$ being minimal. There are now two cases. First, if $\mathtt{X}_1 \cong \mathtt{X}'_1 \cong \mathtt{L}_1 \cong \mathtt{L}'_1$ we are done by induction since $\mathtt{X}/\mathtt{X}_1 \cong \mathtt{X}/\mathtt{X}'_1$ has a shorter filtration with simples factors being either $\mathtt{L}_i$ for $i = 2, \dots, n$, or $\mathtt{L}'_j$ for $j = 2, \dots, n'$, and we can use the induction hypothesis to see that these simples agree up to reordering and isomorphisms. Otherwise, $\mathtt{X}_1 \not\cong \mathtt{X}'_1$ and Schur's lemma implies that $\mathtt{X}_1 \oplus \mathtt{X}'_1$ is a subobject of $\mathtt{X}$ and we can consider $\mathtt{Y} = \mathtt{X}/(\mathtt{X}_1 \oplus \mathtt{X}'_1)$. It is easy to see that $\mathtt{Y}$ has a filtration by simples, say with simple factors $\mathtt{L}^{\mathtt{Y}}_k$, for $k = 1, \dots, r < n$. We then observe that:

- $\mathtt{X}/\mathtt{X}_1$ has a filtration with simple factors $\mathtt{X}'_1$, $\mathtt{L}^{\mathtt{Y}}_k$ for $k = 1, \dots, r$, but also one with the original simple factors $\mathtt{L}_i$ except $\mathtt{X}_1 \cong \mathtt{L}_1$.
- $\mathtt{X}/\mathtt{X}'_1$ has a filtration with simple factors $\mathtt{X}_1$, $\mathtt{L}^{\mathtt{Y}}_k$ for $k = 1, \dots, r$, but also one with the original simple factors $\mathtt{L}'_j$ except $\mathtt{X}'_1 \cong \mathtt{L}'_1$.

By the induction assumption, this means that the collection of simples $\mathtt{X}_1$, $\mathtt{X}'_1$, and $\mathtt{L}_i^{\mathtt{Y}}$ for $i = 1, ..., r$ coincides (up to reordering and isomorphisms) on the one hand with $\mathtt{L}_i$ for $i = 1, ..., n$ and, on the other hand, with $\mathtt{L}'_j$ for $j = 1, ..., n'$. □

6K. **Finiteness assumptions.** For what will follow, we need and want to go to the finite dimensional world:

**Definition 6K.1.** Let $\mathbf{fdVec}_{\mathbb{Z}} \subset \mathbf{Vec}_{\mathbb{Z}}$ be the full subcategory of ***torsion free abelian groups of finite rank***. ◇

Without harm we can think of Definition 6K.1 as being the $\mathbb{Z}$ linear version of $\mathbf{fdVec}_{\Bbbk} \subset \mathbf{Vec}_{\Bbbk}$, and the “fd” refers to finite dimensional: always having an underlying field in mind, we say “finite dimensional” instead of the mouthful “torsion free of finite rank”.

**Definition 6K.2.** Let $\mathbf{C} \in \mathbf{Cat}_{\oplus}$, and assume that $\mathrm{Hom}_{\mathbf{C}}(\mathtt{X}, \mathtt{Y}) \in \mathbf{fdVec}_{\mathbb{Z}}$ for all $\mathtt{X}, \mathtt{Y} \in \mathbf{C}$.

(i) If $\mathbf{C}$ is abelian and any $\mathtt{X} \in \mathbf{C}$ is of finite length, then we call $\mathbf{C}$ ***locally (abelian) finite***.

(ii) If all $\mathtt{X} \in \mathbf{C}$ have a decomposition as in (6J-17) satisfying the Krull–Schmidt theorem, meaning Theorem 6J.19.(ii), then we say that $\mathbf{C}$ is ***locally additively finite***.

Both are examples of ***finiteness assumptions***. ◇

**Example 6K.3.** Here are some prototypical examples:

(a) Not all objects in $\mathbf{Vec}_{\Bbbk}$ have finite length and hom spaces are not finite dimensional, thus $\mathbf{Vec}_{\Bbbk}$ is not locally finite.

(b) The full subcategory $\mathbf{fdVec}_{\Bbbk} \subset \mathbf{Vec}_{\Bbbk}$ is locally finite.

(c) For any group (it may be infinite) G and any $\omega$, the category $\mathbf{Vec}^{\omega}_{\Bbbk\oplus}(\mathrm{G})$ is locally finite, because in the additive closure we only allow finite direct sums.

The difference between (a) and (b) justifies the nomenclature “finiteness assumptions”. ◇

Note that we have

$$\text{locally additively finite} \Leftarrow \text{locally (abelian) finite},$$
$$\text{locally additively finite} \not\Rightarrow \text{locally (abelian) finite},$$

the latter being justified by Example 6K.6. Before we can state it, we need the analog of Definition 6H.4 in this finite setting:

**Definition 6K.4.** An algebra A in $\mathbf{fdVec}_{\mathbb{Z}}$ is called a ***finite dimensional algebra***. The ***category of finite dimensional right*** A ***modules*** for such an algebra is defined to be $\mathbf{fdMod}(\mathrm{A}) = \mathbf{Mod}_{\mathbf{fdVec}_{\mathbb{Z}}}(\mathtt{A})$. We also have full subcategories $\mathbf{fdProj}(\mathrm{A}) \subset \mathbf{fdMod}(\mathrm{A})$ and $\mathbf{fdInj}(\mathrm{A}) \subset \mathbf{fdMod}(\mathrm{A})$ of ***finite dimensional projectives*** and ***finite dimensional injectives***, respectively. ◇

**Definition 6K.5.** For any algebra $\mathtt{A} \in \mathbf{Vec}_{\mathbb{Z}}$ let $\mathbf{fdMod}(\mathrm{A}) = \mathbf{Mod}_{\mathbf{fdVec}_{\mathbb{Z}}}(\mathtt{A}) \subset \mathbf{Mod}_{\mathbf{Vec}_{\mathbb{Z}}}(\mathtt{A})$ denote the corresponding full subcategory of ***finite dimensional modules***. Similarly for finite dimensional projective and injective modules. ◇

**Example 6K.6.** Let us come back to Example 6J.5. The $\mathbb{C}, \mathtt{A} \in \mathbf{fdMod}(\mathrm{A})$, and $\mathbf{fdMod}(\mathrm{A})$ is locally finite. However, $\mathbb{C}$ is neither projective nor injective. Hence, $\mathtt{A}$ does not have any composition series in terms of projectives or injectives. Thus, neither $\mathbf{fdProj}(\mathrm{A})$ nor $\mathbf{fdInj}(\mathrm{A})$ are locally finite, but one can show that both are locally additively finite. ◇

The following are the abelian categories that we will use most of the time.

**Definition 6K.7.** A category $\mathbf{C} \in \mathbf{Cat}_{A}$ is called ***(abelian) finite*** if:

- $\mathbf{C}$ is locally finite;
- the set $\mathrm{Si}(\mathbf{C})$ is finite;
- every simple $\mathtt{L} \in \mathbf{C}$ has a projective cover.

For any such $\mathbf{C}$ we have the full subcategories $\mathbf{fProj}(\mathbf{C})$ and $\mathbf{fInj}(\mathbf{C})$ of ***finite projective*** respectively ***finite injective objects***. We also have the ***category of finite abelian categories*** being the corresponding full subcategory $\mathbf{Cat}_{fA} \subset \mathbf{Cat}_{A}$. ◇

**Example 6K.8.** Back to Example 6K.3:

(a) The abelian category $\mathbf{fdVec}_{\Bbbk}$ is finite.

(b) For any group G and any $\omega$, the abelian category $\mathbf{Vec}^{\omega}_{\Bbbk\oplus}(\mathrm{G})$ is finite if and only if G is finite.

In fact, the slogan is that “Abelian finite = something supported in $\mathbf{fdVec}_{\Bbbk}$”. ◇

We have already seen the explicit description of finite abelian categories, see Theorem 6H.6.(ii). We now sketch a proof.

*Sketch of a proof of Theorem 6H.6.(ii).* Since $\mathrm{Si}(\mathbf{C})$ is finite by assumption, we can number the simples therein $\mathtt{L}_i$ for $i=1,...,n$. Also, by assumption, they have projective covers $\mathtt{P}_i=\mathtt{P}(\mathtt{L}_i)$. Take

$$\mathrm{A}=\mathrm{End}_{\mathbf{C}}\big(\textstyle\bigoplus_{i=1}^n\mathtt{P}_i\big),$$

with $\bigoplus_{i=1}^n\mathtt{P}_i$ usually called a ***projective generator***. Note that A is finite dimensional because the hom spaces are, again by assumption, finite dimensional. Also, $\mathbf{fdMod}(\mathrm{A})$ is finite abelian, by classical representation theory. It is then not hard to see that this is the category we need, *i.e.*

$$\mathrm{Hom}_{\mathbf{C}}\big(\textstyle\bigoplus_{i=1}^n\mathtt{P}_i,_\big)\colon\ \mathbf{C}\xrightarrow{\simeq_e}\mathbf{fdMod}(\mathrm{A})\ .$$

Note hereby that the hom functor is exact since $\bigoplus_{i=1}^n\mathtt{P}_i$ is projective. Finally, $\mathrm{Hom}_{\mathbf{C}}\big(\bigoplus_{i=1}^n\mathtt{P}_i,\mathtt{X}\big)$ for all $\mathtt{X}\in\mathbf{C}$ is a right A module via precomposition. □

*Remark* 6K.9. Let $\mathbf{C}\subset\mathbf{Cat}_{fA}$. Note that the indecomposable projectives in $\mathbf{fProj}(\mathbf{C})$ are the projective covers of the simples in $\mathbf{C}$, while the indecomposable injectives in $\mathbf{fInj}(\mathbf{C})$ are their injective hulls. In particular,

$$\#\mathrm{Si}(\mathbf{C})=\#\mathrm{Pi}(\mathbf{C})=\#\mathrm{Ii}(\mathbf{C}).$$

In fact, in view of the Freyd–Mitchell theorem Theorem 6H.6.(ii), as a right A module

$$\mathrm{A}\cong\textstyle\bigoplus_{i=1}^n\dim(\mathtt{L}_i)\mathtt{P}_i. \tag{6K-10}$$

In turn, Equation 6K-10 gives a good justification for the slightly mysterious definition of projective objects: they are the elemental summands of the algebra presenting $\mathbf{C}$. ◇

**Example 6K.11.** Let us discuss the above proof in two examples.

(a) In $\mathbf{fdVec}_{\mathbb{k}}$ a projective generator is for example $\mathbb{k}$, and $\mathrm{End}_{\mathbf{fdVec}_{\mathbb{k}}}(\mathbb{k})\cong\mathbb{k}$, so that $\mathbf{fdVec}_{\mathbb{k}}\simeq_e\mathbf{fdMod}(\mathbb{k})$. However, $\mathbb{k}\oplus\mathbb{k}$ is also a projective generator and in this case one gets $\mathbf{fdVec}_{\mathbb{k}}\simeq_e\mathbf{fdMod}\big(\mathrm{Mat}_{2\times2}(\mathbb{k})\big)$.

(b) In Example 6J.5 the only simple object is $\mathbb{C}$ itself, and $\mathtt{P}(\mathbb{C})=\mathtt{A}$. Clearly the corresponding algebra $\mathrm{A}=\mathrm{End}_{\mathbf{fdMod}(\mathbb{C}[X]/(X^2))}(\mathtt{P}(\mathbb{C}))$ is isomorphic to $\mathbb{C}[X]/(X^2)$.

These are, of course, rather boring examples as the abelian categories already are of the form $\mathbf{fdMod}(\mathrm{A})$. However, what we want to stress is that the above proof is a generalization of the fact that every monoid M is isomorphic to the monoid $\mathrm{End}_{\bullet\mathrm{M}}(\mathrm{M})$, which we have already seen in the proof of Theorem 2I.5. ◇

Finally, recall the Grothendieck classes, see Definition 1H.15.

**Definition 6K.12.** Let $\mathbf{C}\in\mathbf{Cat}_{fA}$, and let $\mathbf{D}=\mathbf{C}$ or $\mathbf{D}\in\{\mathbf{fProj}(\mathbf{C}),\mathbf{fInj}(\mathbf{C})\}$.

(i) We endow $K_0(\mathbf{C})$ with the structure of an abelian group via

$$\big([\mathtt{Y}]=[\mathtt{X}]+[\mathtt{Z}]\big)\Leftrightarrow\big(\exists\ \mathtt{X}\xrightarrow{\mathrm{i}}\mathtt{Y}\xrightarrow{\mathrm{p}}\mathtt{Z}\ \ \text{SES}\ \big).$$

(ii) We endow $K_0(\mathbf{D})$ with the structure of an abelian group via

$$\big([\mathtt{Y}]=[\mathtt{X}]+[\mathtt{Z}]\big)\Leftrightarrow\big(\mathtt{Y}\cong\mathtt{X}\oplus\mathtt{Z}\big).$$

In order to distinguish the two structures we write $K_0^e(_)$ for the one involving SES and $K_0^{\oplus}(_)$ for the additive version. ◇

The following are easy and omitted.

**Lemma 6K.13.** *Let $\mathbf{C}\in\mathbf{Cat}_{fA}$, and let $\mathbf{D}\in\{\mathbf{Cat}_{fA},\mathbf{fProj}(\mathbf{C}),\mathbf{fInj}(\mathbf{C})\}$. Enumerate the simples in $\mathbf{C}$ or $\mathbf{D}$ by $\mathtt{L}_i$ for $i=1,...,n$, and let $\mathtt{P}_i$ and $\mathtt{I}_i$ for $i=1,...,n$ be their respective projective covers or injective hulls. Then:*

*(i) Definition 6K.12 endows $K_0^e(\mathbf{C})$ and $K_0^{\oplus}(\mathbf{D})$ with the structures of finite dimensional abelian groups.*

*(ii) The set $\mathrm{Si}(\mathbf{C})$ is a basis of $K_0^e(\mathbf{C})$. We have*

$$[\mathtt{X}]=\textstyle\sum_{i=1}^n[\mathtt{X}:\mathtt{L}_i]\cdot[\mathtt{L}_i]\in K_0^e(\mathbf{C}).$$

*(iii) The sets $\mathrm{Pi}(\mathbf{D})$ and $\mathrm{Ii}(\mathbf{D})$ are bases of $K_0^{\oplus}(\mathbf{D})$. We have*

$$[\mathtt{X}]=\textstyle\sum_{i=1}^n(\mathtt{X}:\mathtt{P}_i)\cdot[\mathtt{P}_i]\in K_0^{\oplus}(\mathbf{D}),$$

*and similarly with injectives.* □

**Lemma 6K.14.** *Let $\mathbf{C},\mathbf{C}'\in\mathbf{Cat}_{fA}$, and let $\mathbf{D},\mathbf{D}'\in\{\mathbf{Cat}_{fA},\mathbf{fProj}(\mathbf{C}),\mathbf{fInj}(\mathbf{C})\}$.*

*(i) Any functor* $\mathrm{F} \in \mathbf{Hom}_e(\mathbf{C}, \mathbf{C}')$ *induces a group homomorphism*

$$K_0^e(\mathrm{F})\colon K_0^e(\mathbf{C}) \to K_0^e(\mathbf{C}').$$

*Further, if* $\mathrm{F}$ *is an equivalence, then* $K_0^e(\mathrm{F})$ *is an isomorphism.*

*(ii) Any functor* $\mathrm{F} \in \mathbf{Hom}_{\oplus}(\mathbf{D}, \mathbf{D}')$ *induces a group homomorphism*

$$K_0^{\oplus}(\mathrm{F})\colon K_0^{\oplus}(\mathbf{D}) \to K_0^{\oplus}(\mathbf{D}').$$

*Further, if* $\mathrm{F}$ *is an equivalence, then* $K_0^{\oplus}(\mathrm{F})$ *is an isomorphism.* □

The final definition in this section, Definition 6K.15, is well-defined by Lemma 6K.13.

**Definition 6K.15.** Keeping the notation from Lemma 6K.13, the ***projective*** and ***injective Cartan matrices*** are the $n \times n$ matrices

$$C_p(\mathbf{C}) = \big([\mathtt{P}_i : \mathtt{L}_j]\big)_{i,j=1}^n, \quad C_i(\mathbf{C}) = \big([\mathtt{I}_i : \mathtt{L}_j]\big)_{i,j=1}^n.$$

(These are $n$-by-$n$ matrices with values in $\mathbb{Z}_{\geq 0}$.) ◇

Let us finish this subsection with a bigger example. Before that, let us recall:

*Remark* 6K.16. Let $p' \in \mathbb{Z}_{\geq 0}$ be a prime and $n \in \mathbb{Z}_{>0}$. Recall that there exist a unique, up to isomorphism, finite field $\mathbb{F}_q$ of order $q = (p')^n$ explicitly constructed by:

- If $n = 1$, then $\mathbb{F}_{p'} = \mathbb{Z}/p'\mathbb{Z}$;
- if $n > 1$, then $\mathbb{F}_q = \mathbb{F}_{p'}[X]/(X^q - X)$.

The algebraic closure of $\mathbb{F}_q$ is $\overline{\mathbb{F}}_q = \bigcup_{m\in\mathbb{Z}_{>0}} \mathbb{F}_{q^m}$. (Finite fields can not be algebraically closed by the folk argument: “If $\mathbb{F} = \{z_1, ..., z_r\}$, then $p(X) = 1 + \prod_{i=1}^r (X - z_i)$ has no root in $\mathbb{F}$.”)

Let further $m \in \mathbb{Z}_{>0}$ and consider the polynomial $p(X) = X^m - 1$. Then:

$$p(X) \text{ has } \gcd(m, q-1) \text{ roots in } \mathbb{F}_q. \tag{6K-17}$$

In particular, if $m = p$ is itself a prime, then there are primitive $m$th root of unity in $\overline{\mathbb{F}}_{p'}$ if and only if $p \neq p'$. Explicitly, if $p = 5$, $k \in \mathbb{Z}_{>0}$ and $p' = 5$ or $p' = 7$, then

$$\gcd(5, 5^k - 1) = 1, \quad \gcd(5, 7^k - 1) = \begin{cases} 5 & \text{if } k \equiv 0 \bmod (p-1), \\ 1 & \text{else.} \end{cases}$$

This is easy to generalize. ◇

**Example 6K.18.** Let us consider $\mathrm{A} = \overline{\mathbb{F}}_5[\mathbb{Z}/5\mathbb{Z}]$ and let $\mathbf{C} = \mathbf{fdMod}(\mathrm{A})$.

As already stated, see Remark 6K.9, the sets $\mathrm{Si}(\mathbf{C})$, $\mathrm{Pi}(\mathbf{C})$ and $\mathrm{Ii}(\mathbf{C})$ have all the same size in general, while $\mathrm{In}(\mathbf{C})$ might be bigger. Let us see this explicitly.

For A we can determine a module structure on a $\overline{\mathbb{F}_5}$ vector space by specifying the action of the generator $1 \in \mathbb{Z}/5\mathbb{Z}$ since $\mathbb{Z}/5\mathbb{Z} \cong \langle s \mid s^5 = 1\rangle$ and the isomorphism is given by sending 1 to $s$.

We define five modules

$$\mathtt{Z}_1 = \mathtt{L}_1\colon \begin{pmatrix} 1 \end{pmatrix}, \text{ is simple}$$

$$\mathtt{Z}_2\colon \begin{pmatrix} 1 & 0 \\ 1 & 1 \end{pmatrix}, \text{filtration } \mathtt{0} - \mathtt{L}_1 - \mathtt{Z}_2,$$

$$\mathtt{Z}_3\colon \begin{pmatrix} 1 & 0 & 0 \\ 1 & 1 & 0 \\ 0 & 1 & 1 \end{pmatrix}, \text{filtration } \mathtt{0} - \mathtt{L}_1 - \mathtt{L}_1 - \mathtt{Z}_3,$$

$$\mathtt{Z}_4\colon \begin{pmatrix} 1 & 0 & 0 & 0 \\ 1 & 1 & 0 & 0 \\ 0 & 1 & 1 & 0 \\ 0 & 0 & 1 & 1 \end{pmatrix}, \text{filtration } \mathtt{0} - \mathtt{L}_1 - \mathtt{L}_1 - \mathtt{L}_1 - \mathtt{Z}_4,$$

$$\mathtt{Z}_5 = \mathtt{P}_1\colon \begin{pmatrix} 1 & 0 & 0 & 0 & 0 \\ 1 & 1 & 0 & 0 & 0 \\ 0 & 1 & 1 & 0 & 0 \\ 0 & 0 & 1 & 1 & 0 \\ 0 & 0 & 0 & 1 & 1 \end{pmatrix}, \text{filtration } \mathtt{0} - \mathtt{L}_1 - \mathtt{L}_1 - \mathtt{L}_1 - \mathtt{L}_1 - \mathtt{Z}_5,$$

where we gave the action matrices of 1 and the filtrations by simples, where we give the successive simple quotients. In this case the characteristic 5 version of the Jordan theorem gives

$$\mathrm{Si}(\mathbf{C})=\{\mathtt{L}_1\},\quad \mathrm{Pi}(\mathbf{C})=\{\mathtt{P}_1\}=\mathrm{Ii}(\mathbf{C}),\quad \mathrm{In}(\mathbf{C})=\{\mathtt{Z}_1,\mathtt{Z}_2,\mathtt{Z}_3,\mathtt{Z}_4,\mathtt{Z}_5\},$$
$$K_0^e(\mathbf{C})\cong K_0^{\oplus}\big(\mathbf{fProj}(\mathbf{C})\big)\cong K_0^{\oplus}\big(\mathbf{fInj}(\mathbf{C})\big)\cong\mathbb{Z}.$$

However, the evident group homomorphism

$$K_0^{\oplus}\big(\mathbf{fProj}(\mathbf{C})\big)\to K_0^e(\mathbf{C}),\ [\mathtt{X}]\mapsto[\mathtt{X}]\ \rightsquigarrow\ \mathbb{Z}\to\mathbb{Z},\ 1\mapsto 5,$$

is not a group isomorphism since $[\mathtt{P}_1]=5[\mathtt{L}_1]$ in $K_0^e(\mathbf{C})$, and it corresponds, as indicated, to multiplication by 5. We also have $C_p(\mathbf{C})=(\,5\,)$. Similarly for the injectives. ◇

## 6L. **Exercises.**

*Exercise* 6L.1. Show that the two morphisms

$$\mathrm{e}_+=\tfrac{1}{2}\cdot\big(\,|\ \ |+\times\big),\quad \mathrm{e}_-=\tfrac{1}{2}\cdot\big(\,|\ \ |-\times\big),$$

are orthogonal idempotents in $\mathbf{Br}_{\mathbb{Q}\oplus}$, meaning that $\mathrm{e}_{\pm}^2=\mathrm{e}_{\pm}$ and $\mathrm{e}_{\pm}\mathrm{e}_{\mp}=0$. ◇

*Exercise* 6L.2. Prove Proposition 6F.9 and Lemma 6G.14. ◇

*Exercise* 6L.3. Describe (co)kernels, images, the epic-monic factorizations, simples, projective and injectives in $\mathbf{Vec}_{\Bbbk\oplus}(\mathrm{G})$. Moreover, find a presenting algebra A, *cf.* Theorem 6H.6. ◇

*Exercise* 6L.4. Prove Schur's lemma(s) Lemma 6J.9 and Lemma 6J.10, and find an example for $\Bbbk=\mathbb{Q}$ where (6J-11) does not hold. ◇

*Exercise* 6L.5. For $a,b,c\in\mathbb{C}$ let $\mathrm{A}=\mathbb{C}[X]/(X-a)(X-b)(X-c)$. Consider the cases (a) $a=b=c=0$, (b) $a=b=0,c=2$ and (c) $a=0,b=1,c=2$ and show (*e.g.* via the Chinese reminder theorem, *cf.* Remark 7B.6) that

$$\mathrm{A}\cong\begin{cases}\mathbb{C}[X]/(X^3) & \text{case (a)},\\ \mathbb{C}[X]/(X^2)\oplus\mathbb{C}[X]/(X-2) & \text{case (b)},\\ \mathbb{C}[X]/(X)\oplus\mathbb{C}[X]/(X-1)\oplus\mathbb{C}[X]/(X-2) & \text{case (c)}.\end{cases}$$

(What could the general statement for $\mathbb{C}[X]/\prod_{i=1}^{n}(X-a_i)$ be?) Then compute the Cartan matrix of $\mathbf{fdMod}(\mathrm{A})$ in the above cases. ◇

# 7. Fiat and tensor categories – enrich the concepts from before

Recall that the Grothendieck classes of an additive or abelian category have a commutative addition. So, in some sense, they categorify abelian groups. With a monoidal structure at hand, one additionally gets a multiplication. Thus, a natural question would be:

> What are suitable categorifications of rings or algebras?

Before we begin, it's important to remember that in category theory, morphisms are considered more significant than objects. This remains true for monoidal additive categories, but it takes on a slightly different nuance. The key difference is that, in ordinary categories, objects form a set, meaning that the primary concern is simply counting them. In the monoidal additive context, however, objects form rings, which means they carry more information than just their quantity.

Here is an example of a set $\{1,x,y,z\}$ with two (group) multiplication structures:

(7-1)

| $\cdot$ | 1 | $x$ | $y$ | $z$ |
|---|---|---|---|---|
| 1 | 1 | $x$ | $y$ | $z$ |
| $x$ | $x$ | $y$ | $z$ | 1 |
| $y$ | $y$ | $z$ | 1 | $x$ |
| $z$ | $z$ | 1 | $x$ | $y$ |

,

| $\cdot$ | 1 | $x$ | $y$ | $z$ |
|---|---|---|---|---|
| 1 | 1 | $x$ | $y$ | $z$ |
| $x$ | $x$ | 1 | $z$ | $y$ |
| $y$ | $y$ | $z$ | 1 | $x$ |
| $z$ | $z$ | $y$ | $x$ | 1 |

.

We leave it to the reader to determine the names of the above groups.

*Remark* 7.2. Counting group structures on a set $\{1,x_1,\dots,x_n\}$, of which we think of as a set of colors, is a game in which a) one needs to find a $(n+1)$-by-$(n+1)$ table such that every row and column contain $n+1$ different colors and b) count them up to permutation. For example, in Equation 7-1 every row and column contain four different colors, and these are the only two possibilities that satisfy this pattern. Counting groups is not easy (and quite similar to Sudoku, *cf.* Figure 17), but there are some patterns that can be observed.

For monoids the counting problem is much more involved and not much is known. For starters,

$$\begin{array}{c|cc} \cdot & 1 & x \\ \hline 1 & 1 & x \\ x & x & 1 \end{array}, \qquad \begin{array}{c|cc} + & 1 & x \\ \hline 1 & 1 & x \\ x & x & x \end{array}.$$

are two different monoids on $\{1, x\}$. Note that only one is a group (the left one as the right one does not have different colors in every row and column). The categorical counting problem "How many monoidal categories are there satisfying XYZ?" is akin to the monoid counting problem and very difficult. ◇

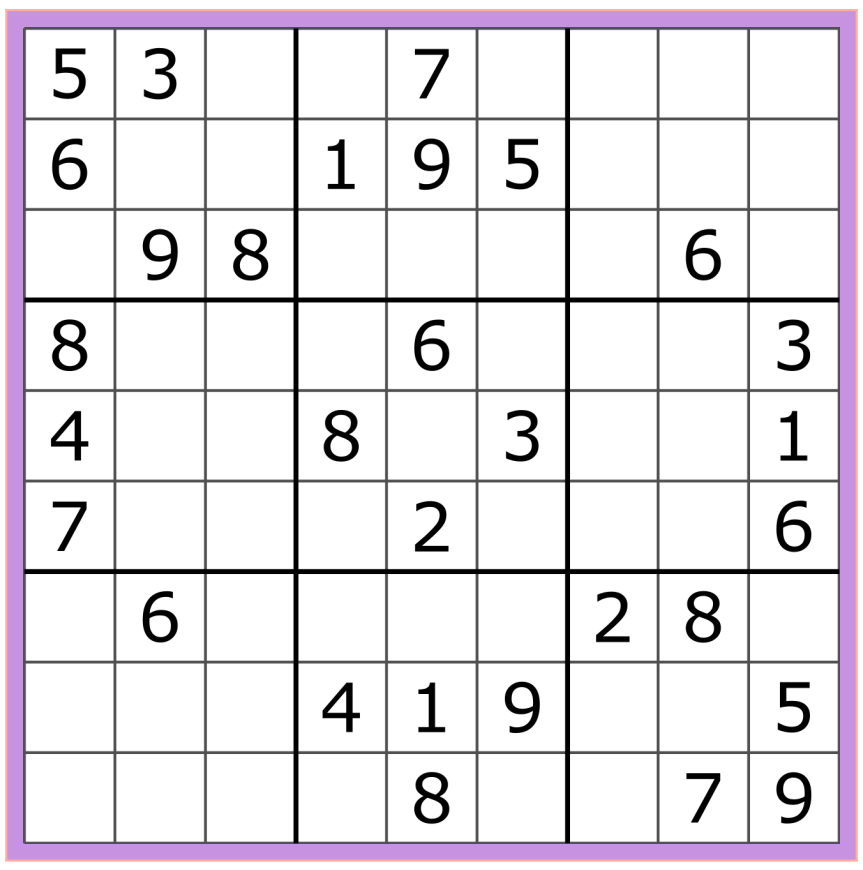

FIGURE 17. In standard Sudoku, the goal is to populate a 9-by-9 grid with digits such that each row, each column, and each of the nine 3-by-3 subgrids contains every digit from 1 to 9 exactly once. Listing all groups (up to isomorphism or permutation) is similar to listing all Sudoku puzzles.

Picture from https://en.wikipedia.org/wiki/Sudoku.

7A. **A word about conventions.** We will maintain the previous conventions and additionally incorporate:

*Convention* 7A.1. From now on we have categories with several structures, and we "stack" the notation; being careful with the hierarchy of the notions. For example, $\mathbf{C} \in \mathbf{RCat}_{\Bbbk A}$ means that $\mathbf{C}$ is a rigid (for which $\mathbf{C}$ needs to be monoidal) $\Bbbk$ linear abelian (for which $\mathbf{C}$ needs to be additive) category. ◇

*Convention* 7A.2. Note that "topological properties" of categories are usually written in front *e.g.* **RCat** means rigid categories, while "algebraic properties" are usually in subscripts, *e.g.* $\mathbf{Cat}_{lA}$ means locally finite abelian categories. ◇

*Convention* 7A.3. We tend to drop the "up to isomorphism" if no confusion can arise. For example, "has one simple" is to be read as "has one simple up to isomorphism". ◇

*Convention* 7A.4. We write $k \cdot \mathtt{X}$ for $\mathtt{X} \oplus \dots \oplus \mathtt{X}$ ($k$ summands). We also use the symbol $\mathtt{X} \Subset \mathtt{Y}$ for "$\mathtt{X}$ is isomorphic to a direct summand of $\mathtt{Y}$". ◇

7B. **"The philosophy of idempotents."** One of the most important operations are the ***idempotent*** operations. These operations can be applied multiple times without changing the result beyond the initial application. Such operations are everywhere, see Figure 18 for a real-world example, although often without us noticing. In mathematics, one should think of them as projections.

Before we can answer the main question of this section, we want to be able to take Grothendieck classes of a larger class of categories. To do this, here is some motivation.

First, let us go back to (6B-3), say for $\mathbf{Vec}_{\Bbbk}$. We write $\mathrm{e}_{\mathtt{X}} = \mathrm{i}_{\mathtt{X}}\mathrm{p}_{\mathtt{X}}$ and $\mathrm{e}_{\mathtt{Y}} = \mathrm{i}_{\mathtt{Y}}\mathrm{p}_{\mathtt{Y}}$. Let us also write $\mathtt{Z} = \mathtt{X} \oplus \mathtt{Y}$. There are now several crucial observations:

- We have $\mathrm{e}_{\mathtt{Y}} = \mathrm{id}_{\mathtt{Z}} - \mathrm{e}_{\mathtt{X}}$, and $\mathrm{Im}(\mathrm{e}_{\mathtt{X}}) \cong \mathtt{X}$ and $\mathrm{Im}(\mathrm{id}_{\mathtt{Z}} - \mathrm{e}_{\mathtt{X}}) \cong \mathtt{Y}$.
- We have idempotency, *i.e.*

(7B-1) $$\mathrm{e}_{\mathtt{X}}^2 = \mathrm{e}_{\mathtt{X}}, \quad (\mathrm{id}_{\mathtt{Z}} - \mathrm{e}_{\mathtt{X}})^2 = \mathrm{id}_{\mathtt{Z}}^2 - 2\mathrm{e}_{\mathtt{X}} + \mathrm{e}_{\mathtt{X}}^2 = \mathrm{id}_{\mathtt{Z}} - \mathrm{e}_{\mathtt{X}}.$$

The property in (7B-1) means that $\mathrm{e}_{\mathtt{X}}$ is an ***idempotent***, and so is $\mathrm{id}_{\mathtt{Z}} - \mathrm{e}_{\mathtt{X}}$. We also have

(7B-2) $$\mathrm{e}_{\mathtt{X}}(\mathrm{id}_{\mathtt{Z}} - \mathrm{e}_{\mathtt{X}}) = \mathrm{e}_{\mathtt{X}} - \mathrm{e}_{\mathtt{X}}^2 = 0 = \mathrm{e}_{\mathtt{X}} - \mathrm{e}_{\mathtt{X}}^2 = (\mathrm{id}_{\mathtt{Z}} - \mathrm{e}_{\mathtt{X}})\mathrm{e}_{\mathtt{X}},$$

(7B-3) $$\mathrm{e}_{\mathtt{X}} + (\mathrm{id}_{\mathtt{Z}} - \mathrm{e}_{\mathtt{X}}) = \mathrm{id}_{\mathtt{Z}},$$

with (7B-2) and (7B-3) being called ***orthogonality*** and ***completeness***, respectively.

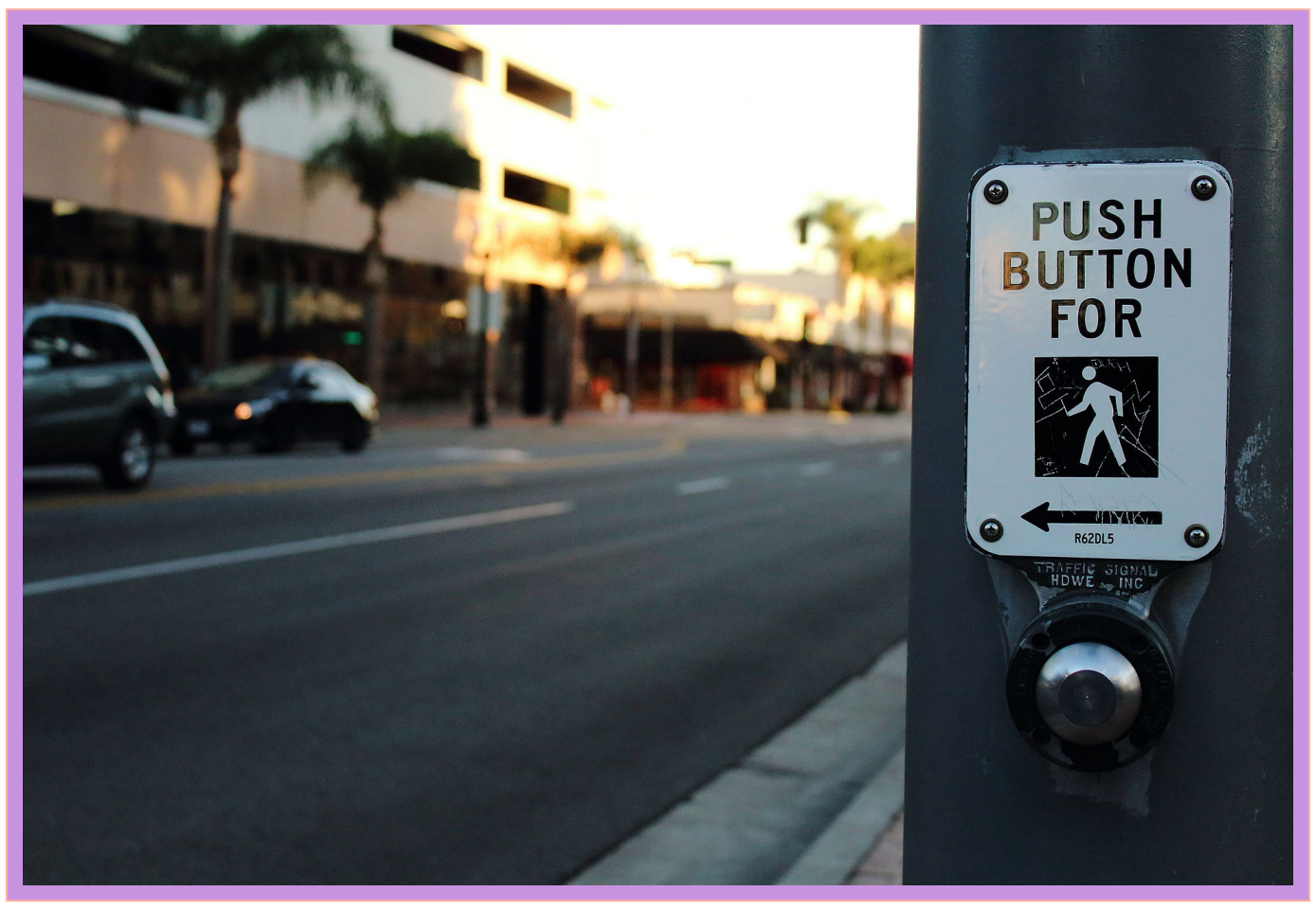

Figure 18. A crosswalk buttons is an idempotent: pressing it once, twice or a million times gives the same result.

Picture from https://en.wikipedia.org/wiki/Idempotence

- We calculate that we have a commuting diagram

(7B-4)

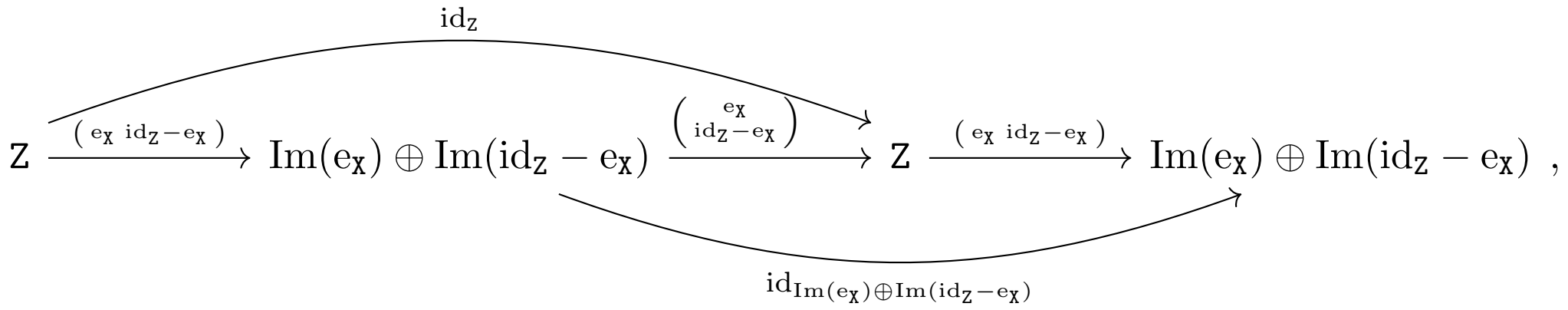

  which implies that $\mathtt{Z} \cong \mathrm{Im}(\mathrm{e}_{\mathtt{X}}) \oplus \mathrm{Im}(\mathrm{id}_{\mathtt{Z}} - \mathrm{e}_{\mathtt{X}})$. The isomorphisms, in the corresponding directions, are the two matrices in (7B-4).
- Note also that the above works both ways: Having a decomposition $\mathtt{Z} \cong \mathtt{X} \oplus \mathtt{Y}$ we get the idempotent $\mathrm{e}_{\mathtt{X}}$ satisfying all the above properties. Conversely, having an idempotent $\mathrm{e}\colon \mathtt{Z} \to \mathtt{Z}$ we get a diagram as in (7B-4).
- The used algebra above is very basic and only uses the existence of images and no other specific properties of $\mathbf{Vec}_{\Bbbk}$.

In $\mathbf{Mat}_{\Bbbk}$ an explicit example is

$$\mathtt{X} = 2, \quad \mathtt{Y} = 1, \quad \mathtt{Z} = 3, \quad \mathrm{e}_{\mathtt{X}} = \begin{pmatrix} 1 & 0 & 0 \\ 0 & 1 & 0 \\ 0 & 0 & 0 \end{pmatrix}, \quad \mathrm{id}_{\mathtt{Z}} - \mathrm{e}_{\mathtt{X}} = \begin{pmatrix} 0 & 0 & 0 \\ 0 & 0 & 0 \\ 0 & 0 & 1 \end{pmatrix}.$$

All of this together is called “The philosophy of idempotents.”, *i.e.*:

Idempotents decompose objects into direct sums.

This is, of course, most useful if the object one might care about does not come directly as $\mathtt{X} \oplus \mathtt{Y}$, but rather in some disguise. Here we do not want to take the ***trivial idempotents*** 0 and $\mathrm{id}_{\mathtt{Z}}$, and idempotents not of this form are called ***non-trivial***.

**Example 7B.5.** Let $\mathrm{A} = \mathbb{C}[X]/(X^2)$ and $\mathrm{B} = \mathbb{C}[X]/(X^2-1)$. We claim that these are quite different algebras in the following sense. An element in either A or B is of the form $a + bX = a \cdot 1 + b \cdot X$, where $a, b \in \mathbb{C}$. We calculate that

$$(a+bX)^2 = a^2 + 2abX + b^2X^2 = \begin{cases} a^2 + 2abX & \text{in A,} \\ a^2 + b^2 + 2abX & \text{in B.} \end{cases}$$

Thus, trying to solve the idempotency equation $(a+bX)^2 = a + bX$ for A gives only the trivial solutions $a = b = 0$ and $a = 1$, $b = 0$, and hence there is no non-trivial idempotent in A. In contrast, in B we get two non-trivial solutions

$$a = b = \tfrac{1}{2} \text{ giving } e_+ = \tfrac{1}{2}(1+X), \qquad a = -b = \tfrac{1}{2} \text{ giving } e_- = \tfrac{1}{2}(1-X),$$

which satisfy (7B-1), (7B-2) and (7B-3). Thus, as $\mathbb{C}$ algebras, we get:

$$\mathbb{C}[X]/(X^2-1) \cong \mathrm{Im}(e_+) \oplus \mathrm{Im}(e_-) \cong \mathbb{C}[X]/(1+X) \oplus \mathbb{C}[X]/(1-X) \cong \mathbb{C} \oplus \mathbb{C},$$

while $\mathbb{C}[X]/(X^2)$ does not decompose further. ◇

*Remark* 7B.6. The isomorphism in Example 7B.5 that decomposes $\mathbb{C}[X]/(X^2-1)$ into two copies of the complex numbers is a special case of the ***Chinese remainder theorem***. This works as follows. For a commutative ring $R$ and two coprime ideals $I$ and $J$ with $I \cap J = 0$ we have an isomorphism

$$f\colon R \xrightarrow{\cong} R/I \times R/J, r \mapsto (r \bmod I, r \bmod J).$$

Moreover, $e_1 = f^{-1}\big((1,0)\big)$ and $e_2 = f^{-1}\big((0,1)\big)$ satisfy (7B-1), (7B-2) and (7B-3). ◇

Group-type-things, *e.g.* groups (but not group algebras) do not have any non-trivial idempotents, as justified by the proof of the following lemma.

**Lemma 7B.7.** *For* $\mathbf{C} \in \mathbf{Cat}_{\mathbb{S}}$*, an idempotent* $\mathrm{e} \in \mathbf{C}$ *is invertible if and only if* $\mathrm{e} = \mathrm{id}$.

*Proof.* An all-time favorite trivial calculation is:

$$(\mathrm{e}^2 = \mathrm{e}) \Rightarrow (\mathrm{e}^{-1}\mathrm{e}^2 = \mathrm{e}^{-1}\mathrm{e}) \Rightarrow (\mathrm{e}^{-1}\mathrm{e}\mathrm{e} = \mathrm{e}^{-1}\mathrm{e}) \Rightarrow (\mathrm{e} = \mathrm{id}),$$

where one assumes that e is invertible. □

For $\mathbf{C} \in \mathbf{Cat}_{\mathbb{S}}$, we say $\mathrm{e} \in \mathbf{C}$ is a ***pseudo-idempotent with eigenvalue*** $\lambda \in \mathbb{S}$ if

$$\mathrm{e}^2 = \lambda\mathrm{e}.$$

**Lemma 7B.8.** *If* $\lambda$ *is invertible in* $\mathbb{S}$*, then* $\frac{1}{\lambda}\mathrm{e}$ *is an idempotent.*

*Proof.* Another all-time favorite trivial calculation: $(\frac{1}{\lambda}\mathrm{e})^2 = \frac{1}{\lambda^2}\mathrm{e}^2 = \frac{1}{\lambda^2}\lambda\mathrm{e} = \frac{1}{\lambda}\mathrm{e}$. □

**Example 7B.9.** Take $G = \mathbb{Z}/3\mathbb{Z} = \langle a | a^3 = 1\rangle$. Then $1 + a + a^2$ is a pseudo-idempotent with eigenvalue 3 in $\mathbb{Z}G$. We can make it an idempotent if we go to a field where 3 is invertible, for example $\mathbb{C}$. Moreover, for $\mathtt{i}$ a square root of $-1$, let $\zeta = \exp(2\pi\mathtt{i}/3) \in \mathbb{C}$ (a complex primitive third root of unity), then

$$e_0 = 1 + a + a^2, \quad e_1 = 1 + \zeta^1 a + \zeta^2 a^2, \quad e_2 = 1 + \zeta^2 a + \zeta^1 a^2,$$

are three pseudo-idempotent with eigenvalue 3, which satisfy (7B-2) and (7B-3) upon division by 3. We get

$$\mathbb{C}G \cong \mathbb{C} \oplus \mathbb{C} \oplus \mathbb{C},$$

by the philosophy of idempotents. ◇

Generalizing Example 7B.9, if $G$ is a finite group, then $\sum_{g\in G} g$ is a pseudo-idempotent with eigenvalue $|G|$. Hence, for all fields where $|G|$ is invertible, we get an idempotent $\frac{1}{|G|}\sum_{g\in G} g$. This gives, for example, $\mathbb{C}G \cong \mathbb{C} \oplus R$ for some "rest" $R$.

## 7C. The idempotent closure.

The above suggests that while images are desirable, they are only necessary for idempotents. So we will artificially add images associated to idempotents in the following way.

**Definition 7C.1.** The ***idempotent closure*** of $\mathbf{C} \in \mathbf{Cat}$, denoted by $\mathbf{C}_{\underline{\in}}$, is the category with

- objects being pairs
$$\mathrm{Ob}(\mathbf{C}_{\underline{\in}}) = \big\{(\mathtt{X}, \mathrm{e}) \mid \mathtt{X} \in \mathbf{C}, \mathrm{e}\colon \mathtt{X} \to \mathtt{X} \text{ idempotent}\big\};$$
- morphisms being $\mathrm{f}\colon (\mathtt{X}, \mathrm{e}) \to (\mathtt{Y}, \mathrm{e}')$ with $(\mathrm{f}\colon \mathtt{X} \to \mathtt{Y}) \in \mathbf{C}$ such that we have a commuting diagram
$$\begin{array}{ccc} \mathtt{Y} & \xrightarrow{\mathrm{e}'} & \mathtt{Y} \\ {\scriptstyle \mathrm{f}}\uparrow & & \uparrow{\scriptstyle \mathrm{f}} \\ \mathtt{X} & \xrightarrow[\mathrm{e}]{} & \mathtt{X} \end{array};$$
- the identities are $\mathrm{id}_{(\mathtt{X},\mathrm{e})} = \mathrm{e}$;
- composition is composition in $\mathbf{C}$.

(Note that being an idempotent makes sense in any category.) ◇

Good news, this works well:

**Proposition 7C.2.** *Let* $\mathbf{C} \in \mathbf{Cat}$*. Then we have:*

*(i)* $\mathbf{C}_{\underline{\in}}$ *is a category.*

*(ii) There exists a well-defined fully faithful functor*
$$\mathrm{K}\colon \mathbf{C} \to \mathbf{C}_{\underline{\in}}, \ \mathtt{X} \mapsto (\mathtt{X}, \mathrm{id}_{\mathtt{X}}), \mathrm{f} \mapsto \mathrm{f}.$$

*(iii) If* $\mathrm{e} \in \mathbf{C}$ *is an idempotent, then it has an image* $\mathrm{Im}(\mathrm{e}) \cong (\mathtt{X}, \mathrm{e})$ *in* $\mathbf{C}_{\underline{\in}}$.

*(iv) We have* $\mathbf{C}_{\underline{\in}} \simeq (\mathbf{C}_{\underline{\in}})_{\underline{\in}}$ *and one can find an equivalence preserving images of idempotents.*

(v) *If* $\mathbf{C} \in \mathbf{Cat}$ *is* $\mathbb{S}$ *linear (or additive or monoidal or rigid or pivotal or braided etc.), then so is* $\mathbf{C}_{\Subset}$ *with its structure induced from* $\mathbf{C}$.

(The notation K comes from the alternative name of $\mathbf{C}_{\Subset}$: it is sometimes called ***Karoubi completion***.)

*Proof.* We only prove (iii), the rest is Exercise 7J.1. To see that $\mathrm{Im}(e) \cong (\mathtt{X}, e)$ we just observe that

$$\begin{array}{ccc} (\mathtt{X}, e) & & \\ {\scriptstyle e}\downarrow & \searrow^{e} & \\ (\mathtt{X}, e) & \underset{e}{\hookrightarrow} & (\mathtt{X}, e) \end{array}$$

commutes, since e is an idempotent and $e = \mathrm{id}_{(\mathtt{X},e)}$. □

We use Proposition 7C.2.(iii) to write Im(e) for the objects of $\mathbf{C}_{\Subset}$, and we also write $\mathtt{X}$ instead of $(\mathtt{X}, \mathrm{id}_{\mathtt{X}})$. Moreover, we call a category $\mathbf{C}$ ***idempotent complete*** if $\mathbf{C} \simeq \mathbf{C}_{\Subset}$. Let us also write $\mathbf{C}_{\Bbbk\Subset} = (\mathbf{C}_{\Bbbk})_{\Subset}$ *etc.*

**Example 7C.3.** The idempotent closures is a technology for non-module-like categories:

(a) We have $\mathbf{Vec}_{\Bbbk} \simeq_{e} \mathbf{Vec}_{\Bbbk\Subset}$, which is thus idempotent complete. The same is true for any abelian category, or any category having images.

(b) Categories of the form $\mathbf{fdMod}(\mathrm{A})$, $\mathbf{fProj}\big(\mathbf{Mod}(\mathrm{A})\big)$ or $\mathbf{fInj}\big(\mathbf{Mod}(\mathrm{A})\big)$ for a finite dimensional algebra A are idempotent complete.

(c) Categories of the form $\mathbf{Vec}^{\omega}_{\Bbbk\oplus}(\mathrm{G})$ are idempotent complete.

On the other hand, set theoretic categories are rarely idempotent complete. ◇

**Example 7C.4.** Diagrammatic categories are almost never idempotent complete and we prefer to think of their idempotent completion as ***coloring with idempotents***: Recall the category **Sym**, see Example 3D.2, and let us make it $\mathbb{Z}[\frac{1}{2}]$ linear additive. Then are orthogonal and complete idempotents in $\mathbf{Sym}_{\mathbb{Z}[\frac{1}{2}]\oplus}$, see also Exercise 6L.1. Thus,

$$\bullet\bullet \cong \mathrm{Im}(e_{+}) \oplus \mathrm{Im}(e_{-}), \quad (\text{in } \mathbf{Sym}_{\mathbb{Z}[\frac{1}{2}]\oplus\Subset}). \tag{7C-5}$$

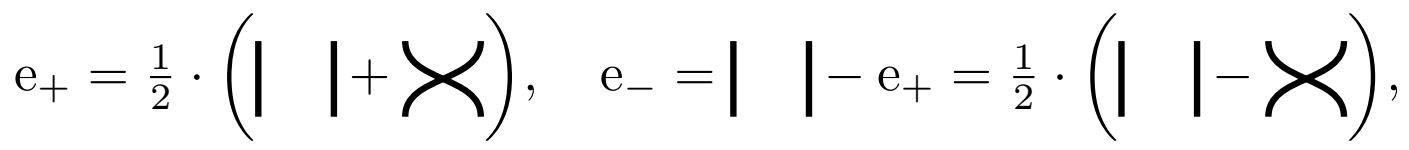

We can think of this as coloring the diagrams with the idempotents $e_{+}$ and $e_{-}$, illustrated say red (and dashed) and green, and (7C-5) becomes

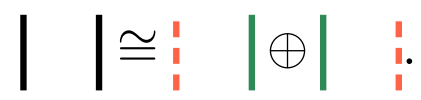

One can show that $\mathbf{Sym}_{\mathbb{Z}[\frac{1}{2}]\oplus\Subset}$ is idempotent complete, but non-abelian. ◇

*Remark* 7C.6. Example 7C.3.(a) and Example 7C.4 show that, for additive categories,

$$\text{idempotent complete} \Leftarrow \text{abelian},$$
$$\text{idempotent complete} \not\Rightarrow \text{abelian}.$$

In fact, nothing is special about $\frac{1}{2}$ in Example 7C.4 and the same results hold over $\mathbb{Z}$, but the idempotent decomposition of $\bullet^{2}$ is different. ◇

Here is the analog of Proposition 6F.10.

**Proposition 7C.7.** *Let* $\mathrm{F} \in \mathbf{Hom}(\mathbf{C}, \mathbf{D})$. *Then there exists an unique* $\mathrm{F}_{\Subset} \in \mathbf{Hom}(\mathbf{C}_{\Subset}, \mathbf{D}_{\Subset})$ *such that we have a commuting diagram*

$$\begin{array}{ccc} \mathbf{C}_{\Subset} & \overset{\exists!}{\underset{\mathrm{F}_{\Subset}}{\dashrightarrow}} & \mathbf{D}_{\Subset} \\ {\scriptstyle \mathrm{K}}\uparrow & & \uparrow{\scriptstyle \mathrm{K}} \\ \mathbf{C} & \underset{\mathrm{F}}{\longrightarrow} & \mathbf{D}. \end{array}$$

*Proof.* The functor $\mathrm{F}_{\Subset}$ is defined by

$$\mathrm{F}_{\Subset} \colon \mathbf{C}_{\Subset} \to \mathbf{D}_{\Subset}\ (\mathtt{X}, e) \mapsto \big(\mathrm{F}(\mathtt{X}), \mathrm{F}(e)\big), f \mapsto \mathrm{F}(f),$$

which satisfies all required properties. □

As for $\mathbb{S}$ linear extensions and additive closures, it can be checked that "all properties that we care about behave nicely with idempotent closures", *e.g.* if $\mathrm{F} \in \mathbf{Hom}_{\otimes}(\mathbf{C}, \mathbf{D})$, then so is $\mathrm{F}_{\Subset} \in \mathbf{Hom}_{\otimes}(\mathbf{C}_{\Subset}, \mathbf{D}_{\Subset})$. In particular, we basically get direct sum decompositions for free:

**Lemma 7C.8.** *Let* $\mathbf{C} \in \mathbf{Cat}_{\oplus}$. *Then*

$$\mathtt{X} \cong \mathrm{Im}(e) \oplus \mathrm{Im}(\mathrm{id}_{\mathtt{X}} - e)$$

*holds in* $\mathbf{C}_{\Subset}$.

*Proof.* The proof is *verbatim* (7B-4). □

**Example 7C.9.** Being idempotent complete is a property, and thus, additive functors are the correct maps between additive idempotent complete categories. Hence, we have the ***category of additive idempotent complete categories*** $\mathbf{Cat}_{\oplus\Subset}$. ◇

7D. **Tensor and fiat categories.** Before going on, we need the additive analog of Definition 6K.7:

**Definition 7D.1.** A category $\mathbf{C} \in \mathbf{Cat}_{\oplus}$ is called ***(additively) finite*** if

- $\mathbf{C}$ is locally additively finite;
- the set $\mathrm{In}(\mathbf{C})$ is finite.

We have the ***category of additively finite categories*** being the corresponding full subcategory $\mathbf{Cat}_{f\oplus} \subset \mathbf{Cat}_{\oplus}$. ◇

The following are additive or abelian categorifications of $\mathbb{S}$ algebras, as we shall see below. (Rigidity is strictly speaking not necessary to categorify algebras, but it makes life easier.) Here "w=weakly", "m=multi" and "l=locally". The study of fiat and finitary categories goes back to [MM11], while tensor categories are around for a long time.

**Definition 7D.2.** A category $\mathbf{C} \in \mathbf{Cat}$ is called ***wml fiat*** (over $\mathbb{S}$) if

- $\mathbf{C} \in \mathbf{RCat}_{\mathbb{S}\oplus\Subset}$;
- $\mathbf{C}$ is locally additively finite in the sense of Definition 6K.2.(b);
- The bifunctor $\otimes$ is $\mathbb{S}$ bilinear.

If additionally

- $\mathbf{C} \in \mathbf{PCat}$, then we drop the "weakly";
- $\mathrm{End}_{\mathbf{C}}(\mathbb{1}) \cong \mathbb{S}$, then we drop the "multi";
- $\mathbf{C}$ is finite in the sense of Definition 7D.1, then we drop the "locally".

(In this definition, everything is additive.) ◇

**Definition 7D.3.** A category $\mathbf{C} \in \mathbf{Cat}$ is called a ***wml tensor category*** (over $\mathbb{S}$) if

- $\mathbf{C} \in \mathbf{RCat}_{\mathbb{S}A}$;
- $\mathbf{C}$ is locally finite in the sense of Definition 6K.2.(b);
- The bifunctor $\otimes$ is $\mathbb{S}$ bilinear.

If additionally

- $\mathbf{C} \in \mathbf{PCat}$, then we drop the "weakly";
- $\mathrm{End}_{\mathbf{C}}(\mathbb{1}) \cong \mathbb{S}$, then we drop the "multi";
- $\mathbf{C}$ is finite in the sense of Definition 6K.7, then we drop the "locally".

(In this definition, everything is abelian.) ◇

Definition 7D.2 and Definition 7D.3 have quite the laundry list of adjectives. The fact to remember is that "fiat" should be thought of as an additive categorification of an algebra, while the "tensor" refers to an abelian categorification. The various adjectives then correspond to finiteness conditions, and other properties algebras can have.

*Remark* 7D.4. The above terminology is not linearly ordered. For example, a wm fiat category cannot be directly compared to an l tensor category. However, we clearly have *e.g.* that w fiat is a stronger notion that wl fiat. ◇

**Example 7D.5.** We have already seen plenty of examples:

(a) The category $\mathbf{fdVec}_{\mathbb{k}}$ is a tensor category. More generally, categories of the form $\mathbf{Vec}^{\omega}_{\mathbb{k}\oplus}(\mathrm{G})$ are l tensor categories, where we can drop the l if and only if G is finite. In both cases, we could have written fiat instead of tensor.

(b) Closures of diagrammatic categories such as $\mathbf{Br}_{\mathbb{k}\Subset}$ can be made into l fiat categories, see Section 7E, but not tensor categories (at least without extra work).

More generally, the "tensor" can be read as vector-space-like, while "fiat" is, in some sense, tailored to work for diagrammatic categories. ◇

**Example 7D.6.** All of these should be thought of as generalizing $\mathbf{Mod}(\mathrm{A})$ for certain "nice" algebras A such as group algebras of finite groups. "Nice" means roughly:

- Bialgebras, see Section 5J, endow $\mathbf{Mod}(\mathrm{A})$ with a monoidal structure;
- Hopf algebras, see Section 5K, provide the duality;
- finite dimensional algebras and finite dimensional modules provide the various finiteness conditions, see Section 7H.

In fact, for $\mathrm{A} = \Bbbk$ the category $\mathbf{fdMod}(\mathrm{A})$ is the fiat and tensor category $\mathbf{fdVec}_{\Bbbk}$. ◇

In other words, Example 7D.6 says that fiat and tensor categories should be thought of as generalizations of $\mathbf{fdMod}\big(\Bbbk[\mathrm{G}]\big)$, where G is a finite group. Explicitly:

**Example 7D.7.** Returning to Example 6K.18, the category $\mathbf{C} = \mathbf{fdMod}\big(\overline{\mathbb{F}}_5[\mathbb{Z}/5\mathbb{Z}]\big)$ is actually pivotal, which we will see completely explicitly in Equation 7H-7 below. Moreover, we have

$$\mathtt{Z}_1\mathtt{P}_1 \cong \mathtt{P}_1, \quad \mathtt{Z}_2\mathtt{P}_1 \cong 2\cdot\mathtt{P}_1, \quad \mathtt{Z}_3\mathtt{P}_1 \cong 3\cdot\mathtt{P}_1, \quad \mathtt{Z}_4\mathtt{P}_1 \cong 4\cdot\mathtt{P}_1, \quad \mathtt{Z}_5\mathtt{P}_1 \cong 5\cdot\mathtt{P}_1, \quad (\mathtt{P}_1)^{\star} \cong \mathtt{P}_1, \tag{7D-8}$$

so $\mathbf{C}' = \mathbf{fdProj}\big(\overline{\mathbb{F}}_5[\mathbb{Z}/5\mathbb{Z}]\big)$ is also pivotal (without monoidal unit). The category $\mathbf{C}$ is a tensor category with one simple and five indecomposables. In contrast, the category $\mathbf{C}'$ is a fiat category without monoidal unit but with a pseudo idempotent instead:

$$\mathtt{P}_1\mathtt{P}_1 \cong 5\cdot\mathtt{P}_1 \quad (\rightsquigarrow e^2 = 5e),$$

and with one indecomposable. Lastly, $\mathbf{C}$ itself is a fiat category, meaning we prioritize indecomposables over simples, and it contains five indecomposables. ◇

Note that we already know the correct functors between fiat respectively tensor categories: such functors should be $\mathbb{S}$ linear additive rigid respectively $\mathbb{S}$ linear exact rigid. Hence:

**Example 7D.9.** We have the ***category of fiat categories*** $\mathbf{Fiat}$, objects being fiat categories and morphisms being $\mathbb{S}$ linear additive rigid functors. We also have the ***category of tensor categories*** $\mathbf{Ten}$, objects being fiat categories and morphisms being $\mathbb{S}$ linear exact rigid functors. Finally, we also have various versions that add the adjectives "weakly", "multi" or "locally". ◇

For a wm fiat category $\mathbf{C}$, by definition, we know that the set $\mathrm{In}(\mathbf{C})$ is finite, so we can enumerate and denote the indecomposables by $\mathtt{Z}_i$ for $i = 1, ..., n$. Similarly, we let $\mathtt{L}_i$ for $i = 1, ..., n$ denote the simples of a wm tensor category.

**Lemma 7D.10.** *Let* $\mathbf{C} \in \mathbf{wmlFiat}$. *Then:*

*(i) The (bi)functors* $_ \otimes _$, $_^{\star}$ *and* ${}^{\star}_$ *are* $\mathbb{S}$ *linear additive.*

*(ii) If* $\mathbf{C} \in \mathbf{wlFiat}$ *is* $\Bbbk$ *linear, then* $\mathbb{1} \in \mathrm{In}(\mathbf{C})$.

*(iii) If* $\mathbf{C} \in \mathbf{wmFiat}$, *then the functor* $_^{\star}$ *induces a bijection*

$$_^{\star}\colon \mathrm{In}(\mathbf{C}) \xrightarrow{\cong} \mathrm{In}(\mathbf{C}). \tag{7D-11}$$

*Similarly for* ${}^{\star}_$.

*Proof.* (i). Since $\Bbbk$ linear implies additive, see Lemma 6E.8, the claim is immediate.

(ii). Note that $\mathrm{End}_{\mathbf{C}}(\mathbb{1}) \cong \Bbbk$ implies that $\mathbb{1}$ is indecomposable as its endomorphism $\Bbbk$ algebra does not have any non-trivial idempotents.

(iii). This follows since the dualities are $\mathbb{S}$ linear additive monoidal functors by (i), so the property of being indecomposable is preserved by them. Moreover, by Proposition 4D.9, they are invertible, thus, induce bijections. □

**Proposition 7D.12.** *Let* $\mathbf{C} \in \mathbf{wmfFiat}$ *or* $\mathbf{C} \in \mathbf{wmfTen}$ *and* $\mathtt{X} \in \mathbf{C}$. *Then* $(\mathtt{X} \otimes _), (_ \otimes \mathtt{X}) \in \mathbf{End}_e(\mathbf{C})$.

*Proof.* Since we have duals, we can use Theorem 4C.8 to see that both functors have right and left adjoints (in the sense of Example 4C.1). It is then not hard to see that such functors preserve the property of being a SES. □

Here is an interesting fact: projective and injective objects form a monoidal ideal (*cf.* (7D-8) for an example: tensoring with a projective gives a projective) in the following sense.

**Proposition 7D.13.** *Let* $\mathbf{C} \in \mathbf{wmfFiat}$ *or* $\mathbf{C} \in \mathbf{wmfTen}$. *Let further* $\mathtt{P} \in \mathbf{Proj}(\mathbf{C})$ *and* $\mathtt{X} \in \mathbf{C}$. *Then* $\mathtt{PX}, \mathtt{XP} \in \mathbf{Proj}(\mathbf{C})$. *Similarly for injective objects.*

*Proof.* Using Theorem 4C.8 we get *e.g.*

$$\mathrm{Hom}_{\mathbf{C}}(\mathtt{PX},\mathtt{Y}) \cong \mathrm{Hom}_{\mathbf{C}}\big(\mathtt{P},\mathtt{Y}({}^{\star}\mathtt{X})\big),$$

which shows, using Proposition 7D.12, that the hom functor for $\mathtt{PX}$ is exact. All other cases that we need to check follow by symmetry. ☐

**Example 7D.14.** Let $G$ be a finite group, and consider $\mathbf{C} = \mathbf{fdMod}\big(\Bbbk[\mathrm{G}]\big)$. Then Proposition 7D.13 recovers the fact that tensoring with a projective $\Bbbk[\mathrm{G}]$ module gives a projective $\Bbbk[\mathrm{G}]$ module.

For a finite monoid M this is not true; see *e.g.* [Ste16, Exercise 17.15] for an explicit example. And, in fact, $\mathbf{fdMod}\big(\Bbbk[\mathrm{M}]\big)$ is monoidal but not pivotal and hence not a tensor category. ◇

7E. **Semisimplicity.** Recall that the elements of, say, an abelian category are the simples, *cf.* Section 6J. The simplest compounds are:

**Definition 7E.1.** An object $\mathtt{X} \in \mathbf{C}$ with $\mathbf{C} \in \mathbf{Cat}_{\oplus\Subset}$ is called ***semisimple*** if

$$\mathtt{X} \cong \mathtt{L}_1 \oplus \dots \oplus \mathtt{L}_r, \quad \text{where } \mathtt{L}_i \in \mathrm{Si}(\mathbf{C}).$$

(Thus, "semisimple = direct sum of simples".) ◇

**Definition 7E.2.** A category $\mathbf{C} \in \mathbf{Cat}_{\oplus\Subset}$ is called ***semisimple*** if all of its objects are semisimple. ◇

**Example 7E.3.** Again, we already know several (non-)examples:

(a) The archetypical example is $\mathbf{fdVec}_{\Bbbk}$, since every finite dimensional vector space $\mathtt{X}$ is isomorphic to $r \cdot \Bbbk$ for some $r \in \mathbb{Z}_{\geq 0}$.

(b) Let G be a finite group. Clearly, all categories of the form $\mathbf{Vec}^{\omega}_{\Bbbk\oplus}(\mathrm{G})$ are semisimple since all objects are direct sums of simple objects, by definition.

(c) Non-examples are the categories $\mathbf{fdMod}(\mathrm{A})$ and $\mathbf{fdMod}\big(\overline{\mathbb{F}}_5[\mathbb{Z}/5\mathbb{Z}]\big)$ from Example 6J.5.(b) and Example 6K.18, respectively. In both cases there is only one simple and its projective cover is not simple, but indecomposable.

We will generalize (c) momentarily. ◇

**Proposition 7E.4.** *Let* A *be a two dimensional* $\mathbb{C}$ *algebra. Then* $\mathbf{fdMod}(\mathrm{A})$ *is semisimple if and only if* $\mathrm{A} \cong \mathbb{C} \oplus \mathbb{C}$.

*Proof.* Any such algebra is of the form $\mathrm{A} = \mathbb{C}[X]/(X-a)(X-b)$ for $a, b \in \mathbb{C}$. For $a \neq b$ the Chinese remainder theorem, as recalled in Remark 7B.6, shows that $\mathrm{A} \cong \mathbb{C} \oplus \mathbb{C}$ and hence $\mathbf{fdMod}(\mathrm{A})$ is semisimple. Similarly, for $a = b$ we get $\mathrm{A} \cong \mathbb{C}[X]/(X^2)$ and we are done by Example 6J.5.(b). ☐

Recalling Schur's lemma Lemma 6J.10, the following says that semisimple categories have controllable hom spaces:

**Proposition 7E.5.** *Let* $\mathbb{K}$ *be algebraically closed, and* $\mathbf{C} \in \mathbf{Cat}_{\oplus\Subset}$ *be locally additively finite. Then* $\mathbf{C}$ *is semisimple if and only if, for any* $\mathtt{X}, \mathtt{Y} \in \mathbf{C}$*, we have an isomorphism*

$$\textstyle\bigoplus_{\mathtt{L}\in\mathrm{Si}(\mathbf{C})} \big(\mathrm{Hom}_{\mathbf{C}}(\mathtt{X},\mathtt{L}) \times \mathrm{Hom}_{\mathbf{C}}(\mathtt{L},\mathtt{Y})\big) \xrightarrow{\cong} \mathrm{Hom}_{\mathbf{C}}(\mathtt{X},\mathtt{Y}), \quad (\mathrm{f},\mathrm{g}) \mapsto \mathrm{gf}. \tag{7E-6}$$

*Proof.* By decomposing $\mathtt{X}$ and $\mathtt{Y}$ into their simple components, one direction is a direct consequence of Lemma 6J.10. To see the converse, note that (7E-6) is equivalent to saying that the finite dimensional $\mathbb{K}$ vector spaces $\mathtt{Z} = \mathrm{Hom}_{\mathbf{C}}(\mathtt{X},\mathtt{L})$ and $\mathtt{Z}^{\star} = \mathrm{Hom}_{\mathbf{C}}(\mathtt{L},\mathtt{X})$ are duals. In particular, if these are non-zero, then the evaluation and coevaluation from Example 4C.3 provide the idempotent $\mathrm{coev}^{\mathtt{Z}}\mathrm{ev}_{\mathtt{Z}} \in \mathtt{Z}\mathtt{Z}^{\star} \cong \mathrm{End}_{\mathbf{C}}(\mathtt{X})$, showing that $\mathtt{L} \Subset \mathtt{X}$. Finally, since $\mathrm{End}_{\mathbf{C}}(\mathtt{X})$ is finite, there are only finitely many $\mathtt{L} \in \mathrm{Si}(\mathbf{C})$ for which $\mathrm{Hom}_{\mathbf{C}}(\mathtt{X},\mathtt{L})$ and $\mathrm{Hom}_{\mathbf{C}}(\mathtt{L},\mathtt{X})$ are non-zero. ☐

**Lemma 7E.7.** *Let* $\mathbf{C} \in \mathbf{Cat}_{\oplus\Subset}$ *be locally additively finite.*

*(i) If* $\mathbf{C}$ *is semisimple, then* $\mathrm{Si}(\mathbf{C}) \subset \mathrm{Pi}(\mathbf{C})$.

*(ii) If* $\mathrm{Si}(\mathbf{C}) \subset \mathrm{Pi}(\mathbf{C})$*, then* $\mathbf{C}$ *is semisimple.*

*(iii) If* $\mathbf{C}$ *is semisimple, then* $\mathrm{Si}(\mathbf{C}) \subset \mathrm{Ii}(\mathbf{C})$.

*(iv) If* $\mathrm{Si}(\mathbf{C}) \subset \mathrm{Ii}(\mathbf{C})$*, then* $\mathbf{C}$ *is semisimple.*

*Proof.* (i). Using Schur's lemma Lemma 6J.9 and semisimplicity, we can fill in the universal diagram as follows. We can only have an epic morphism $\mathrm{p}\colon \mathtt{L} \to \mathtt{X}$ from a simple $\mathtt{L}$ to a non-zero $\mathtt{X}$ if $\mathtt{L} \Subset \mathtt{X}$. Similarly, we can only have an epic morphism $\mathrm{f}\colon \mathtt{Y} \to \mathtt{X}$ from a non-zero $\mathtt{Y}$ to $\mathtt{X}$ if $\mathtt{X} \Subset \mathtt{Y}$. In particular, $\mathtt{L} \Subset \mathtt{Y}$, and we can define the required $\mathrm{u}\colon \mathtt{L} \to \mathtt{Y}$ by the universal property of the direct sum.

(ii). Assume that $\mathtt{X} \cong \mathtt{Z}_1 \oplus \dots \oplus \mathtt{Z}_n$ is a Krull–Schmidt decomposition and that $\mathrm{f}\colon \mathtt{X} \twoheadrightarrow \mathtt{L}$ is a epic morphism to a simple $\mathtt{L}$. Since $\mathtt{L}$ is projective we can use its universal property and get

$$\begin{array}{ccc} & & \mathtt{L} \\ & \exists!\ \swarrow \mathrm{u} & \downarrow \mathrm{id}_{\mathtt{L}} \\ \mathtt{X} & \xrightarrow[\mathrm{f}]{} & \mathtt{L} \end{array}.$$

The morphism fu is an idempotent since $\mathrm{fufu} = \mathrm{fid}_{\mathtt{L}}\mathrm{u}$, so $\mathtt{L} \Subset \mathtt{X}$. This implies that $\mathtt{L} \cong \mathtt{Z}_i$ for some $i$. Now proceed inductively.

(iii)+(iv). From (i) respectively (ii), by symmetry. □

In particular, we get:

**Lemma 7E.8.** *A locally additively finite category* $\mathbf{C} \in \mathbf{Cat}_{\oplus\Subset}$ *is semisimple if and only if* $\mathrm{Si}(\mathbf{C}) = \mathrm{Pi}(\mathbf{C})$ *if and only if* $\mathrm{Si}(\mathbf{C}) = \mathrm{Ii}(\mathbf{C})$. □

Clearly, we have the ***category of semisimple categories*** $\mathbf{Cat}_S \subset \mathbf{Cat}_{\oplus\Subset}$ being the corresponding full subcategory. More surprisingly:

**Theorem 7E.9.** *We have* $\mathbf{Cat}_S \subset \mathbf{Cat}_{lA}$.

In words, semisimple implies locally finite abelian.

*Proof.* For $\mathbf{C} \in \mathbf{Cat}_S$ take $\mathtt{P} = \bigoplus_{\mathtt{L}\in\mathrm{Si}(\mathbf{C})} \mathtt{L}$, which is a projective object by Lemma 7E.8. We consider

$$\mathrm{A} = \mathrm{End}_{\mathbf{C}}(\mathtt{P}).$$

Then we get an exact equivalence

$$\mathrm{Hom}_{\mathbf{C}}(\mathtt{P}, _)\colon \mathbf{C} \xrightarrow{\simeq_e} \mathbf{Mod}(\mathrm{A}).$$

Thus, $\mathbf{C}$ is abelian. That it is also locally finite is a direct consequence of the definition of semisimplicity and Schur's lemma Lemma 6J.10, since each object of $\mathbf{C}$ is a finite direct sum of simples. □

**Example 7E.10.** Example 7E.3 and Theorem 7E.9 immediately imply that $\mathbf{Vec}^{\omega}_{\Bbbk\oplus}(\mathrm{G})$ are all abelian (for G being a finite group). ◇

**Definition 7E.11.** An algebra $\mathtt{A} \in \mathbf{fdVec}_{\Bbbk}$ is called ***semisimple*** if $\mathbf{fdMod}(\mathtt{A})$ is semisimple. ◇

**Example 7E.12.** The classical ***Artin–Wedderburn theorem***, see e.g. [Ben98, Theorem 1.3], says that a $\Bbbk$ algebra $\mathtt{A} \in \mathbf{fdVec}_{\Bbbk}$ is semisimple if and only if $\mathbf{fdMod}(\mathtt{A}) \cong \bigoplus_{i=1}^{r} \mathbf{fdVec}_{\Bbbk}$. Thus, the prototypical examples of semisimple algebras are $\Bbbk$ and direct sums of it. ◇

**Proposition 7E.13.** *A category* $\mathbf{C} \in \mathbf{Cat}_{\Bbbk fA}$ *is semisimple if and only if* $\mathbf{C} \simeq_e \bigoplus_{i=1}^{r} \mathbf{fdVec}_{\Bbbk}$ *for some* $r \in \mathbb{Z}_{\geq 0}$.

*Proof.* That $\bigoplus_{i=1}^{r} \mathbf{fdVec}_{\Bbbk}$ is semisimple is clear. The converse follows from Theorem 6H.6.(ii) and Example 7E.12. □

In words, Proposition 7E.13 suggests that semisimple categories are categorically uninteresting. However, this perspective overlooks additional structure, such as the monoidal structure. (By analogy, in terms of Grothendieck classes, it's like noting that when dealing with a ring, one must not overlook its multiplication.) With this in mind, let us return to fiat and tensor categories. What follows is the categorical counterpart to ***Maschke's theorem***.

**Theorem 7E.14.** *Let* $\mathbf{C} \in \mathbf{wmlFiat}$ *or* $\mathbf{C} \in \mathbf{wmlTen}$. *Then* $\mathbf{C}$ *is semisimple if and only if* $\mathbb{1} \in \mathbf{Proj}(\mathbf{C})$ *if and only if* $\mathbb{1} \in \mathbf{Inj}(\mathbf{C})$.

*Proof.* By combining Proposition 7D.13 and Lemma 7E.8. □

**Example 7E.15.** The classical formulation of Maschke's theorem is the following: "Let G be a finite group of order $m = \#\mathrm{G}$ and $\mathbb{K}$ be algebraically closed. Then $\mathbf{fdMod}(\mathrm{G})$ is semisimple if and only if $\mathrm{char}(\mathbb{K}) \nmid m$.". The original proof of Maschke uses Theorem 7E.14 in the following incarnation: First, observe that $\mathbb{K}[\mathrm{G}] \in \mathbf{fdMod}(\mathrm{G})$ is projective, and so are its direct summands. Further, the sum of all group elements

$$x = \sum_{i=1}^{m} g_i \in \mathbb{K}[\mathrm{G}] \cong \mathrm{End}_{\mathbf{fdMod}(\mathrm{G})}\big(\mathbb{K}[\mathrm{G}]\big), \quad \mathrm{G} = \{g_1, \dots, g_m\},$$

spans a copy of the trivial G module $\mathbb{1}$, which is the monoidal unit in $\mathbf{fdMod}(\mathrm{G})$. Now the crucial calculation:

$$x^2 = m \cdot x.$$

Thus, if $m \neq 0$ in $\mathbb{K}$, then we get an idempotent $\frac{1}{m} \cdot x$ showing that $\mathbb{1} \Subset \mathbb{K}[\mathrm{G}]$. Hence, $\mathbb{1}$ is projective. ◇

**Lemma 7E.16.** *For any* $\Bbbk$ *linear* $\mathbf{C} \in \mathbf{wmfTen}$ *the* $\Bbbk$ *algebra* $\mathrm{End}_{\mathbf{C}}(\mathbb{1})$ *is semisimple.*

*Proof.* We already know that $\mathrm{End}_{\mathbf{C}}(\mathbb{1})$ is a commutative $\Bbbk$ algebra, *cf.* Proposition 2J.4, that is also finite dimensional. By Artin–Wedderburn, it thus remains to show that $f^2 = 0$ implies $f = 0$ for all $f \in \mathrm{End}_{\mathbf{C}}(\mathbb{1})$. So assume that we have such a morphism. We observe that

$$\mathtt{Im}(f)\mathtt{Im}(f) \cong \mathtt{Im}(f^2) \cong \mathtt{Im}(0) \cong 0, \quad \mathtt{Im}(f)\mathtt{Ker}(f) \cong \mathtt{Ker}(f)\mathtt{Im}(f) \cong 0.$$

Thus, $\otimes$ multiplying

$$\mathtt{Ker}(f) \xhookrightarrow{i} \mathbb{1} \xrightarrow{p} \mathtt{Im}(f) \quad \text{SES}$$

with $\mathtt{Im}(f)$ shows, by Proposition 7D.12, that $\mathtt{Im}(f) \cong 0$ and we are done. □

We expect the number of semisimple categories to be small, and we now support this expectation with an analogy. According to [**OEI23**, A038538], the number of semisimple rings $a(n)$ with $n$ elements is given by a sequence so small that its ***asymptotic mean*** is a constant:

$$\lim_{n\to\infty} \frac{1}{n} \sum_{k=1}^{n} a(k) \approx 2.499616.$$

In other words, for any $n \in \mathbb{Z}$, the expected number of semisimple rings is roughly 2.5, regardless of $n$. However, this scarcity does not diminish their intrigue, as illustrated in Figure 19. On the contrary, their simplicity relative to general rings makes them particularly appealing for study.

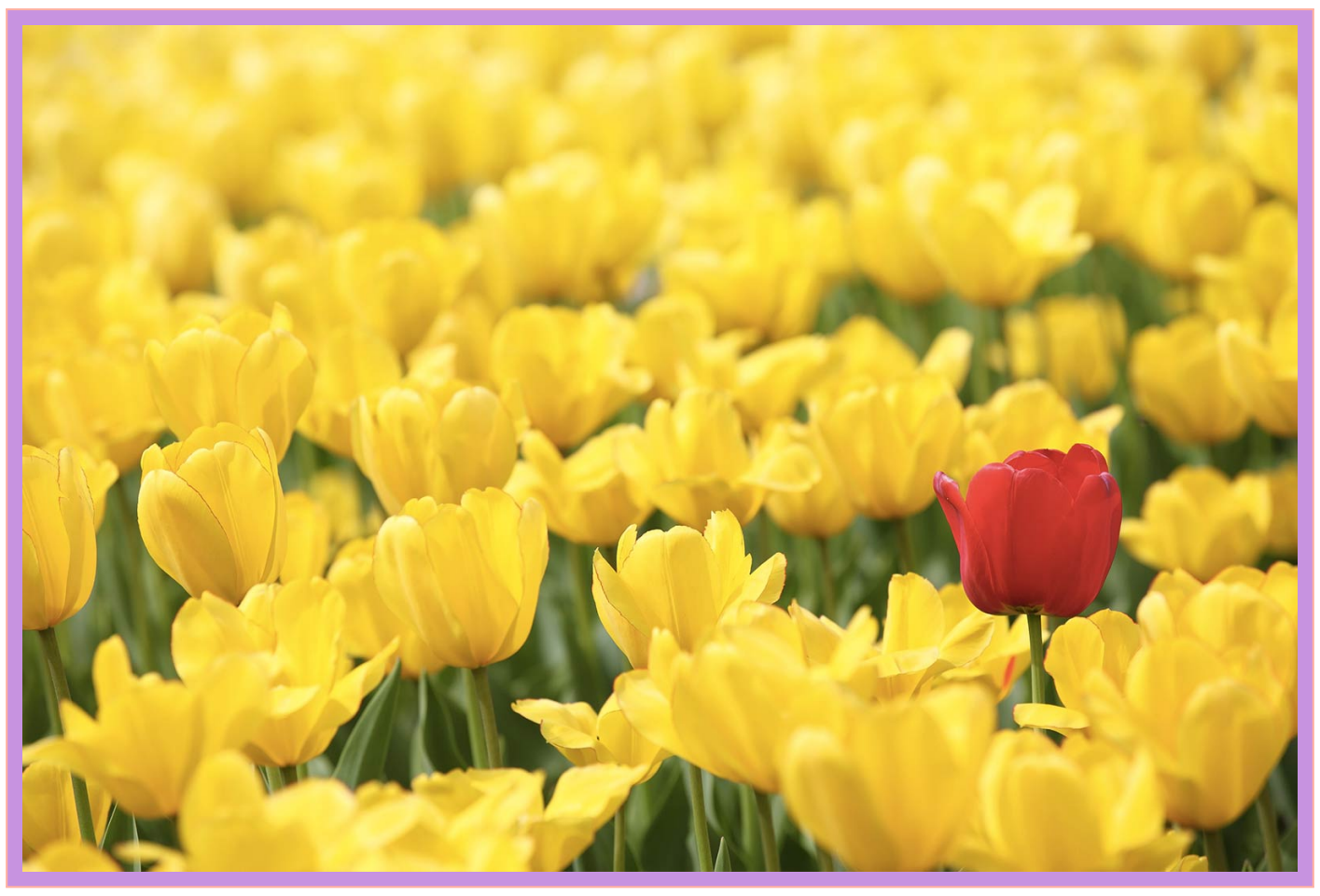

FIGURE 19. Rarity does not inherently equate to dullness; in fact, it can often attract significant interest.

Picture from https://www.chiesi.com/img/aree/sbmsjelqmcmalattie-rare-chiesi-farmaceutici.jpg

7F. **Even more diagrammatics.** Let us revise the categories **TL**, see Example 3D.6, and **Br**, see Example 3D.7.

**Definition 7F.1.** Let $q^{1/2} \in \mathbb{S}^*$. The ***Rumer–Teller–Weyl*** category $\mathbf{TL}^q_{\mathbb{S}\oplus\Subset}$ is the quotient of $\mathbf{TL}_{\mathbb{S}\oplus\Subset}$ by the ***circle removal***

(7F-2) 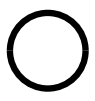 $= -(q + q^{-1})$.

We further endow $\mathbf{TL}^q_{\mathbb{S}\oplus\Subset}$ with the structure of a braided category by extending

(7F-3) 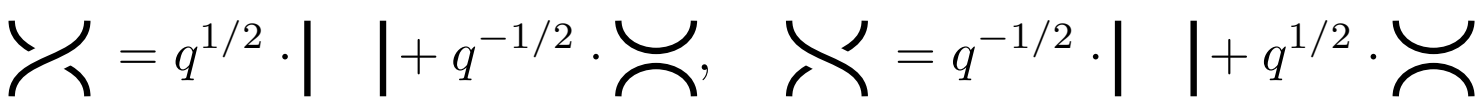

monoidally using Equation 5C-4. ◇

Recall that Equation 5C-4 implies that the overcrossing $\bullet^2\bullet^3 \to \bullet^3\bullet^2$ is

$\beta_{\bullet^2,\bullet^3} =$ 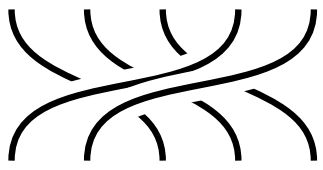 .

Clearly, (7F-3) gives rise to what is known as the ***Kauffman skein relation***:

(7F-4) 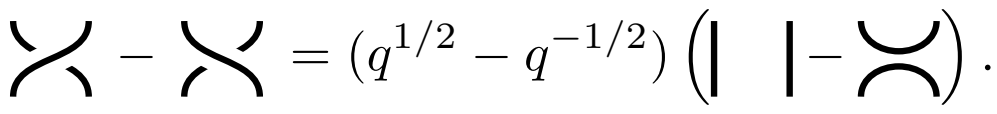

**Lemma 7F.5.** *The category* $\mathbf{TL}^q_{\mathbb{S}\oplus\Subset}$ *is a braided l fiat category (with the braiding in* (7F-3)*).*

*Proof.* First note that (7F-2) and isotopy invariance show that the hom spaces of $\mathbf{TL}^q_{\mathbb{S}\oplus\mathbb{C}}$ are finite dimensional: every diagram can be reduced to a crossingless matching and there are only finitely many of these. The claim that (7F-3) is a braiding is Exercise 7J.4, and the endomorphism of the monoidal unit is

$$\mathrm{End}_{\mathbf{TL}^q_{\mathbb{S}\oplus\mathbb{C}}}(\mathbb{1}) \cong \mathbb{S}, \tag{7F-6}$$

by the circle removal, showing that we can drop the “multi” and the proof is complete. □

A calculation shows that

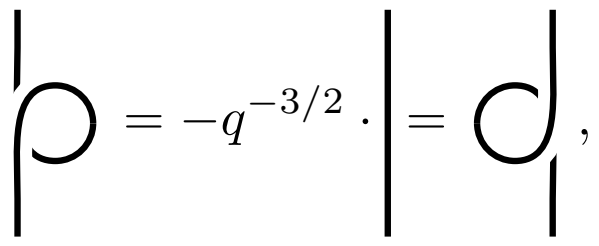

holds in $\mathbf{TL}^q_{\mathbb{S}\oplus\mathbb{C}}$, and thus, $\mathbf{TL}^q_{\mathbb{S}\oplus\mathbb{C}}$ satisfies (5H-3). Hence, we get our first quantum invariant, which is ribbon:

**Proposition 7F.7.** *There exists a well-defined functor*

$$\mathrm{RT}^A_{r=2}\colon \mathbf{1rTan} \to \mathbf{TL}^q_{\mathbb{S}\oplus\mathbb{C}}, \quad \bullet \mapsto \bullet,$$

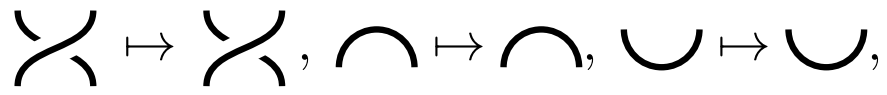

*of braided pivotal categories. Moreover,* $\mathrm{l} \in \mathrm{End}_{\mathbf{1rTan}}(\mathbb{1})$ *is mapped to an element of* $\mathbb{S}$ *under* $\mathrm{RT}^A_{r=2}$.

*Proof.* By construction, there is almost nothing to show: **1rTan** is the free braided pivotal category generated by one self-dual object, and thus there exists the claimed functor by Lemma 7F.5. The final claim follows from Equation 7F-6. □

**Example 7F.8.** The value $\mathrm{RT}^A_{r=2}(\mathrm{l})$ of a ***link*** l, which, by definition, is a morphism $\mathrm{l} \in \mathrm{End}_{\mathbf{1rTan}}(\mathbb{1})$, is an element of $\mathbb{S}$, and an invariant of the link. To be completely explicit, take $\mathbb{S} = \mathbb{Z}[q^{1/2}, q^{-1/2}]$ and $q$ is the corresponding formal variable. Then $\mathrm{RT}^A_{r=2}(\mathrm{l})$ is a (Laurent) polynomial, which is (up to normalization) the so-called ***Jones polynomial***. For instance, take l to be the Hopf link:

$\mathrm{l} =$ 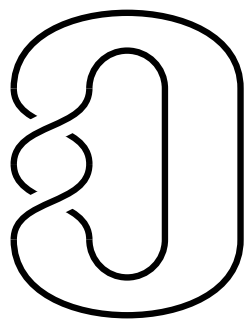 $\in \mathrm{End}_{\mathbf{1rTan}}(\mathbb{1})$.

Then we calculate

$\mathrm{RT}^A_{r=2}(\mathrm{l}) =$

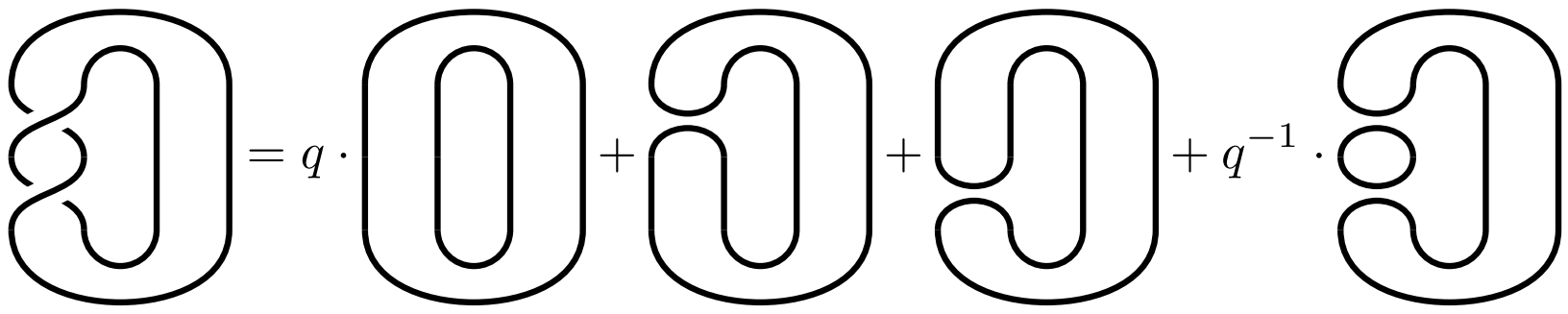

$$= q(q+q^{-1})^2 - 2(q+q^{-1}) + q^{-1}(q+q^{-1})^2 = q^3 + q + q^{-1} + q^{-3},$$

which, up to normalization, is the Jones polynomial of the Hopf link. ◇

Let $\mathrm{TL}_n = \mathrm{End}_{\mathbf{TL}^q_{\mathbb{S}\oplus\mathbb{C}}}(\bullet^n)$, the ***Temperley–Lieb algebra***.

**Lemma 7F.9.** *The Temperley–Lieb algebra admits the following generator–relation presentation. The generators are* $\mathrm{u}_i$ *for* $i \in \{1, \dots, n-1\}$ *and the relations are*

$$\mathrm{u}_i^2 = -(q+q^{-1})\mathrm{u}_i, \quad \mathrm{u}_i\mathrm{u}_j\mathrm{u}_i = \mathrm{u}_i \text{ if } |i-j| = 1, \quad \mathrm{u}_i\mathrm{u}_j = \mathrm{u}_j\mathrm{u}_i \text{ if } |i-j| > 1. \tag{7F-10}$$

*This presentation is meant as an* $\mathbb{S}$ *algebra.*

*Proof.* Let us identify

$\mathrm{u}_i =$ 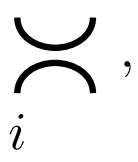 ,

which we wrote in algebraic notation, where the subscript indicates the position, reading left-to-right, of the left strand. With this notation, the relations in Equation 7F-10 hold, so we get a surjection from the Temperley–Lieb algebra onto the algebra generated by the $\mathrm{u}_i$ subject to Equation 7F-10. This surjection is injective is well known, but a bit of a nasty exercise, see [Eas21] for a self-contained argument. □

Similarly:

**Definition 7F.11.** Let $q, a \in \mathbb{S}^*$, with $q \neq q^{-1}$. The quantum Brauer category $\mathbf{rqBr}^{a,q}_{\mathbb{S}\oplus\mathbb{C}}$ is the quotient of $\mathbf{qBr}_{\mathbb{S}\oplus\mathbb{C}}$ by the the ***circle removal***

(7F-12) 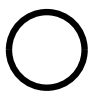 $= \left(\frac{a-a^{-1}}{q-q^{-1}} + 1\right),$

and the ***Kauffman skein*** and ***twist relations***

(7F-13)

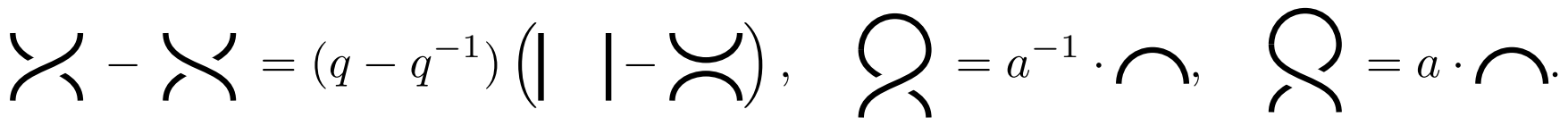

(Note that Equation 7F-4 and Equation 7F-13 look quite alike.) ◇

The following can be proven analogously as Lemma 7F.5 and omitted for the time being.

**Lemma 7F.14.** *The category* $\mathbf{rqBr}^{q}_{\mathbb{S}\oplus\mathbb{C}}$ *is a braided (with the braiding in (7F-13)) l fiat category.* □

The next lemma will be proven in Section 7G below.

**Lemma 7F.15.** *We have*

$$\mathrm{End}_{\mathbf{rqBr}^{q}_{\mathbb{S}\oplus\mathbb{C}}}(\mathbb{1}) \cong \mathbb{S}. \tag{7F-16}$$

Note that we have

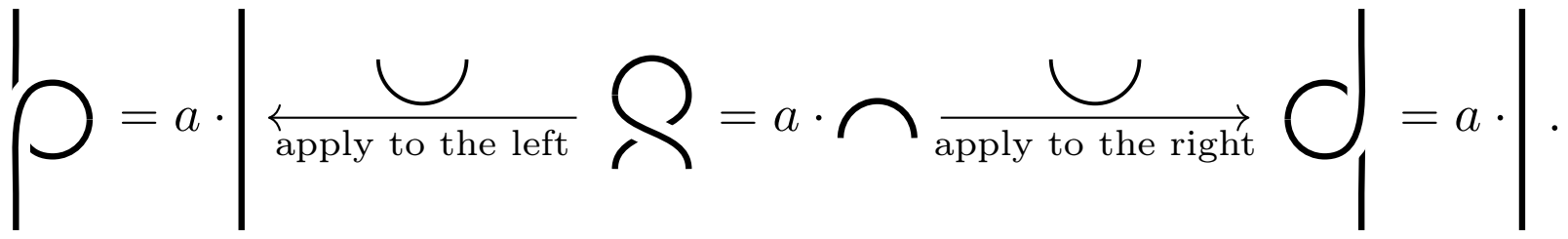

This implies that $\mathbf{rqBr}^{q}_{\mathbb{S}\mathbb{C}}$ is indeed ribbon since it satisfies (5H-3). Hence, we get another quantum invariant, the proof being the same as before (using Equation 7F-16):

**Proposition 7F.17.** *There exists a well-defined functor*

$\mathrm{RT}^{BCD}_{\infty}\colon \mathbf{1rTan} \to \mathbf{rqBr}^{a,q}_{\mathbb{S}\oplus\mathbb{C}}$, 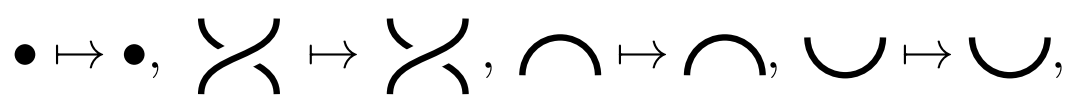

*of braided pivotal categories. Moreover,* $\mathrm{l} \in \mathrm{End}_{\mathbf{1rTan}}(\mathbb{1})$ *is mapped to an element of* $\mathbb{S}$ *under* $\mathrm{RT}^{BCD}_{\infty}$. □

**Example 7F.18.** For the Hopf link as in Example 7F.8 we get

$\mathrm{RT}^{BCD}_{r=\infty}(\mathrm{l}) =$ 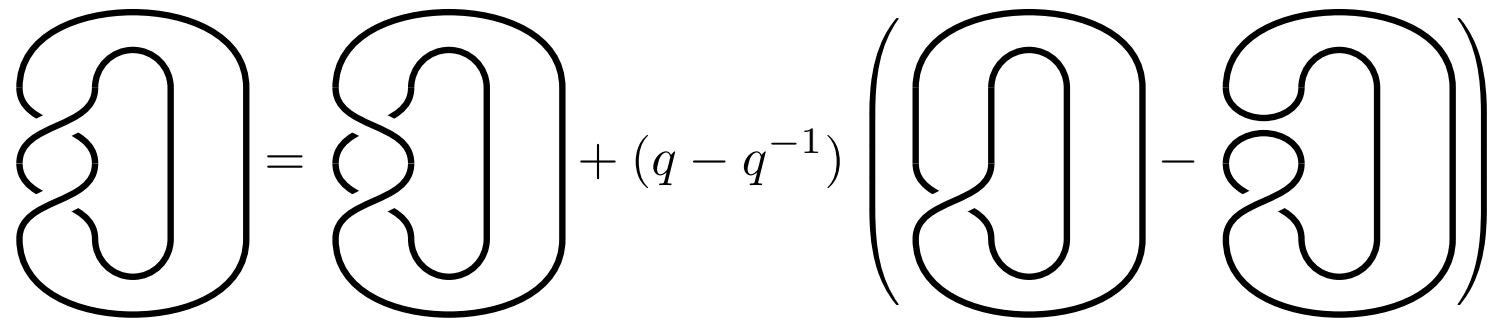

$$= \left(\tfrac{a-a^{-1}}{q-q^{-1}} + 1\right)^2 + (q-q^{-1})(a-a^{-1})\left(\tfrac{a-a^{-1}}{q-q^{-1}} + 1\right)$$
$$\overset{a=q^2}{=} (q+q^{-1})(q^3+1+q^{-3}),$$

where we substituted $a = q^2$ in the last equation to get a nice and short formula. ◇

*Remark* 7F.19. Actually, adding orientations would give quantum invariants of **1Ribbon**. Moreover, one can normalize the two invariants above and similar invariants to get a quantum invariant of **1State**, *i.e.* with honest Reidemeister 1 moves (5H-2). ◇

Let $\mathrm{Br}_n = \mathrm{End}_{\mathbf{rqBr}^{a,q}_{\mathbb{S}\oplus\mathbb{C}}}(\bullet^n)$, the ***quantum Brauer algebra*** or ***BMW algebra***.

**Lemma 7F.20.** *The quantum Brauer algebra admits the following generator–relation presentation. The generators are* $\mathrm{u}_i, \mathrm{b}_i$ *for* $i \in \{1, \ldots, n-1\}$ *and the relations are*

$$\begin{gathered}
\mathrm{b}_i \textit{ is invertible (the inverse is denoted } \mathrm{b}_i^{-1}\textit{)},\\
\mathrm{b}_i\mathrm{b}_j = \mathrm{b}_j\mathrm{b}_i, \textit{if } |i-j| > 1, \quad \mathrm{b}_i\mathrm{b}_{i+1}\mathrm{b}_i = \mathrm{b}_{i+1}\mathrm{b}_i\mathrm{b}_{i+1}, \quad \mathrm{u}_i\mathrm{u}_{i+1}\mathrm{u}_i = \mathrm{u}_i,\\
\mathrm{b}_i - \mathrm{b}_i^{-1} = (q-q^{-1})(1+\mathrm{u}_i), \quad \mathrm{b}_{i+1}\mathrm{b}_i\mathrm{u}_{i+1} = \mathrm{u}_i\mathrm{b}_{i+1}\mathrm{b}_i = \mathrm{u}_i\mathrm{u}_{i+1}, \quad \mathrm{b}_{i+1}\mathrm{u}_i\mathrm{b}_{i+1} = \mathrm{b}_i^{-1}\mathrm{u}_{i+1}\mathrm{b}_i^{-1},\\
\mathrm{b}_{i+1}\mathrm{u}_i\mathrm{u}_{i+1} = \mathrm{b}_i^{-1}\mathrm{u}_{i+1}, \quad \mathrm{u}_{i+1}\mathrm{u}_i\mathrm{b}_{i+1} = \mathrm{u}_{i+1}\mathrm{b}_i^{-1}, \quad \mathrm{b}_i\mathrm{u}_i = \mathrm{u}_i\mathrm{b}_i = a^{-1}\mathrm{u}_i, \quad \mathrm{u}_i\mathrm{b}_{i+1}\mathrm{u}_i = a\mathrm{u}_i.
\end{gathered} \tag{7F-21}$$

*This presentation is meant as an* $\mathbb{S}$ *algebra.*

*Proof.* Using the same notation as in Lemma 7F.9 and

$\mathrm{b}_i =$ 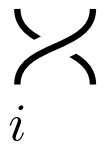 ,

the same arguments as in Lemma 7F.9 work. The final step is however nontrivial, but can be found in [MW89] where they find a basis of the abstract algebra generated as in the statement of the lemma. □

7G. **Why quantum invariants want a nontrivial quantum parameter.** As usual, the following topological statement is neither rigorously formulated nor proven (but is still true).

**Proposition 7G.1.** *Knot theory in 4d is trivial, i.e. every link can be undone.*

*Proof.* Imagine the fourth dimension as the colors of the rainbow, so that red (dotted) and purple (dashed) are far apart, and green (straight) is in the middle. Now:

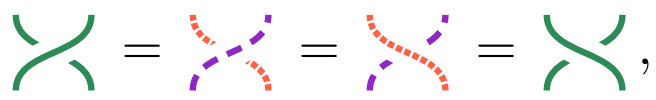

where the color changes happen continuously. Thus, on diagrams we have

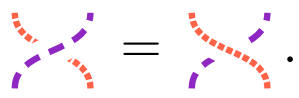

With this relation every knot diagram can be undone. A fast way to see this is to imagine the knot in $\mathbb{R}^3$ but you are allowed to change overcrossings and undercrossings. Start anywhere, denoted • below (and for simplicity in the $x$-$y$ plane), on the knot and travel in some direction. Whenever you visit a crossing that you have not seen before and it is an undercrossing, change it to an overcrossing. In the end you get a helix:

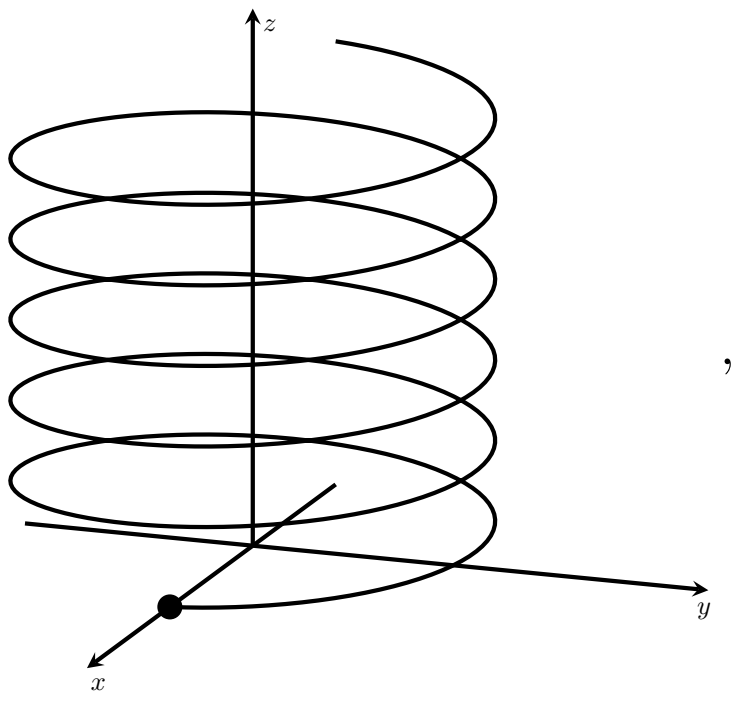

,

and then one travels down again. This is the unknot. The same argument works for links. □

A quantum invariant $\mathrm{Q}\colon \mathbf{1rTan} \to \mathbf{C}$ is ***trivial*** if every link $\mathrm{l} \in \mathrm{End}_{\mathbf{1rTan}}(\mathbb{1})$ is mapped to the same morphism as the unlink with the same number of connected components.

**Theorem 7G.2.** *If a quantum invariant* Q *as defined in Section 5L maps to a symmetric category* **C***, then it is trivial.*

*Proof.* Recall that a braided category is symmetric if

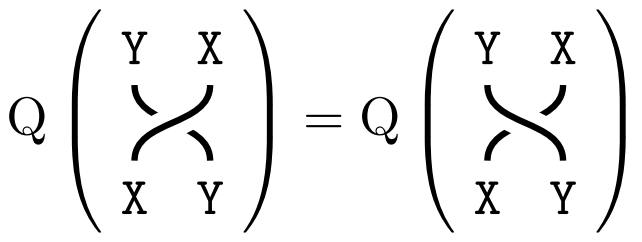

.

Now apply the same arguments as in Proposition 7G.1 □

The case where $q = \pm 1$ is commonly referred to as ***classical***. The following statement asserts that classical invariants are trivial. This point is significant: in a certain sense, quantum topology can be seen as "classical mathematics with a $q$," but the presence of $q$ is, in fact, profoundly important.

**Proposition 7G.3.** *The quantum invariants of Proposition 7F.7 and Proposition 7F.17 are trivial for $q = 1$ and $q = -1$.*

*Proof.* Observe that the two terms on the right-hand side of the Kauffman skein relation (either Equation 7F-4 or Equation 7F-13) can be seen as error terms that prevent the overcrossing from being the undercrossing. However, for $q = \pm 1$ they vanish, so the same arguments as in Proposition 7G.1 apply. □

## 7H. Multiplicative structures on Grothendieck classes.

Let us come back to Definition 6K.12.

**Definition 7H.1.** Let $\mathbf{C} \in \mathbf{Fiat}$ and $\mathbf{D} \in \mathbf{Ten}$. Then we define the ***additive Grothendieck classes*** $K_0^{\oplus}(\mathbf{C})$ of $\mathbf{C}$ respectively the ***SES Grothendieck classes*** $K_0^{e}(\mathbf{D})$ of $\mathbf{D}$ *verbatim* as in Definition 6K.12. ◇

Clearly, we have the analogs of Proposition 4D.15 and Lemma 6K.13:

**Proposition 7H.2.** *Let $\mathbf{C} \in \mathbf{Fiat}$ and $\mathbf{D} \in \mathbf{Ten}$. Then:*

*(i) Definition 7H.1 endows $K_0^{\oplus}(\mathbf{C})$ and $K_0^{e}(\mathbf{D})$ with the structures of finite dimensional abelian groups.*

*(ii) The set $\mathrm{In}(\mathbf{C})$ is a basis of $K_0^{\oplus}(\mathbf{C})$. We have*

$$[\mathtt{X}] = \textstyle\sum_{i=1}^{n} (\mathtt{X} : \mathtt{Z}_i) \cdot [\mathtt{Z}_i] \in K_0^{\oplus}(\mathbf{C}).$$

*(iii) The set $\mathrm{Si}(\mathbf{D})$ is a basis of $K_0^{e}(\mathbf{D})$. We have*

$$[\mathtt{X}] = \textstyle\sum_{i=1}^{n} [\mathtt{X} : \mathtt{L}_i] \cdot [\mathtt{L}_i] \in K_0^{e}(\mathbf{D}).$$

*(iv) For both, $K_0^{\oplus}(\mathbf{C})$ and $K_0^{e}(\mathbf{D})$, the additional structures in Proposition 4D.15 are compatible with the $\mathbb{S}$ linear and additive structures. In particular, $K_0^{\oplus}(\mathbf{C})$ and $K_0^{e}(\mathbf{D})$ are finite dimensional $\mathbb{Z}$ algebras.* □

By Lemma 2H.9 and Lemma 6K.14 we also have:

**Proposition 7H.3.** *Let $\mathbf{C}, \mathbf{C}' \in \mathbf{Fiat}$, and let $\mathbf{D}, \mathbf{D}' \in \mathbf{Ten}$.*

*(i) Any functor $\mathrm{F} \in \mathbf{Hom}_{\Bbbk\oplus\star}(\mathbf{C}, \mathbf{C}')$ induces a $\mathbb{Z}$ algebra homomorphism*

$$K_0^{\oplus}(\mathrm{F})\colon K_0^{\oplus}(\mathbf{C}) \to K_0^{\oplus}(\mathbf{C}'),\ [\mathtt{X}] \mapsto [\mathrm{F}(\mathtt{X})].$$

*Further, if $\mathrm{F}$ is an equivalence, then $K_0^{\oplus}(\mathrm{F})$ is an isomorphism.*

*(ii) Any functor $\mathrm{F} \in \mathbf{Hom}_{\Bbbk e\star}(\mathbf{D}, \mathbf{D}')$ induces a $\mathbb{Z}$ algebra homomorphism*

$$K_0^{e}(\mathrm{F})\colon K_0^{e}(\mathbf{D}) \to K_0^{e}(\mathbf{D}'),\ [\mathtt{X}] \mapsto [\mathrm{F}(\mathtt{X})].$$

*Further, if $\mathrm{F}$ is an equivalence, then $K_0^{e}(\mathrm{F})$ is an isomorphism.* □

This gives a (coarse) numerical invariant:

**Proposition 7H.4.** *Let $\mathbf{C} \in \mathbf{Fiat}$, and let $\mathbf{D} \in \mathbf{Ten}$. The the* **ranks** *$\mathrm{rk}(\mathbf{C})$ and $\mathrm{rk}(\mathbf{D})$, i.e. the dimensions of $K_0^{\oplus}(\mathbf{C})$ and $K_0^{e}(\mathbf{D})$, respectively, are invariants of $\mathbf{C}$ respectively $\mathbf{D}$.* □

*Remark* 7H.5. All of the above finiteness conditions prevent us from running into the ***Eilenberg swindle***: If $\mathtt{X} \cong \mathtt{Y} \oplus \mathtt{Y} \oplus \mathtt{Y} \oplus \mathtt{Y} \oplus \dots$ would be an allowed object, then $\mathtt{X} \oplus \mathtt{Y} \cong \mathtt{X}$ which gives $[\mathtt{Y}] = 0$. This would then hold for any object, since $\mathtt{Y}$ was arbitrary. ◇

**Example 7H.6.** In Example 7D.7, we have isomorphisms of rings

$$[\mathtt{P}_1] \mapsto 5 \rightsquigarrow K_0^{\oplus}(\mathbf{C}') \xrightarrow{\cong} 5\mathbb{Z} \subset \mathbb{Z} \xleftarrow{\cong} K_0^{e}(\mathbf{C}) \leftsquigarrow 1 \leftarrow [\mathtt{L}_1].$$

Further, the endofunctor $_ \otimes \mathtt{P}_1\colon \mathbf{C} \to \mathbf{C}$ gives

$$[_ \otimes \mathtt{P}_1]\colon \mathbb{Z} \to \mathbb{Z}, \quad 1 \mapsto 5.$$

Moreover, $\mathtt{L}_1 = \mathtt{Z}_1$ is the monoidal unit of $\mathbf{C}$, $\mathtt{P}_1 = \mathtt{Z}_5$ a "big pseudo idempotent" and

(7H-7)

| $\otimes$ | $\mathtt{Z}_1$ | $\mathtt{Z}_2$ | $\mathtt{Z}_3$ | $\mathtt{Z}_4$ | $\mathtt{Z}_5$ |
|---|---|---|---|---|---|
| $\mathtt{Z}_1$ | $\mathtt{Z}_1$ | $\mathtt{Z}_2$ | $\mathtt{Z}_3$ | $\mathtt{Z}_4$ | $\mathtt{Z}_5$ |
| $\mathtt{Z}_2$ | $\mathtt{Z}_2$ | $\mathtt{Z}_1 \oplus \mathtt{Z}_3$ | $\mathtt{Z}_2 \oplus \mathtt{Z}_4$ | $\mathtt{Z}_3 \oplus \mathtt{Z}_5$ | $2 \cdot \mathtt{Z}_5$ |
| $\mathtt{Z}_3$ | $\mathtt{Z}_3$ | $\mathtt{Z}_2 \oplus \mathtt{Z}_4$ | $\mathtt{Z}_1 \oplus \mathtt{Z}_3 \oplus \mathtt{Z}_5$ | $\mathtt{Z}_2 \oplus \mathtt{Z}_5 \oplus \mathtt{Z}_5$ | $3 \cdot \mathtt{Z}_5$ |
| $\mathtt{Z}_4$ | $\mathtt{Z}_4$ | $\mathtt{Z}_3 \oplus \mathtt{Z}_5$ | $\mathtt{Z}_2 \oplus \mathtt{Z}_5 \oplus \mathtt{Z}_5$ | $\mathtt{Z}_1 \oplus \mathtt{Z}_5 \oplus \mathtt{Z}_5 \oplus \mathtt{Z}_5$ | $4 \cdot \mathtt{Z}_5$ |
| $\mathtt{Z}_5$ | $\mathtt{Z}_5$ | $2 \cdot \mathtt{Z}_5$ | $3 \cdot \mathtt{Z}_5$ | $4 \cdot \mathtt{Z}_5$ | $5 \cdot \mathtt{Z}_5$ |

(One can check this using the Jordan decomposition over $\overline{\mathbb{F}}_5$.) Thus, we have an isomorphism of rings

$$K_0^{\oplus}(\mathbf{C}) \xrightarrow{\cong} \mathbb{Z}, \quad [\mathtt{Z}_i] \mapsto i.$$

Finally, note that (7H-7) also clearly shows the pivotal structure, since $\mathtt{Z}_1 = \mathbb{1} \in \mathtt{Z}_i\mathtt{Z}_j$ if and only if $i = j$ and $i < 5$. This shows that all indecomposables are self-dual (since $\mathtt{Z}_5$ is the projective cover of $\mathtt{Z}_1$, the monoidal product $\mathtt{Z}_5\mathtt{Z}_5$ has a map to $\mathtt{Z}_1$ regardless whether $\mathtt{Z}_1$ appears as a summand of it). ◇

*Remark* 7H.8. All of the above have appropriate versions in the "weakly" and "multi" setup. ◇

7I. **Finite dimensional algebras in vector spaces.** Here is the main source of examples of tensor categories:

**Theorem 7I.1.** *Let* $\mathtt{A} \in \mathbf{fdVec}_{\mathbb{S}}$ *be a Hopf algebra. Then* $\mathbf{fdMod}(\mathtt{A}) \in \mathbf{mTen}$.

*Proof.* Combining Theorem 6H.6.(ii), showing that $\mathbf{fdMod}(\mathtt{A})$ is finite abelian, and Theorem 5K.6, showing that $\mathbf{fdMod}(\mathtt{A})$ is rigid, and observing that everything is compatible with the $\mathbb{k}$ linear structure. □

So Hopf algebras play a crucial role in the construction of quantum invariants. As an aside, another nice fact about Hopf algebras is that they are "group-like." Let us make this precise. To this end, let $\mathbb{S}\mathbf{Alg}$ denote the ***category of*** $\mathbb{S}$ ***algebras***, objects being $\mathbb{S}$ algebras and morphisms $\mathbb{S}$ algebra homomorphisms. Furthermore, let $\mathbf{Gr} \subset \mathbf{Mon}$ (recall that $\mathbf{Mon}$ is the category of monoids, *cf.* Example 1B.2.(a)) denote the full subcategory whose objects are groups, *i.e.* the ***category of groups***.

**Definition 7I.2.** We call $\mathrm{F} \in \mathbf{Hom}(\mathbb{S}\mathbf{Alg}, \mathbf{Gr})$ ***representable*** if $\mathrm{Forget} \circ \mathrm{F} \in \mathbf{Hom}(\mathbb{S}\mathbf{Alg}, \mathbf{Set})$ is representable in the sense of Example 1H.8. ◇

**Proposition 7I.3.** *Let* $\mathrm{F} \in \mathbf{Hom}(\mathbb{S}\mathbf{Alg}, \mathbf{Gr})$ *be represented by* $\mathtt{A} \in \mathbb{S}\mathbf{Alg}$*. Then:*

*(i) For all* $\mathtt{A} \in \mathbb{S}\mathbf{Alg}$*, the set* $\mathrm{End}_{\mathbb{S}\mathbf{Alg}}(\mathtt{A})$ *has a group structure.*

*(ii) The multiplication and inversion in (i) come from a comultiplication and an antipode on* $\mathtt{A}$*, making* $\mathtt{A}$ *a Hopf algebra.*

For this reason one can say that Hopf algebras are ***cogroup objects in*** $\mathbb{S}\mathbf{Alg}$.

*Proof.* This is Exercise 7J.5. □

7J. **Exercises.**

*Exercise* 7J.1. Prove the missing points in Proposition 7C.2. ◇

*Exercise* 7J.2. Try to make Remark 7D.4 precise by drawing a hierarchy chart and by giving examples whenever Notion A $\not\Rightarrow$ Notion B. ◇

*Exercise* 7J.3. Understand Example 7D.7 and make all claims made in that example precise, *e.g.* the monoidal structure. ◇

*Exercise* 7J.4. Show that (7F-3) defines the structure of a braided category on $\mathbf{TL}^{q}_{\mathbb{S}\oplus\mathbb{C}}$. Compute the quantum invariant $\mathrm{RT}^{A}_{r=2}(_)$ for

$1 =$ 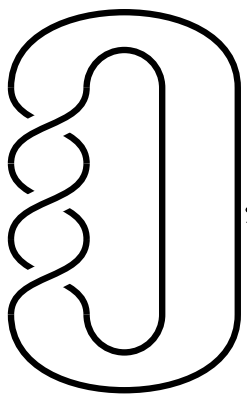 $, 1' =$ 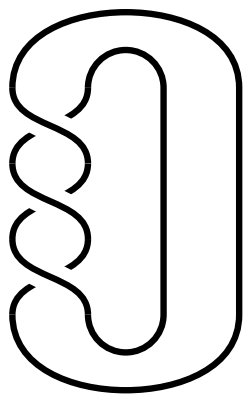 $\in \mathrm{End}_{\mathbf{1Tan}}(\mathbb{1})$.

This is the trefoil knot and its mirror image. Deduce that they are not equivalent. ◇

*Exercise* 7J.5. Prove Proposition 7I.3. Hint: Yoneda. ◇

## 8. Fiat, tensor and fusion categories – definitions and classifications

Fiat and tensor categories serve as a categorification of algebras, and, as we will see, when they are semisimple, they categorify finite groups in a certain sense. When studying finite groups, a natural inclination is to classify them, perhaps by fixing a numerical invariant such as the group's size, known as its order. This leads to statements like, "All finite groups of prime order are cyclic." So:

> Can one hope to classify (semisimple) fiat and tensor categories, maybe after fixing some numerical invariant such as the rank?

The answer will turn out to be "Yes and no." Indeed, that is not surprising: the group-theoretical analog would be to classify simple groups (groups without nontrivial subgroups). When our numerical invariant is a prime, then "All finite groups of prime order are cyclic" directly classifies the simple groups of prime order. Even a tiny step beyond this case and it becomes tricky *cf.* Figure 20.

The answer turns out to be "yes and no." This ambivalence is not surprising: in group theory, an analogous task would be to classify simple groups; those groups without nontrivial proper subgroups. When our numerical invariant is a prime, the statement "All finite groups of prime order are cyclic" directly classifies the simple groups of that order. However, even a small step beyond this case reveals significant complexity; *cf.* Figure 20.

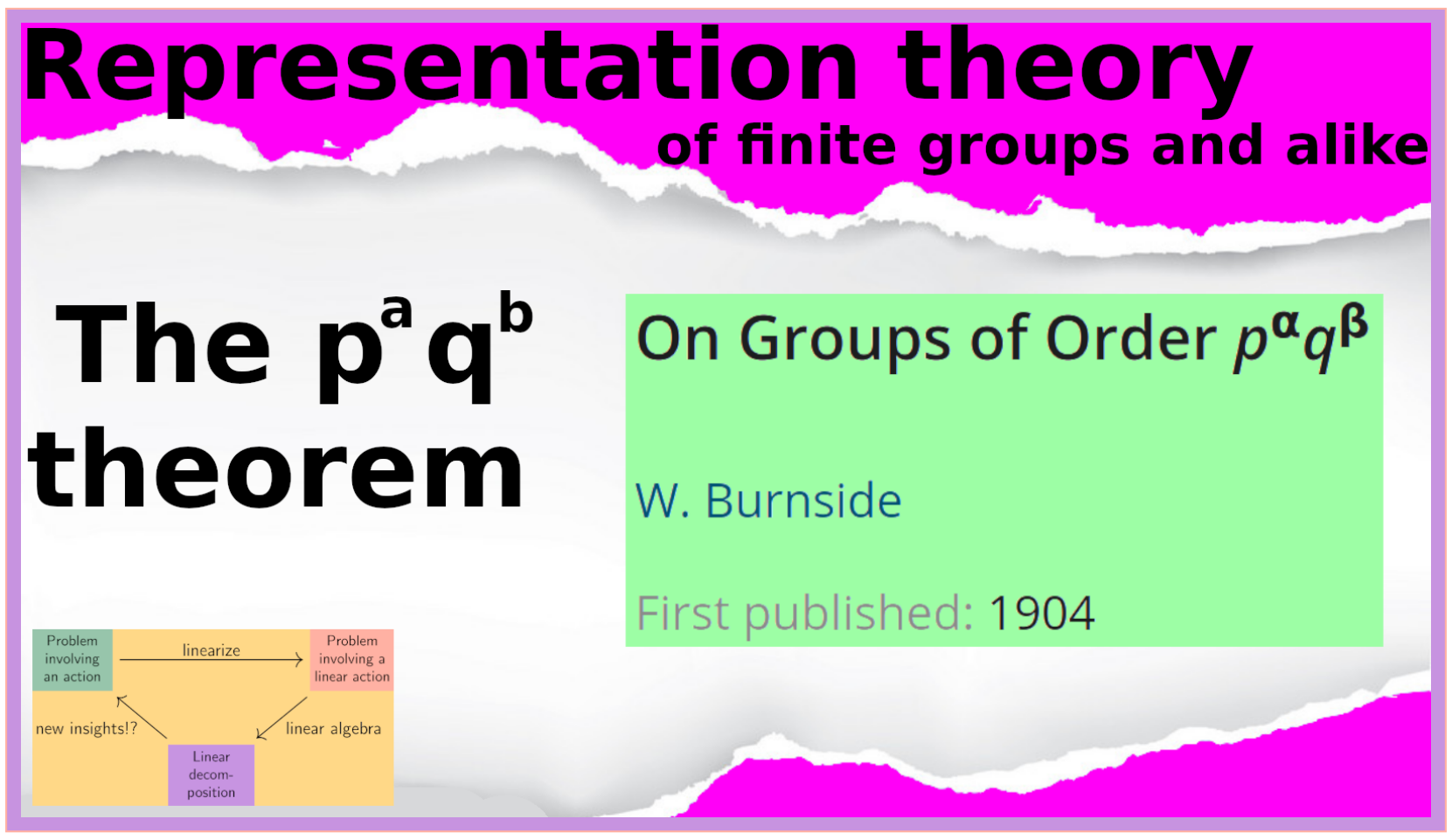

FIGURE 20. Burnsides's $p^a q^b$ theorem implies that noncyclic simple groups have order divisible by at least three primes. In particular, since the is no simple group of order $2 \cdot 3 \cdot 5 = 30$, the smallest possible order for such a group would be $2^2 \cdot 3 \cdot 5 = 60$. (Do you know an example of a simple group of this order?) The proofs of this theorem are all highly nontrivial.
Picture from https://www.youtube.com/watch?v=SJ0A0ewTxg4

8A. **A word about conventions.** Of course, we keep the previous conventions.

*Convention* 8A.1. We will identify directed graphs $\Gamma$ and their adjacency matrices $M$, which we see as matrices with values in $\mathbb{Z}_{\geq 0}$, and we will write $\Gamma$ for both if confusion cannot arise.

Yes, square matrices (with values in $\mathbb{Z}_{\geq 0}$) and graphs are the same.☺

The translation between these two notions is best illustrated in an example:

$$M = \begin{array}{c} \\ v_1 \\ v_2 \end{array} \begin{array}{c} \begin{array}{cc} v_1 & v_2 \end{array} \\ \begin{pmatrix} 0 & 1 \\ 2 & 3 \end{pmatrix} \end{array} \in \mathrm{Mat}_{2\times 2}(\mathbb{Z}_{\geq 0}) \leftrightsquigarrow \Gamma = \; v_1 \underset{2}{\overset{}{\rightleftarrows}} v_2 \circlearrowleft^{3} ,$$

where labels mean parallel edges, with the label 1 being omitted from illustrations. ◇

*Convention* 8A.2. Our notation convention for objects is:

$$\mathtt{X} \leftrightsquigarrow \text{general object}, \quad \mathtt{L} \leftrightsquigarrow \text{simple object}, \quad \mathtt{P} \leftrightsquigarrow \text{projective inde. object},$$
$$\mathtt{I} \leftrightsquigarrow \text{injective inde. object}, \quad \mathtt{Z} \leftrightsquigarrow \text{indecomposable object},$$

(We also tend to use projectives instead of the dual notion of injectives.) Recall also that for semisimple categories

$$\mathtt{X} \text{ is simple} \Leftrightarrow \mathtt{X} \text{ is projective inde.} \Leftrightarrow \mathtt{X} \text{ is injective inde.} \Leftrightarrow \mathtt{X} \text{ is indecomposable},$$

and all of these are interchangeable. Although we are mostly concerned with indecomposable objects, we will use the notation $\mathtt{L}$ in the semisimple case for either of these to stress that this case is easier than in the general situation. ◇

8B. **Representations of groups and Hopf algebras.** We start by discussing a very nicely behaved case in detail: Let $\mathrm{G}$ be a finite group of order $\#\mathrm{G} = m$, and let $p$ be the characteristic of an algebraically closed ground field $\mathbb{K}$. Recall that $\mathbb{K}[\mathrm{G}] = (\mathbb{K}[\mathrm{G}], \mathrm{m}, \mathrm{i}, \mathrm{d}, \mathrm{e}, \mathrm{s})$ is a Hopf algebra in $\mathbf{Vec}_{\mathbb{K}}$. The explicit structure maps are the multiplication $\mathrm{m}$ and unit $\mathrm{i}$ in $\mathbb{K}[\mathrm{G}]$, and

$$\mathrm{d}(g) = g \otimes g, \quad \mathrm{e}(g) = 1, \quad \mathrm{s}(g) = g^{-1}.$$

Thus, $\mathbf{fdMod}\big(\mathbb{K}[\mathrm{G}]\big)$ is $\mathbb{K}$ linear abelian rigid (actually, it is even pivotal). Moreover, by (6K-10) we have

$$\text{(8B-1)} \qquad \textstyle\sum_{i=1}^{n} \dim(\mathtt{L}_i)^2 \leq \sum_{i=1}^{n} \dim(\mathtt{L}_i)\dim(\mathtt{P}_i) = \dim\big(\mathbb{k}[\mathrm{G}]\big) = m.$$

We also know by Lemma 7E.7 that equality holds in (8B-1) if and only if $\mathbf{fdMod}\big(\mathbb{K}[\mathrm{G}]\big)$ is semisimple. The latter, by Maschke's theorem Example 7E.15, happens if and only if $p \nmid m$.

*Remark* 8B.2. In fact, (8B-1) holds for any finite dimensional $\Bbbk$ algebra A, *i.e.*

$$\textstyle\sum_{i=1}^{n}\dim(\mathtt{L}_i)^2 \leq \sum_{i=1}^{n}\dim(\mathtt{L}_i)\dim(\mathtt{P}_i) = \dim(\mathrm{A}),$$

with equality if and only if A is semisimple. ◇

**Example 8B.3.** In a field of characteristic $p$, let us perform the calculation

$$(X-1)^5 = X^5 - 5\cdot X^4 + 10\cdot X^3 - 10\cdot X^2 + 5\cdot X - 1 \overset{p=5}{=} X^5 - 1,$$

which is called ***Freshperson's dream***. (Or my dream: Would it not be great if $(X-1)^n = X^n - 1$? Then many things would be trivial, there may not be any research in mathematics and I could just go to the beach ☺.) Hence, in characteristic 5 there is only the trivial 5th root of unity $\zeta = 1$. In all other cases (for algebraically closed fields), there are five primitive roots of unity $\{1=\zeta^0,\zeta^1,\zeta^2,\zeta^3,\zeta^4\}$, *e.g.* for $\mathbb{K}=\mathbb{C}$ we could let $\zeta = \exp(2\pi i/5)$. (See also (6K-17) and the text below.) Let us now come back to $\mathrm{G}=\mathbb{Z}/5\mathbb{Z}$, see also Example 6K.18, where

$$\begin{pmatrix} \zeta^0 & 0 & 0 & 0 & 0 \\ 0 & \zeta^1 & 0 & 0 & 0 \\ 0 & 0 & \zeta^2 & 0 & 0 \\ 0 & 0 & 0 & \zeta^3 & 0 \\ 0 & 0 & 0 & 0 & \zeta^4 \end{pmatrix} \xleftarrow{p\neq 5} \begin{pmatrix} 0 & 0 & 0 & 0 & 1 \\ 1 & 0 & 0 & 0 & 0 \\ 0 & 1 & 0 & 0 & 0 \\ 0 & 0 & 1 & 0 & 0 \\ 0 & 0 & 0 & 1 & 0 \end{pmatrix} \xrightarrow{p=5} \begin{pmatrix} 1 & 0 & 0 & 0 & 0 \\ 1 & 1 & 0 & 0 & 0 \\ 0 & 1 & 1 & 0 & 0 \\ 0 & 0 & 1 & 1 & 0 \\ 0 & 0 & 0 & 1 & 1 \end{pmatrix},$$

is the matrix for the multiplication action of 1 on $\mathbb{K}[\mathbb{Z}/5\mathbb{Z}]$, which has the characteristic polynomial $X^5-1$. We also gave the corresponding Jordan decompositions. Thus, we get two different cases:

- If $\mathbb{K}$ is not of characteristic 5, then

  $$\mathtt{L}_1 = \mathbb{1}\colon \begin{pmatrix}\zeta^0\end{pmatrix}, \quad \mathtt{L}_2\colon \begin{pmatrix}\zeta^1\end{pmatrix}, \quad \mathtt{L}_3\colon \begin{pmatrix}\zeta^2\end{pmatrix}, \quad \mathtt{L}_4\colon \begin{pmatrix}\zeta^3\end{pmatrix}, \quad \mathtt{L}_5\colon \begin{pmatrix}\zeta^4\end{pmatrix},$$

  defines five simples, all of dimension 1. (They, of course, correspond to the five elements of $\mathbb{Z}/5\mathbb{Z}$, as indicated by $\zeta^i$.) They are also projective; the idempotents that separate them from $\mathbb{K}[\mathbb{Z}/5\mathbb{Z}]$ can be obtained by using the change-of-basis matrices from the action matrix to its Jordan decomposition. Hence, formula (8B-1) takes the form

  $$1^2+1^2+1^2+1^2+1^2=5.$$

- We have already seen the case where $\mathbb{K}$ is of characteristic 5, say $\mathbb{K}=\overline{\mathbb{F}}_5$, in detail, see *e.g.* Example 6K.18. In this case, we have one simple $\mathtt{L}_1=\mathbb{1}$ and its projective cover $\mathtt{P}_1\cong\overline{\mathbb{F}}_5[\mathbb{Z}/5\mathbb{Z}]$. Thus, (8B-1) takes the form

  $$1 = 1^2 \leq 1\cdot 5 = 5.$$

The same works for any prime number instead of five. ◇

Note that (8B-1) implies that the set of simples $\mathrm{Si}\big(\mathbf{fdMod}(\mathbb{K}[\mathrm{G}])\big)$ (or of projectives indecomposables or of injectives indecomposables) of $\mathbf{fdMod}\big(\mathbb{K}[\mathrm{G}]\big)$ is always finite. What about the additive version, *i.e.* what about the set of indecomposables $\mathrm{In}\big(\mathbf{fdMod}(\mathbb{K}[\mathrm{G}])\big)$? We have already seen in the case $\mathrm{G}=\mathbb{Z}/5\mathbb{Z}$ that $\#\mathrm{Si}\big(\mathbf{fdMod}(\mathbb{K}[\mathrm{G}])\big) \leq \#\mathrm{In}\big(\mathbf{fdMod}(\Bbbk[\mathrm{G}])\big)$ with equality if and only if we are in the semisimple situation. Actually, the difference can get arbitrary big:

**Example 8B.4.** Klein's group of order four is $\mathrm{V}_4=\mathbb{Z}/2\mathbb{Z}\times\mathbb{Z}/2\mathbb{Z}=\langle s,t \mid s^2=t^2=1, st=ts\rangle$, with its defining action on the complex plane $\mathbb{C}^2=\{a+ib \mid a,b\in\mathbb{R}\}$ given by reflections:

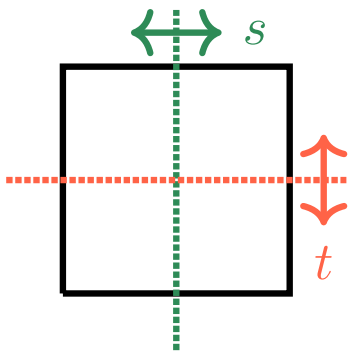

$$s\centerdot(a+ib) = -a+ib, \qquad t\centerdot(a+ib) = a-ib.$$

For $\mathbb{K}$ not of characteristic 2 the category $\mathbf{fdMod}\big(\mathbb{K}[\mathrm{V}_4]\big)$ is semisimple with four simples of dimension one. The case $\mathbb{K}=\overline{\mathbb{F}}_2$ is very different, and we will discuss it now. In this case we have

$$\overline{\mathbb{F}}_2[\mathrm{V}_4] \xrightarrow{\cong} \mathrm{A} = \overline{\mathbb{F}}_2[X,Y]/(X^2,Y^2), \quad s\mapsto X+1, t\mapsto Y+1,$$

by Freshperson's dream.

Let us first discuss the simples and projectives of A. It is easy to see that A has one simple $\mathtt{L}_1$ whose projective cover $\mathtt{P}_1$ is A itself:

(8B-5)
$$\mathtt{L}_1 = \mathbb{1}\colon \ \bullet\ , \qquad \mathtt{P}_1\colon \begin{array}{ccc} & \bullet & \\ {}^{X}\swarrow & & \searrow^{Y} \\ \bullet & & \bullet \\ {}_{Y}\searrow & & \swarrow_{X} \\ & \bullet & \end{array} .$$

Here we use a graph to indicate the modules. This is to be read as follows: the vertices correspond to basis elements while the arrows indicate the non-zero actions of $X$ and $Y$.

In contrast, there are infinitely many indecomposables that are not projective. Here is the list of some of them (see *e.g.* [**Web16**, Section 11.5] for details), using the same notation as in (8B-5):

- For all $2l+3$, where $l \in \mathbb{Z}_{\geq 0}$ (thus, $2\cdot 0+3=3$ is the smallest case), there are two indecomposables $\mathtt{Z}_{2l+1}$ and $\mathtt{Z}^{\star}_{2l+1}$, which are duals:

$$\mathtt{Z}_{2l+1}\colon\ \bullet \xleftarrow{X} \bullet \xrightarrow{Y} \bullet \xleftarrow{X} \bullet \xrightarrow{Y} \bullet \xleftarrow{X} \dots \xrightarrow{Y} \bullet\ ,$$
$$\mathtt{Z}^{\star}_{2l+1}\colon\ \bullet \xrightarrow{Y} \bullet \xleftarrow{X} \bullet \xrightarrow{Y} \bullet \xleftarrow{X} \bullet \xrightarrow{Y} \dots \xleftarrow{X} \bullet\ .$$

  Here and below, the subscript indicates the dimension, *i.e.* the number of vertices.

- For all $2l+2$, where $l \in \mathbb{Z}_{\geq 0}$, there are self-dual indecomposables $\mathtt{Z}_{2l}$ and $\mathtt{Z}'_{2l}$:

$$\mathtt{Z}_{2l}\colon\ \bullet \xleftarrow{X} \bullet \xrightarrow{Y} \bullet \xleftarrow{X} \bullet \xrightarrow{Y} \bullet \xleftarrow{X} \dots \xrightarrow{Y} \bullet \xleftarrow{X} \bullet\ ,$$
$$\mathtt{Z}'_{2l}\colon\ \bullet \xrightarrow{Y} \bullet \xleftarrow{X} \bullet \xrightarrow{Y} \bullet \xleftarrow{X} \dots \xrightarrow{Y} \bullet \xleftarrow{X} \bullet \xrightarrow{Y} \bullet\ .$$

- For $j,l \in \mathbb{Z}_{\geq 0}$ with $j|l$, and an irreducible polynomial $f \in \mathbb{F}_2[Z]$ of degree $l/j$, let $\Theta_{f,j,l} = \sum_{i=0}^{n} \Theta_i Z^i = f(Z)^j$ with $\Theta_n = 1$. For all $2l+2$, where $l \in \mathbb{Z}_{\geq 0}$, and any $\Theta = \Theta_{f,j,l}$ there is another self-dual indecomposable $\mathtt{Z}^{\Theta}_{2l}$:

$$\mathtt{Z}^{\Theta}_{2l}\colon\ \underset{0}{\bullet} \xleftarrow{X} \underset{l+1}{\bullet} \xrightarrow{Y} \underset{1}{\bullet} \xleftarrow{X} \underset{l+2}{\bullet} \xrightarrow{Y} \underset{2}{\bullet} \xleftarrow{X} \underset{\dots}{\dots} \xrightarrow{Y} \underset{l}{\bullet} \xleftarrow{X} \underset{2l+1}{\bullet} \overset{Y}{\underset{\Theta}{\dashrightarrow}}\ ,$$

  where at one end, as indicated, $Y$ acts by $Y(2l+1) = \sum_{i=0}^{n} \Theta_i i$.

  (An explicit example of this family of modules is the case $f(Z) = 1+Z+Z^2$, $j=1$ and $l=2$. Then $\mathtt{Z}^{\Theta}_6$ is six dimensional, and of the form

$$\mathtt{Z}^{\Theta}_{6}\colon\ \underset{0}{\bullet} \xleftarrow{X} \underset{3}{\bullet} \xrightarrow{Y} \underset{1}{\bullet} \xleftarrow{X} \underset{4}{\bullet} \xrightarrow{Y} \underset{2}{\bullet} \xleftarrow{X} \underset{5}{\bullet} \overset{Y}{\underset{\Theta}{\dashrightarrow}}\ ,$$

  and $Y$ acts on the vertex 5 as $Y(5) = 1\cdot 0 + 1\cdot 1 + 1\cdot 2$.)

In fact, this is the so-called ***tame (representation) type***; more later. ◇

The summarized above discussion is as follows.

**Proposition 8B.6.** *Let* G *be a finite group of order* $\#\mathrm{G} = m$*, and let* $p$ *be the characteristic of the algebraically closed ground field* $\mathbb{K}$*. Then:*

*(i) We have* $\mathbf{fdMod}(\mathbb{K}[\mathrm{G}]) \in \mathbf{Cat}_{S}$ *if and only if* $p \nmid m$*.*

*(ii) We have* $\mathbf{fdMod}(\mathbb{K}[\mathrm{G}]) \in \mathbf{Ten}$*.*

*(iii) We have* $\mathbf{fdMod}(\mathbb{K}[\mathrm{G}]) \in \mathbf{Fiat}$ *if and only if (*$p \nmid m$ *or a* $p$ *Sylow subgroup of* G *is cyclic).*

*Proof.* The only things we have not addressed above are the following: First, whether $\mathrm{End}_{\mathbf{fdMod}(\mathbb{K}[\mathrm{G}])}(\mathbb{1}) \cong \mathbb{K}$. However, since $\mathbb{1}$ is the trivial module, Schur's lemma Lemma 6J.10 provides the result. And second, the if and only if condition in (iii), which follows from a classical result, giving an if and only if condition for whether $\#\mathrm{In}(\mathbb{K}[\mathrm{G}]) < \infty$, see [**Hig54**]. □

*Remark* 8B.7. The proof in [**Hig54**] is effective: Let H be a $p$ Sylow subgroup of G. Then [**Hig54**] shows that every indecomposable $\mathbb{K}[\mathrm{G}]$ module $\mathtt{Z}$ can be obtained as a direct summands $\mathtt{Z} \Subset \mathbb{K}[\mathrm{G}] \otimes_{\mathbb{K}[\mathrm{H}]} \mathtt{Z}'$ for some indecomposable $\mathbb{K}[\mathrm{H}]$ module $\mathtt{Z}'$. Because the dimension of $\mathbb{K}[\mathrm{G}] \otimes_{\mathbb{K}[\mathrm{H}]} \mathtt{Z}'$ is $|\mathrm{G}/\mathrm{H}| \dim(\mathtt{Z}')$, there can only be finitely many indecomposables if $\mathbb{K}[\mathrm{H}]$ has only finitely many indecomposables. ◇

**Example 8B.8.** The two cases of G being either $\mathbb{Z}/4\mathbb{Z}$ or Klein's four group $\mathrm{V}_4$ in characteristic 2 are fundamentally different: For $\mathbb{Z}/4\mathbb{Z}$ the unique 2 Sylow subgroup is cyclic and the representation theory of $\overline{\mathbb{F}}_2[\mathbb{Z}/4\mathbb{Z}]$ can be treated verbatim as for $\overline{\mathbb{F}}_5[\mathbb{Z}/5\mathbb{Z}]$, see Example 6K.18. For $\mathrm{V}_4$ the unique 2 Sylow subgroup is not cyclic and $\overline{\mathbb{F}}_2[\mathrm{V}_4]$ has infinitely many indecomposables as listed in Example 8B.4. ◇

Recall that a finite dimensional Hopf algebra A is, by definition, a Hopf algebra $\mathtt{A} = (\mathtt{A}, \mathrm{m}, \mathrm{i}, \mathrm{d}, \mathrm{e}, \mathrm{s})$ in $\mathbf{fdVec}_{\mathbb{S}}$, and $\mathbf{fdMod}(\mathtt{A})$ is its module category in $\mathbf{fdVec}_{\mathbb{S}}$. The version for Hopf algebras, where ***finite (representation) type*** means that $\#\mathrm{In}\big(\mathbf{fdMod}(\mathtt{A})\big) < \infty$, is:

**Proposition 8B.9.** *Let* A *be a finite dimensional Hopf algebra. Then:*

(a) *We have* $\mathbf{fdMod}(\mathtt{A}) \in \mathbf{wTen}$.

(b) *We have* $\mathbf{fdMod}(\mathtt{A}) \in \mathbf{wFiat}$ *if and only if* A *is of finite type.*

*Proof.* All discussion above works for finite dimensional Hopf algebras in general, not just for finite groups: By Theorem 5K.6 we know that $\mathbf{fdMod}(\mathtt{A})$ is rigid, and its also clearly $\mathbb{S}$ linear and abelian with $\mathbb{S}$ bilinear $\otimes$. Moreover, (8B-1) holds for any finite dimensional algebra, so $\mathbf{fdMod}(\mathtt{A})$ is always finite, and additively finite if and only if A is of finite type, by definition. Moreover, A has a trivial module $\mathbb{1}$ obtained by using the counit $\mathrm{e}\colon \mathrm{A} \to \mathbb{S}$, and $\mathbb{1}$ is the monoidal unit of $\mathbf{fdMod}(\mathtt{A})$ giving $\mathrm{End}_{\mathbf{fdMod}(\mathtt{A})}(\mathbb{1}) \cong \mathbb{S}$. □

Thus, fiat and tensor categories can be seen as generalizations of Hopf algebras and their representations. However, fiat categories are much more difficult to find.

To illustrate this, we say a finite dimensional algebra A is of ***wild (representation) type*** if, for any other finite dimensional algebra B, there is an exact functor

$$\mathbf{fdMod}(\mathtt{A}) \to \mathbf{fdMod}(\mathtt{B})$$

preserving indecomposable objects. This should be read as: Classifying indecomposable A representations for a wild algebra implies that we can do the same for any finite dimensional algebra (and, of course, this is not possible). This is similar to NP complete while being of finite type is like a problem in P, *cf.* Figure 21.

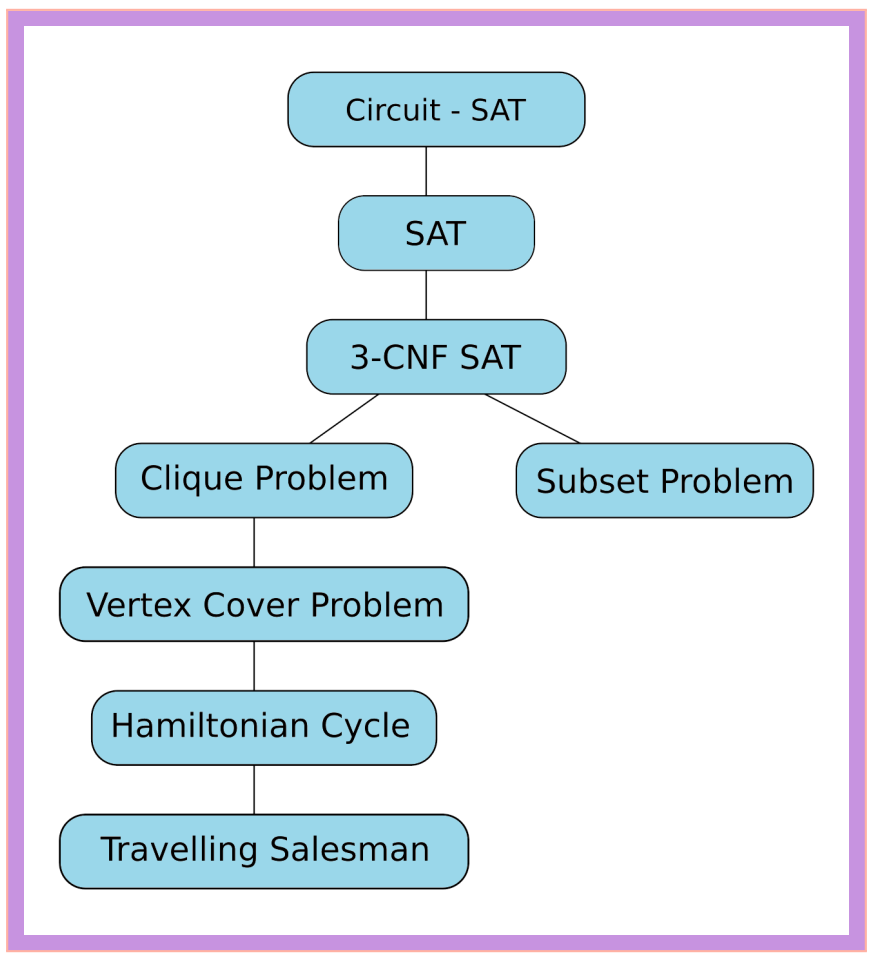

FIGURE 21. Some NP complete problems; they all can be used to simulate one another and also all other NP problems. Under the (almost by everyone accepted) assumption that P is not NP, all of these are "way too difficult" to solve. Wild algebras can be thought of as NP complete problems as they can be used to simulate the classification of indecomposables objects for any finite dimensional algebra.

Picture from https://en.wikipedia.org/wiki/NP-completeness

Let us call a finite dimensional algebra of ***tame (representation) type*** if its not of finite nor of wild type. These are somewhere between P and NP complete problems, which is a repeat of the above analogy. With Proposition 8B.6 in mind, the following statement should be read as "almost everything is wild" when $\mathbf{fdMod}\big(\mathbb{K}[\mathrm{G}]\big) \notin \mathbf{Cat}_{S}$.

**Proposition 8B.10.** *Let* G *be a finite group of order* $\#\mathrm{G} = m$, *and let* $p$ *be the characteristic of the algebraically closed ground field* $\mathbb{K}$ *dividing* $m$. *Then:*

(i) *The algebra* $\mathbb{K}[\mathrm{G}]$ *is of finite type if and only if a* $p$ *Sylow subgroup of* G *is cyclic.*

(ii) *The algebra* $\mathbb{K}[\mathrm{G}]$ *is of tame type if and only if* ($p = 2$ *and a* $p$ *Sylow subgroup of* G *is dihedral, semidihedral or generalized quaternion).*

*Proof.* This can be deduced from [BD77]. □

To see how few groups satisfy, say, the condition in Proposition 8B.10.(i), here is some MAGMA code (see Appendix A for more on MAGMA) that computes the percentage.

```
> // User-defined variables
> max_order := 128;     // Set your desired maximum order here
> prime := 2;          // Set your desired prime here
>
> // Initialize counters
> total_groups := 0;
> cyclic_sylow := 0;
>
> // Loop through all group orders from 1 to max_order
> for n in [1..max_order] do
>     // Get number of groups of order n
>     num_groups := NumberOfSmallGroups(n);
>     total_groups +:= num_groups;
>
>     // Check each group of order n
>     for i in [1..num_groups] do
>         G := SmallGroup(n,i);
>         // Get p-Sylow subgroup for chosen prime
>         S := SylowSubgroup(G, prime);
>         // Check if it's cyclic
>         if IsCyclic(S) then
>             cyclic_sylow +:= 1;
>         end if;
>     end for;
> end for;
>
> // Calculate percentage as a real number
> percentage := (RealField(5)!cyclic_sylow / total_groups) * 100;
>
> // Print results
> total_groups;
> cyclic_sylow;
> percentage;
  ----result----
3596
348
9.6774
```

8C. **Non-negative integral matrices.** The arguably most important numerical invariants associated with a fiat (or tensor category) $\mathbf{C}$ are integral matrices.

**Definition 8C.1.** Let $\mathbf{C} \in \mathbf{wmFiat}$. Then, for $i, j, k \in \{1, ..., n\}$, the ***fusion rules*** and the ***fusion coefficients*** $N_{i,j}^{k} \in \mathbb{Z}_{\geq 0}$ are

$$\mathtt{Z}_i\mathtt{Z}_j \cong \textstyle\bigoplus_{k=1}^{n} N_{i,j}^{k} \cdot \mathtt{Z}_k, \quad \text{where } \mathtt{Z}_l \in \mathrm{In}(\mathbf{C}).$$

(Recall that we write $m \cdot \mathtt{X}$ for $M$ copies $\mathtt{X} \oplus ... \oplus \mathtt{X}$.) ◇

Thus, the fusion coefficients are the structure constants of the $\mathbb{Z}$ algebra $K_0^{\oplus}(\mathbf{C})$. These are most conveniently collected in the ***fusion matrices*** associated to $\mathtt{Z}_j$:

$$K_0^{\oplus}(_ \otimes \mathtt{Z}_j) = M(j) = (N_{i,j}^{k})_{i,k=1}^{n} = \begin{array}{c} \\ \mathtt{Z}_1 \\ \vdots \\ \mathtt{Z}_k \\ \vdots \\ \mathtt{Z}_n \end{array} \begin{array}{c} \mathtt{Z}_1 \ \dots \ \mathtt{Z}_i \ \dots \ \mathtt{Z}_n \\ \begin{pmatrix} & & \\ & & \\ & N_{i,j}^{k} & \\ & & \\ & & \end{pmatrix} \\ \mathtt{Z}_i\mathtt{Z}_j \end{array} \in \mathrm{Mat}_{n \times n}(\mathbb{Z}_{\geq 0}).$$

In words, the fusion matrix $M(j)$ captures the right $\otimes$ action of $\mathtt{Z}_j$ on $\mathbf{C}$. Recall further that we can associate a graph $\Gamma(M)$ with $n$ vertices to each matrix $M \in \mathrm{Mat}_{n\times n}(\mathbb{Z}_{\geq 0})$, see Convention 8A.1, which we identify with $M$. Thus, we have another numerical invariant of fiat categories that captures all the fusion rules:

**Definition 8C.2.** Let $\mathbf{C} \in \mathbf{wmFiat}$. Then, for $i \in \{1,...,n\}$, the ***fusion graphs*** are the directed graphs $\Gamma_i = \Gamma\big(M(i)\big)$, *i.e.* the graphs associated with the fusion matrices. Similarly, the ***fusion graph of*** $\mathtt{X} \in \mathbf{C}$ is the directed graph $\Gamma_{\mathtt{X}}$ associated with the right $\otimes$ action of $\mathtt{X}$ on $\mathbf{C}$. ◇

The following is evident, where the sum of graphs is the graph that one obtains by summing the corresponding matrices, using the identification of these, *cf.* Convention 8A.1.

**Lemma 8C.3.** *Let* $\mathbf{C} \in \mathbf{wmFiat}$. *If* $\mathtt{X} \in \mathbf{C}$ *decomposes as* $\mathtt{X} \cong \bigoplus_{i=1}^{n} (\mathtt{X} : \mathtt{Z}_i) \cdot \mathtt{Z}_i$, *then we have* $\Gamma_{\mathtt{X}} = \sum_{i=1}^{n} (\mathtt{X} : \mathtt{Z}_i) \cdot \Gamma_i$. □

The fusion graphs are invariants:

**Proposition 8C.4.** *Let* $\mathrm{F} \in \mathbf{Hom}_{\Bbbk\oplus\star}(\mathbf{C},\mathbf{D})$ *be an equivalence of categories* $\mathbf{C},\mathbf{D} \in \mathbf{wmFiat}$. *Then, up to reordering, the fusion graphs of* $\mathbf{C}$ *and* $\mathbf{D}$ *are isomorphic as graphs.*

*Proof.* By Proposition 7H.3. □

Note that the fusion graph $\Gamma_{\mathbb{1}}$ associated with the monoidal unit $\mathbb{1}$ is always a completely disconnected graph with one loop per vertex, *e.g.*

$$\Gamma_{\mathbb{1}} = \begin{pmatrix} 1 & 0 & 0 & 0 \\ 0 & 1 & 0 & 0 \\ 0 & 0 & 1 & 0 \\ 0 & 0 & 0 & 1 \end{pmatrix} =$$

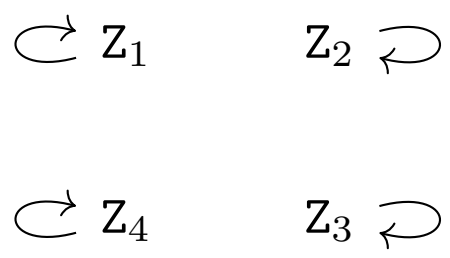

We call these ***trivial fusion graphs***, and all the others ***non-trivial***. Moreover, in tables, we omit the row and column for the fusion rules of $\mathbb{1}$ as they are trivial, see (8C-7).

**Example 8C.5.** Let us consider two examples of semisimple fiat categories and their fusion graphs:

(a) The category $\mathbf{fdMod}\big(\mathbb{C}[\mathbb{Z}/4\mathbb{Z}]\big)$ has simples

$$\mathtt{L}_1 = \mathbb{1}\colon \big(\boxed{1}\big),\quad \mathtt{L}_2\colon \big(\boxed{i}\big),\quad \mathtt{L}_3\colon \big(\boxed{-1}\big),\quad \mathtt{L}_4\colon \big(\boxed{-i}\big),$$

which act on $\mathbf{fdMod}\big(\mathbb{C}[\mathbb{Z}/4\mathbb{Z}]\big)$ as the elements, in order, 0, 1, 2 and 3 in $\mathbb{Z}/4\mathbb{Z}$. Hence, the non-trivial fusion graphs are

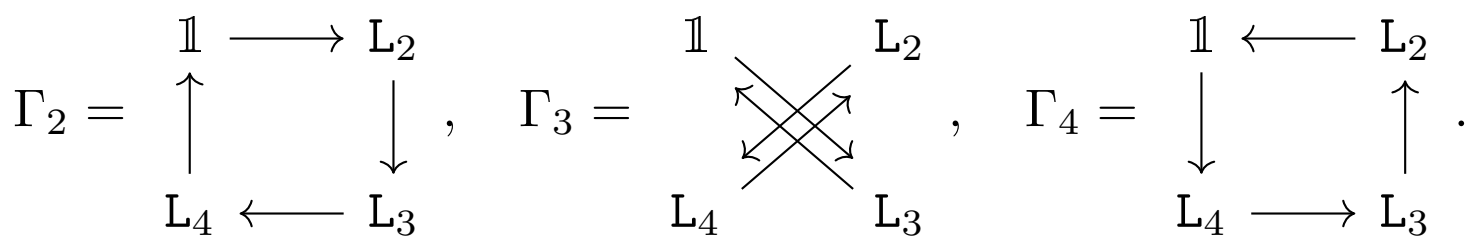

(Careful with the index shift here. For example, $\mathtt{L}_2$ acts as 1 when seen in $\mathbb{Z}/4\mathbb{Z}$.)

(b) The category $\mathbf{fdMod}\big(\mathbb{C}[\mathrm{V}_4]\big)$ (Klein's four group, see Example 8B.4) has simples also corresponding to the elements 1, $s$, $t$ and $st = ts$ in $\mathrm{V}_4$. Hence, the non-trivial fusion graphs are

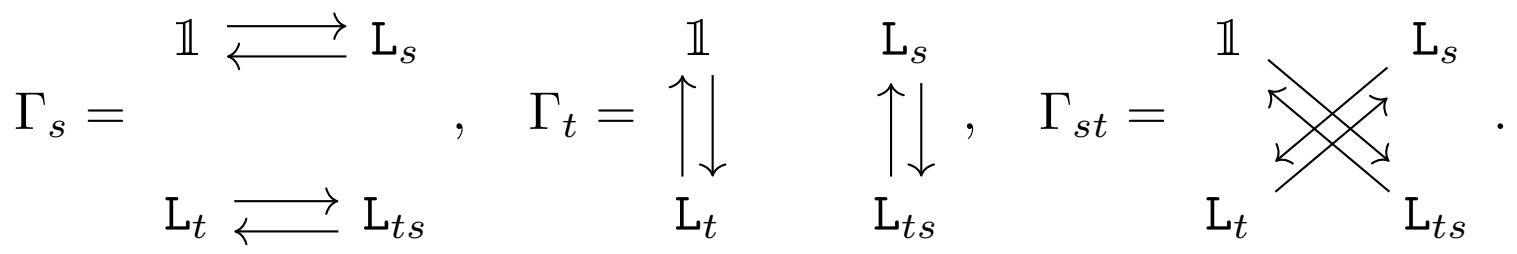

Thus, although $\mathbf{fdMod}\big(\mathbb{C}[\mathbb{Z}/4\mathbb{Z}]\big)$ and $\mathbf{fdMod}\big(\mathbb{C}[\mathrm{V}_4]\big)$ are equivalent as categories, they are not equivalent as fiat categories. ◇

**Example 8C.6.** Let $\mathrm{S}_3$ be the symmetric group of order 6 in three letters. The category $\mathbf{fdMod}\big(\mathbb{C}[\mathrm{S}_3]\big)$ is a semisimple fiat category with simples $\mathtt{L}_1 = \mathbb{1}$, $\mathtt{L}_s$ and $\mathtt{L}_{-1}$ satisfying the following fusion rules:

(8C-7)

| $\otimes$ | $\mathtt{L}_s$ | $\mathtt{L}_{1'}$ |
|---|---|---|
| $\mathtt{L}_s$ | $\mathtt{L}_1 \oplus \mathtt{L}_s \oplus \mathtt{L}_{1'}$ | $\mathtt{L}_s$ |
| $\mathtt{L}_{1'}$ | $\mathtt{L}_s$ | $\mathtt{L}_1$ |

$=$

| $\otimes$ | $\mathtt{L}_1$ | $\mathtt{L}_s$ | $\mathtt{L}_{1'}$ |
|---|---|---|---|
| $\mathtt{L}_1$ | $\mathtt{L}_1$ | $\mathtt{L}_s$ | $\mathtt{L}_{1'}$ |
| $\mathtt{L}_s$ | $\mathtt{L}_s$ | $\mathtt{L}_1 \oplus \mathtt{L}_s \oplus \mathtt{L}_{1'}$ | $\mathtt{L}_s$ |
| $\mathtt{L}_{1'}$ | $\mathtt{L}_{1'}$ | $\mathtt{L}_s$ | $\mathtt{L}_1$ |

Thus, we get the non-trivial fusion graphs

$$\Gamma_s = \begin{pmatrix} 0 & 1 & 0 \\ 1 & 1 & 1 \\ 0 & 1 & 0 \end{pmatrix} =$$

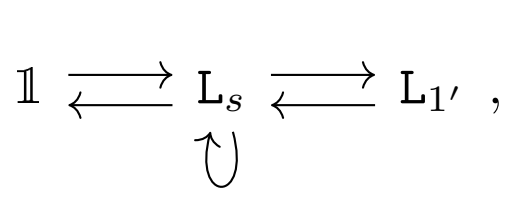

$$\Gamma_{1'} = \begin{pmatrix} 0 & 0 & 1 \\ 0 & 1 & 0 \\ 1 & 0 & 0 \end{pmatrix} =$$

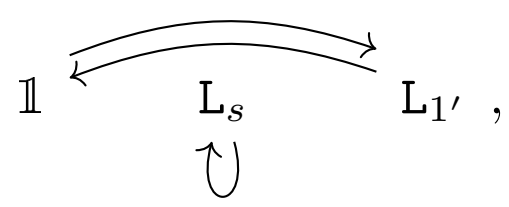

which can be directly read of from Equation 8C-7. ◇

**Example 8C.8.** The fusion graphs are certainly not a complete invariant of fiat categories as they do not involve the morphisms in any form. To be completely explicit, the categories of the form $\mathbf{Vec}^{\omega}_{\Bbbk\oplus}(\mathrm{G})$, for G being a finite group, are fiat categories with the same fusion graphs, independent of $\omega$. ◇

Recall that ***strongly connected*** for graphs means connected as a directed graph.

**Definition 8C.9.** Let $\mathbf{C} \in \mathbf{wmFiat}$. We call $\mathtt{X} \in \mathbf{C}$ a ***fusion generator of* C** if $\Gamma_{\mathtt{X}}$ is strongly connected. In case $\mathbf{C} \in \mathbf{wmFiat}$ has a fusion generator, we call **C** ***transitive***. ◇

**Example 8C.10.** In Example 8C.5.(a) both, $\mathtt{L}_1$ and $\mathtt{L}_3$, are fusion generators, but $\mathtt{L}_2$ is not a fusion generator. In Example 8C.5.(b) none of the simples are fusion generators, but

$$\mathtt{L}_s \oplus \mathtt{L}_t \rightsquigarrow \Gamma_{\mathtt{L}_s \oplus \mathtt{L}_t} = \begin{pmatrix} 0 & 1 & 1 & 0 \\ 1 & 0 & 0 & 1 \\ 1 & 0 & 0 & 1 \\ 0 & 1 & 1 & 0 \end{pmatrix} = \begin{array}{ccc} \mathbb{1} & \rightleftarrows & \mathtt{L}_s \\ \uparrow\downarrow & & \uparrow\downarrow \\ \mathtt{L}_t & \rightleftarrows & \mathtt{L}_{ts} \end{array}$$

is a fusion generator. ◇

The term "generator" is to be understood in this sense:

**Lemma 8C.11.** *Let* $\mathbf{C} \in \mathbf{wmFiat}$ *be transitive,* $\mathtt{X} \in \mathbf{C}$ *be a fusion generator and* $\mathtt{Y} \in \mathbf{C}$ *any object. Then there exists* $k \in \mathbb{Z}_{\geq 0}$ *such that* $\mathtt{Y} \Subset \mathtt{X}^k$.

*Proof.* If $\mathtt{X}$ and $\mathtt{Y}$ are indecomposable, then $k$ can be taken to be the length of a shortest path in $\Gamma_{\mathtt{X}}$ from the vertex corresponding to $\mathtt{X}$ to the vertex corresponding to $\mathtt{Y}$. For general $\mathtt{X}$ and $\mathtt{Y}$ the claim thus follows by additivity. □

For our favorite example, $\mathbf{fdMod}\big(\mathbb{C}[\mathrm{G}]\big)$ for a finite group G, it is easy to decide whether an object is a fusion generator:

**Proposition 8C.12.** *Let* G *be a finite group. Then* $\mathtt{X} \in \mathbf{fdMod}\big(\mathbb{C}[\mathrm{G}]\big)$ *is a fusion generator if and only if* $\dim_{\mathbb{C}} \mathtt{X}$ *appears only once in its character (= row in the character table) if and only if* $\mathtt{X}$ *is faithful.*

*Proof.* The first ⇔ is the Burnside–Brauer–Steinberg theorem, see, for example, [Ste14] for a more general version of this theorem, and a readable account at the same time. For the final ⇔ see Remark AI.1. □

**Example 8C.13.** Let us come back to Example 8C.6. The character table of $\mathrm{S}_3$ is:

```
> CharacterTable(Sym(3));
  ----result----
> Character Table of Group G
>  --------------------------
>
>
> -----------------
> Class |   1  2  3
> Size  |   1  3  2
> Order |   1  2  3
> -----------------
> p  =  2   1  1  3
> p  =  3   1  2  1
> -----------------
> X.1   +   1  1  1
> X.2   +   1 -1  1
> X.3   +   2  0 -1
```

see Appendix A how to generate and read such tables. The first column contains the dimensions, and the order is $\mathtt{L}_1$ (top), $\mathtt{L}_{1'}$ (middle), and $\mathtt{L}_s$ (bottom). So in order to check whether any of these $\mathtt{L}_*$ is a fusion generator, it suffices to check whether the first column entry appears again in its row: if it does, then $\mathtt{L}_*$ is not a fusion generator. Thus, we immediately see that the only fusion generator among these is $\mathtt{L}_s$. ◇

8D. **Perron–Frobenius (PF) theory.** We start with the main player, the ***PF theorem***. To this end, let us look at the example in Figure 22. There is a clear pattern that emerges: There is a leading eigenvalue (which is much bigger than all others) and its eigenvector is nonnegative.

We will now see that this is true in general, which is the content of the PF theorem. Recall that one can associate an oriented and weighted graph, its ***adjacency graph***, to an $n$-by-$n$ matrix $M = (m_{ij})_{1 \leq i,j \leq n} \in \mathrm{Mat}_{n \times n}(\mathbb{R}_{\geq 0})$ as follows:

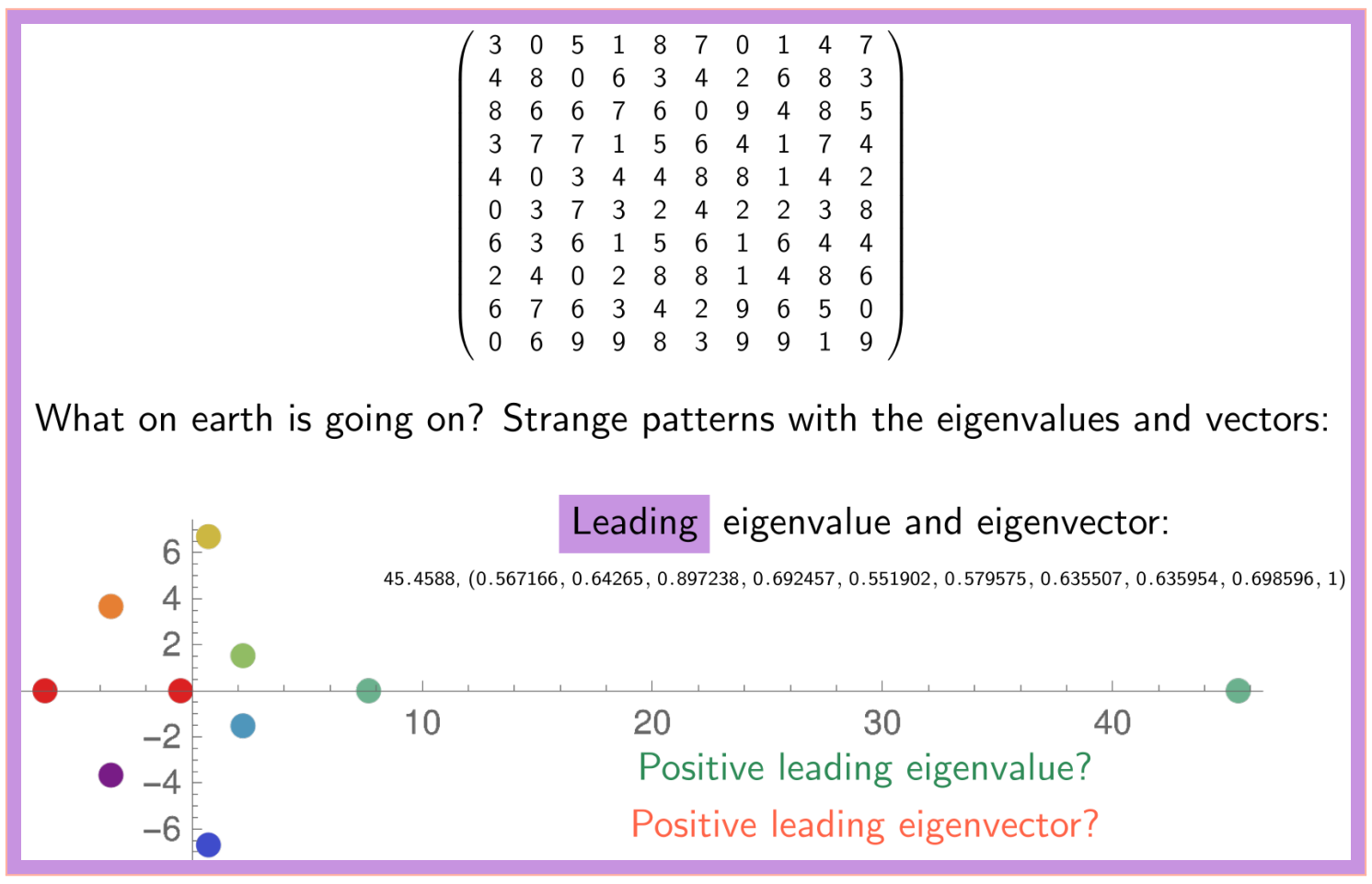

FIGURE 22. A randomly generated 10-by-10 matrix with nonnegative entries and a plot of its eigenvalues and the eigenvector for the rightmost eigenvalue.

Picture from https://www.youtube.com/watch?v=SJ0A0ewTxg4

(i) The vertices are $\{1,\dots,n\}$.

(ii) There is an edge with weight $m_{ij}$ from $i$ to $j$.

This is similar to Convention 8A.1, but now we put weights (edge labels) on the edges.

We call a nonzero matrix $M\in\mathrm{Mat}_{n\times n}(\mathbb{R}_{\geq 0})$ ***strongly connected*** (this is often called irreducible in the literature) if its associated graph is strongly connected. Finally, in this note, a right eigenvector satisfies $Mv=\lambda\cdot v$, and a left eigenvector satisfies $w^{T}M=\lambda\cdot w^{T}$.

**Theorem 8D.1** (PF theorem part I)**.** *Let $M\in\mathrm{Mat}_{n\times n}(\mathbb{R}_{\geq 0})$ be strongly connected.*

**(a)** *$M$ has a* **PF eigenvalue***, that is, $\lambda=\lambda_{pf}(M)\in\mathbb{R}_{>0}$ such that $\lambda\geq|\mu|$ for all other eigenvalues $\mu$. This eigenvalue appears with multiplicity one, and all other eigenvalues with $\lambda=|\mu|$ also appear with multiplicity one.*

**(b)** *There exists $h\in\mathbb{Z}_{\geq 1}$, the* **period***, such that all eigenvalues $\mu$ with $\lambda=|\mu|$ are $\exp(k2\pi i/h)\lambda$ for $k\in\{0,\dots,h-1\}$. We call these* **pseudo dominant eigenvalues***.*

**(c)** *The eigenvectors, left and right, for the PF eigenvalue can be normalized to have entries in $\mathbb{R}_{\geq 0}$.*

*Proof.* A well-known theorem. See for example, Frobenius' paper 92 in Band 3 of [Fro68]. (This is the paper "Über Matrizen aus nicht negativen Elementen".) We sketch a proof of the following variation:

**Lemma 8D.2.** *Let $M\in\mathrm{Mat}_{n\times n}(\mathbb{R}_{\geq 0})$. Then:*

*(i) The matrix $M$ has an eigenvalue $\lambda\in\mathbb{R}_{\geq 0}$ which satisfies*

$$|\lambda|\geq|\mu|,\quad \textit{for all eigenvalues } \mu \textit{ of } M. \tag{8D-3}$$

*Moreover, $M$ has a right eigenvector $v=v_{pf}(M)$ with eigenvalue $\lambda_{pf}(M)$, which can be normalized such that $v\in\mathbb{R}^{n}_{\geq 0}$.*

*(ii) If additionally $M\in\mathrm{Mat}_{n\times n}(\mathbb{R}_{>0})$, then $\lambda\in\mathbb{R}_{>0}$, $v\in\mathbb{R}^{n}_{>0}$ (after normalization), the eigenvalue $\lambda$ is simple, and the inequality in (8D-3) is strict. Further, the eigenvector $v$ is the unique (up to scaling) eigenvector of $M$ with values in $\mathbb{R}_{>0}$.*

$\lambda$ is called the ***PF eigenvalue*** of $M$, and $v$ is called the ***PF eigenvector*** of $M$.

*Proof.* (i)-Existence. The idea is to use Brouwer's fixed-point theorem. Assume that $M$ does not have an eigenvector $v_0(M)\in\mathbb{R}^{n}_{\geq 0}$ of eigenvalue zero. Then, since there are no cancellations due to non-negativity,

$$f\colon\Sigma_n\to\Sigma_n,\quad v\mapsto\frac{Mv}{\sum_{i=1}^{n}(Mv)_i}$$

defines a continuous map from the standard $n$ simplex $\Sigma_n=\{v\in\mathbb{R}^{n}_{\geq 0}\mid\sum_{i=1}^{n}v_i=1\}$ to itself. Thus, Brouwer's fixed-point theorem gives us a fixed point $w$ of $f$, which, by construction, satisfies

$$Mw=\mu w,\quad\text{where } \mu\in\mathbb{R}_{\geq 0}, w\in\mathbb{R}^{n}_{\geq 0}.$$

We can hence define $\lambda_{pf}(M)$ to be the maximal eigenvalue of $M$ having a non-negative eigenvector $v_{pf}(M)$. Since this also works in case $M$ has an eigenvector $v_0(M)\in\mathbb{R}^{n}_{\geq 0}$ of eigenvalue zero, we have now constructed the required eigenvalue and eigenvector, and it remains to show the claimed properties.

(ii)-Positivity. If $M \in \mathrm{Mat}_{n\times n}(\mathbb{R}_{>0})$, then the just constructed eigenvalue $\lambda_{pf}(M)$ and eigenvector $v_{pf}(M)$ are also strictly positive.

(ii)-Simplicity. The eigenvalue $\lambda_{pf}(M)$ is simple: If $w \in \mathbb{R}^n$ is another eigenvector of $\lambda_{pf}(M)$, then define $z = \min\{w_i - v_{pf}(M)_i \mid i = 1, ..., n\}$. Now we observe that $w - z \cdot v_{pf}(M) \in \mathbb{R}^n_{\geq 0}$ is another eigenvector of $M$ with eigenvalue $\lambda_{pf}(M)$. However, this, by positivity, implies that $w - z \cdot v_{pf}(M) = 0$ as it has at least one entry being zero. Thus, $\lambda_{pf}(M)$ is simple.

(ii)-Uniqueness. Assume now that there exists another strictly positive eigenvector $w$ with eigenvalue $\mu$, and let $v_{pf}(M^T)$ denote the PF eigenvector of $v_{pf}(M^T)$. Then $\lambda_{pf}(M)v_{pf}(M^T)w = v_{pf}(M^T)Mw = \mu v_{pf}(M^T)w$, which, by positivity, implies that $\mu = \lambda_{pf}(M)$. Hence, we also have $w = v_{pf}(M^T)$, by simplicity.

(ii)-Inequality. For $w \in \mathbb{C}^n$ let $|w| = \sum_{i=1}^{n} |w_i| v_{pf}(M)_i$. Then one checks that $|Mw| \leq \lambda_{pf}(M)|w|$ and equality holds if and only if all nonzero entries of $w$ have the same argument. Hence, if $w$ is an eigenvector of $M$ with eigenvalue $\mu$, then $|\mu||w| \leq \lambda_{pf}(M)|w|$, which implies that $|\mu| \leq \lambda_{pf}(M)$. Finally, if $|\mu| = \lambda_{pf}(M)$, then all nonzero entries of $w$ have the same arguments and we can normalize $w$ to be strictly positive. Hence, $\mu = \lambda_{pf}(M)$, by uniqueness.

(i)-Rest. Using (ii) and

$$M_N = M + \tfrac{1}{N}\begin{pmatrix} 1 & \dots & 1 \\ \vdots & \ddots & \vdots \\ 1 & \dots & 1 \end{pmatrix} \in \mathrm{Mat}_{n\times n}(\mathbb{R}_{>0}), \qquad \lim_{N\to\infty} M_N = M,$$

this is an easy limit argument. □

We leave it to the reader to fill in the remaining points. □

*Remark* 8D.4. Note that Brouwer's fixed-point theorem does not actually construct its promised fix point, and the proof of the PF theorem above in turn does not give an eigenvalue. If one works over $\mathbb{Q}_{\geq 0}$ instead of $\mathbb{R}_{\geq 0}$, then there is an alternative and constructive proof of the PF theorem that avoids Brouwer's fixed-point theorem. (This is a consequence of Barr's theorem [**Bar74**].) ◇

**Example 8D.5.** Note the difference between strictly positive and positive, and various other properties as *e.g.* "converting to zero" or "nilpotent":

$$M_1 = \begin{pmatrix} \frac{1}{3} & \frac{1}{3} \\ \frac{1}{3} & \frac{1}{3} \end{pmatrix}, \quad M_2 = \begin{pmatrix} \frac{1}{3} & 0 \\ \frac{1}{3} & \frac{1}{3} \end{pmatrix}, \quad M_3 = \begin{pmatrix} 0 & 0 \\ \frac{1}{3} & 0 \end{pmatrix},$$

Let us call them "Case 1", "Case 2" and "Case 3", respectively. Then:

$$\text{Case 1}\colon \begin{cases} \lambda_{pf} = \frac{2}{3}, & v_{pf} = (1,1), \\ \mu_1 = 0, & v_1 = (-1,1), \end{cases} \quad \text{Case 2}\colon \begin{cases} \lambda_{pf} = \frac{1}{3}, \\ v_{pf} = (0,1), \end{cases} \quad \text{Case 3}\colon \begin{cases} \lambda_{pf} = 0, \\ v_{pf} = (0,1). \end{cases}$$

In terms of graphs, Case 1 is strongly connected (in particular, there is a directed cycle), Case 2 is not strongly connected but has a directed cycle, and Case 3 is neither strongly directed nor has a directed cycle. ◇

Fix a function $f\colon \mathbb{Z}_{\geq 0} \to \mathbb{Z}_{\geq 0}$. We say that $f(n)$ converges ***geometrically to*** $a \in \mathbb{R}$ ***with ratio*** $\beta \in [0,1)$ if for all $\gamma \in (\beta, 1)$ we have that $\{(f(n) - a)/\gamma^n\}_{n\in\mathbb{Z}_{\geq 0}}$ is bounded. We, abusing language, will call the infimum of all ratios ***the*** ratio of convergence.

For two matrices of the same size let $\sim_{ew}$ mean that they are asymptotically equal entrywise (note that below $v_i w_i^T$ are matrices). For such matrices we apply the definition of geometric convergence entrywise with the ratio being the maximum of the entrywise ratios. The following accompanies Theorem 8D.1:

**Theorem 8D.6** (PF theorem part II)**.** *Let $M \in \mathrm{Mat}_{n\times n}(\mathbb{R}_{\geq 0})$ be strongly connected, $\lambda$ be its PF eigenvalue and $h$ be its period. Let $\zeta = \exp(2\pi i/h)$. For each $k \in \{0, ..., h-1\}$, choose a left eigenvector $v_k$ and a right eigenvector $w_k$ with eigenvalue $\zeta^k\lambda$, normalized such that $w_k^T v_k = 1$.*

*Then we have:*

$$M^n \sim_{ew} v_0 w_0^T \cdot \lambda^n + v_1 w_1^T \cdot (\zeta\lambda)^n + v_2 w_2^T \cdot (\zeta^2\lambda)^n + \dots + v_{h-1} w_{h-1}^T \cdot (\zeta^{h-1}\lambda)^n.$$

*Moreover, the convergence is geometric with ratio $|\lambda^{sec}/\lambda|$, where $\lambda^{sec}$ is any second largest (in the sense of absolute value) eigenvalue.*

*Proof.* This is known, but proofs are a bit tricky to find in the literature, so we give one. The proof also shows where the vectors $v_i$ and $w_i$ come from.

For any $\mu \in \mathbb{C}$, let $V_\mu$ be the generalized eigenspace of $V = \mathbb{C}^m$ associated with the eigenvalue $\mu$. Then we have

$$V = \bigoplus_{k=0}^{h} V_{\zeta^k\lambda} \oplus \bigoplus_{\mu, |\mu|<\lambda} V_\mu.$$

By Theorem 8D.1, the space $V_{\zeta^k\lambda}$ is the eigenspace associated with the eigenvalue $\zeta^k\lambda$ and $v_k w_k^T$ is the projection onto that subspace.

This implies that we have

$$M^n = v_0 w_0^T \cdot \lambda^n + v_1 w_1^T \cdot (\zeta\lambda)^n + v_2 w_2^T \cdot (\zeta^2\lambda)^n + \ldots + v_{h-1} w_{h-1}^T \cdot (\zeta^{h-1}\lambda)^n + R(n),$$

where $R(n)$ is the multiplication action of $M^n$ onto the rest. Since the eigenvalues $\mu$ of $M$ on the rest satisfies $|\mu| < \lambda$, we have $R(n)/\lambda^n \to_{n\to\infty} 0$ geometrically with ratio $|\lambda^{sec}/\lambda|$. □

For a general matrix $M \in \mathrm{Mat}_{n\times n}(\mathbb{R}_{\geq 0})$ things change, but not too much:

**Theorem 8D.7** (PF theorem part III)**.** *Let $M \in \mathrm{Mat}_{n\times n}(\mathbb{R}_{\geq 0})$.*

**(a)** *$M$ has a* **PF eigenvalue***, that is, $\lambda \in \mathbb{R}_{\geq 0}$ such that $\lambda \geq |\mu|$ for all other eigenvalues $\mu$.*

**(b)** *Let $s$ be the multiplicity of the PF eigenvalue. There exists $(h_1, \ldots, h_s) \in \mathbb{Z}_{\geq 1}^s$, the* **periods***, such that all eigenvalues $\mu$ with $\lambda = |\mu|$ are $\exp(k2\pi i/h_\ell)\lambda$ for $k \in \{0, \ldots, h_\ell - 1\}$, for some period. We call these* **pseudo dominant eigenvalues***.*

**(c)** *The eigenvectors, left and right, for the PF eigenvalues can be normalized to have entries in $\mathbb{R}_{\geq 0}$.*

**(d)** *Let $h = \mathrm{lcm}(h_1, \ldots, h_s)$, and let $\nu$ the maximal dimension of the Jordan blocks of $M$ containing $\lambda$. There exist matrices $S^i(n)$ with polynomial entries of degree $\leq (\nu - 1)$ for $i \in \{0, \ldots, h-1\}$ such that*

$$\lim_{n\to\infty} |\,(M/\lambda)^{hn+i} - S^i(n)| \to 0 \quad \forall i \in \{0, \ldots, h-1\},$$

*and the convergence is geometric with ratio $|\lambda^{sec}/\lambda|^h$. There are also explicit formulas for the matrices $S^i(n)$, see [**Rot81**, Section 5].*

*Proof.* This can be found in [**Rot81**]. See also [**Hog07**, Section I.10] (in the second version) for a useful list of properties of nonnegative matrices. □

8E. **PF theory in categories.** For categories everything is integral, and here is the version of the PF theorem for positive integral matrices, a.k.a. graphs, with slightly better statements:

**Theorem 8E.1.** *Let $\Gamma \in \mathrm{Mat}_{n\times n}(\mathbb{Z}_{\geq 0})$. Then the PF theorem Lemma 8D.2 applies and we have additionally:*

*(i) If $\Gamma$ has a directed cycle, then $\lambda \in \mathbb{R}_{\geq 1}$.*

*(ii) If $\Gamma$ is strongly connected, then $\lambda \in \mathbb{R}_{\geq 1}$, $v \in \mathbb{R}^n_{>0}$ (after normalization), and all $\mu$ satisfying equality in (8D-3) are not in $\mathbb{R}_{>0}$. Further, the only strictly positive eigenvectors of $M$ are of eigenvalue $\lambda_{pf}(M)$.*

*Proof.* Note that $\Gamma^k$ counts the number of paths of length $k$ in $\Gamma$.

(i). As $\Gamma$ has an oriented cycle, $\Gamma^k \neq 0$ for all $k \in \mathbb{Z}_{\geq 0}$. Moreover, since $\Gamma$ has entries from $\mathbb{Z}_{\geq 0}$ we also know that $\lim_{k\to\infty} \Gamma^k \neq 0$. This implies that the PF eigenvalue has to be at least 1.

(ii). We get a strictly positive matrix

$$T = f(\Gamma) = \textstyle\sum_{i=0}^n \Gamma^i \in \mathrm{Mat}_{n\times n}(\mathbb{Z}_{>0}), \quad \text{where } f(X) = \textstyle\sum_{i=0}^n X^i.$$

Hence, we can apply Lemma 8D.2.(ii) to $T$, and simultaneously Lemma 8D.2.(i) to $\Gamma$. These imply that $\sum_{i=1}^n \lambda_{pf}(\Gamma)^i = \lambda_{pf}(T)$, thus $\lambda_{pf}(\Gamma) \in \mathbb{R}_{\geq 1}$, and moreover $v_{pf}(\Gamma) = v_{pf}(T) \in \mathbb{R}^n_{>0}$. The other claims follow by observing that the non-zero roots of $f(X)$ are the $n$th complex roots of unity. □

**Example 8E.2.** Let $\Gamma_1, \Gamma_s \in \mathrm{Mat}_{3\times 3}(\mathbb{Z}_{\geq 0})$ be

$$\Gamma_1 = \begin{pmatrix} 0 & 0 & 1 \\ 1 & 0 & 0 \\ 0 & 1 & 0 \end{pmatrix} =$$

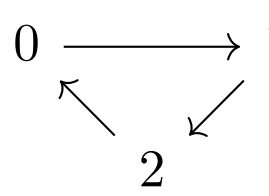

$$, \quad \Gamma_s = \begin{pmatrix} 0 & 1 & 0 \\ 1 & 1 & 1 \\ 0 & 1 & 0 \end{pmatrix} =$$

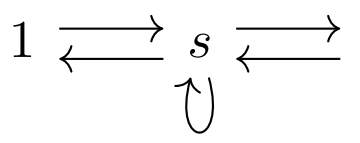

.

(These are action matrices of fusion generators of $\mathbf{fdMod}\big(\mathbb{C}[\mathbb{Z}/3\mathbb{Z}]\big)$ and $\mathbf{fdMod}\big(\mathbb{C}[\mathrm{S}_3]\big)$, respectively.) Let us call them "Case 1" and "Case 2". The eigenvalues and eigenvectors in these two cases are:

$$\text{Case 1:} \begin{cases} \lambda_{pf} = 1, & v_{pf} = (1,1,1), \\ \mu_1 = \tfrac{1}{2}(-1 + i\sqrt{3}), & v_1 = \big(\tfrac{1}{2}(-1 + i\sqrt{3}), \tfrac{1}{2}(-1 - i\sqrt{3}), 1\big), \\ \mu_2 = \tfrac{1}{2}(-1 - i\sqrt{3}), & v_2 = \big(\tfrac{1}{2}(-1 - i\sqrt{3}), \tfrac{1}{2}(-1 + i\sqrt{3}), 1\big), \end{cases}$$

$$\text{Case 2:} \begin{cases} \lambda_{pf} = 2, & v_{pf} = (1,2,1), \\ \mu_1' = -1, & v_1 = (1,-1,1), \\ \mu_2' = 0, & v_2 = (-1,0,1). \end{cases}$$

Note that, in Case 1, $|\mu_1| = |\mu_2| = \lambda_{pf} = 1$, but neither $\mu_1$ nor $\mu_2$ are real numbers. ◇

By Proposition 8C.4, we get the following invariants of fiat categories.

**Definition 8E.3.** Let $\mathbf{C} \in \mathbf{wmFiat}$ and $\mathtt{X} \in \mathbf{C}$. The ***PF dimension of*** $\mathtt{X}$ is $\mathrm{PFdim}(\mathtt{X}) = \lambda_{pf}(\Gamma_{\mathtt{X}})$. The ***PF dimension of*** $\mathbf{C}$ is $\mathrm{PFdim}(\mathbf{C}) = \sum_{i=1}^{n} \mathrm{PFdim}(\mathtt{Z}_i)^2$. ◇

Note that we always have $\mathrm{PFdim}(\mathbb{1}) = 1$. (We will omit this case from the examples.) However, PF dimensions need not to be integers:

**Example 8E.4.** There exists a semisimple fiat category **Fib**, called ***Fibonacci category***, which has two simple objects $\mathtt{L}_1 = \mathbb{1}$ and $\mathtt{L} = \mathtt{L}_2$ with

$$\mathtt{L}^2 \cong \mathbb{1} \oplus \mathtt{L}.$$

Thus, letting $\phi = \frac{1}{2}(1 + \sqrt{5})$ denote the golden ratio, we get

$$\Gamma_{\mathtt{L}} = \begin{pmatrix} 0 & 1 \\ 1 & 1 \end{pmatrix} = \mathbb{1} \rightleftarrows \mathtt{L} \circlearrowleft \Rightarrow \begin{cases} \lambda_{pf}(\Gamma_{\mathtt{L}}) = \mathrm{PFdim}(\mathtt{L}) = \phi, \\ v_{pf}(\Gamma_{\mathtt{L}}) = (1, \phi) = \phi(\phi - 1, 1), \\ \mathrm{PFdim}(\mathbf{Fib}) = 1 + \phi^2 = \frac{1}{2}(5 + \sqrt{5}). \end{cases}$$

We construct this category explicitly in Section 8F below. ◇

**Lemma 8E.5.** *Let* $\mathbf{C} \in \mathbf{wmFiat}$. *Then:*

*(i)* *For* $\mathtt{X}, \mathtt{Y} \in \mathbf{C}$ *we have* $\mathrm{PFdim}(\mathtt{XY}) = \mathrm{PFdim}(\mathtt{X})\mathrm{PFdim}(\mathtt{Y})$.

*(ii)* *If* $\mathtt{X} \in \mathbf{C}$ *is invertible (see Definition 4E.1), then* $\mathrm{PFdim}(\mathtt{X}) = 1$.

*(iii)* *For* $\mathtt{X} \in \mathbf{C}$ *we have* $\mathrm{PFdim}(\mathtt{X}^{\star}) = \mathrm{PFdim}(\mathtt{X}) = \mathrm{PFdim}({}^{\star}\mathtt{X})$. *Moreover, all action matrices,* $\Gamma_{\mathtt{X}^{\star}}$, $\Gamma_{\mathtt{X}}$ *and* $\Gamma_{{}^{\star}\mathtt{X}}$, *agree up to transposition and permutation.*

*(iv)* *The self-dual object* $\mathtt{T} \in \mathbf{C}$, *called the* **total object**, *defined by*

$$\mathtt{T} = \textstyle\sum_{i=1}^{n} \mathtt{Z}_i,$$

*is a fusion generator of* $\mathbf{C}$ *if and only if* $\mathbf{C}$ *is transitive.*

*(v)* *If* $\mathbf{C}$ *is transitive, then there exists a strictly positive virtual object* $\mathtt{R} \in \mathbf{C}$ *(meaning a formal* $\mathbb{R}_{>0}$ *linear combination of indecomposables), called the* **regular object**, *which is the, up to scaling, unique object satisfies the equality*

$$[\mathtt{XR}] = [\mathtt{RX}] = \mathrm{PFdim}(\mathtt{X}) \cdot [\mathtt{R}], \quad (\textit{in } K_0^{\oplus}(\mathbf{C}) \otimes_{\mathbb{Z}} \mathbb{C}),$$

*for all* $\mathtt{X} \in \mathbf{C}$.

*(vi)* *We have* $\mathrm{PFdim}(\mathbf{C}) = \mathrm{PFdim}(\mathtt{R})$.

As we will see in *e.g.* Example 8E.7(i), all of this should be thought of as generalizing very familiar notions from representation theory of groups.

*Proof.* (i). This follows since the PF eigenvalue is multiplicative.

(ii). As an invertible object can not have a nilpotent action matrix this follows from (i) and $\mathrm{PFdim}(\mathtt{X}) \geq 1$, see Theorem 8E.1.(i).

(iii). By Lemma 7D.10.(iii), the functor $_^{\star}$ preserves the property of being indecomposable and induces a bijection as in (7D-11). In other words, duality acts as a permutation on the set of indecomposable objects. Thus, $(\mathtt{ZX})^{\star} \cong (\mathtt{X}^{\star})(\mathtt{Z}^{\star})$ shows that, up to permutation, the action matrix for $\mathtt{X}^{\star}$ is the transpose of the action matrix for $\mathtt{X}$, which implies that $\mathrm{PFdim}(\mathtt{X}^{\star}) = \mathrm{PFdim}(\mathtt{X})$. The other claim follows by symmetry.

(iv). The second claim is clear by additivity, the first, that $\mathtt{T}$ is self-dual, follows from the bijections as in (7D-11).

(v). Since $\mathbf{C}$ is transitive, the total object $\mathtt{T}$ is a fusion generator, see (iv). By Lemma 8D.2.(iii), we can thus take $[\mathtt{R}] = v_{pf}(\Gamma_{\mathtt{T}})$, which is unique up to scaling and strictly positive. Hence, we can interpret $\mathtt{R}$ as a strictly positive sum of indecomposables of $\mathbf{C}$. By this construction and Lemma 8D.2.(iii) it follows that $[\mathtt{XR}]$ and $[\mathtt{RX}]$ must be proportional to $[\mathtt{R}]$. To see this observe that $[\mathtt{XR}]\Gamma_{\mathtt{T}} = \lambda_{pf}(\mathtt{T}) \cdot [\mathtt{XR}] = \Gamma_{\mathtt{T}}[\mathtt{RX}]$. (Note that $\Gamma_{\mathtt{T}}$ is symmetric by (iv).) This implies that $[\mathtt{XR}]$ and $[\mathtt{RX}]$ are strictly positive eigenvectors of $\Gamma_{\mathtt{T}}$, and the claim follows from Lemma 8D.2.(iii).

(vi). Clear by additivity of the PF eigenvalue. □

A crucial feature of PF dimensions is that they come in discrete values:

**Proposition 8E.6.** *Let* $\mathbf{C} \in \mathbf{wmFiat}$ *and* $\mathtt{X} \in \mathbf{C}$. *Then the PF dimensions* $\mathrm{PFdim}(\mathtt{X})$ *and* $\mathrm{PFdim}(\mathbf{C})$ *are algebraic integers, i.e. roots of some* $p \in \mathbb{Z}[X]$, *and* $\geq 1$.

*Proof.* To show that they are algebraic integers we can take $p \in \mathbb{Z}[X]$ to be the characteristic polynomial of $\Gamma_{\mathtt{X}}$, and the claim for $\mathbf{C}$ follows by additivity. For the second claim we observe that $\mathrm{PFdim}(\mathtt{X})^2 = \mathrm{PFdim}(\mathtt{XX}^{\star})$, by Lemma 8E.5.(i), and $\mathrm{PFdim}(\mathtt{XX}^{\star}) \geq 1$ by Theorem 8E.1.(i): $\Gamma_{\mathtt{XX}^{\star}}$ can not be a nilpotent matrix as there should always be a non-degenerate (co)pairing to $\mathbb{1}$. Hence, $\mathrm{PFdim}(\mathtt{X}) \geq 1$ which finishes the proof since, as before, the statement for **C** follows by additivity. □

**Example 8E.7.** Let $\mathrm{S}_3$ be the symmetric group in three letters, which is of order 6, and let $\mathbb{K}$ be algebraically closed. By Proposition 8B.6.(iii) we know that $\mathbf{fdMod}\big(\mathbb{K}[\mathrm{S}_3]\big)$ is fiat, and by Example 7E.15 it is semisimple if and only if the characteristic of $\mathbb{K}$ is not 2 or 3. So we basically have three cases:

(I) The case $\mathbb{K} = \mathbb{C}$, which we already glimpsed upon in Example 8C.6. In this case we get $\mathtt{L}_1 = \mathbb{1}$ and

$$\mathrm{PFdim}(\mathtt{L}_s) = \mathrm{PFdim}\Big(\begin{pmatrix} 0&1&0\\1&1&1\\0&1&0\end{pmatrix}\Big) = 2, \quad \mathrm{PFdim}(\mathtt{L}_{1'}) = \mathrm{PFdim}\Big(\begin{pmatrix} 0&0&1\\0&1&0\\1&0&0\end{pmatrix}\Big) = 1,$$

$$\mathtt{T} = \mathbb{1} \oplus \mathtt{L}_s \oplus \mathtt{L}_{1'}, \quad \Gamma_{\mathtt{T}} = \begin{pmatrix} 1&1&1\\1&3&1\\1&1&1\end{pmatrix} \rightsquigarrow$$

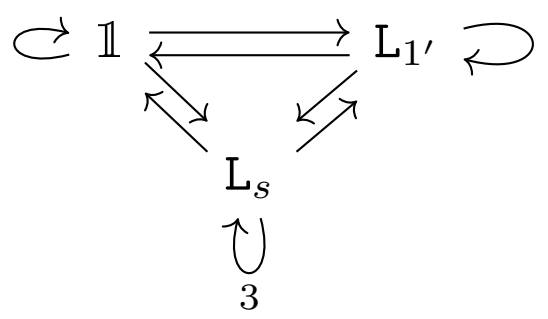

$$\mathtt{R} = \mathbb{1} \oplus 2\cdot\mathtt{L}_s \oplus \mathtt{L}_{1'} \cong \mathbb{C}[\mathrm{S}_3], \quad \mathrm{PFdim}\big(\mathbf{fdMod}\big(\mathbb{C}[\mathrm{S}_3]\big)\big) = 1\cdot 1 + 2\cdot 2 + 1\cdot 1 = 6 = \mathrm{PFdim}(\mathtt{R}).$$

Note that the regular object is the regular representation of $\mathbb{C}[\mathrm{S}_3]$ on itself, hence the name. The PF dimension in this case is the dimension of $\mathbb{C}[\mathrm{S}_3]$, a.k.a. the order of $\mathrm{S}_3$.

(II) The case $\mathbb{K} = \overline{\mathbb{F}}_3$. We want to use the construction in Remark 8B.7. First, we have a 3 Sylow subgroup $\mathbb{Z}/3\mathbb{Z}$ and we want to consider $\overline{\mathbb{F}}_3[\mathbb{Z}/3\mathbb{Z}]$. Similarly as for the case of $\mathbb{Z}/5\mathbb{Z}$, see *e.g.* Example 6K.18, we get that $\overline{\mathbb{F}}_3[\mathbb{Z}/3\mathbb{Z}]$ has three indecomposable modules, given by Jordan blocks for the eigenvalue 1: a $1\times 1$ Jordan block $\mathtt{Z}'_1$, a $2\times 2$ Jordan block $\mathtt{Z}'_2$, and a $3\times 3$ Jordan block $\mathtt{Z}'_3 = \overline{\mathbb{F}}_3[\mathbb{Z}/3\mathbb{Z}]$. Let $z$ be the basis element of $\mathtt{Z}'_1$, and let $1, s$ be the elements of $\mathrm{S}_2 \cong \mathrm{S}_3/(\mathbb{Z}/3\mathbb{Z})$. Then:

$$\overline{\mathbb{F}}_3[\mathrm{S}_3] \otimes_{\overline{\mathbb{F}}_3[\mathbb{Z}/3\mathbb{Z}]} \mathtt{Z}'_1 \cong \mathtt{Z}_1 \oplus \mathtt{Z}_{1'} \cong \mathtt{L}_1 \oplus \mathtt{L}_{1'},$$

for $\mathtt{L}_1 \cong \mathbb{1}$ and $\mathtt{L}_{1'}$ as in (I), which are also the only simples of $\overline{\mathbb{F}}_3[\mathrm{S}_3]$. To see this we simply observe that we can base change

$$\overline{\mathbb{F}}_3[\mathrm{S}_3] \otimes_{\overline{\mathbb{F}}_3[\mathbb{Z}/3\mathbb{Z}]} \mathtt{Z}'_1 = \overline{\mathbb{F}}_3\{1\otimes z, s\otimes z\} = \overline{\mathbb{F}}_3\{\tfrac{1}{2}(1+s)\otimes z, \tfrac{1}{2}(1-s)\otimes z\}.$$

(Also recall the idempotents $\mathrm{e}_{\pm}$ from Example 7C.4.) Moreover, we see analogously that

$$\overline{\mathbb{F}}_3[\mathrm{S}_3] \otimes_{\overline{\mathbb{F}}_3[\mathbb{Z}/3\mathbb{Z}]} \mathtt{Z}'_2 \cong \mathtt{Z}_2 \oplus \mathtt{Z}_{2'}, \quad \overline{\mathbb{F}}_3[\mathrm{S}_3] \otimes_{\overline{\mathbb{F}}_3[\mathbb{Z}/3\mathbb{Z}]} \mathtt{Z}'_3 \cong \mathtt{P}(\mathtt{L}_1) \oplus \mathtt{P}(\mathtt{L}_{1'}) \cong \mathtt{Z}_3 \oplus \mathtt{Z}_{3'},$$

which are all of the indicated dimensions. Summarized:

$$\mathtt{Si}\big(\mathbf{fdMod}\big(\overline{\mathbb{F}}_3[\mathrm{S}_3]\big)\big) = \{\mathtt{Z}_1 \cong \mathbb{1}, \mathtt{Z}_{1'} = \mathtt{L}_{1'}\}, \quad \mathtt{Pi}\big(\mathbf{fdMod}\big(\overline{\mathbb{F}}_3[\mathrm{S}_3]\big)\big) = \{\mathtt{Z}_3 \cong \mathtt{P}_{\mathbb{1}}, \mathtt{Z}_{3'} \cong \mathtt{P}_{1'}\},$$
$$\mathtt{In}\big(\mathbf{fdMod}\big(\overline{\mathbb{F}}_3[\mathrm{S}_3]\big)\big) = \{\mathtt{Z}_1, \mathtt{Z}_{1'}, \mathtt{Z}_2, \mathtt{Z}_{2'}, \mathtt{Z}_3, \mathtt{Z}_{3'}\}, \quad \dim(\mathtt{Z}_i) = i = \dim(\mathtt{Z}_{i'}).$$

To be more explicit, the Jordan–Hölder filtration are of the form

$$\mathtt{Z}_1 \cong \mathbb{1} \text{ is simple}, \quad \mathtt{Z}_{1'} = \mathtt{L}_{1'} \text{ is simple},$$
$$0 - \mathbb{1} - \mathtt{L}_{1'} - \mathtt{Z}_2, \quad 0 - \mathtt{L}_{1'} - \mathbb{1} - \mathtt{Z}_{2'},$$
$$0 - \mathbb{1} - \mathtt{L}_{1'} - \mathbb{1} - \mathtt{Z}_3, \quad 0 - \mathtt{L}_{1'} - \mathbb{1} - \mathtt{L}_{1'} - \mathtt{Z}_{3'}$$

These satisfy the fusion rules:

| $\otimes$ | $\mathtt{Z}_2$ | $\mathtt{P}_{\mathbb{1}}$ | $\mathtt{P}_{1'}$ | $\mathtt{Z}_{2'}$ | $\mathtt{L}_{1'}$ |
|---|---|---|---|---|---|
| $\mathtt{Z}_2$ | $\mathtt{P}_{\mathbb{1}} \oplus \mathtt{L}_{1'}$ | $\mathtt{P}_{\mathbb{1}} \oplus \mathtt{P}_{1'}$ | $\mathtt{P}_{\mathbb{1}} \oplus \mathtt{P}_{1'}$ | $\mathbb{1} \oplus \mathtt{P}_{1'}$ | $\mathtt{Z}_{2'}$ |
| $\mathtt{P}_{\mathbb{1}}$ | $\mathtt{P}_{\mathbb{1}} \oplus \mathtt{P}_{1'}$ | $\mathtt{P}_{\mathbb{1}} \oplus \mathtt{P}_{1'} \oplus \mathtt{P}_{\mathbb{1}}$ | $\mathtt{P}_{1'} \oplus \mathtt{P}_{\mathbb{1}} \oplus \mathtt{P}_{1'}$ | $\mathtt{P}_{\mathbb{1}} \oplus \mathtt{P}_{1'}$ | $\mathtt{P}_{1'}$ |
| $\mathtt{P}_{1'}$ | $\mathtt{P}_{\mathbb{1}} \oplus \mathtt{P}_{1'}$ | $\mathtt{P}_{1'} \oplus \mathtt{P}_{\mathbb{1}} \oplus \mathtt{P}_{1'}$ | $\mathtt{P}_{\mathbb{1}} \oplus \mathtt{P}_{1'} \oplus \mathtt{P}_{\mathbb{1}}$ | $\mathtt{P}_{\mathbb{1}} \oplus \mathtt{P}_{1'}$ | $\mathtt{P}_{\mathbb{1}}$ |
| $\mathtt{Z}_{2'}$ | $\mathbb{1} \oplus \mathtt{P}_{1'}$ | $\mathtt{P}_{\mathbb{1}} \oplus \mathtt{P}_{1'}$ | $\mathtt{P}_{\mathbb{1}} \oplus \mathtt{P}_{1'}$ | $\mathtt{P}_{\mathbb{1}} \oplus \mathtt{L}_{1'}$ | $\mathtt{Z}_2$ |
| $\mathtt{L}_{1'}$ | $\mathtt{Z}_{2'}$ | $\mathtt{P}_{1'}$ | $\mathtt{P}_{\mathbb{1}}$ | $\mathtt{Z}_2$ | $\mathbb{1}$ |

Thus, the action matrices and their PF eigenvalues and eigenvectors are:

$$\Gamma_1 = \begin{pmatrix} 1&0&0&0&0&0\\0&1&0&0&0&0\\0&0&1&0&0&0\\0&0&0&1&0&0\\0&0&0&0&1&0\\0&0&0&0&0&1\end{pmatrix}, \Gamma_{1'} = \begin{pmatrix} 0&0&0&0&0&1\\0&0&0&0&1&0\\0&0&0&1&0&0\\0&0&1&0&0&0\\0&1&0&0&0&0\\1&0&0&0&0&0\end{pmatrix}, \quad \begin{gathered}\lambda_{pf}(\Gamma_1) = 1,\\ v_{pf}(\Gamma_1) = (1,1,1,1,1,1),\\ \lambda_{pf}(\Gamma_{1'}) = 1,\\ v_{pf}(\Gamma_{1'}) = (1,1,1,1,1,1),\end{gathered}$$

$$\Gamma_2 = \begin{pmatrix} 0&0&0&0&1&0\\ 1&0&0&0&0&0\\ 0&1&1&1&0&0\\ 0&0&1&1&1&0\\ 0&0&0&0&0&1\\ 0&1&0&0&0&0 \end{pmatrix}, \Gamma_{2'} = \begin{pmatrix} 0&1&0&0&0&0\\ 0&0&0&0&0&1\\ 0&0&1&1&1&0\\ 0&1&1&1&0&0\\ 1&0&0&0&0&0\\ 0&0&0&0&1&0 \end{pmatrix}, \begin{gathered} \lambda_{pf}(\Gamma_2) = 2, \\ v_{pf}(\Gamma_2) = (0,0,1,1,0,0), \\ \lambda_{pf}(\Gamma_{2'}) = 2, \\ v_{pf}(\Gamma_{2'}) = (0,0,1,1,0,0), \end{gathered}$$

$$\Gamma_3 = \begin{pmatrix} 0&0&0&0&0&0\\ 0&0&0&0&0&0\\ 1&1&2&1&1&0\\ 0&1&1&2&1&1\\ 0&0&0&0&0&0\\ 0&0&0&0&0&0 \end{pmatrix}, \Gamma_{3'} = \begin{pmatrix} 0&0&0&0&0&0\\ 0&0&0&0&0&0\\ 0&1&1&2&1&1\\ 1&1&2&1&1&0\\ 0&0&0&0&0&0\\ 0&0&0&0&0&0 \end{pmatrix}, \begin{gathered} \lambda_{pf}(\Gamma_3) = 3, \\ v_{pf}(\Gamma_3) = (0,0,1,1,0,0), \\ \lambda_{pf}(\Gamma_{3'}) = 3, \\ v_{pf}(\Gamma_{3'}) = (0,0,1,1,0,0). \end{gathered}$$

Note that none of the indecomposables are fusion generators, and only the one dimensional ones are invertible. The total and regular objects are also not fusion generators:

$$\Gamma_{\mathtt{T}} = \begin{pmatrix} 1&1&0&0&1&1\\ 1&1&0&0&1&1\\ 1&3&6&6&3&1\\ 1&3&6&6&3&1\\ 1&1&0&0&1&1\\ 1&1&0&0&1&1 \end{pmatrix}, \Gamma_{\mathtt{R}} = \begin{pmatrix} 1&2&0&0&2&1\\ 2&1&0&0&1&2\\ 3&8&14&14&8&3\\ 3&8&14&14&8&3\\ 2&1&0&0&1&2\\ 1&2&0&0&2&1 \end{pmatrix}, \begin{gathered} \lambda_{pf}(\Gamma_3) = 12, \\ v_{pf}(\Gamma_3) = (0,0,1,1,0,0), \\ \lambda_{pf}(\Gamma_{3'}) = 28, \\ v_{pf}(\Gamma_{3'}) = (0,0,1,1,0,0). \end{gathered}$$

The PF dimension of the category itself is $\mathrm{PFdim}(\mathbf{C}) = \mathrm{PFdim}(\mathtt{R}) = 28$.

(III) The case $\mathbb{K} = \overline{\mathbb{F}}_2$ works *mutatis mutandis* as (II) above. Doing the calculations gives four indecomposable modules $\mathtt{Z}_1 = \mathtt{L}_1 \cong \mathbb{1}$, $\mathtt{Z}_2 = \mathtt{L}_s$, $\mathtt{Z}_3 \cong \mathtt{P}_{\mathbb{1}}$ and $\mathtt{Z}_{3'} \cong \mathtt{P}_s$, which are of the indicated dimensions. The fusion rules are:

| $\otimes$ | $\mathtt{L}_s$ | $\mathtt{P}_{\mathbb{1}}$ | $\mathtt{P}_s$ |
|---|---|---|---|
| $\mathtt{L}_s$ | $\mathtt{L}_{\mathbb{1}} \oplus \mathtt{P}_{\mathbb{1}}$ | $\mathtt{P}_{\mathbb{1}} \oplus \mathtt{P}_s$ | $\mathtt{P}_{\mathbb{1}} \oplus \mathtt{P}_{\mathbb{1}}$ |
| $\mathtt{P}_{\mathbb{1}}$ | $\mathtt{P}_{\mathbb{1}} \oplus \mathtt{P}_s$ | $2\cdot\mathtt{P}_{\mathbb{1}} \oplus \mathtt{P}_s$ | $\mathtt{P}_{\mathbb{1}} \oplus 2\cdot\mathtt{P}_s$ |
| $\mathtt{P}_s$ | $\mathtt{P}_{\mathbb{1}} \oplus \mathtt{P}_{\mathbb{1}}$ | $\mathtt{P}_{\mathbb{1}} \oplus 2\cdot\mathtt{P}_s$ | $2\cdot\mathtt{P}_{\mathbb{1}} \oplus \mathtt{P}_s$ |

To compute the PF dimensions *etc.* is Exercise 8L.4.

It is a fun exercise to use MAGMA to verify all the above, *cf.* Appendix A. ◇

**Example 8E.8.** Let us again consider $\mathrm{S}_3$, but now rather the category $\mathbf{Vec}^{\omega}_{\mathbb{k}\oplus}(\mathrm{S}_3)$, for any 3 cocycle $\omega$. Let $\mathrm{S}_3 = \{1, s, t, ts, st, sts = tst\}$, where, in graphical notation,

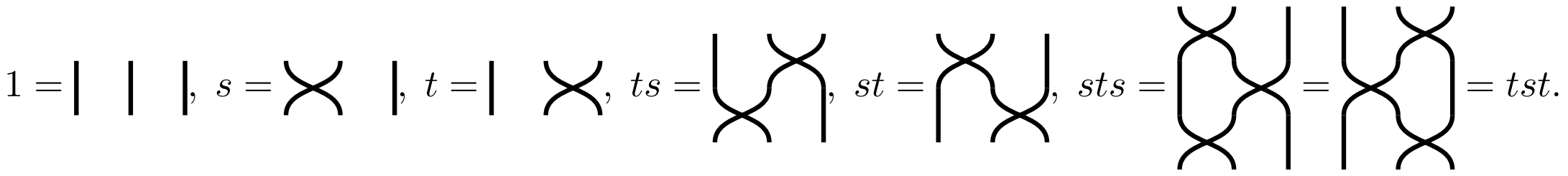

By definition, the fusion rules of $\mathbf{Vec}^{\omega}_{\mathbb{k}\oplus}(\mathrm{S}_3)$ are exactly the multiplication rules of $\mathrm{S}_3$. Thus, the action matrices are just permutation matrices, *e.g.*

$$\Gamma_{st} = \begin{pmatrix} 0&0&0&1&0&0\\ 0&0&0&0&0&1\\ 0&1&0&0&0&0\\ 0&0&0&0&1&0\\ 1&0&0&0&0&0\\ 0&0&1&0&0&0 \end{pmatrix}, \begin{gathered} \lambda_{pf}(\Gamma_{sts}) = 1, \\ v_{pf}(\Gamma_{sts}) = (1,1,1,1,1,1). \end{gathered}$$

Thus, all PF dimensions are 1, and all objects are invertible. Note also that $\mathrm{PFdim}\big(\mathbf{Vec}^{\omega}_{\mathbb{C}\oplus}(\mathrm{S}_3)\big) = 6 = \mathrm{PFdim}\big(\mathbf{fdMod}\big(\mathbb{C}[\mathrm{S}_3]\big)\big)$, the order of $\mathrm{S}_3$. ◇

8F. **Webs and the Fibonacci category.** Let us come back to Example 8E.4. Recall that $\phi = \frac{1}{2}(1+\sqrt{5})$ denotes the golden ratio. Let $\Phi = -\phi^{-1} = -\phi + 1 = \frac{1}{2}(1-\sqrt{5})$, which is the Galois conjugate of $\phi$. Note that these make sense in any algebraically closed field of characteristic not two or five, and we will use the corresponding elements below, using the same notation. In characteristic two the polynomial $x^2 - x - 1$ (the minimal polynomial of $\phi$) is irreducible, and we can let $\phi$ and $\Phi$ denote its roots in $\mathbb{F}_4$.

*Remark* 8F.1. In characteristic five we have $x^2 - x - 1 = (x-3)^2$, and this is the only characteristic where $x^2 - x - 1$ is reducible and has a double root. ◇

**Definition 8F.2.** For an algebraically closed field $\mathbb{K}$ of characteristic not five, let $\mathbf{Fib} = \mathbf{Fib}_{\mathbb{K}\oplus}$ denote the quotient of $\mathbf{Web}_{\mathbb{K}\oplus}$ from Example 3G.2 by the ***circle*** and ***bitri evaluation***, and the ***H=I relation***:

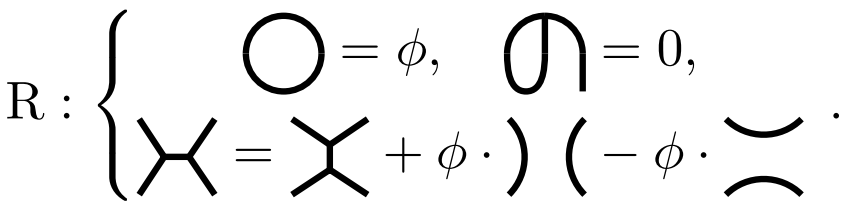

We call **Fib** the ***Fibonacci category***.

For $\mathbb{K}$ of characteristic five we let **Fib** be the full subcategory of $\mathbf{fdMod}\big(\mathbb{K}[\mathbb{Z}/5\mathbb{Z}]\big)/(proj)$ (modulo the ∘-⊗ ideal of projectives) generated by $\mathtt{Z}_3$, *cf.* Example 6K.18. ◇

**Proposition 8F.3.** *The category* **Fib** *is a semisimple fiat category with two simple objects* $\mathtt{L}_1 = \mathbb{1}$ *and* $\mathtt{L} = \mathtt{L}_2$ *with*

$$\mathtt{L}^2 \cong \mathbb{1} \oplus \mathtt{L}.$$

*Proof.* We split the proof into several smaller statements. We rename $\mathtt{L}$ to $\bullet$. Until the very end, we will assume that we are not in characteristic five.

**Lemma 8F.4.** *In* **Fib** *we have*

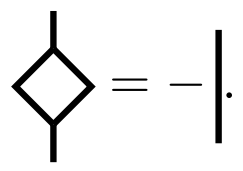

*Proof.* Observe that

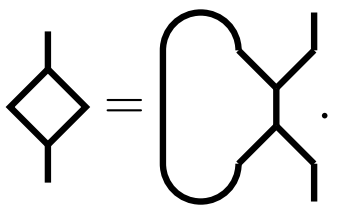

Now apply the H=I relation, followed by the other relations. □

**Lemma 8F.5.** *In* **Fib** *we have*

(8F-6)

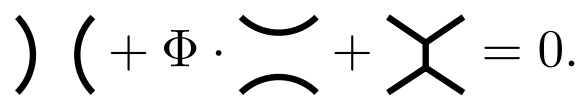

*Proof.* There are four diagrams in $\mathrm{End}_{\mathbf{Fib}}(\bullet^2)$ that one would naturally draw, namely the four diagrams in the H=I relation. The H=I relation itself gives a linear relation between the four diagrams, so three diagrams remain. We choose:

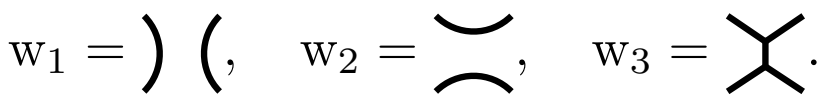

By the defining relations and Lemma 8F.4, their multiplication table is

| | $\mathrm{w}_1$ | $\mathrm{w}_2$ | $\mathrm{w}_3$ |
|---|---|---|---|
| $\mathrm{w}_1$ | $\mathrm{w}_1$ | $\mathrm{w}_2$ | $\mathrm{w}_3$ |
| $\mathrm{w}_2$ | $\mathrm{w}_2$ | $\phi\cdot\mathrm{w}_2$ | $0$ |
| $\mathrm{w}_3$ | $\mathrm{w}_3$ | $0$ | $-\mathrm{w}_3$ |

The traces of the three diagrams are

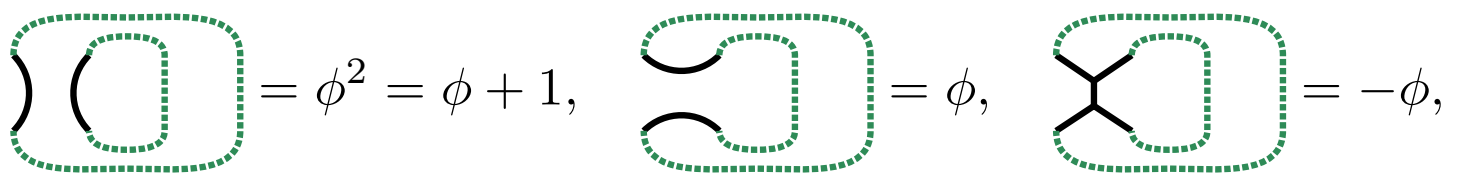

as one can now easily verify. Hence, the trace of the multiplication table (taking traces entry-wise) is

$$P = \begin{pmatrix} \phi+1 & \phi & -\phi \\ \phi & \phi+1 & 0 \\ -\phi & 0 & \phi \end{pmatrix},$$

which has its kernel spanned by $(1, \Phi, 1)$. Thus, we get Equation 8F-6, as desired. □

**Lemma 8F.7.** *In* **Fib***, for each hom space, there is a* $\mathbb{K}$ *spanning set of diagrams without internal edges (edges that do not touch the boundary).*

*Proof.* The H=I relation and Equation 8F-6 imply that we can remove any internal edge. □

**Lemma 8F.8.** *The only idempotents in* $\mathrm{End}_{\mathbf{Fib}}(\bullet^2)$ *are*

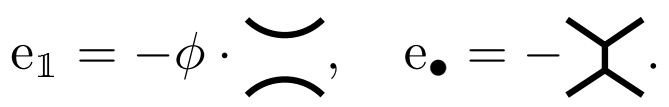

*These form a complete set of orthogonal primitive idempotents.*

*Proof.* The defining relations and Lemma 8F.4 show that these are orthogonal idempotents, while Equation 8F-6 ensures that they add to the identity. Finally, Lemma 8F.7 confirms that there are no others. □

**Lemma 8F.9.** *In* **Fib** *we have*

$$\bullet^2 \cong \mathbb{1} \oplus \bullet.$$

*Hence, there are only two simple objects in* **Fib***, namely* $\mathbb{1}$ *and* $\bullet$*, and the category is semisimple.*

*Proof.* Cutting the idempotents from Lemma 8F.8 in halves gives the required isomorphisms:

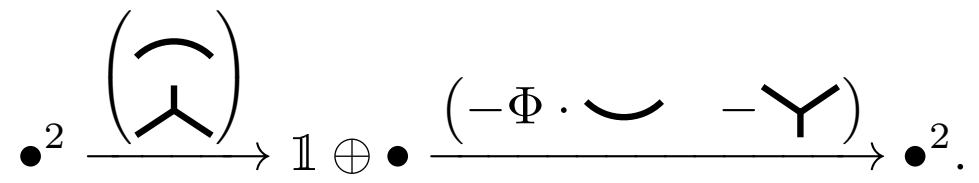

The above verifies that the two columns give isomorphism.

For the final statement, observe that Lemma 8F.7 computes $\mathrm{End}_{\mathbf{Fib}}(\mathbb{1})$ and $\mathrm{End}_{\mathbf{Fib}}(\bullet)$, showing that both objects are simple. The rule $\bullet^2 \cong \mathbb{1} \oplus \bullet$ can then be used to show that **Fib** is semisimple. □

All necessary claims are now easy to verify. Finally, in characteristic five, this follows directly from Example 6K.18 and Equation 7H-7, after field extension from $\overline{\mathbb{F}}_5$ to $\mathbb{K}$. □

8G. **Fusion categories.** Here are the categorical analogs of groups:

**Definition 8G.1.** A semisimple category $\mathbf{C} \in$ **wmFiat** is called a ***weakly multi fusion category***, and a semisimple category $\mathbf{C} \in$ **wFiat** is called a ***weakly fusion category***. If these are also pivotal, we call them ***multi fusion categories*** and ***fusion categories***, respectively. ◇

Without further ado, we get full subcategories *e.g.* of the form **wmFus** ⊂ **wmFiat** (potentially omitting the w and the m), called *e.g.* ***category of weakly multi fusion categories***, equivalence being $\simeq_{\mathbb{S}\oplus\star}$.

**Example 8G.2.** We have already seen the prototypical examples of (weakly) fusion categories, namely $\mathbf{Vec}^{\omega}_{\Bbbk\oplus}(\mathrm{G})$ and $\mathbf{fdMod}\big(\Bbbk[\mathrm{G}]\big)$, where G is a finite group, in both cases, and #G does not divide the characteristic of $\Bbbk$ in the second case, as well as $\mathbf{fdMod}(\mathrm{A})$ where A is a semisimple Hopf algebra. ◇

**Proposition 8G.3.** *Fusion categories can be alternatively defined e.g. as follows:*
*"A semisimple category* $\mathbf{C} \in$ **wmTen** *is called a* **weakly multifusion category***".*

*Proof.* Clear, since semisimple is a stronger notion than abelian, see Theorem 7E.9. □

**Proposition 8G.4.** *Let* $\mathbf{C} \in$ **wmFus***. Then:*

*(i) If* $\mathbf{C} \in$ **wFus** *is* $\Bbbk$ *linear, then* $\mathbb{1} \in \mathrm{Si}(\mathbf{C})$*.*

*(ii) If* $\mathbf{C} \in$ **wFus***, then* $\mathtt{X}^\star \cong {}^\star\mathtt{X}$ *for all* $\mathbf{X} \in \mathbf{C}$*.*

*(iii) The fusion coefficients* $N_{ij}^k$ *are cyclic up to duality, i.e.*

$$\mathtt{L}_i\mathtt{L}_j \cong \textstyle\bigoplus_{k=1}^n N_{i,j}^k \cdot \mathtt{L}_k \Leftrightarrow (\mathtt{L}_k^\star)\mathtt{L}_i \cong \bigoplus_{k=1}^n N_{(k^\star),i}^{j^\star} \cdot \mathtt{L}_j^\star.$$

*(iv) A regular object is*

$$\mathtt{R} = \textstyle\sum_{i=1}^n \mathrm{PFdim}(\mathtt{L}_i) \cdot \mathtt{L}_i.$$

*Proof.* (i). By $\mathrm{End}_{\mathbf{C}}(\mathbb{1}) \cong \Bbbk$, the monoidal unit is indecomposable, hence simple.

(ii). By additivity, it suffices to show $\mathtt{L}^\star \cong {}^\star\mathtt{L}$ for all simples. Here we first note that

$$\big(\mathrm{Hom}_{\mathbf{C}}(\mathbb{1}, \mathtt{L}'\mathtt{L}) \ncong 0\big) \Rightarrow \big(\mathtt{L}' \cong \mathtt{L}^\star\big), \quad \big(\mathrm{Hom}_{\mathbf{C}}(\mathtt{L}'\mathtt{L}, \mathbb{1}) \ncong 0\big) \Rightarrow \big(\mathtt{L}' \cong {}^\star\mathtt{L}\big),$$

by semisimplicity, (i) and Schur's lemma Lemma 6J.9. Moreover, also by semisimplicity,

$$\mathrm{Hom}_{\mathbf{C}}(\mathbb{1}, \mathtt{L}'\mathtt{L}) \cong \mathrm{Hom}_{\mathbf{C}}(\mathtt{L}'\mathtt{L}, \mathbb{1}),$$

which then in turn implies the claim.

(iii). By noting that

$$N_{ij}^k = \dim\big(\mathrm{Hom}_{\mathbf{C}}(\mathtt{L}_i\mathtt{L}_j, \mathtt{L}_k)\big) = \dim\big(\mathrm{Hom}_{\mathbf{C}}((\mathtt{L}_k^\star)\mathtt{L}_i, \mathtt{L}_j^\star)\big) = N_{(k^\star)i}^{j^\star},$$

which holds by Theorem 4C.8 and (ii).

(iv). Using the previous results, we calculate

$$\begin{aligned}\mathtt{L}_i\mathtt{R} &\cong \textstyle\bigoplus_{j=1}^n \mathrm{PFdim}(\mathtt{L}_j) \cdot \mathtt{L}_i\mathtt{L}_j \cong \bigoplus_{j,k=1}^n \mathrm{PFdim}(\mathtt{L}_j) \cdot N_{i,j}^k\mathtt{L}_k \cong \bigoplus_{j,k=1}^n \mathrm{PFdim}(\mathtt{L}_j) \cdot N_{(k^\star)i}^{j^\star}\mathtt{L}_k \\ &\cong \textstyle\bigoplus_{j,k=1}^n \mathrm{PFdim}(\mathtt{L}_j) \cdot N_{(i^\star)k}^{j}\mathtt{L}_k \cong \bigoplus_{k=1}^n \mathrm{PFdim}\big((\mathtt{L}_i^\star)\mathtt{L}_k\big) \cdot \mathtt{L}_k \\ &\cong \mathrm{PFdim}(\mathtt{L}_i^\star) \cdot \big(\textstyle\bigoplus_{k=1}^n \mathrm{PFdim}(\mathtt{L}_k) \cdot \mathtt{L}_k\big) \cong \mathrm{PFdim}(\mathtt{L}_i) \cdot \mathtt{R},\end{aligned}$$

which implies the claim by uniqueness of the regular object. □

**Example 8G.5.** Let G be a finite group, and let $\Bbbk = \mathbb{C}$. Note that for any subgroup $\mathrm{H} \subset \mathrm{G}$ we get an algebra

$$\mathtt{A}_{\mathrm{H}} = \bigoplus_{h \in \mathrm{H}} \mathtt{h} \in \mathbf{Vec}_{\mathbb{C}\oplus}(\mathrm{G}).$$

These are of course the corresponding group algebras and thus

$$\mathbf{Mod}(\mathtt{A}_{\mathrm{H}}) \simeq_{\mathbb{C}\oplus\star} \mathbf{fdMod}\big(\mathbb{C}[\mathrm{H}]\big).$$

The generalization of this construction, as explained in details in [EGNO15, Example 9.7.2], takes an algebra $\mathtt{A}_{\mathrm{H}} \in \mathbf{Vec}^{\omega}_{\mathbb{C}\oplus}(\mathrm{G})$ together with a so-called 2 cochain $\psi$ to twist the multiplication of $\mathtt{A}_{\mathrm{H}}$. The corresponding module categories

$$\mathbf{fdMod}^{\omega,\psi}_{\mathrm{G}}\big(\mathbb{C}[\mathrm{H}]\big) = \mathbf{Mod}(\mathtt{A}_{\mathrm{H}}),$$

$$\text{where } \mathtt{A}_{\mathrm{H}} \in \mathbf{Vec}^{\omega}_{\mathbb{C}\oplus}(\mathrm{G}) \text{ and } \psi \text{ is a 2 cochain with } d_2\psi = \omega|_{\mathrm{H}\times\mathrm{H}\times\mathrm{H}},$$

are all fusion categories, and are sometimes called ***group-like fusion categories***. ◇

*Remark* 8G.6. The construction in Example 8G.5 can also be done for arbitrary $\mathbb{S}$ instead of $\mathbb{C}$, by letting the cochains take values in $\mathbb{S}^*$. ◇

8H. **Verlinde categories – part I.** Let $\mathbb{C}$ be our ground field (for categories and representations). For a (complex) Lie algebra $\mathfrak{g} = \mathfrak{g}(\mathbb{C})$ let $\mathbf{fdMod}(\mathfrak{g})$ denote its category of finite dimensional representations.

**Example 8H.1.** The only example the reader needs to keep in mind is

$$\mathfrak{sl}_2 = \big\{\big(\begin{smallmatrix} a & b \\ c & d \end{smallmatrix}\big) \mid a,b,c,d \in \mathbb{C}, a+d=0\big\}$$

of traceless 2-by-2 matrices. The associated Lie group is

$$\mathrm{SL}_2 = \mathrm{SL}_2(\mathbb{C}) = \big\{\big(\begin{smallmatrix} a & b \\ c & d \end{smallmatrix}\big) \mid a,b,c,d \in \mathbb{C}, ad-bc=1\big\},$$

the ***special linear group*** of 2-by-2 matrices. This is a group so it category of finite dimensional complex representations $\mathbf{fdMod}(\mathrm{SL}_2)$ is symmetric ribbon $\mathbb{C}$ linear, and also semisimple. As a matter of fact,

$$\mathbf{fdMod}(\mathfrak{sl}_2) \simeq_{\mathbb{C}\oplus\star} \mathbf{fdMod}(\mathrm{SL}_2),$$

and the reader unfamiliar with Lie algebras can alternatively think about working with the associated Lie group G. ◇

*Remark* 8H.2. The analogy in Example 8H.1, replacing $\mathbf{fdMod}(\mathfrak{g})$ by $\mathbf{fdMod}(\mathrm{G})$ works well, however with a mild caveat. In general, there is a fully faithful ribbon $\mathbb{C}$ linear functor

$$\mathbf{fdMod}(\mathrm{G}) \to \mathbf{fdMod}(\mathfrak{g}),$$

but the functor needs not to be essentially surjective. ◇

There is a family of fusion categories, which are of paramount importance for the construction of the classical quantum invariants, and also for the theory of fiat and fusion categories. However, they are not easy to construct, and we will postpone a more detailed discussion to the later sections. For now, we just state the theorem:

**Theorem 8H.3.** *We have the following.*

- **(a)** *For any finite dimensional semisimple complex Lie algebra $\mathfrak{g}$, any $k \in \mathbb{Z}_{\geq h}$ (where $h$ is the Coxeter number of $\mathfrak{g}$) and any $q \in \mathbb{C}$ being a primitive $2k$th root of unity there exists a $\mathbb{C}$ linear category $\mathbf{fdMod}^q_k(\mathfrak{g}) \in \mathbf{Fus}$.*
- **(b)** *All of the categories $\mathbf{fdMod}^q_k(\mathfrak{g})$ are ribbon.*
- **(c)** *Let $q = \exp(\pi i/k)$. The limit $k \to \infty$ of $K_0^{\oplus}\big(\mathbf{fdMod}^q_k(\mathfrak{g})\big)$ is $K_0^{\oplus}\big(\mathbf{fdMod}(\mathfrak{g})\big)$.*

The categories of the form $\mathbf{fdMod}^q_k(\mathfrak{g}) \in \mathbf{Fus}$ are called ***Verlinde categories***. By Theorem 8H.3.(c), we can think of them as finite approximations of $\mathbf{fdMod}(\mathfrak{g}$.

*Proof.* We will elaborate later, but the main construction can be found in *e.g.* [And92]. See also [BK01b] for many details on Verlinde categories (denoted $\mathcal{C}(\mathfrak{g},\chi)$ therein.) □

**Example 8H.4.** To be at least a bit more explicit, let $\mathfrak{g} = \mathfrak{sl}_2$ where $h = 2$. Let us explain the fusion rules of the categories $\mathbf{fdMod}^q_k(\mathfrak{sl}_2)$ which only depend on $k$ and not on $q$. So let us fix $k \geq 2$ and choose $q = \exp(\pi i/k)$, for the sake of concreteness. The category $\mathbf{fdMod}^q_k(\mathfrak{sl}_2)$ is ribbon and fusion and has simple objects

$$\mathrm{Si}\big(\mathbf{fdMod}^q_k(\mathfrak{sl}_2)\big) = \{\mathtt{L}_i \mid i = 0,\dots,k-2\}, \quad \mathtt{L}_0 \cong \mathbb{1}, \quad \mathtt{L}_i^\star \cong \mathtt{L}_i.$$

The fusion rules are given by the ***truncated Clebsch–Gordan rule***

$$\mathtt{L}_i\mathtt{L}_j \cong \bigoplus_{l=\max(i+j-k+2,0)}^{\min(i,j)} \mathtt{L}_{i+j-2l}. \tag{8H-5}$$

Let us discuss a few cases for small $k$:

- For $k=2$ we have $\mathbf{fdMod}_2^q(\mathfrak{sl}_2)\simeq_{\mathbb{C}\oplus\star}\mathbf{Vec}_{\mathbb{C}}$. Hence, the fusion generator $\mathbb{1}$ has $\mathrm{PFdim}(\mathbb{1})=1$ giving $\mathrm{PFdim}\big(\mathbf{fdMod}_2^q(\mathfrak{sl}_2)\big)=1$.
- For $k=3$ we have that the fusion generator $\mathtt{L}_1$ satisfies
$$\mathtt{L}_1\mathtt{L}_1\cong\mathbb{1},\quad \Gamma_1=\begin{pmatrix}0&1\\1&0\end{pmatrix}=\mathbb{1}\rightleftarrows\mathtt{L}_1,\quad \mathrm{PFdim}(\Gamma_1)=2\cos(\pi/3)=1.$$
And thus, $\mathrm{PFdim}\big(\mathbf{fdMod}_2^q(\mathfrak{sl}_3)\big)=2$. One can actually show that $\mathbf{fdMod}_k^q(\mathfrak{sl}_2)\simeq_{\mathbb{C}\oplus\star}\mathbf{Vec}_{\mathbb{C}\oplus}(\mathbb{Z}/2\mathbb{Z})$.
- For $k=4$ the fusion rules take the form

| $\otimes$ | $\mathtt{L}_1$ | $\mathtt{L}_2$ |
|---|---|---|
| $\mathtt{L}_1$ | $\mathbb{1}\oplus\mathtt{L}_2$ | $\mathtt{L}_1$ |
| $\mathtt{L}_2$ | $\mathtt{L}_1$ | $\mathbb{1}$ |

.

Thus, we get the fusion graphs
$$\Gamma_1=\begin{pmatrix}0&1&0\\1&0&1\\0&1&0\end{pmatrix}=\mathbb{1}\rightleftarrows\mathtt{L}_1\rightleftarrows\mathtt{L}_2,\quad \mathrm{PFdim}(\Gamma_1)=2\cos(\pi/4)=\sqrt{2},$$
$$\Gamma_2=\begin{pmatrix}0&0&1\\0&1&0\\1&0&0\end{pmatrix}=\mathbb{1}\rightleftarrows\mathtt{L}_2\ \ (\mathtt{L}_1\circlearrowleft),\quad \mathrm{PFdim}(\Gamma_2)=1.$$
Hence, $\mathrm{PFdim}\big(\mathbf{fdMod}_4^q(\mathfrak{sl}_2)\big)=1+2\cos(\pi/4)^2+1=4$.

For general $k\geq 2$ the object $\mathtt{L}_1$ will be a fusion generator of $\mathbf{fdMod}_k^q(\mathfrak{sl}_2)$ linear fusion graph
$$\Gamma_1=\mathbb{1}\rightleftarrows\mathtt{L}_1\rightleftarrows\mathtt{L}_2\rightleftarrows\dots\rightleftarrows\mathtt{L}_k,$$
and $\mathrm{PFdim}(\mathtt{L}_1)=2\cos(\pi/k)$. ◇

**Example 8H.6.** In terms of the group $\mathrm{SL}_2$, *cf.* Example 8H.1, the simple object $\mathtt{L}_i$ in Example 8H.4 corresponds to the $i+1$ dimensional $\mathrm{SL}_2$ representation $\mathrm{Sym}^i\mathbb{C}^2$. If one thinks of the latter as homogeneous degree $i$ polynomials in $X$ and $Y$ with basis $X^iY^0,X^{i-1}Y^1,\dots,X^0Y^i$, then the matrix $\left(\begin{smallmatrix}a&b\\c&d\end{smallmatrix}\right)$ acts as the matrix with the $k$th row given by expanding $(aX+cY)^{i-k+1}(bX+dY)^{k-1}$. Explicitly, if $i=2$, then
$$\left(\begin{smallmatrix}a&b\\c&d\end{smallmatrix}\right)\mapsto\left(\begin{smallmatrix}a^2&ac&c^2\\ab&ad+bc&cd\\b^2&bd&d^2\end{smallmatrix}\right)$$
is the action of $\mathrm{SL}_2$ on $\mathrm{Sym}^2\mathbb{C}^2\cong\mathbb{C}X^2Y^0\oplus\mathbb{C}X^1Y^1\oplus\mathbb{C}X^0Y^2$. ◇

**Example 8H.7.** Using the same notation as in Example 8H.4, it is easy to see from the fusion rules that there is a full subcategory $\mathbf{fdMod}_k^q(\mathfrak{so}_3)\subset\mathbf{fdMod}_k^q(\mathfrak{sl}_2)$ with
$$\big(\mathbf{fdMod}_k^q(\mathfrak{so}_3)\big)=\{\mathtt{L}_i\mid i=0,\dots,k-2,i\text{ even}\},\quad \mathtt{L}_0\cong\mathbb{1},\quad \mathtt{L}_i^\star\cong\mathtt{L}_i,$$
and exactly the same fusion rules, *i.e.* (8H-5), just taking even indexes only. Again, let us discuss some cases in detail:

- The cases $k=2$ and $k=3$ will now collide (on the Grothendieck level) and are as in Example 8H.4.
- The case $k=4$ gives $\mathbf{fdMod}_k^q(\mathfrak{so}_3)\simeq_{\mathbb{C}\oplus\star}\mathbf{Vec}_{\mathbb{C}\oplus}(\mathbb{Z}/2\mathbb{Z})$.
- For $k=5$ one can show that $\mathbf{fdMod}_k^q(\mathfrak{so}_3)\simeq_{\mathbb{C}\oplus\star}\mathbf{Fib}$ (the Fibonacci category, see Example 8E.4) for $q=\exp(\pi i/5)$.
- For $k=6$ the fusion rules take the form

| $\otimes$ | $\mathtt{L}_2$ | $\mathtt{L}_4$ |
|---|---|---|
| $\mathtt{L}_2$ | $\mathbb{1}\oplus\mathtt{L}_2\oplus\mathtt{L}_4$ | $\mathtt{L}_2$ |
| $\mathtt{L}_4$ | $\mathtt{L}_2$ | $\mathbb{1}$ |

,

and we get the fusion graphs
$$\Gamma_2=\begin{pmatrix}0&1&0\\1&1&1\\0&1&0\end{pmatrix}=\mathbb{1}\rightleftarrows\mathtt{L}_2\rightleftarrows\mathtt{L}_4\ \ (\mathtt{L}_2\circlearrowleft),\quad \mathrm{PFdim}(\Gamma_2)=2,$$
$$\Gamma_4=\begin{pmatrix}0&0&1\\0&1&0\\1&0&0\end{pmatrix}=\mathbb{1}\rightleftarrows\mathtt{L}_4\ \ (\mathtt{L}_2\circlearrowleft),\quad \mathrm{PFdim}(\Gamma_1)=1.$$
Hence, $\mathrm{PFdim}\big(\mathbf{fdMod}_6^q(\mathfrak{so}_3)\big)=6$. In fact, $\mathbf{fdMod}_6^q(\mathfrak{so}_3)\simeq_{\mathbb{C}\oplus\star}\mathbf{fdMod}\big(\mathbb{C}[\mathrm{S}_3]\big)$.

- Finally, for $k=7$ we get

| $\otimes$ | $\mathrm{L}_2$ | $\mathrm{L}_4$ |
|---|---|---|
| $\mathrm{L}_2$ | $\mathbb{1}\oplus\mathrm{L}_2\oplus\mathrm{L}_4$ | $\mathrm{L}_2\oplus\mathrm{L}_4$ |
| $\mathrm{L}_4$ | $\mathrm{L}_2\oplus\mathrm{L}_4$ | $\mathbb{1}\oplus\mathrm{L}_2$ |

Thus, the fusion graphs are

$$\Gamma_2=\begin{pmatrix}0&1&0\\1&1&1\\0&1&1\end{pmatrix}=\mathbb{1}\rightleftarrows \underset{\circlearrowright}{\mathrm{L}_2}\rightleftarrows \underset{\circlearrowright}{\mathrm{L}_4},\quad \mathrm{PFdim}(\Gamma_2)=2,$$

$$\Gamma_4=\begin{pmatrix}0&0&1\\0&1&1\\1&1&0\end{pmatrix}=\mathbb{1}\quad \underset{\circlearrowright}{\mathrm{L}_2}\rightleftarrows \mathrm{L}_4\ (\mathbb{1}\rightleftarrows\mathrm{L}_4),\quad \mathrm{PFdim}(\Gamma_1)=1.$$

For general $k=2l\geq 4$ or $k=2l+1\geq 5$ the object $\mathrm{L}_2$ will be a fusion generator of $\mathbf{fdMod}_k^q(\mathfrak{so}_3)$ with fusion graph

$$\Gamma_2=\ \mathbb{1}\rightleftarrows \underset{\circlearrowright}{\mathrm{L}_2}\rightleftarrows \underset{\circlearrowright}{\mathrm{L}_4}\rightleftarrows \underset{\circlearrowright}{\dots}\rightleftarrows \underset{\circlearrowright}{\mathrm{L}_{2l-2}}\rightleftarrows \mathrm{L}_l\qquad \text{if } k=2l\geq 4,$$

$$\Gamma_2=\ \mathbb{1}\rightleftarrows \underset{\circlearrowright}{\mathrm{L}_2}\rightleftarrows \underset{\circlearrowright}{\mathrm{L}_4}\rightleftarrows \underset{\circlearrowright}{\dots}\rightleftarrows \underset{\circlearrowright}{\mathrm{L}_{2l-2}}\rightleftarrows \underset{\circlearrowright}{\mathrm{L}_l}\qquad \text{if } k=2l+1\geq 5.$$

We will come back to this example several times. ◇

We will construct $\mathbf{fdMod}_k^q(\mathfrak{sl}_2)$ and $\mathbf{fdMod}_k^q(\mathfrak{so}_3)$ later one. For now, let us just continue, pretending we have already constructed them. There is almost no redundancy:

**Proposition 8H.8.** *We have*

$$\big(\mathbf{fdMod}_k^q(\mathfrak{sl}_2)\simeq_{\mathbb{C}\oplus\star}\mathbf{fdMod}_{k'}^{q'}(\mathfrak{sl}_2)\big)\Leftrightarrow\big(k=k' \textit{ and } (q=q' \textit{ or } q^{-1}=q')\big),$$
$$\big(\mathbf{fdMod}_k^q(\mathfrak{so}_3)\simeq_{\mathbb{C}\oplus\star}\mathbf{fdMod}_{k'}^{q'}(\mathfrak{so}_3)\big)\Leftrightarrow\big(k=k' \textit{ and } (q=q' \textit{ or } q^{-1}=q')\big).$$

*Proof.* See [**FK93**, Proposition 8.2.3]. □

Proposition 8H.8 identifies the various categories in the following fashion, here $k=5$:

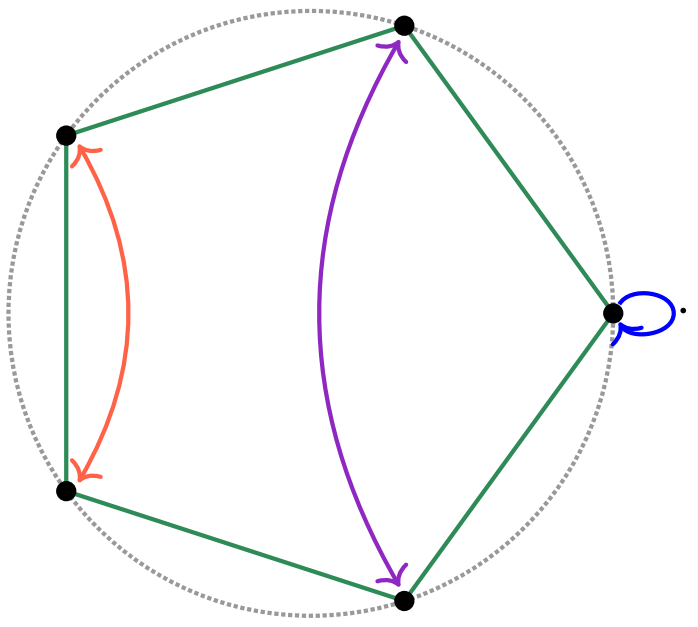

That is, one first chooses $k$, the polygon, and then the isomorphism classes of Verlinde categories are given by the orbits of complex conjugation of the vertices.

**Example 8H.9.** Note that the $k=5$ case in Example 8H.7 has two non-equivalent cases: $\mathbf{fdMod}_k^q(\mathfrak{so}_3)\simeq_{\mathbb{C}\oplus\star}$ **Fib**, which happens for $q=\exp(\pi i/5)$, and another fusion category appearing for $q=\exp(2\pi i/5)$. See also Exercise 8L.2, where we point out that $q^2+1+q^{-2}$ is either $\phi$ or $\Phi$, depending on the choice of $q$. ◇

8I. **Classifying fiat, tensor and fusion categories.** Let us now address some classification problems. Namely, we want to ask (in order) whether one can classify fiat, tensor or fusion categories with a given $K_0^{\oplus}(_)$, with a fixed rank $\mathrm{rk}(_)$ or a fixed PF dimension $\mathrm{PFdim}(_)$.

*Remark* 8I.1. We will be a bit sketchy in this section because we want to state theorems which are easy to understand (and worthwhile to be stated) but sometimes not easy to prove. ◇

The arguable most important theorem in the theory is ***Ocneanu rigidity***, which is a "uniqueness of categorification" type of statement:

**Theorem 8I.2.** *The number of* $\mathbb{k}$ *linear weakly multi fusion categories (up to* $\simeq_{\mathbb{k}\oplus\star}$ *equivalence) with a given* $K_0^{\oplus}(_)$ *is finite.*

*Proof.* The (not easy) proof of this theorem can be found in *e.g.* [**EGNO15**, Theorem 9.1.4]. □

We have already seen two numerical invariants, which only depend only on $K_0^{\oplus}(_)$, of fiat categories: the rank, *cf.* Proposition 7H.4 and the PF dimension, *cf.* Definition 8E.3, and both discrete valued and $\geq 1$. Thus, Theorem 8I.2 motivates the question whether one can classify fiat or fusion categories of a given $K_0^{\oplus}(_)$, of a given rank or of a given PF dimension.

Let us start by fixing $K_0^{\oplus}(_)$.

**Proposition 8I.3.** *Let* $\mathrm{G}$ *be a finite group. If* $\mathbf{C} \in \mathbf{wmFus}$ *is* $\mathbb{C}$ *linear and has* $K_0^{\oplus}(\mathbf{C}) \cong \mathbb{Z}[\mathrm{G}]$ *as* $\mathbb{Z}$ *algebras, then* $\mathbf{C} \simeq_{\mathbb{C}\oplus\star} \mathbf{Vec}^{\omega}_{\mathbb{C}\oplus}(\mathrm{G})$.

*Proof.* By carefully writing down all equations coming from the associativity and unitality constrains, see *e.g.* [**EGNO15**, Proposition 4.10.3] for details. □

For a finite group $\mathrm{G}$ let $\mathrm{TY}_{\mathrm{G}}$ denote the so-called ***Tambara–Yamagami (TY) fusion ring*** given by adjoining a self-dual element $X$ to $\mathbb{Z}[\mathrm{G}]$ satisfying the fusion rules

$$gX = Xg = X, \quad X^2 = \textstyle\sum_{g\in\mathrm{G}} g.$$

(Here we use the terminology from above for the $\mathbb{Z}$ algebra $K_0^{\oplus}(_)$ itself.)

**Example 8I.4.** For $\mathrm{G} = \mathbb{Z}/3\mathbb{Z}$ the fusion rules *etc.* of $\mathrm{TY}_{\mathbb{Z}/3\mathbb{Z}}$ are

| $\otimes$ | 1 | 2 | $X$ |
|---|---|---|---|
| 1 | 2 | 0 | $X$ |
| 2 | 0 | 1 | $X$ |
| $X$ | $X$ | $X$ | $0+1+2$ |

$$\Gamma_X = \begin{pmatrix} 0 & 0 & 0 & 1 \\ 0 & 0 & 0 & 1 \\ 0 & 0 & 0 & 1 \\ 1 & 1 & 1 & 0 \end{pmatrix} = \begin{array}{ccccc} & & 1 & & \\ & & \downarrow\uparrow & & \\ 0 & \rightleftarrows & X & \rightleftarrows & 2 \end{array}, \quad \mathrm{PFdim}(X) = \sqrt{3}.$$

For general $\mathrm{G}$ we have

$$\mathrm{PFdim}(g) = 1, \quad \mathrm{PFdim}(X) = \sqrt{\#\mathrm{G}}, \quad \mathrm{PFdim}(\mathrm{TY}_G) = 2\#\mathrm{G},$$

as one easily checks. ◇

**Proposition 8I.5.** *Let* $\mathrm{G}$ *be a finite group. If* $\mathbf{C} \in \mathbf{wmFus}$ *is* $\mathbb{C}$ *linear and has* $K_0^{\oplus}(\mathbf{C}) \cong \mathrm{TY}_{\mathrm{G}}$ *as* $\mathbb{Z}$ *algebras, then* $\mathrm{G}$ *is abelian. Moreover, for any abelian* $\mathrm{G}$ *there exists a* $\mathbf{C} \in \mathbf{wmFus}$ *with* $K_0^{\oplus}(\mathbf{C}) \cong \mathrm{TY}_{\mathrm{G}}$ *as* $\mathbb{Z}$ *algebras, and such weakly multi fusion categories are parameterized (up to* $\simeq_{\mathbb{C}\oplus\star}$ *equivalence) by symmetric isomorphisms* $\mathrm{G} \xrightarrow{\cong} \mathrm{G}^{\vee}$ *and a choice of sign.*

*Proof.* This is the main result of [**TY98**]. See also Section 8J where we do proof this diagrammatically for $G = \mathbb{Z}/2\mathbb{Z}$. □

Let us continue by fixing the rank:

**Proposition 8I.6.** *We have the following.*

(i) *If* $\mathbf{C} \in \mathbf{wmFiat}$ *is* $\Bbbk$ *linear of rank* $\mathrm{rk}(\mathbf{C}) = 1$*, then* $\mathbf{C} \simeq_{\Bbbk\oplus\star} \mathbf{fdVec}_{\Bbbk}$.

(ii) *If* $\mathbf{D} \in \mathbf{wTen}$ *is* $\Bbbk$ *linear of rank* $\mathrm{rk}(\mathbf{D}) = 1$ *and* $\Bbbk$ *is of characteristic zero, then* $\mathbf{D} \simeq_{\Bbbk e\star} \mathbf{fdVec}_{\Bbbk}$.

*Proof.* (i). Note that any object $\mathtt{X} \in \mathbf{C}$ is a direct sum of the unique indecomposable $\mathtt{Z}$, *i.e.* there exists a $k \in \mathbb{Z}_{\geq 0}$ such that $\mathtt{X} \cong k \cdot \mathtt{Z}$. Hence, we have clearly

$$(m \leq k, l) \Leftrightarrow \begin{array}{ccc} & & k\cdot\mathtt{Z} \\ & \overset{\exists!}{\swarrow} {\scriptstyle \mathrm{u}} & \downarrow {\scriptstyle \mathrm{p}} \\ l\cdot\mathtt{Z} & \xrightarrow{\ \mathrm{f}\ } & m\cdot\mathtt{Z} \end{array},$$

showing that $\mathtt{X}$ is projective. Hence, we are done by *e.g.* Theorem 7E.14 since $\mathbf{C}$ has to be semisimple.

$\mathbb{1} = k \cdot \mathtt{Z}$ for some $k \in \mathbb{Z}_{>0}$, which implies that $\mathbb{1}$ is projective.

(ii). Recall that in each abelian category one can define the abelian group of extensions $\mathrm{Ext}^1_{\mathbf{C}}(\mathtt{X}, \mathtt{Y})$ to be the equivalence class (for an appropriate equivalence) of SES of the form

$$\mathtt{X} \xhookrightarrow{\mathrm{i}} \mathtt{E} \xrightarrow{\mathrm{p}} \mathtt{Y}$$

The SES of this form which split, *i.e.* where $\mathtt{E} \cong \mathtt{X} \oplus \mathtt{Y}$, are trivial in $\mathrm{Ext}^1_{\mathbf{C}}(\mathtt{X}, \mathtt{Y})$.

Back to $\mathbf{D} \in \mathbf{wTen}$, we claim that $\mathrm{Ext}^1_{\mathbf{D}}(\mathbb{1}, \mathbb{1}) = 0$.

To this end, suppose the converse. We want to show that $\mathrm{End}_{\mathbf{D}}(\mathtt{P}_{\mathbb{1}})$ has under this assumption infinitely many modules of dimension one, which is a contradiction since $\mathrm{End}_{\mathbf{D}}(\mathtt{P}_{\mathbb{1}})$ is a finite dimensional $\Bbbk$ algebra. Let $\mathtt{E}$ be a non-trivial extension of $\mathbb{1}$ by itself. Then $\mathrm{Hom}_{\mathbf{D}}(\mathtt{P}_{\mathbb{1}}, \mathtt{E})$ is of dimension two, has a filtration of length 2 with quotients isomorphic to $\mathrm{Hom}_{\mathbf{D}}(\mathtt{P}_{\mathbb{1}},\mathbb{1})$. Note that $\mathrm{End}_{\mathbf{D}}(\mathtt{P}_{\mathbb{1}})$ acts on both, $\mathrm{Hom}_{\mathbf{D}}(\mathtt{P}_{\mathbb{1}},\mathtt{E})$ and $\mathrm{Hom}_{\mathbf{D}}(\mathtt{P}_{\mathbb{1}},\mathbb{1})$, from the right. Thus, taking all of this together and letting $d_0\colon \mathrm{End}_{\mathbf{D}}(\mathtt{P}_{\mathbb{1}})\to \Bbbk$ denote the character obtained from the right action on $\mathrm{Hom}_{\mathbf{D}}(\mathtt{P}_{\mathbb{1}},\mathbb{1})$, we can find a basis of $\mathrm{Hom}_{\mathbf{D}}(\mathtt{P}_{\mathbb{1}},\mathtt{E})$ such that the action matrices of $\mathrm{End}_{\mathbf{D}}(\mathtt{P}_{\mathbb{1}})$ take the form

$$M_a=\begin{pmatrix} d_0(a) & d_1(a)\\ 0 & d_0(a)\end{pmatrix},\begin{pmatrix} d_0(ab) & d_1(ab)\\ 0 & d_0(ab)\end{pmatrix}=M_{ab}=M_aM_b=\begin{pmatrix} d_0(ab) & d_0(a)d_1(b)+d_1(a)d_0(b)\\ 0 & d_0(ab)\end{pmatrix},$$

where $d_1(a)\neq 0$, satisfying $d_1(ab)=d_0(a)d_1(b)+d_1(a)d_0(b)$, as indicated. Similarly, for any $k\in\mathbb{Z}_{>0}$, we get a $2^k$ dimensional $\mathrm{End}_{\mathbf{D}}(\mathtt{P}_{\mathbb{1}})$ module $\mathrm{Hom}_{\mathbf{D}}(\mathtt{P}_{\mathbb{1}},\mathtt{E}^k)$, and one can show that the corresponding $d_k(_)$ satisfies a recursive equality of the form

$$d_k(ab)=\textstyle\sum_{i=0}^{k}\binom{k}{i}d_i(a)d_{k-i}(b).$$

This implies that we can define infinitely many distinct one dimensional $\mathrm{End}_{\mathbf{D}}(\mathtt{P}_{\mathbb{1}})$ modules by the formula

$$e_x(a)=\textstyle\sum_{i=0}^{\infty}\frac{1}{i!}\cdot d_i(a)x^it^i\in\Bbbk((t)).$$

(Here we use that $\Bbbk$ is of characteristic zero because we need $\frac{1}{i!}$, and this formally speaking ends in $\Bbbk((t))$.) A contradiction, and we get $\mathrm{Ext}^1_{\mathbf{D}}(\mathbb{1},\mathbb{1})=0$.

Finally, $\mathrm{Ext}^1_{\mathbf{D}}(\mathbb{1},\mathbb{1})=0$ implies that $\mathbf{D}$ is semisimple, and the claim follows. □

**Example 8I.7.** Note the difference between Proposition 8I.6.(i) and (ii): The first assumes the number of indecomposables to be one, the other assumes the number of simples to be one. In fact, as we have seen in Example 6K.18, there are examples of tensor categories with one simple objects which are not equivalent to $\mathbf{fdVec}_{\Bbbk}$. ◇

Let us try to go to rank 2:

**Proposition 8I.8.** *Let $\mathbf{C}\in\mathbf{wmFus}$ be $\mathbb{C}$ linear of rank $\mathrm{rk}(\mathbf{C})=2$. Then $\mathbf{C}$ is equivalent (as a fusion category) to one of the following cases.*

- $\mathbf{C}\simeq_{\mathbb{C}\oplus\star}(\mathbf{Vec}_{\mathbb{C}}\oplus\mathbf{Vec}_{\mathbb{C}})$.
- $\mathbf{C}\simeq_{\mathbb{C}\oplus\star}\mathbf{Vec}_{\mathbb{C}\oplus}(\mathbb{Z}/2\mathbb{Z})$.
- $\mathbf{C}\simeq_{\mathbb{C}\oplus\star}\mathbf{Vec}^{\omega}_{\mathbb{C}\oplus}(\mathbb{Z}/2\mathbb{Z})$ *for the non-trivial* $\omega\in H^3(\mathbb{Z}/2\mathbb{Z},\mathbb{C}^*)\cong\mathbb{Z}/2\mathbb{Z}$.
- $\mathbf{C}\simeq_{\mathbb{C}\oplus\star}\mathbf{fdMod}^q_5(\mathfrak{so}_3)\simeq_{\mathbb{C}\oplus\star}\mathbf{Fib}$ *for* $q=\exp(\pi i/5)$.
- $\mathbf{C}\simeq_{\mathbb{C}\oplus\star}\mathbf{fdMod}^q_5(\mathfrak{so}_3)$ *for* $q=\exp(2\pi i/5)$.

*Proof.* Let us sketch the proof, and a general proof strategy, details can be found in [Ost03].

We start by observing that, if $\mathbf{C}$ is not transitive, then each simple spans a copy of $\mathbf{Vec}_{\mathbb{C}}$ and we are in the first case. Similarly if $\mathbb{1}$ is not simple.

Thus, we can assume that $\mathbf{C}$ is transitive and $\mathbb{1}\in\mathrm{Si}(\mathbf{C})$. In this case we have another self-dual simple object $\mathtt{L}$ and the fusion rules

$$\mathtt{L}^2\cong m\cdot\mathbb{1}\oplus n\cdot\mathtt{L}. \tag{8I-9}$$

First, the coefficient $m$ of $\mathbb{1}$ has to be one, by rigidity. The main work is now to show that there is no fusion category $\mathbf{C}$ for which $n>2$ in (8I-9). This is non-trivial and needs some clever arguments, and is the main point of [Ost03]:

(8I-10) There is no fusion category with fusion rules as in (8I-9) for $m\neq 1$ and $n>2$.

So let us assume that $n=0$ and $n=1$ are the only possible solutions. In both cases we already know solutions, namely the above listed cases 2 and 3 for $n=0$, where $K_0(\mathbf{C})\cong\mathbb{Z}[\mathbb{Z}/2\mathbb{Z}]$, respectively 4 and 5 for $n=1$. A careful study of the associativity constrains (as we already did for $K_0(\mathbf{C})\cong\mathbb{Z}[\mathbb{Z}/2\mathbb{Z}]$ throughout the previous sections) shows that there can not be other solutions. □

Note that the proof of Proposition 8I.8 had three main features which are part of a general strategy to classify fiat and fusion categories, in increasing difficulty:

- write down the possible solutions on the Grothendieck level, which was (8I-9) above;
- use the categorical properties of $\mathbf{C}$ to rule out cases, which was (8I-10) above;
- in the remaining cases construct the categories and analyze the various categorical constrains to show that one has found all solutions, which was the last step above.

**Example 8I.11.** To rule out the cases $n > 2$ in (8I-9) requires all assumptions. For example, if one drops the assumption on $\mathbf{C}$ to be semisimple, then $m = 0$ and $n$ arbitrary can indeed occur. We have already seen an example, namely $\mathbf{fdProj}\big(\overline{\mathbb{F}}_p(\mathbb{Z}/p\mathbb{Z})\big)$, see *e.g.* (7H-7), where $\mathtt{P}_{\mathbb{1}}\mathtt{P}_{\mathbb{1}} \cong p \cdot \mathtt{P}_{\mathbb{1}}$. (Formally speaking, we would need to adjoin the monoidal unit to $\mathbf{fdProj}\big(\overline{\mathbb{F}}_p(\mathbb{Z}/p\mathbb{Z})\big)$ to make this example solid.) $\diamond$

To continue to try to classify fusion categories by their ranks get tricky, and is doomed to fail from some point on. Let us state the $\mathrm{rk}(\mathbf{C}) = 3$ result, ordered as in the strategy list above:

**Proposition 8I.12.** *Let $\mathbf{C} \in \mathbf{Fus}$ be $\mathbb{C}$ linear of rank $\mathrm{rk}(\mathbf{C}) = 3$. Let $\mathbf{Si}(\mathbf{C}) = \{\mathtt{L}_1 = \mathbb{1}, \mathtt{L}_2, \mathtt{L}_3\}$. Then:*

- *The only possible fusion rules of $\mathbf{C}$ are:*

| $\otimes$ | $\mathtt{L}_2$ | $\mathtt{L}_3$ |
|---|---|---|
| $\mathtt{L}_2$ | $\mathtt{L}_3$ | $\mathbb{1}$ |
| $\mathtt{L}_3$ | $\mathbb{1}$ | $\mathtt{L}_2$ |

,

| $\otimes$ | $\mathtt{L}_2$ | $\mathtt{L}_3$ |
|---|---|---|
| $\mathtt{L}_2$ | $\mathbb{1} \oplus m \cdot \mathtt{L}_2 \oplus k \cdot \mathtt{L}_3$ | $k \cdot \mathtt{L}_2 \oplus l \cdot \mathtt{L}_3$ |
| $\mathtt{L}_3$ | $k \cdot \mathtt{L}_2 \oplus l \cdot \mathtt{L}_3$ | $\mathbb{1} \oplus l \cdot \mathtt{L}_2 \oplus n \cdot \mathtt{L}_3$ |

,

*where $k, l, m, n \in \mathbb{Z}_{\geq 0}$ satisfying $k^2 + l^2 = kn + lm + 1$.*

- *Only the following cases can occur:*

$(A)$:

| $\otimes$ | $\mathtt{L}_2$ | $\mathtt{L}_3$ |
|---|---|---|
| $\mathtt{L}_2$ | $\mathtt{L}_3$ | $\mathbb{1}$ |
| $\mathtt{L}_3$ | $\mathbb{1}$ | $\mathtt{L}_2$ |

, $(B)$:

| $\otimes$ | $\mathtt{L}_2$ | $\mathtt{L}_3$ |
|---|---|---|
| $\mathtt{L}_2$ | $\mathbb{1} \oplus \mathtt{L}_2 \oplus \mathtt{L}_3$ | $\mathtt{L}_2 \oplus \mathtt{L}_3$ |
| $\mathtt{L}_3$ | $\mathtt{L}_2 \oplus \mathtt{L}_3$ | $\mathbb{1} \oplus \mathtt{L}_2$ |

, $(C)$:

| $\otimes$ | $\mathtt{L}_2$ | $\mathtt{L}_3$ |
|---|---|---|
| $\mathtt{L}_2$ | $\mathbb{1} \oplus \mathtt{L}_3$ | $\mathtt{L}_2$ |
| $\mathtt{L}_3$ | $\mathtt{L}_2$ | $\mathbb{1}$ |

,

$(D)$:

| $\otimes$ | $\mathtt{L}_2$ | $\mathtt{L}_3$ |
|---|---|---|
| $\mathtt{L}_2$ | $\mathbb{1} \oplus \mathtt{L}_2 \oplus \mathtt{L}_3$ | $\mathtt{L}_2$ |
| $\mathtt{L}_3$ | $\mathtt{L}_2$ | $\mathbb{1}$ |

, $(E)$:

| $\otimes$ | $\mathtt{L}_2$ | $\mathtt{L}_3$ |
|---|---|---|
| $\mathtt{L}_2$ | $\mathbb{1} \oplus 2 \cdot \mathtt{L}_2 \oplus \mathtt{L}_3$ | $\mathtt{L}_2$ |
| $\mathtt{L}_3$ | $\mathtt{L}_2$ | $\mathbb{1}$ |

.

- *For (A) we have the solutions $\mathbf{Vec}^{\omega}_{\mathbb{C}\oplus}(\mathbb{Z}/3\mathbb{Z})$ (note that $H^3(\mathbb{Z}/3\mathbb{Z}, \mathbb{C}^*) \cong \mathbb{Z}/3\mathbb{Z}$).*
- *For (B) we have the solutions $\mathbf{fdMod}^q_7(\mathfrak{so}_3)$.*
- *For (C) we have the solutions $\mathbf{fdMod}^q_4(\mathfrak{sl}_2)$.*
- *For (D) we have the solutions $\mathbf{fdMod}(\mathrm{S}_3)$ and twists (as in Example 8G.5).*
- *For (E) we have two solutions, a fusion category associated with a subfactor of type $E_6$ or its Galois conjugate. (See e.g. [**HH09**] for the definitions.)*
- *There are no other solutions.*

*Proof.* This is proven in [**Ost15**]. $\square$

Let us have a look now at the PF dimension.

**Theorem 8I.13.** *Let $\mathrm{F} \in \mathbf{Hom}_{\Bbbk\oplus\star}(\mathbf{C}, \mathbf{D})$, where $\mathbf{C}, \mathbf{D} \in \mathbf{wmFiat}$ are $\Bbbk$ linear. Then:*

*(i) If $\mathrm{F}$ is fully faithful, then*

$$\mathrm{PFdim}(\mathbf{C}) \leq \mathrm{PFdim}(\mathbf{D}),$$

*with equality achieved if and only if $\mathrm{F}$ is an equivalence.*

*(ii) If $\mathrm{F}$ is fully faithful, then*

$$\mathrm{PFdim}(\mathbf{C}) \geq \mathrm{PFdim}(\mathbf{D}),$$

*with equality achieved if and only if $\mathrm{F}$ is an equivalence.*

*Proof.* This can be proven *mutatis mutandis* as in [**EGNO15**, Propositions 6.3.3 and 6.3.4]. $\square$

**Proposition 8I.14.** *If $\mathbf{C} \in \mathbf{wFiat}$ is $\Bbbk$ linear of PF dimension $\mathrm{PFdim}(\mathbf{C}) = 1$, then $\mathbf{C} \simeq_{\Bbbk\oplus\star} \mathbf{fdVec}_{\Bbbk}$.*

*Proof.* We already know that $\mathrm{PFdim}(\mathbf{C}) \geq 1$, see Proposition 8E.6. Moreover, there is always a fully faithful functor $\mathrm{F}\colon \mathbf{fdVec}_{\Bbbk} \to \mathbf{C}$ given by $\Bbbk = \mathbb{1} \mapsto \mathbb{1}$. Thus, the claim follows from Theorem 8I.13.(i). $\square$

The analog of "All finite groups of prime order are cyclic." is:

**Proposition 8I.15.** *If $\mathbf{C} \in \mathbf{wFiat}$ is $\mathbb{C}$ linear and satisfies $\mathrm{PFdim}(\mathbf{C}) = p$ for $p \in \mathbb{Z}_{\geq 0}$ being a prime, then $\mathbf{C} \simeq_{\Bbbk\oplus\star} \mathbf{Vec}^{\omega}_{\mathbb{C}\oplus}(\mathbb{Z}/p\mathbb{Z})$.*

*Proof.* See [**ENO05**, Corollary 8.30]. $\square$

An extraordinary fact is that PF dimensions are quantized:

**Proposition 8I.16.** *Let $\mathbf{C} \in \mathbf{wFus}$ be $\Bbbk$ linear, and let $\mathtt{L}_1 \in \mathrm{Si}(\mathbf{C})$ be a fusion generator of PF dimension $\mathrm{PFdim}(\mathtt{L}_1) < 2$. Then:*

*(i) $\mathrm{PFdim}(\mathtt{L}_1) = 2\cos(\pi/k)$ for some $k$.*

(ii) *The fusion graph of* $L_1$ *is one of the following ADE types:*

$$\text{Type } A\colon\ \mathbb{1} \rightleftarrows L_1 \rightleftarrows L_2 \rightleftarrows \dots \rightleftarrows L_k\ ,$$

$$\text{Type } D\colon\ \mathbb{1} \rightleftarrows L_1 \rightleftarrows L_2 \rightleftarrows \dots \rightleftarrows L_k \begin{array}{l} \nearrow\!\!\swarrow\ L_{k+1} \\ \searrow\!\!\nwarrow\ L'_{k+1} \end{array}\ ,$$

(8I-17)
$$\text{Type } E_6\colon\ \begin{array}{c} \qquad\qquad\quad L_5 \\ \qquad\qquad\quad \downarrow\uparrow \\ \mathbb{1} \rightleftarrows L_1 \rightleftarrows L_2 \rightleftarrows L_3 \rightleftarrows L_5 \end{array}\ ,$$

$$\text{Type } E_7\colon\ \begin{array}{c} \qquad\qquad\quad L_6 \\ \qquad\qquad\quad \downarrow\uparrow \\ \mathbb{1} \rightleftarrows L_1 \rightleftarrows L_2 \rightleftarrows L_3 \rightleftarrows L_4 \rightleftarrows L_5 \end{array}\ ,$$

$$\text{Type } E_8\colon\ \begin{array}{c} \qquad\qquad\quad L_7 \\ \qquad\qquad\quad \downarrow\uparrow \\ \mathbb{1} \rightleftarrows L_1 \rightleftarrows L_2 \rightleftarrows L_3 \rightleftarrows L_4 \rightleftarrows L_5 \rightleftarrows L_6 \end{array}\ .$$

(iii) *For all the graphs in* (8I-17) *there exists a fusion category with a fusion generator having the corresponding fusion graph.*

(iv) *In type A the fusion category is of the form* $\mathbf{fdMod}_k^q(\mathfrak{sl}_2)$.

*Proof.* See e.g. [FK93, Chapter 8]. □

8J. **Webs and Tambara–Yamagami fusion categories.** We now prove Proposition 8I.5 for $G = \mathbb{Z}/2\mathbb{Z} = \{0,1\}$ by constructing the corresponding category diagrammatically. Recall the category $\mathbf{Vec}_{\mathbb{C}}^{dia}(\mathbb{Z}/2\mathbb{Z})$ from Example 3G.6, which is a generator-relation category given by

(8J-1)
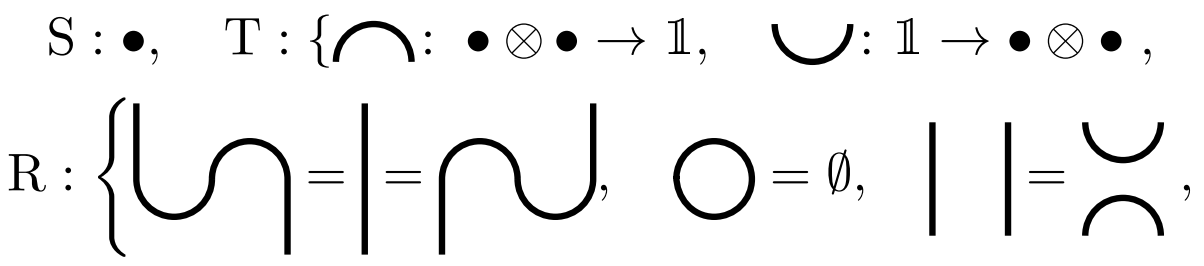

where we identified $\mathbb{1}$ with 0, and $\bullet$ with 1. We extend $\mathbf{Vec}_{\mathbb{C}}^{dia}(\mathbb{Z}/2\mathbb{Z})$ to the category $\mathbf{Vec}_{\mathbb{C}}^{dia+X}(\mathbb{Z}/2\mathbb{Z})$ by adding an extra object and and also additional morphism generators, namely the following (which we add to the generators in Equation 8J-1):

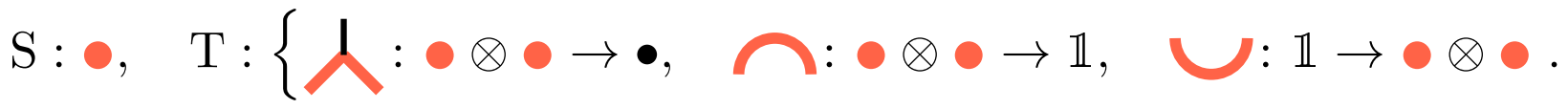

The red bullet $\bullet$ is representing $X$.

*Remark* 8J.2. The reader with a black-and-white copy can distinguish the new object and edges, as they appear slightly larger and thicker, respectively. ◇

**Definition 8J.3.** For an algebraically closed field $\mathbb{K}$ of characteristic not five, let $\mathbf{TY}_{\mathbb{Z}/2\mathbb{Z}} = (\mathbf{TY}_{\mathbb{Z}/2\mathbb{Z}})_{\mathbb{K}\oplus}$ denote the quotient of $\mathbf{Vec}_{\mathbb{C}}^{dia+X}(\mathbb{Z}/2\mathbb{Z})_{\mathbb{K}\oplus}$ by the relations, meaning R, in Example 3G.6 with one edge colored black and two edges colored red around each trivalent vertex, except for ***circle evaluation***, and the ***H=I relation*** which we replace by:

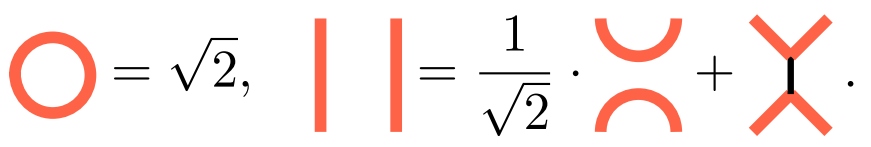

We call $\mathbf{TY}_{\mathbb{Z}/2\mathbb{Z}}$ the ***TY category for*** $\mathbb{Z}/2\mathbb{Z}$. ◇

**Proposition 8J.4.** *The category* $\mathbf{TY}_{\mathbb{Z}/2\mathbb{Z}}$ *is a semisimple fiat category with three simple objects* $L_1 = \mathbb{1}$, $L_2 = \bullet$ *and* $L_3 = \bullet$ *with*

$$L_2^2 \cong \mathbb{1}, \quad L_2L_3 \cong L_3L_2 \cong L_3, \quad L_3^2 \cong \mathbb{1} \oplus L_2.$$

*Proof.* The proof is almost the same as the proof of Proposition 8F.3, so we will only sketch the differences. Working out the details is Exercise 8L.7.

First, the black-only-part can be discussed as in Section 3G. Moreover, we have

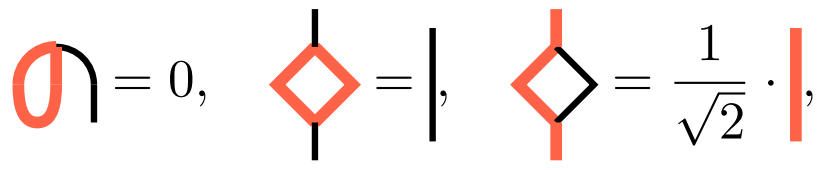

as one easily checks. These relations imply •• ≅ •• ≅ •. Then Lemma 8F.7 holds *verbatim*, while the primitive idempotents in $\mathrm{End}_{\mathbf{TY}_{\mathbb{Z}/2\mathbb{Z}}}(\bullet^2)$ are

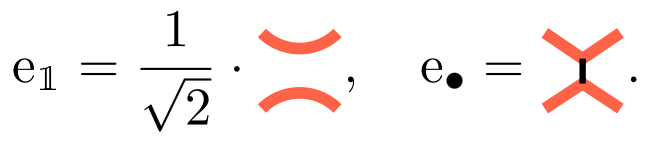

These, together with the diagram

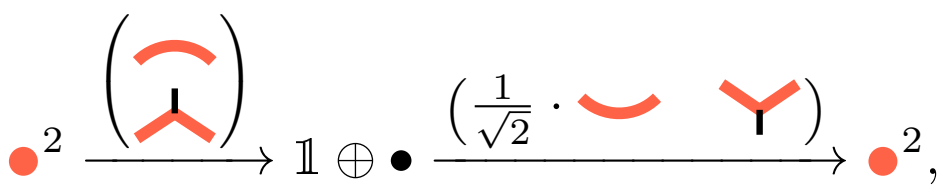

can be used to show that •• ≅ $\mathbb{1} \oplus \bullet$. The final claims can then be tackled as an exercise. □

*Proof of Proposition 8I.5 for* $G = \mathbb{Z}/2\mathbb{Z}$. We first note that using $-\sqrt{2}$ instead of $\sqrt{2}$ in the definition of $\mathbf{TY}_{\mathbb{Z}/2\mathbb{Z}}$ gives another semisimple fiat category $\mathbf{TY}^{-}_{\mathbb{Z}/2\mathbb{Z}}$ satisfying Proposition 8J.4 that is not equivalent (as a fiat category) to $\mathbf{TY}_{\mathbb{Z}/2\mathbb{Z}}$.

We want to show that there are no other categories satisfying Proposition 8J.4. To see this, note that

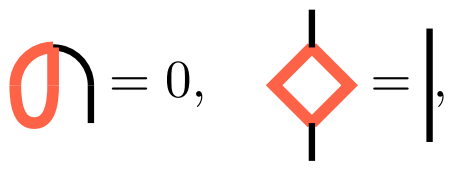

have to hold (the second potentially up to an invertible scalar, but one can divide by that scalar). Now set

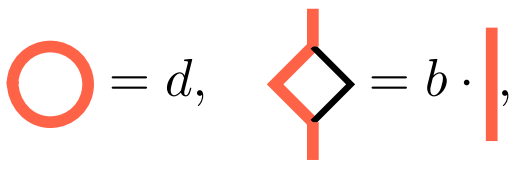

for some variables $d$ and $b$ that we will determine momentarily. Now, as in the proof of Lemma 8F.5, we compute a pairing matrix with

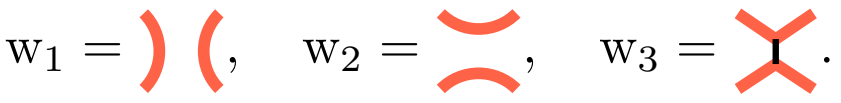

We get the following matrix, that we need of rank two:

$$P = \begin{pmatrix} d^2 & d & bd \\ d & d^2 & 0 \\ bd & 0 & bd \end{pmatrix},$$

and its determinant is $\det P = -bd^3(1 - d^2 + bd)$. Semisimplicity implies $b, d \neq 0$, so we get $b = \frac{d^2-1}{d}$ as the only solution for $\det P = 0$, which also implies $d \neq \pm 1$ as $b \neq 0$. Under these assumptions, the kernel of $P$ is spanned by the vector $(-1, \frac{1}{d}, 1)$ giving the relation

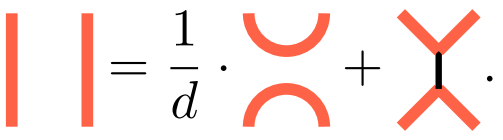

Finally, $b = \frac{1}{d}$ gives

$$\frac{1}{d} = \frac{d^2 - 1}{d},$$

with only two possible solutions: $d = \pm\sqrt{2}$. □

8K. **A pseudo classification – or, summarizing the above.** Let G be a finite group and let us call $\mathbf{Vec}^{\omega}_{\mathbb{C}\oplus}(\mathrm{G})$ for non-trivial $\omega$ a ***twist*** of $\mathbf{Vec}_{\mathbb{C}\oplus}(\mathrm{G})$. Similarly, we have ***twists*** of $\mathbf{fdMod}(\mathrm{G})$, *cf.* Example 8G.5, and we also call the Verlinde categories $\mathbf{fdMod}^q_k(\mathfrak{g})$ for $q \neq \exp(\pm\pi i/k)$ ***twists*** of the standard choice $q = \exp(\pi i/k)$. Then we have the following pseudo classification, motivated by the above.

**Theorem 8K.1** (Pseudo theorem)**.** *All* $\mathbb{C}$ *linear fusion categories are one of the following types:*

*(I) Categories of the form* $\mathbf{Vec}_{\mathbb{C}\oplus}(\mathrm{G})$ *and twists.*

*(II) Categories of the form* $\mathbf{fdMod}(\mathrm{G})$ *and twists.*

*(III) Categories of the form* $\mathbf{fdMod}^q_k(\mathfrak{g})$ *and twists.*

*(IV) Exceptions.* □

The crucial point, which we will explore in the following sections, will be:

> The main source of quantum invariants are the fusion categories of type (III).

The fusion categories of types (I), (II) and (IV) sometimes also give quantum invariants. But it turns out that types (I) and (II) give rather "boring" invariants, while type (IV) remains to be explored further.

## 8L. **Exercises.**

*Exercise* 8L.1. Try to understand the claims in Example 8B.4 and verify as many of them as possible. ◇

*Exercise* 8L.2. Redo the computation in Section 8F for swapped roles of $\phi$ and $\Phi$. Is the resulting fusion category equivalent (as a fusion category) to **Fib**? ◇

*Exercise* 8L.3. Calculate the fusion graphs and PF dimensions of $\mathbf{fdMod}\big(\overline{\mathbb{F}}_5[\mathbb{Z}/5\mathbb{Z}]\big)$ and of $\mathbf{fdMod}\big(\mathbb{C}[S_5]\big)$. The latter is a semisimple fiat category and has the fusion rules $\mathrm{L}_1 \cong \mathbb{1}$ and

| $\otimes$ | $\mathrm{L}_s$ | $\mathrm{L}_b$ | $\mathrm{L}_{s'}$ | $\mathrm{L}_{1'}$ |
|---|---|---|---|---|
| $\mathrm{L}_s$ | $\mathbb{1} \oplus \mathrm{L}_s \oplus \mathrm{L}_b \oplus \mathrm{L}_{s'}$ | $\mathrm{L}_s \oplus \mathrm{L}_{s'}$ | $\mathrm{L}_s \oplus \mathrm{L}_b \oplus \mathrm{L}_{s'} \oplus \mathrm{L}_{1'}$ | $\mathrm{L}_{s'}$ |
| $\mathrm{L}_b$ | $\mathrm{L}_s \oplus \mathrm{L}_{s'}$ | $\mathbb{1} \oplus \mathrm{L}_b \oplus \mathrm{L}_{1'}$ | $\mathrm{L}_s \oplus \mathrm{L}_{s'}$ | $\mathrm{L}_b$ |
| $\mathrm{L}_{s'}$ | $\mathrm{L}_s \oplus \mathrm{L}_b \oplus \mathrm{L}_{s'} \oplus \mathrm{L}_{1'}$ | $\mathrm{L}_s \oplus \mathrm{L}_{s'}$ | $\mathbb{1} \oplus \mathrm{L}_s \oplus \mathrm{L}_b \oplus \mathrm{L}_{s'}$ | $\mathrm{L}_s$ |
| $\mathrm{L}_{1'}$ | $\mathrm{L}_{s'}$ | $\mathrm{L}_b$ | $\mathrm{L}_s$ | $\mathbb{1}$ |

. ◇

*Exercise* 8L.4. Complete the discussion in Example 8E.7. ◇

*Exercise* 8L.5. Verify the calculations in Example 8I.4. ◇

*Exercise* 8L.6. Prove the last claim in Example 8H.4 and Example 8H.7. ◇

*Exercise* 8L.7. Fill in the gaps in Section 8J. ◇

# 9. Fusion and modular categories – definitions and graphical calculus

The question we want to address is:

> Can we separate "topologically boring" fiat categories from "topologically interesting" ones?

What does this mean? Well, think of a quantum invariant Q. The more knots or links it can distinguish, the better *cf.* Figure 23.

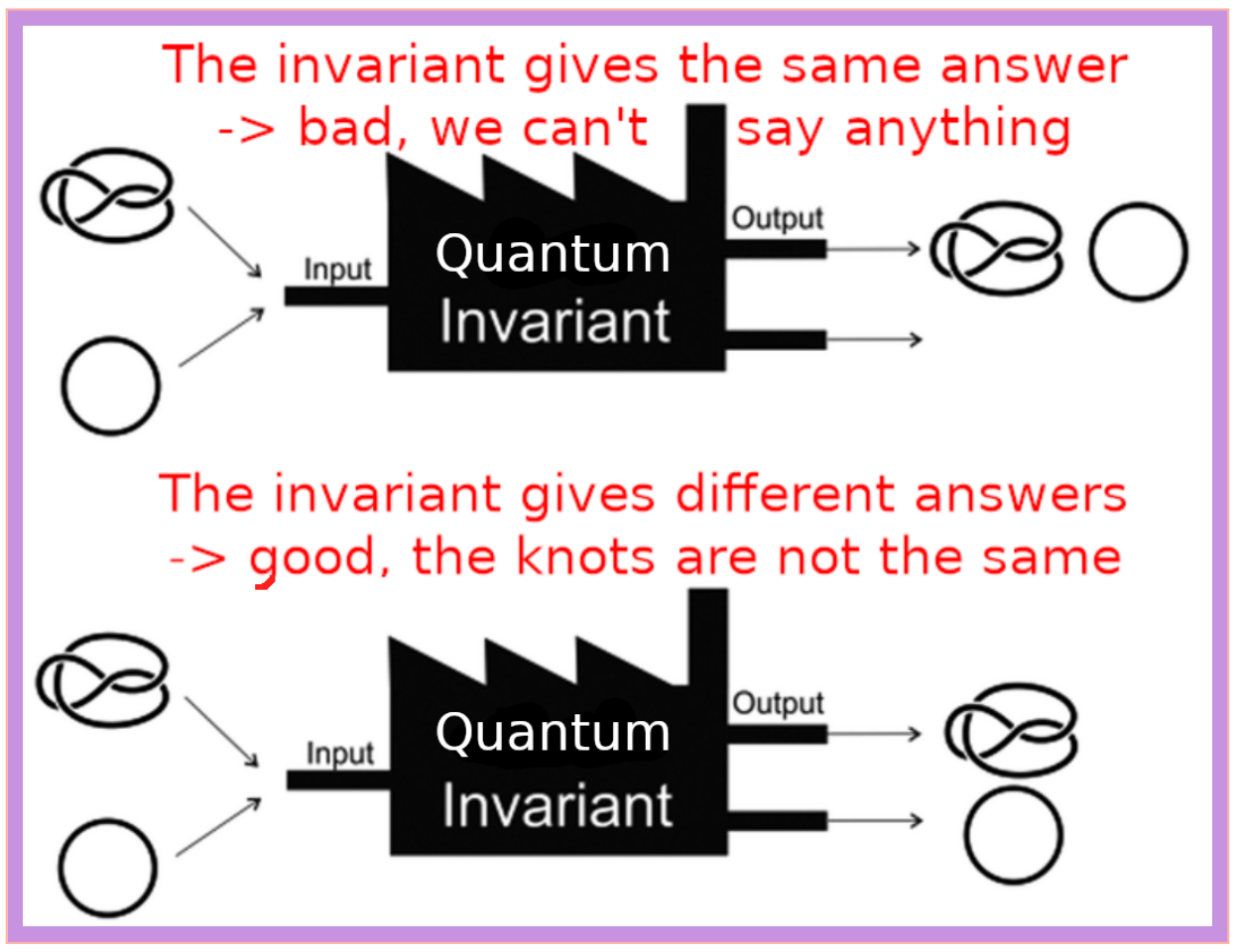

Figure 23. Note that $\mathrm{Q}(\mathrm{L}) \neq \mathrm{Q}(\mathrm{L}')$ implies $\mathrm{L} \neq \mathrm{L}'$, but the converse is not true in general. Hence, if a quantum invariant (or any invariant really, but we have not defined what these are) sends two knots to the same value, then we cannot conclude anything. Thus, a good quantum invariant should take many different values.

Picture from https://dtubbenhauer.com/lecture-geotop-2023.html

The answer will turn out to be "Yes and no". We also address the "topologically boring" and "topologically interesting" via a ***big data approach*** in Section 13 later on.

9A. **A word about conventions.** Of course, we keep the previous conventions.

*Convention* 9A.1. We will revisit several properties which we have seen before and which depend on choices such as being braided. As before we tend to write *e.g.* "ABC is XYZ" instead of the formally correct "there is a choice such that ABC is XYZ" *etc.* ◇

9B. **Hom spaces in fiat, tensor and fusion categories.** Let us start by motivating the diagrammatics which we will see below.

If $\mathbf{C} \in \mathbf{Cat}_{\mathbb{K}S}$ and $\mathrm{Si}(\mathbf{C}) = \{\mathtt{L}_1, ..., \mathtt{L}_m\}$, then Schur's lemma Lemma 6J.10 allows us to compute hom spaces as follows. Let $\mathtt{X}, \mathtt{Y} \in \mathbf{C}$, and decompose them into simples

$$\mathtt{X} \cong \bigoplus_{i=1}^{m} [\mathtt{X} : \mathtt{L}_i] \cdot \mathtt{L}_i, \quad \mathtt{Y} \cong \bigoplus_{i=1}^{m} [\mathtt{Y} : \mathtt{L}_i] \cdot \mathtt{L}_i.$$

Then we have the decomposition and dimension formulas

$$\text{(9B-1)} \qquad \mathrm{Hom}_{\mathbf{C}}(\mathtt{X}, \mathtt{Y}) \cong \bigoplus_{i=1}^{m} \mathrm{Mat}_{[\mathtt{Y}:\mathtt{L}_i] \times [\mathtt{X}:\mathtt{L}_i]}(\mathbb{K}), \quad \dim\big(\mathrm{Hom}_{\mathbf{C}}(\mathtt{X}, \mathtt{Y})\big) = \textstyle\sum_{i=1}^{m} [\mathtt{X} : \mathtt{L}_i][\mathtt{Y} : \mathtt{L}_i].$$

**Example 9B.2.** Assume that $\mathtt{X} \cong 2 \cdot \mathtt{L}_1 \oplus \mathtt{L}_2$ and $\mathtt{Y} \cong \mathtt{L}_1 \oplus \mathtt{L}_3$. Then

$$\mathrm{End}_{\mathbf{C}}(\mathtt{X}) \cong$$
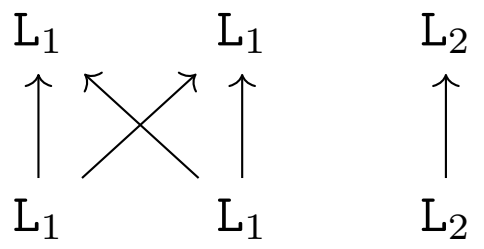
$$\cong \mathrm{Mat}_{2\times 2}(\mathbb{K}) \oplus \mathrm{Mat}_{1\times 1}(\mathbb{K}),$$

$$\mathrm{Hom}_{\mathbf{C}}(\mathtt{X}, \mathtt{Y}) \cong$$
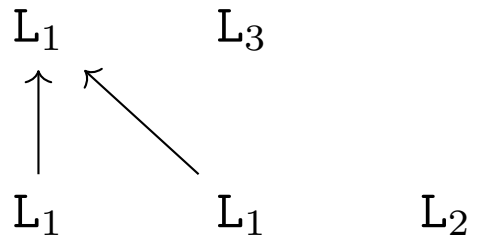
$$\cong \mathrm{Mat}_{1\times 2}(\mathbb{K}),$$

$$\mathrm{End}_{\mathbf{C}}(\mathtt{Y}) \cong$$
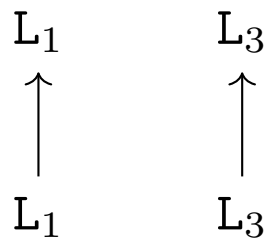
$$\cong \mathrm{Mat}_{1\times 1}(\mathbb{K}) \oplus \mathrm{Mat}_{1\times 1}(\mathbb{K}),$$

illustrates the validity of the formulas in Equation 9B-1, where each arrow represents a, up to scalars unique, basis element of the hom spaces. An object such as $\mathtt{Y}$ is also called ***multiplicity free***, referring to the decomposition of $\mathtt{Y}$ having each simple appear at most once. ◇

Note that Equation 9B-1 fails in the non-semisimple case: For $\mathbf{C} \in \mathbf{Cat}_{\mathbb{K}\oplus\mathbb{E}}$, $\mathrm{Si}(\mathbf{C}) = \{\mathtt{L}_1, ..., \mathtt{L}_m\}$, and $\mathrm{In}(\mathbf{C}) = \{\mathtt{Z}_1, ..., \mathtt{Z}_n\}$, Schur's lemma does not hold between indecomposables and hom spaces need not to be matrix $\mathbb{K}$ algebras, but rather matrix algebras over some local $\mathbb{K}$ algebra.

**Example 9B.3.** Return to Example 6K.18: In $\mathbf{fdMod}\big(\mathbb{F}_5[\mathbb{Z}/5\mathbb{Z}]\big)$ we have seen that $\mathbb{1} = \mathtt{Z}_1$ and its projective cover $\mathtt{P}_{\mathbb{1}} = \mathtt{Z}_5$ are non-isomorphic indecomposables. However, by the definition of the projective cover, the hom space between them is non-zero. ◇

However, we still have ***the idempotent decomposition of*** $\mathrm{id}_{\mathtt{X}}$, *i.e.*

$$\text{(9B-4)} \qquad \begin{gathered} \mathrm{id}_{\mathtt{X}} = \textstyle\sum_{i=1}^{n} \sum_{j=1}^{(\mathtt{X}:\mathtt{Z}_i)} \mathrm{i}_{i_j} \mathrm{p}_{i_j}, \quad \mathrm{p}_{i_k} \mathrm{i}_{j_l} = \delta_{i,j}\delta_{k,l}\mathrm{id}_{\mathtt{Z}_{i_k}}, \\ \mathrm{i}_{i_j} : \mathtt{Z}_{i_j} \hookrightarrow \mathtt{X} \ j\text{th inclusion}, \quad \mathrm{p}_{i_j} : \mathtt{X} \twoheadrightarrow \mathtt{Z}_{i_j} \ j\text{th projection}, \\ \mathrm{id}_{\mathtt{X}} = \textstyle\sum_{i=1}^{n} \mathrm{i}_i \mathrm{p}_i, \quad \mathrm{p}_i \mathrm{i}_j = \delta_{i,j} \mathrm{id}_{(\mathtt{X}:\mathtt{Z}_i)\cdot \mathtt{Z}_i}, \\ \mathrm{i}_i = \textstyle\sum_{j=1}^{(\mathtt{X}:\mathtt{Z}_i)} \mathrm{i}_{i_j} \text{ isotypic inclusion}, \quad \mathrm{p}_i = \sum_{j=1}^{(\mathtt{X}:\mathtt{Z}_i)} \mathrm{p}_{i_j} \text{ isotypic projection.} \end{gathered}$$

(This is just Equation 6B-3, but taking multiplicities into account.) The morphisms $\mathrm{i}_i$ and $\mathrm{p}_i$ are unique up to scaling, and $\mathrm{i}_i$ and $\mathrm{p}_i$ are called the ***isotypic inclusions and projections***, respectively.

**Example 9B.5.** For $\mathrm{End}_{\mathbf{C}}(\mathtt{X})$ as in Example 9B.2, but where the simples are only assumed to be indecomposable, we have

$$\mathrm{id}_{\mathtt{X}} =$$
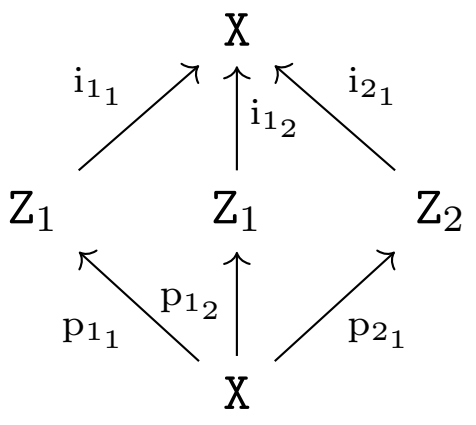

,

$$\mathrm{i}_{1_1}\mathrm{p}_{1_1} \rightsquigarrow \begin{array}{c} \\ \mathtt{Z}_1 \\ \mathtt{Z}_1 \\ \mathtt{Z}_2 \end{array}\!\!\begin{array}{c} \begin{array}{ccc}\mathtt{Z}_1 & \mathtt{Z}_1 & \mathtt{Z}_2\end{array} \\ \begin{pmatrix} 1 & 0 & 0 \\ 0 & 0 & 0 \\ 0 & 0 & 0 \end{pmatrix}\end{array}, \quad \mathrm{i}_{1_2}\mathrm{p}_{1_2} \rightsquigarrow \begin{array}{c} \\ \mathtt{Z}_1 \\ \mathtt{Z}_1 \\ \mathtt{Z}_2 \end{array}\!\!\begin{array}{c} \begin{array}{ccc}\mathtt{Z}_1 & \mathtt{Z}_1 & \mathtt{Z}_2\end{array} \\ \begin{pmatrix} 0 & 0 & 0 \\ 0 & 1 & 0 \\ 0 & 0 & 0 \end{pmatrix}\end{array}, \quad \mathrm{i}_{2_1}\mathrm{p}_{2_1} \rightsquigarrow \begin{array}{c} \\ \mathtt{Z}_1 \\ \mathtt{Z}_1 \\ \mathtt{Z}_2 \end{array}\!\!\begin{array}{c} \begin{array}{ccc}\mathtt{Z}_1 & \mathtt{Z}_1 & \mathtt{Z}_2\end{array} \\ \begin{pmatrix} 0 & 0 & 0 \\ 0 & 0 & 0 \\ 0 & 0 & 1 \end{pmatrix}\end{array},$$

$$\mathrm{i}_{1}\mathrm{p}_{1} \rightsquigarrow \begin{array}{c} \\ \mathtt{Z}_1 \\ \mathtt{Z}_1 \\ \mathtt{Z}_2 \end{array}\!\!\begin{array}{c} \begin{array}{ccc}\mathtt{Z}_1 & \mathtt{Z}_1 & \mathtt{Z}_2\end{array} \\ \begin{pmatrix} 1 & 0 & 0 \\ 0 & 1 & 0 \\ 0 & 0 & 0 \end{pmatrix}\end{array}, \quad \mathrm{i}_{2}\mathrm{p}_{2} \rightsquigarrow \begin{array}{c} \\ \mathtt{Z}_1 \\ \mathtt{Z}_1 \\ \mathtt{Z}_2 \end{array}\!\!\begin{array}{c} \begin{array}{ccc}\mathtt{Z}_1 & \mathtt{Z}_1 & \mathtt{Z}_2\end{array} \\ \begin{pmatrix} 0 & 0 & 0 \\ 0 & 0 & 0 \\ 0 & 0 & 1 \end{pmatrix}\end{array}.$$

Note that the colored zero have to be zero in the semisimple case, by Schur's lemma, but not necessarily in general. ◇

9C. **Feynman diagrams for fiat, tensor and fusion categories.** Let us assume that have have a strict multi fiat category. Then we get, of course, the diagrammatic calculus for pivotal categories as in Section 4G. Additionally, we want to keep track of the morphisms $\mathrm{i}_{i_j}$ and $\mathrm{p}_{i_j}$ from Equation 9B-4, as well as simples. We use the conventions:

(9C-1)

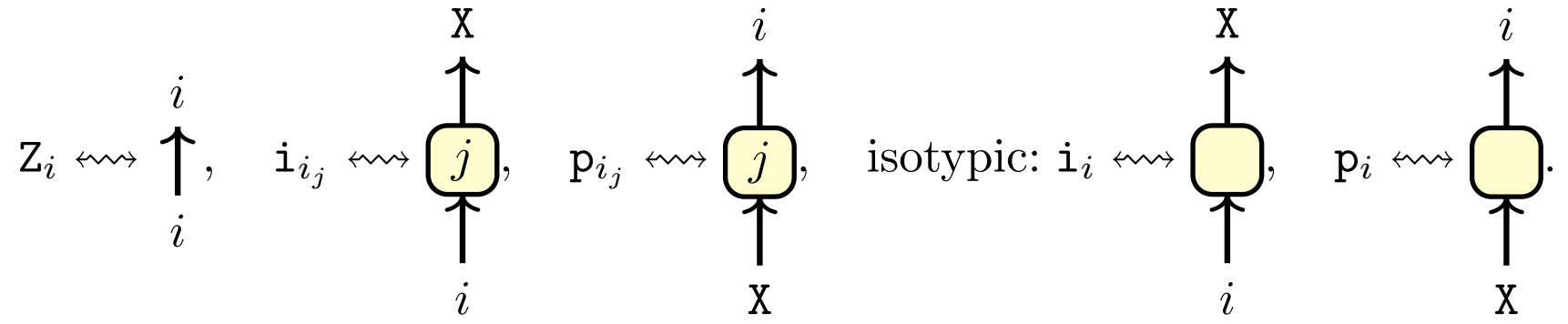

Note that we can distinguish between the inclusions and the projections in Equation 9C-1 by the labeling of the strands, so we can just use colored boxes as indicated. The relations in Equation 9B-4 then are *e.g.*

(9C-2)

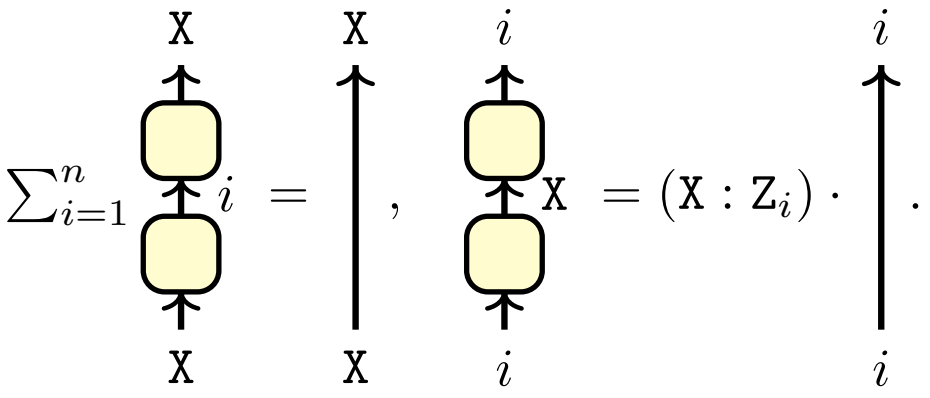

We, of course, still have the topological relations which we have seen, *e.g.* sliding

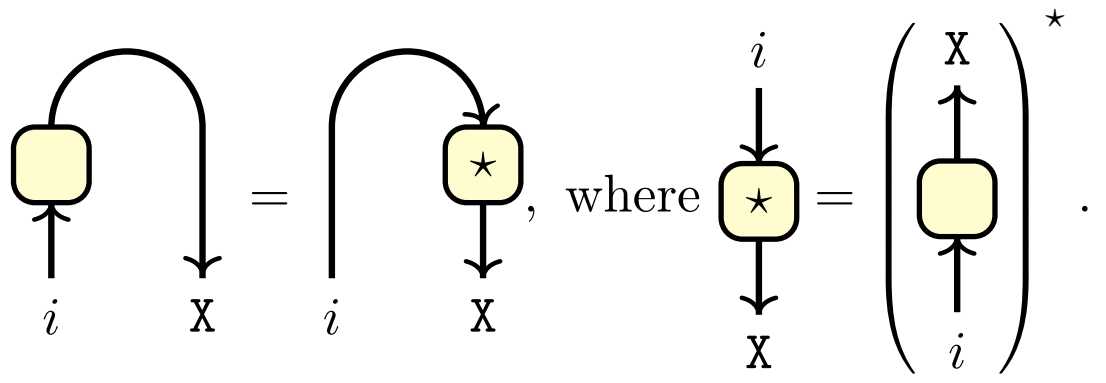

If our category of interest is additionally braided, then we have the power of the Reidemeister calculus, see Section 5F, as well, *e.g.*

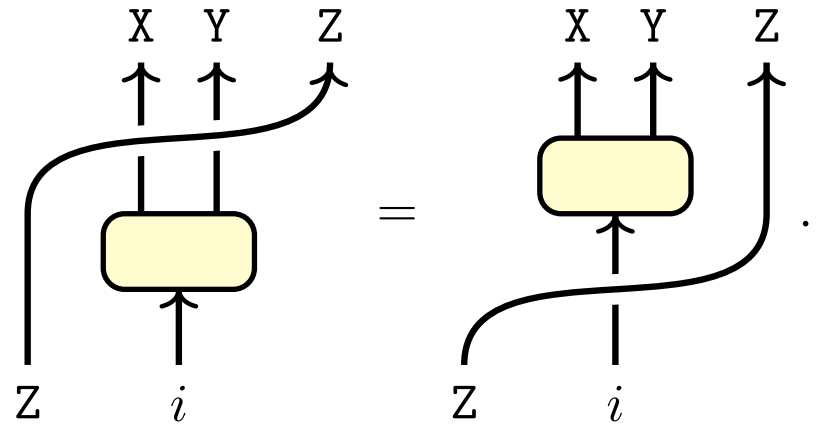

Note that this graphical calculus is, having *e.g.* Theorem 4G.11 established, is automatically consisted as we simply used a special notation for special morphisms.

9D. **Traces and dimensions revisited.** Recall that we had the notions of traces and dimensions for pivotal categories, *cf.* Section 4H. We can say a bit more now. But first an evident lemma:

**Lemma 9D.1.** *Let* $\mathbf{C} \in \mathbf{PCat}_{\oplus}$. *Then*

$$\mathrm{tr}^{\mathbf{C}}(\mathrm{f} \oplus \mathrm{g}) = \mathrm{tr}^{\mathbf{C}}(\mathrm{f}) + \mathrm{tr}^{\mathbf{C}}(\mathrm{g}), \quad {}^{\mathbf{C}}\mathrm{tr}(\mathrm{f} \oplus \mathrm{g}) = {}^{\mathbf{C}}\mathrm{tr}(\mathrm{f}) + {}^{\mathbf{C}}\mathrm{tr}(\mathrm{g}),$$

$$\dim^{\mathbf{C}}(\mathtt{X}\oplus\mathtt{Y})=\dim^{\mathbf{C}}(\mathtt{X})+\dim^{\mathbf{C}}(\mathtt{Y}),\quad {}^{\mathbf{C}}\dim(\mathtt{X}\oplus\mathtt{Y})={}^{\mathbf{C}}\dim(\mathtt{X})+{}^{\mathbf{C}}\dim(\mathtt{Y}),$$

*for all* $\mathtt{X},\mathtt{Y}\in\mathbf{C}$ *and* $(\mathrm{f}\colon\mathtt{X}\to\mathtt{X}),(\mathrm{g}\colon\mathtt{Y}\to\mathtt{Y})\in\mathbf{C}$. □

**Proposition 9D.2.** *Let* $\mathbf{C}\in\mathbf{PCat}_{\mathbb{K}S}$ *and* $\mathtt{L}\in\mathrm{Si}(\mathbf{C})$*. Then we have*

$$\dim^{\mathbf{C}}(\mathtt{L})\neq 0,\quad {}^{\mathbf{C}}\dim(\mathtt{L})\neq 0.$$

*More generally, if* $\mathrm{f}\colon\mathtt{L}\xrightarrow{\cong}\mathtt{L}^{\star\star}$ *is an isomorphism, then*

$$\mathrm{tr}^{\mathbf{C}}(\mathrm{f})\neq 0,\quad {}^{\mathbf{C}}\mathrm{tr}(\mathrm{f})\neq 0.$$

*Proof.* Note that $\dim^{\mathbf{C}}(\mathtt{L})=0$ would contradict Schur's lemma Lemma 6J.10: from $\dim^{\mathbf{C}}(\mathtt{L})=0$ we get

$$\dim\mathrm{End}_{\mathbf{C}}(\mathtt{L})=\dim\mathrm{Hom}_{\mathbf{C}}(\mathtt{L}\mathtt{L}^{\star},\mathbb{1})>1,$$

since we would get a map different from $\mathrm{ev}_{\mathtt{L}}$. Thus, we are done by symmetry as the argument for the traces is exactly the same. □

**Example 9D.3.** Proposition 9D.2 fails in the non-semisimple case, *i.e.* for $\mathtt{Z}\in\mathrm{In}(\mathbf{C})$

$$\dim^{\mathbf{C}}(\mathtt{Z})=0,\quad {}^{\mathbf{C}}\dim(\mathtt{Z})=0,$$

is possible. To give an explicit example, let us consider the Rumer–Teller–Weyl category $\mathbf{TL}^{q}_{\mathbb{C}\Subset}$ as in Definition 7F.1, and let $q=\exp(2\pi i/4)=i\in\mathbb{C}$. Then the circle removal becomes

(9D-4) 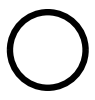 $=0=\dim^{\mathbf{TL}^{q}_{\mathbb{C}\Subset}}(\bullet)$,

and $\bullet$ is, of course, indecomposable. Using Equation 9D-4, we also get an isomorphism of $\mathbb{C}$ algebras

$\mathrm{End}_{\mathbf{TL}^{q}_{\mathbb{C}\Subset}}(\bullet^{2})\xrightarrow{\cong}\mathbb{C}[X]/(X^{2})$, 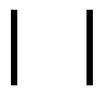 $\mapsto 1$, 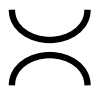 $\mapsto X$.

This implies that $\mathrm{End}_{\mathbf{TL}^{q}_{\mathbb{C}\Subset}}(\bullet^{2})$ is a local $\mathbb{C}$ algebra, and thus, $\bullet^{2}\in\mathrm{In}(\mathbf{TL}^{q}_{\mathbb{C}\Subset})$. We also have

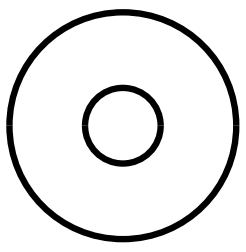 $=0=\dim^{\mathbf{TL}^{q}_{\mathbb{C}\Subset}}(\bullet^{2})$.

More general, one can show that

$$\big(\mathtt{Z}\in\mathrm{In}(\mathbf{TL}^{q}_{\mathbb{C}\Subset})\big)\Rightarrow\big(\dim^{\mathbf{TL}^{q}_{\mathbb{C}\Subset}}(\mathtt{Z})=0\text{ unless }\mathbb{1}=\mathtt{Z}\big).$$

We can therefore say that $\mathbf{TL}^{q}_{\mathbb{C}\Subset}$ is a very degenerate category. ◇

That dimensions of objects are non-zero is, in some sense, a property of semisimple categories:

**Proposition 9D.5.** *Let* $\mathbf{C}\in\mathbf{mlFiat}$*. Then the following are equivalent:*

*(I)* $\mathbf{C}$ *is semisimple;*

*(II)* $\dim^{\mathbf{C}}(\mathtt{P})\neq 0$ *for all* $\mathtt{P}\in\mathrm{Pi}(\mathbf{C})$*;*

*(III)* $\dim^{\mathbf{C}}(\mathtt{P})\neq 0$ *for some* $\mathtt{P}\in\mathbf{Proj}(\mathbf{C})$*;*

*(IV)* ${}^{\mathbf{C}}\dim(\mathtt{P})\neq 0$ *for all* $\mathtt{P}\in\mathrm{Pi}(\mathbf{C})$*;*

*(V)* ${}^{\mathbf{C}}\dim(\mathtt{P})\neq 0$ *for some* $\mathtt{P}\in\mathbf{Proj}(\mathbf{C})$.

*Proof.* (I)⇒(II). By Proposition 9D.2 since all simple objects (and thus, all objects) are projective.

(II)⇒(III). Evident.

(III)⇒(I). Consider

$$\mathbb{1}\xrightarrow{\mathrm{coev}^{\mathtt{P}}}\mathtt{P}({}^{\star}\mathtt{P})\xrightarrow{\cong}\mathtt{P}\mathtt{P}^{\star}\xrightarrow{\mathrm{ev}_{\mathtt{P}}}\mathbb{1}\ ,$$

which calculates $\dim^{\mathbf{C}}(\mathtt{P})$. If this is not zero, then $\mathbb{1}\Subset\mathtt{P}({}^{\star}\mathtt{P})\in\mathbf{Proj}(\mathbf{C})$. Thus, Theorem 7E.14 implies that $\mathbf{C}$ is semisimple.

Finally, by symmetry, (II) and (III) are equivalent to (IV) and (V), so we are done. □

**Lemma 9D.6.** *Let* $\mathbf{C}\in\mathbf{mFiat}$ *and* $\mathtt{X}\in\mathbf{C}$*. Then we have*

$$\dim^{\mathbf{C}}(\mathtt{X})=\textstyle\sum_{i=1}^{n}(\mathtt{X}:\mathtt{Z}_i)\dim^{\mathbf{C}}(\mathtt{L}_i),\quad {}^{\mathbf{C}}\dim(\mathtt{X})=\textstyle\sum_{i=1}^{n}(\mathtt{X}:\mathtt{Z}_i){}^{\mathbf{C}}\dim(\mathtt{L}_i).$$

*Proof.* An easy calculation using Equation 9C-2 and sliding:

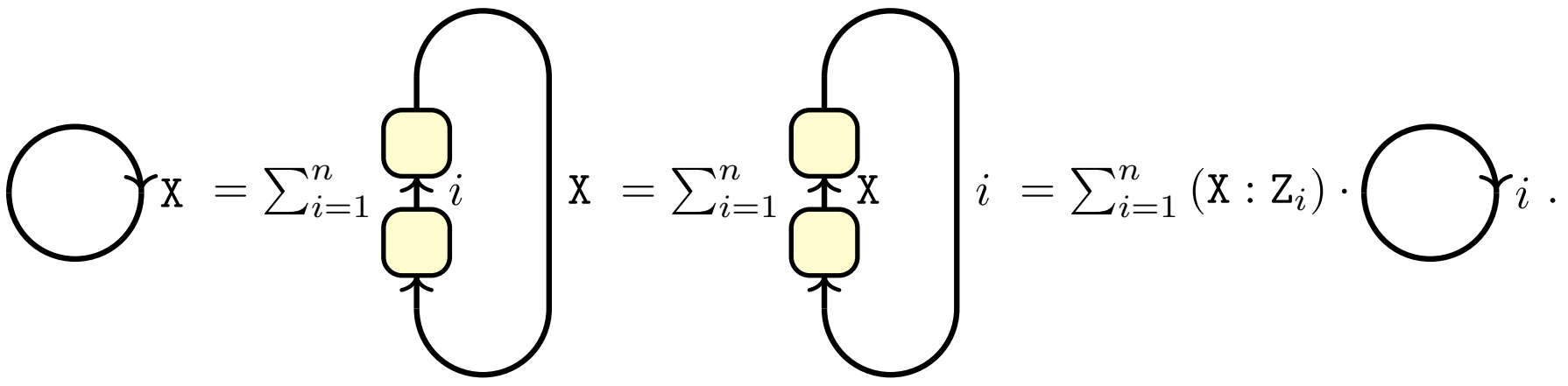

The other case follows by symmetry. □

**Definition 9D.7.** Let $\mathbf{C} \in \mathbf{mFiat}$. Then the ***categorical dimension of*** $\mathbf{C}$ is

$$\dim(\mathbf{C}) = \sum_{i=1}^{n} {}^{\mathbf{C}}\dim(\mathtt{Z}_i)\dim^{\mathbf{C}}(\mathtt{Z}_i).$$

This depends on the choice of a pivotal structure. ◇

Note that, if $\mathbf{C}$ is spherical, then

$$\dim(\mathbf{C}) = \sum_{i=1}^{n} \dim^{\mathbf{C}}(\mathtt{Z}_i)^2. \tag{9D-8}$$

**Example 9D.9.** The categorical dimension generalizes familiar concepts:

(i) For $\mathbf{Vec}_{\Bbbk}$ we get $\dim(\mathbf{Vec}_{\Bbbk}) = 1 = \mathrm{PFdim}(\mathbf{Vec}_{\Bbbk})$.

(ii) For $\mathbf{Vec}^{\omega}_{\mathbb{C}\oplus}(\mathbb{Z}/2\mathbb{Z})$ we have already seen that there were two choices of pivotal structures giving

(9D-10) 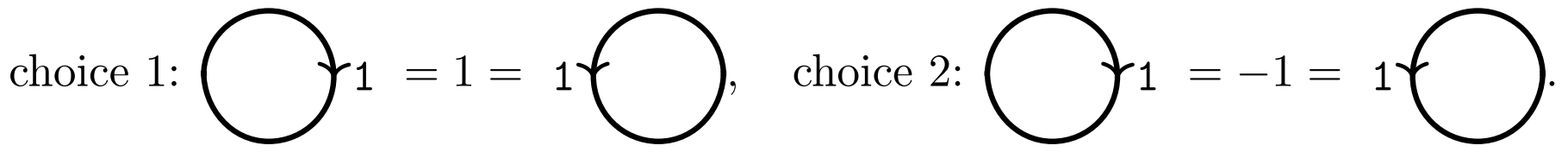

The same work for a trivial $\omega$, *i.e.* there are also two choices of (co)evaluations satisfying Equation 9D-10. Thus, for both choices, we get

$$\mathrm{Dim}\big(\mathbf{Vec}_{\mathbb{C}\oplus}(\mathbb{Z}/2\mathbb{Z})\big) = \mathrm{Dim}\big(\mathbf{Vec}^{\omega}_{\mathbb{C}\oplus}(\mathbb{Z}/2\mathbb{Z})\big) = 2 = \#(\mathbb{Z}/2\mathbb{Z})$$
$$= \mathrm{PFdim}\big(\mathbf{Vec}_{\mathbb{C}\oplus}(\mathbb{Z}/2\mathbb{Z})\big) = \mathrm{PFdim}\big(\mathbf{Vec}^{\omega}_{\mathbb{C}\oplus}(\mathbb{Z}/2\mathbb{Z})\big).$$

(iii) For the non-spherical category such as $\mathbf{Vec}_{\mathbb{C}\oplus}(\mathbb{Z}/3\mathbb{Z})$ with the choice of (co)evaluations from Example 4H.7 we get

$$\mathrm{Dim}\big(\mathbf{Vec}_{\Bbbk\oplus}(\mathbb{Z}/3\mathbb{Z})\big) = 3 = \#(\mathbb{Z}/3\mathbb{Z}) = \mathrm{PFdim}\big(\mathbf{Vec}_{\Bbbk\oplus}(\mathbb{Z}/3\mathbb{Z})\big).$$

In all the above examples, the equality between the categorical dimensions and the PF dimension works without adjustment since $\mathbb{R}_{\geq 0} \subset \mathbb{C}$, the latter being the ground ring for the categorical dimension and the former for the PF dimension. ◇

**Proposition 9D.11.** *Let* $\mathbf{C} \in \mathbf{mFus}$ *be spherical and* $\Bbbk$ *be of characteristic zero. Then* $\dim(\mathbf{C}) \neq 0$. *Moreover, if* $\Bbbk = \mathbb{C}$, *then* $\dim(\mathbf{C}) \geq 1$.

*Proof.* Since $\mathbf{C}$ is spherical, we have Equation 9D-8, which immediately implies the claims. □

**Example 9D.12.** Proposition 9D.11 does not hold in finite characteristic. For example, with reference to Example 9D.9.(b), we have

$$\mathrm{Dim}\big(\mathbf{Vec}_{\mathbb{F}_2\oplus}(\mathbb{Z}/2\mathbb{Z})\big) = 2 = 0 \in \mathrm{End}_{\mathbf{Vec}_{\mathbb{F}_2\oplus}(\mathbb{Z}/2\mathbb{Z})}(\mathbb{1}) \cong \mathbb{F}_2.$$

Note however that

$$\mathrm{PFdim}\big(\mathbf{Vec}_{\mathbb{F}_2\oplus}(\mathbb{Z}/2\mathbb{Z})\big) = 2 \neq 0,$$

since the PF dimension is, by definition, an element in $\mathbb{R}_{\geq 0}$. ◇

Note that Example 9D.12 also illustrates a crucial difference between the categorical dimension and the PF dimension: The first is a categorical notion, meaning it is about morphisms, and lives in $\mathbf{C}$. On the other hand, the PF dimension is a numerical notion, relates to objects and lives in $\mathbb{R}_{\geq 0}$.

9E. **The Alexander–Markov theorem and traces of braids.** Recall that we had the category of braids $\mathbf{qSym}$, see Example 5C.8, which is the free braided category generated by one object. In particular, we get the ***Alexander functor***

(9E-1) $\mathrm{A}\colon \mathbf{qSym} \to \mathbf{oqBr}$, $\bullet \mapsto \bullet$, 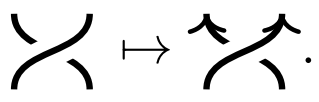

The ***Markov quotient*** of $\mathbf{qSym}$, denoted by $\mathbf{qSym}/\mathbf{MM}$, is the quotient of $\mathbf{qSym}$ by the congruence spanned by ***Markov moves (MM)***. Formally:

**Definition 9E.2.** We let **qSym/MM** $= \langle \mathbf{S}, \mathbf{T} \mid \mathbf{R} \cup \mathbf{MM} \rangle$ with

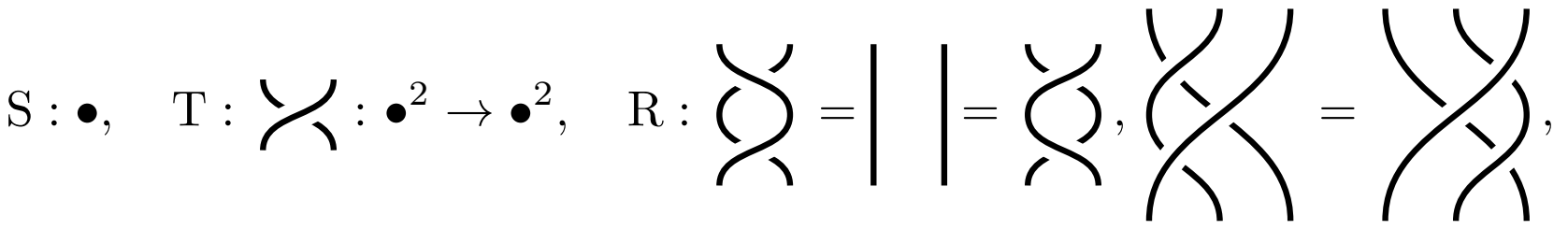

(9E-3)

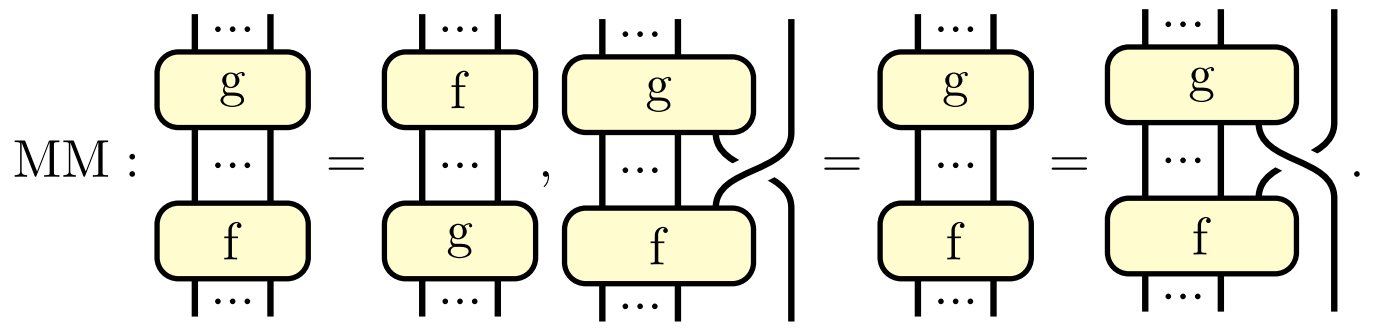

(The Markov moves are imposed for all possible number of strands and all morphisms f and g.) ◇

By definition, we have a full quotient functor

$$\mathrm{M}\colon \mathbf{qSym} \to \mathbf{qSym/MM}, \quad \bullet \mapsto \bullet,$$ 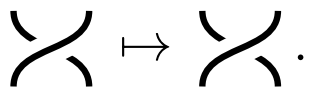

The classical ***Alexander–Markov theorem*** takes now the following form:

**Theorem 9E.4.** *We have the following.*

*(i) The functor* A *from* Equation 9E-1 *is fully faithful.*

*(ii) The functor* A *gives a surjection*

(9E-5) $$\mathrm{A}\colon \textstyle\bigcup_{n\in\mathbb{Z}_{\geq 0}} \mathrm{End}_{\mathbf{qSym}}(\bullet^{n}) \twoheadrightarrow \mathrm{End}_{\mathbf{oqBr}}(\mathbb{1}), \quad \mathrm{f} \mapsto \mathrm{tr}^{\mathbf{oqBr}}\big(\mathrm{A}(\mathrm{f})\big).$$

*(iii) There exists a bijection*

$$\mathrm{AM}^{-1}\colon \textstyle\bigcup_{n\in\mathbb{Z}_{\geq 0}} \mathrm{End}_{\mathbf{qSym/MM}}(\bullet^{n}) \xrightarrow{\cong} \mathrm{End}_{\mathbf{oqBr}}(\mathbb{1}),$$

*making the following diagram (in* **Set***) commutative:*

$$\begin{array}{ccc} \bigcup_{n\in\mathbb{Z}_{\geq 0}} \mathrm{End}_{\mathbf{qSym}}(\bullet^{n}) & \xrightarrow{\quad\mathrm{A}\quad} & \mathrm{End}_{\mathbf{oqBr}}(\mathbb{1}) \\ & \searrow^{\mathrm{M}} \quad\quad \nearrow^{\mathrm{AM}^{-1}}_{\cong} & \\ & \bigcup_{n\in\mathbb{Z}_{\geq 0}} \mathrm{End}_{\mathbf{qSym/MM}}(\bullet^{n}) & \end{array}.$$

*(ii) and (iii) work also for left instead of right traces.*

In words, Theorem 9E.4.(ii) says that every link arises as a closure of a braid, while Theorem 9E.4.(iii) gives a precise condition for when two closures represent the same link.

*Proof.* A proof can be found in *e.g.* [**KT08**, Chapter 2]; we only give a sketch of a proof.

The ***Alexander closure*** A from Equation 9E-5 can be illustrated by "closing a braid to the right":

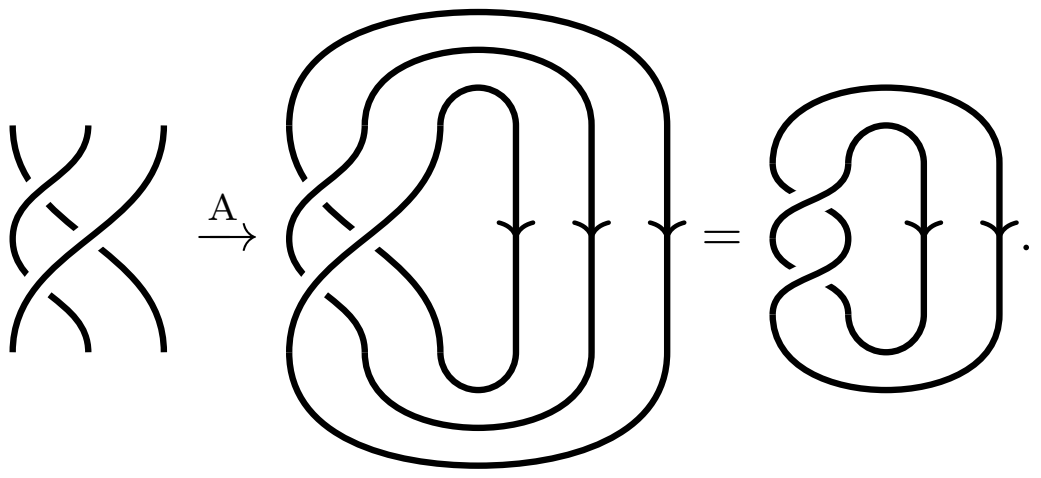

The surprising result of Alexander is then that every link can be arranged such that it has a purely upward-oriented and a purely downward-oriented part, with the latter being trivial. We now sketch a proof (that even gives an algorithm) from [**Vog90**].

The first observation is that Alexander closures can be identified with links on an annulus, *e.g.*:

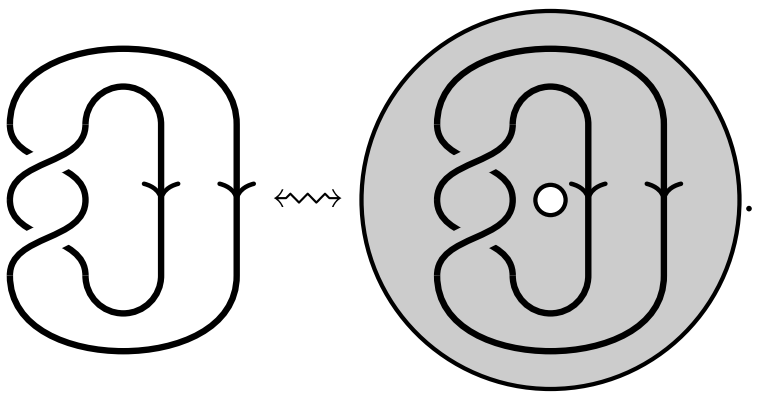

Note the cyclic orientation of such a picture. Now, given the link choose an orientation (any choice works) and get rid of all crossings by

(9E-6)

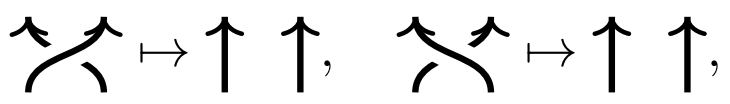

which generates a bunch of circles, but we remember where the crossings were. Now arrange the crossings in cyclic fashion (to show that this always works is the tricky part) a reinsert the crossings. Figure 24 provides an example how to do this on the so-called ***figure eight knot***: only the top (green) bit need to be moved to the bottom (yellow) using Lemma 5I.4.

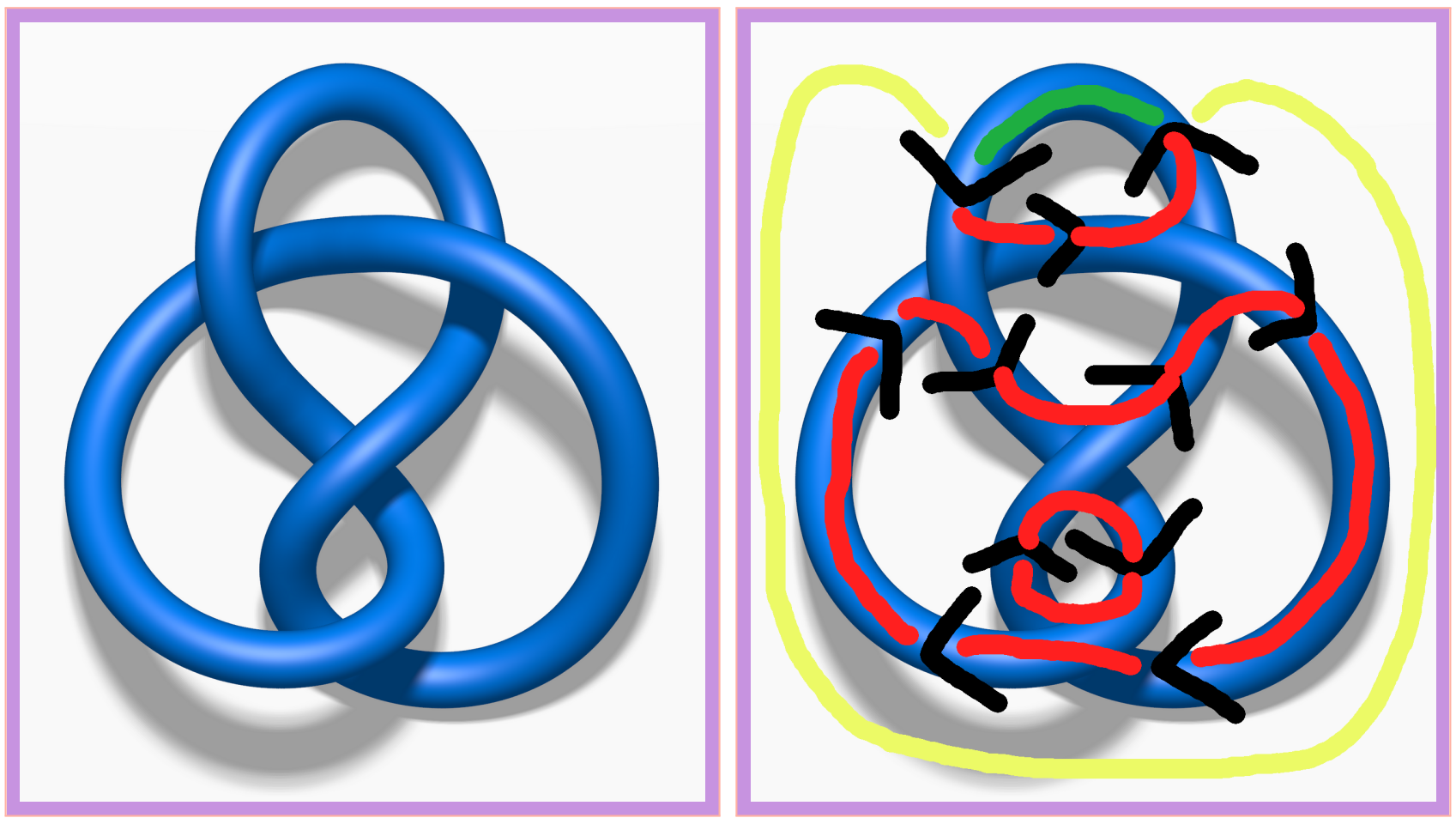

FIGURE 24. The figure eight knot (this picture needs to be rotated 90 degrees when compared to the discussion in this section). We add an orientation and perform Equation 9E-6. The circles are then reorganized using Lemma 5I.4.

Picture from https://en.wikipedia.org/wiki/Figure-eight_knot_(mathematics)

After rearranging the circles, displayed at the left, one then gets

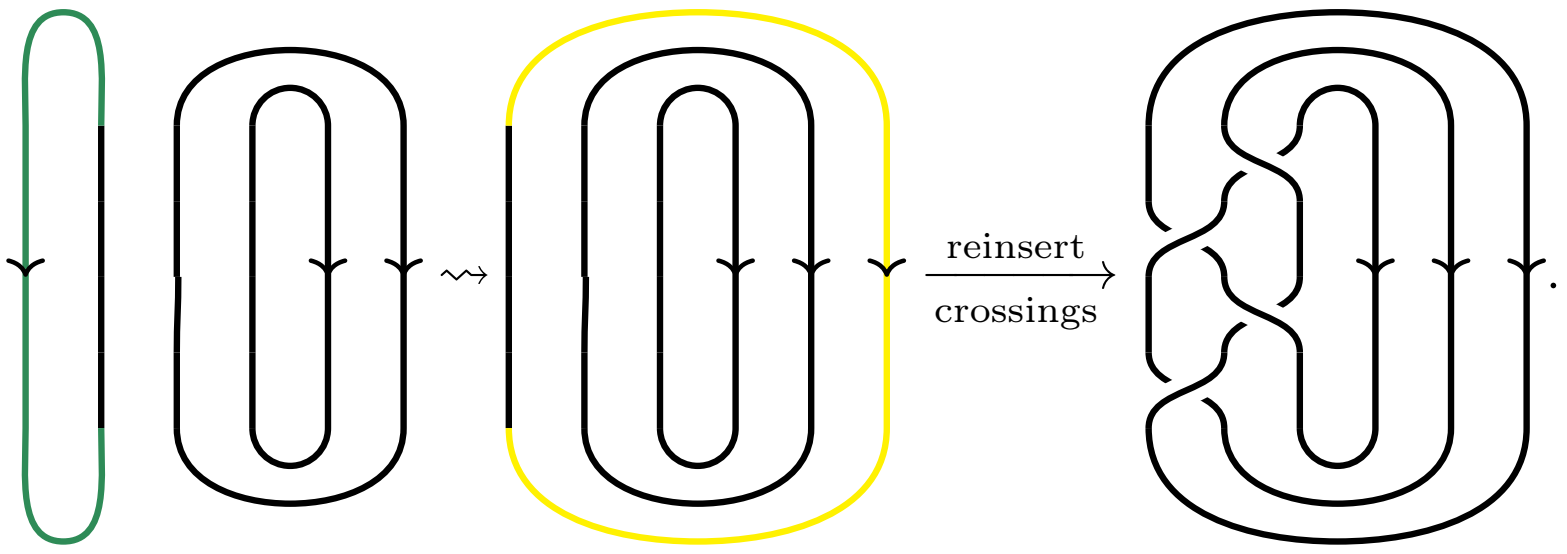

This is a braid presentation for the figure eight knot.

The Markov moves then just take the form of sliding and Reidemeister 1:

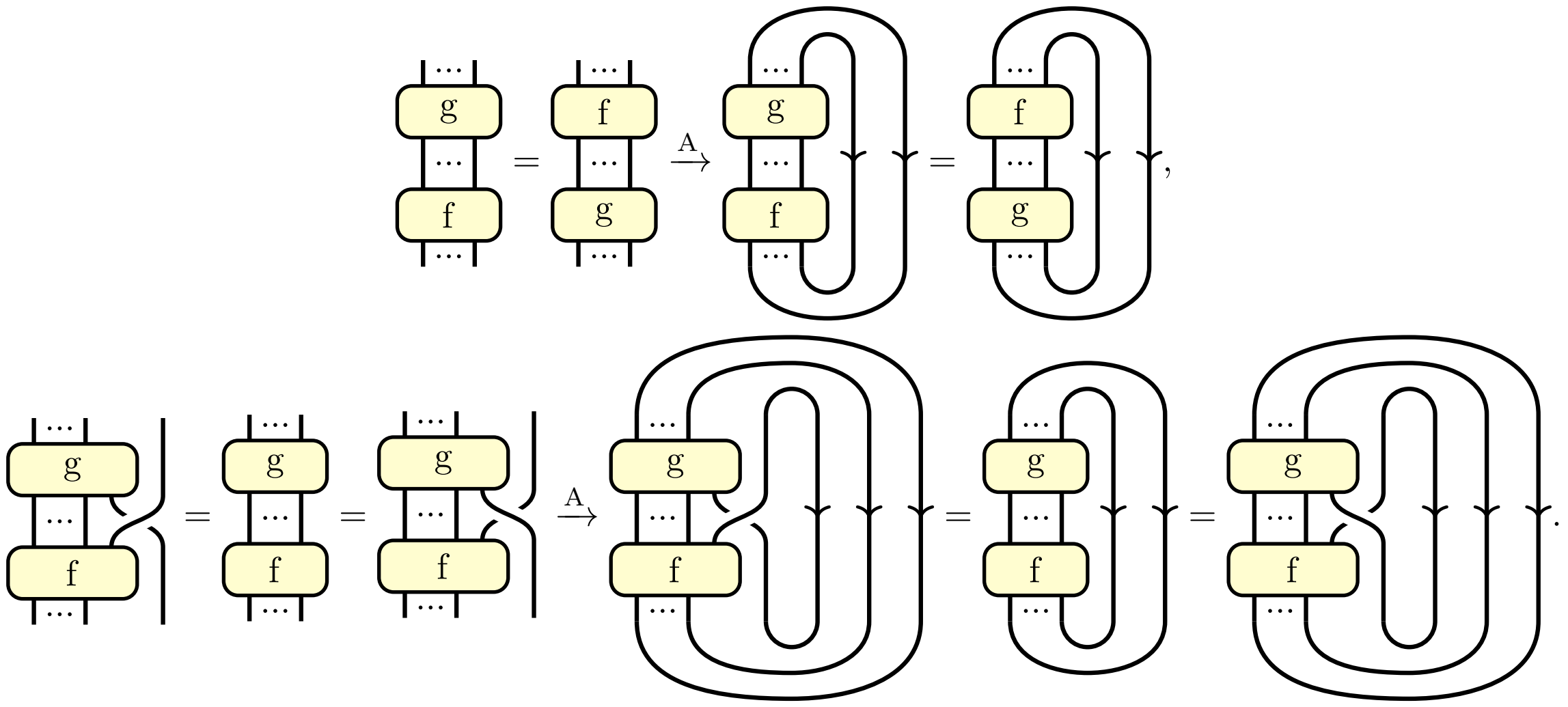

These evidently hold in **oqBr**. (See also Exercise 9M.1.) The point of Theorem 9E.4.(iii) is that these are the only extra relations. □

*Remark* 9E.7. The functor in Equation 9E-1 has, of course, various cousins, *e.g.* we could equally well go to **qBr**. Similarly for the Alexander–Markov theorem Theorem 9E.4, which exists in a variety of flavors. ◇

9F. **Colored braids and links.** The previous section is partially a motivation for the following. More general than Section 9E, if $\mathbf{C} \in \mathbf{BCat}$ is any braided category and $\mathtt{X} \in \mathbf{C}$ is any object, then we get a ***coloring with*** $\mathtt{X}$ ***functor***

$$\mathrm{A}_{\mathtt{X}} \colon \mathbf{qSym} \to \mathbf{C}, \quad \bullet \mapsto \mathtt{X},$$

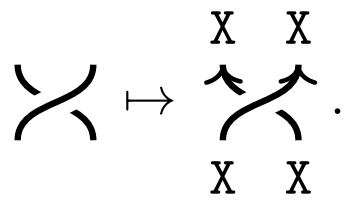

Similarly, if $\mathbf{C} \in \mathbf{BPCat}$ is any braided pivotal category and $\mathtt{X} \in \mathbf{C}$ is any object, then we have a more general ***coloring with*** $\mathtt{X}$ ***functor***

$$\mathrm{A}^{\star}_{\mathtt{X}} \colon \mathbf{oqBr} \to \mathbf{C}, \ \bullet \mapsto \mathtt{X},$$

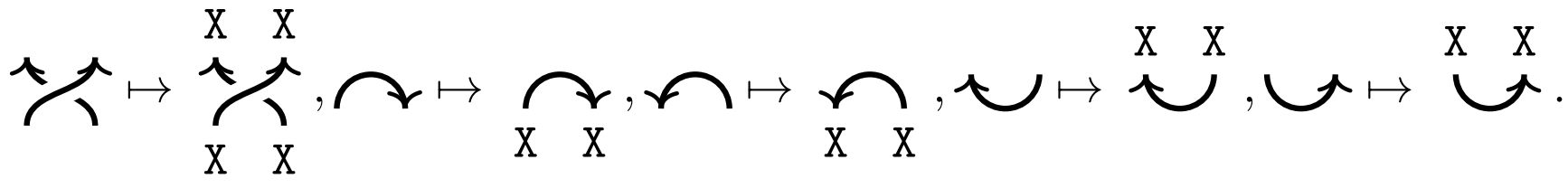

*Remark* 9F.1. The coloring functors as above have the "flaw" that one can only color with one color at a time. This can be corrected by considering the ***category of colored braids*** **cqSym** or the ***colored oriented quantum Brauer category*** **coqBr**. The images of the coloring functors are then called ***colored braids*** or ***colored (oriented) tangles***, respectively. ◇

Our main target for coloring functors (which thus, allow a diagrammatic calculus of colored tangles) are:

**Definition 9F.2.** Let $\mathbf{C} \in \mathbf{mlFiat}$. Then:

- If $\mathbf{C} \in \mathbf{BCat}$, then we call it a ***multi locally bfiat category***. The corresponding category is denoted by **mlBfiat**.
- If $\mathbf{C} \in \mathbf{BCat}_S$, then we call it a ***multi locally bmodular category***. The corresponding category is denoted by **mlBMo**.
- Finally, if these satisfy the ribbon equation Equation 5H-3, then we call them ***rfiat*** or ***rmodular***, with the corresponding notation for the categories.

◇

*Remark* 9F.3. An rmodular category is also called a ***pre-modular category*** in the literature. Note also the hierarchy in the above definitions. ◇

**Example 9F.4.** We have already seen plenty of examples:

(a) Of course, $\mathbf{Vec}_{\mathbb{k}}$ is bmodular.

(b) More generally, $\mathbf{Vec}^{\omega}_{\mathbb{k}\oplus}(\mathrm{G})$ is bmodular if and only if G is a finite abelian group.

(c) For G being a finite group, the category $\mathbf{fdMod}\big(\mathbb{C}[\mathrm{G}]\big)$ is bmodular.

(d) More generally, let $\mathbb{K}$ be algebraically closed and assume that the condition in Proposition 8B.6.(iii) holds. Then $\mathbf{fdMod}\big(\mathbb{K}[\mathrm{G}]\big)$ is bfiat.

These are good examples to keep in mind, but they are all fairly boring as targets for quantum invariants. ◇

**Example 9F.5.** Whether a bfiat or bmodular is ribbon is trickier, as this depends on choices. Let us discuss $\mathbf{Vec}_{\mathbb{C}\oplus}(\mathbb{Z}/3\mathbb{Z})$ in detail, where we recall Example 4H.7 and Lemma 5E.1. In particular, there are several choices of (co)evaluations and braidings given as follows. (Some are equivalent, but let us just list them anyways.) Let

$$d_j(\mathtt{i}) = \zeta^{ij}, \ i,j \in \{0,1,2\},$$

where $\zeta = \exp(2\pi i/3) \in \mathbb{C}$. Then, for $k,l,m \in \{0,1,2\}$,

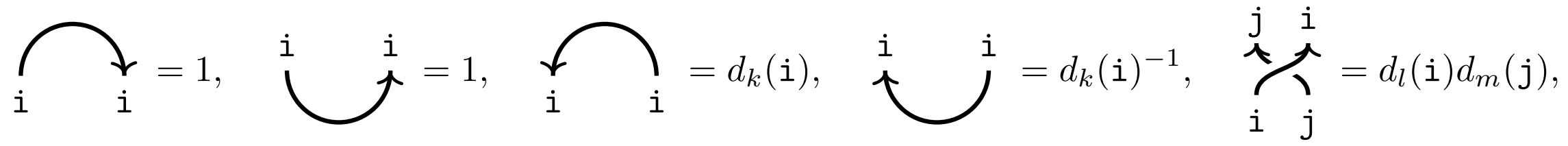

are choices. We then check that

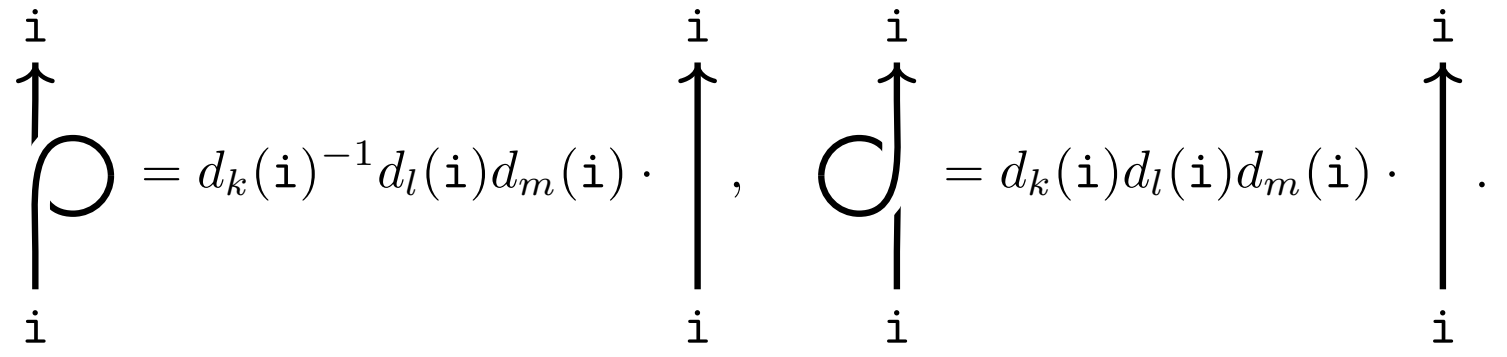

In particular, being ribbon does not depend on the choice of braiding, but is equivalent to $\mathbf{Vec}_{\mathbb{C}\oplus}(\mathbb{Z}/3\mathbb{Z})$ being spherical. ◇

Recall that for a $\mathbf{C} \in \mathbf{BFiat}$ being $\mathbb{S}$ linear we have a finite set of indecomposables $\mathrm{In}(\mathbf{C}) = \{\mathtt{Z}_1, ..., \mathtt{Z}_n\}$ and also $\mathrm{End}_{\mathbf{C}}(\mathbb{1}) \cong \mathbb{S}$, and we can consider the ***colored Hopf braid***

$\mathrm{s}_{ij} = \beta_{\mathtt{Z}_j,\mathtt{Z}_i}\beta_{\mathtt{Z}_i,\mathtt{Z}_j} =$ 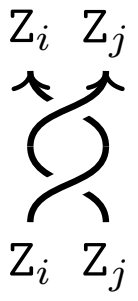

$, \quad \mathrm{tr}^{\mathbf{C}}(\mathrm{s}_{ij}) =$ 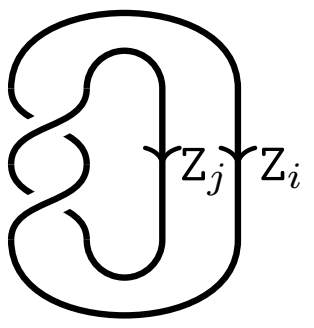

$\in \mathbb{S}.$

These assemble into an important $n \times n$ matrix, called the $S$ ***matrix***

$$S = \big(\mathrm{tr}^{\mathbf{C}}(\mathrm{s}_{ij})\big)_{i,j=1}^{n} \in \mathrm{Mat}_{n\times n}(\mathbb{S}).$$

**Example 9F.6.** If $\mathbf{C}$ has two indecomposable objects, say spinach and orchid (illustrated dashed for readability), then the $S$ matrix is:

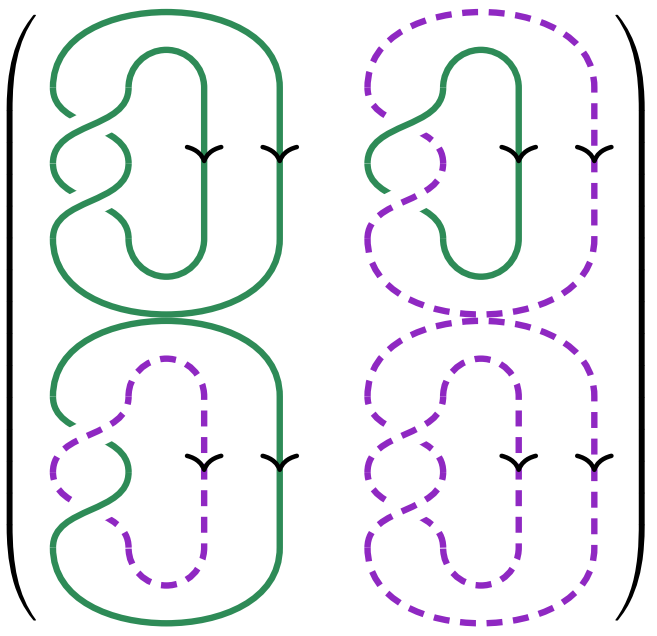
,

which is where the name "coloring" comes from. ◇

Note that, if the braiding is symmetric, then

$\mathrm{tr}^{\mathbf{C}}(\mathrm{s}_{ij}) =$ 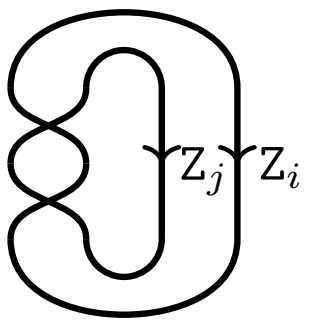

$=$ 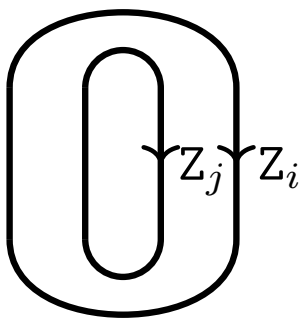

$= \dim^{\mathbf{C}}(\mathtt{Z}_i)\dim^{\mathbf{C}}(\mathtt{Z}_j) \in \mathbb{S}.$

**Example 9F.7.** Let us discuss the case of $\mathbf{Vec}_{\Bbbk\oplus}(\mathrm{G})$ for small G.

(a) Recall from Example 6F.7 that $\mathbf{Vec}_{\mathbb{C}\oplus}(\mathbb{Z}/2\mathbb{Z})$ has two braidings. However, both satisfy $\mathrm{s}_{ij} = \mathrm{id}_{\mathtt{L}_i\mathtt{L}_j}$. So they give the same result for the $S$ matrix. Moreover, recall that $\mathbf{Vec}_{\mathbb{C}\oplus}(\mathbb{Z}/2\mathbb{Z})$ has two choices (co)evaluations, see Example 4H.12. For these we get

$$\text{choice 1: } S = \begin{pmatrix} 1 & 1 \\ 1 & 1 \end{pmatrix}, \quad \det(S) = 0, \quad \text{choice 2: } S = \begin{pmatrix} 1 & -1 \\ -1 & 1 \end{pmatrix}, \quad \det(S) = 0.$$

(b) For $\mathbf{Vec}_{\mathbb{C}\oplus}(\mathbb{Z}/3\mathbb{Z})$, enumerate $\mathrm{Si}\big(\mathbf{Vec}_{\mathbb{C}\oplus}(\mathbb{Z}/3\mathbb{Z})\big) = \{0, 1, 2\}$, let $\zeta = \exp(2\pi i/3)$ and take a braiding such that $s_{ij} = \zeta^{i+j}$. Then, for the standard rigidity structure,

$$S = \begin{pmatrix} 1 & \zeta & \zeta^2 \\ \zeta & \zeta^2 & 1 \\ \zeta^2 & 1 & \zeta \end{pmatrix}, \quad \det(S) = 0.$$

That the $S$ matrices are degenerate is not a coincidence, as we will see. ◇

**Example 9F.8.** Return to $\mathrm{S}_3$, *cf.* Example 8C.6: The category $\mathbf{fdMod}\big(\mathbb{C}[\mathrm{S}_3]\big)$ with the swap map as the braiding and the usual (co)evaluations is rfiat. With this choice the categorical dimension is just the dimension as a $\mathbb{C}$ vector space. Thus, since the braiding is symmetric, we get

$$S = \begin{pmatrix} 1 & 2 & 1 \\ 2 & 4 & 2 \\ 1 & 2 & 1 \end{pmatrix}, \quad \det(S) = 0.$$

As in Example 9F.7, the $S$ matrix is degenerate. ◇

Below we will use the color code from Example 9F.6 to make the equations easier to read.

**Lemma 9F.9.** *For* $\mathbf{C} \in \mathbf{BFiat}$ *and any* $\mathtt{Z}_i \in \mathrm{In}(\mathbf{C})$ *we have*

$$s_{ij}s_{ik} = \dim^{\mathbf{C}}(\mathtt{Z}_i) \cdot \textstyle\sum_{l=1}^{n} N_{jk}^{l} s_{il}.$$

*Proof.* The diagrammatic equation

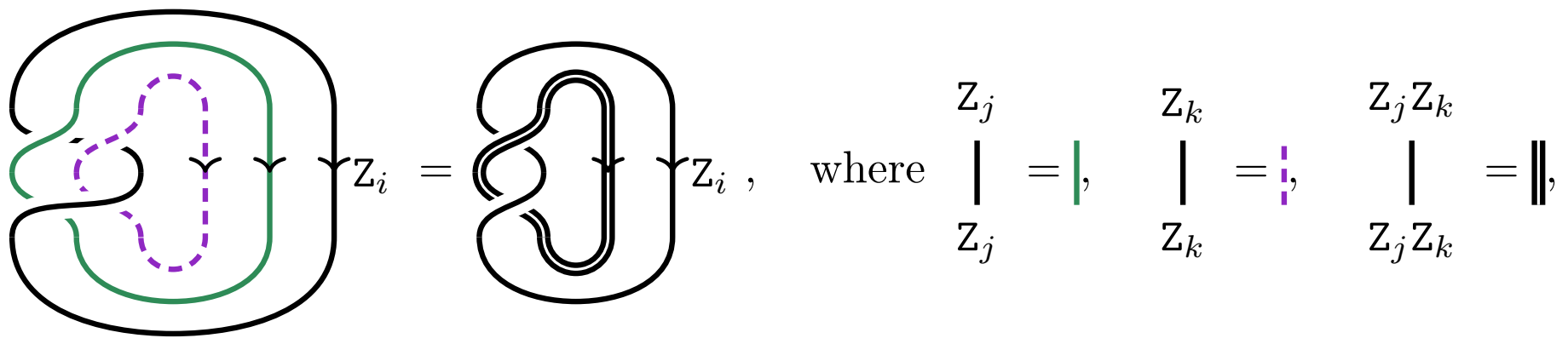

and additivity provide the result. □

**Example 9F.10.** In Example 9F.8 we can easily check that for $i = 2 = j$ and $k = 1$ we have $8 = 2(0 \cdot 2 + 1 \cdot 4 + 0 \cdot 2)$. ◇

**Lemma 9F.11.** *Any* $\mathbf{C} \in \mathbf{BFiat}$ *has a symmetric* $S$ *matrix.*

*Proof.* A Reidemeister-type argument:

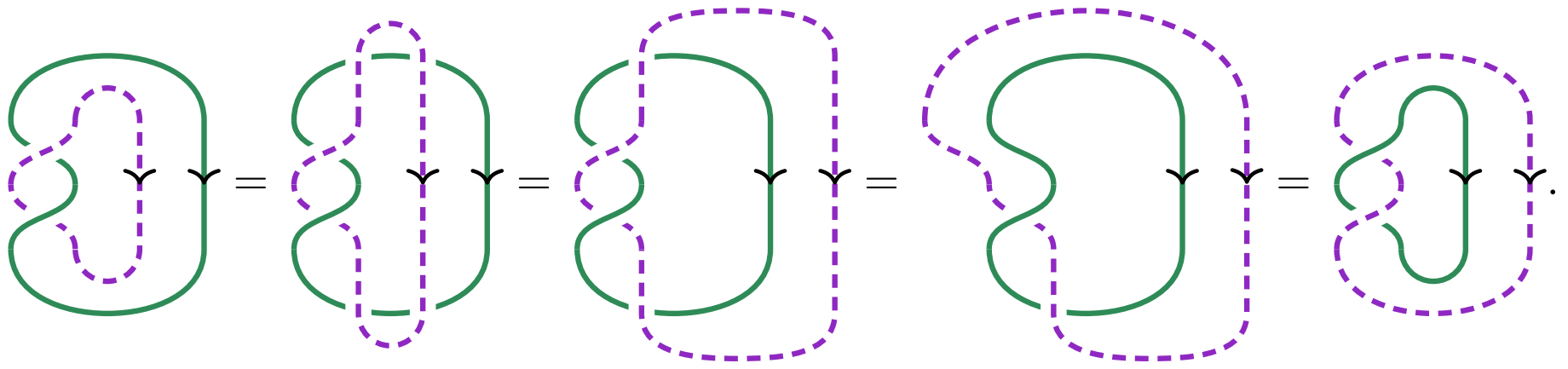

This shows that the $S$ matrix is its own transpose. □

Let us continue this section with an important lemma which should remind us that the Reidemeister 1 move Equation 5H-2 "is not as innocent as it looks".

**Lemma 9F.12.** *Let* $\mathbf{C} \in \mathbf{mlBFiat}$ *be* $\mathbb{K}$ *linear and* $\mathtt{Z} \in \mathrm{In}(\mathbf{C})$ *and* $\mathbb{K}$ *algebraically closed. There exits* $a(\mathtt{Z}) \in \mathbb{K}^*$ *such that*

(9F-13)

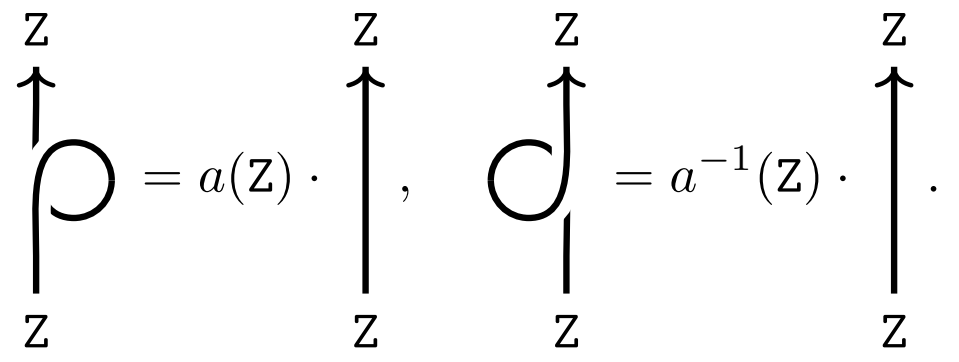

*Proof.* After recalling that the twist is invertible with explicit inverse as given in Lemma 5H.6, this is a direct consequence of Schur's lemma Lemma 6J.10 if $\mathtt{Z}$ is simple. For general $\mathtt{Z}$ we also use Schur's lemma Lemma 6J.9 and additionally observe that the invertible elements in a local ring can be identified with the ground field. □

For $\mathbf{C} \in \mathbf{lBFiat}$ we let

(9F-14)

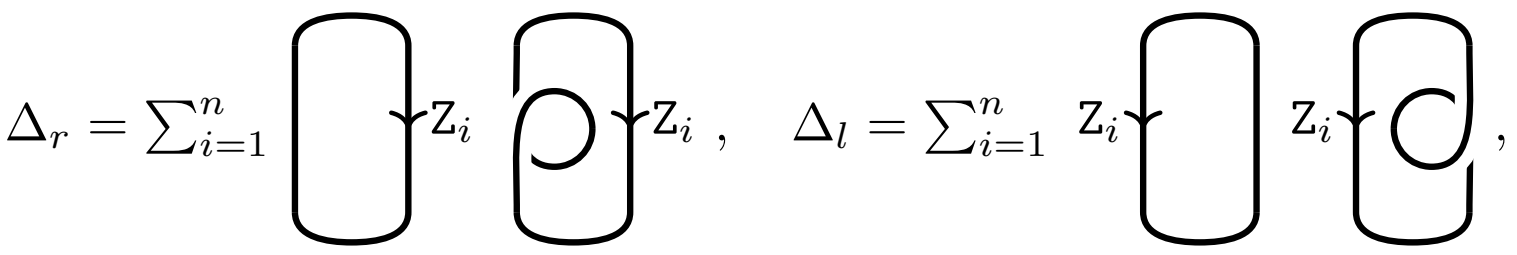

both of which are in $\mathrm{End}_{\mathbf{C}}(\mathbb{1}) \cong \mathbb{S}$.

**Lemma 9F.15.** *Let* $\mathbf{C} \in \mathbf{lBFiat}$ *be* $\mathbb{K}$ *linear. Then:*

*(i) We have* $\Delta_r = \sum_{i=1}^{n} a(\mathtt{Z}_i)\mathrm{dim}^{\mathbf{C}}(\mathtt{Z}_i)^2$ *and* $\Delta_l = \sum_{i=1}^{n} a(\mathtt{Z}_i)^{-1}{}^{\mathbf{C}}\mathrm{dim}(\mathtt{Z}_i)^2$.

*(ii) If* $\mathbf{C}$ *is semisimple, then* $\Delta_r = \sum_{i=1}^{n} a(\mathtt{Z}_i)^{\mathbf{C}}\mathrm{dim}(\mathtt{Z}_i)^2$ *and* $\Delta_l = \sum_{i=1}^{n} a(\mathtt{Z}_i)^{-1}\mathrm{dim}^{\mathbf{C}}(\mathtt{Z}_i)^2$.

*Proof.* An immediate consequence of Lemma 9F.12. □

9G. **Modular categories.** The $S$ matrix is symmetric, but *e.g.* Example 9F.7 shows that it might not be invertible. So:

**Definition 9G.1.** A category $\mathbf{C} \in \mathbf{BFiat}$ with invertible $S$ matrix is called ***mfiat***. If such a $\mathbf{C}$ is additionally semisimple, then $\mathbf{C}$ is called ***modular***. ◇

The corresponding categories are denoted by **MoFiat** and **MoCat**.

**Example 9G.2.** Back to Example 9F.4:

(a) The category $\mathbf{Vec}_{\mathbb{k}}$ is modular.

(b) However, $\mathbf{Vec}^{\omega}_{\mathbb{k}\oplus}(\mathrm{G})$ is rarely modular, *cf.* Example 9F.7. (See also [EGNO15, Example 8.13.5].)

(c) For G being a finite group the category $\mathbf{fdMod}(\mathbb{C}[\mathtt{G}])$ is only modular if G is the trivial group.

(d) More generally, a symmetric category with indecomposable monoidal unit $\mathbb{1}$ cannot be modular if it has more than one indecomposable object.

We will prove the final point as a special case of Proposition 9G.6. ◇

**Example 9G.3.** The Verlinde categories, *cf.* Section 8H, are all modular. This is not trivial, see, for example, [**BK01a**, Section 3.3]. We sketch how this can be proven for the Verlinde category associated to $\mathfrak{sl}_2$ in Section 9H. ◇

Let $\mathbf{C} \in \mathbf{BFiat}$ be $\Bbbk$ linear. Then $\mathtt{X} \in \mathbf{C}$ is ***transparent*** if

(9G-4) $$\beta_{\mathtt{Z},\mathtt{X}}\beta_{\mathtt{X},\mathtt{Z}} = \mathrm{id}_{\mathtt{XZ}}, \quad \beta_{\mathtt{X},\mathtt{Z}}\beta_{\mathtt{Z},\mathtt{X}} = \mathrm{id}_{\mathtt{ZX}}, \quad \text{for all } \mathtt{Z} \in \mathrm{In}(\mathbf{C}).$$

We immediately get:

**Lemma 9G.5.** *Let $\mathbf{C} \in \mathbf{BFiat}$ be $\Bbbk$ linear. Then the monoidal unit $\mathbb{1}$ is transparent.* □

The topological motivation for Definition 9G.1 is the following.

**Proposition 9G.6.** *Let $\mathbf{C} \in \mathbf{BFiat}$ be $\Bbbk$ linear. Then $\mathbf{C} \in \mathbf{MoFiat}$ if and only if $\mathbb{1}$ is the only indecomposable object that is transparent.*

*Proof.* We sketch an argument. Note that Equation 9G-4 for the indecomposable object implies that

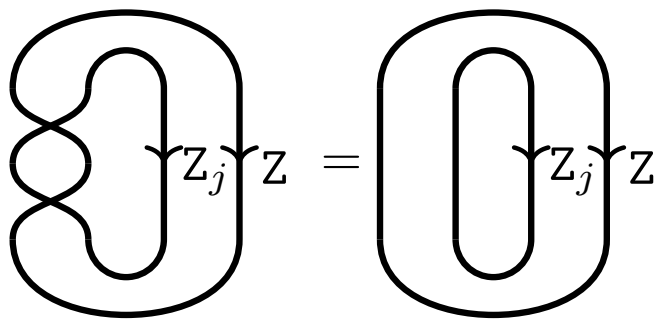

$$= \dim^{\mathbf{C}}(\mathtt{Z})\dim^{\mathbf{C}}(\mathtt{Z}_j).$$

Thus, the column of the $S$ matrix associated to $\mathtt{Z}$ is the (transpose of) $\dim^{\mathbf{C}}(\mathtt{Z})(\dim^{\mathbf{C}}(\mathtt{Z}_1), ..., \dim^{\mathbf{C}}(\mathtt{Z}_n))$. This, after normalization, is $(\dim^{\mathbf{C}}(\mathtt{Z}_1), ..., \dim^{\mathbf{C}}(\mathtt{Z}_n))$. The same is true for any other transparent object, so the $S$ matrix is degenerate as soon as we have two transparent objects. Completing this proof is Exercise 9M.4. □

9H. **More on Rumer–Teller–Weyl categories.** Let us come back to $\mathbf{TL}^q_{\mathbb{S}\oplus\Subset}$ as in Definition 7F.1. This category usually has infinitely many simple and indecomposable objects, but it has nice quotients. Here we use the usual quantum numbers, *i.e.* for $a \in \mathbb{Z}_{\geq 0}$ we let $[0]_q = 0$, $[1]_q = 1$ and

(9H-1) $$[a]_q = q^{a-1} + q^{a-3} + ... + q^{-a+3} + q^{-a+1} \in \mathbb{S}.$$

For $a \in \mathbb{Z}_{<0}$ we let $[a]_q = -[-a]_q$.

**Example 9H.2.** Note that $[a]_q$ depends on the choice of $q$. To be explicit let $\mathbb{S} = \mathbb{C}$ and let $q$ be either 1, $\zeta_2 = i = \exp(2\pi i/4)$, $\zeta_3 = \exp(2\pi i/3)$, $\zeta_4 = \exp(2\pi i/8)$, or $\zeta_5 = \exp(2\pi i/5)$. Then

| | $[1]_q$ | $[2]_q$ | $[3]_q$ | $[4]_q$ | $[5]_q$ | $[6]_q$ | $[7]_q$ | $[8]_q$ |
|---|---|---|---|---|---|---|---|---|
| $q = 1$ | 1 | 2 | 3 | 4 | 5 | 6 | 7 | 8 |
| $q = \zeta_2$ | 1 | 0 | $-1$ | 0 | 1 | 0 | $-1$ | 0 |
| $q = \zeta_3$ | 1 | $-1$ | 0 | 1 | $-1$ | 0 | 1 | $-1$ |
| $q = \zeta_4$ | 1 | $\sqrt{2}$ | 1 | 0 | $-1$ | $-\sqrt{2}$ | $-1$ | 0 |
| $q = \zeta_5$ | 1 | $\frac{1}{2}(-1+\sqrt{5})$ | $\frac{1}{2}(1-\sqrt{5})$ | $-1$ | 0 | 1 | $\frac{1}{2}(-1+\sqrt{5})$ | $\frac{1}{2}(1-\sqrt{5})$ |

(Note the difference between whether the $i$ in $\zeta_i$ is even or odd.) ◇

**Definition 9H.3.** Let $q^2 \in \mathbb{C}^*$ be not a second or third primitive root of unity. For $i = 0, 1, 2, 3$ define the ***ith Jones–Wenzl idempotent*** (JW idempotent for short) $\mathrm{JW}_i \in \mathrm{End}_{\mathbf{TL}^q_{\mathbb{C}\oplus\Subset}}(\bullet^i)$ as follows:

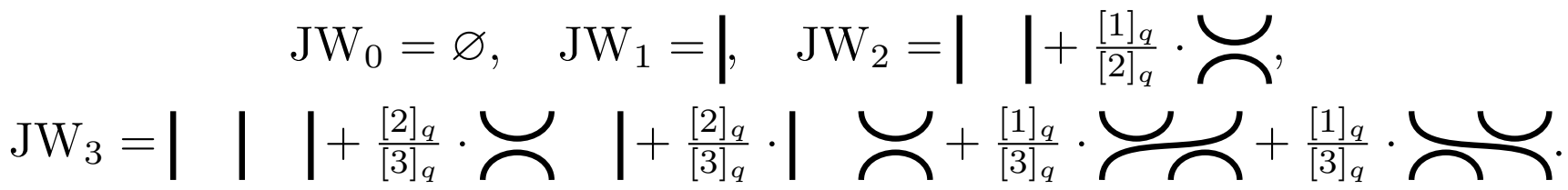

(Note that these are sums, with coefficients, of all possible diagrams with $i$ strands.) ◇

**Example 9H.4.** A calculation shows that we have the traces

$$\mathrm{tr}^{\mathbf{TL}^q_{\mathbb{C}\oplus\Subset}}(\mathrm{JW}_i) = (-1)^{i+1}[i+1]_q.$$

For example, we calculate

$$\mathrm{tr}^{\mathbf{TL}^q_{\mathbb{C}\oplus\Subset}}(\mathrm{JW}_2) =$$

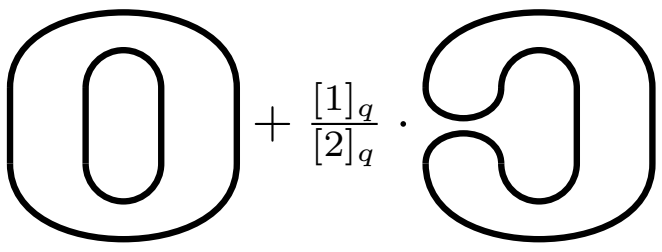

$$= [2]_q^2 - 1 = [3]_q,$$

and

$\mathrm{tr}^{\mathbf{TL}^q_{\mathbb{C}\oplus\mathbb{E}}}(\mathrm{JW}_3) =$

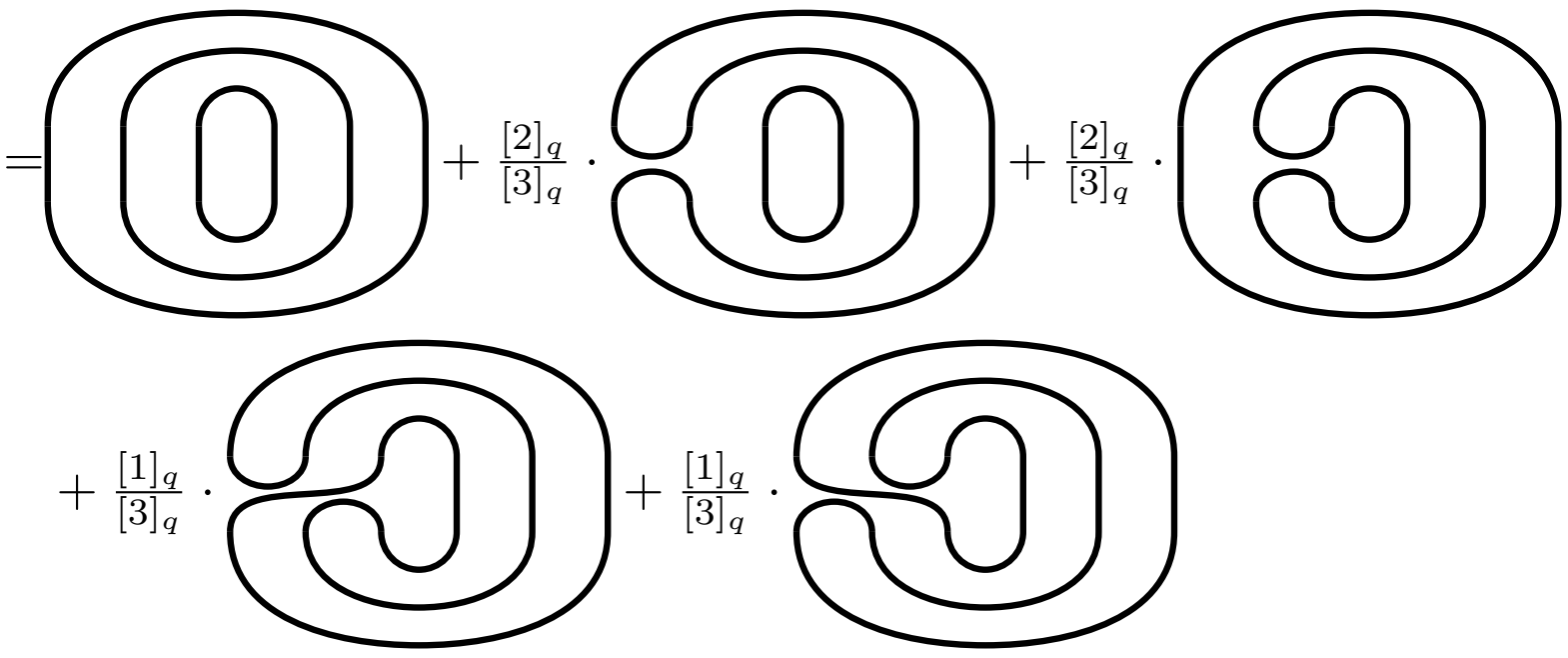

$$= \tfrac{1}{[3]_q}\big( - [2]_q^3[3]_q + [2]_q^3 + [2]_q^3 - [2]_q - [2]_q\big) = -[4]_q.$$

The remaining trace computations are immediate. ◇

**Lemma 9H.5.** *The JW projectors as in Definition 9H.3 are idempotents.*

*Proof.* This is an exercise, *cf.* Exercise 9M.5. □

**Definition 9H.6.** Let $q^2 \in \mathbb{C}^*$ be a primitive forth root of unity. We define the ***level*** $4$ ***semisimplified quotient*** $\mathbf{TL}^4_{\mathbb{C}\oplus\mathbb{E}}$ ***of*** $\mathbf{TL}^q_{\mathbb{C}\oplus\mathbb{E}}$ to be

$$\mathbf{TL}^4_{\mathbb{C}\oplus\mathbb{E}} = \mathbf{TL}^q_{\mathbb{C}\oplus\mathbb{E}}/\langle(\bullet^3, \mathrm{JW}_3)\rangle,$$

where the quotient is given by taking the $\otimes$ ideal generated by the object $(\bullet^3, \mathrm{JW}_3)$. ◇

Recall that we can identify objects of the idempotent completion with Im(e), see Section 7C.

**Proposition 9H.7.** *We have the following.*

*(i)* *We have* $\mathbf{TL}^4_{\mathbb{C}\oplus\mathbb{E}} \in \mathbf{MoCat}$.

*(ii)* *We have* $\mathrm{Si}(\mathbf{TL}^4_{\mathbb{C}\oplus\mathbb{E}}) = \{\mathrm{Im}(\mathrm{JW}_0) = \mathbb{1}, \mathtt{L}_1 = \mathrm{Im}(\mathrm{JW}_1), \mathtt{L}_2 = \mathrm{Im}(\mathrm{JW}_2)\}$.

*(iii)* *The fusion rules are*

| $\otimes$ | $\mathtt{L}_1$ | $\mathtt{L}_2$ |
|---|---|---|
| $\mathtt{L}_1$ | $\mathbb{1} \oplus \mathtt{L}_2$ | $\mathtt{L}_1$ |
| $\mathtt{L}_2$ | $\mathtt{L}_1$ | $\mathbb{1}$ |

*(iv)* *We have*

$$\mathbf{TL}^4_{\mathbb{C}\oplus\mathbb{E}} \simeq_{\mathbb{C}\oplus\star} \mathbf{fdMod}^q_4(\mathfrak{sl}_2).$$

*(v)* *For* $\mathtt{c} \in \{\mathtt{L}_0, \mathtt{L}_1, \mathtt{L}_2\}$ *(read such that the integer c corresponds to the subscript) we have*

(9H-8)

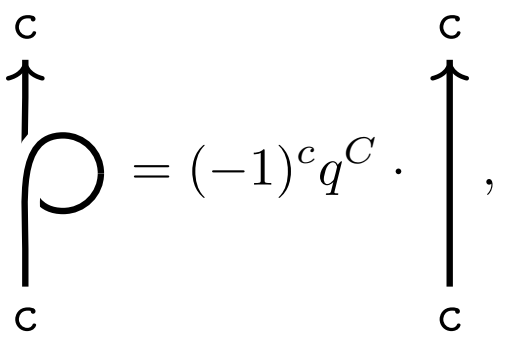

*where* $C = (c^2 + 2c)/2$ *is the* **quadratic Casimir number**.

*Proof.* Let us postpone the proof to a later section, but let us calculate the $S$ matrix of $\mathbf{TL}^4_{\mathbb{C}\oplus\mathbb{E}} \in \mathbf{Mo}$, *i.e.* we need to compute ***colored Jones polynomial*** of the Hopf link. We do this calculation generically, *i.e.* keeping $q$ a formal variable, and specialize $q$ in the end. We also use the diagrammatic notation

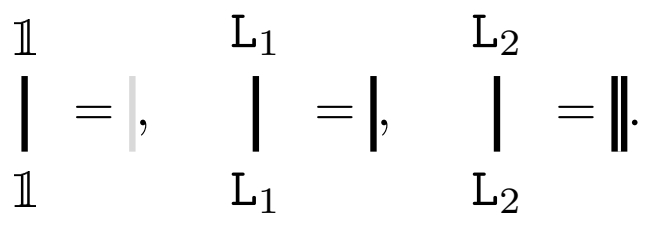

We note that we have the crossing formulas

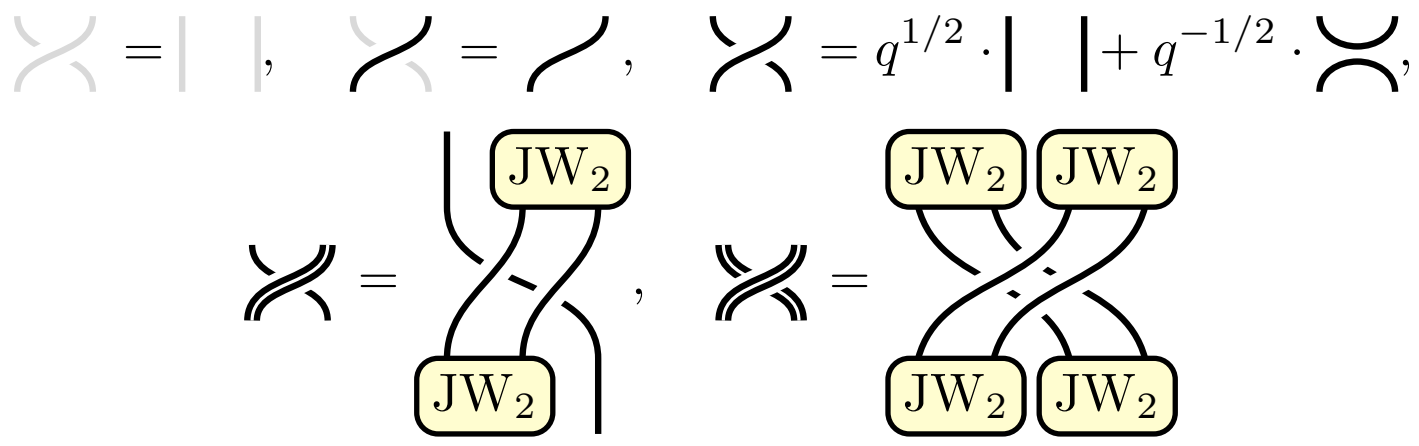

including mirrors. Using these we compute (also with reference to Example 7F.8)

$\mathrm{tr}^{\mathbf{C}}(\mathrm{s}_{\mathbb{1}\mathbb{1}}) = [1]_q$, $\quad \mathrm{tr}^{\mathbf{C}}(\mathrm{s}_{\mathbb{1}1}) = \mathrm{tr}^{\mathbf{C}}(\mathrm{s}_{1\mathbb{1}}) =$ 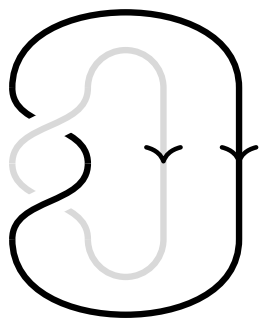 $= -[2]_q$, $\quad \mathrm{tr}^{\mathbf{C}}(\mathrm{s}_{11}) =$ 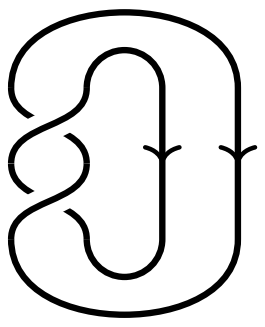 $= [4]_q$,

$\mathrm{tr}^{\mathbf{C}}(\mathrm{s}_{\mathbb{1}2}) = \mathrm{tr}^{\mathbf{C}}(\mathrm{s}_{2\mathbb{1}}) =$ 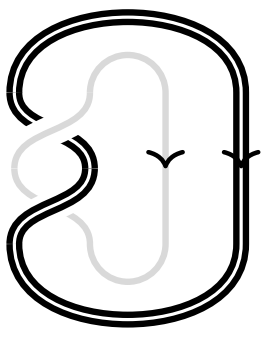 $= [3]_q$, $\quad \mathrm{tr}^{\mathbf{C}}(\mathrm{s}_{12}) = \mathrm{tr}^{\mathbf{C}}(\mathrm{s}_{21}) =$ 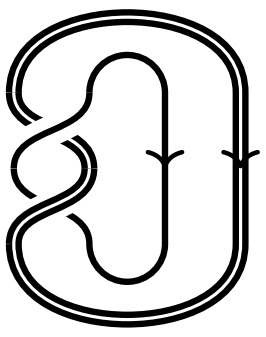 $= -[6]_q$,

$\mathrm{tr}^{\mathbf{C}}(\mathrm{s}_{22}) =$ 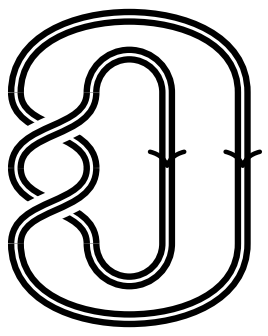 $= [9]_q$.

Thus, we get the $S$ matrix

$$S = \begin{pmatrix} [1]_q & -[2]_q & [3]_q \\ -[2]_q & [4]_q & -[6]_q \\ [3]_q & -[6]_q & [9]_q \end{pmatrix} = \begin{pmatrix} 1 & -\sqrt{2} & 1 \\ -\sqrt{2} & 0 & \sqrt{2} \\ 1 & \sqrt{2} & 1 \end{pmatrix}, \quad \det(S) = 4^3,$$

which shows that $S$ is invertible. A good way to see this directly is to observe that

$$S^2 = \begin{pmatrix} 4 & 0 & 0 \\ 0 & 4 & 0 \\ 0 & 0 & 4 \end{pmatrix}, \quad \det(S) = -8,$$

which then implies that $S$ is invertible with determinant $\sqrt{64}$, up to a sign. □

**Example 9H.9.** The category $\mathbf{TL}^4_{\mathbb{C}\oplus\mathbb{E}}$ from Definition 9H.6 exists more general for any $k \in \mathbb{Z}_{\geq 2}$, more in Section 10, and we one gets

$$\mathbf{TL}^k_{\mathbb{C}\oplus\mathbb{E}} \simeq_{\mathbb{C}\oplus\star} \mathbf{fdMod}^q_k(\mathfrak{sl}_2) \in \mathbf{Mo}.$$

Moreover, the $S$ matrix is

$$S = (\mathrm{tr}^{\mathbf{TL}^k_{\mathbb{C}\oplus\mathbb{E}}}(s_{ij}))^{k-2}_{i,j=0}, \quad \mathrm{tr}^{\mathbf{TL}^k_{\mathbb{C}\oplus\mathbb{E}}}(s_{ij}) = (-1)^{i+j}\big[(i+1)(j+1)\big]_q.$$

That this matrix $S$ is invertible is not difficult but also not trivial, see Proposition 9H.10.

Moreover, the formula in Equation 9H-8 stays exactly what it is. ◇

**Proposition 9H.10.** *The matrix $S$ in Example 9H.9 is invertible.*

*Proof.* Let $q = \exp(\pi i/k) \in \mathbb{C}$ and set $a = \sin(\pi/k)^{-1}$. The first calculation is that

$$[(i+1)(j+1)]_q = a^{-1}\sin(\pi(i+1)(j+1)/k).$$

Hence, if we remove the factor $a^{-1}$ we are left with

$$S' = \big((-1)^{i+j}\sin(\pi(i+1)(j+1)/k)\big)^{k-2}_{i,j=0}.$$

A calculation then yields

$$(S')^2 = \mathrm{diag}(k/2, ..., k/2),$$

which is a diagonal matrix with $k/2$ on the diagonal. Hence, $S$ is invertible. The same argument works for any other primitive $k$th root of unity $q^2$. □

9I. **Modular formulas.** The purpose of this section is to indicate the "Why?" of modular categories, which will be further justified in the upcoming sections.

*Remark* 9I.1. As we will see in Section 9J, modular categories have "good reasons" to have nice number theoretical properties. We are not giving proofs, as this is not our main purpose. (There are more formulas than the ones given below, see *e.g.* [**EGNO15**, Chapter 8].) ◇

We have three important matrices for $\mathbf{C} \in \mathbf{BFiat}$:

- The $S$ ***matrix*** which we have already seen in Section 9F;

- The $T$ ***matrix***

$$t_{ij} = \delta_{i,j} a(\mathtt{Z}_i), \quad T = (t_{ij})_{i,j=1}^{n} \in \mathrm{Mat}_{n\times n}(\mathbb{S}),$$

which is a diagonal matrix having the scalars from Equation 9F-13 on the diagonal;

- The $C$ ***matrix***

$$c_{ij} = \delta_{i,j^\star}, \quad C = (c_{ij})_{i,j=1}^{n} \in \mathrm{Mat}_{n\times n}(\mathbb{S}),$$

which is the $n \times n$ identity matrix if every object is self-dual.

**Example 9I.2.** Say we have two self-dual indecomposable objects $\mathtt{S}$ and $\mathtt{O}$, spinach and orchid, color coded as in Example 9F.6. Then:

$S =$ 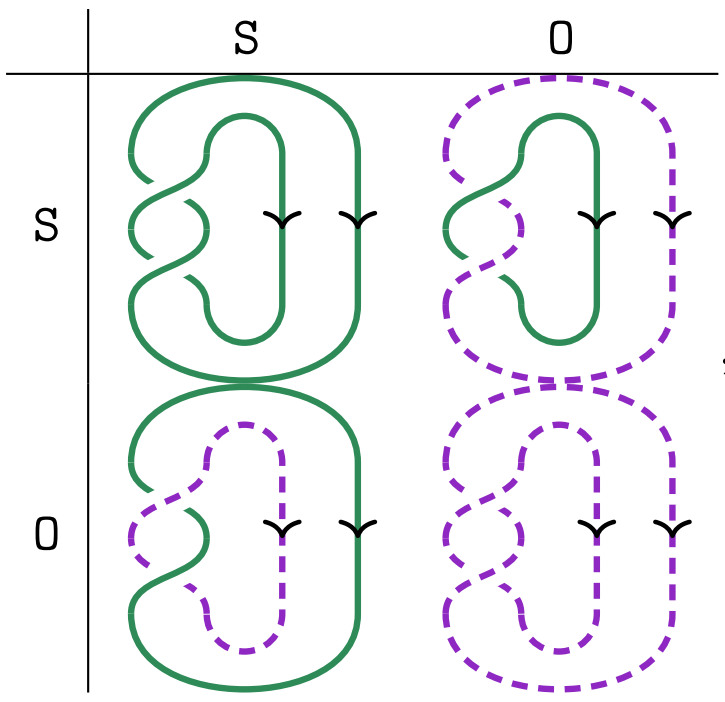

, $T =$ 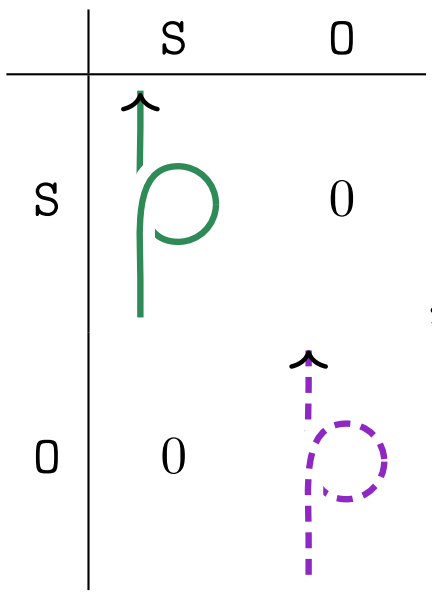

, $C =$

| | S | O |
|---|---|---|
| S | 1 | 0 |
| O | 0 | 1 |

.

The $\Delta_r$, $\Delta_l$, see Equation 9F-14, are sums over pairs of monochromatic circles. ◇

**Proposition 9I.3.** *Let* $\mathbf{C} \in \mathbf{MoCat}$ *be* $\mathbb{k}$ *linear, with* $\mathbb{k}$ *of characteristic zero. Then:*

*(i) We have* $\Delta_r, \Delta_l \in \mathbb{k}^*$.

*(ii) We have* $\dim(\mathbf{C}) = \Delta_r \Delta_l$, *which is non-zero.*

*(iii) We have* $C^2 = \mathrm{id}_n$, *where* $\mathrm{id}_n$ *is the* $n \times n$ *identity matrix.*

*(iv) We have* $S^2 = \dim(\mathbf{C}) \cdot C$.

*(v) We have* $S^4 = \dim(\mathbf{C})^2 \cdot \mathrm{id}_n$.

*(vi) We have* $(ST)^3 = \Delta_r \cdot S^2 = \dim(\mathbf{C})\Delta_r \cdot C$.

*(vii) We have* $TC = CT$.

*Proof.* Let us sketch a proof, see [**TV17**, Section 4.5.2] and [**EGNO15**, Proposition 8.14.2 and Theorem 8.16.1] for details.

*(i)+(ii).* Almost immediate from the definitions, only the non-zero claims need some work. But here we can use Proposition 9D.11 since ribbon categories are spherical by Lemma 5I.4.

*(iii).* Easy.

*(iv)+(vi).* The two non-trivial calculation, but the crucial ones. We omit the proof.

*(v).* From (iii) and (iv).

*(vii).* Easy. □

**Example 9I.4.** Let us come back to Proposition 9H.7 and the calculations therein. We have

$$a(\mathbb{1}) = 1, \quad a(\mathtt{L}_1) = -q^{-3/2}, \quad a(\mathtt{L}_2) = q^{-4} = -1,$$

$$\dim^{\mathbf{C}}(\mathbb{1}) = [1]_q = 1, \quad \dim^{\mathbf{C}}(\mathtt{L}_1) = -[2]_q = -\sqrt{2}, \quad \dim^{\mathbf{C}}(\mathtt{L}_2) = [3]_q = 1,$$

$$\dim(\mathbf{C}) = [1]_q^2 + [2]_q^2 + [3]_q^2 = 4,$$

$$\Delta_r = [1]_q^2 - q^{-3/2}[2]_q^2 + q^{-4}[3]_q^2 = -2\exp(\pi i5/8), \quad \Delta_l = [1]_q^2 - q^{3/2}[2]_q^2 + q^{4}[3]_q^2 = 2\exp(\pi i3/8),$$

$$4 = \big(-2\exp(\pi i5/8)\big)\big(2\exp(\pi i3/8)\big).$$

Moreover, the matrix $C$ is the identity and we have

$$S^2 = \begin{pmatrix} 1 & -\sqrt{2} & 1 \\ -\sqrt{2} & 0 & \sqrt{2} \\ 1 & \sqrt{2} & 1 \end{pmatrix} \cdot \begin{pmatrix} 1 & -\sqrt{2} & 1 \\ -\sqrt{2} & 0 & \sqrt{2} \\ 1 & \sqrt{2} & 1 \end{pmatrix} = \begin{pmatrix} 4 & 0 & 0 \\ 0 & 4 & 0 \\ 0 & 0 & 4 \end{pmatrix}.$$

Finally, we also calculate that

$$(ST)^3 = \begin{pmatrix} -8\exp(\pi i5/8) & 0 & 0 \\ 0 & -8\exp(\pi i5/8) & 0 \\ 0 & 0 & -8\exp(\pi i5/8) \end{pmatrix}.$$

Similar calculations work for the other Verlinde categories. ◇

Recall that $\dim^{\mathbf{C}}(\mathtt{L}) \in \Bbbk^{*}$ if $\mathbf{C}$ is semisimple, see Proposition 9D.2. The Verlinde formula, which is up next, gives us the surprising result that the $S$ matrix is in some sense encoded on the Grothendieck classes:

**Proposition 9I.5.** *Let $\mathbf{C} \in \mathbf{MoCat}$ be $\Bbbk$ linear. Then we have*

$$\dim(\mathbf{C})N_{jk}^{l} = \textstyle\sum_{i=1}^{n} \frac{s_{ij}s_{ik}s_{il^{\star}}}{\dim^{\mathbf{C}}(\mathtt{L}_i)}.$$

*Proof.* Omitted, see [EGNO15, Corollary 8.14.4]. □

**Example 9I.6.** We continue Example 9I.4: The $S$ matrix and the fusion rules are stated in (the proof of) Proposition 9H.7, and we indeed get *e.g.*

$$0 = \frac{-[2]_q[3]_q^2}{[1]_q} + \frac{[4]_q[6]_q^2}{[2]_q} + \frac{-[6]_q[9]_q^2}{[3]_q}, \quad 4 = \frac{[2]_q^2[3]_q}{[1]_q} + \frac{[4]_q^2[6]_q}{[2]_q} + \frac{[6]_q^2[9]_q}{[3]_q},$$

as predicted. ◇

9J. **The modular group.** Let us explain where the name "modular" comes from. To this end, we first recall that there is the ***Möbius group*** given by ***Möbius transformations***, *i.e.*

$$f\colon \mathbb{C} \to \mathbb{C}, \quad f(z) = \frac{az+b}{cz+d}, \quad \text{where } ad-bc \neq 0.$$

Geometrically $f$ can be described by first applying the inverse of the stereographic projection, see Figure 25, then a rigid motion of the sphere, and then applying the stereographic projection again. This usually produces nice pictures; see Figure 26.

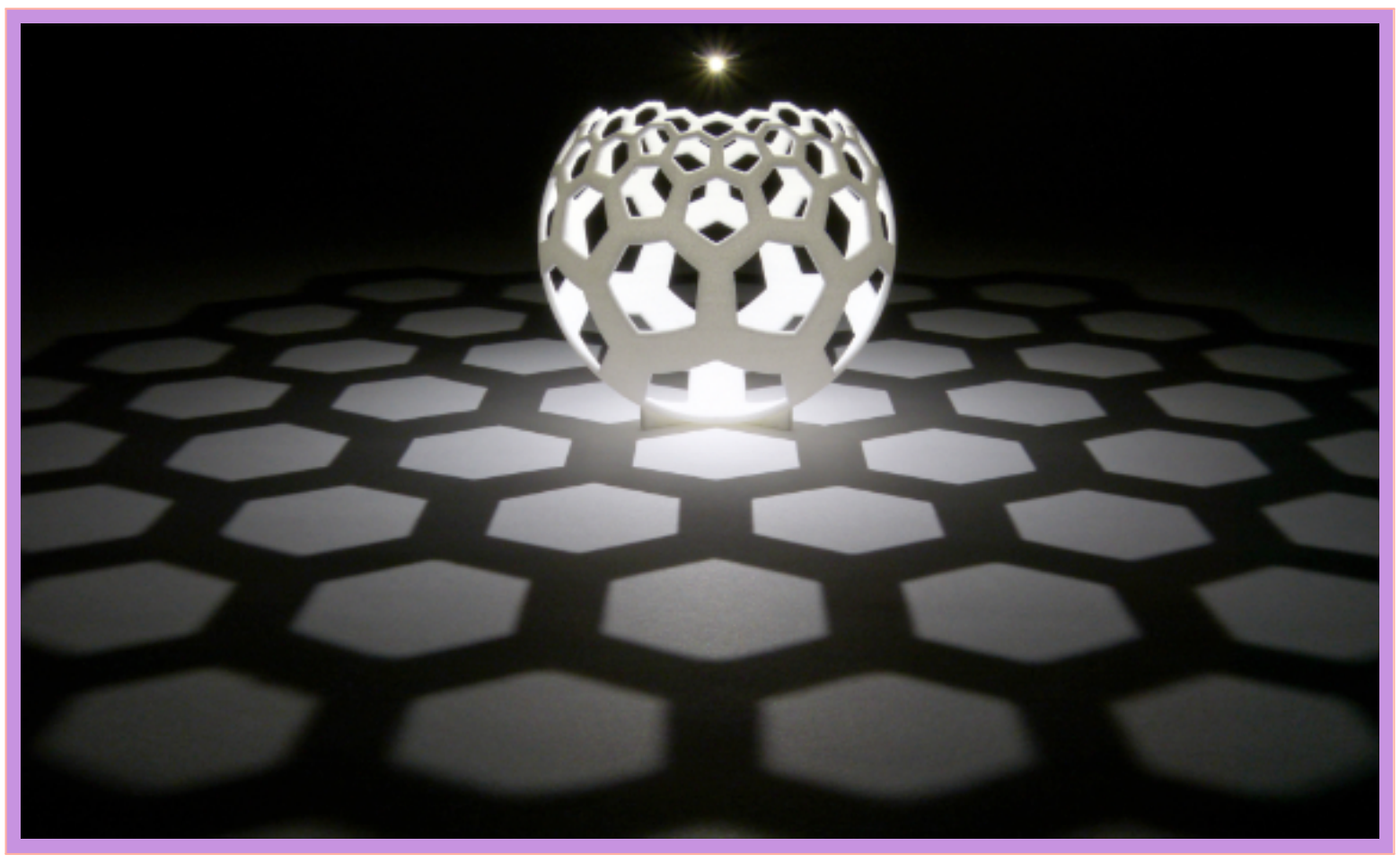

FIGURE 25. The stereographic projection identifies $\mathbb{C} = \mathbb{R}^2$ with the Riemann sphere $S^2 \cong \mathbb{PC}^1$ minus the northpole.

Picture from https://dtubbenhauer.com/lecture-geotop-2023.html

Algebraically speaking the Möbius group is

$$\mathrm{PGL}_2(\mathbb{C}) = \left\{ A = \begin{pmatrix} a & b \\ c & d \end{pmatrix} \middle| A \in \mathrm{Mat}_{2\times 2}(\mathbb{C}), \det(A) = ad-bc \neq 0 \right\} \Big/ (\pm 1),$$

which is the projective linear group of the Riemann sphere $\mathbb{PC}^1$. (Recall that "projective" in this sense should be read as "up to scalars".)

The "algebraic version" of the Möbius group is the ***modular group*** which, depending on the literature, is either $\mathrm{PGL}_2(\mathbb{Z})$ or $\mathrm{PSL}_2(\mathbb{Z})$, and is of crucial importance in *e.g.* number theory. For us the latter is the one we want, and in formulas:

$$\mathrm{PSL}_2(\mathbb{Z}) = \left\{ A = \begin{pmatrix} a & b \\ c & d \end{pmatrix} \middle| A \in \mathrm{Mat}_{2\times 2}(\mathbb{Z}), \det(A) = ad-bc = 1 \right\} \Big/ (\pm 1).$$

We will use either $\mathrm{SL}_2(\mathbb{Z})$ (defined similarly as $\mathrm{PSL}_2(\mathbb{Z})$, but without the quotient) or $\mathrm{PSL}_2(\mathbb{Z})$.

**Lemma 9J.1.** *We have*

$$\mathrm{SL}_2(\mathbb{Z}) \cong \langle S, T \mid S^4 = 1, (ST)^3 = S^2 \rangle,$$
$$\mathrm{PSL}_2(\mathbb{Z}) \cong \langle S, T \mid S^2 = 1, (ST)^3 = 1 \rangle \cong \mathbb{Z}/2\mathbb{Z} * \mathbb{Z}/3\mathbb{Z},$$

*where $S$ and $T$ correspond to the matrices $\left(\begin{smallmatrix} 0 & -1 \\ 1 & 0 \end{smallmatrix}\right)$ and $\left(\begin{smallmatrix} 1 & 1 \\ 0 & 1 \end{smallmatrix}\right)$, respectively.*

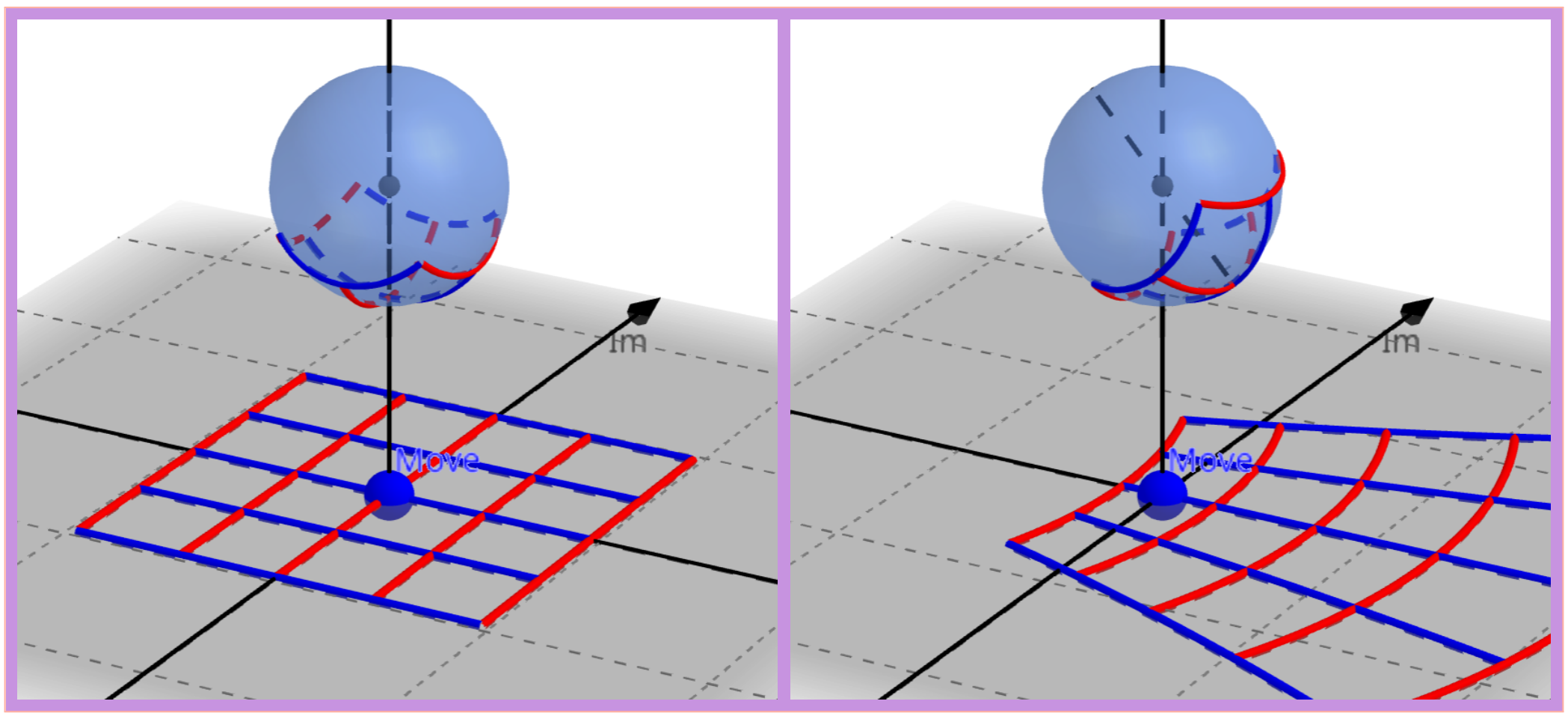

FIGURE 26. An example of a Möbius transformation. Left: One antiprojects the grid onto the sphere. Right: Now rotate the sphere and project back. The grid has now become distorted. Pictures from https://www.geogebra.org/m/GhaSJw3t, due to Juan Carlos Ponce Campuzano.

*Proof.* It is well known that $\mathrm{SL}_2(\mathbb{Z})$ and $\mathrm{PSL}_2(\mathbb{Z})$ have a generator–relation presentation of the form as in the lemma. A nice proof, using trees, can be found in [Ser80]. Note hereby that $S^2$ is $\left(\begin{smallmatrix}-1 & 0\\ 0 & -1\end{smallmatrix}\right)$, which is trivial in $\mathrm{PSL}_2(\mathbb{Z})$. □

Thus, the summary of the above, in particular Proposition 9I.3, is:

**Theorem 9J.2.** *Let* $\mathbf{C} \in \mathbf{MoCat}$ *be* $\Bbbk$ *linear such that* $\sqrt{\dim(\mathbf{C})} \in \Bbbk$. *Then*

$$\mathrm{SL}_2(\mathbb{Z}) \to \mathrm{PSL}_2(\mathbb{Z}) \to \mathrm{End}_{\mathbf{fdVec}_{\Bbbk}}\big(\mathrm{End}_{\mathbf{C}}(\mathbb{1})^n\big), \quad S \mapsto \tfrac{1}{\sqrt{\dim(\mathbf{C})}} \cdot S, T \mapsto T,$$

*defines a projective action of* $\mathrm{SL}_2(\mathbb{Z})$*, and of the modular group if the indecomposable objects of* $\mathbf{C}$ *are self-dual.* □

Let us stress again that "projective" hereby means "up to scalars."

9K. **Classification results.** Theorem 9J.2 indicates that one could try to classify modular categories using representations of $\mathrm{SL}_2(\mathbb{Z})$ and $\mathrm{PSL}_2(\mathbb{Z})$. However, both groups are almost free groups, so it is not hard to map it to other groups, and in particular to produce lots of representations. But there is still a way to produce several examples of $S$ matrices.

For all $N \in \mathbb{Z}_{\geq 2}$ we have a map

$$f_p \colon \mathrm{SL}_2(\mathbb{Z}) \to \mathrm{SL}_2(\mathbb{Z}/N\mathbb{Z}).$$

Since $\mathrm{SL}_2(\mathbb{Z}/N\mathbb{Z})$ is a finite group, its complex representation theory is well understood. Moreover, we have:

**Lemma 9K.1.** *Let* $\rho$ *denote a complex representation of* $\mathrm{SL}_2(\mathbb{Z})$ *on* $\mathbb{C}^r$.

**(a)** *If* $\rho(T)$ *has distinct eigenvalues and is realized by some* $\mathbb{C}$ *linear* $\mathbf{C} \in \mathbf{MoCat}$*, then* $\rho$ *factors through* $\mathrm{SL}_2(\mathbb{Z}/N\mathbb{Z})$ *with* $N$ *being the order of* $\rho(T)$.

**(b)** *If* $\rho$ *is realized by some* $\mathbb{C}$ *linear* $\mathbf{C} \in \mathbf{MoCat}$ *and* $r > 2$*, then* $\rho$ *is not a direct sum of one dimensional representation of* $\mathrm{SL}_2(\mathbb{Z})$.

*Proof.* This is explained in [BNRW16, Section 3.3]. □

Let $\hat{S}$ denote the matrix $1/\sqrt{\dim(\mathbf{C})} \cdot S$ (this is the matrix that actually satisfies $\hat{S}^2$ being trivial). We organize the matrix $\hat{S}$ so that the monoidal unit corresponds to the first entry.

**Theorem 9K.2.** *The possible* 1*-by-*1 *and* 2*-by-*2 *matrices* $\hat{S}$ *are*

$$\hat{S} = \big(1\big), \quad \hat{S} = \begin{pmatrix} 1 & 1 \\ 1 & -1 \end{pmatrix}, \quad \hat{S} = \begin{pmatrix} 1 & -1 \\ -1 & -1 \end{pmatrix}, \quad \hat{S} = \begin{pmatrix} 1 & \Phi \\ \Phi & -1 \end{pmatrix}, \quad \hat{S} = \begin{pmatrix} 1 & \phi \\ \phi & -1 \end{pmatrix},$$

*where* $\phi$ *and* $\Phi$ *are the golden ratio and it Galois conjugate. All of these a realized by some* $\mathbb{C}$ *linear* $\mathbf{C} \in \mathbf{MoCat}$*. Explicitly, and in the same order, the realization categories are:*

$$\mathbf{fdVec}_{\mathbb{C}}, \quad \begin{cases} \mathbf{fdMod}_3^q(\mathfrak{sl}_2), \\ q = \exp(\pi i/3), \end{cases} \quad \begin{cases} \mathbf{fdMod}_3^q(\mathfrak{sl}_2), \\ q = \exp(2\pi i/3), \end{cases} \quad \begin{cases} \mathbf{fdMod}_5^q(\mathfrak{so}_3), \\ q = \exp(\pi i/5), \end{cases} \quad \begin{cases} \mathbf{fdMod}_5^q(\mathfrak{so}_3), \\ q = \exp(2\pi i/5). \end{cases}$$

*Proof.* For the 1-by-1 matrices there is not much to be said: in this case $\mathbf{C}$ has to be $\mathbf{fdVec}_{\mathbb{C}}$ by Proposition 8I.6, while $S$ and $T$ have to be trivial: By Lemma 9K.1, $T$ needs to be of order one, and then $(ST)^3 = S^2$ forces $S$ to be trivial as well.

For the 2-by-2 case we summarize what we can already say without any further computation. We computed the $S$ matrices of $\mathbf{TL}^{k}_{\mathbb{C}\oplus\mathbb{C}} \simeq_{\mathbb{C}\oplus\star} \mathbf{fdMod}^{q}_{k}(\mathfrak{sl}_2) \in \mathbf{Mo}$ in Section 9H, and we can take either this category for $k = 3$ (where $q^2$ is a primitive third root of unity), where it has two indecomposable objects, or for $k = 4$ (where $q^2$ is a primitive forth root of unity) or $k = 5$ (where $q^2$ is a primitive fifth root of unity), but then we need to restrict to the monoidal subcategory generated by $\mathtt{L}_2$, that we call $\mathbf{fdMod}^{q}_{k}(\mathfrak{so}_3)$, so that we still have only two indecomposable objects. Since we have the symmetry that $q \leftrightarrow q^{-1}$ give the same category, the roots of unity to consider are $q = \exp(\pi i/3)$, $q = \exp(2\pi i/3)$, $q = \exp(\pi i/4)$, $q = \exp(3\pi i/4)$, $q = \exp(\pi i/5)$ and $q = \exp(2\pi i/5)$. The $S$ matrices, in the same order as in the previous sentence, are:

$$S = \begin{pmatrix} 1 & -1 \\ -1 & -1 \end{pmatrix}, \quad S = \begin{pmatrix} 1 & 1 \\ 1 & -1 \end{pmatrix}, \quad S = \begin{pmatrix} 1 & 1 \\ 1 & 1 \end{pmatrix}, \quad S = \begin{pmatrix} 1 & 1 \\ 1 & 1 \end{pmatrix}, \quad S = \begin{pmatrix} 1 & \phi \\ \phi & -1 \end{pmatrix}, \quad S = \begin{pmatrix} 1 & \Phi \\ \Phi & -1 \end{pmatrix}.$$

The two middle matrices are degenerate, so we do not need them. Otherwise, we end with the expected classification, so the only question is whether there are more solutions. The $T$ matrix, up to scalars, has to be (again, in the same order):

$$T = \begin{pmatrix} 1 & 0 \\ 0 & -i \end{pmatrix}, \quad T = \begin{pmatrix} 1 & 0 \\ 0 & i \end{pmatrix}, \quad T = \begin{pmatrix} 1 & 0 \\ 0 & \exp(4\pi i/5) \end{pmatrix}, \quad T = \begin{pmatrix} 1 & 0 \\ 0 & \exp(3\pi i/5) \end{pmatrix}.$$

Next, we use MAGMA (see Appendix A) to list all 2 dimensional simple representations of $\mathrm{SL}_2(\mathbb{Z}/p\mathbb{Z})$ when varying the prime $p$.

```
> //CharacterTable(SL(2,2));
> //CharacterTable(SL(2,3));
> CharacterTable(SL(2,5));
> //CharacterTable(SL(2,7));

  ----result----

> Character Table
> ---------------
>
>
> -------------------------------------------
> Class |   1  2  3  4    5    6  7    8    9
> Size  |   1  1 20 30   12   12 20   12   12
> Order |   1  2  3  4    5    5  6   10   10
> -------------------------------------------
> p  =  2   1  1  3  2    6    5  3    6    5
> p  =  3   1  2  1  4    6    5  2    9    8
> p  =  5   1  2  3  4    1    1  7    2    2
> -------------------------------------------
> X.1   +   1  1  1  1    1    1  1    1    1
> X.2   -   2 -2 -1  0   Z1 Z1#2  1  -Z1-Z1#2
> X.3   -   2 -2 -1  0 Z1#2   Z1  1-Z1#2  -Z1
> X.4   +   3  3  0 -1-Z1#2  -Z1  0-Z1#2  -Z1
> X.5   +   3  3  0 -1  -Z1-Z1#2  0  -Z1-Z1#2
> X.6   -   4 -4  1  0   -1   -1 -1    1    1
> X.7   +   4  4  1  0   -1   -1  1   -1   -1
> X.8   +   5  5 -1  1    0    0 -1    0    0
> X.9   -   6 -6  0  0    1    1  0   -1   -1
>
>
> Explanation of Character Value Symbols
> --------------------------------------
>
> # denotes algebraic conjugation, that is,
> #k indicates replacing the root of unity w by w^k
>
> Z1     = (CyclotomicField(5: Sparse := true)) ! [ RationalField() | 0, 0, 1, 1 ]
```

The reader may want to uncomment the other tables, but displaying them here would need too much space; we rather list the result:

| | $p=2$ | $p=3$ | $p=5$ | $p=7$ | $p>7$ |
|---|---|---|---|---|---|
| 2d | 1 | 3 | 2 | None | None |

After $p=5$ there is no other simple representations of dimension 2 and the only simple representations of dimension 1 already appear for $p=2$. This follows from classical character theory as summarized, for example, in [**FH91**, Section II.5].

Let us, for the sake of brevity, only do the $p=5$ case.

```
ComplexField(10)!(CyclotomicField(5: Sparse := true)) ! [ RationalField() | 0, 0, 1, 1 ]

   ----result----

-1.618033989 - 1.164153218E-10*$.1
```

So the `Z1` in the table is $-\phi$, which is the trace of $T=\left(\begin{smallmatrix}1 & 0\\ 0 & \exp(4\pi i/5)\end{smallmatrix}\right)$ after appropriate normalization. One can then check that the two representations `X.2` and `X.3` correspond to the action of $\mathrm{SL}_2(\mathbb{Z})$ coming from $\mathbf{fdMod}_5^q(\mathfrak{so}_3)$ for the two choices of $q^2$ being a primitive fifth root of unity as in the statement. □

9L. **Summary of categories.** Let us summarize the categorical constructions that in the end gave as modular categories.

- "Categorifying sets" ⇝ categories ⇝ access to morphisms.
- "Categorifying monoids" ⇝ monoidal categories ⇝ access to a two dimensional calculus.
- "Categorifying dual vector spaces" ⇝ rigid, pivotal, spherical categories ⇝ access to height operations.
- "Categorifying braid groups" ⇝ braided categories ⇝ access to the Reidemeister calculus.
- "Categorifying abelian groups" ⇝ additive and abelian categories ⇝ access to linear and homological algebra.
- "Categorifying algebras" ⇝ fiat and tensor categories ⇝ access to linear and homological algebra, a two dimensional calculus and height operations.
- "Categorifying semisimple algebras" ⇝ fusion categories ⇝ access to numerical data.
- "Action of the modular group" ⇝ modular categories ⇝ access to number theoretical data.

Here is an informal summary, grouping the categories into vanilla (in white), 2d topological (blueberry), 3d topological (spinach), algebraic (tomato) and all together (orchid):

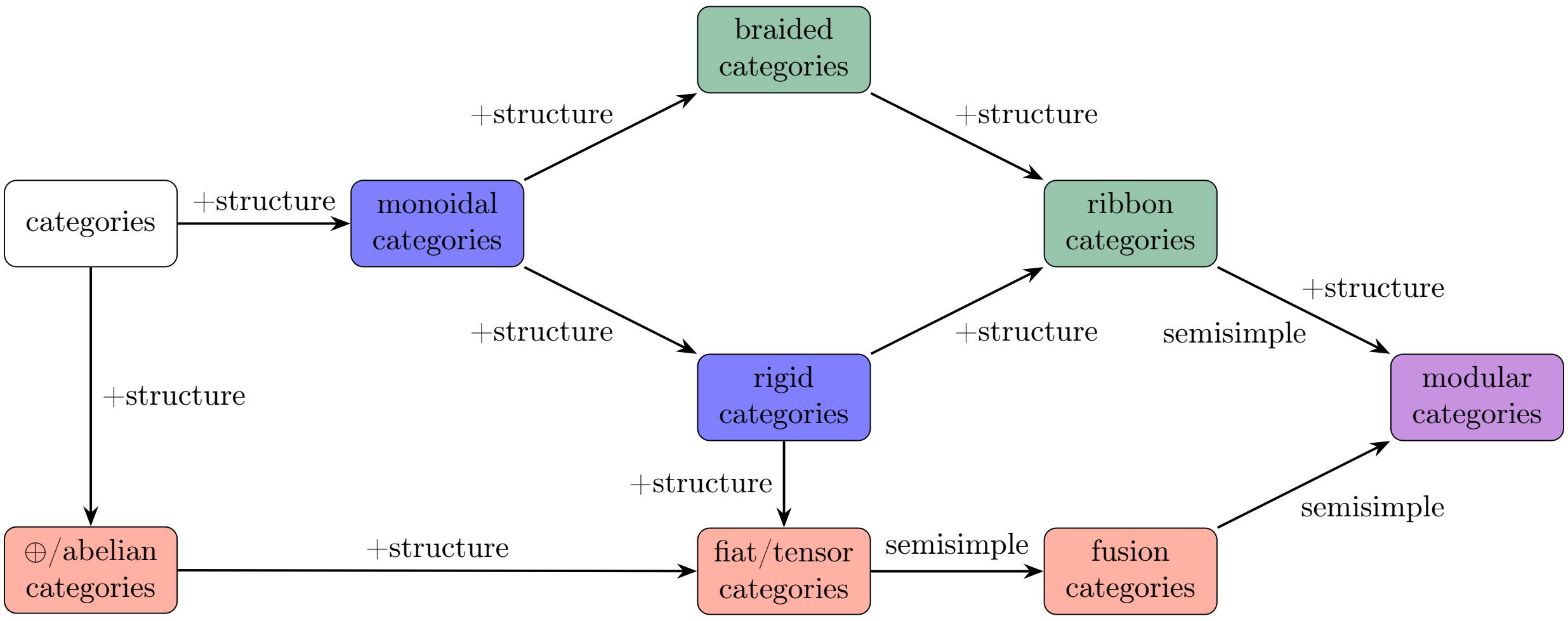

9M. **Exercises.**

*Exercise* 9M.1. Let $\mathrm{l}_i\in\mathrm{End}_{\mathbf{oqBr}}(\mathbb{1})$ for $i=1,2,3,4$ be the Hopf link with various orientations:

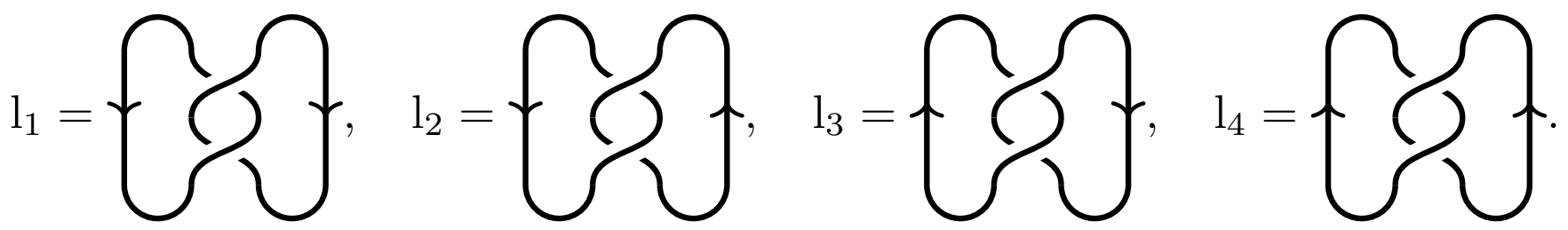

Find $\mathrm{f}_i\in\mathbf{qSym}$ for $i=1,2,3,4$ with $\mathrm{A}(\mathrm{f}_i)=\mathrm{l}_i$ (taking upwards-oriented right traces). Further, prove algebraically that the Markov moves hold after closing in $\mathbf{oqBr}$. ◇

*Exercise* 9M.2. Make Remark 9F.1 precise. For example, what kind of "free as an XYZ" should be satisfied by the category of colored braids? ◇

*Exercise* 9M.3. Let G be a finite group, and consider $\mathbf{fdMod}(\mathbb{C}[G])$ with standard braiding and duality. Show that the $S$ matrix of $\mathbf{fdMod}(\mathbb{C}[G])$ is of rank 1. For which G can $\mathbf{fdMod}(\mathbb{C}[G])$ be modular? ◇

*Exercise* 9M.4. Proof Proposition 9G.6. ◇

*Exercise* 9M.5. Verify as many claims from Section 9G as possible. ◇

## 10. Quantum invariants – a diagrammatic approach

Recall from Section 5L that a quantum invariant is structure preserving functor from a Brauer-type category to a, say, fiat, fusion or modular category.

How to construct quantum invariants?

We address this now for $\mathfrak{sl}_2$. To this end, we will need the ***Catalan numbers*** $C_n$, which is [**OEI23**, A000108] and reads $C_0 = 1, C_1 = 1, C_2 = 2, C_3 = 5, \ldots$. These numbers count many things, including crossingless (perfect) matchings of $\{1, \ldots, 2n\}$, with $n = 0$ just counted as one. The following hopefully explains what a crossingless matching is, where $n = 1$, $n = 2$ and $n = 3$ (so 2, 4 and 6 boundary points):

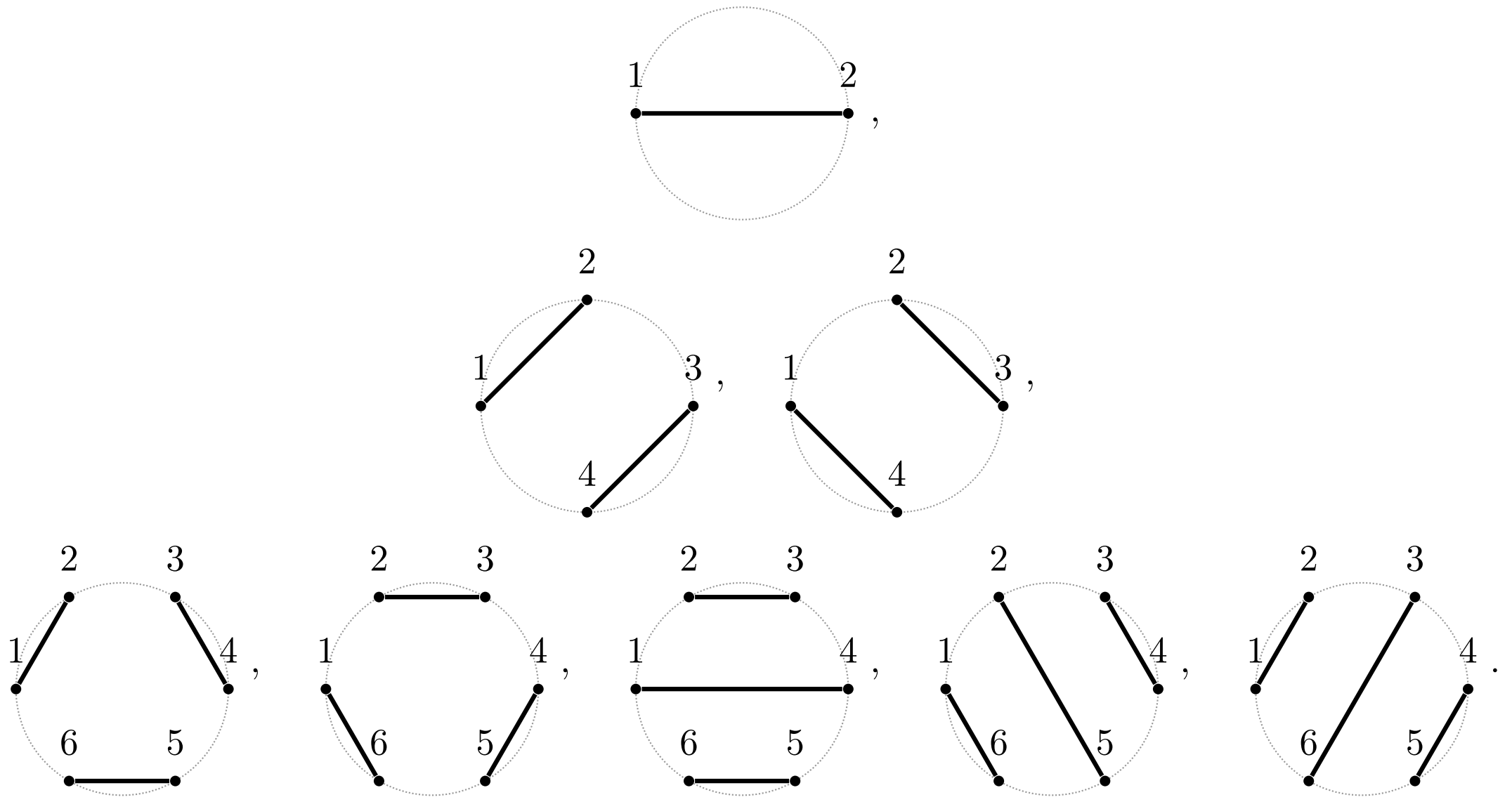

That these numbers are crucial for the study of quantum invariants associated with $\mathfrak{sl}_2$ is well known; see Figure 27 for an early reference.

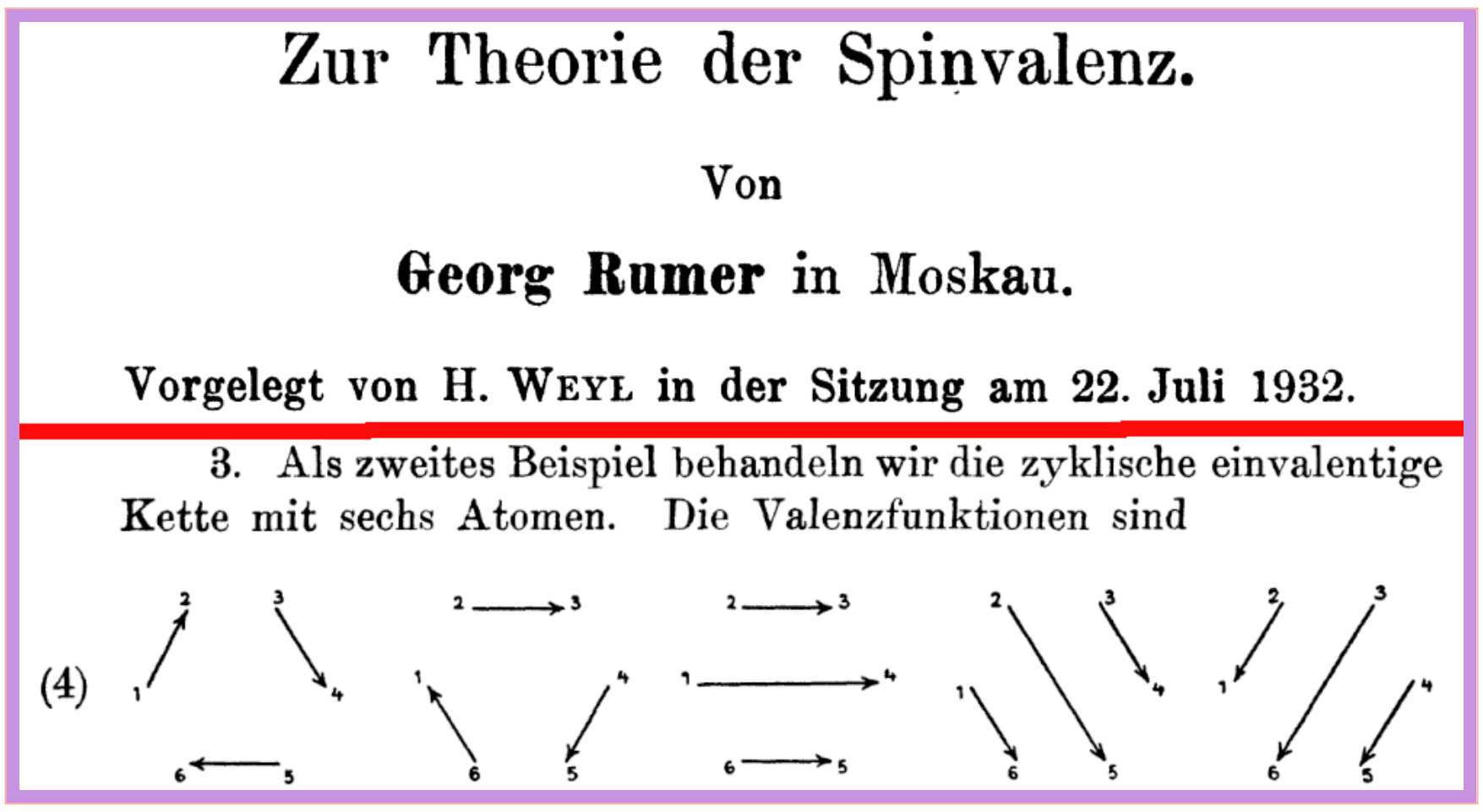

Zur Theorie der Spinvalenz.

Von

Georg Rumer in Moskau.

Vorgelegt von H. Weyl in der Sitzung am 22. Juli 1932.

3. Als zweites Beispiel behandeln wir die zyklische einvalentige Kette mit sechs Atomen. Die Valenzfunktionen sind

(4)

Figure 27. Catalan numbers count linearly independent bonds in one-valed chains with (here) six atoms. (The orientations do not play a role for us.) The picture is from [**Rum32**], and inspired [**RTW32**] which connects all of this with the representation theory of $\mathfrak{sl}_2$.

Picture from https://eudml.org/doc/59385

10A. **A word about conventions.** We need to be careful with the scalars:

*Convention* 10A.1. Recall that $\mathbb{S}$, $\mathbb{k}$ and $\mathbb{K}$ denote a ring, a field and an algebraically closed field, respectively. We further need $\mathbb{A} = \mathbb{Z}[v^{1/2}, v^{-1/2}]$ for $v$ being a formal variable (think: $v$ is generic/transcendental like $v = \pi$, in contrast to $q$ which will always be some specialization to a root of unity like $q = \exp(\pi i/5)$). ◇

*Convention* 10A.2. Since $q = \pm 1 \in \mathbb{S}$ will always behave differently from *e.g.* $q = \exp(2\pi i/l) \in \mathbb{C}$ for $l > 2$, we will not count $q = \pm 1$ as roots of unity. ◇

10B. **A quick count with Catalan numbers.** Before we get started, we recall a standard argument that appears often when working with quantum invariants.

**Lemma 10B.1.** *The Catalan numbers satisfy the following.*

**(a)** $C_0 = 1$ *and* $C_n = \sum_{i=1}^{n} C_{i-1}C_{n-i}$ *for* $n \in \mathbb{Z}_{>0}$.

**(b)** *We have* $C_n = \frac{1}{n+1}\binom{2n}{n}$.

**(c)** *We have* $\lim_{n\to\infty} \sqrt[n]{C_n} = 4$.

*Proof. (a).* Mark the arc from 1 to $2i$. This splits the picture into two sides, one being crossingless matchings of $\{2, ..., 2i-1\}$, counted by $C_{i-1}$, and the other being crossingless matchings of $\{2i+1, ..., 2n\}$, counted by $C_{n-i}$. In pictures:

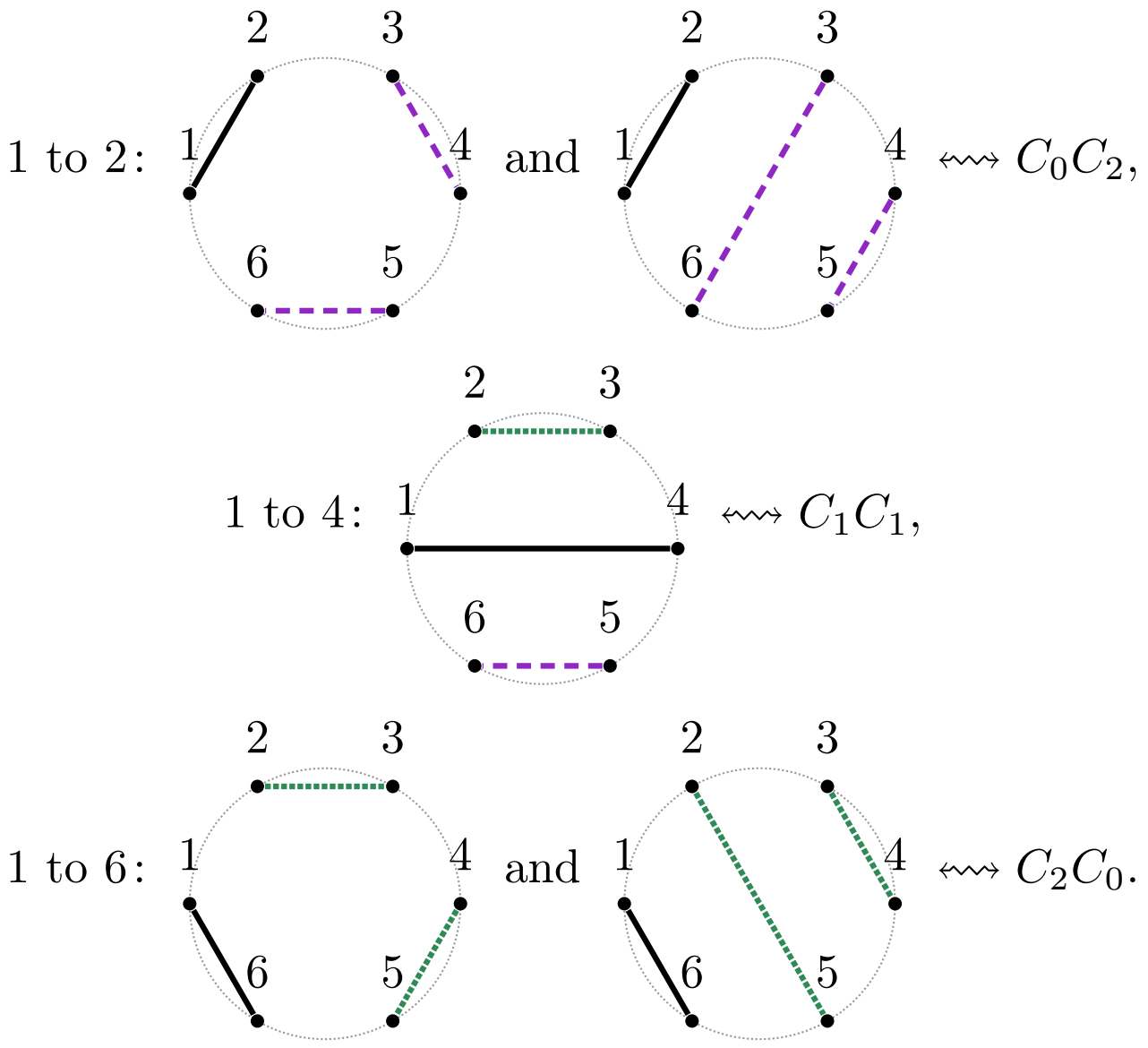

This shows (a).

The next two points follow from standard analysis of generating functions and their growth rates, see [FS09] for a cohesive summary.

*(b).* From (a) we get $c(x) = 1 + xc(x)^2$ for the generating function. Using the quadratic formula gives two possible solutions but only $c(x) = \frac{1}{2x}(1 - \sqrt{1-4x})$ satisfies $\lim_{x\to 0} c(x) = C_0 = 1$. Rewriting $c(x) = \frac{1}{2x}(1 - \sqrt{1-4x})$ using the binomial series gives $c(x) = \sum_{n=0}^{\infty} \frac{1}{n+1}\binom{2n}{n}x^n$, so the result follows.

*(c).* The radius of convergence of the power series $c(x) = \sum_{n=0}^{\infty} \frac{1}{n+1}\binom{2n}{n}x^n = \frac{1}{2x}(1 - \sqrt{1-4x})$ is 1/4, so the overall growth rate of $C_n$ is $4^n$. □

We will return to counting arguments as in Lemma 10B.1 later in Section 12.

10C. **An interlude about specializations.** For any pair $(\mathbb{S}, q)$ of a ring and an element $q^{1/2} \in \mathbb{S}^*$ (we need a square root of the parameters because of the braiding, *cf.* Equation 7F-3), let $\mathbb{A}$ act on $\mathbb{S}$ from the left via

$$\mathbb{A} \circlearrowright \mathbb{S}, \quad 1 \centerdot x = x, v \centerdot x = qx,$$

which makes $\mathbb{S}$ into a left $\mathbb{A}$ module. We thus get a specialization functor

$$\text{(10C-1)} \qquad {}_{-} \otimes_{\mathbb{A}}^{v=q} \mathbb{S} \colon \mathbf{Vec}_{\mathbb{A}} \to \mathbf{Vec}_{\mathbb{S}}, \quad \mathtt{X} \mapsto \mathtt{X} \otimes_{\mathbb{A}}^{v=q} \mathbb{S}.$$

In words, ${}_{-} \otimes_{\mathbb{A}}^{v=q} \mathbb{S}$ extends scalars to $\mathbb{S}$ and substitutes $v = q$. The pair $(\mathbb{S}, q)$ is also called a ***specialization***.

Let $\mathbf{C}_{\mathbb{A}} \in \mathbf{Cat}_{\mathbb{A}}$. Similarly as before, we get a ***category specialized at*** $q$, denoted by $\mathbf{C}_{\mathbb{S}}^{q}$, by extending scalars to $\mathbb{S}$ and substituting $v = q$. Formally:

- First we let $\mathrm{Ob}(\mathbf{C}_{\mathbb{S}}^{q}) = \mathrm{Ob}(\mathbf{C}_{\mathbb{A}})$;

- then we let

$$\mathrm{Hom}_{\mathbf{C}^q_{\mathbb{S}}}(\mathtt{X},\mathtt{Y})=\mathrm{Hom}_{\mathbf{C}_{\mathbb{A}}}(\mathtt{X},\mathtt{Y})\otimes^{v=q}_{\mathbb{A}}\mathbb{S},$$

where $_\otimes^{v=q}_{\mathbb{A}}\mathbb{S}$ is as in Equation 10C-1.

The following is easy, but crucial:

**Lemma 10C.2.** *Let $\mathbf{C}_{\mathbb{A}}$ be as above, and let $(\mathbb{S},q)$ any specialization. If $\mathrm{B}$ is a basis of $\mathrm{Hom}_{\mathbf{C}_{\mathbb{A}}}(\mathtt{X},\mathtt{Y})$, then it is also a basis of $\mathrm{Hom}_{\mathbf{C}^q_{\mathbb{S}}}(\mathtt{X},\mathtt{Y})$.* □

We have a few different looking cases which however will behave grouped as follows.

(I) The ***integral case*** which is either of:
- We stay with $\mathbb{A}$.
- We let $\mathbb{S}=\mathbb{Z}$ and $q=\pm 1$.

(II) The ***generic case*** (or ***generically***) which is either of:
- $\mathbb{S}=\Bbbk$ being a field of characteristic zero and $q=\pm 1$.
- We let $\mathbb{S}=\Bbbk$ be any field and $q\neq\pm 1$ not a root of unity in $\Bbbk$.

(III) The ***complex root of unity case*** where $\mathbb{S}=\Bbbk$ is a field of characteristic $p=0$ and $q\in\Bbbk$ is a root of unity such that $q^2$ is of order $l$.

(IV) The ***mixed root of unity case*** (or ***mixed case***) where $\mathbb{S}=\Bbbk$ is a field of characteristic $p>0$ and $q\in\Bbbk$ is a root of unity such that $q^2$ is of order $l$.

(V) The ***finite characteristic case*** (or ***char $p$ case***) where $\mathbb{S}=\Bbbk$ is a field of characteristic $p>0$ and $q=\pm 1$.

This is not clear at all at this point, but:
- The integral case is universal and treats (II)–(V) at the same time.
- The remaining cases are ordered by difficulty with (II) being the easiest.

**Example 10C.3.** It is a bit confusing, so let us make clear that the generic case includes the choice $\mathbb{S}=\mathbb{C}(q)$, for a formal variable $q$, which is probably the most common ground field in quantum topology and quantum algebra. ◇

The philosophy is that we have a category $\mathbf{C}_{\mathbb{A}}$, defined integrally, with an integral basis and ***integral objects*** ("objects which are always defined"), whose decomposition however depend on the specialization: Usually $\mathbf{C}_{\mathbb{A}}$ has ***pseudo idempotents***, *i.e.* morphisms with

$$\mathrm{e}^2=a\cdot\mathrm{e},\quad a\in\mathbb{A}.$$

As we have already seen in Section 7B, idempotents are very important to understand categories at hand, and they should decompose the integral objects into indecomposables. So we want to divide by $a$ to get an idempotent:

$$\big(\mathrm{e}^2=a\cdot\mathrm{e}\big)\Rightarrow\big((\tfrac{1}{a}\mathrm{e})^2=\tfrac{1}{a}\mathrm{e}\big).$$

So the crucial fact we need is whether the scalar $a$ is invertible, which depends on the choice of specialization. Here is a prototypical example:

**Example 10C.4.** Let us come back to the symmetric group $\mathrm{S}_3$ and let us consider the integral case $\mathbb{Z}[\mathrm{S}_3]$ and its category $\mathbf{C}^1_{\mathbb{Z}}=\mathbf{fdMod}\big(\mathbb{Z}[\mathrm{S}_3]\big)$. In this case an integral object would be $\mathbb{Z}[\mathrm{S}_3]$ itself, which we can always define.

We already know that generically $\mathbf{C}^1_{\mathbb{S}}$ is semisimple, *e.g.* for $\mathbb{S}=\mathbb{C}$ we have

$$\mathbb{C}[\mathrm{S}_3]\cong\mathbb{1}\oplus 2\cdot\mathtt{L}_s\oplus\mathtt{L}_{-1},\tag{10C-5}$$

see Example 8C.6. However, for $\mathbb{S}=\overline{\mathbb{F}}_2$ or $\mathbb{S}=\overline{\mathbb{F}}_3$ this is not the case anymore, see Example 8E.7, and one can see this in $\mathbf{fdMod}\big(\mathbb{Z}[\mathrm{S}_3]\big)$ as follows.

Let $\mathrm{S}_3=\{1,s,t,ts,st,sts=tst\}$, where, in graphical notation,

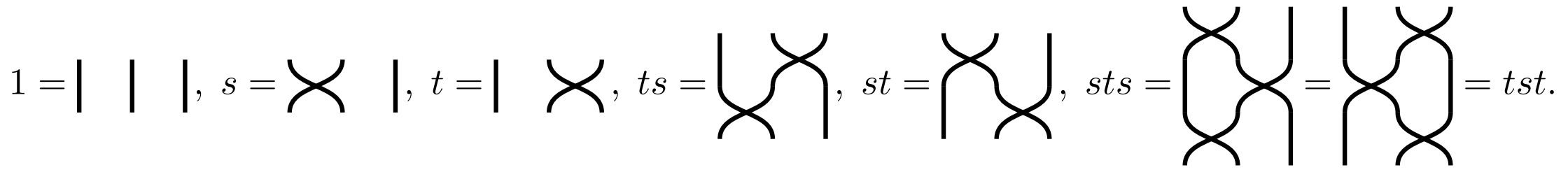

The category $\mathbf{fdMod}\big(\mathbb{Z}[\mathrm{S}_3]\big)$ has the following four pseudo idempotents:

$$\mathrm{e}_{\mathbb{1}}=1+s+t+ts+st+sts,\quad \mathrm{e}_{\mathbb{1}}^2=\boxed{6}\cdot\mathrm{e}_{\mathbb{1}},$$

$$\mathrm{e}_{s,1}=1+s-ts-sts,\quad \mathrm{e}_{s,1}^2=\boxed{3}\cdot\mathrm{e}_{s,1},$$

$$\mathrm{e}_{s,2} = 1 - s - st + sts, \quad \mathrm{e}_{s,2}^2 = \boxed{3} \cdot \mathrm{e}_{s,2},$$
$$\mathrm{e}_{-1} = 1 - s - t + ts + st - sts, \quad \mathrm{e}_{-1}^2 = \boxed{6} \cdot \mathrm{e}_{\mathtt{L}_{-1}},$$
$$\text{orthogonal: } \mathrm{xy} = 0, \quad \text{where } \mathrm{x}, \mathrm{y} \in \{\mathrm{e}_1, \mathrm{e}_{s,1}, \mathrm{e}_{s,2}, \mathrm{e}_{-1}\}, \mathrm{x} \neq \mathrm{y},$$
$$\text{pseudo complete: } \mathrm{e}_1 + 2\mathrm{e}_{s,1} + 2\mathrm{e}_{s,2} + \mathrm{e}_{-1} = 6.$$

We recover the three different cases we were already aware of: In $\mathbb{C}$ we can scale them to be idempotents and we get the decomposition Equation 10C-5. For $\mathbb{S} = \overline{\mathbb{F}}_2$ we can scale the middle two pseudo idempotents to get idempotents, while for $\mathbb{S} = \overline{\mathbb{F}}_3$ no scaling works.

Note that integrally we can not decompose $\mathbb{Z}[\mathrm{S}_3]$. In fact, we get a decomposition into indecomposables depending on the specialization: the generic case is Equation 10C-5, while

$$\overline{\mathbb{F}}_2[\mathrm{S}_3] \cong \mathtt{P}_{\mathbb{1}} \oplus \mathtt{P}_s, \quad \overline{\mathbb{F}}_3[\mathrm{S}_3] \cong \mathtt{P}_{\mathbb{1}} \oplus \mathtt{P}_{\mathbb{1}'},$$

are the decompositions in characteristic 2 and 3, respectively. ◇

In this section we discuss this strategy for the Rumer–Teller–Weyl category, which will ultimately lead to the construction of the Verlinde categories (for $\mathfrak{sl}_2$) and Jones-type invariants.

10D. **An integral basis for the Rumer–Teller–Weyl category.** Let $\mathbf{TL}_{\mathbb{A}}^v$ denote the category as defined in Definition 7F.1, but over the ground ring $\mathbb{A}$ and without taking additive and idempotent closures (for the time being). Recall that $\mathbf{TL}_{\mathbb{A}}^v$ is an $\mathbb{A}$ linear ribbon category.

We let $E = \mathbb{R}$, $E^+ = \mathbb{R}_{\geq 0}$, $X = \mathbb{Z} \subset E$ and $X^+ = \mathbb{Z}_{\geq 0} \subset E$. We also let $\Phi = \{\varepsilon_1 = 1, \varepsilon_{-1} = -1\} \subset E$.

**Definition 10D.1.** A(n integral) ***path*** $\pi$ of length $k$ in $E$ is a word $\pi = \pi_1 \dots \pi_k \in \Phi^k$ of length $k$. Such a path is called ***non-negative*** if $\sum_{i=1}^{j} \pi_i \in X^+$ for all $1 \leq j \leq k$. ◇

**Definition 10D.2.** The ***weight*** of a path $\pi$ is $\lambda(\pi) = \sum_{i=1}^{k} \pi_i \in X$. ◇

**Example 10D.3.** We think of paths as "honest" paths in $E$, starting at 0, using the rules

$\varepsilon_1 \rightsquigarrow \dots$ 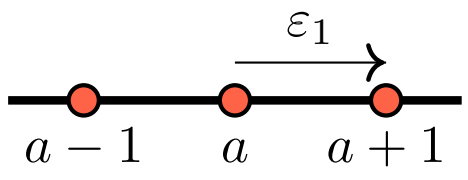 $\dots$, $\varepsilon_{-1} \rightsquigarrow \dots$ 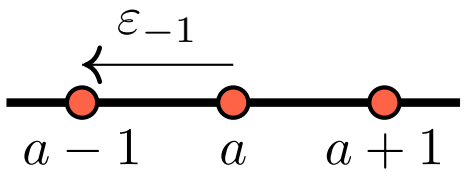 $\dots$ .

Using this interpretation, a path is non-negative if and only if it stays in $E^+$, and the weight is its endpoint in $E^+$. For example,

$\pi = \varepsilon_1\varepsilon_1\varepsilon_{-1}\varepsilon_1 \rightsquigarrow \dots$ 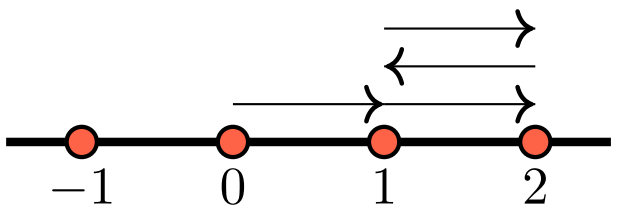 $\dots$

is a non-negative path of length four and weight two. Note also that having a non-negative weight, a.k.a. endpoint, is not enough to be a non-negative path, *e.g.* the following is not a non-negative path, but has non-negative weight:

$\pi = \varepsilon_{-1}\varepsilon_1 \rightsquigarrow \dots$ 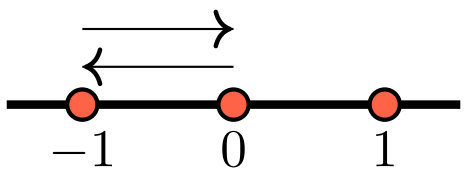 $\dots$ .

Having the path picture in mind is recommended. ◇

**Definition 10D.4.** To $\varepsilon_{\pm 1}$ we associate operators via:

(10D-5) $\varepsilon_1(\mathrm{f})$: 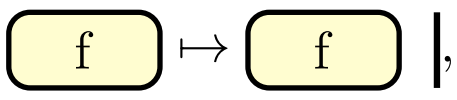, $\varepsilon_{-1}(\mathrm{f})$: 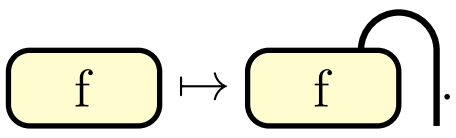.

In words, if we already have a morphism $\mathrm{f} \in \mathbf{TL}_{\mathbb{A}}^v$, then we obtain to new morphism $\mathrm{f}\varepsilon_{\pm 1} \in \mathbf{TL}_{\mathbb{A}}^v$ by either adding a strand or a cap to the right. ◇

**Definition 10D.6.** The ***downward integral ladder*** $\mathrm{d}(\pi) \in \mathbf{TL}_{\mathbb{A}}^v$ associated to a non-negative path $\pi$ is the morphism $\pi(\mathrm{id}_\emptyset)$ obtained by successively using Equation 10D-5. The ***upward integral ladder*** $\mathrm{u}(\pi) \in \mathbf{TL}_{\mathbb{A}}^v$ is the corresponding downward integral ladder horizontally mirrored. ◇

By vertical mirroring it suffices to only calculate down integral ladder, of course, *e.g.*

$\mathrm{d}(\pi) =$ 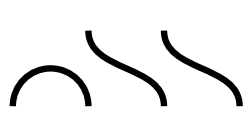 $\Leftrightarrow \mathrm{u}(\pi) =$ 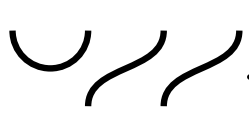.

**Example 10D.7.** One can easily check that

$$\begin{aligned}\lambda = 4\colon &\quad \pi_1 = \varepsilon_1\varepsilon_1\varepsilon_1\varepsilon_1,\\ \lambda = 2\colon &\quad \pi_2 = \varepsilon_1\varepsilon_{-1}\varepsilon_1\varepsilon_1,\quad \pi_3 = \varepsilon_1\varepsilon_1\varepsilon_{-1}\varepsilon_1,\quad \pi_4 = \varepsilon_1\varepsilon_1\varepsilon_1\varepsilon_{-1},\\ \lambda = 0\colon &\quad \pi_5 = \varepsilon_1\varepsilon_{-1}\varepsilon_1\varepsilon_{-1},\quad \pi_6 = \varepsilon_1\varepsilon_1\varepsilon_{-1}\varepsilon_{-1},\end{aligned}$$

are the only non-negative paths of length four. We get:

$\lambda = 4\colon$ $\mathrm{d}(\pi_1) =$ 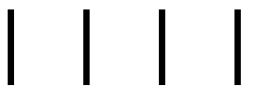

$\lambda = 2\colon$ $\mathrm{d}(\pi_2) =$ 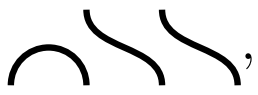, $\mathrm{d}(\pi_3) =$ 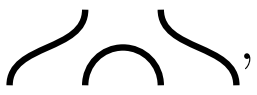, $\mathrm{d}(\pi_4) =$ 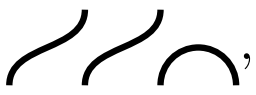,

$\lambda = 0\colon$ $\mathrm{d}(\pi_5) =$ 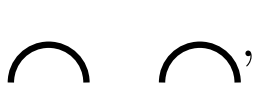, $\mathrm{d}(\pi_6) =$ 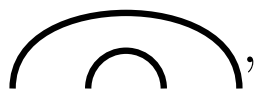,

as the corresponding ladders. ◇

Moreover, we denote by $(\lambda, \pi_d^m, \pi_u^n)$ a triple of a weight $\lambda$, and two non-negative paths $\pi_d^m$ and $\pi_u^n$ of this weight, of length as indicated by the superscripts.

**Definition 10D.8.** The ***integral ladder*** associated to the triple $(\lambda, \pi_d^m, \pi_u^n)$ is the morphism

$$\mathrm{c}^{\lambda}_{\pi_u^n,\pi_d^m} = \mathrm{u}(\pi_u^n)\mathrm{d}(\pi_d^m) \in \mathrm{Hom}_{\mathbf{TL}^v_{\mathbb{A}}}(\bullet^m, \bullet^n).$$

(Essentially, put up on top of down.) ◇

We also write $\mathrm{c}^{\lambda}_{\mathrm{u,d}} = \mathrm{c}^{\lambda}_{\pi_u^n,\pi_d^m}$ *etc.* for simplicity of notation.

**Example 10D.9.** With respect to Example 10D.7: the integral ladders which can be obtained from the diagrams therein are 14 in total:

$\lambda = 4\colon$ $\lambda = 2\colon$ $\lambda = 0\colon$

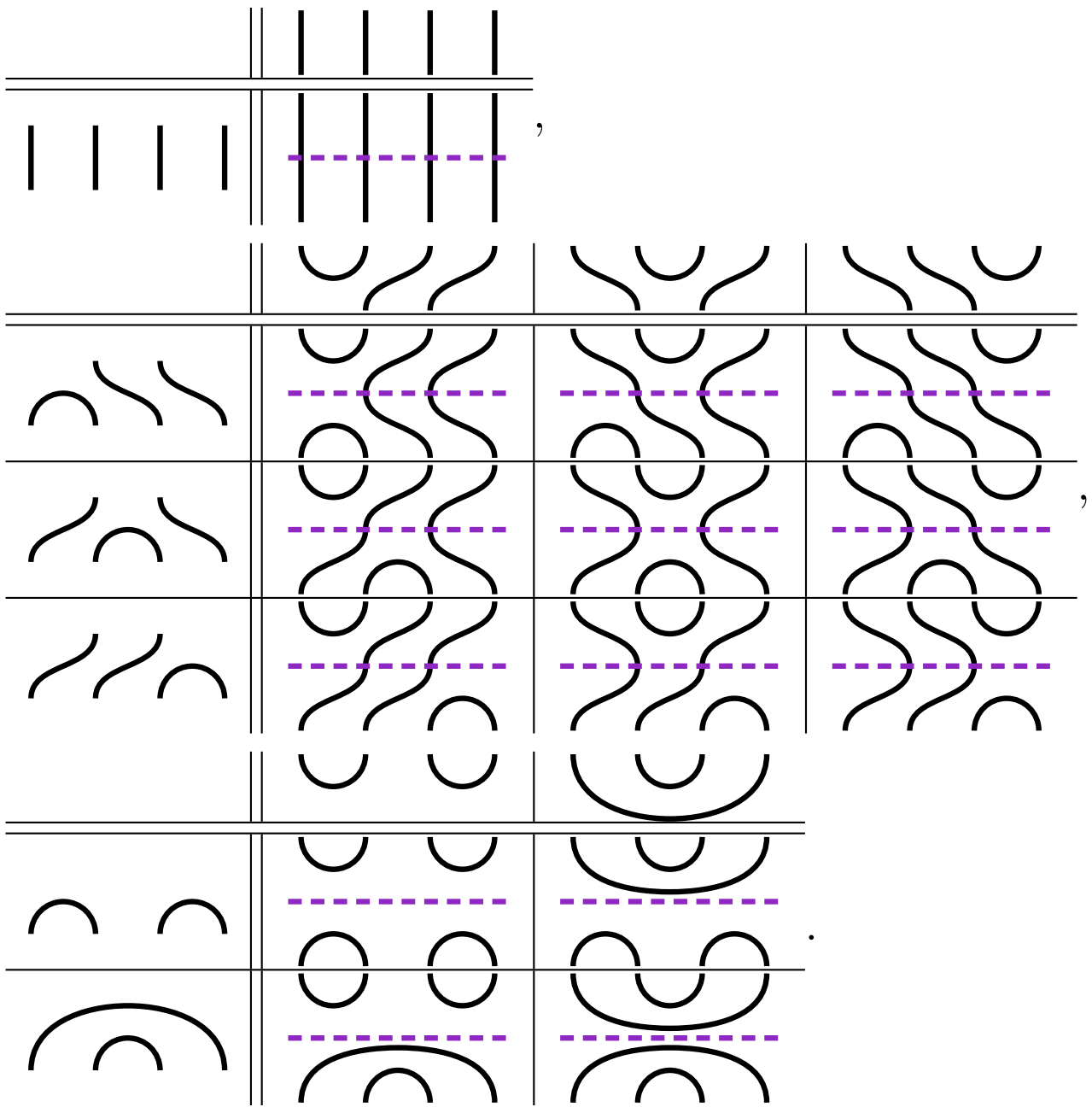

Note that the weight $\lambda$ can be read off in the middle, as indicated. ◇

**Example 10D.10.** In Example 10D.9 we have calculated integral ladder morphisms in $\mathrm{End}_{\mathbf{TL}^v_{\mathbb{A}}}(\bullet^4)$. For, say $\mathrm{Hom}_{\mathbf{TL}^v_{\mathbb{A}}}(\bullet^2, \bullet^4)$ we first observe that we have only two non-negative paths of length 2, and thus, also only two downwards integral ladders:

$\pi_1 = \varepsilon_1\varepsilon_1 \leftrightsquigarrow$ 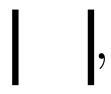, $\pi_2 = \varepsilon_1\varepsilon_{-1} \leftrightsquigarrow$ 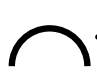.

Thus, we get the following integral ladders in $\mathrm{Hom}_{\mathbf{TL}^v_{\mathbb{A}}}(\bullet^2, \bullet^4)$:

$\lambda = 2\colon$ 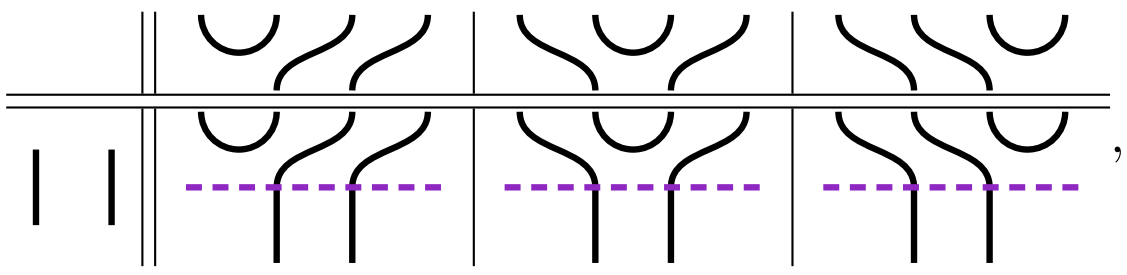, $\lambda = 0\colon$ 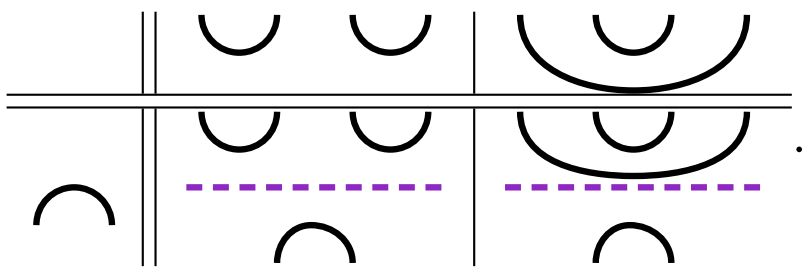.

Thus, there are five morphisms in total. ◇

Recall the Catalan numbers $C_n$ from, *e.g.*, Lemma 10B.1.

**Theorem 10D.11.** *The sets of the form*

$$\mathrm{IL} = \{ c^{\lambda}_{\pi^n_u, \pi^m_d} \mid \lambda \text{ a weight}, \pi^n_u, \pi^m_d \text{ non-negative paths} \} \tag{10D-12}$$

*are bases of* $\mathrm{Hom}_{\mathbf{TL}^v_{\mathbb{A}}}(\bullet^m, \bullet^n)$. *In particular,* $\mathrm{Hom}_{\mathbf{TL}^v_{\mathbb{A}}}(\bullet^m, \bullet^n)$ *is free and of rank* 0*, if* $m+n$ *is odd, and otherwise*

$$\mathrm{rank}_{\mathbb{A}} \mathrm{Hom}_{\mathbf{TL}^v_{\mathbb{A}}}(\bullet^m, \bullet^n) = C_{1/2(m+n)}.$$

*Proof.* In this formulation the crucial observation is [Eli15, Theorem 2.57], showing linear independence. That integral ladders span follows by observing that $\mathrm{Hom}_{\mathbf{TL}^v_{\mathbb{A}}}(\bullet^m, \mathbb{1})$ is clearly spanned by integral ladders which are in one-to-one correspondence with crossingless matchings. To see this, one cuts the circle at one fixed points and bends the pictures over. Here is an example, where we cut at the marked spot:

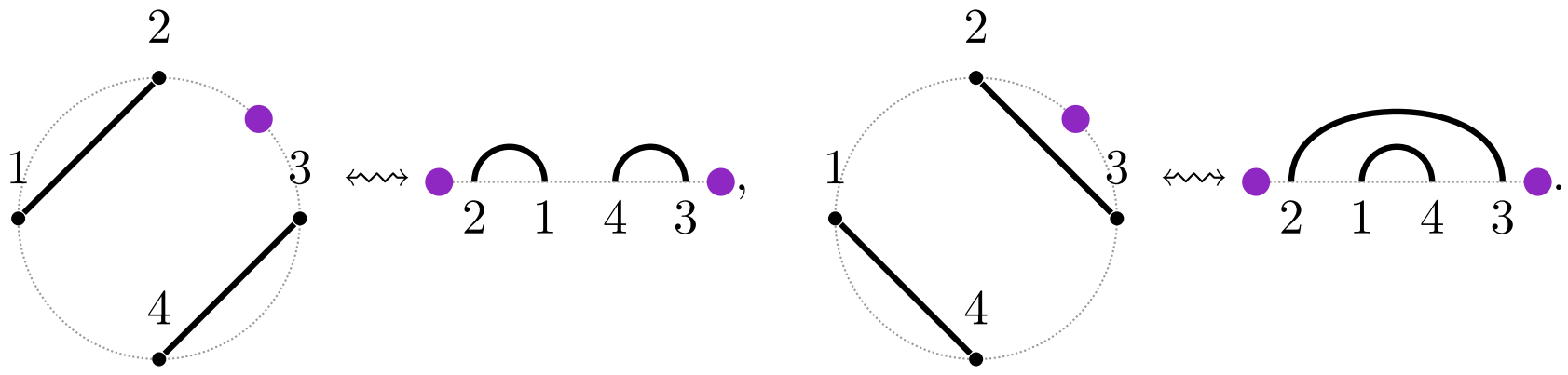

The inverse operation is to identify the two points on the bottom line. This implies the claims since mating preserves this property; see Theorem 4C.8. □

We call IL as in Equation 10D-12 the ***integral ladder basis*** of $\mathbf{TL}^v_{\mathbb{A}}$. Note that this basis is built using a ***bottleneck principle***, and we will also illustrate the basis elements by

(10D-13) 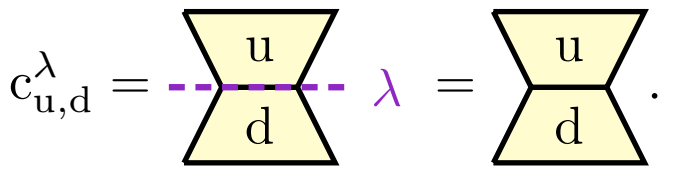

Denote by $\mathbf{I}_{<i}$ the set of morphisms in $\mathbf{TL}^v_{\mathbb{A}}$ which contain at most $i-1$ through strands, *e.g.*

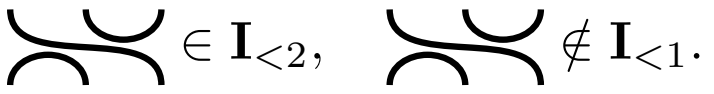

Clearly, $\mathbf{I}_{<j} \subset \mathbf{I}_{<i}$ if $j \leq i$. Moreover:

**Lemma 10D.14.** *The set* $\mathbf{I}_{<i}$ *is an ideal in* $\mathbf{TL}^v_{\mathbb{A}}$*, i.e.*

$$\big(\mathrm{f} \in \mathbf{I}_{<i}, \mathrm{g}, \mathrm{h} \in \mathbf{TL}^v_{\mathbb{A}}\big) \Rightarrow \big(\mathrm{gf}, \mathrm{fh} \in \mathbf{I}_{<i}\big).$$

*Proof.* This is Exercise 10I.1. □

The point is that these "get thinner if we multiply":

**Lemma 10D.15.** *We have*

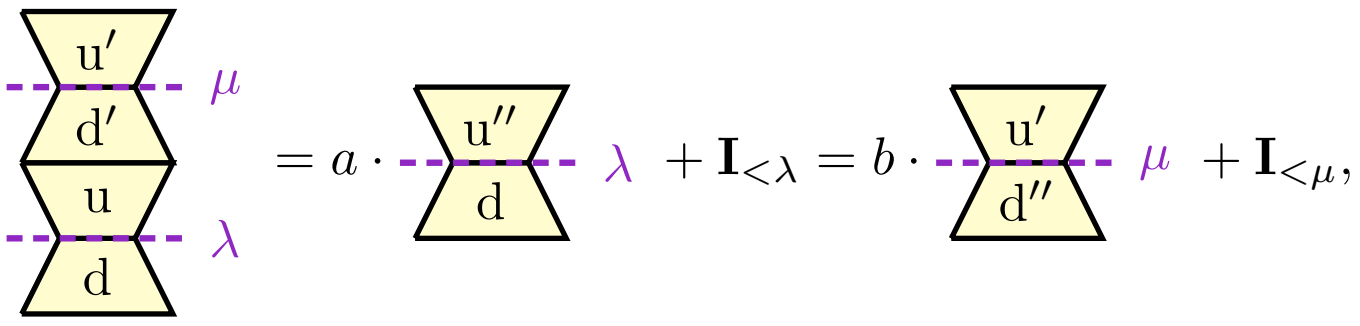

*where the scalars* $a = a(\mathrm{u}, \mathrm{d}'), b = b(\mathrm{u}, \mathrm{d}') \in \mathbb{A}^{\infty}$ *only depend on* u *and* d′.

*Proof.* The integral basis is constructed to get "thinner". Details are supposed to be done in Exercise 10I.3. □

We think of the scalars in Lemma 10D.15 as

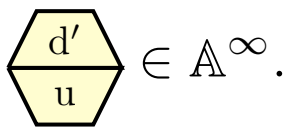

Also, for $\lambda = \mu$ the picture gets a bit easier:

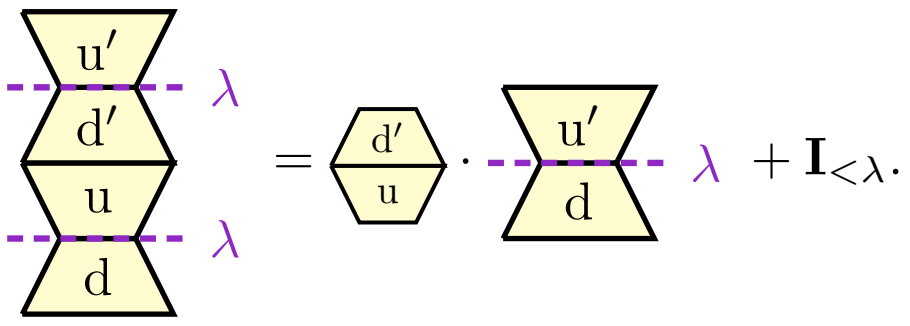

**Example 10D.16.** With the notation as in Example 10D.7 we have *e.g.*

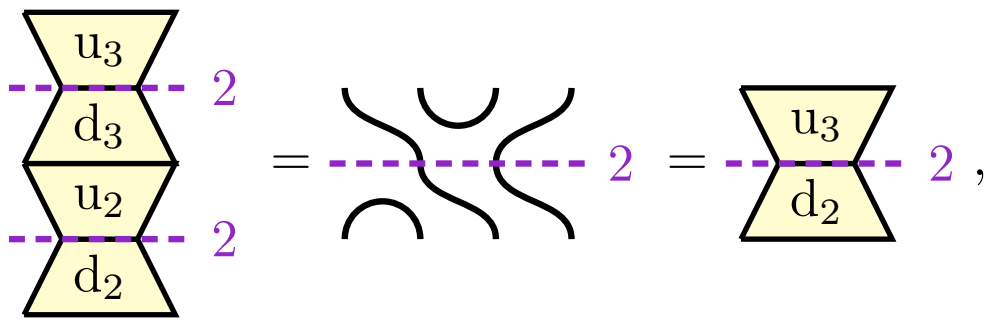

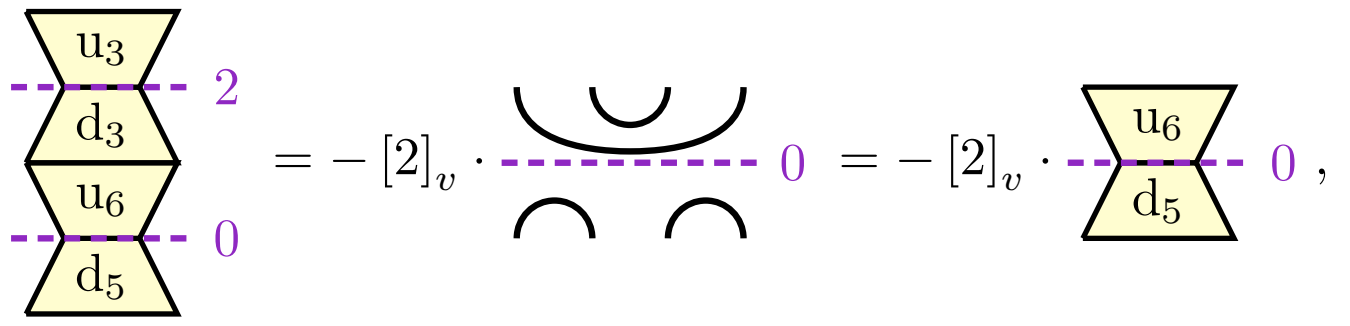

as one easily checks. ◇

10E. **Jones–Wenzl idempotents and their generalizations.** Recall the quantum numbers as in Equation 9H-1. (We will use different subscripts to make it clear whether we work in a specialization or not.) We also need the ***quantum binomials*** for $a \in \mathbb{Z}$ and $b \in \mathbb{Z}_{\geq 0}$. First, by convention, $\begin{bmatrix} a \\ 0 \end{bmatrix}_v = 1$ and otherwise we let

$$\begin{bmatrix} a \\ b \end{bmatrix}_v = \tfrac{[a]_v [a-1]_v \dots [a-b+1]_v}{[b]_v [b-1]_v \dots [1]_v} \in \mathbb{A}.$$

(Note that these are elements of $\mathbb{A}$, as one can check.) Of course, for $q = 1$ the quantum binomial is the usual binomial and for $q = -1$ it is a signed version of the usual binomial. Moreover,

$$\begin{bmatrix} a \\ b \end{bmatrix}_v = 0 \Leftrightarrow a < b,$$

which is however far from being true in specializations.

**Example 10E.1.** Here are some explicit examples for specializations with $p = \mathrm{char}(\mathbb{S})$:

| | $\begin{bmatrix} 8 \\ 0 \end{bmatrix}_q$ | $\begin{bmatrix} 8 \\ 1 \end{bmatrix}_q$ | $\begin{bmatrix} 8 \\ 2 \end{bmatrix}_q$ | $\begin{bmatrix} 8 \\ 3 \end{bmatrix}_q$ | $\begin{bmatrix} 8 \\ 4 \end{bmatrix}_q$ | $\begin{bmatrix} 8 \\ 5 \end{bmatrix}_q$ | $\begin{bmatrix} 8 \\ 6 \end{bmatrix}_q$ | $\begin{bmatrix} 8 \\ 7 \end{bmatrix}_q$ | $\begin{bmatrix} 8 \\ 8 \end{bmatrix}_q$ |
|---|---|---|---|---|---|---|---|---|---|
| $q = 1, p = 0$ | 1 | 8 | 28 | 56 | 70 | 56 | 28 | 8 | 1 |
| $q = -1, p = 0$ | 1 | $-8$ | 28 | $-56$ | 70 | $-56$ | 28 | $-8$ | 1 |
| $q = 2, p = 0$ | 1 | $\frac{21845}{128}$ | $\frac{23859109}{4096}$ | $\frac{1550842085}{32768}$ | $\frac{6221613541}{65536}$ | $\frac{1550842085}{32768}$ | $\frac{23859109}{4096}$ | $\frac{21845}{128}$ | 1 |
| $q = 1, p = 3$ | 1 | 2 | 1 | 2 | 1 | 2 | 1 | 2 | 1 |
| $q = 1, p = 5$ | 1 | 3 | 3 | 1 | 0 | 1 | 3 | 3 | 1 |
| $q = 1, p = 7$ | 1 | 1 | 0 | 0 | 0 | 0 | 0 | 1 | 1 |
| $q = \exp(2\pi i/3), p = 0$ | 1 | $-1$ | 1 | 2 | $-2$ | 2 | 1 | $-1$ | 1 |
| $q = \exp(2\pi i/5), p = 0$ | 1 | $\frac{1}{2}(1-\sqrt{5})$ | $\frac{1}{2}(1-\sqrt{5})$ | 1 | 0 | 1 | $\frac{1}{2}(1-\sqrt{5})$ | $\frac{1}{2}(1-\sqrt{5})$ | 1 |
| $q = \exp(2\pi i/7), p = 0$ | 1 | 1 | 0 | 0 | 0 | 0 | 0 | 1 | 1 |
| $q = 2, p = 13$ | 1 | 4 | 1 | 0 | 0 | 0 | 1 | 4 | 1 |
| $q = 3, p = 13$ | 1 | 12 | 1 | 2 | 11 | 2 | 1 | 12 | 1 |
| $q = 4, p = 13$ | 1 | 1 | 1 | 11 | 11 | 11 | 1 | 1 | 1 |

(Note that the appearing of fractions in the $q = 2$ and $p = 0$ case above is not a contradiction to the claim that $\begin{bmatrix} a \\ b \end{bmatrix}_v \in \mathbb{A}$ since $2^{-k} \in \mathbb{A} \otimes_{\mathbb{A}}^{v=2} \mathbb{S}$.) Let us do two more examples. First $a = 11$:

| | $\begin{bmatrix} 11 \\ 0 \end{bmatrix}_q$ | $\begin{bmatrix} 11 \\ 1 \end{bmatrix}_q$ | $\begin{bmatrix} 11 \\ 2 \end{bmatrix}_q$ | $\begin{bmatrix} 11 \\ 3 \end{bmatrix}_q$ | $\begin{bmatrix} 11 \\ 4 \end{bmatrix}_q$ | $\begin{bmatrix} 11 \\ 5 \end{bmatrix}_q$ | $\begin{bmatrix} 11 \\ 6 \end{bmatrix}_q$ | $\begin{bmatrix} 11 \\ 7 \end{bmatrix}_q$ | $\begin{bmatrix} 11 \\ 8 \end{bmatrix}_q$ | $\begin{bmatrix} 11 \\ 9 \end{bmatrix}_q$ | $\begin{bmatrix} 11 \\ 10 \end{bmatrix}_q$ | $\begin{bmatrix} 11 \\ 11 \end{bmatrix}_q$ |
|---|---|---|---|---|---|---|---|---|---|---|---|---|
| $q = 1, p = 3$ | 1 | 2 | 1 | 0 | 0 | 0 | 0 | 0 | 0 | 1 | 2 | 1 |
| $q = \exp(2\pi i/3), p = 0$ | 1 | $-1$ | 1 | 3 | $-3$ | 3 | 3 | $-3$ | 3 | 1 | $-1$ | 1 |

And finally, $a = 14$:

| | $\begin{bmatrix} 14 \\ 0 \end{bmatrix}_q$ | $\begin{bmatrix} 14 \\ 1 \end{bmatrix}_q$ | $\begin{bmatrix} 14 \\ 2 \end{bmatrix}_q$ | $\begin{bmatrix} 14 \\ 3 \end{bmatrix}_q$ | $\begin{bmatrix} 14 \\ 4 \end{bmatrix}_q$ | $\begin{bmatrix} 14 \\ 5 \end{bmatrix}_q$ | $\begin{bmatrix} 14 \\ 6 \end{bmatrix}_q$ | $\begin{bmatrix} 14 \\ 7 \end{bmatrix}_q$ | $\begin{bmatrix} 14 \\ 8 \end{bmatrix}_q$ | $\begin{bmatrix} 14 \\ 9 \end{bmatrix}_q$ | $\begin{bmatrix} 14 \\ 10 \end{bmatrix}_q$ | $\begin{bmatrix} 14 \\ 11 \end{bmatrix}_q$ | $\begin{bmatrix} 14 \\ 12 \end{bmatrix}_q$ | $\begin{bmatrix} 14 \\ 13 \end{bmatrix}_q$ | $\begin{bmatrix} 14 \\ 14 \end{bmatrix}_q$ |
|---|---|---|---|---|---|---|---|---|---|---|---|---|---|---|---|
| $q = 1, p = 7$ | 1 | 0 | 0 | 0 | 0 | 0 | 0 | 2 | 0 | 0 | 0 | 0 | 0 | 0 | 1 |
| $q = \exp(2\pi i/3), p = 0$ | 1 | $-1$ | 1 | 4 | $-4$ | 4 | 6 | $-6$ | 6 | 4 | $-4$ | 4 | 1 | $-1$ | 1 |
| $q = 2, p = 7$ | 1 | 6 | 1 | 4 | 3 | 4 | 6 | 1 | 6 | 4 | 3 | 4 | 1 | 6 | 1 |

This is best checked with a machine. ◇

For $i \in \mathbb{Z}_{\geq 0}$, let us use the ground ring

$$\mathbb{A}^i = \mathbb{A}\Big[\big(\begin{bmatrix} i \\ j \end{bmatrix}_v\big)^{-1} \mid 0 \leq j \leq i\Big]$$

obtained from $\mathbb{A}$ by formally inverting the quantum binomials. Now we come back to Definition 10E.2, with JW short for Jones–Wenzl:

**Definition 10E.2.** For $i \in \mathbb{Z}_{\geq 0}$ an ***ith JW idempotent*** $\mathrm{e}_i \in \mathrm{End}_{\mathbf{TL}_{\mathbb{A}^i \oplus}}(\bullet^i)$, denoted by,

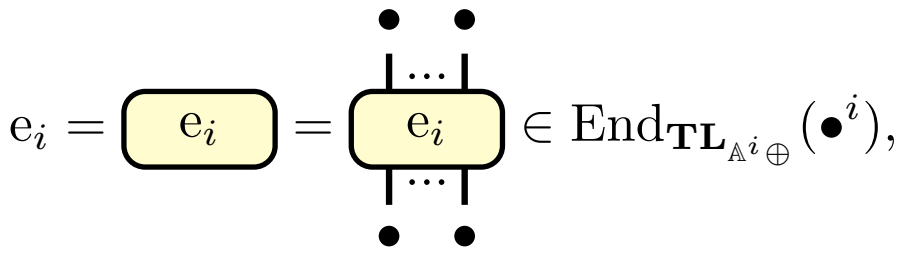

is a morphism satisfying:

- it is an idempotent, *i.e.*

(10E-3) $\mathrm{e}_i^2 = \mathrm{e}_i \rightsquigarrow$ 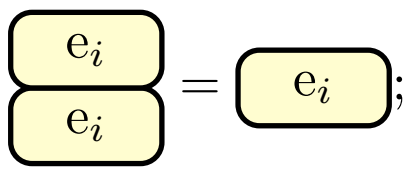

- it annihilates caps and cups, *i.e.*

(10E-4) $\mathbf{I}_{<i}\mathrm{e}_i = 0 = \mathrm{e}_i\mathbf{I}_{<i} \rightsquigarrow$ 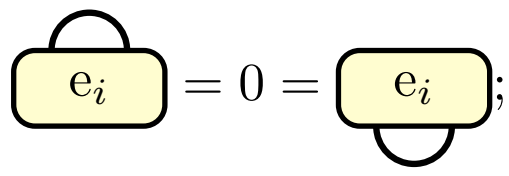

- it contains the identity with coefficient 1, *i.e.*

(10E-5) $(\mathrm{id}_{\bullet^i} - \mathrm{e}_i) \in \mathbf{I}_{<i} \rightsquigarrow$ 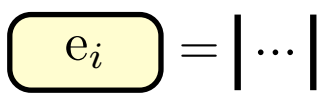 + diagrams with caps and cups.

(That these exists is not immediate from the definition.) ◇

The following is just some algebraic yoga and the crucial point will be the existence of JW idempotents.

**Lemma 10E.6.** *If an ith JW idempotent exists, then it is unique.*

*Proof.* If $\mathrm{e}_i$ and $\mathrm{e}'_i$ are two such idempotents, then Equation 10E-5 implies that $\mathrm{e}_i - \mathrm{e}'_i \in \mathbf{I}_{<i}$. Thus, using the other two defining properties we calculate

$$\mathrm{e}_i - \mathrm{e}_i\mathrm{e}'_i \overset{\text{Equation 10E-3}}{=} \mathrm{e}_i(\mathrm{e}_i - \mathrm{e}'_i) \overset{\text{Equation 10E-4}}{=} 0 \overset{\text{Equation 10E-4}}{=} (\mathrm{e}_i - \mathrm{e}'_i)\mathrm{e}'_i \overset{\text{Equation 10E-3}}{=} \mathrm{e}'_i - \mathrm{e}_i\mathrm{e}'_i,$$

which shows the claim. □

Thus, we will say ***the*** $i$th JW idempotent.

**Proposition 10E.7.** *The ith JW idempotent exists in* $\mathbf{TL}^v_{\mathbb{A}^i}$

*Proof.* We do not know a self-contained proof (*i.e.* using the combinatorics of $\mathbf{TL}^v_{\mathbb{A}^i}$ only) of this fact and refer to [EL17, Theorem A.2]. □

*Remark* 10E.8. As we have seen *e.g.* in Example 10C.4, idempotents tend to have longish expressions. Then same is true for the JW idempotents, see Definition 9H.3 for $\mathrm{e}_2$ and $\mathrm{e}_3$, and the philosophy here would be not to expand them using the recursion from Equation 10E-13 below, but rather the abstract properties. ◇

**Lemma 10E.9.** *JW idempotents satisfy the following.*

*(i) We have* **hom vanishing**, *i.e. for* $0 \leq j \leq i$ *we have*

(10E-10) $$\mathrm{e}_j\mathrm{Hom}_{\mathbf{TL}_{\mathbb{A}^i\oplus}}(\bullet^i, \bullet^j)\mathrm{e}_i = \begin{cases} \mathbb{A}^i\{\mathrm{e}_i\} & \text{if } i{=}j, \\ 0 & \text{else.} \end{cases}$$

*Similarly if* $j \geq i$.

*(ii) We have* **absorption**, *i.e.*

(10E-11) 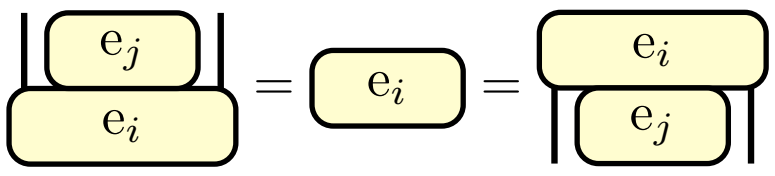 *where* $0 \leq j \leq i$.

*(iii) We have* **partial trace properties**, *i.e.*

(10E-12) 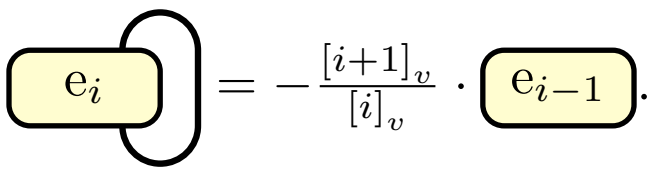

*(iv) The ith JW idempotent satisfies a recursion: First, we have*

$$\mathrm{e}_0 = \varnothing, \quad \mathrm{e}_1 = |.$$

*Then, for* $i \geq 2$, *we have*

(10E-13) 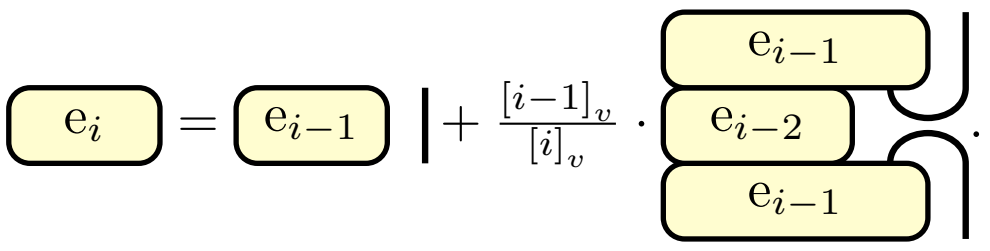

*Proof.* (i). Immediate from the definitions.

(ii). To prove absorption we simply observe that $\mathrm{e}_j = \mathrm{id}_{\bullet^j} + \mathbf{I}_{<j}$, and recall that $\mathrm{e}_i$ annihilates caps and cups Equation 10E-4.

(iii)+(iv). We prove these two claims inductively in tandem. For $i = 0$ or $i = 1$ both claims are clear, so let us suppose that $i \geq 2$ and that (iii) and (iv) hold for all $j < i$. Then define a morphism $\mathrm{e}'_i$ by Equation 10E-13.

Having this expression it is easy to see inductively that the defining properties of an $i$th JW idempotent hold: Equation 10E-5 is clear, while for Equation 10E-4 the crucial calculation is

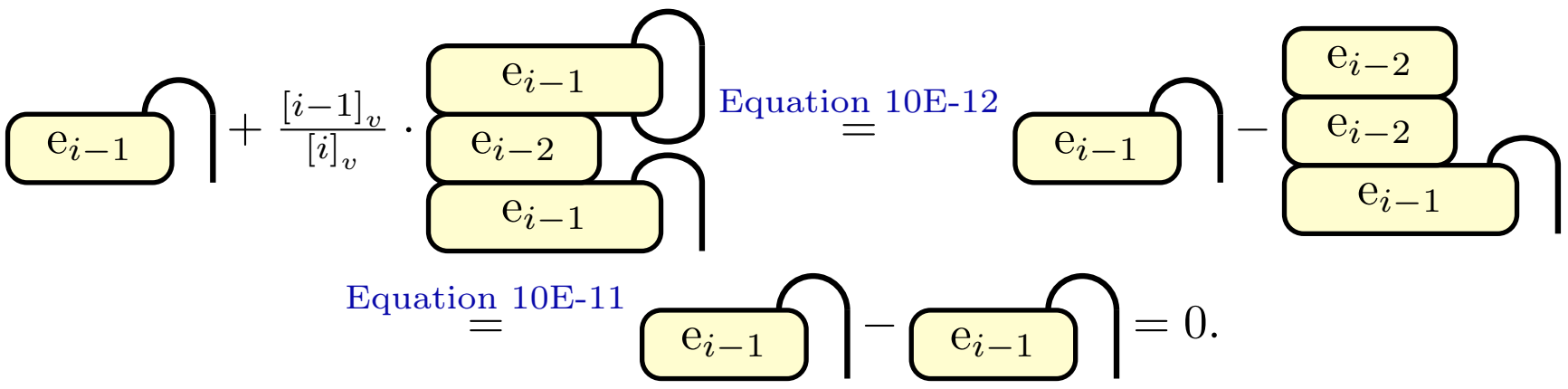

The same calculations shows Equation 10E-3. For partial traces we calculate

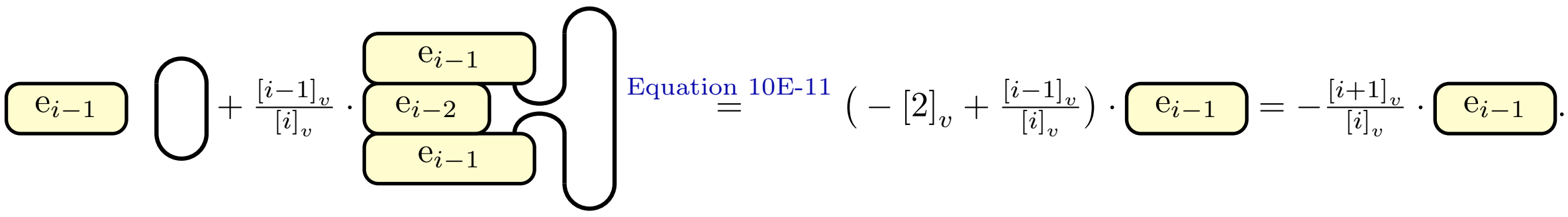

This shows the lemma. □

**Lemma 10E.14.** *For the canonical pivotal structure we have* $\mathrm{tr}^{\mathbf{TL}^{v}_{\mathbb{A}^{i}}}(\mathrm{e}_i) = (-1)^i\,[i+1]_v$.

*Proof.* This is Exercise 10I.3. □

We can also construct a basis using the JW idempotents, which should be compared to the construction of the integral basis from Section 10D (*e.g.* compare Equation 10D-5 and Equation 10E-16). To this end, let us consider the ***generic ground ring***, inverting all quantum binomials,

$$\mathbb{A}^g = \mathbb{A}\Big[\big(\begin{bmatrix} i \\ j \end{bmatrix}_v\big)^{-1} \mid 0 \leq j \leq i, i \in \mathbb{Z}_{\geq 0}\Big],$$

or variations such as $\mathbb{A}^{\leq i}$, having the evident meaning.

**Definition 10E.15.** To $\varepsilon_{\pm 1}$ we associate operators via:

(10E-16) 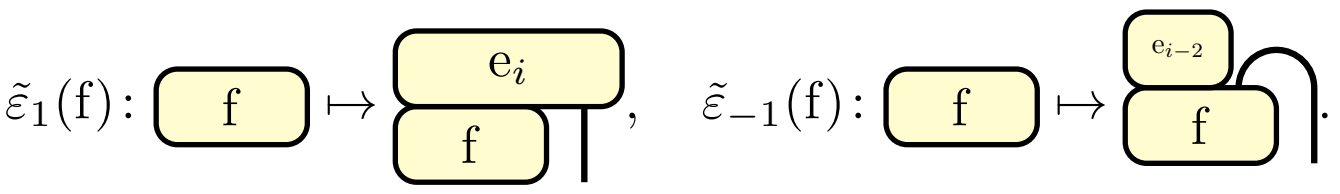

In words, if we already have a morphism $\mathrm{f} \in \mathbf{TL}^{v}_{\mathbb{A}^{\leq i}}$ ending in $\bullet^{i-1}$, then we obtain to new morphism $\mathrm{f}\varepsilon_{\pm 1} \in \mathbf{TL}^{v}_{\mathbb{A}^{\leq i}}$ by either adding a strand or a cap and a JW idempotent. ◇

Copying Section 10D, we obtain ***downward*** $\tilde{\mathrm{d}}(\pi) \in \mathbf{TL}^{v}_{\mathbb{A}^{\leq i}}$ and ***upward Weyl ladders*** $\tilde{\mathrm{u}}(\pi) \in \mathbf{TL}^{v}_{\mathbb{A}^{\leq i}}$, respectively, and also ***Weyl ladders***

$$\tilde{\mathrm{c}}^{\lambda}_{\pi^n_u,\pi^m_d} = \tilde{\mathrm{u}}(\pi^n_u)\tilde{\mathrm{d}}(\pi^m_d) \in \mathrm{Hom}_{\mathbf{TL}^{v}_{\mathbb{A}^{\leq i}}}(\bullet^m, \bullet^n), \quad i = \max\{m, n\},$$

all of which have an associated length *etc.* Not surprisingly, and directly from Proposition 10E.7:

**Proposition 10E.17.** *The morphisms* $\tilde{\mathrm{c}}^{\lambda}_{\pi^n_u,\pi^m_d}$ *exist in* $\mathbf{TL}^{v}_{\mathbb{A}^{\leq i}}$ *for* $i = \max\{m, n\}$. □

**Example 10E.18.** Consider the case of the all non-negative paths $\pi = \varepsilon_1 ... \varepsilon_1$. Then absorption Equation 10E-11 gives inductively

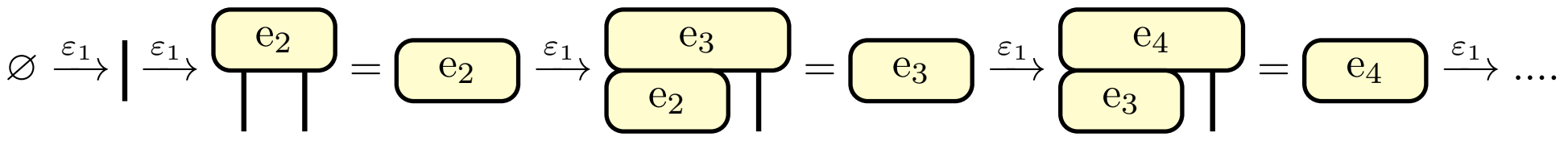

Thus, we get $\tilde{\mathrm{c}}^{\lambda}_{\pi^i,\pi^i} = \tilde{\mathrm{d}}(\pi) = \tilde{\mathrm{u}}(\pi) = \mathrm{e}_i$. ◇

**Example 10E.19.** Let us consider the analog of Example 10D.7, using the same notation. After using absorption we get

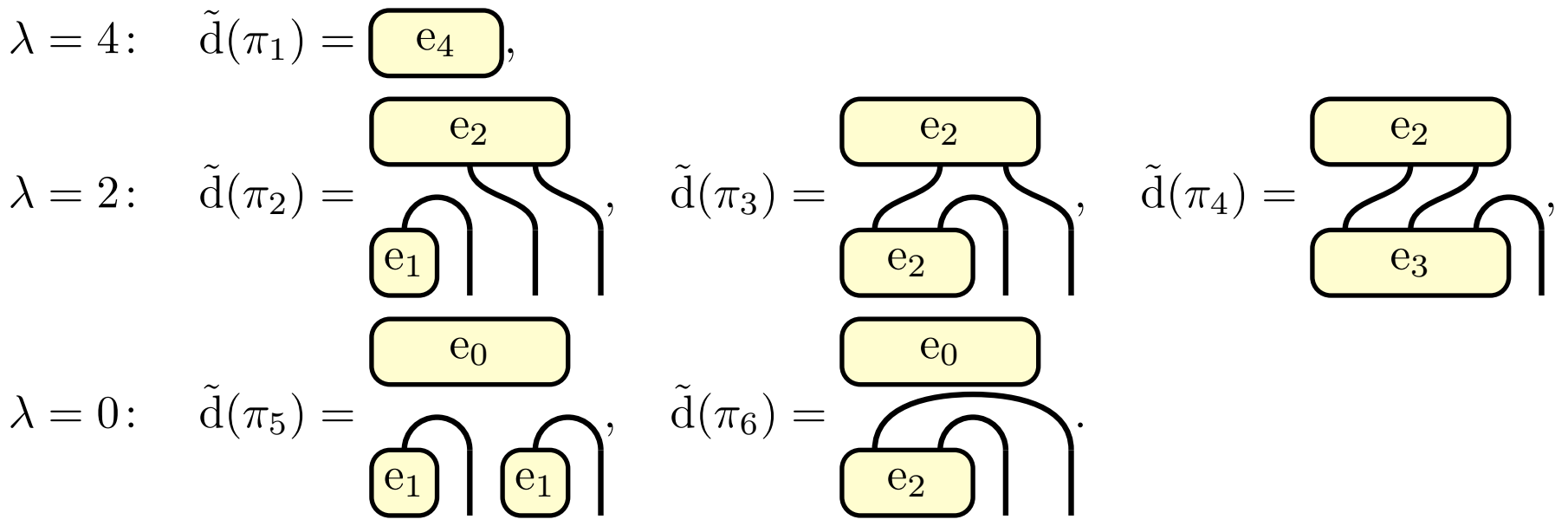

(Of course, the JW idempotents $\mathrm{e}_0$ and $\mathrm{e}_1$ are rather trivial and they are only illustrated to clarify the construction.) ◇

**Theorem 10E.20.** *The sets of the form*

(10E-21) $$\mathrm{WL} = \{\tilde{\mathrm{c}}^{\lambda}_{\pi^n_u,\pi^m_d} \mid \lambda \text{ a weight}, \pi^n_u, \pi^m_d \text{ non-negative paths}\}$$

*are bases of* $\mathrm{Hom}_{\mathbf{TL}^v_{\mathbb{A}^g}}(\bullet^m, \bullet^n)$.

*Proof.* This is (almost) immediate from Theorem 10D.11: Substituting the identity in the coupons of the JW idempotents recovers the integral basis IL. By Equation 10E-5 we thus get an upper triangular change-of-basis matrix between IL and WL. □

Again, we also write *e.g.* $\tilde{\mathrm{c}}^{\lambda}_{\tilde{\mathrm{u}},\tilde{\mathrm{d}}} = \tilde{\mathrm{c}}^{\lambda}_{\pi^n_u,\pi^m_d}$ for simplicity. Moreover, note that these morphisms are constructed using the bottleneck principle as in Equation 10D-13 and we will illustrate these by

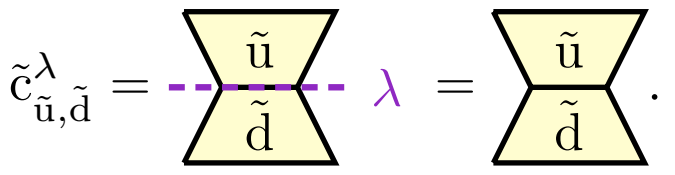

**Lemma 10E.22.** *We have*

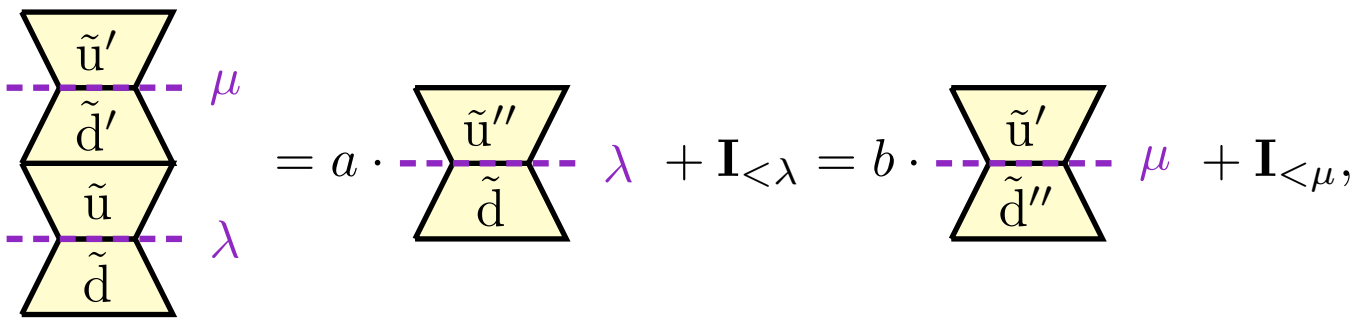

*where the scalars* $a = a(\tilde{\mathrm{u}}, \tilde{\mathrm{d}}'), b = b(\tilde{\mathrm{u}}, \tilde{\mathrm{d}}') \in \mathbb{A}^{\infty}$ *only depend on* $\tilde{\mathrm{u}}$ *and* $\tilde{\mathrm{d}}'$.

*Proof.* This follows using the abstract properties of the JW idempotents. □

**Definition 10E.23.** For $(\lambda, \pi, \pi)$ consisting of a weight and a non-negative path, we define the ***(generalized) JW idempotent*** $\mathrm{e}_\pi$ by

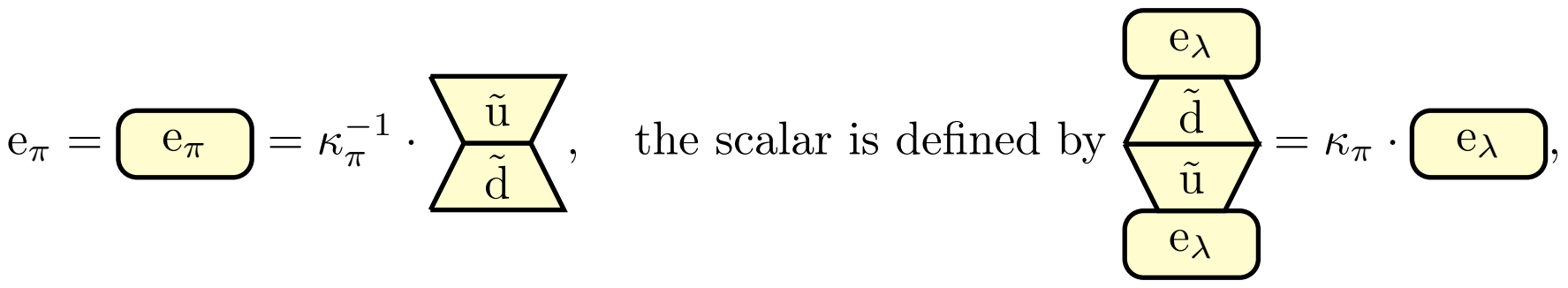

where $\tilde{\mathrm{d}}$ and $\tilde{\mathrm{u}}$ are the downwards and upwards Weyl ladders associated to $\pi$. ◇

**Example 10E.24.** Let us calculate $\kappa^{-1}_{\pi_3}$ for $\pi_3 = \varepsilon_1\varepsilon_1\varepsilon_{-1}\varepsilon_1$:

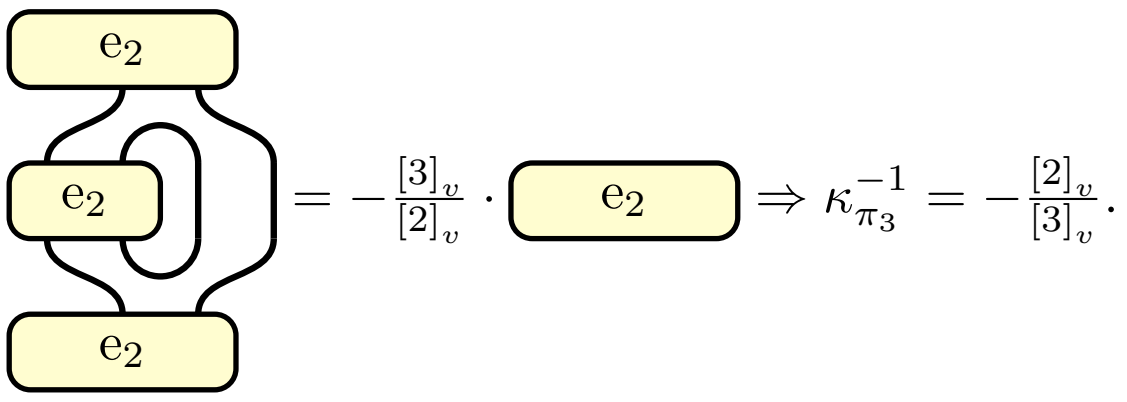

This uses the partial traces Equation 10E-12. ◇

**Theorem 10E.25.** *The generalized JW idempotents are well-defined and the set*

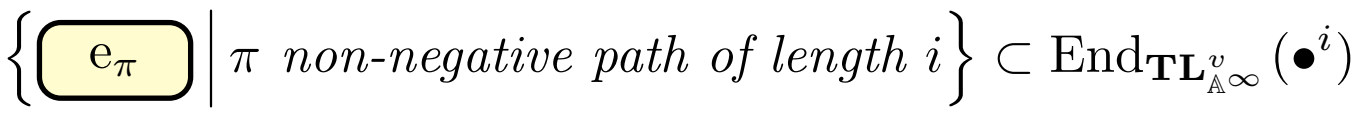

*is a complete set of orthogonal idempotents, i.e.*

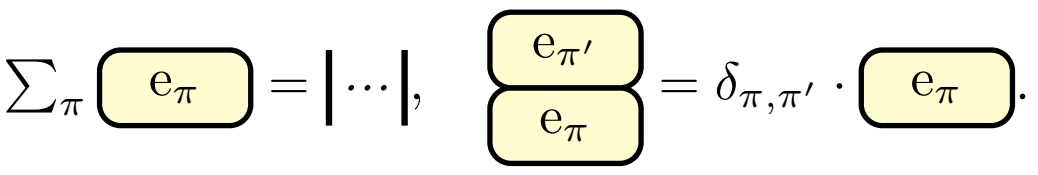

*Proof.* Note first that the scalar $\kappa_\pi$ exists by hom vanishing Equation 10E-10. Moreover, it is easy to see that the scalar $\kappa_\pi$ is an iterative product of partial trace scalars, thus, can be inverted in $\mathbb{A}^{\infty}$. That the generalized JW idempotents are orthogonal idempotents follows from the observation that

(10E-26)
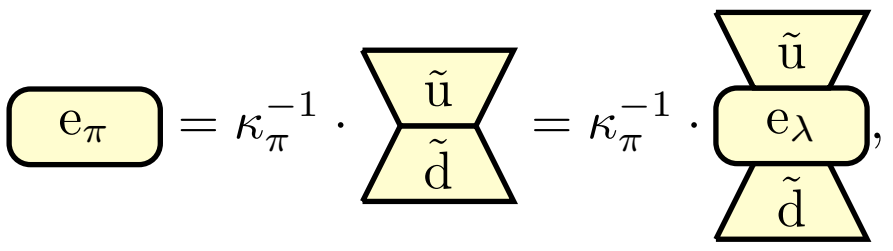

where $\lambda$ is the weight of $\pi$, and the properties of the JW idempotents. Finally, we have

(10E-27)
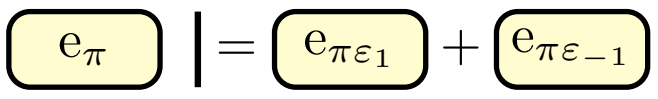

as a consequence of the JW recursion Equation 10E-13, which inductively implies that $\sum_\pi \mathrm{e}_\pi = \mathrm{id}_{\bullet^i}$. □

10F. **The Rumer–Teller–Weyl category – algebra.** Let us further analyze the category $\mathbf{TL}^v_{\mathbb{A}}$ or specializations of it.

**Lemma 10F.1.** *We have the following.*

*(i) We have the decomposition*

$$\bullet^i \cong \textstyle\bigoplus_\pi \mathtt{Im}(\mathrm{e}_\pi) \quad (\textit{in } \mathbf{TL}^v_{\mathbb{A}^\infty \oplus \Subset}).$$

*(ii) The object $\mathtt{Im}(\mathrm{e}_\pi) \in \mathbf{TL}^v_{\mathbb{A}^\infty \oplus \Subset}$ is simple.*

*(iii) We have*

$$\big(\mathtt{Im}(\mathrm{e}_\pi) \cong \mathtt{Im}(\mathrm{e}_{\pi'})\big) \Leftrightarrow \big(\pi \textit{ and } \pi' \textit{ are of the same weight}\big).$$

*(iv) We have*

$$\mathrm{Si}(\mathbf{TL}^v_{\mathbb{A}^\infty \oplus \Subset}) = \mathrm{In}(\mathbf{TL}^v_{\mathbb{A}^\infty \oplus \Subset}) = \big\{\mathtt{Im}(\mathrm{e}_\lambda) \mid \lambda \in \mathbb{Z}_{\geq 0}\big\}.$$

*Proof.* (i)+(ii). These a direct consequences of Theorem 10E.25.

(iii). If $\pi$ and $\pi'$ are not of the same weight, then $\mathrm{Hom}_{\mathbf{TL}^v_{\mathbb{A}^\infty \oplus \Subset}}\big(\mathtt{Im}(\mathrm{e}_\pi), \mathtt{Im}(\mathrm{e}_{\pi'})\big) = 0$ by hom vanishing Equation 10E-10 and Equation 10E-26. Thus, we get $\mathtt{Im}(\mathrm{e}_\pi) \not\cong \mathtt{Im}(\mathrm{e}_{\pi'})$ in this case. For the converse it is enough to consider the case $\mathrm{e}_{\pi'} = \mathrm{e}_\lambda$. Then

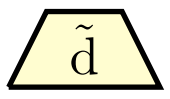
$: \mathtt{Im}(\mathrm{e}_\pi) \to \mathtt{Im}(\mathrm{e}_\lambda)$,
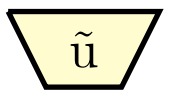
$: \mathtt{Im}(\mathrm{e}_\lambda) \to \mathtt{Im}(\mathrm{e}_\pi)$,

are inverses up to a scalar as Equation 10E-26 shows.

(iv). By (iii) we get that every $\mathtt{Im}(\mathrm{e}_\pi)$ is isomorphic to precisely one $\mathtt{Im}(\mathrm{e}_i)$, while (i) and (ii) show that there are no other simple objects. □

Note that $\mathbf{TL}^q_{\Bbbk \oplus \Subset}$ is always l fiat. Moreover, it is l fusion in exactly the following situation:

**Theorem 10F.2.** *Let $(\Bbbk, q)$ be a specialization. Then $\mathbf{TL}^q_{\Bbbk \oplus \Subset}$ is semisimple if and only if $\mathbb{A}^\infty \subset_{v=q} \Bbbk$ (i.e. all quantum binomials are invertible). Moreover, in the semisimple case we have*

$$\mathrm{Si}(\mathbf{TL}^q_{\Bbbk \oplus \Subset}) = \mathrm{In}(\mathbf{TL}^q_{\Bbbk \oplus \Subset}) = \big\{\mathtt{Im}(\mathrm{e}_\lambda) \mid \lambda \in \mathbb{Z}_{\geq 0}\big\}.$$

*Proof.* If all quantum binomials are invertible, then the specialization $(\Bbbk, q)$ factors through $\mathbb{A}^\infty$ and the claim follows from Lemma 10F.1. On the other hand, if some quantum binomial is not invertible, then there exists some JW idempotent $\mathrm{e}_j$ which is still well-defined, but $\mathrm{e}_{j+1}$ is not. Let $\mathrm{e}_i$ be the minimal such JW idempotent. We claim that $\mathtt{Im}(\mathrm{e}_i \otimes \mathrm{id}_\bullet)$ is indecomposable, but not simple. Indeed,

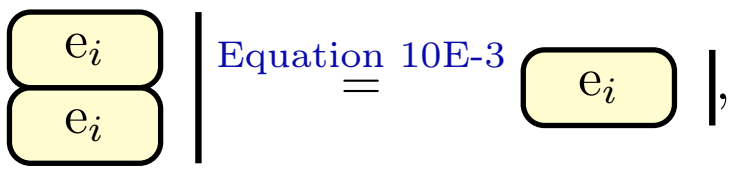

shows that $\mathrm{e}_i \otimes \mathrm{id}_\bullet$ is an idempotent. Moreover, using the standard basis Equation 10D-12 and Equation 10E-10 we have

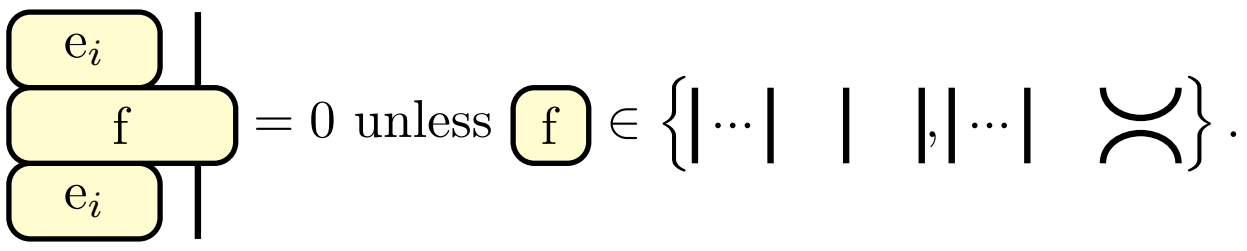

Furthermore, since $\Bbbk$ is a field the quantum binomial $\left[\begin{smallmatrix} i+1 \\ j \end{smallmatrix}\right]_q$ for $0 \leq j \leq i+1$ is only non-invertible if its zero, which, by minimality of $i$, gives $[i+1]_q = 0$. Hence, the calculation shows that the endomorphism ring of $\mathtt{Im}(\mathrm{e}_i \otimes \mathrm{id}_\bullet)$ is $\Bbbk[X]/(X^2)$, which implies that $\mathtt{Im}(\mathrm{id}_\bullet \otimes \mathrm{e}_i)$ is indeed indecomposable, and not simple. □

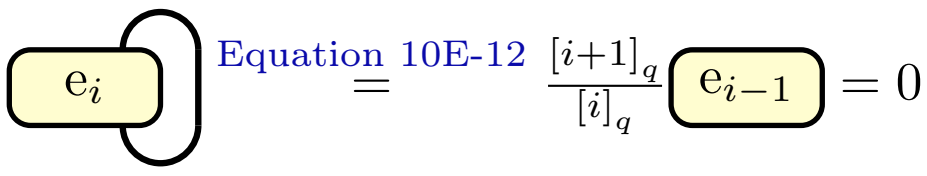

10G. **Some quantum computations.** Let us further study the behavior of quantum numbers. Our main aim is to give "good" conditions for whether the JW idempotents and their generalizations exist, which implies that $\mathbf{TL}^q_{\Bbbk \oplus \Subset}$ is semisimple Theorem 10F.2. For a field $\Bbbk$ this happens if and only if all quantum binomials are non-zero.

**Definition 10G.1.** Define the $q$ ***characteristic*** of a specialization $(\mathbb{S}, q)$ as

$$\mathrm{char}(\mathbb{S}, q) = \min\big\{a \in \mathbb{Z}_{>0} \mid [a]_q = 0\big\},$$

or $\mathrm{char}(\mathbb{S}, q) = \infty$ if $[a]_q \neq 0$ for all $a \in \mathbb{Z}_{>0}$. ◇

There is a well known fact in modular representation theory: it is advantageous to refer to characteristic zero as characteristic infinity, as the behavior in the large prime case is similar to that over the complex numbers. We adopt the same terminology in Definition 10G.1: characteristic $\infty$ is the easiest case.

**Example 10G.2.** For $q = \pm 1$ the $q$ characteristic is the usual characteristic. Here are a few examples of the behavior of the quantum numbers, where $p = \mathrm{char}(\mathbb{S})$.

| | $[0]_q$ | $[1]_q$ | $[2]_q$ | $[3]_q$ | $[4]_q$ | $[5]_q$ | $[6]_q$ | $[7]_q$ | $[8]_q$ |
|---|---|---|---|---|---|---|---|---|---|
| $q = 1, p = 0$ | 0 | 1 | 2 | 3 | 4 | 5 | 6 | 7 | 8 |
| $q = -1, p = 0$ | 0 | 1 | $-2$ | 3 | $-4$ | 5 | $-6$ | 7 | $-8$ |
| $q = 2, p = 0$ | 0 | 1 | $\frac{5}{2}$ | $\frac{21}{4}$ | $\frac{85}{8}$ | $\frac{341}{16}$ | $\frac{1365}{32}$ | $\frac{5461}{64}$ | $\frac{21845}{128}$ |
| $q = 1, p = 3$ | 0 | 1 | 2 | 0 | 1 | 2 | 0 | 1 | 2 |
| $q = 1, p = 5$ | 0 | 1 | 2 | 3 | 4 | 0 | 1 | 2 | 3 |
| $q = 1, p = 7$ | 0 | 1 | 2 | 3 | 4 | 5 | 6 | 0 | 1 |
| $q = 2, p = 13$ | 0 | 1 | 9 | 2 | 9 | 1 | 0 | 12 | 4 |
| $q = 3, p = 13$ | 0 | 1 | 12 | 0 | 1 | 12 | 0 | 1 | 12 |
| $q = 4, p = 13$ | 0 | 1 | 1 | 0 | 12 | 12 | 0 | 1 | 1 |

(The complex root of unity case was already discussed in Example 9H.2, so it is omitted from the above table.) Thus, we have for example $\mathrm{char}(\mathbb{F}_{13}, 2) = 6$. ◇

**Lemma 10G.3.** *Let $a \in \mathbb{Z}_{\neq 0}$. Then:*

*(i) If $\mathrm{char}(\mathbb{S}, q) = 0$, then $[a]_q \in \mathbb{S}$ is non-zero.*

*(ii) If $\mathrm{char}(\mathbb{S}, q) = p > 0$ and $q = \pm 1$, then $[a]_q \in \mathbb{S}$ is non-zero if and only if $p \nmid a$.*

*(iii) If $\mathrm{char}(\mathbb{S}, q) > 0$ and $q \neq \pm 1$, then $[a]_q \in \mathbb{S}$ is non-zero if and only if $q^{2a} \neq 1$.*

*Proof.* If $q = \pm 1$, then we have $[a]_q = \pm a$ and the claims are clear. Note further that $[a]_q = q^a \frac{1-q^{2a}}{q-q^{-1}}$ in case $q \neq \pm 1$. Thus, the roots of $[a]_q$ are exactly the roots of the cyclotomic polynomial $1 - q^{2a}$, which proves the root of unity case. □

For any $c \in \mathbb{Z}_{\geq 0}$ and any $d \in \mathbb{Z}_{\geq 0}$ we use the digits $c_k$ of its $d$-adic expansion:

$$c = [\dots, c_2, c_1, c_0]_d = \textstyle\sum_{k=0}^{\infty} c_k d^k \text{ where } c_k \in \{0, \dots, d-1\}.$$

We also write $*$ for an arbitrary digit.

**Example 10G.4.** Let $q = 1$, so that the $q$ characteristic is the usual characteristic, which we assume to be two. In this case $c = [\dots, c_2, c_1, c_0]_2$ is the binary expansion of $c$, for example, $0 = [0]_2$ $1 = [1]_2$, $2 = [1, 0]_2$, $3 = [1, 1]_2$ *etc.* The ***Lucas' theorem*** (which is valid for any prime) states that

$$\binom{a}{b} = \textstyle\prod_{i=0}^{\infty} \binom{a_i}{b_i},$$

where $a_i$ and $b_i$ are the digits of $a, b$ in binary expansion. In particular,

$$\binom{a}{b} = \begin{cases} 1 & \text{if all } a_i \geq b_i, \\ 0 & \text{otherwise.} \end{cases} \tag{10G-5}$$

Ordering the binomial into a triangle à la Pascal gives ***Sierpinski's triangle*** as in Figure 28. ◇

**Lemma 10G.6.** *Let $a \in \mathbb{Z}_{\geq 0}$. Then:*

*(i) If $\mathrm{char}(\mathbb{S}, q) < a$, then $\begin{bmatrix} a \\ b \end{bmatrix}_q \in \mathbb{S}$ is non-zero for all $0 \leq b \leq a$.*

*(ii) If $\mathrm{char}(\mathbb{S}, q) = \mathrm{char}(\mathbb{S}) = p \geq a$, then $\begin{bmatrix} a \\ b \end{bmatrix}_q \in \mathbb{S}$ is non-zero for all $0 \leq b \leq a$ if and only if*

$$a = [\dots, 0, 0, *, p-1, \dots, p-1]_p.$$

*(iii) If $\mathrm{char}(\mathbb{S}, q) = k \geq a$ and $\mathrm{char}(\mathbb{S}) = 0$, then $\begin{bmatrix} a \\ b \end{bmatrix}_q \in \mathbb{S}$ is non-zero for all $0 \leq b \leq a$ if and only if*

$$a = [\dots, 0, 0, *, \dots, *, k-1]_k.$$

*(iv) If $\mathrm{char}(\mathbb{S}, q) = k \geq a$, $\mathrm{char}(\mathbb{S}) = p > 0$ and $k \neq p$, then $\begin{bmatrix} a \\ b \end{bmatrix}_q \in \mathbb{S}$ is non-zero for all $0 \leq b \leq a$ if and only if*

$$a = [\dots, 0, 0, *, \dots, *, k-1]_k, \quad m = [\dots, 0, 0, *, p-1, \dots, p-1]_p,$$

*where $a = mk + \tilde{a}_0$ for $0 \leq \tilde{a}_0 < k$.*

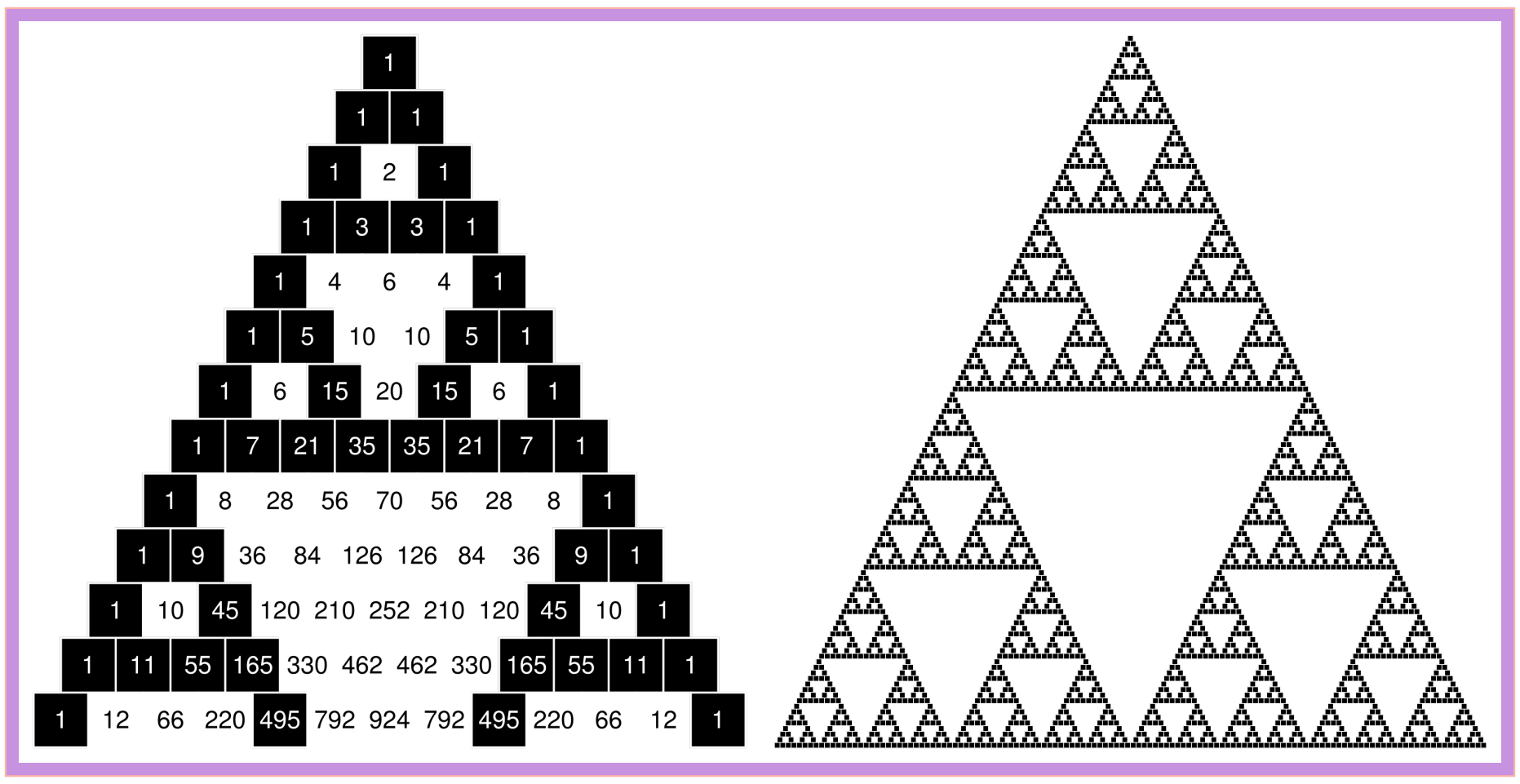

FIGURE 28. Sierpinski's triangle is Pascal's triangle modulo 2, and one gets a fractal type pattern. The fully black rows correspond to $a = 2^k - 1$, since *e.g.* $7 = [1,1,1]_2$ so that the first case in Equation 10G-5 always applies.

Pictures from https://www.youtube.com/watch?v=vWV3yeBD-Gw

*Proof.* In case $\mathrm{char}(\mathbb{S}, q) < a$ the claim is clear, so let us assume that $\mathrm{char}(\mathbb{S}, q) = k \geq a$ and write $a = mk + \tilde{a}_0$ and $b = nk + \tilde{b}_0$ with $0 \leq \tilde{a}_0, \tilde{b}_0 < k$. Recall that then the ***quantum Lucas' theorem*** states that

$$\begin{bmatrix} a \\ b \end{bmatrix}_q = \binom{m}{n} \begin{bmatrix} \tilde{a}_0 \\ \tilde{b}_0 \end{bmatrix}_q, \tag{10G-7}$$

see *e.g.* [Lus10, Lemma 24.1.2]. Note the appearance of the usual binomial: If $\mathrm{char}(\mathbb{S}) = 0$, then this factor is always non-zero and we get the case (iii) in the statement. If however $\mathrm{char}(\mathbb{S}) = p > 0$, then we can apply the classical Lucas' theorem to Equation 10G-7 and get

$$\begin{bmatrix} a \\ b \end{bmatrix}_q = \Big( \textstyle\prod_{i=0}^{\infty} \binom{m_i}{n_i} \Big) \begin{bmatrix} \tilde{a}_0 \\ \tilde{b}_0 \end{bmatrix}_q \overset{p=l}{=} \textstyle\prod_{i=0}^{\infty} \binom{a_i}{b_i} \begin{bmatrix} a_0 \\ b_0 \end{bmatrix}_q .$$

where we distinguish expansion in base $l$ and $p$. □

**Example 10G.8.** If we want to know whether $\begin{bmatrix} a \\ b \end{bmatrix}_q$ is non-zero for all $0 \leq b \leq a$, as in Example 10E.1 for $a = 8$, $a = 11$ or $a = 14$, then:

- Generically this is always the case.
- In char $p$ we would write *e.g.* $8 = [2,2]_3 = [1,3]_5 = [1,1]_7$ which implies that in characteristic 3 the quantum binomials will be non-zero, but it will eventually be zero in characteristic 5 or 7. For $11 = [1,0,2]_3$ the quantum binomial might be zero.
- In the complex root of unity case only the zeroth digit plays a role. For example for $11 = [1,0,2]_3$ and $q = \exp(2\pi i/3)$ the quantum binomial will always be non-zero.
- The mixed case is a mixture of the above two cases. For example, for $q^3 = 1$ and $p = 7$ one would need to expand $a$ in base 3, where only the zeroth digit is important, and then $m$ in base 7, *cf.* $14 = [1,1,2]_3$ and $4 = [4]_7$ in Example 10E.1. Another example where all quantum binomials for $q^3 = 1$ and $p = 7$ are invertible is $146 = 3 \cdot 48 + 2$ since $146 = [1,2,1,0,2]_3$ and $48 = [6,6]_7$.

Its all about vanishing of binomials. ◇

For Lemma 10G.6 and Theorem 10F.2 we immediately get:

**Theorem 10G.9.** *Let $(\Bbbk, q)$ be a specialization. Then $\mathbf{TL}^q_{\Bbbk\oplus\mathbb{C}}$ is semisimple if and only if we are in the generic case. Moreover, in this case we have*

$$\mathrm{Si}(\mathbf{TL}^q_{\Bbbk\oplus\mathbb{C}}) = \mathrm{In}(\mathbf{TL}^q_{\Bbbk\oplus\mathbb{C}}) = \big\{ \mathtt{Im}(\mathrm{e}_\lambda) \mid \lambda \in \mathbb{Z}_{\geq 0} \big\},$$

*as the set of simple objects.* □

10H. **Constructing Verlinde categories.** Let us now finish this section by constructing quantum invariants from (specializations of) $\mathbf{TL}^v_{\mathbb{A}\oplus\mathbb{C}}$.

**Definition 10H.1.** Let $\mathbf{C} \in \mathbf{MCat}$. We call a collection of subspaces

$$\mathbf{I}_{\otimes} = \{ \mathrm{In}(\mathtt{X}, \mathtt{Y}) \subset \mathrm{Hom}_{\mathbf{C}}(\mathtt{X}, \mathtt{Y}) \mid \mathtt{X}, \mathtt{Y} \in \mathbf{C} \}$$

a ***(two-sided) $\otimes$ ideal*** if

- it is closed under vertical composition, *i.e.*

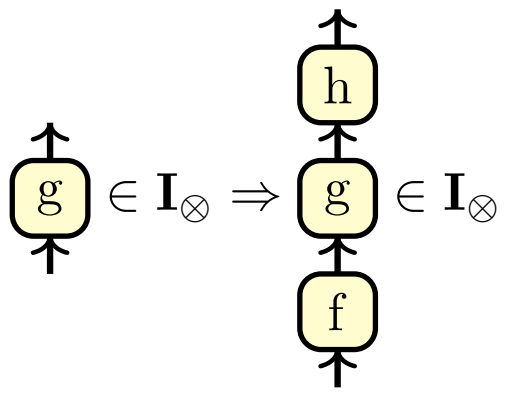

where $\mathrm{f}, \mathrm{h} \in \mathbf{C}$;

- it is closed under horizontal composition, *i.e.*

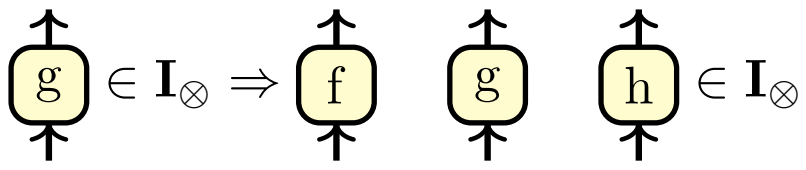

where $\mathrm{f}, \mathrm{h} \in \mathbf{C}$.

(Sometimes such ideals are called $\circ$-$\otimes$ ideals in the literature.) ◇

**Proposition 10H.2.** *Let $\mathbf{C} \in \mathbf{MCat}_{\mathbb{S}}$ and let $\mathbf{I}_{\otimes}$ be a $\otimes$ ideal. Then:*

*(i) The exists a category $\mathbf{C}/\mathbf{I}_{\otimes}$ with*

$$\mathrm{Ob}(\mathbf{C}/\mathbf{I}_{\otimes}) = \mathrm{Ob}(\mathbf{C}), \quad \mathrm{Hom}_{\mathbf{C}/\mathbf{I}_{\otimes}}(\mathtt{X}, \mathtt{Y}) = \mathrm{Hom}_{\mathbf{C}}(\mathtt{X}, \mathtt{Y})/\mathrm{In}(\mathtt{X}, \mathtt{Y}),$$

*and the evident composition.*

*(ii) We have $\mathbf{C}/\mathbf{I}_{\otimes} \in \mathbf{MCat}$ and the identity map on objects and morphisms induces a monoidal and full functor $\mathbf{C} \twoheadrightarrow \mathbf{C}/\mathbf{I}_{\otimes}$.*

*(iii) If $\mathbf{C}$ was braided (or rigid, pivotal, spherical, ribbon), then so is $\mathbf{C}/\mathbf{I}_{\otimes}$.*

*Proof.* This is Exercise 10I.4. □

**Example 10H.3.** Clearly, any $\mathbf{C} \in \mathbf{MCat}_{\mathbb{S}}$ has a trivial $\otimes$ ideal, namely $\mathbf{I}_{\otimes} = \mathbf{C}$. We stress this because of the confusing fact that $\mathbf{C}/\mathbf{C}$ is trivial although $\mathrm{Ob}(\mathbf{C}/\mathbf{C}) = \mathrm{Ob}(\mathbf{C})$. The point is that $\mathrm{Hom}_{\mathbf{C}/\mathbf{C}}(\mathtt{X}, \mathtt{Y}) \cong \{0\}$, and thus all objects are isomorphic. ◇

**Definition 10H.4.** Let $\mathbf{C} \in \mathbf{PCat}_{\mathbb{S}}$. A morphism $\mathrm{f} \in \mathbf{C}$ is called ***right negligible*** if

$$\mathrm{tr}^{\mathbf{C}}(\mathrm{gf}) = 0 \text{ for all } \mathrm{g} \in \mathbf{C},$$

and ***left negligible*** if

$${}^{\mathbf{C}}\mathrm{tr}(\mathrm{gf}) = 0 \text{ for all } \mathrm{g} \in \mathbf{C}.$$

A right and left negligible is called ***negligible***. ◇

For $\mathbf{C} \in \mathbf{PCat}_{\mathbb{S}}$ let $\mathbf{N}_{\mathbf{C}}$ denote the collection of negligible morphisms.

**Proposition 10H.5.** *For any $\mathbf{C} \in \mathbf{PCat}_{\mathbb{S}}$ collection $\mathbf{N}_{\mathbf{C}}$ is a $\otimes$ ideal.*

*Proof.* By definition, the vertical composition of a negligible morphism with any other morphism is negligible. Moreover, up to symmetry,

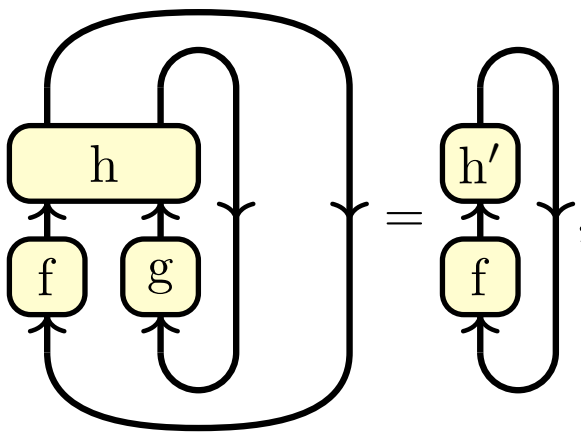

, where $\mathrm{h}' =$

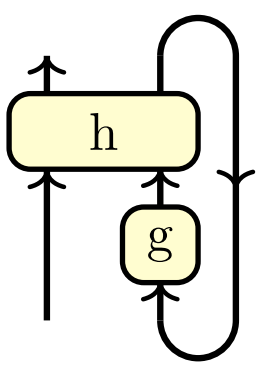

shows the same for the horizontal composition. □

**Definition 10H.6.** Fix the canonical pivotal structure on $\mathbf{TL}_{\mathbb{A}\oplus\Subset}$. For any specialization $(\mathbb{S}, q)$ we call

$$\mathbf{Ver}(\mathbb{S}, q) = \mathbf{TL}^{q}_{\mathbb{S}\oplus\Subset}/\mathbf{N}_{\mathbf{TL}^{q}_{\mathbb{S}\oplus\Subset}}$$

the ***Verlinde category for*** $(\mathbb{S}, q)$. ◇

**Example 10H.7.** In the generic case $\mathbf{Ver}(\Bbbk, q) = \mathbf{TL}^{q}_{\Bbbk\oplus\Subset}$, since $\mathbf{N}_{\mathbf{TL}^{q}_{\Bbbk\oplus\Subset}} = 0$. This follows because we know that the (images of the) JW idempotents are the simple objects in this semisimple category, see Theorem 10G.9, and their traces are non-zero by Lemma 10E.14. ◇

By Proposition 10H.2.(iii) and additivity of categorical traces we immediately see that $\mathbf{Ver}(\mathbb{S}, q) \in \mathbf{lRifiat}$. We get a bit more:

**Proposition 10H.8.** *For any specialization* $(\Bbbk, q)$ *the category* $\mathbf{Ver}(\Bbbk, q)$ *is semisimple, i.e.* $\mathbf{Ver}(\Bbbk, q) \in$ **lRiMo**. *Furthermore, the simple objects of* $\mathbf{Ver}(\Bbbk, q)$ *are the indecomposable objects of* $\mathbf{TL}^{q}_{\Bbbk\oplus\Subset}$ *of non-zero categorical dimension.*

*Proof.* This follows from the following characterization of negligible morphisms. A morphism

$$\mathrm{f} = (f_{i,j})\colon \bigoplus_{i=1}^{k} \mathtt{Z}_i \to \bigoplus_{j=1}^{l} \mathtt{Z}_j$$

between indecomposable objects of $\mathbf{Ver}(\Bbbk, q)$ is negligible if and only if for each $i, j$ either $f_{i,j}$ is not an isomorphism or $\dim^{\mathbf{Ver}(\Bbbk,q)}(\mathtt{Z}_j) = 0$. (This is well-known, see *e.g.* [**EO22**, Lemma 2.2].) □

The Verlinde categories are sometimes even modular, *e.g.* in the complex root of unity case. In general:

The quantum invariants arising from $\mathbf{Ver}(\Bbbk, q)$ are generalized Jones polynomials.

This gives a completely diagrammatic construction of the Jones-type quantum invariants, *i.e.* by coloring strands with (versions of) JW idempotents. Similarly one can construct type BCD versions of these invariants using quantum Brauer categories, or using webs. We explore this in the next section.

## 10I. **Exercises.**

*Exercise* 10I.1. Prove Lemma 10D.14. ◇

*Exercise* 10I.2. Prove Lemma 10D.15. Also try to think what changes in the proof compared to Lemma 10E.22. ◇

*Exercise* 10I.3. Prove Lemma 10E.14. ◇

*Exercise* 10I.4. Prove Proposition 10H.2. ◇

*Exercise* 10I.5. Compute the following quantum invariant.

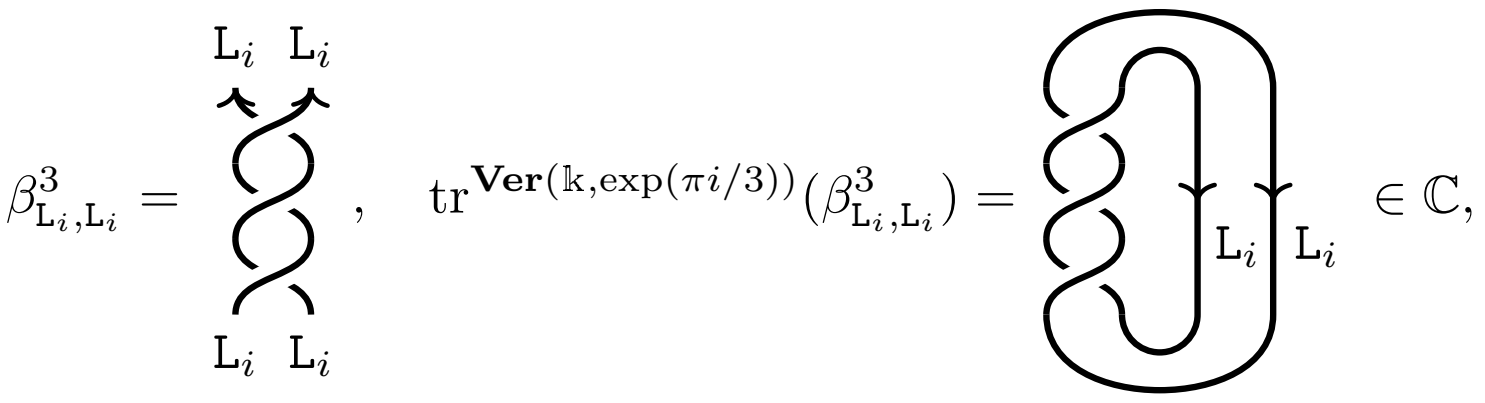

for $i = 0, 1, 2$. (This is the colored Jones polynomial of the trefoil knot.) ◇

# 11. Quantum invariants – a web approach

In Section 10 we have seen quantum invariants constructed from the Rumer–Teller–Weyl category, yielding the Jones polynomial. We have seen that this is $\mathfrak{sl}_2$ representation theory. So we ask:

How to construct quantum invariants beyond the Jones polynomial?

It turns out that we already have a good language at hand: webs. These are diagrammatic objects that originate in work of Tait on the four color theorem, *cf.* Figure 29.

We go into some more details below.

## 11A. **A word about conventions.** Here are our main conventions:

*Convention* 11A.1. We fix $n \in \mathbb{Z}_{\geq 1}$, and a generic pair $(\Bbbk, q)$ (in the sense of Section 10C) such that $q^{\frac{1}{n}} \in \Bbbk$. (We do not really need the $n$th root of $q$ until Section 11G.) Also, let $\mathbb{A} = \mathbb{Z}[q^{\frac{1}{n}}, q^{-\frac{1}{n}}]$ (which is not quite the same as in Convention 10A.1). ◇

*Convention* 11A.2. Since we have a lot of quantum coefficients (such as factorials or binomials) in this section, and they are all evaluated at $q$, we will write e.g. $[k]$ instead of $[k]_q$ to not overload the notation. ◇

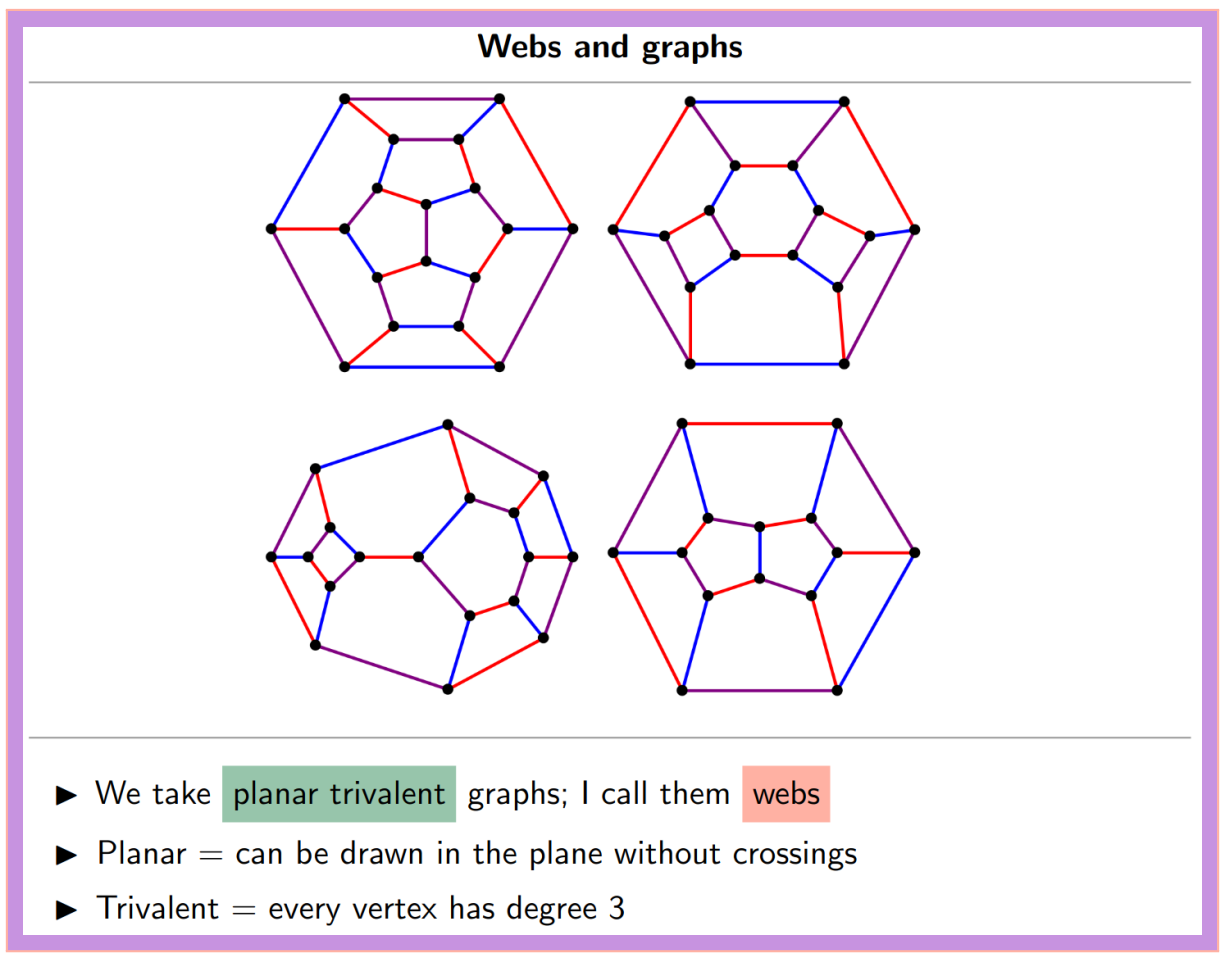

FIGURE 29. Webs are planar trivalent graphs. The name comes from the fact that they tend to look like spider webs.

Picture from https://www.youtube.com/watch?v=-idEadSMnec

11B. **Symmetric GLn webs.** Similarly to Example 3G.2 we define the following. Recall that we think of the generators of **Web** as being bilinear and trilinear forms. Or, alternatively, as in Equation 3G-1 we can rotate these pictures and think about labeled trivalent vertices instead. Now let us add labels and an orientation to the edges:

**Definition 11B.1.** The ***(upward oriented labeled) generic web category without relations*** $\mathbf{lWeb}^{\uparrow}$ is defined as follows. It is monoidally generated by

$$\mathrm{S}: \{k^{\uparrow} \mid k \in \mathbb{Z}_{\geq 0}\},$$

T : { 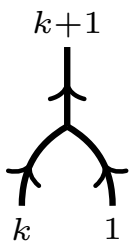 $: k^{\uparrow} \otimes 1^{\uparrow} \to (k+1)^{\uparrow},$ 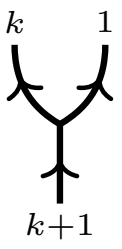 $: (k+1)^{\uparrow} \to k^{\uparrow} \otimes 1^{\uparrow} \mid k \in \mathbb{Z}_{\geq 0}$ },

and the morphism generators are called ***thin merges*** (or $(k,1)$-merges) and ***thin splits*** (or $(k,1)$-splits), respectively. ◇

Here and throughout and as before, we will use the conventions of [LT21] for reading the diagrams, meaning that $v \circ u$ is obtained by gluing $v$ on top of $u$ and $u \otimes v$ is obtained by putting $v$ to the right of $u$, for $u, v \in \mathrm{Hom}_{\mathbf{lWeb}^{\uparrow}}(\vec{k}, \vec{l})$. For example,

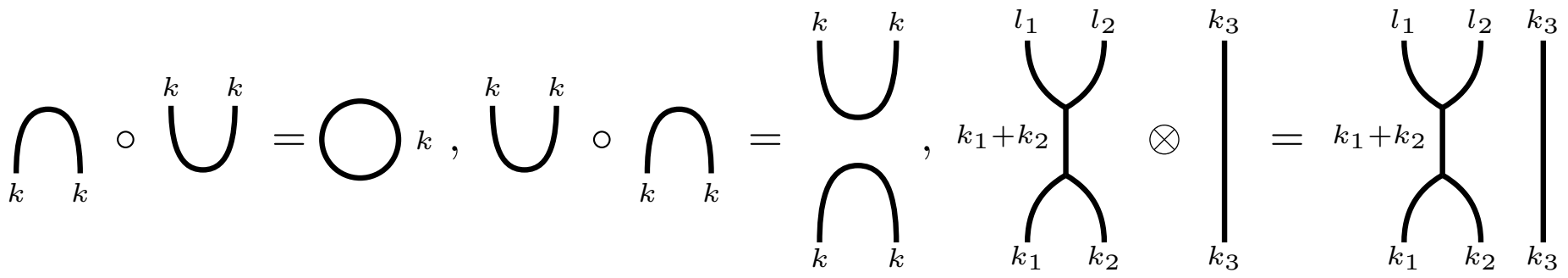

where in the final equation $k_1 + k_2 = l_1 + l_2$. The cup and cap diagrams will be defined later on. Moreover, here we have omitted the arrows, meaning that these conventions hold for all orientations (we elaborate on orientations below).

The next few definitions are morphisms in $\mathbf{lWeb}^{\uparrow}_{\Bbbk\oplus}$.

**Definition 11B.2.** The ***thin overcrossing*** and ***thin undercrossing*** are defined respectively as

(11B-3)

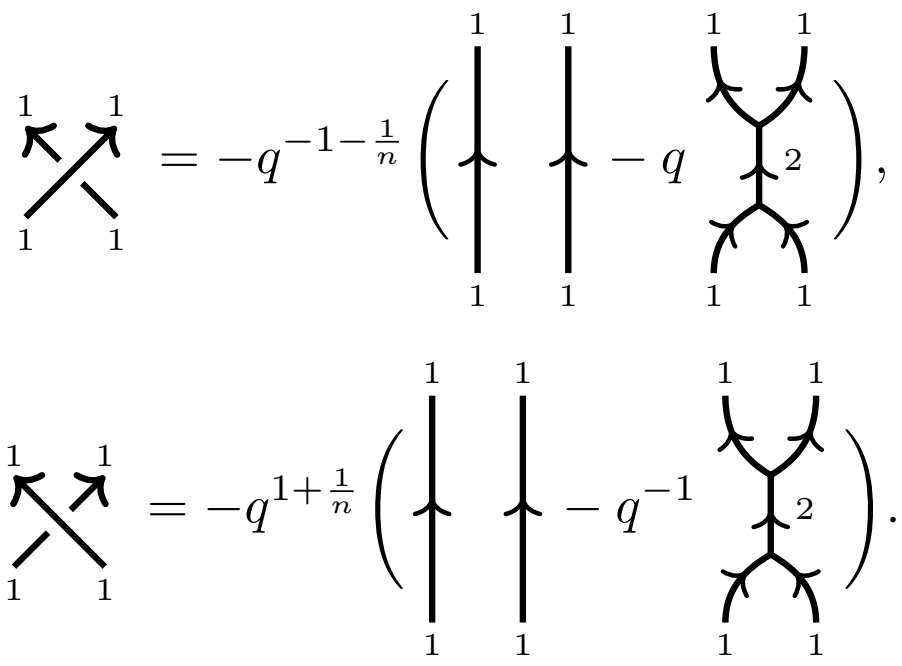

(For now we do not need the $\frac{1}{n}$, and it can be ignored. We will only need it much later and decided to keep it around.) ◇

**Definition 11B.4.** The ***thick overcrossing of type*** $(k, l)$ is defined by "exploding" edges:

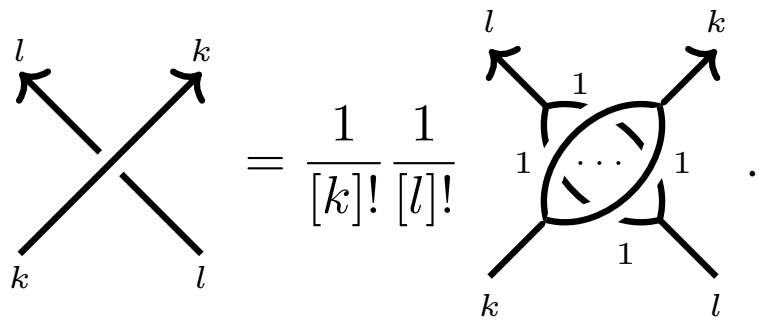

Similarly we define the thick undercrossing. By definition, if $k = 0$ or $l = 0$, then the crossing is the corresponding identity (=empty diagram). ◇

**Definition 11B.5.** We define the $(1, k)$***-merges*** as

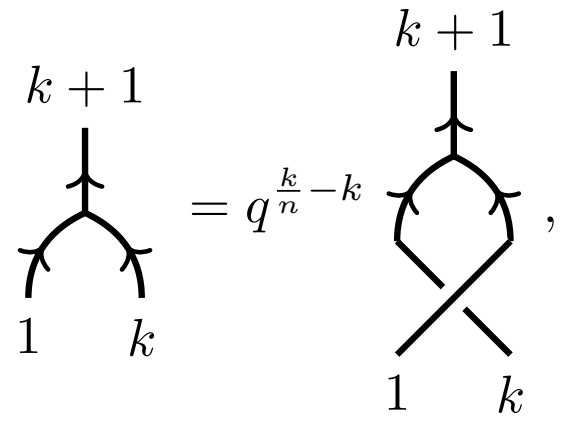

and similarly the $(1, k)$-splits. ◇

Further, we will also need the following.

**Definition 11B.6.** The ***thick merge*** of type $(k, l)$ (or the $(k, l)$-merge) is defined as below

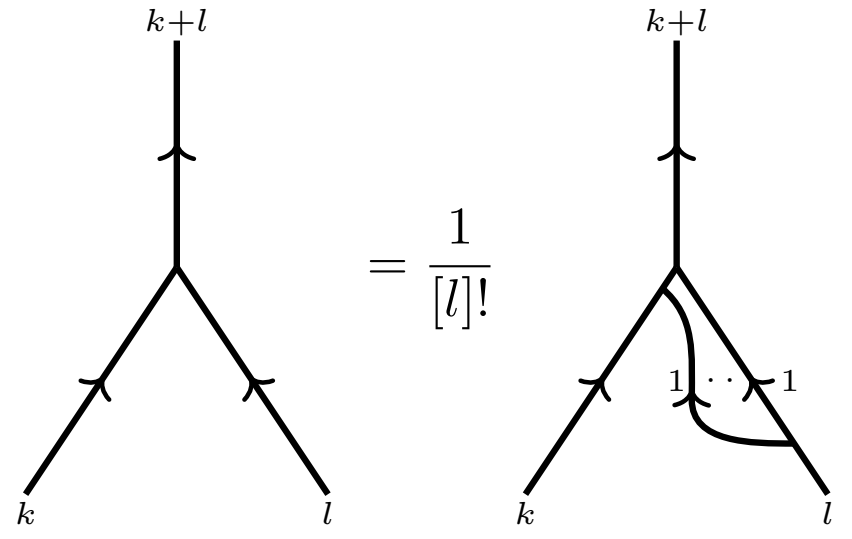

and similarly the $(k, l)$***-split***. ◇

The meticulous reader observes that Definition 11B.4 and Definition 11B.6 are not well-defined since we have not specified the precise form of the webs. However, in Remark 11B.12 we will explain that this is not a problem, so the concerned reader might choose any form of these webs. Similarly for the other web categories that we define later.

**Definition 11B.7.** The $F^{(j)}$ ***and*** $E^{(j)}$***-ladders*** are given by

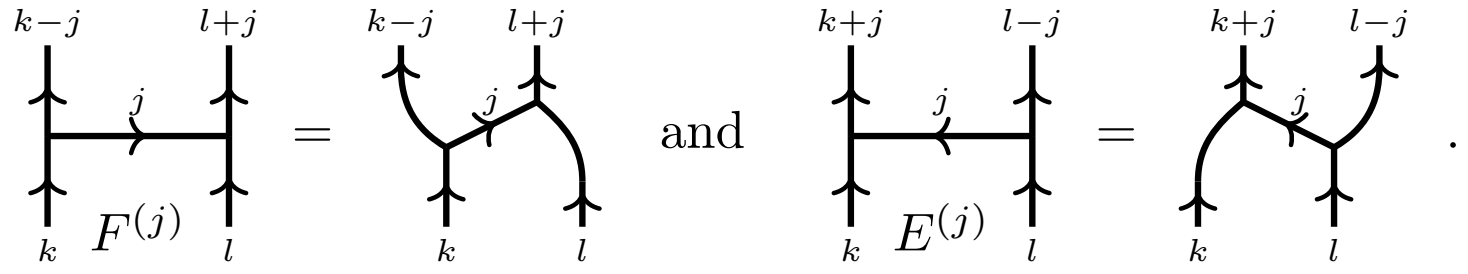

Here the left pictures are shorthand notations for the right pictures. ◇

**Definition 11B.8.** Let $\mathbf{SWeb}^{\uparrow} = \mathbf{SWeb}^{\uparrow}_{\Bbbk\oplus}$ denote the quotient of $\mathbf{lWeb}^{\uparrow}_{\Bbbk\oplus}$ by the following relations, displaying R, that are assumed to hold for all $k, l, h \in \mathbb{Z}_{\geq 0}$:

***Thin associativity and coassociativity***

(11B-9)

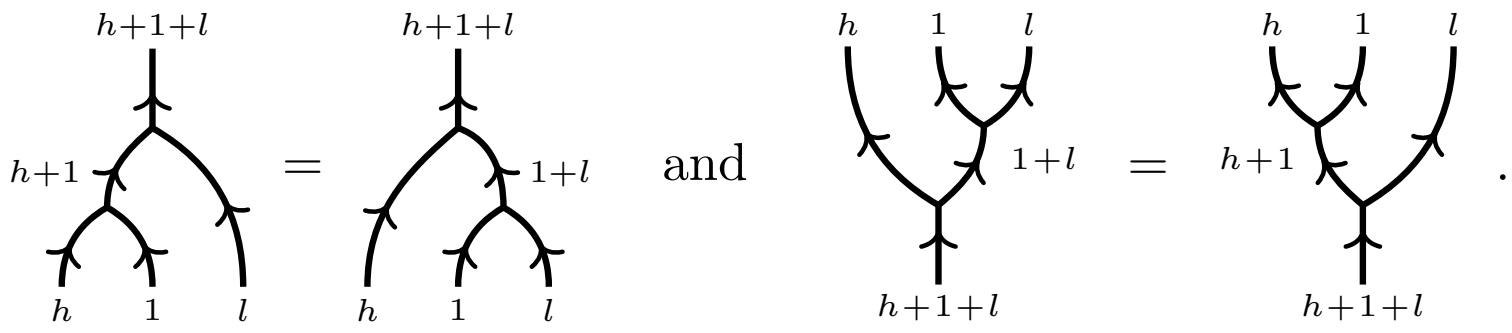

The ***thin square switch***

(11B-10)

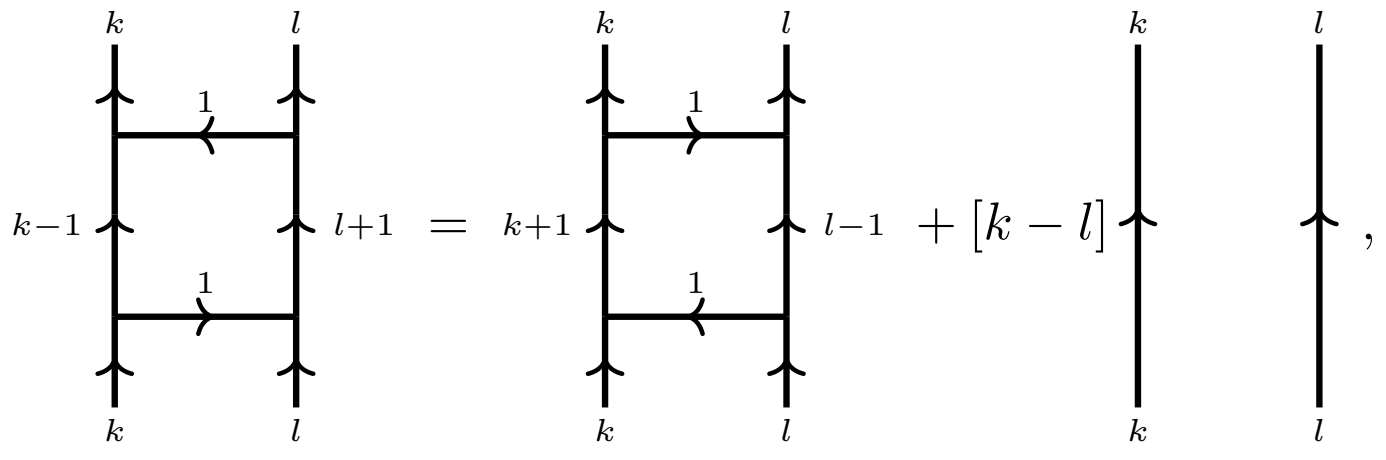

where in both defining relations we used the thick merges, splits and ladders from above.

We call $\mathbf{SWeb}^{\uparrow}$ the ***(upwards oriented labeled) generic web category***. ◇

**Definition 11B.11.** The two types of compositions of merges and splits

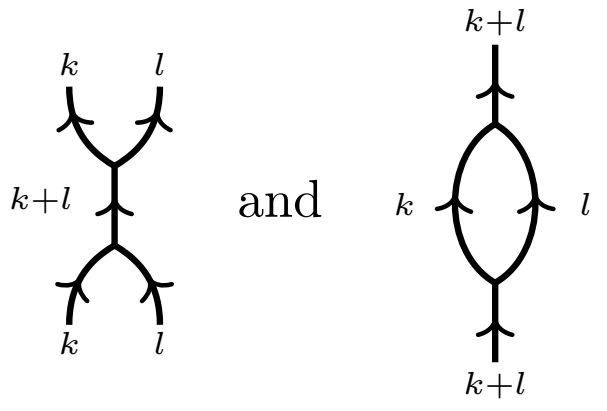

are called the $(k, l)$-***dumbbell*** (or ***dumbbell of thickness*** $k + l$) and the $(k, l)$-***digon***, respectively. ◇

Note that we also illustrate webs with potential zero or negative edges. In these cases, by convention, if a diagram has a 0-labeled edge then we remove that edge and if a diagram has a negative-labeled edge then the entire diagram is 0, e.g. the so-called ***thin digon removal***

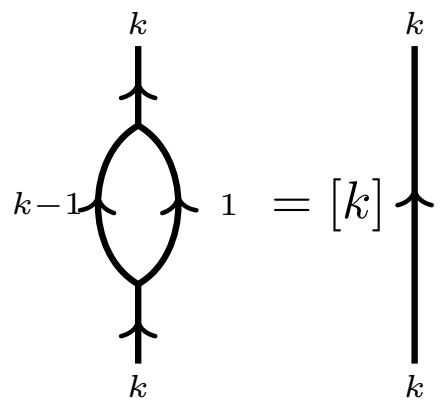

is a special case of the thin square switches Equation 11B-10, namely in is the case for $l = 0$. By Equation 11B-9 the thin digon removal implies

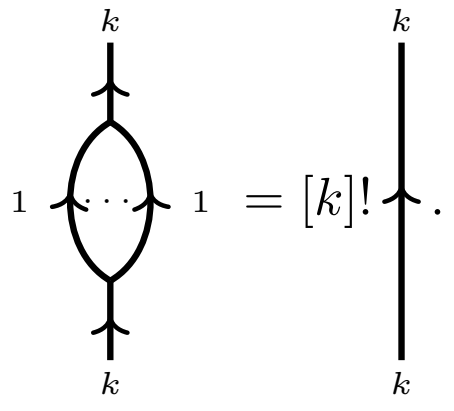

*Remark* 11B.12. Note that since in our ground field is $\Bbbk$ coming from a generic pair, the quantum numbers are invertible, the relation above is an invertible relation, namely we have

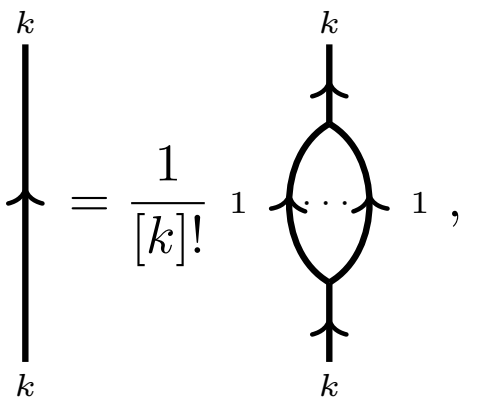

which justifies Definition 11B.4. Moreover, together with the associativity consequence from below

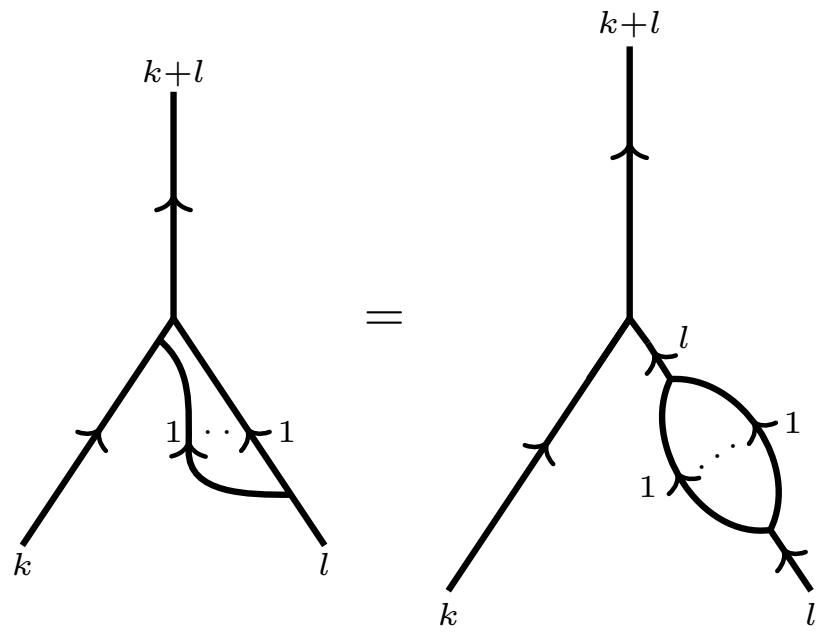

it also justifies Definition 11B.6. ◇

We wanted to define our web category in a way that it is as thin as possible, i.e. with as few $k$-labeled edges as possible, where $k > 1$. We now show that in this thin-defined web category we can deduce thick versions of relations Equation 11B-9 and Equation 11B-10.

**Lemma 11B.13.** *Thick (co)associativity*

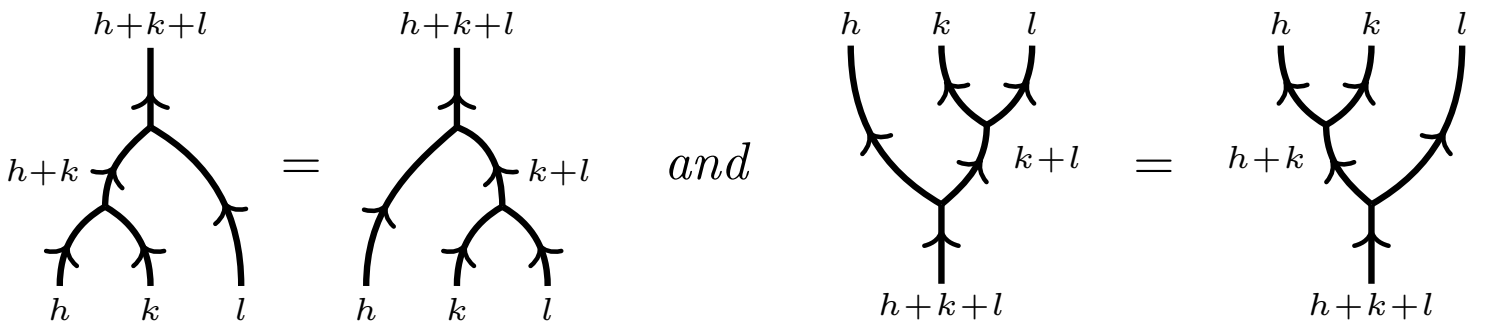

*follows from thin (co)associativity Equation 11B-9 and thin square switch relation Equation 11B-10.*

*Proof.* From the discussion above, we can explode the $k$-labeled edge and proceed like below

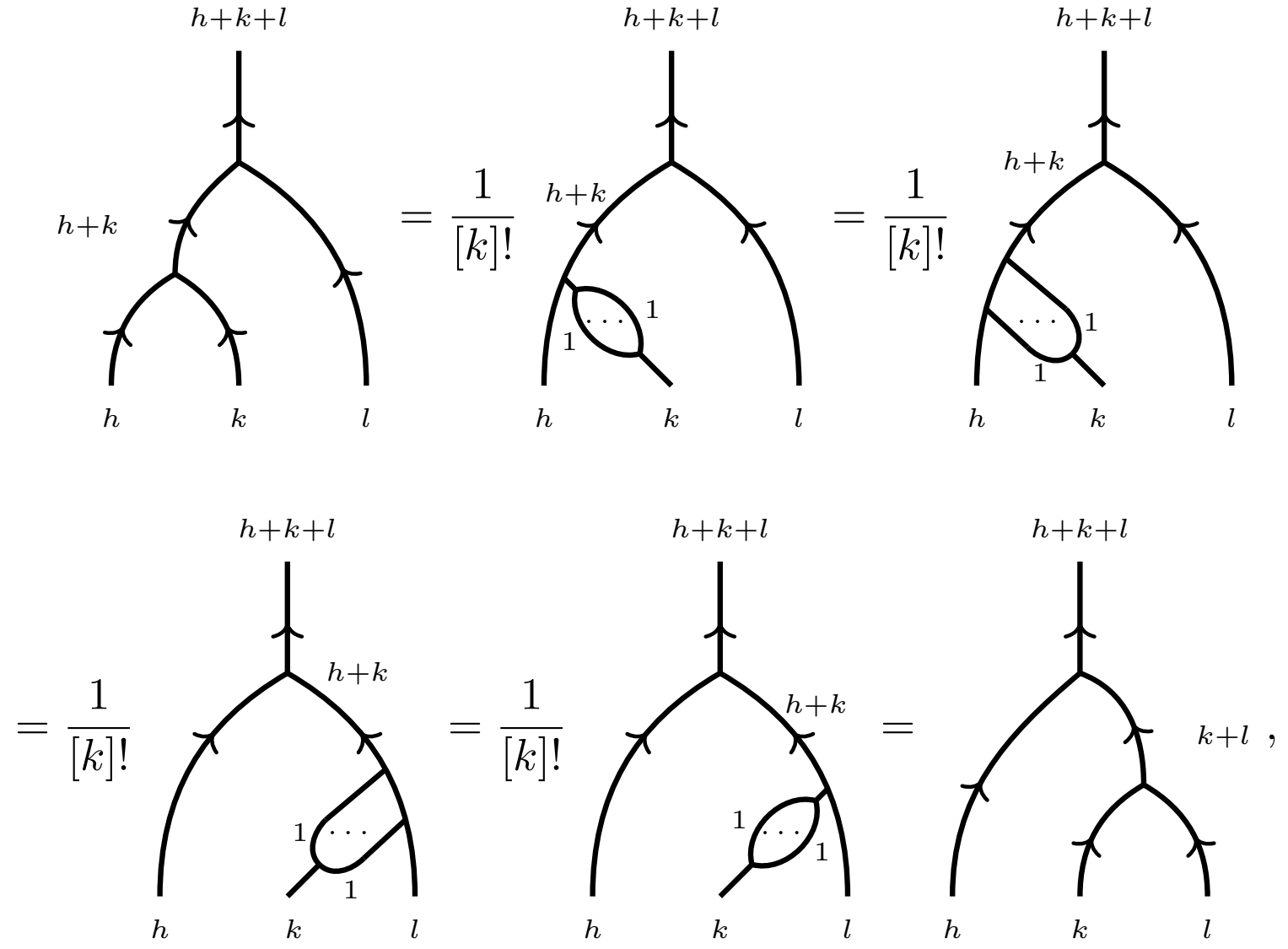

where in the second, third and fourth equality we used relation Equation 11B-9 to pull thin edges along thick ones from one side to the other. Thick coassociativity is proved similarly. □

**Lemma 11B.14.** *The* **thick square switches**

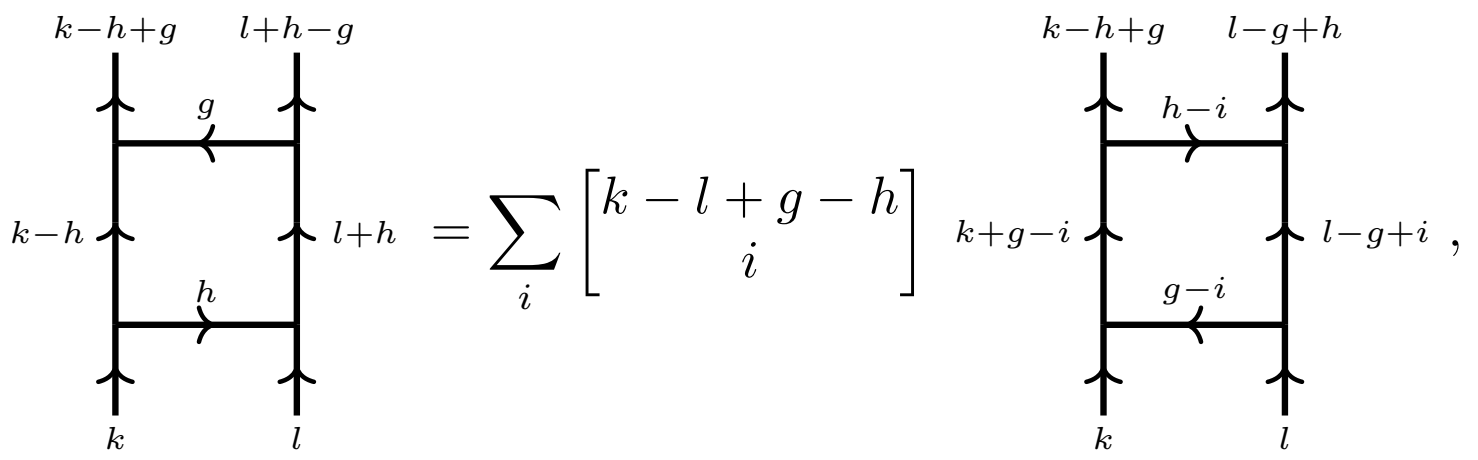

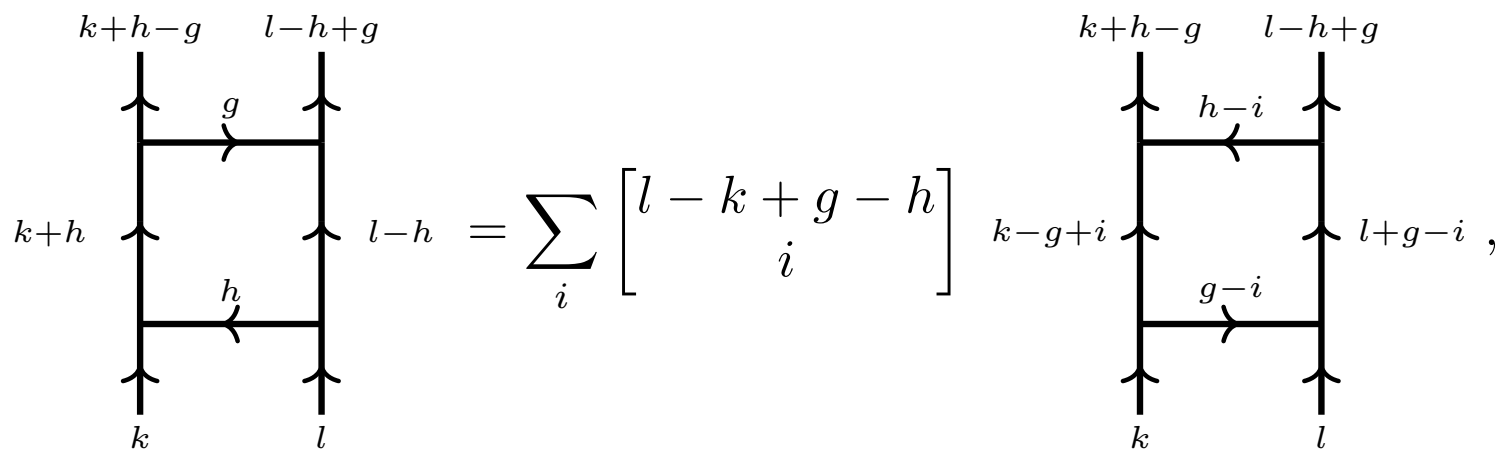

*and* **divided powers collapsing** *(note that the ladder rungs point in the same direction)*

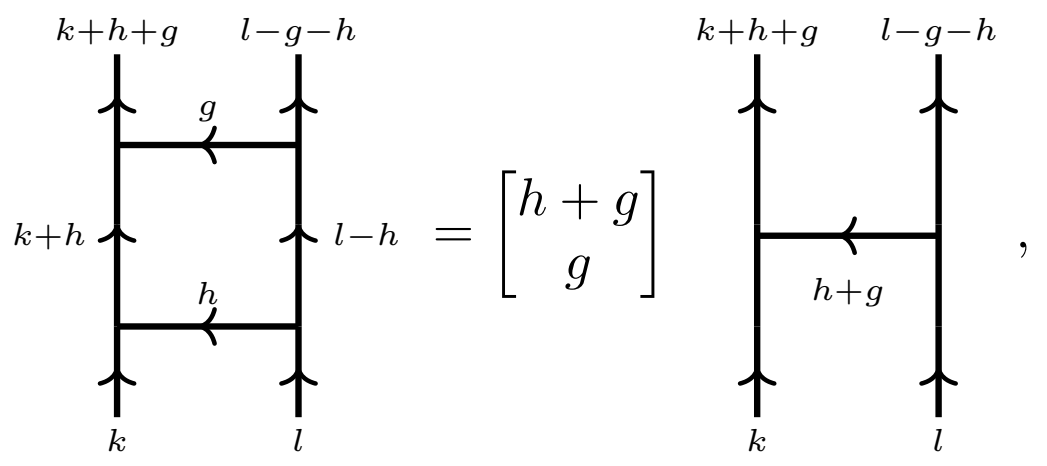

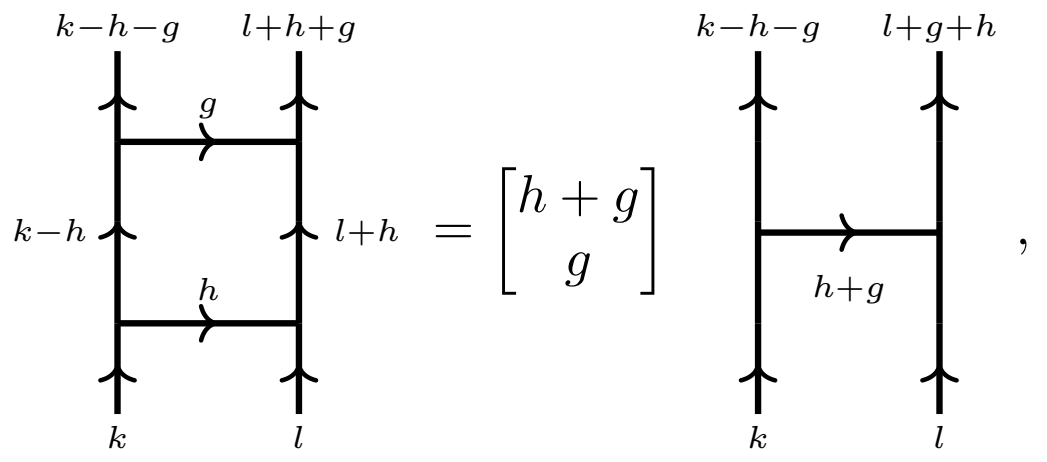

*hold in* $\mathbf{SWeb}^{\uparrow}$.

*Proof.* This is [**ST19**, Remark 2.5] and can be shown by using the definition of $\mathbf{SWeb}^{\uparrow}$ and Lemma 11B.13. □

As a special case of the lemma above, also the ***thick digon removal*** holds, i.e.

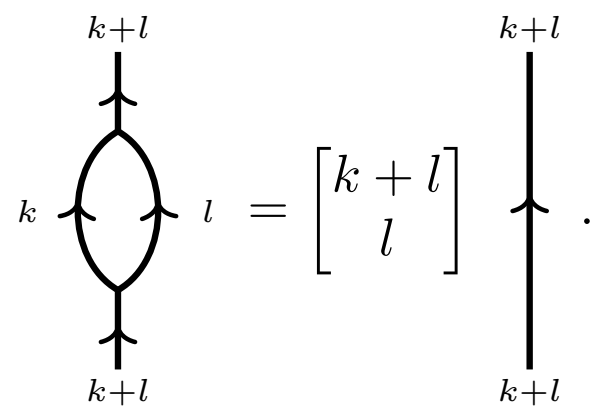

**Lemma 11B.15.** *Generic webs satisfy the so-called* **Serre relation**, *displayed as:*

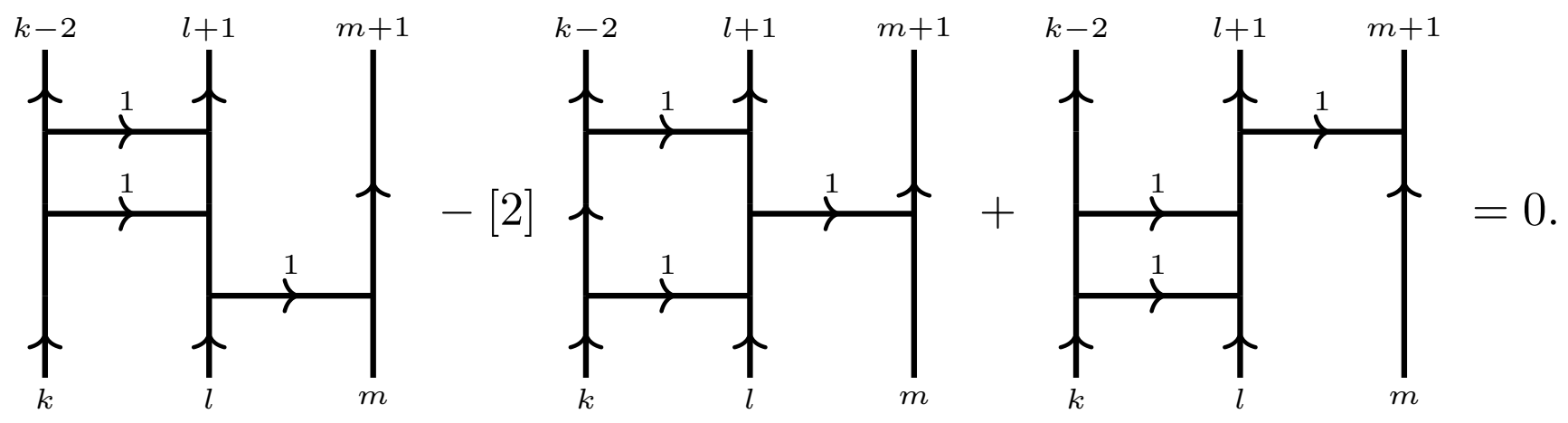

*Proof.* The argument is as in [**CKM14**, Lemma 2.2.1] and follows from applying associativity and the thick square switches to the middle web. Note that the divided powers from the leftmost and rightmost webs collapse, using the previous lemma. □

**Lemma 11B.16.** *The* **higher order Serre relations** *as in e.g. [**Lus10**, Chapter 7] hold in* $\mathbf{SWeb}^{\uparrow}$ *as well.*

*Proof.* For generic pairs the Serre relation implies its higher counterparts. □

Further, we want to emphasize that our thin crossings from Equation 11B-3 satisfy the braiding relations.

**Lemma 11B.17.** *The (thin) Reidemeister 2 and Reidemeister 3 moves hold in* $\mathbf{SWeb}^{\uparrow}$, *i.e.*

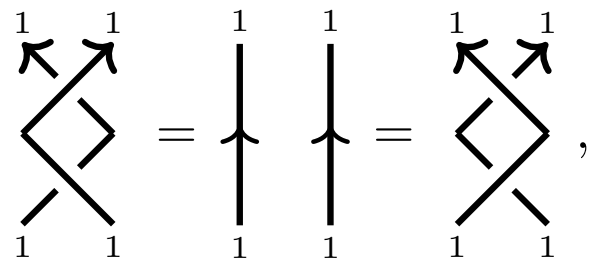

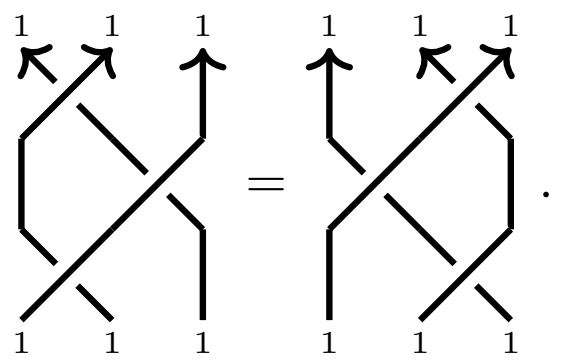

*Proof.* This can be verified via a straightforward calculation. □

We will see later that also the thick version of the previous lemma holds, meaning that $\mathbf{SWeb}^{\uparrow}$ is a braided monoidal category. For now, note that as a direct consequence of the Lemma 11B.17 we have.

**Corollary 11B.18.** *Thin over and undercrossings are inverses of each other.* □

Whenever every object has a dual and we have evaluation and coevaluation maps, we want incorporate the analogous notions into our diagrammatic category. Therefore, we define the following.

**Definition 11B.19.** The ***(upward-downward oriented labeled) generic web category without relations***, denote it by $\mathbf{lWeb}^{\uparrow,\downarrow}$, is the category monoidally generated by $\mathrm{S}: \{k^{\uparrow},\ k_{\downarrow} \mid k \in \mathbb{Z}_{\geq 0}\}$ with $0_{\uparrow} = 0_{\downarrow} = 0$, and T is given by

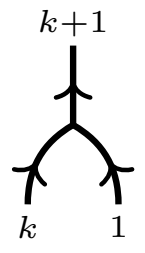 $: k^{\uparrow} \otimes 1^{\uparrow} \to (k+1)^{\uparrow}$, 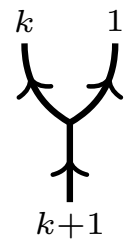 $: (k+1)^{\uparrow} \to k^{\uparrow} \otimes 1^{\uparrow}$,

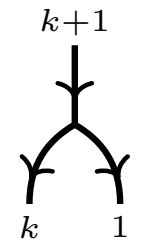 $: k_{\downarrow} \otimes 1_{\downarrow} \to (k+1)_{\downarrow}$, 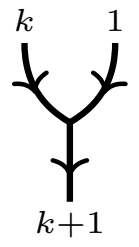 $: (k+1)_{\downarrow} \to k_{\downarrow} \otimes 1_{\downarrow}$,

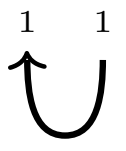 $: 0 \to 1^{\uparrow} \otimes 1_{\downarrow}$, 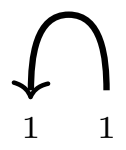 $: 1_{\downarrow} \otimes 1^{\uparrow} \to 0$,

where $k \in \mathbb{Z}_{\geq 0}$. (Note that we only have caps and cups in one direction.) ◇

As before, we want to consider a quotient of $\mathbf{lWeb}^{\uparrow,\downarrow}$ by some relations on the morphisms and in order to give these relations we need to define some other webs, which we will do in $\mathbf{lWeb}^{\uparrow,\downarrow}_{\Bbbk\oplus}$. First, we define upward pointing thin crossings as in Equation 11B-3, thick upward crossings as in Definition 11B.4, the $(1,k)$- merges and splits as in Definition 11B.5 and also the thick merges and splits as in Definition 11B.6. Then we continue as follows.

**Definition 11B.20.** We define the following.
***Leftward (thin) crossings***

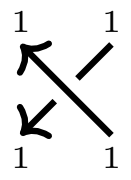 := 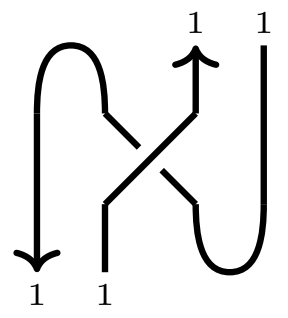 and 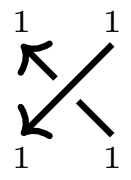 := 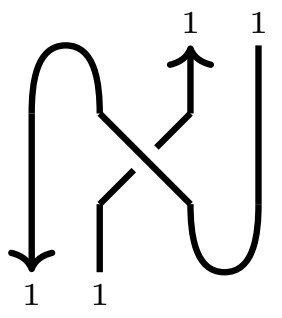 ;

***downward pointing (thin) over and undercrossings***

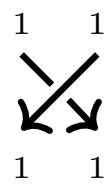 := 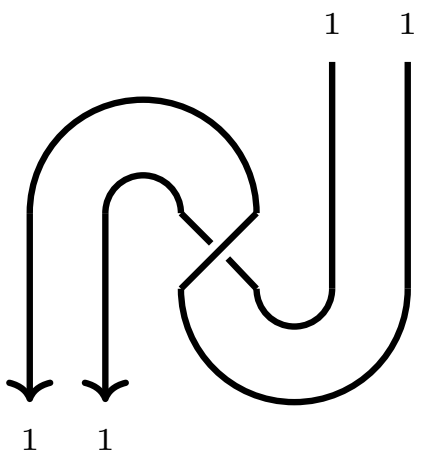 , 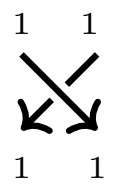 := 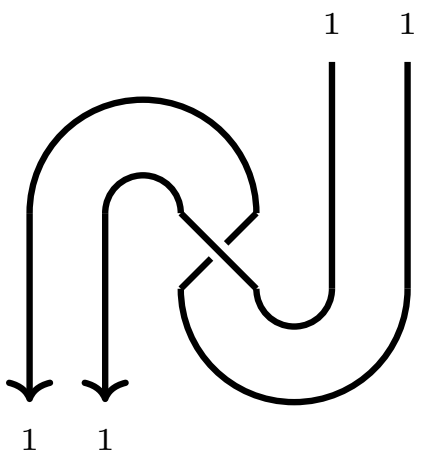 .

(We will also have rightwards crossings later on.) ◇

**Definition 11B.21.** The $k$-***labeled caps*** are defined via explosion

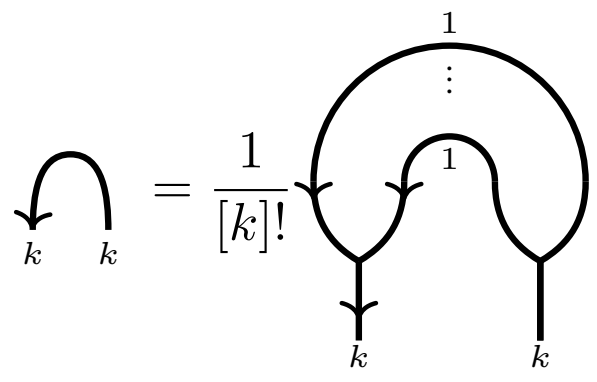

and similarly the $k$-***labeled cups***. ◇

We will later impose that the (thin) leftward crossings are invertible and we will draw their inverses as (thin) ***rightward crossings***, i.e.

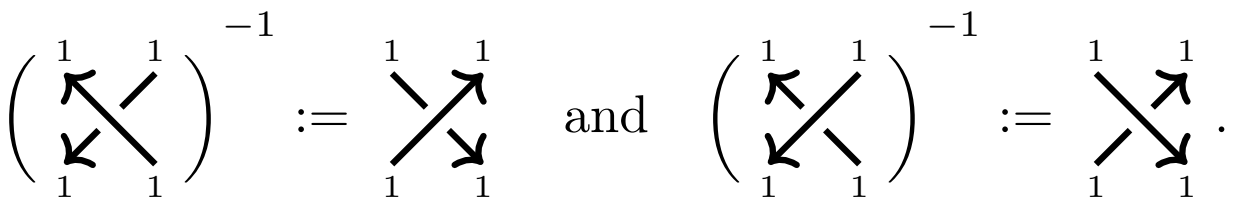

**Definition 11B.22.** The ***thick leftward*** crossings are defined by explosion

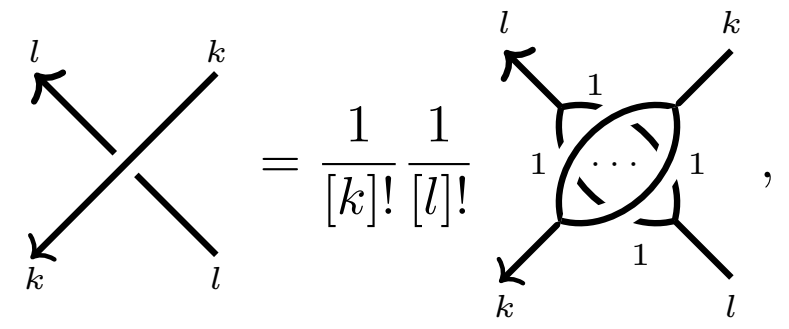

and similarly the ***thick rightward*** crossings. ◇

**Definition 11B.23.** The ***rightward cap and cup*** are defined as

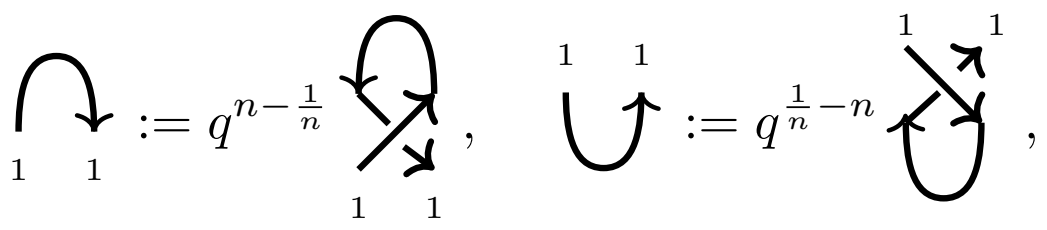

and their thick versions are obtained by explosion. ◇

From now on $n$ actually plays a role.

**Definition 11B.24.** The ***(thin) symmetric (upward-downward pointing) category of*** $\mathfrak{gl}_n$***-webs***, denote it by $\mathbf{sWeb}^{\uparrow,\downarrow}(\mathfrak{gl}_n) = \mathbf{sWeb}^{\uparrow,\downarrow}(\mathfrak{gl}_n)_{\Bbbk\oplus}$, is the quotient category obtained from $\mathbf{lWeb}^{\uparrow,\downarrow}_{\Bbbk\oplus}$ by imposing relations Equation 11B-9 and Equation 11B-10 together with their downward pointing versions, and relations Equation 11B-25-Equation 11B-28 from below, meaning we give R.
***Invertibility of thin leftward crossings***:

(11B-25)

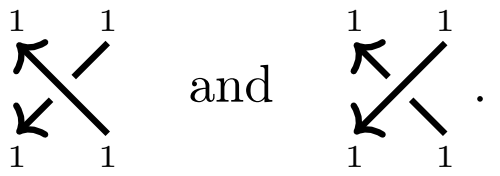

***The thin zigzag relation***:

(11B-26)

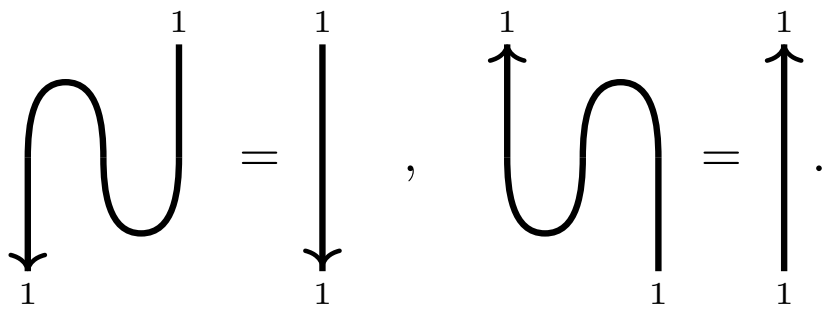

***The merge-split slides***:

(11B-27)

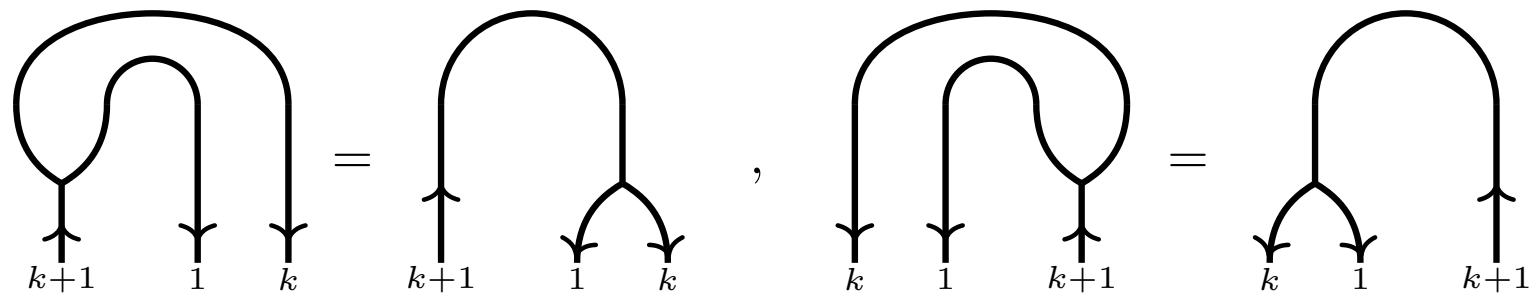

, $k \geq 1$

together with their cup and opposite orientation versions.
Finally,

(11B-28)

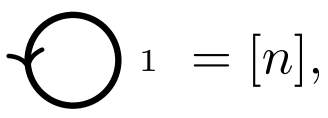

called ***1-labeled circle removal***. ◇

Note that $1_\uparrow$ and $1_\downarrow$ are duals to each other by the thin zigzag relation Equation 11B-26. Later we will show that for every $k \geq 1$, $k_\uparrow$ and $k_\downarrow$ are duals to each other meaning that $\mathbf{SWeb}^{\uparrow,\downarrow}(\mathfrak{gl}_n)$ is endowed with a rigid structure.

*Remark* 11B.29. We could also define $\mathbf{SWeb}^{\uparrow,\downarrow}(\mathfrak{gl}_n)$ using only upward pointing merges and splits but then we would have to include caps and cups of ***any*** thickness as generators. In that case we would define the downward pointing merges and splits like below

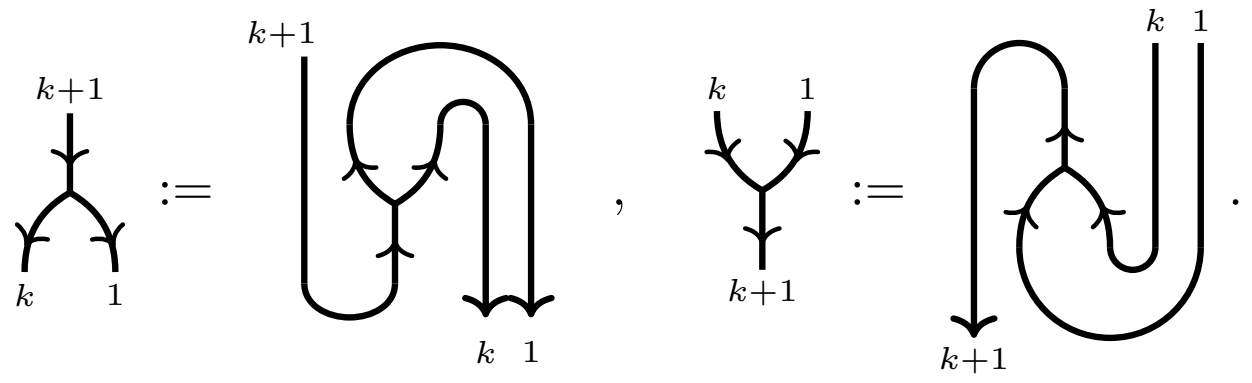

See also Lemma 11C.3. ◇

11C. **Properties of webs.** We now give the first consequences of the defining web relations.

**Lemma 11C.1.** *Relations Equation 11B-26 and Equation 11B-27 imply the thick zigzag relations*

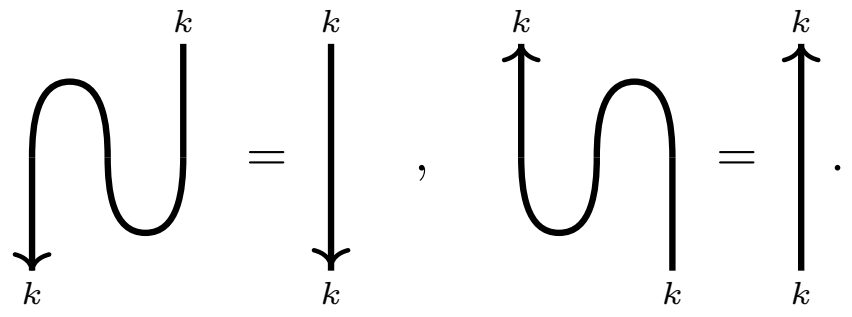

*Proof.* Using Lemma 11B.13 we see that we can adapt the proof of [**RT16**, Lemma 2.21]. Namely, after explosion, using thin zigzag and merge-split slides we have

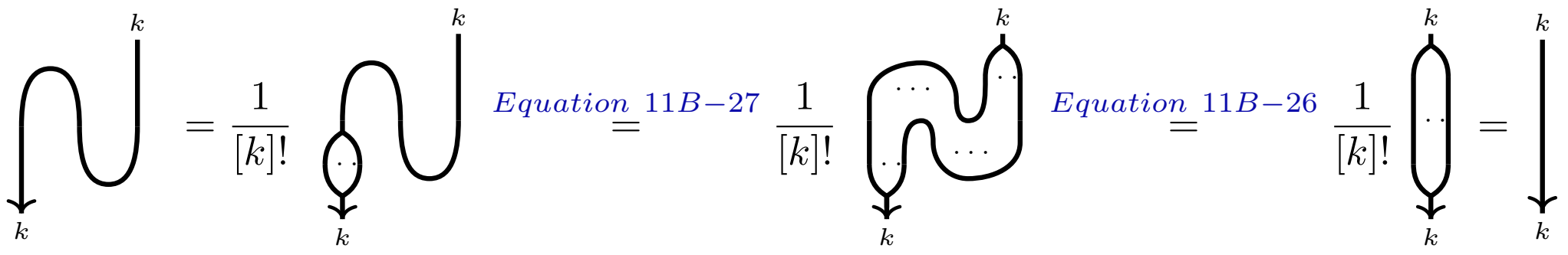

and similarly for the other thick zigzag relation. □

**Lemma 11C.2.** *The thick versions of relation Equation 11B-27 hold, i.e.*

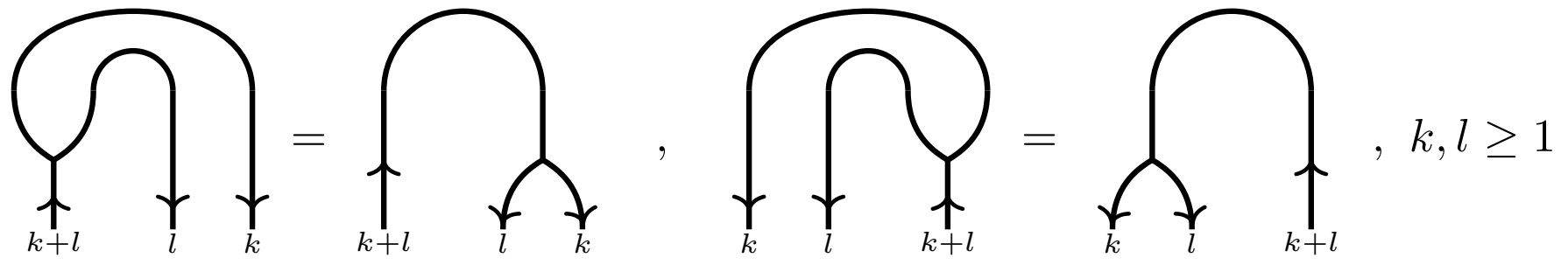

*as well as their cup and opposite orientation versions.*

*Proof.* We prove the first relation via induction on $k+l$. The case $k+l=1$ follows trivially. Assume it holds for thickness smaller or equal to $k+l-1$. We have

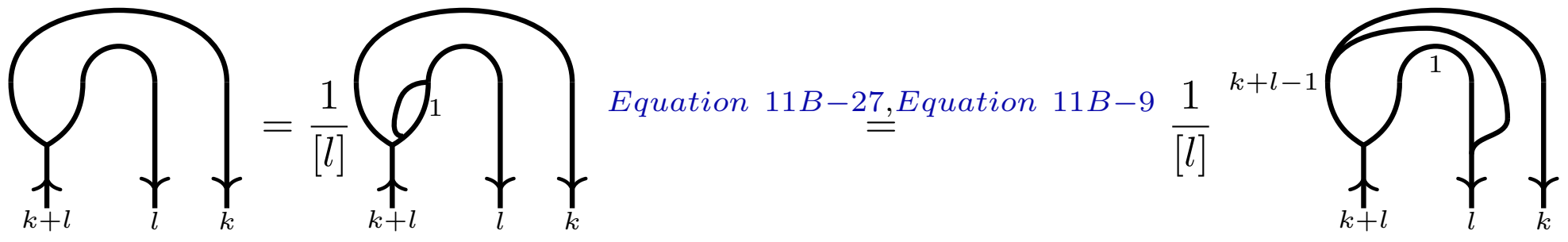

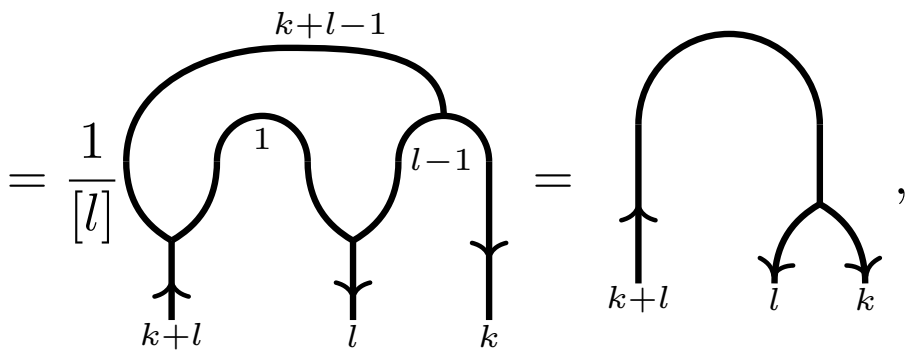

where in the third equality we used the induction hypothesis to slide the $(k, l-1)$ split and in the last equality we used relation Equation 11B-27 for the $(k+l-1, 1)$ split, then associativity and the digon removal. □

**Lemma 11C.3.** *We have the following*

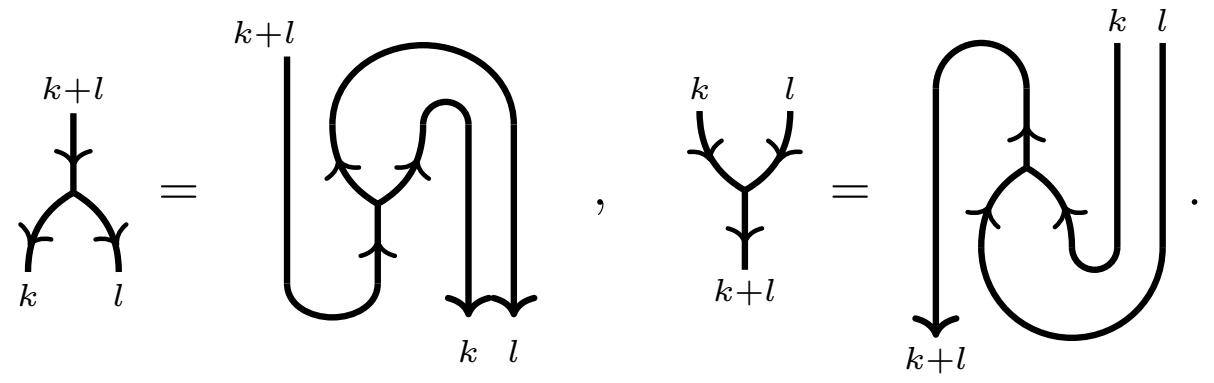

*Proof.* This is a direct consequence of the previous lemma and thick zigzag relation. □

The next lemma shows that in $\mathbf{sWeb}^{\uparrow,\downarrow}(\mathfrak{gl}_n)$ we have a relation akin to the usual planar isotopies, namely a relation of Reidemeister 1 type.

**Lemma 11C.4.** *The following ribbon equations hold. First,*

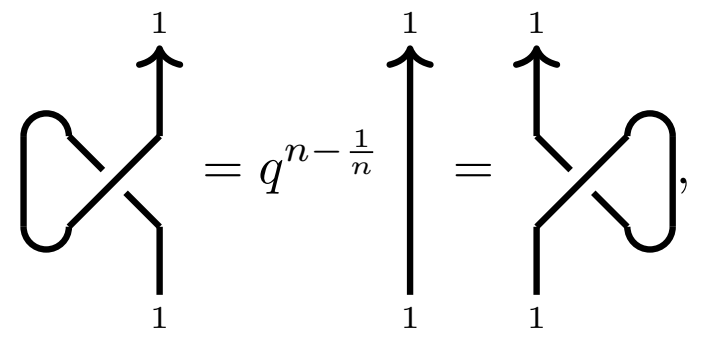

*as well as its downward pointing version, and*

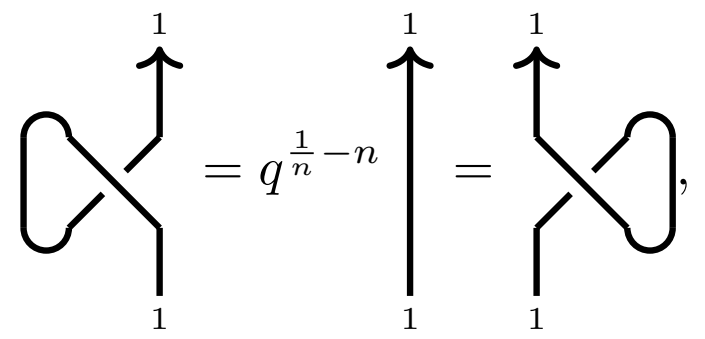

*as well as its downward pointing version.*

*Proof.* From the definition of the rightward cap

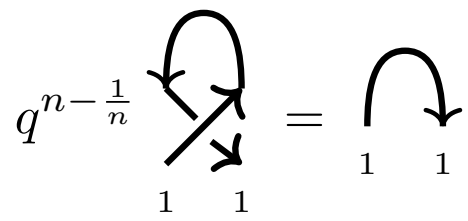

after adding a crossing and using the fact that the right overcrossing is the inverse of the left undercrossing we obtain

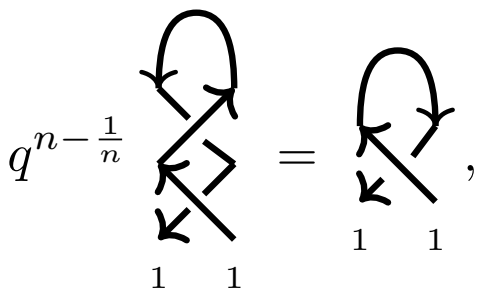

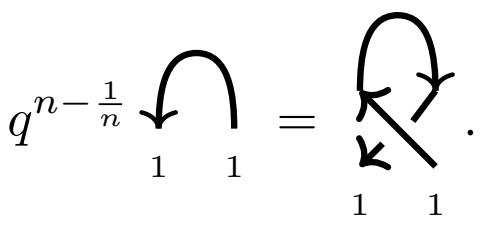

Further, by adding caps on both sides we have

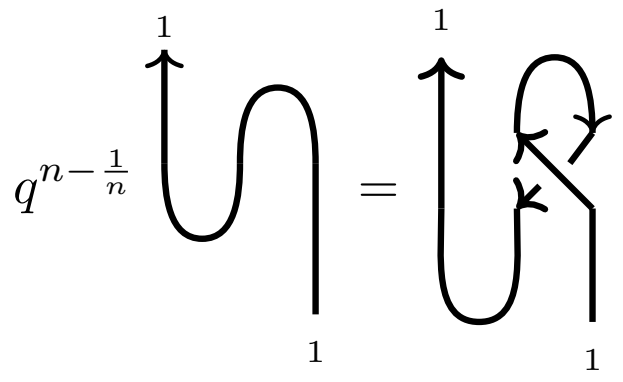

and thus by straightening the left-hand side

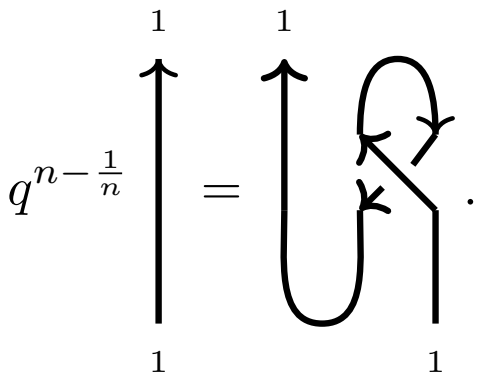

However,

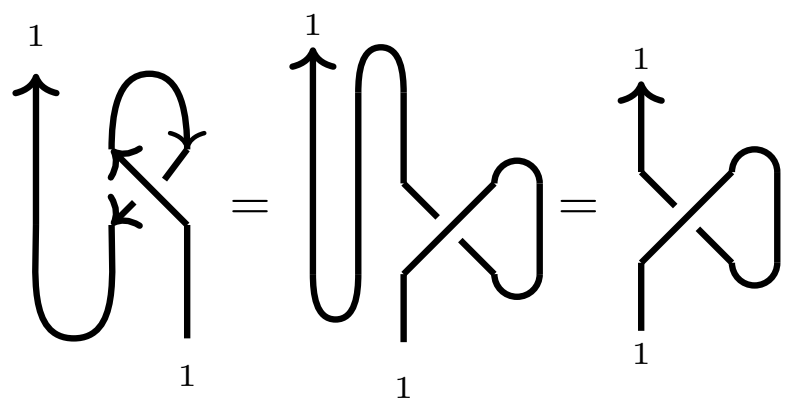

where in the first equality we used the definition of the left crossing and in the second we used the zigzag relation. This means that

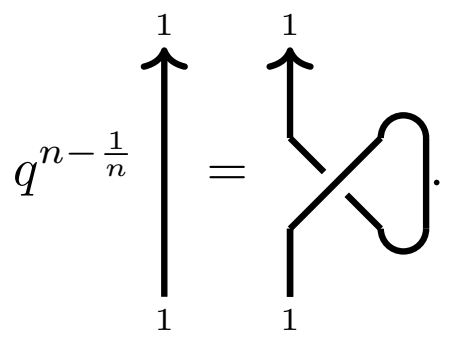

All other equalities can be proved verbatim. ☐

**Lemma 11C.5.** *The 1-circle evaluation with reversed orientation also holds, i.e.*

(11C-6)

*Proof.* This follows from Definition 11B.23, Definition 11B.20, the previous lemma and Equation 11B-28. ☐

**Lemma 11C.7.** *The* **thin sideways dumbbell removal** *holds, i.e.*

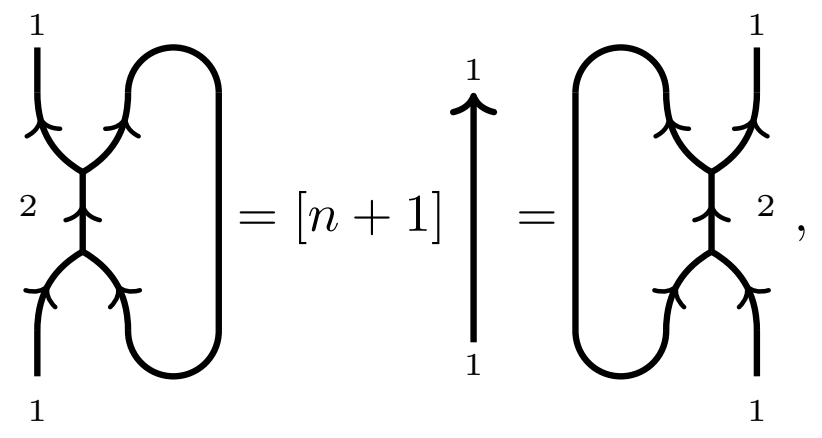

*together with the downward pointing version.*

*Proof.* This follows from Equation 11B-3 and Lemma 11C.4 and Equation 11C-6 as below

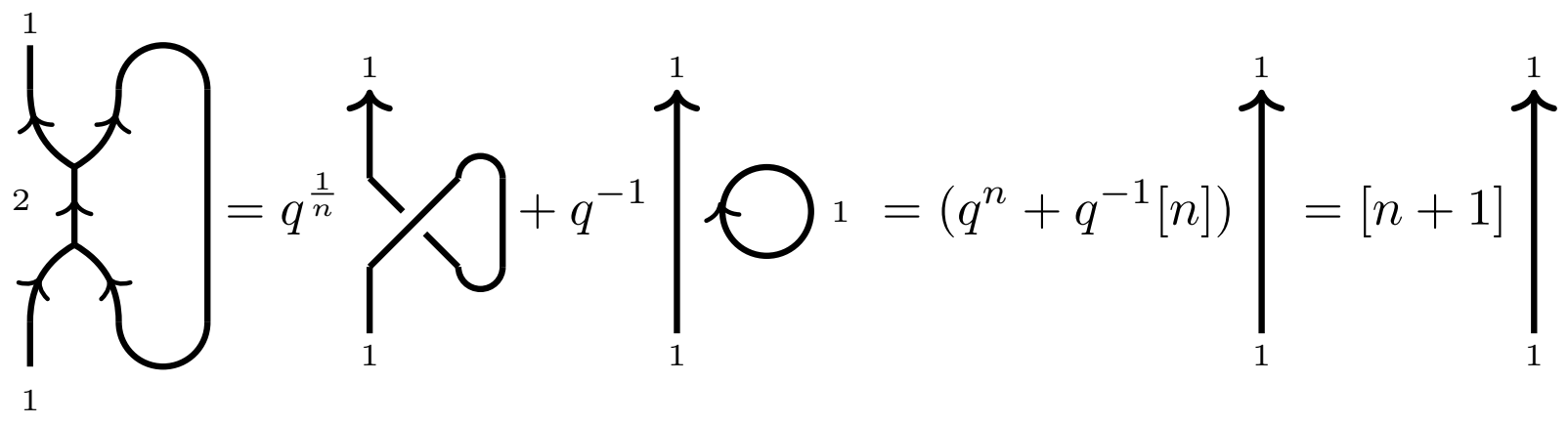

and similarly the downward pointing version. ☐

**Lemma 11C.8.** *We have the following thick circle evaluation formula*

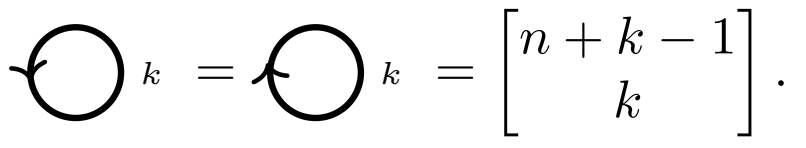

*Proof.* We prove this by induction on $k$. The case $k = 1$ is relation Equation 11B-28. Using the induction hypothesis we have

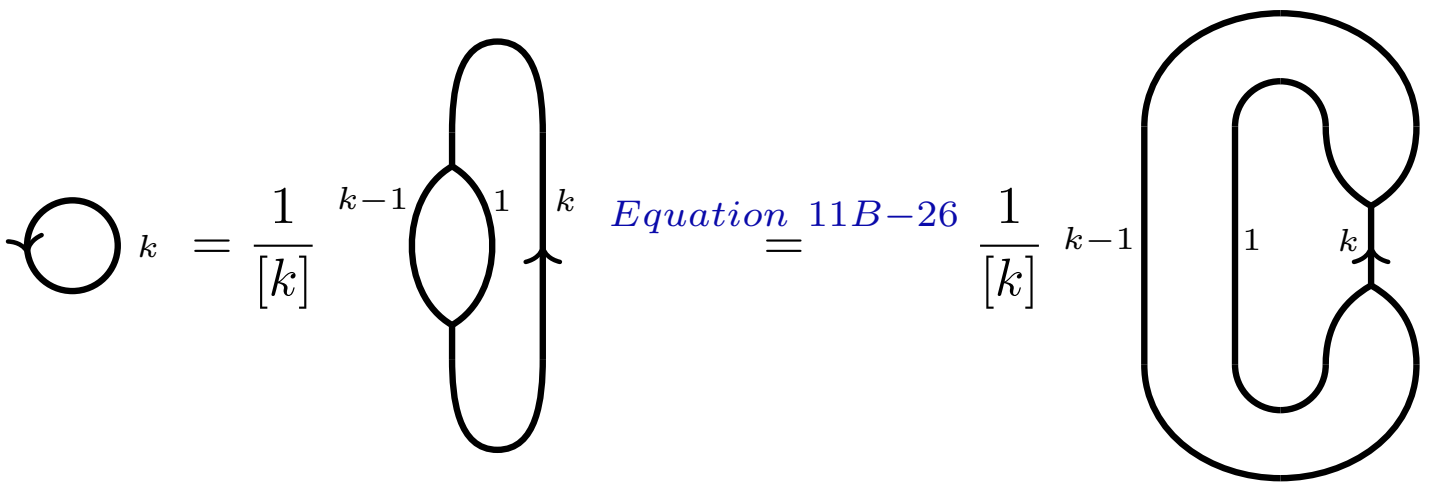

$\overset{\text{Equation 11B−10}}{=} \frac{1}{[k]}$ 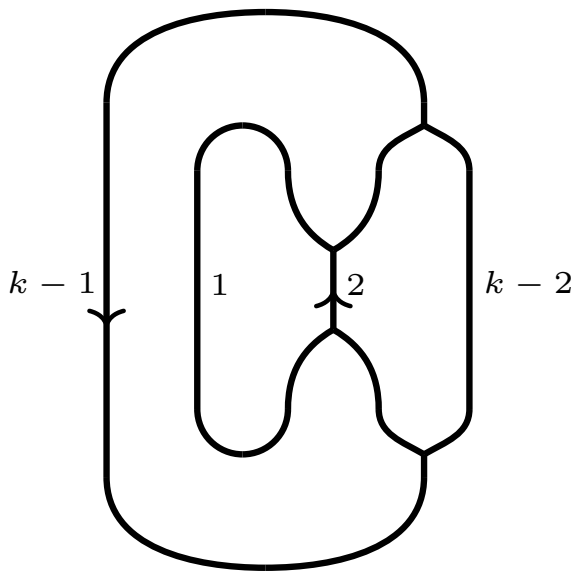

$- \frac{[k-2]}{[k]}$ 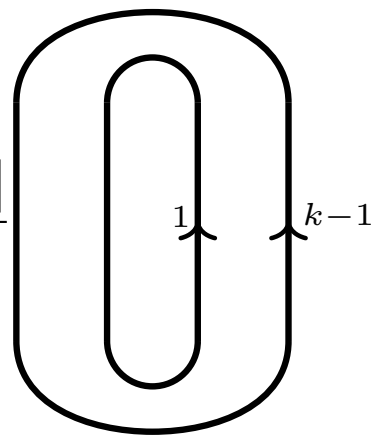

$\overset{\text{Lemma 11C.7}}{=} \frac{1}{[k]}[n+1]$ 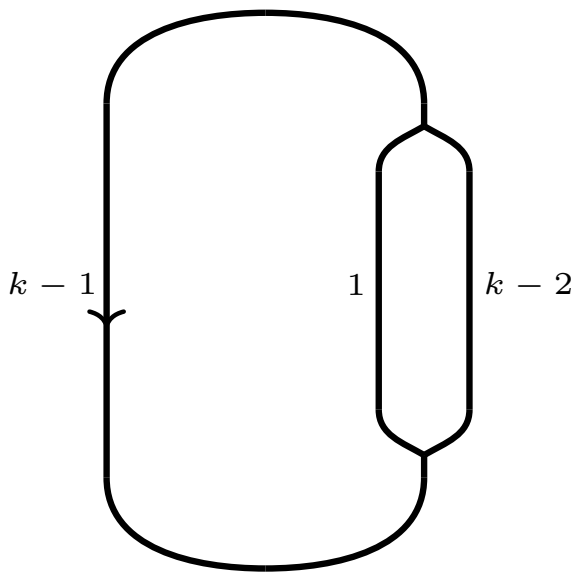

$- \frac{[k-2]}{[k]}[n] \begin{bmatrix} n+k-2 \\ k-1 \end{bmatrix}$

$$= \frac{1}{[k]}[n+1][k-1] \begin{bmatrix} n+k-2 \\ k-1 \end{bmatrix} - \frac{[k-2]}{[k]}[n] \begin{bmatrix} n+k-2 \\ k-1 \end{bmatrix} = \begin{bmatrix} n+k-1 \\ k \end{bmatrix},$$

where the last equality is a straightforward computation. □

**Lemma 11C.9.** *The $(k,l)$-dumbbell removal, i.e.*

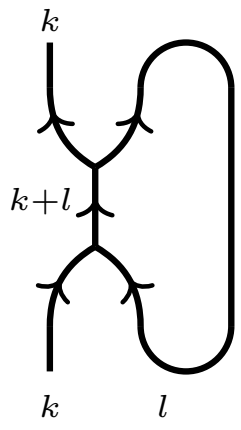

$= \begin{bmatrix} n+k+l-1 \\ l \end{bmatrix}$ 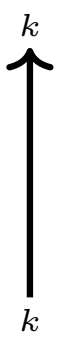

$=$ 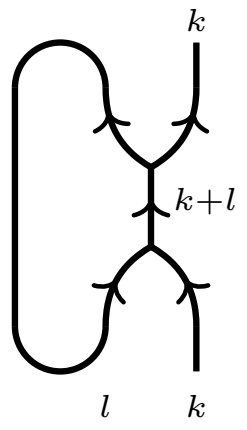

,

*as well as its downward pointing version hold.*

*Proof.* We prove the second equality via induction on $l$. For $l = 1$ we have

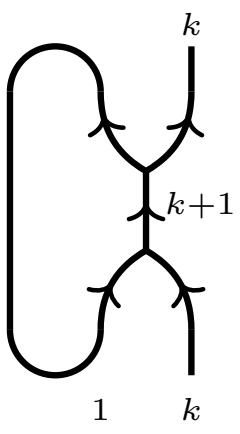

$\overset{\text{Equation 11B−10}}{=}$ 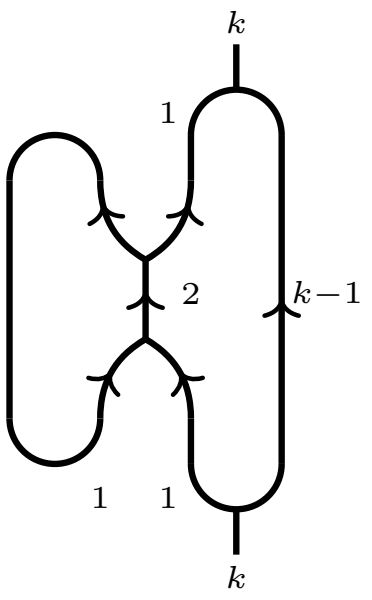

$- [k-2]$ 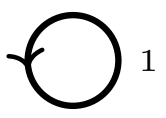

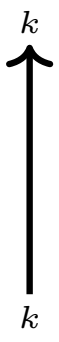

$\overset{\text{Lemma 11C.7}}{=} [n+1]$ 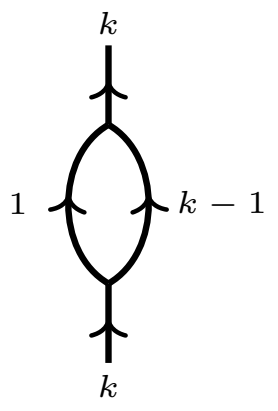

$- [k-2][n]$ 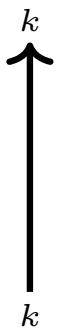

$= ([n+1][k] - [k-2][n])$ 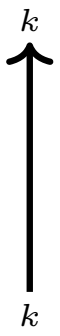

$= [n+k]$ 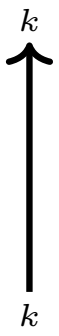

.

Now, assume the relation holds for $(l_1, k)$-dumbbells, where $l_1 \leq l-1$. For the $(l,k)$-dumbbell we have

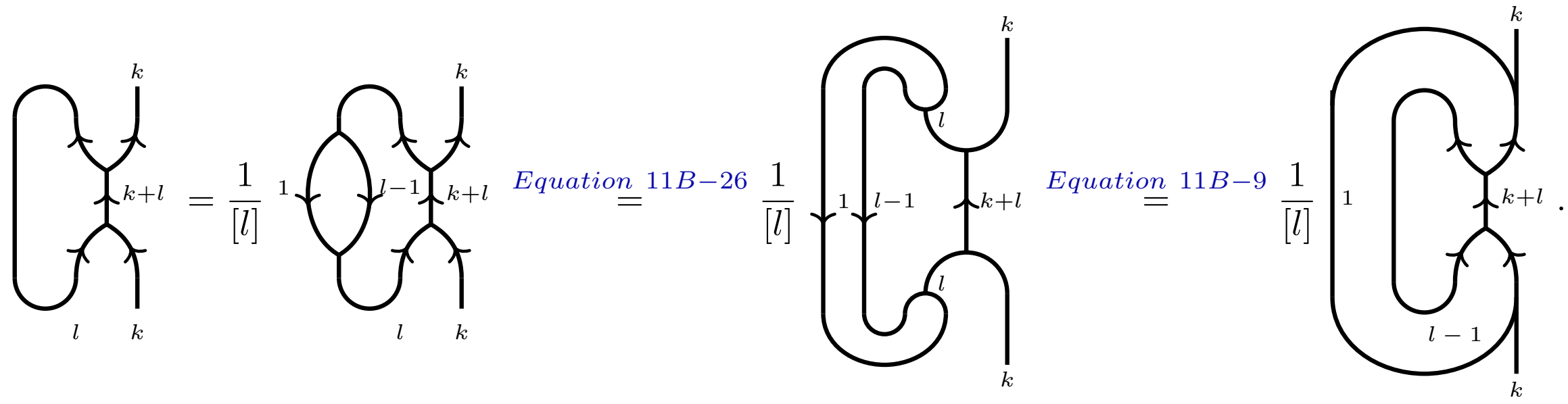

Finally, from the induction hypothesis and the case $l=1$ from above, the later is equal to

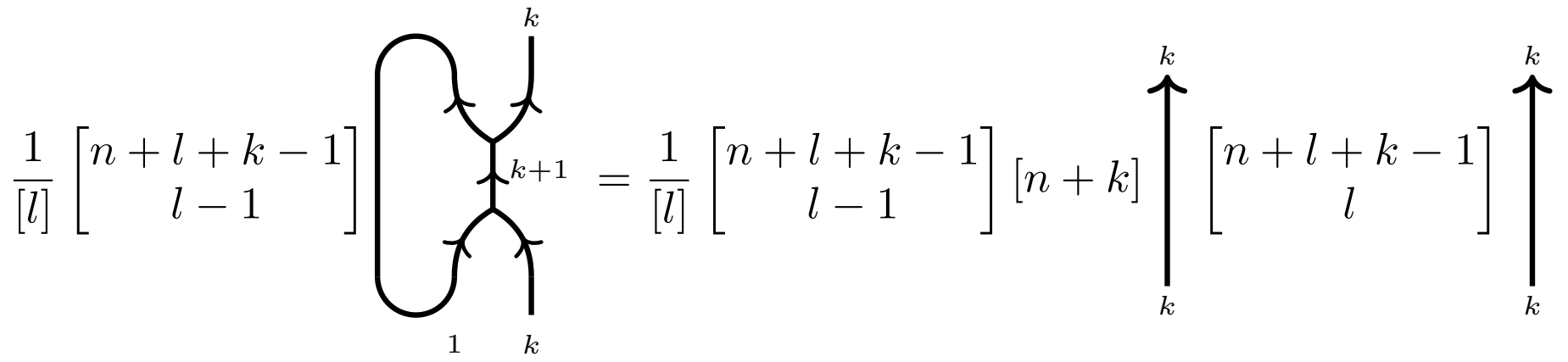

completing the proof. ☐

**Lemma 11C.10.** *The other versions of (thick) zigzag relations also hold:*

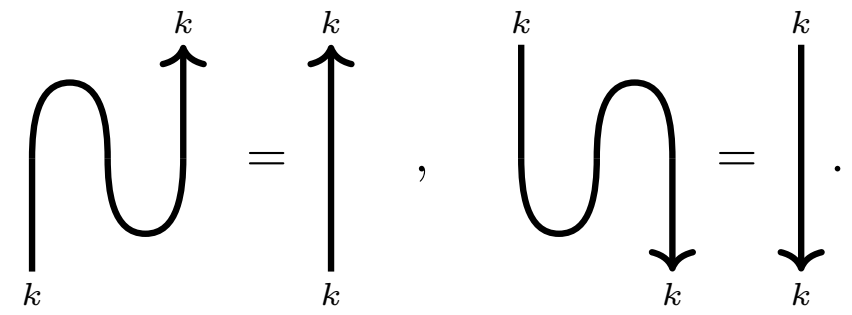

*Proof.* From the previous two lemmas we have

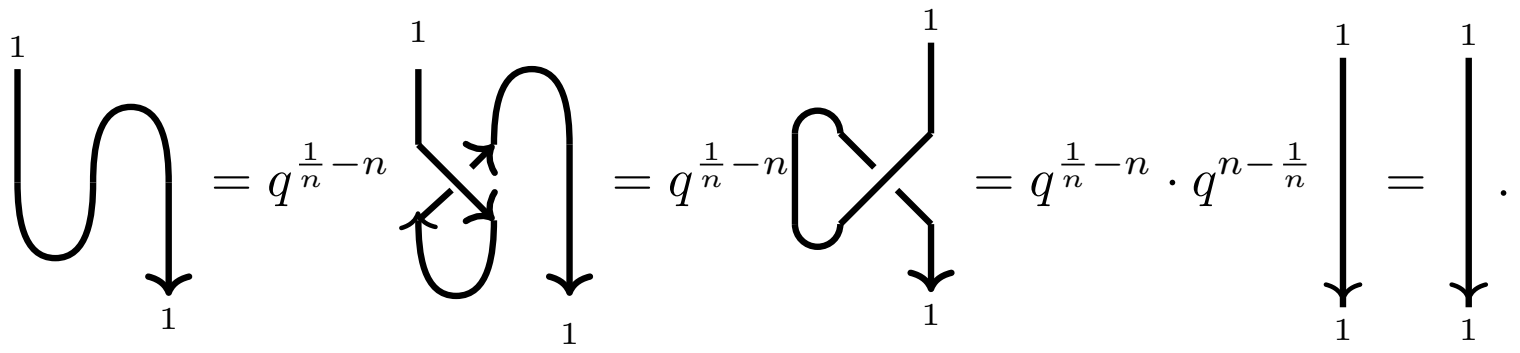

Now we use this together with the digon removal and the relations Equation 11B-26 to obtain the thick version

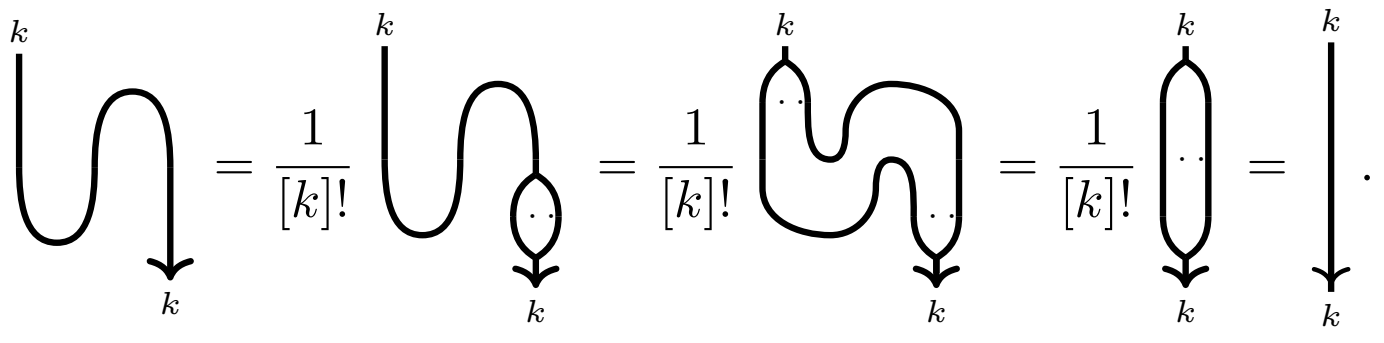

The other equation is proved similarly. ☐

A direct consequence of this lemma and Lemma 11C.1 is the following.

**Corollary 11C.11.** *The category* $\mathbf{sWeb}^{\uparrow,\downarrow}(\mathfrak{gl}_n)$ *can be endowed with a pivotal structure, where* $k_\uparrow$ *and* $k_\downarrow$ *are dual to each other.* ☐

We will use this pivotal structure throughout this section. With this structure, Lemma 11C.8 implies that $k^\uparrow$ and $k^\downarrow$ are both of categorical dimension $\begin{bmatrix} n+k-1 \\ k \end{bmatrix}$.

**Lemma 11C.12.** *The so-called* **pitchfork relations** *hold*

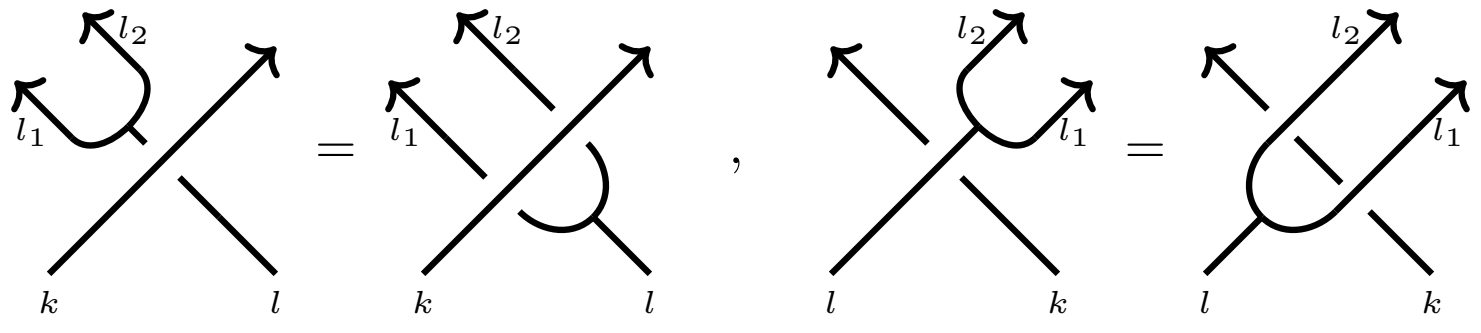

*together with the merge and downwards pointing versions.*

*Proof.* We prove the first relation via induction of $k+l$. For $k+l=1$, this follows trivially. Assume the relation holds for labels whose sum is smaller than $k+l$. We explode the $k$-edge as below

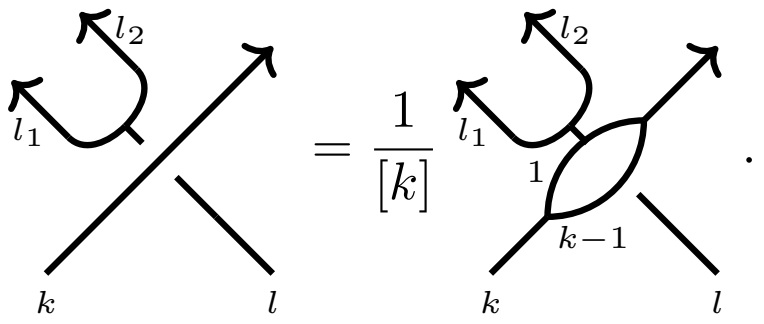

Now, using the induction hypothesis, we slide the edges labeled by 1 and $k-1$ through the $(l_1, l_2)$-split and after removing the digon we obtain

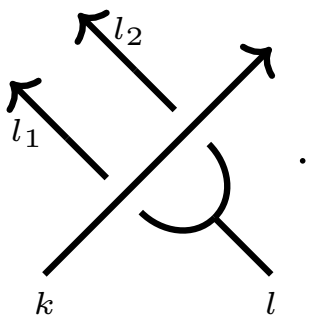

The second relation is shown verbatim. ☐

**Corollary 11C.13.** *Combining the result above for merges and splits, we have that also the sliding of edges through dumbbells hold.* ☐

**Lemma 11C.14.** *The thick versions of Reidemeister 2 and 3 from Lemma 11B.17 hold. Namely,*

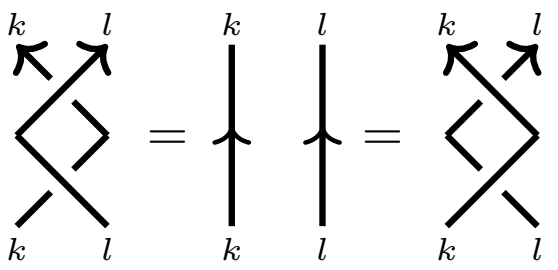

*and*

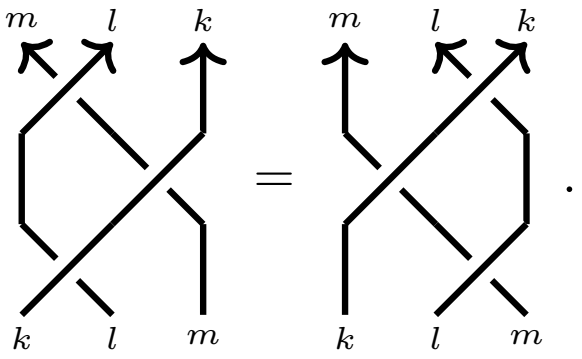

*Proof.* We prove only the first equality in the Reidemeister 2 move. From Definition 11B.4, Corollary 11C.13, Lemma 11B.17 and the digon removal we have

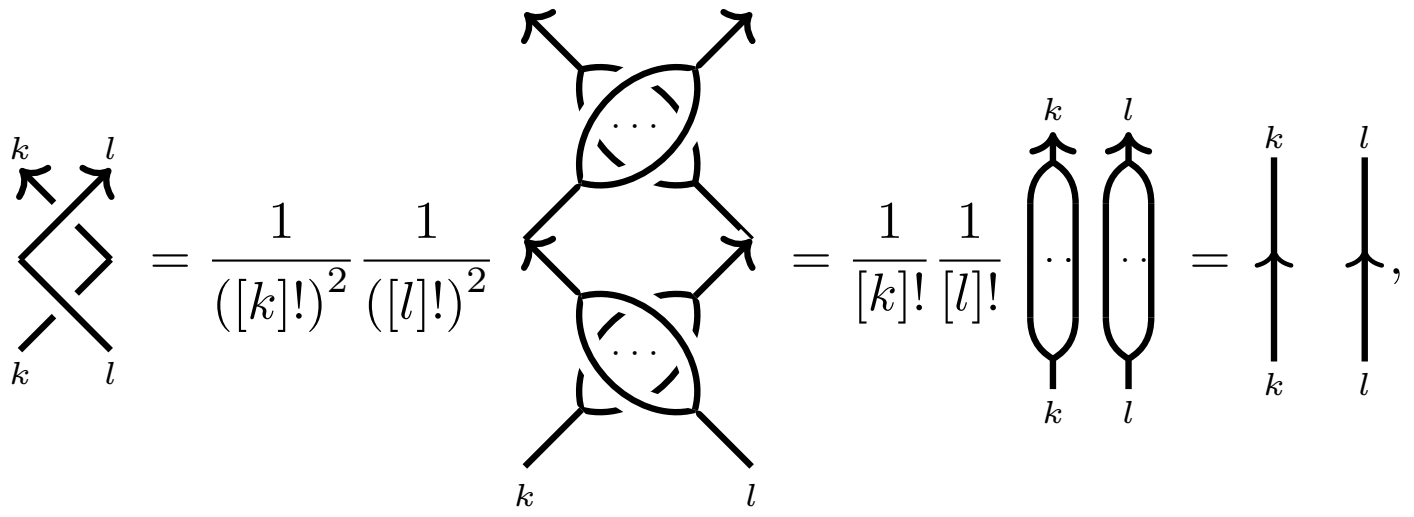

where in the second equation we used Corollary 11C.13 to move the central dumbbells up and then we used the fact that thin over and undercrossings are inverses to each other. ☐

**Lemma 11C.15.** *The sliding of thick edges through thick caps holds, together with its cup, reversed orientation and overcrossing versions*

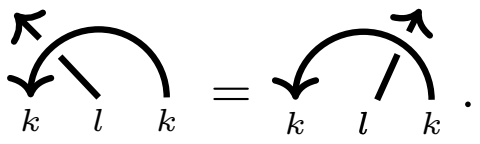

*Proof.* The case $k=l=1$ follows from a direct check using the definition of the leftward crossings. The general case is:

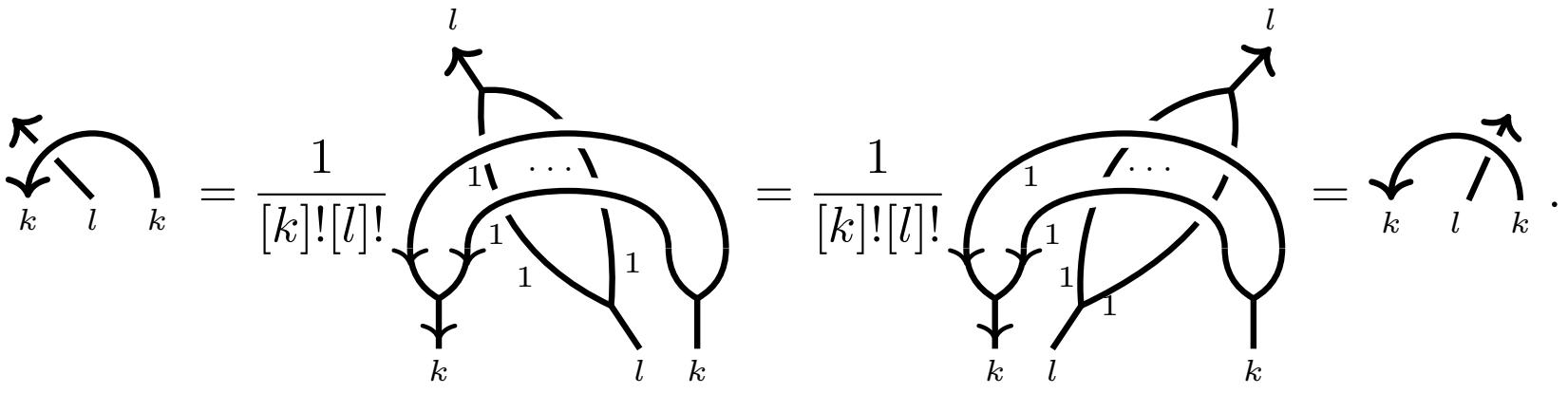

Here we exploded the strings. ☐

Note that from the previous lemmas it follows that $\mathbf{sWeb}^{\uparrow,\downarrow}(\mathfrak{gl}_n)$ is a braided category. In particular we have that the $(l,k)$ overcrossing and the $(k,l)$ undercrossing are inverses to each other. Another direct consequence of the above discussion is the following.

**Corollary 11C.16.** *The leftward and rightward thick crossings are inverse to each other.* □

**Proposition 11C.17.** *For any two objects $\vec{k}$ and $\vec{l}$ of $\mathbf{sWeb}^{\uparrow,\downarrow}(\mathfrak{gl}_n)$ there exist objects $\vec{k}'$ and $\vec{l}'$ (defined in the proof) of $\mathbf{sWeb}^{\uparrow}(\mathfrak{gl}_n)$ such that there is an isomorphism of $\mathbb{k}$-vector spaces*

$$\mathrm{Hom}_{\mathbf{sWeb}^{\uparrow,\downarrow}(\mathfrak{gl}_n)}(\vec{k},\vec{l}) \cong \mathrm{Hom}_{\mathbf{sWeb}^{\uparrow}(\mathfrak{gl}_n)}(\vec{k}',\vec{l}').$$

*Proof.* Let w be a web in $\mathrm{Hom}_{\mathbf{sWeb}^{\uparrow,\downarrow}(\mathfrak{gl}_n)}(\vec{k},\vec{l})$. Invertibility of crossings (Lemma 11C.14 and Corollary 11C.16) and the invertibility of zigzag relations Equation 11B-26 imply that the usual isomorphism in braided pivotal categories holds. In pictures, this isomorphism looks e.g.

$\mathrm{w} \mapsto$ 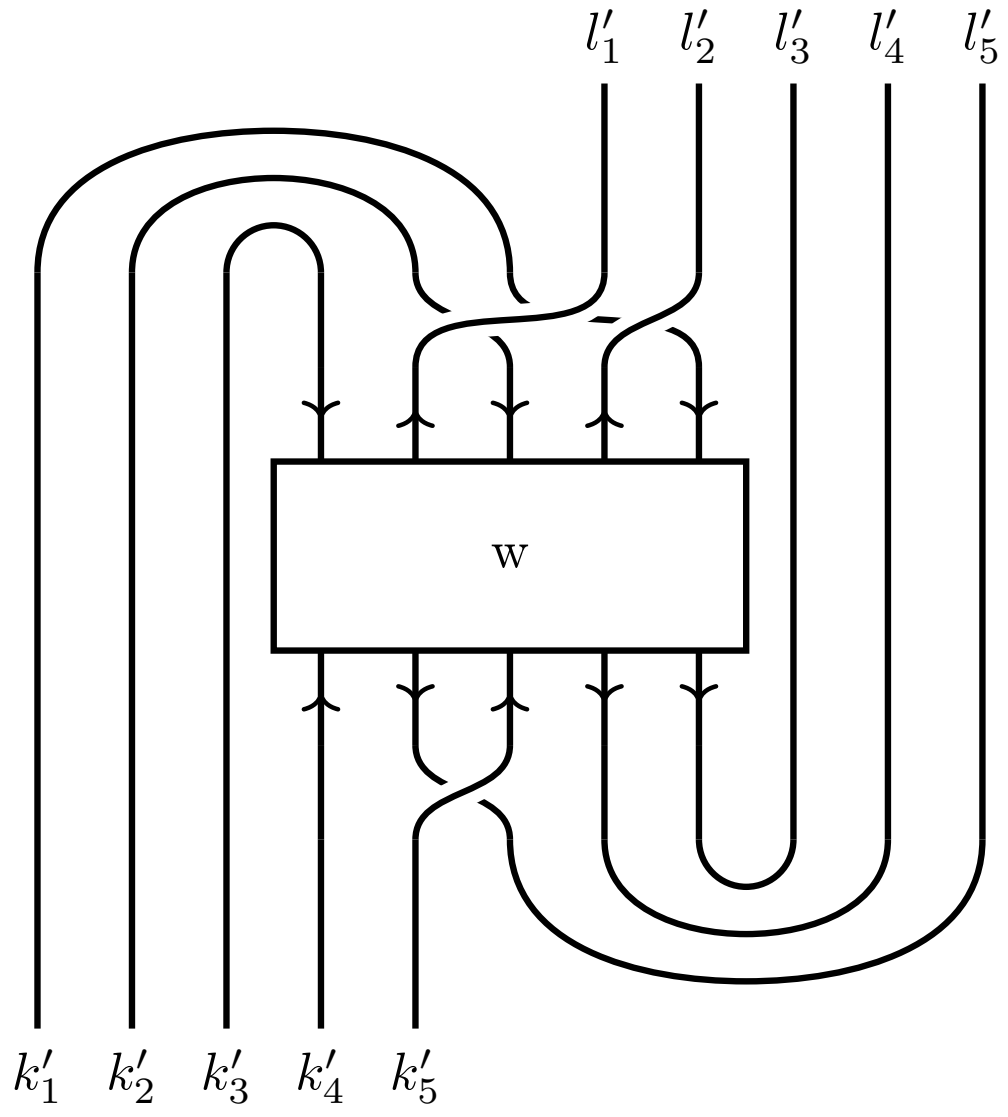 ,

which is evidently an invertible operation. □

*Remark* 11C.18. Note that Proposition 11C.17 only uses invertibility of crossings and the zigzag relations. In particular, it also holds for other web categories with duals and invertible crossings which are going to be discussed in this section, not just $\mathbf{sWeb}^{\uparrow,\downarrow}(\mathfrak{gl}_n)$. ◇

The splits and merges, when composed with the corresponding crossings, are "mirrored" as we will explain in the following lemma.

**Lemma 11C.19.** *We have*

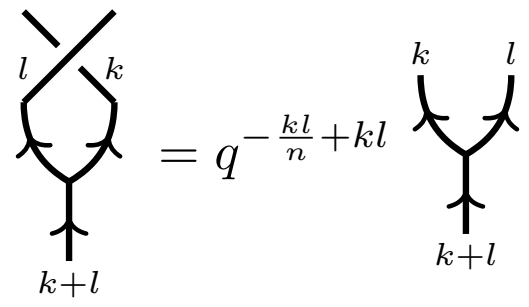

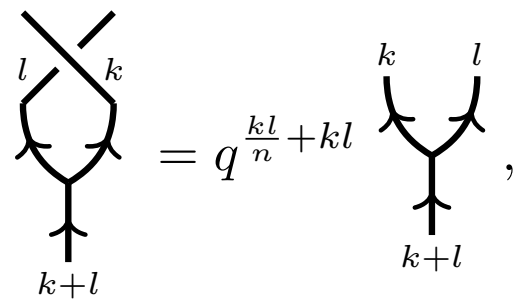

*and merges have the same behavior.*

*Proof.* We prove this via the exploded crossing formula in Definition 11B.4. First, the case $k = l = 1$ can be verified via a direct computation using Equation 11B-3. The relevant diagrammatic calculation is now

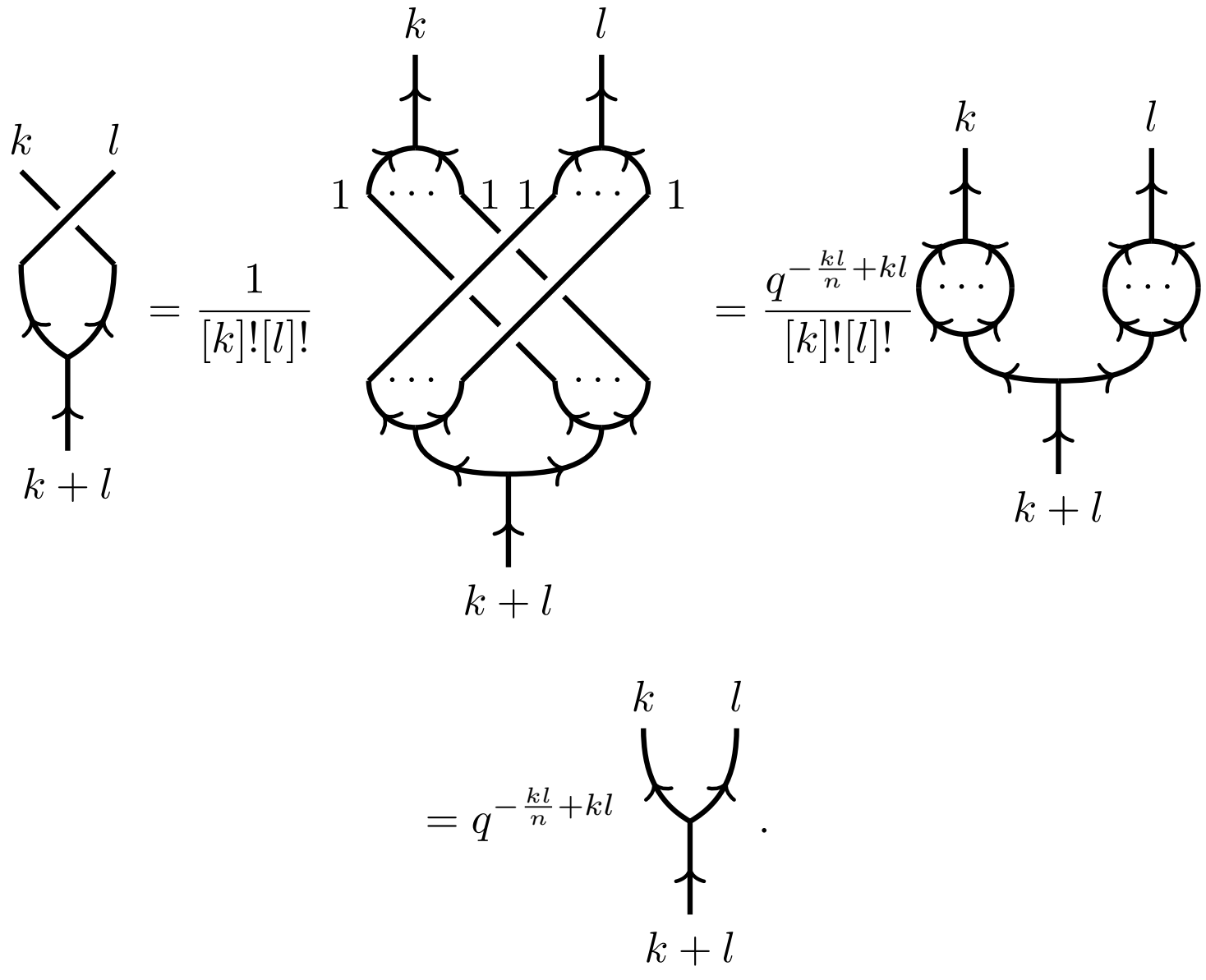

Let us explain the second equality. Successively, using associativity in the bottom bouquet we can always match a crossing with a corresponding split, giving a factor $q^{-\frac{1}{n}+1}$. As there are $kl$ crossings to be absorbed, then the overall factor is $q^{-\frac{kl}{n}+kl}$. All other cases follow mutatis mutandis. □

Further, we have the next two lemmas.

**Lemma 11C.20.** *The following relation of Reidemeister 1 type*

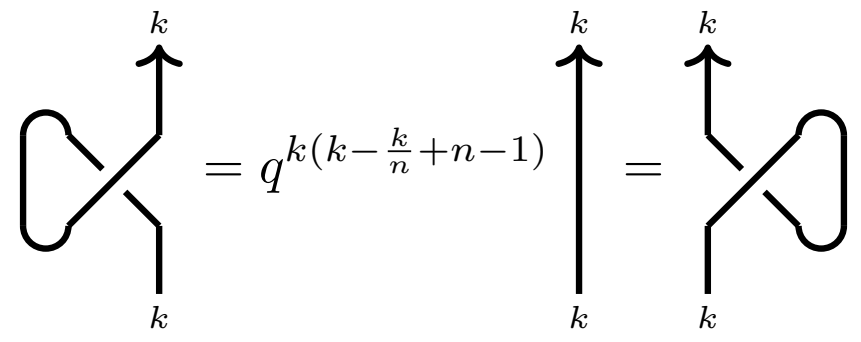

*as well as its downward pointing version and*

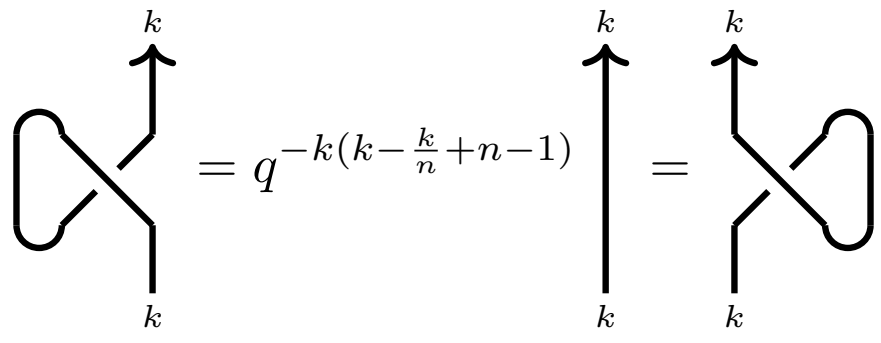

*as well as its downward pointing version hold.*

*Proof.* Using explosion, the previous lemma and the various established topological properties, the main thing to observe is that a Reidemeister 1 picture with $k$ parallel thin strands can be topologically rewritten as follows

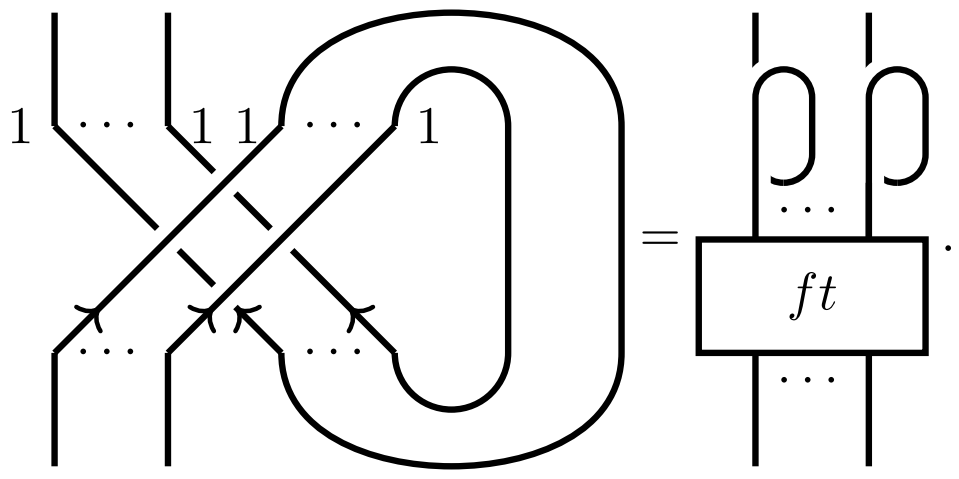

Here $ft$ is the full twist on $k$ strands which has $k(k-1)$ crossings. Now we first need to absorb the full twist using the previous lemma giving a factor $q^{-\frac{k(k-1)}{n}+k(k-1)}$. Then we use the thin Reidemeister 1 move from Lemma 11C.4 to get the factor $q^{k(n-\frac{1}{n})}$. Thus, the overall factor is $q^{k(k-\frac{k}{n}+n-1)}$. All other cases follow mutatis mutandis. □

All the above implies the following.

**Corollary 11C.21.** *The category* $\mathbf{sWeb}^{\uparrow,\downarrow}(\mathfrak{gl}_n)$ *is a ribbon category.* □

**Definition 11C.22.** Let $k \in \mathbb{Z}_{>0}$. The $k$-th ***symmetric*** $\mathfrak{gl}_\mathfrak{n}$***-projector*** $e_k$ is defined via

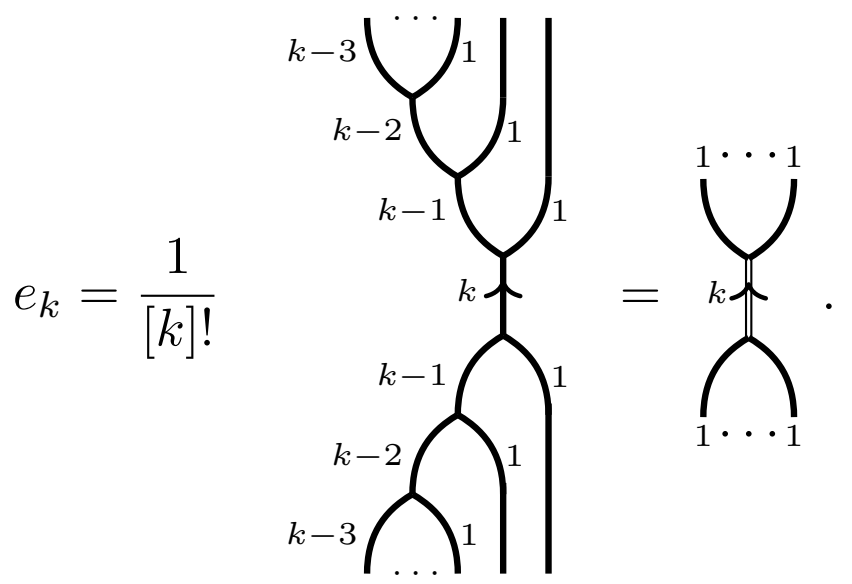

The final, rightmost, diagram means we sometimes draw a very thick edge for this projector. ◇

The following is almost the same as Lemma 10E.9.(iv). However, as a word of warning, the scalars are not the same as in the classical Jones–Wenzl recursion.

**Lemma 11C.23.** *We have:*

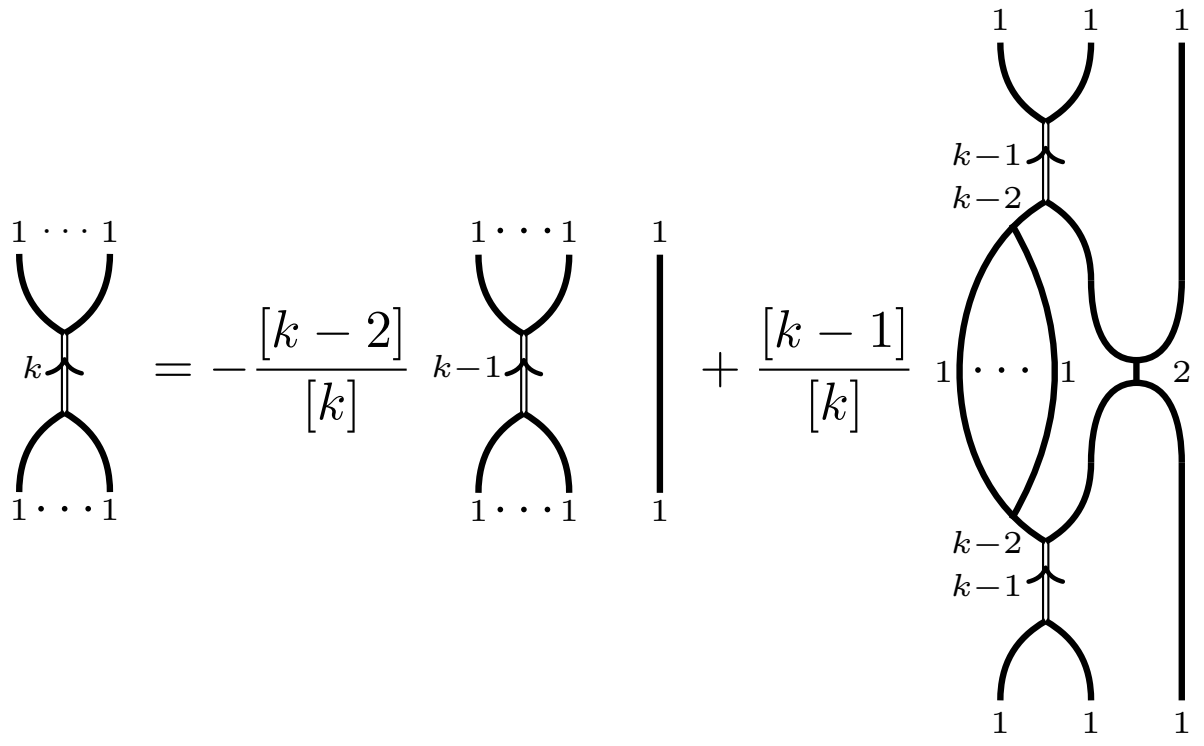

*for all* $k \in \mathbb{Z}_{\geq 3}$.

*Proof.* The proof is similar to that of [RT16, Lemma 2.13]. □

11D. **Webs integrally.** We now point out that webs can be adjusted to work integrally, meaning over $\mathbb{A}$. By specialization, this could provide insights into the representation theory in positive characteristic.

**Definition 11D.1.** The ***thick upward oriented labeled generic web category without relations***, denote by $\mathbf{lWeb}^{\uparrow,\mathrm{int}}$, is monoidally generated as follows. We let $\mathrm{S} = \{k^\uparrow \mid k \in \mathbb{Z}_{\geq 0}\}$, and T is given by the morphisms

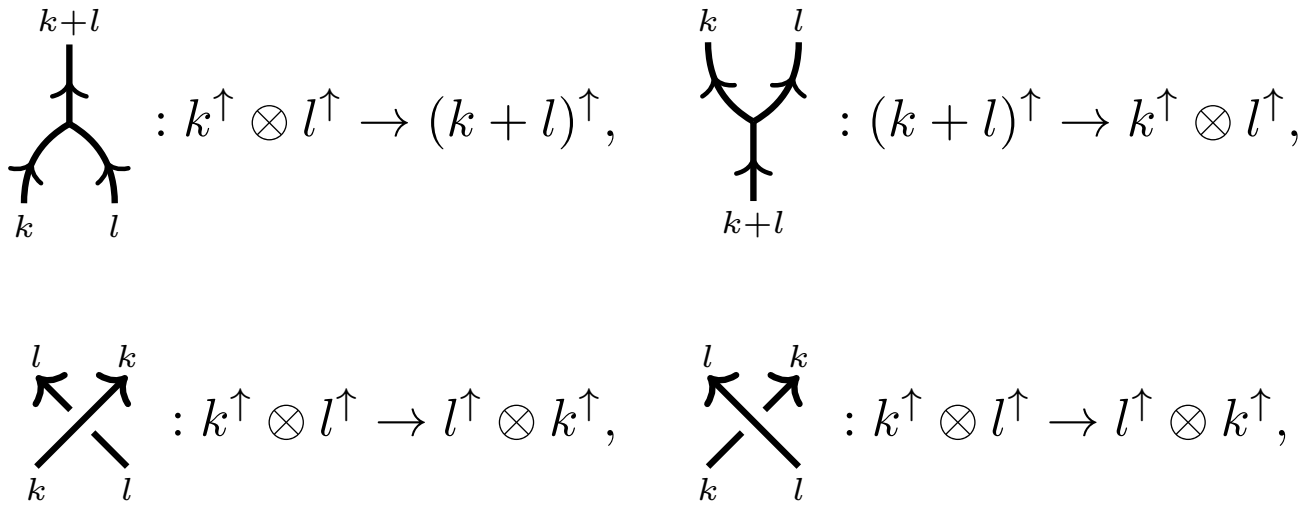

for all $k, l \in \mathbb{Z}_{\geq 0}$. ◇

**Definition 11D.2.** The ***integral symmetric upward oriented labeled generic web category***, denoted by $\mathbf{sWeb}^{\uparrow,\mathrm{int}} = \mathbf{sWeb}^{\uparrow,\mathrm{int}}_{\Bbbk\oplus}$, is the quotient of $\mathbf{lWeb}^{\uparrow,\mathrm{int}}_{\Bbbk\oplus}$ by the following relations, displaying R:

***Associativity and coassociativity***:

(11D-3)

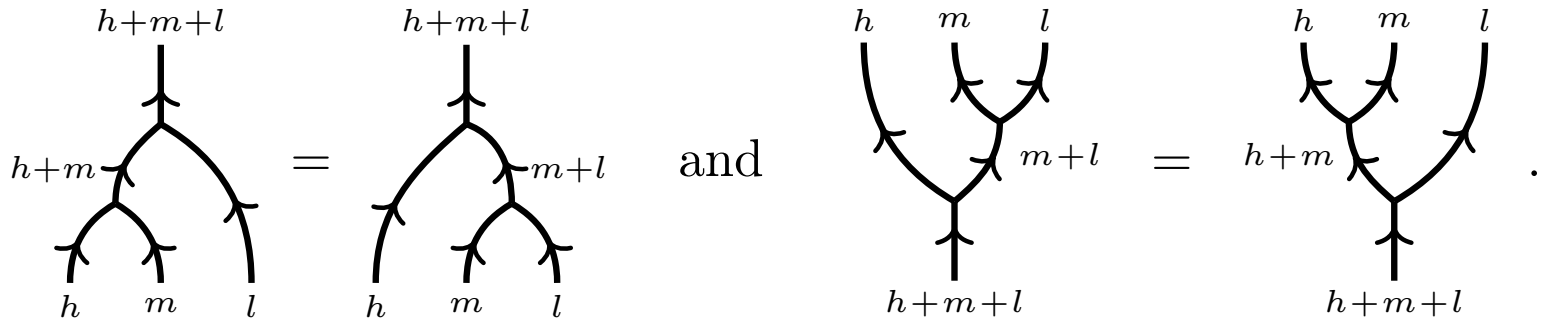

***Digon removal***:

(11D-4)

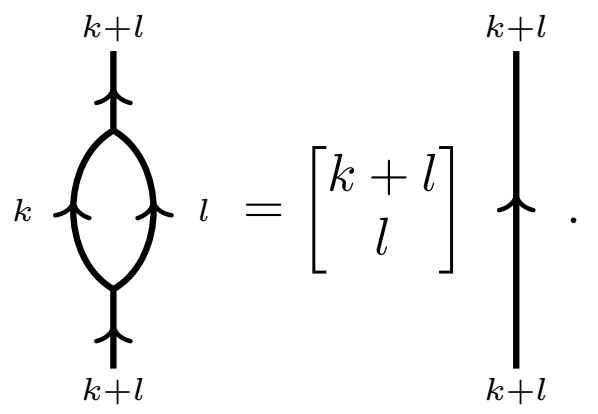

Dumbbell-crossing relations, also called ***Schur relations***:

(11D-5)

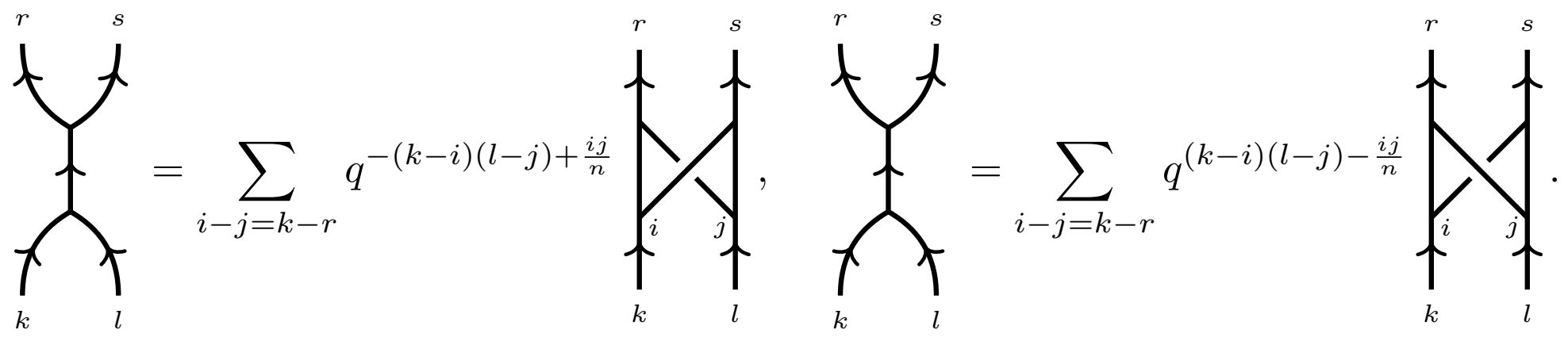

Note that we included relations of arbitrary thickness. ◇

*Remark* 11D.6. For $q = 1$ the relations in Equation 11D-5 are motivated by [**BEAEO20**, (4.23)], or, not diagrammatic, these can be deduced from Green's book on the Schur algebra [**Gre80**] via an interpretation of webs as elements in the Schur algebra. ◇

*Remark* 11D.7. The crossings can be omitted as generators, as we will show in Lemma 11D.9. We only use them above to give a shorter presentation. Note also that we do ***not*** assume that the crossings are invertible. This will be a later consequence. We will actually prove that $\mathbf{sWeb}^{\uparrow,\mathrm{int}}(\mathfrak{gl}_n)$ is a ribbon category. ◇

Now, since we are working over $\mathbb{A}$, the digon removal relation Equation 11D-4 is ***not*** an invertible operation. Consequently we cannot use the idea of exploding edges (from the previous chapters) since this requires the inversion of quantum binomials. This is the reason why in Definition 11D.2 our web generators are thick to begin with.

We point out here that a lot of the properties of generic webs hold also integrally. One of these is the following.

**Lemma 11D.8.** *Thick square switches hold in* $\mathbf{sWeb}^{\uparrow,int}(\mathfrak{gl}_n)$*. Namely,*

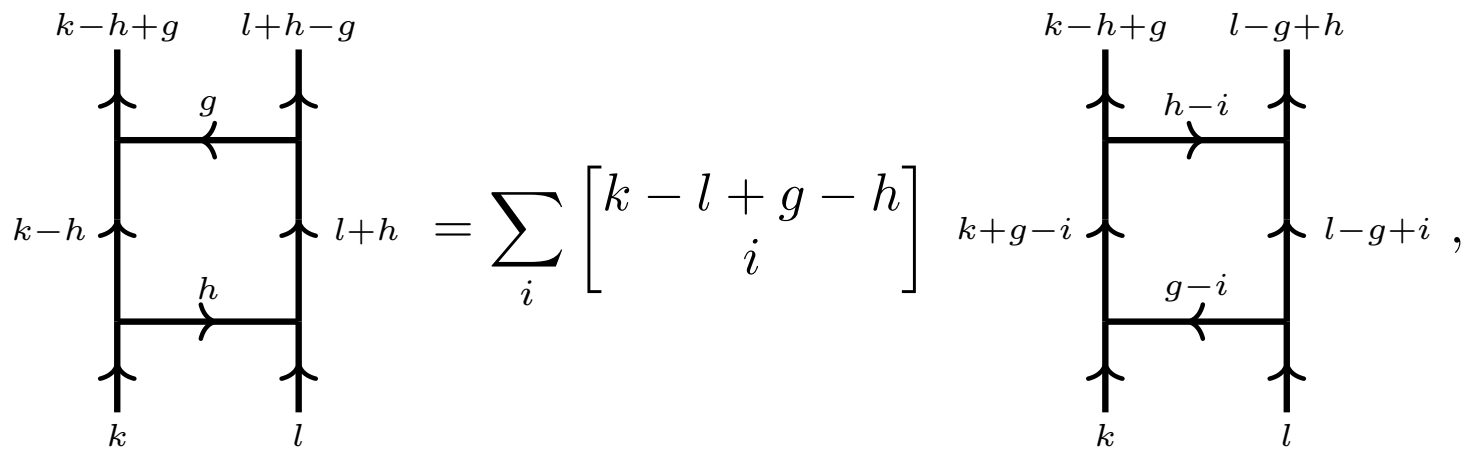

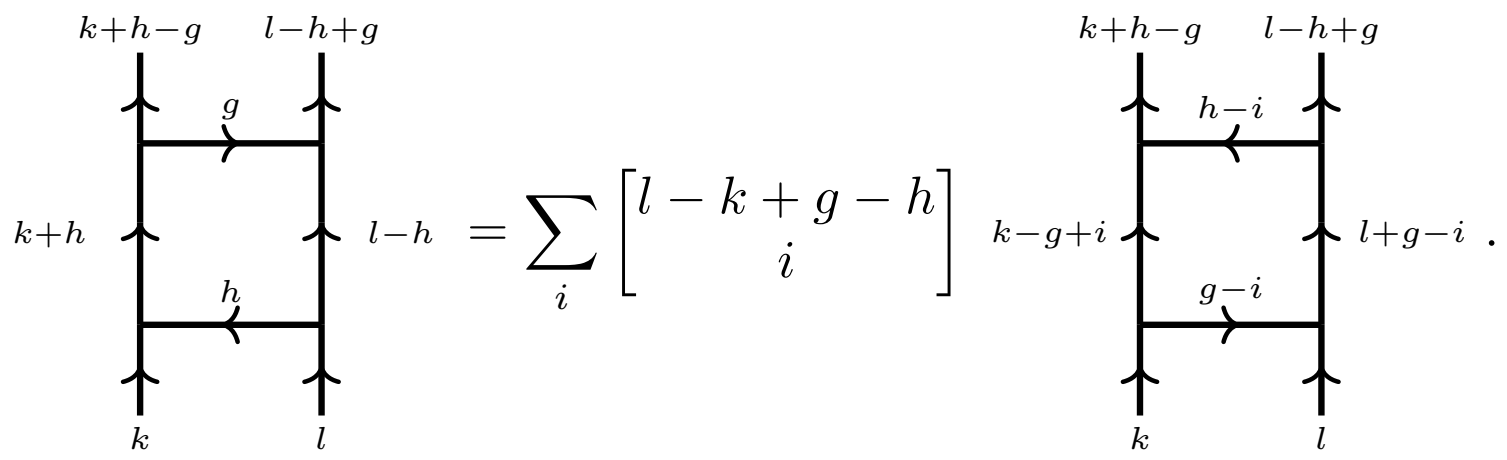

*Proof.* These can be derived from relations Equation 11D-3-Equation 11D-5 and $q$-combinatorics. Indeed, we use relation Equation 11D-5 to write the LHS of the first relation as

$$\sum_{t=\max(0,g-l)}^{\min(g,h)} q^{-t(l-g+t)+\frac{(h-t)(g-t)}{n}}$$

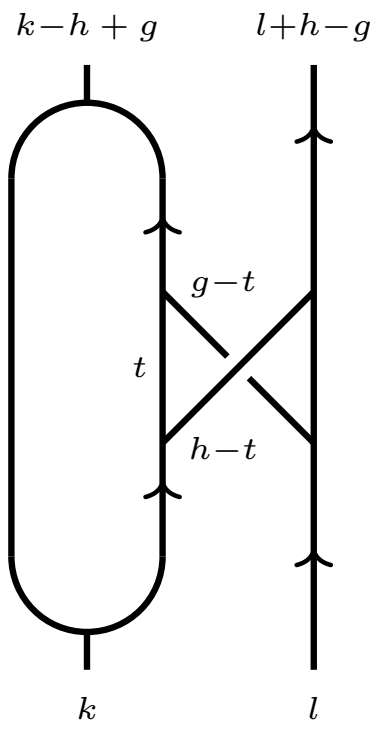

$$\overset{\text{Equation 11D-3 and Equation 11D-4}}{=} \sum_{t=\max(0,g-l)}^{\min(g,h)} q^{-t(l-g+t)+\frac{(g-t)(h-t)}{n}} \begin{bmatrix} k-h+t \\ t \end{bmatrix}$$

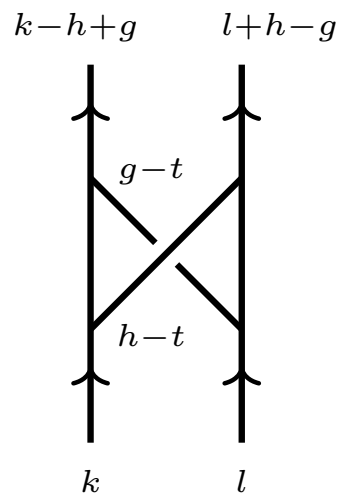

.

Similarly, the RHS is equal to

$$\sum_{i=0}^{\min(g,h)} \begin{bmatrix} k-l+g-h \\ i \end{bmatrix} \sum_{s=0}^{\min(g,h)-i} q^{-s(k-h+i+s)+\frac{(g-i-s)(h-i-s)}{n}}$$

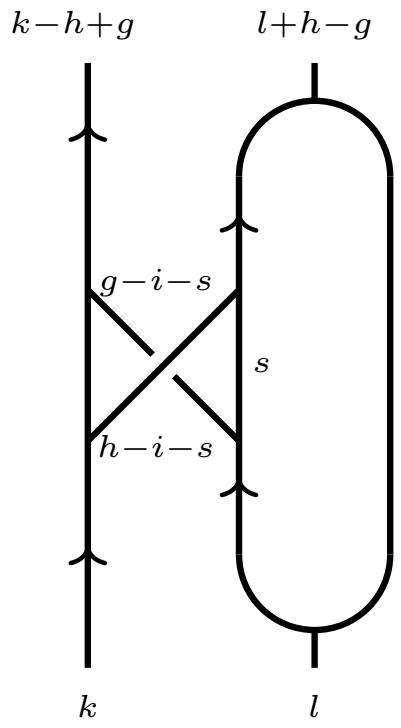

$$= \sum_{i=0}^{\min(g,h)} \begin{bmatrix} k-l+g-h \\ i \end{bmatrix} \sum_{s=0}^{\min(g,h)-i} q^{-s(k-h+i+s)+\frac{(g-i-s)(h-i-s)}{n}} \begin{bmatrix} l-g+i+s \\ s \end{bmatrix}$$

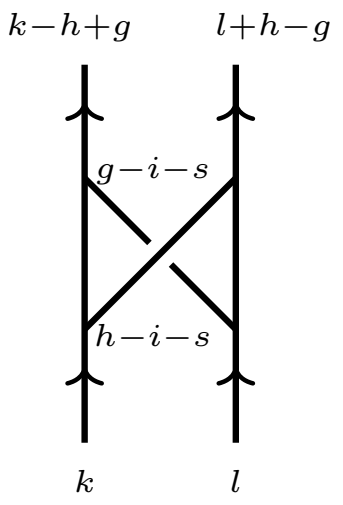

$$= \sum_{i=\max(0,g-l)}^{\min(g,h)} \begin{bmatrix} k-h+i+s \\ i+s \end{bmatrix} q^{-(i+s)(l-g+i+s)+\frac{(g-i-s)(h-i-s)}{n}}$$

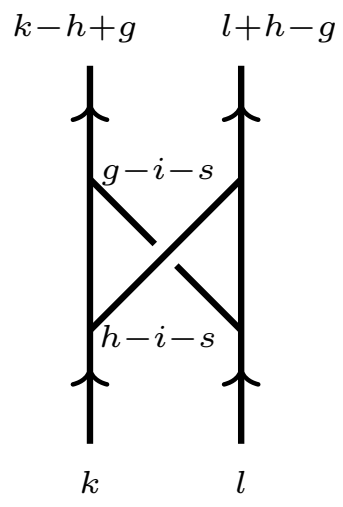

$$\overset{\text{put } i+s=t}{=} \sum_{t=\max(0,g-l)}^{\min(g,h)} q^{-t(l-g+t)+\frac{(g-t)(h-t)}{n}} \begin{bmatrix} k-h+t \\ t \end{bmatrix}$$ 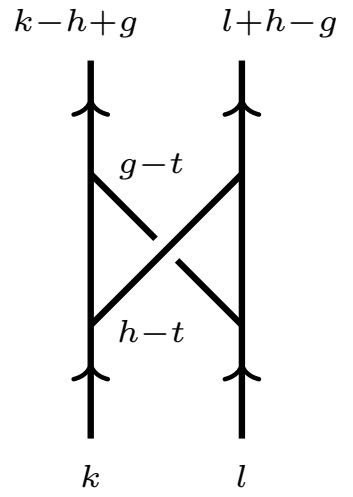 .

The other square switch relation is proved analogously. □

**Lemma 11D.9.** *The crossings from Definition 11D.2 satisfy the following identities*

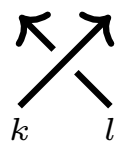 $= q^{-\frac{kl}{n}} \sum_{t=0}^{min(k,l)} (-q)^{-t}$ 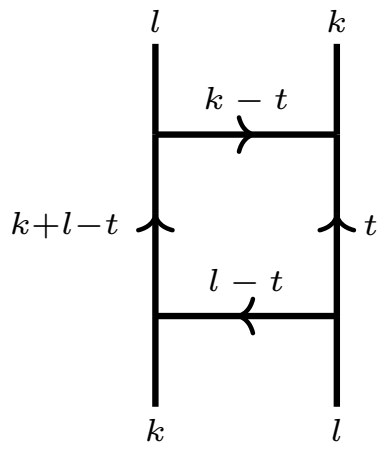 ,

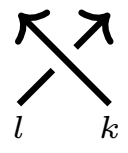 $= q^{\frac{kl}{n}} \sum_{t=0}^{min(k,l)} (-q)^{t}$ 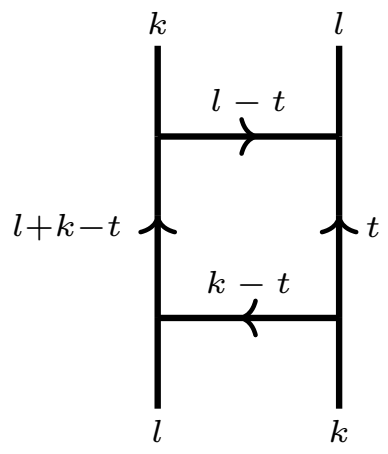 ,

*meaning that they are redundant as generators of* $\mathbf{sWeb}^{\uparrow,int}(\mathfrak{gl}_n)$.

The proof of Lemma 11D.9 is similar to that of Lemma 11D.10. We prove the latter at the end of the section.

**Lemma 11D.10.** *We can deduce Equation 11D-5 from Equation 11D-3, Equation 11D-4 and the thick square switches from Lemma 11D.8.*

*Proof.* We show that Equation 11D-3, Equation 11D-4 and thick square switches from Lemma 11D.8 imply the first equality in Equation 11D-5. The second equality can be proved similarly. At this point, we give the quantum version of [**BEAEO20**, Appendix A]. First, note that up to symmetry, it suffices to prove that the following holds in $\mathbf{sWeb}^{\uparrow}(\mathfrak{gl}_n)$, for $r \geq 0$;

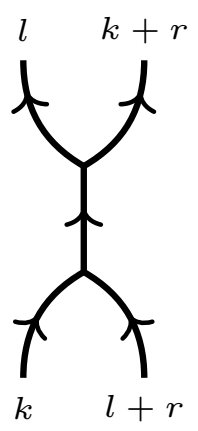 $= \sum_{i-j=k-l} q^{-(k-i)(l+r-j)+\frac{ij}{n}}$ 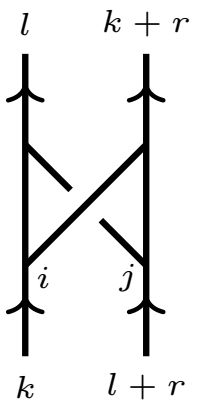 .

For simplicity, we relabel some of the edges above and the final version of relation Equation 11D-5 that we are aiming to prove is

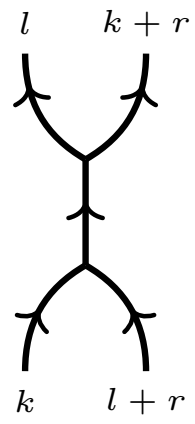 $= \sum_{s=0}^{\min(k,l)} q^{-s(s+r)+\frac{(k-s)(l-s)}{n}}$ 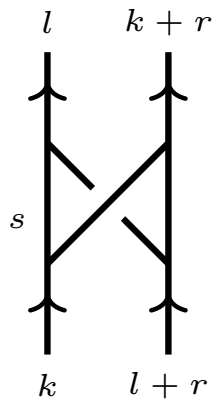 .

Using Lemma 11D.9 (which is proved independently) for the $(k-s, l-s)$ braiding, we have that the right-hand side of the above relation is equal to

$$\sum_{s=0}^{\min(k,l)} q^{-s(s+r)} \sum_{t=0}^{\min(k,l)-s} (-q)^{-t}$$

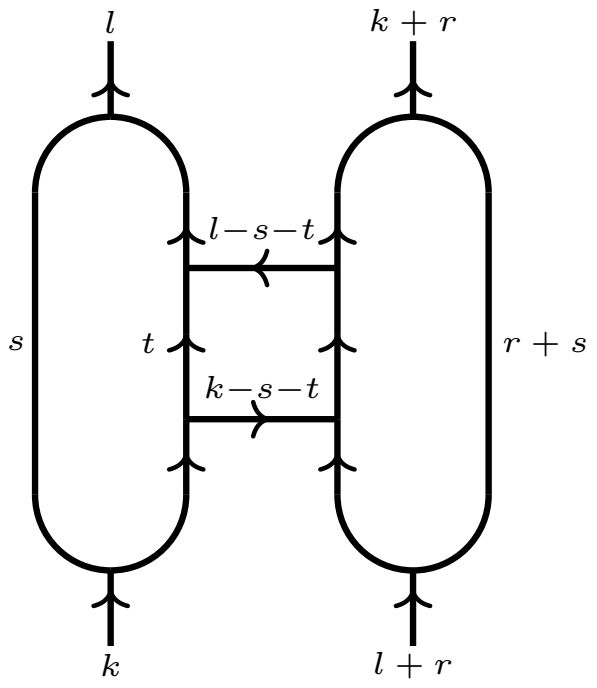

$$\overset{\textit{Equation 11D−3, Equation 11D−4}}{=} \sum_{s=0}^{\min(k,l)} q^{-s(s+r)} \sum_{t=0}^{\min(k,l)-s} (-q)^{-t} \begin{bmatrix} s+t \\ s \end{bmatrix}$$

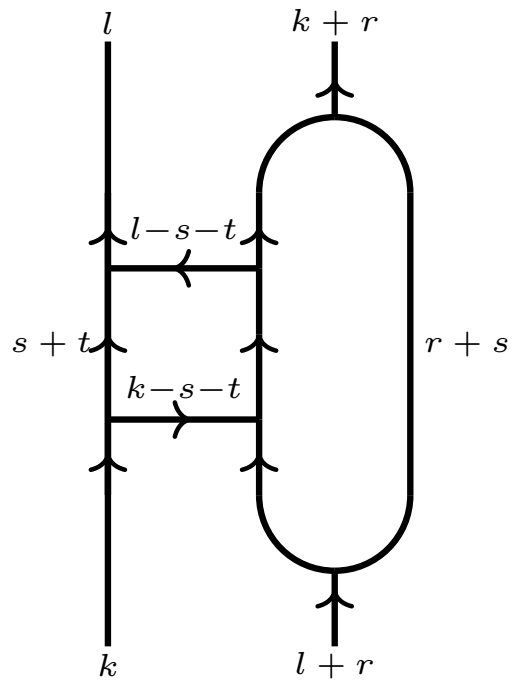

$$\overset{\text{put } s+t=u}{=} \sum_{s=0}^{\min(k,l)} q^{-s(s+r)} \sum_{u=s}^{\min(k,l)} (-q)^{s-u} \begin{bmatrix} u \\ s \end{bmatrix}$$

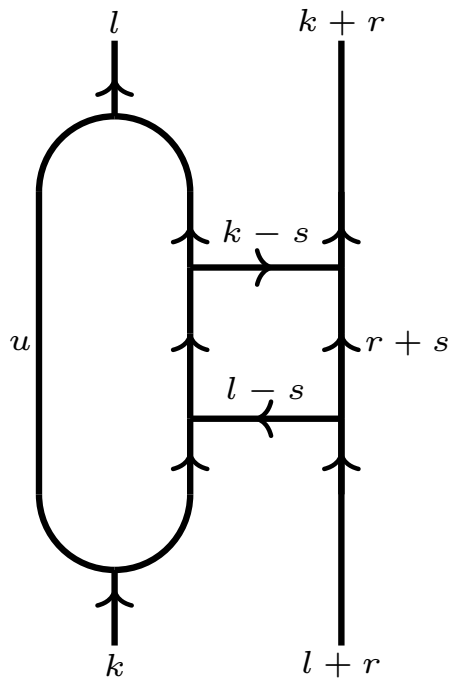

$$\overset{\text{square switch}}{=} \sum_{s=0}^{\min(k,l)} q^{-s(s+r)} \sum_{u=s}^{\min(k,l)} (-q)^{s-u} \begin{bmatrix} u \\ s \end{bmatrix} \sum_{t=u-s}^{\min(k,l)-s} \begin{bmatrix} u+r \\ t \end{bmatrix}$$

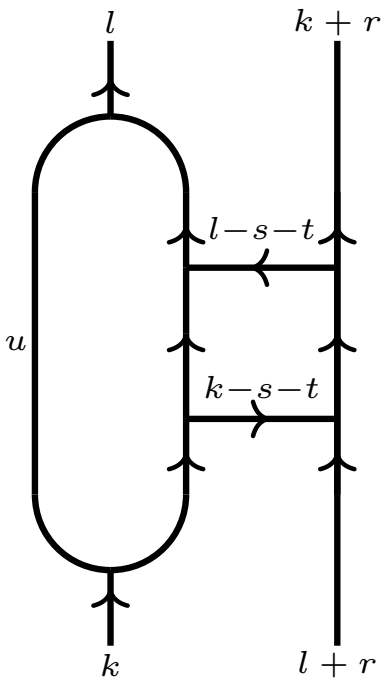

$$\overset{\textit{Equation 11D−3, Equation 11D−4}}{=} \textstyle\sum_{s=0}^{\min(k,l)} q^{-s(s+r)} \sum_{u=s}^{\min(k,l)} (-q)^{s-u} \begin{bmatrix} u \\ s \end{bmatrix} \sum_{t=u-s}^{\min(k,l)-s} \begin{bmatrix} u+r \\ t \end{bmatrix} \begin{bmatrix} s+t \\ u \end{bmatrix}$$

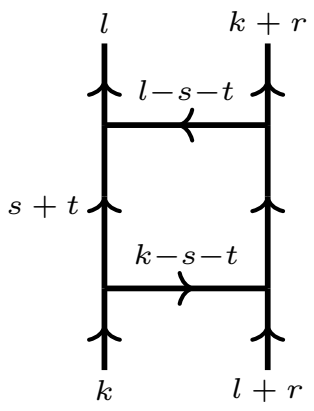

$$\overset{\text{put } v=s+t}{=} \sum_{s=0}^{\min(k,l)} q^{-s(s+r)} \sum_{u=s}^{\min(k,l)} (-1)^{s+u} q^{s-u} \begin{bmatrix} u \\ s \end{bmatrix} \sum_{v=u}^{\min(k,l)} \begin{bmatrix} u+r \\ v-s \end{bmatrix} \begin{bmatrix} v \\ u \end{bmatrix}$$

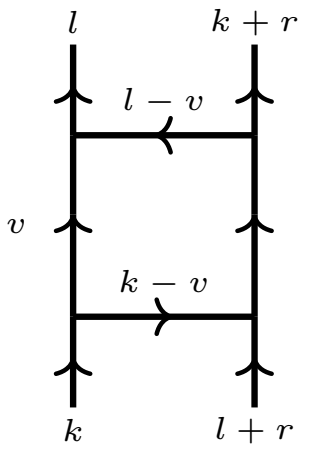

$$= \sum_{s=0}^{\min(k,l)} \sum_{u=s}^{\min(k,l)} \sum_{v=u}^{\min(k,l)} (-1)^{s+u} q^{-s(s+r)} q^{s-u} \begin{bmatrix} u \\ s \end{bmatrix} \begin{bmatrix} u+r \\ v-s \end{bmatrix} \begin{bmatrix} v \\ u \end{bmatrix}$$

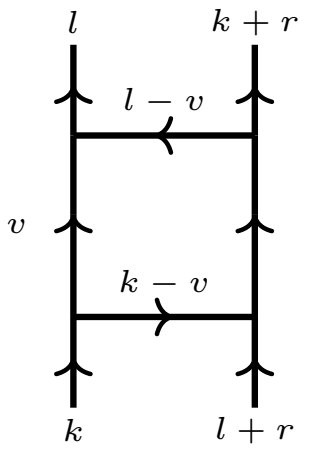

.

Switching the order of summations, while keeping in mind that $(-1)^{s+u} = (-1)^v(-1)^s(-1)^{u-v}$, the last sum is equal to

$$\sum_{v=0}^{\min(k,l)} (-1)^v \sum_{s=0}^{v} (-1)^s \sum_{u=s}^{v} (-1)^{u-v} q^{-s(s+r)+s-u} \begin{bmatrix} v \\ s \end{bmatrix} \begin{bmatrix} u+r \\ u-s, v-u \end{bmatrix}$$

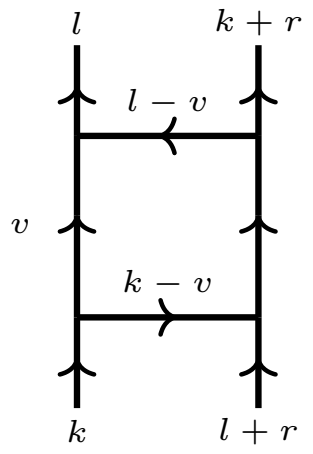

where we used the quantum trinomial coefficient. Further, ignoring the web, we write the last sum as

$$\sum_{v=0}^{\min(k,l)} (-1)^v \sum_{s=0}^{v} (-1)^s q^{vs-s} \begin{bmatrix} v \\ s \end{bmatrix} \sum_{u=s}^{v} (-1)^{u-v} q^{-s(s+r)+s-u-vs+s} \begin{bmatrix} u+r \\ u-s, v-u \end{bmatrix}$$

and this is equal to

$$\sum_{v=0}^{\min(k,l)} (-1)^v \delta_{v,0} \sum_{u=s}^{v} (-1)^{u-v} q^{-s(s+r+v-2)-u} \begin{bmatrix} u+r \\ u-s, v-u \end{bmatrix}.$$

This sum is nontrivial if and only if $v = 0$. In that case also $u = s = 0$ and

$$\sum_{u=s}^{v} (-1)^{u-v} q^{-s(s+r+v-2)-u} \begin{bmatrix} u+r \\ u-s, v-u \end{bmatrix} = 1,$$

meaning that the right-hand side of the initial relation is equal to

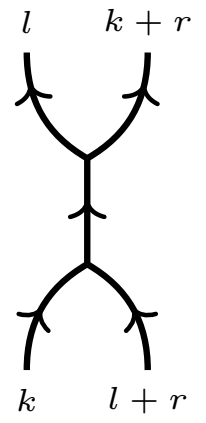

.

This finishes the proof that the relation Equation 11D-5 can be deduced from Equation 11D-3, Equation 11D-4 and thick square switches from Lemma 11D.8. □

We stress that Lemma 11D.8 and Lemma 11D.10 imply that the dumbbell-crossing and the thick square switch are equivalent relations in a certain sense. Let us also stress out that if the underlying field of $\mathbf{sWeb}^{\uparrow,\mathrm{int}}(\mathfrak{gl}_n)$ would be $\Bbbk$ then we could use thinner merges and splits as generators. Then we would be able to define thick merges, splits and braidings via explosion of strings as in the previous sections. Combining this with Lemma 11D.8 and Lemma 11D.10 we have the following.

**Corollary 11D.11.** *If we work over $\Bbbk$ (that is, any field where the quantum binomials are invertible), then $\mathbf{sWeb}^{\uparrow,int}(\mathfrak{gl}_n)$ is equivalent to $\mathbf{sWeb}^{\uparrow}(\mathfrak{gl}_n)$.* □

11E. **Exterior and red-green webs.** From now on, we assume some familiarity with the cited literature. The reader might just want to jump ahead to Section 11G, which is essentially (but not 100%) independent of the following two sections.

All of the above is formulated for symmetric webs, but exterior webs as in [CKM14] or red-green webs as in [TVW17] (red-green paper) work perfectly well using the same strategy. There is one minor difference in the relations: the crossings are defined differently when braiding exterior powers:

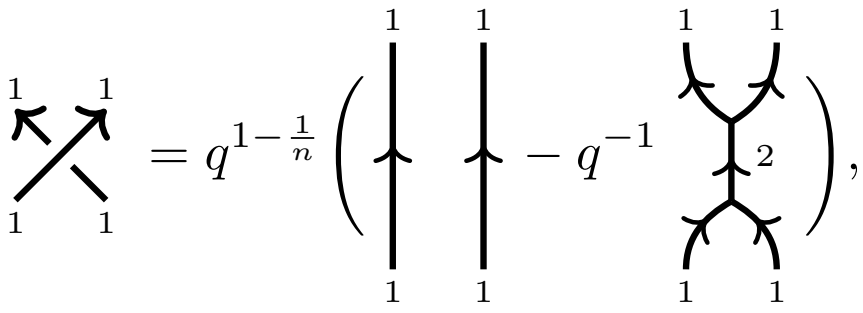

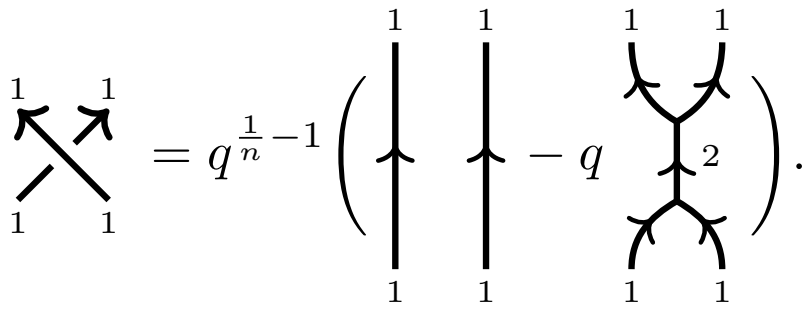

This apparently minor difference has huge impact, in particular, integrally. For example, while exterior webs fit perfectly to tilting modules (see e.g. [Eli15] and [AST18]), symmetric webs do not have any such interpretation known to the author.

More details on minimal and integral exterior and red-green webs can be found in [LT21].

*Remark* 11E.1. From now on we will assume that the reader is reasonably familiar with exterior webs, *cf.* Exercise 11I.2 below. ◇

11F. **A basis for thin webs.** Recall the exterior $\mathfrak{sl}_n$-webs from [CKM14]. We point out here that for $U_q(\mathfrak{gl}_n)$ we have $\bigwedge_q^n(\mathbb{k}_q^n) \not\cong \mathbb{k}_q$, whereas for $U_q(\mathfrak{sl}_n)$ we have $\bigwedge_q^n(\mathbb{k}_q^n) \cong \mathbb{k}_q$. ($\bigwedge_q^n(\mathbb{k}_q^n)$ is the quantum analog of the determinant representation of the general linear group.) So, for exterior $\mathfrak{sl}_n$-webs the $n$-labeled edges are 0-labeled and thus are removed from the illustrations.

In [Eli15], Elias constructed a basis for exterior $\mathfrak{sl}_n$-webs whose elements are called double light ladders. We will be considering only the basis of $\mathrm{End}_{\mathbf{eWeb}_{\uparrow,\downarrow}(\mathfrak{sl}_n)}(1_{\uparrow}^{\otimes k})$, $k \geq 1$. (The notation means exterior webs, see above.) We modify the later by performing an algorithm leading to a basis of $\mathrm{End}_{\mathbf{sWeb}^{\uparrow}(\mathfrak{gl}_n)}(1_{\uparrow}^{\otimes k})$. We color the thick edges of symmetric webs red and those of exterior webs green. These webs are related by the following relation

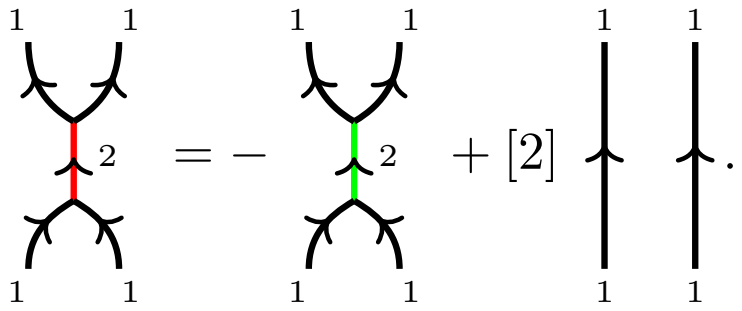

The three definitions from below will be used for both exterior and symmetric webs.

**Definition 11F.1.** A web with 1-labeled boundary edges is called a ***local dumbbell*** if it contains only dumbbells and identity components. ◇

**Lemma 11F.2.** *Let $X$ be a local dumbbell. Then $X$ can be written as a linear combination of webs containing only 2-dumbbells (dumbbells whose thickest label is 2) and identity components.*

*Proof.* This follows from the $\mathfrak{gl}_n$-projector recursion in Lemma 11C.23. □

**Definition 11F.3.** We call the linear combination above the ***thin expression*** of the local dumbbell $X$. ◇

**Definition 11F.4.** We say that a local dumbbell $X$ is ***smaller or equal*** to the local dumbbell $Y$, we write $X \leq Y$, if the total number of the 2-dumbbells appearing in the thin expression of $X$ is smaller or equal to the total number of the 2-dumbbells appearing in the thin expression of $Y$. ◇

Now, let $B_1, ..., B_l$ be the basis elements of $\mathrm{End}_{\mathbf{eWeb}^{\uparrow,\downarrow}(\mathfrak{sl}_n)}(1^{\otimes k})$ from [Eli15]. We perform the following two steps:

**Step 1**. In each of the webs $B_1, ..., B_l$ we change every cup and cap of the form

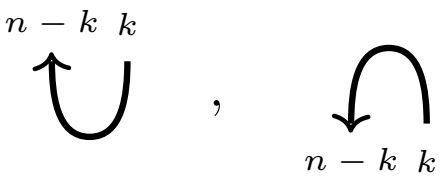

to splits and merges of the form

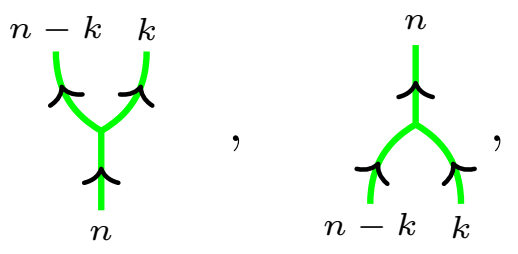

where the $n$-labeled edge has no self-intersections. Similarly we define the operation for the mirrored webs. In order for this operation to be well-defined we take the construction of webs as morphisms in a monoidal category given by generators–relations. Now, whenever this would lead to crossings of webs, we pull the $n$-labeled edges under and outside to the right like below:

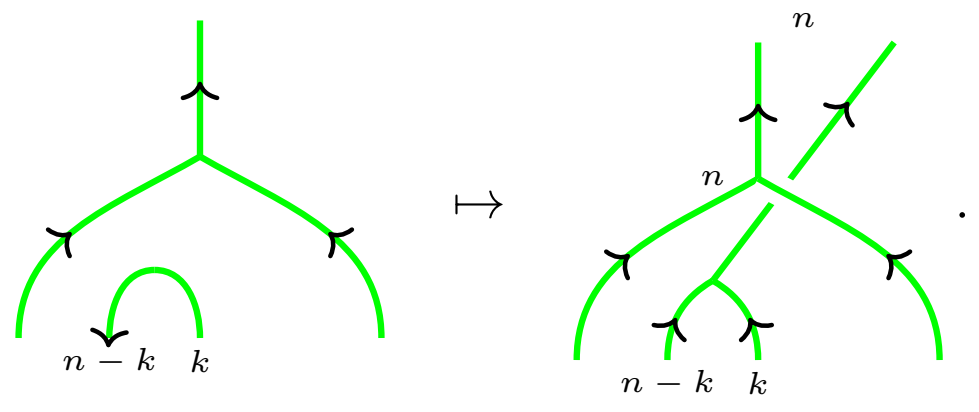

Note that here we are making a choice. The reader might wonder why this is well-defined and what happens if for example we let

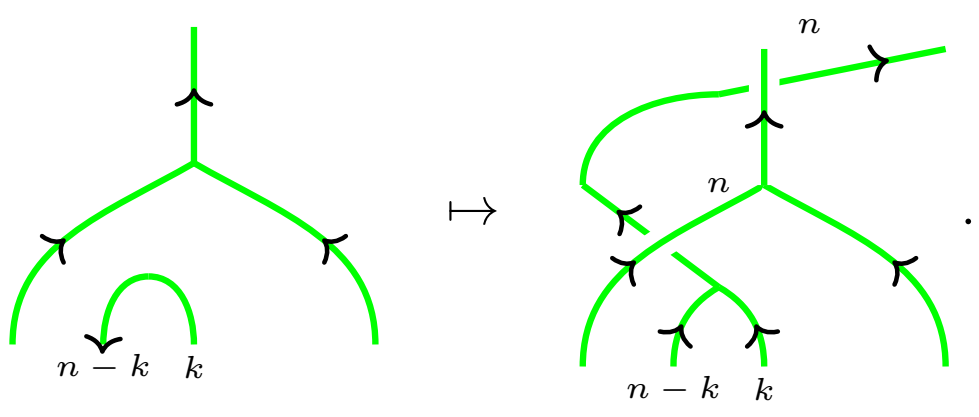

We point out that all such choices would lead to the same web due to the ambient isotopies such as zigzag relations, Reidemeister 2 and 3 moves, pitchfork relations, cup and cap slides etc. Moreover, if there are multiple caps, then we choose an arbitrary order of their end points noting that modulo lower terms.

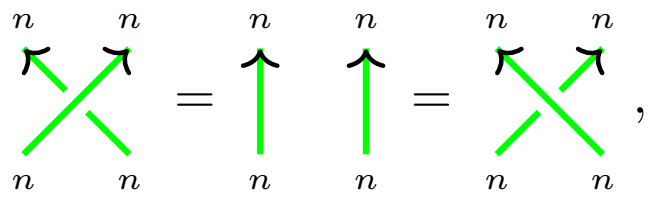

**Step 2**. We change the colors from green to red.

**Lemma 11F.5.** *If we apply Step 1 to the double ladder basis, then we obtain a basis for* $\mathrm{End}_{\mathbf{eWeb}_{\uparrow}(\mathfrak{gl}_n)}(1^{\otimes k})$.

*Proof.* Let $B_1', \dots, B_l'$ be the $\mathfrak{gl}_n$-webs obtained from the $\mathfrak{sl}_n$-webs $B_1, \dots, B_l$ (double ladders) after Step 1. Note that the reason why we are performing Step 1 is because for $\mathfrak{sl}_n$-webs we do not have edges of labels $\geq n$ and exterior $\mathfrak{sl}_n$-webs are a quotient of exterior $\mathfrak{gl}_n$-webs by

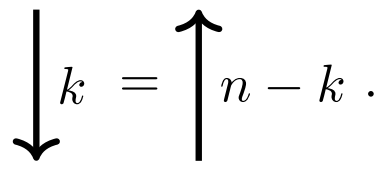

Now, the webs $B_1', \dots, B_l'$ are linearly independent. Indeed, if they were linearly dependent, after passing to $\mathfrak{sl}_n$-webs (by imposing the relation above) we would get that $B_1, \dots, B_l$ are linearly dependent, contradicting the fact that they form a basis. It remains to show that $B_1', \dots, B_l'$ span $\mathrm{End}_{\mathbf{eWeb}_{\uparrow}(\mathfrak{gl}_n)}(1^{\otimes k})$. This is a word for word adaption of the proof of [**Eli15**, Theorem 2.39]. Indeed, Elias points out that his proof would apply to any class of morphism satisfying [**Eli15**, Proposition 2.42]. This finishes the proof. □

**Proposition 11F.6.** *If we apply Step 1 and Step 2 to the the double ladder basis we obtain a basis* $\{\tilde{B}_1, \dots, \tilde{B}_l\}$ *for* $\mathrm{End}_{\mathbf{sWeb}^{\uparrow}(\mathfrak{gl}_n)}(1_{\uparrow}^{\otimes k})$.

*Proof.* Let $\{B_1', \dots, B_l'\}$ be the basis for $\mathrm{End}_{\mathbf{eWeb}^{\uparrow}(\mathfrak{gl}_n)}(1^{\otimes k})$ from the lemma above. We use the following strategy:

(a) First, we replace each $B_i'$ by a linear combination of webs containing only 2-dumbbells and identities locally.
(b) Second, we change the colors of the resulting webs.
(c) In these steps it is crucial that we have an order on webs given by counting 2-dumbbells, so both of the above will produce a triangular change-of-basis matrix.

We now give the details. By construction, each of the webs $B_i'$ is locally built from dumbbells and identities. Now we use the exterior version of $\mathfrak{gl}_n$-projector recursion from Lemma 11C.23 to write every $B_i'$ as a linear combination of webs built from only 2-dumbbells and identities. Without loss of generality we assume that $B_1' \leq \dots \leq B_l'$. Note that there is a biggest web in the $\mathfrak{gl}_n$-projector formula (Lemma 11C.23), namely the rightmost web. After performing all possible recursive steps, we will have one web which is the overall biggest in the thin expression of each $k$-dumbbell and consequently of each $B'_i$. Let $W_i$ be the biggest web appearing in the thin expression of $B'_i$. From the discussion above, the coefficient of $W_i$ in the thin expression of $B'_i$ is invertible. Inductively one concludes that the change-of-basis matrix from $B'_i$ to $W_i$ is triangular with invertible diagonal entries. This means that $W_1, ..., W_l$ are linearly independent and consequently they form a basis for $\mathrm{End}_{\mathbf{eWeb}^{\uparrow}(\mathfrak{gl}_n)}(1^{\otimes k})$. Also, note that $W_1 \leq ... \leq W_l$. After applying Step 2, i.e. changing the colors of the webs $W_1, ..., W_l$, we denote the new webs by $\tilde{W}_1, ..., \tilde{W}_l$. We claim that they are linearly independent. Indeed, using the relation we can write each symmetric web $\tilde{W}_i$ as a linear combination of exterior webs. Recall that by definition, the number of 2-dumbbell components in a web $\tilde{W}_i$ is the same as the number of 2-dumbbell components in $W_i$. So, $\tilde{W}_1 \leq ... \leq \tilde{W}_l$. Now, from the above relation we can write each $\tilde{W}_i$ as a linear combination of exterior basis webs $W_j$, such that $W_j \leq W_i$. Moreover, again from the same relation we have that the coefficient of $W_i$ in this linear combination is $(-1)^a$, where $a$ is the number of 2-dumbbell components of $W_i$. This means that the change-of-basis matrix from symmetric to exterior webs is triangular with determinant $\pm 1$. This concludes the proof that $\tilde{W}_1, ..., \tilde{W}_l$ are linearly independent. Since the algorithm sends exterior web generators to symmetric web generator, then $\{\tilde{W}_1, ..., \tilde{W}_l\}$ is a basis for $\mathrm{End}_{\mathbf{sWeb}^{\uparrow}(\mathfrak{gl}_n)}(1^{\otimes k})$. Finally, note that by construction, $\tilde{W}_i$ is the biggest web appearing in the thin expression of $\tilde{B}_i$. Thus, like before, the change-of-basis matrix from $\tilde{W}_i$ to $\tilde{B}_i$ is triangular with non-zero diagonal entries, meaning that $\{\tilde{B}_1, ..., \tilde{B}_l\}$ is a basis for $\mathrm{End}_{\mathbf{sWeb}^{\uparrow}(\mathfrak{gl}_n)}(1^{\otimes k})$. □

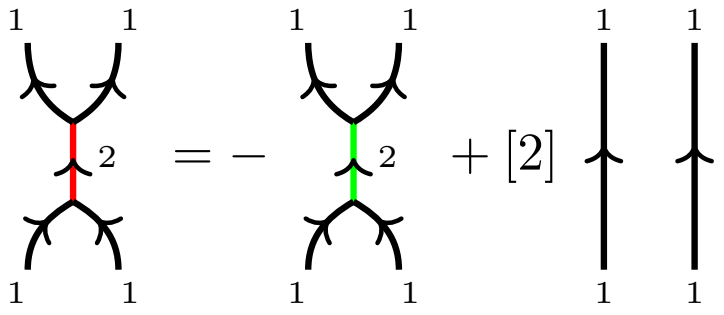

11G. **Equivalence theorem.** We now turn to representation theory and we explain how webs describe the category of finite dimensional modules over the quantum general linear group.

11G.1. *The Lie algebra for the general linear group.* To motivate the below, let us recall that

$$\mathfrak{gl}_m = \mathfrak{gl}_m(\mathbb{C}) = \{(a_{ij})_{i,j=1}^m | a_{ij} \in \mathbb{C}\} = m\text{-by-}m \text{ matrices.}$$
$$\mathfrak{sl}_m = \mathfrak{sl}_m(\mathbb{C}) = \{(a_{ij})_{i,j=1}^m | a_{ij} \in \mathbb{C}, a_{11} + ... + a_{mm} = 0\} = \text{traceless } m\text{-by-}m \text{ matrices.}$$

These are so-called ***Lie algebras*** with Lie bracket operation given by commutators:

$$[A, B] = AB - BA.$$

To be completely explicit, and without loosing much generality, let $m = 3$. Define the following matrices:

$$E_1 = \begin{pmatrix} 0 & 1 & 0 \\ 0 & 0 & 0 \\ 0 & 0 & 0 \end{pmatrix}, \quad E_2 = \begin{pmatrix} 0 & 0 & 0 \\ 0 & 0 & 1 \\ 0 & 0 & 0 \end{pmatrix}, \quad F_1 = \begin{pmatrix} 0 & 0 & 0 \\ 1 & 0 & 0 \\ 0 & 0 & 0 \end{pmatrix}, \quad F_2 = \begin{pmatrix} 0 & 0 & 0 \\ 0 & 0 & 0 \\ 0 & 1 & 0 \end{pmatrix},$$

often called ***Chevalley generators***. Moreover, define

$$\tilde{H}_1 = \begin{pmatrix} 1 & 0 & 0 \\ 0 & 0 & 0 \\ 0 & 0 & 0 \end{pmatrix}, \quad \tilde{H}_2 = \begin{pmatrix} 0 & 0 & 0 \\ 0 & 1 & 0 \\ 0 & 0 & 0 \end{pmatrix}, \quad \tilde{H}_3 = \begin{pmatrix} 0 & 0 & 0 \\ 0 & 0 & 0 \\ 0 & 0 & 1 \end{pmatrix}.$$

Under the Lie bracket, $\tilde{H}_i$, for $i \in \{1, 2, 3\}$, and $E_i, F_i$, for $i \in \{1, 2\}$, generate $\mathfrak{gl}_3$. Indeed, for example,

$$E_{12} = [E_1, E_2] = E_1E_2 - E_2E_1 = \begin{pmatrix} 0 & 0 & 1 \\ 0 & 0 & 0 \\ 0 & 0 & 0 \end{pmatrix}.$$

The Lie algebra $\mathfrak{sl}_3$ is generated by $E_i, F_i$, for $i \in \{1, 2\}$, and $H_1 = \tilde{H}_1 - \tilde{H}_2, H_2 = \tilde{H}_2 - \tilde{H}_3$.

Computing the various Lie brackets between these matrices, for example,

$$[E_1, F_1] = H_1,$$

as one easily checks. All of this is quantized below in the following informal sense. Let $q$ be a formal variable, and let $K_i = q^{H_i}$, viewed formally. Then

$$\lim_{q \to 1} \frac{K_i - K_i^{-1}}{q - q^{-1}} = H_i$$

and we get

$$E_{i_1}F_{i_2} - F_{i_2}E_{i_1} = \delta_{i_1,i_2} \frac{L_{i_1}L_{i_1+1}^{-1} - L_{i_1}^{-1}L_{i_1+1}}{q - q^{-1}} \text{ dequantizes to } [E_i, F_i] = H_i.$$

11G.2. *The quantum general linear group and its idempotented form.* We will use $\epsilon_i = (0, ..., 1, ..., 0) \in \mathbb{Z}^m$, where 1 is in the $i$th coordinate. For $i = 1, ..., m-1$, let $\alpha_i = \epsilon_i - \epsilon_{i+1} = (0, ..., 1, -1, ..., 0) \in \mathbb{Z}^m$,. Moreover, recall that the Euclidean inner product on $\mathbb{Z}^m$ is given by $(\epsilon_i, \epsilon_j) = \delta_{i,j}$.

**Definition 11G.1.** For $m \in \mathbb{Z}_{>1}$, the ***quantum general linear group*** $\mathrm{U}_q(\mathfrak{gl}_m)$ is the associative unital $\Bbbk(q)$-algebra generated by $L_i$ and $L_i^{-1}$, for $i = 1, ..., m$, and $E_i, F_i$, for $i = 1, ..., m-1$, subject to the following relations (for all admissible indices):

$$L_{i_1}L_{i_2} = L_{i_2}L_{i_1},$$
$$L_iL_i^{-1} = L_i^{-1}L_i = 1,$$
$$L_{i_1}E_{i_2} = q^{(\epsilon_{i_1},\alpha_{i_2})}E_{i_2}L_{i_1},$$
$$L_{i_1}F_{i_2} = q^{-(\epsilon_{i_1},\alpha_{i_2})}F_{i_2}L_{i_1},$$
$$E_{i_1}F_{i_2} - F_{i_2}E_{i_1} = \delta_{i_1,i_2}\tfrac{L_{i_1}L_{i_1+1}^{-1} - L_{i_1}^{-1}L_{i_1+1}}{q-q^{-1}}\ ,$$
$$E_{i_1}^2E_{i_2} - [2]E_{i_1}E_{i_2}E_{i_1} + E_{i_2}E_{i_1}^2 = 0, \text{ if } |i_1 - i_2| = 1,\ E_{i_1}E_{i_2} = E_{i_2}E_{i_1}, \text{ else},$$
$$F_{i_1}^2F_{i_2} - [2]F_{i_1}F_{i_2}F_{i_1} + F_{i_2}F_{i_1}^2 = 0, \text{ if } |i_1 - i_2| = 1,\ F_{i_1}F_{i_2} = F_{i_2}F_{i_1}, \text{ else}.$$

The last two relations are called the (quantum) ***Serre relations***. ◇

**Definition 11G.2.** For $m \in \mathbb{Z}_{>1}$ the ***quantum special linear algebra*** $\mathrm{U}_q(\mathfrak{sl}_m)$ is the subalgebra of $\mathrm{U}_q(\mathfrak{gl}_m)$ generated by $E_i, F_i, K_i = L_iL_{i+1}^{-1}$ and $K_i^{-1} = L_{i+1}L_i^{-1}$ for $i = 1, ..., m-1$. ◇

We can endow $\mathrm{U}_q(\mathfrak{gl}_m)$ with a Hopf algebra structure (which will be used throughout this section), where the coproduct $\Delta$ is given by

$$\Delta(E_i) = E_i \otimes L_iL_{i+1}^{-1} + 1 \otimes E_i,$$
$$\Delta(F_i) = F_i \otimes 1 + L_i^{-1}L_{i+1} \otimes F_i \text{ and}$$
$$\Delta(L_i^{\pm 1}) = L_i^{\pm 1} \otimes L_i^{\pm 1}.$$

The antipode and counit are given by

$$S(E_i) = -E_iL_i^{-1}L_{i+1},\ S(F_i) = -L_iL_{i+1}^{-1}F_i,\ S(L_i^{\pm 1}) = L_i^{\mp 1}$$
$$\epsilon(E_i) = \epsilon(F_i) = 0 \text{ and } \epsilon(L_i^{\pm 1}) = 1.$$

The subalgebra $\mathrm{U}_q(\mathfrak{sl}_m)$ inherits the Hopf algebra structure from $\mathrm{U}_q(\mathfrak{gl}_m)$. Recall that from the Hopf algebra structure follows the existence of the trivial representation. It also allows us to extend actions to tensor products of representations. Now, note that the antipode is invertible and the condition below is satisfied.

$$S^2(h) = K_{2\rho}hK_{2\rho}^{-1},$$

for every $h \in \mathrm{U}_q(\mathfrak{gl}_m)$, where $K_{2\rho} = K_1^{m-1}K_2^{m-2}...K_{m-1}$.

The following representation of $\mathrm{U}_q(\mathfrak{gl}_m)$ is of particular interest for us.

**Definition 11G.3.** The ***standard representation*** $\Bbbk_q^m$ of $\mathrm{U}_q(\mathfrak{gl}_m)$ is the $m$-dimensional $\Bbbk(q)$-vector space with basis $b_1, ..., b_m$, where

$$L_i \centerdot b_i = qb_i,\ L_i^{-1} \centerdot b_i = q^{-1}b_i,\ L_i^{\pm 1} \centerdot b_j = b_j,\ j \neq i,$$
$$E_{i-1} \centerdot b_i = b_{i-1},\ E_i \centerdot b_j = 0, \text{ if } i \neq j-1,$$
$$F_i \centerdot b_i = b_{i+1},\ F_i \centerdot b_j = 0, \text{ if } i \neq j.$$

Note that this is the quantum version of the standard (vector) representation of $\mathrm{GL}_m$. ◇

We will later need ***Lusztig's idempotented form*** of $\mathrm{U}_q(\mathfrak{gl}_m)$. The idea is to add to $\mathrm{U}_q(\mathfrak{gl}_m)$ an orthogonal system of idempotents in order to view this algebra as a category. Recall that (integral) $\mathfrak{gl}_m$-weights are tuples $\vec{k} = (k_1, ..., k_m)$ in $\mathbb{Z}^m$. We define the ***divided powers***

$$F_i^{(j)} = \tfrac{F_i^j}{[j]!} \text{ and } E_i^{(j)} = \tfrac{E_i^j}{[j]!}.$$

**Definition 11G.4.** The ***category*** $\dot{\mathrm{U}}_q(\mathfrak{gl}_m)$ is the $\mathbb{Z}[q,q^{-1}]$-linear category with objects the $\mathfrak{gl}_m$-weights $\vec{k}=(k_1,...,k_m)\in\mathbb{Z}^m$. The identity morphism of $\vec{k}$ is denoted by $\mathbf{1}_{\vec{k}}$ and the morphism spaces are denoted by $\mathbf{1}_{\vec{l}}\dot{\mathrm{U}}_q(\mathfrak{gl}_m)\mathbf{1}_{\vec{k}}$. The morphisms are generated by $E_i^{(r)}\mathbf{1}_{\vec{k}}\in\mathbf{1}_{\vec{k}+r\alpha_i}\dot{\mathrm{U}}_q(\mathfrak{gl}_m)\mathbf{1}_{\vec{k}}$ and $F_i^{(r)}\mathbf{1}_{\vec{k}}\in\mathbf{1}_{\vec{k}-r\alpha_i}\dot{\mathrm{U}}_q(\mathfrak{gl}_m)\mathbf{1}_{\vec{k}}$, $r\in\mathbb{N}$ (here $\alpha_i=(0,...,1,-1,...,0)\in\mathbb{Z}^m$ as before). These generating morphisms satisfy the following relations

$$F_{i_1}^{(j_1)}F_{i_2}^{(j_2)}\mathbf{1}_{\vec{k}}=F_{i_2}^{(j_2)}F_{i_1}^{(j_1)}\mathbf{1}_{\vec{k}},\ \text{if } |i_1-i_2|>1,$$

$$F_{i_1}^{(j_1)}E_{i_2}^{(j_2)}\mathbf{1}_{\vec{k}}=E_{i_2}^{(j_2)}F_{i_1}^{(j_1)}\mathbf{1}_{\vec{k}},\ \text{if } |i_1-i_2|>0,$$

(11G-5) $$F_{i_1}F_{i_2}F_{i_1}\mathbf{1}_{\vec{k}}=(F_{i_1}^{(2)}F_{i_2}+F_{i_2}F_{i_1}^{(2)})\mathbf{1}_{\vec{k}},\ \text{if } |i_1-i_2|=1,$$

$$F_i^{(j_1)}F_i^{(j_2)}\mathbf{1}_{\vec{k}}=\begin{bmatrix}j_1+j_2\\ j_1\end{bmatrix}F_i^{(j_1+j_2)}\mathbf{1}_{\vec{k}},$$

$$E_i^{(j_2)}F_i^{(j_1)}\mathbf{1}_{\vec{k}}=\sum_{j'}\begin{bmatrix}k_i-j_1-k_{i+1}+j_2\\ j'\end{bmatrix}F_i^{(j_1-j')}E_i^{(j_2-j')}\mathbf{1}_{\vec{k}},$$

$$E_i^{(j)}\mathbf{1}_{\vec{k}}=\mathbf{1}_{\vec{k}+j\alpha_i}E_i^{(j)},\ F_i^{(j)}\mathbf{1}_{\vec{k}}=\mathbf{1}_{\vec{k}-j\alpha_i}F_i^{(j)},$$

and their analogues for $E_i^{(j)}$ (when the weight is clear from the context we write $E_i^{(j)}$ instead of $E_i^{(j)}\mathbf{1}_{\vec{k}}$ in the first, third and fourth relation). Relation Equation 11G-5 is called the ***Serre relation***. ◇

Note that $\mathrm{Hom}(\vec{k},\vec{l})=0$ unless $\sum k_i=\sum l_i$ by the final relation in Definition 11G.4.

11G.3. *Quantum symmetric Howe functors.* Let $\mathrm{Sym}_q^k\mathbb{k}_q^n$ (or denote just $\mathrm{Sym}_q^k$ for short) be the simple $\mathrm{U}_q(\mathfrak{gl}_m)$-module of highest weight $(k,0,...,0)$. It has a basis given by

$$y_{j_1}\otimes...\otimes y_{j_k},\ \text{for } 1\leq j_1\leq...\leq j_k\leq n,$$

where $\{y_1,...,y_n\}$ is a basis of $\mathbb{k}_q^n$ (see [**BZ08**]) and its dimension is $\binom{n+k-1}{k}$. The following is a consequence of quantum symmetric Howe duality (Theorem 2.6 in [**RT16**]).

**Corollary 11G.6.** *( [**RT16**, Theorem 2.6].) There exists a functor*

$$\Phi_{m,n}:\dot{\mathrm{U}}_q(\mathfrak{gl}_m)\to\mathrm{U}_q(\mathfrak{gl}_n)\text{-}\mathbf{fdMod}_S,$$

*which sends a* $\mathfrak{gl}_m$*-weight* $\vec{k}=(k_1,...,k_m)\in\mathbb{Z}^m_{\geq 0}$ *to the* $U_q(\mathfrak{gl}_n)$*-module* $Sym_q^{k_1}\otimes\cdots\otimes Sym_q^{k_m}$ *and morphisms* $X\in\mathbf{1}_{\vec{l}}\dot{\mathrm{U}}_q(\mathfrak{gl}_m)\mathbf{1}_{\vec{k}}$ *to* $f_{\vec{k}}^{\vec{l}}(X)$*. The* $\mathfrak{gl}_m$*-weights with negative entries are sent to the zero module. This functor is moreover surjective on Hom-spaces as in equation Equation 11G-7:*

(11G-7) $$f_{\vec{k}}^{\vec{l}}:\mathbf{1}_{\vec{l}}\dot{\mathrm{U}}_q(\mathfrak{gl}_m)\mathbf{1}_{\vec{k}}\longrightarrow Hom_{U_q(\mathfrak{gl}_n)}(Sym_q^{k_1}\otimes\cdots\otimes Sym_q^{k_m},Sym_q^{l_1}\otimes\cdots\otimes Sym_q^{l_m})$$

*for any two* $\vec{k},\vec{l}\in\mathbb{Z}^m_{\geq 0}$ *such that* $\sum_{i=0}^m k_i=\sum_{i=0}^m l_i$*, meaning that it is full.* □

By dualising [**RT16**, Theorem 2.6] we obtain:

**Corollary 11G.8.** *There exists a functor*

$$\Phi^*_{m.n}:\dot{\mathrm{U}}_q(\mathfrak{gl}_m)\to\mathrm{U}_q(\mathfrak{gl}_n)\text{-}\mathbf{fdMod}_{S^*},$$

*which sends a* $\mathfrak{gl}_m$*-weight* $\vec{k}=(-k_1,...,-k_m)\in\mathbb{Z}^m_{\leq 0}$ *to the* $U_q(\mathfrak{gl}_n)$*-module* $(Sym_q^{k_m})^*\otimes\cdots\otimes(Sym_q^{k_1})^*$ *and morphisms* $X\in\mathbf{1}_{\vec{l}}\dot{\mathrm{U}}_q(\mathfrak{gl}_m)\mathbf{1}_{\vec{k}}$ *to* $(f_{\vec{k}}^{\vec{l}})^*(X)$*. This functor is surjective on Hom-spaces analogous to those in Equation 11G-7.* □

We refer to $\Phi_{m,n}$ and $\Phi^*_{m,n}$ as the ***q-symmetric Howe functors***.

11G.4. *The functor $\Gamma_{sym}$.* We start by recalling a ladder-type functor which has appeared several times in various forms in the literature. The main idea goes back to [CKM14].

**Lemma 11G.9.** *( [RT16, Lemma 2.17].) For each $m \in \mathbb{Z}_{\geq 0}$, there exists a functor*

$$A_{m,n} \colon \dot{\mathrm{U}}_q(\mathfrak{gl}_m) \to \mathbf{sWeb}_{\uparrow}(\mathfrak{gl}_n)$$

*which sends a $\mathfrak{gl}_m$-weight $\vec{k} \in \mathbb{Z}^m_{\geq 0}$ to the sequence obtained by removing all 0's, and sends all other objects $\vec{k}$ of $\dot{\mathrm{U}}_q(\mathfrak{gl}_m)$ to the zero object. This functor is determined on morphisms by the assignment*

$A_{m,n}(F_i^{(j)}\mathbf{1}_{\vec{k}}) =$ 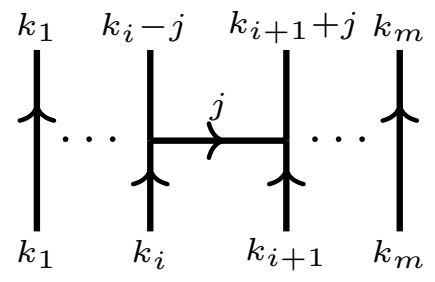

, $A_{m,n}(E_i^{(j)}\mathbf{1}_{\vec{k}}) =$ 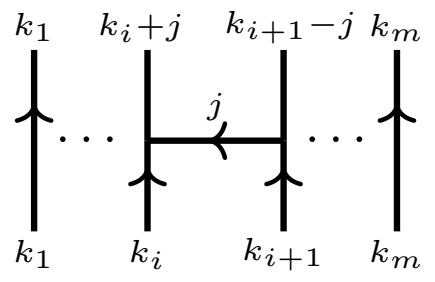

*where we erase any zero labeled edges in the diagrams depicting the images.* □

Let us give some $\mathrm{U}_q(\mathfrak{gl}_n)$-intertwiners, using a thin version of the notations in [RW20, Appendix B].

**Definition 11G.10.** The ***thin projection, inclusion, evaluation and coevaluation maps*** are defined as

$$m_{1,1} : \mathrm{Sym}^1_q \otimes \mathrm{Sym}^1_q \to \mathrm{Sym}^2_q,$$
$$b_i \otimes b_j \mapsto \begin{cases} b_i \otimes b_j, & \text{if} \quad i \leq j, \\ q b_j \otimes b_i, & \text{if} \quad i > j, \end{cases}$$
$$s_{1,1} : \mathrm{Sym}^2_q \to \mathrm{Sym}^1_q \otimes \mathrm{Sym}^1_q,$$
$$b_i \otimes b_j \mapsto \begin{cases} q^{-1} b_i \otimes b_j + b_j \otimes b_i, & \text{if} \quad i < j, \\ [2] b_i \otimes b_i, & \text{if} \quad i = j, \end{cases}$$
$$ev^{\text{left}}_1 : (\mathrm{Sym}^1_q)^* \otimes \mathrm{Sym}^1_q \to \mathbb{k}_q,$$
$$b_i^* \otimes b_j \mapsto b_i^*(b_j),$$
$$coev^{\text{left}}_1 : \mathbb{k}_q \to \mathrm{Sym}^1_q \otimes (\mathrm{Sym}^1_q)^*,$$
$$1 \mapsto \sum_{i=1}^{n} b_i \otimes b_i^*.$$

◇

**Lemma 11G.11.** *The maps $m_{1,1}$, $s_{1,1}$, $ev^{left}_1$, $coev^{left}_1$ are $U_q(\mathfrak{gl}_n)$-intertwiners.*

*Proof.* A calculation. □

Recall now the $q$-symmetric Howe functors from Section 11G.3. We use them to define our main functor.

**Definition 11G.12.** The functor

$$\Gamma_{\text{sym}} : \mathbf{sWeb}^{\uparrow,\downarrow}(\mathfrak{gl}_n) \to \mathrm{U}_q(\mathfrak{gl}_n)\text{-}\mathbf{fdMod}_{S,S^*}$$

is defined as follows.

- On objects $k^{\uparrow}$ and $k_{\downarrow}$, $k \in \mathbb{Z}_{\geq 0}$ we set
$$\Gamma_{\text{sym}}(k^{\uparrow}) = \mathrm{Sym}^k_q \quad \text{and} \quad \Gamma_{\text{sym}}(k_{\downarrow}) = (\mathrm{Sym}^k_q)^*.$$
By convention, we send the zero object to the zero module and the empty tuple to the trivial $\mathrm{U}_q(\mathfrak{gl}_n)$-module.
- On morphisms, we send the generating morphisms of $\mathbf{sWeb}^{\uparrow,\downarrow}(\mathfrak{gl}_n)$ to the following $\mathrm{U}_q(\mathfrak{gl}_n)$-intertwiners and we extend monoidally

$\Gamma_{\text{sym}}\Big($ 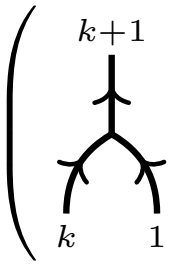

$\Big) = \Phi_{2,n}(E_1^{(1)}\mathbf{1}_{(k,1)})$, $\Gamma_{\text{sym}}\Big($ 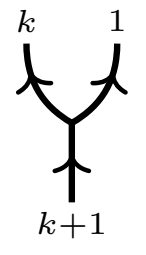

$\Big) = \Phi_{2,n}(F_1^{(1)}\mathbf{1}_{(k+1,0)})$,

$\Gamma_{\text{sym}}\Big($ 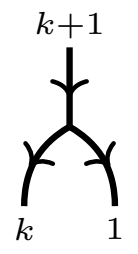

$\Big) = \Phi^*_{2,n}(F_1^{(1)}\mathbf{1}_{(-1,-k)})$, $\Gamma_{\text{sym}}\Big($ 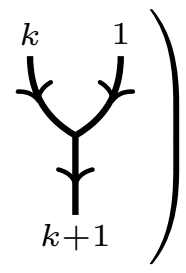

$\Big) = \Phi^*_{2,n}(E_1^{(1)}\mathbf{1}_{(0,-1-k)})$,

$\Gamma_{\text{sym}}\Big($ 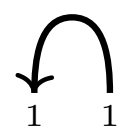

$\Big) = ev^{\text{left}}_1$, $\Gamma_{\text{sym}}\Big($ 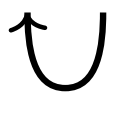

$\Big) = coev^{\text{left}}_1$.

◇

We now show that the functor $\Gamma_{\mathrm{sym}}$ is actually well-defined. Here it will become clear why we constructed the diagrammatic category the way we did.

**Lemma 11G.13.** *The functor*

$$\Gamma_{sym}\colon \mathbf{sWeb}^{\uparrow,\downarrow}(\mathfrak{gl}_n)\to \mathrm{U}_q(\mathfrak{gl}_n)\text{-}\mathbf{fdMod}_{S,S^*}$$

*is well-defined.*

*Proof.* Let us show that the analogues of the relations from Definition 11B.24 hold in $\mathrm{U}_q(\mathfrak{gl}_n)$-**fdMod**. The relations Equation 11B-9 and Equation 11B-10 hold in $\mathrm{U}_q(\mathfrak{gl}_n)$-**fdMod** as a consequence of q-Howe duality, namely of Corollary 11G.6 since they come from relations in $\dot{\mathrm{U}}_q(\mathfrak{gl}_m)$. An easy calculation identifies the morphisms associated to the right-hand sides in Equation 11B-3 as

$$-q^{-1-\frac{1}{n}}\mathrm{id}_1\otimes\mathrm{id}_1+q^{-\frac{1}{n}}s_{1,1}\circ m_{1,1}\text{ and }-q^{-1+\frac{1}{n}}\mathrm{id}_1\otimes\mathrm{id}_1+q^{\frac{1}{n}}s_{1,1}\circ m_{1,1},$$

respectively. These morphisms are invertible, and composing them with evaluation and coevaluation maps does not change that property, meaning that Equation 11B-25 also holds. Now, a simple computation shows that Equation 11B-26 holds in $\mathrm{U}_q(\mathfrak{gl}_n)$-**fdMod**. Indeed, from Definition 11G.10, we have

$$(ev_1^{\mathrm{left}}\circ\mathrm{id}_1)\circ(\mathrm{id}_1\circ coev_1^{\mathrm{left}})(b_j^*\otimes 1)=(ev_1^{\mathrm{left}}\circ\mathrm{id}_1)(b_j^*\otimes\sum_{i=1}^{n}b_i\otimes b_i^*)=\sum_{i=1}^{n}b_j^*(b_i)\otimes b_i^*=1\otimes b_j^*,$$

meaning that $(ev_1^{\mathrm{left}}\circ\mathrm{id}_1)\circ(\mathrm{id}_1\circ coev_1^{\mathrm{left}})=\mathrm{id}_1$, which is exactly the representation theoretical version of the zigzag relation from Equation 11B-26. Showing relation Equation 11B-27 is again a brute force computation which we will omit here, see [LT21] for details. Finally, relation Equation 11B-28 holds because of our definition of $ev_1^{\mathrm{left}}$:

$$ev_1^{\mathrm{left}}\circ coev_1^{\mathrm{right}}(1)=ev_1^{\mathrm{left}}(\sum_{i=1}^{n}q^{-n+2i-1}b_i^*\otimes b_i)=q^{-n+1}+q^{-n+3}+\ldots+q^{n-1}=[n].$$

The proof completes. □

11G.5. *Relating the two web categories.* Now we want to identify $\mathbf{sWeb}^{\uparrow}(\mathfrak{gl}_n)$ as a full subcategory of the category $\mathbf{sWeb}^{\uparrow,\downarrow}(\mathfrak{gl}_n)$. Recall Proposition 11C.17. In particular, from this proposition we can assume that the boundaries of a web in $\mathbf{sWeb}_{\uparrow,\downarrow}(\mathfrak{gl}_n)$ are upward pointing. Below we show that we have an even stronger connection between our two symmetric web categories.

**Theorem 11G.14.** *There exists a fully faithful functor of braided monoidal categories*

$$\mathrm{G}\colon\mathbf{sWeb}^{\uparrow}(\mathfrak{gl}_n)\to\mathbf{sWeb}_{\uparrow,\downarrow}(\mathfrak{gl}_n)$$

*sending objects and morphism from* $\mathbf{sWeb}^{\uparrow}(\mathfrak{gl}_n)$ *to objects and morphisms of the same name in* $\mathbf{sWeb}_{\uparrow,\downarrow}(\mathfrak{gl}_n)$.

*Proof.* We give a line-to-line adaptation of [BDK20, Theorem 7.5.1]. First, it is clear that this functor is well-defined because the defining relations of $\mathbf{sWeb}^{\uparrow}(\mathfrak{gl}_n)$ hold in $\mathbf{sWeb}^{\uparrow,\downarrow}(\mathfrak{gl}_n)$. Also, they have the same monoidal braided structures by construction. Remains to prove fully faithfulness. Let $k^{\uparrow}$, $l^{\uparrow}\in\mathrm{Obj}(\mathbf{sWeb}^{\uparrow}(\mathfrak{gl}_n))$ and $w$ a web in $\mathrm{Hom}_{\mathbf{sWeb}^{\uparrow,\downarrow}(\mathfrak{gl}_n)}(k^{\uparrow},l^{\uparrow})$. We want to show that

$$\mathrm{Hom}_{\mathbf{sWeb}^{\uparrow}(\mathfrak{gl}_n)}(k_{\uparrow},l_{\uparrow})\cong\mathrm{Hom}_{\mathbf{sWeb}^{\uparrow,\downarrow}(\mathfrak{gl}_n)}(\mathrm{G}(k_{\uparrow}),\mathrm{G}(l_{\uparrow})).$$

Let us first show fullness. Note that the bottom and top of a web in $\mathbf{sWeb}_{\uparrow,\downarrow}(\mathfrak{gl}_n)$ are upward pointing but in the middle of it we can have both upward and downward pointing merges, splits, cups, caps and crossings. We use relations Equation 11B-26 and Equation 11B-28 to straighten zigzags and to remove circles. Further, we use Lemma 11C.3 to remove any downward pointing merges and splits. From Lemma 11C.15 and Lemma 11C.12 we can slide crossing, cups and caps through merges and splits putting the former below the later. Proceeding inductively on the total number of crossings, we can therefore write $w$ as a linear combination of diagrams of the form $k^{\uparrow}\to l^{\uparrow}$, where each diagram has an upper part and a lower part. The upper part consists only of upward pointing merges and splits, whereas the lower part consists only of cups, caps and crossings (both upward and downward pointing). The webs of the upper part are clearly morphisms in $\mathbf{sWeb}^{\uparrow}(\mathfrak{gl}_n)$. In the lower part, after removing any circles, the remaining webs are (compositions of) crossings or compositions as in the Reidemeister 1 moves. Note that in the bottom part we have a geometric braid, which can be turned upwards. Now we use the three types of Reidemeister moves, in particular, we use Lemma 11C.20 to write the compositions of crossings with cups and caps as scalar multiples of identities. In a nutshell, every web of form $k^{\uparrow} \to l^{\uparrow}$ in $\mathbf{sWeb}^{\uparrow,\downarrow}(\mathfrak{gl}_n)$ can be written as a linear combination of diagrams of the form $k^{\uparrow} \to l^{\uparrow}$ in $\mathbf{sWeb}^{\uparrow}(\mathfrak{gl}_n)$. This shows fullness of G. Finally, consider the diagram

$$\begin{array}{ccc} \mathbf{sWeb}_{\uparrow}(\mathfrak{gl}_n) & \xrightarrow{\Gamma_{\text{sym}}} & \mathrm{U}_q(\mathfrak{gl}_n)\text{-}\mathbf{fdMod} \\ & {\scriptstyle G}\searrow \quad \nearrow{\scriptstyle \Gamma_{\text{sym}}} & \\ & \mathbf{sWeb}_{\uparrow,\downarrow}(\mathfrak{gl}_n) & \end{array}.$$

It commutes by definition. Moreover, the top functor is faithful from Lemma 11G.17. This then implies that $G$ is faithful. □

11G.6. *Proof that $\Gamma_{sym}$ is an equivalence of categories.* We extend our functor by passing to the additive Karoubi envelopes of our categories. We abuse notation and continue denoting it by $\Gamma_{\text{sym}}$

$$\Gamma_{\text{sym}} : \mathbf{Kar}(\mathbf{sWeb}^{\uparrow,\downarrow}(\mathfrak{gl}_n)) \to \mathrm{U}_q(\mathfrak{gl}_n)\text{-}\mathbf{fdMod}.$$

We are finally ready to prove the following.

**Theorem 11G.15. (Equivalence Theorem for symmetric webs)** *The functor*

$$\Gamma_{\text{sym}} : \mathbf{sWeb}^{\uparrow,\downarrow}(\mathfrak{gl}_n) \to \mathrm{U}_q(\mathfrak{gl}_n)\text{-}\mathbf{fdMod}_{S,S^*}$$

*is fully faithful and*

$$\Gamma_{\text{sym}} : \mathbf{Kar}(\mathbf{sWeb}^{\uparrow,\downarrow}(\mathfrak{gl}_n)) \to \mathrm{U}_q(\mathfrak{gl}_n)\text{-}\mathbf{fdMod}$$

*is an equivalence of ribbon categories.*

*Proof.* We have already proved that $\Gamma_{\text{sym}}$ is well-defined in Lemma 11G.13. Now we show that it is fully faithful and essentially surjective after passing to Karoubi envelopes.
Fully faithful. We have to show that for $\vec{k}, \vec{l} \in \mathbb{Z}^m$, $m > 0$, we have

$$\mathrm{Hom}_{\mathbf{sWeb}^{\uparrow,\downarrow}(\mathfrak{gl}_n)}(\vec{k}, \vec{l}) \cong \mathrm{Hom}_{\mathrm{U}_q(\mathfrak{gl}_n)\text{-}\mathbf{fdMod}}(\Gamma_{\text{sym}}(\vec{k}), \Gamma_{\text{sym}}(\vec{l})).$$

From Proposition 11C.17 and Theorem 11G.14, it suffices to show this for upward pointing webs only and we have proved the later in Lemma 11G.17.
Essentially surjective. This follows from the definition of $\Gamma_{\text{sym}}$ and the fact that every irreducible $\mathrm{U}_q(\mathfrak{gl}_n)$-module appears in a suitable tensor product of symmetric powers of the standard representation and their duals. Indeed, let $V$ be an irreducible $\mathrm{U}_q(\mathfrak{gl}_n)$-module appearing as a direct summand in $\bigotimes_i \mathrm{Sym}_q^{k_i} \otimes \bigotimes_j (\mathrm{Sym}_q^{k_j})^*$. Then $\Gamma_{\text{sym}}(\bigotimes_i k_i{}^{\uparrow} \otimes \bigotimes_j k_{j\downarrow}) = \bigotimes_i \mathrm{Sym}_q^{k_i} \otimes \bigotimes_j (\mathrm{Sym}_q^{k_j})^*$. Since both sides are now idempotent complete, tensor products decompose into irreducible objects and such decompositions are unique up to permutation. Thus, there exists a direct summand of $\bigotimes_i k_i{}^{\uparrow} \otimes \bigotimes_j k_{j\downarrow}$ that is sent to $V$ by $\Gamma_{\text{sym}}$.
This finishes the proof that $\Gamma_{\text{sym}}$ is an equivalence of categories. Finally, after adapting Section 3.1 of [RT16] and Section 6 of [CKM14] to our setting we have that $\Gamma_{\text{sym}}$ is an equivalence of pivotal braided categories. □

Now we give the connection between the three functors of this section. The following is a generalisation of Corollary 2.22 in [RT16].

**Lemma 11G.16.** *The diagram below commutes*

$$\begin{array}{ccc} \dot{\mathrm{U}}_q(\mathfrak{gl}_m) & \xrightarrow{\Phi_{m,n}} & \mathrm{U}_q(\mathfrak{gl}_n)\text{-}\mathbf{fdMod} \\ & {\scriptstyle A_{m,n}}\searrow \quad \nearrow{\scriptstyle \Gamma_{sym}} & \\ & \mathbf{sWeb}_{\uparrow}(\mathfrak{gl}_n) & \end{array}.$$

The following lemma is one of the main ingredients for proving Theorem 11G.15. Note that it concerns upward pointing webs only.

**Lemma 11G.17.** *The functor*

$$\Gamma_{sym} : \mathbf{sWeb}^{\uparrow}(\mathfrak{gl}_n) \to \mathrm{U}_q(\mathfrak{gl}_n)\text{-}\mathbf{fdMod}_S$$

*is fully faithful.*

*Proof.* We have to show that for tuples $k_{\uparrow} = (k_1, ..., k_m)_{\uparrow} \in \mathbb{Z}^m_{\geq 0}$ and $l_{\uparrow} = (l_1, ..., l_{m'})_{\uparrow} \in \mathbb{Z}^{m'}_{\geq 0}$ we have

$$\mathrm{Hom}_{\mathbf{sWeb}^{\uparrow}(\mathfrak{gl}_n)}(k_{\uparrow}, l_{\uparrow}) \cong \mathrm{Hom}_{\mathrm{U}_q(\mathfrak{gl}_n)\text{-}\mathbf{fdMod}}(\Gamma_{\text{sym}}(k_{\uparrow}), \Gamma_{\text{sym}}(l_{\uparrow})).$$

Note that we can assume that $m = m'$, otherwise we add zeros to the shorter tuple. Now, $\Gamma_{\text{sym}}$ by definition comes from a $q$-Howe functor $\Phi_{m,n}$, for $m \in \mathbb{Z}_{\geq 0}$ big enough.

***Surjectivity:*** In the case $\sum_{i=1}^{m} k_i = \sum_{i=1}^{m} l_i$ surjectivity follows from Lemma 11G.16 and the fact that $q$-symmetric Howe functor $\Phi_{m,n}$ is surjective on these Hom-spaces, whereas for $\sum_{i=1}^{m} k_i \neq \sum_{i=1}^{m} l_i$ we have zero on both sides.

***Injectivity:*** By surjectivity and finite-dimensionality of the involved spaces it suffices to show that

$$\dim\mathrm{Hom}_{\mathbf{sWeb}^{\uparrow}(\mathfrak{gl}_n)}(k_{\uparrow}, l_{\uparrow}) = \dim\mathrm{Hom}_{\mathrm{U}_q(\mathfrak{gl}_n)\text{-}\mathbf{fdMod}}(\Gamma_{\mathrm{sym}}(k_{\uparrow}), \Gamma_{\mathrm{sym}}(l_{\uparrow})). \tag{11G-18}$$

We prove this first for $k_{\uparrow} = (1, ..., 1)_{\uparrow} = l_{\uparrow}$. In this case, recall that in $\mathrm{U}_q(\mathfrak{gl}_n)$-**fdMod** we have

$$\mathrm{Sym}_q^1 \Bbbk_q^n \cong \textstyle\bigwedge_q^1 \Bbbk_q^n \cong \Bbbk_q^n.$$

Further, we have

$$\dim\mathrm{End}_{\mathrm{U}_q(\mathfrak{gl}_n)\text{-}\mathbf{fdMod}}(\Bbbk_q^n \otimes ... \otimes \Bbbk_q^n) = \dim\mathrm{End}_{\mathrm{U}_q(\mathfrak{sl}_n)\text{-}\mathbf{fdMod}}(\Bbbk_q^n \otimes ... \otimes \Bbbk_q^n), \tag{11G-19}$$

since every irreducible $\mathrm{U}_q(\mathfrak{sl}_n)$-module comes from a restriction of an irreducible $\mathrm{U}_q(\mathfrak{gl}_n)$-module.

Now, from Proposition 11F.6 we have that

$$\dim\mathrm{End}_{\mathbf{sWeb}_{\uparrow}(\mathfrak{gl}_n)}(1 \otimes ... \otimes 1) = \dim\mathrm{End}_{\mathbf{eWeb}^{\uparrow,\downarrow}(\mathfrak{sl}n)}(1 \otimes ... \otimes 1). \tag{11G-20}$$

Finally, there exists an equivalence of categories (see [CKM14, Theorem 3.3.1])

$$\mathbf{eWeb}^{\uparrow,\downarrow}(\mathfrak{sl}_n) \to \mathrm{U}_q(\mathfrak{sl}_n)\text{-}\mathbf{fdMod}.$$

Combining this equivalence with Equation 11G-19 and Equation 11G-20 shows that the equality in Equation 11G-18 holds for $k_{\uparrow} = (1, ..., 1)_{\uparrow} = l_{\uparrow}$. For the general case, we proceed as in the proof of Theorem 1.10 in [RT16]. Let $k_{\uparrow}, l_{\uparrow}$ be such that $\sum_{i=1}^{m} k_i = \sum_{i=1}^{m} l_i$ and let $u, v \in \mathrm{Hom}_{\mathbf{sWeb}_{\uparrow}}(k_{\uparrow}, l_{\uparrow})$, $u \neq v$. Composing $u$ and $v$ with merges and splits we get

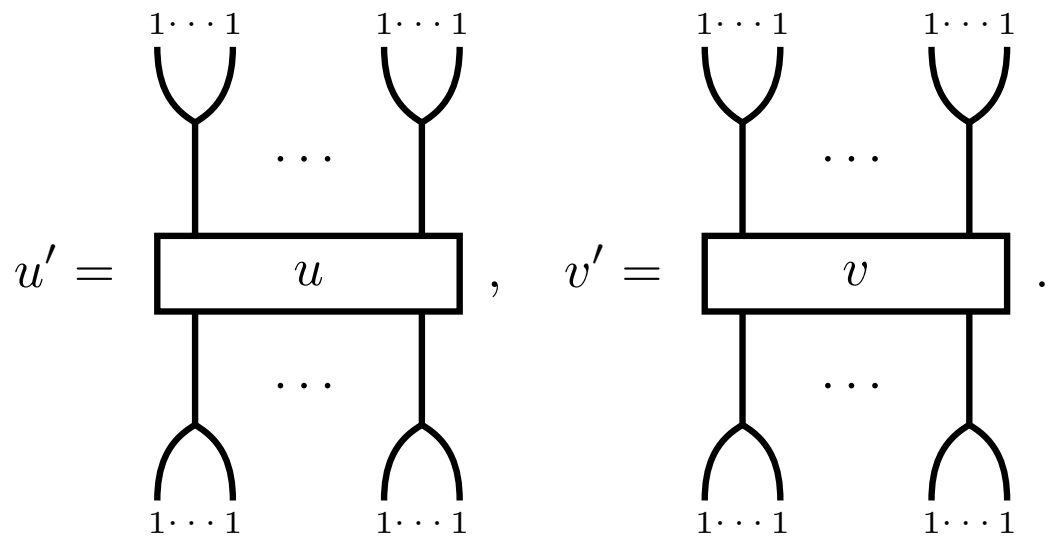

Now, recall that the digon removal is invertible since quantum numbers are invertible (Remark 11B.12) implying that $u' \neq v'$. Finally, from the above argument we have that $\Gamma_{\mathrm{sym}}(u') \neq \Gamma_{\mathrm{sym}}(v')$, showing faithfulness for the general case. □

11H. **Colored Jones polynomials.** We now follow [RT16], and let $n = 2$, and let also $(\Bbbk = \mathbb{C}(q), q)$ be the generic pair we use. We explore how the braiding on $\mathbf{sWeb}^{\uparrow,\downarrow}(\mathfrak{gl}_2)$ can be used to study the colored Jones polynomial of colored, oriented links $L$, which we denote by $J_{\vec{c}}(L_D)$, that we have seen in Section 9H. Here $\vec{c} = (c_1, ..., c_N)$ denotes the colors of the $N$-component, oriented link $L$ and $L_D$ is a colored, oriented diagram for $L$. As before, it can be computed by associating a morphism between trivial representations in $\mathbf{sWeb}^{\uparrow,\downarrow}(\mathfrak{gl}_2)$ to any colored, oriented link diagram $L_D$ of a colored, oriented link $L$ (and rescaling to get an invariant which is not framing-dependent).

This translates to using Definition 11B.4 to view a colored, oriented link diagram for $L_D$ as a morphism in $\mathbf{sWeb}^{\uparrow,\downarrow}(\mathfrak{gl}_2)$, which necessarily evaluates to an element in $\mathbb{C}(q^{\frac{1}{2}})$ (in fact, it is clear from our construction that this always gives an element in $\mathbb{Z}[q^{\frac{1}{2}}, q^{-\frac{1}{2}}]$), and multiplying by a certain normalization factor which can be computed directly from the crossing data of the diagram. For the case of a $c$-colored knot $K$ with diagram $K_D$, this factor is $(-1)^c q^{-C\omega(K_D)}$, where $\omega(K_D)$ is the writhe (difference between the number of over and undercrossings) of $K_D$ and $C = \frac{c^2+2c}{2}$ is the so-called quadratic Casimir number. In general, for a colored, oriented link diagram $L_D$ one normalizes by multiplying by the product of the normalization factors for each of the colored link components.

**Example 11H.1.** As an example, we compute the (1-colored) Jones polynomial of the Hopf link using symmetric webs.

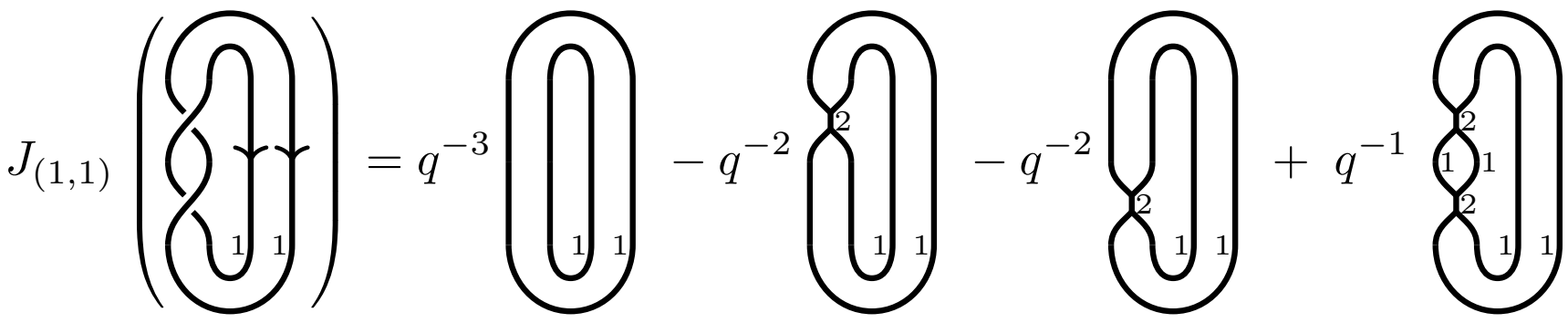

$$= q^{-3}[2]^2 - 2q^{-2}[2][3] + q^{-1}[2]^2[3] = [2](q^2+q^{-2}) = [4].$$

This is what we got in Section 9H as well. ◇

**Example 11H.2.** As another example, we compute the $(2,1)$-colored Jones polynomial of the Hopf link using our approach.

$J_{(2,1)}$ 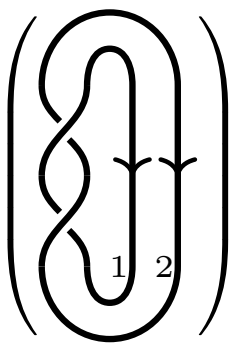

$= q^{-4}$ 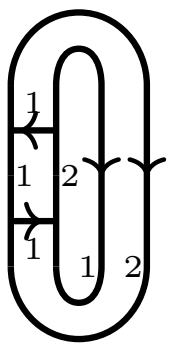

$- q^{-3}$ 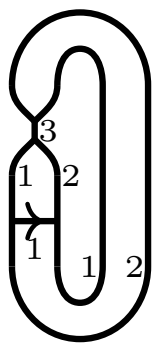

$- q^{-3}$ 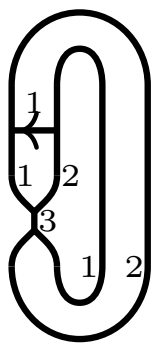

$+ q^{-2}$ 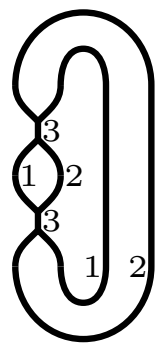

$$= -q^{-4}[2][3]^2 + 2q^{-3}[2][3][4] - q^{-2}[3]^2[4] = -[6].$$

Again, this should be compared with Section 9H. ◇

In both of these examples, we recover the known formula for the colored Jones polynomial of the Hopf link

$J_{(k,l)}$ 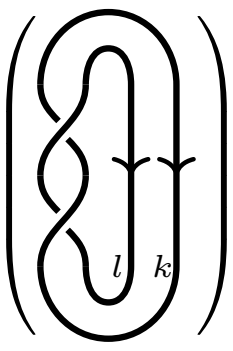

$= (-1)^{k+l}[(k+1)(l+1)]$

for $k,l \in \mathbb{Z}_{\geq 0}$ (the $(-1)^{k+l}$ factor comes from our conventions, and varies in the literature).

## 11I. Exercises.

*Exercise* 11I.1. Find a more general statement of Proposition 11C.17 (beyond web categories). ◇

*Exercise* 11I.2. Redo the first few sections for exterior webs. ◇

*Exercise* 11I.3. Compute the $q \to 1$ limit of $[a] = \frac{q^a - q^{-a}}{q - q^{-1}}$. What are the $q \to 1$ limits of quantum factorials and binomials? ◇

*Exercise* 11I.4. Prove Lemma 11G.11. ◇

*Exercise* 11I.5. Compute, using webs, that $J_{(k,l)} = (-1)^{k+l}[(k+1)(l+1)]$ as claimed above. ◇

*Exercise* 11I.6. Compute, using webs, the value of

$J_{(k)}$ 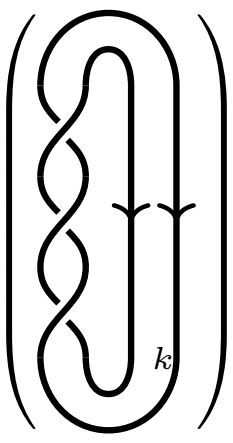

,

for $k \in \mathbb{Z}_{\geq 0}$. This is the $k$-colored trefoil. ◇

# 12. Growth in monoidal categories

The principle of counting is deeply ingrained in human history, evident in practices such as the use of tally marks dating back over 35,000 years. Across civilizations, this desire to count has led to the development of powerful mathematics. In the 21st century, our understanding of counting extends beyond simple enumeration, and is fundamental in all of mathematics and the sciences.

> How to count in monoidal categories?

For sequences $b_n$ arising in many areas of mathematics, the following type of asymptotic behavior (for $n \gg 0$) is ubiquitous, see e.g. [FS09]:

(12-1) $$b_n \sim h(n) \cdot n^a \cdot \beta^n,$$

where $\beta^n$ is the dominating growth, $n^a$ is the subexponential (correction) factor and $h(n)$ is a bounded function (often constant). Remarkably, work in [COT24, CEO24, CEOT24, KST24, LTV23, LTV24a] suggests that the same model Equation 12-1 also describes growth rates in monoidal categories, which we summarize in this section.

What we study in monoidal categories is the following. The classical example of Equation 12-1 is the prime counting function $\pi'_n = \pi_{e^n}$. To see that the dominating growth of $\pi'_n$ is $e^n$ is easy, obtaining the more precise $\pi'_n \sim n^{-1} \cdot e^n$ (implying that $h$ is constant here) is the celebrated prime number theorem, *cf.* Figure 30. A common way of describing how good the approximation Equation 12-1 is studying the behavior of $|b_n - a_n|$, we refer to this question as the ***variance***. The variance for $\pi'_n$ is a direct consequence of the famous and still open Riemann Hypothesis, demonstrating how challenging this problem tends to be.

| $x$ | $\pi(x)$ | $\pi(x) - \frac{x}{\log(x)}$ | $\mathrm{li}(x) - \pi(x)$ | % error | |
|---|---|---|---|---|---|
| | | | | $\frac{x}{\log(x)}$ | $\mathrm{li}(x)$ |
| 10 | 4 | 0 | 2 | 8.22% | 42.606% |
| $10^2$ | 25 | 3 | 5 | 14.06% | 18.597% |
| $10^3$ | 168 | 23 | 10 | 14.85% | 5.561% |
| $10^4$ | 1,229 | 143 | 17 | 12.37% | 1.384% |
| $10^5$ | 9,592 | 906 | 38 | 9.91% | 0.393% |
| $10^6$ | 78,498 | 6,116 | 130 | 8.11% | 0.164% |
| $10^7$ | 664,579 | 44,158 | 339 | 6.87% | 0.051% |
| $10^8$ | 5,761,455 | 332,774 | 754 | 5.94% | 0.013% |
| $10^9$ | 50,847,534 | 2,592,592 | 1,701 | 5.23% | $3.34 \times 10^{-3}$ % |
| $10^{10}$ | 455,052,511 | 20,758,029 | 3,104 | 4.66% | $6.82 \times 10^{-4}$ % |
| $10^{11}$ | 4,118,054,813 | 169,923,159 | 11,588 | 4.21% | $2.81 \times 10^{-4}$ % |
| $10^{12}$ | 37,607,912,018 | 1,416,705,193 | 38,263 | 3.83% | $1.02 \times 10^{-4}$ % |
| $10^{13}$ | 346,065,536,839 | 11,992,858,452 | 108,971 | 3.52% | $3.14 \times 10^{-5}$ % |

FIGURE 30. The prime number theorem: the number of primes $\pi(n)$ is approximately $n/\ln n$: the error goes to zero, and the variance is an order of magnitude smaller than the actual values. (There is also an approximation using the logarithmic integral $li(n)$, but we will ignore this.)
Picture from https://en.wikipedia.org/wiki/Prime_number_theorem

## 12A. **A word about conventions.** We start with the following.

*Convention* 12A.1. We will use sequences $(b_n)_{n \in \mathbb{Z}_{\geq 0}}$ with $b_n \in \mathbb{R}$ and their associated functions $b_n \colon \mathbb{Z}_{\geq 0} \to \mathbb{R}, n \mapsto b_n$ interchangeably. ◇

Throughout, we are interested in ***large $n$ behavior***. We fix some notation that we will use.

*Convention* 12A.2. For functions $f, g \colon \mathbb{Z}_{\geq 0} \to \mathbb{R}_{>0}$ we use:

$$f \in \Omega(g) \Leftrightarrow \exists c > 0, \exists n_0 \text{ such that } \boxed{c \cdot g(n) \leq |f(n)|}, \forall n > n_0,$$

$$f \in O(g) \Leftrightarrow \exists C > 0, \exists n_0 \text{ such that } \boxed{|f(n)| \leq C \cdot g(n)}, \forall n > n_0.$$

We use these, in order, as ***asymptotic lower and upper bounds***. Both together give ***asymptotics***:

$$f \sim g \Leftrightarrow \forall \varepsilon > 0, \exists n_0 \text{ such that } \boxed{\left|\tfrac{f(n)}{g(n)} - 1\right| < \varepsilon}, \forall n > n_0,$$

$$f \in \Theta(g) \Leftrightarrow \big(f \in \Omega(g) \text{ and } f \in O(g)\big).$$

These are almost the same, but $\Theta$ ignores scalars. ◇

## 12B. **Motivation.** Many (almost all?) counting problems are too difficult to solve explicitly, even some that seem easy. But many have nice formulas asymptotically. We have already seen the prime number theorem illustrating this beautifully. Another key example is

$$t_n = \text{number of trees with } n \text{ vertices (up to isomorphism)}, \; n \in \mathbb{Z}_{\geq 0},$$

where we remind the reader that a ***tree*** is a connected graph without cycles. The first few trees are displayed in Figure 31, the sequence $t_n$ is [**OEI23**, A000055].

It seems difficult to get a precise count for $t_n$. This seems surprising since the number of labeled (=colored) trees is just $n^{n-2}$. It turns out that we can still determine its asymptotic growth rate:

$$t_n \sim c \cdot n^{-5/2} \cdot d^n$$

where $c \approx 0.535$ and $d \approx 2.956$ are famous constants given in [**OEI23**, A000055].

The idea how to prove such an asymptotic formula is similar to (and in some sense dual to) how to determine the growth rate of the Fibonacci sequence $(f_n)_{n \in \mathbb{Z}_{\geq 0}} = (1, 1, 2, 3, 5, \ldots)$, and we will focus on this sequence instead. The recursion $f_n = f_{n-1} + f_{n-2}$ can be expressed via the matrix

$$M = \begin{pmatrix} 0 & 1 \\ 1 & 1 \end{pmatrix},$$

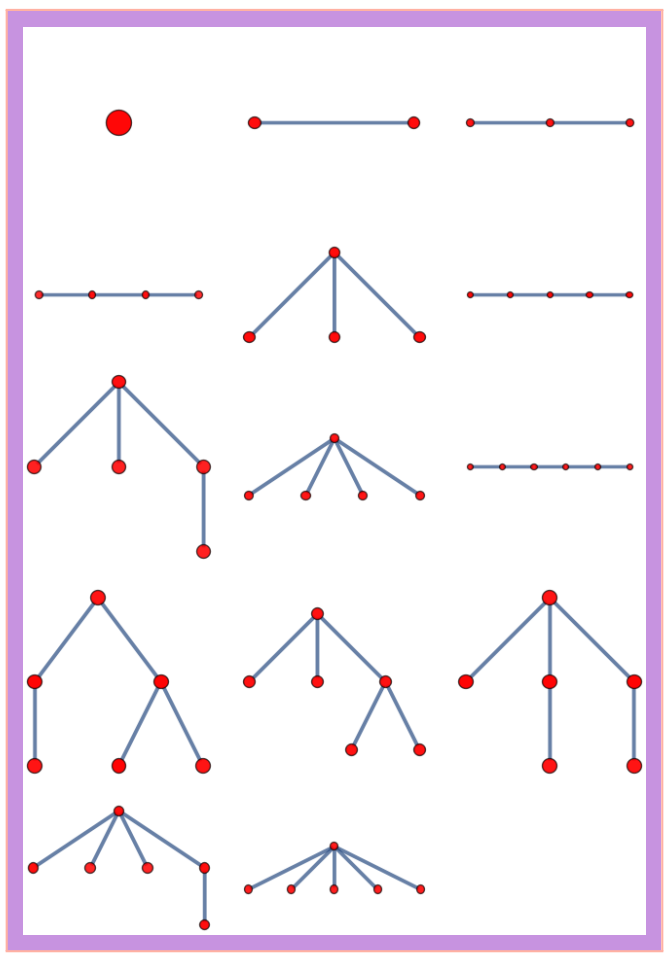

FIGURE 31. The first few trees showing that $(t_n)_{n\in\mathbb{Z}_{\geq 0}} = (1,1,1,2,3,6,...)$.)
Picture created using AI

that we have already seen in Example 8E.4. The PF eigenvalue and vector are

$$\lambda_{pf} = \phi, \quad v_{pf} = (1,\phi).$$

Now we normalize $v_{pf}$ to $v'_{pf} = \frac{1}{\sqrt{1+\phi^2}}$, which has length one. Now we compute (with $T$ meaning transpose)

$$v'_{pf} v^T_{pf} = \begin{pmatrix} \frac{1}{\phi^2+1} & \frac{1}{\sqrt{5}} \\ \frac{1}{\sqrt{5}} & \frac{1}{10}\left(\sqrt{5}+5\right) \end{pmatrix},$$

and the sum of the first column is $\phi/\sqrt{5}$. We get

$$f_n \sim \boxed{\phi/\sqrt{5}} \cdot \boxed{n^0} \cdot \boxed{\phi^n}$$

by the PF theorem, as we will see below. This agrees with the asymptotic one gets from Binet's formula:

$$f_n = \frac{\phi^{n+1} - (-\phi)^{-n-1}}{\sqrt{5}} \sim \frac{\phi^{n+1}}{\sqrt{5}}.$$

12C. **Examples from asymptotic combinatorics.** An involution of $\{1,...,n\}$ is a function $f\colon \{1,...,n\} \to \{1,...,n\}$ such that $f \circ f = id$. In terms of diagrams, this is an endomorphism of $\bullet^n \in \mathbf{Sym}$ as in Example 3D.2 that squares to the identity, e.g. for $n = 3$:

involutions : 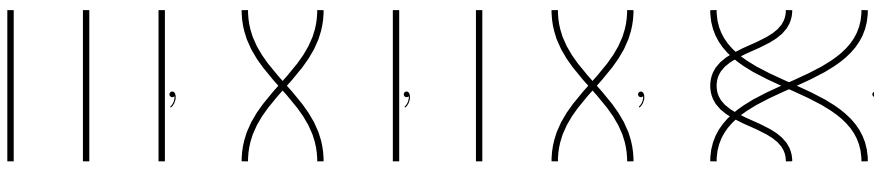

The other two elements in the symmetric group $S_3 = \mathrm{End}_{\mathbf{Sym}}(\bullet^3)$ are not involutions.

Let $s_n$ denote the number of involutions in $S_n$, for example,

$$\{1,1,2,4,10,26,76,232,764,2620,9496\}, \quad s_n \text{ for } n = 0,...,10.$$

A quick search gives that this is [**OEI23**, A000085].

**Lemma 12C.1.** *We have the recursion $s_0 = s_1 = 1$ and $s_n = s_{n-1} + (n-1)s_{n-2}$ for $n \in \mathbb{Z}_{\geq 2}$.*

*Proof.* The following diagrammatic procedure, illustrated for $n = 5$ takes an involution in $S_{n-1}$ and $S_{n-2}$ (seen as a subset of $S_{n-1}$ of diagrams that do not move the rightmost strand), and gives an involution in $S_n$:

adding a strand at the right: 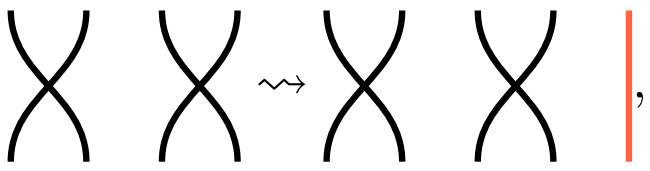 , adding a crossing: 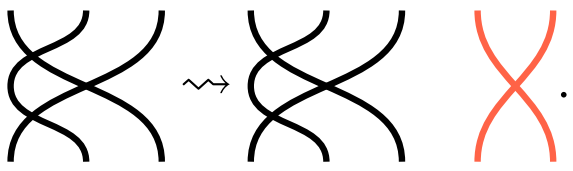 .

On the left-hand side there is no choice, but for the right-hand operation we can decide where to put the left endpoints of the crossings, giving $(n-1)$ options. The recursion thus follows. □

**Lemma 12C.2.** *We have $s_n = \sum_{k=0}^{n/2} n!/((n-2k)!2^k k!)$.*

*Proof.* The formula satisfies the recursion in Lemma 12C.1, so we are done. □

Now, we run the following MAGMA code (as in Appendix A):

```
> // Initialize an empty list to store the total dimensions
> totalDimensions := [];
> // Loop through each symmetric group from S_1 to S_10
> for n in [1..10] do
> // Create the symmetric group S_n
> S := SymmetricGroup(n);
> // Get the character table of S_n
> charTable := CharacterTable(S);
> // Extract the dimensions from the first column of the character table
> dims := [charTable[i, 1] : i in [1..#charTable]];
> // Sum the dimensions
> totalDimension := &+dims;
> // Append the result to the list
> Append(~totalDimensions, totalDimension);
> end for;
> // Output the list of total dimensions
> totalDimensions;
  ----result----
[ 1, 2, 4, 10, 26, 76, 232, 764, 2620, 9496 ]
```

Except that MAGMA does not know $S_0$, this is our sequence $s_n$. We get a surprisingly nice formula (with $e \approx 2.71828$ denoting the usual number):

**Proposition 12C.3.** *Consider* $\mathbf{fdMod}\big(\mathbb{C}[S_n]\big)$.

- **(a)** *The sum of the dimensions of the simple objects is* $s_n = \sum_{k=0}^{n/2} n!/((n-2k)!2^k k!)$.
- **(b)** *We have* $\lim_{n\to\infty} \sqrt[n]{s_n} \sim \lim_{n\to\infty} \sqrt[n]{\sqrt{n!}} \sim \sqrt{n}/\sqrt{e}$.

*Thus, the sum of the dimensions of the simple objects of* $\mathbf{fdMod}\big(\mathbb{C}[S_n]\big)$ *is* $\sum_{simples} \dim \mathrm{L} \approx \sqrt{\# S_n}$.

*Proof.* We only sketch a proof; details are supposed to filled in by the reader, see Exercise 12J.3.

*(a).* We begin with the following fact. We assume the reader has some background in representation theory, in particular, in character theory.

**Lemma 12C.4.** *Let $n(g)$ be the number of square roots of $g \in G$. Then*

$$n(g) = \sum_{\chi} s(\chi)\chi(g),$$

*where the sum runs over all simple characters $\chi$ of $G$, and $s(\chi)$ is the Frobenius–Schur indicator.*

*Proof.* Since $n(g)$ is a class function, it decomposes as a linear combination of simple characters. The coefficient of $\chi$ is given by

$$\langle n, \chi\rangle = \frac{1}{|G|}\sum_{g\in G} n(g)\chi(g) = \frac{1}{|G|}\sum_{g\in G}\sum_{h\in G} \delta_{g,h^2}\chi(g).$$

Interchanging sums gives

$$\langle n, \chi\rangle = \frac{1}{|G|}\sum_{h\in G} \chi(h^2) = s(\chi).$$

The lemma is proved. □

Taking $g = 1$, we get

$$n(1) = \sum_{\chi} s(\chi)\chi(1).$$

If all $s(\chi) = 1$, then $n(1) = \sum_{\chi} \dim L_{\chi}$, where $n(1)$ is the number of involutions in $G$. For $G = S_n$, all indicators are one, so the result follows.

*(b).* Now, the number of elements of the symmetric group is $n!$, which, using Stirling's formula, gives

$$\lim_{n\to\infty} \sqrt[n]{\# S_n} \sim n/e.$$

Since

$$\sum_{simples} (\dim \mathrm{L})^2 = \# S_n,$$

we expect our count to be of the order of the square root of $n!$. Indeed, if $\mathbf{fdMod}\big(\mathbb{C}[S_n]\big)$ would have one simple object, (b) would follow immediately. The other extreme is that all simple objects are of the same dimension $x$, and since the number of simple objects is the number of integer partitions $p(n)$, as in [**OEI23**, A000041], we have

$$n! = \sum_{simples} (\dim \mathtt{L})^2 = p(n)x^2$$

in this case. Using $\lim_{n\to\infty} \sqrt[n]{p(n)} = 1$ we get

$$\lim_{n\to\infty} \sqrt[n]{x^2} = \lim_{n\to\infty} \sqrt[n]{n!/p(n)} \sim n/e,$$

as desired. □

12D. **Growth problems for groups.** Let $G$ be a group. We are interested in studying the growth of the monoidal product in $\mathbf{fdMod}\big(\Bbbk[G]\big)$. We do this as follows.

For any $\mathtt{X} \in \mathbf{fdMod}\big(\Bbbk[G]\big)$ over $\Bbbk$ let $\nu(\mathtt{X}) \in \mathbb{Z}_{\geq 0}$ be an integer such that

$$\mathtt{X} \cong \bigoplus_{i=1}^{\nu(\mathtt{X})} \mathtt{Z}_i,$$

where $\mathtt{Z}_i \in \mathbf{fdMod}\big(\Bbbk[G]\big)$ are indecomposable objects. This number is well defined by the Krull–Schmidt theorem *cf.* Theorem 6J.19. In other words,

$$\nu(\mathtt{X}) = \sum_{\mathtt{Z}\in\mathtt{In}} (\mathtt{X} : \mathtt{Z}).$$

Similarly, define (or redefine) the length $\ell(\mathtt{X})$ by

$$\ell(\mathtt{X}) = \sum_{\mathtt{L}\in\mathtt{Si}} [\mathtt{X} : \mathtt{L}].$$

Now let $\mathtt{V} \in \mathbf{fdMod}\big(\Bbbk[G]\big)$ and define the integer sequences

$$b_n = b_n(\mathtt{V}) := \nu(\mathtt{V}^n) = \text{number of indecomposables summands in } \mathtt{V}^n \text{ counted with multiplicities}, \; n \in \mathbb{Z}_{\geq 0},$$
$$l_n = l_n(\mathtt{V}) := \ell(\mathtt{V}^n) = \text{number of simples factors in } \mathtt{V}^n \text{ counted with multiplicities}, \; n \in \mathbb{Z}_{\geq 0}.$$

Arguably, $b_n$ is the more important number as it measures the "correct" decomposition complexity of $\mathtt{V}^n$ and we will stick with it most of the time.

**Example 12D.1.** For $\mathtt{V} = \mathtt{L}_s \in \mathbf{fdMod}\big(\mathbb{C}[S_3]\big)$ as in Example 8C.6 we get the following count (in this case $b_n = l_n$):

$$\begin{gathered}
\mathtt{V}^0 \cong \mathbb{1} \Rightarrow b_0 = 1,\\
\mathtt{V}^1 = \mathtt{L}_s \Rightarrow b_1 = 1,\\
\mathtt{V}^2 \cong \mathtt{L}_1 \oplus \mathtt{L}_s \oplus \mathtt{L}_{1'} \Rightarrow b_2 = 3,\\
\mathtt{V}^3 \cong \mathtt{L}_1 \oplus 3\mathtt{L}_s \oplus \mathtt{L}_{1'} \Rightarrow b_3 = 5,\\
\mathtt{V}^4 \cong 3\mathtt{L}_1 \oplus 5\mathtt{L}_s \oplus 3\mathtt{L}_{1'} \Rightarrow b_4 = 11,
\end{gathered}$$

and so forth. Precisely,

$$\big(\mathtt{V}^n \cong J(n-1)\mathtt{L}_1 \oplus J(n)\mathtt{L}_s \oplus J(n-1)\mathtt{L}_{1'}\big) \Rightarrow \big(b_n = J(n) + 2J(n-1)\big),$$

where $J(n)$ is the $n$th ***Jacobsthal number*** as, for example, in [**OEI23**, A001045]. From [**OEI23**, A001045] (or many other sources) we get $J(n) =$ nearest integer to $\frac{1}{3}2^n$. We get

$$b_n \sim \frac{2}{3} \cdot n^0 \cdot 2^n.$$

Here, of course, $n^0 = 1$ can be ignored, but we keep it because of Equation 12D-6 below. ◇

**Example 12D.2.** Let $G$ be the dihedral group of order eight. Consider $\mathbf{fdMod}\big(\mathbb{C}[G]\big)$. The group $G$ has one simple representation of dimension two $\mathtt{V}$, let us it for our decomposition problem. Using MAGMA as in Appendix A, we can set this up its character table as follows.

```
> X:=CharacterTable(DihedralGroup(4));
> V:=X[#X];
> for n in [0..10] do
> &+Decomposition(X,V^n);
> end for
  ----result----
```

```
1
1
4
4
16
16
64
64
256
256
1024
```

From there it is easy to guess what the sequence for $b_n = l_n$ is:

$$b_n = 4^{\lfloor n/2 \rfloor}.$$

It is then not difficult to see that

$$b_n \sim \boxed{\tfrac{1}{8}\big(6 + (-1)^n 2\big)} \cdot \boxed{n^0} \cdot \boxed{2^n}.$$

All of this can be verified formally, see Exercise 12J.2. The appearing $f(n) = \frac{1}{8}\big(6 + (-1)^n 2\big)$ is:

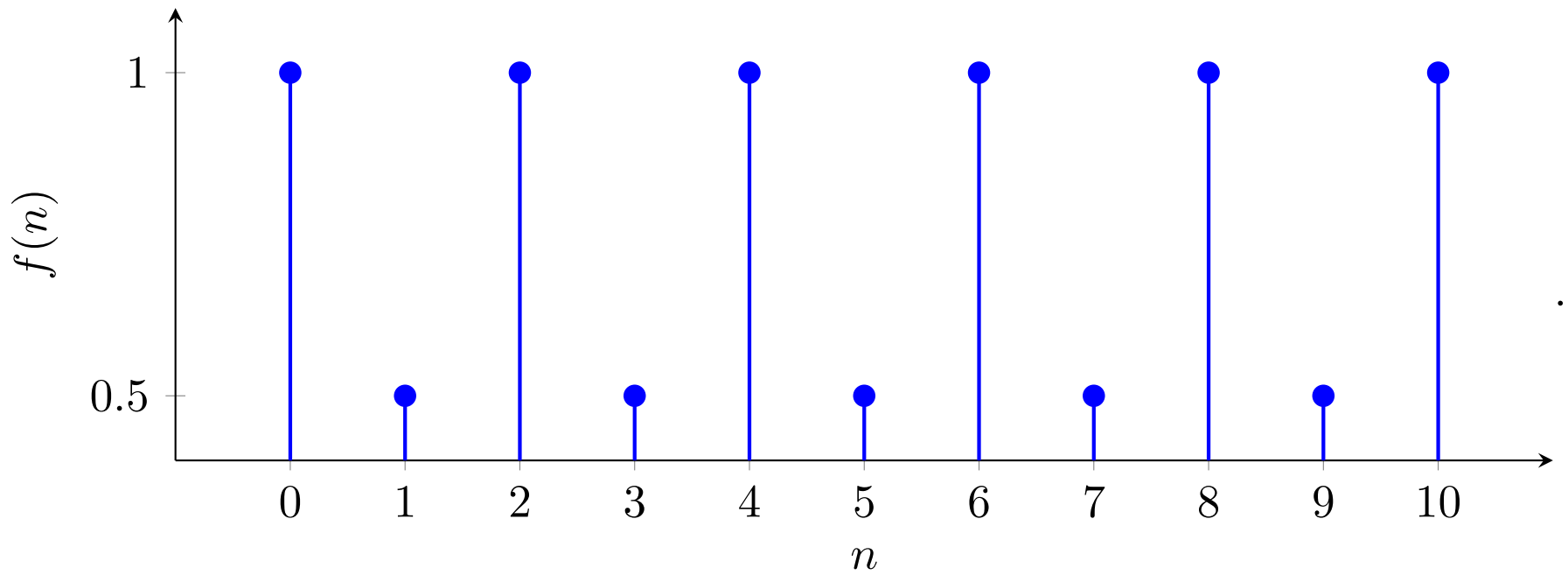

.

It alternates between two constants. ◇

**Example 12D.3.** Take $G = \mathrm{SL}_2(\mathbb{C})$, and we take $\mathtt{V} = \mathbb{C}^2 \in \mathbf{fdMod}\big(\mathbb{C}[G]\big)$, the vector representation. The category $\mathbf{fdMod}\big(\mathbb{C}[G]\big)$ is the $k \to \infty$ limit of the Verlinde category in Example 8H.4, and has the (not truncated) Clebsch–Gordon rules as in Example 8H.4. In particular, the fusion graph for $\mathtt{V} = \mathtt{L}_1$ is the line graph on $\mathbb{Z}_{\geq 0}$:

$$\Gamma_1 = \mathbb{1} \rightleftarrows \mathtt{L}_1 \rightleftarrows \mathtt{L}_2 \rightleftarrows \dots \rightleftarrows \mathtt{L}_k \rightleftarrows \dots .$$

Then the first numbers $b_n$ are:

$$\{1, 1, 2, 3, 6, 10, 20, 35, 70, 126, 252\}, \quad b_n \text{ for } n = 0, \dots, 10.$$

One can then easily guess that these are the ***middle binomial coefficients*** as, for example, in [**OEI23**, A001405]. Assuming that this is true we get

$$b_n \sim \boxed{\big(\tfrac{\pi}{2}\big)^{-1/2}} \cdot \boxed{n^{-1/2}} \cdot \boxed{2^n}.$$

for free (alternatively, from Stirling's formula). Using $\sqrt{2/\pi} \approx 0.798$ we get

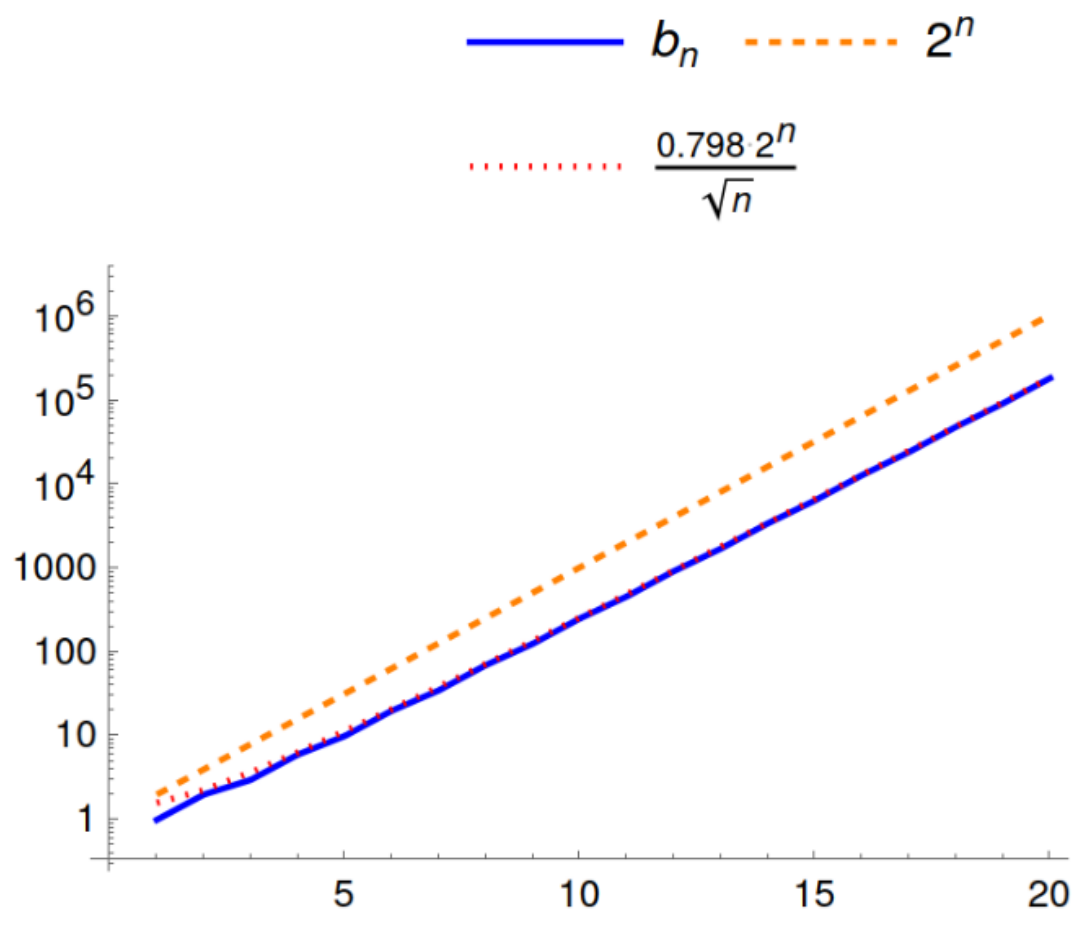

,
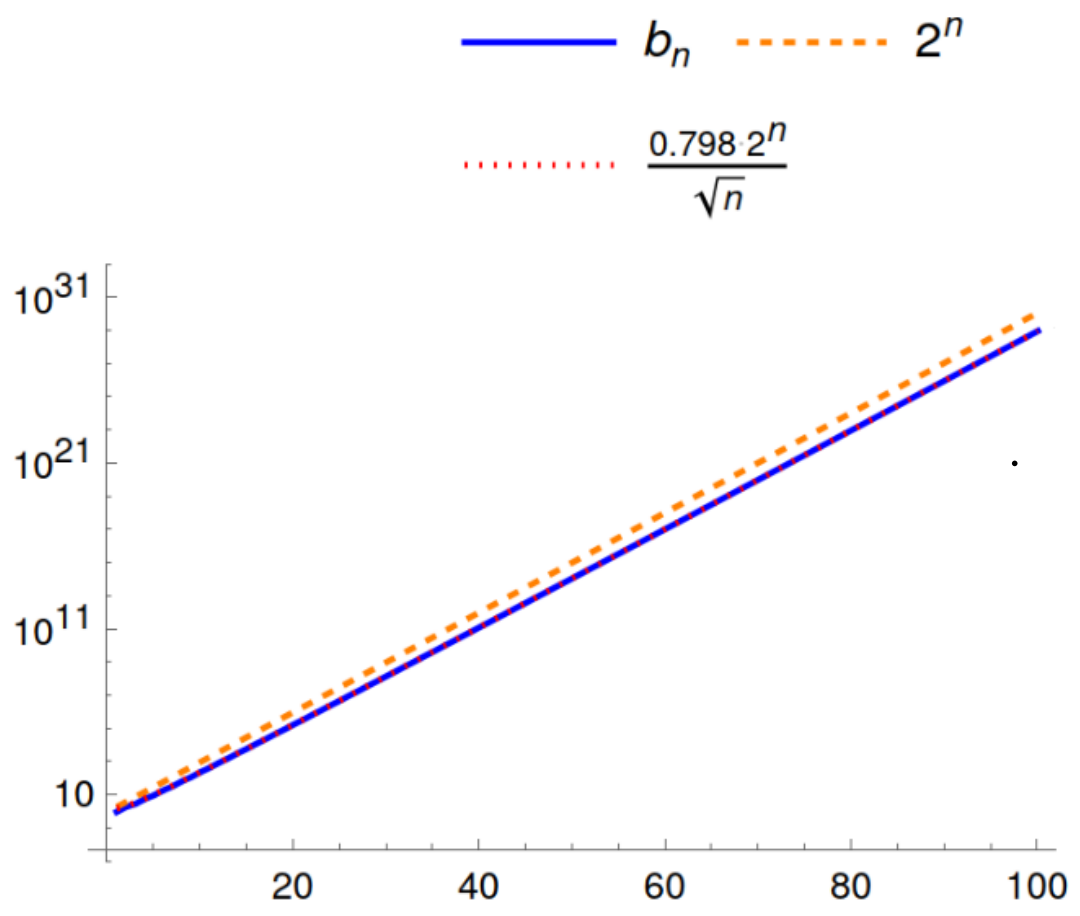

.

as log plots. It remains to argue that $b_n$ is indeed the sequence of middle binomial coefficients. However, that is not too difficult with some Lie theory at hand (the reader unfamiliar might want to have a look at, *e.g.*, [**FH91**] and its discussion of the Lie algebra of $G$). The character of $\mathtt{V}$ is $v^{-1}+v$. Hence, the character of $\mathtt{V}^n$ is $(v^{-1}+v)^n$. Now, the character of $\mathtt{L}_k$ is $[k+1]_v$ (the quantum numbers as in Equation 9H-1), and all of them have either $v^0$ or $v^1$ appearing precisely once, depending on parity. Hence, we are done since $b_n$ is thus the coefficient of either $v^0$ or $v^1$ in $(v^{-1}+v)^n$, depending on parity. ◇

**Example 12D.4.** Generalizing Example 12D.3, $\mathrm{SL}_2(\mathbb{C})$ has a simple representation of dimension $k+1$ for all $k\in\mathbb{Z}_{\geq 0}$, corresponding to $\mathtt{L}_k$ in the Verlinde category. One can check that $b_n$ is given by the ***middle*** $(k+1)$***nomial coefficients*** (the largest coefficient of $(1+x+x^2+...+x^k)^n$) and one gets

$$b_n \sim \boxed{\left(\tfrac{k(k+2)\pi}{6}\right)^{-1/2}} \cdot \boxed{n^{-1/2}} \cdot \boxed{(k+1)^n},$$

directly from well-known results, see [**OEI23**, A002426] for $k=2$. ◇

The study of the growth of $b_n$ and $l_n$ is the main motivation of this section.

**Lemma 12D.5.** *We have* $1\leq b_n\leq l_n\leq (\dim_{\mathbb{k}}\mathtt{V})^n$, *and if* $\mathbf{fdMod}\big(\mathbb{k}[G]\big)$ *is semisimple, then* $b_n=l_n$.

*Proof.* We have $b_n\leq l_n$ with equality in the semisimple case by Lemma 6J.4 and Convention 8A.2 (as recalled therein). Now $\dim_{\mathbb{k}}\mathtt{V}^n=(\dim_{\mathbb{k}}\mathtt{V})^n$, so the upper bound follows, and is achieved only if $\mathtt{V}^n$ is indecomposable. The lower bound comes from the observation that the dimensions of simple factors or indecomposable summands are at least one. □

Motivated by Example 12D.1, Example 12D.2 and Example 12D.3, one could hope that

$$\text{(12D-6)}\qquad b_n\sim \boxed{h(n)}\cdot\boxed{n^{\tau}}\cdot\boxed{\beta^n},\quad \begin{cases} h\colon \mathbb{Z}_{\geq 0}\to\mathbb{R}_{>0} \text{ is a function } \textbf{\textit{bounded away from }} 0,\infty,\\ n^{\tau} \text{ is the } \textbf{\textit{subexponential factor}},\ \tau\in\mathbb{R},\\ \beta^n \text{ is the } \textbf{\textit{exponential factor}},\ \beta\in\mathbb{R}_{\geq 1}.\end{cases}$$

In practice, $h(n)$ is often a constant or alternates between finitely many constants, but sometimes $h(n)$ is more complicated.

We do not know in what generality Equation 12D-6 holds, but expressions of this form are very common in the theory of asymptotics of generating functions; see, for example, [**FS09**] or [**Mis20**] for the relation between counting problems and the analysis of generating functions.

Theorem 12D.7 up next is the $\beta$ ***theorem***, or ***leading growth rate theorem***, and its companion Theorem 12D.8.

**Theorem 12D.7.** *For any group* $G$, *any field* $\mathbb{k}$ *and any* $\mathtt{V}\in\mathbf{fdMod}\big(\mathbb{k}[G]\big)$ *we have*

$$\lim_{n\to\infty}\sqrt[n]{b_n}=\dim_{\mathbb{k}}\mathtt{V}.$$

*That is, if Equation 12D-6 holds, then*

$$\boxed{\beta=\dim_{\mathbb{k}}\mathtt{V}}.$$

For a reductive group, let $R^+$ be the number of positive roots of the associated root system. For example, $R^+=1$ for $G=\mathrm{SL}_2(\mathbb{C})$.

**Theorem 12D.8.** *In the notation of Theorem 12D.7, in the following cases Equation 12D-6 holds:*

**(a)** *When* $G$ *is a finite group. In this case*

$$\boxed{\tau=0},\quad \boxed{h(n)\ \textit{alternates between finitely many constants}}.$$

**(b)** *When* $G$ *is a reductive group over* $\mathbb{k}$, $\mathbb{k}$ *is of characteristic zero and* $\mathtt{V}$ *is over* $\mathbb{k}$. *In this case*

$$\boxed{\tau=-R^+/2},\quad \boxed{h(n)\ \textit{is a constant}}.$$

**(c)** *When* $G=\mathrm{SL}_2(\mathbb{k})$ *for the defining field. If* $\mathbb{k}$ *is finite or of characteristic zero, then we are in one of the previous two cases. Otherwise, let* $p$ *be the characteristic of* $\mathbb{k}$. *Then*

$$\boxed{\tau=\tfrac{1}{2}\big(\log_p(\tfrac{p+1}{2})\big)-1},\quad \boxed{h(n)\ \textit{is quite difficult}},$$

*see [**CEOT24**] or [**Lar24**] for a more precise statement.*

Before we prove these theorems, let us give an example of how a proof could work.

**Example 12D.9.** Return to Example 12D.3. If we are only interested in the leading growth rate we could argue as follows. The maximal dimensional summand of $\mathtt{V}^n$ is the $n$th symmetric power of $\mathbb{C}^2$ (this corresponds to $\mathtt{L}_n$), and thus of dimension $n+1$. So, assuming that all summands of $\mathtt{V}^n$ are maximal dimensional, the ***Jupiter approach*** as in Figure 32, we get $2^n/(n+1)$ as a lower bound for $b_n$. Hence,

$$\big(2^n/(n+1) \leq b_n \leq 2^n\big) \Rightarrow \big(\sqrt[n]{2^n/(n+1)} \leq \sqrt[n]{b_n} \leq \sqrt[n]{2^n}\big) \Rightarrow \big(2\sqrt[n]{1/(n+1)} \leq \sqrt[n]{b_n} \leq 2\big) \Rightarrow$$
$$2 = \lim_{n\to\infty} 2\sqrt[n]{1/(n+1)} \leq \lim_{n\to\infty} \sqrt[n]{b_n} \leq \lim_{n\to\infty} 2 = 2,$$

and we have identified the leading growth rate. ◇

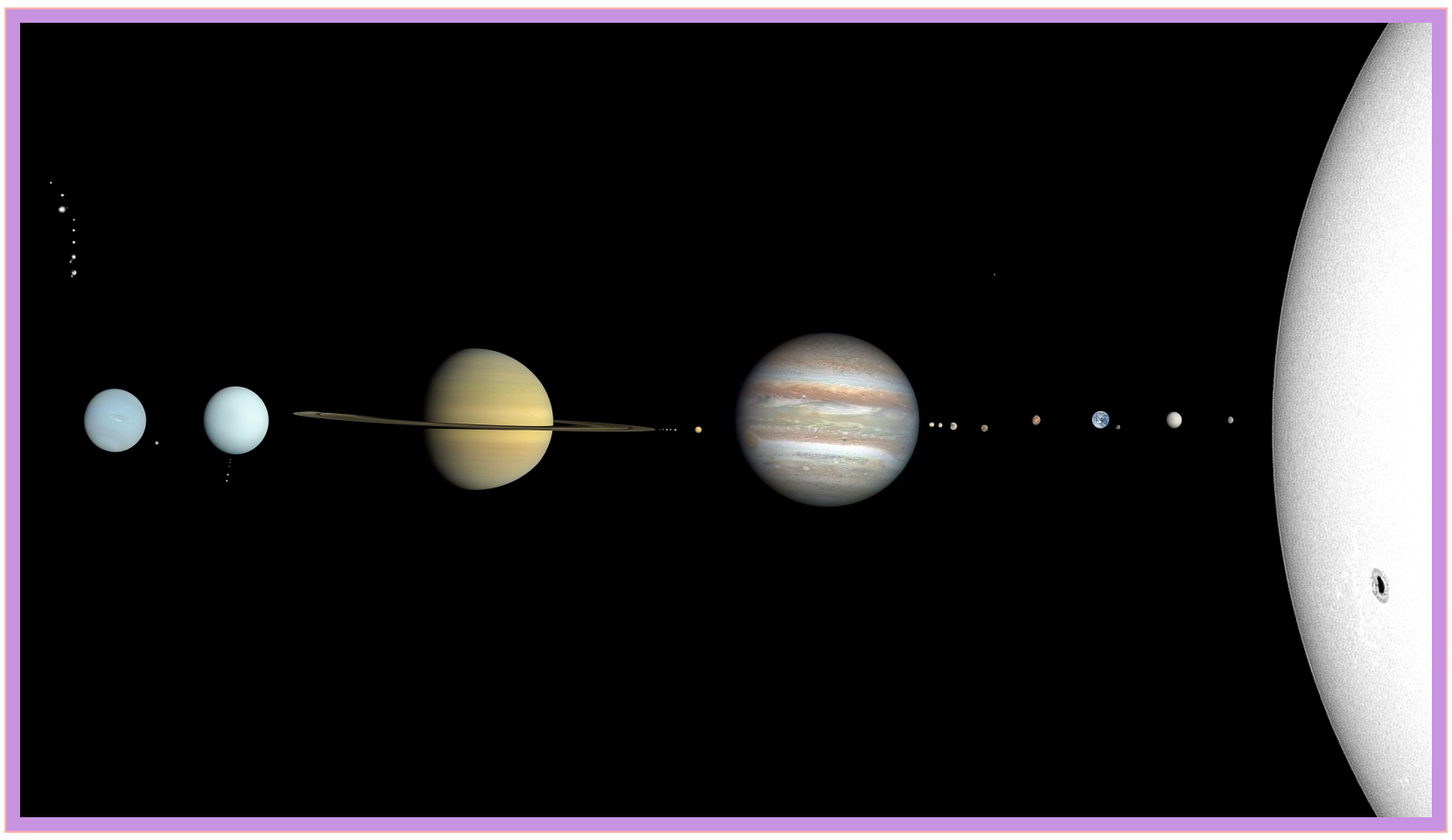

FIGURE 32. Think of our problem of finding $b_n$ as the problem of counting how many planets (the summands) fit into the sun ($\mathtt{V}^n$). To get a precise growth rate of $b_n$ we would need to make a count about the average planet, call it almost-Uranus, and how many fit into the sun, but to get a lower bound it suffices to count how many Jupiters fit into the sun. The sun is so much larger than the planets that this brute force lower bound is actually good enough. For real comparison: $\approx$ 1000 Jupiters fit into the sun, and $\approx$ 4700 almost-Uranuses fit into the sun. (We leave it to the reader to (i) identify the planets in the picture, (ii) double check these numbers generated by ChatGPT.)

Picture from https://en.wikipedia.org/wiki/Solar_System

*Proof of both, Theorem 12D.7 and Theorem 12D.8. The limit.* Let $N = \dim_{\mathbb{k}} \mathtt{V}$. Note that a representation of $G$ on $\mathtt{V}$ is a group homomorphism $G \to \mathrm{GL}(\mathtt{V}) \cong \mathrm{GL}_N(\mathbb{k})$. Thus, the decomposition of $\mathtt{V}^n$ as a $G$-representation is finer than that of $\mathtt{V}^n$ as a $\mathrm{GL}_N(\mathbb{k})$-representation. Hence,

$$b_n(\mathrm{GL}_N) \leq b_n(G).$$

**Lemma 12D.10.** *Theorem 12D.7 is true for $G = \mathrm{GL}_N(\mathbb{k})$, the defining field and $\mathtt{V} = \mathbb{k}^N$ (the vector representation).*

*Proof.* We prove this only for $\mathbb{k} = \mathbb{C}$; the general case is [**COT24**, Proposition 2.2]. For $\mathbb{k} = \mathbb{C}$ ***Weyl's dimension formula*** (see, for example, [**FH91**, Chapter 24]) implies that the maximal dimension of the summands of $\mathtt{V}^n$ is a polynomial $poly(n)$ in $n$. Hence, the Jupiter approach as in Example 12D.9 gives

$$N^n/poly(n) \leq b_n(\mathrm{GL}_N) \leq N^n$$

and the sandwich theorem implies the lemma. □

Now, Lemma 12D.5 and Lemma 12D.10 imply

$$b_n(\mathrm{GL}_N) \leq b_n(G) \leq N^n,$$

and, by sandwiching, we are done.

Thus, by the above discussion the statement follows.

*The list.* We will verify (a) and (b) in this section below. (c) is very difficult to prove, see [**CEOT24**] or [**Lar24**], and we only refer to those papers. □

12E. **Growth problems for monoidal categories.** The below is taken from [LTV24a].

The following mildly generalizes the usual notion of based $\mathbb{R}_{\geq 0}$-algebras as *e.g.* in [KM16, Section 2]. Let $\Bbbk$ be a unital subring of $\mathbb{C}$, for example, $\Bbbk = \mathbb{Z}$. A ***based*** $\mathbb{R}_{\geq 0}$***-algebra*** is a pair $(R, C)$ where $C = (c_i)_{i \in I}$ with $1 \in C$ is a $\Bbbk$-basis of a $\Bbbk$-algebra $R$ such that all structure constants lie in $\mathbb{R}_{\geq 0}$, *i.e.*:

$$c_i c_j = \sum_{k \in I} m_{ij}^k c_k \quad \text{with} \quad m_{ij}^k \in \mathbb{R}_{\geq 0}. \tag{12E-1}$$

Note that $m_{ij}^k \in \Bbbk \cap \mathbb{R}_{\geq 0}$, so that $m_{ij}^k \in \mathbb{Z}_{\geq 0}$ when $\Bbbk = \mathbb{Z}$. We usually only write $R$ for $(R, C)$ if no confusion can arise. Note that the sum in Equation 12E-1 is finite since $C$ is a $\Bbbk$-basis.

*Remark* 12E.2. The reader who wants to consider $\mathbb{R}_{\geq 0}$-algebras without the condition $1 \in C$ can split their growth problem along an idempotent decomposition of the identity, so $1 \in C$ is not a restriction.

By taking a subring of $R$ if necessary, we can, and will, always assume that $I$ is countable and that the identity is the element $c_0$ of the basis $C$. This is no restriction for growth problems as we will see below. ◇

We denote by $(\mathbb{R}_{\geq 0} \cap \Bbbk)C$ the subset of $R$ of finite $(\mathbb{R}_{\geq 0} \cap \Bbbk)$-linear combinations of elements of the basis $C$ with nonnegative real coefficients. Given $c \in (\mathbb{R}_{\geq 0} \cap \Bbbk)C$, we write $cc_j = \sum_{k \in I} m_{c,j}^k c_k$ and we denote $c^n = \sum_{i \in I} m_n^i(c) c_i$. We define the function

$$b^{R,c} \colon \mathbb{Z}_{\geq 0} \to \mathbb{R}_{\geq 0}, n \mapsto b(n) = b^{R,c}(n) := \sum_{i \in I} m_n^i(c).$$

We are interested in the asymptotic behavior of this function. Precisely, the main point of this section is to address, and partially answer, the following questions:

**(a)** What is the ***dominating growth*** of $b_n$?

**(b)** Can we get an ***asymptotic formula*** $a \colon \mathbb{Z}_{\geq 0} \to \mathbb{Z}_{\geq 0}$ expressing $b$, for example,

$$b(n) \sim a(n),$$

where $a$ is "nice", $\sim$ denotes asymptotically equal?

**(c)** Say we have found $a$ as in the previous point. Can we bound the ***variance*** or ***mean absolute difference*** $|b_n - a_n|$, or alternatively the convergence rate of $\lim_{n \to \infty} b(n)/a(n) = 1$?

We refer to **(a)**, **(b)**, and **(c)** as the ***growth problems*** associated with $(R, C)$. We also explore several statements along the same lines, which we will also call growth problems.

*Remark* 12E.3. The order **(a)**, **(b)**, and **(c)** is roughly by difficulty. As an analogy: if $b$ would be the prime counting function, then **(a)** could say that it grows essentially linearly, **(b)** could be the prime number theorem, and **(c)** could be (a consequence of) the Riemann hypothesis. ◇

**Example 12E.4.** If $\mathtt{X}$ and the monoidal unit $\mathbb{1}$ are indecomposable, and we have $\mathtt{X} \otimes \mathtt{X} \cong \mathbb{1} \oplus \mathtt{X}$, then $(b_n)_{n \in \mathbb{Z}_{\geq 0}}$ is the Fibonacci sequence (with the first two terms equal to 1). This happens in the Fibonacci category as in Example 8E.4.

In this case the dominating growth is $\phi^n$ for $\phi \approx 1.618$ the golden ratio. We have $b(n) \sim \frac{\phi}{\sqrt{5}} \cdot \phi^n$ and $|b_n - a_n| < (\phi - 1)^n$ for $\phi - 1 \approx 0.618$.

Note the following: the action matrix of $_ \otimes \mathtt{X}$ acting on $\{\mathbb{1}, \mathtt{X}\}$ is $\left(\begin{smallmatrix} 0 & 1 \\ 1 & 1 \end{smallmatrix}\right)$ has eigenvalues $\phi$ and $-(\phi - 1)$. The largest eigenvalue $\phi$ gives the *dominating* growth, the absolute value of the second largest eigenvalue $\phi - 1$ bounds the variance. ◇

**Example 12E.5.** ***Good examples*** of growth problems (coming from categories) that we can treat nicely are:

**(a)** Let $\mathbf{C} = \langle \mathtt{X} \rangle$ be the additive Krull–Schmidt monoidal category generated by $\mathtt{X}$ meaning that we take direct summands of finite sums of objects $\mathtt{X}^{\otimes n}$, where $n \in \mathbb{Z}_{\geq 0}$ (we call such $\mathtt{X}$ generating objects). This example is good if $\mathbf{C}$ is a finite tensor category over an arbitrary field $\mathbb{F}$ in the sense of [EGNO15, Definition 1.8.5 and Definition 4.1.1] or finite with respect to isomorphism classes of indecomposable objects. This includes:

- **(i)** $(H, \mathtt{X}) = \big(\mathbf{fdMod}(H), \mathtt{X}\big)$ for $H$ a finite dimensional Hopf $\mathbb{F}$-algebra and $\mathtt{X}$ a finite dimensional $H$-representation over $\mathbb{F}$ that generates $\mathbf{fdMod}(H)$. Specifically, $H$ could be the group ring of a finite group and $\mathtt{X}$ a faithful $H$-representation.
- **(ii)** $(\mathbf{SBim}(W), \mathtt{X})$, for a generating object $\mathtt{X}$, or any $\mathtt{X}$ such that some monoidal power of it has the indecomposable Soergel bimodule for the longest word of $W$ appearing as a direct summand (for $\mathbf{SBim}(W)$ we use generating object to refer to either case), see for example [EW16] for details about the category $\mathbf{SBim}(W)$.

**(b)** $(G, \mathtt{X}) = \big(\mathbf{fdMod}(\Bbbk[G]), \mathtt{X}\big)$ where $G$ is a simple reductive group and $\mathtt{X}$ is a $G$-representation, both in characteristic zero. (The assumption on $G$ being simple can be relaxed, see Example 12E.26.)

**(c)** $(\mathbf{C}, \mathtt{X})$ for either of the following: $\mathbf{C} = \mathrm{SL}_{\mathbb{Z}_{\geq 0}}(\mathbb{C})$ and $\mathtt{X} = \mathbb{C}^{\mathbb{Z}_{\geq 0}}$ the defining representation; $\mathbf{C} = \mathit{rep}(S_t, \mathbb{C})$ for $t \in \mathbb{C} \setminus \mathbb{Z}$ and $\mathtt{X}$ the generating object, following the conventions in [CO11], and $\mathbf{C}$ a Delannoy category with $\mathtt{X}$ the (defining) generating object following the conventions in [HSS22] (and [HS22, Theorem 13.2] shows that this is a good example).

We will revisit these several times below. ◇

Let us consider the ***(oriented and weighted) fusion graph*** $\Gamma = \Gamma(c)$ associated to some growth problem $(R, c)$. This graph is defined as follows:

**(i)** The vertices of $\Gamma$ correspond to the basis elements $c_i \in C$ appearing in $c^n$ for some $n$, *i.e.* those $c_i$ for which some $m_n^i \neq 0$.

**(ii)** There is an edge of labeled $m_{c,j}^i$ from the vertex $c_j$ to the vertex $c_i$.

We identify the basis elements $c_i \in C$ with the vertices of $\Gamma$. After fixing some ordering of the vertices, the matrix associated to a fusion graph $\Gamma$ is called the ***action matrix*** $M(\Gamma) = (m_{c,j}^k)_{k,j}$.

**Example 12E.6.** The key example the reader should keep in mind is the growth problem for $\big(\mathrm{SL}_2(\mathbb{C}), \mathbb{C}^2\big)$. In this case the fusion graph is $\mathbb{Z}_{\geq 0}$:

$$\Gamma(\mathrm{SL}_2) = \bullet \rightleftarrows \bullet \rightleftarrows \bullet \rightleftarrows \bullet \rightleftarrows \bullet \rightleftarrows \bullet \rightleftarrows \cdots,$$

(12E-7)

$$M(\mathrm{SL}_2) = \begin{pmatrix} 0 & 1 & 0 & 0 & 0 & \dots \\ 1 & 0 & 1 & 0 & 0 & \dots \\ 0 & 1 & 0 & 1 & 0 & \dots \\ 0 & 0 & 1 & 0 & 1 & \dots \\ 0 & 0 & 0 & 1 & 0 & \dots \\ \vdots & \vdots & \vdots & \vdots & \vdots & \vdots \end{pmatrix}.$$

We will come back to this growth problem several times throughout the paper. ◇

**Example 12E.8.** The Grothendieck group of an additive Krull–Schmidt monoidal category $\mathbf{C}$ is an example of a based $\mathbb{R}_{\geq 0}$-algebra with $\Bbbk = \mathbb{Z}$ with basis given by the classes of indecomposable objects. The growth problems as above can then be reformulated in terms of growth problems for based $\mathbb{R}_{\geq 0}$-algebras, and we do this silently throughout. In this case the fusion graph $\Gamma(\mathtt{X})$ for $\mathtt{X} \in \mathbf{C}$ is the fusion graph on the subcategory generated by $\mathtt{X}^{\otimes n}$ in the usual sense. ◇

**Lemma 12E.9.** *For a fusion graph $\Gamma$ for a growth problem we have:*

**(a)** *$\Gamma$ has countably many vertices.*

**(b)** *There is a path from $c_0$ to any other vertex of $\Gamma$. In particular, $\Gamma$ is connected in the unoriented sense.*

**(c)** *Every vertex of $\Gamma$ is of finite degree.*

*Proof.* This follows because we only consider $c^n$, which, for fixed $n$, has only finitely many nonzero coefficients in terms of $C$. □

The ***naive cutoff*** of $\Gamma$ is the sequence $(\Gamma_k)_{k \in \mathbb{Z}_{\geq 0}}$ of subgraphs of $\Gamma$ of the form

$$(\Gamma_k)_{k \in \mathbb{Z}_{\geq 0}} = \big(\Gamma_0 = \{c_0\} \subset \Gamma_1 \subset \Gamma_2 \subset \dots \subset \Gamma\big) \text{ with } \bigcup_{k \in \mathbb{Z}_{\geq 0}} \Gamma_k = \Gamma,$$

such that $\Gamma_k$ is the induced subgraph of $\Gamma$ that has precisely the vertices that can be reached with a path of length $\leq k$ from $c_0$. Note that Lemma 12E.9.(c) implies that $\Gamma_k$ will be finite.

**Definition 12E.10.** A ***good filtration*** of a growth problem $(R, C)$ is a sequence of (not necessarily associative) $\mathbb{R}_{\geq 0}$-algebras $\big((R_k, C_k)\big)_{k \in \mathbb{Z}_{\geq 0}}$ over $\Bbbk$ such that:

**(a)** $C_k \subset C$ is finite and

$$(C_k)_{k \in \mathbb{Z}_{\geq 0}} = \big(C_0 = \{c_0\} \subset C_1 \subset \dots \subset C_k \subset \dots \subset C\big) \text{ with } \bigcup_{k \in \mathbb{Z}_{\geq 0}} C_k = C.$$

**(b)** $c_0$ is the unit of all $R_k$.

**(c)** The multiplication in $R_k$ is given by taking the product in $(R, C)$ followed by the projection to $(R_k, C_k)$.

The ***naive filtration*** of $(R, C)$ is the good filtration where $C_k$ are the vertices of $\Gamma_k$ for the naive cutoff. ◇

**Example 12E.11.** In Equation 12E-7 the naive cutoff consists of $\Gamma_k$ being a line graph with $k$ vertices. In Equation 12F-16 below, all the $\Gamma_k$ for $k > 1$ will include the vertex labeled $\mathrm{P}_4$. ◇

For a good filtration, in $(R_k, C_k)$ (this is of finite rank over $\Bbbk$) we can define $\mathrm{PFdim}_k c_i$ to be the ***PF eigenvalue*** of the action matrix of $c_i$ with respect to $C_k$.

**Lemma 12E.12.** *We have* $\mathrm{PFdim}_k c_i \leq \mathrm{PFdim}_{k+1} c_i$ *and* $\lim_{k \to \infty} \mathrm{PFdim}_k c_i \in \mathbb{R}_{\geq 0} \cup \{\infty\}$ *is well-defined.*

*Proof.* The first statement follows from the fact that subgraphs of a graphs cannot have larger PF eigenvalues than the original graph. The second statement then follows since the previous point gives us an increasing sequence. □

**Definition 12E.13.** Given $c \in (\mathbb{R}_{\geq 0} \cap \Bbbk)C$, define the (good filtration) ***PF dimension*** as the limit $\mathrm{PFdim}^f c = \lim_{k\to\infty} \mathrm{PFdim}_k c \in \mathbb{R}_{\geq 0} \cup \{\infty\}$. ◇

Note that Definition 12E.13 depends on the choice of a good filtration, see Example 12E.26 below for an almost example (illustrating why one needs to be careful with the choice of the filtration).

**Lemma 12E.14.** *If $|I| < \infty$, then Definition 12E.13 agrees with the usual definition of PF dimension using the largest eigenvalue of the associated matrix (or, equivalently, graph).*

*Proof.* Immediate. □

**Definition 12E.15.** A growth problem $(R, C)$ is called ***PF admissible*** if there exists a good filtration such that:

**(a)** $\mathrm{PFdim}^f$ is superadditive and subpowermultiplicative:

$$\mathrm{PFdim}^f c + \mathrm{PFdim}^f d \leq \mathrm{PFdim}^f (c + d),$$
$$\mathrm{PFdim}^f c^n \leq (\mathrm{PFdim}^f c)^n,$$

for all $c, d \in (\mathbb{R}_{\geq 0} \cap \Bbbk)C$, $n \in \mathbb{Z}_{\geq 0}$.

**(b)** $\mathrm{PFdim}^f c_i \geq 1$ for all $c_i \in C$.

We always associate any such good filtration to a PF admissible growth problem. We use the same terminology for the fusion graphs and action matrices themselves. ◇

*Remark* 12E.16. If $R$ is of finite $\Bbbk$-rank, then two good filtrations will agree on all but finitely many pages of the filtrations. Abusing terminology a bit, we do not need to and will not specify a specific good filtration. ◇

**Example 12E.17.** A growth problem coming from a transitive unital $\mathbb{Z}_{\geq 0}$-ring of finite $\mathbb{Z}$-rank is PF admissible. (This includes decategorifications of fusion categories.) The condition $\mathrm{PFdim}^f c_i \geq 1$ for all $c_i \in C$ holds in this case, see e.g. [**EGNO15**, Proposition 3.3.4]. In fact, this assumption is motivated from these rings. ◇

**Example 12E.18.** By a classical result of Kronecker on algebraic integers, as summarized in e.g. [**EGNO15**, Lemma 3.3.14], every growth problem corresponding to a finite graph satisfies Definition 12E.15.(b) unless the PF dimension is zero. ◇

**Notation 12E.19.** A function $f$ is of ***superexponential growth*** if $f \in \Omega(\gamma^n)$ for all $\gamma \in \mathbb{R}_{>1}$ (in particular, $\sqrt[n]{f(n)}$ is unbounded in this case), and is of ***subexponential growth*** if $f \in O(\gamma^n)$ for all $\gamma \in \mathbb{R}_{>1}$. ◇

**Theorem 12E.20.** *Let $(R, c)$ be a growth problem.*

**(a)** *For a good filtration we have*

$$b(n) \in \Omega(\gamma^n) \text{ for all } \gamma < \mathrm{PFdim}^f c. \tag{12E-21}$$

*Assume now that $(R, c)$ is PF admissible.*

**(b)** *If $\mathrm{PFdim}^f c \neq \infty$, then*

$$b(n) \in O\big((\mathrm{PFdim}^f c)^n\big). \tag{12E-22}$$

**(c)** *We have the* **exponential growth theorem***:*

$$\beta = \lim_{n\to\infty} \sqrt[n]{b_n} = \mathrm{PFdim}^f c.$$

*In particular, $b_n$ grows superexponentially if and only if $\mathrm{PFdim}^f c = \infty$.*

In general the PF dimension depends on the choice of a good filtration. However, if one has two different good filtrations such that the growth problem is PF admissible, then the PF dimension is the same for both of them by Theorem 12E.20.(c).

*Proof of Theorem 12E.20.* As in the proof of [**LTV23**, Theorem 3B.2], the number $b_n$ can be computed as the column sum of the matrix $M(\Gamma)^n$ for the first column corresponding to the identity. (We will use this throughout.) To see this we consider the equation

$$M(\Gamma)c(n-1) = c(n) \xrightarrow{\text{iterate}} M(\Gamma)^n c(0) = c(n),$$

where $c(k) = \big(c_0(k), c_1(k), ...\big) \in \mathbb{R}_{\geq 0}^{\mathbb{Z}_{\geq 0}}$ are vectors such that their $i$th entry is the multiplicity of $c_i$ in $c^k$, and $c(0) = (1, 0, ...)^T$ with the one is in the slot of $c_0 = 1$. Observe $M(\Gamma)^n c(0)$ is the first column of $M(\Gamma)^n$.

*(a).* By an appropriate version of Feteke's subadditive lemma, the sequence $(\sqrt[n]{b_n})_{n\in\mathbb{Z}_{\geq 0}}$ has a limit

$$\lim_{n\to\infty} \sqrt[n]{b_n} = \beta \in \mathbb{R}_{\geq 0} \cup \{\infty\}.$$

Similarly, if we denote by $(b_n^{(k)})_{n\in\mathbb{Z}_{\geq 0}}$ the respective growth problem in the cutoff $R_k$ for $c \in R_k$ and set $b_n^{(k)} = 0$ whenever $c \notin R_k$ (this happens for only finitely many $k \in \mathbb{Z}_{\geq 0}$), then

$$\lim_{n\to\infty} \sqrt[n]{b_n^{(k)}} = \beta_k \in \mathbb{R}_{\geq 0} \cup \{\infty\},$$

exists. Note hereby that $\beta_k = \mathrm{PFdim}_k c$, by [LTV23, Theorem 1] (which uses PF theory of finite graphs). Moreover, the definition of a good filtration Definition 12E.10 implies that $b_n^{(k)} \leq b_n$ which in turn implies that $\beta_k \leq \beta$ for all $k, n \in \mathbb{Z}_{\geq 0}$. Therefore, $\mathrm{PFdim}^f c \leq \beta$ by Definition 12E.13. Hence, $b(n) \in \Omega(\gamma^n)$ for all $\gamma < \mathrm{PFdim}^f c$.

*(b).* We go back to the definition of $b_n$. We have $c^n = \sum_{i\in I} m_n^i(c) c_i$ and we obtain:

$$(\mathrm{PFdim}^f c)^n \geq \mathrm{PFdim}^f c^n \geq \sum_{i\in I} m_n^i(c)\mathrm{PFdim}^f c_i \geq \sum_{i\in I} m_n^i(c) = b_n,$$

the inequalities following from the PF admissibility condition.

*(c).* This follows from the just established Equation 12E-21 and Equation 12E-22. □

**Example 12E.23.** Coming back to Example 12E.6, the cutoff $(\Gamma_i)_{i\in\mathbb{Z}_{\geq 0}}$ of $\Gamma(\mathrm{SL}_2)$ is such that $\Gamma_i$ is the finite subgraph on the first $i$ vertices exemplified by:

$\Gamma_7 =$ 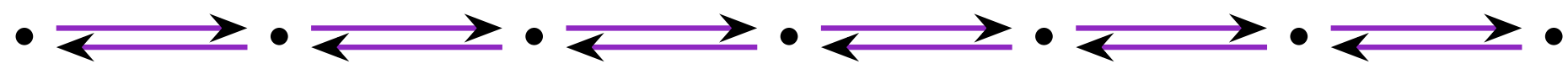 .

We will prove later that this growth problem with the naive truncation is PF admissible, see Proposition 12E.25. It is easy to see that $\mathrm{PFdim}_i\Gamma_i = 2\cos\big(\pi/(i+1)\big)$. Then $\mathrm{PFdim}^f\Gamma(\mathrm{SL}_2) = \lim_{i\to\infty} 2\cos\big(\pi/(i+1)\big) = 2$. Thus, Theorem 12E.20 implies that the dominating growth factor of $b_n$ is $2^n$. Moreover, dividing by the dominating growth factor, the log plots

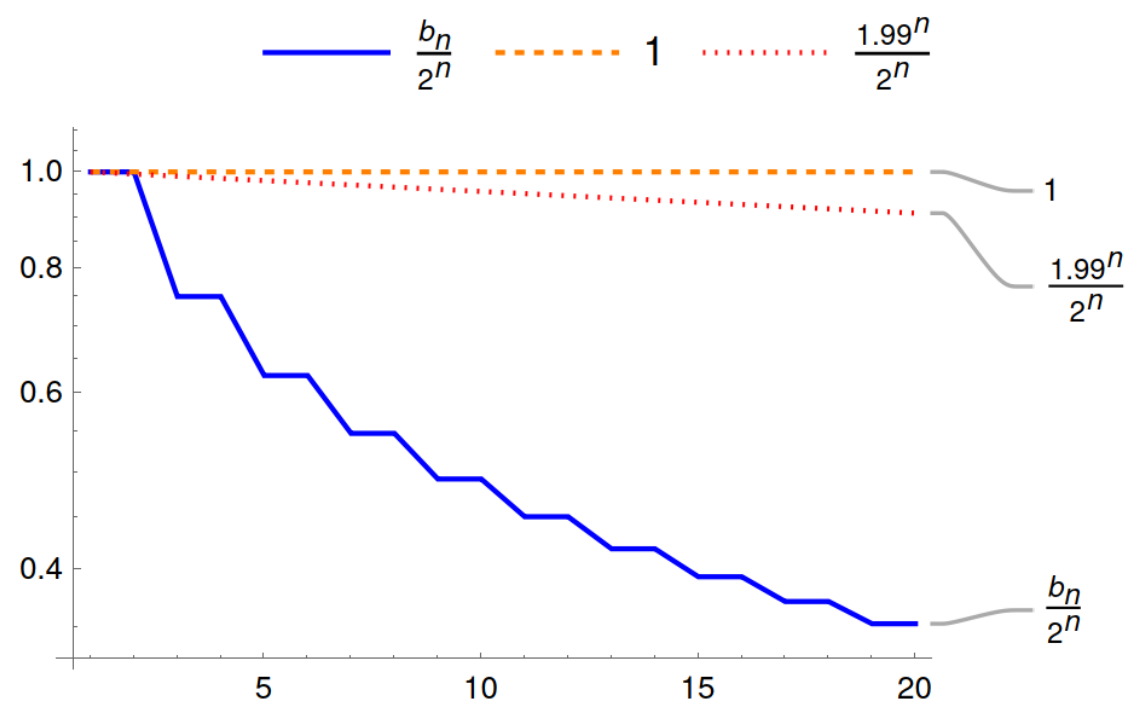

,
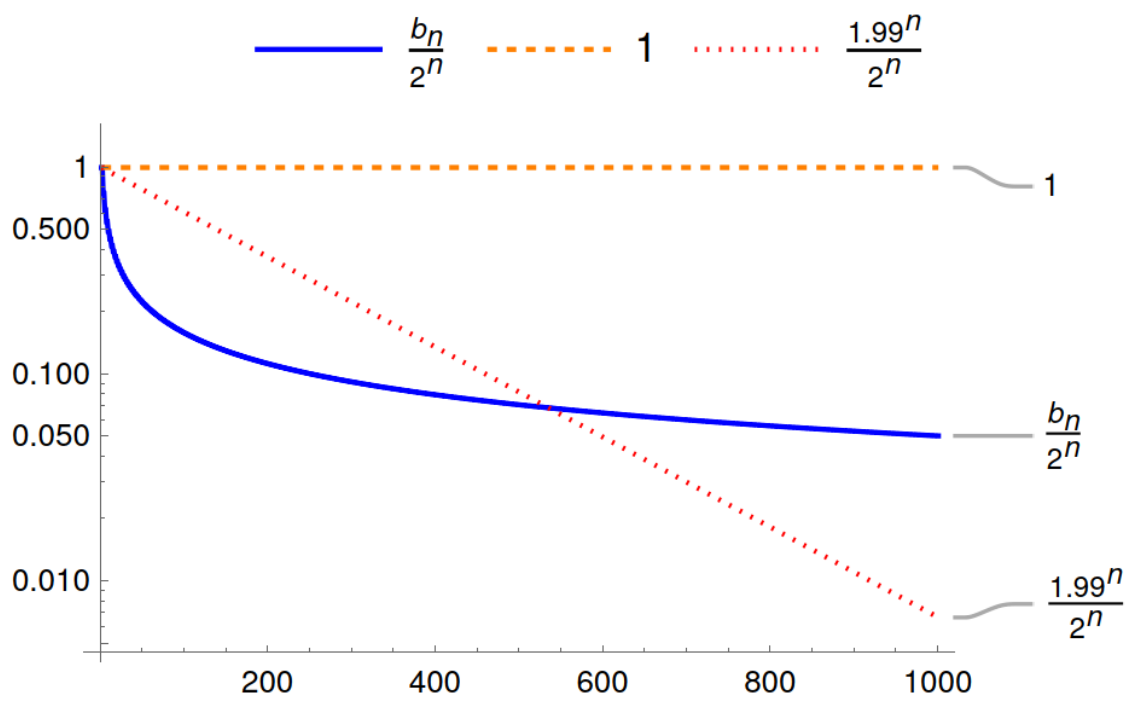

,

shows Equation 12E-21 for $\gamma = 1.99$ and Equation 12E-22. ◇

Next, a "nonexample":

**Example 12E.24.** We give an example where the PF admissibility condition is necessary for Theorem 12E.20. Suppose that our $\mathbb{R}_{\geq 0}$-algebra has a basis $(c_i)_{i\in\mathbb{Z}_{\geq 0}}$ with multiplication $c_ic_j = c_{i+j}$ and that $c = \alpha c_0 + c_1$ for some $\alpha \in \mathbb{R}_{\geq 0}$. Then $c$ satisfies $cc_i = \alpha \cdot c_i + c_{i+1}$ for all $i \in \mathbb{Z}_{\geq 0}$. The action matrix is an infinite Jordan block with $\alpha$-labeled loops and we have

$\Gamma(Jo) =$ 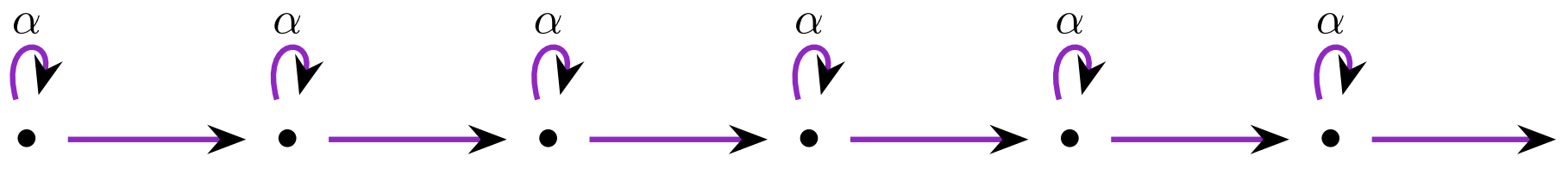

,

$$M(Jo) = \begin{pmatrix} \alpha & 0 & 0 & 0 & 0 & \dots \\ 1 & \alpha & 0 & 0 & 0 & \dots \\ 0 & 1 & \alpha & 0 & 0 & \dots \\ 0 & 0 & 1 & \alpha & 0 & \dots \\ 0 & 0 & 0 & 1 & \alpha & \dots \\ \vdots & \vdots & \vdots & \vdots & \vdots & \vdots \end{pmatrix},$$

as the fusion graph. The naive cutoff is the sequence of finite Jordan blocks, in particular, $\mathrm{PFdim}^f c_i = 0$ for $i > 0$ and $\alpha = 0$, so this example does not satisfy Definition 12E.15.(b). And in fact we get:

$$c^n = \sum_{i=0}^{n} \binom{n}{i} \alpha^{n-i} \cdot c_i,$$

and therefore $b_n = (1+\alpha)^n$ but $\mathrm{PFdim}^f c = \alpha$. ◇

Here are key examples of categories with infinitely many indecomposable objects up to isomorphism that are PF admissible:

**Proposition 12E.25.** *Any growth problem in Example 12E.5.(b) is PF admissible.*

*Proof.* Recall that we only consider simple reductive groups. We will make use of the ***Verlinde categories*** $\mathbf{V}_k(G)$ of level $k$, see for example [**Saw06**] for a nice summary of the main properties we need (when we refer to it we use the arXiv numbering). For example, that the PF dimension is the categorical dimension.

In these semisimple categories the indecomposable (= simple) objects $\overline{L}(\lambda,k)$ are indexed by $\Lambda^k$, which are cutoffs of the dominant integral Weyl chamber, see [**Saw06**, Section 4] for details. Moreover, we have

$$\mathrm{PFdim}_k\overline{L}(\lambda,k)=\dim_{\mathbf{V}_k(G)}\overline{L}(\lambda,k)=\prod_{\beta>0}\frac{q^{\langle\lambda+\rho,\beta\rangle}-q^{-\langle\lambda+\rho,\beta\rangle}}{q^{\langle\rho,\beta\rangle}-q^{-\langle\rho,\beta\rangle}}$$

(the notation is as in [**Saw06**, Section 1]) the so-called ***quantum Weyl dimension formula***. Here $q=q(k)$ is a certain complex root of unity.

The only thing we need to know about this formula is that

$$\lim_{k\to\infty}\mathrm{PFdim}_k\overline{L}(\lambda,k)=\prod_{\beta>0}\frac{\langle\lambda+\rho,\beta\rangle}{\langle\rho,\beta\rangle},$$

which is now the ***classical*** Weyl dimension formula. This is well-known, since $k\to\infty$ corresponds to $q\to 1$.

Let us choose the cutoff such that the vertices of $\Gamma_k$ can be matched with $\Lambda^k$. Using this choice the induced product matches the product in the Verlinde category, by the fusion rules as recalled in [**Saw06**, Section 5]. (It is remarkable that the Verlinde product is a truncation.) Hence, this shows that Definition 12E.15.(b) holds.

Moreover, since the Verlinde categories are braided, it follows that the action matrices $M(c,k)$ of any $\mathbb{R}_{\geq 0}$-linear combination of the $[\overline{L}(\lambda,k)]$ commute with one another, and these matrices can be simultaneously diagonalized. It is easily checked that the vector

$$\sum_{\lambda\in\Lambda^k}\mathrm{PFdim}_k\overline{L}(\lambda,k)\cdot[\overline{L}(\lambda,k)]$$

is an eigenvector of $M(c,k)$. The classical PF theorem implies that the associated eigenvalue of $M(c,k)$ is $\mathrm{PFdim}_kM(c,k)$. It follows that we have subpowermultiplicativity, or more precisely even:

$$\mathrm{PFdim}_k\big(M(c,k)M(c',k)\big)=\mathrm{PFdim}_k\big(M(c,k)\big)\mathrm{PFdim}_k\big(M(c',k)\big).$$

and superadditivity, ore more precisely even:

$$\mathrm{PFdim}_k\big(M(c,k)+M(c',k)\big)=\mathrm{PFdim}_k\big(M(c,k)\big)+\mathrm{PFdim}_k\big(M(c',k)\big).$$

This then implies Definition 12E.15.(a). □

**Example 12E.26.** The argument in the proof of Proposition 12E.25 also works, for example, for the nonsimple group $\mathrm{GL}_n(\mathbb{C})$ and its defining representation. However, the setting we use needs to be adjusted since the induced product is not the product in the Verlinde cutoff.

Explicitly, take $n=2$. Let us denote the vector representation of $\mathrm{GL}_2(\mathbb{C})$ by $L(1,0)$. Then there exists simple $\mathrm{GL}_2(\mathbb{C})$-representations $L(a,b)$ with $a,b\in\mathbb{Z}_{\geq 0},a\geq b$ such that $L(1,0)\otimes L(a,b)\cong L(a+1,b)\oplus L(a,b+1)$, with the convention that $L(a,b)$ is zero if $a,b\in\mathbb{Z}_{\geq 0},a\geq b$ is not satisfied. The fusion graph is then of the form (with the vertices on a grid using $(a,b)$ as coordinates):

$\Gamma(\mathrm{GL}_2) =$ 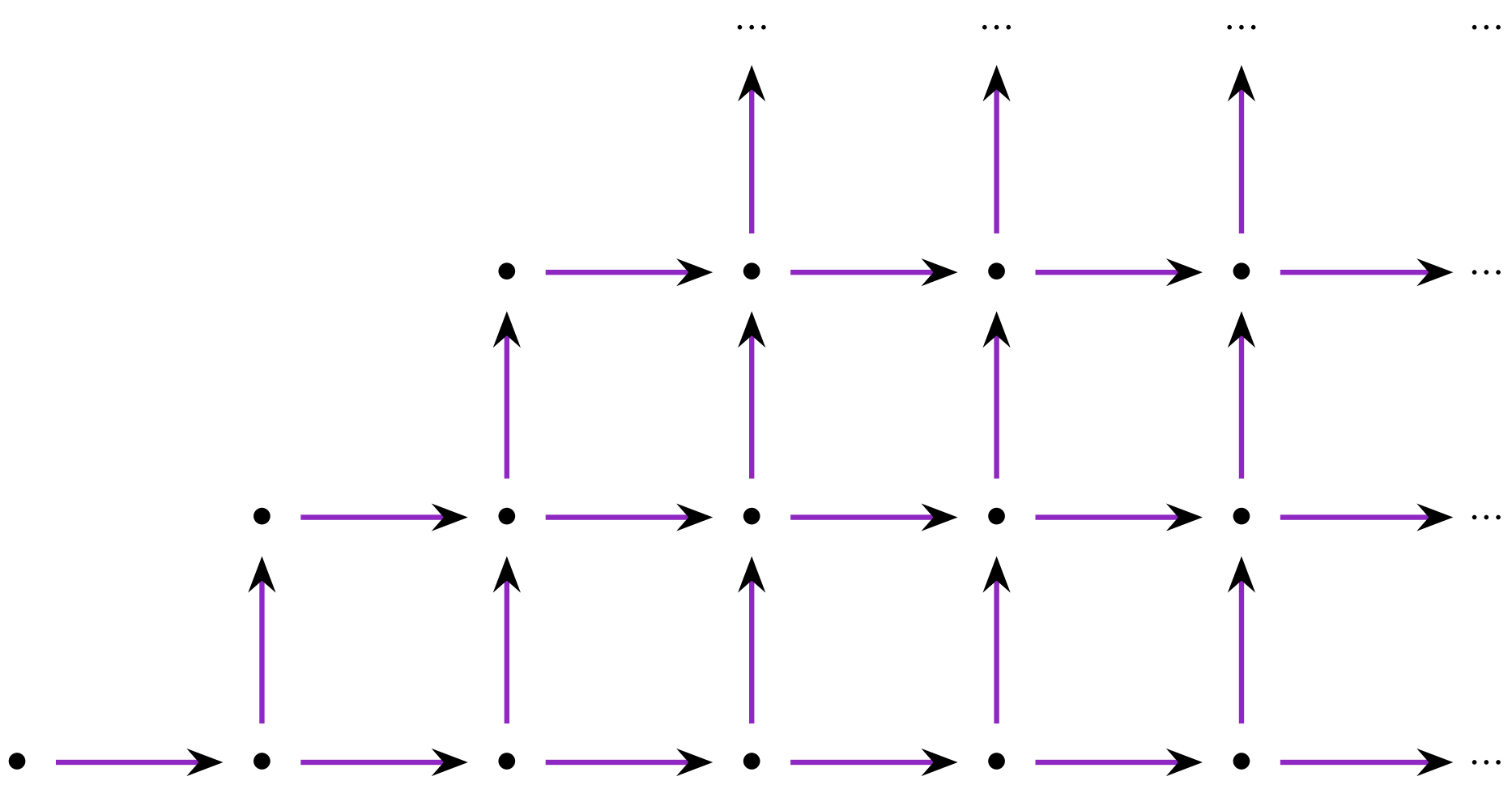 ,

$$M(\mathrm{GL}_2) = \begin{pmatrix} 0 & 0 & 0 & 0 & 0 & \dots \\ J_1 & 0 & 0 & 0 & 0 & \dots \\ 0 & J'_2 & 0 & 0 & 0 & \dots \\ 0 & 0 & J_2 & 0 & 0 & \dots \\ 0 & 0 & 0 & J'_3 & 0 & \dots \\ \vdots & \vdots & \vdots & \vdots & \vdots & \vdots \end{pmatrix},$$

$$\text{with } \begin{cases} J_k = & \text{(lower triangular) Jordan block of size } k \text{ with eigenvalue } 1, \\ J'_k = J_k & \text{without the last column.} \end{cases}$$

The naive filtration therefore gives $\mathrm{PFdim}^f\Gamma(\mathrm{GL}_2) = 0$ since all cutoffs are described by nilpotent matrices, and Definition 12E.15.(b) is not satisfied. In contrast, the Verlinde filtration gives $\mathrm{PFdim}^f\Gamma(\mathrm{GL}_2) = 2$ (as the corresponding graph has backward arrows). ◇

*Remark* 12E.27. Proposition 12E.25 and its proof also generalize [**LSYZ22**, Theorem 1.2]. ◇

**Definition 12E.28.** We say an $\mathbb{R}_{\geq 0}$-algebra is of ***superexponential growth*** if there exists some $c$ such that $b_n$ grows superexponentially.

We say an additive Krull–Schmidt monoidal category is of ***superexponential growth*** if its Grothendieck ring is of superexponential growth in the above sense. ◇

*Remark* 12E.29. Definition 12E.28 runs in parallel to the definition of superexponential growth for abelian monoidal categories using the length of objects. However, note that:

**(a)** An abelian growth problem could be of superexponential growth while its version defined using $b_n$ might not grow superexponentially. We, however, do not know any example where this happens.

**(b)** If $b_n$ grows superexponentially, then the abelian one does so as well.

Hence, the notion from Definition 12E.28 “grows (potentially) slower” than the abelian one. ◇

The next example is simple but crucial to determine categories of superexponential growth.

**Example 12E.30.** Consider the following graph:

$\Gamma(\star_N) =$ 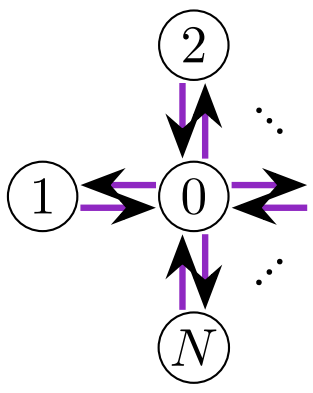

$\cdots$ , $\quad M(\star_N) = \begin{pmatrix} 0 & 1 & 1 & \dots & 1 \\ 1 & 0 & 0 & \dots & 0 \\ 1 & 0 & 0 & \dots & 0 \\ \vdots & \vdots & \vdots & \ddots & \vdots \\ 1 & 0 & 0 & \dots & 0 \end{pmatrix}.$

This is called the ***star graph*** with $N+1$ vertices. It is easy to see that $\mathrm{PFdim}\,\Gamma(\star_N) = \sqrt{N}$.

Let $\Gamma(Y)$ be the Young lattice, considered as an oriented graph by putting an orientation in both directions. The first few layers of this graph are:

$\Gamma(Y) \leftrightsquigarrow$ 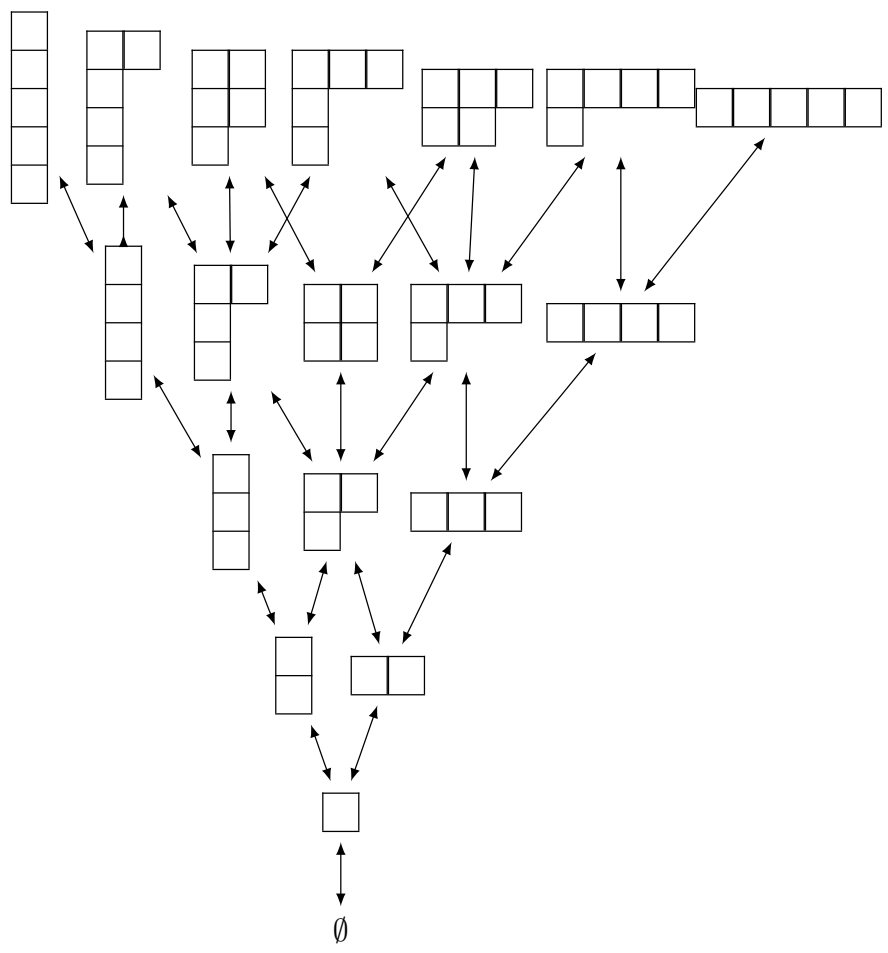 .

We get $\mathrm{PFdim}^f\Gamma(Y) = \infty$ since $\Gamma(Y)$ contains star graphs for arbitrarily large $N$. To see this take the vertex 0 of the star graph to be a staircase partition $(t, t-1, t-2, ..., 1)$ for some large enough $t$. ◇

**Proposition 12E.31.** *The following categories are of superexponential growth.*

- **(a)** *The category* $\mathbf{Rep}\big(\mathrm{SL}_{\mathbb{Z}_{\geq 0}}(\mathbb{C}), \mathbb{C}^{\mathbb{Z}_{\geq 0}}\big) = \mathbf{fdMod}\big(\mathrm{SL}_{\mathbb{Z}_{\geq 0}}(\mathbb{C})\big)$.
- **(b)** *The Deligne category* $\mathbf{Rep}(S_t, \mathbb{C})$*, where* $t \in \mathbb{C} \setminus \mathbb{Z}$.
- **(c)** *The so-called* **Delannoy categories** *from [**HSS22**].*

*Proof.* We can apply Example 12E.30 three times: For $\mathbf{Rep}\big(\mathrm{SL}_{\mathbb{Z}_{\geq 0}}(\mathbb{C}), \mathbb{C}^{\mathbb{Z}_{\geq 0}}\big)$ this follows directly from Theorem 12E.20.(a) and Example 12E.30.

The following argument is based on the description of the fusion graph of $\mathbf{Rep}(S_t, \mathbb{C})$ for $t \in \mathbb{C} \setminus \mathbb{Z}$, which was provided to us by Victor Ostrik via email. We thank Victor for providing us with this description. For $\mathbf{Rep}(S_t, \mathbb{C})$ we recall that the combinatorics of tensoring in is very similar to the combinatorics of tensoring with $\star_N$-representations if $N \gg 0$. The correspondence is as follows: if you have a partition $(p_1, p_2, ...)$ labeling object in the Deligne category, it corresponds to a partition $(N-(p_1+p_2+...), p_1, p_2, ...)$ for $\star_N$. Now tensoring by a 1-node Young diagram object is given by removing one node whenever possible and add it back. This description implies the existence of star graphs, so we are done in the same way as for $\mathbf{Rep}\big(\mathrm{SL}_{\mathbb{Z}_{\geq 0}}(\mathbb{C}), \mathbb{C}^{\mathbb{Z}_{\geq 0}}\big)$.

For the Delannoy categories from [**HSS22**] we can again use star graphs for the fusion graphs computed in [**HSS22**, Section 7]. □

*Remark* 12E.32. By [**GT25b**], "most" of the categories discussed in [**KOK22a**] exhibit superexponential growth, and we suspect that the same is true for the generalizations of the Delannoy categories from [**HSS22**]. ◇

12F. **Recurrent and transient growth problems.** We retain the setup from above. Following the classical theory of random walks, see *e.g.* [**Kit98**] where the reader will also find some standard terminology, we now define ***recurrent and transient growth problems***.

To this end, recall that a vertex $v \in \Gamma$ in a graph with countably many vertices (not just a fusion graph) is ***recurrent*** if the probability of returning to $v$ in the random walk with all edge labels one on $\Gamma$ is 1, and ***transient*** if its not recurrent. If $v \in \Gamma$ is recurrent, then we call it ***positive recurrent*** if the expected amount of time between recurrences is finite, and ***null recurrent*** otherwise.

*Remark* 12F.1. Technically speaking, we are slightly abusing the terminology here. What we refer to as a 'random walk' does not align with the standard definition, since we often assign a weight of one to each edge. In a true random walk, each edge would be given a label such that the sum of all edge weights equals one. ◇

**Lemma 12F.2.** *For a fusion graph* $\Gamma$ *we have:*

- **(a)** *Every vertex* $v$ *of* $\Gamma$ *is either recurrent or transient.*
- **(b)** *Every recurrent vertex of* $\Gamma$ *is either positively recurrent or null recurrent.*
- **(c)** *Let* $C(\Gamma)$ *be a connected component. If* $v \in C(\Gamma)$ *is positively recurrent (or null recurrent or transient), then so is any other* $w \in C(\Gamma)$.

*Proof.* (a) and (b) are true by definition (and only stressed because of the many alternative definitions of recurrent and transient, which are not immediately opposite of each other), and (c) is classical. □

Note that Lemma 12F.2 allows us to consider (positively/null) recurrence and transience for one fixed connected component. We will use this below.

To state the main definition, we recall the notion of PF dimension of irreducible (not necessary finite) $\mathbb{R}_{\geq 0}$-matrices à la Vere-Jones [**VJ67**] (justifying the subscript in the following notation). For such a matrix $M = (m_{ij})_{i,j\in I}$ for a countable set $I$, let

$$\mathrm{PFdim}^{VJ} M = \lim_{n\to\infty} \sqrt[h\cdot n]{m_{ij}^{(h\cdot n)}}, \tag{12F-3}$$

where $m_{ij}^{(n)}$ is the $(i,j)$-entry of $M^n$, and $h \in \mathbb{Z}_{\geq 1}$ is the period. By [**VJ67**, Theorem A] $\mathrm{PFdim}^{VJ} M$ is independent of $i$ and $j$.

Recalling that irreducible matrices correspond to strongly connected graphs, we can define:

**Definition 12F.4.** For a strongly connected graph $\Gamma$ with countably many vertices define the ***PF dimension*** $\mathrm{PFdim}^{VJ}\Gamma$ using Equation 12F-3 with $M$ being the adjacency matrix of $\Gamma$. ◇

*Remark* 12F.5. A generalization of $\mathrm{PFdim}^{VJ}\Gamma$ to graphs that are not necessarily strongly connected was worked out in [**Twe71**]. It however turns out that this generalization is not useful for our purposes. ◇

The following is immediate:

**Lemma 12F.6.** *If $|I| < \infty$, then Definition 12F.4 agrees with the usual definition of PF dimension using the largest eigenvalue of the associated matrix (or, equivalently, graph).* □

*Remark* 12F.7. We use the notation $\mathrm{PFdim}^{VJ}$ to indicate that this PF dimension is different from the one we introduced in ??. They however agree if $|I| < \infty$ by Lemma 12E.14 and Lemma 12F.6. ◇

**Example 12F.8.** We continue with Example 12E.6. In this case we get $\mathrm{PFdim}^{VJ}\Gamma(\mathrm{SL}_2) = 2$ since

$$m_{11}^{(2n)}(\mathrm{SL}_2) = \frac{1}{n+1}\binom{2n}{n} \sim \pi^{-1/2} \cdot n^{-3/2} \cdot 2^{2n}.$$

This equation follows from a standard argument. Note that the period is two, which is why we consider every second value of $m_{11}^{(n)}(\mathrm{SL}_2)$ only. ◇

We also recall the notion of ***final basic classes (FBC)***. Firstly, a ***class*** $C(\Gamma)$ is a strongly connected component of $\Gamma$ and such a class is ***basic*** if:

$$\mathrm{PFdim}^{VJ} C(\Gamma) \geq \mathrm{PFdim}^{VJ} C'(\Gamma) \text{ for all strongly connected components } C'(\Gamma).$$

Finally, a basic class is ***final*** if there is no path to any other basic class. (Note that a FBC does not need to be final in $\Gamma$ itself.)

**Example 12F.9.** For a strongly connected graph $\Gamma$ is the only class, so it is a FBC. In particular, what we will see below is a generalization of the strongly connected situation. ◇

**Example 12F.10.** Consider $G = \mathrm{SL}_2(\bar{\mathbb{F}}_2)$ and take the defining $G$-representation $\mathtt{X} = \bar{\mathbb{F}}_2^2$ over the field $\bar{\mathbb{F}}_2$. The growth problem $(G, \mathtt{X})$ has the associated fusion graph (we illustrate a cutoff):

$\Gamma(\mathrm{SL}_2(\bar{\mathbb{F}}_2)) =$ 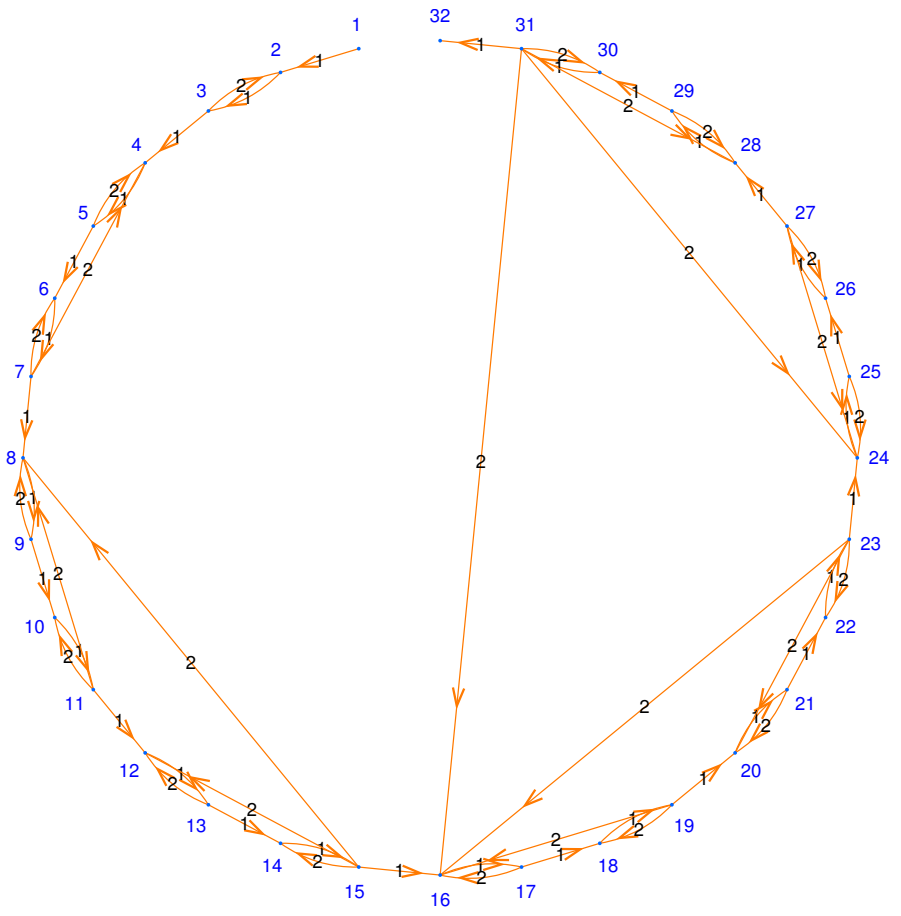 .

The pattern is as follows: The vertices are labeled with $\mathbb{Z}_{\geq 1}$. There is always one forward arrow, from vertex $i$ to vertex $i+1$. The backward arrows are of weight 2 and start at $2^k - 1$ and end at $2^{k-2}$ for all $k \in \mathbb{Z}_{\geq 2}$. See, for example, [**STWZ23**] for details.

The strongly connected components are indexed by $2^k$ for $k \in \mathbb{Z}_{\geq 0}$, with $2^k$ being their minimal vertex. In this case there is no FBC since all strongly connected components have PF dimension strictly smaller than two, but the limit $k \to \infty$ of their PF dimensions is approaching 2.

In a bit more details, one can check (by calculation) that the PF dimension of the cutoffs will be

$$0, \sqrt{2}, \sqrt{2+\sqrt{2}}, \sqrt{2+\sqrt{2+\sqrt{2}}} \text{ etc.}$$

where the jumps happen when the number of vertices increases from $2^k$ to $2^k+1$ in the sequence of naive cutoffs. The limit approaches 2, by classical results about sequences of nested radicals. ◇

The following is a key definition:

**Definition 12F.11.** Let $(R, c)$ be a growth problem, and let FBC refer to its fusion graph $\Gamma$.

**(a)** The growth problem $(R, c)$ is ***recurrent*** if all of its FBCs are recurrent and there exists at least one FBC. The growth problem is ***transient*** otherwise.

**(b)** If $(R, c)$ is recurrent, then we say $(R, c)$ is ***positively recurrent*** if all of its FBCs are positively recurrent, The growth problem is ***null recurrent*** otherwise.

We use the same terminology for graphs in general. ◇

Note that a transient FBC dominates recurrent ones in the sense that the presence of even a single transient FBC renders the growth problem transient, regardless of how many FBCs are recurrent. Similarly, a null recurrent FBC dominates positively recurrent ones.

**Lemma 12F.12.** *Every growth problem is either recurrent or transient. Moreover, every recurrent growth problem is either positively or null recurrent.*

*Proof.* Immediate by definition (but stated for the same reasons as Lemma 12F.2). □

**Example 12F.13.** Take the graph

$\Gamma(\mathbb{Z}) =$ 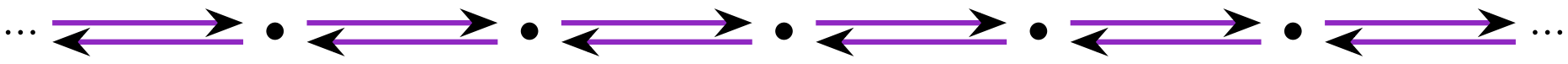 .

The associated random walk is the classical random walk on $\mathbb{Z}$. This is Polya's first example of a recurrent graph [**Pól21**]. Additionally, this example is null recurrent.

In contrast, the random walk on $\Gamma(\mathrm{SL}_2)$ from Equation 12E-7 is transient. Below we will see that this is no coincidence, *cf.* Proposition 12F.23.

The difference becomes evident when one looks at the number of path of length $n$ starting at the origin, and ending at vertex $v$. For $n = 200$, plotting this in $\mathbb{R}^2$ with (end vertex, number of paths) gives:

$\mathbb{Z}$: 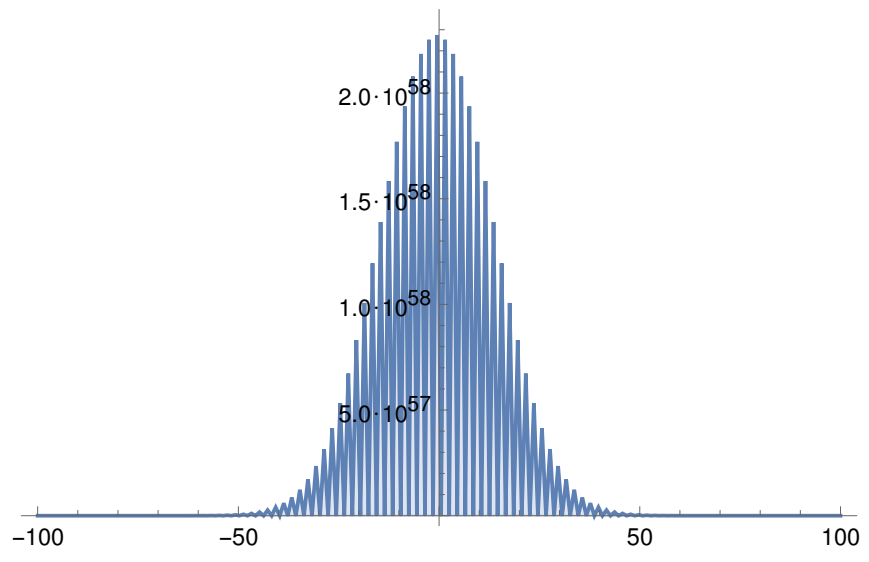

, $\mathbb{Z}_{\geq 0}$: 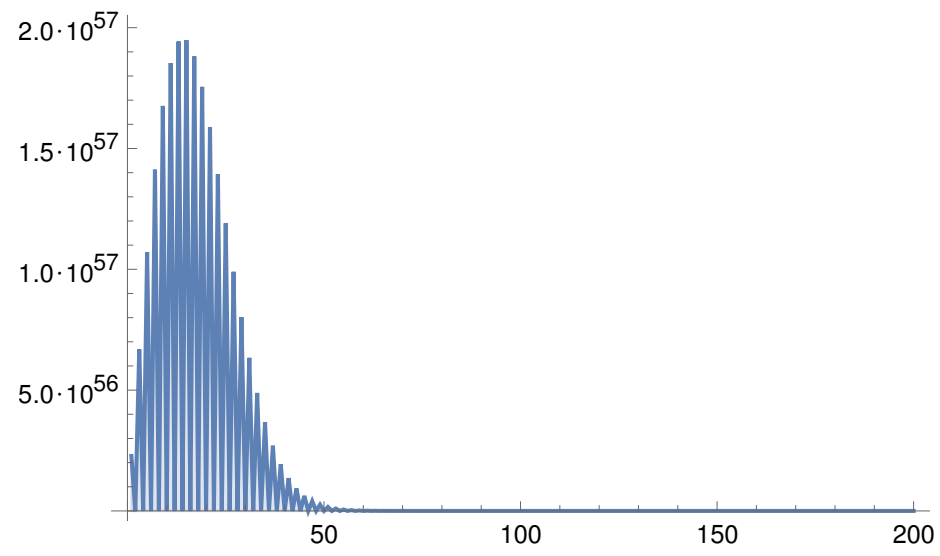

.

Note that for $\Gamma(\mathrm{SL}_2) \cong \mathbb{Z}_{\geq 0}$ the peak of the binomial distribution is roughly at $\sqrt{200} \approx 14.1$. In fact, the endpoints of paths in this case move to infinity, while they stay at the origin for $\mathbb{Z}$.

For completeness: The graph $\Gamma(\mathbb{Z})$ is also the fusion graph for a monoidal category: take $(\mathbf{Rep}(\mathbb{S}^1), \mathtt{X}_1 \oplus \mathtt{X}_{-1})$ with $\mathtt{X}_{\pm 1}$ being the rotation of an irrational angle $\pm\theta$. ◇

**Example 12F.14.** Similarly as for $\mathrm{SL}_2(\mathbb{C})$, the growth problem for $G = \mathrm{SL}_3(\mathbb{C})$ and $\mathtt{X} = \mathbb{C}^3$ with the defining $G$-action is also transient. In this case, the fusion graph is (here a cutoff):

$\Gamma(\mathrm{SL}_3) \rightsquigarrow$ 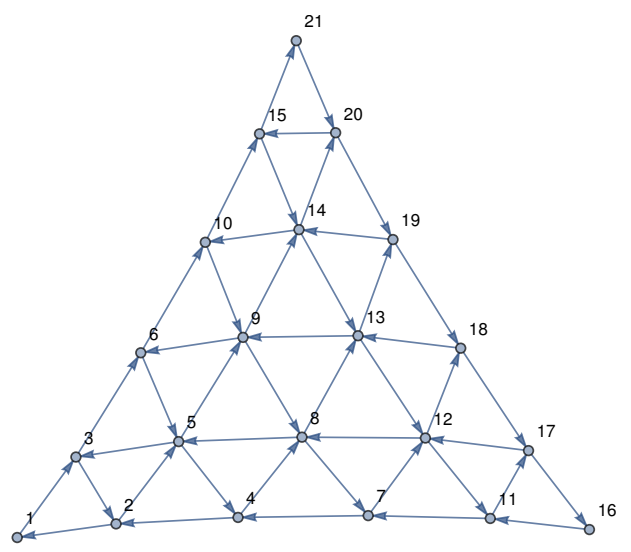

.

And similarly as in Example 12F.13, the number of paths of length $n$ starting at 1 and ending at $v$ moves outwards (illustrated using a large enough cutoff):

$n = 50$: 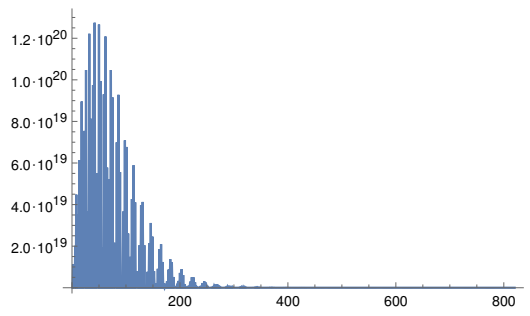

, $n = 100$: 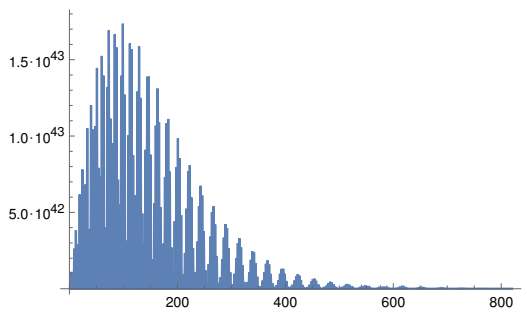

, $n = 200$: 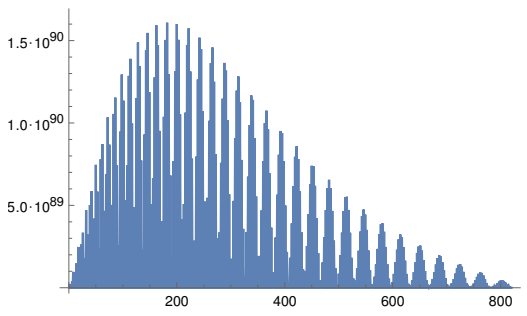

.

In Proposition 12F.23 we will see that this is not a coincidence. ◇

**Example 12F.15.** Every finite growth problem is positively recurrent. More generally, every growth problem where the FBCs are finite is also positively recurrent. Explicitly, let $G = \mathbb{Z}/2\mathbb{Z} \times \mathbb{Z}/2\mathbb{Z}$ be the ***Klein four group***, and we take the ground field $\mathbb{k} = \mathbb{F}_2$. Then one has $\mathbb{k}[G] \cong \mathbb{k}[X, Y]/(X^2, Y^2)$, and using this presentation one can check that

$Z_3 \rightsquigarrow$ 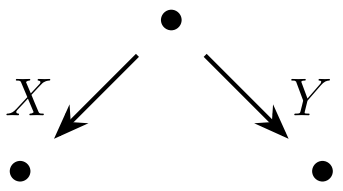

defines an indecomposable $G$-representation of dimension three, call it $V = Z_3$. Moreover, it is also easy to see that $G$ has indecomposable $G$-representations $Z_k$, for odd $k \in \mathbb{Z}_{\geq 1}$, where $Z_k$ is the unique indecomposable summand of $V^{\otimes (k-1)/2}$ that is not projective. Note that $Z_1 \cong L_1$ is the trivial $G$-representation. One can check that $\dim_{\mathbb{k}} Z_k = k$ and that all odd numbers appear.

Let $P_4$ denote the regular $G$-representation. An easy calculation gives us then the following fusion graph.

(12F-16) $\Gamma(Kl) =$ 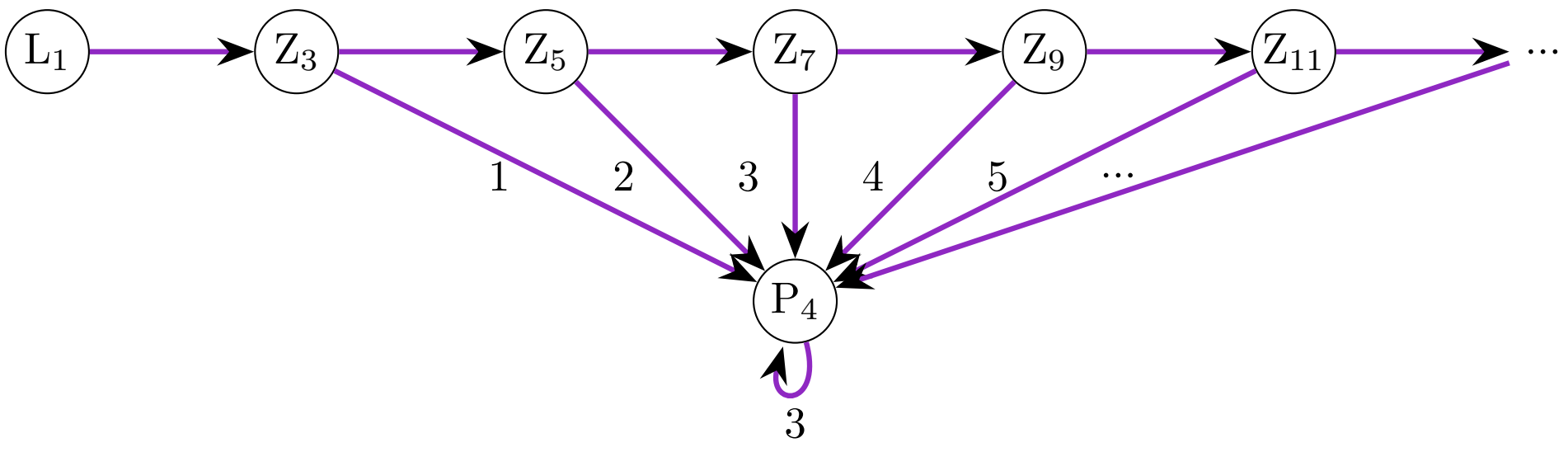

,

$$M(Kl) = \begin{pmatrix} 3 & 0 & 1 & 2 & 3 & 4 & 5 & \dots \\ 0 & 0 & 0 & 0 & 0 & 0 & 0 & \dots \\ 0 & 1 & 0 & 0 & 0 & 0 & 0 & \dots \\ 0 & 0 & 1 & 0 & 0 & 0 & 0 & \dots \\ 0 & 0 & 0 & 1 & 0 & 0 & 0 & \dots \\ 0 & 0 & 0 & 0 & 1 & 0 & 0 & \dots \\ 0 & 0 & 0 & 0 & 0 & 1 & 0 & \dots \\ \vdots & \vdots & \vdots & \vdots & \vdots & \vdots & \vdots & \ddots \end{pmatrix}.$$

The labeled loop at the vertex with label 3 represents three loops.

Let us consider the growth rate problem for $(G, V)$. Then this is positively recurrent since the only FBC is $\{P_4\}$ and therefore a finite graph. ◇

For a more general statement see Proposition 12F.23 below.

**Example 12F.17.** The growth problem in Example 12F.10 is transient since there is no FBC. ◇

Recall that a group $G$ is ***virtually*** $H$ if $H$ is isomorphic to a finite index subgroup of $G$.

**Theorem 12F.18** (Polya's random walk theorem for representations)**.** *Let $\mathbb{k}$ be a field of characteristic zero. For any (finitely presentable) group $G$ and finite dimensional faithful completely reducible $G$-representation $\mathtt{X}$ the following are equivalent:*

**(a)** *The growth problem $(G, \mathtt{X})$ is recurrent.*

**(b)** *The Zariski closure of the image of $G$ in $\mathrm{GL}(\mathtt{X})$ is virtually a torus of rank 0, 1 or 2.*

*Moreover, the following are also equivalent:*

*(i)* *The growth problem $(G, \mathtt{X})$ is positively recurrent.*

*(ii)* *The Zariski closure of the image of $G$ in $\mathrm{GL}(\mathtt{X})$ is virtually a torus of rank 0, i.e. $G$ is finite.*

*If $\mathtt{X}$ is not faithful, then go to the biggest quotient group that acts faithfully and repeat.*

*Proof.* We will make use of the following. We will formulated it for bialgebra, but we use the same notation as for group.

**Lemma 12F.19.** *Let $G$ be a bialgebra, and $H \subset G$ be a subbialgebra. Let $\mathtt{X}$ be a $G$-representation such that:*

**(a)** *The growth problem $(G, \mathtt{X})$ has a strongly connected fusion graph with $\mathrm{PFdim}^{VJ}\mathtt{X} < \infty$.*

**(b)** *The growth problem* $(H, \mathtt{X})$ *obtained via restriction has a strongly connected fusion graph with the same PF dimension* $\mathrm{PFdim}^{VJ}\mathtt{X}$.

*Let* $\mathtt{X} \cong \bigoplus_{i=1}^{j} \mathtt{X}_i$ *as an* $H$*-representation obtained by restriction. If the growth problem* $(H, \mathtt{X}_i)$ *is transient for some* $i \in \{1, ..., j\}$*, then* $(G, \mathtt{X})$ *is also transient.*

*Proof.* We only prove the statement in the aperiodic case, to keep the formulas simple. First note that the transience of $(H, \mathtt{X}_i)$ implies that $(H, \mathtt{X})$ is also transient as there is a subgraph corresponding to $(H, \mathtt{X}_i)$. Now, since $(H, \mathtt{X})$ is transient with strongly connected fusion graph we have

$$\sum_{n \in \mathbb{Z}_{\geq 0}} h_{11}^{(n)} \cdot (\mathrm{PFdim}^{VJ}\mathtt{X})^{-n} < \infty,$$

where $h_{11}^{(n)}$ is the $(1,1)$-entry of the $n$th power of the action matrix corresponding to the fusion graph of $(H, \mathtt{X})$. Using a similar notation for $(G, \mathtt{X})$, we also have

$$\sum_{n \in \mathbb{Z}_{\geq 0}} g_{11}^{(n)} \cdot (\mathrm{PFdim}^{VJ}\mathtt{X})^{-n} \leq \sum_{n \in \mathbb{Z}_{\geq 0}} h_{11}^{(n)} \cdot (\mathrm{PFdim}^{VJ}\mathtt{X})^{-n},$$

by assumption. Taking both together, the claim follows. □

The following is then key:

**Lemma 12F.20.** *Let* $G = \mathrm{SL}_2(\mathbb{C})$ *and* $\mathtt{X}$ *any simple* $G$*-representation of dimension* $\dim_{\mathbb{C}} \mathtt{X} > 1$*. Then the growth problem* $(G, \mathtt{X})$ *is transient. The same is true when replacing* $\mathbb{C}$ *by any characteristic zero field.*

*Proof.* Recall that the fusion rules for $G = \mathrm{SL}_2(\mathbb{C})$ are as follows. We index the simple $G$-representations $L(\lambda)$ by their highest weight $\lambda \in \mathbb{Z}_{\geq 0}$. Then, using the convention that $L(\nu) = 0$ if $\nu \notin \mathbb{Z}_{\geq 0}$, we have

$$L(\lambda) \otimes L(\mu) \cong L(\lambda - \mu) \oplus L(\lambda - \mu + 2) \oplus ... \oplus L(\lambda + \mu - 2) \oplus L(\lambda + \mu)$$

for $\lambda \geq \mu$ and similarly for $\lambda < \mu$.

The fusion graphs are now easily computed. We give a few examples of how their action matrices look like:

$$M(\lambda = 1) = \begin{pmatrix} 0 & 1 & 0 & 0 & 0 & \dots \\ 1 & 0 & 1 & 0 & 0 & \dots \\ 0 & 1 & 0 & 1 & 0 & \dots \\ 0 & 0 & 1 & 0 & 1 & \dots \\ 0 & 0 & 0 & 1 & 0 & \dots \\ \vdots & \vdots & \vdots & \vdots & \vdots & \vdots \end{pmatrix}, \quad M(\lambda = 3) = \begin{pmatrix} 0 & 1 & 0 & 1 & 0 & \dots \\ 1 & 0 & 1 & 0 & 1 & \dots \\ 0 & 1 & 0 & 1 & 0 & \dots \\ 1 & 0 & 1 & 0 & 1 & \dots \\ 0 & 1 & 0 & 1 & 0 & \dots \\ \vdots & \vdots & \vdots & \vdots & \vdots & \vdots \end{pmatrix}, \tag{12F-21}$$

and the number of offdiagonals with 1s keep on increasing with steps of two for $\lambda = 5$, $\lambda = 7$ etc. The case where $\lambda$ is even is similar, but with 1s on the main diagonal except in the top left entry which is still zero. For example:

$$M(\lambda = 2) = \begin{pmatrix} 0 & 1 & 0 & 0 & 0 & \dots \\ 1 & 1 & 1 & 0 & 0 & \dots \\ 0 & 1 & 1 & 1 & 0 & \dots \\ 0 & 0 & 1 & 1 & 1 & \dots \\ 0 & 0 & 0 & 1 & 1 & \dots \\ \vdots & \vdots & \vdots & \vdots & \vdots & \vdots \end{pmatrix}, \quad M(\lambda = 4) = \begin{pmatrix} 0 & 1 & 1 & 0 & 0 & \dots \\ 1 & 1 & 1 & 1 & 0 & \dots \\ 1 & 1 & 1 & 1 & 1 & \dots \\ 0 & 1 & 1 & 1 & 1 & \dots \\ 0 & 0 & 1 & 1 & 1 & \dots \\ \vdots & \vdots & \vdots & \vdots & \vdots & \vdots \end{pmatrix}. \tag{12F-22}$$

Note that the matrices for $\lambda \in 2\mathbb{Z}_{\geq 0}$ are on a different basis on the Grothendieck level (namely on the Grothendieck classes of $L(0)$, $L(2)$, $L(4)$, etc.)

Let us now prove the case where $\lambda = 2$ as an example. In this case $m_{11}^{n}$ is given by A005043 on [**OEI23**]. In particular, $m_{11}^{(n)} \sim \frac{3^{3/2}}{8\sqrt{\pi}} \cdot n^{-3/2} \cdot 3^n$ and this implies that

$$\sum_{n \in \mathbb{Z}_{\geq 0}} m_{11}^{(n)} \cdot 3^{-n} \approx \frac{3^{3/2}}{8\sqrt{\pi}} \sum_{n \in \mathbb{Z}_{\geq 0}} n^{-3/2} < \infty,$$

which is the usual calculation one need to do in order to verify that a problem is transient. (For the reader unfamiliar with random walks, this is one of many possible, and equivalent, definitions of transient, see e.g. [**Kit98**, Chapter 7].)

For general even $\lambda > 0$ we have $m_{11}^{(n)} \sim \text{some scalar} \cdot n^{-3/2} \cdot \lambda^n$, as follows from e.g. [**Bia93**, Theorem 2.2]. So the same calculations as for $\lambda = 2$ works. Moreover, again by e.g. [**Bia93**, Theorem 2.2], when $\lambda$ is odd, then $m_{11}^{(2n)} \sim \text{some scalar} \cdot n^{-3/2} \cdot \lambda^{2n}$, and again the same calculation works. □

We can assume that $G$ is an algebraic group since we can replace it with the Zariski closure of its image in $\mathrm{GL}(\mathtt{X})$. So assume that $G$ is an algebraic group, whose identity component is a reductive group since $\mathtt{X}$ is completely reducible. We can, and will, even replace $G$ by its identity component.

If $G$ is not a torus, then we can apply the Jacobson–Morozov theorem to find a copy of $\mathrm{SL}_2(\mathbb{C})$ in $G$. Then Lemma 12F.20 (see also Example 12F.13) and Lemma 12F.19, whose assumptions are satisfied by [**COT24**, Theorem 1.4], imply that the growth problem $(G, \mathtt{X})$ is transient.

On the other hand, if $G$ is torus, then Polya's classification of recurrent random walks applies, and we get a recurrent growth problem if and only if the rank of the torus is 0, 1 or 2.

The claim about positively recurrent versus null recurrent follows then from Example 12F.13. □

**Proposition 12F.23.** *The following growth problems are positively recurrent.*

- **(a)** *For an arbitrary field, $(G, \mathtt{X})$ for a finite group $G$ and $\mathtt{X}$ any $G$-representation.*
- **(b)** *For an arbitrary field, $(\mathbf{C}, \mathtt{X})$ for a finite tensor category $\mathbf{C}$ and $\mathtt{X} \in \mathbf{C}$ any object.*
- **(c)** *For an arbitrary field, $(\mathbf{SBim}, \mathtt{X})$ for Soergel bimodules $\mathbf{SBim}$ attached to a finite Coxeter group and $\mathtt{X} \in \mathbf{SBim}$ any object.*

*Moreover, the following growth problems are transient.*

- **(i)** *For a field of characteristic zero, $(G, \mathtt{X})$ for a group $G$ and $\mathtt{X}$ any $G$-representation such that the Zariski closure of the image of $G$ in $\mathrm{GL}(\mathtt{X})$ is not virtually a torus.*
- **(ii)** *In the below specified cases, $(U_q, \mathtt{X})$ for a quantum enveloping algebra $U_q = U_q(\mathfrak{g})$ (in the sense of [**Lus90**] or [**APW91**]) and $\mathtt{X}$ any nontrivial tilting $U_q(\mathfrak{g})$-representation (meaning not a direct sum of one dimensional $U_q$-representations). The cases are:*
  - *$\Bbbk$ is arbitrary, and $q$ is not a root of unity.*
  - *$\Bbbk$ is an algebraically closed field of characteristic zero, and $q$ is a root of unity.*
- **(iii)** *For a field of characteristic zero, $(\mathbf{SBim}, \mathtt{X})$ for Soergel bimodules $\mathbf{SBim}$ attached to an affine Weyl group and $\mathtt{X} \in \mathbf{SBim}$ a generating object.*

*Proof.* (a). A special case of (b).

(b). We can and will assume that $\mathtt{X}$ is a generating object: if that is not the case then we would go to a smaller tensor category. By [**CEOT24**, Proposition A.1], every projective indecomposable object of $\mathbf{C}$ is a direct summand of $\mathtt{X}^{\otimes d}$ for some $d \in \mathbb{Z}_{\geq 0}$.

Let us first analyze the subgraph $\Gamma^p$ of the fusion graph $\Gamma$ that contains only vertices for the projective indecomposable objects. We call this the ***projective cell***, and we claim it is strongly connected and basic. To see this we first note that [**CEOT24**, Proposition A.1] also gives us that the regular object appears in some tensor power of $\mathtt{X}$. Then [**EGNO15**, Section 3.3] implies that the projective cell is of maximal PF dimension since going from using Grothendieck classes of simple objects to compute PF dimensions, as in [**EGNO15**], to indecomposable objects can only decrease the PF dimension. That the projective cell is strongly connected also follows from the existence of the regular objects in some tensor power of $\mathtt{X}$.

Moreover, recall that projectives form a $\otimes$-ideal, see e.g. [**EGNO15**, Proposition 4.2.12]. This also implies by [**CEOT24**, Proposition A.1] that no strongly connected component without projective indecomposable objects is final since there is always a path to the projective cell.

Since every finite tensor category has only finitely many indecomposable projective objects, this problem is positively recurrent.

(c). In this case $\Gamma$ is finite.

(i). As in the proof of Theorem 12F.18, the Jacobson–Morozov theorem implies that we will find $\mathrm{SL}_2 \subset \mathrm{GL}(\mathtt{X})$, and Lemma 12F.19 and Lemma 12F.20 then imply the claim.

(ii). If $q$ is not a root of unity, then this problem has the same combinatorics as (i), see for example [**APW91**, Corollary 7.7], so we are done. In the other case there are projective $U_q$-representations and they appear in some tensor power of $\mathtt{X}$:

**Lemma 12F.24.** *Every projective indecomposable tilting $U_q$-representation is a direct summand of $\mathtt{X}^{\otimes d}$ for some $d \in \mathbb{Z}_{\geq 0}$.*

*Proof.* We use two results from [**APW91**, Theorem 9.12]: there are enough projective tilting objects, and every projective tilting object is also injective. With these two facts, the proof of the lemma is, mutatis mutandis, the same as [**CEOT24**, Proof of Proposition A.1]. □

Using Lemma 12F.24, we get the the only final class can be the projective cell, and it remains to argue that this cell is transient. The results in [**And18**, Theorem 3.1 and Remark 2.(2)] imply that the projective cell (called the Steinberg component in [**And18**]) is a copy of the whole category in the semisimple case upon factoring the Steinberg $U_q$-representation. This implies that the projective cell is transient: The majority of paths in the starting category will eventually be in the projective cell, so they move away from the origin. Then the equivalence in [**And18**, Theorem 3.1 and Remark 2.(2)] jumps back to the original problem for which we already know that the majority of paths eventually leave the origin.

(iii). Using the correspondence between the cells from [**Ost97**], very similar as for the quantum group at a complex root of unity. Details are omitted. □

*Remark* 12F.25. One could guess that the analogs of Proposition 12F.23.(i) is also true in positive characteristic, e.g. for $G = \mathrm{SL}_2(\overline{\mathbb{F}}_p)$ this follows from Example 12F.10. Similarly, Proposition 12F.23.(ii) might also be true for quantum groups in the mixed case, and Proposition 12F.23.(iii) might be true for all infinite Coxeter types and regardless of the characteristic of the underlying field. ◇ 12G. **Asymptotics for positively recurrent categories.** Positively recurrent categories are essentially finite with respect to growth problems as justified by the main result of this section. However, we need a special case of positively recurrent growth problems, see Definition 12G.4.

*Convention* 12G.1. In this subsection, all growth problems that we consider are positively recurrent in the sense of Definition 12F.11. ◇

First, let us consider the fusion graph $\Gamma$ associated to some growth problem $(R, C)$. Take the naive cutoff $(\Gamma_k)_{k\in\mathbb{Z}_{\geq 0}}$ and let $\lambda_k = \mathrm{PFdim}\Gamma_k$ (since $\Gamma_k$ is finite, this is the classical PF dimension.) Let $\lambda_k^{sec}$ be any second largest eigenvalue of $\Gamma_k$.

**Lemma 12G.2.** *We have* $\lambda = \lim_{k\to\infty} \lambda_k \in \mathbb{R}_{\geq 0}$ *and* $\lambda^{sec} = \lim_{k\to\infty} \lambda_k^{sec} \in \mathbb{C}$.

*Proof.* Note that $0 \leq \lambda_k \leq \lambda_{k+1}$ so that it remains to argue why $\lambda = \lim_{k\to\infty} \lambda_k \neq \infty$. As usual in the theory, see e.g. [**Kit98**, Section 7], $1/\lambda$ is the radius of convergence of a certain power series, which in turn has a positive radius of convergence on any FBC. Thus, we get $\lambda < \infty$. The second largest eigenvalue exists in $\mathbb{R}_{\geq 0}$ by the same reasoning and $\lambda < \infty$ so $|\lambda^{sec}| < \infty$. □

Using Lemma 12G.2, we will write $\lambda = \lim_{k\to\infty} \lambda_k$ and $\lambda^{sec} = \lim_{k\to\infty} \lambda_k^{sec}$ for a given positively recurrent growth problem.

**Definition 12G.3.** We call $\lambda$ the (positively recurrent) ***PF dimension*** of $(R, c)$, and write $\mathrm{PFdim}^f c = \lambda$. (Here we always use the naive cutoff.) ◇

**Definition 12G.4.** Let $(R, c)$ be a positively recurrent growth problem as defined in Definition 12F.11. We say $(R, c)$ is ***sustainably positively recurrent*** if:

**(a)** There is only one basic class.

**(b)** There exists $\varepsilon > 0$ such that the PF dimensions of the nonbasic classes are smaller than $\lambda - \varepsilon$.

**(c)** We have $(n \mapsto m_{ij}^{(n)}) \in o\big((\lambda + \varepsilon)^n\big)$ for all $\varepsilon \in \mathbb{R}_{>0}$.

We also say positively recurrent and sustainable instead of sustainably positively recurrent. ◇

We will denote the basic class from Definition 12G.4.(a) by $C^{FBC}(\Gamma)$, which necessarily positive recurrent and a final basic class. However, $C^{FBC}(\Gamma)$ is not assumed to be finite or final in $\Gamma$.

**Example 12G.5.** All finite growth problems are sustainably positively recurrent. This includes all positively recurrent growth problems in Theorem 12F.18. ◇

Here are more examples of positively recurrent and sustainable growth problems:

**Proposition 12G.6.** *All the positively recurrent examples in Proposition 12F.23 are sustainable. In all these cases we have* $\mathrm{PFdim}^f \mathtt{X} \in \mathbb{R}_{\geq 1}$, *and* $C^{FBC}(\Gamma)$ *is finite and final in* $\Gamma$.

*Proof.* We have $\mathrm{PFdim}^f \mathtt{X} \in \mathbb{R}_{\geq 1}$ since the fusion graphs have edge weights in $\mathbb{Z}_{\geq 1}$. This follows from the argument in Example 12E.18 together with the observation that every vertex admits a path to the projective cell.

Second, the growth problem in (a) is a special case of the one in (b), while the growth problem in (c) has a finite fusion graph. So it remains to go through in Definition 12G.4 (a)-(c) for case (b), which we do after the following lemma.

**Lemma 12G.7.** *For any growth problem as in Proposition 12F.23.(b) we have*

$$b(n) \sim f(n) \cdot (\mathrm{PFdim}^f \mathtt{X})^n \text{ with } f\colon \mathbb{Z}_{\geq 0} \to (0, 1] \text{ periodic with finite period.}$$

*Proof.* In this case let $l_n = \ell(\mathtt{X}^{\otimes n})$ where $\ell$ is the length. The fusion graph of the growth problem associated to the length is finite and it follows (e.g. by copying [**LTV23**]) that

$$l(n) \sim \tilde{f}(n) \cdot (\mathrm{PFdim}\mathtt{X})^n \text{ with } \tilde{f}\colon \mathbb{Z}_{\geq 0} \to (0, 1] \text{ periodic with finite period,}$$

where $\mathrm{PFdim}\mathtt{X}$ is the usual PF dimension of this finite growth problem. Note that

$$b_n \leq l_n$$

since $b_n$ is defined to count indecomposable summands, while $l_n$ counts simple objects in the Jordan–Hölder filtration. Then [**CEOT24**, Proposition A.1] implies that this $\mathrm{PFdim}^f \mathtt{X}$ is the one from Lemma 12G.2.

As the next step, using the universal grading group of $\mathbf{C}$, one can see that $b(n)/(\mathrm{PFdim}^f \mathtt{X})^n$ has only finitely many limiting points. This implies that the generating function of $b_n$ satisfies the properties necessary to run e.g. [**Mis20**, Section 7.7], and we get $b(hn + s) \sim t_{h,s} \cdot (\mathrm{PFdim}^f \mathtt{X})^{hn}$ for some $t_{h,s} \in (0, 1]$, where $h$ is the number of limiting points and $s \in \{0, ..., h - 1\}$. Collecting the scalars $t_{h,s} \in (0, 1]$ into a piecewise constant function $f\colon \mathbb{Z}_{\geq 0} \to [0, 1)$ shows the claim. □

*Property (a).* As in the proof of Proposition 12F.23, by [**CEOT24**, Proposition A.1] the unique FBC is the projective cell, and its PF dimension is $\mathrm{PFdim}^f \mathtt{X}$. We additionally argue now that all other strongly connected components have PF dimension strictly smaller than $\mathrm{PFdim}^f \mathtt{X}$. Indeed, assuming the contrary leads to a contradiction with Lemma 12G.7. This is easy to see, for example in the toy case

$\Gamma =$ $\mathrm{PFdim}^f \mathtt{X}$ 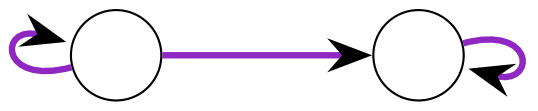 $\mathrm{PFdim}^f \mathtt{X}$ , $\quad M(\Gamma) = \begin{pmatrix} \mathrm{PFdim}^f \mathtt{X} & 0 \\ 1 & \mathrm{PFdim}^f \mathtt{X} \end{pmatrix},$

one already gets $b(n) \sim \frac{n}{\mathrm{PFdim}^f \mathtt{X}} \cdot (\mathrm{PFdim}^f \mathtt{X})^n$. Recalling that all strongly connected components have some path to the projective cell, the general case then reduces to at least this toy case.

*Property (b).* Similarly as in (a), having basic classes $B_n$ with $\lim_{n\to\infty} \mathrm{PFdim}^{VJ} B_n = \mathrm{PFdim}^f \mathtt{X}$ contradicts Lemma 12G.7. For example,

$\Gamma(Kl') =$ 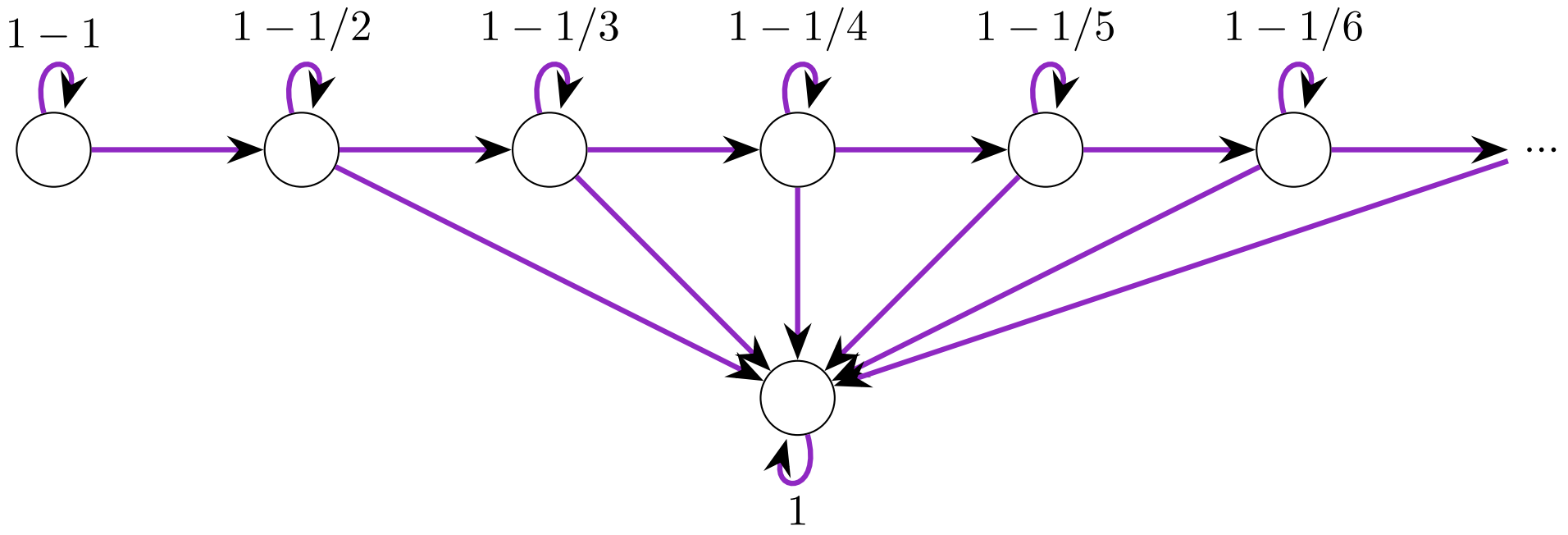

,

has $b(n) \in \Omega(1.5^n)$ while $\mathrm{PFdim}^f \mathtt{X} = 1$.

*Property (c).* As in the points above, nonsustainable would contradict Lemma 12G.7 so we are done. Roughly, consider the following modification of the graph from Equation 12F-16:

$\Gamma(Kl'') =$ 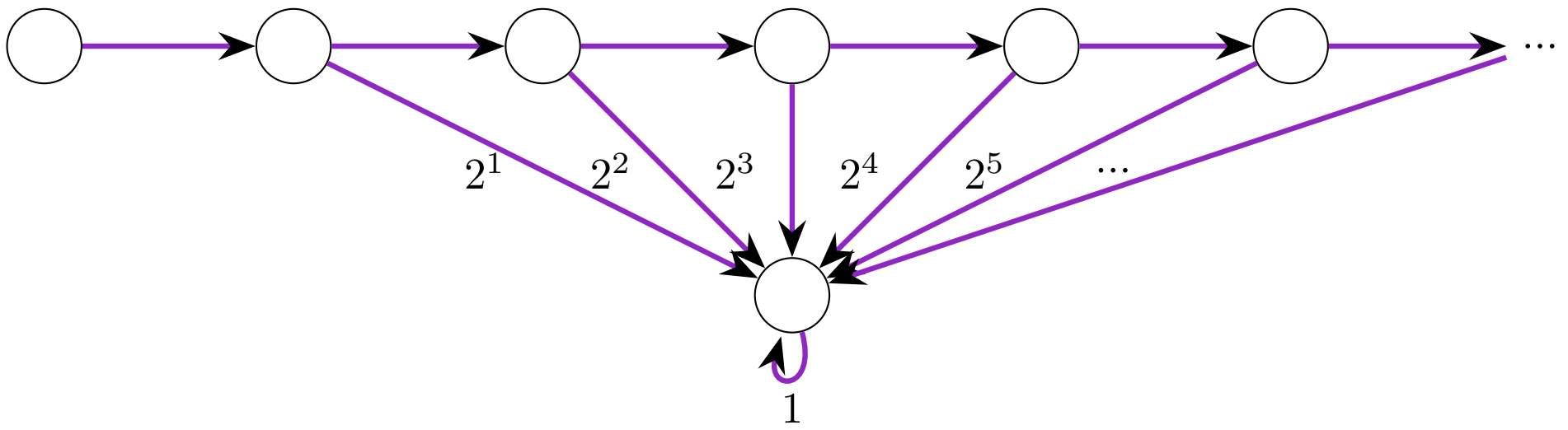

The growth of $b_n$ in this case overshoots $(\mathrm{PFdim}^f \mathtt{X})^n = 1^n$ and we have $b(n) \in \Omega(1.5^n)$.

*Final and finite in* $\Gamma$. We already discussed this above. □

Fix the naive cutoff $\Gamma_k$. Let $h_k$ be the period and $\zeta_k = \exp(2\pi i/h_k)$. Similarly as in [**LTV23**, Section 1], let us denote the right (the one for the left action) and left (the one for the right action) eigenvectors for $\zeta_k^i \lambda_k$ by $v_i^k$ and $w_i^k$, normalized such that $(w_i^k)^T v_i^k = 1$. Let $v_i^k (w_i^k)^T[1]$ denote taking the sum of the first column of the matrix $v_i^k (w_i^k)^T$. Using this, we define

$$a_k(n) = \big(v_0^k (w_0^k)^T[1] \cdot 1 + v_1^k (w_1^k)^T[1] \cdot \zeta^n + v_2^k (w_2^k)^T[1] \cdot (\zeta^2)^n + \dots + v_{h-1}^k (w_{h-1}^k)^T[1] \cdot (\zeta^{h-1})^n\big) \cdot \lambda^n.$$

Similarly, but directly for $\Gamma$, we can define

$$\text{(12G-8)} \quad a(n) = \big(v_0 (w_0)^T[1] \cdot 1 + v_1 (w_1)^T[1] \cdot \zeta^n + v_2 (w_2)^T[1] \cdot (\zeta^2)^n + \dots + v_{h-1} (w_{h-1})^T[1] \cdot (\zeta^{h-1})^n\big) \cdot \lambda^n,$$

where $h$ is the period of the FBC of $\Gamma$ (for strongly connected graphs like the FBC this is defined in e.g. [**Kit98**, Begin of Section 7.1]). Note that, a priori, $a(n)$ might not be a finite expression.

**Theorem 12G.9.** *Assume that the growth problem $(R, C)$ is sustainably positively recurrent. Then:*

**(a)** *We have $a(n) \in \mathbb{R}_{\geq 0}$ and*

$$b(n) \sim a(n).$$

**(b)** *If the FBC is final in $\Gamma$, the limit $\lim_{k\to\infty} a_k(n)$ exists and is equal to $a(n)$.*

**(c)** *If the FBC is finite and final in $\Gamma$, then the convergence of $\lim_{n\to\infty} b_n/a_n = 1$ is geometric with ratio $|\lambda^{sec}/\lambda|$ and $|b_n - a_n| \in O\big((\lambda^{sec})^n + n^d\big)$ for some $d \in \mathbb{R}_{>0}$.*

*Remark* 12G.10. [**LTV23**] states an analog of Theorem 12G.9 for arbitrary finite graphs, but the statement itself is a bit nasty. We decided not to include its positively recurrent analog in this paper. ◇

*Proof of Theorem 12G.9.* We start with the statement about the asymptotic and then proof the other claims that we call “Finite approximation” and “Variance”.

*Asymptotic.* Condition (a) of Definition 12G.4 implies that there is a unique basic class $C^{FBC}(\Gamma)$ (which is a final basic class, of course). Note that $\lambda < \infty$ by Lemma 12G.2. Moreover, $\lambda = 0$ if and only if the growth problem is zero, and hence we can and will assume that $\lambda \in \mathbb{R}_{>0}$.

We borrow from and adjust [**Kit98**, Section 7.1]. There are two differences between our setting and [**Kit98**, Section 7.1] to keep in mind: Firstly, and crucial, we consider graphs that are potentially not strongly connected. Secondly, we can have a period, but that is a rather harmless generalization. We postpone this discussion to the end of the proof and assume first that we are in the aperiodic situation.

Unless stated otherwise, our indexing set below will be the vertices of the fusion graph $\Gamma$ associated to $(R,C)$, and these vertices are often called $i,j$.

**Notation 12G.11.** To have the same conventions as [**Kit98**] in this proof, let $T$ be the transpose of $M$. ◇

We define three generating functions $T$, $L$ and $R$:

$$T_{ij}(z) = \sum_{n\in\mathbb{Z}_{\geq 0}} t_{ij}^{(n)} z^n, \quad L_{ij}(z) = \sum_{n\in\mathbb{Z}_{\geq 0}} l_{ij}(n) z^n, \quad R_{ij}(z) = \sum_{n\in\mathbb{Z}_{\geq 0}} r_{ij}(n) z^n,$$
$$L'_{ij} = \sum_{n\in\mathbb{Z}_{\geq 0}} n l_{ij}(n) z^n, \quad \mu(i) = \lambda^{-1} L'_{ii}(\lambda^{-1}) \text{ (with } \mu(i)=\infty \text{ allowed)},$$

where the coefficients of $L$ and $R$ are defined as follows. The coefficient $l_{ij}(n)$ is the sum of the labels of the paths $i \to j$ of length $n$ that do not return to $i$ with $\leq n-1$ steps. Similarly, $r_{ij}(n)$ is the sum of the labels of the paths $i \to j$ of length $n$ that do not return to $j$ with $\leq n-1$ steps.

*Step 1.* We start with two crucial lemmas. We define $\lambda_{ij} = \limsup_{n\to\infty} \sqrt[n]{t_{ij}^{(n)}}$.

**Lemma 12G.12.** *We have $\lambda_{ij} \leq \lambda$ for all $i,j \in \Gamma$.*

*Proof.* Immediate from $(n \mapsto t_{ij}^{(n)}) \in o\big((\lambda+\varepsilon)^n\big)$ for all $\varepsilon \in \mathbb{R}_{>0}$. □

**Lemma 12G.13.** *For all $i,j \in \Gamma$, the radii of convergences of $T_{ij}(z)$, $L_{ij}(z)$ and $R_{ij}(z)$ are $\geq \lambda^{-1}$.*

*Proof.* The radius of convergence of $T_{ij}(z)$ is $\lambda_{ij}^{-1}$ by the Cauchy–Hadamard theorem, and the result follows from Lemma 12G.12 and the inequalities $l_{ij}(n) \leq t_{ij}^{(n)}$ and $r_{ij}(n) \leq t_{ij}^{(n)}$. □

We are now ready to go through statements in [**Kit98**, Section 7.1]. The references below are all with respect to [**Kit98**].

*Step 2.* Most parts of Lemma 7.1.6 are formal and work verbatim. Some of the used arguments need Lemma 12G.13, but then can be proven mutatis mutandis with it.

*Step 3.* Next, we see a mild modification. Namely, Lemma 7.1.7 for $i,j \in \Gamma$ in the same component works mutatis mutandis since $\lambda$ is the leading eigenvalue so $\lambda^{-1}$ is the minimal radius of convergence, see Lemma 12G.13. Moreover, $L_{ij}(\lambda^{-1})$ is finite in general if there is a path from $j$ to $i$, and vice versa for $R_{ij}(\lambda^{-1})$, by the same arguments as in [**Kit98**, Section 7.1].

*Step 4.* Lemma 7.1.8 remains true with the following changes. If $i \in C^{FBC}(\Gamma)$, then $T_{ij}(\lambda^{-1}) = \infty$ and $L_{ii}(\lambda^{-1}) = 1$ since the FBC is recurrent. We do not need the case where $i \notin C^{FBC}(\Gamma)$.

*Step 5.* Choose some enumeration of the vertices of $\Gamma$. For each $i \in C^{FBC}(\Gamma)$ define a row and a column vector

$$\ell^{(i)} = \begin{pmatrix} L_{i1}(\lambda^{-1}) & L_{i2}(\lambda^{-1}) & \dots & L_{ij}(\lambda^{-1}) & \dots \end{pmatrix},$$
$$r^{(i)} = \begin{pmatrix} R_{1i}(\lambda^{-1}) & R_{2i}(\lambda^{-1}) & \dots & R_{ji}(\lambda^{-1}) & \dots \end{pmatrix}^T.$$

These are well-defined by Step 3.

*Step 6.* Lemma 7.1.9.(i) remains true, that is, the $\ell^{(i)}$ and the $r^{(i)}$ are left and right $\lambda$-eigenvectors of the action matrix $T$:

$$\ell^{(i)} T = \lambda \cdot \ell^{(i)}, \quad T r^{(i)} = \lambda \cdot r^{(i)}.$$

(For the meticulous reader who want to double check that $r^{(i)}$ is an eigenvector we point out that the $R_{ii}$ in Lemma 7.1.6.(viii) should be a $R_{jj}$.)

Choose any $i \in C^{FBC}(\Gamma)$ and define

$$\ell = \ell^{(i)}, \quad r = r^{(i)}.$$

It is easy to see that the choice of the defining $i \in C^{FBC}(\Gamma)$ is not important for what we will do below.

*Step 7.* Lemma 7.1.10 remains true, with the same proof, but needs some adjustments. Firstly, the vectors $x$ and $y$ are assumed to have a nonzero entry corresponding to some $k \in C^{FBC}(\Gamma)$. Then $x$ and $y$ are strictly positive on the entries corresponding to all $k \in C^{FBC}(\Gamma)$. This uses irreducibility of the matrix supported on $C^{FBC}(\Gamma)$. (ii) and (iii) then remain true without further change.

*Step 8.* For Lemma 7.1.11 we replace $xT \leq \lambda x$ and $Ty \leq \lambda y$ by $xT = \lambda x$ and $Ty = \lambda y$ and we keep the same assumptions on $x$ and $y$ from the previous step. We renormalize $x$ and $y$ such that $x_i = 1$ and $y_i = 1$, and the previous step ensure that $x - l^{(i)}$ and $y - r^{(i)}$ are nonnegative and have a zero entry on $C^{FBC}(\Gamma)$. We then obtain that $x - l^{(i)}$ and $y - r^{(i)}$ are supported on the nonbasic classes. We can extract from $x - l^{(i)}$ and $y - r^{(i)}$ eigenvectors of a matrix whose spectrum does not contain $\lambda$, because of the condition (b) in Definition 12G.4. We finally obtain that $x = l^{(i)}$ and $y = r^{(i)}$.

*Step 9.* Since the FBC is positively recurrent, we get $L'_{ii}(\lambda^{-1}) < \infty$ for $i \in C^{FBC}(\Gamma)$. It follows that $l \cdot r = \mu(i) < \infty$ (dot product) by copying the calculation in the proof of Lemma 7.1.14.

*Step 10.* Lemma 7.1.15 remains true (it is independent of the growth problem).

*Step 11.* Theorem 7.1.18 has two cases: for $i \notin C^{FBC}(\Gamma)$ we are in case (i), and for $i \in C^{FBC}(\Gamma)$ we are in case (ii). The justification for $i \notin C^{FBC}(\Gamma)$ is clear, and the proof for $i \in C^{FBC}(\Gamma)$ works mutatis mutandis as in [**Kit98**], using the adjustments made in the steps above.

*Step 12.* Lemma 7.1.19.(i) works verbatim. For Lemma 7.1.19.(ii) there are three cases. If $i \in C^{FBC}(\Gamma)$, then the first equality holds. Moreover, for $j \in C^{FBC}(\Gamma)$ the equality between the first and the last term holds. Thus, for $i, j \in C^{FBC}(\Gamma)$ Lemma 7.1.19.(ii) works in the same way as in [**Kit98**].

*Step 13.* Lemma 7.1.20 works with the extra assumption that $i, k \in C^{FBC}(\Gamma)$ in part (a) of (i), $j, k \in C^{FBC}(\Gamma)$ in part (b) of (i) and $i, j \in C^{FBC}(\Gamma)$ in (ii).

*Step 14.* As in [**Kit98**], steps 1-13 prove [**Kit98**, Theorem 7.1.3.(f)] for sustainably positively recurrent growth problems.

*Step 15.* Now, we can mimic the proof of [**LTV23**, Theorem 7] to obtain the precise asymptotic on $b(n)$.

For periodic matrices evaluate the power series at $\zeta^k \lambda$ for $\zeta$ a root of unity. With this change all the above work mutatis mutandis. See also the final part of [**Kit98**, Section 7.1].

*Finite approximation.* Since the FBC is final, the vector $\ell$ is supported by the FBC. The first column of the matrices in Equation 12G-8 are then supported by the FBC. Hence, [**Kit98**, Theorem 7.1.4] applies and we are done.

*Variance.* As in the previous point, but since the FBC is now assumed to be finite, we can use [**LTV23**, Theorem 1]. □

*Remark* 12G.14. For Krull–Schmidt monoidal categories with finitely many indecomposable objects Theorem 12G.9 recovers the main results of [**LTV23**] (up to the point made in Remark 12G.10). ◇

**Example 12G.15.** Let us come back to Example 12F.15. The random walk on the fusion graph $\Gamma$ in Equation 12F-16 is positively recurrent. Thus, Theorem 12G.9 applies and we get

$$b_n \sim \tfrac{1}{4} \cdot 3^n$$

(note that $|G| = 4$ and $\sum_{L \text{ simples}} \dim_{\mathbb{F}_2} L = 1$) using the finite cutoffs as indicated in Equation 12F-16. The second largest eigenvalue for all finite cutoffs is zero, so the variance is bounded by some polynomial function. Indeed, we get:

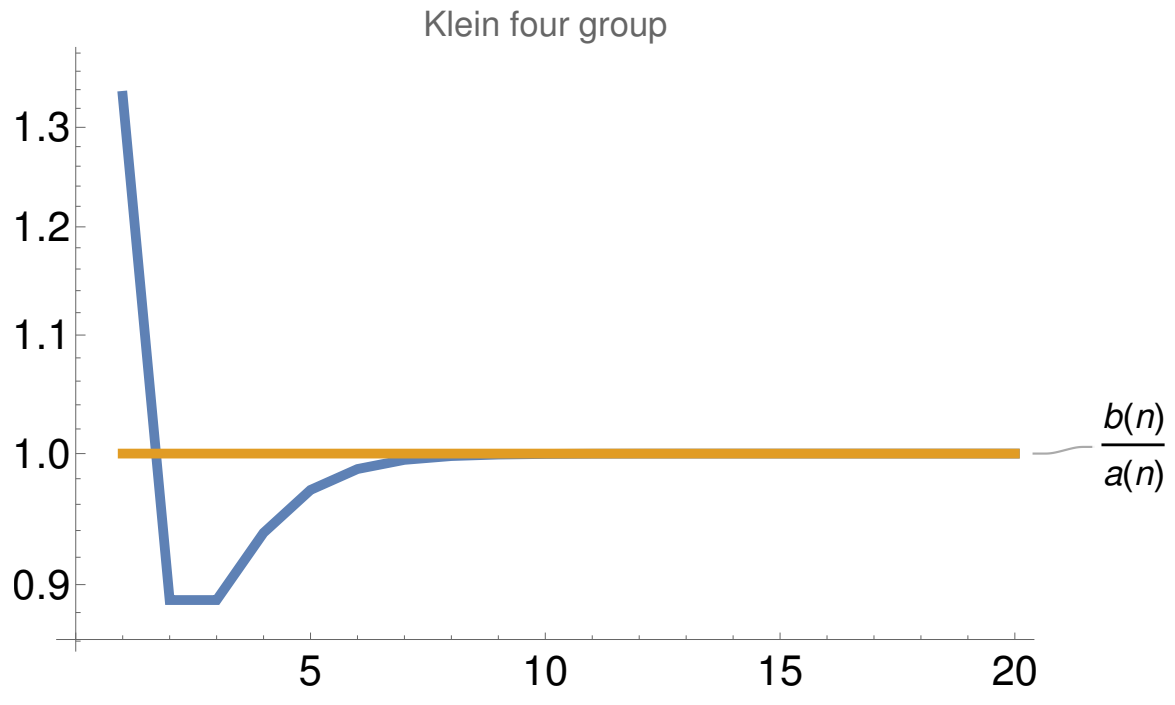

,
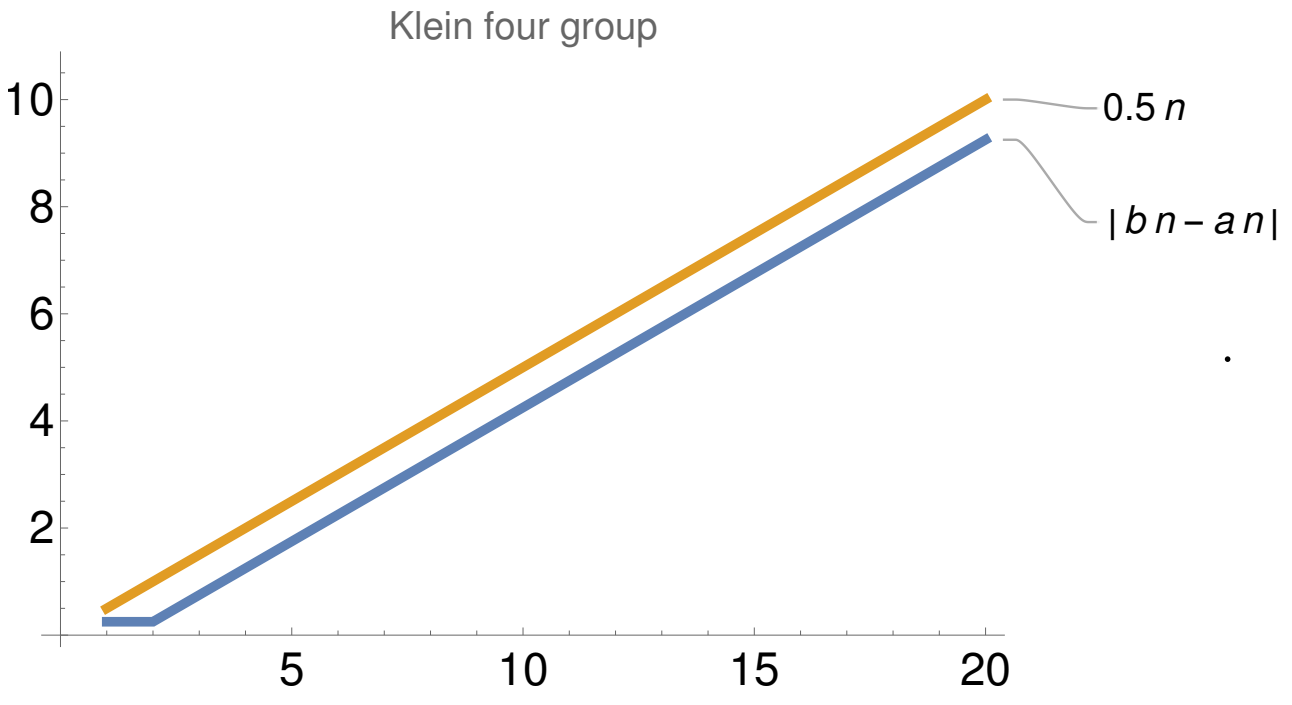

.

The left plot is a logplot, the right is a standard plot. ◇

**Example 12G.16.** Another example from the realm of finite groups is the following. Consider $G = \mathrm{PSL}_2(\mathbb{F}_7)$. Over $\mathbb{F}_2$ the representation theory of $G$ is not semisimple. In this case there are two three dimensional simple representations, and let $\mathtt{X}$ be any of these two. The choice of this example is motivated by [**Cra13**, Theorem 1.4], and is the smallest example on that list.

The fusion graph of the associated growth problem $(G, \mathtt{X})$ can be seen to be (illustrated below with a cutoff)

(12G-17) 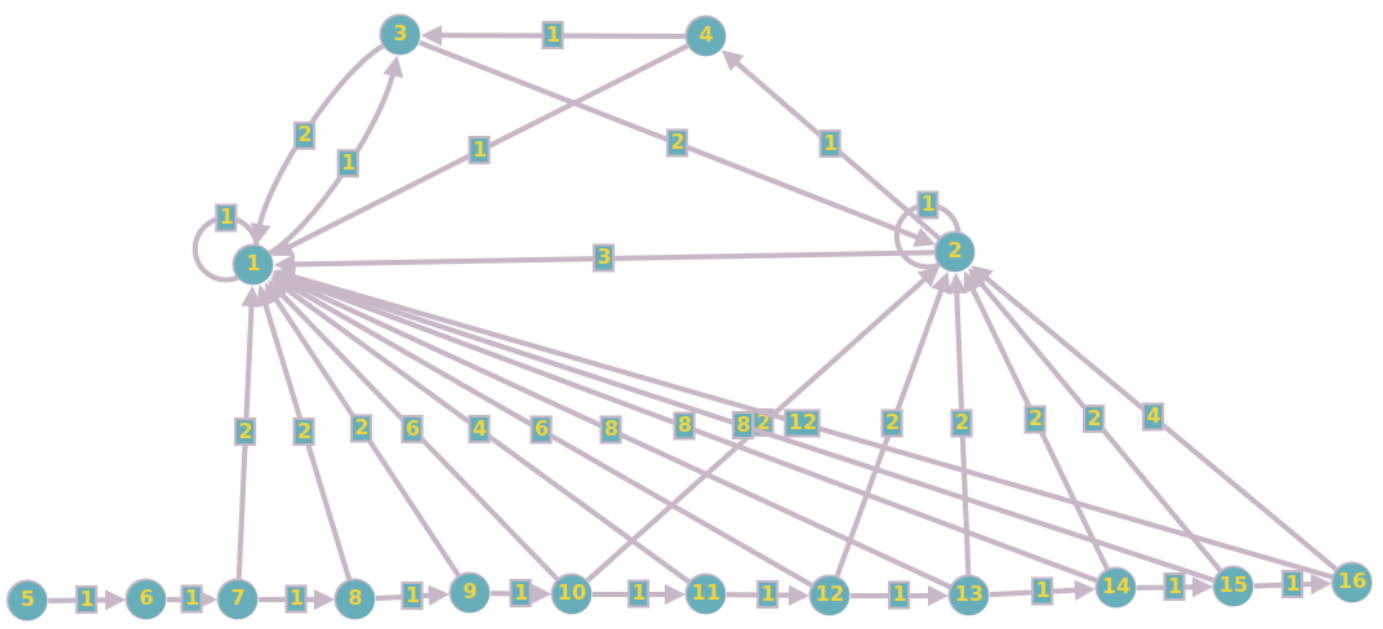,

$$M(PSL2) = \begin{pmatrix} 1&3&2&1&0&0&2&2&2&6&4&6&8&8&8&12 \\ 0&1&2&0&0&0&0&0&0&2&0&2&4&4&4&4 \\ 1&0&0&1&0&0&0&0&0&0&0&0&0&0&0&0 \\ 0&1&0&0&0&0&0&0&0&0&0&0&0&0&0&0 \\ 0&0&0&0&0&0&0&0&0&0&0&0&0&0&0&0 \\ 0&0&0&0&1&0&0&0&0&0&0&0&0&0&0&0 \\ 0&0&0&0&0&1&0&0&0&0&0&0&0&0&0&0 \\ 0&0&0&0&0&0&1&0&0&0&0&0&0&0&0&0 \\ 0&0&0&0&0&0&0&1&0&0&0&0&0&0&0&0 \\ 0&0&0&0&0&0&0&0&1&0&0&0&0&0&0&0 \\ 0&0&0&0&0&0&0&0&0&1&0&0&0&0&0&0 \\ 0&0&0&0&0&0&0&0&0&0&0&0&0&0&0&0 \\ 0&0&0&0&0&0&0&0&0&0&0&1&0&0&0&0 \\ 0&0&0&0&0&0&0&0&0&0&0&0&1&0&0&0 \\ 0&0&0&0&0&0&0&0&0&0&0&0&0&1&0&0 \\ 0&0&0&0&0&0&0&0&0&0&0&0&0&0&1&0 \end{pmatrix}.$$

The projective cell is illustrated at the top and has the adjacency matrix

$$M = \begin{pmatrix} 1&0&1&0 \\ 3&1&0&1 \\ 2&2&0&0 \\ 1&0&1&0 \end{pmatrix}.$$

The dimensions of the projective indecomposables, in the order as in Equation 12G-17, are 8, 16, 16, 8. The eigenvalues of $M$ are $\{3, \frac{1}{2}(-1+i\sqrt{7}), \frac{1}{2}(-1-i\sqrt{7}), 0\}$, and our dominating growth rate is 3.

The dimensions of the indecomposables along the bottom in Equation 12G-17 is

$$1, 3, 9, 11, 17, 35, 25, 43, 49, 51, 57, 75, 65, 83, 89, 91, 97, 115, 105, 123, 129, \dots,$$
$$\text{pattern}: +6, +2, +6, +18, -10, +18, \text{repeat},$$
$$\text{generating function}: \frac{16x^5 + 5x^4 + 7x^3 + 8x^2 + 3x + 1}{(x-1)^2(x^3 + 2x^2 + 2x + 1)},$$

and this pattern continues. From this pattern one gets the whole graph.

Moreover, one can check that $v_0^k (w_0^k)^T[1]$ for the leading eigenvectors converges to $\frac{15}{168} \approx 0.089$ as $k \to \infty$. (Note that $|G| = 168$ and $\sum_{L \text{ simples}} \dim_{\mathbb{F}_2} L = 15$.) Thus, Theorem 12G.9 implies that:

$$b_n \sim \frac{15}{168} \cdot 3^n.$$

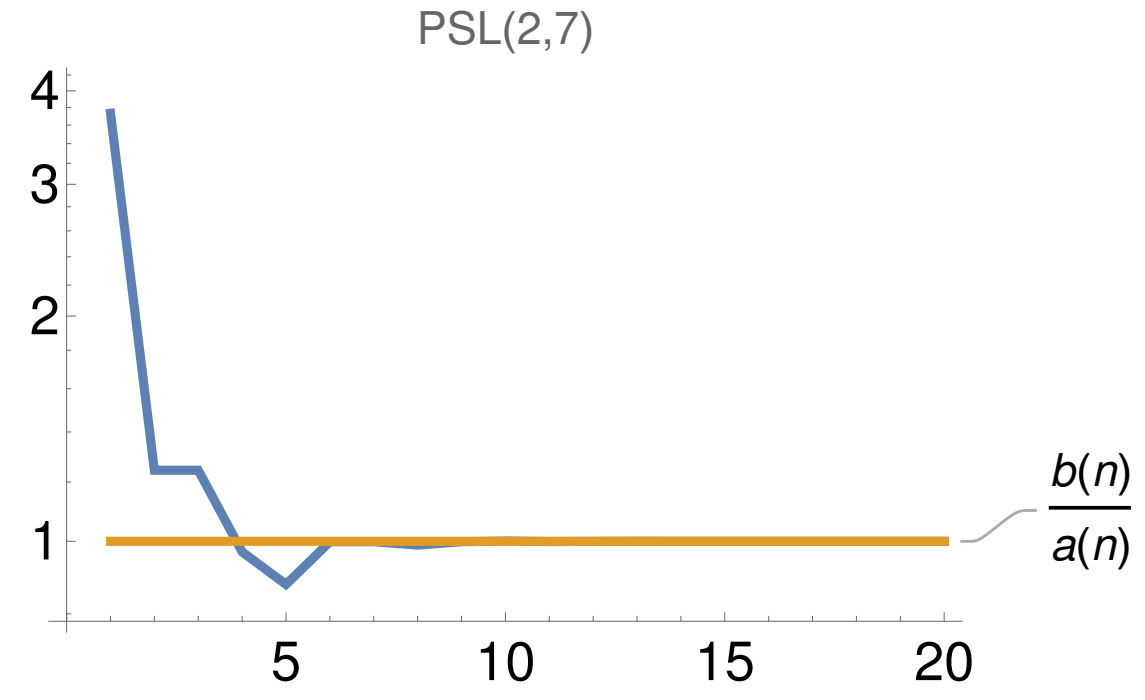

,
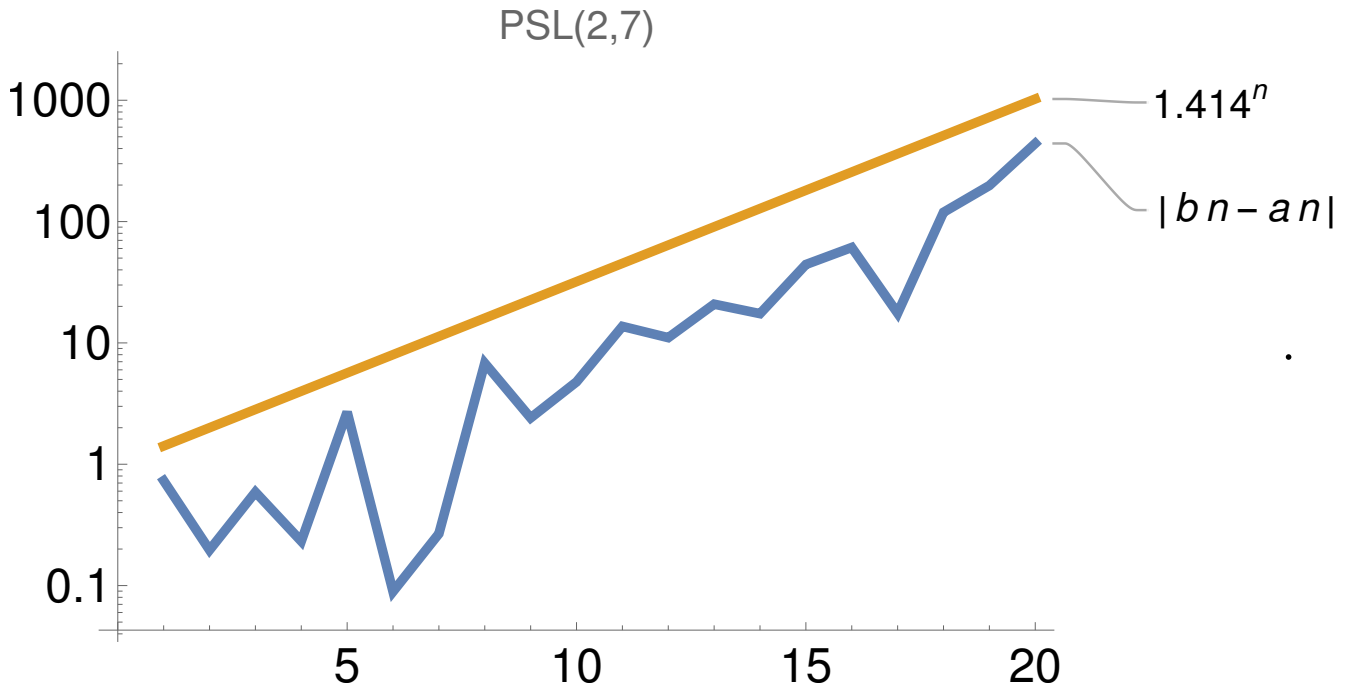

.

Note that $|\lambda^{sec}| = |\frac{1}{2}(-1+i\sqrt{7})| = \sqrt{2} \approx 1.414$. The convergence is geometric with ratio $\sqrt{2}/3 = |\frac{1}{2}(-1+i\sqrt{7})|/3 \approx 0.4714$. ◇

**Example 12G.18.** Let $W$ be a finite Coxeter group and $\mathbf{SBim} = \mathbf{SBim}(W)$ the associated monoidal category of Soergel bimodules. Let $\mathtt{X} \in \mathbf{SBim}$ be some generating object. The growth problem $(\mathbf{SBim}, \mathtt{X})$ is positively recurrent and has the following asymptotics.

In the Grothendieck ring, the group ring $\mathbb{Z}W$, write $[\mathtt{X}] = \sum_{s \in W} m_s \cdot s$ where $m_s \in \mathbb{Z}_{\geq 0}$. One gets

$$b(n) \sim \frac{1}{|W|} \cdot \big(\sum_{s \in W} m_s\big)^n.$$

This follows directly from Theorem 12G.9 since one only has one FBC corresponding to the indecomposable Soergel bimodule associated with the longest word. ◇

*Remark* 12G.19. In Example 12G.18, we lack a general expression for the second largest eigenvalue, except in the dihedral case for $\mathtt{X}$ being the indecomposable object corresponding to the product of the two simple reflections. For this growth problem, the eigenvalue depends on a parity condition related to whether the order of $W$ (always even) is divisible by four. Using, for example, [**Tub22**, Section 3] one can show that the coefficients of the minimal polynomial of the second largest eigenvalue in these cases are the signed versions of the sequences [**OEI23**, A085478] and [**OEI23**, A030528], respectively. Thus, the second largest eigenvalues are closely related, though not identical, to the PF roots of the ***Morgan–Voyce polynomials***. A more general statement would be desirable. ◇

12H. **Schur–Weyl duality, sandwich cellular algebras and growth problems.** For some field $\mathbb{k}$, assume that one has an additive Krull–Schmidt monoidal $\mathbb{k}$-linear category $\mathbf{C}$ with finite-dimensional hom-spaces, and an object $\mathtt{X} \in \mathbf{C}$. For simplicity, assume that $\mathbf{C}$ is semisimple. A version of Schur–Weyl duality implies that ($\mathtt{Y}$ is a simple summand of $\mathtt{X}^{\otimes n}$ that appears with multiplicity $m > 0$ if and only if the semisimple algebra $A_n = \mathrm{End}_{\mathbf{C}}(\mathtt{X}^{\otimes n})$ has a simple representation $L_{\mathtt{Y}}$ of dimension $\dim_{\mathbb{k}} L_{\mathtt{Y}} = m$), and there is a bijection between such $\mathtt{Y}$ and $L_{\mathtt{Y}}$. In particular, we have:

**Lemma 12H.1.** *In the above setting, $A_n$ is semisimple and*

$$b_n = \sum_L \dim_{\mathbb{k}} L, \quad (\textit{sum over simple } A_n\textit{-representations } L).$$

*Proof.* Directly from the above discussion, which, in turn, can be justified as in [**AST18**, Section 4C]. See also [**Erd95**, **Soe99**]. □

We will now assume familiarity with sandwich cellular algebras [**Bro55**, **TV23**, **Tub22**], or at least with some variation of it, most notably, [**FG95**].

**Lemma 12H.2.** *If $A_n$ is a semisimple involutive sandwich cellular algebra with apex set $\mathcal{P}^{ap}$, bottom sets $\mathcal{B}_\lambda$ and sandwiches algebras $\mathcal{H}_\lambda$, then*

$$b_n = \sum_{\lambda \in \mathcal{P}^{ap}} \big(\#\mathcal{B}_\lambda \cdot \sum_L \dim_{\mathbb{k}} L\big), \quad (\textit{inner sum over simple } \mathcal{H}_\lambda\textit{-representations } L).$$

*Proof.* Directly from Lemma 12H.1 and the standard theory of sandwich cellular algebras. □

We will use Lemma 12H.2 silently below.

12I. **Growth problems in diagram categories.** We now summarize [**GT25b**].

We begin with the following table, the meaning of which we will explain shortly:

(12I-1)

| Symbol | Diagrams | X | $\sqrt[n]{b_n} \sim$ |
|---|---|---|---|
| $\mathbf{2Cob}_k$ | 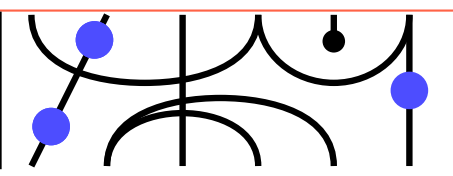 | • | $\frac{2}{e} \cdot \frac{n}{\log 2n/k}$ |

| Symbol | Diagrams | X | $\sqrt[n]{b_n} \sim$ | Symbol | Diagrams | X | $\sqrt[n]{b_n} \sim$ |
|---|---|---|---|---|---|---|---|
| $\mathbf{pPa}_t$ | 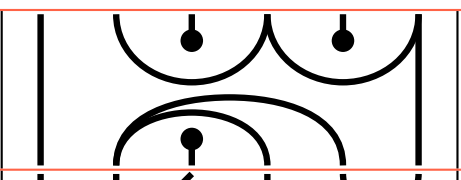 | • | 4 | $\mathbf{Pa}_t$ | 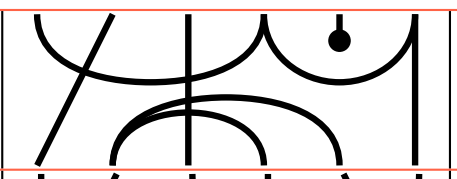 | • | $\frac{2}{e} \cdot \frac{n}{\log 2n}$ |
| $\mathbf{Mo}_t$ | 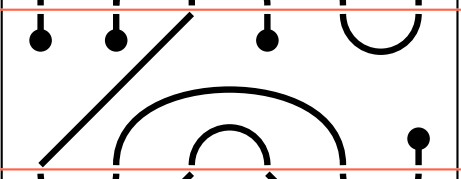 | • | 3 | $\mathbf{RoBr}_t$ | 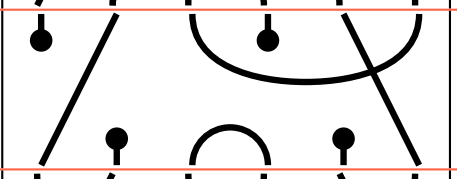 | • | $\frac{\sqrt{2}}{\sqrt{e}} \cdot \sqrt{n}$ |
| $\mathbf{TL}_t$ | 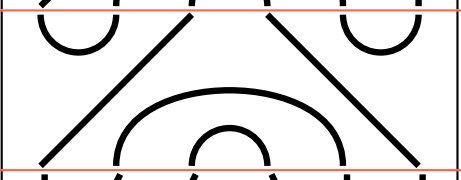 | • | 2 | $\mathbf{Br}_t$ | 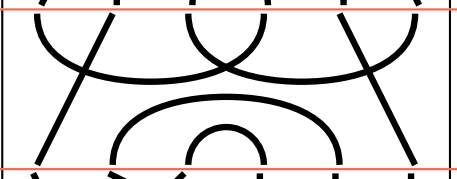 | • | $\frac{\sqrt{2}}{\sqrt{e}} \cdot \sqrt{n}$ |
| $\mathbf{pRo}_t$ | 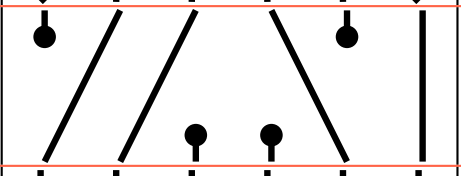 | • | 2 | $\mathbf{Ro}_t$ | 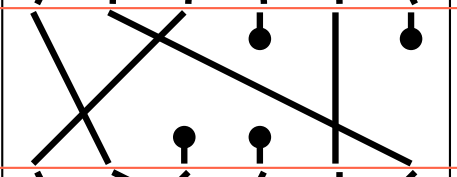 | • | $\frac{\sqrt{1}}{\sqrt{e}} \cdot \sqrt{n}$ |
| $\mathbf{pS}_t$ | 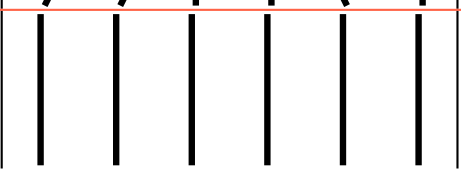 | • | 1 | $\mathbf{S}_t$ | 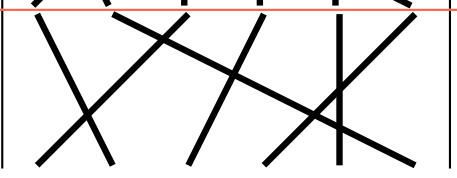 | • | $\frac{\sqrt{1}}{\sqrt{e}} \cdot \sqrt{n}$ |

| Symbol | Diagrams | X | $\sqrt[n]{b_n} \sim$ | Symbol | Diagrams | X | $\sqrt[n]{b_n} \sim$ |
|---|---|---|---|---|---|---|---|
| $\mathbf{oBr}_t$ | 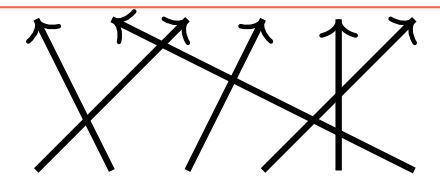 | $\uparrow$ | $\frac{\sqrt{1}}{\sqrt{e}} \cdot \sqrt{n}$ | $\mathbf{oBr}_t$ | 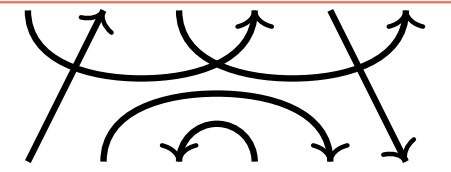 | $\uparrow\downarrow$ | $\frac{2}{e} \cdot n$ |

| Symbol | Diagrams | X | $\sqrt[n]{b_n} \sim$ |
|---|---|---|---|
| $\mathbf{Rep}(\mathrm{GL}_t(\mathbb{F}_q))$ | 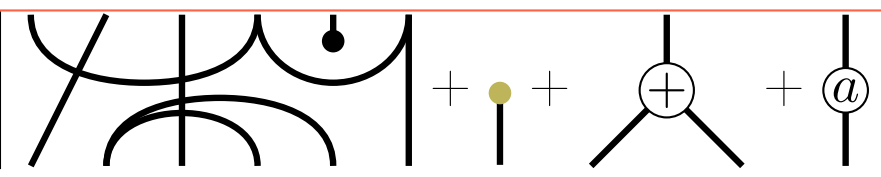 for $a \in \mathbb{F}_q$ | • | $q^{(n+1)/2}$ |

We have seen some of these already:

**(a)** $\mathbf{TL}_t = \mathbf{TL}^q_{\mathbb{S}\oplus\mathbb{E}}$ is the Rumer–Teller–Weyl category from Definition 7F.1, with parameters chosen such that $t = -(q + q^{-1}) \in \mathbb{C}$ is not a root of unity (e.g. transcendental).

**(b)** $\mathbf{Br}_t = \mathbf{TL}^q_{\mathbb{S}\oplus\mathbb{E}}$ is the quantum Brauer category from Definition 7F.11, with parameters chosen such that $q = 1$ and circle value $t \in \mathbb{C}$ is not a root of unity (e.g. transcendental). Note that this does not quite make sense the way we defined the quantum Brauer category, but this is easy to fix, cf. Exercise 12J.4.

**(c)** $\mathbf{S}_t$ is the symmetric group category, and $\mathbf{pS}_t$ is the planar symmetric group category (which is trivial). The parameter $t$ does not play any role for these.

**(d)** We will define $\mathbf{2Cob}_k$ below

**(e)** The reader is supposed to define all the others, except $\mathbf{Rep}(\mathrm{GL}_t(\mathbb{F}_q))$ as this is tricky, see Exercise 12J.4. Everything except $\mathbf{Rep}(\mathrm{GL}_t(\mathbb{F}_q))$ is a subcategory, upon specialization, of $\mathbf{2Cob}_k$.

We now define the ***category of cobordisms***: Here the objects are one-dimensional compact manifolds (circles) and the morphisms are two-dimensional cobordisms (pants), which we will draw using their spines with handles as dots:

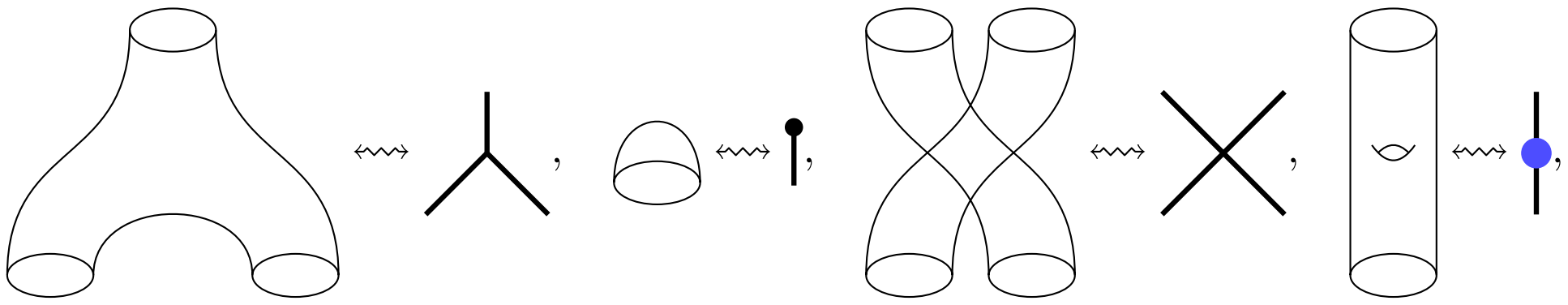 *etc.*

More formally, we define:

**Definition 12I.2.** Let $\mathbf{2Cob}_\infty$ be the monoidal category with $\otimes$-generating object • and generating $\circ$-$\otimes$-generating morphisms

multiplication: 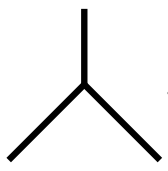, comultiplication: 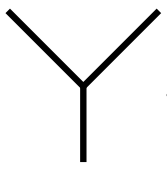, unit: 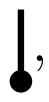, counit: 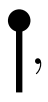, crossing: 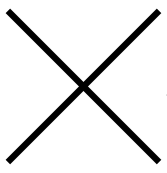,

modulo the $\circ$-$\otimes$-ideal that makes the crossing a symmetry and • a symmetric Frobenius object with the structure maps matching the nomenclature. (This defines S, T and even R; pictures for the relations can be easily drawn or, alternatively, found in e.g. [Koc04].) ◇

The ***diagrammatic antiinvolution*** $*$ flips a cobordism up-side-down. We call a diagram a ***merge diagram*** if it contains only multiplications, counits and a minimal number of crossings. A ***split diagram*** is a $*$-flipped merge diagram. A ***dotted permutation diagram*** contains only dots (=handles) and crossings. Here are some examples (flipping the left illustration gives a split diagram):

Merge diagram: 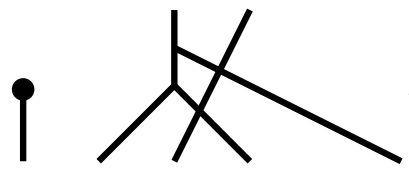 , dotted permutation diagram: 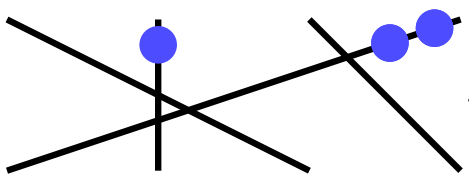 .

Abusing notation, for a field $\Bbbk$ we denote the $\Bbbk$-linear extension of the cobordism category also by $\mathbf{2Cob}_{\infty}$. Recall sandwich cellularity, the picture for $\mathbf{2Cob}_k$ is:

(12I-3) 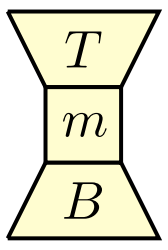

where 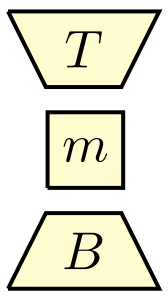

a split diagram, a dotted permutation, a merge diagram. Involution: 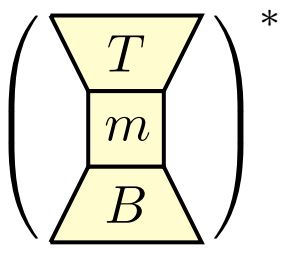

= 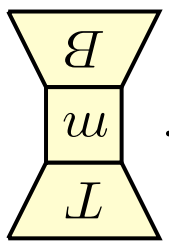 .

We call the existence of a spanning set with a decomposition as above a precell structure and the respective algebras presandwich.

**Lemma 12I.4.** *The endomorphism algebras of* $\mathbf{2Cob}_{\infty}$ *are involutive presandwich cellular with precell structure as in Equation 12I-3.*

*Proof.* Immediate from [**Koc04**, Section 1.4.16]. □

Now fix two polynomials $p, q \in \Bbbk[x]$, and consider the Taylor expansion of $p(x)/q(x) = \sum_{i=0}^{\infty} a_i x^i$, and let $k = \max\{\deg p + 1, \deg q\}$. Let $\mathbf{2Cob}_k$ be the quotient of $\mathbf{2Cob}_{\infty}$ by the $\circ$-$\otimes$-ideal generated by

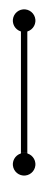 $= a_0$, 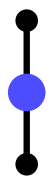 $= a_1$, 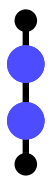 $= a_2$, ....

(The category $\mathbf{2Cob}_k$ actually depends on $p, q$ but we suppress this in the notation.)

**Example 12I.5.** For $p = t$ and $q = 1 - x$ we have $p(x)/q(x) = t(1 + x + x^2 + x^3 + \dots)$ so that all closed surfaces in $\mathbf{2Cob}_k$ evaluate to $t$. We call the resulting category the ***partition category***. ◇

Let us denote the endomorphism monoid of $\bullet^{\otimes n}$ by $\mathbf{2Cob}_k(n)$. Recall the Ariki–Koike algebra (cyclotomic Hecke algebra) as defined in [**AK94**, **BM93**, **Che87**].

**Lemma 12I.6.** *Dotted permutations in* $\mathbf{2Cob}_k(n)$ *span an algebra isomorphic to the Ariki–Koike algebra* $A(n, k)$ *on* $n$ *strands with a cyclotomic relation of degree* $k$ *and trivial quantum parameter. Moreover, Lemma 12I.4 can be refined into a sandwich cell datum with* $A(m, k)$ *for* $m \in \{0, \dots, n\}$ *as the sandwiched algebras.*

*Proof.* This follows from [**KOK22b**], *e.g.* the text around (11) therein, which implies that the handles satisfy a minimal polynomial of degree $k$, and the same arguments as in [**TV23**, Section 6]. □

For the next statement we assume familiarity with the usual tableaux combinatorics, see, for example, [**DJM98**, **Mat99**].

**Proposition 12I.7.** *Let* $\Bbbk$ *be of characteristic* $p$.

- **(a)** *The set of apexes of* $\Bbbk\mathbf{2Cob}_k(n)$ *is* $\{0, \dots, n\}$.
- **(b)** *The finite-dimensional simple* $\Bbbk\mathbf{2Cob}_k(n)$*-modules of apex* $m$*, up to equivalence, are indexed by* $p$*-restricted* $k$*-multipartitions of* $m$.
- **(c)** *The dimensions of the cell representation for a* $p$*-restricted* $k$*-multipartition* $\lambda$ *of* $m$ *is*

$$\#\{\text{merge diagrams with } m \text{ top strands}\} \cdot \#\{\text{standard tableaux of shape } \lambda\}.$$

*Moreover, if* $\Bbbk\mathbf{2Cob}_k(n)$ *is semisimple, then the cell representations are simple.*

*Proof.* Immediate from the standard theory of sandwich cellular algebras as, e.g., in [**Tub22**], Lemma 12I.6 and the cell structure of $A(m, k)$ as, for example, in [**DJM98**, Theorem 3.26]. □

We now need some counting lemmas.

**Lemma 12I.8.** *The number of merge diagrams from* $n$ *bottom strands to* $m$ *top strands is*

$$M_n^m = \sum_{i=m}^{n} \left\{ {n \atop i} \right\} \binom{i}{m},$$

*where the curly brackets denote the Stirling numbers of the second kind.*

*Proof.* A standard count that is independent of $k$, and therefore the same as in the partition category. Details are omitted; however, if the reader encounters difficulties, [**HR15**, Section 4] provides helpful guidance. □

Let $\mathrm{STab}(m,k)$ is the set of standard $k$-multitableaux of $m$ and let $\#\mathrm{STab}(m,k)$ denote its size. Assume from now that $\Bbbk=\mathbb{C}$ and that $\mathbf{2Cob}_k$ is semisimple.

**Lemma 12I.9.** *We have the formula*

$$b_n=\sum_{m=0}^{n} M_n^m\#\mathrm{STab}(m,k).$$

*Proof.* By Proposition 12I.7 and Lemma 12I.8. □

**Lemma 12I.10.** *$b_n$ has exponential generating function* $\exp(\frac{k}{2}\exp(2x)+\exp(x)-\frac{k+2}{2})$.

*Proof.* We first observe that

$$\#\mathrm{STab}(m,k)=k^{\lceil m/2\rceil}\sum_{i\in\mathbb{Z}_{\geq 0}} Bes(m,i)k^{\lfloor m/2\rfloor-i},$$

where we use the Bessel numbers $Bes(m,i)=m!/(i!(n-2i)!2^i)$. Thus, we get

$$b_n=\sum_{m=0}^{n}\sum_{i\in\mathbb{Z}_{\geq 0}} M_n^m k^{\lceil m/2\rceil}Bes(m,i)k^{\lfloor m/2\rfloor-i},$$

and we can use the same calculations as in [**Qua07**, Section 3] (which only uses the exponential generating functions for the Stirling and Bell numbers). □

**Lemma 12I.11.** *We have the asymptotic formula*

$$b_n\sim\frac{\left(\frac{n}{z}\right)^{n+\frac{1}{2}}\exp(\frac{k}{2})\exp(2z+1)\exp(z-n-\frac{k+2}{2})}{\sqrt{2\frac{k}{2}\exp(2z)(2z+1)+\exp(z)(z+1)}},$$

$$z=\frac{W\left(\frac{2n}{k}\right)}{2}-\frac{1}{4\left(\frac{k}{2}\right)n^{1-\frac{1}{2}}\left(W\left(\frac{2n}{k}\right)+1\right)W\left(\frac{2n}{k}\right)^{\frac{1}{2}-2}+\frac{2}{W\left(\frac{2n}{k}\right)}+1},$$

*where $W$ is the Lambert $W$ function.*

*Proof.* Having Lemma 12I.10, the proof of this is automatized, see for example [**GT25a**, **Kot22**]. This works roughly as follows. Let $f$ be the exponential generating function. One then uses Hayman's method [**Hay56**] and computes the asymptotic of $\lim_{x\to\infty} xf'(x)/f(x)$. □

**Theorem 12I.12.** *The formula in Equation 12I-1 holds.*

*Proof.* The proof is also automatized, using Lemma 12I.11 and the code on [**GT25a**]. Essentially, Mathematica has a build in function for this purpose that does exact calculations. (We are referring to 'Asymptotic'; Introduced in 2020 (12.1) | Updated in 2022 (13.2).) □

We will not discuss $\mathbf{Rep}\big(\mathrm{GL}_t(\mathbb{F}_q)\big)$ or the oriented versions here, see [**GT25b**] for details, but for the remaining ones we proceed as follows.

For $\mathbf{Pa}_n$ we simply specialize $k=1$ above. See also Example 12I.5. The corresponding sequence for $b_n$ is [**OEI23**, A002872]. However, following [**OEI23**, A002872] one gets other, slightly nastier, formulas, namely:

$$b_n\sim\left(\frac{2n}{\mathrm{W}(2n)}\right)^n\cdot\exp\left(\frac{n}{\mathrm{W}(2n)}+\left(\frac{2n}{\mathrm{W}(2n)}\right)^{\frac{1}{2}}-n-\frac{7}{4}\right)\Big/\sqrt{1+\mathrm{W}(2n)},$$

$$\sqrt[n]{b_n}\sim\frac{2}{e}\cdot\frac{n}{\log 2n}\sqrt[\log 2n]{e}.$$

Now we consider subalgebras of partition algebras. By an easy (and well known) diagrammatic argument we get $\mathbf{pPa}_n\cong\mathbf{TL}_{2n}$, we already discussed $\mathbf{Pa}_n$ and $\mathbf{pS}_n$ is trivial, so we do not need to address these cases. Let us list the sandwich cellular bases for the remaining diagram categories:

$\mathbf{Mo}_n,\mathbf{TL}_n,\mathbf{pRo}_n$:

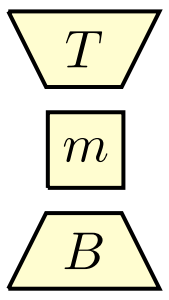

cup and unit diagrams,
an identity,
cap and counit diagram,

$\mathbf{RoBr}_n,\mathbf{Br}_n,\mathbf{Ro}_n$:

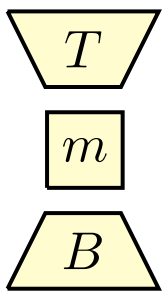

cup and unit diagrams,
a permutation,
cap and counit diagram.

For the symmetric group $\mathbf{S}_n$ the sandwich structure is trivial.

From this one gets explicit formulas, matching the ones in [HR15, Section 4]. Here is the list of remaining sequences:

$$\mathbf{Mo}_n\colon\ [\text{OEI23}, \text{A005773}],\quad \mathbf{RoBr}_n\colon\ [\text{OEI23}, \text{A000898}],$$
$$\mathbf{TL}_n\colon\ [\text{OEI23}, \text{A000984}],\quad \mathbf{Br}_n\colon\ [\text{OEI23}, \text{A047974}],$$
$$\mathbf{pRo}_n\colon\ [\text{OEI23}, \text{A000079}],\quad \mathbf{Ro}_n\colon\ [\text{OEI23}, \text{A005425}],$$
$$\mathbf{S}_n\colon\ [\text{OEI23}, \text{A000085}].$$

Let us just focus on $\mathbf{Br}_n$; the others being similar. In this case the exponential generating function is $\exp(x^2+x)$, and Mathematica gives

$$b_n \sim 2^{\frac{n}{2}-\frac{1}{2}} \exp\left(\sqrt{\frac{n}{2}} - \frac{n}{2} - \frac{1}{8}\right) n^{\frac{n}{2}},$$
$$\sqrt[n]{b_n} \sim \frac{\sqrt{2}}{\sqrt{e}} \cdot \sqrt{n}.$$

This completes the proof that Equation 12I-1 holds.

### 12J. Exercises.

*Exercise* 12J.1. What does the MAGMA code in Section 12C count? Adjust it to do the same for the alternating group. ◇

*Exercise* 12J.2. Prove all statements in Example 12D.2. Note that MAGMA does exact calculations. ◇

*Exercise* 12J.3. Fill in the details in Proposition 12C.3. ◇

*Exercise* 12J.4. How would one define the category $\mathbf{Br}_t$ in Section 12I? Define all the other diagram categories, except $\mathbf{Rep}\big(\mathrm{GL}_t(\mathbb{F}_q)\big)$. ◇

*Exercise* 12J.5. Give more details for the counts leading to Equation 12I-1. ◇

## 13. How good are quantum invariants? – a big data approach

The quantum invariants are one of the central, if not the central, topics of these lecture notes. Is it time to challenge our assumptions?

How well do quantum invariants actually perform?

Story time: "When I first got 1Gbps internet speed in 2024 (which was fast at the time), I immediately tested whether I was actually getting what I paid for. Sure enough, I measured around ~980Mbps at the router and ~800Mbps on my computer. As a rule, I like to double-check things, just to be sure. Or at least, that's what I thought until I attended the outstanding conference 'Diagrammatic Intuition and Deep Learning in Mathematics' (York, July 2024) and Radmila Sazdanovic's talk. During the talk, I realized something crucial: I had been working with quantum invariants for a quarter of a century without ever stopping to ask whether they were truly as reliable as I assumed!"

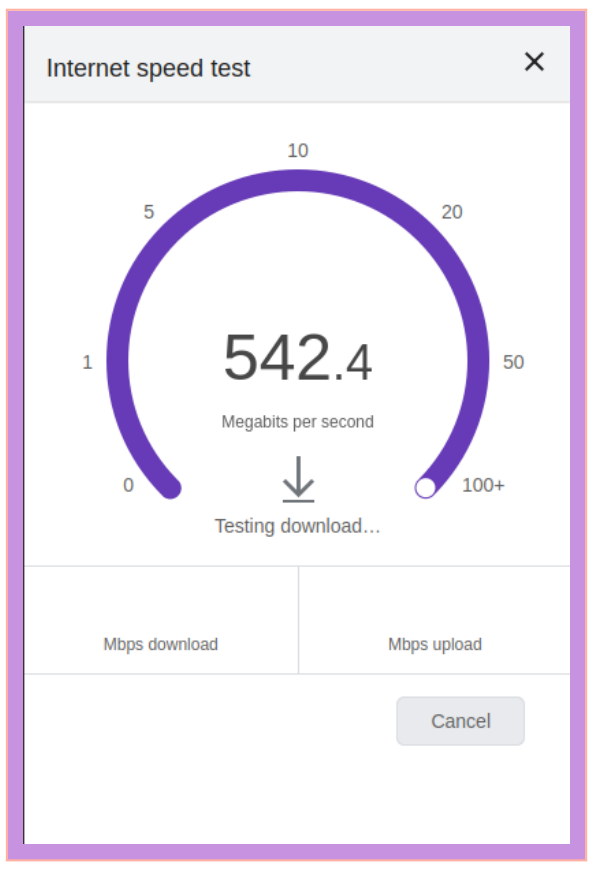

FIGURE 33. Always double-check claims; or, more importantly, be sure of what you should be double-checking.
Screenshot taken 10.Feb.2025.

The story above marks the beginning of the paper [TZ25a], which we will summarize in this section.

13A. **A word about conventions.** Here are some computer-science-style conventions that we will use:

*Convention* 13A.1. In this section, following [**LTV24b**], we have ***conjectures*** and ***speculations***. Conjectures are presented in their standard sense, while speculations refer to preliminary hypotheses that lack full support from the data. We hope that both serve as an inspiration to prove or, equally exciting, disprove the corresponding statements. ◇

*Convention* 13A.2. We collect the various polynomials as ***lists*** (ordered, repetitions allowed, notation [_]). If not stated otherwise, we will identify a polynomial with a vector as follows. We use some upper bound $d$ for the absolute degree of polynomials appearing. Then

$$a_{-d}q^{-d} + ... + a_d q^d \rightsquigarrow [a_{-d}, ..., a_d] \in \mathbb{Z}^{2d+1} \subset \mathbb{R}^{2d+1},$$

where we assign zero to the coefficient of included terms not appearing in the polynomial. In this way, all polynomials are in the same vector space. Similarly, for a two-variable polynomial we get a matrix containing the various coefficients, and we ***flatten*** it into a vector. ◇

*Convention* 13A.3. For a vector $v$, we write $\sum v$ for the sum over all entries of a vector, $\mathtt{max}(v)$ for the maximal entry, $\mathtt{abs}(v)$ for the componentwise absolute value of the vector, $\mathtt{roots}(v)$ for the set of roots of the corresponding normalized polynomial, and we write $\mathtt{spread}(v)$ for the number of entries ignoring padded zeros, *i.e.* zeros at the beginning and end of a vector. This corresponds to the difference between the maximal and minimal degree of a given polynomial. Naturally, `spread` only applies to single variable polynomials. ◇

*Convention* 13A.4. To analyze functions that behave like $a^n$ for some $a \in \mathbb{R}_{>1}$ we use ***successive quotients***. For a function $f\colon \mathbb{Z}_{\geq 0} \to \mathbb{R}_{>0}$, this is the function $n \mapsto \frac{f(n+1)}{f(n)}$, whose value effectively approximates $a$ for $f(n) = a^n$. ◇

*Convention* 13A.5. As a final note, particularly for the data below, we sometimes use a ***log plot*** for displaying data, which means that we use a logarithmic scale on the $y$-axis, and when we show floating point numbers, they are ***truncated*** (*i.e.* a floor) to the appropriate decimal precision. For example, the numbers 0.6669 and 0.6661 are written as 0.666 in three-decimal precision. ◇

13B. **Main questions.** One might ask a naive question:

"How effective are the quantum invariants as tools for distinguishing knots?"

Such a focus underestimates the rich interplay of ideas and deep insights to which these invariants contribute. However, it is precisely this question that forms the focus of this section. It turns out that this may be the wrong question to ask: ***quantum invariants are not expected to be strong knot invariants***. (This seems to be folk knowledge; see, for example, [**Sto04**, **Wes12**] or Theorem 13E.5.) Indeed, they satisfy local relations (*e.g.* skein relations), which make them easy to study but not particularly well suited as invariants.

With this in mind, a better question is:

"How do quantum invariants compare in distinguishing knots?"

This question, along with its variations, forms the central focus of this section.

Two final comments before we start:

(i) The computational complexity of (most) quantum invariants is not polynomial in the number of crossings, but can vary quite a bit; see Section 13D.2 for a more detailed analysis. We will mostly ignore this factor in the comparisons.

(ii) A critical aspect of this study is the scale of the data, *e.g.* in small datasets, there's a risk of overfitting. For us there is no issue as the number of knots grows exponentially in the number of crossings. For example, for $k_n = \#\{\text{prime knots with } n \text{ crossings}\}$, [**ES87**, **Wel92**] proved that $2.68 \leq \liminf_{n\to\infty} \sqrt[n]{k_n} \leq \limsup_{n\to\infty} \sqrt[n]{k_n} \leq 13.5$.

13C. **Background in a nutshell.** We assume that the reader is familiar with the basics of knot theory, who may refer to *e.g.* [**Ada94**].

13C.1. *Comments on the nomenclature.* The quantum invariants that we have seen explicitly so far are the Jones polynomial and its colored versions, see for example Section 11H. Here are a few more, which we will not explain in detail.

**Notation 13C.1.** The quantum invariants of interest in this section are called, for short, ***A2***, ***Alexander*** (A), ***B1***, ***Jones*** (J), ***Khovanov*** (K). All of these are Laurent polynomials in $\mathbb{Z}[q, q^{-1}]$ or $\mathbb{Z}_{\geq 0}[q, q^{-1}, t, t^{-1}]$ (with $q$ and $t$ formal variables), but we, as often in the literature, just call them polynomials. The $t$ is such that the Khovanov homology specialized at $t = -1$ gives the Jones polynomial. We also have the $t = 1$ specialization of Khovanov homology, called ***KhovanovT1*** (KT1), which lives in $\mathbb{Z}_{\geq 0}[q, q^{-1}]$. ***Dequantization*** is the process of setting $q = 1$, which makes all these invariants trivial. ◇

This is how they arise in the notation from the previous sections:

- ▶ We discussed J and B1 in Section 11H, the former one being the invariant for color $c = 1$, and B1 is the invariant for $c = 2$.
- ▶ The invariant A2 can be constructed similarly to J in Section 11H, but using $n = 3$ and $c = 1$.
- ▶ The invariant A arises similarly to J, but with different webs. See [**Sar15**] for details.
- ▶ We will not explain how K is constructed, see [**BN02**] for a nice summery. Take it as a black box; all we need to know is that is is a two variable polynomial $K(q, t)$ such that $J(q) = K(q, -1)$.

All of our main invariants, and all quantum invariants, are associated with a root system (a Dynkin diagram) with simple roots $\alpha_1, ..., \alpha_r$ (nodes) and a dominant weight $a_1\alpha_1 + ... + a_i\alpha_i + ... + a_r\alpha_r$ in $\mathbb{Z}_{\geq 0}\alpha_1 + ... + \mathbb{Z}_{\geq 0}\alpha_r$ (tuples $(a_1, ..., a_r)$). For example, J comes from the root system of type A1 and the dominant weight is just the unique simple root itself. In this notation:

**(a)** ***Coloring*** is the process of replacing $(a_1, ..., a_i, ..., a_r)$ with $(a_1, ..., a_i + 1, ..., a_r)$.

**(b)** ***Rank increase*** is the process of replacing $r$ by $r + 1$ without leaving the type, whenever applicable.

**(c)** ***Categorification*** is the process of adding a superscript, *e.g.* the A1 invariant becomes the A1$^c$ invariant. The mathematics is explained in [**Web17**].

**(d)** ***Leaving the realm of Lie algebras*** is ill-defined. One possibility is taking a Dynkin diagram for a Lie superalgebra as in, for example, [**Kac77**, Proposition 2.5.6].

*Remark* 13C.2. Categorification, for us, includes Khovanov homology but not, for example, ***odd Khovanov homology*** (an alternative categorification of the A1 invariant, see [**ORS13**]) or ***HFK*** (a categorification of the isotropic A1 invariant, see [**OS04**, **Ras03**]). However, our observations hold for odd Khovanov homology as well, see Section 13E, and also most likely for HFK. We expect the overall behavior to be essentially independent of the precise form of categorification that one studies. ◇

In more detail, let us focus on the case where $\mathfrak{g}$ is a simple Lie algebra. Most of the following also work for Lie superalgebras and other objects, but we will not use them much in this paper.

Let $V$ be a simple highest weight representation of $\mathfrak{g}$. The Lie algebra $\mathfrak{g}$ is of ***classical type*** if its Weyl group is of type A, B, C or D. For these, increasing the rank means going to the simple Lie algebra $\mathfrak{h}$ of the same type by adding a vertex as follows.

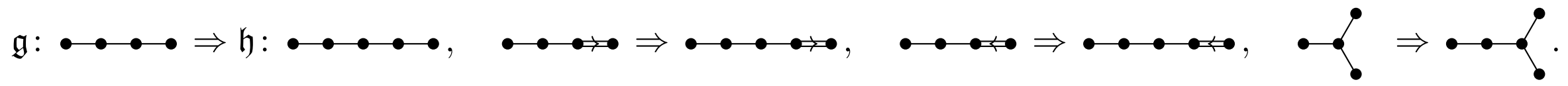

Here, $\mathfrak{g}$ is always on the left, $\mathfrak{h}$ on the right. This induces an injection of highest weights from $\mathfrak{g}$ to $\mathfrak{h}$ by $(a_1, ..., a_r) \mapsto (0, a_1, ..., a_r)$, and it makes sense to take the representation $V'$ associated to the image of $(a_1, ..., a_r)$ for $\mathfrak{h}$. Increasing the rank (rank for short) for an associated $Q_{\mathfrak{g},V,\epsilon}$ means going to $Q_{\mathfrak{h},V',\epsilon}$. Coloring (color for short) means going from $Q_{\mathfrak{g},V,\epsilon}$ to $Q_{\mathfrak{g},W,\epsilon}$, where $V$ is of highest weight $(a_1, ..., a_i, ..., a_r)$ and $W$ of highest weight $(a_1, ..., a_i + 1, ..., a_r)$, and categorification (cat for short) means going from $\epsilon = 0$ to $\epsilon = 1$.

*Remark* 13C.3. Note that $(\mathfrak{sl}_2, \mathrm{Sym}^2\mathbb{C}^2, 0)$, corresponding to $c = 2$ in Section 11H, is the quantum invariant associated with $SO_3$, see *e.g.*, [**Tub23**, Lemma 5A.7] for an explanation. However, we count this as a coloring instead of changing the Lie type. Still, we will call it the B1 invariant. ◇

*Remark* 13C.4. The root system of type $\mathfrak{gl}_{1|1}$ consists of one isotropic root (isotropic roots can only appear for Lie superalgebras), hence the name isotropic A1 invariant. The Alexander polynomial is thus a "special" quantum invariant when compared to the others, but it can also be constructed from an R-matrix; see, for example, [**Sar15**] for a nice summary. ◇

Since quantum invariants behave well with respect to direct sums of representations, restricting to simple $\mathfrak{g}$-representations is a weak restriction.

13C.2. *Prime knots and PD presentations.* We used a ***list of prime knots*** up to 16 crossings that can be found in [**LM25**, **TZ25b**]. The first few are listed in Figure 34 .

From now on, let $n$ be the ***number of crossings***.

*Remark* 13C.5. The sequences $(k_n)_{n=3}^{\infty}$ and $(\sum_{i=3}^{n} k_i)_{n=3}^{\infty}$ are

1,1,2,3,7,21,49,165,552,2176,9988,46972,253293,1388705 , for # of crossings from 3 to 16,

1,2,4,7,14,35,84,249,801,2977,12965,59937,313230,1701935 , for # of crossings from $\leq 3$ to $\leq 16$,

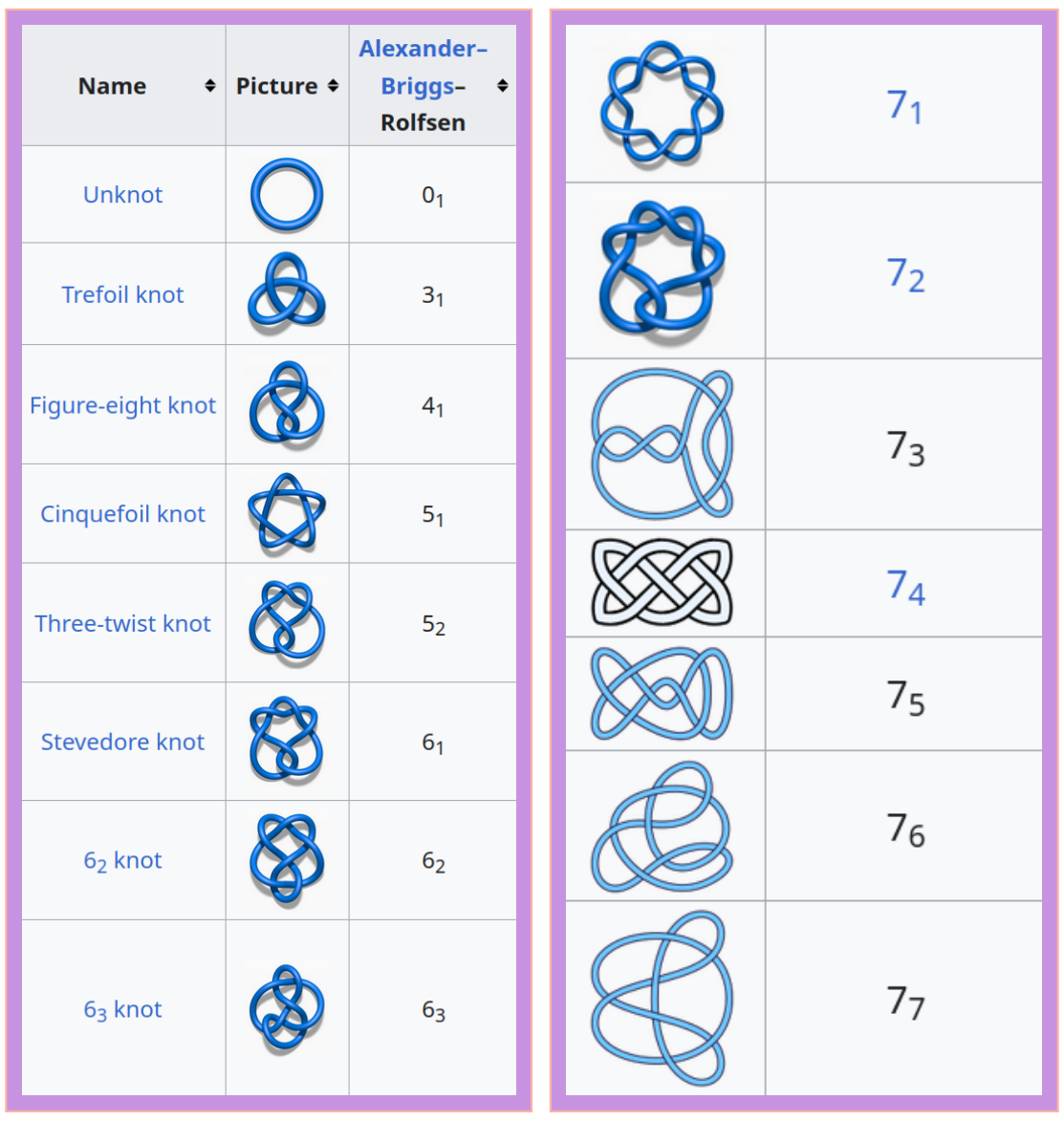

| Name | Picture | Alexander–Briggs–Rolfsen |
|---|---|---|
| Unknot | | $0_1$ |
| Trefoil knot | | $3_1$ |
| Figure-eight knot | | $4_1$ |
| Cinquefoil knot | | $5_1$ |
| Three-twist knot | | $5_2$ |
| Stevedore knot | | $6_1$ |
| $6_2$ knot | | $6_2$ |
| $6_3$ knot | | $6_3$ |

| | |
|---|---|
| | $7_1$ |
| | $7_2$ |
| | $7_3$ |
| | $7_4$ |
| | $7_5$ |
| | $7_6$ |
| | $7_7$ |

FIGURE 34. List of the first few prime knots (not including mirror images).
Picture from https://en.wikipedia.org/wiki/List_of_prime_knots.

see [**OEI23**, A002863]. Starting with the trefoil, every knot is listed explicitly. This is also the list we used for the big data comparison, *i.e.* we do not distinguish a knot from its mirror. ◇

Let us summarize the conventions.

**(a)** The list contains ***planar diagrams*** (PD for short) for each prime knot. Recall that a PD presentation of a knot diagram labels all of its edges with non-repeating numbers $\{1, ..., r\}$, where $r$ is the number of edges. Each edge is then adjacent to two crossings, which induces a labeling of the crossings. We remember the crossings as symbols $X[i, j, k, l]$, where $i$, $j$, $k$ and $l$ are the labels of the edges around that crossing, starting from the incoming lower edge and proceeding counterclockwise. Explicitly, the first two knots on the list and their PD presentations are:

$$PD[X[1,5,2,4], X[3,1,4,6], X[5,3,6,2]] \rightsquigarrow$$

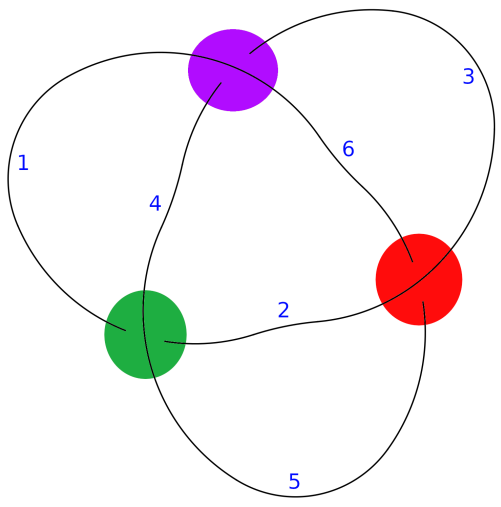

,

$$PD[X[4,2,5,1], X[8,6,1,5], X[6,3,7,4], X[2,7,3,8]] \rightsquigarrow$$

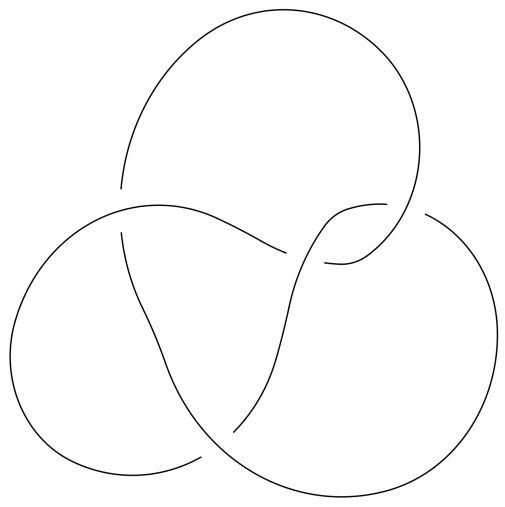

.

We only added the labels and highlighted the crossings for the ***trefoil*** knot at the top. The bottom knot is the ***figure eight knot***.

**(b)** In particular, the mirror image of the trefoil is not on the list. In general, the list identifies knots and their mirror images and has only the PD presentation for one of them.

**(c)** Note that *e.g.*, the Jones polynomial $J$ satisfies $J(K\#L) = J(K)J(L)$, so in some sense the running over prime knots suffices for our purposes.

13D. **Computational aspects.** We now provide a rough overview of the computational complexity associated with quantum invariants.

13D.1. *Computation of quantum invariants.* We now explain the different computational techniques that we have used to compute the quantum invariants.

---

*For A2, A, B1, J; used for A2, B1, J.* The first method is using ***R-matrices***.

We work over $\mathbb{C}(q)$ for a formal variable $q$ that potentially also has roots like $q^{1/k}$. The uncategorified quantum invariants $Q = Q_{\mathfrak{g},V,0}$ that we use, except the Alexander polynomial, can be defined and computed as follows. For a knot, fix a Morse presentation of the knot, arranged vertically.

For simplicity, assume that $V$ is self-dual, which works for B1 and J. In the Morse presentation, we have four basic pieces and an identity, that we name as

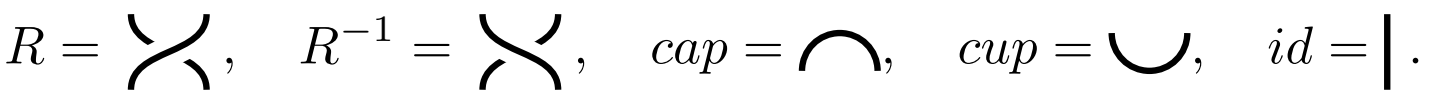

We associate these to linear maps (matrices upon choice of basis) denoted with the same symbols

(13D-1) $R, R^{-1}\colon V_q \otimes V_q \to V_q \otimes V_q, \quad cap\colon V_q \otimes V_q \to \mathbb{C}(q), \quad cup\colon \mathbb{C}(q) \to V_q \otimes V_q, \quad id\colon V_q \to V_q, v \mapsto v.$

Here, $V_q$ is a representation of a quantum group associated with $\mathfrak{g}$ that dequantizes to the $\mathfrak{g}$ representation $V$, and all the above maps are intertwiners for the quantum group that dequantize to, respective to above, the flip maps, the pairing and coparing of $\mathfrak{g}$ representations, and the identity.

Horizontal composition is then the tensor product (Kronecker product upon choice of basis). Write $id \otimes ... \otimes id$ with $k$ factors simply as $k$. Then for the figure eight knot this is could be

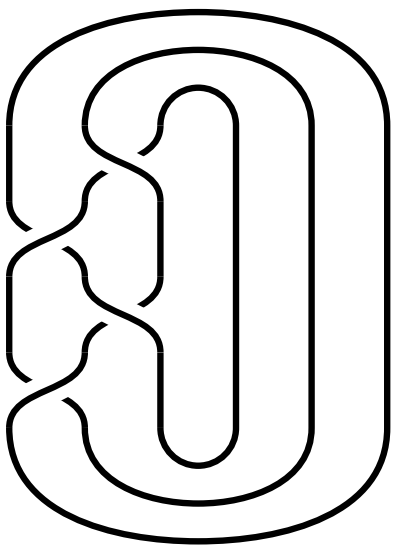

$$\rightsquigarrow \begin{array}{c} 0 \otimes cap \otimes 0 \\ \circ \\ 1 \otimes cap \otimes 1 \\ \circ \\ 2 \otimes cap \otimes 2 \\ \circ \\ 1 \otimes R^{-1} \otimes 3 \\ \circ \\ 0 \otimes R \otimes 4 \\ \circ \\ 1 \otimes R^{-1} \otimes 3 \\ \circ \\ 0 \otimes R \otimes 4 \\ \circ \\ 2 \otimes cup \otimes 2 \\ \circ \\ 1 \otimes cup \otimes 1 \\ \circ \\ 0 \otimes cup \otimes 0 \end{array},$$

with composition, for example from bottom to top.

The maps in Equation 13D-1 are not uniquely determined from what we wrote above. There are some choices involved, but they are mostly irrelevant and just rescale the quantum invariant. To be completely explicit, we used the conventions used in [Atl24]. For example, in coordinates,

$$R = \begin{pmatrix} q^{1/2} & 0 & 0 & 0 \\ 0 & 0 & q & 0 \\ 0 & q & q^{1/2} - q^{3/2} & 0 \\ 0 & 0 & 0 & q^{1/2} \end{pmatrix},$$

is the choice of $R$-matrix for the Jones polynomial.

The Alexander polynomial can be computed using R-matrices, with the slight catch that one needs to renormalize the result appropriately as one gets zero, *cf.* [Sar15], which is a common phenomena when working with isotropic roots.

---

*For A2, A, B1, J, K, KT1; used for B1, K, KT1.* Another way of computing quantum invariants is a ***skein theory approach***.

For Jones this, for example, uses the relations:

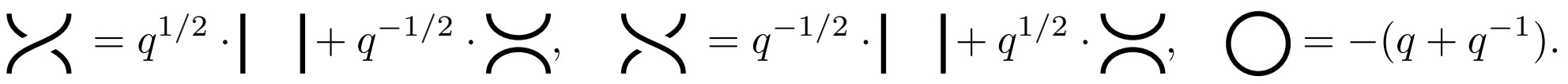

One can then replace all crossings, as usual, and get the desired polynomial. We exemplify this here for the Hopf link:

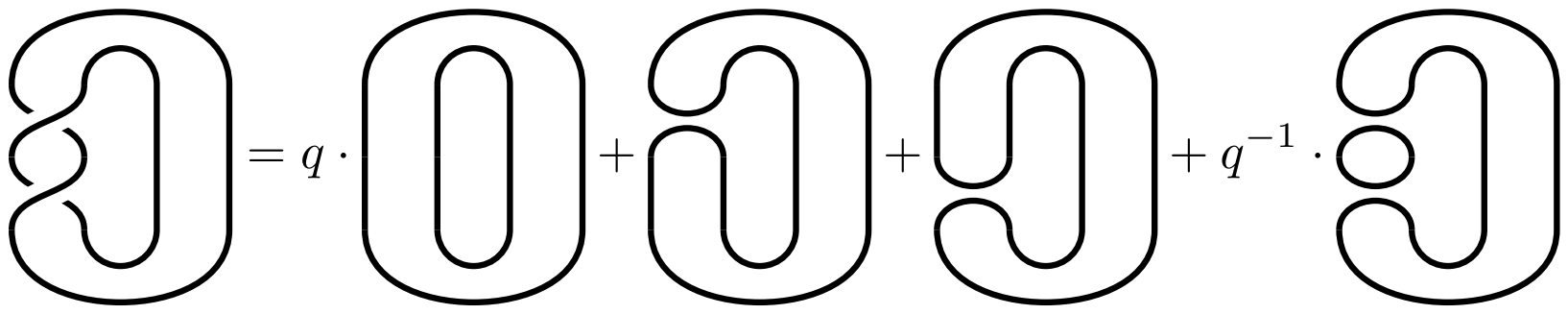

$$= q(q+q^{-1})^2 - 2(q+q^{-1}) + q^{-1}(q+q^{-1})^2 = q^3 + q + q^{-1} + q^{-3}.$$

We never used this for the decategorified quantum invariants, except for B1. We elaborate in [**TZ25b**], for this paper it is enough to know that for B1 we followed the conventions in [**Tub23**, Section 5].

We do not know a general R-matrix strategy to compute categorified quantum invariants. For Khovanov we used [**Atl24**], which is based on [**BN07**] and uses the conventions from [**BN02**]. In particular, it uses a skein theory approach where every crossing is a complex, in a certain category, of the form for a certain differential $d$ and underlined part in homological degree zero. The whole construction is then tensored together, as complexes, over the crossings. The rest of the calculation is similar to the Jones calculation, just in complexes.

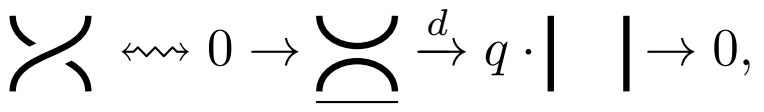

---

*For A; used for A.* The following is the ***determinant approach*** which generally does not work for quantum invariants. It inspired "pseudo quantum invariants" (our choice of name) as in [**BNvdV19**].

The arguably best approach to computing the A invariant is to use the ***Seifert matrix*** $S(K)$ of a knot (which in turn collects the linking numbers of the Seifert circles). Precisely,

$$A(K) = \det(S(K) - qS(K)^T).$$

The example to keep in mind is the trefoil, where the Seifert matrix is $\begin{pmatrix} 1 & 0 \\ -1 & 1 \end{pmatrix}$. Thus:

$$A(PD[X[1,5,2,4], X[3,1,4,6], X[5,3,6,2]]) = \det\left(\begin{pmatrix} 1 & 0 \\ -1 & 1 \end{pmatrix} - q\begin{pmatrix} 1 & -1 \\ 0 & 1 \end{pmatrix}\right) = \det\left(\begin{pmatrix} 1-q & q \\ -1 & 1-q \end{pmatrix}\right) = 1 - q + q^2.$$

This only determines A up to a scalar, and we scale it so that the result is the same when one would use a skein approach, as in [**Atl24**].

13D.2. *Algorithmic complexity.* We will now give a rough estimate of the complexity of the above algorithms. We go through the above list in slightly changed order.

*Remark* 13D.2. We are analyzing specific algorithms, not the problem of computing quantum invariants themselves, which might be very different. We also do not try to give the best upper bounds possible. ◇

---

*Determinant approach.* In this scenario the Seifert matrix $S(K)$ is of size $2g$-by-$2g$, where $g$ is the genus of the diagram of the knot $K$. Note that $g \leq n/2$, so that the Seifert matrix is at most of size $n$-by-$n$. The main computational complexity is then the computation of the determinant, which, using LU decomposition, is in $O(n^3)$. Thus:

$$A \in O(n^3).$$

This is very cheap compared to what is up next.

---

*R-matrix approach.* Let $N = \dim_{\mathbb{C}} V$, and let $p(n)$ denote some polynomial in $n$, and let $m$ be the carving width of the four-valent planar graph associated with a knot. With a bit of care, see [**Mar21**, Theorem 1.1], one can show that $m \in O(\sqrt{n})$ and then

$$Q_{\mathfrak{g},V,0} \in O(p(n)N^{3m/2}) = O(p(n)N^{3\sqrt{n}/2}).$$

Thus, computing quantum invariants using $R$-matrices is superpolynomial in $n$, with the leading factor determined by the dimension of the underlying representation. This is not surprising if one keeps in mind that the main difficulty in this computational approach is not the number of crossings, but rather the number of strands, since this corresponds to tensoring $V$, which, even in decompositions [**COT24**, **KST24**], behaves exponentially in $N$. The algorithm explained in [**Mar21**] is not quite the one used in KnotTheory [**Atl24**], but the one in KnotTheory should have roughly the same complexity (assuming that the algorithm that produces a Morse presentation is optimized; which is not the case for the program used by [**Atl24**]).

*Remark* 13D.3. Finding a presentation of a knot with small carving width is surprisingly inexpensive, namely polynomial in $n$, see [**ST94**]. ◇

*Remark* 13D.4. It is known that computing the Jones polynomial is #P-hard, see, for example, [**Wel93**, Section 6.3]. Thus, under the standard assumptions that complexity classes do not collapse, the R-matrix approach is among the fastest possible. ◇

---

*Skein theory approach.* Let $M$ denote the number of summands in the expression of the crossing, *e.g.* $M = 2$ for the Jones polynomial or Khovanov homology and $M = 3$ for the B1 invariant. As before, let $p(n)$ denote some polynomial in $n$ and $s \in \mathbb{R}_{\geq 1}$ be some scalar. Then

$$Q_{\mathfrak{g},V,\epsilon} \in O(p(n)M^{sn}).$$

In particular, since, to the best of our knowledge, there is no R-matrix approach to Khovanov homology, the computation of it is exponential in the number of crossings.

*Remark* 13D.5. There is evidence that the computation of the categorification should be of the same complexity as the computation of the $\epsilon = 0$ version; see, *e.g.*, [**PS24**]. ◇

---

*Average runtime.* Arguably, the ***worst case*** analysis as above is not the right thing to consider. For example, with the calculations we run it seemed that there are few knots that took a very long time to compute, like 1000 times longer than the others with the same number of crossings. So, the correct measurement might be ***average runtime***.

We sadly do not know of any sources that computed the average runtime (and, indeed, this is usually much harder than the worst case). We here rather list our experimental data for the average runtime (the second plot is a ***box plot***) all of which were run with KnotTheory [**Atl24**] on the same machine: the laptop of the second author. Additionally, we include the data for the Javascript program for B1, run separately on UNSW's Katana servers, which uses the Skein theoretic approach as opposed to the R-matrix approach of the KnotTheory library. The details of these can be found in [**TZ25b**].

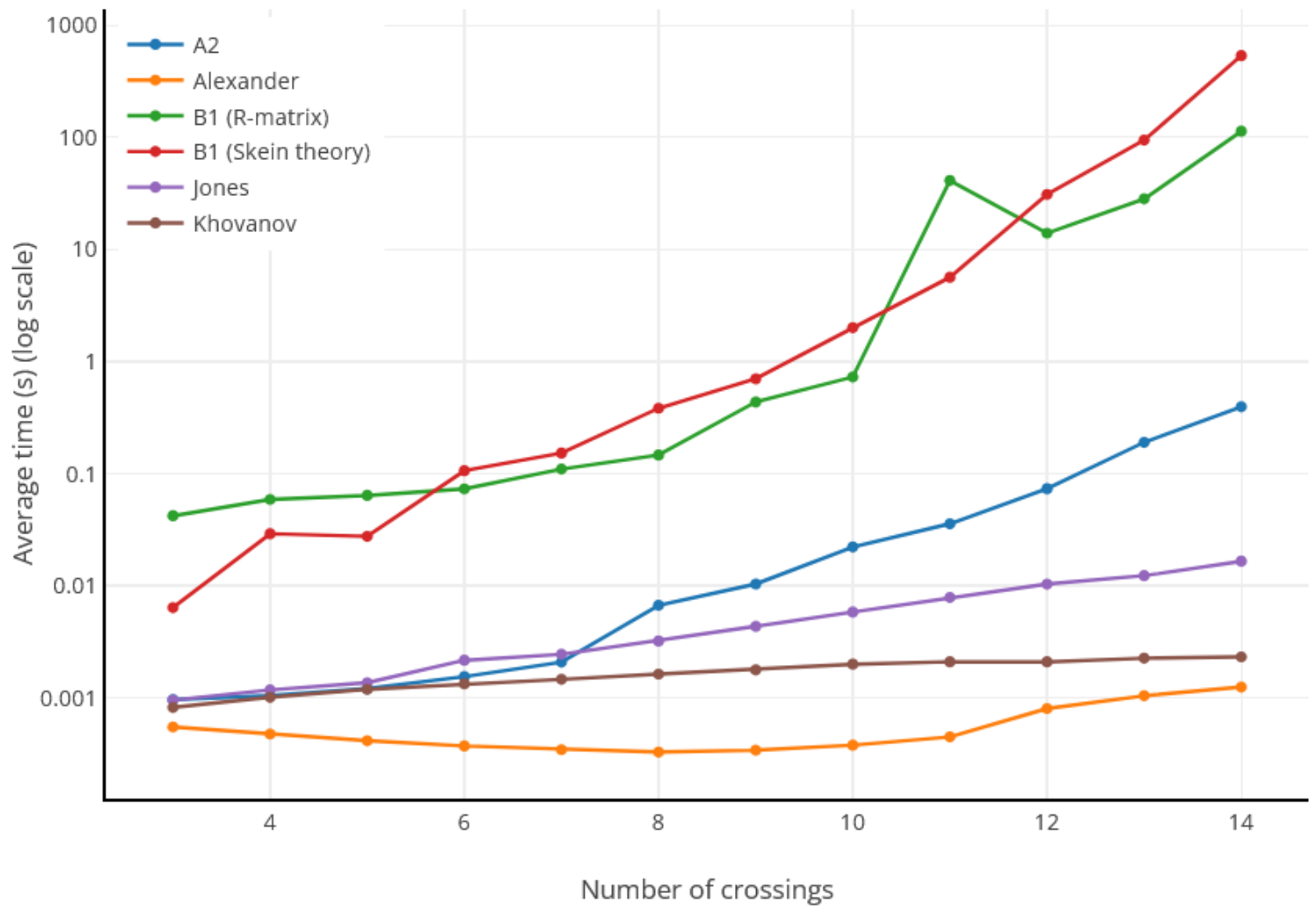

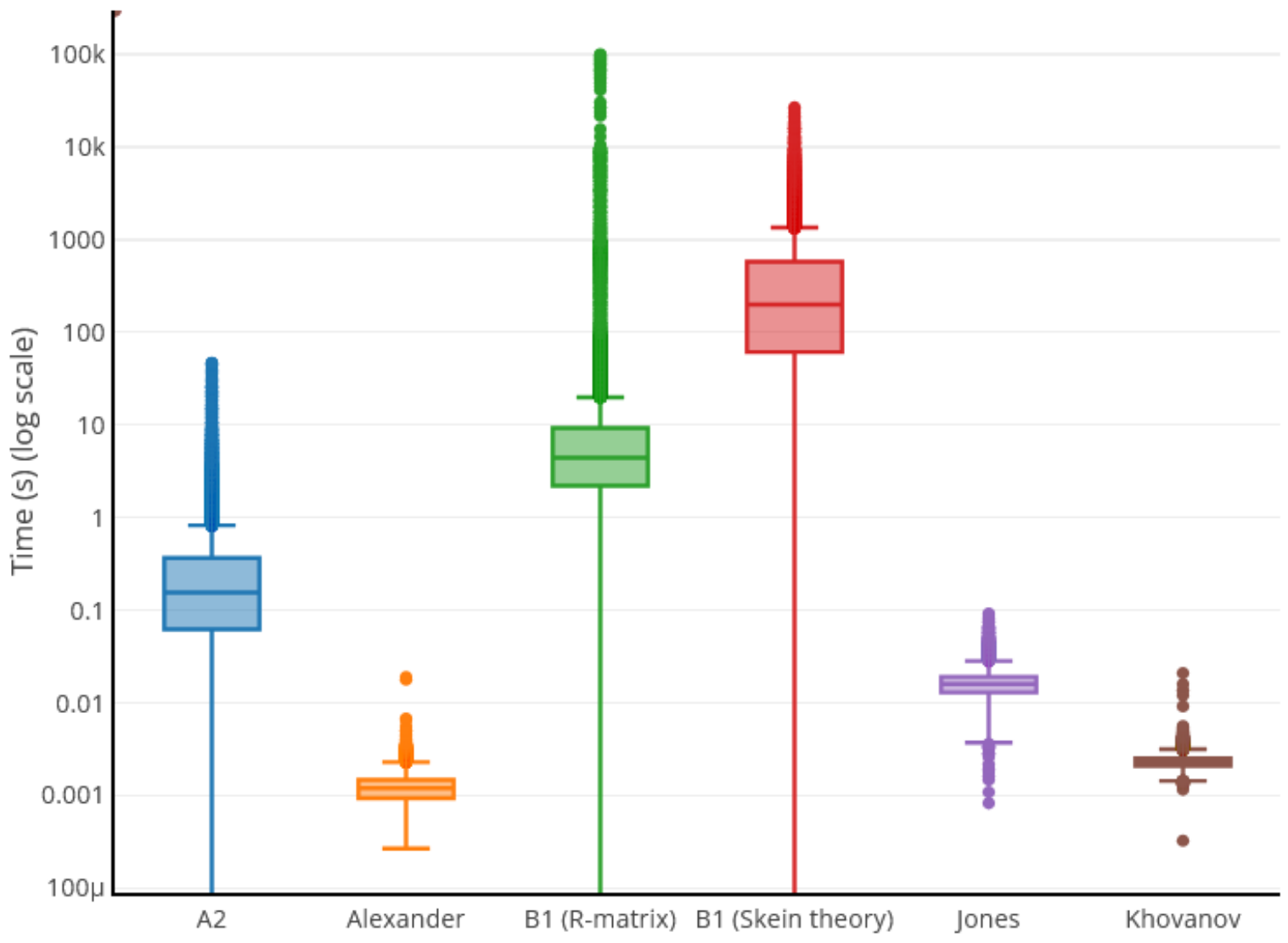

The precise numbers can be found in Figure 35.

We have three comments on this:

(i) The variance is quite high for A2 and B1. Since the runtime of the R-matrix approach has leading factor $N = \dim_{\mathbb{C}} V$, and $N = 3$ for these two invariants, the issue might be an not optimal conversion from a PD presentation to a Morse presentation for some knots. Indeed, there was knot number 775 on the list (11 crossings) where the computation took 21703.436821 seconds, much longer than any other knot. This is one reason we wrote a skein theory program for B1 (that can be found in [**TZ25b**]) and it typically did better on knots that were difficult in the R-matrix approach but was slower on average; see the data above.

(ii) The computation of A was so fast that the program was still loading when the first few knots were computed. Hence, the drop below $n = 8$.

(iii) The computation of K is not run in Mathematica itself, but instead starts a Java program. This could explain why K was computed faster than J.

13E. **Comparison – distinct values.** For the first comparison, we are interested in how many distinct values quantum invariants take on our list of knots, and we measure this as a percentage. In other words, we want to know

$$Q(n)^{\%} = \#\{Q(K) \mid K \in \mathcal{K}_n\}/\#\mathcal{K}_n.$$

These are the ***distinct values*** $Q$ takes.

Recall that A, J, and K stand for Alexander, Jones, and Khovanov, respectively. Let "All" mean that we take all of them together, and we use J+KT1 to take Jones and KT1 together. The data is as follows.

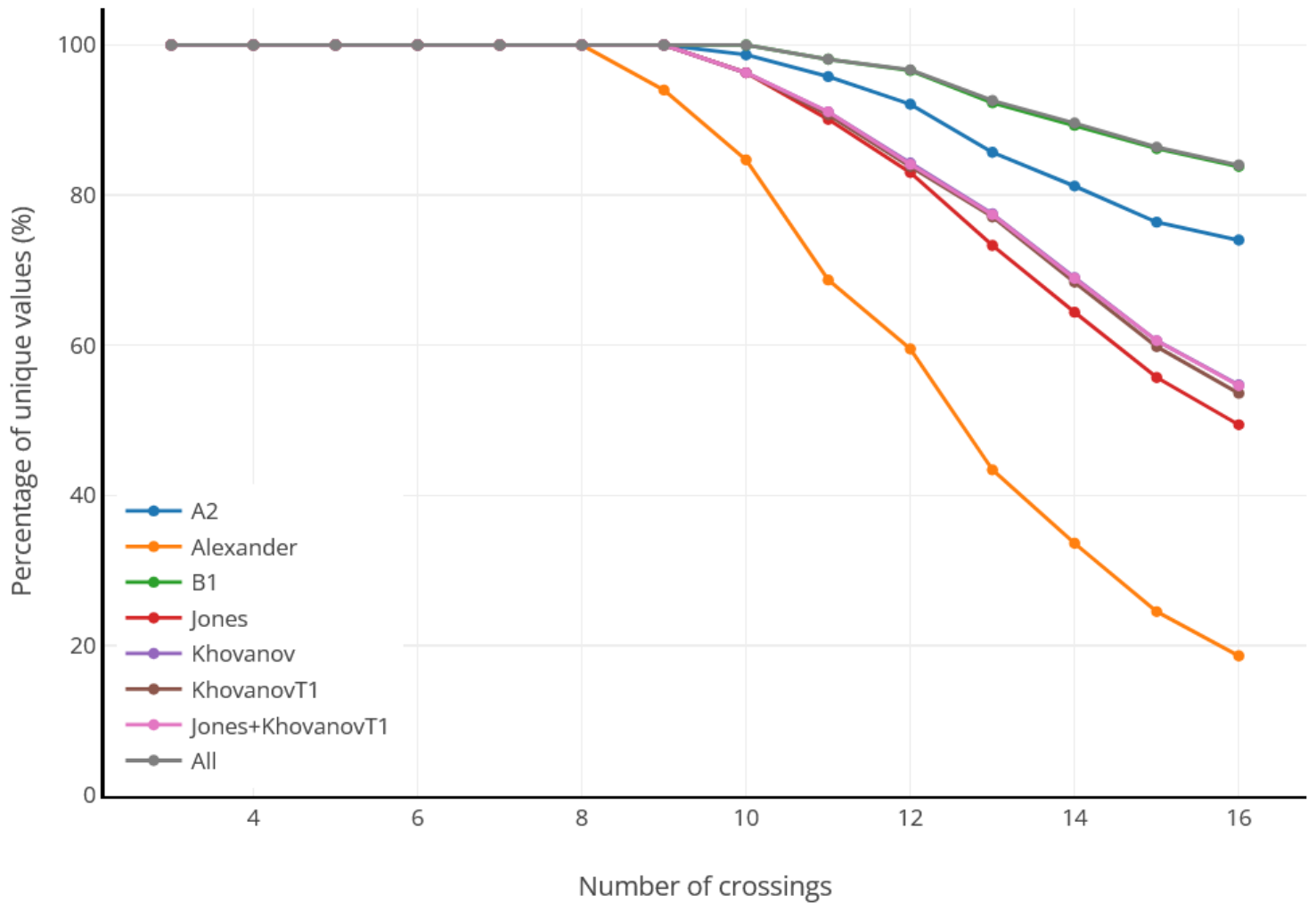

,

with the precise values listed in Figure 36. Here are a few observations.

**(a)** From the data, ***four different classes*** are visible: All/B1, A2, J/K/KT1, and A.

**(b)** They all seem to drop to zero, but the rate of convergence to zero seems to depend on the class.

**(c)** There is essentially no difference between K and J+KT1, so categorification seems to be as strong as two evaluations of it. Additionally, there is no huge difference between K and J or KT1.

**(d)** There is essentially no difference between All and B1; in particular, coloring seems to be the preferred method to distinguish knots.

In general, some ***zero-one law*** should apply, so $Q(n)^{\%}$ should either converge to zero or one. The data suggest that the limit for all of them is zero. From a more detailed analysis of the data we even conclude an ***exponential decay***.

**Conjecture 13E.1.** For $Q \in \{A2, A, B1, K\}$ (and therefore also for $Q = J$ or $Q = KT1$) we have

$$Q(n)^{\%} \in O(\delta^n) \text{ for some } \delta = \delta(Q) \in (0,1).$$

Moreover, we have

$$\delta(B1) > \delta(A2) > \delta(K) > \delta(A).$$

In other words, they all drop to zero exponentially fast, but the exponential factor depends on the class. ◇

**Speculation 13E.2.** For all quantum invariants $Q$ we have

$$Q(n)^{\%} \in O(\delta^n) \text{ for some } \delta = \delta(Q) \in (0,1).$$

In other words, they all drop to zero exponentially fast. ◇

**Speculation 13E.3.** Assume Speculation 13E.2 holds and fix $Q_{\mathfrak{g},V,0}$ for $\mathfrak{g}$ of classical type. Then:

$$\delta(\mathit{color}) > \delta(\mathit{rank}) > \delta(\mathit{cat}).$$

That is, coloring is better than rank increase, which is better than categorification. ◇

*Remark* 13E.4. Detecting certain knots, such as K detects the unknot [**KM11**], seems to be only possible if the knots are very special (unknot, torus knots, *etc.*). ◇

Factoring in computational complexity, one might argue that A is the best invariant. Furthermore, assuming that categorification is actually of the same complexity as its $\epsilon = 0$ counterpart, one might argue that $\mathit{cat} > \mathit{color} > \mathit{rank}$, since coloring and rank increase have a bigger exponential factor in their complexity analysis, *cf.* Section 13D.2. This however needs a deeper analysis of $\delta$ which is beyond this paper.

**Comments on how to potentially prove parts of Conjecture 13E.1.** The rest of the section is dedicated to prove a version of Conjecture 13E.1. For this section, a ***skein relation*** is a relation of the following form. Let $a, b, c$ be elements in a ring, then

$a \cdot$ 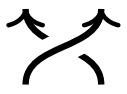 $+ b \cdot$ 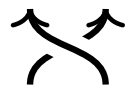 $+ c \cdot$ 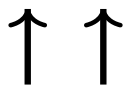 $= 0$, for $a, b$ invertible.

***Multiplicity freeness*** of a quantum invariant $Q_{\mathfrak{g},V,0}$ is the property that $V \otimes V$ is multiplicity free. Examples of quantum invariants that satisfy a skein relation or multiplicity freeness are A2, A, B1, and J.

Let $\mathcal{AL}_n$ be the set of alternating links of $\leq n$ crossings. Similarly to $Q(n)^{\%}$ define

$$Q(n)^{\%}_{AL} = \#\{Q(L) \mid L \in \mathcal{AL}_n\} / \#\mathcal{AL}_n.$$

These are the ***distinct values*** $Q$ takes on ***alternating links***. For the next statement, the reader may want to recall ***Conway mutation*** as, for example, in [**Ada94**, Section 2.3]. Let $K_2$ be Khovanov homology in characteristic 2.

**Theorem 13E.5** (Exponential decay theorem)**.** *For any quantum invariant $Q$ that satisfies a skein relation, is multiplicity free or does not detect Conway mutation we have*

$$Q(n)^{\%}_{AL} \in O(\delta^n) \text{ for some } \delta = \delta(Q) \in (0, 0.996).$$

*This applies to $Q \in \{A2, A, B1, J, K_2\}$.*

*Proof.* Any invariant that satisfies a skein relations or is multiplicity free, this invariant cannot now detect Conway mutation. This is pointed out in [**Weh03**, **Wes12**].

Now, let $\mathcal{ALCM}_n$ denote $\mathcal{AL}_n$ modulo Conway mutation. It is known that

$$\lim_{n\to\infty} \sqrt[n]{\#\mathcal{AL}_n} = \frac{\sqrt{21001} + 101}{40},$$

see [**ST98**, Theorem 1]. By [**Sto04**, Theorem 2], it is further known that

$$\limsup_{n\to\infty} \sqrt[n]{\#\mathcal{ALCM}_n} \leq \frac{\sqrt{24584873929} + 109417}{43334} < \frac{\sqrt{21001} + 101}{40} - 0.004.$$

Thus, the theorem follows for $\delta$ being 0.996 or smaller, and we just have to comment on the various special cases. Over $\mathbb{F}_2$ Khovanov homology fails to detect Conway mutation [**Weh10**], while A satisfies a skein relation, and A2, B1, and J are multiplicity free. □

The reader might have noticed that the proof of Theorem 13E.5 also works for variations of quantum invariants. Explicitly, the HOMFLYPT polynomial satisfies a skein relation, $k$-colored Jones polynomials, for all $k \in \mathbb{Z}_{\geq 0}$, are multiplicity free and odd Khovanov homology does not detect Conway mutation, so they all satisfy the same exponential decay.

*Remark* 13E.6. It seems likely that Theorem 13E.5 also works for generalized mutation as in [**APR89**] and other relations on the set of knots or links. Essentially, "any relation" imposed on the set of knots or links will make the quotient set an order of magnitude smaller, and quantum invariants will satisfy some relation. ◇

*Remark* 13E.7. Khovanov homology and actually KT1 do detect some Conway mutations [**Weh03**], but they still do poorly. So, one should expect that there is another form of "mutation" that is not detected by Khovanov homology and other categorification, and this mutation will eventually dominate the set of knots. In fact, Khovanov homology satisfies a skein relation with an error term [**CGKL21**], which is probably enough to force exponential decay. However, this would require some linearization (or similar) of the set-based proof of Theorem 13E.5. While we do not expect it to be difficult, we do anticipate it being somewhat tedious, which is why we have omitted it from this work. ◇

*Remark* 13E.8. A stronger condition than $V \otimes V$ being multiplicity free is $V$ being strongly multiplicity free, see [**LZ06**, Theorem 3.4] for a classification. Surprisingly, this stronger property, which implies Theorem 13E.5, is also related to faithful representations of the braid group, see [**LTV22**]. This is another example of why braids are considered easier than knots and links. ◇

13F. **Comparison – distinguishing pairs.** In Section 13E we asked how likely one can distinguish a knot from all others. An alternative way to measure how good an invariant is would be to asked how likely one can distinguish two knots from one another. That is, we define a function

$$Q(n)^{\%,\%} = E(\%,\%),$$

where $E(\%,\%)$ is the expected number of times it is necessary to randomly select $K, L \in \mathcal{K}_n$, for $K \neq L$, until $Q(K) = Q(L)$.

*Remark* 13F.1. For the following data we used JavaScript's `Math.random()` to select two distinct knots (pseudo) randomly from the list of all prime knots $\mathcal{K}_n$. In particular, each data point is the average over 100,000 trials and is an approximation of $Q(n)^{\%,\%}$. As noted in ECMAScript specification (2025), `Math.random()` uses an implementation-defined pseudo-random number generation algorithm and thus may differ between machines. By design, the generator cannot be seeded. ◇

With the same notation as in Section 13E, the collected data is as follows. The second plot is a magnification of the tail of the first plot, without a log scale, to accentuate potential classes of invariants.

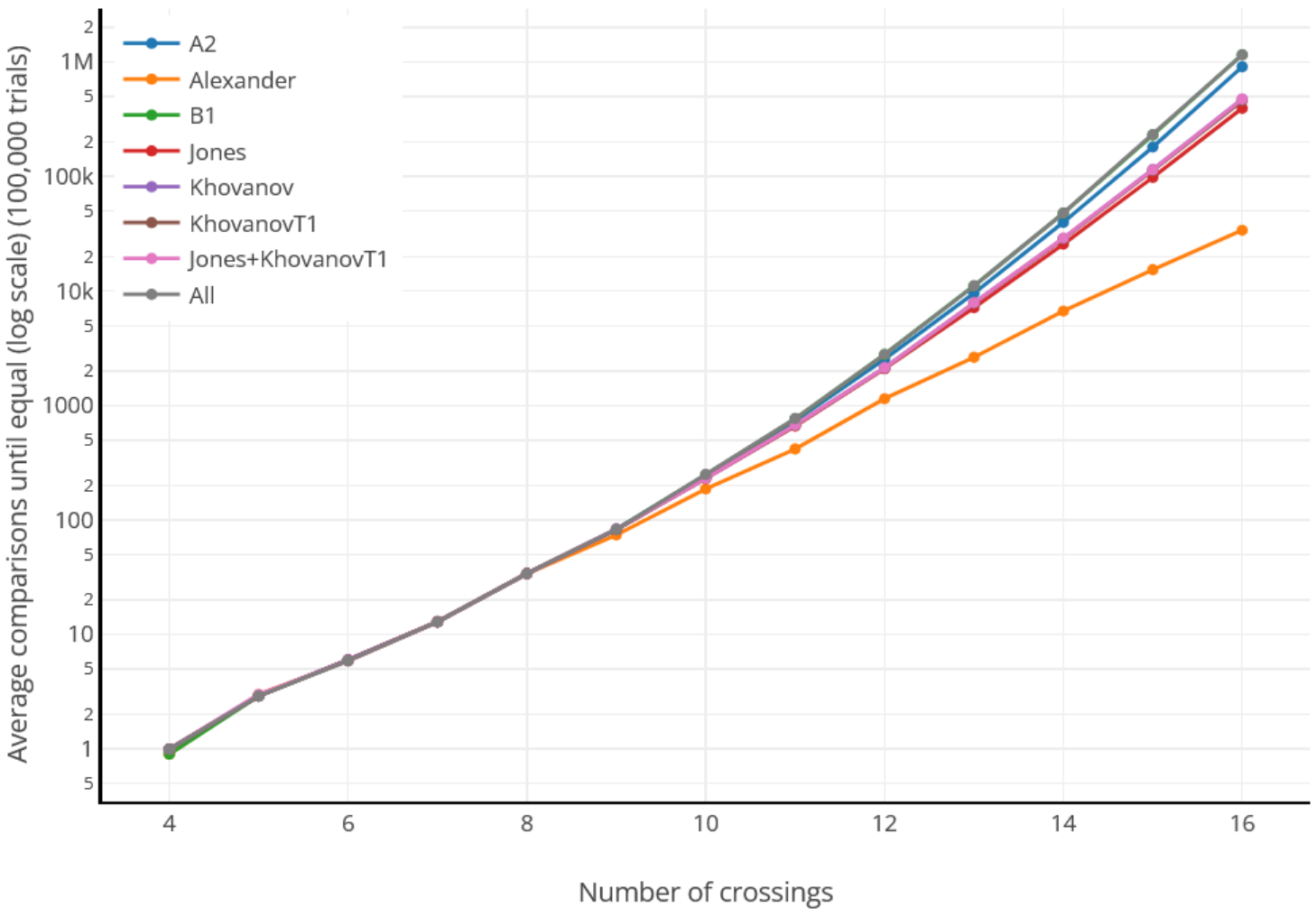

,

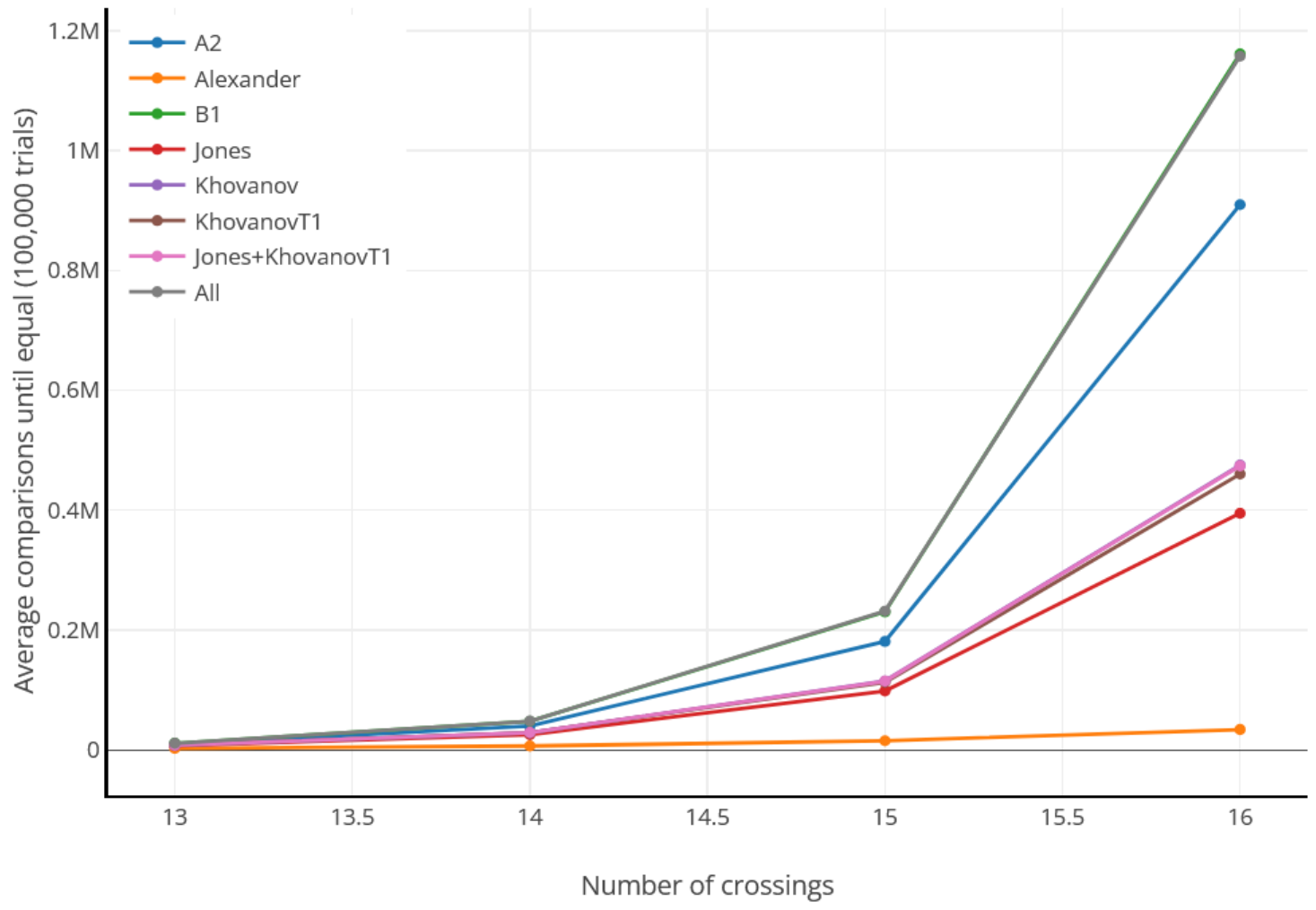

We have listed the numbers in Figure 37. Here are again a few observations:

**(a)** The only change compared to Section 13E is that now the probability of detection goes to one instead of zero.

**(b)** We stress that the preferred method again appears to be coloring.

As before in Section 13E, but slightly adjusted by changing to ***exponential growth***:

**Conjecture 13F.2.** For $Q \in \{A2, A, B1, J, KT1\}$ (and therefore also for $Q = K$) we have

$$Q(n)^{\%,\%} \in \Omega(\gamma^n) \text{ for some } \gamma = \gamma(Q) \in \mathbb{R}_{>1}.$$

Moreover, we have

$$\gamma(B1) > \gamma(A2) > \gamma(K) > \gamma(A).$$

In other words, they all grow exponentially fast, but the exponential factor depends on the class. ◇

We even speculate the following.

**Speculation 13F.3.** For all (nontrivial, e.g. $V$ is of dimension 2 or bigger) quantum invariants $Q$ we have

$$Q(n)^{\%,\%} \in \Omega(\gamma^n) \text{ for some } \gamma = \gamma(Q) \in \mathbb{R}_{>1}.$$

In other words, they all grow exponentially fast. ◇

**Speculation 13F.4.** Assume Speculation 13F.3 holds and fix $Q_{\mathfrak{g},V,0}$ for $\mathfrak{g}$ of classical type. Then:

$$\gamma(color) > \gamma(rank) > \gamma(cat).$$

That is, coloring is better than rank increase is better than categorification. ◇

*Remark* 13F.5. As in Section 13E, after cleaning for computational complexity, one might argue that the ordering in Speculation 13F.4 needs to be adjusted. ◇

In order to prove Conjecture 13F.2, or any of the others, it is likely beneficial to study random knots (see *e.g.* [**EZ17**] for a summary) and their values on quantum invariants. For example, the knots generated by the random model in [**EZHLN16**] seem to be prime with high probability, and it would be interesting to consider whether one can compute, for example, a Jones polynomial in this model.

13G. **Comparison – Big data approach.** The following is inspired by [**LTV24b**]. Fix a quantum invariant $Q = Q_{\mathfrak{g},V,\epsilon}$. Examples of questions one could try to address are the following:

**(a)** What is the asymptotic behavior, for $n \to \infty$, of

$$\mathrm{ev}_n = \mathtt{max}\{\textstyle\sum \mathtt{abs}\big(Q(K)\big) \mid K \in \mathcal{K}_n\}?$$

This corresponds to taking the sum of the coefficients. However, merely doing so would be dequantization, which is trivial, and so instead we take the sum of the absolute values `abs` of the coefficients.

**(b)** Similarly, what is the asymptotic behavior, for $n \to \infty$, of

$$\mathrm{coeff}_n = \mathtt{max}\left\{\mathtt{max}\left(\mathtt{abs}\left(Q(K)\right)\right) \mid K \in \mathcal{K}_n\right\}?$$

That is, we ask how fast the (absolute values of the) coefficients grow.

**(c)** We can ask the same questions as in the previous two points, but for the ***average*** instead of the maximum. That is,

$$\mathrm{ev}_n^{av} = \sum\{\sum \mathtt{abs}\,\big(Q(K)\big) \mid K \in \mathcal{K}_n\}/\#\{Q(K) \mid K \in \mathcal{K}_n\}$$

where we take the average over sums, and

$$\mathrm{coeff}_n^{av} = \mathtt{max}\,\big\{\sum \mathtt{abs}\,\big(Q(K)\big)/\,\mathtt{spread}\,\big(Q(K)\big) \mid K \in \mathcal{K}_n\big\}.$$

where we take the maximum over "average coefficients". Note that for $\mathrm{coeff}_n^{av}$ we do not count padded zeros.

**(d)** One could also ask about maximal and average ***spread*** (maximal minus minimal degree), which we denote as $\mathrm{spread}_n^{av}$ and $\mathrm{spread}_n$.

**(e)** Note that it does not make sense to ask for a spread of a two-variable polynomial, so invariant K is excluded from data that involves spreads or average coefficients.

Now we show the data. First, $\mathrm{coeff}_n$, the maximal coefficient.

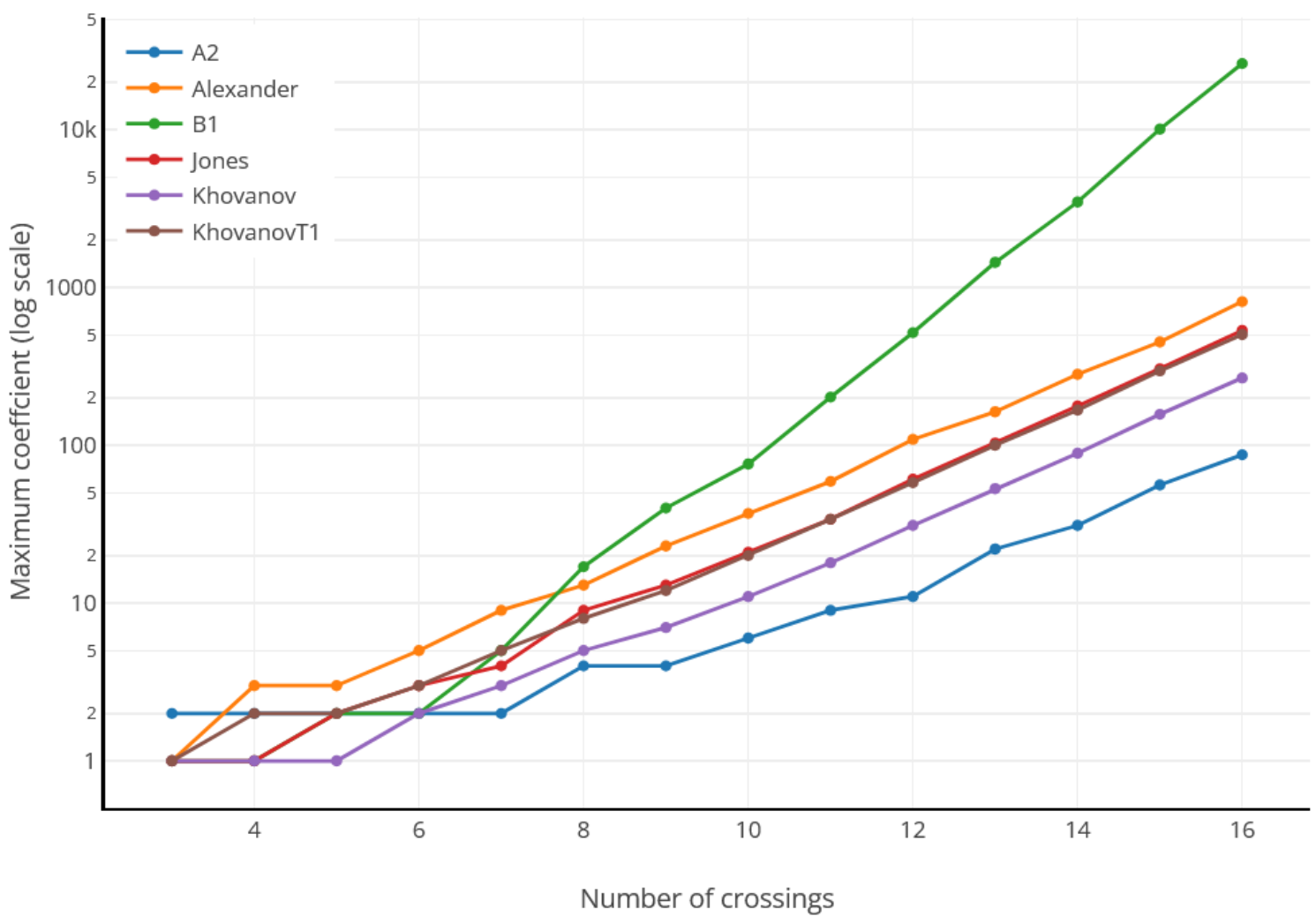

,

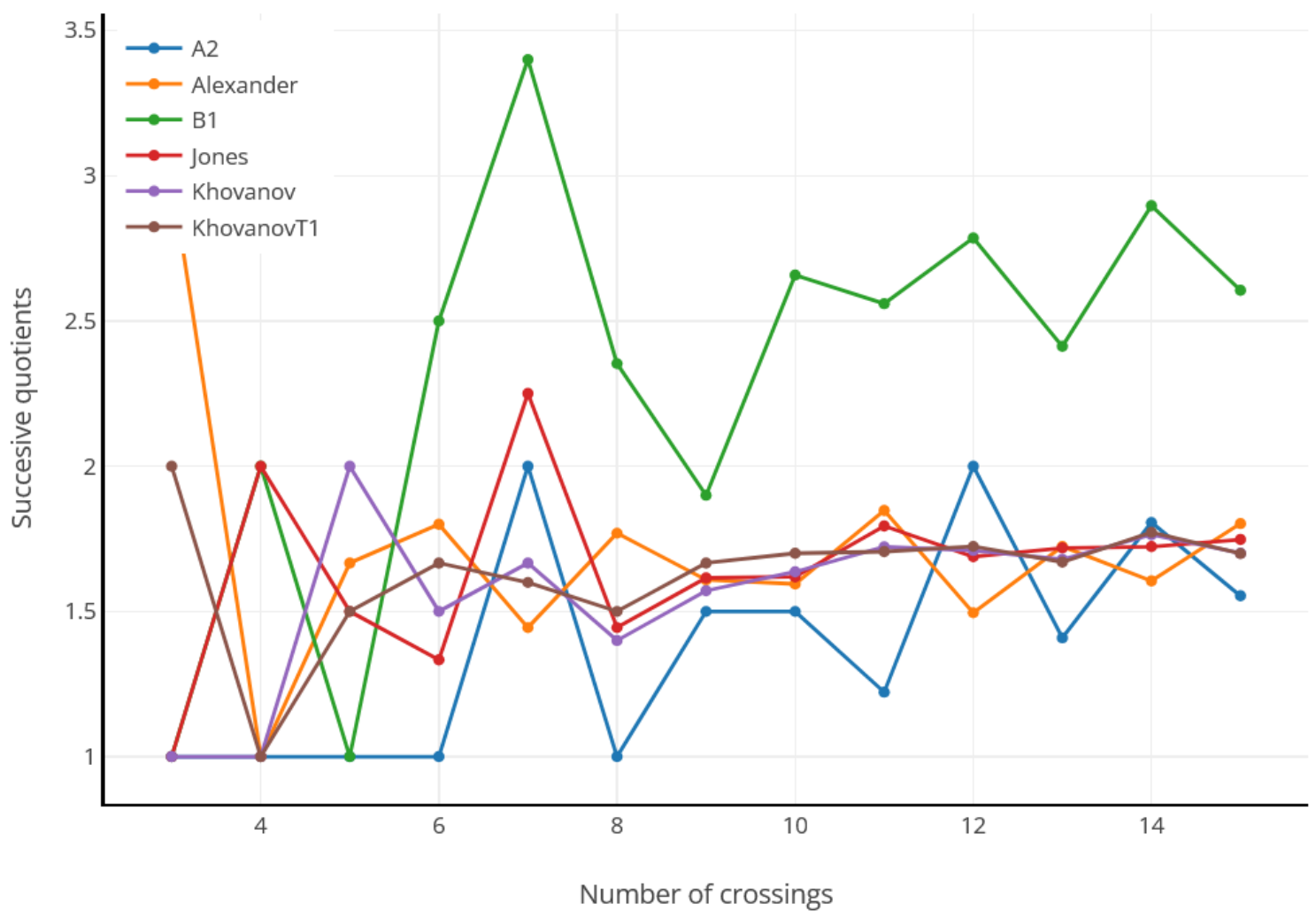

,

and the precise figures can be found in Figure 38. Next, $\mathrm{ev}_n$, the maximal coefficient sum.

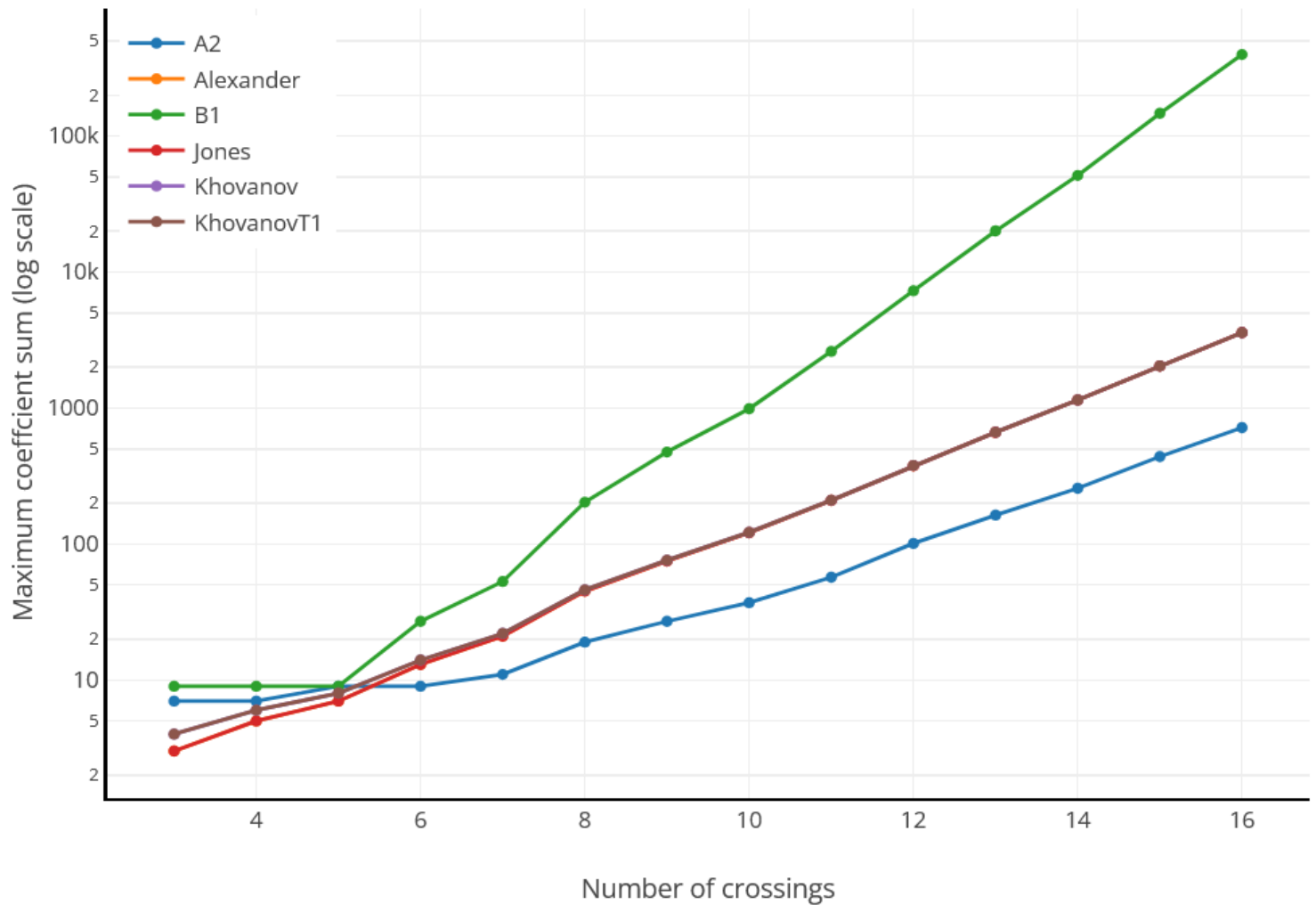

,

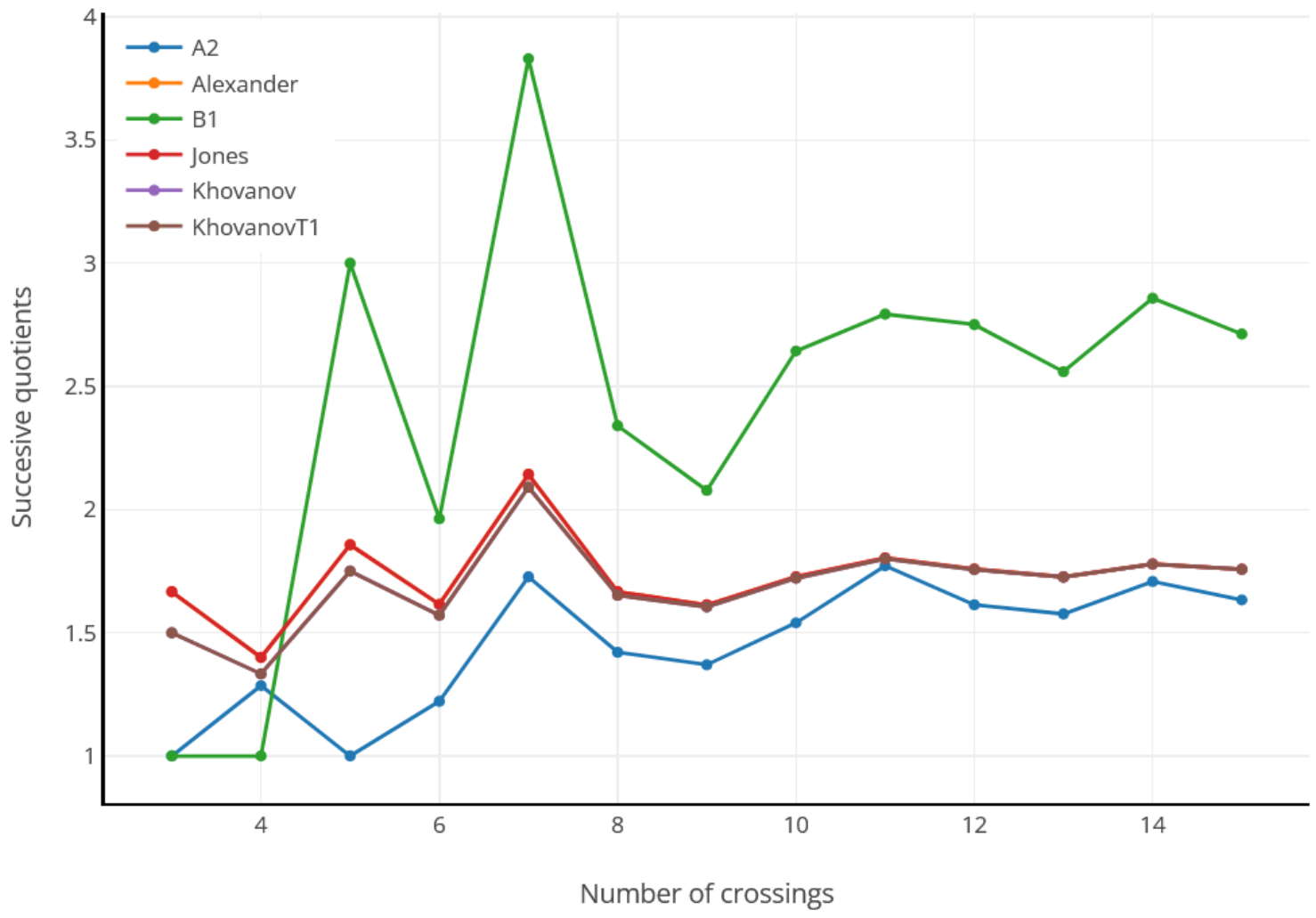

,

with data in Figure 39. Then, taking averages, $\mathrm{ev}_n^{av}$ and $\mathrm{coeff}_n^{av}$.

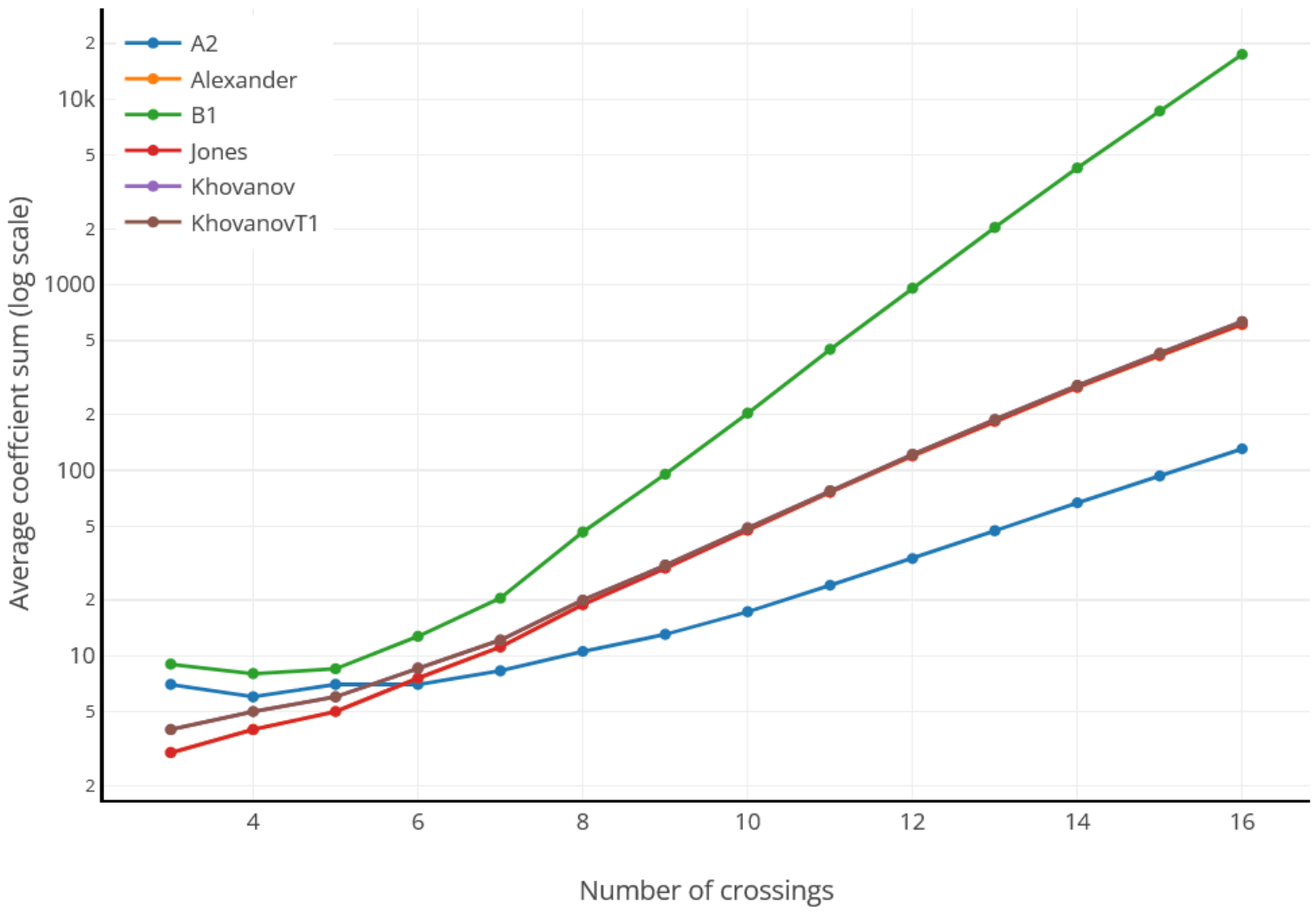

,

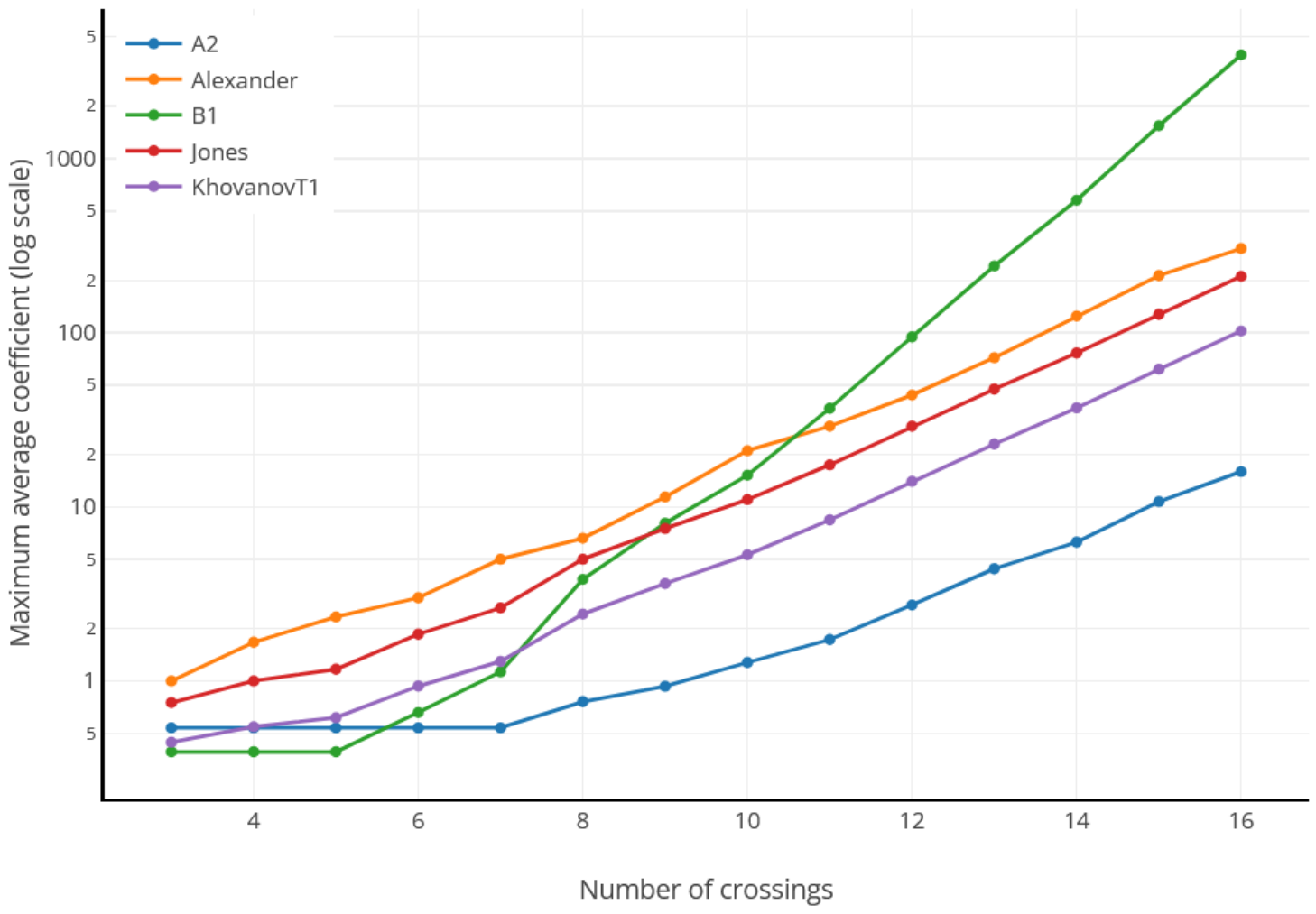

,

with data in Figure 40 and Figure 41. Then we look at spreads, $\mathrm{spread}_n$ and $\mathrm{spread}_n^{av}$.

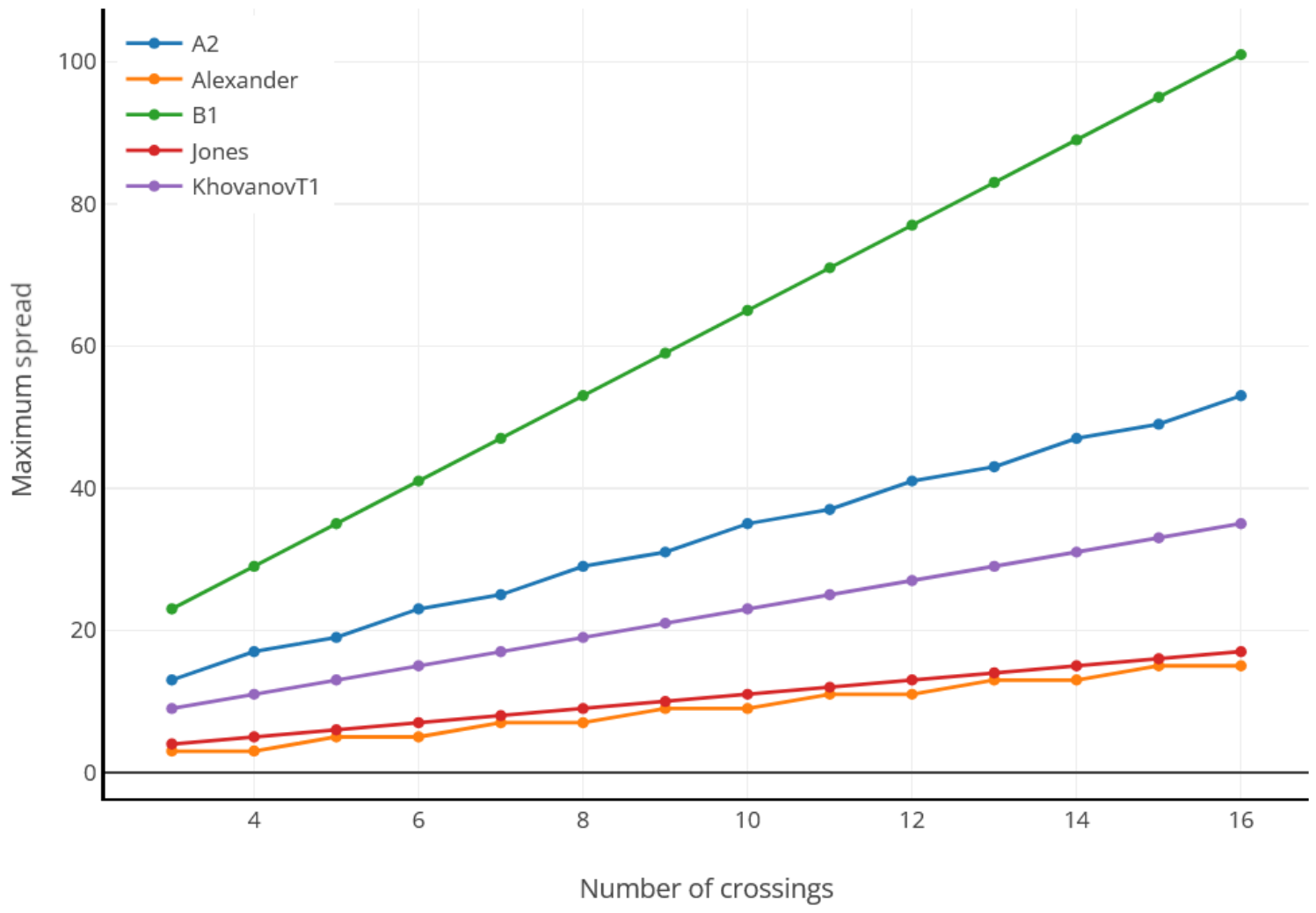

,

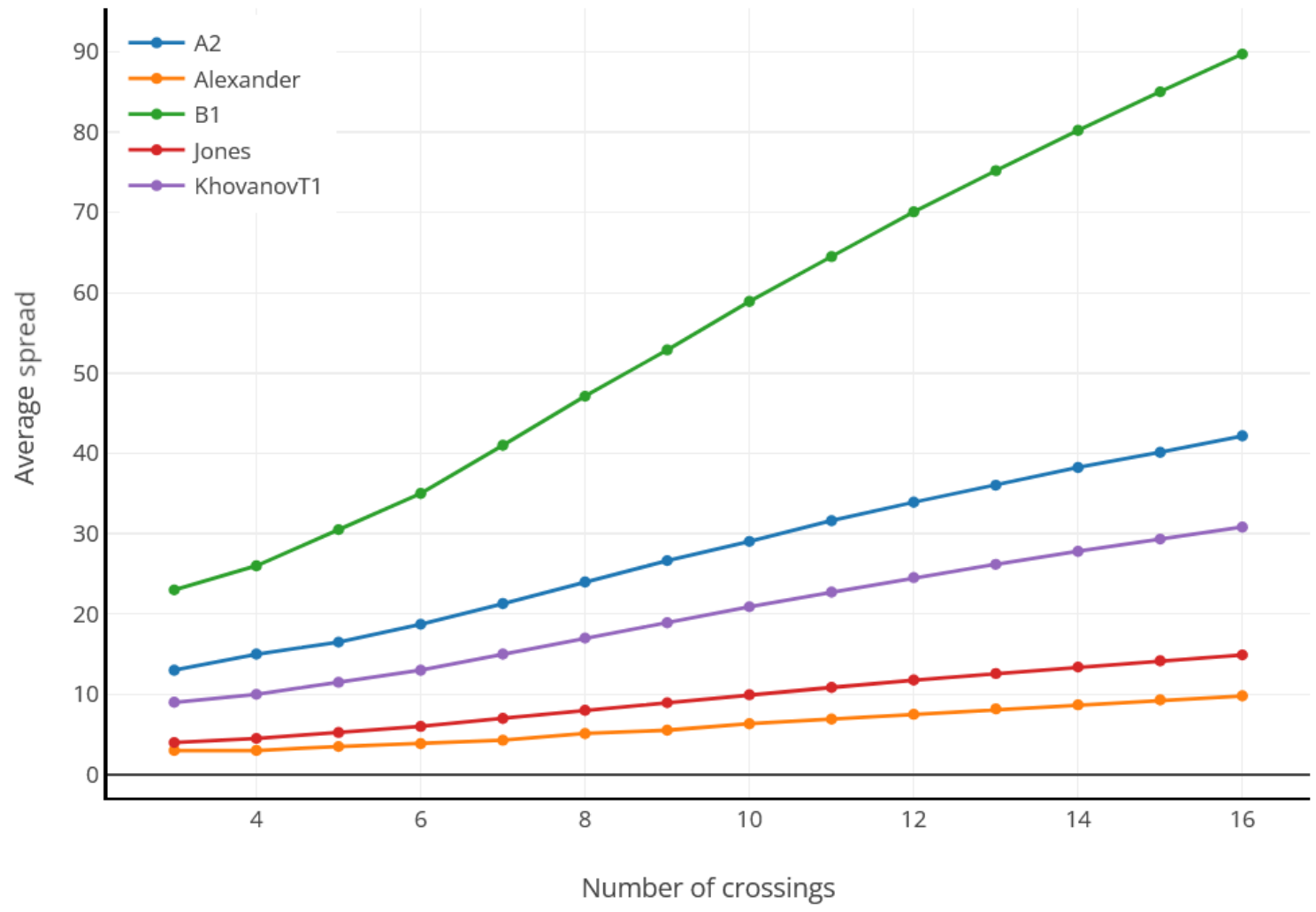

,

where we refer to Figure 42 and Figure 43.

From the data we conjecture:

**Conjecture 13G.1** (Exponential growth)**.** For $Q \in \{A2, A, B1, Kh\}$ (and therefore also for $Q = J$ or $Q = KT1$), we have

$$\mathrm{coeff}_n \in \Omega(\gamma^n) \text{ for some } \gamma \in \mathbb{R}_{>1}.$$

Since $\mathrm{ev}_n \geq \mathrm{coeff}_n$, the same holds for $\mathrm{ev}_n$. We also have

$$\gamma(B1) > \gamma(Q) \text{ for } Q \in \{A2, A, J, KT1\}$$

(and therefore also for $Q = K$). ◇

We now compare pictures of roots, following the ideas in [BCD23, WW01, LTV24b] and the references to various blogs in [BCD23].

The next few pictures show the multiset of roots of the set of $Q \in \{A2, A, B1, J, KT1\}$, with the pictures zoomed in on the right. In formulas, we plot

$$\Big(\mathtt{roots}(p) \mid p \in \{Q(K) \mid K \in \mathcal{K}_{16}\}\Big)$$

where high brightness indicates the high density. Each picture has the same scaling in the axes and are centered on the origin. A faint green circle is drawn to indicate the unit circle. The left pictures are zoomed such that all roots are displayed with nothing more, and the right pictures are zoomed to the same scale as each other.

$A1$ :

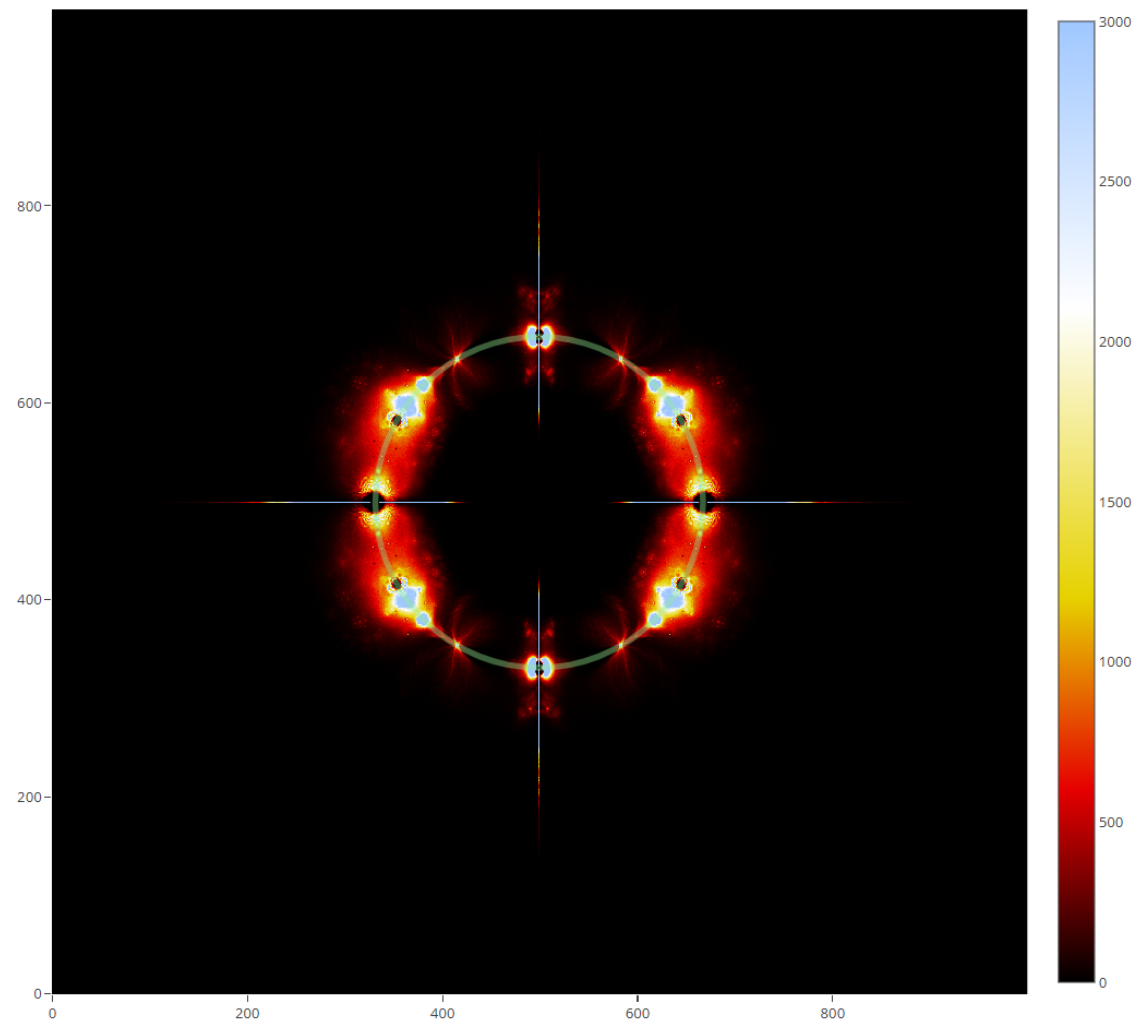

,

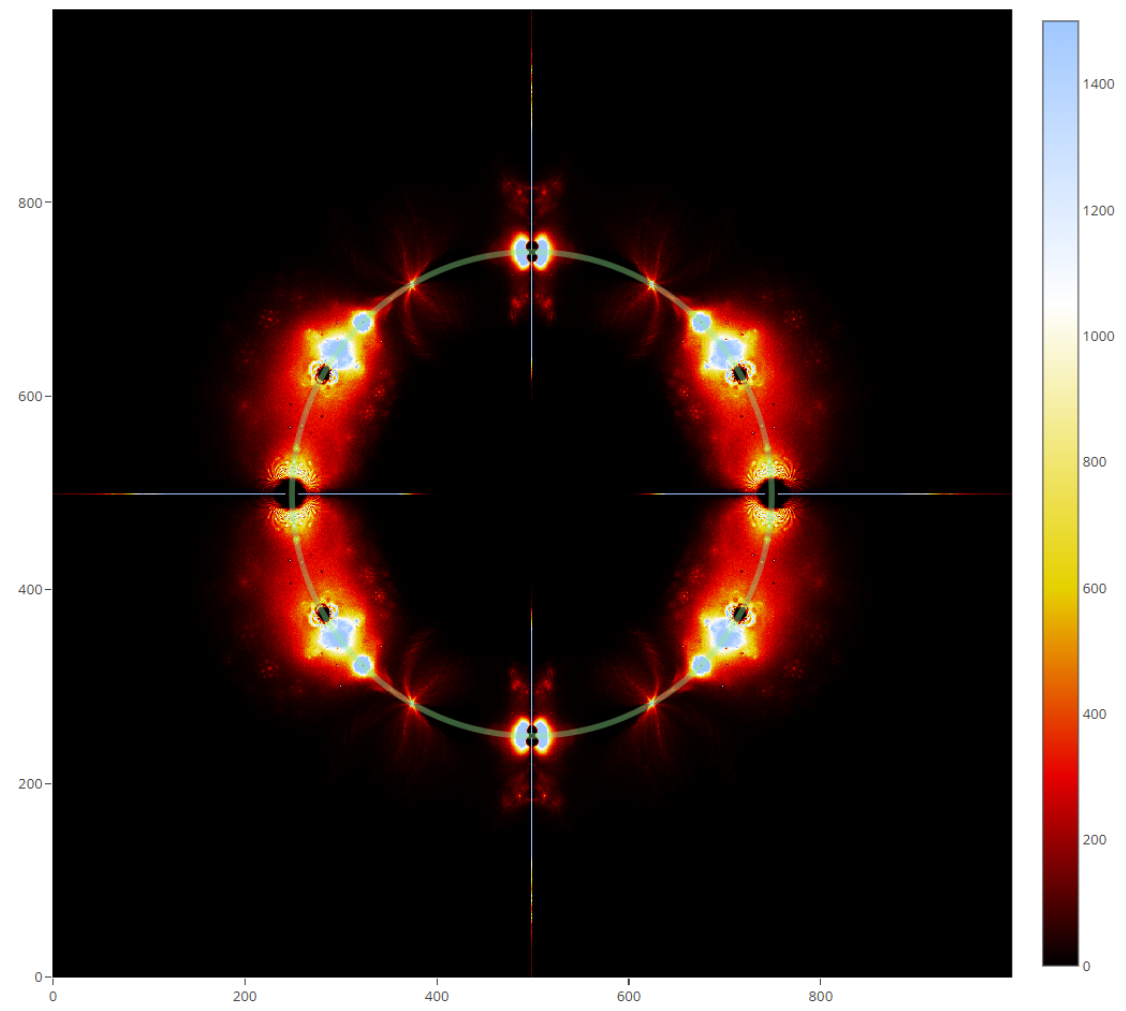

,

$A$ :

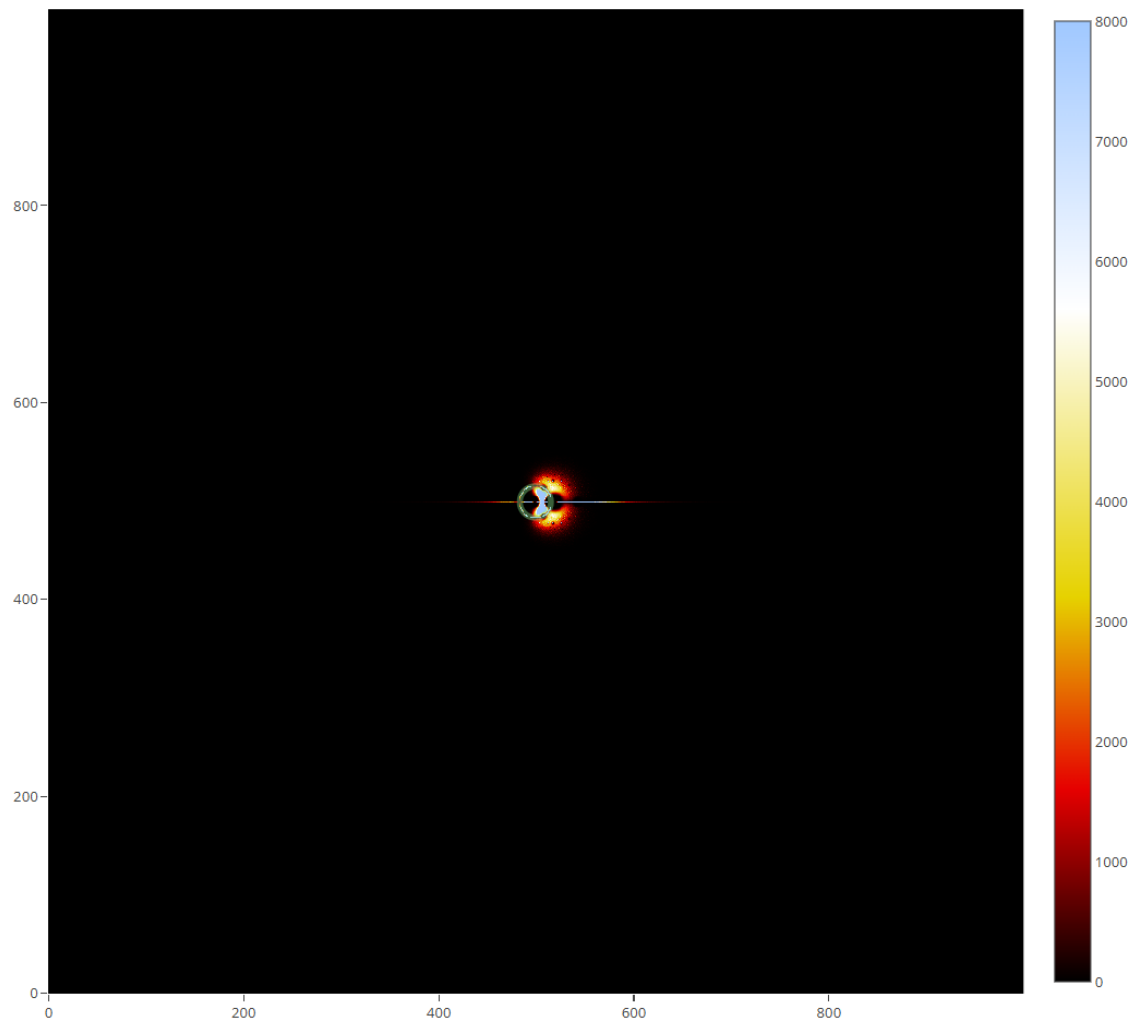

,

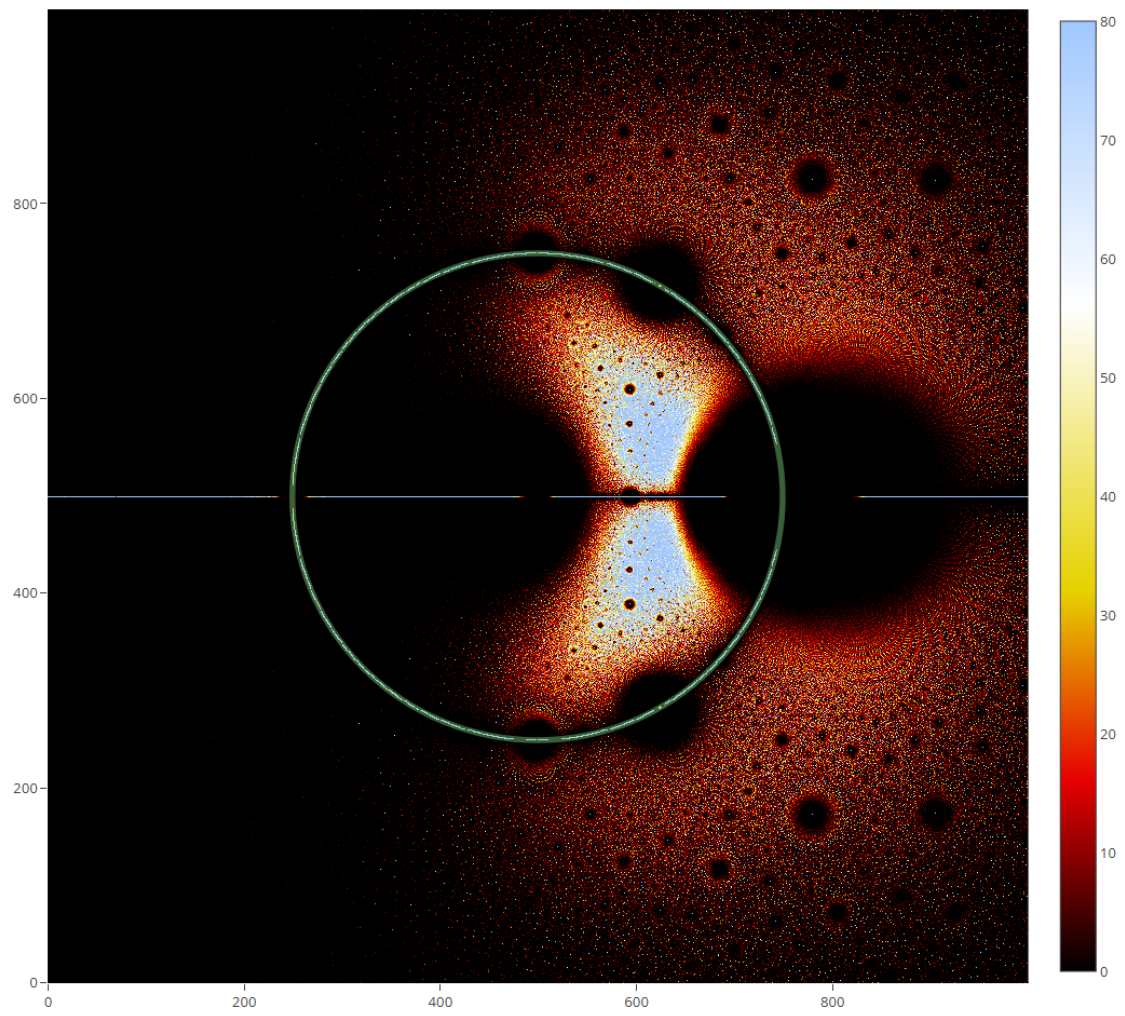

,

$B1:$ 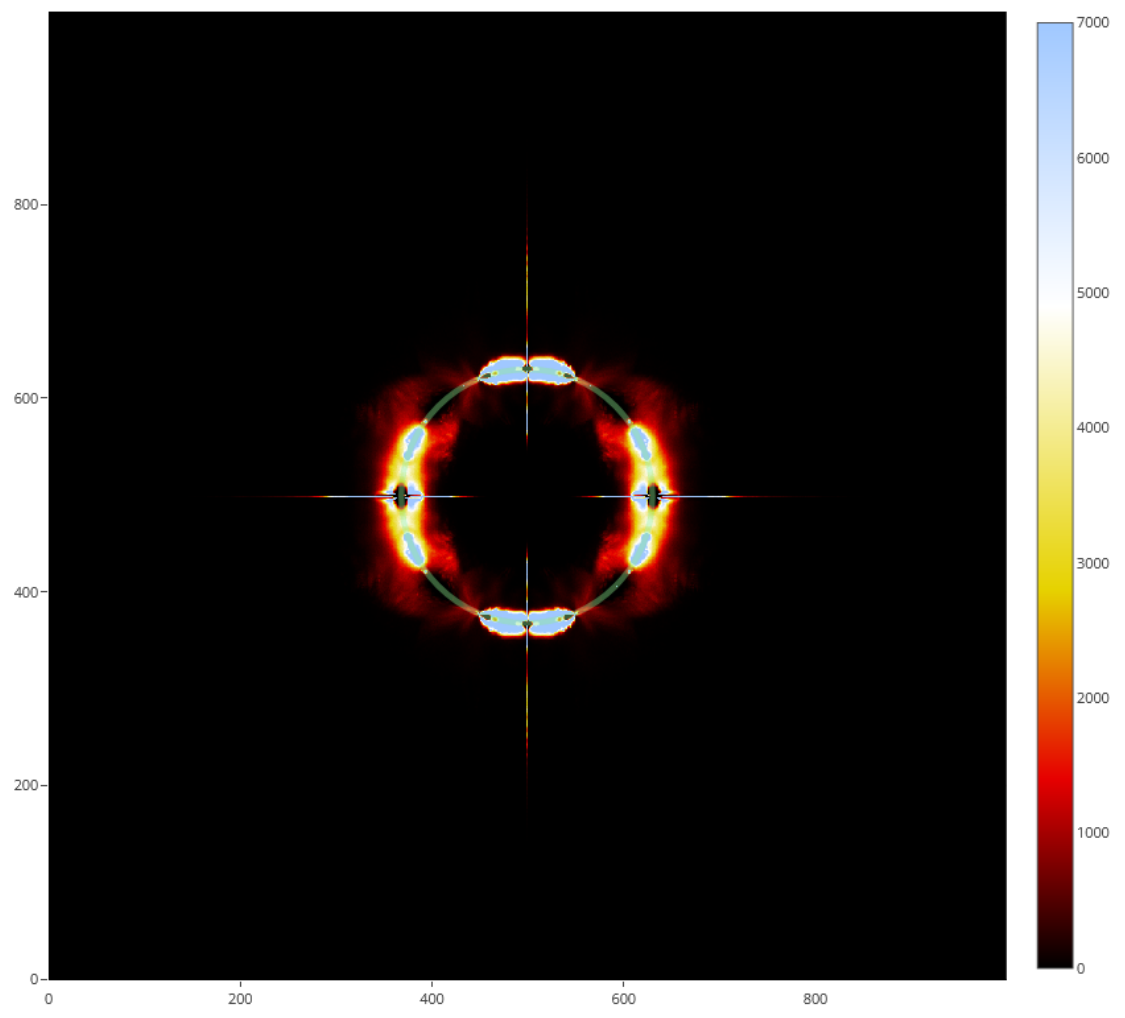 , 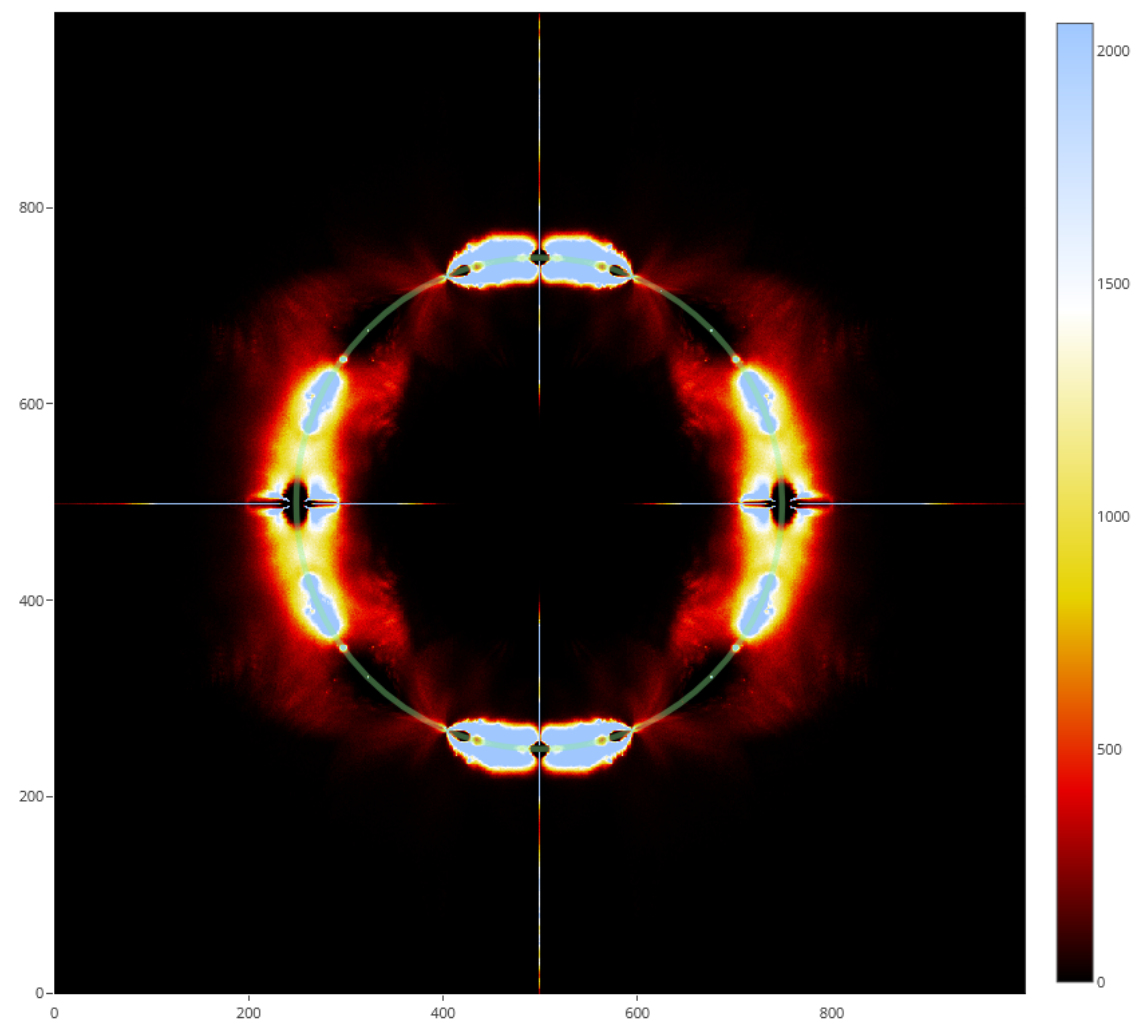 ,

$J:$ 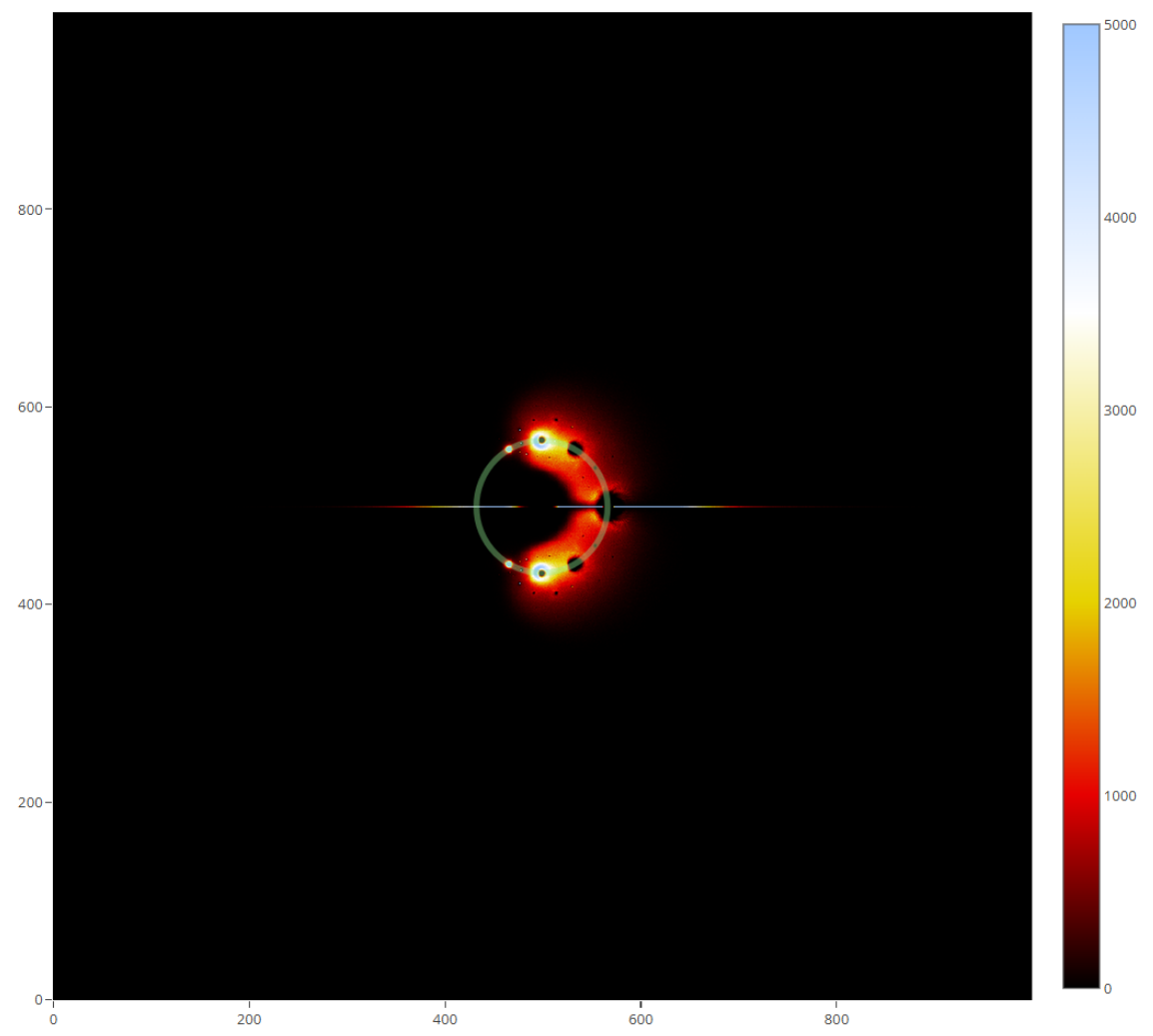 , 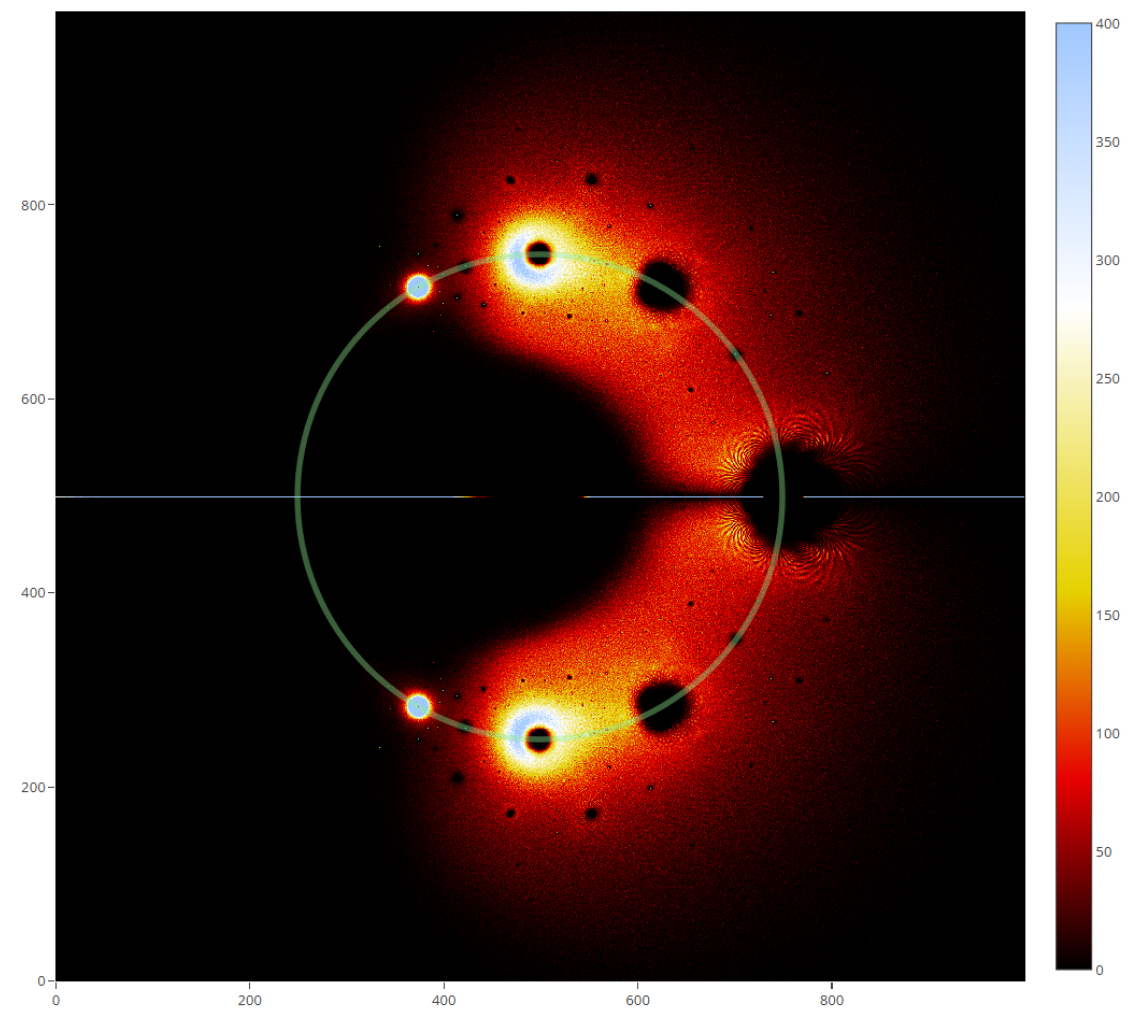 ,

$KT1:$ 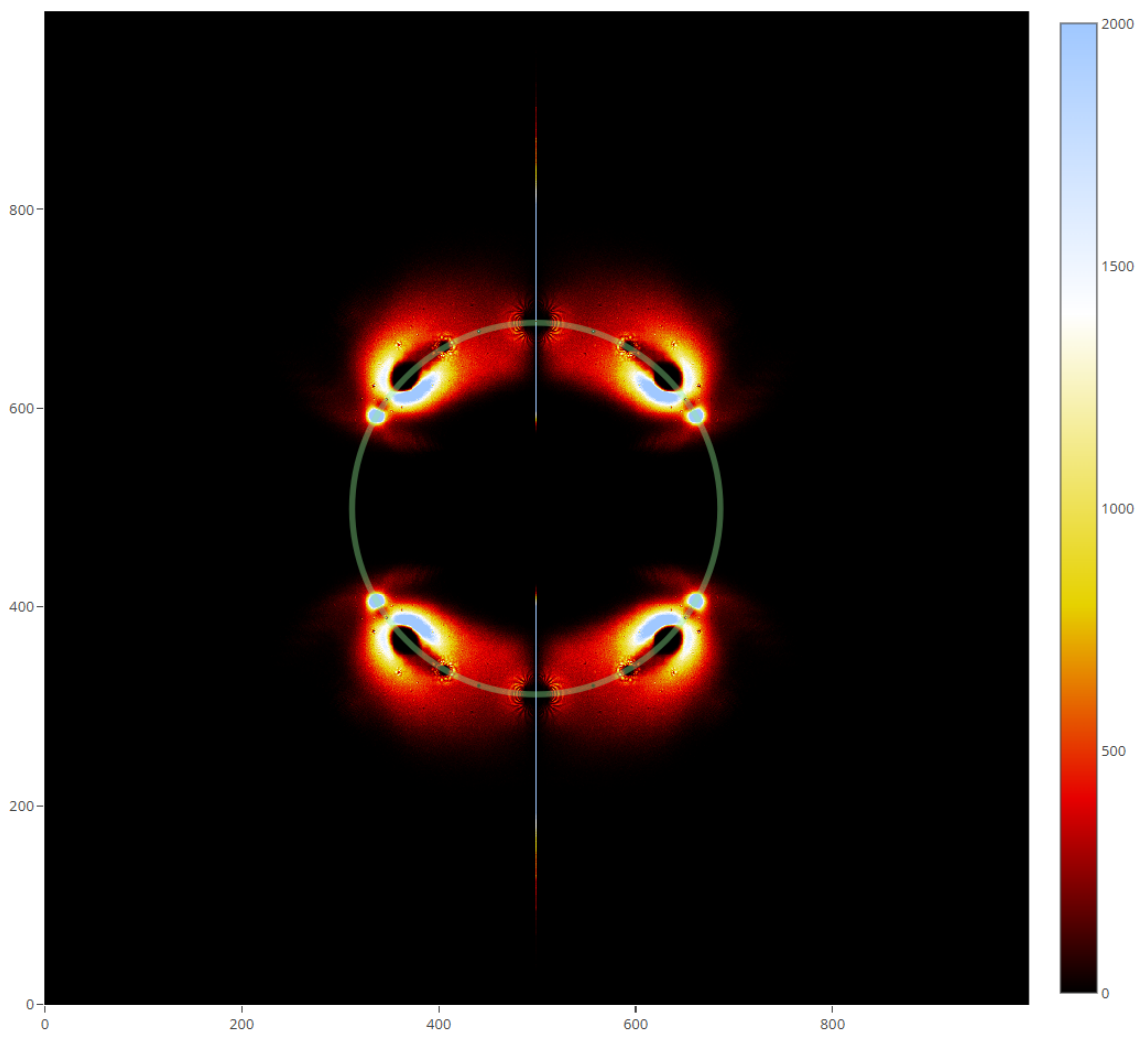 , 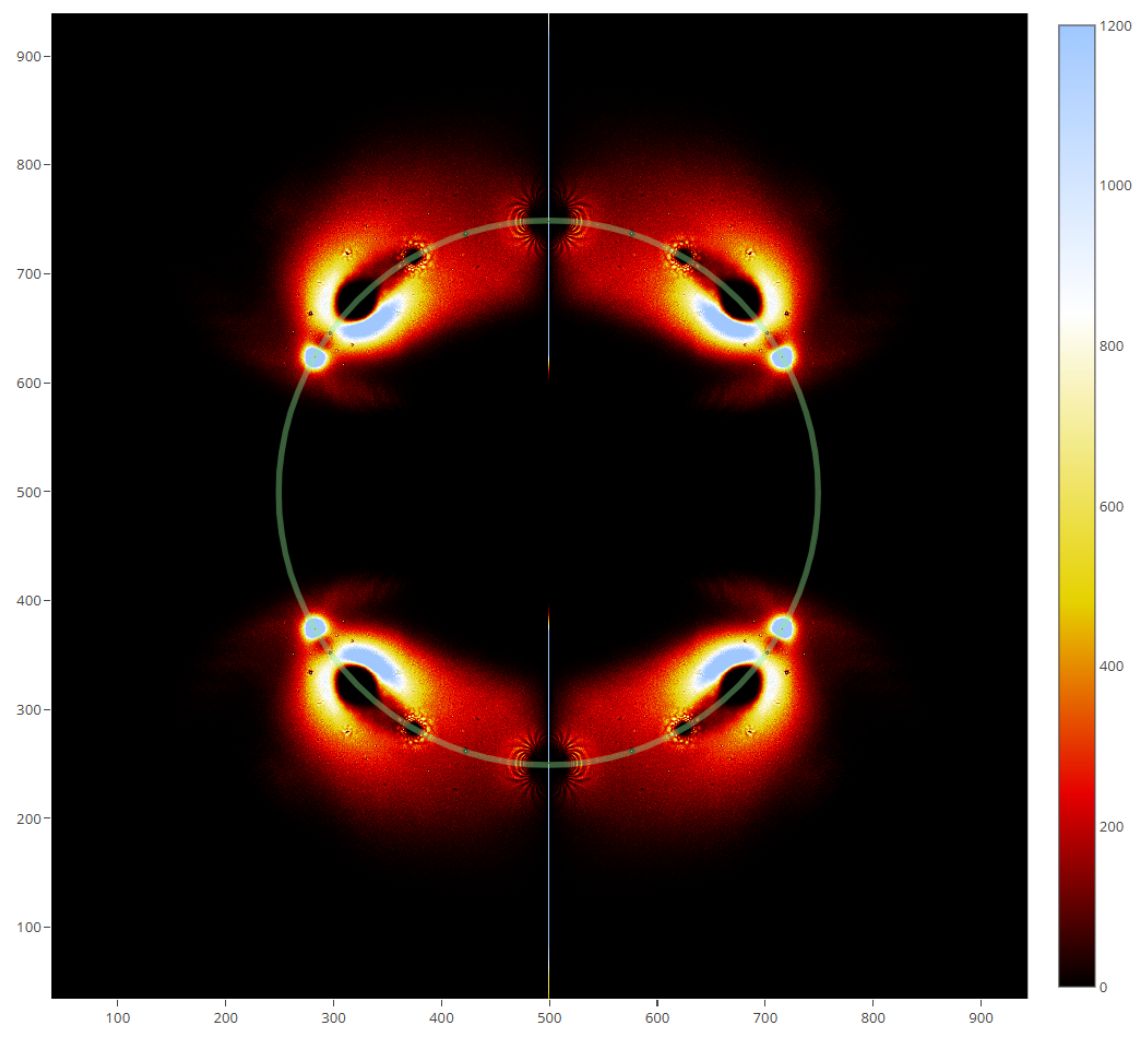 .

Here is the associated data.

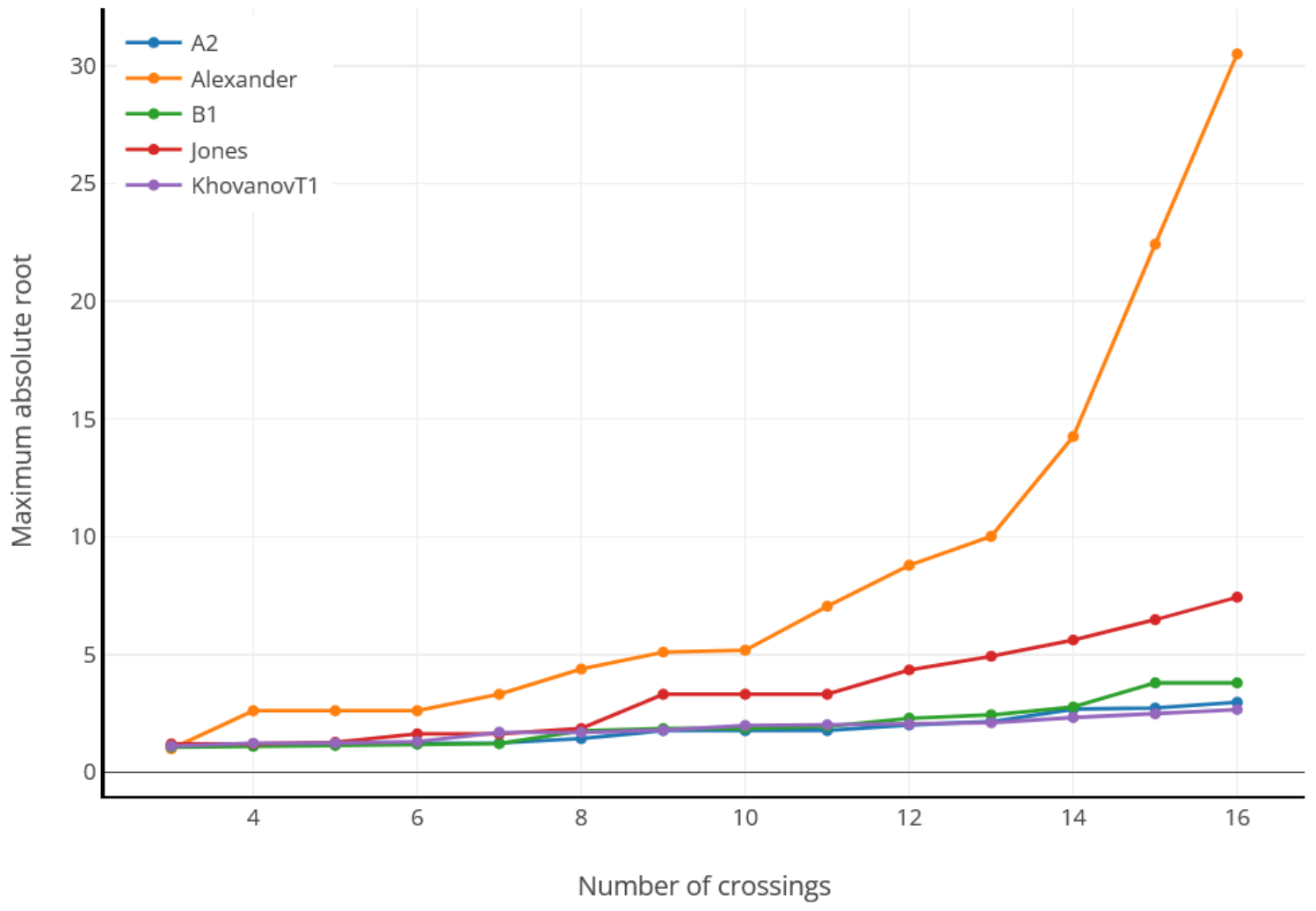

,

see Figure 44 for details. Continued:

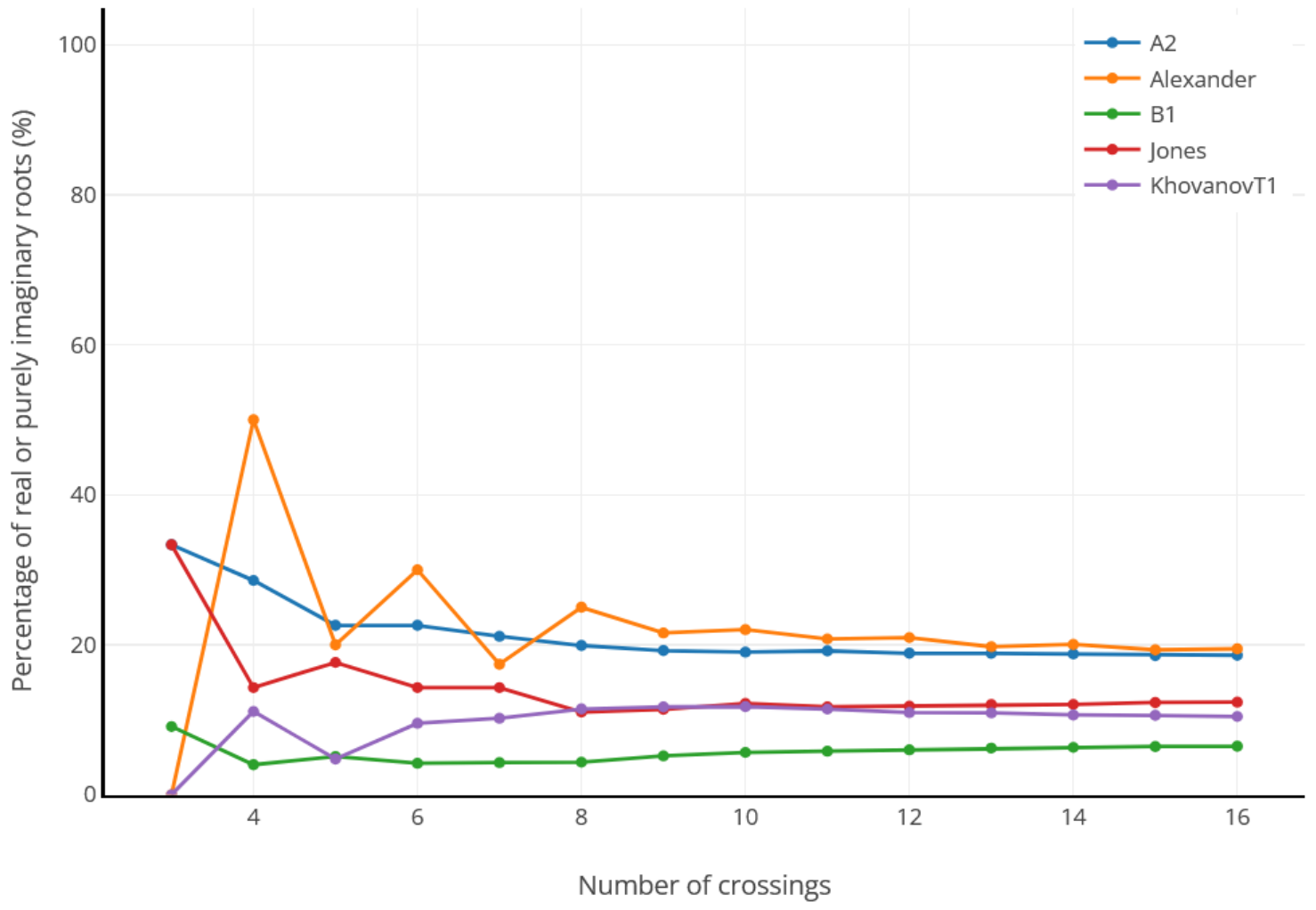

,

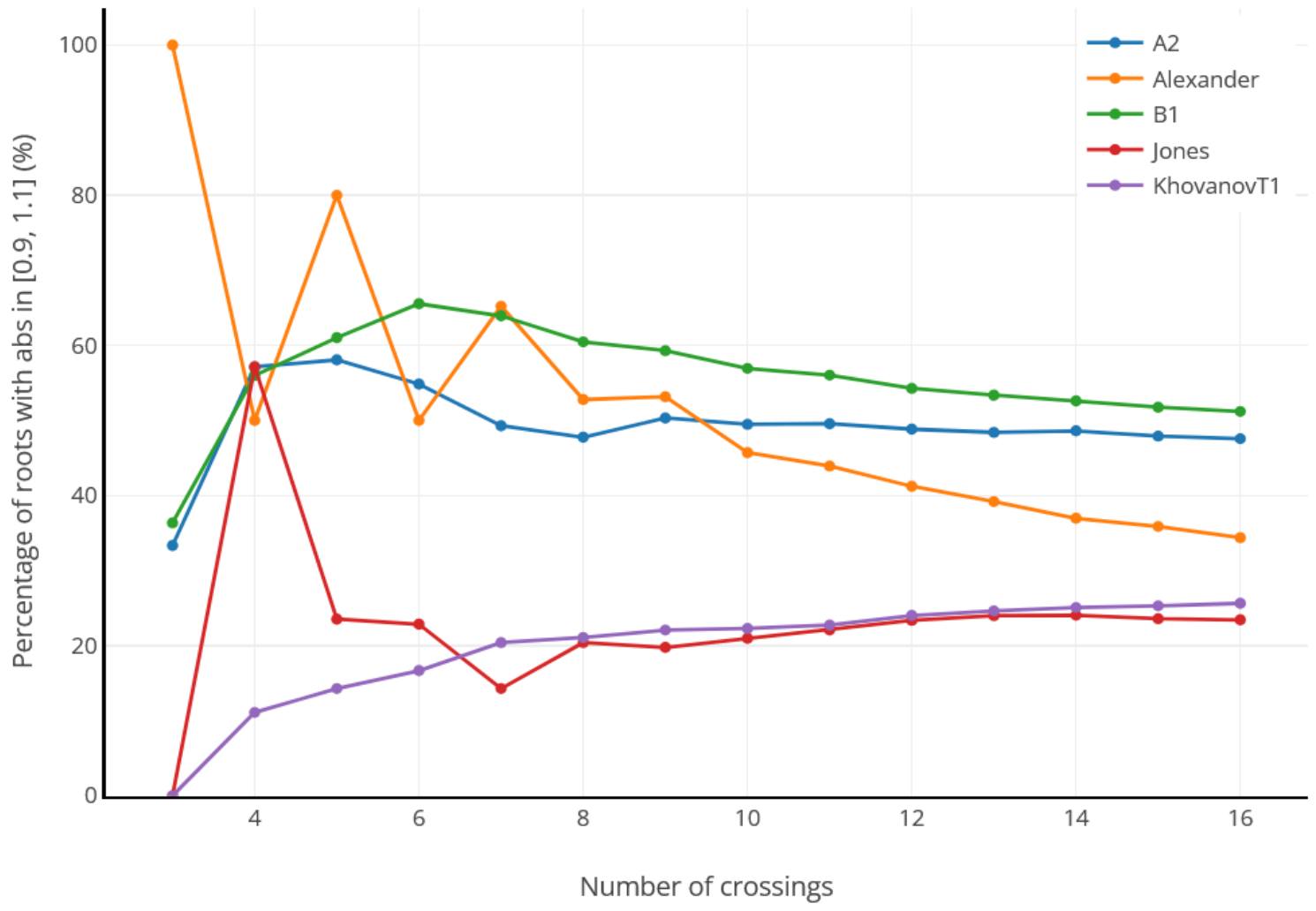

,

which are Figure 45 and Figure 46.

Before we analyze the pictures, we need a bit of terminology.

**Notation 13G.2.** A root is ***purely real or imaginary***, or simply pure, if it is on either of the coordinate axes.

The so-called ***Perron–Frobenius (PF) root*** $\lambda$ is the pure maximal root when taking absolute values; this might not exist, but it almost always exists in our setting. The maximal one, when varying over all polynomials, can be read-off from the left pictures: the edges of the bounding squares in these pictures are of length twice $|\max\{\lambda|\lambda$ is a PF root of $Q(K)$ for $K \in \mathcal{K}_{16}\}|$. ◇

We observe the following.

**(a)** For A2, B1 and KT1 the roots are concentrated around the unit circle. The PF root is fairly small.

**(b)** For A and J the roots are quite spread in the first and fourth quadrants. The PF root is much larger.

For comparison, let us recall a few facts about the distribution of roots of random polynomials. Notable are:

(i) There is an expected tail of pure roots, with one large PF root PFmax(random). In fact, for random polynomials, papers such as [**Gem86**, **Bai97**] implies that the PF root is expected to be very large. That is, letting $M$ be the bound of the coefficients, then, almost surely PFmax(random) $\geq \sqrt{M}$.

(ii) By [**Kac49**], the number of expected real roots of a polynomial of degree $k$ is $2\ln(k)/\pi$. To simplify our calculation, assume that the average degree of our sample is $e^2$ for the usual $e \approx 2.71\ldots$. Then the expected percentage of real roots is $4\pi/e^2 \approx 17.23\%$.

(iii) The clustering of roots around the unit circle is expected for random polynomials; see, for example, [**SV95**]. Fairly explicit formulas for the distribution are known, for example, see [**MBF$^+$97**], but the only thing we notice here is that, in the limit, almost all roots will have an absolute value in $[1-\epsilon, 1+\epsilon]$ for all $\epsilon \in \mathbb{R}_{>0}$.

(iv) We do not know a general statement about holes in the plots of random polynomials, but see [**BCD23**] for similar patterns (this should be true for other integer-valued polynomials as well). The case of coefficients in $\{-1, 0, 1\}$ is addressed, for example, in [**CKW17**].

13H. **Comparison – ballmapper.** Mapper algorithms are fundamental tools in TDA and EDA, first introduced in [**SMC07**], and are used for exploring and visualizing data. These algorithms integrate techniques such as dimensionality reduction, clustering, and graph construction to transform data into a graph. For our purposes, we will use a slight modification known as ball mapper from [**Dło19**], which provides more aesthetically pleasing visualizations. Instead of a formal definition, here is a short summary:

(i) Having a point cloud in $\mathbb{R}^N$ and a fixed $\epsilon \in \mathbb{R}_{>0}$, the ball mapper algorithm creates a graph $G = G(\epsilon)$ by choosing certain anchor points, and then using balls of radius $\epsilon$ around these points.

(ii) The vertices of $G$ are obtained by collapsing all points within one ball into a vertex. The size of the vertex corresponds to the number of points collapsed to that vertex (large = many points).

(iii) The edges of $G$ come from the intersections of the ball, as in this ChatGPT generated picture:

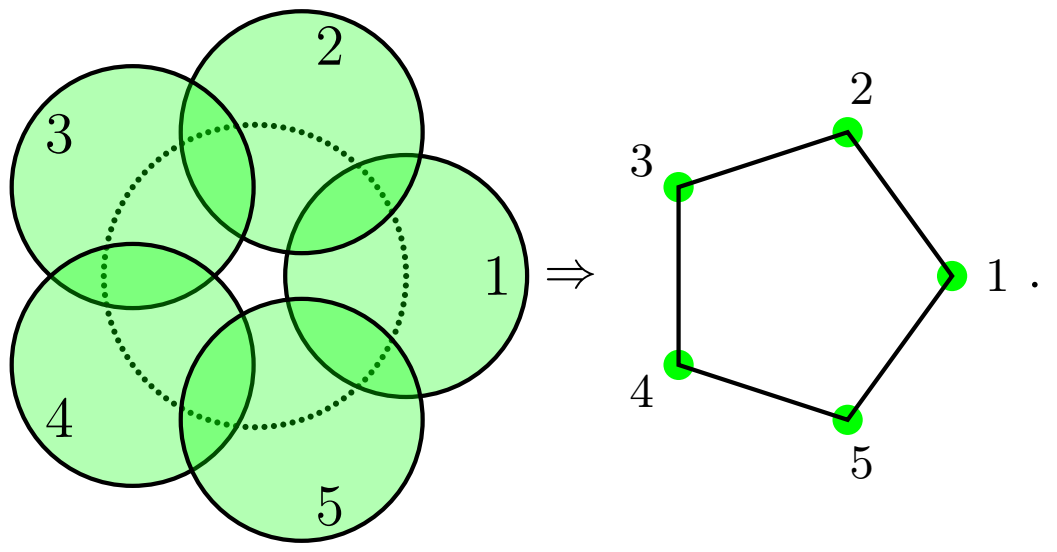

We give only a sample of what can be found in [**TZ25b**]; in particular, the graphs are much more impressive in the interactive plot that can be found in [**TZ25b**]. On that page, we also explain how these were created using blueprint files provided by the webpage of the Dioscuri Centre in Topological Data Analysis.

All plots below are for $n = 15$, and the right picture is a zoomed in version.

A2: 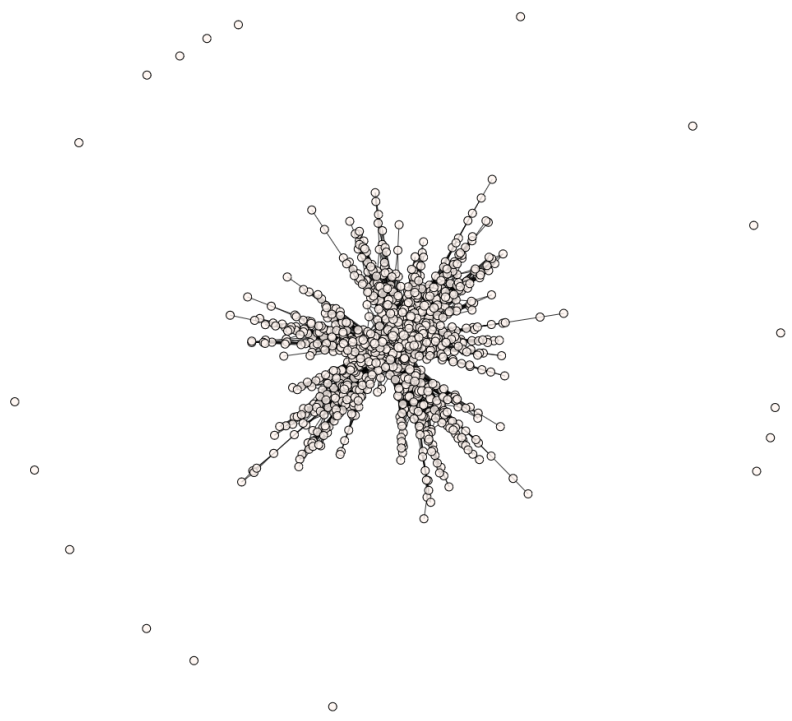, 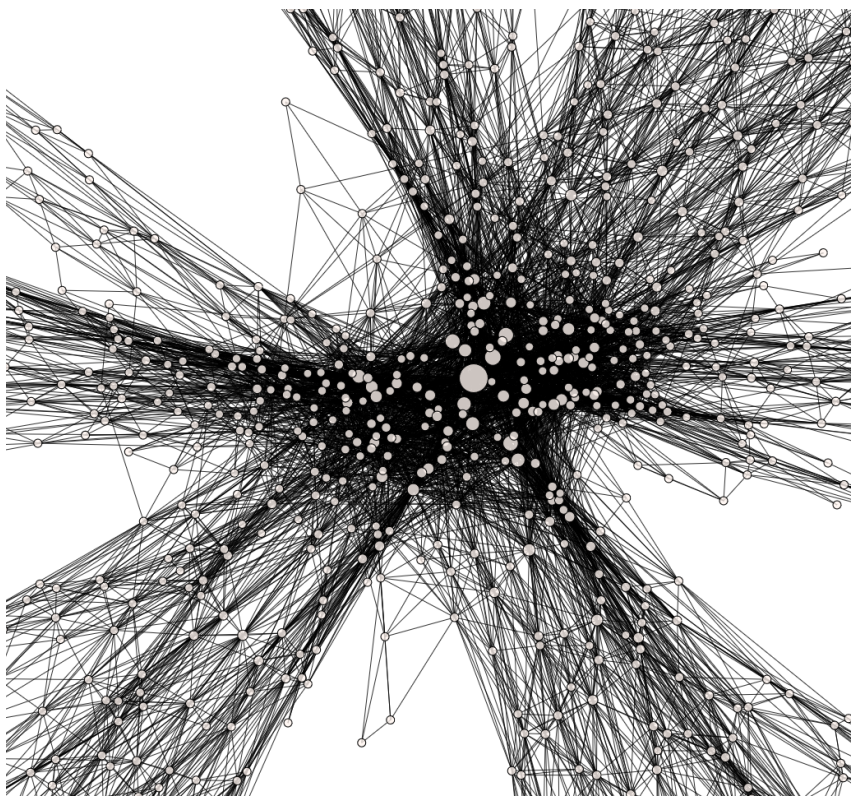,

A: 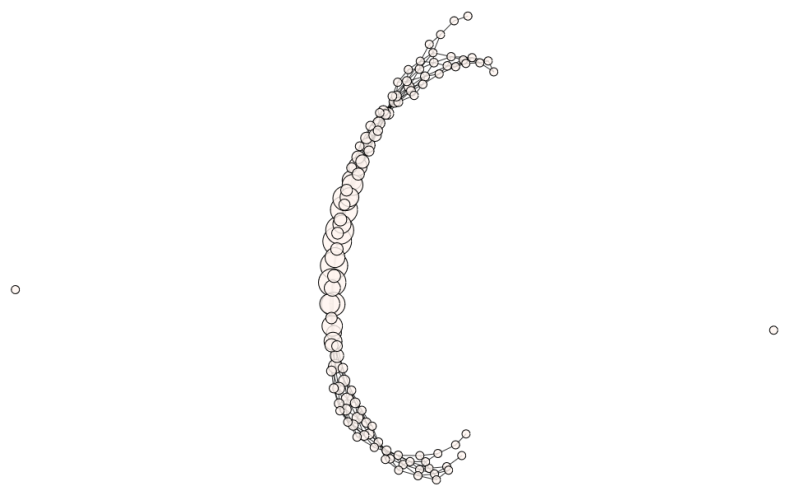, 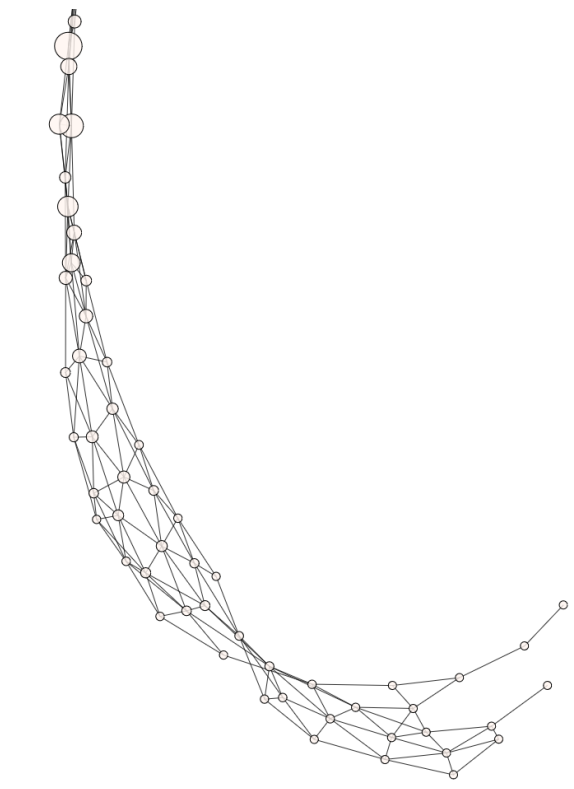,

B1: 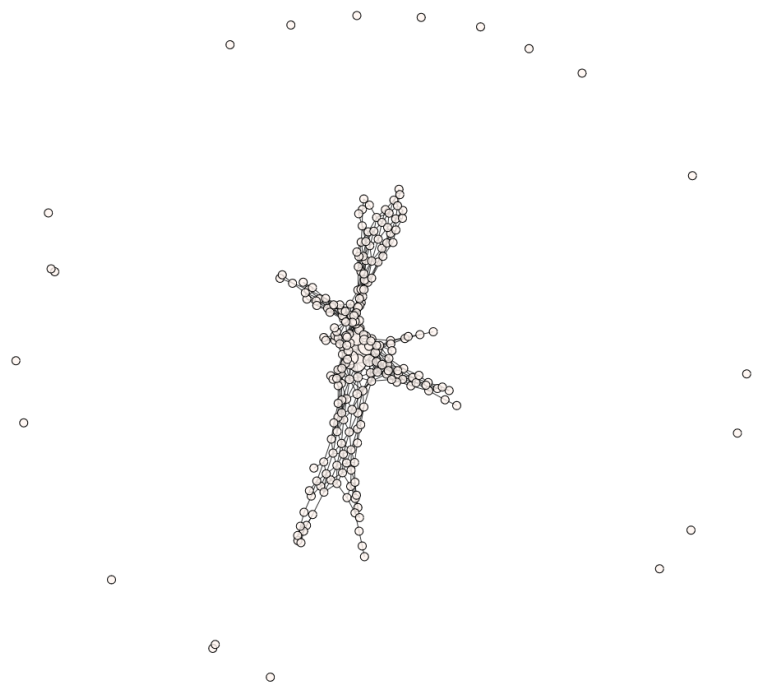, 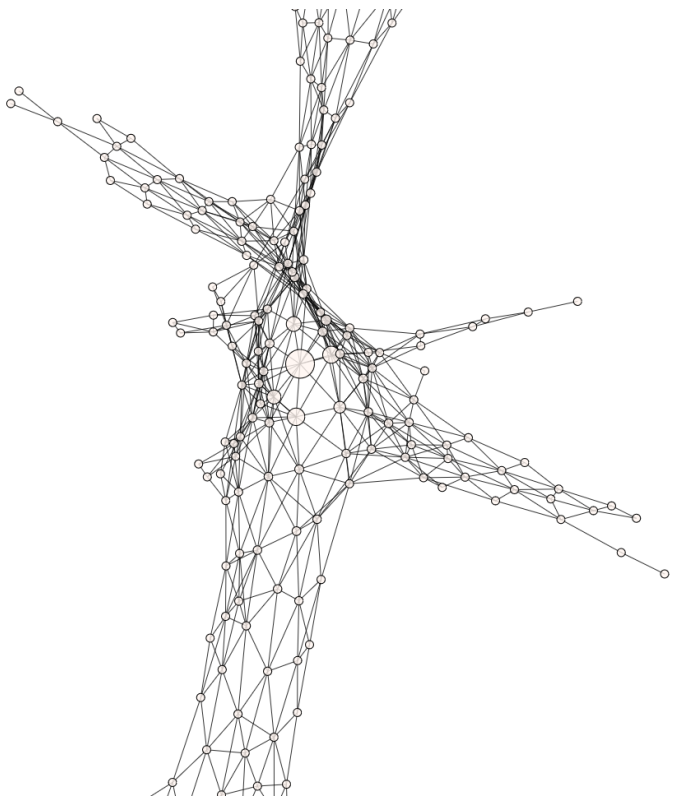,

J: 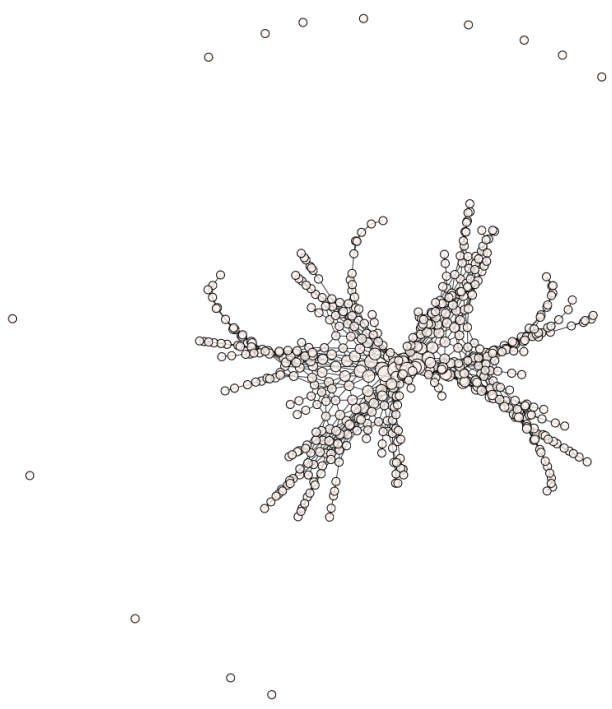, 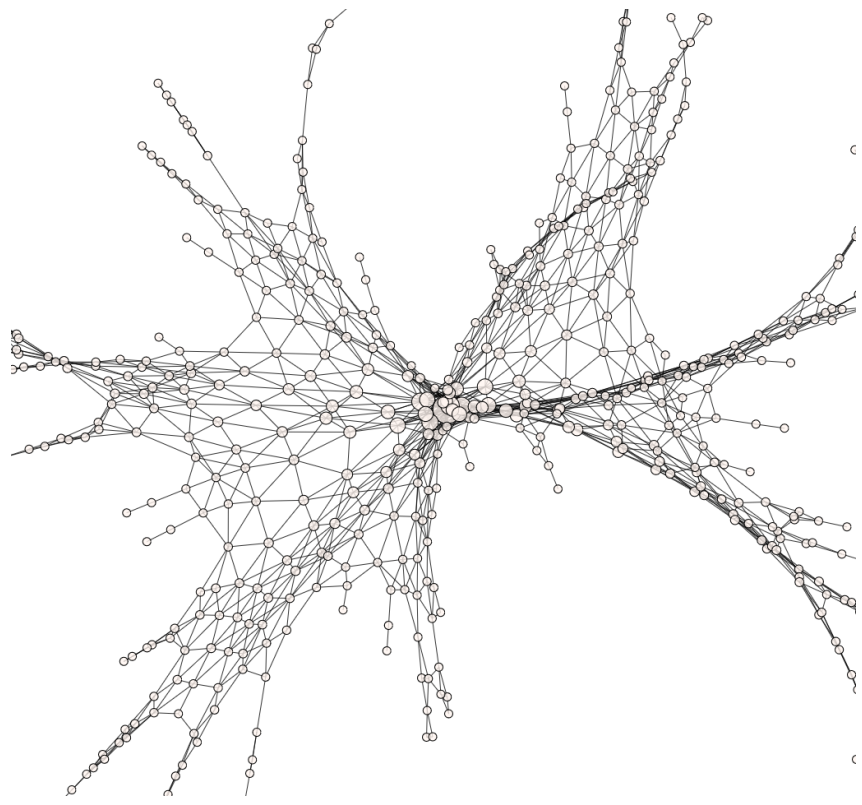,

K: 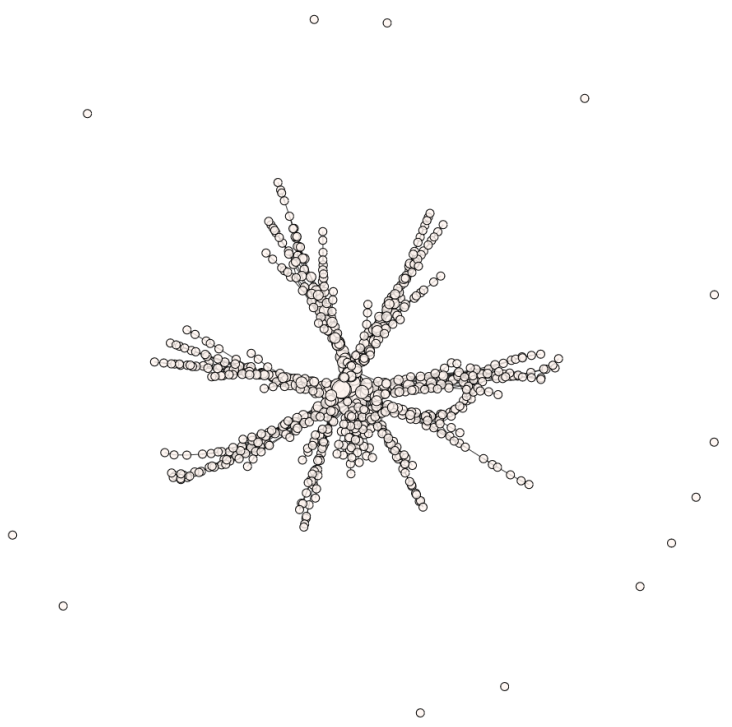 , 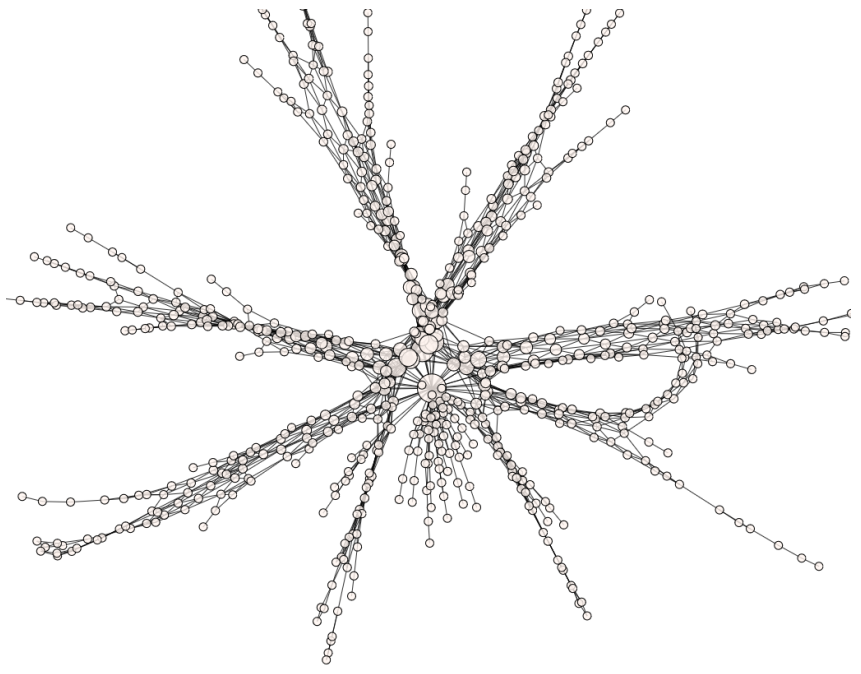 .

These are the ***quantum aliens with their tentacles***.

Additionally, and much nicer in the interactive plot, there is a comparison of the invariants, *e.g.*:

J: 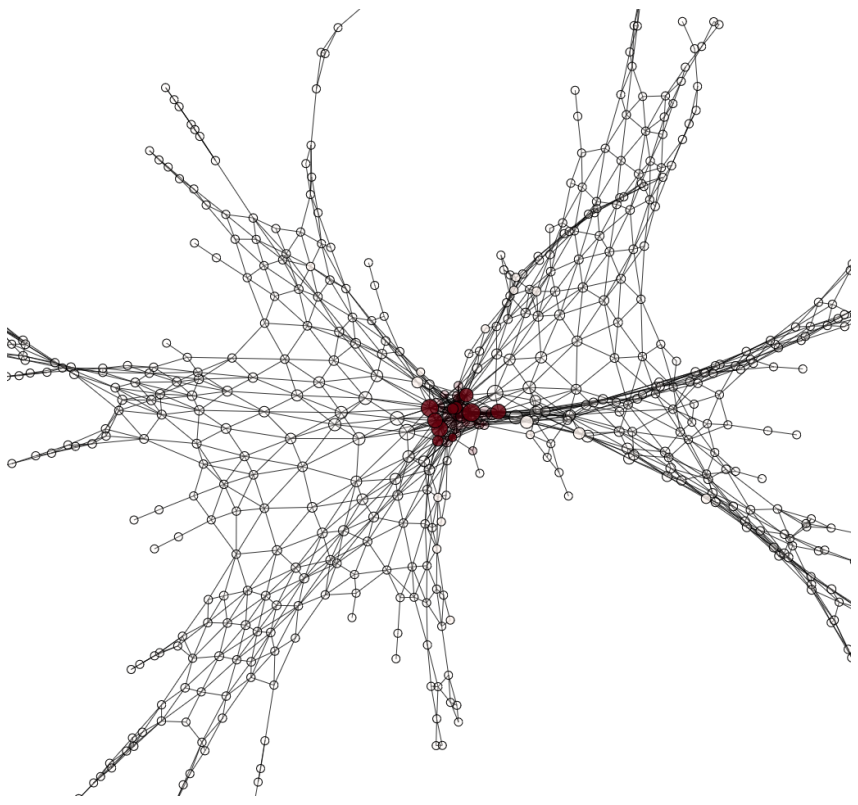 , K: 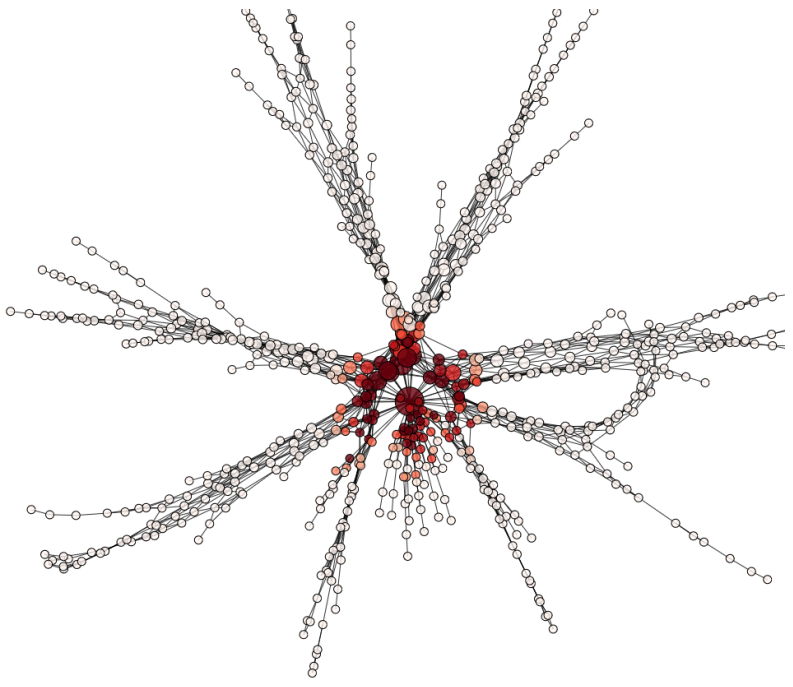 .

This should be read as follows. The selected (and colored) vertices in J appear in the colored range in K.

One observes the following:

- **(a)** All graphs have a dense center. That means that there are a lot of knots with similar quantum invariants.
- **(b)** The vertices on the end of the tentacles are small and mostly very special knots such as alternating or torus knots. These are the knots with very unique quantum invariants.
- **(c)** The only invariant that is different from the others is B1. All other invariants are roughly similarly spread. B1 seems to detect something different from the others, while being more closely packed at the center.

13I. **Copyable data tables.** The title of this section is spot on.

| n | A2 | A | B1 (R) | B1 (Skein) | J | K |
|---|---|---|---|---|---|---|
| 3 | 0.000961 | 0.000548 | 0.04195 | 0.006364 | 0.000952 | 0.000819 |
| 4 | 0.001044 | 0.000474 | 0.058764 | 0.029009 | 0.001174 | 0.001009 |
| 5 | 0.001204 | 0.000413 | 0.063757 | 0.027541 | 0.001359 | 0.001187 |
| 6 | 0.001540 | 0.000371 | 0.073052 | 0.106084 | 0.002156 | 0.001333 |
| 7 | 0.002078 | 0.000344 | 0.109800 | 0.152202 | 0.002442 | 0.001458 |
| 8 | 0.006668 | 0.000327 | 0.146327 | 0.382624 | 0.003204 | 0.001624 |
| 9 | 0.010328 | 0.000340 | 0.434684 | 0.700274 | 0.004329 | 0.001772 |
| 10 | 0.022115 | 0.000377 | 0.727308 | 2.001462 | 0.005811 | 0.001989 |
| 11 | 0.035645 | 0.000447 | 41.247445 | 5.641350 | 0.007834 | 0.002090 |
| 12 | 0.073140 | 0.000800 | 13.954301 | 30.957447 | 0.010320 | 0.002089 |
| 13 | 0.189876 | 0.001039 | 28.231526 | 94.590263 | 0.012261 | 0.002250 |
| 14 | 0.393864 | 0.001244 | 113.210082 | 535.154620 | 0.016506 | 0.002312 |

FIGURE 35. Average time; copyable data.

| n | A2 | A | B1 | J | K | KT1 | J+KT1 | All |
|---|---|---|---|---|---|---|---|---|
| 3 | 100.0 | 100.0 | 100.0 | 100.0 | 100.0 | 100.0 | 100.0 | 100.0 |
| 4 | 100.0 | 100.0 | 100.0 | 100.0 | 100.0 | 100.0 | 100.0 | 100.0 |
| 5 | 100.0 | 100.0 | 100.0 | 100.0 | 100.0 | 100.0 | 100.0 | 100.0 |
| 6 | 100.0 | 100.0 | 100.0 | 100.0 | 100.0 | 100.0 | 100.0 | 100.0 |
| 7 | 100.0 | 100.0 | 100.0 | 100.0 | 100.0 | 100.0 | 100.0 | 100.0 |
| 8 | 100.0 | 100.0 | 100.0 | 100.0 | 100.0 | 100.0 | 100.0 | 100.0 |
| 9 | 100.0 | 94.0 | 100.0 | 100.0 | 100.0 | 100.0 | 100.0 | 100.0 |
| 10 | 98.7 | 84.7 | 100.0 | 96.3 | 96.3 | 96.3 | 96.3 | 100.0 |
| 11 | 95.8 | 68.7 | 98.1 | 90.1 | 91.1 | 90.7 | 91.1 | 98.1 |
| 12 | 92.1 | 59.5 | 96.6 | 83.0 | 84.3 | 83.8 | 84.1 | 96.7 |
| 13 | 85.7 | 43.4 | 92.3 | 73.3 | 77.5 | 77.1 | 77.4 | 92.6 |
| 14 | 81.2 | 33.6 | 89.3 | 64.4 | 69.0 | 68.4 | 68.9 | 89.6 |
| 15 | 76.4 | 24.5 | 86.2 | 55.7 | 60.6 | 59.8 | 60.6 | 86.4 |
| 16 | 74.0 | 18.6 | 83.8 | 49.4 | 54.7 | 53.6 | 54.6 | 84.0 |

FIGURE 36. Percentages of unique values; copyable data.

| n | A2 | A | B1 | J | K | KT1 | J+KT1 | All |
|---|---|---|---|---|---|---|---|---|
| 4 | 1.0 | 0.9 | 0.9 | 1.0 | 1.0 | 1.0 | 1.0 | 1.0 |
| 5 | 2.9 | 3.0 | 2.9 | 2.9 | 2.9 | 2.9 | 3.0 | 2.9 |
| 6 | 6.0 | 5.9 | 6.0 | 6.0 | 6.0 | 5.9 | 5.9 | 5.9 |
| 7 | 13.0 | 12.9 | 12.9 | 13.0 | 12.9 | 12.9 | 12.9 | 12.9 |
| 8 | 34.1 | 33.9 | 34.0 | 34.0 | 33.9 | 34.1 | 34.0 | 34.0 |
| 9 | 83.2 | 74.0 | 82.9 | 82.8 | 83.2 | 82.3 | 83.1 | 82.7 |
| 10 | 241.9 | 187.0 | 247.4 | 231.0 | 232.6 | 231.7 | 230.7 | 250.4 |
| 11 | 738.4 | 417.1 | 764.5 | 662.1 | 673.2 | 671.5 | 676.0 | 768.4 |
| 12 | 2532.7 | 1147.5 | 2791.9 | 2098.2 | 2141.2 | 2117.2 | 2134.0 | 2794.7 |
| 13 | 9579.0 | 2638.0 | 11090.3 | 7145.8 | 7902.2 | 7846.9 | 7916.5 | 11077.2 |
| 14 | 39876.6 | 6679.7 | 47869.6 | 25740.7 | 28845.6 | 28586.1 | 29036.8 | 48079.8 |
| 15 | 181201.6 | 15376.9 | 230546.2 | 98445.1 | 115234.8 | 112986.8 | 114700.2 | 231711.7 |
| 16 | 909889.2 | 34001.2 | 1161532.7 | 394876.7 | 475079.0 | 460340.4 | 473986.2 | 1157489.3 |

FIGURE 37. Average comparisons until equal; copyable data.

| n | A2 | Alexander | B1 | Jones | Khovanov | KhovanovT1 |
|---|---|---|---|---|---|---|
| 3 | 2 | 1 | 1 | 1 | 1 | 1 |
| 4 | 2 | 3 | 1 | 1 | 1 | 2 |
| 5 | 2 | 3 | 2 | 2 | 1 | 2 |
| 6 | 2 | 5 | 2 | 3 | 2 | 3 |
| 7 | 2 | 9 | 5 | 4 | 3 | 5 |
| 8 | 4 | 13 | 17 | 9 | 5 | 8 |
| 9 | 4 | 23 | 40 | 13 | 7 | 12 |
| 10 | 6 | 37 | 76 | 21 | 11 | 20 |
| 11 | 9 | 59 | 202 | 34 | 18 | 34 |
| 12 | 11 | 109 | 517 | 61 | 31 | 58 |
| 13 | 22 | 163 | 1440 | 103 | 53 | 100 |
| 14 | 31 | 281 | 3474 | 177 | 89 | 167 |
| 15 | 56 | 451 | 10063 | 305 | 157 | 296 |
| 16 | 87 | 813 | 26219 | 533 | 267 | 503 |

FIGURE 38. Maximal coefficient; copyable data.

| n | A2 | Alexander | B1 | Jones | Khovanov | KhovanovT1 |
|---|---|---|---|---|---|---|
| 3 | 7 | 3 | 9 | 3 | 4 | 4 |
| 4 | 7 | 5 | 9 | 5 | 6 | 6 |
| 5 | 9 | 7 | 9 | 7 | 8 | 8 |
| 6 | 9 | 13 | 27 | 13 | 14 | 14 |
| 7 | 11 | 21 | 53 | 21 | 22 | 22 |
| 8 | 19 | 45 | 203 | 45 | 46 | 46 |
| 9 | 27 | 75 | 475 | 75 | 76 | 76 |
| 10 | 37 | 121 | 987 | 121 | 122 | 122 |
| 11 | 57 | 209 | 2609 | 209 | 210 | 210 |
| 12 | 101 | 377 | 7287 | 377 | 378 | 378 |
| 13 | 163 | 663 | 20047 | 663 | 664 | 664 |
| 14 | 257 | 1145 | 51311 | 1145 | 1146 | 1146 |
| 15 | 439 | 2037 | 146629 | 2037 | 2038 | 2038 |
| 16 | 717 | 3581 | 397707 | 3581 | 3582 | 3582 |

FIGURE 39. Maximal coefficient sum; copyable data.

| n | A2 | Alexander | B1 | Jones | Khovanov | Khovanovt1 |
|---|---|---|---|---|---|---|
| 3 | 7.0000 | 3.0000 | 9.0000 | 3.0000 | 4.0000 | 4.0000 |
| 4 | 6.0000 | 4.0000 | 8.0000 | 4.0000 | 5.0000 | 5.0000 |
| 5 | 7.0000 | 5.0000 | 8.5000 | 5.0000 | 6.0000 | 6.0000 |
| 6 | 7.0000 | 7.5714 | 12.7142 | 7.5714 | 8.5714 | 8.5714 |
| 7 | 8.2857 | 11.1428 | 20.4285 | 11.1428 | 12.1428 | 12.1428 |
| 8 | 10.5428 | 18.8857 | 46.3714 | 18.8285 | 19.9428 | 19.9428 |
| 9 | 13.0238 | 29.6666 | 95.0476 | 29.6428 | 30.7142 | 30.7142 |
| 10 | 17.2409 | 47.6345 | 202.6867 | 47.5461 | 48.8192 | 48.8192 |
| 11 | 23.9762 | 75.9538 | 447.5218 | 75.9837 | 77.5205 | 77.5205 |
| 12 | 33.3903 | 119.2499 | 955.8572 | 119.3829 | 121.4605 | 121.4605 |
| 13 | 47.1625 | 183.4964 | 2034.3612 | 184.1180 | 187.5322 | 187.5322 |
| 14 | 66.6552 | 278.7659 | 4256.5323 | 280.3088 | 286.0564 | 286.0564 |
| 15 | 92.9985 | 412.9504 | 8620.9249 | 416.3381 | 426.4781 | 426.4781 |
| 16 | 130.3453 | 607.3186 | 17436.7151 | 614.2326 | 632.2439 | 632.2439 |

FIGURE 40. Average coefficient sum; copyable data.

| n | A2 | Alexander | B1 | Jones | Khovanovt1 |
|---|---|---|---|---|---|
| 3 | 0.5384 | 1.0000 | 0.3913 | 0.7500 | 0.4444 |
| 4 | 0.5384 | 1.6666 | 0.3913 | 1.0000 | 0.5454 |
| 5 | 0.5384 | 2.3333 | 0.3913 | 1.1666 | 0.6153 |
| 6 | 0.5384 | 3.0000 | 0.6585 | 1.8571 | 0.9333 |
| 7 | 0.5384 | 5.0000 | 1.1276 | 2.6250 | 1.2941 |
| 8 | 0.7600 | 6.6000 | 3.8301 | 5.0000 | 2.4210 |
| 9 | 0.9310 | 11.4000 | 8.0508 | 7.5000 | 3.6190 |
| 10 | 1.2758 | 21.0000 | 15.1846 | 11.0000 | 5.3043 |
| 11 | 1.7272 | 29.0000 | 36.7464 | 17.4166 | 8.4000 |
| 12 | 2.7297 | 43.8571 | 94.6363 | 29.0000 | 14.0000 |
| 13 | 4.4054 | 71.8571 | 241.5301 | 47.3571 | 22.8965 |
| 14 | 6.2682 | 124.4285 | 576.5280 | 76.3333 | 36.9677 |
| 15 | 10.7073 | 213.0000 | 1543.4631 | 127.3125 | 61.7575 |
| 16 | 15.9333 | 304.1111 | 3937.6930 | 210.6470 | 102.3428 |

FIGURE 41. Maximum average coefficient; copyable data.

| n | A2 | Alexander | B1 | Jones | Khovanovt1 |
|---|---|---|---|---|---|
| 3 | 13 | 3 | 23 | 4 | 9 |
| 4 | 17 | 3 | 29 | 5 | 11 |
| 5 | 19 | 5 | 35 | 6 | 13 |
| 6 | 23 | 5 | 41 | 7 | 15 |
| 7 | 25 | 7 | 47 | 8 | 17 |
| 8 | 29 | 7 | 53 | 9 | 19 |
| 9 | 31 | 9 | 59 | 10 | 21 |
| 10 | 35 | 9 | 65 | 11 | 23 |
| 11 | 37 | 11 | 71 | 12 | 25 |
| 12 | 41 | 11 | 77 | 13 | 27 |
| 13 | 43 | 13 | 83 | 14 | 29 |
| 14 | 47 | 13 | 89 | 15 | 31 |
| 15 | 49 | 15 | 95 | 16 | 33 |
| 16 | 53 | 15 | 101 | 17 | 35 |

FIGURE 42. Maximal spread; copyable data.

| n | A2 | Alexander | B1 | Jones | Khovanovt1 |
|---|---|---|---|---|---|
| 3 | 13.0000 | 3.0000 | 23.0000 | 4.0000 | 9.0000 |
| 4 | 15.0000 | 3.0000 | 26.0000 | 4.5000 | 10.0000 |
| 5 | 16.5000 | 3.5000 | 30.5000 | 5.2500 | 11.5000 |
| 6 | 18.7142 | 3.8571 | 35.0000 | 6.0000 | 13.0000 |
| 7 | 21.2857 | 4.2857 | 41.0000 | 7.0000 | 15.0000 |
| 8 | 23.9714 | 5.1142 | 47.1142 | 8.0000 | 17.0000 |
| 9 | 26.6428 | 5.5238 | 52.8809 | 8.9523 | 18.9047 |
| 10 | 29.0321 | 6.3253 | 58.9036 | 9.9437 | 20.8955 |
| 11 | 31.6117 | 6.9001 | 64.5056 | 10.8514 | 22.7153 |
| 12 | 33.9183 | 7.5072 | 70.0755 | 11.7638 | 24.5391 |
| 13 | 36.0324 | 8.1801 | 75.2020 | 12.5885 | 26.1881 |
| 14 | 38.2418 | 8.6834 | 80.2214 | 13.3948 | 27.8008 |
| 15 | 40.1237 | 9.2937 | 85.0202 | 14.1587 | 29.3301 |
| 16 | 42.1702 | 9.7961 | 89.7204 | 14.9030 | 30.8202 |

FIGURE 43. Average spread; copyable data.

| n | A2 | Alexander | B1 | Jones | Khovanovt1 |
|---|---|---|---|---|---|
| 3 | 1.1278 | 1.0000 | 1.0731 | 1.2106 | 1.1168 |
| 4 | 1.1837 | 2.6180 | 1.1080 | 1.2106 | 1.2406 |
| 5 | 1.1837 | 2.6180 | 1.1388 | 1.2837 | 1.2451 |
| 6 | 1.2196 | 2.6180 | 1.1990 | 1.6355 | 1.3042 |
| 7 | 1.2547 | 3.3165 | 1.2197 | 1.6355 | 1.6993 |
| 8 | 1.4348 | 4.3902 | 1.7655 | 1.8692 | 1.6993 |
| 9 | 1.7765 | 5.1069 | 1.8668 | 3.3168 | 1.7980 |
| 10 | 1.7765 | 5.1903 | 1.8668 | 3.3168 | 1.9942 |
| 11 | 1.7765 | 7.0507 | 1.9518 | 3.3168 | 2.0254 |
| 12 | 2.0162 | 8.7946 | 2.2990 | 4.3519 | 2.0504 |
| 13 | 2.1548 | 10.0233 | 2.4413 | 4.9314 | 2.1074 |
| 14 | 2.6858 | 14.2600 | 2.7764 | 5.6205 | 2.3300 |
| 15 | 2.7322 | 22.4258 | 3.7987 | 6.4831 | 2.4871 |
| 16 | 2.9737 | 30.5071 | 3.7987 | 7.4391 | 2.6661 |

FIGURE 44. Maximal absolute root; copyable data.

| n | A2 | Alexander | B1 | Jones | Khovanovt1 |
|---|---|---|---|---|---|
| 3 | 33.3333 | 0.0000 | 9.0909 | 33.3333 | 0.0000 |
| 4 | 28.5714 | 50.0000 | 4.0000 | 14.2857 | 11.1111 |
| 5 | 22.5806 | 20.0000 | 5.0847 | 17.6470 | 4.7619 |
| 6 | 22.5806 | 30.0000 | 4.2016 | 14.2857 | 9.5238 |
| 7 | 21.1267 | 17.3913 | 4.2857 | 14.2857 | 10.2040 |
| 8 | 19.9004 | 25.0000 | 4.3370 | 11.0204 | 11.4285 |
| 9 | 19.2200 | 21.5789 | 5.1858 | 11.3772 | 11.7021 |
| 10 | 19.0257 | 22.0211 | 5.6457 | 12.1688 | 11.7480 |
| 11 | 19.1843 | 20.7786 | 5.8111 | 11.7222 | 11.4177 |
| 12 | 18.8636 | 20.9580 | 5.9901 | 11.7900 | 10.9652 |
| 13 | 18.8596 | 19.7378 | 6.2311 | 12.0096 | 10.9197 |
| 14 | 18.7889 | 20.0468 | 6.3172 | 12.0427 | 10.6347 |
| 15 | 18.5564 | 19.3164 | 6.4275 | 12.2916 | 10.5783 |
| 16 | 18.5997 | 19.4536 | 6.4558 | 12.3533 | 10.4381 |

FIGURE 45. Percentage of pure roots; copyable data.

| n | A2 | Alexander | B1 | Jones | Khovanovt1 |
|---|---|---|---|---|---|
| 3 | 33.3333 | 100.0000 | 36.3636 | 0.0000 | 0.0000 |
| 4 | 57.1428 | 50.0000 | 56.0000 | 57.1428 | 11.1111 |
| 5 | 58.0645 | 80.0000 | 61.0169 | 23.5294 | 14.2857 |
| 6 | 54.8387 | 50.0000 | 65.5462 | 22.8571 | 16.6666 |
| 7 | 49.2957 | 65.2173 | 63.9285 | 14.2857 | 20.4081 |
| 8 | 47.7611 | 52.7777 | 60.4708 | 20.4081 | 21.0714 |
| 9 | 50.3249 | 53.1578 | 59.2932 | 19.7604 | 22.0744 |
| 10 | 49.4842 | 45.7013 | 56.9288 | 20.9699 | 22.2850 |
| 11 | 49.5595 | 43.9272 | 56.0352 | 22.0884 | 22.7319 |
| 12 | 48.8356 | 41.2244 | 54.2710 | 23.3678 | 24.0139 |
| 13 | 48.4072 | 39.1900 | 53.3642 | 23.9913 | 24.6212 |
| 14 | 48.6044 | 36.9571 | 52.5937 | 24.0480 | 25.0827 |
| 15 | 47.9160 | 35.8780 | 51.7600 | 23.5892 | 25.3112 |
| 16 | 47.5652 | 34.3872 | 51.1914 | 23.4344 | 25.6379 |

FIGURE 46. Percentage of roots in $[0.9, 1.1]$; copyable data.

### 13J. **Exercises.**

*Exercise* 13J.1. Construct the A2 invariant using Section 11H. ◇

*Exercise* 13J.2. Play with the data on [**TZ25b**]. ◇

*Exercise* 13J.3. Above, we wrote “Examples of quantum invariants that satisfy a skein relation or multiplicity freeness are A2, A, B1, and J”. Prove this for A2, B1 and J using the construction in Section 11H (or in Exercise 13J.1). ◇

*Exercise* 13J.4. Prove or disprove any of the stated conjectures ☺ ◇

## APPENDIX A. REPRESENTATION THEORY OF FINITE DIMENSIONAL ALGEBRAS WITH MAGMA

MAGMA is a computer algebra system designed to solve problems in algebra, number theory, geometry, combinatorics and related fields; in particular, representation theory, which is the underlying algebraic part of quantum invariants.

*Remark* A.1.

**(a)** Computer algebra system perform exact calculations; in particular, one can use MAGMA output in papers or theses without loosing the exactness.

**(b)** MAGMA is a noncommercial system, but the costs (such as the preparation of user documentation, the fixing of bugs, and the provision of user support) need to be recovered. So MAGMA is ***non-commercial but not free***, and the distribution is organized on a subscription basis. In order to get MAGMA on your machine use this site: http://magma.maths.usyd.edu.au/magma/ordering/

**(c)** Free, very useful, and completely enough for our calculations, is the ***online calculator*** http://magma.maths.usyd.edu.au/calc/. All the calculations below can be copied and pasted into the online calculator.

We refer to the handbook for more information: [BC23]; a beginner guide from which we stole a lot of the below can be found here [TT23]. ◇

MAGMA has two different ways to study representation theory of finite dimensional algebras: by directly constructing the representations, or using character theory. The latter is for finite groups only, and we will start with it.

AI. **Character theory.** Throughout, we consider finite groups. Recall that the ***character*** of a group representation is a function on the group that associates with each group element $g$ the trace of the corresponding matrix that encodes the action of $g$. The character carries the essential information about the representation in a more condensed form, and the easiest way to study representations of a finite group is to look at their character.

Recall that characters are constant on conjugacy classes, so it is enough to record them in a sequence $(c_1, ..., c_k)$ where $c_i$ is the value of the character on the $i$th conjugacy class.

The first thing to try is:

```
> G:=SymmetricGroup(3);
> X:=CharacterTable(G);
> X
```

----result----

```
> Character Table of Group G
>  --------------------------
>
>
> -----------------
> Class |   1  2  3
> Size  |   1  3  2
> Order |   1  2  3
> -----------------
> p  =  2   1  1  3
> p  =  3   1  2  1
> -----------------
> X.1   +   1  1  1
> X.2   +   1 -1  1
> X.3   +   2  0 -1
```

What do we see? Well, let us ignore finite characteristic and different fields for now, that is let us ignore

```
> -----------------
> p  =  2   1  1  3
> p  =  3   1  2  1
> -----------------
```

and the sign + in the table.

Then what we see is, from top to bottom, a numbering for the conjugacy classes, their sizes, their order (that is, the order of any element in the conjugacy class). Finally, the square matrix of characters, with the rows being the ***character sequences*** as above.

Let us double check the conjugacy classes:

```
> ConjugacyClasses(G);
```

----result----

```
> Conjugacy Classes of group G
> ----------------------------
> [1]     Order 1       Length 1
> Id(G)
>
> [2]     Order 2       Length 3
> (1, 2)
>
> [3]     Order 3       Length 2
> (1, 2, 3)
```

which works out with

```
> -----------------
> Class |   1  2  3
> Size  |   1  3  2
> Order |   1  2  3
> -----------------
```

MAGMA can perform character computations really fast. Here are slightly bigger examples, using the group databases `SmallGroup(n,j)` and `Group(''X'')`:

```
> G := SmallGroup(512,11600);
> time X := CharacterTable(G);
> #X;
> Degree(X[14]);
```

----result----

```
> Time: 0.030
> 44
> 2
```

```
> G := Group(''HS'');
> #G;
> time X := CharacterTable(G);
> #X;
> Degree(X[14]);
```

----result----

```
> Time: 0.160
> 44352000
> 24
> 896
```

Back to `SymmetricGroup(3)`. To get a specific row = character, or even a specific entry, we can use for example:

```
> X[3]; X[3][2];
```

----result----

```
> ( 2, 0, -1 )
> 0
```

Note however that a character is not a sequence. And that is good, because we can add and multiply characters:

```
> Type(X[3]); X[3]+X[2]; X[3]*X[2];
```

----result----

```
> AlgChtrElt
> ( 3, -1, 0 )
> ( 2, 0, -1 )
```

In fact, characters live in the character ring, so there are some other (not quite representation theoretical) operations, e.g.:

```
> X[3]-X[2]; 4*X[3];
> IsCharacter(X[3]-X[2]); IsGeneralizedCharacter(X[3]-X[2]);
```

----result----

```
> ( 1, 1, -2 )
> ( 8, 0, -4 )
> false
> true
```

Some more operations, hopefully self-explanatory, on characters are:

```
> IsIrreducible(X[3]); IsIrreducible(X[3]*X[3]);
> IsFaithful(X[2]); IsFaithful(X[3]);
> IsLinear(X[2]); IsLinear(X[3]);
```

----result----

```
> true
> false
> false
```

```
> true
> true
> false
```

Here linear means one dimensional.

*Remark* AI.1. Let us add the following known but maybe not well-known fact which explains why MAGMA can check faithfulness efficiently from the character table: a finite dimensional complex representation $V$ of a finite group is faithful $\Leftrightarrow$ $\dim V$ appears only in the character entry of the unit. ◇

*Exercise* AI.2. Run the following code

```
> for n in [2..14] do
> G:=SymmetricGroup(n);
> X:=CharacterTable(G);
> M:=[1..#X];
> for i in [1..#X] do
> if(IsFaithful(X[i]) eq true) then M[i]:=1; else M[i]:=0; end if;
> end for;
> M;
> end for;
```

and interpret the output. What happens for `SpecialLinearGroup(2,p)` ($p$ is a prime) instead of the symmetric group?

◇

Also, recall that characters are traces, so we can evaluate them:

```
> G.2 @ X[2]
----result----
> -1
```

Numerical values associated to characters:

```
> Degree(X[3]);
> InnerProduct(X[3],X[3]*X[3]);
> Norm(X[3]*X[3]);
> Indicator(X[3]);
----result----
> 2
> 1
> 3
> 1
```

We will come back to the ***Schur indicator*** `Indicator()` later. Otherwise, the ***norm*** is the inner product with itself, and the ***inner product*** `InnerProduct(x,y);` measures how often $x$ appears in $y$. `Degree` is the dimension.

With the inner product we can decompose characters into ***simple*** characters (simple means irreducible):

```
> G:=SymmetricGroup(3);
> X:=CharacterTable(G);
> y:=X[3]*X[3];
> M:=[1..#X];
> for i in [1..#X] do
> M[i]:=InnerProduct(y,X[i]);
> end for;
> M;
----result----
> [ 1, 1, 1 ]
```

*Exercise* AI.3. Consider following code:

```
> G:=SymmetricGroup(3);
> X:=CharacterTable(G);
> for n in [1..10] do
> y:=X[3]^(n);
> M:=[1..#X];
> for i in [1..#X] do
```

```
> M[i]:=InnerProduct(y,X[i]);
> end for;
> RealField(10)!(&+M)^(1/n);
> end for;
```

----result----

```
> 1.000000000
> 1.732050808
> 1.709975947
> 1.821160287
> 1.838416287
> 1.871736643
> 1.886389086
> 1.901623404
> 1.911688613
> 1.920622757
```

What is the code doing? Modify the code so that it works for power between 90 to 100, and for any representation of `SymmetricGroup(5)`. What is the limit $n \to \infty$? ◇

Let us ***restrict and induce***:

```
> G:=SymmetricGroup(5);
> H:=sub< SymmetricGroup(5) | (1,2), (2,3), (3,4)>;
> X:=CharacterTable(G);
> Y:=CharacterTable(H);
> Induction(Y[3],G);
> Restriction(X[3],H);
```

----result----

```
> ( 10, 0, 2, -2, 0, 0, 0 )
> ( 4, 0, -2, 1, 0 )
```

***Frobenius reciprocity*** (for the symmetric group) is then:

```
> G:=SymmetricGroup(5);
> H:=sub< SymmetricGroup(5) | (1,2), (2,3), (3,4)>;
> K:=sub< G | (1,2), (2,3)>;
> Y:=CharacterTable(H);
> Restriction(Induction(Y[3],G),H);
> Induction(Restriction(Y[3],K),H);
> Restriction(Induction(Y[3],G),H)-Induction(Restriction(Y[3],K),H)-Y[3];
```

----result----

```
> ( 10, 2, 0, -2, 0 )
> ( 8, 0, 0, -1, 0 )
> ( 0, 0, 0, 0, 0 )
```

Restriction and induction can be setup in many ways:

```
> G<x,y>:=PermutationGroup<23|
> [2,1,4,3,5,6,8,7,10,9,11,12,14,13,16,15,17,18,20,19,22,21,23],
> [16,9,1,5,8,22,7,23,21,10,3,2,20,18,17,11,15,6,19,13,12,14,4]>;
> CompositionFactors(G);
> #G;
> Factorization(#G);
> H:=SylowSubgroup(G,2);
> #H;
```

----result----

```
> G
> |  M23
> 1
> 10200960
> [ <2, 7>, <3, 2>, <5, 1>, <7, 1>, <11, 1>, <23, 1> ]
> 128
```

This is the ***Mathieu group*** $M_{23}$, and we took a 2 Sylow subgroup of it. We can now induce and restrict between $G$ and $H$, say the `PermutationCharacter(G,H)` (obtained from the right coset action of $G$ on $H$):

```
> x:=PermutationCharacter(G,H);
> x;
> Restriction(x,H);
> Induction(Restriction(x,H),G);
```

----result----

```
> ( 79695, 735, 0, 19, 0, 0, 0, 0, 1, 0, 0, 0, 0, 0, 0, 0, 0 )
> ( 79695, 735, 735, 735, 735, 735, 735, 735,
> 19, 19, 19, 19, 19, 19, 19, 19, 1 )
> ( 6351293025, 540225, 0, 361, 0, 0, 0, 0, 1, 0, 0, 0, 0, 0, 0, 0, 0 )
```

Finally, let us again say that MAGMA is really fast when it comes to characters. The following computes the character table of `SymmetricGroup(15)` from scratch.

```
> G:=sub< SymmetricGroup(16) | (1,2), (1,2,3,4,5,6,7,8,9,10,11,12,13,14,15)>;
> #G;
> Factorial(15);
> time X:=CharacterTable(G);
```

----result----

```
> 1307674368000
> 1307674368000
> Time: 4.090
```

*Remark* AI.4. Fun fact, the symmetric group `SymmetricGroup(15)` can be generated by two elements: a ***simple transition***, say $(1,2)$, and the ***long cycle***, here $(1,...,15)$. ◇

*Exercise* AI.5. What is the largest $n$ so that MAGMA's online calculator can output the character table of `SymmetricGroup(n)`?

Compare the construction from scratch with `time X:=CharacterTable(SymmetricGroup(n))`. ◇

Let us end with an instance of ***the law of small numbers***: At first sight the characters of the symmetric group seem to have small entries. But that is really only the case for small $n$; here is MAGMA code that computes the average character values of symmetric groups `SymmetricGroup(n)` up to `SymmetricGroup(15)`:

```
> for n in [1..15] do
> G:=SymmetricGroup(n);
> X:=CharacterTable(G);
> Av:=0;
> for i in [1..#X] do
> for j in [1..#X] do
> Av+:=Abs(X[i][j]);
> end for;
> end for;
> RealField(10)!(Av/(#X)^(2));
> end for;
```

----result----

```
> 1.000000000
> 1.000000000
> 1.000000000
> 1.080000000
> 1.244897959
> 1.528925620
> 2.084444444
> 2.950413223
> 4.768888889
> 8.195011338
> 15.95663265
> 31.45977399
> 70.97568866
> 161.9799726
> 402.7822831
```

*Exercise* AI.6. Write MAGMA code that gives the average degree of characters of symmetric groups. Compare to the average character values given above. ◇

AII. **Working with representations.** One can work with ***representations*** instead of characters. The advantage is that this works for any finite dimensional algebra, but at the cost of a slower performance. Let us look at an example to get started:

```
> G:=SymmetricGroup(7);
> X:=CharacterTable(G);
> y:=PermutationCharacter(G);
> y;
> M:=[1..#X];
> for i in [1..#X] do
> M[i]:=InnerProduct(y^5,X[i]);
> end for;
> time M;
```

----result----

```
> ( 7, 5, 1, 3, 4, 1, 3, 1, 2, 0, 2, 1, 0, 0, 0 )
> [ 52, 0, 1, 151, 160, 4, 5, 74, 160, 15, 75, 41, 30, 45, 150 ]
> Time: 0.000
```

```
> G:=SymmetricGroup(7);
> M:=PermutationModule(G,RationalField());
> M;
> time IndecomposableSummands(TensorPower(M,5));
```

----result----

**The computation exceeded the memory limit and so was terminated prematurely.**

```
GModule M of dimension 7 over Rational Field

Current total memory usage: 76.1MB, failed memory request: 1077.6MB
System Error: User memory limit has been reached
```

For the symmetric group, `PermutationModule` and `PermutationCharacter` are the representation respectively character obtained by the action of `SymmetricGroup(n)` on a vector space (using complex numbers for the character and a field that we specify for the representation) of dimension $n$ by permutation of a fixed basis.

*Remark* AII.1. Note that MAGMA actually works with ***modules***, i.e. vector spaces with actions, and not with ***representations***, i.e. the homomorphisms defining the action. We will see that in more details below. In any case, we abuse language and say representations. ◇

*Exercise* AII.2. Compare the third tensor powers, via characters and representations, and their timing. ◇

Ok, computational aspects aside, we can now setup nonsemisimple representations. We start with `PGL(3,4)`, the ***projective linear group*** $PGL(\mathbb{F}_4^3) = GL(\mathbb{F}_4^3)/\text{center}$ over the Galois field `GF(4)` with four elements and look at representations over `GF(3)`:

```
> G:=PGL(3,4);
> M:=PermutationModule(G,GF(3));
> M; IsSemisimple(M);
```

----result----

```
> GModule M of dimension 21 over GF(3)
> false
```

Here we get a 21 dimensional representation using `PermutationModule` because `PGL(3,4)` is realized in MAGMA within `SymmetricGroup(21)`:

```
> PGL(3,4);
```

----result----

```
> Permutation group acting on a set of cardinality 21
> Order = 60480 = 2^6 * 3^3 * 5 * 7
> (4, 14, 21)(5, 15, 17)(7, 10, 20)(8, 9, 12)(11, 13, 19)
> (1, 8, 21, 16, 15, 3, 2)(4, 10, 20, 18, 17, 9, 7)
> (5, 12, 11, 14, 19, 13, 6)
```

To move on, we note that this is indeed not semisimple, so the composition factors are different from the direct summands:

```
> CompositionSeries(M);
> CompositionFactors(M);
> IsIrreducible(M);
> IndecomposableSummands(M);
  ----result----
> [
> GModule of dimension 1 over GF(3),
> GModule of dimension 20 over GF(3),
> GModule M of dimension 21 over GF(3)
> ]
> [
> GModule of dimension 1 over GF(3),
> GModule of dimension 19 over GF(3),
> GModule of dimension 1 over GF(3)
> false
> ]
> [
> GModule M of dimension 21 over GF(3)
> ]
```

Let us test whether the composition factors are the same or not:

```
> X:=CompositionFactors(M); X;
> IsIsomorphic(X[1],X[3]);
> ConstituentsWithMultiplicities(M);
  ----result----
> true
> [
> <GModule of dimension 1 over GF(3), 2>,
> <GModule of dimension 19 over GF(3), 1>
> ]
> [ 1, 2, 1 ]
```

Many other tricks can be played as soon as a representation is setup, e.g.:

```
> IndecomposableSummands(TensorPower(M,2));
> IndecomposableSummands(ExteriorPower(M,2));
> IndecomposableSummands(SymmetricPower(M,2));
  ----result----
> [
> GModule of dimension 21 over GF(3),
> GModule of dimension 21 over GF(3),
> GModule of dimension 84 over GF(3),
> GModule of dimension 126 over GF(3),
> GModule of dimension 189 over GF(3)
> ]
> [
> GModule of dimension 84 over GF(3),
> GModule of dimension 126 over GF(3)
> ]
> [
> GModule of dimension 21 over GF(3),
> GModule of dimension 21 over GF(3),
> GModule of dimension 189 over GF(3)
> ]
```

Or submodules can be constructed:

```
> G:=PGL(3,4);
> M:=PermutationModule(G,GF(3));
> N:=JacobsonRadical(M);
> I:=IndecomposableSummands(TensorProduct(M,N));
```

```
> I; IsIsomorphic(I[1],M);
```

----result----

```
> [
> GModule of dimension 21 over GF(3),
> GModule of dimension 84 over GF(3),
> GModule of dimension 126 over GF(3),
> GModule of dimension 189 over GF(3)
> ]
> true
```

So, since the usual algebra can just be run in MAGMA with only the memory being a problem, the crucial question remains how to construct representations to begin with.

For groups one can construct all simple representation and all projective indecomposables over any finite field in one go, say for the 3-by-3 general linear group with coefficients in `GF(2)`:

```
> G:=PGL(3,2);
> X:=IrreducibleModules(G,GF(2));
> Y:=IrreducibleModules(G,GF(3));
> X; Y;
```

----result----

```
> [
> GModule of dimension 1 over GF(2),
> GModule of dimension 3 over GF(2),
> GModule of dimension 3 over GF(2),
> GModule of dimension 8 over GF(2)
> ]
> [
> GModule of dimension 1 over GF(3),
> GModule of dimension 6 over GF(3),
> GModule of dimension 6 over GF(3),
> GModule of dimension 7 over GF(3)
> ]
```

```
> A:=ProjectiveIndecomposableModules(G,GF(2));
> B:=ProjectiveIndecomposableModules(G,GF(3));
> A; B;
```

----result----

```
> [
> GModule of dimension 8 over GF(2),
> GModule of dimension 16 over GF(2),
> GModule of dimension 16 over GF(2),
> GModule of dimension 8 over GF(2)
> ]
> [
> GModule of dimension 9 over GF(3),
> GModule of dimension 6 over GF(3),
> GModule of dimension 6 over GF(3),
> GModule of dimension 15 over GF(3)
> ]
```

***Projective covers*** can also be constructed easily:

```
> ProjectiveCover(X[1]);
```

----result----

```
> GModule of dimension 8 over GF(2)
>
> [1]
> [1]
> [1]
> [1]
> [1]
> [1]
> [1]
```

```
> [1]
```

Similarly, one can use `InjectiveHull()` to construct ***injective hulls***.

*Remark* AII.3. Note that MAGMA has a bias and likes projective things over injective things. However, for a group this does not make a difference since the projective covers and injective hulls coincide. ◇

Let us look at a slightly bigger example. We consider the cyclic group `CyclicGroup(5)`, this is $\mathbb{Z}/5\mathbb{Z}$, and its representations over `GF(5)`. This group has five indecomposables (called $Z$ below) and only one simple (called $L$ below) representation. The indecomposables are the representations where the generator of $\mathbb{Z}/5\mathbb{Z}$ acts by the following ***Jordan blocks***:

$$Z_1 = L_1 \colon \begin{pmatrix} 1 \end{pmatrix}, \text{ is simple}$$

$$Z_2 \colon \begin{pmatrix} 1 & 0 \\ 1 & 1 \end{pmatrix}, \text{ filtration } 0 - L_1 - Z_2,$$

$$Z_3 \colon \begin{pmatrix} 1 & 0 & 0 \\ 1 & 1 & 0 \\ 0 & 1 & 1 \end{pmatrix}, \text{ filtration } 0 - L_1 - L_1 - Z_3,$$

$$Z_4 \colon \begin{pmatrix} 1 & 0 & 0 & 0 \\ 1 & 1 & 0 & 0 \\ 0 & 1 & 1 & 0 \\ 0 & 0 & 1 & 1 \end{pmatrix}, \text{ filtration } 0 - L_1 - L_1 - L_1 - Z_4,$$

$$Z_5 = P_1 \colon \begin{pmatrix} 1 & 0 & 0 & 0 & 0 \\ 1 & 1 & 0 & 0 & 0 \\ 0 & 1 & 1 & 0 & 0 \\ 0 & 0 & 1 & 1 & 0 \\ 0 & 0 & 0 & 1 & 1 \end{pmatrix}, \text{ filtration } 0 - L_1 - L_1 - L_1 - L_1 - Z_5,$$

We construct the three dimensional $\mathbb{Z}/5\mathbb{Z}$-representation $Z3 = Z_3$ that is indecomposable but not simple:

```
G:=CyclicGroup(5);
Z3:=CompositionSeries(ProjectiveIndecomposables(G,GF(5))[1])[3];
IsIrreducible(Z3);
IndecomposableSummands(Z3);
  ----result----
> false GModule of dimension 1 over GF(5)
> GModule of dimension 2 over GF(5)
> [
> GModule X of dimension 3 over GF(5)
> ]
```

There are additionally $Z1$ and $Z5$ of dimensions one and five, respectively, which are the trivial $\mathbb{Z}/5\mathbb{Z}$-representation and the regular $\mathbb{Z}/5\mathbb{Z}$-representation. We have $Z3 \otimes Z3 \cong Z1 \oplus Z3 \oplus Z5$, with only $Z5$ being projective:

```
> P:=IndecomposableSummands(TensorPower(Z3,2));
> P;
> IsProjective(P[1]);
> IsProjective(P[2]);
> IsProjective(P[3]);
  ----result----
> [
> GModule of dimension 1 over GF(5),
> GModule of dimension 3 over GF(5),
> GModule of dimension 5 over GF(5)
> ]
> false
> false
> true
```

*Exercise* AII.4. Check with MAGMA that $Z3 \otimes Z5 \cong Z5^{\oplus 3}$. ◇

Since $Z1$ is the trivial $\mathbb{Z}/5\mathbb{Z}$-representation, if we annihilate $Z5$, then we have that $Z3^{\otimes 2}$ satisfies $X^2 = 1 + X$ whose roots are the ***golden ratio*** and its Galois conjugate. Thus, tensor powers of $Z3$ give the ***Fibonacci sequence***:

```
> G:=CyclicGroup(5);
> Z3:=CompositionSeries(ProjectiveIndecomposables(G,GF(5))[1])[3];
```

```
> Y:=[Z3];
> n:=10;
> 1;
> for i in [2..n] do
> X:=[];
> for y in Y do
> P:=IndecomposableSummands(TensorProduct(Z3,y));
> for p in P do
> if(Dimension(p) ne 5) then
> X:=Append(X,p);
> end if;
> end for;
> Y:=X;
> end for;
> #X;
> end for;
```

```
----result----
> 1
> 2
> 3
> 5
> 8
> 13
> 21
> 34
> 55
> 89
```

*Exercise* AII.5. Confirm in MAGMA that the $n$th root of the sequence above gives the golden ratio. What happens if one includes $Z5$ and takes the $n$th root?

(Careful: the convergence is slow, and you might not be able to 100% verify what the limit is in the online calculator.) ◇

In general, one can construct representations by specifying the ***acting matrices***. This is in particularly useful when one has a generator-relation presentations. For example:

```
> G<x,y>:=PermutationGroup<11|[1,10,3,11,7,6,5,9,8,2,4],
> [4,5,8,3,6,9,7,1,2,10,11]>;
> F:=GF(3);
> x:=Matrix(F, 5, 5,
> [[0,2,0,1,1],[2,1,1,0,2],[1,1,1,2,2],
> [0,2,2,2,2],[0,2,2,1,0]]);
> y:=Matrix(F, 5, 5,
> [[0,1,1,0,1],[2,0,0,1,0],[0,0,1,2,2],
> [2,1,0,0,0],[0,1,2,2,1]]);
> M:=GModule(G,[x,y]);
> M;
```

```
----result----
> GModule M2 of dimension 5 over GF(3)
```

*Exercise* AII.6. Check what happens if there would be a typo in the action matrices. ◇

We can run the same syntax as above, for example:

```
> IndecomposableSummands(TensorProduct(M,N));
```

```
----result----
> [
> GModule of dimension 10 over GF(3),
> GModule of dimension 15 over GF(3)
> ]
```

For a general finite dimensional algebra, let us setup an algebra given by matrices:

```
> K:=GF(2);
> A := MatrixAlgebra<K, 4 |
```

```
> [1,0,0,0, 0,1,0,0, 0,0,1,0, 0,0,0,1],
> [0,0,0,0, 1,0,0,0, 0,0,0,0, 0,0,1,0],
> [0,0,0,0, 0,0,0,0, 1,0,0,0, 0,1,0,0] >;
> Dimension(A);

  ----result----

> 4;
```

The algebra $A$ above is $A = \mathbb{F}_2[X,Y]/(X^2,Y^2)$, which is the group ring of the ***Klein four group*** $\mathbb{Z}/2\mathbb{Z} \times \mathbb{Z}/2\mathbb{Z}$ in characteristic two. We can now tell MAGMA to work with its ***regular representation***:

```
> V:=RModule(K,4);
> m:=map< CartesianProduct(V, A) -> V | t :-> t[1]*t[2] >;
> M:=Module(A,m);
> IndecomposableSummands(M);

  ----result----

> [
> Right A-module of dimension 4,
> where A is Matrix Algebra of degree 4 with 3
> generators over GF(2)
> ]
```

*Remark* AII.7. With constructions of the form

```
> V:=RModule(K,4);
> m:=map< CartesianProduct(V, A) -> V | t :-> t[1]*t[2] >;
> M:=Module(A,m);
> IndecomposableSummands(M);
```

one can construct any representation as long as one has a good control over the action map encoded by $m$ in this example. ◇

```
> CompositionFactors(M);

  ----result----

> [
> RModule of dimension 1 over GF(2),
> RModule of dimension 1 over GF(2),
> RModule of dimension 1 over GF(2),
> RModule of dimension 1 over GF(2)
> ]
```

*Remark* AII.8. Sometimes one gets bugs. The code

```
> K:=GF(2);
> A := MatrixAlgebra<K, 4 |
> [1,0,0,0, 0,1,0,0, 0,0,1,0, 0,0,0,1],
> [0,0,0,0, 1,0,0,0, 0,0,0,0, 0,0,1,0],
> [0,0,0,0, 0,0,0,0, 1,0,0,0, 0,1,0,0] >;
> V:=RModule(K,4);
> m:=map< CartesianProduct(V, A) -> V | t :-> t[1]*t[2] >;
> M:=Module(A,m);
> CompositionSeries(M);
```

which is not much different from the code above, produced (or produces, depending on the version you use):

```
[
    A-module of dimension 1, where A is

Magma: Internal error

Please mail this entire run [*** WITH THE DETAILS BELOW ***]
    to magma-bugs@maths.usyd.edu.au

You can print the entire input by entering:
%P

Version: 2.28-2
Initial seed: 553935780
Time to this point: 0.01
Memory usage: 32.09MB
Segmentation fault
```

For a different way to setup `CompositionSeries` for the Klein four group see below. ◇

Here are a few constructions to build new representations:

```
> N:=DirectSum([M,M]);
> ActionMatrix(N,A.2);
> W:= ModuleWithBasis([ M.1+M.2+M.3, M.2+M.3, M.3 ]);
> ActionMatrix(W,A.2);
> IndecomposableSummands(W);
```

----result----

```
> [0 0 0 0 0 0 0 0]
> [1 0 0 0 0 0 0 0]
> [0 0 0 0 0 0 0 0]
> [0 0 1 0 0 0 0 0]
> [0 0 0 0 0 0 0 0]
> [0 0 0 0 1 0 0 0]
> [0 0 0 0 0 0 0 0]
> [0 0 0 0 0 0 1 0]
> [1 1 0]
> [1 1 0]
> [0 0 0]
> [
> Right A-module of dimension 3, where A is
> Matrix Algebra of degree 4 and
> dimension 4 with 3 generators over GF(2)
> ]
```

Let us set up another example. Again, take the ***Klein four group*** $\mathbb{Z}/2\mathbb{Z} \times \mathbb{Z}/2\mathbb{Z}$ but setup differently:

```
> G:=SmallGroup(4,2);
> IsCyclic(G);
```

----result----

```
> false
```

Since there are only two groups of order four, and the one we took is not cyclic, it must be the Klein four group. Let us pick an indecomposable three dimensional representation over `GF(2)`:

```
> X:=CompositionSeries(ProjectiveCover(IrreducibleModules(G,GF(2))[1]))[3];
> IndecomposableSummands(X);
```

----result----

```
> [
> GModule X of dimension 3 over GF(2)
> ]
```

Now we take tensor products:

```
> for i in [1..5] do
> IndecomposableSummands(TensorPower(X,i))
> [#IndecomposableSummands(TensorPower(X,i))];
> end for;
```

```
----result----
> GModule X of dimension 3 over GF(2)
> GModule of dimension 5 over GF(2)
> GModule of dimension 7 over GF(2)
> GModule of dimension 9 over GF(2)
> GModule of dimension 11 over GF(2)
```

In each step we get a new indecomposable of the Klein four group! (The dimension verifies that these are really new.) One can in fact check that `SmallGroup(4,2)` has infinitely many indecomposables over `GF(2)`.

To extract the acting matrices, say for the representation of dimension five, we can run the following:

```
> M:=IndecomposableSummands(TensorPower(X,2))
[#IndecomposableSummands(TensorPower(X,2))];
> phi := Representation(M);
> phi(G.1)+phi(G.2);
----result----
> [0 0 1 0 0]
> [0 0 1 0 0]
> [0 0 0 0 0]
> [0 0 0 0 0]
> [1 1 0 0 0]
```

This outputs the sum of the action matrices of $(1,0)$ and $(0,1)$ in $\mathbb{Z}/2\mathbb{Z}\times\mathbb{Z}/2\mathbb{Z}$. Note that the decomposition algorithm underlying `IndecomposableSummands` is randomized and we (can) get different basis, hence action matrices, in every evaluation.

*Exercise* AII.9. Compare the constructions of representations of the Klein four group via `MatrixAlgebra` and `SmallGroup(4,2)` above. What are the differences? ◇

Finally, MAGMA is really awesome when working with ***basic algebras*** (quivers if you want). Here is one example: we setup a "Temperley–Lieb algebra" and compute the dimensions of the projective representations and their composition factors:

```
> K:=GF(2);
> G:=SL(2,K);
> X:=IrreducibleModules(G,K);
> V:=TensorPower(X[2],4);
> A:=EndomorphismAlgebra(V);
> B:=BasicAlgebra(A);
> DimensionsOfProjectiveModules(B);
> DimensionsOfProjectiveModules(B);
> CompositionFactors(ProjectiveModule(B,1));
> CompositionFactors(ProjectiveModule(B,2));
----result----
> [ 1, 2 ]
> [
> AModule of dimension 1 over GF(2)
> ]
> [
> AModule of dimension 1 over GF(2),
> AModule of dimension 1 over GF(2)
> ]
```

AIII. **Representation theory over different fields.** Let us have a look at the character table of $\mathbb{Z}/5\mathbb{Z}$:

```
> X:=CharacterTable(CyclicGroup(5)); X;
----result----
> Character Table
> ---------------
>
>
> -------------------------------
> Class |   1    2    3    4    5
> Size  |   1    1    1    1    1
> Order |   1    5    5    5    5
```

```
> ------------------------------
> p  =  5   1    1    1    1    1
> ------------------------------
> X.1   +   1    1    1    1    1
> X.2   0   1   Z1 Z1#2 Z1#3 Z1#4
> X.3   0   1 Z1#4 Z1#3 Z1#2   Z1
> X.4   0   1 Z1#3   Z1 Z1#4 Z1#2
> X.5   0   1 Z1#2 Z1#4   Z1 Z1#3
>
>
> Explanation of Character Value Symbols
> --------------------------------------
>
> # denotes algebraic conjugation, that is,
> #k indicates replacing the root of unity w by w^k
>
> Z1     = (CyclotomicField(5: Sparse := true)) ! [
> RationalField() | 0, 0, 0, 1 ]
```

One immediately hits the problem of working with roots of unity. Thus, working over $\mathbb{R}$ instead of $\mathbb{C}$ is not completely straightforward. However, the ***Schur indicator*** (or ***Frobenius–Schur indicator***) comes to our rescue. This is the entry $+$, sometimes written 1, and 0 in the table (it could also be $-$ or $-a$, but not in this example), and we can check this as follows:

```
> X:=CharacterTable(CyclicGroup(5));
> for i in [1..#X] do Indicator(X[i]); end for;
----result----
> 1
> 0
> 0
> 0
> 0
```

The indicator gives us all the information we need to work over $\mathbb{R}$:

**(a)** If its 1, then there is nothing to do and the complex simple character is also a real simple character.

**(b)** If its 0, then we need to add two complex simple characters together (one plus its conjugate) to get a real simple character.

**(c)** The case $-1$ is the exercise below.

For example, to get the simple real characters in this case we do:

```
> X:=CharacterTable(CyclicGroup(5));
> IsReal(X[1]);
> IsReal(X[2]); IsReal(X[2]+X[3]);
> IsReal(X[4]); IsReal(X[4]+X[5]);
----result----
> true
> false
> true
> false
> true
```

*Exercise* AIII.1. Use MAGMA to find a group with Schur indicator $-1$. What happens for the representations when you use `IsReal`? ◇

The second thing that comes to mind is character theory over the algebraic closure of the Galois fields. That is easy in MAGMA can controlled by `BrauerCharacterTable(G,p)` with $p$ being the characteristic, e.g.:

```
> G:=SymmetricGroup(4);
> X:=BrauerCharacterTable(G,3);
> X;
> X[3];
> CharacterTable(G);
----result----
```

```
> [
> ( 1, 1, 1, 1 ),
> ( 1, 1, -1, -1 ),
> ( 3, -1, -1, 1 ),
> ( 3, -1, 1, -1 )
> ]
> ( 3, -1, -1, 1 )
>
> Character Table of Group G
> --------------------------
>
>
> -----------------------
> Class |   1  2  3  4  5
> Size  |   1  3  6  8  6
> Order |   1  2  2  3  4
> -----------------------
> p  =  2   1  1  1  4  2
> p  =  3   1  2  3  1  5
> -----------------------
> X.1   +   1  1  1  1  1
> X.2   +   1  1 -1  1 -1
> X.3   +   2  2  0 -1  0
> X.4   +   3 -1 -1  0  1
> X.5   +   3 -1  1  0 -1
```

The Brauer character table is a square matrix of size being the $p$-irregular conjugacy classes. The the example above the 3-irregular conjugacy classes are $1, 2, 3, 5$, while the 3-irregular conjugacy classes are $1, 4, 5$, as indicated by the orders coprime to $p$ in this row:

```
> Order |   1  2  2  3  4
```

*Remark* AIII.2. For soluble groups the Brauer characters can be deduced from the complex character table, and this is what MAGMA does. Without going into details, if one deletes the column four for $p = 3$ one gets:

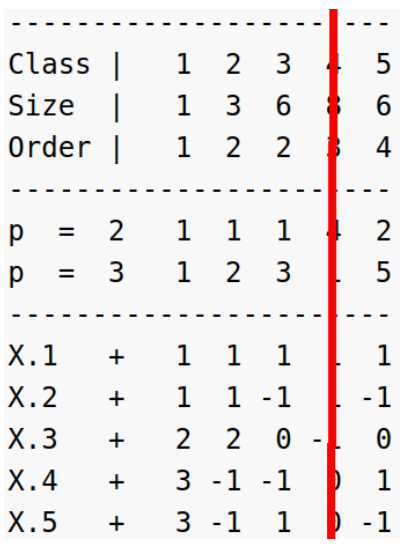

```
---------------------
Class |   1  2  3  4  5
Size  |   1  3  6  8  6
Order |   1  2  2  3  4
---------------------
p  =  2   1  1  1  4  2
p  =  3   1  2  3  1  5
---------------------
X.1   +   1  1  1  1  1
X.2   +   1  1 -1  1 -1
X.3   +   2  2  0 -1  0
X.4   +   3 -1 -1  0  1
X.5   +   3 -1  1  0 -1
```

Now $X.3 = X.1 + X.2$ and that is why this is dropped from the Brauer character table.

For non-soluble groups things are trickier and MAGMA constructs the simple representations. ◇

Essentially the same syntax as for characters works, so we skip that part. But there are also new interesting feature, such as `Blocks` (giving as output the blocks and their ***defect***) and `DefectGroup` (returning the ***defect group*** of a block):

```
> G:=SymmetricGroup(4);
> Y:=CharacterTable(G);
> Blocks(Y,3);
> DefectGroup(Y,Blocks(Y,3)[1],3);
```

```
----result----
```

```
> [
> { 1, 2, 3 },
> { 4 },
> { 5 }
> ]
> [ 1, 0, 0 ]
> Permutation group acting on a set of cardinality 4
> Order = 3
```

```
> (1, 3, 2)
```

These calculations are efficient and fast, for example:

```
> G:=SymmetricGroup(12);
> Y:=CharacterTable(G);
> time IsAbelian(DefectGroup(Y,Blocks(Y,3)[1],3));
  ----result----
> false
> Time: 0.020
```

Another option for characters are ***rational characters***, e.g.:

```
> G:=AlternatingGroup(4);
> CharacterTable(G);
> RationalCharacterTable(G);
  ----result----
> Character Table of Group G
> --------------------------
>
>
> ------------------------
> Class |   1  2    3    4
> Size  |   1  3    4    4
> Order |   1  2    3    3
> ------------------------
> p  =  2   1  1    4    3
> p  =  3   1  2    1    1
> ------------------------
> X.1   +   1  1    1    1
> X.2   0   1  1    J -1-J
> X.3   0   1  1 -1-J    J
> X.4   +   3 -1    0    0
>
>
> Explanation of Character Value Symbols
> --------------------------------------
>
> J = RootOfUnity(3)
>
> [
> ( 1, 1, 1, 1 ),
> ( 2, 2, -1, -1 ),
> ( 3, -1, 0, 0 )
> ]
> [
> [ 1 ],
> [ 2, 3 ],
> [ 4 ]
> ]
```

In this example the characters two and three are not rational, but their direct sum is.

*Remark* AIII.3. In general, the rational characters, are the sums of the Galois orbits on the complex character table. Hence, the story over $\mathbb{Q}$ is similar, but more complicated, than over $\mathbb{R}$. ◇

Working with algebraically closed fields is not easy on a machine. The characters avoid that problem by using a combinatorial approach. For representations we cannot cheat so let us have a look at different fields, starting with algebraically closed ones.

One can construct, e.g., $\overline{\mathbb{Q}}$ or $\overline{\mathbb{F}_p}$, but a lot of operations will not work. For example:

```
> AlgebraicClosure(GF(2));
> AlgebraicClosure(RationalField());
> SpecialLinearGroup(2,AlgebraicClosure(GF(2)));
  ----result----
```

```
> Algebraically closed field with no variables over GF(2)
> Algebraically closed field with no variables over Rational Field
>
> >> SpecialLinearGroup(2,AlgebraicClosure(GF(2)));
> ^
> Runtime error in 'SpecialLinearGroup': Cannot compute generators for matrix
> group
```

In finite characteristic it is better to simulate $\overline{\mathbb{F}_p}$ by making the underlying field big enough. This can be done afterwards, for example by ***base change***:

```
> G:=SpecialLinearGroup(2,4);
> X:=IrreducibleModules(G,GF(2));
> Y:=TensorPower(X[2],2);
> for k in [1..6] do
> #IndecomposableSummands(ChangeRing(Y,GF(2^k)));
> end for;
```

```
----result----
```

```
> 3
> 4
> 3
> 4
> 3
> 4
```

In characteristic zero one can use ***cyclotomic fields*** as approximations. For example:

```
> K:=RationalField();
> F:=CyclotomicField(3);
> G:=AlternatingGroup(4);
> A:=GroupAlgebra(K,G);
> V:=RegularRepresentation(A);
> X, _:=IndecomposableSummands(V);
> Y, _:=IndecomposableSummands(ChangeRing(V,F));
> X; Y;
```

```
----result----
```

```
> [
> Matrix Algebra [ideal of V] of degree 12 and dimension 1 over
> Rational Field,
> Matrix Algebra [ideal of V] of degree 12 and dimension 9 over
> Rational Field,
> Matrix Algebra [ideal of V] of degree 12 and dimension 2 over
> Rational Field
> ]
> [
> Matrix Algebra [ideal] of degree 12 and dimension 1 over F,
> Matrix Algebra [ideal] of degree 12 and dimension 9 over F,
> Matrix Algebra [ideal] of degree 12 and dimension 1 over F,
> Matrix Algebra [ideal] of degree 12 and dimension 1 over F
> ]
```

Compare this to the rational and complex character tables of `AlternatingGroup(4)` above.

*Exercise* AIII.4. The above code gives the decomposition into bimodules, and that is why the $9 = 3 \cdot 3$ is appearing. Ask MAGMA to give the decomposition into representations. For example, you could setup `AlternatingGroup(4)` as a twelve dimensional matrix algebra, and define its regular representation by acting on a twelve dimensional vector space. $\diamond$

D.T.: The University of Sydney, School of Mathematics and Statistics F07, Office Carslaw 827, NSW 2006, Australia, www.dtubbenhauer.com, ORCID 0000-0001-7265-5047

*Email address*: daniel.tubbenhauer@sydney.edu.au